\documentclass[10pt]{report}
\usepackage{amsfonts,amsmath}
\usepackage{color}
\usepackage{epsfig}

\begin{document}

\newcommand{\s}{\mbox{$\mkern 7mu$}}

\newcommand{\ub}{\underline{u}}
\newcommand{\Cb}{\underline{C}}
\newcommand{\Lb}{\underline{L}}
\newcommand{\Nb}{\underline{N}}
\newcommand{\alphab}{\underline{\alpha}}
\newcommand{\phib}{\underline{\phi}}
\newcommand{\Phib}{\underline{\Phi}}
\newcommand{\Pb}{\underline{P}}
\newcommand{\Mb}{\underline{M}}
\newcommand{\rhob}{\underline{\rho}}
\newcommand{\lambdab}{\underline{\lambda}}
\newcommand{\pb}{\underline{p}}
\newcommand{\qb}{\underline{q}}

\newcommand{\hb}{\underline{h}}

\newcommand{\sL}{{\cal L}\mkern-10mu /}
\newcommand{\sPi}{\Pi\mkern-11mu /}
\newcommand{\sF}{F\mkern-11mu /}
\newcommand{\oF}{\overline{F}}
\newcommand{\Fb}{\underline{F}}

\newcommand{\rb}{\underline{r}}
\renewcommand{\sb}{\underline{s}}
\newcommand{\os}{\overline{s}}

\newcommand{\ib}{\underline{i}}
\newcommand{\jb}{\underline{j}}

\newcommand{\sk}{k\mkern-9mu /}
\newcommand{\ok}{\overline{k}}
\newcommand{\kb}{\underline{k}}
\newcommand{\skb}{\underline{\sk}}
\newcommand{\okb}{\underline{\ok}}

\newcommand{\sm}{m\mkern-9mu /}
\newcommand{\om}{\overline{m}}
\newcommand{\mb}{\underline{m}}
\newcommand{\smb}{\underline{\sm}}
\newcommand{\omb}{\underline{\om}}

\newcommand{\sn}{n\mkern-9mu /}
\newcommand{\on}{\overline{n}}
\newcommand{\nb}{\underline{n}}
\newcommand{\snb}{\underline{\sn}}
\newcommand{\onb}{\underline{\on}}

\newcommand{\chib}{\underline{\chi}}
\newcommand{\chih}{\hat{\chi}}
\newcommand{\chibh}{\hat{\chib}}
\newcommand{\etab}{\underline{\eta}}
\newcommand{\psib}{\underline{\psi}}

\newcommand{\tchi}{\tilde{\chi}}
\newcommand{\tchib}{\tilde{\chib}}
\newcommand{\tchih}{\hat{\tilde{\chi}}}
\newcommand{\tchibh}{\hat{\tilde{\chib}}}
\newcommand{\teta}{\tilde{\eta}}
\newcommand{\tetab}{\tilde{\etab}}

\newcommand{\db}{\underline{d}}
\newcommand{\gb}{\underline{g}}
\newcommand{\eb}{\underline{e}}
\newcommand{\lb}{\underline{l}}

\newcommand{\thetab}{\underline{\theta}}
\newcommand{\gammab}{\underline{\gamma}}
\newcommand{\mub}{\underline{\mu}}
\newcommand{\nub}{\underline{\nu}}

\newcommand{\Thetab}{\underline{\Theta}}
\newcommand{\Lambdab}{\underline{\Lambda}}
\newcommand{\Gammab}{\underline{\Gamma}}
\newcommand{\Kb}{\underline{K}}
\newcommand{\Ib}{\underline{I}}
\newcommand{\Wb}{\underline{W}}
\newcommand{\Rb}{\underline{R}}
\newcommand{\Qb}{\underline{Q}}

\newcommand{\Gb}{\underline{G}}

\newcommand{\spi}{\pi\mkern-11mu /}
\newcommand{\sspi}{\spi\mkern-7mu /}

\newcommand{\sJ}{J\mkern-10mu /}
\newcommand{\sV}{V\mkern-13mu /}
\renewcommand{\sp}{p\mkern-9mu /}

\newcommand{\oE}{\overline{E}}
\newcommand{\oD}{\overline{D}}
\newcommand{\oM}{\overline{M}}
\newcommand{\oU}{\overline{U}}
\renewcommand{\oD}{\overline{D}}

\newcommand{\sbeta}{\beta\mkern-9mu /}
\newcommand{\stheta}{\theta\mkern-9mu /}
\newcommand{\sstheta}{{\theta\mkern-10mu /}\mkern-7mu /}
\newcommand{\sxi}{\xi\mkern-9mu /}
\newcommand{\szeta}{\zeta\mkern-9mu /}
\newcommand{\zetab}{\underline{\zeta}}

\renewcommand{\sf}{f\mkern-8mu /}
\newcommand{\of}{\overline{f}}
\newcommand{\fb}{\underline{f}}
\newcommand{\sfb}{\underline{\sf}}
\newcommand{\ofb}{\underline{\of}}

\newcommand{\oA}{\overline{A}}
\newcommand{\Ab}{\underline{A}}
\newcommand{\oAb}{\underline{\oA}}

\newcommand{\oB}{\overline{B}}
\newcommand{\Bb}{\underline{B}}
\newcommand{\oBb}{\underline{\oB}}

\newcommand{\ep}{\epsilon}
\newcommand{\epb}{\underline{\ep}}
\newcommand{\vep}{\varepsilon}
\newcommand{\vepb}{\underline{\vep}}

\newcommand{\svep}{\vep\mkern-9mu /}
\newcommand{\svepb}{\vepb\mkern-9mu /}

\newcommand{\se}{e\mkern-9mu /}
\renewcommand{\oe}{\overline{e}}

\newcommand{\salpha}{\alpha\mkern-9mu /}

\newcommand{\sa}{a\mkern-9mu /}
\newcommand{\sbb}{b\mkern-9mu /}
\newcommand{\sg}{g\mkern-9mu /}
\newcommand{\sh}{h\mkern-9mu /}
\newcommand{\seps}{\epsilon\mkern-8mu /}
\newcommand{\sd}{d\mkern-10mu /}
\newcommand{\sR}{R\mkern-10mu /}
\newcommand{\sP}{P\mkern-10mu /}
\newcommand{\snab}{\nabla\mkern-13mu /}
\newcommand{\sdiv}{\mbox{div}\mkern-19mu /\,\,\,\,}
\newcommand{\scurl}{\mbox{curl}\mkern-19mu /\,\,\,\,}
\newcommand{\slap}{\mbox{$\triangle  \mkern-13mu / \,$}}
\newcommand{\sGamma}{\Gamma\mkern-10mu /}
\newcommand{\sgamma}{\gamma\mkern-9mu /}
\newcommand{\sRic}{\mbox{Ric}\mkern-19mu /\,\,\,\,}

\newcommand{\tgamma}{\tilde{\gamma}}
\newcommand{\ogamma}{\overline{\gamma}}

\newcommand{\sD}{D\mkern-11mu /}

\newcommand{\sB}{B\mkern-11mu /}

\renewcommand{\ss}{s\mkern-8mu /}
\newcommand{\sss}{\ss\mkern-7mu /} 

\newcommand{\cth}{\check{\theta}}
\newcommand{\cthb}{\check{\thetab}}
\newcommand{\cnu}{\check{\nu}}
\newcommand{\cnub}{\check{\nub}}

\newcommand{\cE}{\check{{\cal E}}}
\newcommand{\cEb}{\check{\underline{{\cal E}}}}
\newcommand{\cF}{\check{{\cal F}}}
\newcommand{\cG}{\check{{\cal G}}}

\newcommand{\chF}{\check{F}}
\newcommand{\chG}{\check{G}}
\newcommand{\chFb}{\check{\Fb}}
\newcommand{\chGb}{\check{\Gb}}

\newcommand{\cB}{{\cal B}}
\newcommand{\cBb}{\underline{{\cal B}}}
\newcommand{\cA}{{\cal A}} 

\newcommand{\Sb}{\underline{S}}
\newcommand{\taub}{\underline{\tau}}
\newcommand{\yb}{\underline{y}}
\newcommand{\Pib}{\underline{\Pi}}

\newcommand{\cf}{\check{f}}

\newcommand{\chf}{\check{\hat{f}}}
\newcommand{\cv}{\check{v}}
\newcommand{\cga}{\check{\gamma}}

\newcommand{\ctchi}{\check{\tchi}}
\newcommand{\ctchib}{\check{\tchib}}
\newcommand{\cla}{\check{\lambda}}
\newcommand{\clab}{\check{\lambdab}}
\newcommand{\cmu}{\check{\mu}}
\newcommand{\cmub}{\check{\mub}}

\newcommand{\cxi}{\check{\xi}}
\newcommand{\csxi}{\check{\sxi}}

\newcommand{\Ga}{\Gamma}
\newcommand{\La}{\Lambda}
\newcommand{\oLa}{\overline{\La}}
\newcommand{\Lab}{\underline{\La}}
\newcommand{\oLab}{\underline{\oLa}}
\newcommand{\oX}{\overline{X}}
\newcommand{\Xb}{\underline{X}}
\newcommand{\oXb}{\underline{\oX}}

\newcommand{\xb}{\underline{x}}

\newcommand{\hdl}{\hat{\delta}}
\newcommand{\cdl}{\check{\delta}}
\newcommand{\chdl}{\check{\hat{\delta}}} 

\newcommand{\cO}{{\cal O}}
\newcommand{\ocO}{\overline{\cO}}

\newcommand{\ou}{\overline{u}}
\newcommand{\oub}{\overline{\ub}}

\newcommand{\up}[1]{\stackrel{\circ}{#1}}

\thispagestyle{empty}
\begin{center}
{\Huge \bf \sc The Shock Development Problem}
\end{center}
\bigskip
\begin{center}
{\Large \bf \sc  Demetrios Christodoulou}
\end{center}
\vspace{10mm}
\begin{center}
{\Large May 1, 2017}

\vspace{.10in}

\end{center}
\vfill

\section*{\huge{Contents}}

\vspace{10mm}

\noindent {\bf Prologue}

\vspace{5mm}

\noindent {\bf Chapter 1 : Fluid Mechanics and the Shock Development \\Problem}\\
1.1 The general equations of motion\\
1.2 The irrotational case and the nonlinear wave equation\\
1.3 The non-relativistic limit\\
1.4 The jump conditions\\
1.5 The shock development problem\\
1.6 The restricted shock development problem 

\vspace{5mm}

\noindent {\bf Chapter 2 : The Geometric Construction}\\
2.1 A general construction in Lorentzian geometry\\
2.2 The characteristic system\\ 
2.3 The wave system\\
2.4 Variations by translations and the wave equation for the rectangular \\
components of $\beta$\\ 
2.5 The geometric construction for the shock development problem 

\vspace{5mm}

\noindent {\bf Chapter 3 : The Acoustical Structure Equations}\\ 
3.1 The connection coefficients and the first variation equations\\ 
3.2 The structure functions and the formulas for the torsion forms\\
3.3 The propagation equations for $\lambda$ and $\lambdab$\\
3.4 The second variation and cross variation equations\\
3.5 The case $n=2$\\
3.6 The Codazzi and Gauss equations ($n>2$)

\vspace{5mm}

\noindent {\bf Chapter 4 : The Problem of the Free Boundary}\\
4.1 Analysis of the boundary conditions\\
4.2 The transformation functions and the identification equations\\
4.3 Regularization of the identification equations

\vspace{5mm}

\noindent {\bf Chapter 5 : The Initial and Derived Data}\\
5.1 The propagation equations for $\lambdab$ and $s_{NL}$ on $\underline{{\cal C}}$\\ 
5.2 The propagation equations for the higher order derived data \\
$T^m\lambdab$ and $T^m s_{NL}$, $m\geq 1$, on $\underline{{\cal C}}$\\ 
5.3 The boundary conditions for the higher order derived data  
and the \\ determination of the $T$ derivatives  
of the transformation functions on $\partial_-{\cal B}$

\vspace{5mm}

\noindent {\bf Chapter 6 : The Variation Fields}\\
6.1 The bi-variational stress\\
6.2 The variation fields $V$ and the associated 1-forms $\s^{(V)}\theta^\mu$\\
6.3 The fundamental energy identities\\
6.4 The boundary condition on ${\cal K}$ for the 1-forms $\s^{(V)}\xi$

\vspace{5mm}

\noindent {\bf Chapter 7 : The Multiplier Field}\\ 
7.1 Coercivity at the boundary. The choice of multiplier field\\
7.2 The deformation tensor of the multiplier field. The error integral \\
associated to $\s^{(V)}Q_1$ \\
7.3 The error integral associated to to $\s^{(V)}Q_2$

\vspace{5mm}

\noindent {\bf Chapter 8 : The Commutation Fields}\\
8.1 The commutation fields and the higher order variations\\
8.2 The recursion formulas for the source functions\\
8.3 The deformation tensors of the commutation fields\\
8.4 The principal acoustical error terms

\vspace{5mm}

\noindent {\bf Chapter 9 : The Power Series Approximation Method}\\
9.1 Setup of the truncated power series\\
9.2 Estimates for the quantities by which the $N$th approximants fail to satisfy the characteristic and wave systems\\
9.3 Estimates for the quantities by which the $N$th approximants fail to satisfy the boundary 
conditions\\
9.4 Estimates for the quantities by which the $N$th approximants fail to satisfy the identification equations\\
9.5 Estimates for the quantity by which the $\beta_{\mu,N}$ fail to satisfy the wave equation 
relative to $\tilde{h}_N$ and to $\tilde{h}^\prime_N$\\
9.6 The variation differences $\s^{(m,l)}\check{\dot{\phi}}_\mu$ and the rescaled source differences 
$\s^{(m,l)}\check{\tilde{\rho}}_\mu$\\
9.7 . The difference 1-forms $\s^{(V;m,l)}\check{\xi}$. The difference energies and the difference energy identities

\vspace{5mm}

\noindent {\bf Chapter 10 : The Top Order Acoustical Estimates in the Case $d=2$}\\
10.1 Regularization of the propagation equations for $\tchi$ and $\tchib$\\
10.2 Regularization of the propagation equations for $E^2\lambda$ and $E^2\lambdab$\\
10.3 The structure equations for the $N$th approximants\\
10.4 The propagation equations for $\cth_l$, $\cthb_l$ and for $\cnu_{m,l}$, $\cnub_{m,l}$\\
10.5 Estimates for $\cth_l$\\
10.6 Estimates for $\cnu_{m-1,l+1}$\\
10.7 Estimates for $\cthb_l$ and $\cnub_{m-1,l+1}$ in terms of their boundary values on ${\cal K}$\\
10.8 Boundary conditions on ${\cal K}$ and preliminary estimates for $\cthb_l$ and $\cnub_{m-1,l+1}$ 
on ${\cal K}$

\vspace{5mm} 

\noindent {\bf Chapter 11 : Outline of the Top Order Acoustical Estimates for more than 2 
Spatial Dimensions}\\

\vspace{5mm}

\noindent {\bf Chapter 12 : The Top Order Estimates for the Transformation Functions and the 
Next to Top Order Acoustical Estimates}\\
12.1 The propagation equations for the next to top order acoustical difference quantities 
$(\s^{(n-1)}\ctchi, \s^{(n-1)}\ctchib)$ and $(\s^{(m,n-m)}\cla, \s^{(m,n-m)}\clab)$ : $m=0,...,n$\\
12.2 Estimates for $(\s^{(n-1)}\ctchi, \s^{(n-1)}\ctchib)$ and $(\s^{(0,n)}\cla, \s^{(0,n)}\clab)$\\
12.3 Estimates for $(T\Omega^n\chf, T\Omega^n\cv, T\Omega^n\cga)$\\
12.4 Estimates for $(\s^{(m,n-m)}\cla, \s^{(m,n-m)}\clab)$, for $m=1,...,n$\\
12.5 Estimates for $(T^{m+1}\Omega^{n-m}\chf, T^{m+1}\Omega^{n-m}\cv, T^{m+1}\Omega^{n-m}\cga)$, 
for $m=1,...,n$\\
12.6 Estimates for $(\Omega^{n+1}\chf, \Omega^{n+1}\cv, \Omega^{n+1}\cga)$ 

\vspace{5mm}

\noindent {\bf Chapter 13 : The Top Order Energy Estimates}\\
13.1 Estimates for $\s^{(V;m,n-m)}\check{b}$\\
13.2 The borderline error integrals contributed by $\s^{(V;m,n-m)}\check{Q}_1$, 
$\s^{(V;m,n-m)}\check{Q}_2$\\
13.3 The borderline error integrals associated to $\cth_n$ and to $\cnu_{m-1,n-m+1} : m=1,...,n$\\
13.4 The borderline error integrals associated to $\cthb_n$ and to $\cnub_{m-1,n-m+1} : m=1,...,n$\\
13.5 The top order energy estimates

\vspace{5mm} 

\noindent {\bf Chapter 14 : Lower Order Estimates, Recovery of the Bootstrap Assumptions,  
and Completion of the Argument}\\ 
14.1 The bootstrap assumptions needed\\
14.2 $L^2(S_{\ub,u})$ estimates for $\s^{(n-1)}\ctchi$ and for $\s^{(m,n-m)}\cla \ : \ m=0,...,n-1$\\
14.3 $L^2({\cal K}_\sigma^\tau)$ estimates for the $n$th order acoustical differences\\
14.4 $L^2(S_{\ub,u})$ estimates for the $n$th order variation differences\\
14.5 $L^2(S_{\ub,u})$ estimates for the $n-1$th order acoustical differences\\
14.6 $L^2(S_{\ub,u})$ estimates for all $n$th order derivatives of the $\beta_\mu$\\
14.7 $L^2(S_{\ub,u})$ estimates for $\Omega^{n-1}\log\sh$ and $\Omega^{n-1}b$\\
14.8 Lower order $L^2(S_{\ub,u})$ estimates\\
14.9 Pointwise estimates and recovery of the bootstrap assumptions\\
14.10 Completion of the argument

\vspace{5mm}

\noindent {\bf Bibliography}

\pagebreak 

\section*{\huge{Prologue}}

The subject of this monograph is the shock development problem in fluid mechanics.  
This problem is formulated in the framework of the Eulerian equations of a compressible
perfect fluid as completed by the laws of thermodynamics. These equations express the 
differential conservation laws of mass, momentum and energy and constitute a 
quasilinear hyperbolic 1st order system for the physical variables, that is the fluid velocity and 
the two positive quantities corresponding to a local thermodynamic equilibrium state. 
Smooth initial data for this system of equations leads to the formation of a surface 
in spacetime where the derivatives of the physical quantities with respect to the standard 
rectangular coordinates blow up. Now, there is a mathematical notion of maximal development of 
initial data. As was first shown in the monograph [Ch-S], this maximal development ends 
at a future boundary which consists of a regular part $\underline{{\cal C}}$ and a singular 
part ${\cal B}$ with a common past boundary $\partial_-{\cal B}$, the surface just mentioned. 
A solution of the Eulerian equations in a given spacetime domain defines a cone field on 
this domain, 
the sound cones. This defines a causal structure on the spacetime domain, equivalent to 
a conformal class of Lorentzian metrics, the acoustical causal structure. Relative to this 
structure $\partial_-{\cal B}$ is a spacelike surface, while $\underline{{\cal C}}$ is a null hypersurface. Also ${\cal B}$ is in this sense a null hypersurface, however being singular, 
while its intrinsic geometry is that of a null hypersurface, its extrinsic geometry is that of 
a spacelike hypersurface, for, the past null geodesic cone in the spacetime manifold of a point on 
${\cal B}$ does not intersect ${\cal B}$. The character of ${\cal B}$ and the behavior of the 
the solution at ${\cal B}$ were described in detail in [Ch-S] by means of the introduction 
of a class of coordinates such that the rectangular coordinates as well as the physical 
variables are smooth functions of the new coordinates up to ${\cal B}$, but the Jacobian 
of the transformation to the new coordinates, while strictly negative in the past of ${\cal B}$, vanishes 
at ${\cal B}$ itself, a fact which characterizes ${\cal B}$. Now, the mathematical notion of 
maximal development of initial data, while physically correct up to 
$\underline{{\cal C}}\bigcup\partial_-{\cal B}$ is not physically correct up to ${\cal B}$. 
The problem of the physical continuation of the solution is set up in the Epilogue of [Ch-S] as 
the shock development problem. In this problem one is required to construct a hypersurface 
of discontinuity ${\cal K}$ , the shock hypersurface, lying in the past of ${\cal B}$ but having 
the same past boundary as the latter, namely $\partial_-{\cal B}$, and a solution of the Eulerian 
equations in the spacetime domain bounded in the past by $\underline{{\cal C}}\bigcup{\cal K}$, 
agreeing on $\underline{{\cal C}}$ with the data induced by the maximal development, while having 
jumps across ${\cal K}$ relative to the data induced on ${\cal K}$ by the maximal development, 
these jumps satisfying the jump conditions which follow from the integral form of the mass, momentum 
and energy conservation laws. Moreover, ${\cal K}$ is required to be a spacelike hypersurface 
relative to the acoustical structure corresponding to the prior solution and a timelike hypersurface  relative to the acoustical structure corresponding to the new solution (see Sections 1.4, 1.5). 
Thus, the singular surface $\partial_-{\cal B}$ is the cause generating the shock hypersurface 
${\cal K}$. 

The monograph [Ch-S] actually considered the extension of the Eulerian equations to the framework 
of special relativity. In this framework the underlying geometric structure of the spacetime manifold is 
that of the Minkowski spacetime of special relativity. On the other hand, the underlying 
spacetime structure of the original Eulerian equations, as of all of classical mechanics, is that 
of the Galilean spacetime. The monograph [Ch-Mi] treated the same topics as [Ch-S] reaching similar results in a considerably simpler, self-contained manner. The relationship of the non-relativistic to the 
relativistic theory, how results in the former are deduced as limits of results in the latter, is 
discussed in detail in Sections 1.3 - 1.6. The Galilean structure has a distinguished family of hyperplanes, those of 
absolute simultaneity, while the Minkowski structure is that of a flat Lorentzian manifold, thus 
the latter, with its light cone field, is in a sense more like that corresponding to the sound cone field defined by a solution of the 
Eulerian equations. What we mean by rectangular coordinates in the Minkowskian framework is the 
standard geometric notion and two such systems of coordinates are related by a transformation 
belonging to the Poincar\'{e} group. On the other hand, by rectangular coordinates in the 
Galilean framework we mean a Galilei frame together with rectangular coordinates in Euclidean space 
and two such systems of coordinates are related by the Galilei group, which extends 
the Euclidean group. 

From the mathematical point of view the shock development problem is a free boundary problem, 
with nonlinear conditions at the free boundary ${\cal K}$, for a quasilinear hyperbolic 1st order 
system, with characteristic initial data on $\underline{{\cal C}}$ which are singular, in a 
prescribed manner, at $\partial_-{\cal B}$, the past boundary of $\underline{{\cal C}}$. 
It will be shown that 
the singularity persists, not only as a discontinuity in the physical variables across ${\cal K}$, 
but also as a milder singularity propagating along $\underline{{\cal C}}$. 
While the physical variables and their 1st derivatives extend continuously 
across $\underline{{\cal C}}$, 
the 1st derivatives are only $C^{0,1/2}$ at $\underline{{\cal C}}$ from the point of view of the 
future solution. The spherically symmetric barotropic case of this problem was recently solved in 
[Ch-Li].

What the present monograph solves is not the general shock development problem but what we call 
the restricted shock development problem. This is defined in Section 1.6 and results if we neglect 
the variation along ${\cal K}$ of $\triangle s$, the jump in entropy. As explained there, due to 
the fact that in terms of the coordinate $\tau$ on ${\cal K}$, which is introduced in Section 2.5 
and which vanishes at $\partial_-{\cal B}$, $\triangle s=O(\tau^3)$, the restricted problem retains 
the main difficulty of the general problem, namely the singular behavior as $\tau\rightarrow 0$. 
Moreover, our treatment of the restricted problem is amenable to generalization, being based on 
the 1-form $\beta$ introduced in Section 1.1. In the general case this 1-form is not closed, but the 
2-form $\omega=-d\beta$, the spacetime vorticity 2-form, satisfies the equation
\begin{equation}
i_u\omega=\theta ds
\label{p.1}
\end{equation}
where $u$ is the spacetime fluid velocity, and $\theta$ the temperature. As a consequence 
$\omega$ satisfies the transport equation:
\begin{equation}
{\cal L}_u\omega=d\theta\wedge ds
\label{p.2}
\end{equation}
where the right hand side, like $\omega$ itself, is of differential order 1, the physical variables 
themselves being of differential order 0. The vectorfield $u$ then defines a characteristic field
complementing the sound cone field. The hypersurface generated by the integral curves of $u$ 
initialing at $\partial_-{\cal B}$ divides the domain of the new solution into two subdomains, 
that bounded by $\underline{{\cal C}}$ and the hypersurface in question and that between the hypersurface 
and ${\cal K}$. In the former domain the solution coincides with that of the restricted problem, 
while in the latter domain the spacetime vorticity 2-form $\omega$ does not vanish being determined 
$\triangle s$ and \ref{p.2} through \ref{1.227}, \ref{1.232}. There is then an additional 
3rd order discontinuity in $\omega$ across the dividing hypersurface. 

We shall now describe in brief the mathematical methods introduced in this monograph noting the 
essential points and where the corresponding material is found in the monograph. What follows 
will serve as a guide for the reader. A central role 
is played by the reformulation in Chapter 2 of the Eulerian equations in the domain ${\cal N}$ of the new 
solution. First a diffeomorphism is defined of this domain onto 
\begin{equation}
{\cal R}_{\delta,\delta}=R_{\delta,\delta}\times S^{n-1}
\label{p.a1}
\end{equation}
$n$ being the spatial dimension and:
\begin{equation}
R_{\delta,\delta}=\{(\ub,u)\in\mathbb{R}^2 \ : \ 0\leq\ub\leq u\leq\delta\}
\label{p.3}
\end{equation}
being a domain in $\mathbb{R}^2$, which represents the range in ${\cal N}$ 
of two functions $\ub$ and $u$ the level sets of which are transversal acoustically null hypersurfaces 
denoted by $\Cb_{\ub}$ and $C_u$ respectively, with $\Cb_0=\underline{{\cal C}}$.  In this 
representation the shock hypersuface boundary of ${\cal N}$ is: 
\begin{equation}
{\cal K}^\delta=\{(\tau,\tau) \ : \ \tau\in [0,\delta]\}\times S^{n-1}
\label{p.4}
\end{equation}
and $\partial_-{\cal B}$ is: 
\begin{equation}
\partial_-{\cal B}=(0,0)\times S^{n-1}=S_{0,0}
\label{p.5}
\end{equation}
We denote by $S_{\ub,u}$ 
the surfaces: 
\begin{equation}
S_{\ub,u}=\Cb_{\ub}\bigcap C_u=(\ub,u)\times S^{n-1}
\label{p.6}
\end{equation}
In the reformulation of the Eulerian equations the unknowns constitute a triplet 
$((x^\mu : \mu=0,...,n) , \ b \ , (\beta_\mu : \mu=0,...,n))$, where the $(x^\mu : \mu=0,...,n)$ 
are functions on $R_{\delta,\delta}\times S^{n-1}$ representing rectangular coordinates in 
the corresponding domain in Galilean spacetime in the non-relativistic theory, in Minkowski 
spacetime in the relativistic theory, with $x^0=t$, and the $(\beta_\mu : \mu=0,...,n)$ are 
also functions on $R_{\delta,\delta}\times S^{n-1}$ and represent the rectangular components 
of the 1-form $\beta$. The unknown $b$ is a mapping of $R_{\delta,\delta}$ into the space of 
vectorfields on $S^{n-1}$. The pair $((x^\mu : \mu=0,...,n) , \ b)$ satisfies the 
{\em characteristic system}, a fully nonlinear 1st order system of partial differential equations. 
The $(\beta_\mu : \mu=0,...,n)$ satisfy the {\em wave system}, a quasilinear 1st order system 
of partial differential equations. The two systems are coupled through the 
$(h_{\mu\nu} : \mu, \nu = 0,...,n)$, which represent the rectangular components of the acoustical 
metric, and depend on the $(\beta_\mu : \mu=0,...,n)$. More precisely, 
denoting by $\sd f$ the differential of a function $f$ on the $S_{\ub,u}$, 
the coupling is through 
the functions $N^\mu : \mu=0,...,n$ and $\Nb^\mu : \mu=0,...,n$ defined in terms of the 
$\sd x^\mu : \mu=0,...,n$ pointwise by the conditions:
\begin{equation}
h_{\mu\nu}N^\mu\sd x^\nu=0, \ \ h_{\mu\nu}N^\mu N^\nu=0, \ \ N^0=1
\label{p.7}
\end{equation}
and similarly for the $\Nb^\mu$. By reason of the quadratic nature of the 2nd of these conditions 
a unique pair $(N^\mu : \mu=0,...,n), (\Nb^\mu : \mu=0,...,n))$ up to exchange is pointwise defined 
by $(h_{\mu\nu} : \mu, \nu = 0,...,n)$ and $(\sd x^\mu : \mu=0,...,n)$. The vectorfields $N$, $\Nb$ 
with rectangular components 
$(N^\mu : \mu=0,...,n)$, $(\Nb^\mu : \mu=0,...,n)$ 
are then null normal fields, relative to the acoustical metric $h$, to the surfaces $S_{\ub,u}$. 
The exchange ambiguity is removed by requiring $N$ to be tangential to the $C_u$, 
$\Nb$ to be tangential to the $\Cb_{\ub}$. 
The function
\begin{equation}
c=-\frac{1}{2}h_{\mu\nu}N^\mu\Nb^\nu
\label{p.8}
\end{equation}
is then bounded from below by a positive constant. Defining: 
\begin{equation}
L=\frac{\partial}{\partial\ub}-b, \ \ \Lb=\frac{\partial}{\partial u}+b
\label{p.9}
\end{equation}
and:
\begin{equation}
\rho=Lt, \ \ \rhob=\Lb t
\label{p.10}
\end{equation}
the characteristic system is simply:
\begin{equation}
Lx^\mu=\rho N^\mu,  \ \ \Lb x^\mu=\rhob\Nb^\mu \ \ \ : \ \mu=0,...,n
\label{p.11}
\end{equation} 

What is achieved by the reformulation just described 
is a regularization of the problem. That is, we are now seeking smooth functions on 
$R_{\delta,\delta}\times S^{n-1}$ satisfying the coupled system, the initial data themselves 
being represented by smooth functions on $\Cb_0$. The Jacobian of the transformation 
representing the mapping $((\ub,u),\vartheta)\mapsto (x^\mu((\ub,u),\vartheta) : \mu=0,...,n)$, 
$\vartheta\in S^{n-1}$, is of the form:
\begin{equation}
\frac{\partial(x^0,x^1,...,x^n)}{\partial(\ub,u,\vartheta^1,...,\vartheta^{n-1})}=\rho\rhob d
\label{p.12}
\end{equation}
$(\vartheta^A : A=1,...,n-1)$ being local coordinates on $S^{n-1}$. Here:
\begin{equation}
\rho=Lt, \ \ \rhob=\Lb t, \ \mbox{where} \ 
L=\frac{\partial}{\partial\ub}-b, \ \ \Lb=\frac{\partial}{\partial u}+b
\label{p.13}
\end{equation}
are non-negative functions, the inverse temporal density of the foliation of spacetime by the 
$\Cb_{\ub}$ as measured along the generators of the $C_u$, the inverse temporal density of the 
foliation of spacetime by the $C_u$ as measured along the generators of the $\Cb_{\ub}$, respectively, 
and $d$ is a function bounded from above by a negative constant (see Section 2.2). 
As a consequence, the Jacobian \ref{p.12} vanishes where and only where one of $\rho$, $\rhob$ 
vanishes. 
The function $\rhob$ is given on $\Cb_0$ by the initial data and while positive on 
$\Cb_0\setminus S_{0,0}$, vanishes to 1st order at $S_{0,0}$, the last being a manifestation 
of the singular nature of the surface $\partial_-{\cal B}$. The function $\rho$ on $\Cb_0$ 
represents 1st derived data on $\Cb_0$ (see Chapter 5) and vanishes there to 0th order. 
It turns out that these are the only places where $\rhob$, $\rho$ vanish in 
$R_{\delta,\delta}\times S^{n-1}$. A smooth solution of the coupled characteristic and wave systems 
once obtained then represents a solution of the original Eulerian equations in standard rectangular 
coordinates which is smooth in ${\cal N}\setminus\underline{{\cal C}}$ but singular at 
$\underline{{\cal C}}$ 
with the transversal derivatives of the $\beta_\mu$ being only H\"{o}lder continuous of exponent 
1/2 at  $\underline{{\cal C}}$ and, in addition, a stronger singularity at $\partial_-{\cal B}$,  
namely the blow up of the derivatives of the $\beta_\mu$ at $\partial_-{\cal B}$ in the direction 
tangential to $\underline{{\cal C}}$ but transversal to $\partial_-{\cal B}$. In particular, the 
new solution is smooth at the shock hypersurface ${\cal K}$ except at its past boundary 
$\partial_-{\cal B}$, as is the prior solution which holds in the other side of ${\cal K}$, the 
past side, ${\cal K}\setminus\partial_-{\cal B}$ lying in the interior of the domain of 
the maximal development. 

The wave system consists of the equations 
\begin{equation}
dx^\mu\wedge d\beta_\mu=0, \ \ \ h^{-1}(dx^\mu,d\beta_\mu)=0
\label{p.a2}
\end{equation}
expressed in terms of the representation \ref{p.a1}, that is:
\begin{eqnarray}
&&(Lx^\mu)\Lb\beta_\mu-(\Lb x^\mu)L\beta_\mu=0 \nonumber\\
&&(Lx^\mu)\sd\beta_\mu-(\sd x^\mu)L\beta_\mu=0 \nonumber\\
&&(\Lb x^\mu)\sd\beta_\mu-(\sd x^\mu)\Lb\beta_\mu=0 \nonumber\\
&&\hspace{10mm}\sd x^\mu \wedge \sd\beta_\mu=0, \nonumber\\
&&\frac{1}{2}\left\{(Lx^\mu)\Lb\beta_\mu+(\Lb x^\mu)L\beta_\mu\right\}=a(\sd x^\mu, \sd\beta_\mu)_{\sh}
\label{p.a3}
\end{eqnarray}
where:
\begin{equation}
a=c\rho\rhob
\label{p.a4}
\end{equation}

The characteristic system for the pair $((x^\mu : \mu=0,...,n) , \ b)$ together with the 
$(h_{\mu\nu} : \mu, \nu = 0,...,n)$ which enter this system through the $(N^\mu : \mu=0,...,n)$ 
and the $(\Nb^\mu : \mu=0,...,n)$, manifest a new kind of differential 
geometric structure which involves the interaction of two geometric structures on the same 
underlying manifold, the first of these structures being the background Galilean structure 
in the case of the non-relativistic theory, 
the background Minkowskian structure in the case of the relativistic theory, 
and the other being the Lorentzian geometry deriving from the acoustical metric. As for the 
$(\beta_\mu : \mu=0,...,n)$ of the wave system, this is the set of functions obtained by evaluating 
the 1-form $\beta$ on the set of translation fields of the background structure. Proposition 2.1 
with a translation field substituted for $X$ plays a central role in our approach. The 
proposition asserts that if $X$ is a vectorfield generating an isometries of the background structure 
then:
\begin{equation}
\square_{\tilde{h}}\beta(X)=0
\label{p.14}
\end{equation}
where $\tilde{h}$ is a metric in the conformal class of the acoustical metric $h$. 

The derivatives of $((x^\mu : \mu=0,...,n) , \ b)$ are controlled through the acoustical structure 
equations, the subject of Chapter 3. These are differential consequences of the characteristic system,
bringing out more fully the interaction of the two geometric structures. 
We have the induced metric $\sh$ on the surfaces $S_{\ub,u}$ and the functions: 
\begin{equation}
\lambda=c\rhob, \ \ \ \lambdab=c\rho
\label{p.15}
\end{equation}
(see \ref{p.8}) introduced already in Chapter 2. 
While $\sh$ refers only to the acoustical structure, the functions $\lambda$, $\lambdab$, involve 
the interaction of the two geometric structures. These are acoustical quantities of 0th order. 
The 1st variation equations (Proposition 3.1) 
express $\sL_L\sh$, $\sL_{\Lb}\sh$ in terms of $\chi$, $\chib$, the two 2nd fundamental forms 
of $S_{\ub,u}$. The torsion forms $\eta$, $\etab$ represent the connection in the normal bundle 
of $S_{\ub,u}$ in terms of the vectorfields $L$, $\Lb$ which along $S_{\ub,u}$ constitute 
sections of this bundle. The commutator: 
\begin{equation}
[L,\Lb]=\sL_T b, \ \  T=L+\Lb=\frac{\partial}{\partial\ub}+\frac{\partial}{\partial u}
\label{p.16}
\end{equation}
(see \ref{p.9}) is expressed in terms of $\eta$, $\etab$. While the quantities 
$\chi$, $\chib$, $\eta$, $\etab$ 
refer only to the acoustical structure, the structure equations assume a non-singular form 
only in terms of the  quantities $\tchi$, $\tchib$, $\teta$, $\tetab$ which involve both 
structures. The former are related to the latter as follows: up to an order 1 remainders depending 
only on the $d\beta_\mu : \mu=0,...,n$,  $\chi$ is equal to $\rho\tchi$, 
$\chib$ is equal to $\rhob\tchib$,  $\eta$ is equal to $\rho\teta$, 
$\etab$ is equal to $\rhob\tetab$. In fact, the only order 1 
quantities appearing in the 1st and 3rd remainders are the $L\beta_\mu : \mu=0,...,n$, and 
the only order 1 quantities appearing in the 2nd and 4th remainders are the 
$\Lb\beta_\mu : \mu=0,...,n$. 
Moreover $\teta$, $\tetab$, are expressed in terms of $\sd\lambda$, $\sd\lambdab$ (see Section 3.2). 
Thus $\lambda$, $\lambdab$, $\tchi$, $\tchib$ are the primary acoustical quantities, the first two 
being of 0th order, the second two of 1st order. 

Proposition 3.3, the propagation equations for $\lambda$ and $\lambdab$, plays a central role. 
This proposition expresses $L\lambda$ and $\Lb\lambdab$ in terms of 1st order quantities with 
vanishing 1st order  acoustical part. Proposition 3.4, the 2nd variation equations, also plays a 
central role. This proposition likewise expresses $\sL_L\tchi$ and $\sL_{\Lb}\tchib$ in terms 
of 2nd order quantities with vanishing 2nd order acoustical part. The acoustical structure 
equations are completed by the Codazzi and Gauss equations of Propositions 3.6 and 3.7. 
Denoting by $\sD$ the covariant derivative operator of $(S_{\ub,u},\sh)$, $\sD\tchi$, $\sD\tchib$ 
are 3-covariant tensorfields on each $S_{\ub,u}$ and the Codazzi equations express the 
3-covariant tersorfields obtained from $\sD\tchi$, $\sD\tchib$ by antisymmetrizing with respect 
to the first two entries in terms of 2nd order quantities with vanishing 2nd order acoustical 
part. The Gauss equation likewise expresses the curvature of $(S_{\ub,u},\sh)$ in terms of a  2nd order 
quantity with vanishing 2nd order acoustical part. The Codazzi and Gauss equations trivialize 
in the case of $n=2$ spatial dimensions.

The boundary conditions on ${\cal K}$ are analyzed in Chapter 4. The jumps 
$\triangle\beta_\mu : \mu=0,...,n$ in the rectangular components of the 1-form $\beta$ 
across ${\cal K}$ are subject to two conditions, one of which is linear and the other nonlinear. 
The linear jump condition decomposes into the two conditions:
\begin{equation}
\sd x^\mu \ \triangle\beta_\mu=0, \ \ \ \ T^\mu\triangle\beta_\mu=0
\label{p.17}
\end{equation}
As a consequence of the 1st of \ref{p.17} $\triangle\beta_\mu$ can be expressed as a linear 
combination of the components:
\begin{equation}
\ep=N^\mu\triangle\beta_\mu, \ \ \epb=\Nb^\mu\triangle\beta_\mu 
\label{p.18}
\end{equation}
We denote by $r$ the ratio:
\begin{equation}
r=-\frac{\ep}{\epb}
\label{p.19}
\end{equation}
In  reference to the 2nd of \ref{p.17}, $T^\mu=Tx^\mu$ are the rectangular components of the vectorfield 
$T$ and are given, in view of \ref{p.11}, \ref{p.16} by:
\begin{equation}
T^\mu=\rho N^\mu+\rhob\Nb^\mu 
\label{p.20}
\end{equation}
Then in view of \ref{p.15} the 2nd of \ref{p.18} is equivalent to the following boundary condition 
for $\lambdab$:
\begin{equation}
r\lambdab=\lambda \ \ \mbox{: on ${\cal K}$}
\label{p.21}
\end{equation}
We remark here that the propagation equations for $\lambda$ and for $\tchi$ are supplemented by 
initial conditions on $\Cb_0$, while the propagation equations for $\lambdab$ and for $\tchib$ 
are supplemented by boundary conditions on ${\cal K}$. The boundary condition for $\tchib$ is 
first derived in Chapter 10 in the case of 2 spatial dimensions and afterwards, in Chapter 11, 
in general. The derivation is by applying $\sD$ to the 1st of \ref{p.17}. The boundary condition
for $\tchib$ takes the form of a relation between $r\tchib$ and $\tchi$ on ${\cal K}$ analogous 
to \ref{p.21} (see \ref{11.b6}). The nonlinear jump condition takes the form of a relation between 
$\ep$ and $\epb$ which, in the setting of the shock development problem, is shown 
to be equivalent to: 
\begin{equation}
\ep=-j(\epb)\epb^2, \ \ \mbox{hence} \ \ r=j(\epb)\epb
\label{p.22}
\end{equation}
where $j$ is a smooth function (see Proposition 4.2). 

The jump $\triangle\beta_\mu$ is a function on ${\cal K}$, which 
at a given point on ${\cal K}$ represents the difference of $\beta_\mu$, defined by the new 
solution which holds in the future of ${\cal K}$, at the point, from the corresponding 
quantity, which we denote by $\beta^\prime_\mu$, for the prior solution which holds in the past 
of ${\cal K}$, at the same point in the background spacetime, namely the same point in Galilei 
spacetime in the non-relativistic, the same point in Minkowski spacetime in the relativistic theory. From [Ch-Mi], [Ch-S], the $\beta^\prime_\mu$ are smooth functions of the coordinates 
$(t,u^\prime,\vartheta^\prime)$ where the level sets $C^\prime_{u^\prime}$ are outgoing acoustically 
null hypersurfaces, and $\vartheta^\prime$, the range of which is $S^{n-1}$, is constant along each 
generator of each $C^\prime_{u^\prime}$, and the restriction of $\vartheta^\prime$ to each 
surface $S^\prime_{t,u^\prime}=\Sigma_t\bigcap C^\prime_{u^\prime}$ is a diffeomorphism onto $S^{n-1}$. 
Here we denote by $\Sigma_t$ the level sets of $t$, the hyperplanes of absolute simultaneity 
in the non-relativistic theory, parallel spacelike hyperplanes in Minkowski spacetime in the 
relativistic theory. The spatial rectangular coordinates are likewise 
smooth functions $x^{\prime i}(t,u^\prime,\vartheta^\prime) : i=1,...,n$ on the domain of the maximal development up to 
its boundary. Now according to our construction in Section 2.5, for $c>0$, $C_c$, the $c$-level set 
of $u$ in ${\cal N}$ is the extension across $\underline{{\cal C}}=\Cb_0$ of $C^\prime_c$, the 
$c$-level set of $u^\prime$ in the domain of the prior maximal development. At 
$\partial_-{\cal B}=S_{0,0}$, $u^\prime$ vanishes and $\vartheta^\prime=\vartheta$. On the other hand on 
${\cal K}\setminus S_{0,0}$,  $u^\prime<0$. Therefore $u^\prime$ and $\vartheta^\prime$ induced on 
on ${\cal K}$ by the prior solution are of the form (see \ref{p.4}):
\begin{equation}
u^\prime=w(\tau,\vartheta), \ \ \ \vartheta^\prime=\psi(\tau,\vartheta)
\label{p.23}
\end{equation}
and we have:
\begin{equation}
w(0,\vartheta)=0, \  w(\tau,\vartheta)<0 \ \mbox{: for $\tau>0$}, \ \ \ \ \psi(0,\vartheta)=\vartheta
\label{p.24}
\end{equation}
We denote, in reference to the new solution,
\begin{equation}
f(\tau,\vartheta)=x^0(\tau,\tau,\vartheta), \ \ g^i(\tau,\vartheta)= x^i(\tau,\tau,\vartheta)
 \ : \ i=1,...,n
\label{p.25}
\end{equation}
Then the point $(\tau,\tau,\vartheta)\in{\cal K}$ in terms of the coordinates 
$(\ub,u,\vartheta)$ represents the same point in the background spacetime as the point 
$(t,u^\prime,\vartheta^\prime)$ in terms of the coordinates $(t,u^\prime,\vartheta^\prime)$ if 
with \ref{p.23} we have:
\begin{equation}
x^{\prime i}(f(\tau,\vartheta), w(\tau,\vartheta), \psi(\tau,\vartheta))=g^i(\tau,\vartheta) 
\ : \ i=1,...,n
\label{p.26}
\end{equation}
We call $w$ and $\psi$ {\em transformation functions} and equations \ref{p.26} 
{\em identification equations}. Since they do not involve derivatives of $(w,\psi)$, they are to 
determine $(w,\psi)$ pointwise.  Thinking of $S^{n-1}$ as the unit sphere in Euclidean $n$-dimensional 
space, and denoting by $\exp_{\vartheta}$ the associated exponential map 
$T_\vartheta S^{n-1}\rightarrow S^{n-1}$ we can express the 2nd of \ref{p.23} in the form:
\begin{equation}
\vartheta^\prime=\exp_{\vartheta}(\varphi^\prime) 
\label{p.27}
\end{equation}
for some $\varphi^\prime\in T_{\vartheta}S^{n-1}$, which is unique provided that the distance on 
$S^{n-1}$ of $\vartheta^\prime$ from $\vartheta$ is less than $\pi$, as will be the case in our 
problem. Setting then 
\begin{equation}
F^i((\tau,\vartheta),(u^\prime,\varphi^\prime))=
x^{\prime i}(f(\tau,\vartheta), u^\prime,\exp_{\vartheta}(\varphi^\prime))-g^i(\tau,\vartheta)
\label{p.28}
\end{equation}
the identification equations read:
\begin{equation}
F^i((\tau,\vartheta),(u^\prime,\varphi^\prime))=0 \ : \ i=1,...,n
\label{p.29}
\end{equation}
To arrive at a form of these equations to which the implicit function theorem can be directly 
applied we set: 
\begin{equation}
u^\prime=\tau v, \ \ \ \varphi^\prime=\tau^3\gamma
\label{p.30}
\end{equation} 
and define $\hat{F}^i((\tau,\vartheta),(v,\gamma))$ by:
\begin{equation}
F^i((\tau,\vartheta),(\tau v,\tau^3\gamma))=\tau^3\hat{F}^i((\tau,\vartheta),(v,\gamma))
\label{p.31}
\end{equation}
The identification equations then take their regularized form:
\begin{equation}
\hat{F}^i((\tau,\vartheta),(v,\gamma))=0 \ : \ i=1,...,n
\label{p.32}
\end{equation}
given explicitly by Proposition 4.5. 
The implicit function theorem then applies to determine for each $(\tau,\vartheta)\in {\cal K}^{\delta}$ 
the associated $(v,\gamma)$, $v<0$, $\gamma\in T_{\vartheta}S^{n-1}$. 

On $\Cb_0$ $\vartheta$ is defined by the condition it is constant along the generators of $\Cb_0$ 
and together with the condition that the restriction of $\vartheta$ to $S_{0,0}=\partial_-{\cal B}$ 
is a diffeomorphism of $\partial_-{\cal B}$ onto $S^{n-1}$. We then have:
\begin{equation}
\left. b\right|_{\Cb_0}=0
\label{p.33}
\end{equation}
The initial data for the coupled characteristic and wave systems on $\Cb_0=\underline{{\cal C}}$ for the new solution 
then consist of the pair $((x^\mu : \mu=0,...,n), (\beta_\mu : \mu=0,...,n))$ on $\Cb_0$ which is 
that induced by the prior solution. By virtue of the equations of the characteristic and wave systems 
the transversal derivatives on $\Cb_0$ of 1st order of 
$((x^\mu : \mu=0,...,n), (\beta_\mu : \mu=0,...,n))$, that is 
$((Lx^\mu : \mu=0,...,n), (L\beta_\mu : \mu=0,...,n))$, are directly expressed in terms of the initial 
data, with the exception of the pair $(\lambdab, s_{NL})$, which we call 1st derived data. Here:
\begin{equation}
s_{NL}=N^\mu L\beta_\mu 
\label{p.34}
\end{equation}
According to Proposition 5.1 this pair satisfies along the generators of $\Cb_0$ a linear homogeneous 
system of ordinary differential equations. The initial data for this system are at $S_{0,0}$ and vanish 
as shown in Chapter 4. It then follows that the pair $(\lambdab,s_{NL})$ vanishes everywhere along 
$\Cb_0$. Proceeding to the transversal derivatives on $\Cb_0$ of 2nd order, these are all directly expressed in terms of the preceding with the exception of the pair $(T\lambdab,Ts_{NL})$. This pair 
satisfies along the generators of $\Cb_0$ a linear system of ordinary differential equations, in 
principle inhomogeneous but with the same homogeneous part as that in the case of the 1st derived data. 
Due to the vanishing of the 1st derived data the inhomogeneous terms actually vanish. However, the initial data at $S_{0,0}$ for the pair $(T\lambdab,Ts_{NL})$ do not vanish, $\left.T\lambdab\right|_{S_{0,0}}$ 
being represented by a positive function on $S^{n-1}$. It then follows that $T\lambdab$ has a positive 
minimum on $\Cb_0^{\delta}$ provided that $\delta$ is suitably small depending on the initial data. 
Here we denote by $\Cb_0^{\delta}$ the part of $\Cb_0$ corresponding to $u\leq\delta$, that is, the 
part contained in ${\cal R}_{\delta,\delta}$ (see \ref{p.a1}).  In general, the 
transversal derivatives of order $m+1$, are all directly expressed in terms of the preceding with 
the exception of the pair $(T^m\lambdab,T^ms_{NL})$. This pair 
satisfies along the generators of $\Cb_0$ a linear inhomogeneous system of 
ordinary differential equations with the same homogeneous part as that in the case of the 1st derived data and with inhomogeneous part directly expressed in terms of the preceding. The initial data 
for this system consist of the pair $(T^m\lambdab,T^ms_{NL})$ at $S_{0,0}$, to determine which 
requires determining at the same time $(T^{m-1}v,T^{m-1}\gamma)$, which brings in the regularized 
identification equations as well as the boundary condition \ref{p.21} for $\lambdab$. The problem 
is analyzed in Section 5.3. 

Chapter 6 introduces a new geometric concept, that of a variation field, which plays a central role 
in the present work. Let $X$, $Y$ be arbitrary vectorfields on ${\cal N}$. We define the 
{\em bi-variational stress} associated to the 1-form $\beta$ and to the pair $X$, $Y$, to be the 
$T^1_1$ type tensorfield:
\begin{equation}
\dot{T}=\tilde{h}^{-1}\cdot\dot{T}_\flat
\label{p.35}
\end{equation}
where $\dot{T}_\flat$ is the symmetric 2-covariant tensorfield:
\begin{equation}
\dot{T}_\flat=\frac{1}{2}\left(d\beta(X)\otimes d\beta(Y)+d\beta(Y)\otimes d\beta(X)
-(d\beta(X),d\beta(Y))_h h\right)
\label{p.36}
\end{equation}
which depends only on the conformal class of $h$. We then have the identity:
\begin{equation}
\mbox{div}_{\tilde{h}}\dot{T}=\frac{1}{2}\left(\square_{\tilde{h}}\beta(X)\right)d\beta(Y)
+\frac{1}{2}\left(\square_{\tilde{h}}\beta(Y)\right)d\beta(X)
\label{p.37}
\end{equation}
In particular, if $X$, $Y$ generate isometries of the background structure, then by \ref{p.14}: 
\begin{equation}
\mbox{div}_{\tilde{h}}\dot{T}=0
\label{p.38}
\end{equation}
Setting $X$, $Y$ to be the translation fields:
$$X=\frac{\partial}{\partial x^\mu}, \ \ \ Y=\frac{\partial}{\partial x^\nu}$$
we have: 
$$\beta(X)=\beta_\mu, \ \ \ \beta(Y)=\beta_{\nu}$$
and we denote the corresponding bi-variational stress by $\dot{T}_{\mu\nu}$. The identity \ref{p.37} 
takes in this case the form:
\begin{equation}
\mbox{div}_{\tilde{h}}\dot{T}_{\mu\nu}=\frac{1}{2}(\square_{\tilde{h}}\beta_\mu)d\beta_\nu 
+\frac{1}{2}(\square_{\tilde{h}}\beta_\nu)d\beta_\mu 
\label{p.39}
\end{equation}
The concept of bi-variational stress has been introduced in [Ch-A] in the general context of 
Lagrangian theories of mappings of a manifold ${\cal M}$ into another manifold ${\cal N}$, as 
discussed briefly in Section 6.1. The usefulness of the concept of bi-variational stress in 
the context of a free boundary problem is in conjunction with the concept of variation fields. 
A {\em variation field} is here simply a vectorfield $V$ on ${\cal N}$ which along ${\cal K}$ is 
tangential to ${\cal K}$. This can be expanded in terms of the translation fields 
$\partial/\partial x^\mu : \mu=0,...,n$:
\begin{equation}
V=V^\mu\frac{\partial}{\partial x^\mu}
\label{p.40}
\end{equation}
The coefficients $V^\mu : \mu=0,...,n$ of the expansion are simply the rectangular components of 
$V$. To the variation field $V$ we associate the column of 1-forms:
\begin{equation}
\s^{(V)}\theta^\mu=dV^\mu \ : \ \mu=0,...,n
\label{p.41}
\end{equation}
Note that this depends on the background structure. To a variation field $V$ and to the row of 
functions $(\beta_\mu : \mu=0,...,n)$ we associate the 1-form:
\begin{equation}
\s^{(V)}\xi=V^\mu d\beta_\mu 
\label{p.42}
\end{equation}
To the variation field $V$ is associated the $T^1_1$ type tensorfield:
\begin{equation}
\s^{(V)}S=V^\mu V^\nu\dot{T}_{\mu\nu}
\label{p.43}
\end{equation}
In view of \ref{p.35}, \ref{p.36}, \ref{p.42}, we have:
\begin{equation}
\s^{(V)}S=\tilde{h}^{-1}\cdot\s^{(V)}S_\flat
\label{p.44}
\end{equation}
where $\s^{(V)}S_\flat$ is the symmetric 2-covariant tensorfield:
\begin{equation}
\s^{(V)}S_\flat=\s^{(V)}\xi\otimes\s^{(V)}\xi-\frac{1}{2}(\s^{(V)}\xi,\s^{(V)}\xi)_h h
\label{p.45}
\end{equation}
which depends only on the conformal class of $h$. The identity \ref{p.39} together with the definition 
\ref{p.41} implies the identity:
\begin{eqnarray}
&&\mbox{div}_{\tilde{h}}\s^{(V)}S=(\s^{(V)}\xi,\s^{(V)}\theta^\mu)_{\tilde{h}}d\beta_\mu 
-(\s^{(V)}\xi,d\beta_\mu)_{\tilde{h}}\s^{(V)}\theta^\mu \nonumber\\
&&\hspace{20mm}+(\s^{(V)}\theta^\mu,d\beta_\mu)_{\tilde{h}}\s^{(V)}\xi
+V^\mu(\square_{\tilde{h}}\beta_\mu)\s^{(V)}\xi
\label{p.46}
\end{eqnarray}

A basic requirement on the set of variation fields $V$ is that they span the tangent space to 
${\cal K}$ at each point. The simplest way to achieve this is to choose one of the variation 
fields, which we denote by $Y$, to be at each point of ${\cal N}$ in the linear span of 
$N$ and $\Nb$ and along ${\cal K}$ colinear to $T$, and to choose the 
other variation fields so that at each point of ${\cal N}$ they span the tangent space to 
the surface $S_{\ub,u}$ though that point. We thus set:
\begin{equation}
Y=\gamma N+\ogamma\Nb
\label{p.47}
\end{equation}
In view of \ref{p.21}, the requirement that $Y$ is along ${\cal K}$ colinear to $T$ 
reduces to:
\begin{equation}
\ogamma=r\gamma \ \mbox{: along ${\cal K}$}
\label{p.48}
\end{equation}
The optimal choice is to set:
\begin{equation}
\gamma=1
\label{p.49}
\end{equation}
in which case \ref{p.48} reduces to:
\begin{equation}
\ogamma=r \ \mbox{: along ${\cal K}$}
\label{p.50}
\end{equation}
and to extend $\ogamma$ to ${\cal N}$ by the requirement that it be constant along the integral curves 
of $L$:
\begin{equation}
L\ogamma=0
\label{p.51}
\end{equation}
In 2 spatial dimensions there is an obvious choice of a variation field to complement $Y$, 
namely $E$, the unit tangent field of the curves $S_{\ub,u}$ (with the counterclockwise orientation). 
In higher dimensions, we complement $Y$ with the $(E_{(\mu)} : \mu=0,...,n)$ which are the 
$h$-orthogonal projections to the surfaces $S_{\ub,u}$ of the translation fields 
$(\partial/\partial x^\mu : \mu=0,...,n)$ of the background structure. These are given by 
(see \ref{8.a1}):
\begin{equation}
E_{(\mu)}=h_{\mu\nu}(\sd x^\nu)^\sharp
\label{p.52}
\end{equation}
where, given $\zeta\in T^*_q S_{\ub,u}$ we denote $\zeta^\sharp=\sh_q^{-1}\cdot\zeta\in T_qS_{\ub,u}$, 
considering $\sh_q$ as an isomorphism $\sh_q : T_q S_{\ub,u}\rightarrow T^*_q S_{\ub,u}$. 
The last part of Section 6.2 contains the analysis of the structure forms of the variation fields. 

The fundamental energy identities are discussed in Section 6.3. 
Given a vectorfield $X$, which we call {\em multiplier field}, we consider the vectorfield 
$\s^{(V)}P$ associated to $X$ and to a given variation field $V$ through $\s^{(V)}S$, defined by:
\begin{equation}
\s^{(V)}P=-\s^{(V)}S\cdot X
\label{p.53}
\end{equation}
We call $\s^{(V)}P$ the {\em energy current} associated to $X$ and to $V$. 
Let us denote by $\s^{(V)}Q$ the divergence of $\s^{(V)}P$ with respect to the conformal 
acoustical metric $\tilde{h}=\Omega h$:
\begin{equation}
\mbox{div}_{\tilde{h}} \s^{(V)}P=\s^{(V)}Q
\label{p.54}
\end{equation}
We have:
\begin{equation}
\s^{(V)}Q=\s^{(V)}Q_1+\s^{(V)}Q_2+\s^{(V)}Q_3
\label{p.55}
\end{equation} 
where:
\begin{eqnarray}
&&\s^{(V)}Q_1=-\frac{1}{2}\s^{(V)}S^\sharp\cdot\s^{(X)}\tilde{\pi} \label{p.56}\\
&&\s^{(V)}Q_1=-(\s^{(V)}\xi,\s^{(V)}\theta^\mu)_{\tilde{h}}X\beta_\mu 
+(\s^{(V)}\xi,d\beta_\mu)_{\tilde{h}}\s^{(V)}\theta^\mu(X) \nonumber\\
&&\hspace{15mm}-(\s^{(V)}\theta^\mu,d\beta_\mu)_{\tilde{h}}\s^{(V)}\xi(X) \label{p.57}\\
&&\s^{(V)}Q_3=-\s^{(V)}\xi(X) V^\mu \square_{\tilde{h}}\beta_\mu \label{p.58}
\end{eqnarray}
In \ref{p.56} $\s^{(V)}S^\sharp$ is the symmetric 2-contravariant tensorfield corresponding 
to $\s^{(V)}S$:
\begin{equation}
\s^{(V)}S^\sharp=\s^{(V)}S\cdot\tilde{h}^{-1}
\label{p.59}
\end{equation}
and 
\begin{equation}
\s^{(X)}\tilde{\pi}={\cal L}_X\tilde{h}
\label{p.a5}
\end{equation}
is the {\em deformation tensor} of $X$, the rate of change of the conformal acoustical metric 
with respect to the flow generated by $X$. The concept of a multiplier field goes back to the 
fundamental work of Noether [No] connecting symmetries to conserved quantities.

Integrating \ref{p.54} on a domain in ${\cal R}_{\delta,\delta}$ of the form:
\begin{equation}
{\cal R}_{\ub_1,u_1}=R_{\ub_1,u_1}\times S^{n-1}=\bigcup_{(\ub,u)\in R_{\ub_1,u_1}}S_{\ub,u}
\label{p.60}
\end{equation}
where, with $(\ub_1,u_1)\in R_{\delta,\delta}$ we denote:
\begin{equation}
R_{\ub_1,u_1}=\{(\ub,u) \ : \ u\in[\ub,u_1], \ub\in[0,u_1]\}
\label{p.61}
\end{equation}
we obtain the {\em fundamental energy identity} corresponding to the variation field $V$ and 
to the multiplier field X:
\begin{equation}
\s^{(V)}{\cal E}^{\ub_1}(u_1)+\s^{(V)}\underline{{\cal E}}^{u_1}(\ub_1)+\s^{(V)}{\cal F}^{\ub_1}-\s^{(V)}\underline{{\cal E}}^{u_1}(0)
=\s^{(V)}{\cal G}^{\ub_1,u_1}
\label{p.62}
\end{equation}
Here, $\s^{(V)}{\cal E}^{\ub_1}(u_1)$ and $\s^{(V)}\underline{{\cal E}}^{u_1}(\ub_1)$ are the {\em energies}: 
\begin{eqnarray}
&&\s^{(V)}{\cal E}^{\ub_1}(u_1)=\int_{C_{u_1}^{\ub_1}}\Omega^{(n-1)/2}\s^{(V)}S_\flat(X,L)\nonumber\\
&&\s^{(V)}\underline{{\cal E}}^{u_1}(\ub_1)=\int_{\Cb_{\ub_1}^{u_1}}\Omega^{(n-1)/2}\s^{(V)}S_\flat(X,\Lb)
\label{p.63}
\end{eqnarray}
$\s^{(V)}{\cal F}^{\ub_1}$ is the {\em flux}:
\begin{equation}
\s^{(V)}{\cal F}^{\ub_1}=\int_{{\cal K}^{\ub_1}}\Omega^{(n-1)/2}\s^{(V)}S_\flat(X,M)
\label{p.64}
\end{equation}
where 
\begin{equation}
M=\Lb-L
\label{p.65}
\end{equation}
is a normal to ${\cal K}$ pointing to the interior of ${\cal K}$. In \ref{p.63} we denote by 
$C_{u_1}^{\ub_1}$ the part of $C_{u_1}$ corresponding to $\ub\leq\ub_1$ and by $\Cb_{\ub_1}^{u_1}$ 
the part of $\Cb_{\ub_1}$ corresponding to $u\leq u_1$. The right hand side of \ref{p.62} 
is the {\em error integral}: 
\begin{equation}
\s^{(V)}{\cal G}^{\ub_1,u_1}=\int_{{\cal R}_{\ub_1,u_1}}2a\Omega^{(n+1)/2}\s^{(V)}Q
\label{p.66}
\end{equation}
The energies are positive semi-definite if the multiplier field $X$ acoustically timelike 
future-directed. In Chapter 7 the conditions on $X$ are investigated which make the flux coercive 
when $\s^{(V)}\xi$ satisfies on ${\cal K}$ the boundary condition to be presently discussed. 

The boundary condition for $\s^{(V)}\xi$ on ${\cal K}$ appears as Proposition 6.2 in the form:
\begin{equation}
\s^{(V)}\xi_+(A_+)=\s^{(V)}\xi_-(A_-)
\label{p.67}
\end{equation}
where the subscripts $+$ and $-$ denote the future and past sides of ${\cal K}$ respectively 
and $A_\pm$ are the vectorfields:
\begin{equation}
A_\pm=\frac{\triangle I}{\delta}-K_\pm
\label{p.68}
\end{equation} 
Here $I$ is the particle current 
\begin{equation}
I=n u 
\label{p.69}
\end{equation}
where $n$ is the rest mass density and $u$ the spacetime fluid velocity 
(see \ref{1.12} in regard to the relativistic theory, \ref{1.110} in regard to the non-relativistic 
theory). We define $\zeta$ to be the covectorfield along ${\cal K}$ such that at each $p\in{\cal K}$ 
the null space of $\zeta_p$ is $T_p{\cal K}$, $\zeta_p(U)>0$ if the vector $U$ points to the interior 
of ${\cal N}$, $\zeta$ being normalized in the relativistic theory to be of unit magnitude 
with respect to the Minkowski metric  and in the non-relativistic theory by the condition that 
$\overline{\zeta}$, the restriction of $\zeta$ to the $\Sigma_t$, to be of unit magnitude with 
respect to the Euclidean metric. The nonlinear jump condition can then be stated in the form:
\begin{equation}
\zeta\cdot\triangle I=0
\label{p.70}
\end{equation}
that is, the vectorfield $\triangle I$ along ${\cal K}$ is tangential to ${\cal K}$, while the 
linear jump conditions \ref{p.17} take the form:
\begin{equation}
\triangle\beta=\delta \zeta
\label{p.71}
\end{equation}
for some function $\delta$ on ${\cal K}$. This clarifies the meaning of the 1st term on the right in \ref{p.68}. As for the 2nd term, $K_\pm$ it is a normal, relative to the acoustical metric, vectorfield 
along the future and past sides of ${\cal K}$ respectively, given by:
\begin{equation}
K_\pm^\mu=G_\pm(h_\pm^{-1})^{\mu\nu}\zeta_\nu 
\label{p.72}
\end{equation}
where $G$ is in the relativistic theory the thermodynamic function \ref{1.72} and reduces to $n$ in the 
non-relativistic theory. In Section 13.1 we revisit the proof of Proposition 6.2, clarifying things 
further. Proposition 6.3 plays a central role in the analysis of the coercivity of the flux integrant 
in Chapter 7. This proposition states that the 1st term on the right in \ref{p.68}, a vectorfield 
along ${\cal K}$ tangential to ${\cal K}$, is timelike future-directed with respect to the acoustical 
metric defined by the future solution. The vectorfields $A_\pm$ are singular at $S_{0,0}$. 
In the following we drop the subscript + in reference to quantities defined on ${\cal K}$ by the future 
solution. 
To express the boundary condition \ref{p.67} in a form which is regular at $S_{0,0}$ we multiply 
$A$, $A_-$ by $\kappa G^{-1}$, defining:
\begin{equation}
B=\kappa G^{-1}A, \ \ B_-=\kappa G^{-1}A_-
\label{p.73}
\end{equation}
where $\kappa$ is the positive function defined by:
\begin{equation}
\kappa K=\frac{1}{2}G M
\label{p.74}
\end{equation}
noting that by \ref{p.72} $K$ is colinear and in the same sense as $M$ (see \ref{p.65}). The 
boundary condition \ref{p.67} is then equivalent to its regularized form:
\begin{equation}
\s^{(V)}\xi(B)=\s^{(V)}\xi_-(B_-)
\label{p.75}
\end{equation}
and we have:
\begin{equation}
B=B_{||}+B_{\bot}, \ \ B_-=B_{||}+B_{-\bot}
\label{p.76}
\end{equation}
where $B_{||}$ is tangential to ${\cal K}$, 
\begin{equation}
B_{\bot}=-\frac{1}{2}M
\label{p.77}
\end{equation}
is an exterior normal to ${\cal K}$ relative to the acoustical metric defined by the future 
solution, and $B_{-\bot}$ is a normal to ${\cal K}$ in the same sense relative to the acoustical 
metric defined by the past solution. 
The form of the vectorfields $B_{||}$, $B_{-\bot}$ in a neighborhood of $S_{0,0}$ in ${\cal K}$ 
is analyzed in the last part of Section 6.4.

Section 7.1 determines the necessary and sufficient conditions on the multiplier field $X$ for the 
flux integrant (see \ref{p.64}) to be coercive. This rests on the fundamental work of 
G$\stackrel{\circ}{\mbox{a}}$rding 
[Ga] who first showed in connection with the initial - boundary value problem for the wave equation 
$\square_g\phi=0$ on a Lorentzian manifold $({\cal M},g)$ with timelike boundary ${\cal K}$ that 
the boundary condition of prescribing $B\phi$ on ${\cal K}$, where $B$ is a vectorfield along 
${\cal K}$, is well posed if $B$ is of the form \ref{p.76} with $B_{||}$ tangential to ${\cal K}$ 
and timelike future-directed while $B_{\bot}$ is an exterior normal to ${\cal K}$, for in this 
case an appropriate energy inequality can be derived. Here, in  Proposition 7.1 we give a simple 
geometric characterization of the set of multiplier fields making the flux integrant coercive. 
We then show that a suitable choice for the multiplier field is:
\begin{equation}
X=3L+\Lb \ \ \mbox{: on ${\cal N}$}
\label{p.78}
\end{equation}
With this choice there is a constant $C^\prime$ such that 
\begin{equation}
\s^{(V)}{\cal F}^{\prime\ub_1}=\s^{(V)}{\cal F}^{\ub_1}+2C^\prime\int_{{\cal K}^{\ub_1}}\Omega^{(n-1)/2}
(\s^{(V)}b)^2
\label{p.79}
\end{equation}
is positive-definite (see \ref{9.a34}). Here:
\begin{equation}
\s^{(V)}b=\s^{(V)}\xi_-(B_-)=B_-^\mu V\beta_{-\mu}
\label{p.80}
\end{equation}
(see \ref{13.28}). Adding the 2nd term on the right in \ref{p.79} to both sides of the energy 
identity \ref{p.62}, the last takes the form:
\begin{eqnarray}
&&\s^{(V)}{\cal E}^{\ub_1}(u_1)+\s^{(V)}\underline{{\cal E}}^{u_1}(\ub_1)+\s^{(V)}{\cal F}^{\prime\ub_1}\nonumber\\
&&\hspace{18mm}=\s^{(V)}\underline{{\cal E}}^{u_1}(0)+\s^{(V)}{\cal G}^{\ub_1,u_1}
+2C^\prime\int_{{\cal K}^{\ub_1}}\Omega^{(n-1)/2}(\s^{(V)}b)^2 \nonumber\\
&&\label{p.81}
\end{eqnarray}
The error integral \ref{p.66} decomposes into:
\begin{equation}
\s^{(V)}{\cal G}^{\ub_1,u_1}=\s^{(V)}{\cal G}_1^{\ub_1,u_1}+\s^{(V)}{\cal G}_2^{\ub_1,u_1}
+\s^{(V)}{\cal G}_3^{\ub_1,u_1}
\label{p.82}
\end{equation}
according to the decomposition \ref{p.54} of $\s^{(V)}Q$. 

The deformation tensor $\s^{(X)}\tilde{\pi}$ of the multiplier field is analyzed in Section 7.2. 
The error integral $\s^{(V)}{\cal G}_1^{\ub_1,u_1}$ is then estimated. In the estimates \ref{7.122}, 
\ref{7.124} singular integrals first appear. But, as we shall see below, this is only the tip of the 
iceberg. In Section 7.3 the error integral $\s^{(V)}{\cal G}_2^{\ub_1,u_1}$ is estimated using the results 
of Section 6.2 on the structure forms of the variation fields. 

The {\em commutation fields } which are used to control the higher order analogues of the 
functions $\beta_\mu : \mu=0,....,n$ are defined in Section 7.1. Commutation fields where first 
introduced by Klainerman in his derivation in [Kl] of the decay properties of the solutions 
of the wave equation in Minkowski spacetime using the fact that this equation is 
invariant under the Poincar\'{e} group, the isometry group of Minkowski spacetime, 
the commutation fields being the vectorfields generating the group action. The scope of 
commutation fields was substantially extended in [Ch-Kl] where the problem of the stability 
of the Minkowski metric in the context of the vacuum Einstein equations of general relativity 
was solved. In that work, while the metric arising as the development of  
general asymptotically flat initial data does not possess a non-trivial isometry group, nevertheless 
a large enough subgroup of the scale extended Poincar\'{e} group was found and an action of 
this subgroup which approximates that of an isometry in the sense that the deformation tensors 
of the vectorfields generating this action, the commutation fields, are appropriately bounded with 
decay. In the present monograph, denoting (see \ref{p.4}), for $\sigma\in[0,\delta]$, 
\begin{equation}
{\cal K}^\delta_\sigma=\{(\tau,\sigma) \ : \ \tau\in [0,\delta-\sigma]\}\times S^{n-1}
\label{p.83}
\end{equation}
(note that ${\cal K}^\delta_0={\cal K}^\delta$) we require that at each point $q\in{\cal N}$, 
$q\in {\cal K}^\delta_\sigma$, the set of commutation fields $C$ 
to span $T_q{\cal K}^\delta_\sigma$. As first of the commutation fields we take the vectorfield $T$. 
The remaining commutation fields are then required to span the tangent space to the $S_{\ub,u}$ at 
each point. In $n=2$ spatial dimensions we choose $E$ to complement $T$ as a commutation field. 
For $n>2$ we choose the $E_{(\mu)} : \mu=0,...,n$ (see \ref{p.51}) to complement $T$. Thus 
$E$ for $n=2$ and the $E_{(\mu)}$ for $n>2$ play a dual role being commutation fields as well as 
variation fields. However, what characterizes the action of a variation field $V$ is the corresponding 
structure form $\s^{(V)}\theta$, what characterizes the action of a commutation field $C$ is the 
corresponding deformation tensor $\s^{(C)}\tilde{\pi}={\cal L}_C\tilde{h}$. 

The commutation fields generate higher order analogues of the functions $\beta_\mu : \mu=0,...,n$. 
At order $m+l$ we have, for $n=2$:
\begin{equation}
\s^{(m,l)}\beta_\mu=E^l T^m\beta_\mu 
\label{p.84}
\end{equation}
and for $n>2$:
\begin{equation}
\s^{(m,\nu_1...\nu_l)}\beta_\mu=E_{(\nu_l)} ... E_{(\nu_1)} T^m\beta_\mu 
\label{p.85}
\end{equation}
To these and to the variation field $V$ there correspond higher order analogues of the 
1-form $\s^{(V)}\xi$, namely, for $n=2$:
\begin{equation}
\s^{(V;m,l)}\xi=V^\mu d\s^{(m,l)}\beta_\mu
\label{p.86}
\end{equation}
and for $n>2$:
\begin{equation}
\s^{(V;m,\nu_1...\nu_l)}\xi=V^\mu d\s^{(m,\nu_1...\nu_l)}\beta_\mu
\label{p.87}
\end{equation} 
The preceding identities \ref{p.39}, \ref{p.46}, \ref{p.53} - \ref{p.58}, \ref{p.62}, \ref{p.81} - 
\ref{p.82} which refer to $\beta_\mu$ and to 
$\s^{(V)}\xi$ all 
hold with these higher order analogues in the role of $\beta_\mu$ and $\s^{(V)}\xi$ respectively, 
the boundary condition for $\s^{(V;m,l)}\xi$ (case $n=2$) being of the form:
\begin{equation}
\s^{(V;m,l)}\xi(B)=\s^{(V;m,l)}b
\label{p.88}
\end{equation}
where $\s^{(V;m,l)}b$ is analyzed in Section 13.1. Similarly, for $n>2$, $\s^{(V;m,\nu_1...\nu_l)}\xi$ 
is of the form:
\begin{equation}
\s^{(V;m,\nu_1...\nu_l)}\xi(B)=\s^{(V;m,\nu_1...\nu_l)}b
\label{p.89}
\end{equation}
 However while for the original 
$\beta_\mu$ we have $\square_{\tilde{h}}\beta_\mu=0$, hence the error term $\s^{(V)}Q_3$ (see 
\ref{p.58}) vanishes, this is no longer true for the higher order analogues. Instead we have:
\begin{eqnarray}
&&\Omega a\square_{\tilde{h}}\s^{(m,l)}\beta_\mu=\s^{(m,l)}\tilde{\rho}_\mu \ \ : \ \mbox{in the case 
$n=2$} \label{p.90}\\
&&\Omega a\square_{\tilde{h}}\s^{(m,\nu_1...\nu_l)}\beta_\mu=\s^{(m,\nu_1...\nu_l)}\tilde{\rho}_\mu \ \ 
: \ \mbox{for $n>2$} \label{p.91}
\end{eqnarray}
The $\s^{(m,l)}\tilde{\rho}_\mu$, $\s^{(m,\nu_1...\nu_l)}\tilde{\rho}_\mu$, which we call (rescaled) 
source functions, obey certain recursion formulas, deduced in Section 8.2, 
which determine them for all $m$ and $l$. 
In the last section of Chapter 8, Section 8.4, the error terms at order $m+l$  are discerned 
which contain the acoustical quantities of highest order, $m+l+1$. These are contained in the 
error integral $\s^{(V;m,l)}{\cal G}_3^{\ub_1,u_1}$ (case $n=2$), 
$\s^{(V;m,\nu_1...\nu_l)}{\cal G}_3^{\ub_1,u_1}$ (case $n>2$), which by \ref{p.66} and \ref{p.57} 
is given by:
\begin{eqnarray}
&&\s^{(V;m,l)}{\cal G}_3^{\ub_1,u_1}=-\int_{{\cal R}_{\ub_1,u_1}}2\Omega^{(n-1)/2}\s^{(V;m,l)}\xi(X)
V^\mu\s^{(m,l)}\tilde{\rho}_\mu \nonumber\\
&&\hspace{2mm}(\mbox{case $n=2$})\label{p.92}
\end{eqnarray}
\begin{eqnarray}
&&\s^{(V;m,\nu_1...\nu_l)}{\cal G}_3^{\ub_1,u_1}=-\int_{{\cal R}_{\ub_1,u_1}}2\Omega^{(n-1)/2}\s^{(V;m,\nu_1...\nu_l)}\xi(X)
V^\mu\s^{(m,\nu_1...\nu_l)}\tilde{\rho}_\mu \nonumber\\
&&\hspace{12mm}(\mbox{for $n>2$})\label{p.93}
\end{eqnarray} 
In the case $n=2$ the leading terms in $\s^{(m,l)}\tilde{\rho}$ involving the acoustical quantities 
of highest order are (see \ref{8.158}, \ref{8.159}):
\begin{eqnarray}
&&\mbox{for $m=0$:}\hspace{7mm}\frac{1}{2}\rho(\Lb\beta_\mu)E^l\tchi+\frac{1}{2}\rhob(L\beta_\mu)E^l\tchib 
\label{p.94}\\
&&\mbox{for $m\geq 1$:}\hspace{7mm}\rho(\Lb\beta_\mu)E^l T^{m-1}E^2\lambda
+\rhob(L\beta_\mu)E^l T^{m-1}E^2\lambdab \label{p.95}
\end{eqnarray}
In the case $n>2$ the leading terms in $\s^{(m,\nu_1...\nu_l)}\tilde{\rho}$ involving the acoustical quantities 
of highest order are (see \ref{8.160}, \ref{8.161}):
\begin{eqnarray}
&&\mbox{for $m=0$:}\hspace{7mm}\frac{1}{2}\rho(\Lb\beta_\mu)E_{(\nu_l)}...E_{(\nu_1)}\mbox{tr}\tchi
+\frac{1}{2}\rhob(L\beta_\mu)E_{(\nu_l)}...E_{(\nu_1)}\mbox{tr}\tchib\nonumber\\
&&\hspace{22mm}+\frac{ah_{\nu_1,\kappa}}{2c}(\sd\beta_\mu)^\sharp\cdot
\left(\Nb^{\kappa}\sd(E_{(\nu_l)}...E_{(\nu_2)}\mbox{tr}\tchi
+N^{\kappa}\sd(E_{(\nu_l)}...E_{(\nu_2)}\mbox{tr}\tchib\right) \nonumber\\
&&\label{p.96}\\
&&\mbox{for $m\geq 1$:}\hspace{7mm}\rho(\Lb\beta_\mu)E_{(\nu_l)}...E_{(\nu_1)}T^{m-1}\slap\lambda+\rhob(L\beta_\mu)E_{(\nu_l)}...E_{(\nu_1)}T^{m-1}\slap\lambdab \nonumber\\
&&\label{p.97}
\end{eqnarray}

We now come to the main analytic method introduced in this monograph. To motivate the introduction 
of this method, we need to first discuss the difficulties encountered. After this discussion it will 
become evident that the new analytic method is a natural way to overcome these difficulties. 
In Chapters 1 - 8 which we have just reviewed, we denote by $n$ the spatial dimension, a notation 
which we have followed in the above review. However in the remainder of the monograph, Chapters 
9 - 14, which we shall presently review, we designate by $n$ the top order of the 
$\s^{(m,l)}\beta_\mu$, that is $m+l=n$, and we denote the spatial dimension by $d$. We shall now 
follow the latter notation.

The difficulties arise in estimating the contribution of the terms involving the top order 
(order $n+1$) acoustical quantities, \ref{p.94}, \ref{p.95} in the case $d=2$, \ref{p.96}, \ref{p.97} 
for $d>2$, to the error integral $\s^{(V;m,l)}{\cal G}_3^{\ub_1,u_1}$, 
$\s^{(V;m,\nu_1...\nu_l)}{\cal G}_3^{\ub_1,u_1}$. 
To understand the origin of these difficulties 
it is advantageous to go directly to Chapter 11 where the approach applicable to any $d\geq 2$ is 
laid, an approach which simplifies in the case $d=2$ where the detailed estimates are deduced in 
Chapter 10. In regard to \ref{p.96} we must derive appropriate estimates for 
\begin{equation}
\sd(E_{(\nu_{l-1})}...E_{(\nu_1)}\mbox{tr}\tchi), \ 
\sd(E_{(\nu_{l-1})}...E_{(\nu_1)}\mbox{tr}\tchib) \  : \  \nu_1,...,\nu_{l-1}=0,...,d
\label{p.98}
\end{equation} 
In regard to \ref{p.97} we must derive appropriate estimates for
\begin{equation} 
E_{(\nu_l)}...E_{(\nu_1)}T^{m-1}\slap\lambda, \ 
E_{(\nu_l)}...E_{\nu_1)}T^{m-1}\slap\lambdab \ : \  \nu_1,...,\nu_l=0,...,d. 
\label{p.99}
\end{equation}

Now, as mentioned above, Proposition 3.4, the 2nd variation equations, 
express $\sL_L\tchi$ and $\sL_{\Lb}\tchib$ 
in terms of 2nd order quantities with vanishing 2nd order acoustical part. These imply 
expressions for $L\mbox{tr}\tchi$ and $\Lb\mbox{tr}\tchib$ again in terms of 2nd order quantities 
with vanishing 2nd order acoustical part. To be able then to estimate $\mbox{tr}\tchi$, 
$\mbox{tr}\tchib$ in terms of 1st order quantities, so that we can estimate \ref{p.98} in 
terms of quantities of the top order $l+1=n+1$, we must express the principal (2nd order) part of 
the expressions for $L\mbox{tr}\tchi$ and $\Lb\mbox{tr}\tchib$ in the form 
$-L\hat{f}$ and $-\Lb\hat{\fb}$ respectively, up to lower order terms, with $\hat{f}$ and $\hat{\fb}$ 
being quantities of 1st order. That this possible follows from the fact that the quantities 
(see \ref{11.14}):
\begin{equation}
M=\frac{1}{2}\beta_N^2(a\slap H-L(\Lb H)), \ \ \Mb=\frac{1}{2}\beta_{\Nb}^2(a\slap H-\Lb(LH))
\label{p.100}
\end{equation}
with $H$ a function of the $\beta_\mu : \mu=0,...,n$, are actually 1st order quantities. This 
fact is a direct consequence of the equation $\square_{\tilde{h}}\beta_\mu=0$. The functions 
$\hat{f}$, $\hat{\fb}$ each contain a singular term with coefficient $\lambda^{-1}$, $\lambdab^{-1}$ 
respectively. The functions $f=\lambda\hat{f}$, $\fb=\lambdab\hat{\fb}$ are then regular, and 
transferring the corresponding terms to the left hand side, we obtain propagation equations 
for the quantities:
\begin{equation}
\theta=\lambda\mbox{tr}\tchi+f, \ \ \thetab=\lambdab\mbox{tr}\tchib+\fb
\label{p.101}
\end{equation}
of the form:
\begin{equation}
L\theta=R, \ \ \Lb\theta=\Rb
\label{p.102}
\end{equation}
where $R$, $\Rb$ are again quantities of order 1, their 1st order acoustical parts being given by 
Proposition 11.1. The leading terms in $R$, $\Rb$ from the point of view of behavior as we approach 
the singularity at $\ub=0$ (that is $\Cb_0$) and the stronger singularity at $u=0$ (that is $S_{0,0}$) 
are the terms:
\begin{equation}
2(L\lambda)\mbox{tr}\tchi=2\lambda^{-1}(L\lambda)(\theta-f) \ : \mbox{in $R$}, \ \ \ 
2(\Lb\lambdab)\mbox{tr}\tchib=2\lambdab^{-1}(\Lb\lambdab)(\thetab-\fb) \ : \mbox{in $\Rb$}
\label{p.103}
\end{equation}
To estimate \ref{p.98} we introduce the quantities, of order $l+1=n+1$,
\begin{eqnarray}
&&\s^{(\nu_1...\nu_{l-1})}\theta_l=\lambda \sd\left(E_{(\nu_{l-1})} . . . E_{(\nu_1)}
\mbox{tr}\tchi\right)+\sd\left(E_{(\nu_{l-1})} . . . E_{(\nu_1)}f\right)\nonumber\\
&&\s^{(\nu_1...\nu_{l-1})}\thetab_l=\lambdab \sd\left(E_{(\nu_{l-1})} . . . E_{(\nu_1)}
\mbox{tr}\tchib\right)+\sd\left(E_{(\nu_{l-1})} . . . E_{(\nu_1)}\fb\right)\nonumber\\
&&\label{p.104}
\end{eqnarray}
These quantities then satisfy propagation 
equations of the form:
\begin{eqnarray}
&&\sL_L\s^{(\nu_1...\nu_{l-1})}\theta_l=\s^{(\nu_1...\nu_{l-1})}R_l\nonumber\\
&&\sL_L\s^{(\nu_1...\nu_{l-1})}\thetab_l=\s^{(\nu_1...\nu_{l-1})}\Rb_l \label{p.105}
\end{eqnarray}
where $\s^{(\nu_1...\nu_{l-1})}R_l$, $\s^{(\nu_1...\nu_{l-1})}\Rb_l$ 
are likewise also of order $l+1=n+1$ and their leading terms from the point of view of behavior as we 
approach the singularities are: 
\begin{eqnarray}
&&2(L\lambda)\sd(E_{(\nu_{l-1})} . . . E_{(\nu_1)}\mbox{tr}\tchi)=
2\lambda^{-1}(L\lambda)\left(\s^{(\nu_1...\nu_{l-1})}\theta_l
-\sd(E_{(\nu_{l-1})} . . . E_{(\nu_1)}f)\right), \nonumber\\
&&2(\Lb\lambdab)\sd(E_{(\nu_{l-1})} . . . E_{(\nu_1)}\mbox{tr}\tchib)=
2\lambdab^{-1}(\Lb\lambdab)\left(\s^{(\nu_1...\nu_{l-1})}\thetab_l
-\sd(E_{(\nu_{l-1})} . . . E_{(\nu_1)}\fb)\right)\nonumber\\
&&\label{p.106}
\end{eqnarray}

As mentioned above, Proposition 3.3, the propagation equations for $\lambda$, $\lambdab$, 
express $L\lambda$ and $\Lb\lambdab$ 
in terms of 1st order quantities with vanishing 1st order acoustical part. These imply 
expressions for $L\slap\lambda$ and $\Lb\slap\lambdab$ in terms of 3rd order quantities 
with vanishing 3rd order acoustical part. To be able then to estimate $\slap\lambda$, 
$\slap\lambdab$ in terms of 2nd order quantities, so that we can estimate \ref{p.99} in 
terms of quantities of the top order $m+l+1=n+1$, we must express the principal (3rd order) part of 
the expressions for $L\slap\lambda$ and $\Lb\slap\lambdab$ in the form 
$L\hat{j}$ and $\Lb\hat{\jb}$ respectively, up to lower order terms, with $\hat{j}$ and $\hat{\jb}$ 
being quantities of 2nd order. That this possible follows again from the fact that the quantities 
$M$ and $\Mb$ defined by \ref{p.100} are actually 1st order quantities. The functions 
$\hat{j}$, $\hat{\jb}$ each contain a singular term with coefficient $\lambda^{-1}$, $\lambdab^{-1}$ 
respectively. The functions $j=\lambda\hat{j}$, $\jb=\lambdab\hat{\jb}$ are then regular, and 
transferring the corresponding terms to the left hand side, we obtain propagation equations 
for the quantities:
\begin{equation}
\nu=\lambda\slap\lambda-j, \ \ \nub=\lambdab\slap\lambdab-\jb
\label{p.107}
\end{equation}
of the form:
\begin{equation}
L(\nu-\tau)=I, \ \ \Lb(\nub-\taub)=\Ib
\label{p.108}
\end{equation}
where $\tau$, $\taub$ are quantities of order 1, while $I$, $\Ib$ are quantities of order 2, their 2nd  order acoustical parts being given by Proposition 11.2. The leading terms in $I$, $\Ib$ from the point of view of behavior as we approach 
the singularity at $\ub=0$ (that is $\Cb_0$) and the stronger singularity at $u=0$ (that is $S_{0,0}$) 
are the terms:
\begin{equation}
2(L\lambda)\slap\lambda=2\lambda^{-1}(L\lambda)(\nu+j) \ : \mbox{in $I$}, \ \ \ 
2(\Lb\lambdab)\slap\lambdab=2\lambdab^{-1}(\Lb\lambdab)(\nub+\jb) \ : \mbox{in $\Ib$}
\label{p.109}
\end{equation}
To estimate \ref{p.99} we introduce the quantities, of order $m+l+1=n+1$,
\begin{eqnarray}
&&\s^{(\nu_1 ... \nu_l)}\nu_{m-1,l+1}=\lambda E_{(\nu_l)} ... E_{(\nu_1)}T^{m-1}\slap\lambda
-E_{(\nu_l)} ... E_{(\nu_1)}T^{m-1}j\nonumber\\
&&\s^{(\nu_1 ... \nu_l)}\nub_{m-1,l+1}=\lambdab E_{(\nu_l)} ... E_{(\nu_1)}T^{m-1}\slap\lambdab
-E_{(\nu_l)} ... E_{(\nu_1)}T^{m-1}\jb \nonumber\\
&&\label{p.110}
\end{eqnarray}
These quantities then satisfy propagation equations of the form:
\begin{eqnarray}
&&L\left(\s^{(\nu_1 ... \nu_l)}\nu_{m-1,l+1}-\s^{(\nu_1 ... \nu_l)}\tau_{m-1,l+1}\right)
=\s^{(\nu_1 ... \nu_l)}I_{m-1,l+1}\nonumber\\
&&\Lb\left(\s^{(\nu_1 ... \nu_l)}\nub_{m-1,l+1}-\s^{(\nu_1 ... \nu_l)}\taub_{m-1,l+1}\right)
=\s^{(\nu_1 ... \nu_l)}\Ib_{m-1,l+1}\nonumber\\
\label{p.111}
\end{eqnarray}
where $\s^{(\nu_1 ... \nu_l)}\tau_{m-1,l+1}$, $\s^{(\nu_1 ... \nu_l)}\taub_{m-1,l+1}$ are quantities of order $m+l=n$, while $\s^{(\nu_1 ... \nu_l)}I_{m-1,l+1}$, $\s^{(\nu_1 ... \nu_l)}\Ib_{m-1,l+1}$ 
are quantities of order $m+l+1=n+1$ and their leading terms from the point of view of behavior as we 
approach the singularities are:  
\begin{eqnarray}
&&2(L\lambda)\sd(E_{(\nu_l)} . . . E_{(\nu_1)}T^{m-1}\slap\lambda) \nonumber\\
&&\hspace{20mm}=2\lambda^{-1}(L\lambda)\left(\s^{(\nu_1 ... \nu_l)}\nu_{m-1,l+1}+E_{(\nu_l)} ... E_{(\nu_1)}T^{m-1}j\right), \nonumber\\
&&2(\Lb\lambdab)\sd(E_{(\nu_l)} . . . E_{(\nu_1)}T^{m-1}\slap\lambdab) \nonumber\\
&&\hspace{20mm}=2\lambdab^{-1}(\Lb\lambdab)\left(\s^{(\nu_1 ... \nu_l)}\nub_{m-1,l+1}
+E_{(\nu_l)} ... E_{(\nu_1)}T^{m-1}\jb\right) \nonumber\\
&&\label{p.112}
\end{eqnarray}

Now, in accordance with the discussion following \ref{p.12}, we have:
\begin{equation}
\lambdab\sim \ub, \ \ \ \lambda\sim u^2
\label{p.113}
\end{equation}
where we denote by $\sim$ equality up to positive multiplicative constants. The propagation 
equations for $\lambda$, $\lambdab$ of Proposition 3.3 then imply:
\begin{equation}
L\lambda=O(\ub), \ \ \ \Lb\lambdab=O(\ub)
\label{p.114}
\end{equation}
Then, integrating the 1st of \ref{p.105} from $\Cb_0$, the contribution of the term 
$$-2\lambda^{-1}(L\lambda)\sd(E_{(\nu_{l-1})}...E_{(\nu_1)}f)$$
from the 1st of \ref{p.106} to $\|\s^{(\nu_1...\nu_{l-1})}\theta_l\|_{L^2(S_{\ub_1,u})}$ is bounded by:
\begin{equation}
Cu^{-2}\int_0^{\ub_1}\|\sd(E_{(\nu_{l-1})}...E_{(\nu_1)}f)\|_{L^2(S_{\ub,u})}\ub d\ub
\label{p.115}
\end{equation}
Now from the expression for the function $f$ of Proposition 11.1 we see that the term 
$$-\frac{1}{2}\beta_N^2\Lb H$$
in $f$ makes the leading contribution, and this contribution to \ref{p.115} is bounded, up to lower 
order terms, by:
\begin{equation}
Cu^{-2}\int_0^{\ub_1}\sum_{\nu_l}\|\Nb^\mu \Lb\s^{(0,\nu_1...\nu_{l-1}\nu_l)}\beta_\mu\|
_{L^2(S_{\ub,u})}\ub d\ub
\label{p.116}
\end{equation}
(compare with \ref{10.370}). From the 2nd of \ref{p.63} and from \ref{p.45} with 
$\s^{(0,\nu_1...\nu_l)}\beta_\mu$, $\s^{(V;0,\nu_1...\nu_l)}\xi$ in the roles of $\beta_\mu$, 
$\s^{(V)}\xi$, in view of \ref{p.78}, we have:
\begin{equation}
\s^{(V;0,\nu_1...\nu_l)}\underline{{\cal E}}^{u_1}(\ub)
=\int_{\Cb_{\ub}^{u_1}}\Omega^{(d-1)/2}\left((\s^{(V;0,\nu_1...\nu_l)}\xi_{\Lb})^2
+3a|\s^{(V;0,\nu_1...\nu_l)}\sxi|^2\right)
\label{p.117}
\end{equation}
($\s^{(V)}\sxi$ is the 1-form $\s^{(V)}\xi$ induced on the $S_{\ub,u}$). Now, the quantity 
\begin{equation}
\|\Nb^\mu \Lb\s^{(0,\nu_1...\nu_l)}\beta_\mu\|_{L^2(S_{\ub,u})}
\label{p.118}
\end{equation}
in the integrant in \ref{p.116} can only be estimated through 
\begin{equation}
\|\s^{(Y;0,\nu_1...\nu_l)}\xi_{\Lb}\|_{L^2(S_{\ub,u})}
\label{p.119}
\end{equation}
Comparing with \ref{p.47}, in view of the fact that by \ref{p.49} - \ref{p.51} and \ref{p.113} we have:
\begin{equation}
\ogamma\sim u
\label{p.120}
\end{equation}
we conclude that \ref{p.118} can only be bounded in terms of:
\begin{equation}
Cu^{-1}\|\s^{(Y;0,\nu_1...\nu_l)}\xi_{\Lb}\|_{L^2(S_{\ub,u})}
\label{p.121}
\end{equation}
therefore \ref{p.116} can only be bounded in terms of:
\begin{equation}
Cu^{-3}\sum_{\nu_l}\int_0^{\ub_1}\|\s^{(Y;0,\nu_1...\nu_{l-1}\nu_l)}\xi_{\Lb}\|_{L^2(S_{\ub,u})}\ub d\ub
\label{p.122}
\end{equation}
This bounds the leading contribution to $\|\s^{(\nu_1...\nu_{l-1})}\theta_l\|_{L^2(S_{\ub_1,u})}$. 
The corresponding contribution to $\|\s^{(\nu_1...\nu_{l-1})}\theta_l\|_{L^2(\Cb_{\ub_1}^{u_1})}$ is 
then bounded by:
\begin{equation}
C\sum_{\nu_l}\left\{\int_{\ub_1}^{u_1}\left(u^{-3}\int_0^{\ub_1}\|\s^{(Y;0,\nu_1...\nu_{l-1}\nu_l)}\xi_{\Lb}\|_{L^2(S_{\ub,u})}\ub d\ub\right)^2 du\right\}^{1/2}
\label{p.123}
\end{equation}
To bound the integral in the square root we must use the Schwartz inequality to replace it by 
(1/3 times):
\begin{eqnarray}
&&\ub_1^3\int_{\ub_1}^{u_1}u^{-6}\left(\int_0^{\ub_1}
\|\s^{(Y;0,\nu_1...\nu_{l-1}\nu_l)}\xi_{\Lb}\|^2_{L^2(S_{\ub,u})}d\ub\right)du\nonumber\\
&&=\ub_1^3\int_0^{\ub_1}\left(\int_{\ub_1}^{u_1} u^{-6}
\|\s^{(Y;0,\nu_1...\nu_{l-1}\nu_l)}\xi_{\Lb}\|^2_{L^2(S_{\ub,u})}du\right)d\ub\nonumber\\
&&\leq \ub_1^{-3}\int_0^{\ub_1}\|\s^{(Y;0,\nu_1...\nu_{l-1}\nu_l)}\xi_{\Lb}\|^2_{L^2(\Cb_{\ub}^{u_1})}
d\ub\nonumber\\
&&\leq C\ub_1^{-3}\int_0^{\ub_1}\s^{(Y;0,\nu_1...\nu_{l-1}\nu_l)}\underline{{\cal E}}^{u_1}(\ub)d\ub 
\label{p.124}
\end{eqnarray}
Then \ref{p.123} is bounded by:
\begin{equation}
C\sum_{\nu_l}\ub_1^{-3/2}\left\{\int_0^{\ub_1}\s^{(Y;0,\nu_1...\nu_{l-1}\nu_l)}\underline{{\cal E}}^{u_1}(\ub)d\ub\right\}^{1/2}
\label{p.125}
\end{equation}
In regard to \ref{p.93} with $m=0$, $V=Y$, the factor $Y^\mu\s^{(0,\nu_1...\nu_l)}\tilde{\rho}_\mu$ 
in the integrant, contains the term 
\begin{equation}
\frac{1}{2}\rho(Y^\mu\Lb\beta_\mu)E_{(\nu_l)}...E_{(\nu_1)}\mbox{tr}\tchi
\label{p.126}
\end{equation}
contributed by the 1st term in \ref{p.96}. Since by \ref{p.47}, \ref{p.49}:
\begin{equation}
Y^\mu\Lb\beta_\mu=N^\mu\Lb\beta_\mu+\ogamma\Nb^\mu\Lb\beta_\mu=\lambda\mbox{tr}\sss+\ogamma s_{\Nb\Lb}
\label{p.127}
\end{equation}
where:
\begin{equation}
\sss=\sd x^\mu\otimes\sd\beta_\mu \ \ \mbox{hence} \ \ \mbox{tr}\sss=(\sd x^\mu,\sd\beta_\mu)_{\sh}
\label{p.128}
\end{equation}
and (see \ref{p.34}):
\begin{equation}
s_{\Nb\Lb}=\Nb^\mu\Lb\beta_\mu 
\label{p.129}
\end{equation}
\ref{p.126} contributes through \ref{p.104} the term:
\begin{equation}
\frac{1}{2}s_{\Nb\Lb}\rho\ogamma\lambda^{-1}E_{(\nu_l)}\cdot\s^{(\nu_1...\nu_{l-1})}\theta_l
\label{p.130}
\end{equation}
In view of \ref{p.113}, \ref{p.120}, this term contributes to 
$\s^{(Y;0,\nu_1...\nu_l)}{\cal G}_3^{\ub_1,u_1}$ (see \ref{p.93}) through the $\Lb$ component of $X$ 
(see \ref{p.78}) a term which can only be bounded in terms of:
\begin{eqnarray}
&&C\int_{{\cal R}_{\ub_1,u_1}}\ub u^{-1}|\s^{(Y;0,\nu_1...,\nu_l)}\xi_{\Lb}||E_{\nu_l}\cdot\s^{(\nu_1...\nu_{l-1})}\theta_l|\nonumber\\
&&\leq C\int_0^{\ub_1}\|\s^{(Y;0,\nu_1...\nu_l)}\xi_{\Lb}\|_{L^2(\Cb_{\ub}^{u_1})}
\|E_{(\nu_l)}\cdot\s^{(\nu_1...\nu_{l-1})}\theta_l\|_{L^2(\Cb_{\ub}^{u_1})}d\ub 
\label{p.131}
\end{eqnarray}
Substituting for $\|\s^{(\nu_1...\nu_{l-1})}\theta_l\|_{L^2(\Cb_{\ub}^{u_1})}$ the leading 
contribution \ref{p.125} we see that \ref{p.131} can in turn only be bounded in terms of:
\begin{equation}
C\sum_\mu\int_0^{\ub_1}\left(\s^{(Y;0,\nu_1...\nu_l)}\underline{{\cal E}}^{u_1}(\ub)\right)^{1/2}
\left\{\frac{1}{\ub}\int_0^{\ub}\s^{(Y;0,\nu_1...\nu_{l-1}\mu)}\underline{{\cal E}}^{u_1}(\ub^\prime)
d\ub^\prime\right\}^{1/2}\frac{d\ub}{\ub}
\label{p.132}
\end{equation}
a singular integral. 

However, this is the lesser of the difficulties one faces in regard to \ref{p.93} with $m=0$. 
The greater difficulty arises in estimating the contribution through the factor 
$V^\mu\s^{(0,\nu_1...\nu_l)}\tilde{\rho}_\mu$ in the integrant of the terms
\begin{eqnarray}
&&\frac{1}{2}\rhob(V^\mu L\beta_\mu)E_{(\nu_l)}...E_{(\nu_1)}\mbox{tr}\tchib \nonumber\\
&&+\frac{a}{2c}h_{\nu_1,\kappa}N^\kappa V^\mu(\sd\beta_\mu)^\sharp\cdot
\sd(E_{(\nu_l)}...E_{(\nu_2)}\mbox{tr}\tchib) 
\label{p.133}
\end{eqnarray}
contributed by the terms involving $\mbox{tr}\tchib$ in \ref{p.96}. For, the 2nd of \ref{p.105} must be integrated from 
${\cal K}$ where a boundary condition holds which equates 
$r\sd (E_{(\nu_{l-1})}...E_{(\nu_1)}\mbox{tr}\tchib)$ to 
$\sd(E_{(\nu_{l-1})}...E_{(\nu_1)}\mbox{tr}\tchi)$
to leading terms (see paragraph following \ref{p.21}). In view of the definitions \ref{p.109} and 
\ref{p.21} this takes the form of a relation between $r^2\s^{(\nu_1...\nu_{l-1})}\thetab_l$ 
and $\s^{(\nu_1...\nu_{l-1})}\theta_l$ on ${\cal K}$. Integrating the 2nd of \ref{p.105} from 
${\cal K}$, the contribution of the term 
$$-2\lambdab^{-1}(\Lb\lambdab)\sd(E_{(\nu_{l-1})}...E_{(\nu_1)}\fb)$$
from the 2nd of \ref{p.106} to 
$\|\s^{(\nu_1...\nu_{l-1})}\thetab_l\|_{L^2(S_{\ub,u_1})}$ is, in view of \ref{p.113}, \ref{p.114}, 
bounded by:
\begin{equation}
C\int_{\ub}^{u_1}\|\sd(E_{(\nu_{l-1})}...E_{(\nu_1)}\fb)\|_{L^2(S_{\ub,u})}du
\label{p.134}
\end{equation}
(compare with \ref{p.115}). This causes no difficulty. The difficulty arises from the boundary term 
$\|\s^{(\nu_1...\nu_{l-1})}\thetab_l\|_{L^2(S_{\ub,\ub})}$. Taking account of the fact that 
$r\sim\tau$ along ${\cal K}$, we multiply the inequality for 
$\|\s^{(\nu_1...\nu_{l-1})}\thetab_l\|_{L^2(S_{\ub,u_1})}$  by $\ub^2$ and take the $L^2$ norm with 
respect to $\ub$ on $[0,\ub_1]$ to seek a bound for 
$\|\ub^2\s^{(\nu_1...\nu_{l-1})}\thetab_l\|_{L^2(C_{u_1}^{\ub_1})}$. The contribution of the 
boundary term to this is bounded by a constant multiple of 
$\|r^2\s^{(\nu_1...\nu_{l-1})}\thetab_l\|_{L^2({\cal K}^{\ub_1})}$, which by the above is bounded 
in terms of:
\begin{equation}
\|\s^{(\nu_1...\nu_{l-1})}\theta_l\|_{L^2({\cal K}^{\ub_1})}
\label{p.135}
\end{equation}
To estimate this we must revisit the argument leading from \ref{p.122} to \ref{p.125}. Here we are 
to set $\ub_1=u$ in \ref{p.122} to obtain the leading contribution to 
$\|\s^{(\nu_1...\nu_{l-1})}\theta_l\|_{L^2(S_{u,u})}$. The corresponding contribution to 
$\|\s^{(\nu_1...\nu_{l-1})}\theta_l\|_{L^2({\cal K}^{\ub_1})}$ is then bounded by:
\begin{equation}
C\sum_{\nu_l}\left\{\int_0^{\ub_1}\left(u^{-3}\int_0^u\|\s^{(Y;0,\nu_1...\nu_{l-1}\nu_l)}\xi_{\Lb}
\|_{L^2(S_{\ub,u})}\ub d\ub\right)^2 du\right\}^{1/2}
\label{p.136}
\end{equation}
The integral in the square root is now bounded by (1/3 times):
\begin{eqnarray}
&&\int_{\ub_1}^{u_1}u^{-3}\left(\int_0^{\ub_1}
\|\s^{(Y;0,\nu_1...\nu_{l-1}\nu_l)}\xi_{\Lb}\|^2_{L^2(S_{\ub,u})}d\ub\right)du\nonumber\\
&&=\int_0^{\ub_1}\left(\int_{\ub_1}^{u_1} u^{-3}
\|\s^{(Y;0,\nu_1...\nu_{l-1}\nu_l)}\xi_{\Lb}\|^2_{L^2(S_{\ub,u})}du\right)d\ub\nonumber\\
&&\leq \int_0^{\ub_1}\ub^{-3}\|\s^{(Y;0,\nu_1...\nu_{l-1}\nu_l)}\xi_{\Lb}\|^2_{L^2(\Cb_{\ub}^{u_1})}
d\ub\nonumber\\
&&\leq C\int_0^{\ub_1}\ub^{-3}\s^{(Y;0,\nu_1...\nu_{l-1}\nu_l)}\underline{{\cal E}}^{u_1}(\ub)d\ub 
\label{p.137}
\end{eqnarray}
and \ref{p.136} is bounded by:
\begin{equation}
C\sum_{\nu_l}\left\{\int_0^{\ub_1}\s^{(Y;0,\nu_1...\nu_{l-1}\nu_l)}\underline{{\cal E}}^{u_1}(\ub)
\frac{d\ub}{\ub^3}\right\}^{1/2}
\label{p.138}
\end{equation}
This is not only singular but, placing a lower cutoff $\varepsilon>0$ on $\ub$, blows up like 
$\varepsilon^{-1}$ as $\varepsilon\rightarrow 0$. Moreover, this only bounds 
$\|\ub^2\s^{(\nu_1...\nu_{l-1})}\thetab_l\|_{L^2(C_{u_1}^{\ub_1})}$ rather than 
$\|\s^{(\nu_1...\nu_{l-1})}\thetab_l\|_{L^2(C_{u_1}^{\ub_1})}$. Nevertheless, from the point of 
view of scaling, the singular bound \ref{p.138} is similar to the bound \ref{p.125}. Also, 
in view of \ref{p.113}, \ref{p.114}, the contribution to \ref{p.133} of 
$\s^{(\nu_2...\nu_l)}\thetab_l$ is pointwise bounded by:
\begin{equation}
Cu^2|\s^{(\nu_2...\nu_l)}\thetab_l|
\label{p.139}
\end{equation}
as compared with \ref{p.130} which is pointwise bounded by:
\begin{equation}
C\ub u^{-1}|\s^{(\nu_1...\nu_{l-1})}\theta_l|
\label{p.140}
\end{equation}
As a consequence the two contributions to $\s^{(Y;0,\nu_1...\nu_l)}{\cal G}_3^{\ub_1,u_1}$ are similar 
from the point of view of scaling. 

Similar results are obtained in regard to the error integral $\s^{(V;m,\nu_1...\nu_l)}
{\cal G}_3^{\ub_1,u_1}$ for $m\geq 1$, with $E_{(\nu_l)}...E_{(\nu_1)}T^{m-1}\slap\lambda$, 
$E_{(\nu_l)}...E_{(\nu_1)}T^{m-1}\slap\lambdab$ playing the roles of 
$\sd(E_{(\nu_{l-1})}...E_{(\nu_1)}\mbox{tr}\tchi)$, $\sd(E_{(\nu_{l-1})}...E_{(\nu_1)}\mbox{tr}\tchib)$ 
respectively. (Compare \ref{p.97} with \ref{p.96}.) In fact, integrating the 1st of \ref{p.111} from 
$\Cb_0$, the contribution of the term 
$$2\lambda^{-1}(L\lambda)E_{(\nu_l)}...E_{(\nu_1}T^{m-1}j$$
from the 1st of \ref{p.112} to $\|\s^{(\nu_1...\nu_l)}\nu_{m-1,l+1}\|_{L^2(S_{\ub_1,u})}$ is 
bounded by:
\begin{equation}
Cu^{-2}\int_0^{\ub_1}\|E_{(\nu_l)}...E_{(\nu_1)}T^{m-1}j\|_{L^2(S_{\ub,u})}\ub d\ub
\label{p.141}
\end{equation}
From the expression for the function $j$ of Proposition 11.2 we see that the term 
$$\frac{1}{4}\beta_N^2\Lb^2 H$$
makes the leading contribution. More precisely, writing $\Lb^2 H=\Lb TH-\Lb LH$ the leading 
contribution is that of
$$\frac{1}{4}\beta_N^2\Lb TH$$
and this contribution to \ref{p.141} is bounded, up to lower order terms, by:
\begin{equation}
Cu^{-2}\int_0^{\ub_1}\|\Nb^\mu\Lb\s^{(m,\nu_1...\nu_l)}\beta_{\mu}\|_{L^2(S_{\ub,u})}\ub d\ub
\label{p.142}
\end{equation}
Using the fact that 
\begin{equation}
\s^{(V;m,\nu_1...\nu_l)}\underline{{\cal E}}^{u_1}(\ub)
=\int_{\Cb_{\ub}^{u_1}}\Omega^{(d-1)/2}\left((\s^{(V;m,\nu_1...\nu_l)}\xi_{\Lb})^2
+3a|\s^{(V;m,\nu_1...\nu_l)}\sxi|^2\right)
\label{p.143}
\end{equation}
we deduce, in analogy with \ref{p.122}, that the leading contribution to 
$\|\s^{(\nu_1...\nu_l)}\nu_{m-1,l+1}\|_{L^2(S_{\ub_1,u})}$ is bounded by:
\begin{equation}
Cu^{-3}\int_0^{\ub_1}\|\s^{(Y;m,\nu_1...\nu_l)}\xi_{\Lb}\|_{L^2(S_{\ub,u})}\ub d\ub
\label{p.144}
\end{equation}
By an argument analogous to that leading from \ref{p.122} to \ref{p.125} we then deduce that the 
leading contribution to $\|\s^{(\nu_1...\nu_l)}\nu_{m-1,l+1}\|_{L^2(\Cb_{\ub_1}^{u_1})}$ is bounded by: 
\begin{equation}
C\ub_1^{-3/2}\left\{\int_0^{\ub_1}\s^{(Y;m,\nu_1...\nu_l)}\underline{{\cal E}}^{u_1}(\ub)d\ub\right\}^{1/2}
\label{p.145}
\end{equation}
Then in analogy with \ref{p.132} the corresponding contribution to the error integral 
$\s^{(Y;m,\nu_1...\nu_l)}{\cal G}_3^{\ub_1,u_1}$ is bounded in terms of the mildly singular 
integral:
\begin{equation}
C\int_0^{\ub_1}\left(\s^{(Y;m,\nu_1...\nu_l)}\underline{{\cal E}}^{u_1}(\ub)\right)^{1/2}
\left\{\frac{1}{\ub}\int_0^{\ub}\s^{(Y;m,\nu_1...\nu_l)}\underline{{\cal E}}^{u_1}(\ub^\prime)
d\ub^\prime\right\}^{1/2}\frac{d\ub}{\ub}
\label{p.146}
\end{equation}

On the other hand, the contribution through the factor $V^\mu \s^{(m,\nu_1...\nu_l)}\tilde{\rho}_\mu$ 
in the integrant in \ref{p.93} for $m\geq 1$, of the term
\begin{equation}
\rhob(V^\mu L\beta_\mu)E_{(\nu_l)}...E_{(\nu_1)}T^{m-1}\slap\lambdab
\label{p.147}
\end{equation}
contributed by the term involving $\slap\lambdab$ in \ref{p.97}, can only be bounded in terms of a 
severely singular integral. This is because the 2nd of \ref{p.111} must be integrated from ${\cal K}$ 
where a boundary condition which equates $rE_{(\nu_l)}...E_{(\nu_1)}T^{m-1}\slap\lambdab$ to 
$E_{(\nu_l)}...E_{(\nu_1)}T^{m-1}\slap\lambda$ holds, a differential consequence of the 
boundary condition \ref{p.21}. In view of the definitions \ref{p.110} and \ref{p.21} this takes the 
form of a relation between $r^2\s^{(\nu_1...\nu_l)}\nub_{m-1,l+1}$ and 
$\s^{(\nu_1...\nu_l)}\nu_{m-1,l+1}$ on ${\cal K}$.  The contribution of the 
boundary term to $\|\ub^2\s^{(\nu_1...\nu_l)}\nub_{m-1,l+1}\|_{L^2(C_{u_1}^{\ub_1})}$ is 
then bounded by a constant multiple of 
$\|r^2\s^{(\nu_1...\nu_l)}\nub_{m-1,l+1}|_{L^2({\cal K}^{\ub_1})}$, which in view of this relation 
is bounded in terms of:
\begin{equation}
\|\s^{(\nu_1...\nu_l)}\nu_{m-1,l+1}\|_{L^2({\cal K}^{\ub_1})}
\label{p.148}
\end{equation}
Proceeding as in the argument leading to \ref{p.138} we conclude that the leading contribution to 
\ref{p.148} is only bounded in terms of the severely singular integral:
\begin{equation}
C\left\{\int_0^{\ub_1}\s^{(Y;m,\nu_1...\nu_{l-1}\nu_l)}\underline{{\cal E}}^{u_1}(\ub)
\frac{d\ub}{\ub^3}\right\}^{1/2}
\label{p.149}
\end{equation}

Additional difficulties of the same kind as those discussed above arise 
indirectly in estimating the contribution of the terms involving the top order (order $n+1$) derivatives of the transformation functions to the integral which appears as the last term on the right 
in \ref{p.81} when $\s^{(V)}b$ is replaced by $\s^{(V;m,l)}b$, $\s^{(V;m,\nu_1...\nu_l)}b$  
according to \ref{p.88}, \ref{p.89}.  

The new analytic method is designed to overcome the difficulties due to the appearance of the 
singular integrals. The starting point is the observation that, in view of the fact that these 
integrals are borderline from the point of view of scaling, they would become regular borderline 
integrals if the energies $\s^{(V;m,\nu_1....\nu_l)}{\cal E}^{\ub}(u)$, 
$\s^{(V;m,\nu_1...\nu_l)}\underline{{\cal E}}^{u}(\ub)$ had a growth from the singularities at 
$\ub=0$ ($\Cb_0$) and  $u=0$ ($S_{0,0}$) like $\ub^{2a_m}u^{2b_m}$ for some sufficiently large 
exponents $a_m$ and $b_m$, and, moreover, the flux $\s^{(V;m,\nu_1...\nu_l)}{\cal F}^{\prime \tau}$ 
had a growth from the singularity at $\tau=0$ ($S_{0,0}$) on ${\cal K}$ like $\tau^{2c_m}$, 
where $c_m=a_m+b_m$. This observation seems at first sight irrelevant since the required growth properties cannot hold for these quantities. However, given that the initial data of the problem 
are expressed as smooth functions of $(u,\vartheta)$, as we have seen above in our review of Chapter 5, 
the derived data, that is the $T$ derivatives on $\Cb_0$ of up to any desired order $N$ of the 
unknowns $((x^\mu : \mu=0,...,n) , \ b \ , (\beta_\mu : \mu=0,...,n))$ in the characteristic and 
wave systems, are determined as smooth functions of $(u,\vartheta)$. The $T$ derivatives 
of the transformation functions $(v,\gamma)$ on $S_{0,0}$ of up to order $N-2$ 
are also determined at the same time as smooth functions of $\vartheta$. Therefore we can define 
the $N$th approximants $((x_N^\mu : \mu=0,...,n) , \ b_N \ , (\beta_{\mu,N} : \mu=0,...,n))$ 
as the corresponding $N$th degree polynomials in $\tau=\ub$ with coefficients which are known 
smooth functions of $(\sigma=u-\ub,\vartheta)$. Similarly we can define the $N$th approximants 
$(v_N,\gamma_N)$ as the corresponding $N-2$ th degree polynomials in $\tau$ with coefficients 
which are known smooth functions of $\vartheta$. If the $N$ approximants, so defined, are inserted 
into the equations of the characteristic and wave systems, these equations fail to be satisfied 
by errors which are known smooth functions of $(\ub,u,\vartheta)$ and whose derivatives up to 
order $n$ are bounded by a known constant times $\tau^{N-n+k}$ where $k$ is a fixed integer depending 
on which equation we are considering. Similarly, inserting the $N$ approximants into the boundary 
conditions and into the identification equations, these fail to be satisfied by errors which are 
known smooth functions of $(\tau,\vartheta)$ and whose derivatives up to order $n$ are likewise 
bounded by a known constant times $\tau^{N-n+k}$ where $k$ is a fixed integer depending 
on which equation we are considering. Moreover, in connection with the equation 
$\square_{\tilde{h}}\beta_\mu=0$ satisfied by an actual solution, the corresponding $N$th 
approximant quantity $\square_{\tilde{h}_N}\beta_{\mu,N}$ is a known smooth function of 
of $(\ub,u,\vartheta)$ and whose derivatives up to 
order $n$ are bounded by a known constant times $\tau^{N-n+k}$ where $k$ is a fixed integer. 
The above are discussed in detail in Chapter 9, in the case $d=2$, the extension to the case $d>2$ 
being straightforward. Note that in the above construction we have, in terms of 
$(\tau=\ub, \sigma=u-\ub,\vartheta)$ coordinates:
\begin{equation}
T_N=\frac{\partial}{\partial\tau}=T, \ \ \ \Omega_N=\frac{\partial}{\partial\vartheta}=\Omega
\label{p.150}
\end{equation}
independently of the approximation. On the other hand, $L_N$, $\Lb_N$, $E_N$ depend on the 
approximation, the first two through $b_N$, and the last through $\sh_N$, where:
\begin{equation}
\sh_N=h_{\mu\nu,N}(\Omega x_N^\mu)(\Omega x_N^\nu)
\label{p.151}
\end{equation}
We then define the difference quantities:
\begin{equation}
\s^{(m,l)}\check{\beta}_\mu=E^l T^m \beta_\mu-E_N^l T^m \beta_{\mu,N}
\label{p.152}
\end{equation}
\begin{equation}
\s^{(V;m,l)}\cxi=V^\mu d\s^{(m,l)}\check{\beta}_\mu 
\label{p.153}
\end{equation}
We also define $\s^{(V;m,l)}\check{S}$ as in \ref{p.44}, \ref{p.45} with $\s^{(V;m,l)}\cxi$ 
in the role of $\s^{(V)}\xi$, that is:
\begin{equation}
\s^{(V;m,l)}\check{S}=\tilde{h}^{-1}\cdot\s^{(V;m,l)}\check{S}_\flat
\label{p.154}
\end{equation}
where: 
\begin{equation}
\s^{(V;m,l)}\check{S}_\flat=\s^{(V;m,l)}\cxi\otimes\s^{(V;m,l)}\cxi
-\frac{1}{2}(\s^{(V;m,l)}\cxi,\s^{(V;m,l)}\cxi)_h h
\label{p.155}
\end{equation}
Defining moreover, in analogy with \ref{p.53}:
\begin{equation}
\s^{(V;m,l)}\check{P}=-\s^{(V;m,l)}\check{S}\cdot X
\label{p.156}
\end{equation}
we have, in analogy with \ref{p.54} - \ref{p.58}:
\begin{equation}
\mbox{div}_{\tilde{h}} \s^{(V;m,l)}\check{P}=\s^{(V;m,l)}\check{Q}
\label{p.157}
\end{equation}
where:
\begin{equation}
\s^{(V;m,l)}\check{Q}=\s^{(V;m,l)}\check{Q}_1+\s^{(V;m,l)}\check{Q}_2+\s^{(V;m,l)}\check{Q}_3
\label{p.158}
\end{equation} 
with:
\begin{eqnarray}
&&\s^{(V;m,l)}\check{Q}_1=-\frac{1}{2}\s^{(V;m,l)}\check{S}^\sharp\cdot\s^{(X)}\tilde{\pi} 
\label{p.159}\\
&&\s^{(V;m,l)}\check{Q}_1=-(\s^{(V;m,l)}\cxi,\s^{(V)}\theta^\mu)_{\tilde{h}}
X\s^{(m,l)}\check{\beta}_\mu 
+(\s^{(V;m,l)}\cxi,d\s^{(m,l)}\check{\beta}_\mu)_{\tilde{h}}\s^{(V)}\theta^\mu(X) \nonumber\\
&&\hspace{15mm}-(\s^{(V)}\theta^\mu,d\s^{(m,l)}\check{\beta}_\mu)_{\tilde{h}}\s^{(V;m,l)}\cxi(X) \label{p.160}\\
&&\s^{(V)}\check{Q}_3=-\s^{(V;m,l)}\cxi(X) V^\mu \square_{\tilde{h}}\s^{(m,l)}\beta_\mu \label{p.161}
\end{eqnarray}
We then deduce, in analogy with \ref{p.81} the $(m,l)$ {\em difference energy identity}:
\begin{eqnarray}
&&\s^{(V;m,l)}\cE^{\ub_1}(u_1)+\s^{(V;m,l)}\cEb^{u_1}(\ub_1)+\s^{(V;m,l)}\cF^{\prime\ub_1}\nonumber\\
&&\hspace{18mm}=\s^{(V;m,l)}\cEb^{u_1}(0)+\s^{(V;m,l)}\cG^{\ub_1,u_1}
+2C^\prime\int_{{\cal K}^{\ub_1}}\Omega^{1/2}(\s^{(V;m,l)}b)^2 \nonumber\\
&&\label{p.162}
\end{eqnarray}
(here $d=2$), where (see \ref{p.63}): 
\begin{eqnarray}
&&\s^{(V;m,l)}\cE^{\ub_1}(u_1)=\int_{C_{u_1}^{\ub_1}}\Omega^{1/2}\s^{(V;m,l)}\check{S}_\flat(X,L)
\nonumber\\
&&\hspace{24mm}=\int_{C_{u_1}^{\ub_1}}\Omega^{1/2}\left(3(\s^{(V;m,l)}\cxi_L)^2+a(\s^{V;m,l)}\csxi)^2\right)\nonumber\\
&&\s^{(V;m,l)}\cEb^{u_1}(\ub_1)=\int_{\Cb_{\ub_1}^{u_1}}\Omega^{1/2}\s^{(V;m,l)}\check{S}_\flat(X,\Lb)
\nonumber\\
&&\hspace{24mm}=\int_{\Cb_{\ub_1}^{u_1}}\Omega^{1/2}\left(3a(\s^{(V;m,l)}\csxi)^2+(\s^{(V;m,l)}\cxi_{\Lb})^2\right)\label{p.163}
\end{eqnarray}
are the $(m,l)$ {\em difference energies}. Also, in \ref{p.162}, 
\begin{equation}
\s^{(V;m,l)}\cF^{\prime\ub_1}=\s^{(V;m,l)}\cF^{\ub_1}+2C^\prime\int_{{\cal K}^{\ub_1}}\Omega^{1/2}
(\s^{(V;m,l)}\check{b})^2
\label{p.164}
\end{equation}
is positive-definite (see \ref{9.a34}), 
$\s^{(V;m,l)}\cF^{\ub_1}$ being the $(m,l)$ {\em difference flux}:
\begin{equation}
\s^{(V;m,l)}\cF^{\ub_1}=\int_{{\cal K}^{\ub_1}}\Omega^{1/2}\s^{(V;m,l)}\check{S}_\flat(X,M)
\label{p.165}
\end{equation}
and $\s^{(V;m,l)}\check{b}$ representing the boundary values on ${\cal K}$ of $\s^{(V;m,l)}\cxi(B)$ 
(see \ref{13.1}). Finally, in \ref{p.162}, 
\begin{equation}
\s^{(V;m,l)}\check{{\cal G}}^{\ub_1,u_1}:=\int_{{\cal R}_{\ub_1,u_1}}2a\Omega^{3/2}
\s^{(V;m,l)}\check{Q}
\label{p.166}
\end{equation}
is the $(m,l)$ {\em difference error integral}. 

Now, from the definition \ref{p.152} and the preceding discussion it follows that for any solution 
of problem the functions $\s^{(m,l)}\check{\beta}_\mu$ vanish with all their $T$ derivatives up to order $n$ 
on $\Cb_0$ if we choose $N\geq n$. We then have:
\begin{equation}
\s^{(V;m,l)}\cEb^{u_1}(0) \ : \ \mbox{for all $m=0,...,n$}
\label{p.167}
\end{equation}
therefore the $(m,l)$ {\em difference energy identity} simplifies to:
\begin{eqnarray}
&&\s^{(V;m,l)}\cE^{\ub_1}(u_1)+\s^{(V;m,l)}\cEb^{u_1}(\ub_1)+\s^{(V;m,l)}\cF^{\prime\ub_1}\nonumber\\
&&\hspace{23mm}=\s^{(V;m,l)}\cG^{\ub_1,u_1}
+2C^\prime\int_{{\cal K}^{\ub_1}}\Omega^{1/2}(\s^{(V;m,l)}b)^2 \nonumber\\
&&\label{p.168}
\end{eqnarray}
Expecting that the growth of $\s^{(V;m,l)}\cE^{\ub}(u)$, $\s^{(V;m,l)}\cEb^{u}(\ub)$ is like 
$\ub^{2a_m}u^{2b_m}$ and that the growth of $\s^{(V;m,l)}\cF^{\prime\tau}$ is like $\tau^{2c_m}$, 
we define the weighted quantities:
\begin{equation}
\s^{(V;m,l)}\cB(\ub_1,u_1)=\sup_{(\ub,u)\in R_{\ub_1,u_1}}\ub^{-2a_m}u^{-2b_m}
\s^{(V;m,l)}\cE^{\ub}(u)
\label{p.169}
\end{equation}
\begin{equation}
\s^{(V;m,l)}\cBb(\ub_1,u_1)=\sup_{(\ub,u)\in R_{\ub_1,u_1}}\ub^{-2a_m}u^{-2b_m}
\s^{(V;m,l)}\cEb^u(\ub)
\label{p.170}
\end{equation}
and:
\begin{equation}
\s^{(V;m,l)}\cA(\tau_1)=\sup_{\tau\in[0,\tau_1]}\tau^{-2(a_m+b_m)}\s^{(V;m,l)}\cF^{\prime\tau}
\label{p.171}
\end{equation}
the exponents $a_m$, $b_m$ being non-negative real numbers which are non-increasing with $m$. 
Of course the above definitions do not make sense unless we already know that the quantities 
$\s^{(V;m,l)}\cE^{\ub}(u)$, $\s^{(V;m,l)}\cEb^{u}(\ub)$, $\s^{(V;m,l)}\cF^{\prime\tau}$ have 
the appropriate growth properties. Making this assumption would introduce a vicious circle into  
the argument, so this is not what we do. What we actually do is presented in Section 14.10. 
There we regularize the problem by giving initial data not on $\Cb_0$ but on $\Cb_{\tau_0}$ 
for $\tau_0>0$ but not exceeding a certain fixed positive number which is much smaller than $\delta$. 
The initial data on $\Cb_{\tau_0}$ 
is modeled after the restriction to $\Cb_{\tau_0}$of the $N$th approximate solution, the difference being 
bounded by a fixed constant times $\tau_0^{N-1}$. Similarly considering the $m$th derived data 
on $\Cb_{\tau_0}$, $m=1,...,n+1$ we show that the difference from the corresponding 
$N$th approximants on $\Cb_{\tau_0}$ is bounded by a fixed constant times $\tau_0^{N-1-m}$. 
The $(m,l)$ difference energy identity now refers to the domain ${\cal R}_{\ub_1,u_1,\tau_0}$ in 
${\cal N}$ which corresponds  to the domain:
\begin{equation}
R_{\ub_1,u_1,\tau_0}=\{(\ub,u) \ : \ u\in[\ub,u_1], \ \ub\in[\tau_0,\ub_1]\}
\label{p.172}
\end{equation}
in $\mathbb{R}^2$, so the 1st term on the right in \ref{p.162} is replaced by:
\begin{equation}
\s^{(V;m,l)}\cEb^{u_1}(\tau_0)\leq C\tau_0^{2(N-1-m)}
\label{p.173}
\end{equation}
Moreover, in reference to \ref{p.172}, if $\ub_1-\tau_0$ is suitably small then given a solution 
defined in ${\cal R}_{\ub_1,u_1,\tau_0}$ of the problem with initial data on $\Cb_{\tau_0}$, the quantities $\s^{(V;m,l)}\cEb^{u}(\ub)$, $\s^{(V;m,l)}\cE^{\ub}(u)$, and $\s^{(V;m,l)}\cF$, for 
$(\ub,u)\in R_{\ub_1,u_1,\tau_0}$, are likewise bounded by a constant times 
$\tau_0^{2(N-1-m)}$. Therefore replacing the supremum over $R_{\ub_1,u_1}$, the supremum over 
$[0,\tau_1]$, in the definitions \ref{p.170}, \ref{p.171}, respectively, and taking $N\geq m+1+c_m$, 
everything now makes sense. In fact, we take $N>m+\frac{5}{2}+c_m$, in which case the modifications 
in the resulting estimates for the quantities \ref{p.170}, \ref{p.171}, tend to 0 as 
$\tau_0\rightarrow 0$. The estimates in the preceding chapters, namely Chapters 10, 12, and 13, 
concluding with the top order energy estimates in Chapter 13, are thus derived with the foreknowledge 
that such modifications involving a very small positive $\tau_0$ will be made which will tend to 0 
as $\tau_0\rightarrow 0$. We have chosen not to introduce such a small positive $\tau_0$ from 
the beginning, because that would have made the exposition more complicated lengthening  
unnecessarily the monograph. 

The argument relies on the derivation of energy estimates of the 
top order $m+l=n$ only. Since these are derived in Chapter 13 before the regularization of the problem by
the introduction of the small positive $\tau_0$, they are to be thought of as {\em a priori estimates}. 
In our exposition we at first ignore all but the principal terms involved. Thus in Chapters 10 and 11, 
the aim of which is to derive estimates for the top order acoustical difference quantities corresponding to \ref{p.104}, \ref{p.110}:
\begin{eqnarray}
&&\hspace{8mm}\s^{(\nu_1,...,\nu_{l-1})}\cth_l=\s^{(\nu_1...\nu_{l-1})}\theta_l-\s^{(\nu_1...\nu_{l-1})}\theta_{l,N}\nonumber\\
&&\hspace{8mm}\s^{(\nu_1,...,\nu_{l-1})}\cthb_l=\s^{(\nu_1...\nu_{l-1})}\thetab_l-\s^{(\nu_1...\nu_{l-1})}\thetab_{l,N}\nonumber\\
&&\s^{(\nu_1...\nu_l)}\cnu_{m-1,l+1}=\s^{(\nu_1...\nu_l)}\nu_{m-1,l+1}
-\s^{(\nu_1...\nu_l)}\nu_{m-1,l+1,N}\nonumber\\
&&\s^{(\nu_1...\nu_l)}\cnub_{m-1,l+1}=\s^{(\nu_1...\nu_l)}\nub_{m-1,l+1}
-\s^{(\nu_1...\nu_l)}\nub_{m-1,l+1,N}
\label{p.174}
\end{eqnarray}
for $d>2$, 
\begin{eqnarray}
&&\hspace{30mm}\cth_l=\theta_l-\theta_{l,N}, \ \ \ \cthb_l=\thetab_l-\thetab_{l,N}\nonumber\\
&&\cnu_{m-1,l+1}=\nu_{m-1,l+1}-\nu_{m-1,l+1,N}, \ \ \cnub_{m-1,l+1}=\nub_{m-1,l+1}-\nub_{m-1,l+1,N}
\nonumber\\
&&\label{p.175}
\end{eqnarray}
in the case $d=2$, we ignore all but the top order terms, namely those of order $n+1$. 
An exception occurs in regard to the terms in the boundary estimates for $\cnub_{m-1,l+1}$ 
involving the transformation function differences $\chf$, $\cv$, $\cga$ (see Proposition 10.9), 
where the terms of orders 2, 1, 3 lower must be kept because they enter the estimates multiplied 
by correspondingly lower powers of $\tau$. The only essential change when passing from the case $d=2$ 
to $d>2$, is that in higher dimensions we must also derive $L^2(S_{\ub,u})$ estimates for 
$\lambda\sD\s^{(\nu_1...\nu_{l-1})}\ctchi_{l-1}$, $\lambdab\sD\s^{(\nu_1...\nu_{l-1})}\ctchib_{l-1}$ 
in terms of $L^2(S_{\ub,u})$ estimates for $\s^{(\nu_1,...,\nu_{l-1})}\cth_l$, 
$\s^{(\nu_1,...,\nu_{l-1})}\cthb_l$ respectively, where:
\begin{eqnarray}
&&\s^{(\nu_1 ... \nu_{l-1})}\ctchi_{l-1}=\sL_{E_{(\nu_{l-1})}} ... \sL_{E_{(\nu_1)}}\tchi-
\sL_{E_{(\nu_{l-1}),N}} ... \sL_{E_{(\nu_1),N}}\tchi_N\nonumber\\
&&\s^{(\nu_1 ... \nu_{l-1})}\ctchib_{l-1}=\sL_{E_{(\nu_{l-1})}} ... \sL_{E_{(\nu_1)}}\tchib-
\sL_{E_{(\nu_{l-1}),N}} ... \sL_{E_{(\nu_1),N}}\tchib_N\nonumber\\
&&\label{p.176}
\end{eqnarray}
The required estimates are derived using elliptic theory in connection with the Codazzi equations
on the $S_{\ub,u}$. We must also derive $L^2(S_{\ub,u})$ estimates for 
$\lambda\sD^2\s^{(\nu_1...\nu_l)}\cla_{m-1,l}$, $\lambdab\sD^2\s^{(\nu_1...\nu_l)}\clab_{m-1,l}$ 
in terms of $L^2(S_{\ub,u})$ estimates for $\s^{(\nu_1...\nu_l)}\cnu_{m-1,l+1}$, 
$\s^{(\nu_1...\nu_l)}\cnu_{m-1,l+1}$ respectively, where:
\begin{eqnarray}
&&\s^{(\nu_1 ... \nu_l)}\cla_{m-1,l}=E_{(\nu_l)} ... E_{(\nu_1)}T^{m-1}\lambda
-E_{(\nu_l),N} ... E_{(\nu_1),N}T^{m-1}\lambda_N\nonumber\\
&&\s^{(\nu_1 ... \nu_l)}\clab_{m-1,l}=E_{(\nu_l)} ... E_{(\nu_1)}T^{m-1}\lambdab
-E_{(\nu_l),N} ... E_{(\nu_1),N}T^{m-1}\lambdab_N\nonumber\\
&&\label{p.177}
\end{eqnarray}
The required estimates are derived using elliptic theory for the Laplacian on $S_{\ub,u}$. 

The aim of Chapter 12 is to derive estimates for the derivatives of top order, $n+1$, of the transformation function differences $\chf$, $\cv$, $\cga$, and at the same time more precise 
estimates for the next to top order acoustical differences 
$\s^{(l-1)}\tchi, \s^{(l-1)}\tchib \ : \ l=n$, and 
$\s^{(m,l)}\cla, \s^{(m,l)}\clab \ : \ m+l=n$, in terms of the top order difference energies and fluxes. 
Here all terms of order up to $n-1$ are ignored and of the terms of order $n$ only those involving 
acoustical quantities of order $n$, namely the quantities being estimated, are considered. 

In Chapter 13, which contains the derivation of the top order energy estimates, we again ignore all 
but the top order terms. 

The estimates of Chapters 10 - 13 require taking the exponents $a_m$, $b_m$, $c_m$ to be 
suitably large. This is in accordance with our preceding heuristic discussion. Moreover the 
estimates of Chapters 12 and 13, in particular the top order energy estimates, require that 
$\delta$ does not exceed a positive constant which is independent of $m$ and $n$. 

In the above only the treatment of the principal terms is shown, and these are estimated using 
only the fundamental bootstrap assumptions. The full treatment, which includes all the lower order 
terms, uses the complete set of {\em bootstrap assumptions} stated in Section 14.1. In the same 
section it is shown how all the lower order terms are treated. We return to the treatment of the 
lower order terms in Section 14.8 where it is shown that the lower order terms are treated in an optimal 
manner if we choose all the exponents $a_m$ equal and all the exponents $b_m$ equal:
\begin{equation}
a_m=a, \ \ \ b_m=b \ \ \ : \ m=0,...,n
\label{p.178}
\end{equation}
and this optimal choice is taken from that point on. In Section 14.4 $L^2(S_{\ub,u})$ estimates 
for the $\s^{(m,l)}\beta_\mu : m+l=n$ are deduced. In Section 14.5 $L^2(S^{\ub,u})$ estimates 
for the $n-1$th order acoustical differences are deduced. In Section 14.6 $L^2(S_{\ub,u})$ estimates 
for all $n$th order derivatives are deduced. In Section 14.9 pointwise estimates are deduced and 
the bootstrap assumptions are recovered as strict inequalities. This recovery however requires 
a further smallness condition on $\delta$ which now depends on $n$. The derivation of pointwise estimates from 
$L^2(S_{\ub,u})$ estimates of derivatives intrinsic to the $S_{\ub,u}$ requires $k_0$ derivatives 
where $k_0$ is the smallest integer greater than $(d-1)/2$, so in the case $d=2$ we have $k_0=1$. 
We then first set:
\begin{equation}
n=n_0:=2k_0+2
\label{p.179}
\end{equation}
The smallness condition on $\delta$ needed to recover the bootstrap assumptions is then that 
corresponding to $n=n_0$. Given then any $n>n_0$, the nonlinear argument having closed at order $n_0$, 
the bootstrap assumptions are no longer needed, therefore no new smallness conditions on $\delta$ 
are required to proceed inductively to orders $n_0+1,...,n$. We remark that $n_0$ can be lowered 
to:
\begin{equation}
n_0=k_0+1
\label{p.180}
\end{equation}
if our treatment of the lower order terms, which uses only $L^\infty(S_{\ub,u})$ and $L^2(S_{\ub,u})$ estimates, 
is refined to make full use of the Sobolev inequalities on the $S_{\ub,u}$. 

In Section 14.1 we first regularize the problem by giving the initial data on $\Cb_{\tau_0}$  
as discussed above.
We then apply a continuity argument to establish the existence of a solution to this 
regularized problem defined on the whole of ${\cal R}_{\delta,\delta,\tau_0}$. This argument 
relies on the solution by Majda and Thomann [Ma-Th] of the restricted local shock continuation 
problem. This is the problem of continuing locally in time a solution displaying a shock discontinuity 
initially. Majda had in fact earlier addressed in his fundamental works [Ma1], [Ma2] the general shock development problem. 
Our continuity argument relies on [Ma-Th] at two points. First, to construct a solution on 
${\cal R}_{\tau_0+\vep,\tau_0+\vep, \tau_0}$ for some $\vep>0$, which could be 
very small, even much smaller than $\tau_0$. Second, at the end of the continuity argument, starting 
with a solution defined on a maximal ${\cal R}_{\ub_*,\ub_*,\tau_0}$, 
to extend a solution which has previously been shown to be defined on 
${\cal R}_{\ou_*,\ou_*+\vep,\tau_0}$ 
to the domain corresponding to the triangle: 
$$\{(\ub,u) \ : \ u\in[\ub,\ub+\vep_*], \ \ub\in[\ou_*,\ou_*+\vep_*]\}$$
for some $\vep_*>0$, not exceeding $\vep$.  Since the union contains 
${\cal R}_{\ub_*+\vep_*,\ub_*+\vep_*,\tau_0}$ we then arrive at a contradiction to the maximality of 
$\ou_*$, unless $\ou_*=\delta$. This double use of a local existence theorem is typical of 
continuity arguments. 

After obtaining a solution on ${\cal R}_{\delta,\delta,\tau_0}$ satisfying the appropriate estimates 
we take $\tau_0$ to be any member of a sequence $(\tau_{0,m} \ :\ m=M, M+1, M+2, \ . \ . \ . )$ converging to 0, 
and pass to the limit in a subsequence to obtain the solution to our problem. The theorem which 
concludes this monograph is stated at the end. 

We conclude this prologue with a few historical remarks to acknowledge the fact that the 
methods of this monograph have their roots in the past and that at best we simply develop things 
a little further in a particular direction. The method of continuity, using a local existence theorem 
together with a priori bounds, derived on the basis of bootstrap assumptions which are recovered in 
the course of the argument, originated in the field of ordinary differential equations, in particular 
the classical problem of the motion of $N$ point masses under their mutual gravitational attraction 
formulated by Newton in his {\em Principia} [Ne]. This corresponds to a Hamiltonian system with Hamiltonian 
function $H=K+V$ where
$$K=\sum_{\alpha=1}^N \frac{|p_\alpha|^2}{2M_\alpha}$$
is the kinetic energy, $p_\alpha$, being the momentum  of the mass labeled $\alpha$, 
a vector in Euclidean 3 dimensional space $\mathbb{E}^3$ at the point $x_\alpha\in\mathbb{E}^3$, the position of the mass $\alpha$, and $M_\alpha>0$ its mass. Also, 
$$V=-\frac{1}{2}\sum_{\alpha\neq\beta}\frac{GM_\alpha M_\beta}{r_{\alpha\beta}}$$
is the potential energy, $r_{\alpha\beta}$ being the distance between the masses $\alpha$ and $\beta$ 
and $G$ being Newton's gravitational constant. The potential energy is defined when the $x_\alpha$ 
are distinct. If it is regularized by replacing $r_{\alpha\beta}$ by $\sqrt{r^2_{\alpha\beta}+\vep^2}$ 
the simplest application of the method of continuity, without bootstrap assumptions, shows that 
for any initial condition we have a solution defined for all future time. For the actual problem 
the same method gives a basic theorem (see [Si-Mo]) stating that for any initial condition either 
we have a solution defined for all future time, or the maximal existence interval is $[0,t_*)$, 
in which case, with $\rho=\min_{\alpha\neq\beta}r_{\alpha\beta}$ the minimal mutual distance, 
we have $\rho(t)\rightarrow 0$ as $t\rightarrow t_*$, so in a sense at time $t_*$ we have a collision. 
To this day it is still not known whether in the $6N$ dimensional space of initial conditions 
the subset leading to a collision has zero measure or positive measure. Incidentally, the 
{\em Principia} is the work where power series were first introduced. 

Another example of the method of continuity in conjunction with a local existence theorem and 
a priori bounds is from the field of elliptic partial differential equations. This is the 
non-parametric minimal hypersurface problem: Given a domain in a hyperplane in $\mathbb{E}^{n+1}$ 
and a graph over the boundary, to find an extension to a graph over the interior which is of minimal 
$n$ dimensional area. The graph is that of the height function above the hyperplane. We can assume
that the mean value of the height function on the boundary vanishes. To apply the method of 
continuity we multiply the boundary height function, the data of the problem, by $t\in[0,1]$. 
Then for $t=0$ we have the trivial solution where the graph coincides with the hyperplane domain. 
For any $t_0\in[0,1]$ for which we have a regular solution, the implicit function theorem plays the role of a local existence theorem giving us a solution for $t$ suitably close to $t_0$. The problem 
then reduces to that of deriving the appropriate a priori bounds. For once these have been established,  
the method of continuity yields a solution for all $t\in[0,1]$, in particular for $t=1$, which is the 
original problem. The required a priori bounds where established in two steps. The first step 
was the derivation of gradient bounds on the height function and was done by Bernstein in 
his fundamental works [Be1], [Be2], [Be3]. The second step was the derivation of H\"{o}lder estimates 
for the gradient. This step was done half a century later, independently by De Giorgi [DG] and by
Nash [Na]. While De Giorgi's work was in the framework of minimal surface theory, Nash's motivation 
was actually a very different problem, the study of the evolution of a viscous, heat conducting 
compressible fluid, which is why he addressed an analogous problem for parabolic equations 
and deduced the result for elliptic equations as a corollary in the time independent case. 
The complete result for the original problem, which is that mean convexity of the boundary of 
the domain is both necessary and sufficient for the problem to be solvable for all data, was deduced 
later by Jenkins and Serrin [Je-Se]. 

A last example of the method of continuity, which uses a local existence theorem 
and a priori bounds which are derived on the basis of bootstrap assumptions recovered in 
the course of the argument, is from the field of hyperbolic partial differential 
equations, therefore closer to the topic of the present monograph. This is the work [Ch-Kl],  
mentioned above, on the problem of the stability 
of the Minkowski metric. In this case the local existence theorem had been established in 
the basic work of Choquet-Bruhat [CB]. We then consider the least upper bound $t_*$ of the set 
of all positive times $t_1$ for which we have a solution on the spacetime slab corresponding 
to the time interval $[0,t_1]$ satisfying the bootstrap assumptions. So if $t_*$ is finite, 
while the solution extends to the spacetime slab corresponding to $[0,t_*]$, the inequalities 
corresponding to the bootstrap assumptions are saturated on that slab. The a priori bounds  
are then derived through a construction on the slab in question which proceeds from the future 
boundary of the slab, which corresponds to the time $t_*$. This is because the derivation of the 
a priori bounds uses approximate symmetries, which a posteriori are shown to become exact in the limit 
$t_*\rightarrow\infty$. The a priori bounds show that the inequalities corresponding to the 
bootstrap assumptions are in fact not saturated on the slab in question. Then after another application 
of the local existence theorem we arrive at a contradiction to the maximality of $t_*$. 

\vspace{2.5mm}

\noindent{\bf Acknowledgement :} This work was supported by ERC Advanced Grant 246574, 
{\em Partial Differential Equations of Classical Physics.}

\pagebreak

\chapter{Fluid Mechanics and the Shock Development Problem}
\section{The General Equations of Motion}

The mechanics of a perfect fluid is described in the framework of the Minkowski spacetime of 
special relativity by a future-directed unit time-like vectorfield $u$, the fluid {\em spacetime 
velocity} and two positive functions $n$ and $s$, the {\em number of particles per unit volume} 
(in the local rest frame of the fluid) and the {\em entropy per particle}. In terms of a system 
of rectangular coordinates $(x^0,x^1,x^2,x^3)$, with $x^0$ a time coordinate and $(x^1,x^2,x^3)$ spatial coordinates, choosing, as is customary in the relativistic framework, the relation of the units of temporal to spatial lengths so as to set the universal constant $c$ represented by the speed 
of light in vacuum equal to 1, the components $g_{\mu\nu}$ of the Minkowski metric are given by: 
\begin{equation}
g_{00}=-1, \ g_{11}=g_{22}=g_{33}=1, \ \  g_{\mu\nu}=0 \ : \ \mbox{if $\mu\neq\nu$}
\label{1.1}
\end{equation}
The conditions on the components $u^\mu$ of the spacetime velocity are then:
\begin{equation}
g_{\mu\nu}u^\mu u^\nu=-1, \ \  u^0>0
\label{1.2}
\end{equation}
Here, and throughout this monograph, we follow the summation convention, according to which repeated 
upper and lower indices are summed over their range. 

The mechanical properties of a perfect fluid are specified once we give the {\em equation of state}, 
which expresses the relativistic {\em energy density} (in the local rest frame of the fluid) $\rho$ 
as a function of $n$ and $s$:
\begin{equation}
\rho=\rho(n,s)
\label{1.3}
\end{equation}
According to the laws of thermodynamics, the {\em pressure} $p$ and the {\em temperature} $\theta$ 
are then given by:
\begin{equation}
p=n\frac{\partial\rho}{\partial n}-\rho, \ \ \theta=\frac{1}{n}\frac{\partial\rho}{\partial s}
\label{1.4}
\end{equation}
In terms of the {\em volume per particle} 
\begin{equation}
V=\frac{1}{n}
\label{1.5}
\end{equation}
and the relativistic {\em energy per particle}
\begin{equation}
e=\rho V,
\label{1.6}
\end{equation}
these relations take the familiar form: 
\begin{equation}
e=e(V,s)
\label{1.7}
\end{equation}
\begin{equation}
de=-pdV+\theta ds
\label{1.8}
\end{equation}
We note that the relativistic energy per particle contains the rest mass contribution, $mc^2$ in 
conventional units, $m$ being the particle rest mass. Under ordinary circumstances this is in fact the dominant contribution to $e$. The corresponding contribution to $\rho$ is $nmc^2$ in conventional units, $nm$ being the rest mass density. When seeking to deduce the equation of state \ref{1.3} or 
\ref{1.7} for a given substance in a given phase from microscopic physics through statistical mechanics, by ``particle" we actually mean a microscopic constituent particle which has a 
definite rest mass $m$  when considered in isolation in free space. In fluid mechanics however we take the macroscopic point of view, and what we 
call ``particle" is an arbitrarily chosen but fixed number of microscopic constituent particles, therefore $m$ is the corresponding multiple of the microscopic rest mass. We may thus fix $m$ to be the unit of mass, so we can conveniently set $m=1$. Then all ``extensive" quantities, like $V$, $s$, 
and $e$, which we have referred to as quantities ``per particle" are actually quantities ``per unit 
rest mass", or ``specific" quantities in the customary terminology. Thus $V$, $s$, and $e$, are the 
{\em specific volume}, {\em specific entropy}, and {\em specific energy}, respectively, and 
$n$ coincides with the {\em rest mass density}. 

The positive function $h$ defined by: 
\begin{equation}
h=\frac{(\rho+p)}{n}
\label{1.9}
\end{equation}
or equivalently 
\begin{equation}
h=e+pV
\label{1.10}
\end{equation}
is the {\em enthalpy per particle} or {\em specific enthalpy}. By virtue of equations \ref{1.3},  
\ref{1.4} or \ref{1.7}, \ref{1.8}, $h$ can be considered to be a function of $p$ and $s$, and its differential is given by:
\begin{equation}
dh=Vdp+\theta ds
\label{1.11}
\end{equation}
We may in fact use $p$ and $s$ instead of $V$ and $s$ as the basic thermodynamic variables. 

The {\em particle current} is the vectorfield $I$ whose components are given by:
\begin{equation}
I^\mu=n u^\mu
\label{1.12}
\end{equation}
The {\em energy-momentum-stress} tensor is the symmetric 2-contravariant tensorfield $T$ whose 
components are:
\begin{equation}
T^{\mu\nu}=(\rho+p)u^\mu u^\nu+p(g^{-1})^{\mu\nu}
\label{1.13}
\end{equation}
Here $(g^{-1})^{\mu\nu}$ are the components of the reciprocal Minkowski metric, given in rectangular 
coordinates by:
\begin{equation}
(g^{-1})^{00}=-1, \ (g^{-1})^{11}=(g^{-1})^{22}=(g^{-1})^{33}=1, \ \ (g^{-1})^{\mu\nu}=0 \ : \ 
\mbox{if $\mu\neq\nu$}
\label{1.14}
\end{equation}
The {\em equations of motion} of a perfect fluid are the {\em differential particle conservation law} 
or {\em equation of continuity}:
\begin{equation}
\nabla_\mu I^\mu=0
\label{1.15}
\end{equation}
together with the {\em differential energy-momentum conservation law}:
\begin{equation}
\nabla_\nu T^{\mu\nu}=0
\label{1.16}
\end{equation}
In \ref{1.15}, \ref{1.16}, and throughout this monograph, $\nabla$ denotes the covariant derivative 
operator of the Minkowski metric $g$, and equations \ref{1.15}, \ref{1.16}, like the definitions 
\ref{1.12}, \ref{1.13}, are valid in any coordinate system. Nevertheless, Minkowski spacetime having 
the structure of an affine space, there are preferred coordinate systems, the linear coordinates, in which the metric components are constants and the operator $\nabla$ reduces to 
\begin{equation}
\partial_\mu=\frac{\partial}{\partial x^\mu}
\label{1.17}
\end{equation}
applied to the components of tensorfields in these coordinates. A particular case of linear 
coordinates are rectangular coordinates, in which the Minkowski metric components take the form \ref{1.1}. 

Taking the component of equation \ref{1.16} along $u$ by contracting with 
$$u_\mu=g_{\mu\nu}u^\nu$$
yields the equation:
\begin{equation}
u^\mu\partial_\mu\rho+(\rho+p)\nabla_\mu u^\mu=0
\label{1.18}
\end{equation}
Substituting for $\nabla_\mu u^\mu$ from equation \ref{1.15} brings this to the form:
\begin{equation}
u^\mu\partial_\mu\rho=\frac{(\rho+p)}{n}u^\mu\partial_\mu n
\label{1.19}
\end{equation}
On the other hand by virtue of the equation of state \ref{1.3} and the definitions \ref{1.4},
\begin{equation}
\partial_\mu\rho=\frac{(\rho+p)}{n}\partial_\mu n+n\theta\partial_\mu s
\label{1.20}
\end{equation}
Comparing \ref{1.19} and \ref{1.20} we conclude that:
\begin{equation}
u^\mu\partial_\mu s=0
\label{1.21}
\end{equation}
This is the {\em adiabatic condition}, namely the condition that the entropy per particle is constant along the 
{\em flow lines}, that is the integral curves of $u$. It holds as long as we are dealing with a solution of the equations of motion in the classical sense, that is the components $u^\mu \ : \ 
\mu=0,1,2,3$ of $u$, $p$, and $s$ are $C^1$ functions of the rectangular coordinates. A portion of a 
fluid is called {\em isentropic} if $s$ is constant throughout this portion. The adiabatic condition  implies that if a portion of the fluid is isentropic at some time then the same portion, as defined 
by the flow of $u$, is isentropic for all later time. This holds as long as the solution remains 
$C^1$. 

Let us denote by $\Pi$ the $T^1_1$-type tensorfield which is at each point $x$ in spacetime the 
operator of projection to $\Sigma_x$, the {\em local simultaneous space of the fluid} at $x$, 
namely the orthogonal complement in the tangent space at $x$ of the linear span of $u$. The 
components of $\Pi$ are given by:
\begin{equation}
\Pi^\mu_\nu=\delta^\mu_\nu+u^\mu_\nu
\label{1.22}
\end{equation}
where $\delta^\mu_\nu$ is the Kronecker symbol. The projection of \ref{1.16} to $\Sigma_x$ at 
each point reads:
\begin{equation}
(\rho+p)u^\nu\nabla_\nu u^\mu+\Pi^{\mu\nu}\partial_\nu p=0
\label{1.23}
\end{equation}
where
\begin{equation}
\Pi^{\mu\nu}=\Pi^\mu_\lambda(g^{-1})^{\lambda\nu}=(g^{-1})^{\mu\nu}+u^\mu u^\nu
\label{1.24}
\end{equation} 

The sound speed $\eta$ is in conventional units defined by: 
\begin{equation}
\left(\frac{dp}{d\rho}\right)_s=\frac{\eta^2}{c^2}
\label{1.25}
\end{equation}
{\em a fundamental thermodynamic assumption being that the left hand side}, 
a dimensionless quantity, {\em is positive}. Then $\eta$ is defined to be positive. Another 
condition on $\eta$ in the framework of special relativity is that $\eta<c$, namely that the sound 
speed is less than the speed of light in vacuum. In units where $c=1$ this condition reads $\eta<1$. 

By virtue of \ref{1.21} and the definition \ref{1.25} we have:
\begin{equation}
u^\mu\partial_\mu\rho=u^\mu\left(\left(\frac{d\rho}{dp}\right)_s\partial_\mu p+
\left(\frac{d\rho}{ds}\right)_p\partial_\mu s\right)
=\frac{u^\mu}{\eta^2}\partial_\mu p
\label{1.26}
\end{equation}
The equations of motion of a perfect fluid are then seen to be equivalent to the following first 
order quasilinear system of partial differential equations in rectangular coordinates for the unknowns $u^\mu \ : \ \mu=0,1,2,3$, $p$, and $s$, :
\begin{eqnarray}
(\rho+p)u^\nu\partial_\nu u^\mu+\Pi^{\mu\nu}\partial_\nu p=0&\nonumber\\
u^\mu\partial_\mu p+\eta^2(\rho+p)\partial_\mu u^\mu=0&\nonumber\\
u^\mu\partial_\mu s=0&\label{1.27}
\end{eqnarray}

Consider the transformation: 
\begin{equation}
p\mapsto p+b, \ \ \mbox{: $b$  a constant}
\label{1.b1}
\end{equation}
keeping $s$ and $h$ unchanged. Then (see \ref{1.5}, \ref{1.11}) $n$ is unchanged, but (see \ref{1.9}) $\rho$ transforms according to:
\begin{equation}
\rho\mapsto \rho-b
\label{1.b2}
\end{equation}
so that the sum $\rho+p$ is unchanged. The particle current \ref{1.12} is unchanged, but the energy-momentum-stress tensor \ref{1.13} transforms according to:
\begin{equation}
T^{\mu\nu}\mapsto T^{\mu\nu}+b(g^{-1})^{\mu\nu}
\label{1.b3}
\end{equation}
The equations of motion \ref{1.15}, \ref{1.16}, or \ref{1.27} are unaffected.

Let $(u,p,s)$ be a given solution of the equations of motion \ref{1.27} and let 
$$\{(u_\lambda,p_\lambda,s_\lambda) \ : \ \lambda\in I\},$$
$I$ an open interval of the real line containing 0, be a differentiable 1-parameter family of solutions such that:
$$(u_0,p_0,s_0)=(u,p,s).$$
Then
\begin{equation}
(\dot{u},\dot{p},\dot{s})=((du_\lambda/d\lambda)_{\lambda=0},(dp_\lambda/d\lambda)_{\lambda=0},(ds_\lambda/d\lambda)_{\lambda=0})
\label{1.28}
\end{equation}
is a {\em variation} of $(u,p,s)$ through solutions. Note that the constraint \ref{1.2} on $u$ implies the following constraint on $\dot{u}$:
\begin{equation}
u_\mu\dot{u}^\mu=0
\label{1.29}
\end{equation}
Differentiating the equations of motion \ref{1.27} with respect to $\lambda$ at $\lambda=0$ we obtain 
the {\em equations of variation}:
\begin{eqnarray}
(\rho+p)u^\nu\partial_\nu\dot{u}^\mu+\Pi^{\mu\nu}\partial_\nu\dot{p} &=&-(\rho+p)\dot{u}^\nu\partial_\nu u^\mu
+\frac{(\dot{\rho}+\dot{p})}{(\rho+p)}\Pi^{\mu\nu}\partial_\nu p-\dot{\Pi}^{\mu\nu}\partial_\nu p
\nonumber\\
u^\mu\partial_\mu\dot{p}+q\partial_\mu\dot{u}^\mu &=&-\dot{u}^\mu\partial_\mu p-\dot{q}\partial_\mu u^\mu\nonumber\\
u^\mu\partial_\mu\dot{s} &=&-\dot{u}^\mu\partial_\mu s \label{1.30}
\end{eqnarray}
Here $q$ is the function:
\begin{equation}
q=\eta^2(\rho+p)
\label{1.31}
\end{equation}
and for any function of the thermodynamic variables alone:
\begin{equation}
\dot{f}=\left(\frac{df}{dp}\right)_s\dot{p}+\left(\frac{df}{ds}\right)_p\dot{s}
\label{1.32}
\end{equation}
Also:
\begin{equation}
\dot{\Pi}^{\mu\nu}=\dot{u}^\mu u^\nu+u^\mu\dot{u}^\nu
\label{1.33}
\end{equation}
In equations \ref{1.30} we have placed on the left the principal terms, which are linear in the first derivatives of the variation, and on the right the lower order terms, which are linear in the variation itself, with coefficients depending linearly on the first derivatives of the background solution. 

To the system \ref{1.30} and the variation $(\dot{u},\dot{p},\dot{s})$ is associated the {\em energy 
variation current} $\dot{J}$ and the {\em entropy variation current} $\dot{K}$. These are vectorfields 
with components:
\begin{equation}
\dot{J}^\mu=(\rho+p)u^\mu g_{\nu\lambda}\dot{u}^\nu\dot{u}^\lambda+2\dot{u}^\mu\dot{p}
+\frac{1}{\eta^2(\rho+p)}u^\mu\dot{p}^2
\label{1.34}
\end{equation}
and
\begin{equation}
\dot{K}^\mu=u^\mu\dot{s}^2
\label{1.35}
\end{equation}
respectively. Thus $\dot{J}$ is a quadratic form in $(\dot{u},\dot{p})$, while $\dot{K}$ is 
a quadratic form in $\dot{s}$. Let us calculate the divergences of these vectorfields. Substituting from the equations of variation \ref{1.30}, the principal terms cancel and we find:
\begin{eqnarray}
\partial_\mu\dot{J}^\mu&=&-2(\rho+p)\dot{u}^\mu\dot{u}^\nu\partial_\mu u_\nu 
+\frac{2}{(\rho+p)}\left[\left(\frac{d\rho}{ds}\right)_p\dot{s}+\dot{p}\right]\dot{u}^\mu\partial_\mu p\nonumber\\
&+&\left\{\eta^2(\rho+p)g_{\nu\lambda}\dot{u}^\nu\dot{u}^\lambda
+\frac{1}{q}\left[-2\dot{p}\dot{q}+\dot{p}^2\left(1+\left(\frac{dq}{dp}\right)_s\right)\right]\right\}\partial_\mu u^\mu\nonumber\\
&\,&\label{1.36}
\end{eqnarray}
and
\begin{equation}
\partial_\mu\dot{K}^\mu=-2\dot{s}\dot{u}^\mu\partial_\mu s+\dot{s}^2\partial_\mu u^\mu 
\label{1.37}
\end{equation}
The right hand sides of \ref{1.36} and \ref{1.37} are quadratic in the variation with coefficients 
which depend linearly on the first derivatives of the background solution. 

Note that
\begin{equation}
-u_\mu\dot{J}^\mu=(\rho+p)g_{\nu\lambda}\dot{u}^\nu\dot{u}^\lambda+\frac{1}{\eta^2(\rho+p)}\dot{p}^2
\label{1.38}
\end{equation}
and
\begin{equation}
-u_\mu\dot{K}^\mu=\dot{s}^2
\label{1.39}
\end{equation}
In view of the fact that $\dot{u}$ is subject to the constraint \ref{1.29}, so that 
$$g_{\nu\lambda}\dot{u}^\nu\dot{u}^\lambda=\Pi_{\nu\lambda}\dot{u}^\nu\dot{u}^\lambda$$
(where $\Pi_{\nu\lambda}=g_{\nu\mu}\Pi^\mu_\lambda$), \ref{1.38} is a positive definite quadratic 
form in $(\dot{u},\dot{p})$. Also, \ref{1.39} is trivially a positive definite quadratic form in $\dot{s}$. Consider for any covector $\xi$ in the cotangent space at a given point the quadratic 
forms $\xi_\mu\dot{J}^\mu$, $\xi_\mu\dot{K}^\mu$. For any pair of variations 
$(\dot{u},\dot{p},\dot{s})$, $(\dot{u}^\prime,\dot{p}^\prime,\dot{s}^\prime)$ the corresponding symmetric bilinear forms are:
\begin{eqnarray}
&&\xi_\mu\dot{J}^\mu((\dot{u},\dot{p}),(\dot{u}^\prime,\dot{p}^\prime))\label{1.40}\\
&& \ \ \ =\xi_\mu\left\{(\rho+p)u^\mu g_{\nu\lambda}\dot{u}^\nu\dot{u}^{\prime\lambda}
+\dot{u}^\mu\dot{p}^\prime+\dot{p}\dot{u}^{\prime\mu}+
\frac{1}{\eta^2(\rho+p)}u^\mu\dot{p}\dot{p}^\prime\right\}\nonumber
\end{eqnarray}
and
\begin{equation}
\xi_\mu\dot{K}^\mu(\dot{s},\dot{s}^\prime)=\xi_\mu u^\mu \dot{s}\dot{s}^\prime
\label{1.41}
\end{equation}
Consider now the set of all covectors $\xi$ at $x$ such that the symmetric bilinear forms \ref{1.40} 
and \ref{1.41} are degenerate, that is, there is a non-zero variation 
$(\dot{u},\dot{p},\dot{s})$ such that for all variations $(\dot{u}^\prime,\dot{p}^\prime,\dot{s}^\prime)$ we have:
\begin{equation}
\xi_\mu\dot{J}^\mu((\dot{u},\dot{p}),(\dot{u}^\prime,\dot{p}^\prime))=0
\label{1.42}
\end{equation}
and:
\begin{equation}
\xi_\mu\dot{K}^\mu(\dot{s},\dot{s}^\prime)=0
\label{1.43}
\end{equation}
This defines the {\em characteristic subset}  of the cotangent space at $x$. 
Taking into account the constraint \ref{1.29} we see that \ref{1.42} is equivalent to the linear 
system:
\begin{eqnarray}
\Pi^\lambda_\nu\xi_\lambda\dot{p}+(\rho+p)(\xi_\mu u^\mu)\dot{u}_\nu=&0\nonumber\\
\frac{(\xi_\mu u^\mu)\dot{p}}{\eta^2(\rho+p)}+\xi_\mu\dot{u}^\mu=&0\label{1.44}
\end{eqnarray}
while \ref{1.43} is equivalent to:
\begin{equation}
(\xi_\mu u^\mu)\dot{s}=0
\label{1.45}
\end{equation}
It follows that either:
\begin{equation}
\xi_\mu u^\mu=0
\label{1.46}
\end{equation}
in which case:
\begin{equation}
\dot{p}=0 \ \ \mbox{and} \ \ \xi_\mu\dot{u}^\mu=0
\label{1.47}
\end{equation}
or:
\begin{equation}
(h^{-1})^{\mu\nu}\xi_\mu\xi_\nu=0
\label{1.48}
\end{equation}
in which case:
\begin{equation}
\dot{s}=0 \ \ \mbox{and} \ \ \dot{u}^\nu=-\frac{\Pi^{\nu\lambda}\xi_\lambda\dot{p}}
{(\rho+p)(\xi_\mu u^\mu)}
\label{1.49}
\end{equation}
In \ref{1.48} $h^{-1}$ is the {\em reciprocal acoustical metric}, given by:
\begin{equation}
(h^{-1})^{\mu\nu}=\Pi^{\mu\nu}-\frac{1}{\eta^2}u^\mu u^\nu=(g^{-1})^{\mu\nu}
-\left(\frac{1}{\eta^2}-1\right)u^\mu u^\nu
\label{1.50}
\end{equation}
It is the reciprocal of the {\em acoustical metric} $h$, given by:
\begin{equation}
h_{\mu\nu}=g_{\mu\nu}+(1-\eta^2)u_\mu u_\nu
\label{1.51}
\end{equation}
a Lorentzian metric on the spacetime manifold. The subset of the cotangent space at $x$ defined by the 
condition \ref{1.46} is a hyperplane $P^*_x$ while the subset of the cotangent space at $x$ defined by 
the condition \ref{1.48} is a cone $C^*_x$. The dual to $P^*_x$ characteristic subset of the tangent 
space at $x$ is the linear span of $u$, while the dual to $C^*_x$ characteristic subset $C_x$ is the 
set of all non-zero vectors $X$ at $x$ satisfying:
\begin{equation}
h_{\mu\nu}X^\mu X^\nu=0
\label{1.52}
\end{equation}
This is the {\em sound cone} at $x$. The requirement that $\eta<1$ is the requirement that at each 
point $x$ the sound cone \ref{1.52} lies within the light cone. It is equivalent to the requirement 
that in the cotangent space at each point $x$ the cone $(g^{-1})^{\mu\nu}\xi_\mu\xi_\nu=0$ lies within 
the cone $C^*_x$. The last is in turn equivalent to the condition that a covector $\xi$ at $x$ whose null space is a space-like, relative to $g$, hyperplane $H_x$,  and which has positive evaluation in the future of $H_x$, belongs to the interior of the positive component of $C^*_x$. 
It follows from the above that the set of all covectors $\xi$ at $x$ such that the 
quadratic forms $\xi(\dot{J})$, $\xi(\dot{K})$ are both positive definite is the interior of the positive component of $C^*_x$. 

We introduce the 1-form $\beta$, given by:
\begin{equation}
\beta_\mu=-h u_\mu 
\label{1.53}
\end{equation}
The 1-form $\beta$ has at each point $x$ in spacetime positive evaluation on the future light cone 
in the tangent space at $x$ and its interior. Moreover, by the constraint \ref{1.2}, we can express:
\begin{equation}
h=\sqrt{-(g^{-1})^{\mu\nu}\beta_\mu\beta_\nu}
\label{1.54}
\end{equation}
By \ref{1.9} and the first of equations \ref{1.27}, ${\cal L}_u\beta$, the Lie derivative with respect to $u$ of 
$\beta$, is given by:
\begin{equation}
({\cal L}_u\beta)_\mu=\frac{1}{n}\partial_\mu p+u_\mu u^\nu\left(\frac{1}{n}\partial_\nu p-\partial_\nu h\right)
\label{1.55}
\end{equation}
Now by \ref{1.5} and \ref{1.11} the expression in the last parenthesis is equal to $-\theta\partial_\nu s$, therefore by virtue of the third of equations \ref{1.27} the last term vanishes. Taking again account of \ref{1.11} then yields:
\begin{equation}
{\cal L}_u\beta=dh-\theta ds
\label{1.56}
\end{equation}

The {\em vorticity 2-form} $\omega$ is minus the exterior derivative of the 1-form $\beta$:
\begin{equation}
\omega=-d\beta
\label{1.57}
\end{equation}

Let us recall the following general relation between the exterior derivative and the Lie derivative with respect to a vectorfield $X$, as applied to an exterior differential form $\vartheta$ of any rank:
\begin{equation}
{\cal L}_X\vartheta=i_X d\vartheta+di_X\vartheta
\label{1.58}
\end{equation}
where $i_X$ denotes contraction on the left by $X$. Taking $X=u$ and $\vartheta=\beta$ we have from 
\ref{1.53}:
\begin{equation}
i_u\beta=h
\label{1.59}
\end{equation}
Comparing then \ref{1.49} in this case with \ref{1.56} we conclude that:
\begin{equation}
i_u\omega=\theta ds
\label{1.60}
\end{equation}
This equation is equivalent to the first and third of equations \ref{1.27}, more precisely the 
$u$ component of \ref{1.60} is equivalent to the third of \ref{1.27} while the projection of \ref{1.60} orthogonal to $u$ is equivalent, modulo the third of \ref{1.27}, to the first of \ref{1.27}. It follows that, for $C^1$ solutions the equations of motion \ref{1.15}, \ref{1.16} 
are equivalent to equation \ref{1.60} together with \ref{1.15}. Moreover, we can take $(h,s)$ as the basic thermodynamic variables and consider the equation of state as giving $p$ as a function 
of $h$ and $s$, in which case \ref{1.11} takes the form (see \ref{1.5}):
\begin{equation}
p=p(h,s), \ \ \ dp=n(dh-\theta ds)
\label{1.d1}
\end{equation}
Since $h$ is expressed  in terms of $\beta$ by \ref{1.54}, the unknowns in this formulation are $\beta$ and $s$. Also, introducing the thermodynamic function: 
\begin{equation}
G=\frac{n}{h}
\label{1.72}
\end{equation}
we can by \ref{1.53} express \ref{1.15} in the form:
\begin{equation}
\nabla_\mu(G\beta^\mu)=0 \ \ \  \mbox{where} \ \beta^\mu=(g^{-1})^{\mu\nu}\beta_\nu
\label{1.d2}
\end{equation}
so for $C^1$ solutions the equations of motion \ref{1.15}, \ref{1.16} are equivalent to equations 
\ref{1.60}, \ref{1.d2} for $\beta$ and $s$. Note that:
\begin{equation}
\frac{h}{G}\left(\frac{dG}{dh}\right)_s=\frac{h}{n}\left(\frac{dn}{dh}\right)_s-1
\label{1.d3}
\end{equation}
On the other hand, by \ref{1.9} $\rho=nh-p$, hence:
$$\left(\frac{d\rho}{dh}\right)_s=h\left(\frac{dn}{dh}\right)_s$$
and, since also $(dp/dh)_s=n$, we obtain:
\begin{equation}
\frac{h}{n}\left(\frac{dn}{dh}\right)_s=\frac{(d\rho/dh)_s}{(dp/dh)_s}=\frac{1}{\eta^2}
\label{1.d4}
\end{equation}
Comparing with \ref{1.d3} we conclude that:
\begin{equation}
\frac{h}{G}\left(\frac{dG}{dh}\right)_s=\frac{1}{\eta^2}-1
\label{1.d5}
\end{equation}
This implies that:
\begin{eqnarray}
\partial_\mu G&=&\left(\frac{dG}{dh}\right)_s\partial_\mu h+
\left(\frac{dG}{ds}\right)_h\partial_\mu s\nonumber\\
&=&-F\beta^\nu\nabla_\mu\beta_\nu+\left(\frac{dG}{ds}\right)_h\partial_\mu s
\label{1.d6}
\end{eqnarray}
where $F$ is the function:
\begin{equation}
F=\frac{1}{hG}\left(\frac{dG}{dh}\right)_s=\frac{1}{h^2}\left(\frac{1}{\eta^2}-1\right)
\label{1.d7}
\end{equation}
by \ref{1.d5}. In view of the $u$ component of \ref{1.60}, which is equivalent to the 3rd of 
\ref{1.27}, and \ref{1.53}, it then follows that:
\begin{equation}
\beta^\mu\partial_\mu G=-F\beta^\mu\beta^\nu\nabla_\mu\beta_\nu
\label{1.d8}
\end{equation}
We thus see that equation \ref{1.d2} can be written in the form:
\begin{equation}
(h^{-1})^{\mu\nu}\nabla_\mu\beta_\nu=0
\label{1.d9}
\end{equation}
the reciprocal acoustical metric \ref{1.50}, being, in view of \ref{1.53} and \ref{1.d7}, given by:
\begin{equation}
(h^{-1})^{\mu\nu}=(g^{-1})^{\mu\nu}-F\beta^\mu\beta^\nu
\label{1.d10}
\end{equation}
Note that the acoustical metric \ref{1.51} itself is given by:
\begin{equation}
h_{\mu\nu}=g_{\mu\nu}+H\beta_\mu\beta_\nu
\label{1.d11}
\end{equation}
where $H$ is the function:
\begin{equation}
H=\frac{F}{1+h^2 F}
\label{1.d12}
\end{equation}
and we have:
\begin{equation}
1-h^2 H=\eta^2=\frac{1}{1+h^2 F}
\label{1.d13}
\end{equation}

If we think of the vorticity 2-form $\omega$ as analogous to the electromagnetic field, equation 
\ref{1.60} states that the ``electric" part of $\omega$ is equal to the heat differential $\theta ds$. The ``magnetic" part of $\omega$ is given by the vectorfield:
\begin{equation}
\varpi^\mu=\frac{1}{2}(\epsilon^{-1})^{\mu\alpha\beta\gamma}u_\alpha\omega_{\beta\gamma}
\label{1.61}
\end{equation}
The electric-magnetic decomposition refers to the local rest frame of the fluid. In \ref{1.61} 
$\epsilon^{-1}$ is the reciprocal volume form of $g$, that is the volume form induced by $g$ in the 
cotangent space at each point. We denote by $\epsilon$ the volume form of $g$. Thus, if $(E_0,E_1,E_2,E_3)$ is a positive basis for the tangent space at $x$ which is orthonormal with 
respect to $g$ and $(\vartheta_0,\vartheta_1,\vartheta_2,\vartheta_3)$ is the dual basis for the 
cotangent space at $x$, then:
\begin{equation}
\epsilon^{-1}(\vartheta_0,\vartheta_1,\vartheta_2,\vartheta_3)=\epsilon(E_0,E_1,E_2,E_3)=1
\label{1.62}
\end{equation}
The components of $\epsilon$ and $\epsilon^{-1}$ are given in an arbitrary system of coordinates by:
\begin{equation}
\epsilon_{\alpha\beta\gamma\delta}=\sqrt{-\mbox{det}g}[\alpha\beta\gamma\delta], \ \ \ 
(\epsilon^{-1})^{\alpha\beta\gamma\delta}=\frac{[\alpha\beta\gamma\delta]}{\sqrt{-\mbox{det}g}}
\label{1.63}
\end{equation}
where $[\alpha\beta\gamma\delta]$ is the 4-dimensional fully antisymmetric symbol. In rectangular 
coordinates \ref{1.63} reduces to:
\begin{equation}
\epsilon_{\alpha\beta\gamma\delta}=(\epsilon^{-1})^{\alpha\beta\gamma\delta}=[\alpha\beta\gamma\delta]
\label{1.64}
\end{equation}
The vectorfield $\varpi$ is the obstruction to integrability of the distribution of local simultaneous spaces $\{\Sigma_x\}$. We call $\varpi$ {\em vorticity vector}. At each point $x$, the vector $\varpi(x)$ belongs to $\Sigma_x$. 

In the case of an isentropic fluid portion equation \ref{1.60} reduces to:
\begin{equation}
i_u\omega=0
\label{1.65}
\end{equation}
Since by the definition \ref{1.57} we have, in general:
\begin{equation}
d\omega=0
\label{1.66}
\end{equation}
it follows, taking in \ref{1.58} $X=u$ and $\vartheta=\omega$, that in an isentropic fluid portion 
we have:
\begin{equation}
{\cal L}_u\omega=0
\label{1.67}
\end{equation}
that is, in an isentropic fluid portion the vorticity 2-form is Lie transported along the 
flow. A portion of the fluid is called {\em irrotational} if $\omega$ vanishes throughout this 
portion. Remark that by \ref{1.60} a fluid portion which is irrotational in this sense is also isentropic. 

\vspace{5mm}

\section{The Irrotational Case and the Nonlinear Wave Equation}

Consider now a fluid portion which is irrotational at a given time. Then according to the 
preceding the same portion, as defined by the flow of $u$, remains irrotational for all later time, 
as long as the solution remains $C^1$ in terms of the basic variables $u$, $p$ and $s$. If the portion is simply connected we can then introduce a function $\phi$, determined up to an additive constant, 
such that
\begin{equation}
\beta=d\phi
\label{1.68}
\end{equation}
It follows from the definition of $\beta$ that the derivative of $\phi$ along any future-directed 
non-space-like vector is positive. By \ref{1.53} and \ref{1.68} we can express:
\begin{equation}
h=\sqrt{-(g^{-1})^{\mu\nu}\partial_\mu\phi\partial_\nu\phi}, \ \ \ u^\mu=-\frac{\partial^\mu\phi}{h} 
\ \ \ (\partial^\mu=(g^{-1})^{\mu\nu}\partial_\nu)
\label{1.69}
\end{equation}
With \ref{1.68} $\omega$ vanishes so with constant $s$ equation \ref{1.60} is identically satisfied. 
It follows that the whole content of the equations of motion is in this case contained in  
equation \ref{1.d2}, which takes the form of a nonlinear wave equation for $\phi$:
\begin{equation}
\nabla_\mu(G\partial^\mu\phi)=0
\label{1.70}
\end{equation}
$G$ being the thermodynamic function \ref{1.72}, which here, since $s$ is constant, is a function of $h$ alone:
\begin{equation}
G=G(h)
\label{1.71}
\end{equation}
Likewise $p$ is a function of $h$ alone:
\begin{equation}
p=p(h), \ \ \ \frac{dp}{dh}=n
\label{1.73}
\end{equation}
the first being an expression of the equation of state at the given value of $s$ (a given 
adiabat). Equation \ref{1.70} is the Euler-Lagrange equation corresponding to the Lagrangian density
\begin{equation}
L=p(h)
\label{1.74}
\end{equation}
with $h$ {\em defined} by the first of \ref{1.69}. Remark that $p$ is a convex function of $h$:
\begin{equation}
\frac{d^2 p}{dh^2}>0
\label{1.75}
\end{equation}
For, 
\begin{equation}
\frac{d^2 p}{dh^2}=\frac{dn}{dh}=\frac{1}{\eta^2}\frac{n}{h}
\label{1.76}
\end{equation}
by \ref{1.d4}. 

The nonlinear wave equation \ref{1.70} can be written in rectangular coordinates in the form: 
\begin{equation}
(h^{-1})^{\mu\nu}\frac{\partial^2\phi}{\partial x^\mu\partial x^\nu}=0
\label{1.80}
\end{equation} 
which corresponds to \ref{1.d9}. The reciprocal acoustical metric \ref{1.d10}, and the acoustical metric \ref{1.d11} take in the irrotational case the form:
\begin{equation}
(h^{-1})^{\mu\nu}=(g^{-1})^{\mu\nu}-F\partial^\mu\phi\partial^\nu\phi
\label{1.81}
\end{equation}
and
\begin{equation}
h_{\mu\nu}=g_{\mu\nu}+H\partial_\mu\phi\partial_\nu\phi
\label{1.88}
\end{equation}
where the functions $F$ and $H$  are here functions of $h$ alone, $s$ being constant.           

The acoustical metric \ref{1.88} is:
\begin{equation}
h=g+Hd\phi\otimes d\phi
\label{1.a1}
\end{equation}
The case in which, along the adiabat corresponding to the given value of $s$, 
\begin{equation}
H=H_0 \ \ \mbox{: a positive constant}
\label{1.a2}
\end{equation}
plays a special role, as we shall discuss in Section 1.4. In this case the acoustical metric 
is the metric induced on the graph:
\begin{equation}
y=\sqrt{H_0}\phi(x) \ : \ x\in \mathbb{M}^3
\label{1.a3}
\end{equation}
in Minkowski spacetime $\mathbb{M}^4$ of 4 spatial dimensions over the standard Minkowski spacetime 
$\mathbb{M}^3$ of 3 spatial dimensions, the extra spatial dimension being $y$. The metric of $\mathbb{M}^4$ is:
\begin{equation}
g+dy\otimes dy
\label{1.a4}
\end{equation}
$g$ denoting as usual the metric of $\mathbb{M}^3$. By \ref{1.d13} the sound speed is in this case 
given by:
\begin{equation}
\eta=\sqrt{1-H_0 h^2}, \ \ \mbox{so we have} \ \ 0<h<\frac{1}{\sqrt{H_0}}.
\label{1.a5}
\end{equation}
Integrating the general equation \ref{1.d4} with $\eta^2$ given by \ref{1.a5} we obtain:
\begin{equation}
n=k\frac{h}{\sqrt{1-H_0 h^2}}, \ \ \mbox{$k$ : a positive constant.}
\label{1.a6}
\end{equation}
Integrating then the general equation \ref{1.73} with $n$ given by \ref{1.a6} yields:
\begin{equation}
p=p_0-\frac{k}{H_0}\sqrt{1-H_0 h^2}
\label{1.a7}
\end{equation}
where $p_0$ is an arbitrary constant. 
Here $p_0$ is the limit of $p$ as $h\rightarrow 1/\sqrt{H_0}$. 
Let $a$ be the positive constant:
\begin{equation}
a=\frac{k}{H_0}
\label{1.a8}
\end{equation}
By adding a suitable constant to $p$ (see discussion following \ref{1.b1}) we may set:
\begin{equation}
p_0=a
\label{1.a9}
\end{equation} 
Then in the limit $h\rightarrow 0$ we have $\eta\rightarrow 1$, $n\rightarrow 0$, $p\rightarrow 0$, 
while in the limit $h\rightarrow 1/\sqrt{H_0}$ we have $\eta\rightarrow 0$, $n\rightarrow\infty$, 
$p\rightarrow a$. Also, since in general:
\begin{equation}
\rho+p=nh
\label{1.a10}
\end{equation}
$p$ is expressed in this special case as a function of $\rho$ by:
\begin{equation}
p=\frac{\rho}{1+(\rho/a)}
\label{1.a11}
\end{equation}
$p$ increasing from 0 to $a$ as $\rho$ increases from 0 to $\infty$. Moreover, $p$ is 
given in terms of $\eta$ by:
\begin{equation}
p=a(1-\eta)
\label{1.a12}
\end{equation}
Now, by the general formula \ref{1.51} we have, in rectangular coordinates, 
\begin{equation}
\mbox{det}h=-\eta^2
\label{1.93}
\end{equation}
This implies that the volume elements of the metrics $h$ and $g$ are related by:
\begin{equation}
d\mu_h=\eta d\mu_g
\label{1.a13}
\end{equation}
The action corresponding to the Lagrangian density \ref{1.74} and to a spacetime domain ${\cal D}$ 
being in general 
\begin{equation}
\int_{{\cal D}}pd\mu_g
\label{1.a14}
\end{equation}
this reduces by \ref{1.a13}, \ref{1.93} in the special case to:
\begin{equation}
a(\mbox{Vol}_g({\cal D})-\mbox{Vol}_h({\cal D}))
\label{1.a15}
\end{equation}
that is, in view of the above, up to the multiplicative constant $a$, the volume of ${\cal D}$ minus 
the volume of the graph over ${\cal D}$ in $\mathbb{M}^4$. It follows that the nonlinear wave equation 
\ref{1.70} reduces in the special case to the {\em minimal hypersurface equation}. 

Let $\phi$ be a given solution of the nonlinear wave equation \ref{1.70} and let $\{\phi_\lambda \ : \ 
\lambda\in I\}$, $I$ an open interval of the real line containing 0, be a differentiable 1-parameter 
family of solutions such that 
$$\phi_0=\phi.$$
Then 
\begin{equation}
\dot{\phi}=\left(\frac{d\phi_\lambda}{d\lambda}\right)_{\lambda=0}
\label{1.91}
\end{equation}
is a {\em variation} of $\phi$ through solutions. Differentiating the equation, in rectangular coordinates, 
$$\partial_\mu(G(h_\lambda)\partial^\mu\phi_\lambda)=0, \ \ \ 
h_\lambda=\sqrt{-\partial^\mu\phi_\lambda\partial_\mu\phi_\lambda}$$
with respect to $\lambda$ at $\lambda=0$ we obtain the {\em equation of variation}:
\begin{equation}
\partial_\mu(G(h)(h^{-1})^{\mu\nu}\partial_\nu\dot{\phi})=0
\label{1.92}
\end{equation}
Here $h^{-1}$ is the reciprocal acoustical metric \ref{1.81} corresponding to the solution $\phi$. 

Consider now the {\em conformal acoustical metric}:
\begin{equation}
\tilde{h}=\Omega h
\label{1.94}
\end{equation}
The conformal factor $\Omega$ is to be appropriately chosen below. The metrics $\tilde{h}$ and $h$ being conformal, they define the same null cone at each point, the sound cone. In the physical case of 3 spatial dimensions we have:
\begin{equation}
\sqrt{-\mbox{det}\tilde{h}}=\Omega^2\sqrt{-\mbox{det}h}=\Omega^2\eta
\label{1.95}
\end{equation}
(see \ref{1.93}). Since also
$$(\tilde{h}^{-1})^{\mu\nu}=\Omega^{-1}(h^{-1})^{\mu\nu}$$
it follows that:
$$\sqrt{-\mbox{det}\tilde{h}}(\tilde{h}^{-1})^{\mu\nu}=\Omega\eta(h^{-1})^{\mu\nu}$$
Comparing with \ref{1.92} we see that, choosing
\begin{equation}
\Omega=\frac{G}{\eta}
\label{1.96}
\end{equation}
the equation of variation \ref{1.92} becomes 
$$\partial_\mu(\sqrt{-\mbox{det}\tilde{h}}(\tilde{h}^{-1})^{\mu\nu}\partial_\nu\dot{\phi})=0$$
This is simply the linear wave equation corresponding to the metric $\tilde{h}$:
\begin{equation}
\square_{\tilde{h}}\dot{\phi}=0
\label{1.97}
\end{equation}
and in this form it is valid in any system of coordinates. 

To equation \ref{1.97} and to the variation $\dot{\phi}$ is associated the {\em variational stress}, 
$\dot{T}$, a $T^1_1$-type tensorfield given in arbitrary coordinates by:
\begin{equation}
\dot{T}^\mu_\nu=(\tilde{h}^{-1})^{\mu\kappa}\partial_\kappa\dot{\phi}\partial_\nu\dot{\phi}-\frac{1}{2}\delta^\mu_\nu(\tilde{h}^{-1})^{\kappa\lambda}\partial_\kappa\dot{\phi}\partial_\lambda\dot{\phi}
\label{1.98}
\end{equation}
Thus $\dot{T}$ is a quadratic form in $d\dot{\phi}$. 
Denoting by $\tilde{D}$ the covariant derivative operator of the conformal acoustical metric $\tilde{h}$, the following identity holds:
\begin{equation}
\tilde{D}_\mu\dot{T}^\mu_\nu=\partial_\nu\dot{\phi}\square_{\tilde{h}}\dot{\phi}
\label{1.99}
\end{equation}
hence as a consequence of equation \ref{1.97} we have:
\begin{equation}
\tilde{D}_\mu\dot{T}^\mu_\nu=0
\label{1.100}
\end{equation} 
If $(X,\xi)$ is a pair of a vector $X$ and a covector $\xi$ at $x$, 
$$-\dot{T}^\mu_\nu(x)X^\nu\xi_\mu,$$
a quadratic form in $d\dot{\phi}(x)$, is positive definite, if and only if $X$ belongs to the interior 
of the positive component of $C_x$ and $\xi$ belongs to the interior of the positive component of 
$C^*_x$, or else $X$ belongs to the interior of the negative component of $C_x$ and $\xi$ belongs to 
the interior of the negative component of $C^*_x$. 

\vspace{5mm}

\section{The Non-Relativistic Limit}

We shall now show how the non-relativistic theory is obtained by taking the limit $c\rightarrow\infty$. What this means is that time intervals and lengths are measured in conventional units and we consider the limit where the speed of light in vacuum $c$ in 
these units is considered to be unboundedly large. To be specific, if $(x^0; x^i \ : \ i=1,2,3)$ are rectangular coordinates in Minkowski spacetime as above, we set:
\begin{equation}
x^0=ct
\label{1.101}
\end{equation}
and consider the limit $c\rightarrow \infty$. 
In this limit, the affine structure of spacetime 
is preserved, but the closed exterior of the light cone  at each spacetime  point $x$ flattens out to become a spatial hyperplane of absolute simultaneity $\Sigma_t$. So the limit spacetime, the Galilei 
spacetime of classical mechanics, is endowed, like the Minkowski spacetime, with an affine structure, 
but there is in addition a canonical projection $\pi$ to the 1-dimensional affine space of time, taking the spacetime point $x$ to $t$. The $\Sigma_t=\pi^{-1}(t)$ constitute a family of parallel hyperplanes with respect to the affine spacetime structure. The Euclidean induced metric on each $\Sigma_t$ is preserved in the limit, but this only measures the arc lengths of curves which lie on a given $\Sigma_t$. Only such curves remain space-like in the limit. Any curve transversal to the $\Sigma_t$ is time-like in the limit and the temporal length of any arc of such a curve is in the limit simply the difference of the values of $t$ at its end points. Taking any line transversal to the 
$\Sigma_t$ we consider the family of all lines parallel to the given one, relative to the affine spacetime structure. Such a family defines an isometry of each $\Sigma_t$ onto every other $\Sigma_t$. 
These families are called {\em Galilei frames} and correspond to inertial systems of observers. 
A given Galilei frame makes Galilei spacetime into a product of the 1-dimensional affine space of  time, with a unit, which, picking an origin, we represent by $\mathbb{R}$, with 3-dimensional Euclidean space $\mathbb{E}^3$, which represents any one of the $\Sigma_t$. Although Galilei spacetime 
is not endowed with a metric, it is endowed with a volume element, given by:
\begin{equation}
d\mu=dt\wedge d\overline{\mu}
\label{1.c1}
\end{equation}
where $d\overline{\mu}$ is the volume element of $\mathbb{E}^3$:
\begin{equation}
d\overline{\mu}=dx^1\wedge dx^2\wedge dx^3 \ \ \ \mbox {: in rectangular coordinates on $\mathbb{E}^3$.}
\label{1.c2}
\end{equation}
The volume element \ref{1.c1} is independent of the choice of Galilei frame. Because of the fact that 
Galilei spacetime is not endowed with a metric, there is no natural isomorphism of the cotangent space 
onto the tangent space at each point. However, given a covector $\xi$ at a spacetime point $x$, we can 
consider $\overline{\xi}$, the restriction of $\xi$ to vectors in $\Sigma_t$, the hyperplane of absolute simultaneity through $x$. Since $\Sigma_t$ has Euclidean structure, there is a vector 
$\overline{\xi}^\sharp$ tangent to $\Sigma_t$ at $x$ corresponding to $\overline{\xi}$. However the linear map $\xi\mapsto\overline{\xi}^\sharp$ is not onto, and its kernel is the 1-dimensional space of covectors at $x$ whose null space is $\Sigma_t$. 

Coming now to fluid mechanics, the components $u^\mu$ of the relativistic spacetime fluid velocity 
are given in any rectangular coordinate system by:
\begin{equation}
u^0=\frac{1}{\sqrt{1-|v|^2/c^2}}; \ \ \ u^i=\frac{v^i/c}{\sqrt{1-|v|^2/c^2}} \ : \ i=1,2,3
\label{1.102}
\end{equation}
where 
\begin{equation}
v^i=c\frac{u^i}{u^0}
\label{1.103}
\end{equation}
are the components of the spatial fluid velocity in conventional units. Note that the componets $u^\mu$ are dimensionless quantities while the componets $v^i$ have the dimensions of velocity. 
Here and in the following the Latin indices range in $\{1,2,3\}$ (the Greek indices range in 
$\{0,1,2,3\}$). In the limit $c\rightarrow\infty$ the spacetime velocity vectorfield
$$cu^\mu\frac{\partial}{\partial x^\mu}$$
becomes:
\begin{equation}
u=\frac{\partial}{\partial t}+v^i\frac{\partial}{\partial x^i}
\label{1.104}
\end{equation}

The concepts of volume per particle $V$, hence also number of particles per unit volume $n$, entropy per particle $s$, pressure $p$, and temperature $\theta$, are the same in the non-relativistic theory, 
however the non-relativistic energy per particle does not contain the rest mass contribution. We therefore write:
\begin{equation}
e=mc^2+e^\prime
\label{1.105}
\end{equation}
What has a non-relativistic limit is the remainder $e^\prime$. For the same reason (see \ref{1.10}) 
the non-relativistic enthalpy per particle does not contain the rest mass contribution, therefore 
writing:
\begin{equation}
h=mc^2+h^\prime
\label{1.106}
\end{equation}
what has a non-relativistic limit is the remainder $h^\prime$. As in the relativistic theory we can 
choose the particle mass $m$ to be 1, so all quantities ``per particle", like $v$, $s$, $e^\prime$, 
and $h^\prime$, are quantities per unit mass or ``specific" quantities, and $n$ is the {\em mass density}. In particular $e^\prime$ 
is called {\em specific internal energy}. Since the subtraction of a constant does not affect the 
differential, the analogues of equations \ref{1.8}, \ref{1.11}, namely:
\begin{equation}
de^\prime=-pdV+\theta ds, \ \ \ dh^\prime=Vdp+\theta ds
\label{1.107}
\end{equation}
hold. Also, we have:
\begin{equation}
\rho=nc^2+\rho^\prime
\label{1.108}
\end{equation}
where
\begin{equation}
\rho^\prime=ne^\prime
\label{1.109}
\end{equation}
is the {\em internal energy per unit volume}. 

We remark that in the non-relativistic theory the specific internal energy, hence also the specific 
enthalpy, is only defined up to an additive constant. So the transformation 
\begin{equation}
e^\prime\mapsto e^\prime+a, \ \ h^\prime\mapsto h^\prime+a, \ \ \ \mbox{ $a$ : a constant}
\label{1.b4}
\end{equation}
which makes the internal energy per unit volume transform according to:
\begin{equation}
\rho^\prime\mapsto \rho^\prime+na
\label{1.b5}
\end{equation}
should not affect the equations of motion.

The non-relativistic particle current is the limit of the spacetime vectorfield
$$cI^\mu\frac{\partial}{\partial x^\mu}$$
which is simply:
\begin{equation}
I=nu
\label{1.110}
\end{equation}
with $u$ given by \ref{1.104}. The non-relativistic limit of equation \ref{1.15} is then seen to be the 
{\em classical equation of continuity}: 
\begin{equation}
\frac{\partial n}{\partial t}+\frac{\partial(nv^i)}{\partial x^i}=0
\label{1.111}
\end{equation} 

In regard to the components \ref{1.13} of the energy-momentum-stress tensor, substituting \ref{1.102},  
\ref{1.108} we obtain:
\begin{eqnarray}
&T^{00}=nc^2+n|v|^2+\rho^\prime+O(c^{-2})\nonumber\\
&T^{0i}=ncv^i+c^{-1}(n|v|^2+\rho^\prime+p)v^i+O(c^{-3})\nonumber\\
&T^{ij}=nv^iv^j+p\delta^{ij}+O(c^{-2})
\label{1.112}
\end{eqnarray}
The non-relativistic limit of the spatial components of \ref{1.16} is the seen to be the 
{\em classical differential momentum conservation law}:
\begin{equation}
\frac{\partial(nv^i)}{\partial t}+\frac{\partial(nv^iv^j+p\delta^{ij})}{\partial x^j}=0
\label{1.113}
\end{equation}
The 0-component of \ref{1.16}, multiplied by $c^{-1}$, reproduces in the limit $c\rightarrow\infty$ 
equation \ref{1.111}. To obtain something new, we consider instead the following linear combination
of \ref{1.15} with the 0-component of \ref{1.16}:
\begin{equation}
\frac{\partial}{\partial x^\nu}(T^{0\nu}-c^2 I^\nu)=0
\label{1.114}
\end{equation}
Substituting \ref{1.102} in \ref{1.12} we obtain:
\begin{eqnarray}
&c^2 I^0=nc^2+\frac{1}{2}n|v|^2+O(c^{-2})\nonumber\\
&c^2 I^i=nv^i c+\frac{1}{2}c^{-1}n|v|^2 v^i+O(c^{-3})
\label{1.115}
\end{eqnarray}
Combining with the first two of \ref{1.112} gives:
\begin{eqnarray}
&&T^{00}-c^2 I^0=\frac{1}{2}n|v|^2+\rho^\prime+O(c^{-2})\nonumber\\
&&T^{0i}-c^2 I^i=c^{-1}\left(\frac{1}{2}n|v|^2+\rho^\prime+p\right)v^i+O(c^{-3})
\label{1.116}
\end{eqnarray}
Therefore multiplying \ref{1.114} by $c$ and taking the limit $c\rightarrow\infty$ yields the 
{\em classical differential energy conservation law}:
\begin{equation}
\frac{\partial}{\partial t}\left(\frac{1}{2}n|v|^2+\rho^\prime\right)+
\frac{\partial}{\partial x^i}\left[\left(\frac{1}{2}n|v|^2+\rho^\prime+p\right)v^i\right]=0
\label{1.117}
\end{equation}
Equations \ref{1.111}, \ref{1.113}, \ref{1.117} constitute the {\em classical equations of motion} 
for a perfect fluid. We see that by virtue of the equation of continuity \ref{1.111} the energy 
equation \ref{1.117} is unaffected by the transformation \ref{1.b5}. Moreover, as in the relativistic theory, the transformation 
\begin{equation}
p\mapsto p+b, \ \ \mbox{: $b$  a constant}
\label{1.b6}
\end{equation}
keeping $s$ and $h$ unchanged, leaves $n$ unchanged, but, since $\rho^\prime+p=nh$, transforms $\rho^\prime$ according to 
\begin{equation}
\rho^\prime\mapsto \rho^\prime-b
\label{1.b7}
\end{equation}
This transformation similarly leaves the equations of motion unaffected.

The argument leading from the relativistic equations of motion \ref{1.15}, \ref{1.16}, for $C^1$ 
solutions to the adiabatic condition in the form \ref{1.21}, leads similarly for $C^1$ solutions 
of the classical equations of motion to the classical form of the adiabatic condition:
\begin{equation}
\frac{\partial s}{\partial t}+v^i\frac{\partial s}{\partial x^i}=0
\label{1.118}
\end{equation}

In view of \ref{1.108} the non-relativistic limit of the definition \ref{1.25} of the sound speed $\eta$ (a formula already in conventional units) is:
\begin{equation}
\left(\frac{dp}{dn}\right)_s=\eta^2
\label{1.119}
\end{equation}

From \ref{1.51}, together with \ref{1.101}, \ref{1.102}, we see that in the non-relativistic limit 
the acoustical metric $h$ becomes:
\begin{equation}
h=-\eta^2 dt\otimes dt+(dx^i-v^idt)\otimes(dx^i-v^i dt)
\label{1.120}
\end{equation}
This is a Lorentzian metric on the Galilei spacetime manifold the null cones of which are the 
classical sound cones. Remark that the metric induced by $h$ on the hyperplanes of absolute simultaneity $\Sigma_t$ coincides with the Euclidean metric. In the non-relativistic limit the reciprocal acoustical metric \ref{1.50} becomes: 
\begin{equation}
h^{-1}=-\frac{1}{\eta^2}u\otimes u
+\frac{\partial}{\partial x^i}\otimes\frac{\partial}{\partial x^i}
\label{1.121}
\end{equation}
with $u$ given by \ref{1.104}. This is the reciprocal of \ref{1.120}. 

If we denote in the framework of the non-relativistic theory: 
\begin{equation}
x^0=t
\label{1.144}
\end{equation}
(contrast with \ref{1.101}), then $(x^0,x^1,x^2,x^3)$, with $(x^1,x^2,x^3)$ rectangular coordinates 
in Euclidean space, will be referred to as {\em Galilean coordinates} in spacetime. In such coordinates the components of the acoustical metric \ref{1.120} and those of the reciprocal acoustical metric \ref{1.121} read:
\begin{equation}
h_{00}=-\eta^2+|v|^2, \ \ \ h_{0i}=-v^i, \ \ \ h_{ij}=\delta_{ij}
\label{1.145}
\end{equation}
\begin{equation}
(h^{-1})^{00}=-\eta^{-2}, \ \ \ (h^{-1})^{0i}=-\eta^{-2} v^i, \ \ \ (h^{-1})^{ij}=\delta^{ij}-\eta^{-2}v^i v^j
\label{1.146}
\end{equation}
Note that because of the fact that $x^0$, as defined by \ref{1.144}, has the physical dimensions of 
time while the $x^i \ : \ i=1,2,3$ have the physical dimensions of length, the different components of 
spacetime tensors have different physical dimensions. 

In proceeding to consider the non-relativistic analogue of the 1-form $\beta$, we first remark that 
$\beta$ is actually defined by \ref{1.53} only if the units are chosen so that $c$ is set equal to 1. 
In general, the rectangular components of $\beta$ should be dimensionless quantities. Therefore \ref{1.53} in arbitrary units should read:
\begin{equation}
\beta_\mu=-c^{-2}hu_\mu 
\label{1.122}
\end{equation}
Then equations in equations \ref{1.54}-\ref{1.56} $\beta$ should be replaced by $c^2\beta$ when these equations are written in arbitrary units. The definition \ref{1.57} is unaffected because $\omega$ 
is similarly to be replaced by $c^2\omega$, but the last replacement affects equation \ref{1.59}. 
 
Substituting \ref{1.102} and \ref{1.106} (with $m=1$) in \ref{1.122} we find, for the components of the 
1-form $\beta$, 
\begin{eqnarray}
&&\beta_0=1+c^{-2}\left(h^\prime+\frac{1}{2}|v|^2\right)+O(c^{-4})\nonumber\\
&&\beta_i=-\frac{v^i}{c}+O(c^{-3})
\label{1.123}
\end{eqnarray}
Consider the spacetime 1-form:
\begin{equation}
\beta^\prime=c\beta-c^2 dt
\label{1.124}
\end{equation}
(see \ref{1.101}). Then as $c\rightarrow\infty$ this 1-form becomes the spacetime 1-form:
\begin{equation}
\beta^\prime=\left(h^\prime+\frac{1}{2}|v|^2\right)dt-v^idx^i
\label{1.125}
\end{equation}
In analogy with \ref{1.57} let us define in the non-relativistic limit the {\em spacetime vorticity 2-form} $\omega$ by:
\begin{equation}
\omega=-d\beta^\prime
\label{1.126}
\end{equation}
Substituting \ref{1.125} we find:
\begin{equation}
\omega=\left(\frac{\partial v^i}{\partial t}+\frac{\partial}{\partial x^i}\left(h^\prime+\frac{1}{2}|v|^2\right)\right)dt\wedge dx^i+\frac{1}{2}\omega_{ij}dx^i\wedge dx^j
\label{1.127}
\end{equation}
where 
\begin{equation}
\frac{1}{2}\omega_{ij}dx^i\wedge dx^j:=\overline{\omega}, \ \ \ \omega_{ij}=\frac{\partial v^j}{\partial x^i}-\frac{\partial v^i}{\partial x^j}
\label{1.128}
\end{equation}
is the {\em classical vorticity 2-form} which is a spatial 2-form. In view of \ref{1.124} the non-relativistic spacetime vorticity 2-form $\omega$ is 
the limit as $c\rightarrow\infty$ of $c\omega$ with $\omega$ the relativistic vorticity 2-form. 
Taking then the limit as $c\rightarrow\infty$ of the relativistic equation \ref{1.60} yields 
the equation
\begin{equation}
i_u\omega=\theta ds
\label{1.132}
\end{equation}
in the non-relativistic framework, which is identical in form to the relativistic equation.  
Equation \ref{1.132} can be directly verified. 
Also, taking the limit of the relativistic equation 
\ref{1.56} we obtain in the equation:
\begin{equation}
{\cal L}_u\beta^\prime=d\left(h^\prime-\frac{1}{2}|v|^2\right)-\theta ds
\label{1.d14}
\end{equation}
in the non-relativistic framework, an equation which can also be directly verified. In fact, since 
by \ref{1.125} and \ref{1.104} we have: 
\begin{equation}
i_u\beta^\prime=h^\prime-\frac{1}{2}|v|^2
\label{1.d15}
\end{equation} 
the argument used to deduce \ref{1.60} from \ref{1.56} in the relativistic framework, applies 
in the non-relativistic framework as well to deduce \ref{1.132} from \ref{1.d14} . 
The $u$ component of \ref{1.132} is \ref{1.118} while the restriction of \ref{1.132} to $\Sigma_t$ 
is equivalent, modulo \ref{1.111} to \ref{1.113}. It thus follows that for $C^1$ solutions the 
classical equations of motion are equivalent to \ref{1.132} together with \ref{1.111}. 
Moreover, as in the relativistic theory we can take $(h^\prime,s)$ as the basic thermodynamic variables and consider the equation of state as giving $p$ as a function 
of $h^\prime$ and $s$, in which case the 2nd of \ref{1.107} takes the form:
\begin{equation}
p=p(h^\prime,s), \ \ \ dp=n(dh^\prime-\theta ds)
\label{1.d16}
\end{equation}
Then by \ref{1.119} we have: 
\begin{equation}
\left(\frac{dn}{dh^\prime}\right)_s=\frac{n}{h^\prime}
\label{1.d17}
\end{equation}
Since
$$dn=\left(\frac{dn}{dh^\prime}\right)_s dh^\prime+\left(\frac{dn}{ds}\right)_{h^\prime} ds,$$
it then follows by \ref{1.118} that:
\begin{equation}
\frac{\partial n}{\partial t}+v^i\frac{\partial n}{\partial x^i}=
\frac{n}{\eta^2}\left(\frac{\partial h^\prime}{\partial t}
+v^i\frac{\partial h^\prime}{\partial x^i}\right)
\label{1.d18}
\end{equation}
Therefore equation \ref{1.111} takes the form:
\begin{equation}
\frac{1}{\eta^2}\left(\frac{\partial h^\prime}{\partial t}
+v^i\frac{\partial h^\prime}{\partial x^i}\right)+\frac{\partial v^i}{\partial x^i}=0
\label{1.d19}
\end{equation} 
Substituting from \ref{1.125}:
\begin{equation}
h^\prime=\beta^\prime_0-\frac{1}{2}\beta^\prime_i\beta^\prime_i, \ \ \ v^i=-\beta^\prime_i
\label{1.d20}
\end{equation}
and comparing with \ref{1.146} we see that equation \ref{1.d19}, that is, the equation of continuity \ref{1.111}, reduces to:
\begin{equation}
(h^{-1})^{\mu\nu}\partial_\mu\beta^\prime_\nu=0
\label{1.d21}
\end{equation}
(Galilean coordinates). 

To find the non-relativistic analogue of the function $H$ we first write down the definition \ref{1.d13} in arbitrary units as:
\begin{equation}
1-\frac{h^2}{c^4}H=\frac{\eta^2}{c^2}
\label{1.a16}
\end{equation}
In view of \ref{1.121} (where we have set $m=1$) we write:
\begin{equation}
H=1+c^{-2}H^\prime
\label{1.a17}
\end{equation}
Taking then the limit $c\rightarrow\infty$ we see that $H^\prime$ becomes:
\begin{equation}
H^\prime=-\eta^2-2h^\prime
\label{1.a18}
\end{equation}
This is the function which in the non-relativistic theory corresponds to $H$. Note that the 
relativistic function $H$ is dimensionless while the non-relativistic function $H^\prime$ has the 
physical dimensions of $(\mbox{length}/\mbox{time})^2$. In terms of 1-form $\beta^\prime$ and the function $H^\prime$ the acoustical metric \ref{1.120} takes the form:
\begin{equation}
h=H^\prime dt\otimes dt+dt\otimes\beta^\prime+\beta^\prime\otimes dt+dx^i\otimes dx^i
\label{1.d22}
\end{equation}
its components in Galilean coordinates being given by:
\begin{equation}
h_{00}=H^\prime+2\beta^\prime_0, \ \ \ h_{0i}=\beta^\prime_i, \ \ \ h_{ij}=\delta_{ij}
\label{1.d23}
\end{equation}

The limit of the spacetime vectorfield 
$$c\varpi^\mu\frac{\partial}{\partial x^\mu}$$
where $\varpi$ is the relativistic vorticity vector \ref{1.61}, is the {\em classical vorticity 
vector}:
\begin{equation}
\varpi=\varpi^k\frac{\partial}{\partial x^k}, \ \ \ \varpi^k=\frac{1}{2}(\epsilon^{-1})^{kij}\omega_{ij}
\label{1.133}
\end{equation}
(with $\omega_{ij}$ given by \ref{1.128}), a spatial vectorfield. Here $\epsilon^{-1}$ is the reciprocal volume form of the Euclidean metric, 
\begin{equation}
(\epsilon^{-1})^{kij}=[kij] \ : \ \mbox{the fully antisymmetric 3-dimensional symbol}
\label{1.134}
\end{equation}
in rectangular coordinates. We remark that the local simultaneous spaces of the fluid of the 
relativistic theory having collapsed in the non-relativistic limit to the hyperplanes of absolute 
simultaneity $\Sigma_t$, the vorticity vector no longer appears as an obstruction to integrability. 

The irrotational case is characterized as in the relativistic theory by the vanishing of the spacetime vorticity. Then, assuming as above that the fluid portion under consideration is 
simply connected, we can introduce a function $\phi^\prime$, determined up to an additive constant, 
such that:
\begin{equation}
\beta^\prime=d\phi^\prime
\label{1.135}
\end{equation}
Comparing with \ref{1.125} we see that in the irrotational case:
\begin{equation}
h^\prime=\frac{\partial\phi^\prime}{\partial t}-\frac{1}{2}|\nabla \phi^\prime|^2, \ \ \ 
v^i=-\frac{\partial\phi^\prime}{\partial x^i}
\label{1.136}
\end{equation}
where $|\nabla\phi^\prime|^2$ denotes the squared magnitude of the spatial gradient of $\phi^\prime$:
\begin{equation}
|\nabla \phi^\prime|^2=\frac{\partial\phi^\prime}{\partial x^i}\frac{\partial\phi^\prime}{\partial x^i}
\label{1.137}
\end{equation}
By reason of the first of these equations the function $\phi^\prime$ is called {\em velocity potential}. Note that according to the above $\phi^\prime$ is obtained from the function $\phi$ 
of the relativistic theory as  the limit of $c(\phi-ct)$ as $c\rightarrow\infty$. Note also that 
the relativistic function $\phi$ has the physical dimensions of length while the non-relativistic function $\phi^\prime$ has the physical dimensions of $(\mbox{length})^2/\mbox{time}$. Remark also 
that in the relativistic theory since the derivatives of $\phi$ in rectangular coordinates are dimensionless quantities, equations \ref{1.69} in arbitrary units read:
\begin{equation}
\frac{h}{c^2}=\sqrt{-(g^{-1})^{\mu\nu}\partial_\mu\phi\partial_\nu\phi}, \ \ \ 
u^\mu=-\frac{c^2\partial^\mu\phi}{h}
\label{1.138}
\end{equation}

With \ref{1.135} $\omega$ vanishes so with constant $s$ equation \ref{1.132} is identically satisfied. It follows that, as in the relativistic theory, in the irrotational case the whole content 
of the equations of motion is contained in the equation of continuity \ref{1.111}, which, in view 
of \ref{1.136} and the fact that $n$ is a function of $h^\prime$, takes the form of a nonlinear wave 
equation for $\phi^\prime$. As in the relativistic theory, this equation is the Euler-Lagrange equation corresponding to the Lagrangian density:
\begin{equation}
L=p(h^\prime)
\label{1.139}
\end{equation}
with $h^\prime$ {\em defined} by the 1st of \ref{1.136}. For, $p=p(h^\prime)$ expresses the equation of state along the adiabat corresponding to the given value of $s$, and by the 2nd of \ref{1.107} 
we have:
\begin{equation}
\frac{dp}{dh^\prime}=n
\label{1.140}
\end{equation}
Again $p$ is a convex function of $h^\prime$ and by \ref{1.119}, \ref{1.d17} we now have:
\begin{equation}
\frac{d^2 p}{dh^{\prime 2}}=\frac{dn}{dh^\prime}=\frac{n}{\eta^2}
\label{1.141}
\end{equation}
The nonlinear wave equation for $\phi^\prime$ is equation \ref{1.d21}: 
\begin{equation}
(h^{-1})^{\mu\nu}\frac{\partial^2\phi^\prime}{\partial x^\mu\partial x^\nu}=0
\label{1.149}
\end{equation}
This is similar in form to the relativistic equation \ref{1.80}. Explicitly, equation \ref{1.149} 
reads:
\begin{equation}
-\frac{1}{\eta^2}\left(\frac{\partial^2\phi^\prime}{\partial t^2}
-2\frac{\partial\phi^\prime}{\partial x^i}\frac{\partial^2\phi^\prime}{\partial t\partial x^i}
+\frac{\partial\phi^\prime}{\partial x^i}\frac{\partial\phi^\prime}{\partial x^j}
\frac{\partial^2\phi^\prime}{\partial x^i\partial x^j}\right)+\triangle\phi^\prime=0
\label{1.142}
\end{equation}
where $\triangle$ is the standard Laplacian associated to the Euclidean metric:
\begin{equation}
\triangle=\frac{\partial}{\partial x^i}\frac{\partial}{\partial x^i}
\label{1.143}
\end{equation}

From \ref{1.d22} in the irrotational case the acoustical metric takes the form: 
\begin{equation}
h=H^\prime dt\otimes dt+dt\otimes d\phi^\prime+d\phi^\prime\otimes dt+dx^i\otimes dx^i
\label{1.a19}
\end{equation}
As in the relativistic theory, the case in which, along the adiabat corresponding to the given value of $s$,
\begin{equation}
H^\prime=H^\prime_0 \ \ \mbox{: a constant}
\label{1.a20}
\end{equation}
plays a special role. In this case the acoustical metric is the metric induced on the graph:
\begin{equation}
y^\prime=\frac{1}{2}H^\prime_0 t+\phi^\prime(t,\overline{x}), \ \ : \ t\in\mathbb{R}, \  \overline{x}\in\mathbb{E}^3
\label{1.a21}
\end{equation}
over the hyperplane $y^\prime=0$ in the manifold $\mathbb{R}^5$ with the metric:
\begin{equation}
dt\otimes dy^\prime+dy^\prime\otimes dt+dx^i\otimes dx^i
\label{1.a22}
\end{equation}
This is isometric to $\mathbb{M}^4$, but the hyperplanes of constant $t$ and the hyperplanes of constant $y^\prime$ are now null hyperplanes. To understand why this is so, we recall that 
the velocity potential $\phi^\prime$ of the non-relativistic theory is the limit as $c\rightarrow\infty$ of the function $\phi^\prime$ defined in the framework of the relativistic theory 
by setting 
\begin{equation}
\phi=ct+c^{-1}\phi^\prime
\label{1.a23}
\end{equation}
Going then back to the metric \ref{1.a4} of $\mathbb{M}^4$, we correspondingly set:
\begin{equation}
y=ct+c^{-1}y^\prime
\label{1.a24}
\end{equation}
Recalling also that
$$g=-c^2 dt\otimes dt+dx^i\otimes dx^i$$
we see that the limit $c\rightarrow\infty$ of the metric \ref{1.a4} is indeed the metric \ref{1.a22}. 

By \ref{1.a18}, in the case \ref{1.a20} the sound speed is given by:
\begin{equation}
\eta=\sqrt{-H^\prime_0-2h^\prime},\ \ \mbox{so we have} \ \ h^\prime<-\frac{1}{2}H^\prime_0. 
\label{1.a25}
\end{equation}
Integrating the general equation, from \ref{1.141},
\begin{equation}
\frac{d\log n}{dh^\prime}=\frac{1}{\eta^2}
\label{1.a26}
\end{equation}
with $\eta^2$ given by \ref{1.a25} we obtain:
\begin{equation}
n=\frac{k}{\sqrt{-H^\prime_0-2h^\prime}}, \ \ \mbox{$k$ : a positive constant.}
\label{1.a27}
\end{equation}
Integrating then the general equation \ref{1.140} with $n$ given by \ref{1.a27} yields:
\begin{equation}
p=p_0-k\sqrt{-H^\prime_0-2h^\prime}
\label{1.a28}
\end{equation}
where $p_0$ is an arbitrary constant. Comparing with \ref{1.a27} and with \ref{1.a25} we conclude 
that:
$$p=p_0-\frac{k^2}{n}$$
or:
\begin{equation}
p=p_0-k^2 V
\label{1.a29}
\end{equation}
and:
\begin{equation}
p=p_0-k\eta
\label{1.a30}
\end{equation}
In the limit $h^\prime\rightarrow -\infty$ we have $\eta\rightarrow\infty$, $V\rightarrow\infty$, 
$p\rightarrow -\infty$, while in the limit $h^\prime\rightarrow -H^\prime_0/2$ we have 
$\eta\rightarrow 0$, $V\rightarrow 0$, $p\rightarrow p_0$. Now, the action corresponding to the 
Lagrangian density \ref{1.139} and to a spacetime domain ${\cal D}$ 
is in general 
\begin{equation}
\int_{{\cal D}}pd\mu
\label{1.a31}
\end{equation} 
where $d\mu$ is the volume element \ref{1.c1}. 
From \ref{1.120} we have, in Galilean coordinates, 
\begin{equation}
\mbox{det}h=-\eta^2
\label{1.176}
\end{equation}
(compare with \ref{1.93}). It follows that the volume element of the acoustical metric is given 
in general by:
\begin{equation}
d\mu_h=\eta d\mu
\label{1.a32}
\end{equation}
Thus, by \ref{1.a30}, \ref{1.a32}, the action \ref{1.a31} reduces in the special case to:
\begin{equation}
p_0\mbox{Vol}({\cal D})-k\mbox{Vol}_h({\cal D})
\label{1.a33}
\end{equation}
that is, up to an additive constant proportional to the volume of ${\cal D}$, minus a positive 
multiplicative constant time 
the volume of the graph over ${\cal D}$ considered as a domain in the hyperplane $y^\prime=0$ in $\mathbb{M}^4$ with the metric \ref{1.a22}. It follows that the nonlinear wave equation 
\ref{1.142} of the non-relativistic theory reduces in the special case to the {\em minimal hypersurface equation}, just as in the relativistic theory. 

Going now back to the general non-relativistic equations of motion \ref{1.111}, \ref{1.113}, \ref{1.117}, taking as the unknowns $v^i \ : \ i=1,2,3$, the components of the spatial velocity vectorfield 
\begin{equation}
v=v^i\frac{\partial}{\partial x^i}
\label{1.150}
\end{equation}
the pressure $p$ and the specific entropy $s$, we may write them, for $C^1$ solutions, in the form: 
\begin{eqnarray}
n\left(\frac{\partial v^i}{\partial t}+v^j\frac{\partial v^i}{\partial x^j}\right)+\frac{\partial p}{\partial x^i}=0&\nonumber\\
\frac{\partial p}{\partial t}+v^i\frac{\partial p}{\partial x^i}+\eta^2 n\frac{\partial v^i}{\partial x^i}=0&\nonumber\\
\frac{\partial s}{\partial t}+v^i\frac{\partial s}{\partial x^i}=0&
\label{1.151}
\end{eqnarray}
which is similar to the form \ref{1.27} of relativistic equations of motion. The corresponding 
equations of variation are:
\begin{eqnarray}
n\left(\frac{\partial\dot{v}^i}{\partial t}+v^j\frac{\partial\dot{v}^i}{\partial x^j}\right)
+\frac{\partial\dot{p}}{\partial x^i}&=&\frac{\dot{n}}{n}\frac{\partial p}{\partial x^i} -n\dot{v}^j\frac{\partial v^i}{\partial x^j}\nonumber\\
\frac{\partial\dot{p}}{\partial t}+v^i\frac{\partial\dot{p}}{\partial x^i}+q\frac{\partial\dot{v}^i}{\partial x^i}&=&-\dot{v}^i\frac{\partial p}{\partial x^i}-
\dot{q}\frac{\partial v^i}{\partial x^i}\nonumber\\
\frac{\partial\dot{s}}{\partial t}+v^i\frac{\partial\dot{s}}{\partial x^i}
&=&-\dot{v}^i\frac{\partial s}{\partial x^i}
\label{1.152}
\end{eqnarray}
where we now denote by $q$ the thermodynamic function:
\begin{equation}
q=\eta^2 n
\label{1.153}
\end{equation}
(see \ref{1.32}). To the system \ref{1.152} and the variation $(\dot{v},\dot{p},\dot{s})$, is associated the {\em energy variation current}:
\begin{equation}
\dot{J}=nu|\dot{v}|^2+2\dot{v}\dot{p}+\frac{u}{\eta^2 n}\dot{p}^2
\label{1.154}
\end{equation}
and the {\em entropy variation current}:
\begin{equation}
\dot{K}=u\dot{s}^2
\label{1.155}
\end{equation}
both spacetime vectorfields, $\dot{J}$ a quadratic form in $(\dot{v},\dot{p})$ and $\dot{K}$ a 
quadratic form in $\dot{s}$. The divergences of these spacetime vectorfields are given by (Galilean 
coordinates):
\begin{eqnarray}
\frac{\partial\dot{J}^\mu}{\partial x^\mu}&=&-2n\dot{v}^i\dot{v}^j\frac{\partial v^j}{\partial x^i}+\frac{2}{n}\left(\frac{dn}{ds}\right)_p\dot{s}\dot{v}^i\frac{\partial p}{\partial x^i}\nonumber\\
&&+\frac{1}{q}\left[-2\dot{p}\dot{q}+\dot{p}^2\left(1+\left(\frac{dq}{dp}\right)_s\right)\right]\frac{\partial v^i}{\partial x^i}
\label{1.156}\\
\frac{\partial\dot{K}^\mu}{\partial x^\mu}&=&-2\dot{s}\dot{v}^i\frac{\partial s}{\partial x^i}+\frac{\partial v^i}{\partial x^i}\dot{s}^2
\label{1.157}
\end{eqnarray} 

Consider for any covector $\xi$ in the cotangent space at a given spacetime point the quadratic 
forms $\xi(\dot{J})$, $\xi(\dot{K})$. For any pair of variations 
$(\dot{v},\dot{p},\dot{s})$, $(\dot{v}^\prime,\dot{p}^\prime,\dot{s}^\prime)$ the corresponding symmetric bilinear forms are:
\begin{eqnarray}
&&\xi(\dot{J})((\dot{v},\dot{p}),(\dot{v}^\prime,\dot{p}^\prime))\label{1.158}\\
&& \ \ \ =n\xi(u)<\dot{v},\dot{v}^\prime>
+\xi(\dot{v})\dot{p}^\prime+\xi(\dot{v}^\prime)\dot{p}+
\frac{\xi(u)}{\eta^2 n}\dot{p}\dot{p}^\prime\nonumber
\end{eqnarray}
and
\begin{equation}
\xi(\dot{K})(\dot{s},\dot{s}^\prime)=\xi(u)\dot{s}\dot{s}^\prime
\label{1.159}
\end{equation}
Here $< \ , \ >$ denotes the Euclidean inner product. 
Consider now the set of all covectors $\xi$ at a given spacetime point $x$ such that the symmetric bilinear forms \ref{1.158} 
and \ref{1.159} are degenerate, that is, there is a non-zero variation 
$(\dot{v},\dot{p},\dot{s})$ such that for all variations $(\dot{v}^\prime,\dot{p}^\prime,\dot{s}^\prime)$ we have:
\begin{equation}
\xi(\dot{J})((\dot{v},\dot{p}),(\dot{v}^\prime,\dot{p}^\prime))=0
\label{1.160}
\end{equation}
and:
\begin{equation}
\xi(\dot{K})(\dot{s},\dot{s}^\prime)=0
\label{1.161}
\end{equation}
This defines, as in the relativistic theory, the {\em characteristic subset}  of the cotangent space at spacetime point $x$. Denoting by $\overline{\xi}$ the restriction of $\xi$ to $\Sigma_t$ (the 
hyperplane of absolute simultaneity through $x$) and by $\overline{\xi}^\sharp$ the corresponding 
Euclidean vector, we see that \ref{1.160} is equivalent to the linear system:
\begin{eqnarray}
n\xi(u)\dot{v}+\dot{p}\overline{\xi}^\sharp=0&\nonumber\\
\frac{\xi(u)\dot{p}}{\eta^2 n}+\overline{\xi}(\dot{v})=0&
\label{1.162}
\end{eqnarray}
while \ref{1.161} is equivalent to:
\begin{equation}
\xi(u)\dot{s}=0
\label{1.163}
\end{equation}
It follows that either:
\begin{equation}
\xi(u)=0
\label{1.164}
\end{equation}
in which case:
\begin{equation}
\dot{p}=0 \ \ \mbox{and} \ \ \overline{\xi}(\dot{v})=0
\label{1.165}
\end{equation}
or:
\begin{equation}
-\eta^{-2}(\xi(u))^2+|\overline{\xi}|^2:=(h^{-1})^{\mu\nu}\xi_\mu\xi_\nu=0
\label{1.166}
\end{equation}
in which case:
\begin{equation}
\dot{s}=0 \ \ \mbox{and} \ \ \dot{v}=-\frac{\overline{\xi}^\sharp\dot{p}}
{n\xi(u)}
\label{1.167}
\end{equation}
In \ref{1.166} $h^{-1}$ is the reciprocal acoustical metric \ref{1.121}. 

As in the relativistic theory, the subset of the cotangent space at the spacetime point $x$ 
defined by the 
condition \ref{1.164} is a hyperplane $P^*_x$ while the subset of the cotangent space at $x$ defined by 
the condition \ref{1.165} is a cone $C^*_x$. The dual to $P^*_x$ characteristic subset of the tangent 
space at $x$ is the linear span of $u$, while the dual to $C^*_x$ characteristic subset $C_x$ is the 
set of all non-zero vectors $X$ at $x$ satisfying:
\begin{equation}
h_{\mu\nu}X^\mu X^\nu=0
\label{1.168}
\end{equation}
with $h$ the acoustical metric \ref{1.120}, the {\em sound cone} at $x$. A covector $\xi$ at a spacetime point $x$ whose null space is the hyperplane $\Sigma_t$ of absolute simultaneity through $x$, and which has positive evaluation in the future of $\Sigma_t$, belongs to the interior of the positive component of $C^*_x$. 
It follows from the above that the set of all covectors $\xi$ at $x$ such that the 
quadratic forms $\xi(\dot{J})$, $\xi(\dot{K})$ are both positive definite is the interior of the positive component of $C^*_x$. 

We now consider the irrotational case in the non-relativistic theory, where with the definitions \ref{1.136} the equations of motion reduce to \ref{1.111}, that is:
\begin{eqnarray}
&&\frac{\partial n}{\partial t}-\frac{\partial}{\partial x^i}\left(n\frac{\partial\phi^\prime}{\partial x^i}\right)=0,\label{1.169}\\
\mbox{with} && n=n(h^\prime), \ \ \ h^\prime=\frac{\partial \phi^\prime}{\partial t}
-\frac{1}{2}|\nabla\phi^\prime|^2.\nonumber
\end{eqnarray}
The corresponding equation of variation is:
\begin{equation}
\frac{\partial\dot{n}}{\partial t}-\frac{\partial}{\partial x^i}\left(\dot{n}\frac{\partial\phi^\prime}{\partial x^i}+n\frac{\partial\dot{\phi}^\prime}{\partial x^i}\right)=0
\label{1.171}
\end{equation}
We have:
\begin{equation}
\dot{n}=\frac{n}{\eta^2}\dot{h}^\prime, \ \ \ 
\dot{h}^\prime=\frac{\partial\dot{\phi}^\prime}{\partial t}-\frac{\partial\phi^\prime}{\partial x^i}\frac{\partial\dot{\phi}^\prime}{\partial x^i}=u\dot{\phi}^\prime
\label{1.172}
\end{equation}
(see \ref{1.104}), hence equation \ref{1.171} takes the form:
\begin{equation}
-\frac{\partial}{\partial t}\left(\frac{n}{\eta^2}u\dot{\phi}^\prime\right)
+\frac{\partial}{\partial x^i}\left(\frac{n}{\eta^2}\frac{\partial\phi^\prime}{\partial x^i}u\dot{\phi}^\prime+n\frac{\partial\dot{\phi}^\prime}{\partial x^i}\right)=0
\label{1.173}
\end{equation}
From \ref{1.121} and \ref{1.104} we have, in Galilean coordinates, 
\begin{eqnarray*}
(h^{-1})^{\mu\nu}\frac{\partial\dot{\phi}^\prime}{\partial x^\nu}\frac{\partial}{\partial x^\mu}
 &=&-\frac{1}{\eta^2}(u\dot{\phi}^\prime)u+\frac{\partial\dot{\phi}^\prime}{\partial x^i}\frac{\partial}{\partial x^i}\\
&=&-\frac{1}{\eta^2}(u\dot{\phi}^\prime)\frac{\partial}{\partial t}+\left(\frac{1}{\eta^2}\frac{\partial\phi^\prime}{\partial x^i}u\dot{\phi}^\prime+\frac{\partial\dot{\phi}^\prime}{\partial x^i}\right)\frac{\partial}{\partial x^i}
\end{eqnarray*}
Comparing with \ref{1.173} we see that the equation of variation is:
\begin{equation}
\frac{\partial}{\partial x^\mu}\left(n(h^{-1})^{\mu\nu}\frac{\partial\dot{\phi}^\prime}{\partial x^\nu}\right)=0
\label{1.175}
\end{equation} 
Consider now the {\em conformal acoustical metric}:
\begin{equation}
\tilde{h}=\Omega h
\label{1.177}
\end{equation}
Since in the physical case of 3 spatial dimensions we have:
$$\sqrt{-\mbox{det}\tilde{h}}=\Omega^2\sqrt{-\mbox{det}h}=\Omega^2\eta,$$
(see \ref{1.176}) it follows that:
$$\sqrt{-\mbox{det}\tilde{h}}(\tilde{h}^{-1})^{\mu\nu}=\Omega\eta(h^{-1})^{\mu\nu}.$$
Comparing with \ref{1.175} we see that, choosing
\begin{equation}
\Omega=\frac{n}{\eta}
\label{1.178}
\end{equation}
the equation of variation \ref{1.175} becomes 
$$\frac{\partial}{\partial x^\mu}\left(\sqrt{-\mbox{det}\tilde{h}}(\tilde{h}^{-1})^{\mu\nu}\frac{\partial}{\partial x^\nu}\dot{\phi}^\prime\right)=0$$
which is simply the linear wave equation corresponding to the metric $\tilde{h}$:
\begin{equation}
\square_{\tilde{h}}\dot{\phi}^\prime=0
\label{1.179}
\end{equation}
and in this form it is valid in any system of coordinates. This equation can also be more simply 
derived from the corresponding relativistic equation \ref{1.97} by taking the limit $c\rightarrow\infty$.

We have arrived at a structure which is identical in form to that of the relativistic theory. 
To equation \ref{1.179} and to the variation $\dot{\phi}^\prime$ we then associate the 
{\em variational stress} 
$\dot{T}$, a quadratic form in $d\dot{\phi}^\prime$ at each spacetime point with values in the $T^1_1$-type tensors at the point, given in arbitrary coordinates by a formula identical to \ref{1.98}, but with $\dot{\phi}^\prime$ in the role of $\dot{\phi}$. The statements 
made earlier in regard to the variational stress in the relativistic framework hold then equally well 
in the non-relativistic framework. 

\vspace{5mm}

\section{The Jump Conditions}

Let ${\cal K}$ be a hypersurface in Minkowski spacetime which is time-like  relative to the Minkowski metric $g$, with a neighborhood ${\cal U}$ such that the components $I^\mu$ of the particle current and the components $T^{\mu\nu}$ of the energy-momentum-stress tensor in rectangular coordinates 
are $C^1$ functions of the coordinates in the closure of each connected component of the complement of ${\cal K}$ in ${\cal U}$ but are not continuous across ${\cal K}$. We assume that for each point $x\in{\cal K}$, 
the flow line in ${\cal U}$ initiating at $x$ is contained in one component ${\cal U}_+$ of 
${\cal U}$ which we accordingly call {\em future component}, while the flow line in ${\cal U}$ terminating at $x$ is contained in the other component ${\cal U}_-$ of ${\cal U}$ which we accordingly 
call {\em past component}. The past boundary of ${\cal U}_+$ represents the {\em future side} of ${\cal K}$, while the future boundary of ${\cal U}_-$ represents the {\em past side} of ${\cal K}$. 
We thus have a classical $C^1$ solution of the equations of motion \ref{1.15}, \ref{1.16} in 
${\cal U}_-$ and in ${\cal U}_+$. We call the solution in ${\cal U}_-$ {\em past solution} and 
the solution in ${\cal U}_+$ {\em future solution}. For any fluid quantity $q$ we denote by $q_-$ and 
$q_+$ the corresponding functions induced on ${\cal K}$ by the past and future solutions respectively. 
Thus $q_-$ and $q_+$ are functions on ${\cal K}$ representing $q$ on the past and future sides of ${\cal K}$ respectively. We denote by
\begin{equation} 
\triangle q=q_+-q_-
\label{1.180}
\end{equation}
the function on ${\cal K}$ which represents the jump in $q$ across ${\cal K}$ in going from the past side to the future side. 

Consider the 3-form $I^*$  dual to the particle current vectorfield $I$, that is, in components 
relative to arbitrary coordinates, 
\begin{equation}
I^*_{\alpha\beta\gamma}=I^\mu\epsilon_{\mu\alpha\beta\gamma}
\label{1.181}
\end{equation}
where $\epsilon_{\mu\alpha\beta\gamma}$ are the components of the volume form of the Minkowski 
metric $g$ (see \ref{1.63}). Then the differential particle conservation law \ref{1.15} is equivalent to:
\begin{equation}
dI^*=0
\label{1.182}
\end{equation}

In regard to the energy-momentum-stress tensor, given a vectorfield $X$, let $P_X$ be the vectorfield 
given in terms of components relative to arbitrary coordinates by:
\begin{equation}
P_X^\mu=-T^{\mu\nu}g_{\nu\lambda}X^\lambda
\label{1.183}
\end{equation}
Then by the differential energy-momentum conservation law, in view of the symmetry of $T^{\mu\nu}$, we have:
\begin{equation}
\nabla_\mu P_X^\mu=-\frac{1}{2}T^{\mu\nu}\left({\cal L}_X g\right)_{\mu\nu}
\label{1.184}
\end{equation}
where ${\cal L}_X g$ is the Lie derivative of $g$ with respect to $X$. Thus if $X$ is a Killing field 
of $g$ we have:
\begin{equation}
\nabla_\mu P_X^\mu=0
\label{1.185}
\end{equation}
This holds in particular if $X$ is a {\em translation field}, so $X$ is in rectangular coordinates 
given by:
\begin{equation} 
X=k^\alpha\frac{\partial}{\partial x^\alpha} \ \ \ \mbox{where the coefficients $k^\alpha \ : \ \alpha=0,1,2,3$ are constants}
\label{1.186}
\end{equation}
Conversely, \ref{1.185} holding for all translation fields $X$ is equivalent to \ref{1.16}. Consider then, for each translation field $X$ the 3-form $P^*_X$ dual to the vectorfield $P_X$:
\begin{equation}
(P^*_X)_{\alpha\beta\gamma}=P^\mu_X\epsilon_{\mu\alpha\beta\gamma}
\label{1.187}
\end{equation}
Then \ref{1.185} is equivalent to:
\begin{equation}
dP^*_X=0
\label{1.188}
\end{equation}
Thus, \ref{1.188} holding for all translation fields $X$ is equivalent to the differential energy-momentum conservation law \ref{1.16}. 

Integrating \ref{1.182} and \ref{1.188} on a spacetime domain ${\cal V}$ we obtain, by 
Stokes' theorem, 
\begin{equation}
\int_{\partial{\cal V}}I^*=0
\label{1.189}
\end{equation}
and 
\begin{equation}
\int_{\partial{\cal V}}P^*_X=0
\label{1.190}
\end{equation}
While the {\em differential conservation laws} \ref{1.182} and \ref{1.188} require $I^*$ and $P^*_X$ to be $C^1$ to be classically interpreted, the {\em integral conservation laws} \ref{1.189} and \ref{1.190} only require  
integrability of $I^*$ and $P^*_X$ on hypersurfaces to be meaningful. Therefore, in the present context, we stipulate that \ref{1.189} and \ref{1.190} hold not only when ${\cal V}$ is wholly contained either in 
${\cal U}_-$ or in ${\cal U}_+$ but also when ${\cal V}$ is an arbitrary subdomain of ${\cal U}$ 
which may intersect both. 

Consider now an arbitrary point $x\in{\cal K}$ and let ${\cal W}_x\subset{\cal U}$ be a neighborhood 
of $x$ in spacetime. We denote by ${\cal S}_x$ the corresponding neighborhood ${\cal S}_x={\cal K}\bigcap{\cal W}_x$ of $x$ in ${\cal K}$. Let $Y$ be a vectorfield without critical points on some larger neighborhood 
${\cal W}^\prime_x\supset {\cal W}_x$ and transversal to ${\cal K}$. Let $L_\delta(y)$ denote the 
segment of the integral curve of $Y$ through $y\in {\cal S}_x$ corresponding to the parameter interval 
$(-\delta,\delta)$:
$$L_\delta(y)=\{ f_s(y) \ : \ s\in (-\delta,\delta) \}$$
where $f_s$ is the local 1-parameter group of diffeomorphisms generated by $Y$. We then define the 
neighborhood ${\cal V}_{x,\delta}$ of $x$ in spacetime by:
$${\cal V}_{x,\delta}=\bigcup_{y\in {\cal S}_x}L_\delta(y)$$
Equations \ref{1.189}, \ref{1.190} with ${\cal V}_{x\delta}$ in the role of ${\cal V}$ read:
\begin{eqnarray}
\int_{\partial{\cal V}_{x,\delta}}I^*=0&\label{1.191}\\
\int_{\partial{\cal V}_{x,\delta}}P^*_X=0\label{1.192}
\end{eqnarray}
Now the boundary of ${\cal V}_{x,\delta}$ consists of the hypersurfaces:
$${\cal S}_{x,-\delta}=\{ f_{-\delta}(y) \ : \ y\in {\cal S}_x \}, \ \ \ 
{\cal S}_{x,\delta}=\{ f_\delta(y) \ : \ y\in {\cal S}_x \}$$
together with the lateral hypersurface:
$$\bigcup_{y\in\partial {\cal S}_x}L_\delta(y)$$
Since this lateral component is bounded in measure by a constant multiple of $\delta$, taking the limit $\delta\rightarrow 0$ in \ref{1.191}, \ref{1.192} we obtain:
\begin{eqnarray}
\int_{{\cal S}_x}\triangle I^*=0&\label{1.193}\\
\int_{{\cal S}_x}\triangle P^*_X=0&\label{1.194}
\end{eqnarray}
These being valid for any neighborhood ${\cal S}_x$ of $x$ in ${\cal K}$ implies that the corresponding 3-forms induced on ${\cal K}$ from the 2-sides coincide at $x$, or, equivalently, that:
\begin{equation}
\xi_\mu\triangle I^\mu=0, \  \  \ \xi_\mu\triangle P^\mu_X=0
\label{1.195}
\end{equation}
where $\xi$ is a covector at $x$ so that the null space of $\xi$ is $T_x{\cal K}$, the tangent space to ${\cal K}$ at $x$. We may orient $\xi$ so that:
\begin{equation}
\mbox{$\xi(V)>0$ for any vector $V$ at $x$ pointing to ${\cal U}_+$}
\label{1.196}
\end{equation}
Moreover, we can normalize $\xi$ to be of unit magnitude with respect to the Minkowski metric $g$:
\begin{equation}
(g^{-1})^{\mu\nu}\xi_\mu\xi_\nu=1
\label{1.197}
\end{equation}
Then the vector $\xi^\sharp$ corresponding to the covector $\xi$ through $g$, that is, the vector with 
componets $(g^{-1})^{\mu\nu}\xi_\nu$, is the unit normal with respect to $g$ to ${\cal K}$ at $x$, pointing to ${\cal U}_+$. This is a space-like relative to $g$ vector. The 2nd of \ref{1.195} holding for all translation fields $X$ is 
equivalent to the condition:
\begin{equation}
\xi_\mu\triangle T^{\mu\nu}=0
\label{1.198}
\end{equation}
The 1st of \ref{1.195} together with \ref{1.198} constitute the {\em jump conditions} across 
${\cal K}$. As we have seen, these follow from the integral form \ref{1.189}, \ref{1.190} of the 
particle and energy-momentum conservation laws. Conversely, a $C^1$ solution of the differential 
conservation laws \ref{1.15}, \ref{1.16} in the closure of ${\cal U}_-$ and of ${\cal U}_+$ 
which satisfies across ${\cal K}$ the jump conditions, satisfies the integral 
conservation laws in ${\cal U}$. 

Let us denote:
\begin{equation}
u_{\bot\pm}=\xi(u_{\pm})
\label{1.199}
\end{equation}
According to the above we have:
\begin{equation}
u_{\bot\pm}>0
\label{1.200}
\end{equation}
The jump condition 1st of \ref{1.195} reads:
\begin{equation}
n_+ u_{\bot+}=n_- u_{\bot-}:=f
\label{1.201}
\end{equation}
The positive function $f$ on ${\cal K}$ is the {\em particle flux} through ${\cal K}$. 
In terms of the specific volume $V$ \ref{1.201} reads:
\begin{equation}
u_{\bot-}=fV_-, \ \ \ u_{\bot+}=fV_+, \ \ \ f>0
\label{1.202}
\end{equation}
In view of the definition \ref{1.13}, the jump condition \ref{1.198} reads:
\begin{equation}
(\rho_++p_+)u_{\bot+}u_++p_+\xi^\sharp=(\rho_-+p_-)u_{\bot-}u_-+p_-\xi^\sharp
\label{1.203}
\end{equation}
Therefore at each point $x$ of ${\cal K}$ the vectors $u_-$, $u_+$, and $\xi^\sharp$, the normal  (relative to $g$) vector to ${\cal K}$ at $x$ pointing to ${\cal U}_+$, all three lie in the same time-like (relative to $g$) plane $P_x$. Substituting \ref{1.201} in \ref{1.203} the latter reduces in view of \ref{1.9} to:
\begin{equation}
fh_+u_++p_+\xi^\sharp=fh_-u_-+p_-\xi^\sharp
\label{1.204}
\end{equation}
We may also write \ref{1.204} using the 1-form $\beta$ defined by \ref{1.53} as:
\begin{equation}
f\beta_+-p_+\xi=f\beta_--p_-\xi
\label{1.205}
\end{equation}
Evaluating $\xi$ on each side of \ref{1.204} and substituting, in view of \ref{1.200}, from \ref{1.202}, we obtain:
\begin{equation}
p_+-p_-=-f^2(h_+V_+-h_-V_-)
\label{1.206}
\end{equation}
On the other hand, taking the $g$-inner product of each side of \ref{1.204} with itself we obtain:
\begin{equation}
p_+^2-p_-^2=f^2(h_+^2-h_-^2-2h_+V_+p_++2h_-V_-p_-)
\label{1.207}
\end{equation}
Comparing \ref{1.206} and \ref{1.207} we conclude that:
\begin{equation}
h_+^2-h_-^2=(h_+V_++h_-V_-)(p_+-p_-)
\label{1.208}
\end{equation}
This is the relativistic {\em Hugoniot relation}. 

We now consider $\triangle s$, the jump in specific entropy across ${\cal K}$. We shall relate this to 
$\triangle p$, the jump in pressure, assuming that both are small. We may then expand $\triangle h$, 
the jump in specific enthalpy, in powers of $\triangle p$ and $\triangle s$. By equation \ref{1.11} 
we have:
\begin{eqnarray}
&&\triangle h=V_-\triangle p+\frac{1}{2}\left(\frac{\partial V}{\partial p}\right)_-(\triangle p)^2
+\frac{1}{6}\left(\frac{\partial^2 V}{\partial p^2}\right)_-(\triangle p)^3+\theta_-\triangle s\nonumber\\
&&\hspace{6mm}+O((\triangle p)^4)
+O(\triangle p \triangle s)+O((\triangle s)^2)
\label{1.209}
\end{eqnarray}
Hence:
\begin{eqnarray}
&&h_+^2-h_-^2=2h_-V_-\triangle p+\left[h_-\left(\frac{\partial V}{\partial p}\right)_-+V_-^2\right](\triangle p)^2\nonumber\\
&&\hspace{10mm}+\left[\frac{h_-}{3}\left(\frac{\partial^2 V}{\partial p^2}\right)_-+V_-\left(\frac{\partial V}{\partial p}\right)_-\right](\triangle p)^3+2h_-\theta_-\triangle s\nonumber\\
&&\hspace{10mm}+O((\triangle p)^4)+O(\triangle p \triangle s)+O((\triangle s)^2)
\label{1.210}
\end{eqnarray}
Also,
\begin{equation}
\triangle V=\left(\frac{\partial V}{\partial p}\right)_-\triangle p+\frac{1}{2}
\left(\frac{\partial^2 V}{\partial p^2}\right)_-(\triangle p)^2+O((\triangle p)^3)+O(\triangle s)
\label{1.211}
\end{equation}
hence:
\begin{eqnarray}
&&(h_+V_++h_-V_-)(p_+-p_-)=2h_-V_-+\left[h_-\left(\frac{\partial V}{\partial p}\right)_-+V_-^2\right]
(\triangle p)^2\nonumber\\
&&\hspace{25mm}+\left[\frac{h_-}{2}\left(\frac{\partial^2 V}{\partial p^2}\right)_-+\frac{3V_-}{2}
\left(\frac{\partial V}{\partial p}\right)_-\right](\triangle p)^3\nonumber\\
&&\hspace{25mm}+O((\triangle p)^4)+O(\triangle s\triangle p)
\label{1.212}
\end{eqnarray}
Comparing \ref{1.210}, \ref{1.212} with the Hugoniot relation \ref{1.208} we conclude that:
\begin{equation}
\triangle s=\frac{1}{12h_-\theta_-}\left[h_-\left(\frac{\partial^2 V}{\partial p^2}\right)_-
+3V_-\left(\frac{\partial V}{\partial p}\right)_-\right](\triangle p)^3+O((\triangle p)^4)
\label{1.213}
\end{equation}

We shall now relate the coefficient in square brackets of the $(\triangle p)^3$ term in \ref{1.213} to the function 
$(dH/dh)_s$
induced on ${\cal K}$ by the past solution, $H$ being the function defined by \ref{1.d13}. We have: 
$$\left(\frac{dH}{dh}\right)_s=\frac{(dH/dp)_s}{(dh/dp)_s}=\frac{1}{V}\left(\frac{dH}{dp}\right)_s$$
hence by \ref{1.d13}:
\begin{equation}
Vh^2\left(\frac{dH}{dh}\right)_s=-\left(\frac{d(\eta^2)}{dp}\right)_s-\frac{2V}{h}(1-\eta^2)
\label{1.214}
\end{equation}
By \ref{1.d4}:
\begin{equation}
\frac{1}{\eta^2}=-\frac{h}{V}\left(\frac{dV}{dh}\right)_s=-\frac{h}{V^2}\left(\frac{dV}{dp}\right)_s
\label{1.215}
\end{equation}
Substituting \ref{1.215} and its derivative with respect to $p$ at constant $s$ in \ref{1.214} we find:
\begin{equation}
-\frac{V^3 h^2}{\eta^4}\left(\frac{dH}{dh}\right)_s=
h\left(\frac{d^2 V}{dp^2}\right)_s+3V\left(\frac{dV}{dp}\right)_s
\label{1.216}
\end{equation}

Assuming that the past solution is irrotational, so 
\begin{equation}
\omega_-=0
\label{1.217}
\end{equation}
shall now determine $\triangle\omega$, the jump in vorticity across ${\cal K}$, relating this to $\triangle s$, the jump in specific entropy. 
Let us denote by $\beta_{||-}$ and $\beta_{||+}$ the 1-forms $\beta$ induced on ${\cal K}$ by the 
past solution and the future solution respectively:
\begin{equation}
\beta_{||\pm}=\left.\beta_{\pm}\right|_{T{\cal K}}
\label{1.219}
\end{equation}
Restricting \ref{1.205} to $T{\cal K}$ we obtain:
\begin{equation}
\beta_{||+}=\beta_{||-}
\label{1.220}
\end{equation}
Let us denote by $\omega_{||-}$ and $\omega_{||+}$ the vorticity 2-forms induced on ${\cal K}$ by 
the past solution and the future solution respectively:
\begin{equation}
\omega_{||\pm}=\left.\omega_{\pm}\right|_{T{\cal K}}
\label{1.221}
\end{equation}
We then have:
\begin{equation}
\omega_{||\pm}=-d_{||}\beta_{||\pm}
\label{1.222}
\end{equation}
where $d_{||}$ denotes the exterior derivative on ${\cal K}$. Then \ref{1.220} implies:
\begin{equation}
\omega_{||-}=\omega_{||+}
\label{1.223}
\end{equation}
therefore the assumption \ref{1.217} implies:
\begin{equation}
\omega_{||+}=0
\label{1.224}
\end{equation}
The remaining component of the vorticity along ${\cal K}$ is $\omega_{\bot-}$, corresponding to 
the past solution, $\omega_{\bot+}$, corresponding to the future solution, where:
\begin{equation}
\omega_{\bot\pm}=\left.i_{\xi^\sharp}\omega_{\pm}\right|_{T{\cal K}}
\label{1.225}
\end{equation}
The assumption \ref{1.217} implies that also $\omega_{\bot-}=0$. In general, if we extend 
$\omega_{||\pm}$ and $\omega_{\bot\pm}$ to the tangent bundle to the spacetime manifold over ${\cal K}$ by the condition that they vanish if one of their entries is $\xi^\sharp$, we can express:
\begin{equation}
\omega_{\pm}=\omega_{||\pm}+\xi\wedge\omega_{\bot\pm}
\label{1.226}
\end{equation}
Here by \ref{1.224} we have: 
\begin{equation}
\omega_+=\xi\wedge\omega_{\bot+}
\label{1.227}
\end{equation}
Substituting this in \ref{1.60} and restricting the resulting equation to $T{\cal K}$ yields: 
\begin{equation}
\omega_{\bot+}=\frac{\theta_+}{u_{\bot+}}d_{||}s_+
\label{1.231}
\end{equation}
Moreover, since $s_-$ is constant, the past solution being isentropic, we may replace $s_+$ on the right by the jump $\triangle s$:
\begin{equation}
\omega_{\bot+}=\frac{\theta_+}{u_{\bot+}}d_{||}\triangle s
\label{1.232}
\end{equation}
Consider finally the vorticity vectors $\varpi_-$ and $\varpi_+$ on ${\cal K}$ corresponding to the 
past and future solutions respectively. By reason of the assumption \ref{1.217} we have $\varpi_-=0$, 
while by \ref{1.227} and \ref{1.232},
\begin{equation}
\varpi_+^\mu=\frac{\theta_+}{u_{\bot+}}(\epsilon^{-1})^{\mu\alpha\beta\gamma}u_{+\alpha}\xi_\beta\partial_\gamma\triangle s
\label{1.233}
\end{equation}
Let at each $x\in {\cal K}$, $\Pi_x$ be the space-like plane which is $g$-orthogonal complement of the time-like plane $P_x$. 
We have:
\begin{equation}
\Pi_x=\Sigma_{x-}\bigcap T_x{\cal K}=\Sigma_{x+}\bigcap T_x{\cal K}
\label{1.234}
\end{equation}
where $\Sigma_{x-}$ and $\Sigma_{x+}$ are the simultaneous spaces of the fluid at $x$ corresponding to 
the past and future solutions respectively. Let us denote by $\sPi$ the $g$-orthogonal projection to 
$\Pi_x$ at each $x\in{\cal K}$, and for any function $f$ defined on ${\cal K}$, let us denote by 
$\sd f$ its derivative tangentially to $\Pi_x$, at each $x\in{\cal K}$. Then only $\sd\triangle s$ 
contributes to $\varpi_+$. Denoting by $\sd\triangle s^\sharp$ the corresponding, through 
$\left. g\right|_{\Pi_x}$, vector, element of $\Pi_x$, \ref{1.233} takes the form:
\begin{equation}
\varpi_+=\frac{\theta_+}{u_{\bot+}}\s^*\sd\triangle s^\sharp
\label{1.235}
\end{equation}
where for $V\in\Pi_x$ we denote by $\s^*V$ the result of rotating $V$ counterclockwise, relative to 
$\xi^\sharp$, by a right angle. 

We now come to the {\em determinism conditions}. These are two conditions. The {\em first determinism condition} is that the hypersurface of discontinuity ${\cal K}$ is {\em acoustically space-like} 
relative to the past solution. That is, at each $x\in {\cal K}$, the tangent space $T_x{\cal K}$ 
intersects the closure of the interior of $C_{x-}$, the sound cone at $x$ corresponding to the past 
solution, only at the origin. The {\em second determinism condition} is that ${\cal K}$ is 
{\em acoustically time-like} relative to the future solution. That is, at each $x\in {\cal K}$, 
$T_x{\cal K}$ intersects the interior of $C_{x+}$, the sound cone at $x$ corresponding to the future 
solution. The first determinism condition has the effect that the physical conditions along 
${\cal K}$ corresponding to the past solution are determined by the initial conditions which 
determine the past solution. The second determinism condition has the effect that the knowledge 
which derives from the jump conditions of the  physical conditions along ${\cal K}$ corresponding 
to the future solution, when complemented by the initial conditions, determines the future solution. 
The two determinism conditions together with the jump conditions determine the complete solution  
including the hypersurface of discontinuity  itself. Note that since at each spacetime point the light cone contains the sound cone, the second determinism condition implies our initial assumption that 
${\cal K}$ is time-like relative to the Minkowski metric $g$. 

In terms of $h_-^{-1}$, the reciprocal acoustical metric along ${\cal K}$ corresponding to the past 
solution, the first determinism condition reads:
\begin{equation}
(h_-^{-1})^{\mu\nu}\xi_\mu\xi_\nu<0
\label{1.236}
\end{equation}
In terms of $h_+^{-1}$, the reciprocal acoustical metric along ${\cal K}$ corresponding to the future 
solution, the second determinism condition reads:
\begin{equation}
(h_+^{-1})^{\mu\nu}\xi_\mu\xi_\nu>0
\label{1.237}
\end{equation}
From \ref{1.50} and \ref{1.199} we obtain:
\begin{equation}
(h_{\pm}^{-1})^{\mu\nu}\xi_\mu\xi_\nu=1-\left(\frac{1}{\eta_{\pm}^2}-1\right)(u_{\bot\pm})^2
\label{1.238}
\end{equation}
Therefore the determinism conditions \ref{1.236}, \ref{1.237} are:
\begin{equation}
u_{\bot-}>\frac{\eta_-}{\sqrt{1-\eta_-^2}}, \ \ \ 
u_{\bot+}<\frac{\eta_+}{\sqrt{1-\eta_+^2}}
\label{1.239}
\end{equation}
Substituting from \ref{1.202} these become:
\begin{equation}
f>\frac{\eta_-/V_-}{\sqrt{1-\eta_-^2}}, \ \ \ 
f<\frac{\eta_+/V_+}{\sqrt{1-\eta_+^2}}
\label{1.240}
\end{equation}
We conclude that when the jump conditions are satisfied the determinism conditions reduce to:
\begin{equation}
\frac{\eta_-/V_-}{\sqrt{1-\eta_-^2}}<\frac{\eta_+/V_+}{\sqrt{1-\eta_+^2}}
\label{1.241}
\end{equation}

Defining the quantity:
\begin{equation}
w=\left(\frac{1}{\eta^2}-1\right)V^2
\label{1.242}
\end{equation}
the condition \ref{1.241} is seen to be equivalent to:
\begin{equation}
\triangle w<0
\label{1.243}
\end{equation}
Taking as in the preceding $p$ and $s$ as the basic thermodynamic variables, we have:
\begin{equation}
\frac{1}{\eta^2}=\frac{\partial\rho}{\partial p}=-\frac{h}{V^2}\frac{\partial V}{\partial p}
\label{1.244}
\end{equation}
hence $w$ is given by:
\begin{equation}
w=-h\frac{\partial V}{\partial p}-V^2
\label{1.245}
\end{equation}
We then obtain:
\begin{equation}
\frac{\partial w}{\partial p}=-h\frac{\partial^2 V}{\partial p^2}-3V\frac{\partial V}{\partial p}
\label{1.246}
\end{equation}
In view of the fact that by \ref{1.213} $\triangle s=O(\triangle p)^3)$, it follows that:
\begin{equation}
\triangle w=-\left[h_-\left(\frac{\partial^2 V}{\partial p^2}\right)_-+
3V_-\left(\frac{\partial V}{\partial p}\right)_-\right]\triangle p+O((\triangle p)^2)
\label{1.247}
\end{equation}
Remark that the coefficient in square brackets of the $\triangle p$ term in \ref{1.247} coincides 
with the coefficient in square brackets of the $(\triangle p)^3$ term in \ref{1.213}. Therefore, at 
least for suitably small $\triangle p$, condition \ref{1.243} is equivalent to the {\em entropy 
condition}:
\begin{equation}
\triangle s>0
\label{1.248}
\end{equation}
The fact that the requirement of determinism and that of the second law of thermodynamics coincide 
in the present context recalls the interpretation of entropy increase as loss of information. 
Note also that in view of \ref{1.216} this requirement implies that $\triangle p$ is positive 
or negative according as to whether $(dH/dh)_{s-}$ is negative or positive. 

We shall now discuss briefly the jump conditions in the non-relativistic theory. These are obtained 
in a straightforward manner by taking the limit $c\rightarrow\infty$ as explained in the previous 
section. In particular, in the non-relativistic theory we have the particle current vectorfield 
\ref{1.110} limit of the corresponding vectorfield of the relativistic theory, the {\em energy current} vectorfield:
\begin{equation}
E=\left(\frac{1}{2}n|v|^2+\rho^\prime\right)\frac{\partial}{\partial t}
+\left(\frac{1}{2}n|v|^2+\rho^\prime+p\right)v^i\frac{\partial}{\partial x^i} 
\label{1.250}
\end{equation}
(see \ref{1.117}) limit of the vectorfield 
$$T^{0\mu}-c^2 I^\mu$$
of the relativistic theory as explained in the previous section (see \ref{1.114}), and, for each 
$i=1,2,3$, the momentum current vectorfield:
\begin{equation}
\overline{P}^i=nv^i\frac{\partial}{\partial t}+(nv^iv^j+p\delta^{ij})\frac{\partial}{\partial x^j}
\label{1.251}
\end{equation}
(see \ref{1.113}). These vectorfields are obtained by taking the limit of the vectorfield $P_X$ of the relativistic theory, given by \ref{1.183}, in which we set $X$ to be each of the spatial translation fields: 
\begin{equation}
\overline{X}=\frac{\partial}{\partial x^i} \ \ : \ i=1,2,3
\label{1.252}
\end{equation}
corresponding to a rectangular coordinate system in $\mathbb{E}^3$. The duals of the above five spacetime vectorfields relative to the spacetime volume form \ref{1.c1} are the 3-forms: 
$$I^*, \ \ E^*, \ \ \overline{P}^{i*}\ : \ i=1,2,3$$
and in the non-relativistic theory the following integral conservation laws are stipulated:
\begin{equation}
\int_{\partial{\cal V}}I^*=0,
\label{1.253}
\end{equation}
\begin{equation}
\int_{\partial{\cal V}}E^*=0
\label{1.254}
\end{equation}
and
\begin{equation}
\int_{\partial{\cal V}}\overline{P}^{i*}=0 \ \ : \ i=1,2,3
\label{1.255}
\end{equation}
the conservation laws of mass, energy, and momentum respectively. As in the relativistic theory 
the integral conservation laws yield the jump conditions:
\begin{equation}
\xi_\mu\triangle I^\mu=0, \ \ \ \xi_\mu\triangle E^\mu=0, \ \ \ 
\xi_\mu\triangle(\overline{P}^i)^\mu=0 \ : \ i=1,2,3
\label{1.256}
\end{equation}
(Galilean coordinates). As in the relativistic theory, at each $x\in{\cal K}$, $\xi$ is a covector at 
$x$ so that the null space of $\xi$ is $T_x{\cal K}$.  We may 
again orient $\xi$ according to \ref{1.196}. However we cannot normalize $\xi$ as in \ref{1.197}, 
because of the absence of a metric on Galilei spacetime. Instead, let $x$ belong to the hyperplane 
of absolute simultaneity $\Sigma_t$. The hypersurface ${\cal K}$ being transversal to the $\Sigma_t$, 
the intersection $\Sigma_t\bigcap {\cal K}$ is a surface ${\cal K}_t$ and 
$\Sigma_t\bigcap T_x{\cal K}$ is the tangent plane to ${\cal K}_t$ at $x$. If we denote by $\overline{\xi}$ the restriction of $\xi$ to $\Sigma_t$, then $\overline{\xi}^\sharp$, the corresponding Euclidean vector, is a normal in $\Sigma_t$ to ${\cal K}_t$ at $x$, pointing 
to the interior of $\Sigma_t\bigcap{\cal U}_+$. Here we are appealing to the Euclidean structure 
of $\Sigma_t$. We now normalize $\xi$ by the condition that $\overline{\xi}^\sharp$ is of unit magnitude, or:
\begin{equation}
\overline{\xi}(\overline{\xi}^\sharp)=1
\label{1.257}
\end{equation}
We then have:
\begin{equation}
\xi=-\nu dt+\overline{\xi}
\label{1.258}
\end{equation}
where $\nu$ is the velocity at $x$ of ${\cal K}_t$ in $\Sigma_t$ in the direction of the unit normal
$\overline{\xi}^\sharp$. 

In view of \ref{1.104} and \ref{1.258} the definitions \ref{1.199} become, in the non-relativistic theory,
\begin{equation}
u_{\bot\pm}=v_{\bot\pm}-\nu
\label{1.259}
\end{equation}
where:
\begin{equation}
v_{\bot\pm}=\overline{\xi}(v_\pm)
\label{1.260}
\end{equation}
The condition \ref{1.200} again holds and the jump condition 1st of \ref{1.256} again takes the form 
\ref{1.201}, with $u_\pm$ now given by \ref{1.259}, so the {\em particle flux} $f$ is defined in 
the same way, and in terms of the specific volume $V$ we have, as in \ref{1.202},
\begin{equation}
u_{\bot-}=fV_-, \ \ \ u_{\bot+}=fV_+, \ \ \ f>0
\label{1.261}
\end{equation}
In view of the definition \ref{1.251} and \ref{1.201} the last of the jump conditions \ref{1.256} read:
\begin{equation}
fv_++p_+\overline{\xi}^\sharp=fv_-+p_-\overline{\xi}^\sharp
\label{1.262}
\end{equation}
Let us denote by $v_{||\pm}$ the Euclidean projection of $v_{\pm}$ to the plane orthogonal to $\overline{\xi}^\sharp$. Then \ref{1.263} implies:
\begin{equation}
v_{||+}=v_{||-}
\label{1.263}
\end{equation}
Evaluating $\overline{\xi}$ on each side of \ref{1.262}, taking into account the fact that by \ref{1.259} and \ref{1.261}:
\begin{equation}
v_{\bot+}-v_{\bot-}=u_{\bot+}-u_{\bot-}=f(V_+-V_-)
\label{1.264}
\end{equation}
we obtain, in view of \ref{1.257},
\begin{equation}
p_+-p_-=-f(v_{\bot+}-v_{\bot-})=-f^2(V_+-V_-)
\label{1.265}
\end{equation}
In view of the definition \ref{1.250}, \ref{1.201}, and \ref{1.109}, 
the jump condition 2nd of \ref{1.256} reads: 
\begin{equation}
f\left[\frac{1}{2}\left(|v_+|^2-|v_-|^2\right)+e^\prime_+-e^\prime_-\right]
+v_{\bot+}p_+-v_{\bot-}p_-=0
\label{1.266}
\end{equation}
By \ref{1.263} we have:
$$|v_+|^2-|v_-|^2=v_{\bot+}^2-v_{\bot-}^2$$
Writing
$$v_{\bot+}p_+-v_{\bot-}p_-=\frac{1}{2}(v_{\bot+}+v_{\bot-})(p_+-p_-)+(v_{\bot+}-v_{\bot-})(p_+-p_-)$$
we then see that:
\begin{eqnarray*}
&&\frac{1}{2}f\left(|v_+|^2-|v_-|^2\right)+v_{\bot+}p_+-v_{\bot-}p_-=\\
&&\hspace{20mm}\frac{1}{2}(v_{\bot+}+v_{\bot-})\left[f(v_{\bot+}-v_{\bot-})+p_+-p_-\right]\\
&&\hspace{20mm}+\frac{1}{2}(v_{\bot+}-v_{\bot-})(p_++p_-)
\end{eqnarray*}
By \ref{1.264} and \ref{1.265} the first term on the right vanishes, therefore, again by \ref{1.264}, 
this reduces to:
$$f(V_+-V_-)(p_++p_-)$$
Substituting then in \ref{1.266} yields the classical {\em Hugoniot relation:}
\begin{equation}
e^\prime_+-e^\prime_-+\frac{1}{2}(p_++p_-)(V_+-V_-)=0
\label{1.267}
\end{equation}
or, in terms of the enthalpy:
\begin{equation}
h^\prime_+-h^\prime_-=\frac{1}{2}(V_++V_-)(p_+-p_-)
\label{1.268}
\end{equation}
The same is immediately derived from the relativistic relation \ref{1.208} substituting \ref{1.106} 
and taking the limit $c\rightarrow\infty$. 

Going back to condition \ref{1.266} and taking into 
account the fact that by \ref{1.259}, \ref{1.261}, 
$$v_{\bot\pm}-fV_\pm=\nu,$$
we see that this condition takes in terms of the enthalpy the form:
\begin{equation}
f\left[h^\prime_+-h^\prime_-+\frac{1}{2}\left(|v_+|^2-|v_-|^2\right)\right]
+\nu (p_+-p_-)=0
\label{1.269}
\end{equation}
Comparing \ref{1.269} and \ref{1.262} with \ref{1.258} and \ref{1.125} we see that the 2nd and last of the jump conditions \ref{1.256} can be combined to the following condition in terms of the 1-form 
$\beta^\prime$:
\begin{equation}
f\beta^\prime_+-p_+\xi=f\beta^\prime_--p_-\xi
\label{1.270}
\end{equation}
This is identical in form to the relativistic condition \ref{1.205}. Proceeding then in the same 
way as in the relativistic theory, with the non-relativistic quantities $\beta^\prime$ and 
$\omega$ in the role of the relativistic quantities $\beta$ and $\omega$, placing also 
the vector $\overline{\xi}^\sharp$ in the role of the vector $\xi^\sharp$ of the relativistic theory, and with 
the non-relativistic equations \ref{1.126} and \ref{1.132} in place of the relativistic equations 
\ref{1.57} and \ref{1.60}, we deduce, assuming that the past solution is irrotational, 
the equations:
\begin{equation}
\omega_{||+}=0
\label{1.271}
\end{equation}
\begin{equation}
\omega_{\bot+}=\frac{\theta_+}{u_{\bot+}}d_{||}\triangle s
\label{1.272}
\end{equation}
which are identical in form to the corresponding relativistic equations, but with the 1-form $\omega_{\bot+}$ being here defined by: 
\begin{equation}
\omega_{\bot+}=\left.i_{\overline{\xi}^\sharp}\omega_{+}\right|_{T{\cal K}}
\label{1.273}
\end{equation}
(compare with \ref{1.225}). Then $\overline{\omega}_+$, the classical vorticity 2-form of the future solution at the surface ${\cal K}_t$, is given by:
\begin{equation}
\overline{\omega}_+=\frac{\theta_+}{u_{\bot+}}\overline{\xi}\wedge\overline{d}_{||}\triangle s
\label{1.274}
\end{equation}
where $\overline{d}_{||}$ denotes exterior derivative on ${\cal K}_t$. The corresponding vorticity 
vector is: 
\begin{equation}
\varpi_+=\frac{\theta_+}{u_{\bot+}}\s^*\overline{d}_{||}\triangle s^\sharp
\label{1.275}
\end{equation}
where by $\overline{d}_{||}\triangle s^\sharp$ we denote the corresponding tangent vector to the surface ${\cal K}_t$ and by $\s^*$ the rotation in $T_x{\cal K}_t$ by a right angle, counterclockwise 
relative to the normal $\overline{\xi}^\sharp$.

From \ref{1.268} in analogy with the relativistic theory, 
or directly from \ref{1.213} taking the limit $c\rightarrow\infty$, we deduce in the non-relativistic theory the relation:
\begin{equation}
\triangle s=\frac{1}{12\theta_-}\left(\frac{\partial^2 V}{\partial p^2}\right)_-(\triangle p)^3+O((\triangle p)^4)
\label{1.276}
\end{equation}
Moreover, we can relate the coefficient of the $(\triangle p)^3$ term in \ref{1.269} to the function 
$(dH^\prime/dh^\prime)_s$
induced on ${\cal K}$ by the past solution, $H^\prime$ being the function defined by \ref{1.a18}. 
We find: 
\begin{equation}
-\frac{V^3}{\eta^4}\left(\frac{dH^\prime}{dh^\prime}\right)_s=
\left(\frac{d^2 V}{dp^2}\right)_s
\label{1.277}
\end{equation} 

The {\em determinism conditions} in the non-relativistic theory are entirely analogous to those of the 
relativistic theory. We thus have, as in the relativistic theory, the 1st and 2nd determinism 
conditions, expressed by \ref{1.236} and \ref{1.237} respectively. However, the reciprocal 
acoustical metric being given in the non-relativistic theory by \ref{1.121}, we have, by \ref{1.258} 
and \ref{1.199}, 
\begin{equation}
(h^{-1}_{\pm})^{\mu\nu}\xi_\mu\xi_\nu=-\frac{1}{\eta_\pm^2}u_{\bot\pm}^2+1
\label{1.278}
\end{equation}
Therefore the determinism conditions become:
\begin{equation}
u_{\bot-}>\eta_-,  \ \ \ u_{\bot+}<\eta_+
\label{1.279}
\end{equation}
or, substituting from \ref{1.261},
\begin{equation}
f>\frac{\eta_-}{V_-}, \ \ \ f<\frac{\eta_+}{V_+}
\label{1.280}
\end{equation}
We conclude that when the jump conditions are satisfied the determinism conditions reduce to:
\begin{equation}
\frac{\eta_-}{V_-}<\frac{\eta_+}{V_+}
\label{1.281}
\end{equation}
Defining in the non-relativistic framework, taking the limit of \ref{1.242}, the quantity:
\begin{equation}
w=\frac{V^2}{\eta^2}
\label{1.282}
\end{equation}
the condition \ref{1.281} is equivalent to:
\begin{equation}
\triangle w<0
\label{1.283}
\end{equation}
Since in the non-relativistic theory we have
$$\frac{1}{\eta^2}=\frac{\partial n}{\partial p}=-\frac{1}{V^2}\frac{\partial V}{\partial p}$$
$w$ is simply:
\begin{equation}
w=-\frac{\partial V}{\partial p}
\label{1.284}
\end{equation}
In view of the fact that by \ref{1.276} $\triangle s=O((\triangle p)^3)$ it follows that:
\begin{equation}
\triangle w=-\left(\frac{\partial^2 V}{\partial p^2}\right)_-\triangle p
+O((\triangle p)^2)
\label{1.285}
\end{equation}
Comparing with \ref{1.276} we conclude that, at 
least for suitably small $\triangle p$, condition \ref{1.283} is equivalent to the {\em entropy 
condition}:
\begin{equation}
\triangle s>0
\label{1.286}
\end{equation}
Note also that in view of \ref{1.277} this requirement implies that $\triangle p$ is positive 
or negative according as to whether $(dH^\prime/dh^\prime)_{s-}$ is negative or positive. These results are completely analogous to those of the relativistic theory. 

\vspace{5mm}

\section{The Shock Development Problem}

In the framework of special relativity, we choose a time function $t$ in Minkowski spacetime, equal to 
the coordinate $x^0$ of a rectangular coordinate system and we denote by $\Sigma_t$ an arbitrary level set of the function $t$. So the $\Sigma_t$ are parallel spacelike hyperplanes with respect to the Minkowski metric $g$.  Initial data for the equations of motion \ref{1.27} is given on a domain 
$\overline{{\cal D}}_0$ in the hyperplane $\Sigma_0$, which may be the whole of $\Sigma_0$. 
It consists of 
the specification of the triplet $(u,p,s)$ on $\overline{{\cal D}}_0$ with $(p,s)$ and the rectangular components of $u$ as smooth functions of the rectangular coordinates. To any given initial data set there corresponds a unique 
{\em maximal classical development} of the equations of motion \ref{1.27}. 
In the context of fluid mechanics in special 
relativity, the notion of maximal classical development of given initial data is the following. Given an 
initial data set, the local existence theorem asserts the existence of a 
{\em classical development} of this data, 
namely of a domain ${\cal D}$ in Minkowski spacetime whose past boundary is $\overline{{\cal D}}_0$, 
and of a solution $(u,p,s)$ defined in ${\cal D}$, with $(p,s)$ and the rectangular components of $u$ 
being smooth functions of the rectangular coordinates, taking the given data at the past boundary, such that the 
following condition holds: If we consider any point $p\in {\cal D}$ and any curve issuing at $p$ with 
the property that its tangent vector at any point $q$ belongs to the interior or the boundary of the 
past component of the sound cone at $q$, then the curve terminates at a point of 
$\overline{{\cal D}}_0$. The local uniqueness theorem asserts that two classical developments of the same 
initial data set, with domains ${\cal D}_1$ and ${\cal D}_2$, coincide in 
${\cal D}_1\bigcap {\cal D}_2$. It follows that the union of all classical developments of a given initial 
data set is itself a classical development, the unique maximal classical development of the initial data set. 

In the monograph [Ch-S], we considered initial data on the whole of $\Sigma_0$ for the 
general equations of motion \ref{1.27}, which outside a sphere coincide with the data corresponding 
to a constant state. That is, on $\Sigma_0$ outside the sphere, $p$ and $s$ are constant, 
$(p,s)=(p_0,s_0)$, while $u$ coincides with the future-directed unit normal to $\Sigma_0$. Then, 
given any $\varepsilon\in(0,1]$, we showed that under a suitable restriction on the size of the departure of the initial data from those of the constant state, the solution can be controlled for 
a time interval the ratio of which to the radius of the initial sphere is at least $1/\varepsilon\eta_0$, where $\eta_0=\eta(p_0,s_0)$ is the sound speed in the surrounding constant state. At the end of this time interval an annular region has formed, bounded by concentric spheres 
with the ratio of the radius of the inner sphere to the radius of the outer sphere being at most 
$\varepsilon$, where the flow is irrotational, the constant state holding outside the outer sphere. 
Taking then as new initial hyperplane the hyperplane corresponding the end point of the said time interval, we studied the maximal classical development of the restriction of the data induced on 
this hyperplane to the exterior of the inner sphere. Thus, the global aspects of the work [Ch-S] 
pertain to the irrotational case. The study of the general equations of motion for the initial time  interval was in order to show that the appearance of a singular boundary of the maximal development 
is implied under very general conditions on the initial data. 

To state the results of [Ch-S] in regard to the boundary of the maximal classical development we 
must introduce the {\em acoustical structure}. To be specific, with $\Sigma_0$ the 0-level set of the 
time function $t$ as redefined to correspond to the end point of the initial time interval, 
let $S_{0,0}$ be the outer sphere mentioned above outside which the constant state holds. We then 
define the function $u$ in $\Sigma_0\setminus O$, $O$ being the center of $S_{0,0}$, as minus the 
signed Euclidean distance from $S_{0,0}$. We denote by $S_{0,u}$ the level sets of $u$ on 
$\Sigma_0\setminus O$. So these are spheres centered at $O$, the inner sphere mentioned above being 
$S_{0,u_0}$ for some $u_0>0$. What is considered in [Ch-S] is the restriction of the 
initial data induced on $\Sigma_0$ to $\Sigma_0^{u_0}$, the exterior of $S_{0,u_0}$, the 
non-trivial data corresponding to the range $[0,u_0]$ of $u$, while the negative range of $u$ 
corresponds to the surrounding constant state. What the monograph [Ch-S] studies is the maximal 
classical development of this data on $\Sigma_0^{u_0}$. The function $u$ is extended to the domain of the maximal 
classical development by the condition that its level sets are outgoing null hypersurfaces relative 
to the acoustical metric $h$. We call $u$, extended in this way, {\em acoustical function} and we 
denote by $C_u$ an arbitrary level set of $u$. So each $C_u$ is an outgoing {\em characteristic 
hypersurface}, and is generated by {\em bi-characteristics}, namely null geodesics of $h$. In 
[Ch-S] we defined $L$ to be the tangent vectorfield to these null geodesic generators, parametrized 
not affinely but by $t$. We also defined the surfaces $S_{t,u}$ to be $\Sigma_t\bigcap C_u$. 
Moreover, we defined the vectorfield $T$ to be tangential to the $\Sigma_t$ and so that the flow 
generated by $T$ on each $\Sigma_t$ is the normal, relative to the induced on $\Sigma_t$ acoustical 
metric $\overline{h}$, flow of the foliation of $\Sigma_t$ by the surfaces $S_{t,u}$. So $T$ is the 
tangent vectorfield to the normal curves parametrized by $u$. 

What plays a central role in [Ch-S] is the density of the foliation of spacetime by the $C_u$. 
One measure of this density is the {\em inverse spatial density}, that is, the inverse density of 
the foliation of each spatial hyperplane $\Sigma_t$ by the surafces $S_{t,u}$. This is simply the 
magnitude $\kappa$ of the vectorfield $T$ with respect to $\overline{h}$. Another is the inverse 
{\em temporal density}, namely the function $\mu$, given in arbitrary coordinates by:
\begin{equation}
\frac{1}{\mu}=-(h^{-1})^{\mu\nu}\partial_\mu t\partial_\nu u
\label{1.287}
\end{equation}
The two measures are related by:
\begin{equation}
\mu=\alpha\kappa
\label{1.288}
\end{equation}
where $\alpha$ is the inverse density, with respect to the acoustical metric $h$, of the foliation 
of spacetime by the hyperplanes $\Sigma_t$:
\begin{equation}
\frac{1}{\alpha}=\sqrt{-(h^{-1})^{\mu\nu}\partial_\mu t\partial_\nu t}
\label{1.289}
\end{equation}
(the corresponding inverse density with respect to the Minkowski metric $g$ being of course 1). 
As was shown in [Ch-S], the function $\alpha$ is bounded above and below by positive constants, 
consequently $\mu$ and $\kappa$ are equivalent measures of the density of the foliation of spacetime 
by the $C_u$. 

The domain of the maximal classical development being a domain in Minkowski spacetime, which by a 
choice of rectangular coordinates is identified with $\mathbb{R}^4$, inherits the subset topology 
and the standard differential structure induced by the rectangular coordinates $x^\alpha \ : \ 
\alpha=0,1,2,3$. Choosing an acoustical function $u$ as above, we introduce 
{\em acoustical coordinates}  
$(t,u,\vartheta)$, $\vartheta\in S^2$, the coordinate lines corresponding to a given value of $u$ 
and to constant values of $\vartheta$ being the generators of the $C_u$. Given a choice of $u$  
this property does not determine the $S^2$ valued function $\vartheta$ uniquely. To do so we must 
specify $\vartheta$ on a hypersurface transversal to the $C_u$. For example, we may stipulate that 
on $\Sigma_0^{u_0}$, $\vartheta$ is constant along the normal curves to the foliation of 
$\Sigma_0^{u_0}$ by the surfaces $S_{0,u}$. Otherwise $\vartheta$ is determined only up to transformations of the form:
\begin{equation}
\vartheta\mapsto f(u,\vartheta)
\label{1.290}
\end{equation}
where $f$ is a smooth function. In terms of acoustical coordinates the vectorfields $L$ and $T$ are 
given by:
\begin{equation}
L=\frac{\partial}{\partial t}, \ \ \ T=\frac{\partial}{\partial u}-\Xi
\label{1.291}
\end{equation}
where $\Xi$ is a vectorfield which is tangential to the surfaces $S_{t,u}$. This can be thought of 
as a vectorfield on $S^2$ depending on $(t,u)$. The specification of $\vartheta$ on a hypersurface 
transversal to the $C_u$ corresponds to the specification of $\Xi$ on such a hypersurface. 
In particular, the example of specifying $\vartheta$ to be constant along the normal curves to the 
foliation of $\Sigma_0^{u_0}$ by the $S_{0,u}$ corresponds to setting $\Xi=0$ on $\Sigma_0^{u_0}$. 
If local coordinates $(\vartheta ^A \ : \ A=1,2)$ are chosen on $S^2$, 
\begin{equation}
\left(\frac{\partial}{\partial\vartheta^A} \ : \ A=1,2\right)
\label{1.292}
\end{equation}
being a coordinate frame field for the $S_{t,u}$, we can express:
\begin{equation}
\Xi=\Xi^A\frac{\partial}{\partial\vartheta^A}
\label{1.293}
\end{equation}
Moreover, with 
\begin{equation}
\sh_{AB}=\sh\left(\frac{\partial}{\partial\vartheta^A},\frac{\partial}{\partial\vartheta^B}\right)
\label{1.294}
\end{equation}
being the components of the metric $\sh$ induced by $h$ on the surfaces $S_{t,u}$ in the frame field 
\ref{1.292}, the acoustical metric $h$ takes in the acoustical coordinate system 
$(t,u, \vartheta^A:A=1,2)$ the form:
\begin{equation}
h=-\mu(dt\otimes du+du\otimes dt)+\kappa^2 du\otimes du+
\sh_{AB}(d\vartheta^A+\Xi^A du)\otimes(d\vartheta^B+\Xi^B du)
\label{1.295}
\end{equation}

As was shown in [Ch-S], the rectangular coordinates $(x^\alpha)$ are 
smooth functions of the acoustical coordinates $(t,u,\vartheta)$ and the Jacobian of the 
transformation is, up to a multiplicative factor which is bounded above and below by positive 
constants, the inverse temporal density function $\mu$. The acoustical coordinates induce another 
differential structure on the same underlying topological manifold. However since $\mu$ is positive 
in the interior of the domain of the maximal classical development, 
the two differential structures coincide in this interior. 
The main theorem of [Ch-S] asserts that relative to the differential structure induced by the 
acoustical coordinates the maximal classical development extends smoothly to the boundary of its 
domain. This boundary contains however a singular part ${\cal B}$ where the function $\mu$ 
vanishes. So the rectangular coordinates themselves extend smoothly to the boundary but the Jacobian 
vanishes on the singular part of the boundary. The mapping from acoustical to rectangular coordinates 
has a continuous but not differentiable inverse at ${\cal B}$. As a result, the two differential 
structures no longer coincide when the singular boundary ${\cal B}$ is included. With respect to 
the standard differential structure the solution is continuous but not differentiable at ${\cal B}$, 
the derivative, in terms of rectangular coordinates, $\hat{T}^\mu\hat{T}^\nu\partial_\mu\beta_\nu$ 
blowing up as we approach ${\cal B}$. Here $\hat{T}=\kappa^{-1}T$ is the vectorfield of unit 
magnitude, with respect to $\overline{h}$, corresponding to $T$, and $\beta$ is the 1-form \ref{1.53}. 
With respect to the standard differential structure, the acoustical metric $h$ is everywhere in the 
closure of the domain of the maximal classical development non-degenerate and continuous, but it is 
not differentiable at ${\cal B}$, while with respect to the differential structure induced by the 
acoustical coordinates $h$ is everywhere smooth, but it is degenerate at ${\cal B}$, reducing along 
${\cal B}$ to (see \ref{1.295}):
\begin{equation}
h=\sh_{AB}(d\vartheta^A+\Xi^A du)\otimes(d\vartheta^B+\Xi^B du)
\label{1.296}
\end{equation}

According to the results of [Ch-S] the boundary of the domain of the classical maximal development 
consists of a regular part $\underline{{\cal C}}$ together with the singular part ${\cal B}$. 
Each of these is not necessarily connected. Each connected component of $\underline{{\cal C}}$ is an 
incoming acoustically null hypersurface which is smooth relative to both differential structures and 
has a singular past boundary which coincides with the past boundary of an associated component of 
${\cal B}$. The union of these singular past boundaries we denote by $\partial_-{\cal B}$. 

The singular part ${\cal B}$ can be described in acoustical coordinates as the graph of a smooth 
function:
\begin{equation}
{\cal B}=\{(t_*(u,\vartheta),u,\vartheta) \ : \ (u,\vartheta)\in{\cal D}\}
\label{1.297}
\end{equation}
where ${\cal D}$ is a domain in $[0,u_0]\times S^2$, each component of $\overline{{\cal B}}$ being a graph 
over the corresponding component of $\overline{{\cal D}}$. In fact, each component of 
$\overline{{\cal B}}$ is a smooth hypersurface relative to the standard as well as the acoustical differential structure. Using $(u,\vartheta)$ as coordinates on ${\cal B}$, 
the induced metric $h_*$ on ${\cal B}$ is expressed as (see \ref{1.296}):
\begin{equation}
h_*=(\sh_*)_{AB}(d\vartheta^A+\Xi_*^A du)\otimes(d\vartheta^B+\Xi_*^B du)
\label{1.298}
\end{equation}
where:
\begin{equation}
(\sh_*)_{AB}(u,\vartheta)=\sh_{AB}(t_*(u,\vartheta),u,\vartheta)
\label{1.299}
\end{equation}
are the components of a positive-definite metric on $S^2$ depending on $u$ and 
\begin{equation}
\Xi_*^A(u,\vartheta)=\Xi^A(t_*(u,\vartheta),u,\vartheta)
\label{1.300}
\end{equation}
are the components of a vectorfield on $S^2$ depending on $u$. The functions $(\sh_*)_{AB}$, 
$\Xi_*^A$ are all smooth. The induced metric $h_*$ on ${\cal B}$ being a smooth degenerate metric, 
the singular boundary ${\cal B}$ is {\em from the point of view of its intrinsic geometry similar to 
a regular null hypersurface in a regular spacetime}. The vectorfield $V$ given in $(u,\vartheta)$ coordinates on ${\cal B}$ by:
\begin{equation}
V=\frac{\partial}{\partial u}-\Xi_*
\label{1.301}
\end{equation}
is a nowhere vanishing null vectorfield on ${\cal B}$. Its integral curves, which we call {\em invariant curves} have zero arc length with respect to the acoustical metric and are parametrized 
by $u$. The singular boundary ${\cal B}$ is ruled by these curves. The past boundary 
$\partial_-{\cal B}$ of ${\cal B}$ is the set of past end points of the invariant curves. 

According to the results of [Ch-S], we have:
\begin{equation}
\frac{\partial\mu}{\partial t}<0 \ \ \mbox{: on $\overline{{\cal B}}$}
\label{1.302}
\end{equation}
Consequently {\em relative to the acoustical differential structure} the hypersurface ${\cal B}$ is transversal to the $C_u$. We may then adapt the coordinates $(u,\vartheta)$ on ${\cal B}$ so that 
the coordinate lines corresponding to constant values of $\vartheta$ coincide with the invariant curves. This choice is equivalent to the condition:
\begin{equation}
\Xi_*=0
\label{1.303}
\end{equation}
We call the corresponding acoustical coordinates {\em canonical}. In such coordinates the induced metric (see \ref{1.298}) takes the form, simply:
\begin{equation}
h_*=(\sh_*)_{AB}d\vartheta^A\otimes d\vartheta^B
\label{1.304}
\end{equation}
Moreover, by virtue of condition \ref{1.303}, we have, for $t\leq t_*(u,\vartheta)$:
\begin{equation}
\Xi^A(t,u,\vartheta)=-\int_t^{t_*(u,\vartheta)}\left(\frac{\partial\Xi^A}{\partial t}\right)(t^\prime,u,\vartheta)dt^\prime
\label{1.305}
\end{equation}
Since also:
\begin{equation}
\mu(t,u,\vartheta)=-\int_t^{t_*(u,\vartheta)}\left(\frac{\partial\mu}{\partial t}\right)(t^\prime,u,\vartheta)dt^\prime
\label{1.306}
\end{equation}
it follows that the functions:
\begin{equation}
\hat{\Xi}^A=\mu^{-1}\Xi^A
\label{1.307}
\end{equation}
which a priori are defined only for $t<t_*(u,\vartheta)$, are actually given by:
\begin{equation}
\hat{\Xi}^A(t,u,\vartheta)=\frac{\mbox{mean value on $[t,t_*(u,\vartheta)]$ of 
$(\partial\Xi^A/\partial t)(\cdot,u,\vartheta)$}}
{\mbox{mean value on $[t,t_*(u,\vartheta)]$ of $(\partial\mu/\partial t)(\cdot,u,\vartheta)$}}
\label{1.308}
\end{equation}
hence extend smoothly to $t=t_*(u,\vartheta)$, that is, to ${\cal B}$. These functions are the 
components of the vectorfield:
\begin{equation}
\hat{\Xi}=\hat{\Xi}^A\frac{\partial}{\partial\vartheta^A}=\mu^{-1}\Xi
\label{1.309}
\end{equation}
which is tangential to the $S_{t,u}$.

According to the results of [Ch-S], the singular boundary ${\cal B}$ is {\em from the extrinsic point 
of view}, that is, from the point of view of how it is embedded in the acoustical spacetime {\em a 
spacelike hypersurface}. This means that at each point $q\in{\cal B}$, the past null geodesic conoid 
of $q$ relative to the acoustical metric does not intersect ${\cal B}$. This is connected with the fact that in canonical acoustical 
coordinates
\begin{equation}
\frac{\partial\mu}{\partial u}<0 \ \ \mbox{: on ${\cal B}$}
\label{1.310}
\end{equation}
The connection is through the fact that: 
\begin{equation}
\mu h^{-1}(d\mu,d\mu)\rightarrow -2\frac{\partial\mu}{\partial t}\frac{\partial\mu}{\partial u}
\label{1.311}
\end{equation}
in canonical acoustical coordinates, as we approach ${\cal B}$. On the other hand, in contrast with 
\ref{1.302} we have:
\begin{equation}
\frac{\partial\mu}{\partial u}=0 \ \ \mbox{: at $\partial_-{\cal B}$}
\label{1.312}
\end{equation}
which implies that ${\cal B}$ {\em becomes acoustically null along its past boundary} 
$\partial_-{\cal B}$. Since $\overline{{\cal B}}$ is where $\mu$ vanishes, we have:
\begin{equation}
\mu(t_*(u,\vartheta),u,\vartheta)=0 \ : \ \mbox{for all $(u,\vartheta)\in\overline{{\cal D}}$}
\label{1.313}
\end{equation}
Differentiating this equation implicitly with respect to $u$, which, since we are using canonical 
acoustical coordinates, means differentiating along the invariant curves, we obtain:
\begin{equation}
\left(\frac{\partial t_*}{\partial u}\right)(u,\vartheta)=
-\left(\frac{\partial\mu/\partial u}{\partial\mu/\partial t}\right)(t_*(u,\vartheta),u,\vartheta) 
\ \ : \forall (u,\vartheta)\in\overline{{\cal D}}
\label{1.314}
\end{equation}
Then in view of \ref{1.302}, the properties \ref{1.310} and \ref{1.312} are equivalent to:
\begin{equation}
\left(\frac{\partial t_*}{\partial u}\right)(u,\vartheta)<0 \ \ : \ \forall (u,\vartheta)\in{\cal D}
\label{1.315}
\end{equation}
and
\begin{equation}
\left(\frac{\partial t_*}{\partial u}\right)(u,\vartheta)=0 \ \ : \ 
\forall (u,\vartheta)\in\partial_-{\cal D}
\label{1.316}
\end{equation}
respectively, $\partial_-{\cal B}$ being a graph over $\partial_-{\cal D}$. Moreover, in regard to 
$\partial_-{\cal B}$, in addition to \ref{1.312} we have (canonical acoustical coordinates):
\begin{equation}
\frac{\partial^2\mu}{\partial u^2}>0 \ \ \mbox{: at $\partial_-{\cal B}$}
\label{1.317}
\end{equation}
Also, differentiating again \ref{1.314} with respect to $u$ and evaluating the result on 
$\partial_-{\cal D}$ we obtain:
\begin{equation}
\left(\frac{\partial^2 t_*}{\partial u^2}\right)(u,\vartheta)=
-\left(\frac{\partial^2\mu/\partial u^2}{\partial\mu/\partial t}\right)(t_*(u,\vartheta),u,\vartheta) 
\ \ : \forall (u,\vartheta)\in\partial_-{\cal D}
\label{1.318}
\end{equation}
hence in view of \ref{1.302} the property \ref{1.317} is equivalent to:
\begin{equation}
\left(\frac{\partial^2 t_*}{\partial u^2}\right)(u,\vartheta)>0 \ \ : \ 
\forall (u,\vartheta)\in\partial_-{\cal D}
\label{1.319}
\end{equation}
Thus, $t$ as a function on each invariant curve reaches a minimum at its past end point on  $\partial_-{\cal B}$. Moreover, each component of $\partial_-{\cal B}$ is an acoustically 
spacelike surface which is smooth relative to both differential structures. 

We remark that the tangent field $V$ to the invariant curves (see \ref{1.301}) expressed as a vectorfield in spacetime along ${\cal B}$ in canonical acoustical coordinates takes the form:
\begin{equation}
V=\frac{\partial t_*}{\partial u}\frac{\partial}{\partial t} +\frac{\partial}{\partial u}
\label{1.320}
\end{equation}
The fact that relative to the standard differential structure the vectorfield $\partial/\partial u$ vanishes along ${\cal B}$ implies that at any point 
$q=(t_*(u,\vartheta),u,\vartheta)\in{\cal B}$ the invariant curve through $q$ is tangent to the 
generator of $C_u$ corresponding to the value of $\vartheta$, relative to the standard differential structure. 

The last part of the work [Ch-S] contains a result which is important for reaching an understanding.  
This is the {\em trichotomy theorem}. According to this theorem, for each point 
$q\in\overline{\cal B}$, the intersection of the past null geodesic conoid of $q$ (relative to the 
acoustical metric) with any $\Sigma_t$ in the past of $q$ splits into three parts, the parts 
corresponding to the outgoing and to the incoming sets of null geodesics ending at $q$ being 
embedded disks with a common boundary, an embedded circle, which corresponds to the set of the 
remaining null geodesics ending at $q$. {\em All outgoing null geodesics ending at $q$ have the same 
tangent vector at $q$}. This vector is then an invariant null vector associated to the singular point 
$q$ and defines the tangent line at $q$ to the invariant curve through $q$. This is the reason why 
the freedom in the choice of acoustical function does not matter in the end. While each acoustical 
function is defined by its restriction to $\Sigma_0$, in [Ch-S] the change is investigated, 
from one acoustical function $u$ to another $u^\prime$, such that the null geodesic generators 
of $C_u$, the level sets of $u$, and of $C^\prime_{u^\prime}$, the level sets of $u^\prime$, 
which have future end points on the singular boundary of the domain of the maximal classical 
development, both belong to the class of outgoing null geodesics ending at these points. This change 
can be viewed as being generated by a coordinate transformation on $\overline{{\cal B}}$:
\begin{equation}
(u,\vartheta)\mapsto (u^\prime,\vartheta^\prime)=(f(u,\vartheta),g(\vartheta))
\label{1.321}
\end{equation}
preserving the invariant curves and their orientation, so that $\partial f/\partial u>0$ and 
$g$ is orientation-preserving. It was then shown that the foliations corresponding to the two 
different families of outgoing acoustically null hypersurfaces, the $C_u$ and the 
$C^\prime_{u^\prime}$, degenerate in precisely the same way on the same singular boundary 
$\overline{{\cal B}}$. 

Now, as stated in the epilogue of [Ch-S], the notion of maximal classical development, on which the 
whole work [Ch-S] is based, while perfectly reasonable from the mathematical point of view, and 
also the correct notion from the physical point of view up to 
$\underline{{\cal C}}\bigcup\partial_-{\cal B}$, {\em it is not the correct notion from the 
physical point of view up to ${\cal B}$}. Let us consider a given component of $\partial_-{\cal B}$ 
and the corresponding components of $\underline{{\cal C}}$ and ${\cal B}$. {\em The actual physical 
problem associated to the given component of $\partial_-{\cal B}$ is the following:}

\vspace{5mm}

Find a hypersurface ${\cal K}$ in Minkowski spacetime, lying in the past of the component of 
${\cal B}$ and with the same past boundary, namely the component of $\partial_-{\cal B}$, and the 
same tangent hyperplane at each point along this boundary, and a solution $(u,p,s)$ of the equations 
of motion in the domain in Minkowski spacetime bounded in the past by $\underline{\cal C}$ and 
${\cal K}$, such that on $\underline{\cal C}$ $(u,p,s)$ coincides with the data induced by the maximal 
classical development, while across ${\cal K}$ there are jumps relative to the data induced by the 
maximal classical development, jumps which satisfy the jump conditions, discussed in the previous 
section, following from the integral form of the conservation laws. Moreover, the hypersurface 
${\cal K}$ is to be space-like relative to the acoustical metric induced by the maximal classical 
solution, which holds in the past of ${\cal K}$, and time-like relative to the new solution, which 
holds in the future of ${\cal K}$. The last conditions are the determinism conditions, also discussed 
in the previous section. 

\vspace{5mm}

We call this the {\em shock development problem}. The shock is the hypersurface of discontinuity 
${\cal K}$. The singular acoustically space-like surface $\partial_-{\cal B}$ is the 
origin of the shock, hence may be thought of as its cause. Since what is sought in the problem 
is the hypersurface ${\cal K}$ as well as the new solution in the domain bounded in the past by 
$\underline{{\cal C}}$ and ${\cal K}$, and data is given along $\underline{{\cal C}}$ making 
$\underline{{\cal C}}$ an incoming acoustically null hypersurface, while nonlinear boundary 
conditions are imposed at ${\cal K}$ and ${\cal K}$ is to be an acoustically time-like hypersurface 
relative to the new solution, the shock development problem is a characteristic initial - free 
boundary value problem for the 1st order quasilinear hyperbolic system \ref{1.27} with non-linear 
boundary conditions at the free boundary ${\cal K}$. Moreover, and {\em this is the main difficulty of 
the problem}, the initial conditions on the characteristic hypersurface $\underline{{\cal C}}$ become 
singular at its past boundary $\partial_-{\cal B}$, while the boundary conditions at ${\cal K}$ are 
stated in terms of the prior solution which becomes singular at $\partial_-{\cal B}$ the last being 
at the same time the past boundary of ${\cal K}$. 

In the non-relativistic framework the shock development problem is stated in the same way as above. 
In the role of Minkowski spacetime we have of course the Galilei spacetime. The $\Sigma_t$ are then 
the hyperplanes of absolute simultaneity, and the function $t$, unique up to an additive constant 
is absolute time. In the role of the relativistic equations of motion \ref{1.27} we have the 
equations of motion of the non-relativistic theory, equations \ref{1.151}, the unknowns here being 
$(v,p,s)$. In the monograph [Ch-Mi] we considered these equations in the irrotational case, where 
the equations reduce to the nonlinear wave equation \ref{1.142}. We considered initial data for this 
equation outside a sphere such that the data outside a larger concentric sphere corresponds to that of a constant state. Since in the non-relativistic framework the specific enthalpy $h^\prime$ is 
only defined up to an additive constant (see \ref{1.b4}), we may set $h^\prime=0$ in the constant 
state. Hence $\phi^\prime=0$ in the constant state and the initial data for \ref{1.142} simply 
vanishes outside the larger sphere. What we studied in [Ch-Mi] is again the maximal classical 
development of the data in the exterior of the inner sphere in $\Sigma_0$. An acoustical function 
$u$ and acoustical coordinates $(t,u,\vartheta)$ are defined in exactly the same way as in the relativistic theory. So are the vectorfields $L$ and $T$, the inverse spatial density $\kappa$ 
and the inverse temporal density $\mu$ of the foliation of spacetime by the $C_u$. The last two 
are again related by \ref{1.288} with $\alpha$ defined by \ref{1.289}, but since $h^{-1}$ is given 
in the non-relativistic theory by \ref{1.121}, we now have:
$$-(h^{-1})^{\mu\nu}\partial_\mu t\partial_\nu t=\eta^{-2}$$
so in the non-relativistic framework it holds simply:
\begin{equation}
\alpha=\eta
\label{1.322}
\end{equation}
It is similarly shown that $\eta$ is bounded above and below by positive constants. The treatment 
in [Ch-Mi] is self-contained, making no reference to the relativistic theory. It is considerably 
simpler than the corresponding relativistic study, due chiefly to the following fact: the function 
$t$, which plays an important role in the study being one of the acoustical coordinates, has level 
sets $\Sigma_t$ whose induced acoustical metric $\overline{h}$ coincides in the non-relativistic 
theory with the Euclidean metric, as is evident from \ref{1.120}. The results of [Ch-Mi] in regard 
to the boundary of the maximal classical development are identical in form to those of [Ch-S] 
discussed above. In the non-relativistic framework, the jump conditions referred to in the 
statement of the shock development problem are of course those of the non-relativistic theory, 
which have also been discussed in the previous section. 

\vspace{5mm} 

\section{The Restricted Shock Development Problem}

As discussed in the previous section, we are considering the case in which the maximal classical 
solution is irrotational. In this case, from Section 1.4, in the relativistic theory the new solution 
which we are seeking in the shock development problem has a vorticity 2-form which satisfies on 
${\cal K}$ the boundary conditions \ref{1.224} and \ref{1.232}, and is given along ${\cal K}$ by 
\ref{1.227}. By \ref{1.220}, the 1-form $\beta_{||+}$ induced on ${\cal K}$ by the 
new solution is given by: 
\begin{equation}
\beta_{||+}=d_{||}\phi_-
\label{1.323}
\end{equation}
Although the 2-form $\omega_{||+}$ induced on ${\cal K}$ by the new solution vanishes, the new 
solution is not irrotational due to the non-vanishing of the component $\omega_{\bot+}$, given 
by \ref{1.232}. As equation \ref{1.232} shows, $\omega_{\bot+}$ is generated by the variation 
of $\triangle s$ along ${\cal K}$. The portion of the fluid which does not cross the shock 
hypersurface ${\cal K}$, that is, the set of flow lines in the new solution which have past 
end points on $\underline{{\cal C}}$ rather than on ${\cal K}$, remains isentropic and therefore 
by \ref{1.67} irrotational. Now, $\triangle s$ is given in terms of $\triangle p$ for small 
values of the jumps by \ref{1.213}. At $\partial_-{\cal B}$ the jumps vanish, as there continuity 
holds. Since according to \ref{1.213} $\triangle s=O((\triangle p)^3)$, in a neighborhood of 
$\partial_-{\cal B}$, where, as noted in the previous section, the main difficulty of the shock development problem lies, $\triangle s$ may be considered to be a higher order quantity. We thus 
retain the main difficulty of the problem  if we neglect it altogether. This leads to the 
{\em restricted shock development problem}, the subject of the present monograph. In this problem  
also the new solution is irrotational, that is the 1-form $\beta$ corresponding to the new solution 
is again of the form \ref{1.68} for a new function $\phi$, and $s$ has for the new solution the 
same constant value as in the prior maximal classical solution. By the condition \ref{1.323} and 
continuity at $\partial_-{\cal B}$, $\phi$ is continuous across ${\cal K}$:
\begin{equation}
\triangle\phi=0
\label{1.324}
\end{equation}
Conversely, since the new function $\phi$ is to be continuously differentiable up to the boundary, 
\ref{1.324} implies \ref{1.323}. Since according to the actual shock development problem 
$\beta$ and $s$ must be continuous across $\underline{{\cal C}}$, in the restricted shock development 
problem $\phi$ as well as $d\phi$ are to be continuous across $\underline{{\cal C}}$. The only 
remaining jump condition is then \ref{1.201}, which is the 1st of \ref{1.195} and is associated to 
equation \ref{1.15}, which is equation \ref{1.d2}. By \ref{1.53} and \ref{1.72}:
\begin{equation}
n_\pm u_{\bot\pm}=-\xi_\mu G_\pm\beta_\pm^\mu
\label{1.325}
\end{equation}
so the jump condition \ref{1.201} reads:
\begin{equation}
\xi_\mu\triangle(G\beta^\mu)=0
\label{1.326}
\end{equation}
Now, the condition \ref{1.323}, which as we have seen is equivalent to \ref{1.324}, is: 
\begin{equation}
\triangle\beta_{||}=0
\label{1.327}
\end{equation}
This means:
\begin{equation}
X^\mu\triangle\beta_\mu=0 \ : \ \forall X\in T{\cal K}
\label{1.328}
\end{equation}
Recalling the definition of the covector $\xi$, this means that at each point on ${\cal K}$ the 
covectors $\xi$ and $\triangle\beta$ are colinear. Therefore the condition \ref{1.326} can be expressed in 
the form:
\begin{equation}
(\triangle\beta_\mu)\triangle(G\beta^\mu)=0
\label{1.329}
\end{equation}
To summarize, in the restricted shock development problem we have two jump conditions across 
the shock hypersurface ${\cal K}$: the {\em linear jump condition} \ref{1.328}, and {\em the nonlinear 
jump condition} \ref{1.329}. These equations have been written in terms of components in an 
arbitrary frame field, in particular components relative to a rectangular coordinate system. 
The restricted shock development problem is formulated in a manner very similar to that of the actual 
shock development problem:

\vspace{5mm}

Find a hypersurface ${\cal K}$ in Minkowski spacetime, lying in the past of the component of 
${\cal B}$ and with the same past boundary, namely the component of $\partial_-{\cal B}$, and the 
same tangent hyperplane at each point along this boundary, and a solution $\phi$ of the nonlinear 
wave equation \ref{1.70} in the domain in Minkowski spacetime bounded in the past by 
$\underline{\cal C}$ and ${\cal K}$, such that on $\underline{\cal C}$, $\beta=d\phi$ coincides with 
the data induced by the maximal classical development, while across ${\cal K}$ there are jumps 
relative to the data induced by the maximal classical development, jumps which satisfy the jump 
conditions \ref{1.328} and \ref{1.329}. Moreover, the hypersurface ${\cal K}$ is to be space-like 
relative to the acoustical metric induced by the maximal classical solution, which holds in the past 
of ${\cal K}$, and time-like relative to the new solution, which holds in the future of ${\cal K}$. 

\vspace{5mm}

In the framework of the non-relativistic theory we have the condition: 
\begin{equation}
\triangle\beta^\prime_{||}=0
\label{1.330}
\end{equation}
which is equivalent to: 
\begin{equation}
\triangle\phi^\prime=0
\label{1.331}
\end{equation}
These are identical in form to \ref{1.327} and \ref{1.324} respectively. We also have the jump 
condition 1st of \ref{1.256} which is associated to equation \ref{1.111}. This condition reads: 
\begin{equation}
\xi_\mu\triangle I^\mu=0
\label{1.332}
\end{equation}
where $I^\mu$ are the components of the non-relativistic particle current \ref{1.110}  
\begin{equation}
I^0=n, \ \ \ I^i=-n\beta^\prime_i
\label{1.333}
\end{equation}
in Galilean coordinates. Condition \ref{1.330} is the {\em linear jump condition}:
\begin{equation}
X^\mu\triangle\beta_\mu=0 \ : \ \forall X\in T{\cal K}
\label{1.334}
\end{equation}
which is identical in form to \ref{1.328}, allowing us to express condition \ref{1.332} 
as the {\em nonlinear jump condition}: 
\begin{equation}
(\triangle\beta^\prime_\mu)\triangle I^\mu=0
\label{1.335}
\end{equation}
(Galilean coordinates). The same can also be obtained from the relativistic condition \ref{1.329} by 
taking the limit $c\rightarrow\infty$. 
The restricted shock development problem is then stated in the same way as 
in the relativistic theory, with Galilei spacetime in the role of Minkowski spacetime, 
$\phi^\prime$ and $\beta^\prime$ in the role of $\phi$ and $\beta$, the non-relativistic 
nonlinear wave equation \ref{1.169} in the role of the relativistic nonlinear wave equation 
\ref{1.70}, and the jump conditions \ref{1.334} and \ref{1.335} in the role of the jump 
conditions \ref{1.328} and \ref{1.329}. 

\pagebreak
\chapter{The Geometric Construction}
\section{A General Construction in Lorentzian Geometry}

We begin with a general construction in Lorentzian geometry. To attain a full understanding of the 
mathematical structure we consider the general case of $n$ spatial dimensions. So, we have an 
$n+1$ dimensional oriented manifold ${\cal M}$ equipped with a Lorentzian metric $h$, which is also 
time-oriented. In the realization 
of the structure in fluid mechanics, $h$ shall be the acoustical metric. Let $u$ and $\ub$ be 
functions defined in a domain in ${\cal M}$, increasing toward the future, such that the $C_u$, 
the level sets of $u$, are outgoing 
null hypersurfaces and the $\Cb_{\ub}$, the level sets of $\ub$ are incoming null hypersurfaces, 
relative to $h$. The functions $u$ and $\ub$ are solutions of the 
{\em eikonal equation}:
\begin{equation}
h^{-1}(du,du)=0, \ \ \ h^{-1}(d\ub,d\ub)=0
\label{2.1}
\end{equation}
Each $\Cb_{\ub}$ intersects each $C_u$ transversally. We denote by $S_{\ub,u}$ the surfaces of 
intersection:
\begin{equation}
S_{\ub,u}=\Cb_{\ub}\bigcap C_u
\label{2.2}
\end{equation}
These are space-like surfaces relative to $h$, that is, the metric $\sh$ induced by $h$ on the 
$S_{\ub,u}$ is positive-definite. We shall take the $S_{\ub,u}$ to be diffeomorphic to 
$S^{n-1}$. Alternatively, we could take the $S_{\ub,u}$ to be diffeomorphic to $\mathbb{R}^{n-1}$, 
however, in view of the fact that $\mathbb{R}^{n-1}$ is diffeomorphic to $S^{n-1}$ minus a point 
and the construction in connection with the shock development problem can be localized 
by reason of the domain of dependence property, there is no loss of generality 
in restricting ourselves to the first alternative. We have:
\begin{equation}
C_u=\bigcup_{\ub\in\underline{I}_u}S_{\ub,u}, \ \ \ \Cb_{\ub}=\bigcup_{u\in I_{\ub}}S_{{\ub},u}
\label{2.3}
\end{equation}
where each $\underline{I}_u$  and each $I_{\ub}$ is an interval in $\mathbb{R}$. We thus have 
a domain $D\subset\mathbb{R}^2$, such that:
\begin{equation}
\underline{I}_u=\{\ub \ : \ (\ub,u)\in D\}, \ \ \ I_{\ub}=\{u \ : \ (\ub,u)\in D\}
\label{2.4}
\end{equation}
and we are considering the domain ${\cal M}_D$ in ${\cal M}$ given by:
\begin{equation}
{\cal M}_D=\bigcup_{(\ub,u)\in D}S_{\ub,u}
\label{2.5}
\end{equation}
Thus, each $C_u$ is foliated by the 1-parameter family of sections 
$\{S_{\ub,u} \ : \ \ub\in\underline{I}_u\}$, 
each $\Cb_{\ub}$ is foliated by the 1-parameter family of sections 
$\{S_{\ub,u} \ : \ u\in I_{\ub}\}$, and the spacetime domain 
${\cal M}_D$ is foliated by the 2-parameter family $\{S_{\ub,u} \ : \ (\ub,u)\in D\}$. 

Each $C_u$ being a null hypersurface, it is generated by null geodesic segments. Let $L$ be the 
tangent field to these null geodesics, parametrized not affinely but by the function $\ub$. Then 
$L$ is a future-directed null vectorfield defined in ${\cal M}_D$ such that
\begin{equation}
Lu=0, \ \ \ L\ub=1, \ \ \ \mbox{$L$ is  $h$-orthogonal to the $S_{\ub,u}$}
\label{2.6}
\end{equation}
and
\begin{equation}
D_L L=\omega L
\label{2.7}
\end{equation}
for some function $\omega$. Here $D$ is the covariant derivative operator associated to 
the metric $h$. Similarly, each $\Cb_{\ub}$ being a null hypersurface, it is generated by null 
geodesic segments. Let $\Lb$ be the tangent field to these null geodesics, parametrized not affinely 
but by the function $u$. Then $\Lb$ is a future-directed null vectorfield defined in ${\cal M}_D$ 
such that
\begin{equation}
\Lb\ub=0, \ \ \ \Lb u=1, \ \ \ \mbox{$\Lb$ is $h$-orthogonal to the $S_{\ub,u}$}
\label{2.8}
\end{equation}
and
\begin{equation}
D_{\Lb}\Lb=\underline{\omega}\Lb
\label{2.9}
\end{equation}
for some function $\underline{\omega}$. 

The vectorfields $L$, $\Lb$ being both future-directed null, there is a positive function $a$ such 
that 
\begin{equation}
h(L,\Lb)=-2a
\label{2.10}
\end{equation}
Consider the vectorfields:
\begin{equation}
L^\prime=a^{-1}L, \ \ \ \Lb^\prime=a^{-1}\Lb
\label{2.11}
\end{equation}
We have:
\begin{equation}
h(L^\prime,\Lb)=h(\Lb^\prime,L)=-2
\label{2.12}
\end{equation}
We shall presently show that:
\begin{equation}
L^\prime=-2(du)^\sharp, \ \ \ \Lb^\prime=-2(d\ub)^\sharp
\label{2.13}
\end{equation}
This means that for any vector $V$ we have:
\begin{equation}
h(L^\prime,V)=-2Vu, \ \ \ h(\Lb^\prime,V)=-2V\ub
\label{2.14}
\end{equation}
In fact, it suffices to verify the cases $V=L$, $V=\Lb$, and $V\in TS_{\ub,u}$. In the case $V=L$ 
we have, by \ref{2.6}, \ref{2.11}, \ref{2.12},
$$h(L^\prime,L)=0 \ \ \mbox{and} \ \ -2Lu=0,$$
$$h(\Lb^\prime,L)=-2 \ \ \mbox{and} \ \ -2L\ub=-2$$
In the case $V=\Lb$ we have, by \ref{2.8}, \ref{2.11}, \ref{2.12},
$$h(L^\prime,\Lb)=-2 \ \ \mbox{and} \ \ -2\Lb u=-2,$$
$$h(\Lb^\prime,\Lb)=0, \ \ \mbox{and} \ \ -2\Lb\ub=0$$
Finally, in the case $V\in TS_{\ub,u}$ we have, by \ref{2.6}, \ref{2.8}, and \ref{2.11}:
$$h(L^\prime,V)=0=-2Vu, \ \ \ h(\Lb^\prime,V)=0=-2V\ub$$
By a well known proposition in Lorentzian geometry, \ref{2.13} imply that the vectorfields 
$L^\prime$, $\Lb^\prime$ are geodesic:
\begin{equation}
D_{L^\prime}L^\prime=0, \ \ \ D_{\Lb^\prime}\Lb^\prime=0
\label{2.15}
\end{equation}
Thus, the integral curves of $L^\prime$ are the null geodesic generators of the $C_u$ affinely 
parametrized, and the integral curves of $\Lb^\prime$ are the null geodesic generators of the 
$\Cb_{\ub}$ affinely parametrized. Comparing \ref{2.15} with \ref{2.7}, \ref{2.9}, and with 
\ref{2.11} we see that the functions $\omega$, $\underline{\omega}$ are given by:
\begin{equation}
\omega=a^{-1}La, \ \ \ \underline{\omega}=a^{-1}\Lb a
\label{2.16}
\end{equation}

Let us denote by $Z$ the commutator:
\begin{equation}
Z=[\Lb,L]
\label{2.17}
\end{equation}
From \ref{2.6}, \ref{2.8} we see that $Zu=Z\ub=0$. Therefore the vectorfield $Z$ is tangential to 
the surfaces $S_{\ub,u}$. Let $\Phi_{\sb}$ be the flow generated by $L$ and $\Phib_s$ the flow 
generated by $\Lb$. So if $p\in S_{\ub,u}$, then $\sb\mapsto\Phi_{\sb}(p)$ is the generator of $C_u$ 
through $p$, and $s\mapsto\Phib_s(p)$ is the generator of $\Cb_{\ub}$ through $p$. 
Thus $\Phi_{\sb}$ preserves the $C_u$, while $\Phib_s$ preserves the 
$\Cb_{\ub}$. Moreover, $\left.\Phi_{\sb}\right|_{S_{\ub,u}}$ is a diffeomorphism of  
$S_{\ub,u}$ onto $S_{\ub+\sb,u}$, while $\left.\Phib_s\right|_{S_{\ub,u}}$ is a diffeomorphism of 
$S_{\ub,u}$ onto $S_{\ub,u+s}$. The non-vanishing of $Z$ implies that the two flows do not commute, 
the mapping $\Phib_s^{-1}\circ\Phi_{\sb}^{-1}\circ\Phib_s\circ\Phi_{\sb}$ preserving the $S_{\ub,u}$ 
but differing from the identity. 

For each $(\ub,u)\in D$ we can define a diffeomorphism:
\begin{equation}
\Phi_{(\ub,u)} : S_{\ub,u}\rightarrow S^{n-1}
\label{2.18}
\end{equation} 
such that these diffeomorphisms are related according to:
\begin{equation}
\Phi_{(\ub+\sb,u)}=\Phi_{(\ub,u)}\circ\Phi_{\sb}^{-1}
\label{2.19}
\end{equation}
Thus the family of diffeomorphisms $\{\Phi_{(\ub,u)} \ : \ (\ub,u)\in D\}$ respects the flow of $L$:
$$\Phi_{(\ub+\sb,u)}(\Phi_{\sb}(p))=\Phi_{(\ub,u)}(p) \ \ : \ \forall p\in S_{\ub,u}$$ 
That is, each generator of each $C_u$ is mapped to a single point of $S^{n-1}$. We then have 
a diffeomorphism:
\begin{equation}
{\cal M}_D \rightarrow D\times S^{n-1} \ \ \mbox{by:} \ \ p\in S_{\ub,u}\mapsto 
((\ub,u)\in D,\vartheta=\Phi_{(\ub,u)}(p)\in S^{n-1})
\label{2.20}
\end{equation} 
In terms of this representation of ${\cal M}_D$, the vectorfields $L$ and $\Lb$ are expressed by:
\begin{equation}
L=\frac{\partial}{\partial\ub}, \ \ \ \Lb=\frac{\partial}{\partial u}+b
\label{2.21}
\end{equation}
where $b$ is a vectorfield on $S^{n-1}$ depending on $(\ub,u)$. Also $Z$ is represented as a 
vectorfield on $S^{n-1}$ depending on $(\ub,u)$ and by \ref{2.17} we have:
\begin{equation}
\frac{\partial b}{\partial\ub}=-Z
\label{2.22}
\end{equation}

Similarly, for each $(\ub,u)\in D$ we can define a diffeomorphism:
\begin{equation}
\Phib_{(\ub,u)} : S_{\ub,u}\rightarrow S^{n-1}
\label{2.23}
\end{equation}
such that these diffeomorphisms are related according to:
\begin{equation}
\Phib_{(\ub,u+s)}=\Phib_{(\ub,u)}\circ\Phib_s^{-1}
\label{2.24}
\end{equation}
Thus the family of diffeomorphisms $\{\Phib_{(\ub,u)} \ : \ (\ub,u)\in D\}$ respects the flow of $\Lb$: 
$$\Phib_{(\ub,u+s)}(\Phib_s(p))=\Phib_{(\ub,u)}(p) \ \ : \ \forall p\in S_{\ub,u}$$
That is, each generator of each $\Cb_{\ub}$ is mapped to a single point of $S^{n-1}$. We then have 
a diffeomorphism:
\begin{equation}
{\cal M}_D \rightarrow D\times S^{n-1} \ \ \mbox{by:} \ \ p\in S_{\ub,u}\mapsto 
((\ub,u)\in D,\vartheta=\Phib_{(\ub,u)}(p)\in S^{n-1})
\label{2.25}
\end{equation}
In terms of this representation of ${\cal M}_D$, the vectorfields $L$ and $\Lb$ are expressed by:
\begin{equation}
L=\frac{\partial}{\partial\ub}-b, \ \ \ \Lb=\frac{\partial}{\partial u}
\label{2.26}
\end{equation}
where again $b$ is a vectorfield on $S^{n-1}$ depending on $(\ub,u)$. Also $Z$ is represented as a 
vectorfield on $S^{n-1}$ depending on $(\ub,u)$ and by \ref{2.17} we have:
\begin{equation}
\frac{\partial b}{\partial u}=-Z
\label{2.27}
\end{equation}

Each of the representations \ref{2.20}, \ref{2.25} is appropriate in a particular context. However, 
as we shall see, the representation which is generally appropriate in the framework of a free boundary 
problem, such as the shock development problem, is one adapted to the flow $F_\tau$ of the 
future-directed timelike vectorfield
\begin{equation}
T=L+\Lb
\label{2.28}
\end{equation}
We have:
\begin{equation}
Tu=T\ub=1, \ \  \mbox{and $T$ is  $h$-orthogonal to the $S_{\ub,u}$}
\label{2.29}
\end{equation}
Then $\left.F_\tau\right|_{S_{\ub,u}}$ is a diffeomorphism of $S_{\ub,u}$ onto $S_{\ub+\tau,u+\tau}$. 
Let us then define, for each $(\ub,u)\in D$ a diffeomorphism:
\begin{equation}
F_{(\ub,u)} : S_{\ub,u}\rightarrow S^{n-1}
\label{2.30}
\end{equation}
such that these diffeomorphisms are related according to:
\begin{equation}
F_{(\ub+\tau,u+\tau)}=F_{(\ub,u)}\circ F_\tau^{-1}
\label{2.31}
\end{equation}
Thus the family of diffeomorphism $\{F_{(\ub,u)} \ : \ (\ub,u)\in D\}$ respects the flow of $T$: 
$$F_{(\ub+\tau,u+\tau)}(F_\tau(p))=F_{(\ub,u)}(p) \ \ : \ \forall p\in S_{\ub,u}$$
That is, each integral curve of $T$ is mapped to a single point of $S^{n-1}$. We then have 
a diffeomorphism:
\begin{equation}
{\cal M}_D \rightarrow D\times S^{n-1} \ \ \mbox{by:} \ \ p\in S_{\ub,u}\mapsto 
((\ub,u)\in D,\vartheta=F_{(\ub,u)}(p)\in S^{n-1})
\label{2.32}
\end{equation}
In terms of this representation of ${\cal M}_D$ the vectorfields $L$, $\Lb$ are expressed by:
\begin{equation}
L=\frac{\partial}{\partial\ub}-b, \ \ \ \Lb=\frac{\partial}{\partial u}+b
\label{2.33}
\end{equation}
and the vectorfield $T$ simply by:
\begin{equation}
T=\frac{\partial}{\partial\ub}+\frac{\partial}{\partial u}
\label{2.34}
\end{equation}
In \ref{2.33} $b$ is again a vectorfield on $S^{n-1}$ depending on $(\ub,u)$. 
Likewise $Z$ is represented as a vectorfield on $S^{n-1}$ depending on $(\ub,u)$ 
and by \ref{2.17} we have:
\begin{equation}
\frac{\partial b}{\partial\ub}+\frac{\partial b}{\partial u}=-Z
\label{2.35}
\end{equation}

Let $(\vartheta^A : A=1,...,n-1)$ be local coordinates on $S^{n-1}$, defined in a domain 
$U\subset S^{n-1}$. Then we have a system of local coordinates $(\ub,u, \vartheta^A:A=1,...,n-1)$ in the inverse image of $D\times U$ in ${\cal M}_D$ by the diffeomorphism \ref{2.32}. Let us denote 
by $\Omega_A : A=1,...,n-1$ the vectorfields, in terms of these coordinates,
\begin{equation}
\Omega_A=\frac{\partial}{\partial\vartheta^A} \ : \ A=1,...,n-1
\label{2.36}
\end{equation}
which are tangential to the surfaces $S_{\ub,u}$. At each $S_{\ub,u}$ these constitute a local frame 
field for $S_{\ub,u}$. The components
\begin{equation}
\sh_{AB}=h(\Omega_A,\Omega_B)
\label{2.37}
\end{equation}
of the induced metric relative to this frame field constitute a positive-definite matrix 
at each point. 

Also, $(L,\Lb,\Omega^A:A=1,...,n-1)$ is a local frame field for ${\cal M}_D$ in the inverse image of 
$D\times U$ by \ref{2.32}. By \ref{2.10} and \ref{2.37}, the reciprocal metric $h^{-1}$ is 
expressed in terms of this frame field as:
\begin{equation}
h^{-1}=-\frac{1}{2a}(L\otimes\Lb+\Lb\otimes L)+(\sh^{-1})^{AB}\Omega_A\otimes\Omega_B
\label{2.38}
\end{equation}
According to \ref{2.6} and \ref{2.8} we have:
\begin{equation}
d\ub\cdot L=1, \ \ \ d\ub\cdot\Lb=0, \ \ \ d\ub\cdot\Omega_A=0 \ : \ A=1,...,n-1
\label{2.39}
\end{equation}
and:
\begin{equation}
du\cdot L=0, \ \ \ du\cdot\Lb=1, \ \ \ du\cdot\Omega_A=0 \ : \ A=1,...,n-1
\label{2.40}
\end{equation}
Let $\theta^A : A=1,...,n-1$ be the 1-forms:
\begin{equation}
\theta^A=d\vartheta^A+b^A(d\ub-du)
\label{2.41}
\end{equation}
where $b^A$ are the components of the $S_{\ub,u}$ -tangential vectorfield $b$ in the frame field 
$(\Omega_A:A=1,...,n-1)$:
\begin{equation}
b=b^A\Omega_A
\label{2.42}
\end{equation}
Then by \ref{2.33} and \ref{2.36} we have:
\begin{equation}
\theta^A\cdot L=0, \ \ \ \theta^A\cdot\Lb=0, \ \ \ \theta^A\cdot\Omega_B=\delta^A_B \ : \ B=1,...,n-1
\label{2.43}
\end{equation}
From \ref{2.39}, \ref{2.41}, and \ref{2.43} we see that $(d\ub,du,\theta^A:A=1,...,n-1)$ is the 
co-frame field dual to $(L,\Lb,\Omega_A:A=1,...,n-1)$. The metric $h$ is then expressed in terms 
of this co-frame field as:
\begin{equation}
h=-2a(d\ub\otimes du+du\otimes d\ub)+\sh_{AB}\theta^A\otimes\theta^B
\label{2.44}
\end{equation}
This expression shows in particular that the function $a$ is the {\em inverse density of the double 
null foliation}, that is the foliation of ${\cal M}_D$ by the $C_u$, $\Cb_{\ub}$. 

\vspace{5mm}

\section{The Characteristic System}

The mathematical structure we have discussed in the previous section does not 
depend on any underlying structure on the differentiable manifold ${\cal M}$. From this point on 
however we bring in the underlying structure of ${\cal M}$, namely that of $\mathbb{M}^n$, 
Minkowski spacetime of $n$ spatial dimensions in framework of the relativistic theory, or that of 
$\mathbb{G}^n$, Galilei spacetime of $n$ spatial dimensions in the non-relativistic framework. 
In both cases there is an underlying affine structure so there are preferred systems of linear 
coordinates, more precisely rectangular coordinates $(x^\mu:\mu=0,1,...,n)$ in the relativistic 
case, Galilean coordinates $(x^0=t, x^i:i=1,...,n)$ in the non-relativistic case. Moreover we are 
now addressing fluid mechanics, so $h$ is the acoustical metric. We assume that the time orientation relative to $h$, mentioned above, agrees with the time orientation 
of the underlying $\mathbb{M}^n$, or $\mathbb{G}^n$. That is, at each spacetime point $p$, the 
positive cone of $h$ at $p$, the future sound cone at $p$, is contained in the future light cone at $p$ in the relativistic case, 
the future of $\Sigma_t$, the hyperplane of absolute simultaneity through $p$, 
in the non-relativistic case. In the present context we call $\ub$ and $u$ {\em acoustical functions} 
and the coordinates $(\ub,u,\vartheta)$ as in \ref{2.32} {\em acoustical coordinates}.  

As we have seen in the previous chapter, the components $h_{\mu\nu}$ of the acoustical metric in 
the preferred coordinates, being given in the relativistic case by \ref{1.d11} (with $H$ a given 
smooth function of the relativistic enthalpy $h$ which is itself given by \ref{1.54}) and in the 
non-relativistic case by \ref{1.d23} (with $H^\prime$ a given smooth function of the non-relativistic 
enthalpy $h^\prime$ which is itself given by the 1st of \ref{1.d20}), are given smooth functions 
of the components $\beta_\mu$, or of the components $\beta^\prime_\mu$, of the 1-form $\beta$, or of 
the 1-form $\beta^\prime$, in the preferred coordinates, accordingly. {\em In the present 
section we shall take the point of view that components $\beta_\mu$, or $\beta^\prime_\mu$, therefore 
also the components $h_{\mu\nu}$, in the preferred coordinates, are known smooth functions of the acoustical coordinates 
$(\ub,u,\vartheta)$, and what we are seeking are the conditions that must be satisfied by the 
preferred coordinates $x^\mu$ (rectangular or Galilean according as to whether we are treating the relativistic or 
the non-relativistic theory) as functions of the acoustical coordinates 
$(\ub,u,\vartheta)$.}

Consider then, for a fixed value of $(\ub,u)\in D$, the surface $S_{\ub,u}$. This surface is described 
as an embedding 
\begin{equation}
S^{n-1}\rightarrow\mathbb{M}^n, \ \mbox{or} \ S^{n-1}\rightarrow\mathbb{G}^n
\label{2.45}
\end{equation}
accordingly, expressed by giving the preferred coordinates $x^\mu$ as functions on $S^{n-1}$ at 
the given value $(\ub,u)\in D$:
\begin{equation}
x^\mu=x^\mu(\ub,u,\vartheta) \ : \ \vartheta\in S^{n-1} 
\label{2.46}
\end{equation}
In a domain $U\in S^{n-1}$ these are functions of the local coordinates $(\vartheta^A:A=1,...,n-1)$ 
defined in $U$. Then
\begin{equation}
\Omega_A^\mu=\Omega_A x^\mu=\frac{\partial x^\mu}{\partial\vartheta^A} \ \ : \ A=1,...,n-1
\label{2.47}
\end{equation}
are the components of the vectorfields $\Omega_A$ along $S_{\ub,u}$ in the preferred coordinates. 

Consider now at a given point $p\in S_{\ub,u}$ the $h$- orthogonal complement of $T_p S_{\ub,u}$. 
This is a 2-dimensional plane $P_p$ consisting of all the vectors $M$ at $p$ which satisfy:
\begin{equation}
h(M,\Omega_A)=0 \ : \ \forall A=1,...,n-1
\label{2.48}
\end{equation}
or, in terms of the components in the preferred coordinates:
\begin{equation}
h_{\mu\nu}M^\mu\Omega_A^\nu=0 \ : \ \forall A=1,....,n-1
\label{2.49}
\end{equation}
Now, the $n-1$ dimensional plane $T_p S_{\ub,u}$ being space-like relative to $h$, 
the 2-dimensional plane $P_p$ is time-like 
relative to $h$. Therefore there are exactly two future-directed null rays in $P_p$. These are 
the vectors $M\in P_p$ which satisfy:
\begin{equation}
h(M,M)=0, \ \ M^0>0
\label{2.50}
\end{equation}
or, in terms of the components in the preferred coordinates:
\begin{equation}
h_{\mu\nu}M^\mu M^\nu=0, \ \ M^0>0
\label{2.51}
\end{equation}
Since given any element of one of the two rays the whole ray is obtained by multiplying the given 
element by all the positive real numbers, we may represent the ray by the element which 
satisfies:
\begin{equation}
M^0=1
\label{2.52}
\end{equation}
Thus, the $n$-tuplet $\overline{M}=(M^1,...,M^n)$ is subject to the conditions \ref{2.49} which read:
\begin{equation}
h_{00}\Omega_A^0+h_{0i}\Omega_A^i+(h_{0i}\Omega_A^0+h_{ij}\Omega_A^j)M^i=0 \ \ : \ A=1,...,n-1
\label{2.53}
\end{equation}
as well as to the condition, 1st of \ref{2.51}, which reads:
\begin{equation}
h_{00}+2h_{0i}M^i+h_{ij}M^i M^j=0
\label{2.54}
\end{equation}
Equations \ref{2.53} constitute a system of $n-1$ inhomogeneous linear equations for the $n$ 
unknowns $M^i:i=1,...,n$. This defines a line, not through the origin, in $\mathbb{R}^n$:
\begin{equation}
\overline{M}=\alpha s+\beta
\label{2.55}
\end{equation}
$\beta\in\mathbb{R}^n$ being a particular solution and $\alpha s$, for a particular 
$\alpha\in\mathbb{R}^n$ and any $s\in\mathbb{R}$ being the general solution of the corresponding 
homogeneous system. Substituting then \ref{2.55} in equation \ref{2.54} we obtain a quadratic 
equation for $s$ whose discriminant is: 
\begin{equation}
\mbox{det}\sh, \ \ \mbox{where $\sh_{AB}$ is the matrix:} \ \ 
\sh_{AB}=h(\Omega_A,\Omega_B)=h_{\mu\nu}\Omega_A^\mu\Omega_B^\nu
\label{2.56}
\end{equation}
This is a homogeneous polynomial in the $\Omega_A^\mu=\partial x^\mu/\partial\vartheta^A$ of 
degree $2(n-1)$ and represents the square of the area element of $T_p S_{\ub,u}$. We then obtain 
two solutions $\overline{M}_\pm$ representing the two future-directed null rays in $P_p$, of the 
form:
\begin{equation}
\overline{M}_\pm=\frac{P\pm R\sqrt{\mbox{det}\sh}}{Q}
\label{2.57}
\end{equation}
where $P$ and $Q$ are respectively $\mathbb{R}^n$ and $\mathbb{R}$ valued homogeneous polynomials in the $\Omega_A^\mu$ of degree $2(n-1)$  and $R$ is an $\mathbb{R}^n$ valued homogeneous polynomial 
in the $\Omega_A^\mu$ of degree $n-1$. In the case $n=2$ in particular, $P^i$ is the homogeneous 
quadratic polynomial
$$P^i=\Omega^0\Omega^i-(h^{-1})^{0i}h(\Omega,\Omega),$$
$Q$ is the homogeneous quadratic polynomial
$$Q=-(h^{-1})^{00}\overline{h}^{-1}(\overline{\theta},\overline{\theta})$$
where $\overline{h}$, with components $h_{ij}$, is the restriction of $h$ to the $\Sigma_t$, and 
$\overline{\theta}$, with components $\theta_i$, is the restriction of the 1-form $\theta$ with components $\theta_\mu=h_{\mu\nu}\Omega^\nu$ to the $\Sigma_t$. Also, $R^i$ is the homogeneous 
linear polynomial
$$R^i=[ij]\eta^{-1}\theta_j$$
where $[ij]$ is the fully antisymmetric 2-dimensional symbol. We denote, in general, 
\begin{equation}
N^\mu=M_+^{\mu}, \ \ \ \Nb^\mu=M_-^{\mu}, \ \ \ \mbox{so} \ \ N^0=\Nb^0=1
\label{2.58}
\end{equation}
These are the components of $N$ and $\Nb$:
\begin{equation}
N=N^\mu\frac{\partial}{\partial x^\mu}, \ \ \ \Nb=\Nb^\mu\frac{\partial}{\partial x^\mu}
\label{2.59}
\end{equation}
the future-directed outward and inward null normals to $S_{\ub,u}$ at $p$ 
relative to $h$ respectively, satisfying the normalization condition $Nx^0=\Nb x^0=1$. Note 
that in the relativistic theory the vectors $N$ and $\Nb$ themselves have the physical dimensions 
of inverse length and their components are dimensionless, while in the non-relativistic theory 
the vectors $N$ and $\Nb$ themselves have the physical dimensions of inverse time, the 
time components are of course dimensionless, but the spatial components $N^i$ and 
$\Nb^i$ have the physical dimensions of length/time, that is of velocity. 

Now the vectorfields $L$ and $\Lb$, given in the acoustical coordinates by \ref{2.33} (see also 
\ref{2.36}, \ref{2.42}), also define future-directed outward and inward null normals to the 
$S_{\ub,u}$ relative to $h$ respectively. It follows that there are positive functions $\rho$ and 
$\rhob$ such that:
\begin{equation}
L=\rho N, \ \ \ \Lb=\rhob\Nb
\label{2.60}
\end{equation}
In view of the normalization condition $N^0=\Nb^0=1$, the functions $\rho$ and $\rhob$ are actually 
{\em defined} by the time components of equations \ref{2.60}:
\begin{equation}
\rho=L^0, \ \ \ \rhob=\Lb^0
\label{2.61}
\end{equation}
Applying then \ref{2.60} to the $x^\mu$ as functions of the acoustical coordinates we obtain, in 
view of \ref{2.33}, \ref{2.36}, \ref{2.42}, the system: 
\begin{eqnarray}
&&\frac{\partial x^i}{\partial\ub}-b^A\frac{\partial x^i}{\partial\vartheta^A}
=N^i\left(\frac{\partial t}{\partial\ub}-b^A\frac{\partial t}{\partial\vartheta^A}\right) 
\ \ : \ i=1,...,n
\nonumber\\
&&\frac{\partial x^i}{\partial u}+b^A\frac{\partial x^i}{\partial\vartheta^A}
=\Nb^i\left(\frac{\partial t}{\partial u}+b^A\frac{\partial t}{\partial\vartheta^A}\right) 
\ \ : \ i=1,...,n
\label{2.62}
\end{eqnarray}
Here we denote $x^0=t$ in the relativistic as well as in the non-relativistic theory, as in the former 
units where the speed of light in vacuum is equal to 1 have been chosen. This is a system of 
$2n$ equations for the $2n$ unknown functions of the acoustical coordinates, namely the function 
$t$, the $n$ functions $x^i:i=1,...,n$, and the $n-1$ functions $b^A:A=1,...,n-1$. By the discussion 
above, the $N^i:i=1,...,n$ and $\Nb^i:i=1,...,n$ are homogeneous degree 0 algebraic functions of 
the partial derivatives $\partial t/\partial\vartheta^A, \partial x^i/\partial\vartheta^A : A=1,...,n-1; \ i=1,...,n$. Thus \ref{2.62} constitutes a fully nonlinear 1st order system. We 
call this the {\em characteristic system}. From this system we shall derive in the next chapter 
the {\em acoustical structure equations} through which we shall be able to control the derivatives 
of the $x^\mu$ with respect to the $(\ub,u,\vartheta)$ of all orders.

Now, in the representation \ref{2.32} $C_u$ is the manifold:
\begin{equation}
C_u=\{(\ub,u,\vartheta) \ : \ \ub\in\Ib_u, \ \vartheta\in S^{n-1}\}
\label{2.63}
\end{equation}
and $\Cb_{\ub}$ is the manifold:
\begin{equation}
\Cb_{\ub}=\{(\ub,u,\vartheta) \ : \ u\in I_{\ub}, \ \vartheta\in S^{n-1}\}
\label{2.64}
\end{equation}
Then $(L,\Omega^A:A=1,...,n-1)$ is a frame field for $C_u$ and $(\Lb,\Omega^A:A=1,...,n-1)$ is a 
frame field for $\Cb_{\ub}$. Let $m$ be the induced metric on $C_u$ and $\mb$ the induced metric 
on $\Cb_{\ub}$. We then have:
\begin{eqnarray}
&&m(L,L)=h_{\mu\nu}L^\mu L^\nu=\rho^2 h_{\mu\nu}N^\mu N^\nu=0,\nonumber\\
&&m(L,\Omega_A)=h_{\mu\nu}L^\mu\Omega_A^\nu=\rho h_{\mu\nu}N^\mu\Omega_A^\nu=0
\label{2.65}
\end{eqnarray}
by the 1st of \ref{2.62}, which by \ref{2.33}, \ref{2.36}, \ref{2.42}, \ref{2.60}, \ref{2.61} is simply $$L^\mu=\rho N^\mu,$$
and the conditions defining $N^\mu$, while
\begin{equation}
m(\Omega_A,\Omega_B)=h_{\mu\nu}\Omega_A^\mu\Omega_B^\nu=\sh_{AB}
\label{2.66}
\end{equation}
is positive-definite. Therefore the induced metric $m$ on $C_u$ is degenerate with $L$ the null tangential 
vectorfield, that is $C_u$ is a null hypersurface relative to $h$ and its generators are the integral curves of $L$. 
Similarly, we have:
\begin{eqnarray}
&&\mb(\Lb,\Lb)=h_{\mu\nu}\Lb^\mu \Lb^\nu=\rhob^2 h_{\mu\nu}\Nb^\mu\Nb^\nu=0,\nonumber\\ 
&&\mb(\Lb,\Omega_A)=h_{\mu\nu}\Lb^\mu\Omega_A^\nu=\rhob h_{\mu\nu}\Nb^\mu\Omega_A^\nu=0
\label{2.67}
\end{eqnarray}
by the 2nd of \ref{2.62}, which by \ref{2.33}, \ref{2.36}, \ref{2.42}, \ref{2.60}, \ref{2.61} 
is simply $$\Lb^\mu=\rhob\Nb^\mu,$$
and the conditions defining $\Nb^\mu$, while
\begin{equation}
\mb(\Omega_A,\Omega_B)=h_{\mu\nu}\Omega_A^\mu\Omega_B^\nu=\sh_{AB}
\label{2.68}
\end{equation}
is positive-definite. Therefore the induced metric $\mb$ on $\Cb_{\ub}$ is degenerate with $\Lb$ the 
null tangential vectorfield, that is $\Cb_{\ub}$ is a null hypersurface relative to $h$ and its 
generators are the integral curves of $\Lb$. Thus, the 1st subsystem of the characteristic system 
\ref{2.62} expresses the fact that the $C_u$ are null hypersurfaces relative to $h$, while the 
2nd subsystem of the characteristic system \ref{2.62} expresses the fact that also the $\Cb_{\ub}$ are 
null hypersurfaces relative to $h$. 

In regard to the meaning of the quantities $\rho$ and $\rhob$ we note that since
$$Lt=\rho \ \ \mbox{while} \ \ L\ub=1$$
along a generator $\Gamma$ of a $C_u$, integral curve of $L$, we have:
\begin{equation}
\rho=\left(\frac{dt}{d\ub}\right)_{\Gamma}
\label{2.69}
\end{equation}
Thus {\em $\rho$ is the inverse temporal density of the foliation of spacetime by the $\Cb_{\ub}$ 
as measured along the generators of the $C_u$}. Similarly, since 
$$\Lb t=\rhob \ \ \mbox{while} \ \ \Lb u=1$$
along a generator $\Gammab$ of a $\Cb_{\ub}$, integral curve of $\Lb$, we have:
\begin{equation}
\rhob=\left(\frac{dt}{du}\right)_{\Gammab}
\label{2.70}
\end{equation}
Thus {\em $\rhob$ is the inverse temporal density of the foliation of spacetime by the $C_u$ 
as measured along the generators of a $\Cb_{\ub}$}. 

Now, the vectorfields $N$ and $\Nb$ being both null and future-directed there is a positive function 
$c$ such that:
\begin{equation}
h(N,\Nb)=-2c
\label{2.71}
\end{equation}
Comparing with \ref{2.10} and \ref{2.60} we see that the function $a$ is expressed as:
\begin{equation}
a=c\rho\rhob
\label{2.72}
\end{equation}
In view of \ref{2.38} we can express the components of the reciprocal acoustical metric in the form:
\begin{equation}
(h^{-1})^{\mu\nu}=-\frac{1}{2c}(N^\mu\Nb^\nu+\Nb^\mu N^\nu)+(\sh^{-1})^{AB}\Omega_A^\mu\Omega_B^\nu
\label{2.73}
\end{equation}

As we shall see in the next chapter it is the functions 
\begin{equation}
\lambda=c\rhob \ \ \mbox{and} \ \ \lambdab=c\rho
\label{2.74}
\end{equation}
rather than $\rhob$ and $\rho$ which satisfy appropriate propagation equations along the generators 
of the $C_u$ and $\Cb_{\ub}$. 
Recalling the function $\mu$ defined by \ref{1.287}, which plays a leading role in [Ch-S] and [Ch-Mi], 
by \ref{2.11}, \ref{2.13} and \ref{2.72} this function is in the present context given by:
\begin{equation}
\mu=2\lambda
\label{2.75}
\end{equation}
Similarly we can define here in analogy with \ref{1.287} the function $\mub$ by:
\begin{equation}
\frac{1}{\mub}=-(h^{-1})^{\mu\nu}\partial_\mu t\partial_\nu\ub
\label{2.76}
\end{equation}
Then from \ref{2.11}, \ref{2.13} and \ref{2.72} we have:
\begin{equation}
\mub=2\lambdab
\label{2.77}
\end{equation}

Consider finally the Jacobian of the transformation from acoustical to preferred (that is, rectangular 
in the relativistic theory, Galilean in the non-relativistic theory) coordinates:
\begin{equation}
\frac{\partial(x^0,x^1,...,x^n)}{\partial(\ub,u,\vartheta^1,...,\vartheta^{n-1})}=
\left|\begin{array}{cccc}
\partial x^0/\partial\ub&\partial x^1/\partial\ub&\cdot\cdot\cdot 
&\partial x^n/\partial\ub\\
\partial x^0/\partial u&\partial x^1/\partial u&\cdot\cdot\cdot 
&\partial x^n/\partial u\\
\partial x^0/\partial\vartheta^1&\partial x^1/\partial\vartheta^1&\cdot\cdot\cdot 
&\partial x^n/\partial\vartheta^1\\
\cdot &\cdot &\cdot\cdot\cdot &\cdot\\
\cdot &\cdot &\cdot\cdot\cdot &\cdot\\
\partial x^0/\partial\vartheta^{n-1}&\partial  x^1/\partial\vartheta^{n-1}&\cdot\cdot\cdot 
&\partial x^n/\partial\vartheta^{n-1}
\end{array}\right|
\label{2.78}
\end{equation}
The last is equal to:
\begin{equation}
\left|\begin{array}{cccc}
L^0&L^1&\cdot\cdot\cdot&L^n\\
\Lb^0&\Lb^1&\cdot\cdot\cdot&\Lb^n\\
\Omega_1^0&\Omega_1^1&\cdot\cdot\cdot&\Omega_1^n\\
\cdot &\cdot &\cdot\cdot\cdot &\cdot\\
\cdot &\cdot &\cdot\cdot\cdot &\cdot\\
\Omega_{n-1}^0&\Omega_{n-1}^1&\cdot\cdot\cdot&\Omega_{n-1}^n
\end{array}\right|
=\rho\rhob d
\label{2.79}
\end{equation}
where:
\begin{equation}
d=\left|\begin{array}{cccc}
1&N^1&\cdot\cdot\cdot&N^n\\
1&\Nb^1&\cdot\cdot\cdot&\Nb^n\\
\Omega_1^0&\Omega_1^1&\cdot\cdot\cdot&\Omega_1^n\\
\cdot &\cdot &\cdot\cdot\cdot &\cdot\\
\cdot &\cdot &\cdot\cdot\cdot &\cdot\\
\Omega_{n-1}^0&\Omega_{n-1}^1&\cdot\cdot\cdot&\Omega_{n-1}^n
\end{array}\right|
\label{2.80}
\end{equation}
The function $c$ being bounded from below by a positive constant and $\sh_{AB}$ being uniformly 
positive-definite imply that $d$ is bounded from above by a negative constant. Then the Jacobian 
\ref{2.78} being 
\begin{equation}
\frac{\partial(x^0,x^1,...,x^n)}{\partial(\ub,u,\vartheta^1,...,\vartheta^{n-1})}=\rho\rhob d
\label{2.81}
\end{equation}
vanishes where and only where one of $\rho$, $\rhob$ vanishes. 

We finally remark in regard to the differential order of the different quantities 
appearing in the characteristic system \ref{2.62}. We consider the 1st derivatives of the $x^\mu$, 
the components $b^A$,  
and the components $\beta_\mu$ to be of order 0. Then the components $h_{\mu\nu}$ are of order 0 
and so are the $\Omega_A^\mu$, as well as $N^\mu$, $\Nb^\mu$, $c$, $\rho$, $\rhob$, $\lambda$, 
$\lambdab$, $a$, and $\sh_{AB}$. The characteristic system \ref{2.62} itself is a fully nonlinear system 
of order 0 for the 1st derivatives of the $x^\mu$ and for the $b^A$.

\vspace{5mm}

\section{The Wave System}

As we have remarked in the previous section, when considering the characteristic system we think 
in the relativistic framework of the components $\beta_\mu$ of the 1-form $\beta$ in the rectangular 
coordinates, and in the non-relativistic framework of the components $\beta^\prime_\mu$ of the 
1-form $\beta^\prime$ in the Galilean coordinates, as given functions of the acoustical coordinates. 
Consequently the characteristic system must be supplemented by the system of differential conditions 
to which the components $\beta_\mu$ or $\beta^\prime_\mu$ accordingly are subject to as functions 
of the acoustical coordinates. 

In the relativistic theory the 1-form $\beta$ is subject to equation \ref{1.d9} which is valid 
in any system of coordinates, together with the condition $\omega=0$. Now, equation \ref{1.d9} 
is in particular valid in rectangular coordinates in which the equation reads:
\begin{equation}
(h^{-1})^{\mu\nu}\frac{\partial\beta_\nu}{\partial x^\mu}=0
\label{2.82}
\end{equation}
If we think of the {\em rectangular components} $\beta_\nu$ as being functions on the spacetime 
manifold which may be locally expressed in any local coordinate system whatever, we can write, 
in terms of $d\beta_\nu$, the differentials of these functions, 
\begin{equation}
\frac{\partial\beta_\nu}{\partial x^\mu}=d\beta_\nu\cdot\frac{\partial}{\partial x^\mu}
\label{2.83}
\end{equation}
Thinking similarly of the {\em rectangular coordinates} $x^\nu$ as being functions on the spacetime 
manifold which may be locally expressed in any local coordinate system whatever, we have, in terms 
of $dx^\nu$, the differentials of these functions,
\begin{equation}
dx^\nu\cdot\frac{\partial}{\partial x^\lambda}=\delta^\nu_\lambda
\label{2.84}
\end{equation}
Moreover, 
\begin{equation}
(h^{-1})^{\mu\nu}=(h^{-1})^{\mu\lambda}\delta^\nu_\lambda \ \ \mbox{and} \ \ 
h^{-1}=(h^{-1})^{\mu\lambda}\frac{\partial}{\partial x^\mu}\otimes\frac{\partial}{\partial x^\lambda}
\label{2.85}
\end{equation}
In view of \ref{2.83}, \ref{2.84}, \ref{2.85}, equation \ref{2.82} can be expressed in the form:
\begin{equation}
h^{-1}(d\beta_\nu,dx^\nu)=0
\label{2.86}
\end{equation}
where the rectangular components $\beta_\nu$ as well as the rectangular coordinates $x^\nu$ are 
considered as functions on the spacetime manifold which may be locally expressed in any local coordinate system whatever. 

Since the 1-form $\beta$ is expressed in terms of its rectangular components and the differentials 
of the rectangular coordinates by:
\begin{equation}
\beta=\beta_\nu dx^\nu
\label{2.87}
\end{equation}
the condition $\omega=0$, that is $d\beta=0$, reads:
\begin{equation}
d\beta_\nu\wedge dx^\nu=0
\label{2.88}
\end{equation}

Equations \ref{2.86}, \ref{2.88} together constitute the general form of what we call the 
{\em wave system}, in the relativistic theory. 

In the non-relativistic theory the 1-form $\beta^\prime$ is subject to equation \ref{1.d21} which is 
valid in Galilean coordinates. In exactly the same way as above we show that this is equivalent 
to the equation
\begin{equation}
h^{-1}(d\beta^\prime_\nu,dx^\nu)=0
\label{2.89}
\end{equation}
where the components $\beta^\prime_\nu$ of $\beta^\prime$ in Galilean coordinates as well as the 
Galilean coordinates $x^\nu$ themselves are considered as functions on the spacetime manifold 
which may be locally expressed in any local coordinate system whatever. Also the condition 
$\omega^\prime=0$, that is $d\beta^\prime=0$, of the non-relativistic theory, is:
\begin{equation}
d\beta^\prime_\nu\wedge dx^\nu=0
\label{2.90}
\end{equation}
Equations \ref{2.89}, \ref{2.90} together constitute the general form of the {\em non-relativistic 
wave system}, identical in form to the corresponding relativistic system. This being the case, we 
shall only deal explicitly with the relativistic system from this point on. 

The specific form of the wave system refers to expressing $x^\mu$ and $\beta_\mu$ in acoustical 
coordinates. We substitute in \ref{2.86} the expression \ref{2.38} for $h^{-1}$ to obtain the specific form of this equation:
\begin{equation}
\frac{1}{2}\{(Lx^\mu)(\Lb\beta_\mu)+(\Lb x^\mu)(L\beta_\mu)\}
-a(\sh^{-1})^{AB}(\Omega_A x^\mu)(\Omega_B\beta_\mu)=0
\label{2.91}
\end{equation}
Also, we evaluate \ref{2.88} on the pair $(L,\Lb)$, the pairs $(L,\Omega_A)$, the pairs 
$(\Lb,\Omega_A)$, as well as the pairs $(\Omega_A,\Omega_B)$, to obtain the specific form of this 
equation:
\begin{eqnarray}
&&(Lx^\mu)(\Lb\beta_\mu)-(\Lb x^\mu)(L\beta_\mu)=0\nonumber\\
&&(Lx^\mu)(\Omega_A\beta_\mu)-(\Omega_A x^\mu)(L\beta_\mu)=0 \ \ : \ A=1,...,n-1\nonumber\\
&&(\Lb x^\mu)(\Omega_A\beta_\mu)-(\Omega_A x^\mu)(\Lb\beta_\mu)=0 \ \ : \ A=1,...,n-1\nonumber\\
&&(\Omega_A x^\mu)(\Omega_B\beta_\mu)-(\Omega_B x^\mu)(\Omega_A\beta_\mu)=0 \ \ 
: \ A,B=1,...,n-1 
\label{2.92}
\end{eqnarray}
$n(n+1)/2$ algebraically independent equations in all. Note that by the 1st of these equations, we have, in regard 
to the 1st term in \ref{2.91},
\begin{equation}
\frac{1}{2}\{(Lx^\mu)(\Lb\beta_\mu)+(\Lb x^\mu)(L\beta_\mu)\}=(Lx^\mu)(\Lb\beta_\mu)
=(\Lb x^\mu)(L\beta_\mu)
\label{2.93}
\end{equation}
Equations \ref{2.91} and \ref{2.92} in which \ref{2.33}, \ref{2.36}, that is:
\begin{equation}
L=\frac{\partial}{\partial\ub}-b^A\frac{\partial}{\partial\vartheta^A}, \ \ \ 
\Lb=\frac{\partial}{\partial u}+b^A\frac{\partial}{\partial\vartheta^A}, \ \ \ 
\Omega_A=\frac{\partial}{\partial\vartheta^A}
\label{2.94}
\end{equation}
are substituted, constitute the specific form of the wave system. If we take the point of view, 
complementary to that of the previous section, that the $x^\mu$ and $b^A$ are known smooth functions of 
the acoustical coordinates, the wave system constitutes a 1st order quasilinear system for the 
$\beta_\mu$ as functions of the acoustical coordinates, $a$ being defined by \ref{2.10}, that is:
\begin{equation}
a=-\frac{1}{2}h_{\mu\nu}(L x^\mu)(\Lb x^\nu)
\label{2.95}
\end{equation}
and $\sh_{AB}$ by \ref{2.37}, that is:
\begin{equation}
\sh_{AB}=h_{\mu\nu}(\Omega_A x^\mu)(\Omega_B x^\nu)
\label{2.96}
\end{equation}
recalling that the components $h_{\mu\nu}$ are expressed in terms of the $\beta_\mu$. 

We finally remark in regard to differential order that the wave system \ref{2.91}, \ref{2.92} is a 1st order quasilinear system for the $\beta_\mu$. 

\vspace{5mm}

\section{Variations by Translations and the Wave Equation for the Rectangular Components of $\beta$}

In the previous two sections, to show the general nature of the structures introduced, not only did 
we consider the general case of $n$ spatial dimensions, but also we treated simultaneously 
the relativistic and non-relativistic theories. In the present section we consider the variations 
through solutions which are generated from a given solution of the equations of motion by applying 
a 1-parameter subgroup of the group of automorphisms of the underlying spacetime structure. In 
the case of the relativistic theory this group is the isometry group 
of Minkowski spacetime, the Poincar\'{e} group, extended by the 1-parameter group of scalings. 
In the non-relativistic theory the group is the Euclidean group, extended by the 
time translations, the Galilei group, which corresponds to the transformations between Galilei 
frames, and the scaling group. In the non-relativistic theory the scaling group is 2-dimensional, 
as lengths and time intervals can be independently scaled, however the equation of state, which 
expresses $p$ as a function of $h^\prime$ ($s$ is a constant), in general reduces this to a 1-parameter group. 
Homogeneous equations of state constitute an exception, in which the full 2-dimensional scaling 
symmetry is retained (see Chapter 1 of [Ch-Mi]). 

In particular, in the irrotational case of the relativistic theory, if $\phi$ is a solution of the 
nonlinear wave equation \ref{1.70} and $f_\lambda$ is a 1-parameter group of isometries of the 
Minkowski metric $g$, then 
\begin{equation}
\phi_\lambda=\phi\circ f_\lambda
\label{2.97}
\end{equation}
is a 1-parameter family of solutions and (see \ref{1.91})
\begin{equation}
\dot{\phi}=\left(\frac{d\phi_\lambda}{d\lambda}\right)_{\lambda=0}=X\phi
\label{2.98}
\end{equation}
is a variation through solutions. Here $X$ is the vectorfield generating the 1-parameter group 
$f_\lambda$. The isometry group of Minkowski spacetime consists of the spacetime translations, 
generated by the vectorfields:
\begin{equation}
X_{(\mu)}=\frac{\partial}{\partial x^{\mu}} \ \mbox{: in terms of rectangular coordinates}
\label{2.99}
\end{equation}
and the spacetime rotations or Lorentz transformations, generated by the vectorfields:
\begin{equation}
X_{(\mu\nu)}=x_\mu\frac{\partial}{\partial x^\nu}-x_\nu\frac{\partial}{\partial x^\mu}, 
\ \mu\neq\nu \mbox{: in terms of rectangular coordinates}
\label{2.100}
\end{equation}
Here:
\begin{equation}
x_\mu=g_{\mu\nu}x^\nu \ \mbox{: in terms of rectangular coordinates}
\label{2.101}
\end{equation}
The 1-parameter group of scalings from an origin is given in terms of any system of linear coordinates by:
\begin{equation}
f_\lambda(x)=e^\lambda x
\label{2.102}
\end{equation}
however in this case the 1-parameter family of solutions is:
\begin{equation}
\phi_\lambda=e^{-\lambda}\phi\circ f_\lambda
\label{2.103}
\end{equation}
If $h_\lambda$ is the enthalpy associated to $\phi_\lambda$ we then have:
\begin{equation}
h_\lambda=h\circ f_\lambda
\label{2.104}
\end{equation}
The vectorfield $S$ generating the the scaling group is, in terms of any system of linear coordinates, 
\begin{equation}
S=x^\mu\frac{\partial}{\partial x^\mu}
\label{2.105}
\end{equation}
and the corresponding variation through solutions is: 
\begin{equation}
\dot{\phi}=\left(\frac{d\phi_\lambda}{d\lambda}\right)_{\lambda=0}=S\phi-\phi
\label{2.106}
\end{equation}

As we have seen in Chapter 1, in the irrotational case of the relativistic theory a variation 
through solutions $\dot{\phi}$ of the function $\phi$ satisfies the wave equation \ref{1.96} 
associated to the conformal acoustical metric \ref{1.94}. The conformal factor $\Omega$ is actually 
given by \ref{1.96} only in the case of 3 spatial dimensions. This is because equation \ref{1.95} 
holds as it stands only in 3 spatial dimensions. The general form of this equation is:
\begin{equation}
\sqrt{-\mbox{det}\tilde{h}}=\Omega^{(n+1)/2}\sqrt{-\mbox{det}h}=\Omega^{(n+1)/2}\eta
\label{2.107}
\end{equation}
As a consequence $\Omega$ is in general given by:
\begin{equation}
\Omega^{(n-1)/2}=\frac{G}{\eta}, \ \ \mbox{that is:} \ \ \Omega=\left(\frac{G}{\eta}\right)^{2/(n-1)}
\label{2.108}
\end{equation}

In particular, the variations generated by the translations \ref{2.99} are the rectangular components 
$\beta_\mu$ of the 1-form $\beta$, for:
\begin{equation}
X_{(\mu)}\phi=\beta(X_{(\mu)})=\beta_\mu
\label{2.109}
\end{equation}
We conclude that the functions $\beta_\mu$ satisfy the wave equation associated to the conformal acoustical 
metric:
\begin{equation}
\square_{\tilde{h}}\beta_\mu=0
\label{2.110}
\end{equation}
In fact, since for any vectorfield $X$ we have
$$X\phi=\beta(X),$$
the functions $\beta(X_{(\mu\nu)})$, with $X_{(\mu\nu)}$ the generators of spacetime 
rotations \ref{2.100}, also satisfy the same equation. 

In the non-relativistic theory, the Euclidean group extended by the time translations acts on 
Galilei spacetime in a manner similar to that of an isometry group, so if $f_\lambda$ is the 1-parameter subgroup generated by the vectorfield $X$ and $\phi^\prime$ is a solution of \ref{1.169} 
then 
\begin{equation}
\phi^\prime_\lambda=\phi^\prime\circ f_\lambda
\label{2.111}
\end{equation}
is a 1-parameter family of solutions of the same equation and 
\begin{equation}
\dot{\phi^\prime}=\left(\frac{d\phi^\prime_\lambda}{d\lambda}\right)_{\lambda=0}=X\phi^\prime
\label{2.112}
\end{equation}
is a variation through solutions. The spacetime translations are again generated by the vectorfields 
$X_{(\mu)}$ given in terms of Galilean coordinates by \ref{2.99}, and the spatial rotations 
are given by \ref{2.100} with the indices $\mu, \nu$ restricted to $i,j=1,...,n$, that is:
\begin{equation}
X_{(ij)}=x^i\frac{\partial}{\partial x^j}-x^j\frac{\partial}{\partial x^i}, \ i\neq j 
\mbox{: in terms of Galilean coordinates} 
\label{2.113}
\end{equation}

What is more complicated in the non-relativistic framework is the action of the Galilei group 
on a solution $\phi^\prime$ of \ref{1.169}. If we denote 
$\overline{x}=(x^1,..., x^n)\in\mathbb{R}^n$, the 
Galilei group is the group:
\begin{equation}
f^\prime_k=(t,\overline{x})\mapsto (t,\overline{x}+kt) \ : \ k\in\mathbb{R}^n
\label{2.114}
\end{equation}
As a set of transformations, it is the limit as $c\rightarrow\infty$ of the boosts 
\begin{equation}
f_k=(t,\overline{x})\mapsto\left(\frac{t+c^{-2}<k,\overline{x}>}{\sqrt{1-c^{-2}|k|^2}},
\frac{\overline{x}+kt}{\sqrt{1-c^{-2}|k|^2}}\right) \ : \ k\in\mathbb{R}^n
\label{2.115}
\end{equation}
which are transformations belonging to the Lorentz group of Minkowski spacetime. Since $\phi^\prime$ 
is obtained as the limit as $c\rightarrow\infty$ of $c(\phi-ct)$, to find how $\phi^\prime$ 
transforms under $f^\prime_k$ we consider the limit as $c\rightarrow\infty$ of
$$c(\phi\circ f_k-ct)$$
This yields the conclusion that under  $f^\prime_k$,  $\phi^\prime$ transforms to $\phi^\prime_k$, where: 
\begin{equation}
\phi^\prime_k(t,\overline{x})=\frac{1}{2}|k|^2 t+<k,\overline{x}>+\phi^\prime(t,\overline{x}+kt)
\label{2.116}
\end{equation}
the last term being $\phi^\prime\circ f_k$ at $(t,\overline{x})$. Then $v_k$ and $h^\prime_k$ 
the spatial velocity and enthalpy corresponding to $\phi^\prime_k$ are given by:
\begin{equation}
v_k=v\circ f^\prime_k-k, \ \ \ h^\prime_k=h^\prime\circ f^\prime_k
\label{2.117}
\end{equation}
It is straightforward to check that for each $k\in \mathbb{R}^n$, $\phi^\prime_k$ is a solution 
of \ref{1.169}. Setting $k=\lambda e_i$ where $e_i$ is the $i$th basis vector in $\mathbb{R}^n$, 
the corresponding variation through solutions is: 
\begin{equation}
x^i+t\frac{\partial\phi^\prime}{\partial x^i}
\label{2.118}
\end{equation}

As already mentioned above, in the non-relativistic theory the scaling group is a priori 2-dimensional 
since lengths and time intervals can be independently scaled. So, we may define:
\begin{equation}
f_{(\mu,\lambda)} : (t,\overline{x})\mapsto (e^\mu t,e^\lambda\overline{x})
\label{2.119}
\end{equation}
Since $\phi^\prime$ has the physical dimensions of $(\mbox{length})^2/\mbox{time}$, under 
$f_{(\mu,\lambda)}$, $\phi^\prime$ transforms to $\phi^\prime_{(\mu,\lambda)}$, where: 
\begin{equation}
\phi^\prime_{(\mu,\lambda)}=e^{\mu-2\lambda}\phi^\prime\circ f_{(\mu,\lambda)}
\label{2.120}
\end{equation}
Then $v_{\mu,\lambda}$ and $h^\prime_{(\mu,\lambda)}$ the spatial velocity and enthalpy corresponding 
to $\phi^\prime_{(\mu,\lambda)}$ are given by:
\begin{equation}
v_{(\mu,\lambda)}=e^{\mu-\lambda}v\circ f_{(\mu,\lambda)}, \ \ \ 
h^\prime_{(\mu,\lambda)}=e^{2(\mu-\lambda)}h^\prime\circ f_{(\mu,\lambda)}
\label{2.121}
\end{equation}
Because of the transformation law of the enthalpy, unless the equation of state gives $p$ as a homogeneous function of $h^\prime$, $\phi^\prime_{(\mu,\lambda)}$ is a solution of \ref{1.169} only 
when $\mu=\lambda$. Thus in general the equation of state reduces the 2-parameter scaling group 
to the 1-parameter subgroup
\begin{equation}
f_\lambda=f_{(\lambda,\lambda)}
\label{2.122}
\end{equation}
under which $\phi^\prime$ transforms to $\phi^\prime_\lambda$, where:
\begin{equation}
\phi^\prime_\lambda=e^{-\lambda}\phi^\prime\circ f_\lambda
\label{2.123}
\end{equation}
The spatial velocity and enthalpy associated to $\phi^\prime_\lambda$ are:
\begin{equation}
v_\lambda=v\circ f_\lambda, \ \ \ h^\prime_\lambda=h^\prime\circ f_\lambda
\label{2.124}
\end{equation}
and the corresponding variation through solutions is:
\begin{equation}
S\phi^\prime-\phi^\prime
\label{2.125}
\end{equation}
where $S$ is given in Galilean coordinates by \ref{2.105}. 

As we have seen in Chapter 1, in the irrotational case of the non-relativistic theory a variation 
through solutions $\dot{\phi}^\prime$ of the function $\phi^\prime$ satisfies the wave equation \ref{1.179} 
associated to the conformal acoustical metric \ref{1.177}, with the conformal factor $\Omega$ 
given in the case of 3 spatial dimensions by \ref{1.178}. This is because equation \ref{1.95} 
holds as it stands only in 3 spatial dimensions. Since in the general case of $n$ spatial dimensions 
we have, as in the relativistic theory (see \ref{2.107}),
$$\sqrt{-\mbox{det}\tilde{h}}=\Omega^{(n+1)/2}\sqrt{-\mbox{det}h}=\Omega^{(n+1)/2}\eta$$
$\Omega$ is in the non-relativistic theory in $n$ spatial dimensions given by:
\begin{equation}
\Omega=\left(\frac{n}{\eta}\right)^{2/(n-1)}
\label{2.126}
\end{equation}

In particular, the variations generated by the spacetime translations  are the components 
$\beta^\prime_\mu$ of the 1-form $\beta^\prime$ in Galilean coordinates, for:
\begin{equation}
X_{(\mu)}\phi^\prime=\beta^\prime(X_{(\mu)})=\beta^\prime_\mu
\label{2.127}
\end{equation}
We conclude that in analogy with the relativistic theory the functions $\beta^\prime_\mu$ 
in the non-relativistic theory satisfy the wave equation associated to the conformal acoustical 
metric \ref{1.177} with the  conformal factor \ref{2.126}:
\begin{equation}
\square_{\tilde{h}}\beta^\prime_\mu=0
\label{2.128}
\end{equation}
In fact, since for any vectorfield $X$ we have
$$X\phi^\prime=\beta^\prime(X),$$
the functions $\beta(X_{(ij)})$, with $X_{(ij)}$ the generators of spatial 
rotations \ref{2.113}, also satisfy the same equation. 

Returning now to the relativistic theory, to elucidate things further, we replace the underlying 
Minkowski spacetime $\mathbb{M}^n$, by an arbitrary $n+1$-dimensional Lorentzian manifold 
$({\cal M},g)$. We denote from this point on by $\sigma$ the squared enthalpy, given in general 
by:
\begin{equation}
\sigma=-g^{-1}(\beta,\beta)
\label{2.129}
\end{equation}
and consider $G$, $F$, and $H$ to be functions of $\sigma$. Thus, $s$ here being a constant, 
equation \ref{1.d5} defines the sound speed $\eta$, itself a function of $\sigma$, by:
\begin{equation}
\frac{2\sigma}{G}\frac{dG}{d\sigma}=\frac{1}{\eta^2}-1
\label{2.130}
\end{equation}
Also, the functions $F$ and $H$ defined by \ref{1.d7} and \ref{1.d12}, that is:
\begin{equation}
F=\frac{2}{G}\frac{dG}{d\sigma}, \ \ \ H=\frac{F}{1+\sigma F}
\label{2.131}
\end{equation}
satisfy (see \ref{1.d13}):
\begin{equation}
1-\sigma H=\eta^2=\frac{1}{1+\sigma F}
\label{2.132}
\end{equation}
The acoustical metric is in general given by \ref{1.d11}:
\begin{equation}
h=g+H\beta\otimes\beta
\label{2.133}
\end{equation}
and the reciprocal acoustical metric by \ref{1.d10}. 
We shall establish the following proposition without referring at all to the function $\phi$. 

\vspace{5mm}

\noindent {\bf Proposition 2.1} \ \ \ Let $({\cal M},g)$ be a $n+1$ dimensional Lorentzian manifold 
and $\beta$ a closed 1-form on ${\cal M}$, such that with 
$$\sigma=-g^{-1}(\beta,\beta)>0$$
and $G>0$ a given function of $\sigma$ such that
$$\frac{2}{G}\frac{dG}{d\sigma}>0,$$
we have, in terms of arbitrary local coordinates, 
$$\nabla_\mu (G\beta^\mu)=0, \ \mbox{where} \ \beta^\mu=(g^{-1})^{\mu\nu}\beta_\nu.$$
Let $X$ be a vectorfield generating a 1-parameter group of isometries of $({\cal M},g)$. Then the 
function $\beta(X)$ satisfies the equation:
$$\square_{\tilde{h}}\beta(X)=0$$
with $\tilde{h}=\Omega h$, 
$$h=g+H\beta\otimes\beta, \ \ \ \Omega=\left(\frac{G}{\eta}\right)^{2/(n-1)},$$
$\eta>0$, being defined by: 
$$\frac{1}{\eta^2}=1+\frac{2}{G}\frac{dG}{d\sigma},$$
so also $\eta<1$, and $H>0$ by:
$$\sigma H=1-\eta^2.$$
so $h$ and $\tilde{h}$ are Lorentzian metrics on ${\cal M}$.

\vspace{2.5mm}

\noindent {\em Proof:} Let $A$ be the vectorfield, with components, in arbitrary local coordinates, 
\begin{equation}
A^\mu=G(g^{-1})^{\mu\nu}\beta_\nu
\label{2.134}
\end{equation}
so $A$ is divergence-free relative to $g$:
\begin{equation}
\nabla_\mu A^\mu=0
\label{2.135}
\end{equation}
Let $\alpha$ be the $n$-form corresponding to $A$ through the volume $n+1$ form $\epsilon$ of 
$({\cal M},g)$:
\begin{equation}
\alpha=i_A\epsilon
\label{2.136}
\end{equation}
Thus, the $n$-form $\alpha$ is the Hodge dual of the 1-form $G\beta$ relative to $g$:
\begin{equation}
\alpha=\s^*(G\beta)
\label{2.137}
\end{equation}
In terms of $\alpha$, equation \ref{2.135} reads, simply:
\begin{equation}
d\alpha=0
\label{2.138}
\end{equation}
$d$ being the exterior derivative. Since for any vectorfield $X$, ${\cal L}_X$ commutes with $d$, 
this implies:
\begin{equation}
d{\cal L}_X\alpha=0
\label{2.139}
\end{equation}
for any vectorfield $X$. If $X$ is a Killing field of $g$, that is $X$ generates isometries of 
$({\cal M},g)$, then ${\cal L}_X$ also commutes with the Hodge dual, hence:
\begin{equation}
{\cal L}_X\alpha=\s^*({\cal L}_X(G\beta))
\label{2.140}
\end{equation}
We have:
\begin{equation}
{\cal L}_X(G\beta)=(XG)\beta+G{\cal L}_X\beta
\label{2.141}
\end{equation}
Now, in view of the general relation \ref{1.58}, by virtue of the fact that $d\beta=0$ we have:
\begin{equation}
{\cal L}_X\beta=d\beta(X)
\label{2.142}
\end{equation}
Next, by the 1st of \ref{2.131}, 
\begin{equation}
XG=\frac{dG}{d\sigma}X\sigma=\frac{1}{2}GFX\sigma
\label{2.143}
\end{equation}
From the definition \ref{2.129} of $\sigma$, the vectorfield $X$ being a Killing field of $g$ we have:
\begin{equation}
X\sigma=-2g^{-1}(\beta,{\cal L}_X\beta)
\label{2.144}
\end{equation}
Therefore in view of \ref{2.142} we obtain:
\begin{equation}
XG=-GFg^{-1}(\beta,d\beta(X))
\label{2.145}
\end{equation}
Substituting \ref{2.142} and \ref{2.145} in \ref{2.141} we obtain:
\begin{equation}
{\cal L}_X(G\beta)=G(d\beta(X)-Fg^{-1}(\beta,d\beta(X))\beta)
\label{2.146}
\end{equation}
Now, $h^{-1}$ takes a 1-form $\theta$ to the vectorfield:
\begin{equation}
h^{-1}\cdot\theta=g^{-1}\cdot(\theta-Fg^{-1}(\beta,\theta)\beta)
\label{2.147}
\end{equation}
Taking $\theta=d\beta(X)$ and comparing with \ref{2.146} we see that: 
\begin{equation}
g^{-1}\cdot{\cal L}_X(G\beta)=GB
\label{2.148}
\end{equation}
where $B$ is the vectorfield:
\begin{equation}
B=h^{-1}\cdot d\beta(X)
\label{2.149}
\end{equation}
In view of \ref{2.140} it follows that:
\begin{equation}
{\cal L}_X\alpha=Gi_B\epsilon
\label{2.150}
\end{equation}
(compare with \ref{2.136}). Let now $\Sigma$ be the null space of $\beta$ at a point. Then 
$\Sigma$ is space-like relative to $g$ and by the definition of $h$ we have:
\begin{equation}
h|_\Sigma=g|_\Sigma
\label{2.151}
\end{equation}
More generally, if $Y\in\Sigma$ and $V$ is an arbitrary vector at the same point, it holds:
\begin{equation}
h(V,Y)=g(V,Y)
\label{2.152}
\end{equation}
Let $u$ be vector at the same point which is $g$-orthogonal to $\Sigma$. Then $u$ is time-like 
relative to $g$. Taking $u$ to be of unit magnitude with respect to $g$, and taking $(E_i:1=1,...,n)$ 
to be a positive orthonormal basis for $\Sigma$, then $(u,E_1,...,E_n)$ is a positive orthonormal 
basis with respect to $g$ for the tangent space at the point. Therefore we can express at this point 
$g^{-1}$ in 
terms of this basis as:
\begin{equation}
g^{-1}=-u\otimes u+E_i\otimes E_i
\label{2.153}
\end{equation}
It then follows from the definition of $\sigma$, in view of the fact that $\beta(E_i)=0:i=1,...,n$, 
that:
\begin{equation}
\sigma=(\beta(u))^2
\label{2.154}
\end{equation}
hence, by the definition of $h$, while $g(u,u)=-1$, we have:
\begin{equation}
h(u,u)=-1+\sigma H=-\eta^2
\label{2.155}
\end{equation}
Thus $u$ is also $h$-orthogonal to $\Sigma$, but its magnitude with respect to $h$ is $\eta$. 
Therefore the corresponding positive orthonormal basis with respect to $h$ is 
$(\eta^{-1}u,E_1,...,E_n)$. Consequently the volume form $e$ associated to $h$ is related to 
the volume form $\epsilon$ associated to $g$ by:
\begin{equation}
e=\eta\epsilon
\label{2.156}
\end{equation}
Then \ref{2.150} can be written in the form:
\begin{equation}
{\cal L}_X\alpha=\frac{G}{\eta}i_B e
\label{2.157}
\end{equation}
Let $\tilde{B}$ be the vectorfield analogous to $B$ (see \ref{2.149}), but defined with the 
conformal metric $\tilde{h}=\Omega h$ in the role of $h$:
\begin{equation}
\tilde{B}=\tilde{h}^{-1}\cdot d\beta(X)
\label{2.158}
\end{equation}
Then obviously:
\begin{equation}
\tilde{B}=\Omega^{-1} B
\label{2.159}
\end{equation}
On the other hand, the volume form $\tilde{e}$ associated to $\tilde{h}$ is related to the volume 
form $e$ associated to $h$ by:
\begin{equation}
\tilde{e}=\Omega^{(n+1)/2} e
\label{2.160}
\end{equation}
It follows that \ref{2.157} can be written in the form:
$${\cal L}_X\alpha=\frac{G}{\eta}\Omega^{-(n-1)/2}i_{\tilde{B}}\tilde{e}$$
But, in view of the definition of the conformal factor $\Omega$, this is simply:
\begin{equation}
{\cal L}_X\alpha=i_{\tilde{B}}\tilde{e}
\label{2.161}
\end{equation}
Moreover, in view of the definition \ref{2.158}, the right hand side is the Hodge dual of the 
1-form $d\beta(X)$ relative to the metric $\tilde{h}$. Denoting this Hodge dual with a $\star$, 
\ref{2.161} is simply:
\begin{equation}
{\cal L}_X\alpha=\s^\star(d\beta(X))
\label{2.162}
\end{equation}
Therefore by equation \ref{2.139} we have:
\begin{equation}
d\s^\star(d\beta(X))=0
\label{2.163}
\end{equation}
which is equivalent to:
\begin{equation}
\square_{\tilde{h}}\beta(X)=0
\label{2.164}
\end{equation}
and the proposition is proved. 

\vspace{5mm}

\section{The Geometric Construction for the Shock Development Problem}

We now return to the formulation of the shock development problem in Section 1.5 and that of the 
restricted shock development problem in Section 1.6. As stated there, the problem is associated to a given 
component of $\partial_-{\cal B}$ and the associated components of ${\cal B}$ and of 
$\underline{{\cal C}}$. From this point, by $\partial_-{\cal B}$, ${\cal B}$ and 
$\underline{{\cal C}}$ we shall mean the components in question. We denote by ${\cal M}$ the domain 
of the maximal classical development and by ${\cal N}$ the domain of the new solution sought. 
Thus, the past boundary of ${\cal N}$ is $\underline{{\cal C}}\bigcup{\cal K}$, ${\cal K}$ being 
the shock hypersurface, which is also sought in the problem. What is known in regard to ${\cal K}$ is 
its past boundary, which coincides with that of ${\cal B}$:
\begin{equation}
\partial_-{\cal K}=\partial_-{\cal B}
\label{2.165}
\end{equation}
Moreover we know that ${\cal K}$ is to precede ${\cal B}$ except at its past boundary. In the part of 
${\cal M}$ lying in the past of ${\cal K}$ the solution given by the maximal classical development 
holds. 

In the physical case of 3 spatial dimensions ${\cal B}$ is represented in acoustical coordinates 
$(t,u,\vartheta)$ in the sense of [Ch-S] by \ref{1.297}. Also, in canonical acoustical coordinates 
in the sense of [Ch-S], where the invariant curves of ${\cal B}$ are the curves of constant 
$\vartheta$ on ${\cal B}$, by \ref{1.315}, \ref{1.316}, \ref{1.319}, the function $u$ attains its maximum $u_M$ 
along an invariant curve at its past end point on $\partial_-{\cal B}$. Thus $\partial_-{\cal B}$ 
is represented in canonical acoustical coordinates by:
\begin{equation}
\partial_-{\cal B}=\{(t_*(u_M(\vartheta),\vartheta),u_M(\vartheta),\vartheta) \ : \ \vartheta\in U\}
\label{2.166}
\end{equation}
and ${\cal B}$ is represented in the same coordinates by:
\begin{equation}
{\cal B}=\{(t_*(u,\vartheta),u,\vartheta) \ : \ u<u_M(\vartheta),\ \vartheta\in U\}
\label{2.167}
\end{equation}
where $U$ is $S^2$ or a domain in $S^2$, $u_M$ being a smooth function on $U$. 
As the problem can be localized 
in regard to $\partial_-{\cal B}$ by reason of the domain of dependence property, 
there is no loss of generality in restricting ourselves to 
the case that $U$ is $S^2$, that is, in taking $\partial_-{\cal B}$ to be diffeomorphic to $S^2$. 
Performing a suitable transformation of the form \ref{1.321}, thereby changing to 
a different acoustical function $u^\prime$, $\partial_-{\cal B}$ can then be represented simply as:
\begin{equation}
\partial_-{\cal B}=\{(t_*(0,\vartheta^\prime),0,\vartheta^\prime) \ : \ \vartheta^\prime\in S^2\}
\label{2.168}
\end{equation}
and ${\cal B}$ as:
\begin{equation}
{\cal B}=\{(t_*(u^\prime,\vartheta^\prime) \ : \ u^\prime<0,\ \vartheta^\prime\in S^2\}
\label{2.169}
\end{equation}
In fact it suffices to choose the functions $f$ and $g$ in \ref{1.321} so that:
\begin{equation}
f(u_M(\vartheta),\vartheta)=0, \ \ \ g(\vartheta)=\vartheta
\label{2.170}
\end{equation}

From this point, whenever necessary to avoid confusion, we shall denote quantities referring to 
the maximal classical development with a prime. We first extend the acoustical function $u^\prime$ 
to ${\cal N}$ by the function $u$, defined by the condition:
\begin{equation}
u|_{\underline{{\cal C}}}=u^\prime_{\underline{{\cal C}}}
\label{2.171}
\end{equation}
together with the condition that the level sets $C_u$ of $u$ are outgoing acoustically null hypersurfaces 
in ${\cal N}$. The boundary hypersurface ${\cal K}$ being acoustically time-like relative to 
${\cal N}$, each $C_u$, level set of $u$, intersects ${\cal K}$ in an acoustically space-like 
section $S^*_u$, diffeomorphic to $S^2$. Also, since at any point in ${\cal N}$, $du\cdot V>0$ 
for any vector $V$ belonging to the interior of the future sound cone at the point, 
$u|_{\cal K}$ is a time function relative to the Lorentzian, induced from ${\cal N}$, metric on 
${\cal K}$. Thus $u|_{\cal K}$ is increasing toward the future. On the other hand, ${\cal K}$ is acoustically space-like relative to ${\cal M}$, so the metric 
induced on ${\cal K}$ by the maximal classical development is positive definite, and 
$u^\prime|_{\cal K}$ is decreasing outward. The functions $u^\prime$ and $u$ both vanish at 
$\partial_-{\cal B}$. Thus, if we view $u$ as being the extension of $u^\prime$ to ${\cal N}$, 
this acoustical function is discontinuous across ${\cal K}$, suffering a positive jump as we cross 
${\cal K}$ from the past to the future. Moreover, the size of the jump increases toward the future 
along ${\cal K}$. 

Next we define the conjugate acoustical function $\ub$ in ${\cal N}$, by the condition:
\begin{equation}
\ub|_{{\cal K}}=u|_{{\cal K}}
\label{2.172}
\end{equation}
together with the condition that the level sets $\Cb_{\ub}$ of $\ub$ are incoming acoustically 
null hypersurfaces in ${\cal N}$. We have:
\begin{equation}
\Cb_0=\underline{{\cal C}}
\label{2.173}
\end{equation}
Denoting by $S_{\ub,u}$ the surfaces of intersection $\Cb_{\ub}\bigcap C_u$, as in \ref{2.2}, we have:
\begin{equation}
S_{\tau,\tau}=S^*_\tau \ : \ \tau\geq 0, \ \ \ S_{0,0}=S^*_0=\partial_-{\cal B}
\label{2.174}
\end{equation}

Thus in ${\cal N}$ we have the structures discussed in the previous sections of the present 
chapter. Moreover, the hypotheses on which our construction relies are of a qualitative nature 
and will hold in any number $n$ of spatial dimensions greater than 1. We shall consider the 
general case of $n$ spatial dimensions when discussing the main features of the mathematical 
structure, the aim being a fuller understanding. 

The $S^{n-1}$ valued function $\vartheta$ is defined in ${\cal N}$ as follows. 
It is first defined on $\partial_-{\cal B}$ by the requirement:
\begin{equation}
\vartheta|_{\partial_-{\cal B}}=\vartheta^\prime|_{\partial_-{\cal B}}
\label{2.175}
\end{equation}
the $S^{n-1}$ valued function $\vartheta^\prime$ being that defined in ${\cal M}$, as discussed 
above. It is then defined on $\underline{{\cal C}}=\Cb_0$ by the condition that it is constant  
along the generators of $\Cb_0$. Finally it is defined in ${\cal N}$, in accordance with the 
prescription of Section 2.1, by the condition that it is constant along the integral curves of 
$T$. Note that ${\cal K}$ is itself generated by integral curves of $T$. This is the reason why 
we follow the said prescription. 

In terms of the representation \ref{2.32}, ${\cal N}$ corresponds to the domain:
\begin{equation}
D=\{(\ub,u) \ : \ 0\leq\ub\leq u\}
\label{2.176}
\end{equation}
and its boundary $\underline{{\cal C}}\bigcup{\cal K}$ to: 
\begin{equation}
\partial D=\{(0,u) \ : \ u\geq 0\}\bigcup\{(\tau,\tau) \ : \ \tau\geq 0\}
\label{2.177}
\end{equation}
with $\underline{{\cal C}}$ corresponding to $\{(0,u) \ : \ u\geq 0\}$ and ${\cal K}$ to 
$\{(\tau,\tau) \ : \ \tau\geq 0\}$, and their common past boundary, $\partial_-{\cal B}$, to 
$(0,0)$. 

\pagebreak

\chapter{The Acoustical Structure Equations}
\section{The Connection Coefficients 
and the First Variation Equations}

In the present section we shall define quantities and shall deduce certain 
equations in regard to these quantities, which, as in Section 2.1, depend only on 
$({\cal N},h)$ as a Lorentzian manifold and not on the underlying Minkowskian structure. 
We denote by $D$ the covariant derivative operator associated to $({\cal N},h)$. 

First some basic concepts and notation. Let $\xi$ a 1-form on ${\cal N}$ such that
\begin{equation}
\xi(L)=\xi(\Lb)=0
\label{3.1}
\end{equation}
We then say that $\xi$ is an $S$ 1-form. Conversely, given for each $(\ub,u)\in D$ a 1-form $\xi$ 
intrinsic to the surface $S_{\ub,u}$, then along each $S_{\ub,u}$, $\xi$ extends to the tangent bundle 
of ${\cal N}$ over $S_{\ub,u}$ by the conditions \ref{3.1}. We then obtain an $S$ 1-form on 
${\cal N}$. Thus an $S$ 1-form may be thought of as the specification of a 1-form intrinsic to 
$S_{\ub,u}$ for each $(\ub,u)$. Similar considerations apply to $p$-covariant tensorfields. Thus, 
a $p$-covariant $S$ tensorfield is a $p$-covariant tensorfield which vanishes if either $L$ or $\Lb$ 
is inserted in one of its $p$ entries, and can be thought of as the specification of a $p$-covariant 
tensorfield intrinsic to $S_{\ub,u}$ for each $(\ub,u)$. Next, an $S$ vectorfield on ${\cal N}$ 
is a vectorfield $X$ on ${\cal N}$ such that $X$ is at each point $x\in{\cal N}$ tangential to the surface $S_{\ub,u}$ through $x$. Thus an $S$ vectorfield may be thought of as the specification of 
a vectorfield intrinsic to $S_{\ub,u}$ for each $(\ub,u)$. Similarly, a $q$-contravariant $S$ 
tensorfield  on ${\cal N}$ is a $q$-contravariant  tensorfield $W$ on ${\cal N}$ such that at 
each point $x\in{\cal N}$ $W$ at $x$ belongs to $\otimes^q T_x S_{\ub,u}$, and can be thought of as 
a $q$-contravariant vectorfield intrinsic to $S_{\ub,u}$ for each $(\ub,u)$. Finally, a type $T^q_p$ 
$S$ tensorfield $\theta$ on ${\cal N}$ is a type $T^q_p$ tensorfield on ${\cal N}$  such that at each 
$x\in{\cal N}$ and each $X_1,...,X_p\in T_x{\cal N}$ we have 
$\theta(X_1,...,X_p)\in\otimes^q T_xS_{\ub,u}$ and $\theta(X_1,...,X_p)=0$ if one of 
$X_1,...,X_p$ is either $L$ or $\Lb$. Thus a general type $T^q_p$ $S$ tensorfield may be thought of 
as the specification of a type $T^q_p$ tensorfield intrinsic to $S_{\ub,u}$ for each $(\ub,u)$. 

Let $\xi$ be an $S$ 1-form. We define $\sL_L\xi$ to be, for each $(\ub,u)\in D$, the restriction to 
$TS_{\ub,u}$ of ${\cal L}_L\xi$, the Lie derivative of $\xi$ with respect to $L$, a notion intrinsic 
to each hypersurface $C_u$. Then $\sL_L\xi$ is an $S$ 1-form as well. Note that in any case we have 
$({\cal L}_L\xi)(L)=0$. Similarly, we define $\sL_{\Lb}\xi$ to be, for each $(\ub,u)\in D$, the 
restriction to $TS_{\ub,u}$ of ${\cal L}_{\Lb}\xi$, the Lie derivative of $\xi$ with respect to $\Lb$, 
a notion intrinsic to each hypersurface $\Cb_{\ub}$. Then $\sL_{\Lb}\xi$ is an $S$ 1-form as well. 
Note that in any case we have $({\cal L}_{\Lb}\xi)(\Lb)=0$. For any $p$-covariant $S$ tensorfield 
$\xi$, the derivatives $\sL_L\xi$ and $\sL_{\Lb}\xi$ are similarly defined. 

Consider next the case of an $S$ vectorfield $X$. 

\vspace{5mm}

\noindent {\bf Lemma 3.1} \ \ \ Let $X$ be a vectorfield defined along a given $C_u$ and tangential to 
its $S_{\ub,u}$ sections. Then the vectorfield $[L,X]$ defined along $C_u$ is also tangential to its 
$S_{\ub,u}$ sections. Also, let $X$ be a vectorfield defined along a given $\Cb_{\ub}$ and tangential 
to its $S_{\ub,u}$ sections. Then the vectorfield $[\Lb,X]$ defined along $\Cb_{\ub}$ is also 
tangential to its $S_{\ub,u}$ sections.

\vspace{2.5mm}

\noindent {\em Proof:} The first statement is obvious from the facts that $[L,X]u=L(Xu)-X(Lu)=0$, 
$[L,X]\ub=L(X\ub)-X(L\ub)=0$. The second statement is likewise obvious from the facts that 
$[\Lb,X]u=\Lb(Xu)-X(\Lb u)=0$, $[\Lb,X]\ub=\Lb(X\ub)-X(\Lb\ub)=0$. 

\vspace{2.5mm}

If $Y$ is an $S$ vectorfield, then according to the above lemma ${\cal L}_L Y=[L,Y]$ and 
${\cal L}_{\Lb}Y=[\Lb,Y]$ are also $S$ vectorfields. We thus define the derivatives $\sL_L Y$ and 
$\sL_{\Lb}Y$ to be simply ${\cal L}_L Y$ and ${\cal L}_{\Lb}Y$ respectively. The case of any  
$q$-contravariant tensorfield $W$ in the role of $Y$ is formally identical to the case of an $S$ 
vectorfield, such a tensorfield being expressible as the sum of tensor products of $S$ vectorfields. 
Thus the derivatives $\sL_L W$ and $\sL_{\Lb}W$ are defined to be simply ${\cal L}_L W$ and 
${\cal L}_{\Lb}W$ respectively, which are themselves $q$-contravariant $S$ tensorfields. 

Consider finally the general case of an arbitrary type $T^q_p$ $S$ tensorfield $\theta$. We may 
consider, in accordance with the above, $\theta$ as being, at each surface $S_{\ub,u}$ and at each 
point $x\in S_{\ub,u}$, a $p$-linear form in $T_x S_{\ub,u}$ with values in 
$\otimes^q T_x S_{\ub,u}$. Then $\sL_L\theta$ is defined by considering $\theta$ on each $C_u$ 
extended to $TC_u$ according to the condition that it vanishes if one of the entries is $L$ and 
setting $\sL_L\theta$ equal to the restriction to the $TS_{\ub,u}$ of the usual Lie derivative with 
respect to $L$ of this extension. Similarly, $\sL_{\Lb}\theta$ is defined by considering $\theta$ on 
each $\Cb_{\ub}$ extended to $T\Cb_{\ub}$ according to the condition that it vanishes if one of the 
entries is $\Lb$ and setting $\sL_{\Lb}\theta$ equal to the restriction to the $TS_{\ub,u}$ of the 
usual Lie derivative with respect to $\Lb$ of this extension. Then $\sL_L\theta$ and $\sL_{\Lb}\theta$ 
are themselves type $T^q_p$ $S$ tensorfields. 

For any function $f$ defined on ${\cal N}$ we denote by $\sd f$ the $S$ 1-form which is 
the differential of the restriction to each $S_{\ub,u}$ of $f$. 

\vspace{2.5mm}

\noindent {\bf Lemma 3.2} \ \ \  For any function $f$ defined on ${\cal N}$ we have
$$\sL_L\sd f=\sd Lf, \ \ \ \sL_{\Lb}\sd f=\sd\Lb f$$

\vspace{2.5mm}

\noindent {\em Proof:} Consider the first equality. If we evaluate each side on $L$ or $\Lb$, then 
both sides vanish by definition. On the other hand, if we evaluate on an $S$ vectorfield $X$, the 
left hand side is 
$$({\cal L}_L\sd f)(X)=L((\sd f)(X))-(\sd f)([L,X])=L(Xf)-[L,X]f=X(Lf)$$
and the right hand side is
$$(\sd Lf)(X)=X(Lf)$$
as well. This establishes the first equality. The second equality is established in a similar manner.

\vspace{2.5mm}

Let now $X$ be an $S$ vectorfield and $\xi$ an $S$ 1-form. Then $\sL_X\xi$ is defined on each 
$S_{\ub,u}$ as the usual Lie derivative of $\xi$ as a 1-form on $S_{\ub,u}$ with respect to $X$ as a 
vectorfield on $S_{\ub,u}$, a notion intrinsic to $S_{\ub,u}$, which makes no reference to the 
ambient spacetime. If $\xi$ is any $p$-covariant $S$ tensorfield, $\sL_X\xi$ is similarly defined. 
Consider next the case of another $S$ vectorfield $Y$. Then $\sL_X Y$ is defined to be simply 
${\cal L}_X Y=[X,Y]$, which is itself an $S$ vectorfield. The case of any $q$-contravariant $S$ 
tensorfield $W$ in the role of $Y$ is formally identical to the case of an $S$ vectorfield, such a 
tensorfield being expressible as the sum of tensor products of $S$ vectorfields. Thus $\sL_X W$ is 
defined to be simply ${\cal L}_X W$, which is itself a $q$-contravariant $S$ tensorfield. Consider 
finally the general case of an arbitrary type $T^q_p$ $S$ tensorfield $\theta$. Then $\sL_X\theta$ is 
again defined on each $S_{\ub,u}$ as the usual Lie derivative of $\theta$ as a type $T^q_p$ 
tensorfield on $S_{\ub,u}$ with respect to $X$ as a vectorfield on $S_{\ub,u}$, a notion intrinsic to 
$S_{\ub,u}$, which makes no reference to the ambient spacetime. 

Consider now a given $S_{\ub,u}$ as a section of $C_u$. The {\em 2nd fundamental form} of 
$S_{\ub,u}$ relative to $C_u$ is the 2-covariant tensorfield $\chi$ on $S_{\ub,u}$ defined by setting
at each $p\in S_{\ub,u}$ and for any $X,Y\in T_p S_{\ub,u}$:
\begin{equation}
\chi(X,Y)=h(D_X L,Y)
\label{3.2}
\end{equation}
Similarly, we may consider $S_{\ub,u}$ as a section of $\Cb_{\ub}$. The {\em 2nd fundamental form} of 
$S_{\ub,u}$ relative to $\Cb_{\ub}$ is the 2-covariant tensorfield $\chib$ on $S_{\ub,u}$ defined by 
setting at each $p\in S_{\ub,u}$ and for any $X,Y\in T_p S_{\ub,u}$:
\begin{equation}
\chib(X,Y)=h(D_X\Lb,Y)
\label{3.3}
\end{equation}
The 2nd fundamental forms $\chi$ and $\chib$ are then defined as 2-covariant $S$ tensorfields on 
${\cal N}$. These tensorfields are symmetric. For if $X,Y$ is any pair of $S$ vectorfields, then 
since $h(X,L)=h(Y,L)=0$, we have:
$$\chi(X,Y)-\chi(Y,X)=-h(L,D_X Y)+h(L,D_Y X)=-h(L,[X,Y])=0$$
$[X,Y]$ being itself an $S$ vectorfield. Similarly, since $h(X,\Lb)=h(Y,\Lb)=0$, we have:
$$\chib(X,Y)-\chib(Y,X)=-h(\Lb,D_X Y)+h(\Lb,D_Y X)=-h(\Lb,[X,Y])=0$$

The induced metric $\sh$ on the $S_{\ub,u}$ is a positive-definite symmetric 2-covariant $S$ 
tensorfield. 

\vspace{2.5mm}

\noindent{\bf Proposition 3.1} \ \ \ The {\em first variational formulas:}
$$\sL_L\sh=2\chi, \ \ \ \sL_{\Lb}\sh=2\chib$$
hold. 

\vspace{2.5mm}

\noindent{\em Proof:} Let $X,Y$ be a pair of $S$ vectorfields. Then by Lemma 3.1,
\begin{eqnarray*}
&&(\sL_L\sh)(X,Y)=({\cal L}_L\sh)(X,Y)=L(\sh(X,Y))-\sh([L,X],Y)-\sh(X,[L,Y])\\
&&L(h(X,Y))-h([L,X],Y)-h(X,[L,Y])=h(D_X L,Y)+h(X,D_Y L)\\
&&=2\chi(X,Y)
\end{eqnarray*}
Similarly,
\begin{eqnarray*}
&&(\sL_{\Lb}\sh)(X,Y)-({\cal L}_{\Lb})(X,Y)=\Lb(\sh(X,Y))-\sh([\Lb,X],Y)-\sh(X,[\Lb,Y])\\
&&\Lb(h(X,Y))-h([\Lb,X],Y)-h(X,[\Lb,Y])=h(D_X\Lb,Y)+h(X,D_Y\Lb)\\
&&=2\chib(X,Y)
\end{eqnarray*}

\vspace{2.5mm}

In terms of the frame field $\Omega_A:A=1,...,n-1$ the first variational formulas of the above proposition read (see \ref{2.33}, \ref{2.36}, \ref{2.37}):
\begin{eqnarray}
&&\frac{\partial\sh_{AB}}{\partial\ub}-(\sL_b\sh)_{AB}=2\chi_{AB}\label{3.4}\\
&&\frac{\partial\sh_{AB}}{\partial u}+(\sL_b\sh)_{AB}=2\chib_{AB}\label{3.5}
\end{eqnarray}
where:
\begin{equation}
\chi_{AB}=\chi(\Omega_A,\Omega_B), \ \ \ \chib_{AB}=\chib(\Omega_A,\Omega_B)
\label{3.6}
\end{equation}

Consider a given surface $S_{\ub,u}$. At each point $p\in S_{\ub,u}$, the normal (relative to $h$) 
plane $P_p$ to $T_p S_{\ub,u}$ is a 2-dimensional time-like plane. The pair of future-directed 
null vectors $(L^\prime_p,\Lb_p)$ is a basis for this plane, which is normalized as 
$h(L^\prime_p,\Lb_p)=-2$ (see \ref{2.12}). Thus $(L^\prime,\Lb)$ along $S_{\ub,u}$ 
constitute normalized null basis sections 
for the normal bundle of $S_{\ub,u}$. These are the natural basis sections if we view $S_{\ub,u}$ as 
being a section of $C_u$. The {\em torsion} of 
$S_{\ub,u}$ relative to $C_u$ is the 1-form $\eta$ on $S_{\ub,u}$ defined by setting
at each $p\in S_{\ub,u}$ and for any $X\in T_p S_{\ub,u}$:
\begin{equation}
\eta(X)=\frac{a}{2}h(D_X L^\prime,\Lb)
\label{3.7}
\end{equation}
The pair of future-directed null vectors $(\Lb^\prime_p,L_p)$ at $p$ is another 
basis for $P_p$, which is also normalized as $h(\Lb^\prime_p,L_p)=-2$ (see \ref{2.12}). 
Thus $(\Lb^\prime,L)$ along $S_{\ub,u}$ also constitute normalized null basis sections for the 
normal bundle of $S_{\ub,u}$. These are the natural basis sections if we view $S_{\ub,u}$ as being a 
section of $\Cb_{\ub}$. The {\em torsion} of 
$S_{\ub,u}$ relative to $\Cb_{\ub}$ is the 1-form $\etab$ on $S_{\ub,u}$ defined by setting
at each $p\in S_{\ub,u}$ and for any $X\in T_p S_{\ub,u}$:
\begin{equation}
\etab(X)=\frac{a}{2}h(D_X\Lb^\prime,L)
\label{3.8}
\end{equation}
The torsions $\eta$ and $\etab$ are then defined as $S$ 1-forms on ${\cal N}$. On each $S_{\ub,u}$,  
they represent the connection of the normal bundle of $S_{\ub,u}$. 

By \ref{2.11}, \ref{2.12} we can rewrite \ref{3.7}, \ref{3.8} in the form:
\begin{eqnarray}
&&\eta(X)=\frac{1}{2}h(D_X L,\Lb)+Xa\label{3.9}\\
&&\etab(X)=\frac{1}{2}h(D_X\Lb,L)+Xa\label{3.10}
\end{eqnarray}
Since $h(D_X L,\Lb)+h(D_X\Lb,L)=X(h(L,\Lb))=-2Xa$ by \ref{2.10}, adding these equations we obtain:
$$\eta(X)+\etab(X)=Xa$$
or:
\begin{equation}
\eta+\etab=\sd a
\label{3.11}
\end{equation}

Consider next $D_{\Lb}L$. Since $L$ is null, we have $h(D_{\Lb} L,L)=0$. Also, since 
$h(L,\Lb^\prime)=0$, we have: 
$$h(D_{\Lb}L,\Lb^\prime)=-h(L,D_{\Lb}\Lb^\prime)=0$$
by \ref{2.15}. We conclude that $D_{\Lb}L$ is an $S$ vectorfield. Let then $X$ be an arbitrary 
$S$ vectorfield an consider $h(D_{\Lb}L,X)$. Substituting (see \ref{2.11}) $L=aL^\prime$ we obtain:
\begin{equation}
h(D_{\Lb}L,X)=ah(D_{\Lb}L^\prime,X)
\label{3.12}
\end{equation}
Now, by \ref{2.13}:
$$h\cdot L^\prime=-2du, \ \ \mbox{hence} \ \ h\cdot DL^\prime=-2Ddu$$
and $Ddu$, the Hessian of $u$, is symmetric. Thus, for any pair of vectorfields $V_1,V_2$ on 
${\cal N}$ we have:
\begin{equation}
h(D_{V_1}L^\prime,V_2)=h(D_{V_2}L^\prime,V_1)
\label{3.13}
\end{equation}
Applying this in the case $V_1=\Lb$, $V_2=X$, yields:
$$ah(D_{\Lb}L^\prime,X)=ah(D_X L^\prime, \Lb)=2\eta$$
in view of the definition \ref{3.7}. Substituting in \ref{3.12} we obtain:
$$h(D_{\Lb}L,X)=2\eta(X)$$
We then conclude that:
\begin{equation}
D_{\Lb}L=2\eta^\sharp
\label{3.14}
\end{equation}
where $\eta^\sharp$ is the $S$ vectorfield corresponding to the $S$ 1-form $\eta$ though the induced 
metric $\sh$, that is: $\eta=\sh\cdot\eta^\sharp$. Similarly, we deduce:
\begin{equation}
D_L\Lb=2\etab^\sharp
\label{3.15}
\end{equation}
Equations \ref{3.14} and \ref{3.15} give, in view of the definition \ref{2.17}:
\begin{equation}
Z=2(\eta^\sharp-\etab^\sharp)
\label{3.16}
\end{equation}
Equations \ref{3.11} and \ref{3.16} allow us to express:
\begin{eqnarray}
&&\eta=\frac{1}{2}\sd a+\frac{1}{4}\sh\cdot Z\label{3.17}\\
&&\etab=\frac{1}{2}\sd a-\frac{1}{4}\sh\cdot Z\label{3.18}
\end{eqnarray}
Comparing \ref{3.9}, \ref{3.10} with \ref{3.17}, \ref{3.18}, we see that:
\begin{equation}
h(D_X L,\Lb)=-2\etab(X), \ \ \ h(D_X\Lb,L)=-2\eta(X)
\label{3.19}
\end{equation}
Together with \ref{3.2}, \ref{3.3}, these allow us to express, for any $S$ vectorfield $X$:
\begin{eqnarray}
&&D_X L=\chi^\sharp\cdot X+a^{-1}\etab(X)L\label{3.20}\\
&&D_X\Lb=\chib^\sharp\cdot X+a^{-1}\eta(X)\Lb\label{3.21}
\end{eqnarray}
Here $\chi^\sharp$ and $\chib^\sharp$ are the $T^1_1$-type $S$ tensorfields corresponding through 
$\sh$ to $\chi$ and $\chib$ respectively, that is, for any pair $X$, $Y$ of $S$ vectorfields 
we have:
\begin{equation}
\sh(X,\chi^\sharp\cdot Y)=\chi(X,Y), \ \ \ \sh(X,\chib^\sharp\cdot Y)=\chib(X,Y)
\label{3.22}
\end{equation}

We combine the results \ref{2.7}, \ref{2.9}, \ref{2.16}, and \ref{3.14}, \ref{3.15} into the 
following table:
\begin{equation}
\begin{array}{llll}
&D_L L=a^{-1}(La)L &\s &D_{\Lb}\Lb=a^{-1}(\Lb a)\Lb\\
&D_{\Lb} L=2\eta^\sharp &\s &D_L\Lb=2\etab^\sharp\\
&D_X L=\chi^\sharp\cdot X+a^{-1}\etab(X)L &\s &D_X\Lb=\chib^\sharp\cdot X+a^{-1}\eta(X)\Lb
\end{array}
\label{3.23}
\end{equation}
the $S$ 1-forms $\eta$ and $\etab$ being given by \ref{3.17} and \ref{3.18}. Here $X$ is an 
arbitrary $S$ vectorfield.

Since $D_L X-D_X L=[L,X]$, $D_{\Lb}X-D_X\Lb=[\Lb,X]$, equations 
\ref{3.20}, \ref{3.21} allow us to express $D_L X$ and $D_{\Lb}X$ for any $S$ vectorfield $X$ in 
terms of the commutators $[L,X]$ and $[\Lb,X]$, which by Lemma 3.1 are $S$ vectorfields. 

Consider finally $D_X Y$, $X$, $Y$ being a pair of $S$ vectorfields. Denoting by $\sD$ the covariant 
derivative operator on the surfaces $S_{\ub,u}$ with respect to the induced metric $\sh$, we have:
\begin{equation}
\sD_X Y=\Pi\cdot D_X Y
\label{3.24}
\end{equation}
$\Pi$ being the $h$-orthogonal projection to the surfaces $S_{\ub,u}$. For any vectorfield 
$V$ on ${\cal N}$ we have:
\begin{equation}
\Pi\cdot V=V+\frac{1}{2a}\left(h(V,\Lb)L+h(V,L)\Lb\right)
\label{3.25}
\end{equation}
Setting then $V=D_X Y$, since
$$h(D_X Y,L)=-\chi(X,Y), \ \ \ h(D_X Y,\Lb)=-\chib(X,Y),$$
we conclude that:
\begin{equation}
D_X Y=\sD_X Y+\frac{1}{2a}\left(\chib(X,Y)L+\chi(X,Y)\Lb\right)
\label{3.26}
\end{equation}

If we substitute for $X$ in the table \ref{3.23} and for $X$, $Y$ in \ref{3.26} each one of 
the vectorfields $\Omega_A:A=1,...,n-1$ (see \ref{2.36}) we shall obtain the connection 
coefficients of the frame field $(L,\Lb,\Omega^A:A=1,...,n-1)$ for ${\cal N}$. 

In the following we call {\em conjugation} the formal operation of exchanging $u$ and $L$ with $\ub$ 
and $\Lb$ respectively. We call two entities conjugate if the definition of one is obtained from the definition of the other by conjugation. Thus $\chib$ and $\etab$ are the conjugates of $\chi$ and 
$\eta$ respectively, the function $a$ is self-conjugate and the conjugate of $Z$ is $-Z$. Also, 
the operator $\sL_{\Lb}$ is the conjugate of the operator $\sL_L$. We likewise call two formulas 
{\em conjugate} if one is obtained from the other by conjugation. Thus the first variational formulas of Proposition 3.1 form a conjugate pair, $\sh$ being self-conjugate. Similarly, equations \ref{3.17}, 
\ref{3.18} for a conjugate pair. 

\vspace{5mm}

\section{The Structure Functions and the Formulas for the Torsion Forms}

We now bring in the underlying Minkowskian structure (Galilean structure in the non-relativistic 
framework). So we have an underlying affine structure and preferred systems of linear coordinates 
$(x^\mu:\mu=0,...,n)$. 
If $U$ and $V$ are vectorfields and $\nabla$ the covariant derivative operator defined by this 
affine structure, then if $V^\mu=Vx^\mu:\mu=0,...,n$ are the components of $V$ in such coordinates, so 
that 
$$V=V^\mu\frac{\partial}{\partial x^\mu}$$
then 
\begin{equation}
\nabla_U V=(UV^\mu)\frac{\partial}{\partial x^\mu}
\label{3.27}
\end{equation}
so the components of the vectorfield $\nabla_U V$ in the same coordinates are 
$UV^\mu:\mu=0,...,n$. Now, the vectorfields $N,\Nb,\Omega_A:A=1,...,n-1$ constitute a frame field 
for ${\cal N}$. For any pair of vectorfields $U$, $V$ on ${\cal N}$ we can expand the 
vectorfield $\nabla_U V$ in this frame field. In view of \ref{3.27} this means that there are 
functions $\sf^A$, $f$, and $\of$ such that:
\begin{equation}
UV^\mu=\sf^A\Omega_A^\mu+fN^\mu+\of\Nb^\mu \ \ : \ \mu=0,...,n
\label{3.28}
\end{equation}
This holds in particular in rectangular coordinates (Galilean coordinates in the non-relativistic 
framework). 

Placing then each of the $\Omega_A:A=1,...,n-1$ in the role of $U$ and $N$ in the role of $V$ we have: 
\begin{equation}
\Omega_A N^\mu=\sk^B_A\Omega_B^\mu+k_A N^\mu+\ok_A\Nb^\mu
\label{3.29}
\end{equation}
Similarly, placing $\Nb$ in the role of $V$ we have:
\begin{equation}
\Omega_A\Nb^\mu=\skb^B_A\Omega_B^\mu+\kb_A N^\mu+\okb_A\Nb^\mu
\label{3.30}
\end{equation}
The functions $k_A$, $\ok_A$, $\kb_A$, $\okb_A$ are the components of the $S$ 1-forms 
$k$, $\ok$, $\kb$, $\okb$ respectively, while the functions $\sk^B_A$, $\skb^B_A$ are the 
components of the $T^1_1$-type $S$ tensorfields $\sk$, $\skb$ respectively, in the 
$\Omega_A:A=1,...,n-1$ frame field for the $S_{\ub,u}$. 

Placing $L$ in the role of $U$ and $N$ in the role of $V$ in \ref{3.28} we have:
\begin{equation}
LN^\mu=\sm^A\Omega_A^\mu+mN^\mu+\om\Nb^\mu
\label{3.31}
\end{equation}
Similarly, placing $\Lb$ in the role of $U$ and $\Nb$ in the role of $V$ we have:
\begin{equation}
\Lb\Nb^\mu=\smb^A\Omega_A^\mu+\mb N^\mu+\omb\Nb^\mu
\label{3.32}
\end{equation}
The functions $\sm^A$, $\smb^A$ are the components of the $S$ vectorfields $\sm$, $\smb$ respectively, 
in the $\Omega_A:A=1,...,n-1$ frame field for the $S_{\ub,u}$. 

Placing $\Lb$ in the role of $U$ and $N$ in the role of $V$ in \ref{3.28} we have:
\begin{equation}
\Lb N^\mu=\sn^A\Omega_A^\mu+nN^\mu+\on\Nb^\mu
\label{3.33}
\end{equation}
Similarly, placing $L$ in the role of $U$ and $\Nb$ in the role of $V$ we have:
\begin{equation}
L\Nb^\mu=\snb^A\Omega_A^\mu+\nb N^\mu+\onb\Nb^\mu
\label{3.34}
\end{equation}
The functions $\sn^A$, $\snb^A$ are the components of the $S$ vectorfields $\sn$, $\snb$ respectively, 
in the $\Omega_A:A=1,...,n-1$ frame field for the $S_{\ub,u}$.

The operation of conjugation exchanges $N$ and $\Nb$. The conjugate of $\rho$ is $\rhob$ (see 
\ref{2.60}, \ref{2.61}), the function $c$ (see \ref{2.71} is self-conjugate, and the conjugate 
of $\lambda$ is $\lambdab$ (see \ref{2.74}). Then the conjugate of $\sk$ is $\skb$, while the 
conjugate of $k$ is $\okb$ and the conjugate of $\ok$ is $\kb$. Similarly, the conjugate of 
$\sm$ is $\smb$, while the conjugate of $m$ is $\omb$ and the conjugate of $\om$ is $\mb$. 
Also, the conjugate of $\sn$ is $\snb$, while the conjugate of $n$ is $\onb$ and the conjugate of $\on$ is $\nb$. 

We call the coefficients in the expansions \ref{3.29} - \ref{3.34} {\em structure functions}. 
Assuming that the $x^\mu$ and the $\beta_\mu$ are smooth functions of the acoustical coordinates 
$(\ub,u,\vartheta)$, so are the $\Omega_A^\mu$ and, from Section 2.2, so are the $N^\mu$ and 
$\Nb^\mu$. Assuming also that the vectorfield $b$ is smooth in acoustical coordinates, it then 
follows that the left hand sides of \ref{3.29} - \ref{3.34} are smooth functions of the acoustical 
coordinates. Therefore so are the structure functions. 

\vspace{2.5mm}

Let us define the functions: 
\begin{equation}
\beta_N=\beta_\mu N^\mu, \ \ \ \beta_{\Nb}=\beta_\mu\Nb^\mu, \ \ \ \sbeta_A=\beta_\mu\Omega_A^\mu 
\label{3.35}
\end{equation}
These are the components of the 1-form $\beta$ in the $N,\Nb,\Omega_A:A=1,...,n-1$ frame. 
In particular $\sbeta_A$ are the components in the $\Omega_A:A=1,...,n-1$ frame of the $S$ 1-form 
$\sbeta$ induced by $\beta$ on the $S_{\ub,u}$. The $\beta_N$, $\beta_{\Nb}$, $\sbeta_A$ are all 
smooth functions of the acoustical coordinates. 

We denote by $s$ the 2-covariant tensorfield:
\begin{equation}
s=dx^\mu\otimes d\beta_\mu
\label{3.36}
\end{equation}
By virtue of equation \ref{2.88} of the wave system the tensorfield $S$ is symmetric. If $U$, $V$ 
is any pair of vectorfields on ${\cal N}$ we have:
\begin{equation}
s(U,V)=U^\mu V\beta_\mu=V^\mu U\beta_\mu=s(V,U)
\label{3.37}
\end{equation}
In particular, we denote:
\begin{equation}
(\ss_N)_A=s(N,\Omega_A)=N^\mu\Omega_A\beta_\mu, \ \ \ 
(\ss_{\Nb})_A=s(\Nb,\Omega_A)=\Nb^\mu\Omega_A\beta_\mu
\label{3.38}
\end{equation}
and:
\begin{equation}
\sss_{AB}=s(\Omega_A,\Omega_B)=\Omega_A^\mu\Omega_B\beta_\mu
\label{3.39}
\end{equation}
Then $(\ss_N)_A$, $(\ss_{\Nb})_A$ are the components of the $S$ 1-forms $\ss_N$ and $\ss_{\Nb}$ 
respectively, and $\sss_{AB}$ are the components of the symmetric 2-covariant $S$ tensorfield 
$\sss$ induced by $s$ on the $S_{\ub,u}$, in the $\Omega_A:A=1,...,n-1$ frame field for the 
$S_{\ub,u}$. By virue of our assumptions above, the  $(\ss_N)_A$, $(\ss_{\Nb})_A$, $\sss_{AB}$ are 
all smooth functions of the acoustical coordinates. 

We have (see \ref{2.73}):
\begin{equation}
h^{-1}=-\frac{1}{2c}\left(N\otimes\Nb+\Nb\otimes N\right)+(\sh^{-1})^{AB}\Omega_A\otimes\Omega_B
\label{3.40}
\end{equation}
It follows that, denoting 
\begin{equation}
s_{N\Nb}=s(N,\Nb)
\label{3.41}
\end{equation}
equation \ref{2.86} of the wave system takes the form:
\begin{equation}
s_{N\Nb}=c\mbox{tr}\sss,
\label{3.42}
\end{equation}
the trace being taken relative to $\sh$. By virtue of our assumptions above, the function $c$ 
and the components $\sh_{AB}$ are smooth functions of the acoustical coordinates, therefore so is 
$s_{N\Nb}$. Finally, the components
\begin{equation}
s_{NL}=N^\mu L\beta_\mu, \ \ \ s_{\Nb\Lb}=\Nb^\mu \Lb\beta_\mu
\label{3.43}
\end{equation}
are also smooth functions of the acoustical coordinates. 

Consider now equations \ref{3.29} - \ref{3.34}. Let us denote by $\pi$ the $S$ 1-form
\begin{equation}
\pi=\sd t
\label{3.44}
\end{equation}
Its components in the in the $\Omega_A:A=1,...,n-1$ frame field for the $S_{\ub,u}$ are
\begin{equation}
\pi_A=\Omega_A t=\Omega_A^0
\label{3.45}
\end{equation}
Since $N^0=\Nb^0=1$, the $\mu=0$ components of equations \ref{3.29} - \ref{3.33} read, 
in view of the definition \ref{3.45}, 
\begin{eqnarray}
\sk^B_A\pi_B+k_A+\ok_A=0, \ && \ \skb^B_A\pi_B+\kb_A+\okb_A=0\label{3.46}\\
\sm^A\pi_A+m+\om=0, \ && \ \smb^A\pi_A+\mb+\omb=0\label{3.47}\\
\sn^A\pi_A+n+\on=0, \ && \ \snb^A\pi_A+\nb+\onb=0\label{3.48}
\end{eqnarray} 

Consider $\ok_A$. Since according to \ref{2.71}:
\begin{equation}
h_{\mu\nu}\Nb^\mu N^\nu=-2c
\label{3.49}
\end{equation}
\ref{3.29} implies:
$$\ok_A=-\frac{1}{2c}h_{\mu\nu}(\Omega_A N^\mu)N^\nu$$
Since $h_{\mu\nu}N^\mu N^\nu=0$ this reduces to:
\begin{equation}
\ok_A=\frac{1}{4c}(\Omega_A h_{\mu\nu})N^\mu N^\nu
\label{3.50}
\end{equation}
Recalling \ref{1.d11}, that is:
\begin{equation}
h_{\mu\nu}=g_{\mu\nu}+H\beta_\mu\beta_\nu
\label{3.51}
\end{equation}
we can express in terms of the definitions \ref{3.35}, \ref{3.38}:
\begin{equation}
(\Omega_A h_{\mu\nu})N^\mu N^\nu=\beta_N^2\sd_A H+2H\beta_N(\ss_N)_A
\label{3.52}
\end{equation}
We then conclude that:
\begin{equation}
\ok_A=\frac{1}{4c}\left(\beta_N^2\sd_A H+2H\beta_N(\ss_N)_A\right)
\label{3.53}
\end{equation}
and by conjugation:
\begin{equation}
\kb_A=\frac{1}{4c}\left(\beta_{\Nb}^2\sd_A H+2H\beta_{\Nb}(\ss_{\Nb})_A\right)
\label{3.54}
\end{equation}

Consider $\om$. By \ref{3.49} the expansion \ref{3.31} implies:
$$\om=-\frac{1}{2c}h_{\mu\nu}(LN^\mu)N^\nu$$
which, since $h_{\mu\nu}N^\mu N^\nu=0$, reduces to:
\begin{equation}
\om=\frac{1}{4c}(Lh_{\mu\nu})N^\mu N^\nu
\label{3.55}
\end{equation}
In view of \ref{3.51}, we can express in terms of the definitions \ref{3.35}, \ref{3.43}:
\begin{equation}
(Lh_{\mu\nu})N^\mu N^\nu=\beta_N^2 LH+2H\beta_N s_{NL}
\label{3.56}
\end{equation}
We then conclude that:
\begin{equation}
\om=\frac{1}{4c}\left(\beta_N^2 LH+2H\beta_N s_{NL}\right)
\label{3.57}
\end{equation}
and by conjugation:
\begin{equation}
\mb=\frac{1}{4c}\left(\beta_{\Nb}^2 \Lb H+2H\beta_{\Nb} s_{\Nb\Lb}\right)
\label{3.58}
\end{equation}

Consider $\on$. By \ref{3.49} the expansion \ref{3.33} implies:
$$\on=-\frac{1}{2c}h_{\mu\nu}(\Lb N^\mu)N^\nu$$
Since $h_{\mu\nu}N^\mu N^\nu=0$, this reduces to:
\begin{equation}
\on=\frac{1}{4c}(\Lb h_{\mu\nu})N^\mu N^\nu
\label{3.59}
\end{equation}
In view of \ref{3.51} we can express in terms of the definitions \ref{3.35}, \ref{3.41}:
\begin{equation}
(\Lb h_{\mu\nu})N^\mu N^\nu=\beta_N^2\Lb H+2\rhob H\beta_N s_{N\Nb}
\label{3.60}
\end{equation}
We then conclude that:
\begin{equation}
\on=\frac{1}{4c}\left(\beta_N^2\Lb H+2\rhob H\beta_N s_{N\Nb}\right)
\label{3.61}
\end{equation}
and by conjugation:
\begin{equation}
\nb=\frac{1}{4c}\left(\beta_{\Nb}^2 LH+2\rho H\beta_{\Nb} s_{N\Nb}\right)
\label{3.62}
\end{equation}

Consider next $\sm$. By \ref{3.49} the expansion \ref{3.31} implies:
\begin{equation}
\sh_{AB}\sm^B=h_{\mu\nu}(LN^\mu)\Omega_A^\nu
\label{3.63}
\end{equation}
Since $h_{\mu\nu}N^\mu\Omega_A^\nu=0$ this becomes:
\begin{equation}
\sh_{AB}\sm^B=-(Lh_{\mu\nu})N^\mu\Omega_A^\nu-h_{\mu\nu}N^\mu L\Omega_A^\nu
\label{3.64}
\end{equation}
In regard to the 2nd term on the right, since $[L,\Omega_A]$ is an $S$ vectorfield, we have:
$$h_{\mu\nu}N^\mu[L,\Omega_A]^\nu=0$$
Therefore:
\begin{eqnarray}
&&h_{\mu\nu}N^\mu L\Omega_A^\nu=h_{\mu\nu}N^\mu\Omega_A L^\nu
=h_{\mu\nu}N^\mu\Omega_A(\rho N^\nu)\nonumber\\
&&\hspace{19mm}=\rho h_{\mu\nu}N^\mu\Omega_A N^\nu=-\frac{1}{2}\rho(\Omega_A h_{\mu\nu})N^\mu N^\nu
\label{3.65}
\end{eqnarray}
in view of the fact that $h_{\mu\nu}N^\mu N^\nu=0$. Substituting in \ref{3.64} we obtain:
\begin{equation}
\sh_{AB}\sm^B=-(Lh_{\mu\nu})N^\mu\Omega_A^\nu+\frac{1}{2}\rho(\Omega_A h_{\mu\nu})N^\mu N^\nu
\label{3.66}
\end{equation}
By \ref{3.51} we can express in terms of the definitions \ref{3.35}, \ref{3.38}, \ref{3.43}:
\begin{equation}
(Lh_{\mu\nu})N^\mu\Omega_A^\nu=\beta_N\sbeta_A LH+H(\sbeta_A s_{NL}+\rho\beta_N(\ss_N)_A)
\label{3.67}
\end{equation}
Substituting this together with \ref{3.52} in \ref{3.66} we conclude that:
\begin{equation}
\sh_{AB}\sm^B=-\beta_N\sbeta_A LH+\frac{1}{2}\rho\beta_N^2(\sd H)_A-H\sbeta_A s_{NL}
\label{3.68}
\end{equation}
and by conjugation:
\begin{equation}
\sh_{AB}\smb^B=-\beta_{\Nb}\sbeta_A\Lb H+\frac{1}{2}\rhob\beta_{\Nb}^2(\sd H)_A-H\sbeta_A s_{\Nb\Lb}
\label{3.69}
\end{equation}

Consider now the structure functions $\sk$, $\skb$. We shall relate these to the connection 
coefficients $\chi$ and $\chib$. To find this relation we remark that if $U$ and $V$ are arbitrary 
vectorfields on ${\cal N}$, then $D_U V$ is expressed in arbitrary local coordinates, in particular 
in linear coordinates of the background affine structure by:
\begin{equation}
D_U V=(D_U V)^\mu\frac{\partial}{\partial x^\mu}, \  \mbox{where} \ 
(D_U V)^\mu=UV^\mu+\Gamma^\mu_{\nu\lambda} U^\nu V^\lambda
\label{3.70}
\end{equation}
where $\Gamma^\mu_{\nu\lambda}$ are the connection coefficients associated to the metric $h$ in 
the said coordinates. Hence if $W$ is another arbitrary vectorfield on ${\cal N}$, we have:
\begin{equation}
h(D_U V,W)=h_{\mu\kappa}(UV^\mu)W^\kappa+\Gamma_{\nu\lambda\kappa}U^\nu V^\lambda W^\kappa
\label{3.71}
\end{equation}
where:
\begin{equation}
\Gamma_{\nu\lambda\kappa}=h_{\mu\kappa}\Gamma^\mu_{\nu\lambda}
=\frac{1}{2}
(\partial_\nu h_{\lambda\kappa}+\partial_\lambda h_{\nu\kappa}-\partial_\kappa h_{\nu\lambda})
\label{3.72}
\end{equation}
Substituting for $h$ from \ref{3.51} and recalling the definition \ref{3.36} we find that in 
rectangular coordinates:
\begin{equation}
\Gamma_{\nu\lambda\kappa}=\frac{1}{2}(\beta_\lambda\beta_\kappa\partial_\nu H
+\beta_\nu\beta_\kappa\partial_\lambda H-\beta_\nu\beta_\lambda\partial_\kappa H)
+H\beta_\kappa s_{\nu\lambda}
\label{3.73}
\end{equation}
Therefore, with $U, V, W$ an arbitrary triplet of vectorfields on ${\cal N}$ we have:
\begin{eqnarray}
&&\Gamma_{\nu\lambda\kappa}U^\nu V^\lambda W^\kappa=
\frac{1}{2}\left\{\beta(V)\beta(W) UH+\beta(U)\beta(W) VH-\beta(U)\beta(V) WH\right\}\nonumber\\
&&\hspace{24mm}+H\beta(W)s(U,V)
\label{3.74}
\end{eqnarray}
Let us set $U=\Omega_A$, $V=L$, $W=\Omega_B$. Then the left hand side of \ref{3.71} is, by \ref{3.2}, 
\ref{3.6},
$$h(D_{\Omega_A}L,\Omega_B)=\chi_{AB}$$
On the other hand, since $L^\mu=\rho N^\mu$, by \ref{3.29} the 1st term on the right hand side of 
\ref{3.71} is:
$$h_{\mu\kappa}(\Omega_A L^\mu)\Omega_B^\kappa=\rho\sk^C_A\sh_{CB}$$
and by \ref{3.74} and the definitions \ref{3.35} and \ref{3.38} the 2nd term on the right hand side 
of \ref{3.71} is:
\begin{eqnarray*}
&&\Gamma_{\nu\lambda\kappa}\Omega_A^\nu L^\lambda\Omega_B^\kappa=\frac{1}{2}
\left\{\rho\beta_N(\sbeta_B(\sd H)_A-\sbeta_A(\sd H)_B)+\sbeta_A\sbeta_B LH\right\}\\
&&\hspace{24mm}+\rho H\sbeta_B(\ss)_N)_A
\end{eqnarray*}
Defining then the symmetric 2-covariant $S$ tensorfield $\tchi$ by:
\begin{equation}
\chi=\rho\tchi+\frac{1}{2}(LH)\sbeta\otimes\sbeta 
\label{3.75}
\end{equation}
equation \ref{3.71} takes in the case $U=\Omega_A$, $V=L$, $W=\Omega_B$ the form:
$$\sk^C_A\sh_{CB}=\tchi_{AB}+\frac{1}{2}\beta_N(\sbeta_A(\sd H)_B-\sbeta_B(\sd H)_A)
-H\sbeta_B(\ss_N)_A$$
or:
\begin{equation}
\sk\cdot\sh=\tchi+\frac{1}{2}\beta_N\sbeta\wedge\sd H-H\ss_N\otimes\sbeta
\label{3.76}
\end{equation}
The conjugate equation is:
\begin{equation}
\skb\cdot\sh=\tchib+\frac{1}{2}\beta_{\Nb}\sbeta\wedge\sd H-H\ss_{\Nb}\otimes\sbeta
\label{3.77}
\end{equation}
where $\tchib$ is the symmetric 2-covariant $S$ tensorfield defined by:
\begin{equation}
\chib=\rhob\tchib+\frac{1}{2}(\Lb H)\sbeta\otimes\sbeta 
\label{3.78}
\end{equation}
We remark that since the $\sk^B_A$, $\skb^B_A$ are smooth functions of the acoustical coordinates 
$(\ub,u,\vartheta)$, so are the $\tchi_{AB}$, $\tchib_{AB}$. 

Consider finally the structure functions $\sn^A$, $\snb^A$. 
The remaining structure functions $k_A$, $\okb_A$, $m$, $\omb$, $n$, $\onb$, are determined through 
equations \ref{3.46}, \ref{3.47}, \ref{3.48}. 

We shall first relate $\sn$ and $\snb$ to the connection 
coefficients $\eta$ and $\etab$. To find this relation we set $U=\Lb$, $V=L$, $W=\Omega_A$ in 
\ref{3.71}. Then the left hand side of \ref{3.71} is, by \ref{3.14}, 
$$h(D_{\Lb}L,\Omega_A)=2\eta_A$$
On the other hand, since $L^\mu=\rho N^\mu$, by \ref{3.33} the 1st term on the right hand side of 
\ref{3.71} is:
$$h_{\mu\kappa}(\Lb L^\mu)\Omega_A^\kappa=\rho\sn^B\sh_{AB}$$
and by \ref{3.74} and the definitions \ref{3.35}, \ref{3.38}, \ref{3.42}  the 2nd term on the right hand side 
of \ref{3.71} is:
\begin{eqnarray*}
&&\Gamma_{\nu\lambda\kappa}\Lb^\nu L^\lambda\Omega_A^\kappa=
\frac{1}{2}\left\{\rho\beta_N\sbeta_A\Lb H+\rhob\beta_{\Nb}\sbeta_A LH-\rho\rhob\beta_N\beta_{\Nb}
(\sd H)_A\right\}\\
&&\hspace{24mm}+\rho\rhob H\sbeta_A s_{N\Nb}
\end{eqnarray*}
Defining then the $S$ 1-form $\teta$ by:
\begin{equation}
\eta=\rho\teta+\frac{1}{4}\rhob\beta_{\Nb}(LH)\sbeta
\label{3.79}
\end{equation}
equation \ref{3.71} takes in the case $U=\Lb$, $V=L$, $W=\Omega_A$ the form:
\begin{equation}
\sn^B\sh_{AB}=2\teta_A-\frac{1}{2}\beta_N\sbeta_A\Lb H+\frac{1}{2}\rhob(\beta_N\beta_{\Nb}(\sd H)_A
-2H\sbeta_A s_{N\Nb})
\label{3.80}
\end{equation}
The conjugate equation is:
\begin{equation}
\snb^B\sh_{AB}=2\tetab_A-\frac{1}{2}\beta_{\Nb}\sbeta_A LH+\frac{1}{2}\rho(\beta_N\beta_{\Nb}(\sd H)_A
-2H\sbeta_A s_{N\Nb})
\label{3.81}
\end{equation}
where $\tetab$ is the $S$ 1-form defined by:
\begin{equation}
\etab=\rhob\tetab+\frac{1}{4}\rho\beta_N(\Lb H)\sbeta
\label{3.82}
\end{equation}
We remark that since the $\sn^A$, $\snb^A$ are smooth functions of the acoustical coordinates 
$(\ub,u,\vartheta)$, so are the $\teta_A$, $\tetab_A$. 

We shall now relate $\sn$ and $\snb$ to $\sd\rhob$ and $\sd\rho$. By \ref{3.33} and the fact 
that $h_{\mu\nu}N^\mu\Omega_A^\nu=0$ we have:
\begin{eqnarray}
&&\sn^B\sh_{AB}=h_{\mu\nu}(\Lb N^\mu)\Omega_A^\nu \label{3.a1}\\
&&\hspace{12mm}=-(\Lb h_{\mu\nu})N^\mu\Omega_A^\nu-h_{\mu\nu}N^\mu(\Lb\Omega_A^\nu)\nonumber
\end{eqnarray}
By Lemma 3.1 $[\Lb,\Omega_A]$ is an $S$ vectorfield. It follows that 
\begin{eqnarray}
&&h_{\mu\nu}N^\mu(\Lb\Omega_A^\nu)=h_{\mu\nu}N^\mu(\Omega_A\Lb^\nu)=
h_{\mu\nu}N^\mu\Omega_A(\rhob\Nb^\nu)\nonumber\\
&&\hspace{22mm}=-2c\Omega_A\rhob+\rhob h_{\mu\nu}N^\mu(\Omega\Nb^\nu)\nonumber\\
&&\hspace{22mm}=-2c\left((\sd\rhob)_A+\okb_A\rhob\right)\label{3.a2}
\end{eqnarray}
Substituting this together with 
\begin{equation}
(\Lb h_{\mu\nu})N^\mu\Omega_A^\nu=\beta_N\sbeta_A\Lb H
+\rhob H(\beta_N(\ss_{\Nb})_A+\sbeta_A s_{N\Nb})
\label{3.a3}
\end{equation}
in \ref{3.a1} we obtain:
\begin{equation}
\sn^B\sh_{AB}=2c\left((\sd\rhob)_A+\okb_A\rhob\right)-\beta_N\sbeta_A\Lb H
-\rhob H(\beta_N(\ss_{\Nb})_A+\sbeta_A s_{N\Nb})
\label{3.a4}
\end{equation}
The conjugate equation is:
\begin{equation}
\snb^B\sh_{AB}=2c\left((\sd\rho)_A+k_A\rho\right)-\beta_{\Nb}\sbeta_A LH
-\rho H(\beta_{\Nb}(\ss_N)_A+\sbeta_A s_{N\Nb})
\label{3.a5}
\end{equation}

\vspace{5mm}

Consider now the positive function $c$ defined by \ref{2.71}:
\begin{equation}
c=-\frac{1}{2}h_{\mu\nu}N^\mu\Nb^\nu
\label{3.83}
\end{equation}
In view of \ref{3.29} - \ref{3.34} the 1st derivatives of $c$ are determined through the structure 
functions. We have:
\begin{eqnarray*}
&&\Omega_A c=-\frac{1}{2}(\Omega_A h_{\mu\nu})N^\mu\Nb^\nu
-\frac{1}{2}h_{\mu\nu}(\Omega_A N^\mu)\Nb^\nu-\frac{1}{2}h_{\mu\nu}N^\mu(\Omega_A\Nb^\nu)\\
&&\hspace{8mm}=-\frac{1}{2}(\Omega_A h_{\mu\nu})N^\mu \Nb^\nu+c(k_A+\okb_A)
\end{eqnarray*}
which gives:
\begin{equation}
(\sd c)_A=-\frac{1}{2}\left(\beta_N\beta_{\Nb}(\sd H)_A+H(\beta_N(\ss_{\Nb})_A
+\beta_{\Nb}(\ss_N)_A)\right)+c(k_A+\okb_A)
\label{3.84}
\end{equation}
We have:
\begin{eqnarray*}
&&Lc=-\frac{1}{2}(Lh_{\mu\nu})N^\mu\Nb^\nu-\frac{1}{2}h_{\mu\nu}(LN^\mu)\Nb^\nu
-\frac{1}{2}h_{\mu\nu}N^\mu(L\Nb^\nu)\\
&&\hspace{6mm}=-\frac{1}{2}(Lh_{\mu\nu})N^\mu\Nb^\nu+c(m+\onb)
\end{eqnarray*}
which gives:
\begin{equation}
Lc=-\frac{1}{2}\left\{\beta_N\beta_{\Nb}LH+H(\rho\beta_N s_{N\Nb}+\beta_{\Nb}s_{NL})\right\}
+c(m+\onb)
\label{3.85}
\end{equation}
The conjugate equation is:
\begin{equation}
\Lb c=-\frac{1}{2}\left\{\beta_N\beta_{\Nb}\Lb H
+H(\rhob\beta_{\Nb}s_{N\Nb}+\beta_N s_{\Nb\Lb})\right\}+c(\omb+n)
\label{3.86}
\end{equation}

\vspace{5mm}

Using \ref{3.84} we can rewrite equations \ref{3.a4} and \ref{3.a5} in terms of $\sd\lambda$ and 
$\sd\lambdab$ in the form:
\begin{equation}
\sn^B\sh_{AB}=2(\sd\lambda)_A-2k_A\lambda+\rhob\beta_N\beta_{\Nb}(\sd H)_A-\beta_N\sbeta_A\Lb H 
+\rhob H(\beta_{\Nb}\ss_N-\sbeta s_{N\Nb})
\label{3.a6}
\end{equation}
\begin{equation}
\snb^B\sh_{AB}=2(\sd\lambdab)_A-2\okb_A\lambdab+\rho\beta_N\beta_{\Nb}(\sd H)_A-\beta_{\Nb}\sbeta_A LH
+\rho H(\beta_N\ss_{\Nb}-\sbeta s_{N\Nb})
\label{3.a7}
\end{equation}
Comparing finally \ref{3.a6} and \ref{3.a7} with \ref{3.80} and \ref{3.81} respectively, we 
arrive at the following proposition.

\vspace{2.5mm}

\noindent{\bf Proposition 3.2} \ \ \ The torsion forms $\teta$, $\tetab$ are given in terms of 
$\sd\lambda$, $\sd\lambdab$ by:
\begin{eqnarray*}
&&\teta=\sd\lambda
+\left(\frac{1}{4c}\left(\beta_N\beta_{\Nb}\sd H+2H\beta_{\Nb}\ss_N\right)-k\right)\lambda
-\frac{1}{4}\beta_N\sbeta\Lb H\\
&&\tetab=\sd\lambdab
+\left(\frac{1}{4c}\left(\beta_N\beta_{\Nb}\sd H+2H\beta_N \ss_{\Nb}\right)-\okb\right)\lambdab
-\frac{1}{4}\beta_{\Nb}\sbeta LH
\end{eqnarray*}
with $k$ given in terms of $\ok$, $\sk$ and $\okb$ given in terms of $\kb$, $\skb$ by \ref{3.46}, 
while $\ok$, $\kb$ are given by \ref{3.53}, \ref{3.54} and $\sk$, $\skb$ are given in terms of 
the 2nd fundamental forms $\tchi$, $\tchib$ by \ref{3.76}, \ref{3.77}.

\vspace{2.5mm}

Recalling the remarks at the end of Section 2.2, we call {\em acoustical quantities} the quantities 
of a given order which represent properties of the foliation of ${\cal N}$ by  the level sets of the acoustical functions $u$ and $\ub$, of that order. 
At order 0 the independent acoustical quantities are the $\sd x^\mu$, which, in view of the 
fact that: 
\begin{equation}
\sh=h_{\mu\nu}\sd x^\mu\otimes\sd x^\nu
\label{3.120}
\end{equation}
determine $\sh$, as well as, according to Section 2.2, $N^\mu$ and $\Nb^\mu$, therefore also $c$, 
the $S$ vectorfield $b$, and the quantities $\rho$ and $\rhob$, or equivalently the quantities 
$\lambdab$ and $\lambda$. At order 1, by virtue of the above proposition, we may take $\tchi$, 
$\tchib$ and the 1st derivatives $\sd\lambda$, $\Lb\lambda$ and $\sd\lambdab$, $L\lambdab$ of $\lambda$ and $\lambdab$ respectively to be 
the independent acoustical quantities. As we shall see in Proposition 3.3, the 1st order quantities $L\lambda$ and $\Lb\lambdab$ are not acoustical, as they depend on acoustical quantities of 0th order only. 

\vspace{2.5mm}

In regard to the structure functions as defined by \ref{3.29} - \ref{3.34} we remark that what is missing is 
$$\Omega_A \Omega_B^\mu=\frac{\partial^2 x^\mu}{\partial\vartheta^A\partial\vartheta^B}$$
However this does not have geometrical meaning. What does have geometrical meaning, and shall 
presently be determined, are the components 
\begin{equation}
(\sD^2 x^\mu)_{AB}=(\sD^2 x^\mu)(\Omega_A,\Omega_B)
\label{3.87}
\end{equation}
of the covariant Hessian of the restriction of the rectangular coordinate functions $x^\mu$ 
to the surfaces $S_{\ub,u}$ with respect to the induced metric $\sh$. 

The corresponding components 
of the covariant Hessian of the $x^\mu$ in ${\cal N}$ are:
\begin{equation}
(D^2 x^\mu)_{AB}=(D^2 x^\mu)(\Omega_A,\Omega_B)
\label{3.88}
\end{equation}
We have:
\begin{equation}
(\sD^2 x^\mu)_{AB}=\Omega_A(\Omega_B x^\mu)-(\sD_{\Omega_A}\Omega_B)x^\mu
\label{3.89}
\end{equation}
and:
\begin{equation}
(D^2 x^\mu)_{AB}=\Omega_A(\Omega_B x^\mu)-(D_{\Omega_A}\Omega_B) x^\mu
\label{3.90}
\end{equation}
Subtracting \ref{3.90} from \ref{3.89} gives:
\begin{equation}
(\sD^2 x^\mu)_{AB}-(D^2 x^\mu)_{AB}=(D_{\Omega_A}\Omega_B-\sD_{\Omega_A}\Omega_B) x^\mu
\label{3.91}
\end{equation}
On the left hand side, $D^2 x^\mu$ is readily obtained in rectangular coordinates:
\begin{equation}
(D^2 x^\mu)_{\nu\lambda}=\frac{\partial^2 x^\mu}{\partial x^\nu\partial x^\lambda}
-\Gamma^\kappa_{\nu\lambda}\frac{\partial x^\mu}{\partial x^\kappa}=-\Gamma^\mu_{\nu\lambda}
\label{3.92}
\end{equation}
Hence, setting
\begin{equation}
\gamma^\mu_{AB}=\Gamma^\mu_{\nu\lambda}\Omega_A^\nu\Omega_B^\lambda,
\label{3.93}
\end{equation}
we have:
\begin{equation}
(D^2 x^\mu)_{AB}=-\gamma^\mu_{AB}
\label{3.94}
\end{equation}
We can expand the vectorfield 
$$\gamma_{AB}=\gamma^\mu_{AB}\frac{\partial}{\partial x^\mu}$$ 
in terms of the frame field $\Omega_A:A=1,...,n-1$, $N$, $\Nb$:
\begin{equation}
\gamma^\mu_{AB}=\sgamma^C_{AB}\Omega_C^\mu+\gamma^N_{AB}N^\mu+\gamma^{\Nb}_{AB}\Nb^\mu
\label{3.95}
\end{equation}
The coefficients of the expansion are determined by:
\begin{equation}
h(\gamma_{AB},\Omega_C)=\sgamma_{AB}^D\sh_{CD}, \ \  
h(\gamma_{AB},\Nb)=-2c\gamma^N_{AB}, \ \  
h(\gamma_{AB},N)=-2c\gamma^{\Nb}_{AB}
\label{3.96}
\end{equation}
In view of \ref{3.93} we have:
\begin{eqnarray}
&&h(\gamma_{AB},\Omega_C)=\Gamma_{\nu\lambda\kappa}\Omega_A^\nu\Omega_B^\lambda\Omega_C^\kappa
\nonumber\\
&&h(\gamma_{AB},\Lb)=\Gamma_{\nu\lambda\kappa}\Omega_A^\nu\Omega_B^\lambda\Lb^\kappa\nonumber\\
&&h(\gamma_{AB}, L)=\Gamma_{\nu\lambda\kappa}\Omega_A^\nu\Omega_B^\lambda L^\kappa
\label{3.97}
\end{eqnarray}
Through the formula \ref{3.74} we then obtain:
\begin{equation}
\sgamma^D_{AB}\sh_{CD}=\frac{1}{2}\left(\sbeta_B\sbeta_C(\sd H)_A+\sbeta_A\sbeta_C(\sd H)_B
-\sbeta_A\sbeta_B(\sd H)_C\right)+H\sbeta_C(\sss)_{AB}
\label{3.98}
\end{equation}
Moreover, setting:
\begin{eqnarray}
&&\gamma^N_{AB}=\tgamma^N_{AB}+\frac{1}{4c\rhob}\sbeta_A\sbeta_B \Lb H\nonumber\\ 
&&\gamma^{\Nb}_{AB}=\tgamma^{\Nb}_{AB}+\frac{1}{4c\rho}\sbeta_A\sbeta_B LH
\label{3.99}
\end{eqnarray}
we find:
\begin{eqnarray}
&&\tgamma^N_{AB}=-\frac{1}{4c}\beta_{\Nb}\left(\sbeta_B(\sd H)_A+\sbeta_A(\sd H)_B
+2H(\sss)_{AB}\right)\nonumber\\
&&\tgamma^{\Nb}_{AB}=-\frac{1}{4c}\beta_N\left(\sbeta_B(\sd H)_A+\sbeta_A(\sd H)_B
+2H(\sss)_{AB}\right)
\label{3.100}
\end{eqnarray}
We remark that the $\tgamma^N_{AB}$, $\tgamma^{\Nb}_{AB}$ are, like the $\sgamma^C_{AB}$, 
smooth functions of the acoustical coordinates $(\ub,u,\vartheta)$ even where the product 
$\rho\rhob$ vanishes. Defining then $\tgamma^\mu_{AB}$ by:
\begin{equation}
\tgamma^\mu_{AB}=\sgamma^C_{AB}\Omega_C^\mu +\tgamma^N_{AB} N^\mu +\tgamma^{\Nb}_{AB}\Nb^\mu
\label{3.101}
\end{equation}
we have:
\begin{equation}
\gamma^\mu_{AB}=\tgamma^\mu_{AB}+\frac{1}{4c}
\left(\frac{(\Lb H)}{\rhob}N^\mu +\frac{(LH)}{\rho}\Nb^\mu\right)\sbeta_A\sbeta_B
\label{3.102}
\end{equation}
and the functions $\tgamma^\mu_{AB}$ are smooth functions of the acoustical coordinates 
$(\ub,u,\vartheta)$ even where the product $\rho\rhob$ vanishes.

Going back to \ref{3.91}, by \ref{3.26} the right hand side is given by (see \ref{2.72}):
\begin{equation}
(D_{\Omega_A}\Omega_B-\sD_{\Omega_A}\Omega_B) x^\mu=
\frac{1}{2c}\left(\frac{\chib_{AB}}{\rhob}N^\mu+\frac{\chi_{AB}}{\rho}\Nb^\mu\right)
\label{3.103}
\end{equation}
Substituting also \ref{3.94}, equation \ref{3.91} gives:
\begin{equation}
(\sD^2 x^\mu)_{AB}=\frac{1}{2c}\left(\frac{\chib_{AB}}{\rhob}N^\mu+\frac{\chi_{AB}}{\rho}\Nb^\mu\right)
-\gamma^\mu_{AB}
\label{3.104}
\end{equation}
Substituting \ref{3.75}, \ref{3.78}, and \ref{3.102}, the terms which are singular where the 
product $\rho\rhob$ vanishes cancel, and we obtain:
\begin{equation}
(\sD^2 x^\mu)_{AB}=\frac{1}{2c}(\tchib_{AB}N^\mu+\tchi_{AB}\Nb^\mu)-\tgamma^\mu_{AB}
\label{3.105}
\end{equation}

\vspace{5mm}

\section{The Propagation Equations for $\lambda$ and $\lambdab$}

Let us now return to the characteristic system \ref{2.62}. In view of \ref{2.33}, \ref{2.36}, \ref{2.42}, 
\ref{2.61}, this system is:
\begin{equation}
Lx^\mu=\rho N^\mu, \ \ \ \Lb x^\mu=\rhob \Nb^\mu \ \ \ : \ \mu=0,...,n
\label{3.106}
\end{equation}
the $\mu=0$ components being the definitions \ref{2.61}. Applying $\Lb$ to the first and $L$ to the 
second and subtracting, we obtain, since by \ref{2.17} and the fact that $Z$ is an $S$ vectorfield 
we have
\begin{equation}
[\Lb,L]x^\mu=Zx^\mu=Z^A\Omega_A^\mu,
\label{3.107}
\end{equation}
the equation:
\begin{equation}
Z^A\Omega_A^\mu = \Lb(\rho N^\mu)-L(\rhob \Nb^\mu)
\label{3.108}
\end{equation}
Substituting for $\Lb N^\mu$, $L\Nb^\mu$ from \ref{3.33}, \ref{3.34} then yields the equation:
\begin{eqnarray}
&&(\rho\sn^A-\rhob\snb^A-Z^A)\Omega_A^\mu\label{3.109}\\
&&+(\Lb\rho+n\rho-\nb\rhob)N^\mu-(L\rhob+\onb\rhob-\on\rho)\Nb^\mu=0\nonumber
\end{eqnarray}
This is equivalent to the three the coefficients of $\Omega_A^\mu$, $N^\mu$, $\Nb^\mu$  equal to zero. 
The first equation is 
\begin{equation}
Z^A=\rho\sn^A-\rhob\snb^A
\label{3.110}
\end{equation}
substituting for $\sn$, $\snb$ in terms of $\teta$, $\tetab$ from \ref{3.80}, \ref{3.81} 
we recover the formula \ref{3.16}. So the first equation contains no new information. The second and 
third equations read:
\begin{eqnarray}
&&\Lb\rho+n\rho-\nb\rhob=0\label{3.111}\\
&&L\rhob+\onb\rhob-\on\rho=0\label{3.112}
\end{eqnarray}
Multiplying by $c$ these can be written in terms of $\lambda$, $\lambdab$ in the form:
\begin{eqnarray}
&&L\lambda-\rhob Lc+\onb\lambda-\on\lambdab=0\label{3.113}\\
&&\Lb\lambdab-\rho\Lb c+n\lambdab-\nb\lambda=0\label{3.114}
\end{eqnarray}
Substituting then for $Lc$, $\Lb c$ from \ref{3.85}, \ref{3.86} and for $\on$, $\nb$ from 
\ref{3.61}, \ref{3.62} we arrive at the following proposition.

\vspace{2.5mm}

\noindent{\bf Proposition 3.3} \ \ \ The inverse density functions $\lambda$, $\lambdab$ satisfy 
the propagation equations:
$$L\lambda=p\lambda+q\lambdab, \ \ \ \Lb\lambdab=\pb\lambdab+\qb\lambda$$
where 
\begin{eqnarray*}
&&p=m-\frac{1}{2c}\left(\beta_N\beta_{\Nb} LH+H\beta_{\Nb}s_{NL}\right), \ \ \ 
q=\frac{1}{4c}\beta_N^2\Lb H\\
&&\pb=\omb-\frac{1}{2c}\left(\beta_N\beta_{\Nb}\Lb H+H\beta_N s_{\Nb\Lb}\right), \ \ \ 
\qb=\frac{1}{4c}\beta_{\Nb}^2 LH
\end{eqnarray*}
with $m$ given in terms of $\om$, $\sm$ and $\omb$ given in terms of $\mb$, $\smb$ by \ref{3.47}, 
while $\om$, $\mb$ are given by \ref{3.57}, \ref{3.58} and $\sm$, $\smb$ by \ref{3.68}, \ref{3.69}. 

\vspace{2.5mm}

We remark that the functions $p$, $q$, $\pb$, $\qb$ are of order 1 but depend on acoustical quantities 
of only the 0th order. 

\vspace{5mm}

\section{The Second Variation and Cross Variation Equations}

By \ref{3.29} we have:
$$\sk^C_A\sh_{CB}=h_{\mu\nu}(\Omega_A N^\mu)\Omega_B^\nu$$
This expresses in components in the $\Omega_A:A=1,...,n-1$ frame field for the $S_{\ub,u}$ the 
equation:
\begin{equation}
\sk\cdot\sh=h_{\mu\nu}\sd N^\mu\otimes\sd x^\nu
\label{3.115}
\end{equation}
Applying $\sL_L$ to both sides of this equation we obtain, in view of Lemma 3.2, 
\begin{equation}
\sL_L(\sk\cdot\sh)=h_{\mu\nu}\sd(LN^\mu)\otimes\sd x^\nu+h_{\mu\nu}\sd N^\mu\otimes\sd L^\nu
+(Lh_{\mu\nu})\sd N^\mu\otimes\sd x^\nu
\label{3.116}
\end{equation}
Consider the 1st term on the right. According to \ref{3.31} we have:
\begin{equation}
LN^\mu=\sm\cdot\sd x^\mu+mN^\mu+\om\Nb^\mu
\label{3.117}
\end{equation}
Hence:
\begin{eqnarray}
&&\sd(LN^\mu)=\sD\sm\cdot\sd x^\mu+\sm\cdot\sD^2 x^\mu+m\sd N^\mu+\om\sd\Nb^\mu\nonumber\\
&&\hspace{13mm}+N^\mu\sd m+\Nb^\mu\sd\om\label{3.118}
\end{eqnarray}
In view of \ref{3.115}, the conjugate equation
\begin{equation}
\skb\cdot\sh=h_{\mu\nu}\sd\Nb^\mu\otimes\sd x^\nu
\label{3.119}
\end{equation}
\ref{3.120}, and the fact that by \ref{3.104}, \ref{3.102} we have
\begin{equation}
h_{\mu\nu}\sD^2 x^\mu\otimes\sd x^\nu=-\sgamma\cdot\sh
\label{3.121}
\end{equation}
where $\sgamma$ is the type $T^1_2$ $S$ tensorfield with components $\sgamma^C_{AB}$ (symmetric in 
the lower indices) given by \ref{3.98},
we then obtain, in regard to the 1st term on the right in \ref{3.116}, 
\begin{equation}
h_{\mu\nu}\sd(LN^\mu)\otimes\sd x^\nu=(\sD\sm-\sm\cdot\sgamma+m\sk+\om\skb)\cdot\sh
\label{3.122}
\end{equation}
Consider next the 2nd term on the right in \ref{3.116}. This term is
\begin{equation}
h_{\mu\nu}\sd N^\mu\otimes\sd(\rho N^\nu)=h_{\mu\nu}N^\nu\sd N^\mu\otimes\sd\rho
+\rho h_{\mu\nu}\sd N^\mu\otimes\sd N^\nu
\label{3.123}
\end{equation}
According to the expansion \ref{3.29}:
\begin{equation}
\sd N^\mu=\sk^A\Omega_A^\mu+k N^\mu+\ok\Nb^\mu
\label{3.124}
\end{equation}
where we consider the components $\sk^A:A=1,...,n-1$ of the type $T^1_1$ $S$ tensorfield $\sk$ 
as $S$ 1-forms. Using this representation we obtain, in regard to the two terms on the right in 
\ref{3.123}, 
\begin{equation}
h_{\mu\nu}N^\nu\sd N^\mu=-2c\ok
\label{3.125}
\end{equation}
\begin{equation}
h_{\mu\nu}\sd N^\mu\otimes\sd N^\nu=\sk\times\sk-2c(k\otimes\ok+\ok\otimes k)
\label{3.126}
\end{equation}
where for any pair of $T^1_1$ type $S$ tensorfields $\theta$, $\xi$ we denote:
\begin{equation}
\theta\times\xi=\sh_{AB}\theta^A\otimes\xi^B
\label{3.127}
\end{equation}
a 2-covariant $S$ tensorfield. Finally, in regard to the 3rd term on the right in \ref{3.116},  
substituting 
$$Lh_{\mu\nu}=\beta_\mu\beta_\nu LH+H(\beta_\nu L\beta_\mu+\beta_\mu L\beta_\nu)$$
and the expansion \ref{3.124}, we find that this term is:
\begin{eqnarray}
&&(Lh_{\mu\nu})\sd N^\mu\otimes\sd x^\nu=\left(\sk\cdot\sbeta
+\beta_N k+\beta_{\Nb}\ok\right)\otimes\left((LH)\sbeta+\rho H\ss_N\right)\nonumber\\
&&\hspace{28mm}+H\left(\rho\sk\cdot\ss_N+s_{NL}k+\rho s_{N\Nb}\ok\right)\otimes\sbeta\label{3.128}
\end{eqnarray}
Combining the above results we arrive at the following proposition. 

\vspace{2.5mm}

\noindent{\bf Proposition 3.4} \ \ \ The 2-covariant $S$ tensorfield $\sk\cdot\sh$ satisfies 
the following propagation equation:
\begin{eqnarray*}
&&\sL_L(\sk\cdot\sh)=(\sD\sm-\sm\cdot\sgamma+m\sk+\om\skb)\cdot\sh\\
&&\hspace{15mm}-2c\ok\otimes\sd\rho+\rho\left(\sk\times\sk-2c(k\otimes\ok+\ok\otimes k)\right)\\
&&\hspace{15mm}+\left(\sk\cdot\sbeta
+\beta_N k+\beta_{\Nb}\ok\right)\otimes\left((LH)\sbeta+\rho H\ss_N\right)\\
&&\hspace{15mm}+H\left(\rho\sk\cdot\ss_N+s_{NL}k+\rho s_{N\Nb}\ok\right)\otimes\sbeta
\end{eqnarray*}
Similarly, the 2-covariant $S$ tensorfield $\skb\cdot\sh$ satisfies the conjugate propagation 
equation:
\begin{eqnarray*}
&&\sL_{\Lb}(\skb\cdot\sh)=(\sD\smb-\smb\cdot\sgamma+\omb\skb+\mb\sk)\cdot\sh\\
&&\hspace{15mm}-2c\kb\otimes\sd\rhob+\rhob\left(\skb\times\skb-2c(\okb\otimes\kb+\kb\otimes \okb)\right)\\
&&\hspace{15mm}+\left(\skb\cdot\sbeta
+\beta_{\Nb}\okb+\beta_N\kb\right)\otimes\left((\Lb H)\sbeta+\rhob H\ss_{\Nb}\right)\\
&&\hspace{15mm}+H\left(\rhob\skb\cdot\ss_{\Nb}+s_{\Nb\Lb}\okb+\rhob s_{N\Nb}\kb\right)\otimes\sbeta
\end{eqnarray*}
In view of equations \ref{3.76}, \ref{3.77} the above constititute the 
{\em second variational formulas}, which govern the evolution of $\tchi$ along the generators of 
the $C_u$ and the evolution of $\tchib$ along the generators of the $\Cb_{\ub}$. 

\vspace{2.5mm}

The principal terms on the right hand sides of the above second variation equations are the order 2 
terms $\sD(\sm\cdot \sh)$ and $\sD(\smb\cdot\sh)$ respectively. We see from \ref{3.68}, \ref{3.69} 
that the 1st order quantities $\sm\cdot\sh$ and $\smb\cdot\sh$ depend on acoustical quantities of the 
0th order only. Therefore there are no principal acoustical terms on the right hand sides of the 
second variation equations. These right hand sides depend on acoustical quantities of the 1st order 
only. 

\vspace{2.5mm}

Applying $\sL_{\Lb}$ to both sides of equation \ref{3.115} we obtain, in view of Lemma 3.2, 
\begin{equation}
\sL_{\Lb}(\sk\cdot\sh)=h_{\mu\nu}\sd(\Lb N^\mu)\otimes\sd x^\nu+h_{\mu\nu}\sd N^\mu\otimes\sd\Lb^\nu
+(\Lb h_{\mu\nu})\sd N^\mu\otimes\sd x^\nu
\label{3.129}
\end{equation}
Consider the 1st term on the right. According to \ref{3.33} we have:
\begin{equation}
\Lb N^\mu=\sn\cdot\sd x^\mu+nN^\mu+\on\Nb^\mu
\label{3.130}
\end{equation}
Hence:
\begin{eqnarray}
&&\sd(LN^\mu)=\sD\sn\cdot\sd x^\mu+\sn\cdot\sD^2 x^\mu+n\sd N^\mu+\on\sd\Nb^\mu\nonumber\\
&&\hspace{13mm}+N^\mu\sd n+\Nb^\mu\sd\on\label{3.131}
\end{eqnarray}
In view of \ref{3.115}, \ref{3.119}, \ref{3.120}, and \ref{3.121}, 
we then obtain, in regard to the 1st term on the right in \ref{3.129}, 
\begin{equation}
h_{\mu\nu}\sd(\Lb N^\mu)\otimes\sd x^\nu=(\sD\sn-\sn\cdot\sgamma+n\sk+\on\skb)\cdot\sh
\label{3.132}
\end{equation}
Consider next the 2nd term on the right in \ref{3.129}. This term is
\begin{equation}
h_{\mu\nu}\sd N^\mu\otimes\sd(\rhob\Nb^\nu)=h_{\mu\nu}\Nb^\nu\sd N^\mu\otimes\sd\rhob
+\rhob h_{\mu\nu}\sd N^\mu\otimes\sd\Nb^\nu
\label{3.133}
\end{equation}
According to the expansion \ref{3.30}:
\begin{equation}
\sd\Nb^\mu=\skb^A\Omega_A^\mu+\kb N^\mu+\okb\Nb^\mu
\label{3.134}
\end{equation}
where we consider the components $\skb^A:A=1,...,n-1$ of the type $T^1_1$ $S$ tensorfield $\skb$ 
as $S$ 1-forms. Using the representations \ref{3.124} and \ref{3.134} we obtain, in regard to the two 
terms on the right in \ref{3.133}, 
\begin{equation}
h_{\mu\nu}\Nb^\nu\sd N^\mu=-2ck
\label{3.135}
\end{equation}
\begin{equation}
h_{\mu\nu}\sd N^\mu\otimes\sd\Nb^\nu=\sk\times\skb-2c(k\otimes\okb+\ok\otimes\kb)
\label{3.136}
\end{equation}
Finally, in regard to the 3rd term on the right in \ref{3.129},  
substituting 
$$\Lb h_{\mu\nu}=\beta_\mu\beta_\nu\Lb H+H(\beta_\nu\Lb\beta_\mu+\beta_\mu\Lb\beta_\nu)$$
and the expansion \ref{3.124}, we find that this term is:
\begin{eqnarray}
&&(\Lb h_{\mu\nu})\sd N^\mu\otimes\sd x^\nu=\left(\sk\cdot\sbeta
+\beta_N k+\beta_{\Nb}\ok\right)\otimes\left((\Lb H)\sbeta+\rhob H\ss_{\Nb}\right)\nonumber\\
&&\hspace{28mm}+H\left(\rhob\sk\cdot\ss_{\Nb}+\rhob s_{N\Nb}k+s_{\Nb\Lb}\ok\right)\otimes\sbeta\label{3.137}
\end{eqnarray}
Combining the above results we arrive at the following proposition. 

\vspace{2.5mm}

\noindent{\bf Proposition 3.5} \ \ \ The 2-covariant $S$ tensorfield $\sk\cdot\sh$ satisfies 
the following equation:
\begin{eqnarray*}
&&\sL_{\Lb}(\sk\cdot\sh)=(\sD\sn-\sn\cdot\sgamma+n\sk+\on\skb)\cdot\sh\\
&&\hspace{15mm}-2ck\otimes\sd\rhob+\rhob\left(\sk\times\skb-2c(k\otimes\okb+\ok\otimes\kb)\right)\\
&&\hspace{15mm}+\left(\sk\cdot\sbeta
+\beta_N k+\beta_{\Nb}\ok\right)\otimes\left((\Lb H)\sbeta+\rhob H\ss_{\Nb}\right)\\
&&\hspace{15mm}+H\left(\rhob\sk\cdot\ss_{\Nb}+\rhob s_{N\Nb}k+s_{\Nb\Lb}\ok\right)\otimes\sbeta
\end{eqnarray*}
Similarly, the 2-covariant $S$ tensorfield $\skb\cdot\sh$ satisfies the conjugate equation:
\begin{eqnarray*}
&&\sL_L(\skb\cdot\sh)=(\sD\snb-\snb\cdot\sgamma+\onb\skb+\nb\sk)\cdot\sh\\
&&\hspace{15mm}-2c\okb\otimes\sd\rho+\rho\left(\skb\times\sk-2c(\okb\otimes k+\kb\otimes\ok)\right)\\
&&\hspace{15mm}+\left(\skb\cdot\sbeta
+\beta_{\Nb}\okb+\beta_N\kb\right)\otimes\left((LH)\sbeta+\rho H\ss_N\right)\\
&&\hspace{15mm}+H\left(\rho\skb\cdot\ss_N+\rho s_{N\Nb}\okb+s_{NL}\kb\right)\otimes\sbeta
\end{eqnarray*}
In view of equations \ref{3.76}, \ref{3.77} the above equations  
govern the evolution of $\tchi$ along the generators of 
the $\Cb_{\ub}$ and the evolution of $\tchib$ along the generators of the $\Cb_{\ub}$. For this 
reason we call them {\em cross variation equations}. 

\vspace{2.5mm}

The principal terms on the right hand sides of the above cross variation equations are the 2nd order 
terms $\sD(\sn\cdot \sh)$ and $\sD(\snb\cdot\sh)$ respectively. We see from \ref{3.a6}, \ref{3.a7} 
that the 1st order quantities $\sn\cdot\sh$ and $\snb\cdot\sh$ do contain a 1st order acoustical 
terms. Denoting by $[\s]_{P.A.}$ the principal acoustical part, namely for a quantity of order $l$ 
the part containing acoustical quantities of order $l$, we have:
$$[\sn\cdot\sh]_{P.A.}=2[\sd\lambda-k\lambda]_{P.A.}, \ \ \ 
[\snb\cdot\sh]_{P.A.}=2[\sd\lambdab-\okb\lambdab]_{P.A.}$$
and from \ref{3.46}, together with \ref{3.53}, \ref{3.54} and \ref{3.76}, \ref{3.77},
$$[k]_{P.A.}=-[\sk\cdot\pi]_{P.A.}=-\tchi^\sharp\cdot\pi, \ \ \ 
[\okb]_{P.A.}=-[\skb\cdot\pi]_{P.A.}=-\tchib^\sharp\cdot\pi$$
Here for any 2-covariant $S$ tensorfield $\theta$ we denote by $\theta^\sharp$ the $T^1_1$ type 
$S$ tensorfield such that for any pair of $S$ vectorfields $X$, $Y$ we have:
\begin{equation}
h(\theta^\sharp\cdot X,Y)=\theta(X,Y)
\label{3.138}
\end{equation}
We then have:
\begin{equation}
[\sn\cdot\sh]_{P.A.}=2(\sd\lambda+\lambda\tchi^\sharp\cdot\pi), \ \ \ 
[\snb\cdot\sh]_{P.A.}=2(\sd\lambdab+\lambdab\tchib^\sharp\cdot\pi)
\label{3.139}
\end{equation}
It follows that:
\begin{equation}
[\sD(\sn\cdot\sh]_{P.A.}=2\left(\sD^2\lambda+\lambda(\sD\tchi^\sharp)\cdot\pi\right), \ \ \ 
[\sD(\snb\cdot\sh]_{P.A.}=2\left(\sD^2\lambdab+\lambdab(\sD\tchib^\sharp)\cdot\pi\right)
\label{3.140}
\end{equation}
These are then the principal acoustical parts of the right hand sides of the cross 
variation equations.

\vspace{5mm}

\section{The Case $n=2$}

In the case $n=2$ the $S_{\ub,u}$ which in the general case we call surfaces are in this case curves. 
Since the index $A$ takes only the value 1, we drop it altogether and write for example \ref{2.94}, 
which gives the basic vectorfields, as: 
\begin{equation}
L=\frac{\partial}{\partial\ub}-b\frac{\partial}{\partial\vartheta}, \ \ \ 
\Lb=\frac{\partial}{\partial u}+b\frac{\partial}{\partial\vartheta}, \ \ \ 
\Omega=\frac{\partial}{\partial\vartheta}
\label{3.a8}
\end{equation}
Equation \ref{2.96} becomes:
\begin{equation}
\sh=h_{\mu\nu}(\Omega x^\mu)(\Omega x^\nu)
\label{3.a9}
\end{equation}
and gives the squared magnitude of the vectorfield $\Omega$. There is a vectorfield 
colinear and in the same sense as $\Omega$, which is of unit magnitude, namely:
\begin{equation}
E=\frac{1}{\sqrt{\sh}}\Omega 
\label{3.a10}
\end{equation}
the unit tangent field of the oriented curves $S_{\ub,u}$. This is a canonical frame field for 
the $S_{\ub,u}$, and any $S$ tensorfield reduces to a scalar, the single component of the tensorfield 
with respect to this frame field. Setting $X=E$ in \ref{3.23}, we obtain: 
\begin{equation}
\begin{array}{llll}
&D_L L=a^{-1}(La)L &\s &D_{\Lb}\Lb=a^{-1}(\Lb a)\Lb\\
&D_{\Lb} L=2\eta E &\s &D_L\Lb=2\etab E\\
&D_E L=\chi E+a^{-1}\etab L &\s &D_E\Lb=\chib E+a^{-1}\eta\Lb
\end{array}
\label{3.a11}
\end{equation}
Also, since for any function $f$ on ${\cal N}$ the $S$ 1-form $\sd f$ is replaced by the function
$Ef$, equation \ref{3.11} becomes:
\begin{equation}
\eta+\etab=Ea
\label{3.a12}
\end{equation}
Moreover, since: 
$$h(D_L E,L)=-h(E,D_L L)=0, \ \ \ h(D_{\Lb} E,\Lb)=-h(E,D_{\Lb},\Lb)=0,$$
$$h(D_L E,\Lb)=-h(E,D_L\Lb)=-2\etab, \ \ \ h(D_{\Lb}E,L)=-h(E,D_{\Lb}L)=-2\eta,$$
$$h(D_E E,L)=-h(E,D_E L)=-\chi, \ \ \ h(D_E E,\Lb)=-h(E,D_E\Lb)=-\chib$$
and:
$$h(D_L E,E)=h(D_{\Lb}E,E)=h(D_E E,E)=0$$
we have:
\begin{eqnarray}
&&D_L E=a^{-1}\etab L \ \ \ \ \ D_{\Lb}E=a^{-1}\eta\Lb \nonumber\\
&&\hspace{7mm} D_E E=(2a)^{-1}(\chib L+\chi\Lb)
\label{3.a13}
\end{eqnarray}
The table \ref{3.a11} together with \ref{3.a13} give the connection 
coefficients of $({\cal N},h)$ relative to the frame field $(L,\Lb,E)$. They imply the commutation 
relations:
\begin{eqnarray}
&&[L,E]=-\chi E \ \ \ \ \ [\Lb,E]=-\chib E \nonumber\\
&&\hspace{7mm}[\Lb,L]=Z=2(\eta-\etab)E
\label{3.a14}
\end{eqnarray}

In the case $n=2$, \ref{3.29} - \ref{3.34} take the form:
\begin{eqnarray}
&&EN^\mu=\sk E^\mu+kN^\mu+\ok\Nb^\mu\nonumber\\
&&E\Nb^\mu=\skb E^\mu+\kb N^\mu+\okb\Nb^\mu\nonumber\\
&&LN^\mu=\sm E^\mu+mN^\mu+\om\Nb^\mu\nonumber\\
&&\Lb\Nb^\mu=\smb E^\mu+\mb N^\mu+\omb\Nb^\mu\nonumber\\
&&\Lb N^\mu=\sn E^\mu+n N^\mu+\on\Nb^\mu\nonumber\\
&&L\Nb^\mu=\snb E^\mu+\nb N^\mu+\onb\Nb^\mu 
\label{3.a15}
\end{eqnarray}
all the coefficients being now scalars. 

In the case $n=2$, $\sbeta$ (see \ref{3.35}) is a scalar:
\begin{equation}
\sbeta=\beta_\mu E^\mu 
\label{3.a16}
\end{equation}
Also, $\ss_N$, $\ss_{\Nb}$, and $\sss$ (see \ref{3.38}, \ref{3.39}) are scalars:
\begin{equation}
\ss_N=N^\mu E\beta_\mu, \ \ \ \ss_{\Nb}=\Nb^\mu E\beta_\mu, \ \ \ \sss=E^\mu E\beta_\mu 
\label{3.a17}
\end{equation}
and equation \ref{3.42} takes the form:
\begin{equation}
s_{N\Nb}=c\sss
\label{3.a18}
\end{equation}
Also, the 1-form $\pi$ (see \ref{3.44}, \ref{3.45}) is a scalar:
\begin{equation}
\pi=Et=E^0
\label{3.a19}
\end{equation}

In the case $n=2$ equations \ref{3.75}, \ref{3.78} take the form:
\begin{eqnarray}
&&\chi=\rho\tchi+\frac{1}{2}\sbeta^2 LH\nonumber\\
&&\chib=\rhob\tchib+\frac{1}{2}\sbeta^2\Lb H
\label{3.a20}
\end{eqnarray}
Also, the relations \ref{3.76}, \ref{3.77} reduce to:
\begin{eqnarray}
&&\sk=\tchi-H\sbeta\ss_N\nonumber\\
&&\skb=\tchib-H\sbeta\ss_{\Nb}
\label{3.a21}
\end{eqnarray}
The relations \ref{3.79}, \ref{3.82} remain unaltered in form:
\begin{eqnarray}
&&\eta=\rho\teta+\frac{1}{4}\rhob\beta_{\Nb}\sbeta LH\nonumber\\
&&\etab=\rhob\tetab+\frac{1}{4}\rho\beta_N\sbeta\Lb H
\label{3.a22}
\end{eqnarray}
while equation \ref{3.84} takes the form:
\begin{equation}
Ec=-\frac{1}{2}\left(\beta_N\beta_{\Nb}EH+H(\beta_N\ss_{\Nb}+\beta_{\Nb}\ss_N)\right)+c(k+\okb)
\label{3.a23}
\end{equation}
The relations \ref{3.87}, \ref{3.88} reduce to:
\begin{eqnarray}
&&\sn=2E\lambda-2k\lambda+\rhob\beta_N\beta_{\Nb}EH-\beta_N\sbeta\Lb H\nonumber\\
&&\snb=2E\lambdab-2\okb\lambdab+\rho\beta_N\beta_{\Nb}EH-\beta_{\Nb}\sbeta LH
\label{3.a24}
\end{eqnarray}
Finally, the equations of proposition 3.2 take in the case $n=2$ the form:
\begin{eqnarray}
&&\teta=E\lambda+\left(\frac{1}{4c}(\beta_N\beta_{\Nb}EH+2H\beta_{\Nb}\ss_N)-k\right)\lambda
-\frac{1}{4}\beta_N\sbeta\Lb H\nonumber\\
&&\tetab=E\lambdab+\left(\frac{1}{4c}(\beta_N\beta_{\Nb}EH+2H\beta_N\ss_{\Nb})-\okb\right)\lambdab
-\frac{1}{4}\beta_{\Nb}\sbeta LH
\label{3.a25}
\end{eqnarray}

Since $h(D_E E,E)=0$,by \ref{3.24} we have:
\begin{equation}
\sD_E E=0
\label{3.a26}
\end{equation}
It follows that for any function $f$ on ${\cal N}$, $\sD^2 f$ reduces to 
\begin{equation}
(\sD^2 f)=E(Ef)
\label{3.a27}
\end{equation}
In particular this holds for the rectangular coordinate functions $x^\mu$. We then obtain from 
\ref{3.105}:
\begin{equation}
EE^\mu=\frac{1}{2c}(\tchib N^\mu+\tchi\Nb^\mu)-\tgamma^\mu 
\label{3.a28}
\end{equation}
where by \ref{3.98}, \ref{3.100}, \ref{3.101}:
\begin{equation}
\tgamma^\mu=\sgamma E^\mu+\tgamma^N N^\mu+\tgamma^{\Nb}\Nb^\mu 
\label{3.a29}
\end{equation}
with:
\begin{eqnarray}
&&\sgamma=\frac{1}{2}\sbeta^2 EH+H\sbeta\sss\nonumber\\
&&\tgamma^N=-\frac{1}{2c}\beta_{\Nb}(\sbeta EH+H\sss)\nonumber\\
&&\tgamma^{\Nb}=-\frac{1}{2c}\beta_N(\sbeta EH+H\sss)
\label{3.a30}
\end{eqnarray}
We can write \ref{3.a28} in the form:
\begin{equation}
EE^\mu=\se E^\mu+eN^\mu+\oe\Nb^\mu 
\label{3.a31}
\end{equation}
where:
\begin{eqnarray}
&&\se=-\frac{1}{2}\sbeta^2 EH-H\sbeta\sss\nonumber\\
&&e=\frac{1}{2c}\left(\tchib+\beta_{\Nb}(\sbeta EH+H\sss)\right)\nonumber\\
&&\oe=\frac{1}{2c}\left(\tchi+\beta_N(\sbeta EH+H\sss)\right)
\label{3.a32}
\end{eqnarray}

By the 1st of \ref{3.a14}:
$$LE^\mu-EL^\mu=-\chi E^\mu$$
Substituting $L^\mu=\rho N^\mu$ and expressing $EN^\mu$ by the 1st of \ref{3.a15} we then obtain:
\begin{equation}
LE^\mu=\sf E^\mu+fN^\mu+\of\Nb^\mu 
\label{3.a33}
\end{equation}
where:
\begin{equation}
\sf=\rho\sk-\chi, \ \ \ f=E\rho+\rho k, \ \ \ \of=\rho\ok
\label{3.a34}
\end{equation}
Substituting in the first of these the 1st of \ref{3.a20}, \ref{3.a21} yields:
\begin{equation}
\sf=-\frac{1}{2}\sbeta^2 LH-\rho H\sbeta\ss_N
\label{3.a35}
\end{equation}
In the second of \ref{3.a34} we write $E\rho=E(c^{-1}\lambdab)$ and use \ref{3.a23} to express 
$Ec$. This gives:
$$f=c^{-1}(E\lambdab-\okb\lambdab)+(2c)^{-1}\rho\left(\beta_N\beta_{\Nb}EH+H(\beta_N\ss_{\Nb}
+\beta_{\Nb}\ss_N)\right)$$
By the 2nd of \ref{3.46}, substituting also \ref{3.a21} and \ref{3.54}, the first term on the right here is:
$$\frac{1}{c}\left(E\lambdab+\pi\lambdab(\tchib-H\sbeta\ss_{\Nb})\right)+
\frac{\rho}{4c}\left(\beta_{\Nb}^2 EH+2H\beta_{\Nb}\ss_{\Nb}\right) $$
We thus obtain:
\begin{eqnarray}
&&f=\frac{1}{c}\left(E\lambdab+\pi\lambdab\tchib\right)\label{3.a36}\\
&&\hspace{5mm}+\frac{\rho}{2c}\left\{\frac{1}{2}\beta_{\Nb}\left(\beta_{\Nb}+2\beta_N\right)EH
+H\left((\beta_N+\beta_{\Nb}-2c\pi\sbeta)\ss_{\Nb}+\beta_{\Nb}\ss_N\right)\right\}\nonumber
\end{eqnarray}
Also, in view of \ref{3.53} the third of \ref{3.a34} gives:
\begin{equation}
\of=\frac{\rho}{4c}\left(\beta_N^2 EH+2H\beta_N\ss_N\right)
\label{3.a37}
\end{equation}
The conjugate of equation \ref{3.a33} is:
\begin{equation}
\Lb E^\mu=\sfb E^\mu+\fb N^\mu+\ofb\Nb^\mu 
\label{3.a38}
\end{equation}
with the coefficients given by the conjugates of \ref{3.a35}, \ref{3.a36}, \ref{3.a37}, that is:
\begin{equation}
\sfb==-\frac{1}{2}\sbeta^2\Lb H-\rhob H\sbeta\ss_{\Nb}
\label{3.a39}
\end{equation}
\begin{equation}
\fb=\frac{\rhob}{4c}\left(\beta_{\Nb}^2 EH+2H\beta_{\Nb}\ss_{\Nb}\right)
\label{3.a40}
\end{equation}
\begin{eqnarray}
&&\ofb=\frac{1}{c}\left(E\lambda+\pi\lambda\tchi\right)\label{3.a41}\\
&&\hspace{5mm}+\frac{\rhob}{2c}\left\{\frac{1}{2}\beta_N\left(\beta_N+2\beta_{\Nb}\right)EH
+H\left((\beta_{\Nb}+\beta_N-2c\pi\sbeta)\ss_N+\beta_N\ss_{\Nb}\right)\right\}\nonumber
\end{eqnarray}

We turn to the 2nd variation equations of Proposition 3.4 in the case $n=2$. Let $\xi$ be an arbitrary 
2-covariant $S$ tensorfield. In the case $n=2$ this is represented by the scalar $\xi(E,E)$. 
To avoid confusion let us denote the last by $\sxi$. We then have:
\begin{eqnarray}
&&(\sL_L\xi)(E,E)=L\sxi-\xi([L,E],E)-\xi(E,[L,E])\nonumber\\
&&\hspace{20mm}=L\sxi+2\chi\sxi
\label{3.a42}
\end{eqnarray}
by the first of the commutation relations \ref{3.a14}. Similarly, we have:
\begin{eqnarray}
&&(\sL_{\Lb}\xi)(E,E)=\Lb\sxi-\xi([\Lb,E],E)-\xi(E,[\Lb,E])\nonumber\\
&&\hspace{20mm}=\Lb\sxi+2\chib\sxi
\label{3.a43}
\end{eqnarray}
by the second of the commutation relations \ref{3.a14}. Consequently, the 2nd variation equations 
of Proposition 3.4 reduce in the case $n=2$ to: 
\begin{eqnarray}
&&L\sk+2\chi\sk=E\sm-\sm\sgamma+m\sk+\om\skb-2c\ok E\rho\nonumber\\
&&\hspace{18mm}+\rho(\sk^2-4ck\ok)+(\sbeta\sk+\beta_N k+\beta_{\Nb}\ok)(\sbeta LH+\rho H\ss_N)\nonumber\\
&&\hspace{18mm}+H\sbeta(\rho\ss_N\sk+s_{NL}k+\rho s_{N\Nb}\ok)
\label{3.a44}
\end{eqnarray}
\begin{eqnarray}
&&\Lb\skb+2\chib\skb=E\smb-\smb\sgamma+\omb\skb+\mb\sk-2c\kb E\rhob\nonumber\\
&&\hspace{18mm}+\rhob(\skb^2-4c\okb\kb)+(\sbeta\skb+\beta_{\Nb}\okb+\beta_N\kb)(\sbeta\Lb H+\rhob H\ss_{\Nb})\nonumber\\
&&\hspace{18mm}+H\sbeta(\rhob\ss_{\Nb}\skb+s_{\Nb\Lb}\okb+\rhob s_{N\Nb}\kb)
\label{3.a45}
\end{eqnarray}

Finally, the cross variation equations of Proposition 3.5 reduce in the case $n=2$ to:
\begin{eqnarray}
&&\Lb\sk+2\chib\sk=E\sn-\sn\sgamma+n\sk+\on\skb-2ck E\rhob\nonumber\\
&&\hspace{18mm}+\rhob(\sk\skb-2c(k\okb+\ok\kb))
+(\sbeta\sk+\beta_N k+\beta_{\Nb}\ok)(\sbeta\Lb H+\rhob H\ss_{\Nb})\nonumber\\
&&\hspace{18mm}+H\sbeta(\rhob\ss_{\Nb}\sk+\rhob s_{N\Nb}k+s_{\Nb\Lb}\ok)
\label{3.a46}
\end{eqnarray}
\begin{eqnarray}
&&L\skb+2\chi\skb=E\snb-\snb\sgamma+\onb\skb+\nb\sk-2c\okb E\rho\nonumber\\
&&\hspace{18mm}+\rho(\skb\sk-2c(\okb k+\kb\ok))
+(\sbeta\skb+\beta_{\Nb}\okb+\beta_N\kb)(\sbeta LH+\rho H\ss_N)\nonumber\\
&&\hspace{18mm}+H\sbeta(\rho\ss_N\skb+\rho s_{N\Nb}\okb+s_{NL}\kb)
\label{3.a47}
\end{eqnarray}

\vspace{5mm}

\section{The Codazzi and Gauss Equations ($n>2$)}

In the case of $n=2$ spatial dimensions the acoustical structure equations which we have developed 
so far form a complete system, the $S_{\ub,u}$ being in this case diffeomorphic to $S^1$. However for 
$n>2$ spatial dimensions the above must be supplemented by the Codazzi and Gauss equations of the 
embedding of the $S_{\ub,u}$, which are diffeomorphic to $S^{n-1}$, in the spacetime $({\cal N},h)$. 

We shall first derive the Codazzi equations. Consider the 2-covariant 
$S$ tensorfield $\sk\cdot\sh$, that is, in components:
\begin{equation}
(\sk\cdot\sh)_{AB}=\sk_A^C\sh_{CB}:=\sk_{AB}
\label{3.141}
\end{equation}
From \ref{3.29} we have:
\begin{equation}
\sk_{AB}=h_{\mu\nu}(\Omega_A N^\mu)\Omega_B^\nu
\label{3.142}
\end{equation}
Then:
\begin{equation}
(\sD(\sk\cdot\sh))_{ABC}=\Omega_A(\sk_{BC})-(\sk\cdot\sh)(\sD_{\Omega_A}\Omega_B,\Omega_C)
-(\sk\cdot\sh)(\Omega_B,\sD_{\Omega_A}\Omega_C)
\label{3.143}
\end{equation}
Now, we have:
\begin{eqnarray}
&&\Omega_A(\sk_{BC})=h_{\mu\nu}(\Omega_A\Omega_B N^\mu)\Omega_C^\nu
+h_{\mu\nu}(\Omega_B N^\mu)(\Omega_A\Omega_C^\nu)\nonumber\\
&&\hspace{15mm}+(\Omega_A h_{\mu\nu})(\Omega_B N^\mu)\Omega_C^\nu
\label{3.144}
\end{eqnarray}
Substituting this in \ref{3.143}, interchanging $A$ and $B$, and subtracting, noting that 
$[\Omega_A,\Omega_B]=0$, we obtain:
\begin{eqnarray}
&&(\sD(\sk\cdot\sh))_{ABC}-(\sD(\sk\cdot\sh))_{BAC}\label{3.145}\\
&&\hspace{15mm}=h_{\mu\nu}(\Omega_B N^\mu)(\Omega_A\Omega_C^\nu)
-h_{\mu\nu}(\Omega_A N^\mu)(\Omega_B\Omega_C^\nu)\nonumber\\
&&\hspace{15mm}+\left\{(\Omega_A h_{\mu\nu})(\Omega_B N^\mu)
-(\Omega_B h_{\mu\nu})(\Omega_A N^\mu)\right\}\Omega_C^\nu\nonumber\\
&&\hspace{15mm}-\sk_{BD}(\sD_{\Omega_A}\Omega_C)^D+\sk_{AD}(\sD_{\Omega_B}\Omega_C)^D\nonumber
\end{eqnarray}
where we denote by $(\sD_{\Omega_A}\Omega_B)^C$ the components of the $S$ vectorfield 
$\sD_{\Omega_A}\Omega_B$ in the $\Omega_1,...,\Omega_{n-1}$ frame:
\begin{equation}
\sD_{\Omega_A}\Omega_B=(\sD_{\Omega_A}\Omega_B)^C\Omega_C
\label{3.146}
\end{equation}
that is, the connection coefficients associated to the induced metric $\sh$ and the said frame field. 
In the first two terms on the right in \ref{3.145} we substitute the expansion \ref{3.29} to deduce 
that the contribution of these terms is:
\begin{eqnarray}
&&\sk^D_B h_{\mu\nu}\Omega_D^\mu(\Omega_A\Omega_C^\nu)
-\sk^D_A h_{\mu\nu}\Omega_D^\mu(\Omega_B\Omega_C^\nu)\nonumber\\
&&+k_B h_{\mu\nu}N^\mu(\Omega_A\Omega_C^\nu)-k_A h_{\mu\nu} N^\mu(\Omega_B\Omega_C^\nu)\nonumber\\
&&+\ok_B h_{\mu\nu}\Nb^\mu(\Omega_A\Omega_C^\nu)-\ok_A h_{\mu\nu}\Nb^\mu(\Omega_B\Omega_C^\nu)
\label{3.147}
\end{eqnarray}
Now, by \ref{3.89}, \ref{3.105}, \ref{3.101}, and \ref{3.100}, defining 
\begin{eqnarray}
&&\psi_{AB}=\tchi_{AB}-2c\tgamma^{\Nb}_{AB}
=\tchi_{AB}+\frac{1}{2}\beta_N\left(\sbeta_A(\sd H)_B+\sbeta_B(\sd H)_A+2H\sss_{AB}\right)\nonumber\\
&&\psib_{AB}=\tchib_{AB}-2c\tgamma^N_{AB}
=\tchib_{AB}+\frac{1}{2}\beta_{\Nb}\left(\sbeta_A(\sd H)_B+\sbeta_B(\sd H)_A+2H\sss_{AB}\right)\nonumber\\
&&
\label{3.148}
\end{eqnarray}
and noting the fact that 
$$\nabla_{\Omega_A}\Omega^B=(\Omega_A\Omega_B^\mu)\frac{\partial}{\partial x^\mu}$$
we obtain:
\begin{equation}
\nabla_{\Omega_A}\Omega_B=\sD_{\Omega_A}\Omega_B-\sgamma_{AB}^C\Omega_C
+\frac{1}{2c}(\psib_{AB}N+\psi_{AB}\Nb)
\label{3.149}
\end{equation}
Substituting this in \ref{3.147}, the contribution of the 1st term in \ref{3.149} is
$$\sk_{BD}(\sD_{\Omega_A}\Omega_C)^D-\sk_{AD}(\sD_{\Omega_B}\Omega_C)^D$$
which cancels the last two terms in \ref{3.145}, the contribution of the 2nd term in \ref{3.149} is
$$-\sk_{BD}\sgamma^D_{AC}+\sk_{AD}\sgamma^D_{BC}$$
and the contribution of the last term in \ref{3.149} is:
$$-\ok_B\psib_{AC}+\ok_A\psib_{BC}-k_B\psi_{AC}+k_A\psi_{BC}$$
What then remains to be considered in \ref{3.145} is the term:
$$\left\{(\Omega_A h_{\mu\nu})(\Omega_B N^\mu)
-(\Omega_B h_{\mu\nu})(\Omega_A N^\mu)\right\}\Omega_C^\nu$$
Substituting the expansion \ref{3.29}, we express this term as in the form:
\begin{eqnarray}
&&(\Omega_A h_{\mu\nu})\sk^D_B\Omega_D^\mu\Omega_C^\nu
-(\Omega_B h_{\mu\nu})\sk^D_A\Omega_D^\mu\Omega_C^\nu\nonumber\\
&&+(\Omega_A h_{\mu\nu})k_B N^\mu\Omega_C^\nu-(\Omega_B h_{\mu\nu})k_A N^\mu\Omega_C^\nu\nonumber\\
&&+(\Omega_A h_{\mu\nu})\ok_B\Nb^\mu\Omega_C^\nu
-(\Omega_B h_{\mu\nu})\ok_A\Nb^\mu\Omega_C^\nu\label{3.150}\\
&&=\left[\sbeta_C\sbeta_D(\sd H)_A+H(\sbeta_C\sss_{AD}+\sbeta_D\sss_{AC})\right]\sk^D_B\nonumber\\
&&\hspace{18mm}-\left[\sbeta_C\sbeta_D(\sd H)_B+H(\sbeta_C\sss_{BD}+\sbeta_D\sss_{BC})\right]\sk^D_A\nonumber\\
&&+\left[\beta_N\sbeta_C(\sd H)_A+H(\sbeta_C(\ss_N)_A+\beta_N\sss_{AC})\right]k_B\nonumber\\
&&\hspace{18mm}-\left[\beta_N\sbeta_C(\sd H)_B+H(\sbeta_C(\ss_N)_B+\beta_N\sss_{BC})\right]k_A\nonumber\\
&&+\left[\beta_{\Nb}\sbeta_C(\sd H)_A+H(\sbeta_C(\ss_{\Nb})_A+\beta_{\Nb}\sss_{AC})\right]\ok_B\nonumber\\
&&\hspace{18mm}-\left[\beta_{\Nb}\sbeta_C(\sd H)_B+H(\sbeta_C(\ss_{\Nb})_B+\beta_{\Nb}\sss_{BC})\right]\ok_A\nonumber
\end{eqnarray}
We thus arrive at the following proposition.

\vspace{2.5mm}

\noindent{\bf Proposition 3.6} \ \ \ The 2nd fundamental form $\tchi$ satisfies through 
\ref{3.76} the 
Codazzi equation: 
\begin{eqnarray*}
&&(\sD(\sk\cdot\sh))_{ABC}-(\sD(\sk\cdot\sh))_{BAC}=\\
&&\hspace{5mm}\left[\frac{1}{2}\left(\sbeta_A\sbeta_C(\sd H)_D+\sbeta_D\sbeta_C(\sd H)_A
-\sbeta_A\sbeta_D(\sd H)_C\right)+H\sbeta_C\sss_{AD}\right]\sk^D_B\\
&&\hspace{5mm}-\left[\frac{1}{2}\left(\sbeta_B\sbeta_C(\sd H)_D+\sbeta_D\sbeta_C(\sd H)_B
-\sbeta_B\sbeta_D(\sd H)_C\right)+H\sbeta_C\sss_{BD}\right]\sk^D_A\\
&&\hspace{5mm}+k_A\left((\sk\cdot\sh)_{BC}-\frac{1}{2}H\beta_N\sss_{BC}\right)
-k_B\left((\sk\cdot\sh)_{AC}-\frac{1}{2}H\beta_N\sss_{AC}\right)\\
&&\hspace{5mm}+\ok_A\left((\skb\cdot\sh)_{BC}-\frac{1}{2}H\beta_{\Nb}\sss_{BC}\right)
-\ok_B\left((\skb\cdot\sh)_{AC}-\frac{1}{2}H\beta_{\Nb}\sss_{AC}\right)
\end{eqnarray*}
Similarly, the conjugate 2nd fundamental form $\tchib$ satisfies through \ref{3.77} the conjugate 
Codazzi equation:
\begin{eqnarray*}
&&(\sD(\skb\cdot\sh))_{ABC}-(\sD(\skb\cdot\sh))_{BAC}=\\
&&\hspace{5mm}\left[\frac{1}{2}\left(\sbeta_A\sbeta_C(\sd H)_D+\sbeta_D\sbeta_C(\sd H)_A
-\sbeta_A\sbeta_D(\sd H)_C\right)+H\sbeta_C\sss_{AD}\right]\skb^D_B\\
&&\hspace{5mm}-\left[\frac{1}{2}\left(\sbeta_B\sbeta_C(\sd H)_D+\sbeta_D\sbeta_C(\sd H)_B
-\sbeta_B\sbeta_D(\sd H)_C\right)+H\sbeta_C\sss_{BD}\right]\skb^D_A\\
&&\hspace{5mm}+\okb_A\left((\skb\cdot\sh)_{BC}-\frac{1}{2}H\beta_{\Nb}\sss_{BC}\right)
-\okb_B\left((\skb\cdot\sh)_{AC}-\frac{1}{2}H\beta_{\Nb}\sss_{AC}\right)\\
&&\hspace{5mm}+\kb_A\left((\sk\cdot\sh)_{BC}-\frac{1}{2}H\beta_N\sss_{BC}\right)
-\kb_B\left((\sk\cdot\sh)_{AC}-\frac{1}{2}H\beta_N\sss_{AC}\right)
\end{eqnarray*}

\vspace{2.5mm}

Note that there are no principal terms on the right hand sides of the Codazzi equations as given by the above proposition. 

\vspace{5mm}

We proceed to the derivation of the Gauss equation. Here we must bring in the curvature $R$ of the 
acoustical spacetime $({\cal N},h)$. In [Ch-S] it was shown that this {\em acoustical curvature} 
decomposes into:
\begin{equation}
R_{\mu\nu\kappa\lambda}=P_{\mu\nu\kappa\lambda}+N_{\mu\nu\kappa\lambda}
\label{3.151}
\end{equation}
where $P_{\mu\nu\kappa\lambda}$ is the principal part, given by:
\begin{equation}
P_{\mu\nu\kappa\lambda}=\frac{1}{2}\left((D^2 H)_{\kappa\nu}\beta_\lambda\beta_\mu
+(D^2 H)_{\lambda\mu}\beta_\kappa\beta_\nu-(D^2 H)_{\kappa\mu}\beta_\lambda\beta_\nu
-(D^2 H)_{\lambda\nu}\beta_\kappa\beta_\mu\right)
\label{3.152}
\end{equation}
and $N_{\mu\nu\kappa\lambda}$ is the non-principal part, given by:
\begin{equation}
N_{\mu\nu\kappa\lambda}=-\eta^{-2}HA_{\mu\nu\kappa\lambda}-\frac{1}{2}\eta^{-2}B_{\mu\nu\kappa\lambda}
-\frac{1}{4}\sigma\eta^{-2}C_{\mu\nu\kappa\lambda}
\label{3.153}
\end{equation}
Here, $\eta$ is the sound speed and:
\begin{eqnarray}
&&A_{\mu\nu\kappa\lambda}=s_{\mu\lambda}s_{\nu\kappa}-s_{\mu\kappa}s_{\nu\lambda}\label{3.154}\\
&&B_{\mu\nu\kappa\lambda}=(dH)_\mu\xi_{\nu\kappa\lambda}-(dH)_\nu\xi_{\mu\kappa\lambda}
+(dH)_\kappa\xi_{\lambda\mu\nu}-(dH)_\lambda\xi_{\kappa\mu\nu}\label{3.155}\\
&&C_{\mu\nu\kappa\lambda}=((dH)_\mu\beta_\nu-(dH)_\nu\beta_\mu)
((dH)_\kappa\beta_\lambda-(dH)_\lambda\beta_\kappa))\label{3.156}
\end{eqnarray}
Moreover, in \ref{3.155}:
\begin{equation}
\xi_{\mu\kappa\lambda}=s_{\kappa\mu}\beta_\lambda-s_{\lambda\mu}\beta_\kappa
\label{3.157}
\end{equation}

Let now $\sR$ be the curvature of the $(S_{\ub,u}, \sh)$ and let $W$, $X$, $Y$, $Z$ be $S$ 
vectorfields. We have:
\begin{eqnarray}
&&\sR(W,X,Y,Z)=\sh(W,\sR(X,Y)Z)\label{3.158}\\
&&\hspace{22mm}=h(W,\sD_X\sD_Y Z-\sD_Y\sD_X Z-\sD_{[X,Y]}Z)\nonumber\\
&&\hspace{22mm}=h(W,D_X\sD_Y Z-D_Y\sD_X Z-D_{[X,Y]}Z)\nonumber
\end{eqnarray}
By \ref{3.26},
\begin{eqnarray*}
&&D_X\sD_Y Z=D_X D_Y Z-\frac{1}{2a}\left(\chi(Y,Z)D_X\Lb+\chib(Y,Z)D_X L\right)\\
&&\hspace{24mm}-\frac{\Lb}{2}X(a^{-1}\chi(Y,Z))-\frac{L}{2} X(a^{-1}\chib(Y,Z))
\end{eqnarray*}
hence:
\begin{equation}
h(W,D_X\sD_Y Z)=h(W,D_X D_Y Z)-\frac{1}{2a}(\chi(Y,Z)\chib(X,W)+\chib(Y,Z)\chi(X,W))
\label{3.159}
\end{equation}
Substituting this and a similar expression with $X$ and $Y$ interchanged in \ref{3.158} and noting 
that 
$$h(W, D_X D_Y Z-D_Y D_X Z-D_{[X,Y]}Z)=h(W,R(X,Y)Z)=R(W,Z,X,Y)$$
we obtain the Gauss equation:
\begin{eqnarray}
&&\sR(W,Z,X,Y)=R(W,Z,X,Y)-\frac{1}{2a}\left((\chi(X,W)\chib(Y,Z)+\chib(X,W)\chi(Y,Z)\right.\nonumber\\
&&\hspace{45mm}-\left.\chi(Y,W)\chib(X,Z)-\chib(Y,W)\chi(X,Z)\right)
\label{3.160}
\end{eqnarray}
Setting $W=\Omega_A$, $Z=\Omega_B$, $X=\Omega_C$, $Y=\Omega_D$ and denoting
$$\sR_{ABCD}=\sR(\Omega_A,\Omega_B,\Omega_C,\Omega_D), \ \ \ 
R_{ABCD}=R(\Omega_A,\Omega_B,\Omega_C,\Omega_D)$$
this reads:
\begin{equation}
\sR_{ABCD}=R_{ABCD}-\frac{1}{2a}\left(\chi_{AC}\chib_{BD}+\chib_{AC}\chi_{BD}
-\chi_{AD}\chib_{BC}-\chib_{AD}\chi_{BC}\right)
\label{3.161}
\end{equation}
We are now to substitute for 
$$R_{ABCD}=R_{\mu\nu\kappa\lambda}\Omega_A^\mu\Omega_B^\nu\Omega_C^\kappa\Omega_D^\lambda$$
from the formulas \ref{3.151} - \ref{3.157}. We have
$$(D^2 H)_{AB}=\Omega_A(\Omega_B H)-(D_{\Omega_A}\Omega_B)H$$
and
$$(\sD^2 H)_{AB}=\Omega_A(\Omega_B H)-(\sD_{\Omega_A}\Omega_B)H$$
therefore:
\begin{equation}
(D^2 H)_{AB}-(\sD^2 H)_{AB}=-(D_{\Omega_A}\Omega_B-\sD_{\Omega_A}\Omega_B)H
\label{3.162}
\end{equation}
On the other hand, by \ref{3.26}:
\begin{equation}
D_{\Omega_A}\Omega_B-\sD_{\Omega_A}\Omega_B=\frac{1}{2a}(\chib_{AB}L+\chi_{AB}\Lb)
\label{3.163}
\end{equation}
Consequently:
\begin{equation}
(D^2 H)_{AB}=(\sD^2 H)_{AB}-\frac{1}{2a}(\chib_{AB}LH+\chi_{AB}\Lb H)
\label{3.164}
\end{equation}
Substituting this in 
$P_{ABCD}=P_{\mu\nu\kappa\lambda}\Omega_A^\mu\Omega_B^\nu\Omega_C^\kappa\Omega_D^\lambda$ 
as given by \ref{3.152} we obtain:
\begin{eqnarray}
&&P_{ABCD}=\tilde{P}_{ABCD}\label{3.165}\\
&&\hspace{12mm}-\frac{1}{4a}\left(\chib_{CB}\sbeta_D\sbeta_A+\chib_{DA}\sbeta_C\sbeta_B
-\chib_{CA}\sbeta_D\sbeta_B-\chib_{DB}\sbeta_C\sbeta_A\right)LH\nonumber\\
&&\hspace{12mm}-\frac{1}{4a}\left(\chi_{CB}\sbeta_D\sbeta_A+\chi_{DA}\sbeta_C\sbeta_B
-\chi_{CA}\sbeta_D\sbeta_B-\chi_{DB}\sbeta_C\sbeta_A\right)\Lb H\nonumber
\end{eqnarray}
where:
\begin{eqnarray}
&&\tilde{P}_{ABCD}=\label{3.166}\\
&&\frac{1}{2}\left((\sD^2 H)_{CB}\sbeta_D\sbeta_A+(\sD^2 H)_{DA}\sbeta_C\sbeta_B
-(\sD^2 H)_{CA}\sbeta_D\sbeta_B-(\sD^2 H)_{DB}\sbeta_C\sbeta_A\right)\nonumber
\end{eqnarray}
Then in regard to the right hand side of \ref{3.161}, substituting for $\chi$ and $\chib$ in terms 
of $\tchi$ and $\tchib$ from \ref{3.75} and \ref{3.78} respectively the terms which are singular 
when the product $\rho\rhob$ vanishes cancel and we obtain:
\begin{eqnarray}
&&P_{ABCD}-\frac{1}{2a}\left(\chi_{AC}\chib_{BD}+\chib_{AC}\chi_{BD}
-\chi_{AD}\chib_{BC}-\chib_{AD}\chi_{BC}\right)\label{3.167}\\
&&=\tilde{P}_{ABCD}-\frac{1}{2c}\left(\tchi_{AC}\tchib_{BD}+\tchib_{AC}\tchi_{BD}
-\tchi_{AD}\tchib_{BC}-\tchib_{AD}\tchi_{BC}\right)\nonumber
\end{eqnarray}
Also, $N_{ABCD}=N_{\mu\nu\kappa\lambda}\Omega_A^\mu\Omega_B^\nu\Omega_C^\kappa\Omega_D^\lambda$ 
is given by \ref{3.153} - \ref{3.157}. Thus \ref{3.161} yields the following proposition. 

\vspace{2.5mm}

\noindent{\bf Proposition 3.7} \ \ \ The Gauss equation of the embedding of $S_{\ub,u}$ in 
the acoustical spacetime $({\cal N}, h)$ takes the form:
\begin{eqnarray*}
&&\sR_{ABCD}=\tilde{P}_{ABCD}-\frac{1}{2c}\left(\tchi_{AC}\tchib_{BD}+\tchib_{AC}\tchi_{BD}
-\tchi_{AD}\chib_{BC}-\tchib_{AD}\tchi_{BC}\right)\\
&&\hspace{15mm}+N_{ABCD}
\end{eqnarray*}
where $\tilde{P}_{ABCD}$ is the principal part, given by:
\begin{eqnarray*}
&&\tilde{P}_{ABCD}=\frac{1}{2}\left((\sD^2 H)_{CB}\sbeta_D\sbeta_A+(\sD^2 H)_{DA}\sbeta_C\sbeta_B
\right.\\
&&\hspace{22mm}-\left.(\sD^2 H)_{CA}\sbeta_D\sbeta_B-(\sD^2 H)_{DB}\sbeta_C\sbeta_A\right)
\end{eqnarray*}
while $N_{ABCD}$ is given by:
$$N_{ABCD}=-\eta^{-2}HA_{ABCD}-\frac{1}{2}\eta^{-2}B_{ABCD}-\frac{1}{4}\sigma\eta^{-2}C_{ABCD}$$
with
$$A_{ABCD}=\sss_{AD}\sss_{BC}-\sss_{AC}\sss_{BD},$$
$$B_{ABCD}=(\sd H)_A\xi_{BCD}-(\sd H)_B\xi_{ACD}+(\sd H)_C\xi_{DAB}-(\sd H)_D\xi_{CAB},$$
$$C_{ABCD}=((\sd H)_A\sbeta_B-(\sd H)_B\sbeta_A)((\sd H)_C\sbeta_D-(\sd H)_D\sbeta_C),$$
and 
$$\xi_{ABC}=\sss_{BA}\sbeta_C-\sss_{CA}\sbeta_B.$$

\vspace{2.5mm}

We remark that the Codazzi equations as well the 2nd variation and cross variation equations 
can also be derived by appealing to the acoustical curvature. In particular, the Codazzi equations 
for $\tchi$ and $\tchib$ can be derived by considering the acoustical curvature components 
$R(Z,L,X,Y)$ and $R(Z,\Lb,X,Y)$ respectively, with $Z$, $X$, $Y$ being arbitrary $S$ vectorfields. 
The 2nd variation equations for $\tchi$ and $\tchib$ can be derived by considering the acoustical curvature components $R(X,L,Y,L)$ and $R(X,\Lb,Y,\Lb)$ respectively, with $X$, $Y$ being arbitrary 
$S$ vectorfields. Also, the cross variation equations can be derived by considering the acoustical 
curvature components $R(X,L,Y,\Lb)$ with $X$, $Y$ being arbitrary $S$ vectorfields. However, these 
acoustical curvature components, like the component $R(W,Z,X,Y)$ which enters the Gauss equation 
(see \ref{3.160}) are singular where the product $\rho\rhob$, equivalently $a$, vanishes, the 
acoustical metric $h$ being degenerate there (see \ref{2.44}). The fact that these acoustical 
structure equations are actually regular where $a$ vanishes will appear in such a derivation as the 
result of delicate  cancellations. In contrast, the manner in which we have established Propositions 
3.4 to 3.6 above, which avoids considering the acoustical curvature, is straightforward. 
On the other hand in deriving the Gauss equation in Proposition 3.7 appealing to the acoustical 
curvature cannot be avoided. 

\pagebreak

\chapter{The Problem of the Free Boundary}

\section{Analysis of the Boundary Conditions}

We now consider the two boundary conditions on ${\cal K}$ stated in Section 1.6 : the linear 
jump condition \ref{1.328} and the nonlinear jump condition \ref{1.329}. Since along ${\cal K}$ 
the vectorfields $\Omega_A:A=1,...,n-1$ (see \ref{2.36}) together with the vectorfield $T$ 
(see \ref{2.28}) constitute a frame field for ${\cal K}$, the linear condition \ref{1.328} 
is equivalent to the conditions: 
\begin{equation}
\Omega_A^\mu\triangle\beta_\mu=0 \ : \ A=1,...,n-1
\label{4.1}
\end{equation}
together with the condition:
\begin{equation}
T^\mu\triangle\beta_\mu=0
\label{4.2}
\end{equation}
Let us denote:
\begin{equation}
\ep=N^\mu\triangle\beta_\mu, \ \ \ \epb=\Nb^\mu\triangle\beta_\mu
\label{4.3}
\end{equation}
Then according to \ref{4.1} the vectorfield along ${\cal K}$ with rectangular components 
$(h^{-1})^{\mu\nu}\triangle\beta_\nu$ is given by:
\begin{equation}
(h^{-1})^{\mu\nu}\triangle\beta_\nu=-\frac{1}{2c}(\epb N^\mu+\ep\Nb^\mu)
\label{4.4}
\end{equation}
Also, since $T^\mu=L^\mu+\Lb^\mu=\rho N^\mu+\rhob\Nb^\mu$, \ref{4.2} reads:
\begin{equation}
\rho\ep+\rhob\epb=0 \ \ \mbox{or:} \ \ r:=\frac{\rhob}{\rho}=-\frac{\ep}{\epb}
\label{4.5}
\end{equation}
which in terms of the definitions \ref{2.74} is expressed as: 
\begin{equation}
r:=\frac{\lambda}{\lambdab}=-\frac{\ep}{\epb}
\label{4.6}
\end{equation}
Since the functions $\rho$, $\rhob$ are both positive, it follows that $\ep$, $\epb$ have 
opposite signs. Moreover, since ${\cal K}$ must have along its past boundary, which 
coincides with $\partial_-{\cal B}$, the past boundary of ${\cal B}$, the same tangent hyperplane 
at each point as ${\cal B}$, and ${\cal B}$ is acoustically null outgoing along $\partial_-{\cal B}$, 
it follows that the ratio $r=\rhob/\rho$ tends to zero as we approach $\partial_-{\cal B}$ along 
${\cal K}$. Hence by the condition \ref{4.5} the ratio $\ep/\epb$ also tends to zero as we 
approach $\partial_-{\cal B}$ along ${\cal K}$. 

We turn to the nonlinear jump condition \ref{1.329}. We have:
\begin{eqnarray*}
&&\triangle(G\beta^\mu)=G_+\beta_+^\mu-G_-\beta_-^\mu\\
&&\hspace{14mm}=(\triangle G)\beta^\mu+(G-\triangle G)\triangle\beta^\mu
\end{eqnarray*}
where for quantities defined on ${\cal K}$ by the future solution we omit the subscript $+$ when 
there can be no misunderstanding. We then obtain:
\begin{equation}
\triangle(G\beta^\mu)\triangle\beta_\mu=\nu\triangle G+\delta^2 (G-\triangle G)
\label{4.7}
\end{equation}
the quantities $\nu$ and $\delta^2$ being defined by:
\begin{equation}
\nu=\beta^\mu\triangle\beta_\mu, \ \ \ \delta^2=\triangle\beta^\mu\triangle\beta_\mu
\label{4.8}
\end{equation}
The quantity $\delta$ itself shall be defined in the sequel. Recall from \ref{1.d2} that 
$\beta^\mu$ stands for  $(g^{-1})^{\mu\nu}\beta_\nu$. By \ref{1.d10}, \ref{1.d13}, 
\begin{equation}
(h^{-1})^{\mu\nu}\beta_\nu=\beta^\mu(1+\sigma F)=\eta^{-2}\beta^\mu \ \ 
\mbox{or} \ \ \beta^\mu=\eta^2(h^{-1})^{\mu\nu}\beta_\nu
\label{4.9}
\end{equation}
In view of the expression \ref{4.4} we then obtain:
\begin{eqnarray}
&&\nu=\eta^2(h^{-1})^{\mu\nu}\beta_\mu\triangle\beta_\nu\nonumber\\
&&\hspace{3mm}=-\frac{\eta^2}{2c}(\beta_N\epb+\beta_{\Nb}\ep)
\label{4.10}
\end{eqnarray}
Let us also define the quantity:
\begin{equation}
\mu=(h^{-1})^{\mu\nu}\triangle\beta_\mu\triangle\beta_\nu
\label{4.11}
\end{equation}
Substituting from the expression \ref{4.4} we obtain:
\begin{equation}
\mu=-\frac{1}{c}\epb\ep
\label{4.12}
\end{equation}
Since $\ep$, $\epb$ have opposite signs, $\mu$ is a positive quantity. By \ref{1.d10}:
\begin{eqnarray}
&&\delta^2=((h^{-1})^{\mu\nu}+F\beta^\mu\beta^\nu)\triangle\beta_\mu\triangle\beta_\nu\nonumber\\
&&\hspace{4mm}=\mu+F\nu^2
\label{4.13}
\end{eqnarray}
therefore $\delta^2$ is a fortiori a positive quantity.

We have: 
\begin{eqnarray*}
&&\triangle\sigma=-\beta_+^\mu\beta_{+\mu}+\beta_-^\mu\beta_{-\mu}=-\beta^\mu\beta_\mu+(\beta^\mu-\triangle\beta^\mu)(\beta_\mu-\triangle\beta_\mu)\\
&&\hspace{6mm}=-2\beta^\mu\triangle\beta_\mu+\triangle\beta^\mu\triangle\beta_\mu
\end{eqnarray*}
that is, in terms of the definitions \ref{4.8}, \ref{4.11}, 
\begin{equation}
\triangle\sigma=-2\nu+\delta^2=-2\nu+\mu+F\nu^2
\label{4.14}
\end{equation}
We have:
\begin{eqnarray}
&&\triangle G=G_+-G_-=-\int_0^1\frac{d}{ds}G(\sigma-s\triangle\sigma)ds\nonumber\\
&&\hspace{7mm}=K\triangle\sigma
\label{4.15}
\end{eqnarray}
where, denoting $G^\prime=dG/d\sigma$, 
\begin{equation}
K(\sigma,\triangle\sigma)=\int_0^1 G^\prime(\sigma-s\triangle\sigma)ds
\label{4.16}
\end{equation}
is a known smooth function of $(\sigma,\triangle\sigma)$. 

Substituting \ref{4.13} and \ref{4.15}, in which \ref{4.14} is substituted, in \ref{4.7}, 
we conclude that:
\begin{equation}
\triangle(G\beta^\mu)\triangle\beta_\mu=Q
\label{4.17}
\end{equation}
where $Q$ is a smooth function of $(\sigma,\mu,\nu)$, given by:
\begin{equation}
Q=K(-2\nu+\mu+F\nu^2)(\nu-\mu-F\nu^2)+G(\mu+F\nu^2)
\label{4.18}
\end{equation}
Here $K$ stands for $K(\sigma,\triangle\sigma)$, with $\triangle\sigma$ given by \ref{4.14}, so this 
is a known smooth function of $(\sigma,\mu,\nu)$. From \ref{4.16} we have, denoting 
$G^{\prime\prime}=d^2 G/d\sigma^2$, 
\begin{equation}
K(\sigma,\triangle\sigma)=G^\prime(\sigma)-\frac{1}{2}G^{\prime\prime}(\sigma)\triangle\sigma
+(\triangle\sigma)^2\tilde{K}_2(\sigma,\triangle\sigma)
\label{4.19}
\end{equation}
where $\tilde{K}_2$ is a smooth function of $(\sigma,\triangle\sigma)$. According to \ref{1.d7}:
\begin{equation}
G^\prime=\frac{1}{2}GF \ \ \mbox{hence} \ \ 
G^{\prime\prime}=\frac{1}{2}G\left(F^\prime+\frac{1}{2}F^2\right)
\label{4.20}
\end{equation}
where we denote $F^\prime=dF/d\sigma$. 

We have:
\begin{equation}
Q=Q_0+\mu Q_1+\mu^2\tilde{Q}_2
\label{4.21}
\end{equation}
where 
\begin{equation}
Q_0=\left. Q\right|_{\mu=0}, \ \ \ Q_1=\left.\frac{\partial Q}{\partial\mu}\right|_{\mu=0}
\label{4.22}
\end{equation}
and $\tilde{Q}_2$ is a smooth function of $(\sigma,\mu,\nu)$. 

From \ref{4.18},
\begin{equation}
Q_0=K_0(-2\nu+F\nu^2)(\nu-F\nu^2)+GF\nu^2
\label{4.23}
\end{equation}
where 
\begin{equation}
K_0=K(\sigma,\left.\triangle\sigma\right|_{\mu=0})
\label{4.24}
\end{equation}
with (see \ref{4.14})
\begin{equation}
\left.\triangle\sigma\right|_{\mu=0}=-2\nu+F\nu^2
\label{4.25}
\end{equation} 
From \ref{4.19} we then obtain:
\begin{equation}
K_0(\sigma,\nu)=G^\prime(\sigma)+G^{\prime\prime}(\sigma)\nu+\nu^2\tilde{K}_{0,2}(\sigma,\nu)
\label{4.26}
\end{equation}
where $\tilde{K}_{0,2}$ is a smooth function of $(\sigma,\nu)$. Substituting \ref{4.26} in \ref{4.23}, 
by \ref{4.20} the terms proportional to $\nu^2$ cancel while the terms proportional to $\nu^3$ are
$$(3G^\prime F-2G^{\prime\prime})\nu^3=-G(F^\prime-F^2)\nu^3$$
Now by \ref{1.d12}, \ref{1.d13}:
\begin{equation}
F^\prime-F^2=\eta^{-4}H^\prime
\label{4.27}
\end{equation}
where we denote $H^\prime=dH/d\sigma$. Thus the terms in \ref{4.23} proportional to $\nu^3$ 
are
$$-\eta^{-4}GH^\prime\nu^3$$
and we conclude that:
\begin{equation}
Q_0=\nu^3\hat{Q}_0 \ \ \mbox{where} \ \ \hat{Q}_0=-\eta^{-4}GH^\prime+\nu\tilde{Q}_{0,4}
\label{4.28}
\end{equation}
where $\tilde{Q}_{0,4}$ is a smooth function of $(\sigma,\nu)$. 

From \ref{4.18} (see \ref{4.22}):
\begin{equation}
Q_1=G+(3\nu-2F\nu^2)K_0-K_1(2\nu-F\nu^2)(\nu-F\nu^2)
\label{4.29}
\end{equation}
where we denote
$$K_1=\left.\frac{\partial K}{\partial\mu}\right|_{\mu=0}$$
By \ref{4.16}, in view of the fact that by \ref{4.14} $\partial\triangle\sigma/\partial\mu=1$, 
we have:
\begin{equation}
K_1=\frac{\partial K}{\partial\triangle\sigma}(\sigma,\left.\triangle\sigma\right|_{\mu=0})
=-\int_0^1 G^{\prime\prime}(\sigma-s\left.\triangle\sigma\right|_{\mu=0})sds
\label{4.30}
\end{equation}
This is a smooth function of $(\sigma,\nu)$. 

Substituting \ref{4.28} and \ref{4.29} in \ref{4.21} and setting: 
\begin{equation}
\mu=m\nu^3
\label{4.31}
\end{equation}
we see that the nonlinear jump condition $Q=0$ (see \ref{4.17}) becomes:
\begin{equation}
\hat{Q}_0(\nu)+mQ_1(\nu)+m^2\nu^3\tilde{Q}_2(m\nu^3,\nu)=0
\label{4.32}
\end{equation}
the dependence of the coefficients $\hat{Q}_0$, $Q_1$, $\tilde{Q}_2$ on $\sigma$ being understood. 
At $\nu=0$ this reduces to:
\begin{equation}
\hat{Q}_0(0)+mQ_1(0)=0
\label{4.33}
\end{equation}
Since by \ref{4.28} and \ref{4.29}:
\begin{equation}
\hat{Q}_0(0)=-\eta^{-4}GH^\prime, \ \ \ Q_1(0)=G
\label{4.34}
\end{equation}
this yields:
\begin{equation}
\left. m\right|_{\nu=0}:=m_0=\eta^{-4}H^\prime
\label{4.35}
\end{equation}
Applying the implicit function theorem to \ref{4.32} we then obtain the following proposition.

\vspace{2.5mm}

\noindent{\bf Proposition 4.1} \ \ \ There is an open $\nu$-interval containing 0 and 
a smooth function $m(\nu)$ defined in this interval with $m(0)=m_0$ such that \ref{4.32} 
has in the said interval the unique solution $m=m(\nu)$. 

\vspace{2.5mm} 

Writing, in accordance with \ref{4.5} or \ref{4.6},
\begin{equation}
\ep=-r\epb
\label{4.36}
\end{equation}
we express \ref{4.10} and \ref{4.12} in the form:
\begin{equation}
\nu=-\frac{\eta^2}{2c}(\beta_N-r\beta_{\Nb})\epb, \ \ \ \mu=\frac{1}{c}r\epb^2
\label{4.37}
\end{equation}
Substituting these in the equation
$$\mu-m(\sigma,\nu)\nu^3=0$$
this equation takes the form:
\begin{equation}
r+\frac{\eta^6}{8c^2}(\beta_N-r\beta_{\Nb})^3 \epb
m\left(\sigma,-\frac{\eta^2}{2c}(\beta_N-r\beta_{\Nb})\epb\right)=0
\label{4.38}
\end{equation}
Together with the fact that $r$ tends to 0 as we approach $\partial_-{\cal B}$ along ${\cal K}$, 
this implies that the ratio $r/\epb$ tends to $j_0$ as we approach $\partial_-{\cal B}$ along 
${\cal K}$, where:
\begin{equation}
j_0=-\frac{\eta^6\beta_N^3}{8c^2}m_0=-\frac{\eta^2\beta_N^3}{8c^2}H^\prime
\label{4.39}
\end{equation}
(see \ref{4.35}). Setting then: 
\begin{equation}
r=s\epb
\label{4.40}
\end{equation}
equation \ref{4.38} becomes:
\begin{equation}
s+\frac{\eta^6}{8c^2}(\beta_N-s\epb\beta_{\Nb})^3 
m\left(\sigma,-\frac{\eta^2}{2c}(\beta_N-s\epb\beta_{\Nb})\epb\right)=0
\label{4.41}
\end{equation}
This equation is of the form:
\begin{equation}
s-i(\kappa,s\epb,\epb)=0
\label{4.42}
\end{equation}
where $i$ is a smooth function of $(\kappa,s\epb,\epb)$ with $\kappa$ standing for the quadruplet:
\begin{equation}
\kappa=(\sigma,c,\beta_N,\beta_{\Nb})
\label{4.43}
\end{equation}
At $\epb=0$ we have:
\begin{equation}
i(\kappa,0,0)=j_0(\kappa)
\label{4.44}
\end{equation}
Applying then the implicit function theorem to \ref{4.42} we obtain the following proposition.

\vspace{2.5mm}

\noindent{\bf Proposition 4.2} \ \ \ There is an open $\epb$-interval containing 0 and 
a smooth function $j(\kappa,\epb)$ defined in this interval with $j(\kappa,0)=j_0(\kappa)$ 
such that \ref{4.42} has in the said interval the unique solution 
$$s=j(\kappa,\epb)$$ 
Thus, the nonlinear jump condition \ref{1.329} determines $\ep$ as a function of $\epb$ 
according to
$$\ep=-j(\kappa,\epb)\epb^2$$

\vspace{2.5mm} 

We now come to the quantity $\delta$ whose square is defined by the 2nd of \ref{4.8}. We recall from 
Section 1.4 the covector field $\xi$ defined along ${\cal K}$ such that at each point on ${\cal K}$ 
the null space of $\xi$ is the tangent space to ${\cal K}$ at the point and $\xi(V)>0$ for any vector 
$V$ at the point which points to the interior of ${\cal N}$, the domain of the new solution. Moreover 
$\xi$ was normalized to be of unit magnitude with respect to the Minkowski metric $g$ 
(see \ref{1.197}). Then $\xi^\sharp$, the vectorfield along ${\cal K}$ corresponding to $\xi$ 
through $g$, is the unit normal to ${\cal K}$ with respect to $g$ pointing to the interior of 
${\cal N}$, a space-like relative to $g$ vectorfield. Denoting then 
\begin{equation}
\beta_{\bot\pm}=\beta_\pm\cdot\xi^\sharp
\label{4.45}
\end{equation}
we define:
\begin{equation}
\delta=\triangle\beta_\bot
\label{4.46}
\end{equation}
Recalling from Section 1.6 that at each point on ${\cal K}$ the covectors $\xi$ and $\triangle\beta$ 
are colinear, we have $\triangle\beta=(\triangle\beta_\bot)\xi$, that is:
\begin{equation}
\triangle\beta=\delta\xi
\label{4.47}
\end{equation}
Hence $\delta^2$ is indeed given by the 2nd of \ref{4.8}. 

From \ref{1.53}, \ref{1.199} we have:
\begin{equation}
\beta_{\bot\pm}=-h_\pm u_{\bot\pm}
\label{4.48}
\end{equation}
hence by \ref{1.202}:
\begin{equation}
\beta_{\bot\pm}=-fh_\pm V_\pm, \ \ f>0
\label{4.49}
\end{equation}
Therefore:
\begin{equation}
\delta=-f\triangle(hV)
\label{4.50}
\end{equation}
Differentiating with respect to $p$ at constant $s$ we obtain:
$$\frac{d}{dp}(hV)=V^2+h\frac{dV}{dp}$$
and by \ref{1.215}:
$$h\frac{dV}{dp}=-\frac{V^2}{\eta^2}$$
Hence:
\begin{equation}
\frac{d}{dp}(hV)=-V^2\left(\frac{1}{\eta^2}-1\right)<0
\label{4.51}
\end{equation}
It follows that the sign of $\triangle(hV)$ is opposite to that of $\triangle p$. Consequently, 
{\em the sign of $\delta$ is the same as that of $\triangle p$}, positive in the 
{\em compressive} case: $\triangle p>0$, and negative in the {\em depressive} case: $\triangle p<0$. 
Recall here from Section 1.4 that $\triangle p$ is positive or negative according as to whether 
$H^\prime_-$ is negative or positive. 

\vspace{5mm}

\section{The Transformation Functions and the Identification Equations}

Recall from Section 2.5 that the free boundary ${\cal K}$ is represented from the point of view of 
${\cal N}$, the domain of the new solution by: 
\begin{equation}
{\cal K}=\{(\tau,\tau,\vartheta) \  : \tau\geq 0, \vartheta\in S^{n-1}\}
\label{4.52}
\end{equation}
Thus a given point of ${\cal K}$ corresponds to a given pair $(\tau,\vartheta)$ where $\tau\geq 0$ 
and $\vartheta\in S^{n-1}$. Now ${\cal K}\setminus\partial_-{\cal K}$ is contained 
in ${\cal M}$, the domain of the prior maximal classical development, 
and $\partial_-{\cal K}$, the past boundary of ${\cal K}$, coincides with 
$\partial_-{\cal B}$, the past boundary of ${\cal B}$, which is part of the boundary of ${\cal M}$. The point on ${\cal K}\setminus\partial_-{\cal K}$ 
which is represented as $(\tau,\tau,\vartheta)$, $\tau>0$, $\vartheta\in S^{n-1}$ 
from the point of view of the new solution, must be identified with a point in ${\cal M}$. 
Let $(t,u^\prime,\vartheta^\prime)$ be the coordinates of 
this point from the point of view of the prior solution in the acoustical coordinate system discussed 
in Sections 1.5 and 2.5. Now the spatial rectangular coordinates $x^{\prime i}$ have been determined 
by the prior solution as functions of $(t,u^\prime,\vartheta^\prime)$ in ${\cal M}$. 
Suppose then that the rectangular coordinates $x^\mu$ are also determined in 
${\cal N}$, as functions of the acoustical coordinates $(\ub,u,\vartheta)$. Consider the 
restrictions of these functions to the boundary ${\cal K}$:
\begin{equation}
f(\tau,\vartheta)=x^0(\tau,\tau,\vartheta), \ \ \ g^i(\tau,\vartheta)=x^i(\tau,\tau,\vartheta) \ : \ 
i=1,..., n
\label{4.53}
\end{equation}
Now {\em the rectangular coordinates of points which are identified must coincide}. Therefore, we 
must have:
\begin{equation}
t=f(\tau,\vartheta), \ \ \ x^{\prime i}(t,u^\prime,\vartheta^\prime)=g^i(\tau,\vartheta) 
\ : i=1,...,n
\label{4.54}
\end{equation}
That is, representing ${\cal K}$ and $\partial_-{\cal K}$ by:
\begin{equation}
{\cal K}=\{(\tau,\vartheta) \ : \tau\geq 0, \vartheta\in S^{n-1}\}, \ \ \ 
\partial_-{\cal K}=\{(0,\vartheta) : \vartheta\in S^{n-1}\},
\label{4.55}
\end{equation}
there must be a smooth function $w$ on ${\cal K}$ and a smooth mapping $\psi$ of ${\cal K}$ into 
$S^{n-1}$ such that with 
\begin{equation}
u^\prime=w(\tau,\vartheta), \ \ \ \vartheta^\prime=\psi(\tau,\vartheta)
\label{4.56}
\end{equation}
and $t$ given by the 1st of \ref{4.54} the $n$ equations constituting the 2nd of \ref{4.54} hold. 
In other words, we have:
\begin{equation}
x^{\prime i}(f(\tau,\vartheta),w(\tau,\vartheta),\psi(\tau,\vartheta))=g^i(\tau,\vartheta) 
\ : i=1,...,n \ \ : \ \forall (\tau,\vartheta)\in{\cal K} 
\label{4.57}
\end{equation}
Moreover, at each $\tau\geq 0$, $\psi|_{S^*_\tau}$ is a diffeomorphism of 
\begin{equation}
S^*_\tau=\{\tau\}\times S^{n-1}
\label{4.58}
\end{equation}
onto $S^{n-1}$, and, according to the discussion in Section 2.5, $w$ is a negative function decreasing 
with $\tau$ on ${\cal K}\setminus\partial_-{\cal K}$ and vanishes on $\partial_-{\cal K}=S^*_0$. 
Also, $\psi|_{S^*_0}$ is the identity. Since $\partial_-{\cal K}$ coincides with 
$\partial_-{\cal B}$ which is known from the prior solution, the values of the functions $f$ and $g^i$ on $\partial_-{\cal K}$ are known, and we have:
\begin{equation}
g^i(0,\vartheta)=x^{\prime i}(f(0,\vartheta),0,\vartheta) \ : \ \forall\vartheta\in S^{n-1}
\label{4.59}
\end{equation}
We call $w$ and $\psi$ {\em transformation functions} and equations \ref{4.57} {\em identification 
equations}. The identification equations constitute a system of $n$ equations in $n$ unknowns, namely 
$w$ and the $n-1$ functions representing the mapping $\psi$ in local coordinates on $S^{n-1}$. 

Now, in the case of $n=2$ spatial dimensions, choosing the unit of length so that the perimeter of 
$\partial_-{\cal B}$ is set equal to $2\pi$, we achieve $\partial_-{\cal B}$ to be {\em isometric}, 
not merely diffeomorphic, to the standard $S^1$. The universal cover of $S^1$ being $\mathbb{R}^1$ 
which has affine structure, we can define the notion of a {\em displacement} on $S^1$. That is to say, 
$\vartheta$ may be set equal to the arc length (relative to the induced metric $\sh$) on 
$\partial_-{\cal B}$ from an abritrary origin and a displacement is meant in terms of arc length. 
We can then express $\vartheta^\prime$ as a sum:
\begin{equation}
\vartheta^\prime=\vartheta+\varphi^\prime
\label{4.60}
\end{equation}
Then with \ref{4.56}, the identification equations \ref{4.57} take the form:
\begin{equation}
F^i((\tau,\vartheta),(u^\prime,\varphi^\prime))=0 \ : \ i=1,2
\label{4.61}
\end{equation}
where:
\begin{equation}
F^i((\tau,\vartheta),(u^\prime,\varphi^\prime))=
x^{\prime i}(f(\tau,\vartheta),u^\prime,\vartheta+\varphi^\prime)-g^i(\tau,\vartheta)
\label{4.62}
\end{equation}

For $n>2$ spatial dimensions, even though we may choose the unit of length so that the area 
(volume) of 
$\partial_-{\cal B}$ is set equal to that of the standard $S^{n-1}$, we cannot achieve 
$\partial_-{\cal B}$ to be isometric to the standard $S^{n-1}$, as in general 
$(\partial_-{\cal B}, \sh)$ is a Riemannian manifold of non-constant curvature. 
Moreover, there is no affine structure associated to $S^{n-1}$. 
As a consequence, the notion of displacement on $\partial_-{\cal B}$ cannot be canonically defined. 
Nevertheless this notion can be defined relative to a given diffeomorphism of $\partial_-{\cal B}$ 
onto $S^{n-1}$, and a given atlas for $S^{n-1}$, for example a stereographic atlas. Such an atlas 
consists of two charts $\chi$ and $\tilde{\chi}$, with domains ${\cal U}$ and $\tilde{{\cal U}}$ respectively. Let ${\cal V}$ and $\tilde{{\cal V}}$ be the ranges of $\chi$ and $\tilde{\chi}$ 
respectively. These are domains in $\mathbb{R}^{n-1}$, which we can take to be convex. If 
$\vartheta\in{\cal U}$ then $\vartheta$ is represented in the chart $\chi$ by 
$$\chi(\vartheta)=(\vartheta^1,...,\vartheta^{n-1})\in{\cal V}$$
and if $\vartheta\in\tilde{{\cal U}}$ then $\vartheta$ is represented in the chart $\tilde{\chi}$ by 
$$\tilde{\chi}(\vartheta)=(\tilde{\vartheta}^1,...,\tilde{\vartheta}^{n-1})\in\tilde{{\cal V}}$$
And if $\vartheta\in{\cal U}\bigcap\tilde{{\cal U}}$, there are smooth functions 
$$f^A(\vartheta^1,...,\vartheta^{n-1}) \ : \ A=1,...,n-1 \ \ \mbox{on} \ 
\chi({\cal U}\bigcap\tilde{{\cal U}})\subset{\cal V}$$
such that:
\begin{equation}
\tilde{\vartheta}^A=f^A(\vartheta^1,...,\vartheta^{n-1}) \ : \ A=1,...,n-1
\label{4.63}
\end{equation}
In fact $f=(f^1,...,f^{n-1})$ is a diffeomorphism of 
$\chi({\cal U}\bigcap\tilde{{\cal U}})\subset{\cal V}$ onto \\
$\tilde{\chi}({\cal U}\bigcap\tilde{{\cal U}})\subset\tilde{{\cal V}}$. 
Then with $\vartheta, \vartheta^\prime \in{\cal U}$, $\vartheta^\prime$ is represented in the 
chart $\chi$ by 
$$(\vartheta^{\prime A}=\vartheta^A+\varphi^{\prime A} \ : \ A=1,...,n-1)$$
Likewise, with $\vartheta, \vartheta^\prime \in\tilde{{\cal U}}$, $\vartheta^\prime$ is represented 
in the chart $\tilde{\chi}$ by 
$$(\tilde{}\vartheta^{\prime A}=\tilde{\vartheta}^A+\tilde{\varphi}^{\prime A} \ : \ A=1,...,n-1)$$
Moreover, by \ref{4.63}, if $\vartheta, \vartheta^\prime \in {\cal U}\bigcap\tilde{{\cal U}}$, we have:
$$\tilde{\vartheta}^A+\tilde{\varphi}^{\prime A}=f^A(\vartheta^1+\varphi^{\prime 1},...,\vartheta^{n-1}+\varphi^{\prime n-1}) 
\ : \ A=1,...,n-1$$
Subtracting \ref{4.63} we then obtain:
\begin{equation}
\tilde{\varphi}^{\prime A}=
f^A(\vartheta^1+\varphi^{\prime 1},...,\vartheta^{n-1}+\varphi^{\prime n-1})
-f^A(\vartheta^1,...,\vartheta^{n-1}) \ : \ A=1,...,n-1
\label{4.64}
\end{equation}
The $f^A \ :A=1,...,n-1$ being smooth functions we have a Taylor expansion of the 
$\tilde{\varphi}^{\prime A} \ :A=1,...,n-1$ in $(\varphi^{\prime 1},...,\varphi^{\prime n-1})$,  
with remainder, to arbitrary order, the coefficients being smooth functions of 
$(\vartheta^1,...,\vartheta^{n-1})$. We can then think of a {\em displacement} $\varphi^\prime$ 
as being represented in the chart $\chi$ by the components 
$(\varphi^{\prime 1},...,\varphi^{\prime n-1})$ and in the chart $\tilde{\chi}$ by the components 
$(\tilde{\varphi}^{\prime 1},...,\tilde{\varphi}^{\prime n-1})$ and transforming under the change 
of charts according to \ref{4.64}. And we can think of the point $\vartheta^\prime$ as resulting 
from the point $\vartheta$ by the displacement $\varphi^\prime$ and write:
\begin{equation}
\vartheta^\prime=\vartheta+\varphi^\prime
\label{4.65}
\end{equation}
Setting then as in \ref{4.62}:
\begin{equation}
F^i((\tau,\vartheta),(u^\prime,\varphi^\prime))=
x^{\prime i}(f(\tau,\vartheta),u^\prime,\vartheta+\varphi^\prime)-g^i(\tau,\vartheta)
\label{4.66}
\end{equation}
the identification equations take in general the form:
\begin{equation}
F^i((\tau,\vartheta),(u^\prime,\varphi^\prime))=0 \ : \ i=1,...,n
\label{4.67}
\end{equation}
In treating the case $n>2$ we proceed from this point in a given chart. 

Consider the difference:
\begin{eqnarray*}
&&F^i((\tau,\vartheta),(u^\prime,\varphi^\prime))-F^i((\tau,\vartheta),(u^\prime,0))\\
&&\hspace{10mm}=x^{\prime i}(f(\tau,\vartheta),u^\prime,\vartheta+\varphi^\prime)
-x^{\prime i}(f(\tau,\vartheta),u^\prime,\vartheta)
\end{eqnarray*}
This can be written as:
$$\int_0^1\frac{d}{ds}x^{\prime i}(f(\tau,\vartheta),u^\prime,\vartheta+s\varphi^\prime)ds 
=\varphi^{\prime A}\int_0^1\Omega^{\prime i}_A
(f(\tau,\vartheta),u^\prime,\vartheta+s\varphi^\prime)ds$$
Note that the displacement $s\varphi^\prime$ only has meaning in a given chart. Here, 
\begin{equation}
\Omega^{\prime i}_A=\frac{\partial x^{\prime i}}{\partial\vartheta^{\prime A}}
\label{4.68}
\end{equation}
defined in ${\cal M}$ by the prior solution. Defining then:
\begin{equation}
G_A^i((\tau,\vartheta),(u^\prime,\varphi^\prime))=
\int_0^1\Omega^{\prime i}_A(f(\tau,\vartheta),u^\prime,\vartheta+s\varphi^\prime)ds
\label{4.69}
\end{equation}
we have:
\begin{equation}
F^i((\tau,\vartheta),(u^\prime,\varphi^\prime))-F^i((\tau,\vartheta),(u^\prime,0))
=\varphi^{\prime A}G_A^i((\tau,\vartheta),(u^\prime,\varphi^\prime))
\label{4.70}
\end{equation}
Furthermore, setting:
\begin{eqnarray}
&&H^i((\tau,\vartheta),u^\prime)=F^i((\tau,\vartheta),(u^\prime,0))\nonumber\\
&&\hspace{20mm}=x^{\prime i}(f(\tau,\vartheta),u^\prime,\vartheta)-g^i(\tau,\vartheta)
\label{4.71}
\end{eqnarray}
we have the expression:
\begin{equation}
F^i((\tau,\vartheta),(u^\prime,\varphi^\prime))=H^i((\tau,\vartheta),u^\prime)
+\varphi^{\prime A}G_A^i((\tau,\vartheta),(u^\prime, \varphi^\prime))
\label{4.72}
\end{equation}

At $\tau=0$ the solution of the identification equations \ref{4.67} is:
\begin{equation}
u^\prime=0, \ \ \ \varphi^{\prime A}=0 \ :A=1,...,n-1 \ \ \mbox{: at $\tau=0$}
\label{4.73}
\end{equation}
(see \ref{4.59}) and  \ref{4.69} reduces to:
\begin{equation}
G_A^i((0,\vartheta),(0,0))=\Omega^{\prime i}_A(f(0,\vartheta),0,\vartheta)
=\left.\Omega^{\prime i}_A\right|_{\partial_-{\cal B}}(\vartheta)
\label{4.74}
\end{equation}

As discussed in Section 2.5 the acoustical function $u$ defined in ${\cal N}$ extends 
across $\underline{{\cal C}}$ the acoustical function $u^\prime$ defined in ${\cal M}$  
(see \ref{2.171}). The vectorfield $L^\prime$ defined by the prior solution 
being the tangent field to the generators of the level sets of $u^\prime$ parametrized by 
$t=x^0$, we see that along $\underline{{\cal C}}=\Cb_0$, $L^\prime$ must coincide with 
$N$, which has the same properties relative to $u$:
\begin{equation}
L^\prime=N \ \mbox{: along $\Cb_0$}
\label{4.75}
\end{equation}
Now, as discussed in Section 1.5, the coordinate $\vartheta^\prime$ is canonically defined 
in the region in ${\cal M}$ which 
corresponds to the union of the generators of the $C^\prime_{u^\prime}$, the level sets of $u^\prime$, 
which intersect ${\cal B}$, by the condition that $\vartheta^\prime$ is constant along the invariant 
curves of ${\cal B}$. In the region in ${\cal M}$ in question we have 
$u^\prime\leq 0$. The coordinate $\vartheta^\prime$ may be extended to the region in ${\cal M}$ 
where $u^\prime>0$, the future boundary of which is $\underline{{\cal C}}$, by the condition that 
$\vartheta^\prime$ is constant along the generators of $\underline{{\cal C}}=\Cb_0$.

Consider a point $p\in\Cb_0$. To $p$ there corresponds a definite value of $t$ of the rectangular coordinate $x^0$ and a definite common value of $u=u^\prime$. Thus $p\in S_{0,u}=\Cb_0\bigcap C_u$
and at the same time $p\in S^\prime_{t,u^\prime}=\Sigma_t\bigcap C^\prime_{u^\prime}$. At $p$ we have 
the frame $(N,\Nb,\Omega_A:A=1,...,n-1)$ which refers to the new solution, as well as the frame 
$(L^\prime,\hat{T}^\prime,\Omega^\prime_A:A=1,...,n-1)$ which refers to the prior solution. Here
\begin{equation}
\hat{T}^\prime=\kappa^{\prime -1}T^\prime
\label{4.76}
\end{equation}
is the unit vector corresponding to $T^\prime$ (see Section 1.5). We shall find the relation between 
the two frames. First, we have already noted \ref{4.75}. Concerning the relation between the 
$\Omega_A$ and the $\Omega^\prime_A$ we note, referring to Section 2.5, that:
\begin{equation}
\vartheta=\vartheta^\prime \ \mbox{: along $\Cb_0$}
\label{4.77}
\end{equation}
The rectangular components of the $\Omega_A$ are:
\begin{equation}
\Omega_A^\mu=\frac{\partial x^\mu}{\partial\vartheta^A} \ \mbox{: at constant $(\ub,u)$} \ \ ; \ 
\mu=0,...,n
\label{4.78}
\end{equation}
while the rectangular components of the $\Omega^\prime_A$ are:
\begin{equation}
\Omega^{\prime 0}_A=0, \ \ \ 
\Omega^{\prime i}_A=\frac{\partial x^{\prime i}}{\partial\vartheta^{\prime A}} \ 
\mbox{: at constant $(t,u^\prime)$} \ \ ; \ i=1,...,n
\label{4.79}
\end{equation}
Since the rectangular coordinate functions coincide along $\underline{{\cal C}}$, that is 
$$x^0=t \ \ \mbox{and} \ \ x^i=x^{\prime i}\ :i=1,...,n \ \ \mbox{along $\Cb_0$}$$
and $u^\prime=u$ while $\vartheta^{\prime A}=\vartheta^A$ along $\Cb_0$, it follows that, 
along $\Cb_0$:
\begin{equation}
\frac{\partial x^i}{\partial\vartheta^A}=\frac{\partial x^{\prime i}}{\partial\vartheta^{\prime A}}
+\frac{\partial x^{\prime i}}{\partial t}\frac{\partial t}{\partial\vartheta^A}
\label{4.80}
\end{equation}
or, in view of the fact that
\begin{equation}
\frac{\partial x^{\prime i}}{\partial t}=L^{\prime i}, 
\label{4.81}
\end{equation}
which along $\Cb_0$ is equal to $N^i$, 
\begin{equation}
\Omega^\prime_A=\Omega_A-\pi_A N \ \mbox{: along $\Cb_0$}
\label{4.82}
\end{equation}
using the notation \ref{3.45}. We finally consider $\hat{T}^\prime$. This is characterized by the 
conditions that it is tangential to the $\Sigma_t$, $h$-orthogonal to the $\Omega^\prime_A:A=1,...,n-1$, of unit magnitude and inward pointing. Remark that the rectangular components of 
the 1-form $\beta$, hence also of the acoustical metric $h$, are continuous across 
$\underline{{\cal C}}=\Cb_0$. Expanding $\hat{T}^\prime$ along $\Cb_0$ in the 
$(N,\Nb,\Omega_A:A=1,...,n-1)$ frame:
\begin{equation}
\hat{T}^\prime=\sa^A\Omega_A+a N+\overline{a}\Nb
\label{4.83}
\end{equation}
we have $\hat{T}^\prime x^0=0$ which reads:
\begin{equation}
\sa^A\pi_A+a+\overline{a}=0
\label{4.84}
\end{equation}
We also have $h(\hat{T}^\prime,\Omega^\prime_A)=0$ which by \ref{4.82} reads:
\begin{equation}
\sh_{AB}\sa^B+2c\pi_A\overline{a}=0
\label{4.85}
\end{equation}
Lastly we have $h(\hat{T}^\prime,\hat{T}^\prime)=1$, which reads:
\begin{equation}
\sh_{AB}\sa^A\sa^B-4ca\overline{a}=1
\label{4.86}
\end{equation}
Solving \ref{4.85} for $\sa^A$ we obtain:
\begin{equation}
\sa^A=-2c\overline{a}(\sh^{-1})^{AB}\pi_B
\label{4.87}
\end{equation}
Substituting in \ref{4.84} and solving for $a$ we obtain:
\begin{equation}
a=-(1-2c(\sh^{-1})^{AB}\pi_A\pi_B)\overline{a}
\label{4.88}
\end{equation}
Substituting then \ref{4.87} and \ref{4.88} in \ref{4.86} the last becomes:
\begin{equation}
4c(1-c(\sh^{-1})^{AB}\pi_A\pi_B)\overline{a}^2=1
\label{4.89}
\end{equation}
Now by the expansion \ref{2.73} of the components of the reciprocal acoustical metric we have, 
in view of the definition \ref{3.45} and the fact that $N^0=\Nb^0=1$, 
\begin{equation}
\frac{1}{c}-(\sh^{-1})^{AB}\pi_A\pi_B=-(h^{-1})^{00}=1+F\beta_0^2
\label{4.90}
\end{equation}
The right hand side is the function denoted by $\alpha^{-2}$ in [Ch-S] (see \ref{1.289}). Thus, 
\begin{equation}
\frac{1}{c}-(\sh^{-1})^{AB}\pi_A\pi_B=\frac{1}{\alpha^2}
\label{4.91}
\end{equation}
and \ref{4.89} becomes:
$$\frac{4c^2}{\alpha^2}\overline{a}^2=1$$
Finally the condition that $\hat{T}^\prime$ is inward pointing means $\overline{a}>0$. 
We then obtain:
\begin{equation}
\overline{a}=\frac{\alpha}{2c}
\label{4.92}
\end{equation}
which, together with \ref{4.87} and \ref{4.88} determine the coefficients of the expansion \ref{4.83}. 
Substituting for $\Omega_A$ in terms of $\Omega^\prime_A$ from \ref{4.82} this expansion with the 
coefficients so determined takes the form:
\begin{equation}
\frac{\hat{T}^\prime}{\alpha}+(\sh^{-1})^{AB}\pi_B\Omega^\prime_A=\frac{\Nb-N}{2c}
\label{4.93}
\end{equation}
Remark that here both sides are tangential to the $\Sigma_t$. Summarizing, we have established the 
following proposition.

\vspace{2.5mm}

\noindent{\bf Proposition 4.3} \ \ \ The frames $(N,\Nb,\Omega_A:A=1,...,n-1)$ and 
$(L^\prime,\hat{T}^\prime,\Omega^\prime_A:A=1,...,n-1)$ are related along $\underline{{\cal C}}$ 
by:
$$L^\prime=N, \ \ \Omega^\prime_A=\Omega_A-\pi_A N, \ \ 
\frac{\hat{T}^\prime}{\alpha}+(\sh^{-1})^{AB}\pi_B\Omega^\prime_A=\frac{\Nb-N}{2c}$$

\vspace{2.5mm}

Now, in the context of [Ch-S] we have the incoming acoustically null vectorfield 
\begin{equation}
\Lb^\prime=\frac{1}{2}\alpha^{-1}\kappa^\prime L^\prime+T^\prime
\label{4.94}
\end{equation}
which is $h$-orthogonal to the surfaces $S^\prime_{t,u^\prime}=\Sigma_t\bigcap C^\prime_{u^\prime}$ 
and satisfies
\begin{equation}
\Lb^\prime u^\prime=1
\label{4.95}
\end{equation}
Substituting for $T^\prime$ in terms of $\hat{T}^\prime$ from \ref{4.76} in \ref{4.94} we obtain:
\begin{equation}
\Lb^\prime=\kappa^\prime\hat{\Lb^\prime} \ \ \mbox{where} \ \ 
\hat{\Lb^\prime}=\frac{1}{2}\alpha^{-1}L^\prime+\hat{T}^\prime
\label{4.96}
\end{equation}
Along $\Cb_0$ we may substitute for $\hat{T}^\prime$ from Proposition 4.3 in $\hat{\Lb^\prime}$ 
to obtain, along $\Cb_0$:
\begin{equation}
\hat{\Lb^\prime}=-\alpha(\sh^{-1})^{AB}\pi_B\Omega^\prime_A+
\frac{1}{2}\left(\frac{1}{\alpha}-\frac{\alpha}{c}\right)L^\prime+\frac{\alpha}{2c}\Nb
\label{4.97}
\end{equation}
Since $u^\prime$ coincides with $u$ along $\Cb_0$, we have:
\begin{equation}
\Nb u^\prime=\Nb u=\frac{\Lb u}{\rhob}=\frac{1}{\rhob} \ \ \mbox{: along $\Cb_0$}
\label{4.98}
\end{equation}
Then, since $\Omega^\prime_A u^\prime=L^\prime u^\prime=0$, applying \ref{4.97} to $u^\prime$ 
we obtain:
\begin{equation}
\hat{\Lb^\prime} u^\prime=\frac{\alpha}{2\lambda} \ \ \mbox{: along $\Cb_0$}
\label{4.99}
\end{equation}
hence, in view of the fact that $\alpha\kappa^\prime=\mu^\prime$ (see \ref{1.288}), 
\begin{equation}
\Lb^\prime u^\prime=\frac{\mu^\prime}{2\lambda} \ \ \mbox{: along $\Cb_0$}
\label{4.100}
\end{equation}
Comparing with \ref{4.95} yields the following (compare with \ref{2.75}).

\vspace{2.5mm}

\noindent{\bf Proposition 4.4} \ \ \ Along $\underline{{\cal C}}$ the following relation holds:
$$\lambda=\frac{1}{2}\mu^\prime$$
the right hand side defined by the prior solution and the left hand side by the new solution. 
Consequently, the derivatives of $\lambda$ along $\Cb_0$ can all be obtained from the corresponding 
derivatives of $\mu^\prime$ determined by the prior solution.

\vspace{2.5mm}

In view of Proposition 4.4, multiplying by $\mu^\prime$ the expression for $(\Nb-N)/2c$ 
along $\Cb_0$ given by Proposition 4.3 yields the relation:
\begin{equation}
T^\prime+\mu^\prime\left((\sh^{-1})^{AB}\pi_B\Omega^\prime_A+\frac{L^\prime}{2c}\right)=\Lb 
\ \ \mbox{: along $\Cb_0$}
\label{4.101}
\end{equation}

Now, since, according to the results of [Ch-S] quoted in Section 1.5, $\mu^\prime$ vanishes on the 
whole of ${\cal B}$, in particular at $\partial_-{\cal B}$, and by \ref{1.307} and the 2nd of 
\ref{1.291}:
\begin{equation}
T^\prime=\frac{\partial}{\partial u^\prime}
-\mu^\prime\hat{\Xi}^A\frac{\partial}{\partial\vartheta^{\prime A}}
\label{4.102}
\end{equation}
and, moreover,\ref{1.312} holds, we in fact have:
\begin{equation}
\mu^\prime|_{\partial_-{\cal B}}=T^\prime\mu^\prime|_{\partial_-{\cal B}}=0
\label{4.103}
\end{equation}
This implies:
$$T^{\prime 2}\mu^\prime|_{\partial_-{\cal B}}=
\left.\frac{\partial^2\mu^\prime}{\partial u^{\prime 2}}\right|_{\partial_-{\cal B}}$$
Hence, according to \ref{1.317} the function $k$ on $\partial_-{\cal B}$ defined by:
\begin{equation}
k=T^{\prime 2}\mu^\prime|_{\partial_-{\cal B}}
\label{4.104}
\end{equation}
is {\em strictly positive}. 

By Proposition 4.4 and \ref{4.102} along $\Cb_0$ we have:
\begin{equation}
T^\prime\mu^\prime+\mu^\prime\left((\sh^{-1})^{AB}\pi_B\Omega^\prime_A\mu^\prime
+\frac{L^\prime\mu^\prime}{2c}\right)=2\Lb\lambda
\label{4.105}
\end{equation}
In particular at $\partial_-{\cal B}$ by \ref{4.103} it holds:
\begin{equation}
\lambda|_{\partial_-{\cal B}}=\Lb\lambda|_{\partial_-{\cal B}}=0
\label{4.106}
\end{equation}
Moreover, applying \ref{4.102} again and evaluating the result at $\partial_-{\cal B}$ we obtain, 
in terms of the definition \ref{4.104}:
\begin{equation}
\Lb^2\lambda|_{\partial_-{\cal B}}=\frac{1}{2}k
\label{4.107}
\end{equation}

According to \ref{1.302}, $\partial\mu^\prime/\partial t<0$ on the whole of ${\cal B}$, in particular 
at $\partial_-{\cal B}$. Thus, in view of the 1st of \ref{1.291}:
\begin{equation}
L^\prime=\frac{\partial}{\partial t}
\label{4.108}
\end{equation}
the function $l$ on $\partial_-{\cal B}$ defined by:
\begin{equation}
l=L^\prime\mu^\prime|_{\partial_-{\cal B}}
\label{4.109}
\end{equation}
is {\em strictly negative}. Moreover, from the propagation equation of [Ch-S] for $\mu^\prime$ 
along the integral curves of $L^\prime$ we see that on the whole of ${\cal B}$, where $\mu^\prime$ 
vanishes, we have: 
$$L^\prime\mu^\prime=\frac{1}{2}\beta_{L^\prime}^2 T^\prime H$$
which, evaluated at $\partial_-{\cal B}$ gives:
\begin{equation}
l=\frac{1}{2}\beta_N^2 T^\prime H|_{\partial_-{\cal B}}
\label{4.110}
\end{equation}
The rectangular components $\beta_\mu$ of the 1-form $\beta$ being 
continuous across $\underline{{\cal C}}=\Cb_0$ their derivatives along $\Cb_0$ are also continuous. 
Applying then \ref{4.102} to $H$ and evaluating the result at $\partial_-{\cal B}$ we conclude that:
\begin{equation}
T^\prime H|_{\partial_-{\cal B}}=\Lb H|_{\partial_-{\cal B}}
\label{4.111}
\end{equation}
Comparing with \ref{4.110} we see that:
\begin{equation}
l=\frac{1}{2}\beta_N^2\Lb H|_{\partial_-{\cal B}}
\label{4.112}
\end{equation}

Consider now the boundary condition \ref{4.6}, that is:
\begin{equation}
r\lambdab=\lambda \ \ \mbox{: along ${\cal K}$}
\label{4.113}
\end{equation}
Applying the vectorfield $T=L+\Lb$, which is tangential to ${\cal K}$, to this we obtain:
\begin{equation}
(Tr)\lambdab+rT\lambdab=T\lambda \ \ \mbox{: along ${\cal K}$}
\label{4.114}
\end{equation}
Let us evaluate this at $\partial_-{\cal B}$. By \ref{4.40} and Proposition 4.2 $r$ vanishes 
at $\partial_-{\cal B}$, as the jump $\epb$ vanishes there. Moreover, by \ref{4.106} 
$T\lambda=L\lambda$ at $\partial_-{\cal B}$. The last is given by the propagation equation 
for $\lambda$ of Proposition 3.3, which at $\partial_-{\cal B}$ reduces to 
$L\lambda=q\lambdab$, in view of the vanishing of $\lambda$ there. Therefore \ref{4.114} reduces at 
$\partial_-{\cal B}$ to:
\begin{equation}
(Tr-q)\lambdab=0 \ \ \mbox{: at $\partial_-{\cal B}$}
\label{4.115}
\end{equation}
Since $r$ vanishes at $\partial_-{\cal B}$ and is positive on ${\cal K}\setminus\partial_-{\cal K}$,
it follows that 
\begin{equation}
Tr\geq 0 \ \ \mbox{: at $\partial_-{\cal B}$}
\label{4.116}
\end{equation}
On the other hand, according to Proposition 3.3:
$$q=\frac{1}{4c}\beta_N^2\Lb H$$
Comparing with \ref{4.112} we see that:
\begin{equation}
cq|_{\partial_-{\cal B}}=\frac{1}{2}l<0
\label{4.117}
\end{equation}
We then conclude through \ref{4.115} that:
\begin{equation}
\lambdab=0 \ \ \mbox{: at $\partial_-{\cal B}$}
\label{4.118}
\end{equation}

Consider now the functions $f$ and $g^i:i=1,...,n$ (see \ref{4.53}). Since 
$$T=L+\Lb=\rho N+\rhob \Nb,$$
and on ${\cal K}$ :
$$T=\frac{\partial}{\partial\tau},$$
we have:
\begin{eqnarray}
&&\frac{\partial f}{\partial\tau}(\tau,\vartheta)=(\rho+\rhob)(\tau,\tau,\vartheta)\nonumber\\
&&\frac{\partial g^i}{\partial\tau}(\tau,\vartheta)=
(\rho N^i+\rhob\Nb^i)(\tau,\tau,\vartheta) \ : \ i=1,...,n
\label{4.119}
\end{eqnarray}
The vanishing of $\lambdab$ and $\lambda$ at $\partial_-{\cal B}$, just demonstrated, is equivalent 
to the vanishing of $\rho$ and $\rhob$ there (see \ref{2.74}). As $\tau=0$ corresponds to 
$\partial_-{\cal K}=\partial_-{\cal B}$, it follows, evaluating \ref{4.119} at $\tau=0$, that:
\begin{equation}
\left.\frac{\partial f}{\partial\tau}\right|_{\tau=0}=0, \ \ \ 
\left.\frac{\partial g^i}{\partial\tau}\right|_{\tau=0}=0 \ : \ i=1,...,n
\label{4.120}
\end{equation}

Consider the functions $H^i:i=1,...,n$ defined by \ref{4.71}. In view of the fact that 
(see \ref{4.108})
\begin{equation}
\frac{\partial x^{\prime i}}{\partial t}=L^{\prime i}
\label{4.121}
\end{equation}
we have:
\begin{equation}
\frac{\partial H^i}{\partial\tau}((\tau,\vartheta),u^\prime)=
L^{\prime i}(f(\tau,\vartheta),u^\prime,\vartheta)\frac{\partial f}{\partial\tau}(\tau,\vartheta)
-\frac{\partial g^i}{\partial\tau}(\tau,\vartheta)
\label{4.122}
\end{equation}
By \ref{4.120} we then obtain:
\begin{equation}
\frac{\partial H^i}{\partial\tau}((0,\vartheta),0)=0
\label{4.123}
\end{equation}
By \ref{4.102} and \ref{4.76}):
\begin{equation}
\frac{\partial x^{\prime i}}{\partial u^\prime}=\mu^\prime S^i
\label{4.124}
\end{equation}
where $S^i$ are the components of the vectorfield:
\begin{equation}
S=\alpha^{-1}\hat{T}^\prime+\hat{\Xi}^A\Omega^\prime_A
\label{4.125}
\end{equation}
which is tangential to the $\Sigma_t$ and is defined by the prior solution. 
We then have:
\begin{equation}
\frac{\partial H^i}{\partial u^\prime}((\tau,\vartheta),u^\prime)=
(\mu^\prime S^i)(f(\tau,\vartheta),u^\prime,\vartheta)
\label{4.126}
\end{equation}
By virtue of the fact that $\mu^\prime$ vanishes at $\partial_-{\cal B}$ we then obtain:
\begin{equation}
\frac{\partial H^i}{\partial u^\prime}((0,\vartheta),0)=0
\label{4.127}
\end{equation}

Next, differentiating \ref{4.122} with respect to $\tau$ we obtain:
\begin{eqnarray}
&&\frac{\partial^2 H^i}{\partial\tau^2}((\tau,\vartheta),u^\prime)=
L^i(f(\tau,\vartheta),u^\prime,\vartheta)
\frac{\partial^2 f}{\partial\tau^2}(\tau,\vartheta)\nonumber\\
&&\hspace{25mm}+(L^\prime L^{\prime i}(f(\tau,\vartheta),u^\prime,\vartheta)
\left(\frac{\partial f}{\partial\tau}\right)^2(\tau,\vartheta)\nonumber\\
&&\hspace{25mm}-\frac{\partial^2 g^i}{\partial\tau^2}(\tau,\vartheta)
\label{4.128}
\end{eqnarray}
Applying $T$ to \ref{4.119} we obtain:
\begin{eqnarray}
&&\frac{\partial^2 f}{\partial\tau^2}(\tau,\vartheta)
=(T\rho+T\rhob)(\tau,\tau,\vartheta)\nonumber\\
&&\frac{\partial^2 g^i}{\partial\tau^2}(\tau,\vartheta)=
(N^i T\rho+\Nb^iT\rhob+\rho TN^i+\rhob T\Nb^i)(\tau,\tau,\vartheta)
\label{4.129}
\end{eqnarray}
In view of the vanishing of $\lambda$ and $\lambdab$ at $\partial_-{\cal B}$ (see \ref{4.106}, 
\ref{4.118}) the propagation equations for $\lambda$ and $\lambdab$ of Proposition 3.3 evaluated 
at $\partial_-{\cal B}$ give:
\begin{equation}
L\lambda=\Lb\lambdab=0 \ \ \mbox{: at $\partial_-{\cal B}$}
\label{4.130}
\end{equation}
Together with \ref{4.106} and \ref{4.118} these imply, in view of \ref{2.74} and the fact that 
$T=L+\Lb$, that \ref{4.129} reduce at $\tau=0$ to:
\begin{equation}
\left.\frac{\partial^2 f}{\partial\tau^2}\right|_{\tau=0}=L\rho|_{\partial_-{\cal B}}, \ \ \ 
\left.\frac{\partial^2 g^i}{\partial\tau^2}\right|_{\tau=0}=N^i L\rho|_{\partial_-{\cal B}} \ : \ 
i=1,...,n
\label{4.131}
\end{equation}
Since by the 1st of \ref{4.120}, the formula \ref{4.128} reduces at $(\tau,u^\prime)=(0,0)$ to
$$\frac{\partial^2 H^i}{\partial\tau^2}((0,\vartheta),0)=L^{\prime i}(f(0,\vartheta),0,\vartheta)
\frac{\partial^2 f}{\partial\tau^2}(0,\vartheta)-\frac{\partial^2 g^i}{\partial\tau^2}(0,\vartheta)$$
recalling from Proposition 4.3 that $N^i=L^{\prime i}$ along $\Cb_0$, in particular at 
$\partial_-{\cal B}$, we then conclude that:
\begin{equation}
\frac{\partial^2 H^i}{\partial\tau^2}((0,\vartheta),0)=0
\label{4.132}
\end{equation}
Differentiating \ref{4.126} with respect to $\tau$ we obtain:
\begin{equation}
\frac{\partial^2 H^i}{\partial\tau\partial u^\prime}((\tau,\vartheta),u^\prime)
=\frac{\partial(\mu^\prime S^i)}{\partial t}(f(\tau,\vartheta),u^\prime,\vartheta)
\frac{\partial f}{\partial\tau}(\tau,\vartheta)
\label{4.133}
\end{equation}
Then by the 1st of \ref{4.120}:
\begin{equation}
\frac{\partial^2 H^i}{\partial\tau\partial u^\prime}((0,\vartheta),0)=0
\label{4.134}
\end{equation}
Differentiating \ref{4.126} with respect to $u^\prime$ we obtain:
\begin{equation}
\frac{\partial^2 H^i}{\partial u^{\prime 2}}((\tau,\vartheta),u^\prime)=
\frac{\partial(\mu^\prime S^i)}{\partial u^\prime}(f(\tau,\vartheta),u^\prime,\vartheta)
\label{4.135}
\end{equation}
Since according to \ref{4.103} (see \ref{4.102}):
\begin{equation}
\mu^\prime|_{\partial_-{\cal B}}=
\left.\frac{\partial\mu^\prime}{\partial u^\prime}\right|_{\partial_-{\cal B}}=0
\label{4.136}
\end{equation}
it follows that:
\begin{equation}
\frac{\partial^2 H^i}{\partial u^{\prime 2}}((0,\vartheta),0)=0
\label{4.137}
\end{equation}

Recalling also that the functions $H^i$, being defined by \ref{4.71}, vanish at 
$(\tau,u^\prime)=(0,0)$ by \ref{4.59}, the results \ref{4.123}, \ref{4.127}, and 
\ref{4.132}, \ref{4.134}, \ref{4.137}, imply that the Taylor expansion of 
$H^i((\tau,\vartheta),u^\prime)$ in $(\tau^\prime,u^\prime)$ at fixed $\vartheta$ begins with 
cubic terms. By \ref{4.73}, the functions $w$ and $\psi$ which according to \ref{4.56} 
express the solution of the identification equations \ref{4.67} satisfy:
\begin{equation}
w(0,\vartheta)=0, \ \ \ \psi(0,\vartheta)=\vartheta
\label{4.138}
\end{equation}
These functions being assumed smooth, the first implies that for the solution of the identification 
equations $u^\prime=O(\tau)$, hence $H^i=O(\tau^3)$. Then by \ref{4.72}, \ref{4.74}, and the fact 
that at each point $p$ in the closure of ${\cal M}$, in particular at each $p\in\partial_-{\cal B}$, 
the vectors $\Omega^\prime_A: A=1,...,n-1$ constitute a basis for $T_p S^\prime_{t,u^\prime}$, 
hence are linearly independent, we conclude that for the solution of the identification equations 
$\varphi^{\prime A}=O(\tau^3)$. Thus, in the sense of the given chart (see \ref{4.65}) we have:
\begin{equation}
\frac{\partial\psi^A}{\partial\tau}(0,\vartheta)=
\frac{\partial^2\psi^A}{\partial\tau^2}(0,\vartheta)=0
\label{4.139}
\end{equation}
Let us then set:
\begin{equation}
\frac{\partial w}{\partial\tau}(0,\vartheta)=v_0(\vartheta), \ \ \ 
\frac{\partial^3\psi^A}{\partial\tau^3}(0,\vartheta)=6\gamma_0^A(\vartheta)
\label{4.140}
\end{equation}
{\em The function $v_0$ must be non-positive}, $w$ being negative for $\tau>0$ and decreasing 
with $\tau$. For the solution of the identification equations we then have:
\begin{equation}
u^\prime=v_0(\vartheta)\tau+O(\tau^2), \ \ \ 
\varphi^{\prime A}=\gamma_0^A(\vartheta)\tau^3+O(\tau^4)
\label{4.141}
\end{equation}

Going now back to \ref{4.114} and applying $T$ again we obtain:
\begin{equation}
2(Tr)T\lambdab+rT^2\lambdab+\lambdab T^2 r=T^2\lambda \ \ \mbox{: along ${\cal K}$}
\label{4.142}
\end{equation}
At $\partial_-{\cal B}$ this reduces to:
\begin{equation}
2(Tr)T\lambdab=T^2\lambda \ \ \mbox{: at $\partial_-{\cal B}$}
\label{4.143}
\end{equation}
Now, we have:
$$T^2\lambda=TL\lambda+T\Lb\lambda$$
Applying $T$ to the propagation equation for $\lambda$ of Proposition 3.3, we obtain, since 
$T\lambda$ vanishes at $\partial_-{\cal B}$,
$$TL\lambda=qT\lambdab \ \ \mbox{: at $\partial_-{\cal B}$}$$
Also, we have:
$$T\Lb\lambda=\Lb^2\lambda+L\Lb\lambda$$
By \ref{2.17} $L\Lb\lambda=\Lb L\lambda-Z\cdot\sd\lambda$ and $\sd\lambda$ vanishes at 
$\partial_-{\cal B}=S_{0,0}$, while applying $\Lb$ to the propagation equation for $\lambda$ of 
Proposition 3.3 we see that by the 2nd of \ref{4.130} and \ref{4.106} $\Lb L\lambda$ vanishes 
at $\partial_-{\cal B}$. Combining, and recalling \ref{4.107}, we obtain:
\begin{equation}
T^2\lambda=qT\lambdab+\frac{1}{2}k \ \ \mbox{: at $\partial_-{\cal B}$}
\label{4.144}
\end{equation}
Inserting this in \ref{4.143} yields the equation:
\begin{equation}
(2Tr-q)T\lambdab=\frac{1}{2}k \ \ \mbox{: at $\partial_-{\cal B}$}
\label{4.145}
\end{equation}

To determine then $T\lambdab$ at $\partial_-{\cal B}$, we must determine $Tr$ there. By \ref{4.40} 
and Proposition 4.2 we have:
\begin{equation}
Tr|_{\partial_-{\cal B}}=s_0 T\epb|_{\partial_-{\cal B}}
\label{4.146}
\end{equation}
where 
\begin{equation}
s_0=\left.s\right|_{\partial_-{\cal B}}=j_0(\kappa_0), \ \ \ 
\kappa_0=\left.\kappa\right|_{\partial_-{\cal B}}
\label{4.a1}
\end{equation}
($\kappa$ being the quadruplet \ref{4.43}).
From \ref{4.3} we have, since the jumps $\triangle\beta_\mu$ vanish at $\partial_-{\cal B}$,
\begin{equation}
T\epb=\Nb^\mu T\triangle\beta_\mu \ \ \mbox{: at $\partial_-{\cal B}$}
\label{4.147}
\end{equation}
Let us denote from now on by $\beta^\prime_\mu$ the rectangular components of the 1-form $\beta$ 
corresponding to the prior solution. These are smooth functions of $(t,u^\prime,\vartheta^\prime)$. 
We can then write:
\begin{equation}
\triangle\beta_\mu=\left.\beta_\mu\right|_{{\cal K}}-\left.\beta^\prime_\mu\right|_{{\cal K}}
\label{4.148}
\end{equation}
or, in view of \ref{4.56}, in terms of the $(\tau,\vartheta)$ coordinates on ${\cal K}$, 
\begin{equation}
\triangle\beta_\mu(\tau,\vartheta)=\beta_\mu(\tau,\tau,\vartheta)
-\beta^\prime_\mu(f(\tau,\vartheta),w(\tau,\vartheta),\psi(\tau,\vartheta))
\label{4.149}
\end{equation}
Applying $T=\partial/\partial\tau$ we then obtain:
\begin{eqnarray}
&&(T\triangle\beta_\mu)(\tau,\vartheta)=(T\beta_\mu)(\tau,\tau,\vartheta)
-\left(\frac{\partial\beta^\prime_\mu}{\partial t}\right)
(f(\tau,\vartheta),w(\tau,\vartheta),\psi(\tau,\vartheta))
\frac{\partial f}{\partial\tau}(\tau,\vartheta)\nonumber\\
&&\hspace{25mm}-\left(\frac{\partial\beta^\prime_\mu}{\partial u^\prime}\right)
(f(\tau,\vartheta),w(\tau,\vartheta),\psi(\tau,\vartheta))
\frac{\partial w}{\partial\tau}(\tau,\vartheta)\nonumber\\
&&\hspace{25mm}-\left(\frac{\partial\beta^\prime_\mu}{\partial\vartheta^{\prime A}}\right)
(f(\tau,\vartheta),w(\tau,\vartheta),\psi(\tau,\vartheta))
\frac{\partial\psi^A}{\partial\tau}(\tau,\vartheta)
\label{4.150}
\end{eqnarray}
By the 1st of \ref{4.120} and by \ref{4.139}, \ref{4.140} this reduces at $\tau=0$, that is at 
$\partial_-{\cal B}$ to:
\begin{equation}
T\triangle\beta_\mu|_{\partial_-{\cal B}}=\left.\left(T\beta_\mu-
v_0\frac{\partial\beta^\prime_\mu}{\partial u^\prime}\right)\right|_{\partial_-{\cal B}}
\label{4.151}
\end{equation}
Now, in reference to \ref{4.147} we have (see \ref{3.42}):
$$\Nb^\mu L\beta_\mu=\rho s_{N\Nb}=0 \ \ \mbox{: at $\partial_-{\cal B}$}$$
hence:
\begin{equation}
\Nb^\mu T\beta_\mu=s_{\Nb\Lb} \ \ \mbox{: at $\partial_-{\cal B}$}
\label{4.152}
\end{equation}
Also (see \ref{4.102}, \ref{4.103}),
$$\frac{\partial\beta^\prime_\mu}{\partial u^\prime}=T^\prime\beta^\prime_\mu=\Lb\beta^\prime_\mu \ \ 
\mbox{: at $\partial_-{\cal B}$}$$
the last equality by virtue of \ref{4.101}. Since $\beta^\prime_\mu=\beta_\mu$ along 
$\underline{{\cal C}}$,  we have $\Lb\beta^\prime_\mu=\Lb\beta_\mu$ along $\Cb_0$. 
Therefore:
\begin{equation}
\Nb^\mu\frac{\partial\beta^\prime_\mu}{\partial u^\prime}=s_{\Nb\Lb} \ \ 
\mbox{: at $\partial_-{\cal B}$}
\label{4.153}
\end{equation}
Combining \ref{4.152} and \ref{4.153} we conclude through \ref{4.151} and \ref{4.147} that:
\begin{equation}
T\epb=(1-v_0)s_{\Nb\Lb} \ \ \mbox{: at $\partial_-{\cal B}$}
\label{4.154}
\end{equation}
Substituting this as well as \ref{4.39} in \ref{4.146} yields:
\begin{equation}
Tr=-\frac{(1-v_0)}{8}\frac{\eta^2}{c^2}\beta_N^3\frac{dH}{d\sigma}s_{\Nb\Lb} \ \ 
\mbox{: at $\partial_-{\cal B}$}
\label{4.155}
\end{equation}
We shall presently relate the coefficient of $1-v_0$ here to the negative function $l$ 
(see \ref{4.112}). We have:
\begin{equation}
\Lb H=-2\frac{dH}{d\sigma}\beta^\mu\Lb\beta_\mu \ \ \mbox{where} \ \ 
\beta^\mu=(g^{-1})^{\mu\nu}\beta_\nu
\label{4.156}
\end{equation}
Now $\beta^\mu$ being the rectangular components of a vectorfield which we may denote by 
$\beta^\sharp$, we can expand this vectorfield in the frame field $(N,\Nb,\Omega_A:A=1,...,n-1)$:
\begin{equation}
\beta^\sharp=\sbb^A\Omega_A+bN+\overline{b}\Nb
\label{4.157}
\end{equation}
If $X$ is an arbitrary vectorfield, we have 
$$h(\beta^\sharp,X)=h_{\mu\nu}\beta^\mu X^\nu$$
and from \ref{1.d11}:
\begin{equation}
h_{\mu\nu}\beta^\mu=(g_{\mu\nu}+H\beta_\mu\beta_\nu)\beta^\mu=\beta_\nu-H\sigma\beta_\nu=\eta^2\beta_\nu
\label{4.158}
\end{equation}
(see \ref{1.d13}). Hence, with $X$ an arbitrary vectorfield,  
\begin{equation}
h(\beta^\sharp,X)=\eta^2\beta\cdot X
\label{4.159}
\end{equation}
Substituting then the expansion \ref{4.157} and setting $X$ equal to $\Omega_A:A=1,...,n-1$, $\Nb$, 
$N$, successively, we determine the coefficients to be:
\begin{equation}
\sbb^A=\eta^2(\sh^{-1})^{AB}\sbeta_B \ : \ A=1,...,n-1, \ \ \ b=-\frac{\eta^2}{2c}\beta_{\Nb}, 
\ \ \ \overline{b}=-\frac{\eta^2}{2c}\beta_N
\label{4.160}
\end{equation}
We conclude that the following expansion holds:
\begin{equation}
\beta^\mu=\eta^2\left((\sh^{-1})^{AB}\sbeta_B\Omega_A^\mu-\frac{1}{2c}\beta_{\Nb}N^\mu
-\frac{1}{2c}\beta_N\Nb^\mu\right)
\label{4.161}
\end{equation}
Substituting this in \ref{4.156} yields:
\begin{equation}
\Lb H=-2\eta^2\frac{dH}{d\sigma}\left[\rhob\left((\sh^{-1})^{AB}\sbeta_B(\ss_{\Nb})_A
-\frac{1}{2c}\beta_{\Nb}s_{N\Nb}\right)-\frac{1}{2c}\beta_N s_{\Nb\Lb}\right]
\label{4.162}
\end{equation}
At $\partial_-{\cal B}$ this reduces to:
\begin{equation}
\Lb H=\frac{\eta^2}{c}\beta_N\frac{dH}{d\sigma}s_{\Nb\Lb} \ \ \mbox{: at $\partial_-{\cal B}$}
\label{4.163}
\end{equation}
Comparing with \ref{4.155} we see that:
\begin{equation}
Tr=-\frac{(1-v_0)}{8c}\beta_N^2\Lb H \ \ \mbox{: at $\partial_-{\cal B}$}
\label{4.164}
\end{equation}
or, comparing with \ref{4.112},
\begin{equation}
Tr=-\frac{(1-v_0)}{4c}l \ \ \mbox{: at $\partial_-{\cal B}$}
\label{4.165}
\end{equation}
Remark that the right hand side is positive at it should be. Comparing finally with \ref{4.117} 
we see that:
\begin{equation}
Tr=-\frac{(1-v_0)}{2}q \ \ \mbox{: at $\partial_-{\cal B}$}
\label{4.166}
\end{equation}

Substituting the result \ref{4.166} in \ref{4.145} yields:
\begin{equation}
T\lambdab=-\frac{k}{2(2-v_0)q}=-\frac{ck}{(2-v_0)l} \ \ \mbox{: at $\partial_-{\cal B}$}
\label{4.167}
\end{equation}
Also, substituting this in \ref{4.144} yields:
\begin{equation}
T^2\lambda=\frac{(1-v_0)k}{2(2-v_0)} \ \ \mbox{: at $\partial_-{\cal B}$}
\label{4.168}
\end{equation}

By Proposition 4.2 and the fact that the jump $\epb$ vanishes at $\partial_-{\cal B}$ we 
obtain:
\begin{equation}
T\ep=0 \ \ \mbox{: at $\partial_-{\cal B}$}
\label{4.169}
\end{equation}
Recalling the definitions \ref{4.3} this is $N^\mu T\triangle\beta_\mu=0$ : at $\partial_-{\cal B}$, 
or, in terms of \ref{4.148}, 
\begin{equation}
N^\mu T\beta_\mu=N^\mu T\beta^\prime_\mu \ \ \mbox{: at $\partial_-{\cal B}$}
\label{4.170}
\end{equation}
Now, we have $N^\mu\Lb\beta_\mu=\rhob s_{N\Nb}$ which vanishes at $\partial_-{\cal B}$. Therefore 
\ref{4.170} reduces to:
\begin{equation}
s_{NL}|_{\partial_-{\cal B}}=N^\mu T\beta^\prime_\mu|_{\partial_-{\cal B}}
\label{4.171}
\end{equation}
Along ${\cal K}$ we have:
\begin{equation}
\left.T\beta^\prime_\mu\right|_{{\cal K}}=
\left.\frac{\partial\beta^\prime_\mu}{\partial t}\right|_{{\cal K}}\frac{\partial f}{\partial\tau}+
\left.\frac{\partial\beta^\prime_\mu}{\partial u^\prime}\right|_{{\cal K}}
\frac{\partial w}{\partial\tau}+
\left.\frac{\partial\beta^\prime_\mu}{\partial\vartheta^{\prime A}}\right|_{{\cal K}}
\frac{\partial\psi^A}{\partial\tau}
\label{4.172}
\end{equation}
By the 1st of \ref{4.120}, \ref{4.139}, \ref{4.140} this reduces at $\partial_-{\cal B}$ to:
\begin{equation}
\left.T\beta^\prime_\mu\right|_{\partial_-{\cal B}}=v_0 
\left.T^\prime\beta^\prime_\mu\right|_{\partial_-{\cal B}}
\label{4.173}
\end{equation}
Therefore:
\begin{equation}
N^\mu T\beta^\prime_\mu|_{\partial_-{\cal B}}=v_0 
L^{\prime\mu}T^\prime\beta^\prime_\mu|_{\partial_-{\cal B}}
=v_0 T^{\prime\mu}L^\prime\beta^\prime_\mu|_{\partial_-{\cal B}}=0
\label{4.174}
\end{equation}
and by \ref{4.171} we conclude that:
\begin{equation}
s_{NL}|_{\partial_-{\cal B}}=0
\label{4.175}
\end{equation}
According to the conjugate of \ref{4.162}:
\begin{equation}
LH=-2\eta^2\frac{dH}{d\sigma}\left[\rho\left((\sh^{-1})^{AB}\sbeta_B(\ss_N)_A
-\frac{1}{2c}\beta_N s_{N\Nb}\right)-\frac{1}{2c}\beta_{\Nb} s_{NL}\right]
\label{4.176}
\end{equation}
Then by virtue of \ref{4.175} we have:
\begin{equation}
LH|_{\partial_-{\cal B}}=0
\label{4.177}
\end{equation}

We proceed to calculate the 3rd derivatives of the functions $H^i$ with respect to $(\tau,u^\prime)$ 
at $(\tau,u^\prime)=(0,0)$. Differentiating \ref{4.128} with respect to $\tau$ and evaluating the 
result at $(\tau,u^\prime)=(0,0)$ we obtain, in view of the 1st of \ref{4.120}, 
\begin{equation}
\left.\frac{\partial^3 H^i}{\partial\tau^3}\right|_{(\tau,u^\prime)=(0,0)}
=N^i\left.\frac{\partial^3 f}{\partial\tau^3}\right|_{\tau=0}
-\left.\frac{\partial^3 g^i}{\partial\tau^3}\right|_{\tau=0}
\label{4.178}
\end{equation}
Applying $T$ to \ref{4.129} we obtain:
\begin{eqnarray}
&&\frac{\partial^3 f}{\partial\tau^3}=T^2\rho+T^2\rhob\nonumber\\
&&\frac{\partial^3 g^i}{\partial\tau^3}=N^i T^2\rho+\Nb^i T^2\rhob
+2(TN^i)T\rho+2(T\Nb^i)T\rhob\nonumber\\
&&\hspace{12mm}+\rho T^2 N^i+\rhob T^2\Nb^i
\label{4.179}
\end{eqnarray}
At $\tau=0$ the last reduces to:
\begin{equation}
\left.\frac{\partial^3 g^i}{\partial\tau^3}\right|_{\tau=0}=N^2 T^2\rho+\Nb^i T^2\rhob
+2(TN^i)T\rho
\label{4.180}
\end{equation}
in view of the fact that $T\rhob$ vanishes at $\partial_-{\cal B}$. Substituting in \ref{4.176} 
we obtain:
\begin{equation}
\left.\frac{\partial^3 H^i}{\partial\tau^3}\right|_{(\tau,u^\prime)=(0,0)}
=\left.\left((N^i-\Nb^i)T^2\rhob-2(TN^i)T\rho\right)\right|_{\partial_-{\cal B}}
\label{4.181}
\end{equation}
Here we must determine $TN^i$ at $\partial_-{\cal B}$. From the expansion \ref{3.31} and from 
\ref{3.57}, \ref{3.68}, and the 1st of \ref{3.47}, we see that by virtue of \ref{4.175}, \ref{4.177}:
\begin{equation}
LN^i|_{\partial_-{\cal B}}=0
\label{4.182}
\end{equation}
As for $\Lb N^i$ we have the expansion \ref{3.33}. By Proposition 4.3, along $\Cb_0$ this 
takes the form (see also 1st of \ref{3.48}):
\begin{equation}
\Lb N^i=\sn^A\Omega^{\prime i}_A+\on(\Nb^i-N^i) \ \ \mbox{: along $\Cb_0$}
\label{4.183}
\end{equation}
Now, $\sn^A$ is given by \ref{3.a6} which at $\partial_-{\cal B}$ reduces to:
\begin{equation}
\sn^B\sh_{AB}=-\beta_N\sbeta_A\Lb H \ \ \mbox{: at $\partial_-{\cal B}$}
\label{4.184}
\end{equation}
Also, $\on$ is given by \ref{3.61} which at $\partial_-{\cal B}$ reduces to:
\begin{equation}
\on=\frac{1}{4c}\beta_N^2\Lb H \ \ \mbox{: at $\partial_-{\cal B}$}
\label{4.185}
\end{equation}
Defining then:
\begin{equation}
m^A=\left.(\sh^{-1})^{AB}\frac{\sbeta_B}{\beta_N}\right|_{\partial_-{\cal B}} \ : \ A=1,...,n-1
\label{4.186}
\end{equation}
and comparing with \ref{4.112} we see that:
\begin{equation}
\left.\sn^A\right|_{\partial_-{\cal B}}=-2lm^A, \ \ \ 
\left.c\on\right|_{\partial_-{\cal B}}=\frac{l}{2}
\label{4.187}
\end{equation}
Then in view of \ref{4.182}, \ref{4.183} we obtain:
\begin{equation}
TN^i|_{\partial_-{\cal B}}=
l\left.\left(-2m^A\Omega^{\prime i}_A+\frac{\Nb^i-N^i}{2c}\right)\right|_{\partial_-{\cal B}}
\label{4.188}
\end{equation}
Substituting this together with \ref{4.167}, \ref{4.168} in \ref{4.181} we obtain:
\begin{equation}
\left.\frac{\partial^3 H^i}{\partial\tau^3}\right|_{(\tau,u^\prime)=(0,0)}=
\frac{k}{(2-v_0)}\left.\left((1+v_0)\frac{(\Nb^i-N^i)}{2c}-4m^A\Omega^{\prime i}_A\right)
\right|_{\partial_-{\cal B}}
\label{4.189}
\end{equation}

Now, by Proposition 4.3:
\begin{equation}
\left.\frac{(\Nb^i-N^i)}{2c}\right|_{\partial_-{\cal B}}=
\left.\left(\frac{\hat{T}^{\prime i}}{\alpha}+(\sh^{-1})^{AB}\pi_B\Omega^{\prime i}_A\right)
\right|_{\partial_-{\cal B}}
\label{4.190}
\end{equation}
We shall express the right hand side in terms of:
\begin{equation}
S^i_0:=\left.S^i\right|_{\partial_-{\cal B}} \ \ \mbox{and the} \ \ 
\Omega^{\prime i}_{A 0}:=\left.\Omega^{\prime i}_A\right|_{\partial_-{\cal B}} \ : A=1,...,n-1
\label{4.191}
\end{equation}
where $S$ is the vectorfield defined by \ref{4.125}. For this purpose we must determine 
$\hat{\Xi}^A|_{\partial_-{\cal B}}$. By \ref{1.308} we have
\begin{equation}
\left.\hat{\Xi}^A\right|_{\partial_-{\cal B}}=\frac{\partial\Xi^A/\partial t|_{\partial_-{\cal B}}}
{\partial\mu^\prime/\partial t|_{\partial_-{\cal B}}}
\label{4.192}
\end{equation}
in $(t,u^\prime,\vartheta^\prime)$ coordinates on ${\cal M}$. By \ref{4.109} the denominator  
here is $l$. As for the numerator, by \ref{1.291} it is $-\Lambda^A|_{\partial_-{\cal B}}$, where 
$\Lambda=\Lambda^A\Omega^\prime_A$ is the commutator $\Lambda=[L^\prime,T^\prime]$. To determine 
this we shall refer to results from Chapter 3 of [Ch-S]. We have (see page 62 of [Ch-S]):
\begin{equation} 
\Lambda^A=-(\sh^{-1})^{AB}(2\zeta^\prime_B+\Omega_B^\prime\mu^\prime)
\label{4.193}
\end{equation}
On ${\cal B}$, $\zeta^\prime_A$ reduces to (see page 76 of [Ch-S]):
\begin{equation}
\zeta^\prime_A=-\frac{1}{2}(\beta^\prime\cdot L^\prime)(\beta^\prime\cdot\Omega^\prime_A)T^\prime H
\ \ \mbox{: on ${\cal B}$}
\label{4.194}
\end{equation}
This holds in particular at $\partial_-{\cal B}$, where by Proposition 4.3:
$$\beta^\prime\cdot L^\prime=\beta_N, \ \ \ \beta^\prime\cdot\Omega^\prime_A=\sbeta_A-\pi_A\beta_N$$
Comparing with \ref{4.186} and \ref{4.110} we then obtain:
\begin{equation}
\zeta^\prime_A|_{\partial_-{\cal B}}=l(\pi_A-m_A) 
\label{4.195}
\end{equation}
(where $m_A=\sh_{AB}m^B$). Moreover, since $\mu^\prime$ vanishes at $\partial_-{\cal B}$, 
$\Omega_A\mu^\prime$ likewise vanishes there, hence by Proposition 4.3:
\begin{equation}
\Omega^\prime_A\mu^\prime=-\pi_A L^\prime\mu^\prime=-\pi_A l \ \ \mbox{: at $\partial_-{\cal B}$}
\label{4.196}
\end{equation}
Substituting the results \ref{4.195}, \ref{4.196} in \ref{4.193} yields:
\begin{equation}
\Lambda^A|_{\partial_-{\cal B}}=-l(\pi^A-m^A)
\label{4.197}
\end{equation}
(where $\pi^A=(\sh^{-1})^{AB}\pi_B$). As a consequence we obtain through \ref{4.192}:
\begin{equation}
\left.\hat{\Xi}^A\right|_{\partial_-{\cal B}}=\pi^A-2m^A
\label{4.198}
\end{equation}
Comparing then with \ref{4.125} we see that the right hand side of \ref{4.190} is equal to:
$$S^i_0+2m^A\Omega^{\prime i}_{A 0}$$
and \ref{4.190} becomes:
\begin{equation}
\left.\frac{(\Nb^i-N^i)}{2c}\right|_{\partial_-{\cal B}}=S^i_0+2m^A\Omega^{\prime i}_{A 0}
\label{4.199}
\end{equation}
Substituting this finally in \ref{4.189} yields:
\begin{equation} 
\left.\frac{\partial^3 H^i}{\partial\tau^3}\right|_{(\tau,u^\prime)=(0,0)}=
\frac{k}{(2-v_0)}\left((1+v_0)S^i_0-2(1-v_0)m^A\Omega^{\prime i}_{A 0}\right)
\label{4.200}
\end{equation}

Next, differentiating \ref{4.133} with respect to $\tau$ we obtain:
\begin{eqnarray}
&&\frac{\partial^3 H^i}{\partial\tau^2\partial u^\prime}((\tau,\vartheta),u^\prime)=
\frac{\partial(\mu^\prime S^i)}{\partial t}(f(\tau,\vartheta,u^\prime,\vartheta)
\frac{\partial^2 f}{\partial\tau^2}(\tau,\vartheta)\nonumber\\
&&\hspace{21mm}+\frac{\partial^2(\mu^\prime S^i)}{\partial t^2}(f(\tau,\vartheta,u^\prime,\vartheta)
\left(\frac{\partial f}{\partial\tau}\right)^2(\tau,\vartheta)
\label{4.201}
\end{eqnarray}
By the 1st of \ref{4.120} and the fact that $\mu^\prime$ vanishes at $\partial_-{\cal B}$ this 
reduces at $(\tau,u^\prime)=(0,0)$ to:
\begin{equation}
\left.\frac{\partial^3 H^i}{\partial\tau^2\partial u^\prime}\right|_{(\tau,u^\prime)=(0,0)}=
l S^i_0\left.\frac{\partial^2 f}{\partial\tau^2}\right|_{\tau=0}
\label{4.202}
\end{equation}
recalling \ref{4.109}. We have (see 1st of \ref{4.129}):
\begin{equation}
\left.\frac{\partial^2 f}{\partial\tau^2}\right|_{\tau=0}=T\rho|_{\partial_-{\cal B}}
=-\frac{k}{(2-v_0)l}
\label{4.203}
\end{equation}
by \ref{4.167}. Substituting then gives:
\begin{equation}
\left.\frac{\partial^3 H^i}{\partial\tau^2\partial u^\prime}\right|_{(\tau,u^\prime)=(0,0)}=
-\frac{k S^i_0}{(2-v_0)}
\label{4.204}
\end{equation}

Differentiating \ref{4.135} with respect to $\tau$ we obtain:
\begin{equation}
\frac{\partial^3 H^i}{\partial\tau\partial u^{\prime 2}}((\tau,\vartheta),u^\prime)=
\frac{\partial^2(\mu^\prime S^i)}{\partial t\partial u^\prime}(f(\tau,\vartheta),u^\prime,\vartheta)
\frac{\partial f}{\partial\tau}(\tau,\vartheta)
\label{4.205}
\end{equation}
In view of the 1st of \ref{4.120} this reduces at $(\tau,\vartheta)=(0,0)$ to:
\begin{equation}
\left.\frac{\partial^3 H^i}{\partial\tau\partial u^{\prime 2}}\right|_{(\tau,u^\prime)=(0,0)}=0
\label{4.206}
\end{equation}

Finally, differentiating \ref{4.135} with respect to $u^\prime$ we obtain:
\begin{equation}
\frac{\partial^3 H^i}{\partial u^{\prime 3}}((\tau,\vartheta),u^\prime)=
\frac{\partial^2(\mu^\prime S^i)}{\partial u^{\prime 2}}(f(\tau,\vartheta),u^\prime,\vartheta)
\label{4.207}
\end{equation}
Recalling \ref{4.103}, \ref{4.104} this reduces at $(\tau,u^\prime)=(0,0)$ to:
\begin{equation}
\left.\frac{\partial^3 H^i}{\partial u^{\prime 3}}\right|_{(\tau,u^\prime)=(0,0)}=k S^i_0
\label{4.208}
\end{equation}

Recalling the fact that for the solution of the identification equations $u^\prime=O(\tau)$, the 
results \ref{4.200}, \ref{4.204}, \ref{4.206}, \ref{4.208} imply that for the solution of the 
identification equations:
\begin{eqnarray}
&&H^i=\frac{k S^i_0}{6(2-v_0)}\left\{(1+v_0)\tau^3-3\tau^2 u^\prime
+(2-v_0)u^{\prime 3}\right\}\nonumber\\
&&\hspace{15mm}-\frac{(1-v_0)km^A\Omega^{\prime i}_{A 0}}{3(2-v_0)}\tau^3 + O(\tau^4)
\label{4.209}
\end{eqnarray}
Recalling also that for the solution of the identification equations $\varphi^{\prime A}=O(\tau^3)$, by \ref{4.69}, \ref{4.74}, for the solution of the identification equations we have:
\begin{equation}
G^i_A=\Omega^{\prime i}_{A 0}+O(\tau)
\label{4.210}
\end{equation}
Therefore in view of \ref{4.72}, according to the identification equations \ref{4.67} we have:
\begin{eqnarray}
&&\frac{k S^i_0}{6(2-v_0)}\left\{(1+v_0)\tau^3-3\tau^2 u^\prime
+(2-v_0)u^{\prime 3}\right\}\nonumber\\
&&-\frac{(1-v_0)km^A\Omega^{\prime i}_{A 0}}{3(2-v_0)}\tau^3
+\varphi^{\prime A}\Omega^{\prime i}_{A 0}=O(\tau^4)
\label{4.211}
\end{eqnarray}
Substituting \ref{4.141}, multiplying by $\tau^{-3}$, and taking the limit $\tau\rightarrow 0$, 
we deduce: 
\begin{equation}
\frac{k(1-v_0)^3(1+v_0)}{6(2-v_0)}S^i_0+
\left(\gamma_0^A-\frac{km^A(1-v_0)}{3(2-v_0)}\right)\Omega^{\prime i}_{A 0}=0
\label{4.212}
\end{equation}
Since the vectors $S_0$, $\Omega^\prime_A : A=1,...,n-1$ are linearly independent, spanning 
$\Sigma_t$ at each point on $\partial_-{\cal B}$, the coefficients of $S^i_0$ and 
$\Omega^{\prime i}_{A 0}$ must separately vanish. Since $v_0\leq 0$ the first condition is 
equivalent to:
\begin{equation}
v_0=-1
\label{4.213}
\end{equation}
Substituting this in the second condition yields:
\begin{equation}
\gamma_0^A=\frac{2km^A}{9} \ \ : \ A=1,...,n-1
\label{4.214}
\end{equation}
In conclusion, we have established the following proposition. 

\vspace{2.5mm}

\noindent{\bf Proposition 4.4} \ \ \ The solution of the identification equations \ref{4.67} is of the 
form:
$$u^\prime=v\tau, \ \ \ \varphi^{\prime A}=\gamma^A\tau^3$$
with:
$$\left. v\right|_{\tau=0}=-1, \ \ \ \left.\gamma\right|_{\tau=0}=\frac{2km^A}{9}$$

\vspace{2.5mm}

Substituting \ref{4.213} in \ref{4.165} we obtain:
\begin{equation}
Tr|_{\partial_-{\cal B}}=-\frac{l}{2c_0} \ \ \mbox{where we denote: $c_0=c|_{\partial_-{\cal B}}$}
\label{4.215}
\end{equation}
Substituting \ref{4.213} in \ref{4.167} and \ref{4.168} we obtain:
\begin{equation}
T\lambdab|_{\partial_-{\cal B}}=-\frac{c_0 k}{3 l}, \ \ \ 
T^2\lambda|_{\partial_-{\cal B}}=\frac{k}{3}
\label{4.216}
\end{equation}
Substituting finally \ref{4.213} in \ref{4.203} we obtain:
\begin{equation}
\left.\frac{\partial^2 f}{\partial\tau^2}\right|_{\tau=0}=-\frac{k}{3l}
\label{4.217}
\end{equation}
In view of the 1st of \ref{4.120}, this implies:
\begin{equation}
f(\tau,\vartheta)=f(0,\vartheta)-\left(\frac{k}{6l}\right)(\vartheta)\tau^2+O(\tau^3)
\label{4.218}
\end{equation}
On the other hand, by \ref{1.316} and \ref{1.318}:
\begin{equation}
t_*(u^\prime,\vartheta^\prime)=t_*(0,\vartheta^\prime)-
\left(\frac{k}{2l}\right)(\vartheta^\prime)u^{\prime 2}+O(u^{\prime 3})
\label{4.219}
\end{equation}
Since $f(0,\vartheta)=t_*(0,\vartheta)$ and by Proposition 4.4
$$u^\prime=-\tau+O(\tau^2), \ \ \ \vartheta^\prime=\vartheta+O(\tau^3)$$
comparing \ref{4.218} with \ref{4.219} we see that, at least in a neighborhood of ${\cal B}$, 
${\cal K}$ does indeed lie in the past of ${\cal B}$, in accordance with the conditions of the 
shock development problem.

\vspace{5mm}

\section{Regularization of the Identification Equations}

Let us consider again the functions $H^i: i=1,...,n$ defined by \ref{4.71}. Recalling that the
\begin{equation}
N_0^i(\vartheta)=L^{\prime i}_0(\vartheta)=\frac{\partial x^{\prime i}}{\partial t}
(f(0,\vartheta),0,\vartheta)
\label{4.220}
\end{equation}
are known smooth functions on $S^{n-1}$ defined by the prior solution, and that by \ref{4.59} the 
$g^i(0,\vartheta)$ are likewise known smooth functions on $S^{n-1}$ defined by the prior solution, 
let us define the functions:
\begin{equation}
\triangle^i(s,u^\prime,\vartheta)=x^{\prime i}(f(0,\vartheta)+s,u^\prime,\vartheta)
-g^i(0,\vartheta)-N_0^i(\vartheta)s
\label{4.221}
\end{equation}
So these are known smooth functions of their arguments, defined by the prior solution. Setting 
\begin{equation}
\check{f}(\tau,\vartheta)=f(\tau,\vartheta)-f(0,\vartheta), \ \ \ 
\check{g}^i(\tau,\vartheta)=g^i(\tau,\vartheta)-g^i(0,\vartheta)
\label{4.222}
\end{equation}
let us also define the functions:
\begin{equation}
\delta^i(\tau,\vartheta)=\check{g}^i(\tau,\vartheta)-N^i_0(\vartheta)\check{f}(\tau,\vartheta)
\label{4.223}
\end{equation}
So these are a priori unknown functions, depending on the new solution. Comparing \ref{4.222}, 
\ref{4.223} with \ref{4.71} we see that:
\begin{equation}
H^i((\tau,\vartheta),u^\prime)=\triangle^i(\check{f}(\tau,\vartheta),u^\prime,\vartheta)
-\delta^i(\tau,\vartheta)
\label{4.224}
\end{equation}

Consider the functions $\triangle^i(s,u^\prime,\vartheta)$. By \ref{4.59},
\begin{equation}
\triangle^i(0,0,\vartheta)=0
\label{4.225}
\end{equation}
and by \ref{4.220},
\begin{equation}
\frac{\partial\triangle^i}{\partial s}(0,0,\vartheta)=0
\label{4.226}
\end{equation}
According to \ref{4.124} we have:
\begin{equation}
\frac{\partial\triangle^i}{\partial u^\prime}(s,u^\prime,\vartheta)=
(\mu^\prime S^i)(f(0,\vartheta)+s,u^\prime,\vartheta)
\label{4.227}
\end{equation}
In particular, by \ref{4.103},
\begin{equation}
\frac{\partial\triangle^i}{\partial u^\prime}(0,0,\vartheta)=0
\label{4.228}
\end{equation}
Differentiating \ref{4.227} with respect to $s$ and evaluating the result at $(s,u^\prime)=(0,0)$ 
we obtain, in view of \ref{4.109},
\begin{equation}
\frac{\partial^2\triangle^i}{\partial s\partial u^\prime}(0,0,\vartheta)=l(\vartheta)S^i_0(\vartheta)
\label{4.229}
\end{equation}
Differentiating \ref{4.227} with respect to $u^\prime$ and evaluating the result at 
$(s,u^\prime)=(0,0)$ we obtain, by \ref{4.103},
\begin{equation}
\frac{\partial^2\triangle^i}{\partial u^{\prime 2}}(0,0,\vartheta)=0
\label{4.230}
\end{equation}
Let us denote:
\begin{equation}
\frac{\partial^2\triangle^i}{\partial s^2}(0,0,\vartheta)=a_0^i(\vartheta)
\label{4.231}
\end{equation}
These are known smooth functions on $S^{n-1}$, determined by the prior solution. The above results 
\ref{4.225}, \ref{4.226}, \ref{4.228}, \ref{4.229}, \ref{4.230}, \ref{4.231}, imply that:
\begin{equation}
\triangle^i(s,u^\prime,\vartheta)=\frac{1}{2}a_0^i(\vartheta)s^2+l(\vartheta)S^i_0(\vartheta)su^\prime
+\triangle^{\prime i}(s,u^\prime,\vartheta)
\label{4.232}
\end{equation}
where the $\triangle^{\prime i}(s,u^\prime,\vartheta)$ are known smooth functions of their arguments 
whose Taylor expansion in $(s,u^\prime)$ begins with cubic terms. Furthermore, differentiating 
\ref{4.227} twice with respect to $u^\prime$ and evaluating the result at $(s,u^\prime)=(0,0)$ we 
obtain, by \ref{4.103} and \ref{4.104},
\begin{equation}
\frac{\partial^3\triangle^i}{\partial u^{\prime 3}}(0,0,\vartheta)=k(\vartheta)S^i_0(\vartheta)
\label{4.233}
\end{equation}
Let us denote:
\begin{eqnarray}
&&\frac{\partial^3\triangle^i}{\partial s\partial u^{\prime 2}}(0,0,\vartheta)=a_1^i(\vartheta)
\label{4.234}\\
&&\frac{\partial^3\triangle^i}{\partial s^2\partial u^\prime}(0,0,\vartheta)=a_2^i(\vartheta)
\label{4.235}\\
&&\frac{\partial^3\triangle^i}{\partial s^3}(0,0,\vartheta)=a_3^i(\vartheta)
\label{4.236}
\end{eqnarray}
These are all known smooth functions on $S^{n-1}$, determined by the prior solution. Then \ref{4.232} 
is precised to:
\begin{eqnarray}
&&\triangle^i(s,u^\prime,\vartheta)
=\frac{1}{2}a_0^i(\vartheta)s^2+l(\vartheta)S^i_0(\vartheta)su^\prime
\nonumber\\
&&\hspace{15mm}+\frac{1}{6}a_3^i(\vartheta)s^3+\frac{1}{2}a_2^i(\vartheta)s^2 u^\prime
+\frac{1}{2}a_1^i(\vartheta)s u^{\prime 2}+\frac{1}{6}k(\vartheta)S_0^i(\vartheta)u^{\prime 3}
\nonumber\\
&&\hspace{20mm}+\triangle^{\prime\prime i}(s,u^\prime,\vartheta)
\label{4.237}
\end{eqnarray}
where the $\triangle^{\prime\prime i}(s,u^\prime,\vartheta)$ are known smooth functions of 
their arguments whose Taylor expansion in $(s,u^\prime)$ begins with quartic terms.

We shall repeatedly apply the following standard lemma.

\vspace{2.5mm}

\noindent{\bf Lemma 4.1} \ \ \ Let $f$ be a smooth function defined in an open set $U$ in 
$\mathbb{R}^n$ which is star-shaped relative to the origin. Suppose that $f$ and its partial 
derivatives of order up to $k-1$ vanish at the origin. Then there are smooth functions 
$g_{i_1...i_k}$ defined in $U$ such that
$$f(x)=\sum_{i_1,...,i_k}x^{i_1}\cdot\cdot\cdot x^{i_k} g_{i_1...i_k}(x)$$
The functions $g_{i_1...i_k}$ are in fact given by
$$g_{i_1...i_k}(x)=\int_0^1\left\{\int_0^{t_1}\cdot\cdot\cdot\int_0^{t_{k-1}}
\frac{\partial^k f}{\partial x^{i_1}...\partial x^{i_k}}(t_k x)dt_k...dt_2\right\}dt_1$$

\vspace{2.5mm}

We apply the above lemma to the functions $\triangle^{\prime\prime i}$ to conclude that there 
are known smooth functions $A_l^i \ : \ l=0,1,2,3,4$ of $(s,u^\prime,\vartheta)$ such that:
\begin{eqnarray}
&&\triangle^{\prime\prime i}(s,u^\prime,\vartheta)=u^{\prime 4}A_0^i(s,u^\prime,\vartheta)\nonumber\\
&&\hspace{20mm}+s u^{\prime 3}A_1^i(s,u^\prime,\vartheta)
+s^2 u^{\prime 2}A_2^i(s,u^\prime,\vartheta)\nonumber\\
&&\hspace{20mm}+s^3 u^\prime A_3^i(s,u^\prime,\vartheta)+s^4 A_4^i(s,u^\prime,\vartheta)
\label{4.238}
\end{eqnarray}

In view of the fact that:
\begin{equation}
\check{f}(0,\vartheta)=\frac{\partial\check{f}}{\partial\tau}(0,\vartheta)=0
\label{4.239}
\end{equation}
we can set:
\begin{equation}
\check{f}(\tau,\vartheta)=\tau^2\hat{f}(\tau,\vartheta)
\label{4.240}
\end{equation}
In fact by \ref{4.218} we have:
\begin{equation}
\hat{f}(0,\vartheta)=-\frac{k(\vartheta)}{6l(\vartheta)}
\label{4.241}
\end{equation}
Setting then (see \ref{4.224}), according to Proposition 4.4:
$$s=\check{f}(\tau,\vartheta)=\tau^2\hat{f}(\tau,\vartheta), \ \ \ u^\prime=\tau v$$
in \ref{4.237}, \ref{4.238} we obtain:
\begin{eqnarray}
&&\triangle^i(\check{f}(\tau,\vartheta),\tau v,\vartheta)=
\tau^3 S_0^i(\vartheta)v\left[l(\vartheta)\hat{f}(\tau,\vartheta)+\frac{1}{6}k(\vartheta)v^2\right]
\label{4.242}\\
&&\hspace{10mm}+\tau^4\left\{\frac{1}{2}\left(a_0^i(\vartheta)
\hat{f}(\tau,\vartheta)+a_1^i(\vartheta)v^2\right)\hat{f}(\tau,\vartheta)
+\frac{1}{2}\tau a_2^i(\vartheta)v\hat{f}^2(\tau,\vartheta)\right.\nonumber\\
&&\hspace{10mm}+\frac{1}{6}\tau^2 a_3^i(\vartheta)\hat{f}^3(\tau,\vartheta)
+v^4 A_0^i(\tau^2\hat{f}(\tau,\vartheta),\tau v,\vartheta)\nonumber\\
&&\hspace{7.5mm}+\tau v^3\hat{f}(\tau,\vartheta)
A_1^i(\tau^2\hat{f}(\tau,\vartheta),\tau v,\vartheta)+\tau^2 v^2\hat{f}^2(\tau,\vartheta)
A_2^i(\tau^2\hat{f}(\tau,\vartheta),\tau v,\vartheta)\nonumber\\
&&\hspace{7.5mm}\left.+\tau^3 v\hat{f}^3(\tau,\vartheta)
A_3^i(\tau^2\hat{f}(\tau,\vartheta),\tau v,\vartheta)+\tau^4\hat{f}^4(\tau,\vartheta)
A_4^i(\tau^2\hat{f}(\tau,\vartheta),\tau v,\vartheta)\right\}\nonumber
\end{eqnarray}

Consider next the functions $\delta^i$ defined by \ref{4.223}. By \ref{4.120}, \ref{4.131} we have:
\begin{equation}
\delta^i(0,\vartheta)=\frac{\partial\delta^i}{\partial\tau}(0,\vartheta)
=\frac{\partial^2\delta^i}{\partial\tau^2}(0,\vartheta)=0
\label{4.243}
\end{equation}
We can therefore set:
\begin{equation}
\delta^i(\tau,\vartheta)=\tau^3\hat{\delta}^i(\tau,\vartheta)
\label{4.244}
\end{equation}
In fact, by \ref{4.178}:
\begin{eqnarray}
&&\frac{\partial^3\delta^i}{\partial\tau^3}(0,\vartheta)=
-\frac{\partial^3 H^i}{\partial\tau^3}((0,\vartheta),0)\nonumber\\
&&\hspace{17mm}=\frac{4k(\vartheta)m^A(\vartheta)}{3}\Omega^{\prime i}_{A 0}(\vartheta)
\label{4.245}
\end{eqnarray}
by \ref{4.200} and \ref{4.213}. Hence, comparing with \ref{4.214},
\begin{equation}
\hat{\delta}^i(0,\vartheta)=\gamma^A_0(\vartheta)\Omega^{\prime i}_{A 0}(\vartheta)
\label{4.246}
\end{equation}

Consider now the functions $G^i_A((\tau,\vartheta),(u^\prime,\varphi^\prime))$ defined by 
\ref{4.69}:
\begin{equation}
G_A^i((\tau,\vartheta),(u^\prime,\varphi^\prime))=
\int_0^1\Omega^{\prime i}_A(f(\tau,\vartheta),u^\prime,\vartheta+r\varphi^\prime)dr
\label{4.247}
\end{equation}
We have:
\begin{equation}
G_A^i((\tau,\vartheta),(u^\prime,0))=\Omega^{\prime i}_A(f(\tau,\vartheta),u^\prime,\vartheta)
\label{4.248}
\end{equation}
In particular:
\begin{equation}
G_A^i((0,\vartheta),(0,0))=\Omega^{\prime i}_{A 0}(\vartheta)
\label{4.249}
\end{equation}
Let us define the functions:
\begin{equation}
\Gamma_A^i(s,u^\prime,\varphi^\prime,\vartheta)=
\int_0^1\Omega^{\prime i}_A(f(0,\vartheta)+s,u^\prime,\vartheta+r\varphi^\prime)dr
-\Omega^{\prime i}_{A 0}(\vartheta)
\label{4.250}
\end{equation}
These are known smooth functions of their arguments, determined by the prior solution, 
vanishing at $(s,u^\prime,\varphi^\prime)=(0,0,0)$. We can then apply Lemma 4.1 to conclude that 
there are known smooth functions $\Theta_{l A}^i \ : \ l=0,1$ and $\Phi_{AB}^i \ : \ B=1,...,n-1$ 
of $(s,u^\prime,\varphi^\prime,\vartheta)$ such that:
\begin{eqnarray}
&&\Gamma_A^i(s,u^\prime,\varphi^\prime,\vartheta)=
s\Theta_{0 A}^i(s,u^\prime,\varphi^\prime,\vartheta)
+u^\prime\Theta_{1 A}^i(s,u^\prime,\varphi^\prime,\vartheta)\nonumber\\
&&\hspace{25mm}+\varphi^{\prime B}\Phi_{AB}^i(s,u^\prime,\varphi^\prime,\vartheta)
\label{4.251}
\end{eqnarray}
Comparing \ref{4.247} with \ref{4.250} we see that:
\begin{equation}
G_A^i((\tau,\vartheta),(u^\prime,\varphi^\prime))=\Omega^{\prime i}_{A 0}(\vartheta)
+\Gamma_A^i(\check{f}(\tau,\vartheta),u^\prime,\varphi^\prime,\vartheta)
\label{4.252}
\end{equation}
Substituting then \ref{4.240} and
$$u^\prime=\tau v, \ \ \ \varphi^{\prime A}=\tau^3\gamma^A$$
according to Proposition 4.4, we obtain:
\begin{eqnarray}
&&G_A^i((\tau,\vartheta),(\tau v,\tau^3\gamma))=\Omega^{\prime i}_{A 0}(\vartheta)
+\tau\left\{v\Theta_{1 A}^i(\tau^2\hat{f}(\tau,\vartheta),\tau v,\tau^3\gamma,\vartheta)\right.\nonumber\\
&&\hspace{35mm}+\tau\hat{f}(\tau,\vartheta)
\Theta_{0 A}^i(\tau^2\hat{f}(\tau,\vartheta),\tau v,\tau^3\gamma,\vartheta)\nonumber\\
&&\hspace{35mm}+\left.\tau^2\gamma^B\Phi_{AB}^i(\tau^2\hat{f}(\tau,\vartheta),\tau v,\tau^3\gamma,\vartheta)\right\}
\label{4.253}
\end{eqnarray}

In view of the expressions \ref{4.72}, \ref{4.224}, \ref{4.242}, \ref{4.244}, \ref{4.253} we arrive 
at the following proposition. 

\vspace{2.5mm}

\noindent{\bf Proposition 4.5} \ \ \ We have:
$$F^i((\tau,\vartheta),(\tau v,\tau^3\gamma))=\tau^3\hat{F}^i((\tau,\vartheta),(v,\gamma))$$
where:
\begin{eqnarray*}
&&\hat{F}^i((\tau,\vartheta),(v,\gamma))=S_0^i(\vartheta)v\left[l(\vartheta)\hat{f}(\tau,\vartheta)
+\frac{1}{6}k(\vartheta)v^2\right]\\
&&\hspace{26mm}+\Omega^{\prime i}_{A 0}(\vartheta)\gamma^A-\hat{\delta}^i(\tau,\vartheta)
+\tau E^i((\tau,\vartheta),(v,\gamma))
\end{eqnarray*}
the functions $E^i$ being given by:
\begin{eqnarray*}
&&E^i((\tau,\vartheta),(v,\gamma))=\frac{1}{2}\left(a_0^i(\vartheta)\hat{f}(\tau,\vartheta)
+a_1^i(\vartheta)v^2\right)\hat{f}(\tau,\vartheta)\\
&&\hspace{26mm}+\frac{1}{2}\tau a_2^i(\vartheta)v\hat{f}^2(\tau,\vartheta)
+\frac{1}{6}\tau^2 a_3^i(\vartheta)\hat{f}^3(\tau,\vartheta)\\
&&\hspace{26mm}+v^4 A_0^i(\tau^2\hat{f}(\tau,\vartheta),\tau v,\vartheta)
+\tau v^3\hat{f}(\tau,\vartheta)A_1^i(\tau^2\hat{f}(\tau,\vartheta),\tau v,\vartheta)\\
&&\hspace{13mm}+\tau^2 v^2\hat{f}^2(\tau,\vartheta)
A_2^i(\tau^2\hat{f}(\tau,\vartheta),\tau v,\vartheta)
+\tau^3 v\hat{f}^3(\tau,\vartheta) A_3^i(\tau^2\hat{f}(\tau,\vartheta),\tau v,\vartheta)\\
&&\hspace{35mm}+\tau^4\hat{f}^4(\tau,\vartheta) A_4^i(\tau^2\hat{f}(\tau,\vartheta),\tau v,\vartheta)\\
&&\hspace{13mm}+\gamma^A\left\{ v\Theta_{1 A}^i(\tau^2\hat{f}(\tau,\vartheta),\tau v,\tau^3\gamma,\vartheta)+\tau\hat{f}(\tau,\vartheta)\Theta_{0 A}^i(\tau^2\hat{f}(\tau,\vartheta),\tau v, 
\tau^3\gamma,\vartheta)\right.\\
&&\hspace{30mm}+\left.\tau^2\gamma^B\Phi_{AB}^i(\tau^2\hat{f}(\tau,\vartheta),\tau v,\tau^3\gamma,
\vartheta)\right\}
\end{eqnarray*}
where the $a_l^i \ : \ l=0,1,2,3$ are known smooth functions of their argument, the 
$A_l^i \ : \ l=0,1,2,3,4$ are known smooth functions of their three arguments, and the 
$\Theta_{l A}^i \ : \ l=0,1$ and $\Phi_{AB}^i$ are known smooth functions of their four 
arguments, all determined by the prior solution. The identification equations \ref{4.67}, 
with 
$$u^\prime=\tau v, \ \ \ \varphi^\prime=\tau^3\gamma,$$
in accordance with Proposition 4.4, are then equivalent to their regularized form:
$$\hat{F}^i((\tau,\vartheta),(v,\gamma))=0$$

\vspace{2.5mm}

Note that by \ref{4.241}, \ref{4.246} at $(\tau,\vartheta)=(0,\vartheta)$ the solution of 
the regularized identification equations is 
$(v,\gamma)=(-1,\gamma_0(\vartheta))$. 

\pagebreak

\chapter{The Initial and Derived Data}

\section{The Propagation Equations for $\lambdab$ and $s_{NL}$ \\on $\underline{{\cal C}}$}

Let us consider the initial data on $\underline{{\cal C}}=\Cb_0$ for the characteristic system \ref{2.62} and for the wave system 
\ref{2.91}, \ref{2.92}. The unknowns in the characteristic system are the $n+1$ rectangular 
coordinate functions $x^\mu:\mu=0,...,n$ and the $n-1$ components $b^A:A=1,...,n-1$ of the $S$ 
vectorfield $b$. The unknowns in the wave system are the $n+1$ rectangular components 
$\beta_\mu:\mu=0,...,n$ of the 1-form $\beta$. The initial data for the characteristic system 
consist of the specification of the $x^\mu$ on $\Cb_0$, that is of the $x^\mu$ as functions of 
$(u,\vartheta)$ at $\ub=0$. According to the construction in Section 2.5, $\vartheta$ is constant 
along the generators of $\Cb_0$. These being integral curves of $\Lb$ and $\Lb$ being given by 
the 2nd of \ref{2.94}, this means that
\begin{equation}
\left.b^A\right|_{\Cb_0}=0 \ : \ A=1,...,n-1
\label{5.1}
\end{equation}
The 2nd of the characteristic system \ref{2.62} then reduces on $\Cb_0$ to:
\begin{equation}
\frac{\partial x^i}{\partial u}=\Nb^i\frac{\partial t}{\partial u} \ : \ i=1,...,n \ \ \ 
\mbox{: on $\Cb_0$}
\label{5.2}
\end{equation}
The initial data for the characteristic system determine the rectangular components (see \ref{2.47})
\begin{equation}
\Omega_A^\mu=\frac{\partial x^\mu}{\partial\vartheta^A} \ \ : \ A=1,...,n-1
\label{5.3}
\end{equation}
of the vectorfields $\Omega_A:A=1,...,n-1$ along $\Cb_0$. 
The initial data for the wave system consist of the specification of the $\beta_\mu$ on $\Cb_0$, 
that is of the $\beta_\mu$ as functions of $(u,\vartheta)$ at $\ub=0$. These determine the rectangular 
components $h_{\mu\nu}$ of the acoustical metric along $\Cb_0$, as (see \ref{1.d11})
\begin{equation}
h_{\mu\nu}=g_{\mu\nu}+H\beta_\mu\beta_\nu, \ \ H=H(\sigma), 
\ \ \sigma=-(g^{-1})^{\mu\nu}\beta_\mu\beta_\nu
\label{5.4}
\end{equation}
According to Section 2.2, the $N^\mu$, $\Nb^\mu$, the rectangular components of the vectorfields 
$N$, $\Nb$, are then determined along $\Cb_0$. Consequently, equations \ref{5.2} constitute $n$ 
constraints on the initial data for the characteristic system. We can thus say that in regard to  
the initial data for the characteristic system there is only one free function of $(u,\vartheta)$, 
and we can think of this function as representing $t$ on $\Cb_0$. Furthermore, the components of the 
induced metric (see \ref{2.56})
\begin{equation}
\sh_{AB}=h_{\mu\nu}\Omega_A^\mu\Omega_B^\nu
\label{5.5}
\end{equation}
on the $S_{0,u}$ sections of $\Cb_0$, are determined, as is the positive function (see \ref{2.71})
\begin{equation}
c=-\frac{1}{2}h_{\mu\nu}N^\mu\Nb^\nu
\label{5.6}
\end{equation}
on $\Cb_0$. Also, $\rhob=\Lb t$ is determined:
\begin{equation}
\rhob=\frac{\partial t}{\partial u} \ \mbox{: on $\Cb_0$}
\label{5.7}
\end{equation}
hence so is $\lambda=c\rhob$ (see \ref{2.74}). The last two of \ref{2.92} along $\Cb_0$ read:
\begin{eqnarray}
&&\rhob\Nb^\mu\frac{\partial\beta_\mu}{\partial\vartheta^A}-
\Omega_A^\mu\frac{\partial\beta_\mu}{\partial u}=0 \ : \ A=1,...,n-1 \nonumber\\
&&\Omega_A^\mu\frac{\partial\beta_\mu}{\partial\vartheta^B}-
\Omega_B^\mu\frac{\partial\beta_\mu}{\partial\vartheta^A}=0 \ : \ A,B=1,...,n-1
\label{5.8}
\end{eqnarray}
Moreover, by the 1st of \ref{2.92}, and the fact that (see \ref{2.72}) $a=c\rho\rhob$, equation 
\ref{2.91} takes along $\Cb_0$ the form:
\begin{equation}
N^\mu\frac{\partial\beta_\mu}{\partial u}-
c\rhob(\sh^{-1})^{AB}\Omega_A^\mu\frac{\partial\beta_\mu}{\partial\vartheta^B}=0
\label{5.9}
\end{equation}
Equations \ref{5.8} and \ref{5.9} constitute $n(n-1)/2$ plus 1 algebraically independent 
constraints on the initial data for the wave system. However only $n$ of these constraints are 
differentially independent. We can thus say that also in regard to  the initial data for the 
wave system there is only one free function of $(u,\vartheta)$, and we can think of this function 
as representing $\phi$ on $\Cb_0$. For our problem the initial data on $\Cb_0=\underline{{\cal C}}$ 
is induced by the prior solution and the constraints are automatically satisfied. 

According to the above, of the 1st derivatives of the $x^\mu$ along $\Cb_0$ what is not yet 
determined  is the transversal derivative
\begin{equation}
\frac{\partial t}{\partial\ub}=\rho \ \ \mbox{: along $\Cb_0$}
\label{5.10}
\end{equation}
or equivalently the function $\lambdab=c\rho$. Once this is determined the transversal derivatives 
$\partial x^i/\partial\ub \ : \ i=1,...,n$ would be determined through the 1st of the 
characteristic system \ref{2.62}, which reduces along $\Cb_0$ to:
\begin{equation}
\frac{\partial x^i}{\partial\ub}=N^i\frac{\partial t}{\partial\ub} \ : \ i=1,...,n \ \ \ 
\mbox{: on $\Cb_0$}
\label{5.11}
\end{equation}
Also, of the 1st derivatives of the $\beta_\mu$ along $\Cb_0$, these being expressed by the 
components of the tensorfield $s$ along $\Cb_0$ (see \ref{3.36} - \ref{3.43}), what is not yet 
determined is $s_{NL}$, the transversal component of the transversal derivative. 

Now $\lambdab$ satisfies along the generators of the $\Cb_{\ub}$, in particular the generators 
of $\Cb_0$, the propagation equation given by Proposition 3.3. On $\Cb_0$ this reads:
\begin{equation}
\frac{\partial\lambdab}{\partial u}=\pb\lambdab+\qb\lambda
\label{5.12}
\end{equation}
where:
\begin{equation}
\pb=\omb-\frac{1}{2c}\left(\beta_N\beta_{\Nb}\Lb H+H\beta_N s_{\Nb\Lb}\right), \ \ \
\qb=\frac{1}{4c}\beta^2_{\Nb}LH
\label{5.13}
\end{equation}
Here it is important to distinguish the quantities which are already determined by the data 
induced on $\Cb_0$ by the prior solution according to the above discussion, from the 
quantities not yet determined as they depend on the pair $(\lambdab,s_{NL})$. Recall that 
$\omb$ is given in terms of $\mb$, $\smb$ by \ref{3.47} (recall \ref{3.45}) and $\mb$, $\smb$ 
are given by \ref{3.58}, \ref{3.69}. We see that these are already determined, hence so is 
$\omb$ and $\pb$. On the other hand $\qb$ is not yet determined. We have:
\begin{equation}
LH=\frac{dH}{d\sigma}L\sigma, \ \ \ L\sigma=-2\beta^\mu L\beta_\mu
\label{5.14}
\end{equation}
Substituting for $\beta^\mu$ from \ref{4.161} we obtain:
\begin{equation}
L\sigma=\frac{\eta^2}{c}\left[\left(\beta_N\mbox{tr}\sss
-2(\sh^{-1})^{AB}\sbeta_B(\ss_N)_A\right)\lambdab+\beta_{\Nb}s_{NL}\right]
\label{5.15}
\end{equation}
This is a linear combination of $\lambdab$ and $s_{NL}$, hence so is $LH$ and $\qb$. 
We conclude that the right hand side of \ref{5.12} is a linear combination of $\lambdab$ and 
$s_{NL}$, that is this equation is of the form:
\begin{equation}
\frac{\partial\lambdab}{\partial u}=a_{11}\lambdab+a_{12}s_{NL} \ \ \mbox{: on $\Cb_0$}
\label{5.16}
\end{equation}
where the coefficients $a_{11}$, $a_{12}$ are smooth functions on $\Cb_0$ determined by 
the data induced by the prior solution. 

Before we proceed to derive the propagation equation for $s_{NL}$ along the generators of 
$\Cb_0$, we first consider on $\Cb_0$ the 1st order acoustical quantities $\tchi$ , $\tchib$ 
which are related to the 2nd fundamental forms $\chi$, $\chib$ by \ref{3.75}, \ref{3.78} 
respectively, as well as the $S$ 1-forms $\teta$, $\tetab$ which are related to the torsion 
forms $\eta$, $\etab$ by \ref{3.79}, \ref{3.82} respectively. From the definitions of the structure 
functions $\sk_A^B$, $k_A$, $\ok_A$ through the expansion \ref{3.29} and the definition 
of the structure functions $\skb_A^B$, $\kb_A$, $\okb_A$ through the expansion \ref{3.30}, 
we conclude, in view of the fact that the $N^\mu$, $\Nb^\mu$ are known smooth functions of 
$(u,\vartheta)$ on $\Cb_0$, that so are these structure functions. This is true in particular of 
$\sk_A^B$, $\skb_A^B$. Then the relations \ref{3.76}, \ref{3.77} show that $\tchi$, $\tchib$ are on 
$\Cb_0$ known smooth $S$ tensorfields.  Consequently, by the relation \ref{3.78} $\chib$ is on 
$\Cb_0$ a known smooth $S$ tensorfield, while by the relation \ref{3.75} and the 
expression \ref{5.15} $\chi$ is on $\Cb_0$ a linear combination of $\lambdab$ and $s_{NL}$ with 
coefficients which are known smooth $S$ tensorfields. Proposition 3.2 shows that on $\Cb_0$ 
$\teta$ is a known smooth $S$ 1-form, while $\tetab-\sd\lambdab$ is a linear combination of 
$\lambdab$ and $s_{NL}$ with coefficients which are known smooth $S$ 1-forms. Then by the relation 
\ref{3.79} and the expression \ref{5.15} $\eta$ is on $\Cb_0$ a linear combination of $\lambdab$ 
and $s_{NL}$ with coefficients which are known smooth $S$ 1-forms, while by the relation \ref{3.82} 
$\etab-\rhob\sd\lambdab$ is on $\Cb_0$ such a linear combination. 

We now proceed to derive the propagation equation for $s_{NL}$ along the generators of $\Cb_0$. 
We recall from Chapter 2 that the rectangular components $\beta_\mu$ of the 1-form $\beta$ satisfy 
the wave equation associated to the conformal acoustical metric:
\begin{equation}
\square_{\tilde{h}}\beta_\mu=0 \ \ \ \mbox{with} \ \ \tilde{h}=\Omega h, \ 
\Omega=\left(\frac{G}{\eta}\right)^{2/(n-1)}
\label{5.17}
\end{equation}
(see \ref{2.108}, \ref{2.110}). Now, for any function $f$ on the $n+1$ dimensional spacetime manifold 
${\cal N}$ we have:
\begin{equation}
\square_{\tilde{h}}f=\Omega^{-1}\left(\square_h f+\frac{(n-1)}{2}h^{-1}(d\log\Omega,df)\right)
\label{5.18}
\end{equation}
Applying this to the functions $\beta_\mu$, we conclude that by virtue of \ref{5.17}:
\begin{equation}
\square_h\beta_\mu+\frac{(n-1)}{2}h^{-1}(d\log\Omega,d\beta_\mu)=0
\label{5.19}
\end{equation}
In regard to the 2nd term, since $\Omega=\Omega(\sigma)$ we have:
$$d\log\Omega=\frac{d\log\Omega}{d\sigma}d\sigma$$
hence, substituting the expansion \ref{2.38} of $h^{-1}$, 
\begin{equation}
h^{-1}(d\log\Omega,d\beta_\mu)=-\frac{1}{2a}\left((\Lb\sigma)L\beta_\mu+(L\sigma)\Lb\beta_\mu\right)
+(\sh^{-1})^{AB}(\Omega_A\sigma)\Omega_B\beta_\mu
\label{5.20}
\end{equation}

Now, for any function $f$ on ${\cal N}$, $\square_h f$ is the trace relative to $h$ of $D^2 f=Ddf$, 
the covariant relative to $h$ Hessian of $f$:
\begin{equation}
\square_h f=\mbox{tr}(h^{-1}\cdot D^2 f)
\label{5.21}
\end{equation}
Substituting again the expansion \ref{2.38} of $h^{-1}$ and noting that $D^2 f$ is symmetric 
we obtain:
\begin{equation}
\square_h f=-a^{-1}(D^2 f)(\Lb,L)+(\sh^{-1})^{AB}(D^2 f)(\Omega_A,\Omega_B)
\label{5.22}
\end{equation}
Now, for any pair $X,Y$ of vectorfields on ${\cal N}$ we have:
\begin{equation}
(D^2 f)(X,Y)=X(Yf)-(D_X Y)f
\label{5.23}
\end{equation}
Applying this taking $X=\Lb$, $Y=L$ we obtain, in view of \ref{3.14},
\begin{equation}
(D^2 f)(\Lb,L)=\Lb(Lf)-2\eta^\sharp\cdot\sd f
\label{5.24}
\end{equation}
Applying \ref{5.23} taking $X=\Omega_A$, $Y=\Omega_B$ we obtain:
\begin{equation}
(D^2 f)(\Omega_A,\Omega_B)=\Omega_A(\Omega_B f)-(D_{\Omega_A}\Omega_B)f
\label{5.25}
\end{equation}
On the other hand a formula similar to \ref{5.23} holds for $\sD^2 f$, the symmetric 2-covariant 
$S$ tensorfield which is the covariant Hessian of the restriction of $f$ to each of the $S_{\ub,u}$ 
relative to the induced metric $\sh$: 
\begin{equation}
(\sD^2 f)(X,Y)=X(Yf)-(\sD_X Y)f
\label{5.26}
\end{equation}
where now $X,Y$ are $S$ vectorfields. Taking then $X=\Omega_A$, $Y=\Omega_B$ we obtain:
\begin{equation}
(\sD^2 f)(\Omega_A,\Omega_B)=\Omega_A(\Omega_B f)-(\sD_{\Omega_A}\Omega_B)f
\label{5.27}
\end{equation}
Subtracting \ref{5.27} from \ref{5.25} and appealing to \ref{3.26} then yields:
\begin{equation}
(D^2 f)(\Omega_A,\Omega_B)=(\sD^2 f)(\Omega_A,\Omega_B)
-\frac{1}{2a}\left(\chib_{AB}Lf+\chi_{AB}\Lb f\right)
\label{5.28}
\end{equation}
Substituting finally \ref{5.24} and \ref{5.28} in \ref{5.22} we obtain the following expression:
\begin{equation}
\square_h f=\frac{1}{a}\left\{-\Lb(Lf)+2(\sh^{-1})^{AB}\eta_B\Omega_A f
-\frac{1}{2}\left(\mbox{tr}\chi\Lb f+\mbox{tr}\chib Lf\right)\right\}+\slap f
\label{5.29}
\end{equation}
where $\slap f=\mbox{tr}(\sh^{-1}\cdot\sD^2 f)$ is the Laplacian of the restriction of $f$ to each of 
the $S_{\ub,u}$ relative to the induced metric $\sh$. 

Substituting the expression \ref{5.29} with $\beta_\mu$ in the role of $f$ in \ref{5.19}, substituting also \ref{5.20} and multiplying the resulting equation by $-a=-c^{-1}\lambda\lambdab$ yields the 
equation:
\begin{eqnarray}
&&\Lb(L\beta_\mu)-2(\sh^{-1})^{AB}\eta_B\Omega_A\beta_\mu
+\frac{1}{2}\left(\mbox{tr}\chi\Lb\beta_\mu+\mbox{tr}\chib L\beta_\mu\right)
-c^{-1}\lambda\lambdab\slap\beta_\mu\nonumber\\
&&+\frac{(n-1)}{2}\frac{d\log\Omega}{d\sigma}\left[\frac{1}{2}\left((\Lb\sigma)L\beta_\mu+
(L\sigma)\Lb\beta_\mu\right)-c^{-1}\lambda\lambdab
(\sh^{-1})^{AB}(\Omega_B\sigma)\Omega_A\beta_\mu\right]\nonumber\\
&&\hspace{90mm}=0
\label{5.30}
\end{eqnarray}
This equation holds in particular on $\Cb_0$, where by \ref{5.15} and the preceding discussion 
$\eta_B$, $\mbox{tr}\chi$, and $L\sigma$ are linear combinations of $\lambdab$ and $s_{NL}$ with 
coefficients which are known smooth functions of $(u,\vartheta)$ defined by the data induced on 
$\Cb_0$ by the prior solution. We see then that on $\Cb_0$ \ref{5.30} is of the form:
\begin{equation}
\Lb(L\beta_\mu)+f_\mu\lambdab+e_\mu s_{NL}+gL\beta_\mu=0
\label{5.31}
\end{equation}
where the coefficients $f_\mu$, $e_\mu$, and $g$ are known smooth functions of $(u,\vartheta)$ on 
$\Cb_0$. Let us multiply this equation by $N^\mu$. Since
$$N^\mu\Lb(L\beta_\mu)=\Lb s_{NL}-(\Lb^\mu)L\beta_\mu,$$
substituting the expansion \ref{3.33} we conclude that $s_{NL}$ satisfies on $\Cb_0$ the equation:
\begin{equation}
\frac{\partial s_{NL}}{\partial u}=a_{21}\lambdab+a_{22}s_{NL} \ \ \mbox{: on $\Cb_0$}
\label{5.32}
\end{equation}
where the coefficients $a_{21}$, $a_{22}$ are given by:
\begin{eqnarray}
&&a_{21}=c^{-1}\sn^A(\ss_N)_A+\mbox{tr}\sss\on-N^\mu f_\mu\nonumber\\
&&a_{22}=n-g-N^\mu e_\mu
\label{5.33}
\end{eqnarray}
(see \ref{3.42}). By \ref{3.a6}, \ref{3.61} and the 1st of \ref{3.48} these are smooth functions on $\Cb_0$ defined by the data induced by the prior solution. Recalling from Chapter 4 the results 
\ref{4.118}, \ref{4.175} we arrive at the following proposition. 

\vspace{2.5mm}

\noindent{\bf Proposition 5.1} \ \ \ The pair of functions $(\lambdab,s_{NL})$ satisfy along the 
generators of $\Cb_0$ the linear homogeneous system of ordinary differential equations:
\begin{eqnarray*}
&&\frac{\partial\lambdab}{\partial u}=a_{11}\lambdab+a_{12}s_{NL}\\
&&\frac{\partial s_{NL}}{\partial u}=a_{21}\lambdab+a_{22}s_{NL}
\end{eqnarray*}
where the coefficients constitute a matrix function
$$\left(
\begin{array}{lll}
a_{11}&a_{12}\\
a_{21}&a_{22}
\end{array}
\right)$$
which is a known smooth matrix function of $(u,\vartheta)$ defined by the data induced on $\Cb_0$ by 
the prior solution. Since a given generator corresponds to a given $\vartheta\in S^{n-1}$ and 
initiates at $u=0$ at the corresponding point on $S_{0,0}=\partial_-{\cal B}$, by virtue of 
\ref{4.118}, \ref{4.175} the initial data for the pair $(\lambdab,s_{NL})$ vanish. {\em Consequently,  
this pair of functions vanishes everywhere along $\Cb_0$.}

\vspace{2.5mm}

We obtain a number of vanishing results on $\Cb_0$ as corollaries. 
First, $L\beta_\mu$ vanish on $\Cb_0$ as follows by considering the evaluation on $N^\mu$, $\Nb^\mu$, 
and the $\Omega_A^\mu:A=1,...,n-1$. In fact, on any $\Cb_{\ub}$ we can express:
\begin{equation}
L\beta_\mu=c^{-1}h_{\mu\nu}\left[\left((\sh^{-1})^{AB}(\ss_N)_A\Omega_B^\nu-\frac{1}{2}\mbox{tr}\sss N^\nu\right)\lambdab-\frac{1}{2}\Nb^\nu s_{NL}\right]
\label{5.a1}
\end{equation}
Next, the quantities which were shown to be linear combinations 
of $\lambdab$ and $s_{NL}$ on $\Cb_0$, such as $\chi$, $\eta$, $\tetab-\sd\lambdab$, 
$\etab-\rhob\sd\lambdab$, all vanish everywhere along $\Cb_0$. But since $\lambdab$ vanishes on 
$\Cb_0$, so do $\tetab$ and $\etab$. Then in view of \ref{3.16} $Z$ likewise vanishes on $\Cb_0$:
\begin{equation}
\left.Z\right|_{\Cb_0}=0
\label{5.34}
\end{equation} 
which by \ref{2.35} means:
\begin{equation}
\left.\frac{\partial b^A}{\partial\ub}\right|_{\Cb_0}=0 \ : \ A=1,...,n-1
\label{5.35}
\end{equation}
Next by \ref{3.57}, \ref{3.68}, and the 1st of \ref{3.47}, $\sm$, $\om$, $m$ all vanish on $\Cb_0$. 
Then by the expansion \ref{3.31} so does $LN^\mu$:
\begin{equation}
\left.LN^\mu\right|_{\Cb_0}=0
\label{5.36}
\end{equation}
It then follows through Proposition 3.3 that 
$p$ and $L\lambda$ vanish on $\Cb_0$:
\begin{equation}
\left.L\lambda\right|_{\Cb_0}=0
\label{5.37}
\end{equation}
Also, by \ref{3.a7}, \ref{3.62}, and the 2nd of \ref{3.48}, 
$\snb$, $\nb$, $\onb$ all vanish on $\Cb_0$. Then by the expansion \ref{3.34} so does $L\Nb^\mu$:
\begin{equation}
\left.L\Nb^\mu\right|_{\Cb_0}=0
\label{5.38}
\end{equation} 
It follows by \ref{3.85} that $Lc$ vanishes on $\Cb_0$:
\begin{equation}
\left.Lc\right|_{\Cb_0}=0
\label{5.39}
\end{equation}
Next, by \ref{3.44} and Lemma 3.2 the $S$ 
1-form $\sL_L\pi$ vanishes on $\Cb_0$:
\begin{equation}
\left.\sL_L\pi\right|_{\Cb_0}=0
\label{5.40}
\end{equation}
By Proposition 3.4 and 3.5 the 2-covariant $S$ tensorfields 
$\sL_L(\sk\cdot\sh)$ and $\sL_L(\skb\cdot\sh)$ vanish on $\Cb_0$. By the relations 
\ref{3.76}, \ref{3.77}, in view of what has already been shown, to demonstrate the vanishing of 
symmetric 2-covariant $S$ tensorfields $\sL_L\tchi$ and $\sL_L\tchib$ on $\Cb_0$, it suffices to 
show that the $S$ 1-forms $\sL_L\sbeta$, $\sL_L\ss_N$, $\sL_L\ss_{\Nb}$ all vanish on $\Cb_0$. 
In the case of the first we have:
\begin{eqnarray*}
&&(\sL_L\sbeta)(\Omega_A)=({\cal L}_L\sbeta)(\Omega_A)=L(\beta(\Omega_A))-\beta([L,\Omega_A])\\
&&\hspace{12mm}=L(\beta_\mu\Omega_A^\mu)-\beta_\mu(L\Omega_A^\mu-\Omega_A L^\mu)
=\beta_\mu\Omega_A L^\mu\\
&&\hspace{12mm}=\beta_\mu\Omega_A(\rho N^\mu)=0 \ \ \ \mbox{: on $\Cb_0$}
\end{eqnarray*}
hence: 
\begin{equation}
\left.\sL_L\sbeta\right|_{\Cb_0}=0
\label{5.41}
\end{equation}
by the vanishing of $L\beta_\mu$ and $\rho$ on $\Cb_0$. In the case of the second we have:
\begin{eqnarray*}
&&(\sL_L\ss_N)(\Omega_A)=L(\ss_N(\Omega_A))-\ss_N([L,\Omega_A])\\
&&\hspace{12mm}=L(N^\mu\Omega_A\beta_\mu)-N^\mu(L\Omega_A-\Omega_A L)\beta_\mu\\
&&\hspace{12mm}=N^\mu L\Omega_A\beta_\mu-N^\mu(L\Omega_A-\Omega_A L)\beta_\mu\\
&&\hspace{12mm}=N^\mu\Omega_A L\beta_\mu=0 \ \ \ \mbox{: on $\Cb_0$}
\end{eqnarray*}
hence:
\begin{equation}
\left.\sL_L\ss_N\right|_{\Cb_0}=0
\label{5.42}
\end{equation}
by the vanishing of $LN^\mu$ and $L\beta_\mu$ on $\Cb_0$. Similarly, in the case of the third 
we have:
\begin{eqnarray*}
&&(\sL_L\ss_{\Nb})(\Omega_A)=L(\ss_{\Nb}(\Omega_A))-\ss_{\Nb}([L,\Omega_A])\\
&&\hspace{12mm}=L(\Nb^\mu\Omega_A\beta_\mu)-\Nb^\mu(L\Omega_A-\Omega_A L)\beta_\mu\\
&&\hspace{12mm}=\Nb^\mu L\Omega_A\beta_\mu-\Nb^\mu(L\Omega_A-\Omega_A L)\beta_\mu\\
&&\hspace{12mm}=\Nb^\mu\Omega_A L\beta_\mu=0 \ \ \ \mbox{: on $\Cb_0$}
\end{eqnarray*}
hence:
\begin{equation}
\left.\sL_L\ss_N\right|_{\Cb_0}=0
\label{5.43}
\end{equation}
by the vanishing of $L\Nb^\mu$ and $L\beta_\mu$ on $\Cb_0$. We conclude that indeed $\sL_L\tchi$ and  
$\sL_L\tchib$ vanish everywhere along $\Cb_0$: 
\begin{equation}
\left.\sL_L\tchi\right|_{\Cb_0}=\left.\sL_L\tchib\right|_{\Cb_0}=0
\label{5.44}
\end{equation}
From the 2nd of these, together with \ref{5.34}, \ref{5.37}, \ref{5.39}, \ref{5.41}, in view of 
the relation \ref{3.78} we deduce:
\begin{equation}
\left.\sL_L\chib\right|_{\Cb_0}=0
\label{5.45}
\end{equation}
We also note that
\begin{equation}
\left.L\beta_N\right|_{\Cb_0}=\left.L\beta_{\Nb}\right|_{\Cb_0}=0
\label{5.46}
\end{equation}
Moreover, we have
\begin{eqnarray*}
&&(\sL_L\sss)(\Omega_A,\Omega_B)=L(\sss(\Omega_A,\Omega_B))-\sss([L,\Omega_A],\Omega_B)
-\sss(\Omega_A,[L,\Omega_B])\\
&&\hspace{12mm}=L(\Omega_A^\mu\Omega_B\beta_\mu)-(L\Omega_A^\mu-\Omega_A L^\mu)\Omega_B\beta_\mu
-\sss(\Omega_A,[L,\Omega_B])\\
&&\hspace{12mm}=(\Omega_A L^\mu)(\Omega_B\beta_\mu)+\Omega_A^\mu L\Omega_B\beta_\mu-
\Omega_A^\mu[L,\Omega_B]\beta_\mu\\
&&\hspace{12mm}=(\Omega_A(\rho N^\mu))(\Omega_B\beta_\mu)+\Omega_A^\mu\Omega_B(L\beta_\mu)
\end{eqnarray*}
Since both $\rho$ and $L\beta_\mu$ vanish on $\Cb_0$, we conclude that:
\begin{equation}
\left.\sL_L\sss\right|_{\Cb_0}=0
\label{5.47}
\end{equation}
Also, since 
\begin{eqnarray*}
&&Ls_{\Nb\Lb}=\Nb^\mu L\Lb\beta_\mu+(L\Nb^\mu)\Lb\beta_\mu\\
&&\hspace{10mm}=\Nb^\mu\Lb L\beta_\mu-\ss_N\cdot Z+(L\Nb^\mu)\Lb\beta_\mu
\end{eqnarray*}
the above imply:
\begin{equation}
\left.Ls_{\Nb\Lb}\right|_{\Cb_0}=0
\label{5.48}
\end{equation}

Summarizing the above vanishing results on $\Cb_0$ we can say that the quantities which represent 
1st derivatives transversal to the $\Cb_{\ub}$ can be expressed in terms of $\lambdab$, $s_{NL}$ 
and their derivatives tangential to the $\Cb_{\ub}$, hence all these quantities vanish on $\Cb_0$ by 
virtue of Proposition 5.1. 

We conclude the present section with a formula for $\sL_T Z$ at $\partial_-{\cal B}$. By \ref{5.34} 
this coincides with $\sL_L Z$ at $\partial_-{\cal B}$. The vanishing of $\chi$ on $\Cb_0$ 
is equivalent, according to Proposition 3.1, to the vanishing of $\sL_L\sh$ on $\Cb_0$. Then by 
\ref{3.16}:
\begin{equation}
\sL_L Z=2\left((\sL_L\eta)^\sharp-(\sL_L\etab)^\sharp\right) \ \ \mbox{: on $\Cb_0$}
\label{5.a2}
\end{equation}
Besides $\rho$, by \ref{5.37}, \ref{5.39} also $L\rhob$ vanishes on $\Cb_0$. Moreover $\rho$ itself 
vanishes at $\partial_-{\cal B}$. It then follows from \ref{3.79}, \ref{3.82} that at 
$\partial_-{\cal B}$ we have:
$$\sL_L\eta=(L\rho)\teta, \ \ \ \sL_L\etab=\frac{1}{4}(L\rho)\beta_N\sbeta\Lb H$$ 
Moreover, by Proposition 3.2,
$$\teta=-\frac{1}{4}\beta_N\sbeta\Lb H \ \ \mbox{: at $\partial_-{\cal B}$}$$
Therefore:
\begin{equation}
\left.\sL_L\etab\right|_{\partial_-{\cal B}}=-\left.\sL_L\eta\right|_{\partial_-{\cal B}}
=\left.\frac{1}{4}(L\rho)\beta_N\sbeta\Lb H\right|_{\partial_-{\cal B}}
\label{5.a3}
\end{equation}
Since on $\Cb_0$ $L\rho=T\rho$ and by the 1st of \ref{4.216}
$$\left.T\rho\right|_{\partial_-{\cal B}}=-\frac{k}{3l},$$
in view of \ref{4.112} and \ref{4.186}, we have: 
\begin{equation}
\left.(\sL_L\etab)^\sharp\right|_{\partial_-{\cal B}}
=-\left(\sL_L\eta)^\sharp\right|_{\partial_-{\cal B}}
=-\frac{k}{6}m
\label{5.a4}
\end{equation}
We then conclude that
\begin{equation}
\left.\sL_T Z\right|_{\partial_-{\cal B}}=\frac{2k}{3}m
\label{5.a5}
\end{equation}
a relation which, together with \ref{4.214}, elucidates the role of the vectorfield $m$ on 
$\partial_-{\cal B}$.

\vspace{5mm}

\section{The Propagation Equations for the Higher Order Derived Data 
$T^m\lambdab$ and $T^m s_{NL}$, $m\geq 1$, on $\underline{{\cal C}}$} 

Recall that $\lambdab$ satisfies along the generators of each $\Cb_{\ub}$ the propagation equation 
of Proposition 3.3. Following the argument leading to \ref{5.16} we similarly deduce along the 
generators of each $\Cb_{\ub}$ the propagation equation
\begin{equation}
\Lb\lambdab=a_{11}\lambdab+a_{12}s_{NL}
\label{5.49}
\end{equation}
where the coefficient functions $a_{11}$, $a_{12}$ are given on any $\Cb_{\ub}$ by:
\begin{eqnarray}
&&a_{11}=\pb+\frac{\lambda\eta^2}{4c^2}\frac{dH}{d\sigma}\beta_{\Nb}^2\left(\beta_N\mbox{tr}\sss
-2\sh^{-1}(\sbeta,\ss_N)\right)\nonumber\\
&&a_{12}=\frac{\lambda\eta^2}{4c^2}\frac{dH}{d\sigma}\beta_{\Nb}^3
\label{5.50}
\end{eqnarray}
(here $\eta$ is the sound speed). Also, following the argument leading to \ref{5.32} we similarly 
deduce along the generators of each $\Cb_{\ub}$ the propagation equation
\begin{equation}
\Lb s_{NL}=a_{21}\lambdab+a_{22}s_{NL}
\label{5.51}
\end{equation}
where the coefficient functions $a_{21}$, $a_{22}$ are given on any $\Cb_{\ub}$ by \ref{5.33}, 
that is:
\begin{eqnarray}
&&a_{21}=c^{-1}\ss_N\cdot\sn+\mbox{tr}\sss\on-f_N\nonumber\\
&&a_{22}=n-g-e_N
\label{5.52}
\end{eqnarray}
where we denote $f_N=N^\mu f_\mu$, $e_N=N^\mu e_\mu$. Following the argument leading to \ref{5.31} 
we in fact find:
\begin{equation}
g=\frac{1}{2}\mbox{tr}\chib+\frac{(n-1)}{4}\frac{d\log\Omega}{d\sigma}\Lb\sigma
\label{5.53}
\end{equation}
\begin{eqnarray}
&&f_N=-c^{-1}\lambda N^\mu\slap\beta_\mu-2c^{-1}\sh^{-1}(\teta,\ss_N)\nonumber\\
&&\hspace{5mm}+\frac{1}{2}c^{-1}\lambda\mbox{tr}\tchi\mbox{tr}\sss
-\frac{(n-1)}{2}c^{-1}\lambda\frac{d\log\Omega}{d\sigma}\sh^{-1}(\sd\sigma,\ss_N)\nonumber\\
&&\hspace{5mm}+\frac{1}{4}c^{-1}\lambda\eta^2\left(\beta_N\mbox{tr}\sss-2\sh^{-1}(\sbeta,\ss_N)\right)
\left[-2c^{-1}\beta_{\Nb}\frac{dH}{d\sigma}\sh^{-1}(\sbeta,\ss_N)\right.\nonumber\\
&&\hspace{40mm}+\left.\left(|\sbeta|^2\frac{dH}{d\sigma}+(n-1)\frac{d\log\Omega}{d\sigma}\right)
\mbox{tr}\sss\right]
\label{5.54}
\end{eqnarray}
\begin{eqnarray}
&&e_N=\frac{1}{4}c^{-1}\lambda\eta^2\beta_{\Nb}\left[-2c^{-1}\beta_{\Nb}\frac{dH}{d\sigma}\sh^{-1}(\sbeta,\ss_N)\right.\nonumber\\
&&\hspace{30mm}+\left.\left(|\sbeta|^2\frac{dH}{d\sigma}+(n-1)\frac{d\log\Omega}{d\sigma}\right)
\mbox{tr}\sss\right]
\label{5.55}
\end{eqnarray}
(where $|\sbeta|^2=\sh^{-1}(\sbeta,\sbeta)$). 

In regard to the 1st of \ref{5.50}, from Proposition 3.3:
\begin{equation}
\pb=\omb-\frac{1}{2c}\left(\beta_N\beta_{\Nb}\Lb H+H\beta_N s_{\Nb\Lb}\right)
\label{5.56}
\end{equation}
and $\omb$ is given through $\smb$, $\mb$ by the 2nd of \ref{3.47}, while $\mb$ is given by \ref{3.58} 
and $\smb$ by \ref{3.69}. We see that $\omb$, hence also $\pb$, involves only derivatives tangential 
to the $\Cb_{\ub}$. We then conclude through \ref{5.50} that {\em the coefficients $a_{11}$, $a_{12}$ 
involve only derivatives tangential to the $\Cb_{\ub}$}. 

In regard to \ref{5.52}, consider first $\sn$, $\on$, and $n$. Since according to the expansion 
\ref{3.29} we have:
\begin{equation}
k=-\frac{1}{2c}h_{\mu\nu}\Nb^\nu\sd N^\mu
\label{5.57}
\end{equation}
the $S$ 1-from $k$ involves only derivatives tangential to the $\Cb_{\ub}$. Then by the formula 
\ref{3.a6} (see also \ref{3.42}) so does the $S$ vectorfield $\sn$. By \ref{3.61} this is also true for the function 
$\on$. Then by the 1st of \ref{3.48} the function $n$ likewise involves only derivatives tangential 
to the $\Cb_{\ub}$. Consider next $g$, $f_N$, and $e_N$, given by \ref{5.53}, \ref{5.54}, and 
\ref{5.55}, respectively. The first obviously involves only derivatives tangential to the $\Cb_{\ub}$. 
Showing that this is true also of the second and the third, reduces to showing that $\tchi$ and 
$\teta$ involve only derivatives tangential to the $\Cb_{\ub}$. That this is true for $\tchi$ follows 
from the relation \ref{3.76} and the fact that, in view of the definition of $\sk$ through the 
expansion \ref{3.29}, $\sk$ involves only derivatives tangential to the $\Cb_{\ub}$. That it is also 
true for $\teta$ follows from the relation \ref{3.80} and the corresponding fact about $\sn$ already 
established. Thus $g$, $f_N$, and $e_N$ involve only derivatives tangential to the $\Cb_{\ub}$. 
We then conclude through \ref{5.52} that {\em the coefficients $a_{21}$, $a_{22}$ involve 
only derivatives tangential to the $\Cb_{\ub}$.} 

In view of the definition of $Z$ as the commutator $[\Lb,L]$ according to \ref{2.17}, we have:
\begin{equation}
[\Lb,T]=Z
\label{5.58}
\end{equation}
Applying then $T$ to both sides of equations \ref{5.49}, \ref{5.51} we obtain the following 
propagation equations for the pair $T\lambdab$, $Ts_{NL}$ along the generators of each $\Cb_{\ub}$:
\begin{eqnarray}
&&\Lb T\lambdab=a_{11}T\lambdab+a_{12}Ts_{NL}+\stackrel{(1)}{b_1}\nonumber\\
&&\Lb Ts_{NL}=a_{21}T\lambdab+a_{22}Ts_{NL}+\stackrel{(1)}{b_2}\label{5.59}
\end{eqnarray}
a linear system similar to \ref{5.49}, \ref{5.51}, which is however inhomogeneous, the 
inhomogeneous terms being given by:
\begin{eqnarray}
&&\stackrel{(1)}{b_1}=Z\lambdab+(Ta_{11})\lambdab+(Ta_{12})s_{NL}\nonumber\\
&&\stackrel{(1)}{b_2}=Zs_{NL}+(Ta_{21})\lambdab+(Ta_{22})s_{NL}\label{5.60}
\end{eqnarray}
Now, since the coefficients $a_{11}$, $a_{12}$, $a_{21}$, $a_{22}$ involve only derivatives 
tangential to the $\Cb_{\ub}$, their $T$-derivatives $Ta_{11}$, $Ta_{12}$, $Ta_{21}$, $Ta_{22}$ 
involve only 1st derivatives transversal to the $\Cb_{ub}$, which can be expressed in terms of 
$\lambdab$, $s_{NL}$ and their derivatives tangential to the $\Cb_{\ub}$. We conclude that 
$\stackrel{(1)}{b_1}$, $\stackrel{(1)}{b_2}$, the inhomogeneous terms in \ref{5.59}, involve only 
$\lambdab$, $s_{NL}$ and their derivatives tangential to the $\Cb_{\ub}$. In particular on $\Cb_0$, 
by virtue of Proposition 5.1 we have:
\begin{equation}
\stackrel{(1)}{b_1}=\stackrel{(1)}{b_2}=0 \ \mbox{: on $\Cb_0$}
\label{5.61}
\end{equation}
and in view of \ref{5.1} the system reduces along the generators of $\Cb_0$ to the homogeneous 
linear system of ordinary differential equations
\begin{eqnarray}
&&\frac{\partial T\lambdab}{\partial u}=a_{11}T\lambdab+a_{12}Ts_{NL}\nonumber\\
&&\frac{\partial Ts_{NL}}{\partial u}=a_{21}T\lambdab+a_{22}Ts_{NL}
\label{5.62}
\end{eqnarray}

Proceeding in this way, that is applying $T$ repeatedly to both sides of \ref{5.49}, \ref{5.51} we 
obtain an inhomogenenous linear system of propagation equations for the pair $T^m \lambdab$, 
$T^m s_{NL}$, $m=1,2,3,...$:
\begin{eqnarray}
&&\Lb T^m \lambdab=a_{11}T^m \lambdab+a_{12}T^m s_{NL}+\stackrel{(m)}{b_1}\nonumber\\
&&\Lb T^m s_{NL}=a_{21}T^m \lambdab+a_{22}T^m s_{NL}+\stackrel{(m)}{b_2}
\label{5.63}
\end{eqnarray}
where the inhomogeneous terms are determined by the recursion relations:
\begin{eqnarray}
&&\stackrel{(m)}{b_1}=T\stackrel{(m-1)}{b_1}+ZT^{m-1}\lambdab+(Ta_{11})T^{m-1}\lambdab
+(Ta_{12})T^{m-1}s_{NL}\nonumber\\
&&\stackrel{(m)}{b_2}=T\stackrel{(m-1)}{b_2}+ZT^{m-1}s_{NL}+(Ta_{21})T^{m-1}\lambdab
+(Ta_{22})T^{m-1}s_{NL}
\label{5.64}
\end{eqnarray}
together with the conditions:
\begin{equation}
\stackrel{(0)}{b_1}=\stackrel{(0)}{b_2}=0
\label{5.65}
\end{equation}
It follows by induction that $\stackrel{(m)}{b_1}$, $\stackrel{(m)}{b_2}$ depend on 
$T^l \lambdab$, $T^l s_{NL}$ for $l=0,...,m-1$ and their derivatives tangential to the $\Cb_{\ub}$. 
On $\Cb_0$, in view of \ref{5.1}, the system \ref{5.63} reduces  to the following 
inhomogeneous linear system of ordinary differential equations along the generators of $\Cb_0$:
\begin{eqnarray}
&&\frac{\partial T^m\lambdab}{\partial u}=a_{11}T^m\lambdab+a_{12}T^m s_{NL}
+\stackrel{(m)}{b_1}\nonumber\\
&&\frac{\partial T^m s_{NL}}{\partial u}=a_{21}T^m\lambdab+a_{22}T^m s_{NL}
+\stackrel{(m)}{b_2}
\label{5.66}
\end{eqnarray}
where the inhomogeneous terms $\stackrel{(m)}{b_1}$, $\stackrel{(m)}{b_2}$ depend on the previously 
determined quantities 
$T^l \lambdab$, $T^l s_{NL}$ for $l=0,...,m-1$ and their derivatives tangential to $\Cb_0$. The 
system \ref{5.66} then allows us to determine recursively the higher order derived data $T^m\lambdab$, 
$T^m s_{NL}$, $m=1,2,3,...$ on $\Cb_0$, once boundary conditions for these on 
$S_{0,0}=\partial_-{\cal B}$ have been provided. 

\vspace{5mm}

\section{The Boundary Conditions for the Higher Order Derived Data  
and the Determination of the $T$ Derivatives  
of the Transformation Functions on $\partial_-{\cal B}$}

We proceed to deduce the boundary conditions for the 2nd derived data $T\lambdab$, $Ts_{NL}$ on 
$\partial_-{\cal B}$. In the case of $T\lambdab$ we have already deduced the result, 1st of 
\ref{4.216}:
\begin{equation}
\left.T\lambdab\right|_{\partial_-{\cal B}}=-\frac{c_0 k}{3l}
\label{5.67}
\end{equation}
What then remains to be determined at this level is $Ts_{NT}$ at $\partial_-{\cal B}$. For, the 
difference of $s_{NL}$ from $s_{NT}$ is (see \ref{3.42}):
\begin{equation}
s_{N\Lb}=\lambda\mbox{tr}\sss
\label{5.68}
\end{equation}
which vanishes together with its $T$ derivative at $\partial_-{\cal B}$. Along ${\cal K}$ we have:
\begin{equation}
s_{NT}=N^\mu T\beta_\mu=N^\mu T\beta^\prime_\mu+N^\mu T\triangle\beta_\mu
\label{5.69}
\end{equation}
Applying $T$ to this we obtain, along ${\cal K}$,
\begin{equation}
Ts_{NT}=N^\mu T^2\beta^\prime_\mu+N^\mu T^2\triangle\beta_\mu+
(TN^\mu)T\beta^\prime_\mu+(TN^\nu)T\triangle\beta_\mu
\label{5.70}
\end{equation}
From the expression \ref{4.4} for $\triangle\beta_\mu$, noting that the jumps $\epb$ and $\ep$ 
vanish at $\partial_-{\cal B}$ and by Proposition 4.2 $T\ep$ vanishes there as well, we obtain:
\begin{equation}
T\triangle\beta_\mu=-\frac{1}{2c} h_{\mu\nu}N^\nu T\epb \ \ \mbox{: at $\partial_-{\cal B}$}
\label{5.71}
\end{equation}
Since $LN^\mu$ vanishes along $\Cb_0$, 
$$TN^\mu=\Lb N^\mu=\sn^A\Omega_A^\mu+nN^\mu+\on\Nb^\mu \ \ \mbox{: along $\Cb_0$}$$
according to the expansion \ref{3.33}. Then by \ref{5.71}:
\begin{equation}
(TN^\mu)T\triangle\beta_\mu=\on T\epb \ \ \mbox{: at $\partial_-{\cal B}$}
\label{5.72}
\end{equation}
Substituting for $\on$ at $\partial_-{\cal B}$ from the 2nd of \ref{4.187} and noting that by 
\ref{4.154} and Proposition 4.4 we have:
\begin{equation}
\left.T\epb\right|_{\partial_-{\cal B}}=2\left.s_{\Nb\Lb}\right|_{\partial_-{\cal B}}
\label{5.73}
\end{equation}
we obtain:
\begin{equation}
\left.(TN^\mu)T\triangle\beta_\mu\right|_{\partial_-{\cal B}}
=\frac{l}{c_0}\left.s_{\Nb\Lb}\right|_{\partial_-{\cal B}}
\label{5.74}
\end{equation}
By \ref{4.172} and the fact that by the 1st of \ref{4.120} and by Proposition 4.4:
\begin{equation}
\left.\frac{\partial f}{\partial\tau}\right|_{\tau=0}=0, \ \ \ 
\left.\frac{\partial w}{\partial\tau}\right|_{\tau=0}=-1, \ \ \ 
\left.\frac{\partial\psi^A}{\partial\tau}\right|_{\tau=0}=
\left.\frac{\partial^2\psi^A}{\partial\tau^2}\right|_{\tau=0}=0 \ : \ A=1,...,n-1
\label{5.75}
\end{equation}
we obtain:
\begin{equation}
\left.T\beta^\prime_\mu\right|_{\partial_-{\cal B}}=
-\left.T^\prime\beta^\prime_\mu\right|_{\partial_-{\cal B}}=
-\left.\Lb\beta_\mu\right|_{\partial_-{\cal B}}
\label{5.76}
\end{equation}
the last equality by virtue of the relation \ref{4.101}, the continuity of the tangential 
derivatives across $\Cb_0$, and the fact that $\mu^\prime$ vanishes at $\partial_-{\cal B}$. 
It follows that:
\begin{equation}
\left.(TN^\mu)T\beta^\prime_\mu\right|_{\partial_-{\cal B}}=
-\left.\on s_{\Nb\Lb}\right|_{\partial_-{\cal B}}
=-\frac{l}{2c_0}\left.s_{\Nb\Lb}\right|_{\partial_-{\cal B}}
\label{5.77}
\end{equation}
We turn to the first two terms on the right in \ref{5.70}. 
From the expression \ref{4.4} for $\triangle\beta_\mu$, noting again that the jumps $\epb$ and $\ep$ 
vanish at $\partial_-{\cal B}$ and by Proposition 4.2 $T\ep$ vanishes there as well, we deduce:
\begin{eqnarray}
&&T^2\triangle\beta_\mu=-\frac{1}{2}c^{-1}h_{\mu\nu}N^\nu T^2\epb-T(c^{-1}h_{\mu\nu}N^\nu)T\epb
\nonumber\\
&&\hspace{15mm}-\frac{1}{2}c^{-1}h_{\mu\nu}\Nb^\nu T^2\ep \ \ \ 
\mbox{: at $\partial_-{\cal B}$}
\label{5.78}
\end{eqnarray}
Therefore:
\begin{equation}
N^\mu T^2\triangle\beta_\mu=-N^\mu T(c^{-1}h_{\mu\nu}N^\nu)T\epb+T^2\ep \ \ \ 
\mbox{: at $\partial_-{\cal B}$}
\label{5.79}
\end{equation}
Now, since $h_{\mu\nu}N^\mu N^\nu=0$, we have:
$$-N^\mu T(c^{-1}h_{\mu\nu}N^\nu)=c^{-1}h_{\mu\nu}N^\nu TN^\mu$$
Substituting for $TN^\mu$ from above yields:
$$\left. c^{-1}h_{\mu\nu}N^\nu TN^\mu\right|_{\partial_-{\cal B}}=
-2\left.\on\right|_{\partial_-{\cal B}}=-\frac{l}{c_0}$$
by the 2nd of \ref{4.187}. Hence \ref{5.79} becomes:
\begin{equation}
\left.N^\mu T^2\triangle\beta_\mu\right|_{\partial_-{\cal B}}
=-\frac{l}{c_0}\left.T\epb\right|_{\partial_-{\cal B}}
+\left.T^2\ep\right|_{\partial_-{\cal B}}
\label{5.80}
\end{equation}
By \ref{4.36} and the fact that both $r$ and $\epb$ vanish at $\partial_-{\cal B}$ we have:
\begin{equation}
T^2\ep=-2(Tr)T\epb \ \ \mbox{: at $\partial_-{\cal B}$}
\label{5.81}
\end{equation}
Then in view of \ref{4.215} we obtain simply:
\begin{equation}
\left.N^\mu T^2\triangle\beta_\mu\right|_{\partial_-{\cal B}}=0
\label{5.82}
\end{equation}
The first term on the right in \ref{5.70} remains to be considered. Applying $T$ to \ref{4.172} we 
have, along ${\cal K}$:
\begin{eqnarray}
&&\left.T^2\beta^\prime_\mu\right|_{{\cal K}}=
\left.\frac{\partial\beta^\prime_\mu}{\partial t}\right|_{{\cal K}}\frac{\partial^2 f}{\partial\tau^2}
+\left.\frac{\partial\beta^\prime_\mu}{\partial u^\prime}\right|_{{\cal K}}\frac{\partial^2 w}{\partial\tau^2}
+\left.\frac{\partial\beta^\prime_\mu}{\partial\vartheta^{\prime A}}\right|_{{\cal K}}\frac{\partial^2\psi^A}{\partial\tau^2}\nonumber\\
&&\hspace{12mm}+\left(\left.\frac{\partial^2\beta^\prime_\mu}{\partial t^2}\right|_{{\cal K}}\frac{\partial f}{\partial\tau}+
\left.\frac{\partial^2\beta^\prime_\mu}{\partial t\partial u^\prime}\right|_{{\cal K}}\frac{\partial w}{\partial\tau}+
\left.\frac{\partial^2\beta^\prime_\mu}{\partial t\partial\vartheta^{\prime A}}\right|_{{\cal K}}
\frac{\partial\psi^A}{\partial\tau}\right)\frac{\partial f}{\partial\tau}\nonumber\\
&&\hspace{12mm}+\left(\left.\frac{\partial^2\beta^\prime_\mu}{\partial u^\prime\partial t}\right|_{{\cal K}}\frac{\partial f}{\partial\tau}
+\left.\frac{\partial^2\beta^\prime_\mu}{\partial u^{\prime 2}}\right|_{{\cal K}}\frac{\partial w}{\partial\tau}+
\left.\frac{\partial^2\beta^\prime_\mu}{\partial u^\prime\partial\vartheta^{\prime A}}\right|_{{\cal K}}\frac{\partial\psi^A}{\partial\tau}\right)\frac{\partial w}{\partial\tau}\nonumber\\
&&\hspace{12mm}+\left(\left.\frac{\partial^2\beta^\prime_\mu}{\partial\vartheta^{\prime A}\partial t}\right|_{{\cal K}}\frac{\partial f}{\partial\tau}
+\left.\frac{\partial^2\beta^\prime_\mu}{\partial\vartheta^{\prime A}\partial u^\prime}\right|_{{\cal K}}\frac{\partial w}{\partial\tau}+
\left.\frac{\partial^2\beta^\prime_\mu}{\partial\vartheta^{\prime A}\partial\vartheta^{\prime B}}\right|_{{\cal K}}\frac{\partial\psi^B}{\partial\tau}\right)\frac{\partial\psi^A}{\partial\tau}
\nonumber\\
&&\hspace{50mm}
\label{5.83}
\end{eqnarray}
By Proposition 4.4:
\begin{equation}
\left.\frac{\partial^2 w}{\partial\tau^2}\right|_{\tau=0}=
2\left.\frac{\partial v}{\partial\tau}\right|_{\tau=0}
\label{5.84}
\end{equation}
Recalling also \ref{4.217} and \ref{5.76} we see that \ref{5.83} reduces at $\partial_-{\cal B}$ to:
\begin{equation}
\left.T^2\beta^\prime_\mu\right|_{\partial_-{\cal B}}=
\left.\left(2\frac{\partial v}{\partial\tau}\frac{\partial\beta^\prime_\mu}{\partial u^\prime}
-\frac{k}{3l}\frac{\partial\beta^\prime_\mu}{\partial t}
+\frac{\partial^2\beta^\prime_\mu}{\partial u^{\prime 2}}\right)\right|_{\tau=0}
\label{5.85}
\end{equation}
Now, when considering $\left.N^\mu T^2\beta^\prime_\mu\right|_{\partial_-{\cal B}}$ the first term in 
parenthesis on the right does not contribute, for, as in \ref{4.174} we have:
\begin{equation}
\left.N^\mu\frac{\partial\beta^\prime_\mu}{\partial u^\prime}\right|_{\partial_-{\cal B}}
=\left. L^{\prime\mu}T^\prime\beta^\prime_\mu\right|_{\partial_-{\cal B}}
=\left. T^{\prime\mu}L^\prime\beta^\prime_\mu\right|_{\partial_-{\cal B}}=0
\label{5.86}
\end{equation}
Therefore we obtain:
\begin{equation}
\left.N^\mu T^2\beta^\prime_\mu\right|_{\partial_-{\cal B}}=
-\frac{k}{3l}\left.L^{\prime\mu}L^\prime\beta^\prime_\mu\right|_{\partial_-{\cal B}}
+\left.L^{\prime\mu}T^{\prime 2}\beta^\prime_\mu\right|_{\partial_-{\cal B}}:=\xi
\label{5.87}
\end{equation}
a known smooth function on $\partial_-{\cal B}$ defined by the prior solution. Combining the results 
\ref{5.74}, \ref{5.77}, \ref{5.82}, and \ref{5.87}, we conclude through \ref{5.70} that:
\begin{equation}
\left.Ts_{NT}\right|_{\partial_-{\cal B}}=\frac{l}{2c_0}\left.s_{\Nb\Lb}\right|_{\partial_-{\cal B}}+\xi
\label{5.88}
\end{equation}
This completes the derivation of the boundary conditions for the 2nd derived data on 
$\partial_-{\cal B}$. 

We proceed to the derivation of the boundary conditions for the 3rd and higher order derived data 
$T^m\lambdab$, $T^m s_{NL}$, $m\geq 2$ on $\partial_-{\cal B}$. By virtue of \ref{5.68} we may 
replace the 2nd of these by $T^m s_{NT}$, because the difference involves at $\partial_-{\cal B}$ only 
the derived data of order up to $m-1$. We shall see that the determination of $T^m\lambdab$ on 
$\partial_-{\cal B}$ is coupled to the determination of $\partial^{m-1}v/\partial\tau^{m-1}$ at 
$\tau=0$ through the regularized identification equations of Proposition 4.5. As this coupling first 
manifests itself in the case $m=2$, namely boundary conditions for the 3rd derived data, 
we shall first consider this case before proceeding to the general higher order case. 

Applying $T^3$ to the boundary condition \ref{4.113} and evaluating the result at $\partial_-{\cal B}$ 
we obtain, in view of the fact that both $r$ and $\lambdab$ vanish at $\partial_-{\cal B}$, 
\begin{equation}
3(Tr)T^2\lambdab+3(T\lambdab)T^2 r=T^3\lambda \ \ \mbox{: at $\partial_-{\cal B}$}
\label{5.89}
\end{equation}
In regard to the right hand side at $\partial_-{\cal B}$, $\Lb^3\lambda$ is given along $\Cb_0$ by the 
data induced by the prior solution. Also, by \ref{5.38} $\Lb^2 L\lambda$ vanishes along $\Cb_0$. 
Moreover, in view of \ref{5.35}, $L\Lb\lambda=\Lb L\lambda=0$ along $\Cb_0$, hence also 
$\Lb L\Lb\lambda$ vanishes along $\Cb_0$ and $L\Lb^2\lambda=\Lb L\Lb\lambda=0$ along $\Cb_0$ as well. 
Next, applying $L$ to the propagation equation for $\lambda$ of Proposition 3.3 and evaluating the 
result along $\Cb_0$ we obtain, since $p$ and $\lambdab$ vanish along $\Cb_0$,
\begin{equation}
L^2\lambda=(Lp)\lambda+qL\lambdab \ \ \mbox{: along $\Cb_0$}
\label{5.90}
\end{equation}
Applying $\Lb$ to this and evaluating the result at $\partial_-{\cal B}$ gives, since $\lambda$ and 
$\Lb\lambda$ vanish at $\partial_-{\cal B}$,
\begin{equation}
\Lb L^2\lambda=\Lb(qL\lambdab) \ \ \mbox{: at $\partial_-{\cal B}$}
\label{5.91}
\end{equation}
Here the right hand side is already determined along $\Cb_0$, being expressed in terms of 2nd 
derived data. 
Moreover, in view of \ref{5.36} along $\Cb_0$ we have $L\Lb L\lambda=\Lb L^2\lambda$. Finally, 
we have:
$$L^2\Lb\lambda=L(\Lb L\lambda-Z\lambda)=L\Lb L\lambda-\sd\lambda\cdot\sL_L Z$$
which at $\partial_-{\cal B}$ reduces to:
\begin{equation}
L^2\Lb\lambda=L\Lb L\lambda \ \ \mbox{: at $\partial_-{\cal B}$}
\label{5.92}
\end{equation}
As a consequence of the above we can write: 
\begin{equation}
T^3\lambda=L^3\lambda+\xi \ \ \mbox{: at $\partial_-{\cal B}$}
\label{5.93}
\end{equation}
where from now on we denote by $\xi$ various known smooth functions on $\partial_-{\cal B}$ 
(represented as smooth functions on $S^{n-1}$). Applying $L^2$ to the propagation equation for 
$\lambda$ of Proposition 3.3 and evaluating the result at $\partial_-{\cal B}$ we obtain, in view of 
the fact that $p$, $L\lambda$ as well as $Lq$ vanish along $\Cb_0$, 
\begin{equation}
L^3\lambda=qL^2\lambdab \ \ \mbox{: at $\partial_-{\cal B}$}
\label{5.94}
\end{equation}
Since along $\Cb_0$ we have $T^2\lambdab=L^2\lambdab$, we then obtain from \ref{5.93}:
\begin{equation}
T^3\lambda=qT^2\lambdab+\xi \ \ \mbox{: at $\partial_-{\cal B}$}
\label{5.95}
\end{equation}
Substituting this in \ref{5.89} then yields:
\begin{equation}
(3Tr-q)T^2\lambdab+3(T\lambdab)T^2 r=\xi \ \ \mbox{: at $\partial_-{\cal B}$}
\label{5.96}
\end{equation}
By \ref{4.166}, since $v_0=-1$ (\ref{4.213}), we have:
\begin{equation}
Tr=-q=-\frac{l}{2c_0} \ \ \mbox{: at $\partial_-{\cal B}$}
\label{5.97}
\end{equation}
according to \ref{4.215}. Also, $T\lambdab$ at $\partial_-{\cal B}$ is given by the 
1st of \ref{4.216}. Thus, what remains to be determined is $T^2 r$ at $\partial_-{\cal B}$. 

Applying $T^2$ to \ref{4.40} and evaluating the result at $\partial_-{\cal B}$ we obtain:
\begin{equation}
T^2 r=s_0 T^2\epb+2s_1 T\epb \ \ \mbox{: at $\partial_-{\cal B}$, where 
$s_1=\left.Ts\right|_{\partial_-{\cal B}}$}
\label{5.98}
\end{equation}
By \ref{4.154}, since $v_0=-1$, we have:
\begin{equation}
T\epb=2s_{\Nb\Lb} \ \ \mbox{: at $\partial_-{\cal B}$}
\label{5.99}
\end{equation}
By Proposition 4.2 $s_1$ is likewise a known smooth function on $\partial_-{\cal B}$, 
$\left.T\kappa\right|_{\partial_-{\cal B}}$ being a quadruplet of known smooth functions on 
$\partial_-{\cal B}$. Therefore we can write \ref{5.98} as:
\begin{equation}
T^2 r=s_0 T^2\epb+\xi \ \ \mbox{: at $\partial_-{\cal B}$}
\label{5.100}
\end{equation}
Now, since $T^2\beta_\mu$ is determined along $\Cb_0$ in terms of the 2nd derived data, in view 
of \ref{5.85}, we have:
\begin{equation}
\left.T^2\triangle\beta_\mu\right|_{\partial_-{\cal B}}=
-2\left.\frac{\partial v}{\partial\tau}\frac{\partial\beta^\prime_\mu}{\partial u^\prime}
\right|_{\tau=0}+\xi_\mu=-2\left.\frac{\partial v}{\partial\tau}\right|_{\tau=0}\left.\Lb\beta_\mu\right|_{\partial_-{\cal B}}+\xi_\mu
\label{5.101}
\end{equation}
(see \ref{5.76}) $\xi_\mu$ denoting again known smooth functions on $\partial_-{\cal B}$. Hence:
\begin{equation}
\left.\Nb^\mu T^2\triangle\beta_\mu\right|_{\partial_-{\cal B}}
=-2\left.\frac{\partial v}{\partial\tau}\right|_{\tau=0}\left.s_{\Nb\Lb}\right|_{\partial_-{\cal B}}
+\xi
\label{5.102}
\end{equation}
It follows that:
\begin{equation}
\left.T^2\epb\right|_{\partial_-{\cal B}}=-2\left.\frac{\partial v}{\partial\tau}\right|_{\tau=0}
\left.s_{\Nb\Lb}\right|_{\partial_-{\cal B}}+\xi
\label{5.103}
\end{equation}
Substituting this in \ref{5.100} and recalling from \ref{4.146}, \ref{5.97}, \ref{5.99} that
\begin{equation}
2s_0\left.s_{\Nb\Lb}\right|_{\partial_-{\cal B}}=-\frac{l}{2c_0}
\label{5.104}
\end{equation}
we obtain:
\begin{equation}
\left.T^2 r\right|_{\partial_-{\cal B}}=\frac{l}{2c_0}\left.\frac{\partial v}{\partial\tau}
\right|_{\tau=0}+\xi
\label{5.105}
\end{equation}
This is to be substituted in \ref{5.96}. In view of \ref{5.97}, taking also into account the 1st of 
\ref{4.216} we then deduce:
\begin{equation}
\left.T^2\lambdab\right|_{\partial_-{\cal B}}=-\frac{c_0 k}{4l}\left.\frac{\partial v}{\partial\tau}
\right|_{\tau=0}+\xi
\label{5.106}
\end{equation}

To obtain an equation for $\left.(\partial v/\partial\tau)\right|_{\tau=0}$ we must appeal to 
the regularized identification equations of Proposition 4.5. This requires that we express 
$\left.(\partial\hat{f}/\partial\tau)\right|_{\tau=0}$ and 
$\left.(\partial\hat{\delta}^i/\partial\tau)\right|_{\tau=0}$ in terms of 
$\left.(\partial v/\partial\tau)\right|_{\tau=0}$. 
We have: 
\begin{equation}
\left.\frac{\partial\hat{f}}{\partial\tau}\right|_{\tau=0}=
\frac{1}{6}\left.\frac{\partial^3 f}{\partial\tau^3}\right|_{\tau=0}=
\frac{1}{6}\left.T^3 f\right|_{\partial_-{\cal B}}
\label{5.107}
\end{equation}
According to the 1st of \ref{4.179}:
$$T^3 f=T^2\rho+T^2\rhob$$
Since $\rho$ and $Z=[\Lb,L]$ vanish along $\Cb_0$, $\Lb^2\rho$, $L\Lb\rho$ and $\Lb L\rho$ likewise 
vanish along $\Cb_0$, hence
$$T^2\rho=L^2\rho \ \ \mbox{along $\Cb_0$}$$
Also, since $Lc$ and $L\lambda$ vanish along $\Cb_0$ (see \ref{5.37}, \ref{5.39}) so does $L\rhob$. It then follows that 
$\Lb L\rhob$ and $L\Lb\rhob$ likewise vanish along $\Cb_0$, hence
$$T^2\rhob=L^2\rhob+\Lb^2\rhob \ \ \mbox{along $\Cb_0$}$$
Consequently,
\begin{equation}
T^3 f=L^2\rho+L^2\rhob+\Lb^2\rhob \ \ \mbox{along $\Cb_0$}
\label{5.108}
\end{equation}
Moreover, since $\rho$ and $Lc$ vanish along $\Cb_0$, we have:
$$L^2\lambdab=cL^2\rho \ \ \mbox{along $\Cb_0$}$$
and since $\rhob$ vanishes at $\partial_-{\cal B}$:
$$L^2\lambda=cL^2\rhob \ \ \mbox{at $\partial_-{\cal B}$}$$
Also, since $\rhob$, $\Lb\rhob$ both vanish at $\partial_-{\cal B}$, we have:
$$\Lb^2\lambda=c\Lb^2\rhob \ \ \mbox{at $\partial_-{\cal B}$}$$
Consequently, we can write \ref{5.108} in the form:
\begin{equation}
\left.cT^3 f\right|_{\partial_-{\cal B}}=
\left.L^2\lambdab\right|_{\partial_-{\cal B}}
+\left.L^2\lambda\right|_{\partial_-{\cal B}}
+\left.\Lb^2\lambdab\right|_{\partial_-{\cal B}}
\label{5.109}
\end{equation}
Moreover, in analogy with the above we have:
$$L^2\lambdab=T^2\lambdab \ \ \mbox{along $\Cb_0$},$$
by \ref{5.90} $L^2\lambda=qT\lambdab$ at $\partial_-{\cal B}$, and $\Lb^2\lambda$ is given along 
$\Cb_0$ by the data induced on $\Cb_0$ by the prior solution. Therefore \ref{5.109} takes the form:
\begin{equation}
\left.cT^3 f\right|_{\partial_-{\cal B}}=\left.T^2\lambdab\right|_{\partial_-{\cal B}}+\xi
\label{5.110}
\end{equation}
Substituting from \ref{5.106} we obtain:
\begin{equation}
\left.T^3 f\right|_{\partial_-{\cal B}}=-\frac{k}{4l}\left.\frac{\partial v}{\partial\tau}
\right|_{\tau=0}+\xi
\label{5.111}
\end{equation}
which by \ref{5.107} yields:
\begin{equation}
\left.\frac{\partial\hat{f}}{\partial\tau}\right|_{\tau=0}=-\frac{k}{24l}\left.
\frac{\partial v}{\partial\tau}\right|_{\tau=0}+\xi
\label{5.112}
\end{equation}

We turn to $\left.(\partial\hat{\delta}^i/\partial\tau)\right|_{\tau=0}$. From \ref{4.244} we have:
\begin{equation}
\left.\frac{\partial\hat{\delta}^i}{\partial\tau}\right|_{\tau=0}=
\frac{1}{24}\left.\frac{\partial^4\delta^i}{\partial\tau^4}\right|_{\tau=0}=
\frac{1}{24}\left.T^4\delta^i\right|_{\partial_-{\cal B}}
\label{5.113}
\end{equation}
and from \ref{4.223}:
\begin{equation}
\left.T^4\delta^i\right|_{\partial_-{\cal B}}=\left.T^4 g^i\right|_{\partial_-{\cal B}}
-\left.N^i T^4 f\right|_{\partial_-{\cal B}}
\label{5.114}
\end{equation}
From \ref{4.179} we deduce:
$$T^4 g^i-N^i T^4 f=3(TN^i)T^2\rho+(\Nb^i-N^i)T^3\rhob+\xi^i \ \ \mbox{: at $\partial_-{\cal B}$}$$
or:
\begin{equation}
c(T^4 g^i-N^i T^4 f)=3(TN^i)T^2\lambdab+(\Nb^i-N^i)T^3\lambda+\xi^i \ \ 
\mbox{: at $\partial_-{\cal B}$}
\label{5.115}
\end{equation}
where the $\xi^i$ are known smooth functions on $\partial_-{\cal B}$. Now by \ref{4.188} and 
\ref{4.199}:
\begin{equation}
\left.TN^i\right|_{\partial_-{\cal B}}=lS^i_0
\label{5.116}
\end{equation}
Taking into account the fact that by \ref{5.95}, \ref{5.97}:
\begin{equation}
\left.T^3\lambda\right|_{\partial_-{\cal B}}=\frac{l}{2c_0}\left.T^2\lambdab
\right|_{\partial_-{\cal B}}+\xi
\label{5.117}
\end{equation}
as well as \ref{4.199}, we then obtain:
\begin{equation}
\left.c(T^4 g^i-N^i T^4 f)\right|_{\partial_-{\cal B}}
=4l\left(S^i_0+\frac{1}{2}m^A\Omega^{\prime i}_{A 0}\right)\left.T^2\lambdab\right|_{\partial_-{\cal B}}+\xi^i
\label{5.118}
\end{equation}
Substituting for $\left.T^2\lambdab\right|_{\partial_-{\cal B}}$ from \ref{5.106} and substituting 
the result in \ref{5.114} and \ref{5.113} then yields:
\begin{equation}
\left.\frac{\partial\hat{\delta}^i}{\partial\tau}\right|_{\tau=0}=
-\frac{k}{24}\left(S^i_0+\frac{1}{2}m^A\Omega^{\prime i}_{A 0}\right)
\left.\frac{\partial v}{\partial\tau}\right|_{\tau=0}+\xi^i
\label{5.119}
\end{equation}

To determine $\left.(\partial v/\partial\tau)\right|_{\tau=0}$ we substitute 
$$v=v(\tau,\vartheta), \ \ \ \gamma=\gamma(\tau,\vartheta)$$
in the regularized identification equations 
$$\hat{F}^i((\tau,\vartheta),(v,\gamma))=0$$
of Proposition 4.5 and differentiate implicitly with respect to $\tau$ to obtain:
\begin{equation}
\frac{\partial\hat{F}^i}{\partial v}\frac{\partial v}{\partial\tau}+
\frac{\partial\hat{F}^i}{\partial\gamma^A}\frac{\partial\gamma^A}{\partial\tau}+
\frac{\partial\hat{F}^i}{\partial\tau}=0
\label{5.120}
\end{equation}
This is to be evaluated at $\tau=0$ using the result of Proposition 4.4:
$$v(0,\vartheta)=v_0(\vartheta)=-1, \ \ \ \gamma^A(0,\vartheta)=\gamma^A_0(\vartheta)=
\frac{2k(\vartheta)m^A(\vartheta)}{9}$$
From the expression for $\hat{F}^i((\tau,\vartheta),(v,\gamma))$ of Proposition 4.5 we obtain:
\begin{eqnarray}
&&\frac{\partial\hat{F}^i}{\partial v}((0,\vartheta),(v_0(\vartheta),\gamma_0(\vartheta)))
=S^i_0(\vartheta)\left[l(\vartheta)\hat{f}(0,\vartheta+\frac{1}{2}k(\vartheta)\right]\nonumber\\
&&\hspace{40mm}=\frac{k(\vartheta)}{3}S_0^i(\vartheta)
\label{5.121}
\end{eqnarray}
the 2nd equality by virtue of \ref{4.241}. Also:
\begin{equation}
\frac{\partial\hat{F}^i}{\partial\gamma^A}((0,\vartheta),(v_0(\vartheta),\gamma_0(\vartheta)))
=\Omega^{\prime i}_{A 0}(\vartheta)
\label{5.122}
\end{equation}
Moreover, we obtain:
\begin{eqnarray}
&&\frac{\partial\hat{F}^i}{\partial\tau}((0,\vartheta),(v_0(\vartheta),\gamma_0(\vartheta)))
=-S^i_0(\vartheta)l(\vartheta)\frac{\partial\hat{f}}{\partial\tau}(0,\vartheta)
-\frac{\partial\hat{\delta}^i}{\partial\tau}(0,\vartheta)\nonumber\\
&&\hspace{42mm}+E^i((0,\vartheta),(v_0(\vartheta),\gamma_0(\vartheta)))
\label{5.123}
\end{eqnarray}
where:
\begin{eqnarray}
&&E^i((0,\vartheta),(v_0(\vartheta),\gamma_0(\vartheta)))=
\frac{1}{2}(a_0^i(\vartheta)\hat{f}(0,\vartheta)+a_1^i(\vartheta))\hat{f}(0,\vartheta)\nonumber\\
&&\hspace{38mm}+A_0^i(0,0,\vartheta)-\gamma_0^A(\vartheta)\Theta^i_{1 A}(0,0,0,\vartheta)
\label{5.124}
\end{eqnarray}
are known smooth functions on $S^{n-1}$. Substituting \ref{5.112} and \ref{5.119} in \ref{5.123} 
yields:
\begin{equation}
\frac{\partial\hat{F}^i}{\partial\tau}((0,\vartheta),(v_0(\vartheta),\gamma_0(\vartheta)))
=\frac{k(\vartheta)}{12}\left(S^i_0(\vartheta)+\frac{1}{4}m^A(\vartheta)
\Omega^{\prime i}_{A 0}(\vartheta)\right)
\frac{\partial v}{\partial\tau}(0,\vartheta)+\xi^i(\vartheta)
\label{5.125}
\end{equation}
Substituting this together with \ref{5.121}, \ref{5.122}, in \ref{5.120} evaluated at $\tau=0$,  
results in the equations:
\begin{equation}
\frac{k}{12}\left(5S^i_0+\frac{1}{4}m^A\Omega^{\prime i}_{A 0}\right)
\left.\frac{\partial v}{\partial\tau}\right|_{\tau=0}
+\Omega^{\prime i}_{A 0}\left.\frac{\partial\gamma^A}{\partial\tau}\right|_{\tau=0}+\xi^i=0
\label{5.126}
\end{equation}
Now, at each point on $\partial_-{\cal B}$ the spatial vectors
$$S_0=S_0^i\frac{\partial}{\partial x^i}, \ \ \ 
\Omega^\prime_{A 0}=\Omega_{A 0}^{\prime i}\frac{\partial}{\partial x^i} \ : \ A=1,...,n-1$$
constitute a basis for the spatial hyperplane $\Sigma_t$ through the point. These vectors are 
therefore linearly independent. Equations \ref{5.126} constitute a spatial vector equation. The 
$S_0$-component of this equation determines $\left.(\partial v/\partial\tau)\right|_{\tau=0}$, 
after which the $\Omega_{A 0}$-components determine 
$\left.(\partial\gamma^A/\partial\tau)\right|_{\tau=0}$, in terms of known smooth functions 
on $S^{n-1}$, that is on $\partial_-{\cal B}$. Substituting the first in \ref{5.106} then so 
determines $\left.T^2\lambdab\right|_{\partial_-{\cal B}}$. 

To complete the derivation of the boundary conditions for the 3rd derived data on 
$\partial_-{\cal B}$, what remains to be done is to determine 
$\left.T^2s_{NT}\right|_{\partial_-{\cal B}}$. By \ref{5.103} 
$\left.T^2\epb\right|_{\partial_-{\cal B}}$ is now a known smooth function on $\partial_-{\cal B}$, 
and by \ref{5.101} so is $\left.T^2\triangle\beta_\mu\right|_{\partial_-{\cal B}}$. By \ref{5.85} 
$\left.T^2\beta^\prime_\mu\right|_{\partial_-{\cal B}}$ is likewise already determined. Applying then 
$T$ to \ref{5.70} and evaluating the result at $\partial_-{\cal B}$ we obtain:
\begin{equation}
\left.T^2 s_{NT}\right|_{\partial_-{\cal B}}=
\left.N^\mu T^3\beta^\prime_\mu\right|_{\partial_-{\cal B}}+
\left.N^\mu T^3\triangle\beta_\mu\right|_{\partial_-{\cal B}}+\xi
\label{5.127}
\end{equation}
Now from \ref{4.4} we have, at $\partial_-{\cal B}$, 
\begin{equation}
\left.T^3\triangle\beta_\mu\right|_{\partial_-{\cal B}}=
\left.-\frac{h_{\mu\nu}}{2c}\left(N^\nu T^3\epb+\Nb^\nu T^3\ep\right)\right|_{\partial_-{\cal B}}
+\xi_\mu
\label{5.128}
\end{equation}
Since by Proposition 4.2 $\ep=-j\epb^2$, we have (see \ref{4.a1}):
\begin{equation}
\left.T^3\ep\right|_{\partial_-{\cal B}}=-6j_0(\kappa_0)\left.T\epb\right|_{\partial_-{\cal B}}
\left.T^2\epb\right|_{\partial_-{\cal B}}-6(\left.T\epb\right|_{\partial_-{\cal B}})^2
\left.Tj\right|_{\partial_-{\cal B}}
\label{5.129}
\end{equation}
these are now known smooth functions on $\partial_-{\cal B}$. Therefore \ref{5.128} reduces to:
\begin{equation}
\left.T^3\triangle\beta_\mu\right|_{\partial_-{\cal B}}=
\left.-\frac{h_{\mu\nu}}{2c}N^\nu T^3\epb\right|_{\partial_-{\cal B}}+\xi_\mu
\label{5.130}
\end{equation}
Multiplying by $\left.N^\mu\right|_{\partial_-{\cal B}}$ then gives simply:
\begin{equation}
\left. N^\mu T^3\triangle\beta_\mu\right|_{\partial_-{\cal B}}=\xi
\label{5.131}
\end{equation}
a known smooth function on $\partial_-{\cal B}$. We turn to the 1st term on the right in \ref{5.127}. 
Applying $T$ to \ref{5.83} and evaluating the result at $\partial_-{\cal B}$, taking into account the  
fact that by \ref{5.110} $\left.(\partial^3 f/\partial\tau^3)\right|_{\tau=0}$ is now known, while 
$$\left.\frac{\partial^3 w}{\partial\tau^3}\right|_{\tau=0}=
\left.3\frac{\partial^2 v}{\partial\tau^2}\right|_{\tau=0}, \ \ \ 
\left.\frac{\partial\psi^A}{\partial\tau^3}\right|_{\tau=0}=6\gamma_0^A,$$
we obtain:
\begin{equation}
\left.T^3\beta^\prime_\mu\right|_{\partial_-{\cal B}}=
3\left.\frac{\partial^2 v}{\partial\tau^2}\right|_{\tau=0}
\left.\frac{\partial\beta^\prime_\mu}{\partial u^\prime}\right|_{\tau=0}+\xi_\mu
\label{5.132}
\end{equation}
In view of \ref{5.86}, multiplying by $\left.N^\mu\right|_{\partial_-{\cal B}}$ then gives simply:
\begin{equation}
\left. N^\mu T^3\beta^\prime_\mu\right|_{\partial_-{\cal B}}=\xi
\label{5.133}
\end{equation}
a known smooth function on $\partial_-{\cal B}$. In view of the results \ref{5.131} and \ref{5.133} 
we conclude through \ref{5.127} that $\left.T^2 s_{NT}\right|_{\partial_-{\cal B}}$ is now a known 
smooth function on $\partial_-{\cal B}$. This completes the derivation of the boundary conditions 
for the 3rd derived data on $\partial_-{\cal B}$. 

We proceed to the derivation of the boundary conditions at $\partial_-{\cal B}$ for the higher order 
derived data. At the $m+1$ th step we are to determine $T^m\lambdab$ and $T^m s_{NT}$ at 
$\partial_-{\cal B}$. This requires determining at the same time 
$\partial^{m-1} v/\partial\tau^{m-1}$, $\partial^{m-1}\gamma^A/\partial\tau^{m-1}$ at $\tau=0$, 
which brings in the regularized identification equations (Proposition 4.5) and involves 
$\partial^{m-1}\hat{f}/\partial\tau^{m-1}$, $\partial^{m-1}\hat{\delta}^i/\partial\tau^{m-1}$ 
at $\tau=0$. These in turn depend on $T^m\lambdab$ at $\partial_-{\cal B}$. 

The argument being inductive, at the $m+1$ th step we may assume that $T^k\lambdab$ and $T^k s_{NT}$ 
at $\partial_-{\cal B}$ for $k=0,...,m-1$ have already been determined and are represented as smooth 
functions on $S^{n-1}$. Also, $\partial^{k-1} v/\partial\tau^{k-1}$, 
$\partial^{k-1}\gamma^A/\partial\tau^{k-1}$ at $\tau=0$ for $k=1,...,m-1$ have likewise already been 
determined. Moreover, $\partial^{k-1}\hat{f}/\partial\tau^{k-1}$, 
$\partial^{k-1}\hat{\delta}^i/\partial\tau^{k-1}$ at $\tau=0$ for $k=1,...,m-1$ have similarly been 
determined as well (inductive hypothesis). 

Since the coefficients of the $\tau^{k-1}$ term in the Taylor expansions of $\hat{f}$ and 
$\hat{\delta}^i$ with respect to $\tau$ are 
$$\frac{1}{(k-1)!}\left.\frac{\partial^{k-1}\hat{f}}{\partial\tau^{k-1}}\right|_{\tau=0} 
\ \ \mbox{and} \ \ \frac{1}{(k-1)!}\left.\frac{\partial^{k-1}\hat{\delta}^i}{\partial\tau^{k-1}}\right|_{\tau=0}$$
respectively, and these are at the same time the coefficients of the $\tau^{k+1}$ term and the 
$\tau^{k+2}$ term in the Taylor expansions of $f$ and $\delta^i$ with respect to $\tau$ respectively, 
that is
$$\frac{1}{(k+1)!}\left.\frac{\partial^{k+1} f}{\partial\tau^{k+1}}\right|_{\tau=0} \ \ 
\mbox{and} \ \ \frac{1}{(k+2)!}\left.\frac{\partial^{k+2}\delta^i}{\partial\tau^{k+2}}\right|_{\tau=0}$$
it follows that: 
\begin{eqnarray}
&&\left.\frac{\partial^{k-1}\hat{f}}{\partial\tau^{k-1}}\right|_{\tau=0}=
\frac{1}{k(k+1)}\left.\frac{\partial^{k+1}f}{\partial\tau^{k+1}}\right|_{\tau=0}
\label{5.134}\\
&&\left.\frac{\partial^{k-1}\hat{\delta}^i}{\partial\tau^{k-1}}\right|_{\tau=0}=
\frac{1}{k(k+1)(k+2)}\left.\frac{\partial^{k+2}\delta^i}{\partial\tau^{k+2}}\right|_{\tau=0}
\label{5.135}
\end{eqnarray}
Thus, by the inductive hypothesis, at the $m+1$ th step $\partial^{k+1} f/\partial\tau^{k+1}$, 
$\partial^{k+2}\delta^i/\partial\tau^{k+2}$ at $\tau=0$ for $k=1,...,m-1$ are known smooth functions 
on $S^{n-1}$. 

Applying now $T^{m-1}$ to \ref{4.172} and evaluating the result at $\partial_-{\cal B}$ we obtain, 
by virtue of the inductive hypothesis, 
\begin{equation}
\left.T^m\beta^\prime_\mu\right|_{\partial_-{\cal B}}=
m\left.\frac{\partial^{m-1} v}{\partial\tau^{m-1}}\right|_{\tau=0}
\left.\Lb\beta_\mu\right|_{\partial_-{\cal B}}+\xi_\mu
\label{5.136}
\end{equation}
(see \ref{5.76}) where we again denote by $\xi$ various known smooth functions on $S^{n-1}$. 
Since $\left.T^m\beta_\mu\right|_{\partial_-{\cal B}}$ are by the inductive hypothesis known smooth 
functions on $S^{n-1}$ being determined by the up to the $m$ th order derived data, it follows that:
\begin{equation}
\left.T^m\triangle\beta_\mu\right|_{\partial_-{\cal B}}=
-m\left.\frac{\partial^{m-1} v}{\partial\tau^{m-1}}\right|_{\tau=0}
\left.\Lb\beta_\mu\right|_{\partial_-{\cal B}}+\xi_\mu
\label{5.137}
\end{equation}
Since moreover $\left.T^k\Nb^\mu\right|_{\partial_-{\cal B}}$ for $k=0,...,m$ are likewise known 
smooth functions on $S^{n-1}$, we conclude that:
\begin{equation}
\left.T^m\epb\right|_{\partial_-{\cal B}}=-m\left.\frac{\partial^{m-1} v}{\partial\tau^{m-1}}\right|_{\tau=0}\left.s_{\Nb\Lb}\right|_{\partial_-{\cal B}}+\xi
\label{5.138}
\end{equation}
Applying $T^m$ to \ref{4.40} and evaluating the result at $\partial_-{\cal B}$ we obtain, again by 
virtue of the inductive hypothesis, 
\begin{equation}
\left.T^m r\right|_{\partial_-{\cal B}}=s_0\left.T^m\epb\right|_{\partial_-{\cal B}}+\xi
\label{5.139}
\end{equation}
Substituting \ref{5.138} and taking into account \ref{5.104} then yields:
\begin{equation}
\left.T^m r\right|_{\partial_-{\cal B}}=
\frac{l}{4c_0}m\left.\frac{\partial^{m-1} v}{\partial\tau^{m-1}}\right|_{\tau=0}+\xi
\label{5.140}
\end{equation}

We now apply $T^{m+1}$ to the boundary condition \ref{4.113} and evaluate the result at 
$\partial_-{\cal B}$. Since by virtue of the inductive hypothesis 
$\left.T^k\lambdab\right|_{\partial_-{\cal B}}$, $\left.T^{k+1}\lambda\right|_{\partial_-{\cal B}}$ 
and $\left.T^k r\right|_{\partial_-{\cal B}}$ for $k=1,...,m-1$ are already known smooth functions, 
this gives:
\begin{equation}
(m+1)(Tr)T^m\lambdab+(m+1)(T\lambdab)T^m r=T^{m+1}\lambda+\xi \ \ \mbox{: at $\partial_-{\cal B}$}
\label{5.141}
\end{equation}
Now, by virtue of the inductive hypothesis the derived data on $\Cb_0$ of order up to $m$ are  
determined by solving the system \ref{5.66}, with $k=0,...,m$ in the role of $m$, with the  
boundary conditions on $S_{0,0}=\partial_-{\cal B}$ provided. Thus $T^k\lambda$ and $\Lb T^k\lambda$ 
are determined along $\Cb_0$ for $k=0,...,m$. Since 
$$\Lb T^k\lambda-T^k\Lb\lambda=\sum_{l=0}^{k-1} T^l(ZT^{k-1-l}\lambda)$$
and $(\sL_T)^l Z$ on $\Cb_0$ for $l=0,...,k-1$ is determined by the derived data on $\Cb_0$ of order 
up to $k$, it follows that $T^k\Lb\lambda$ is determined on $\Cb_0$ for $k=0,...,m$. 
In particular $\left.T^m\Lb\lambda\right|_{\partial_-{\cal B}}$ is a known smooth function. 
Therefore we can write:
\begin{equation}
\left.T^{m+1}\lambda\right|_{\partial_-{\cal B}}=\left.T^m L\lambda\right|_{\partial_-{\cal B}}+\xi
\label{5.142}
\end{equation}
Recall from Proposition 3.3 that
$$L\lambda=p\lambda+q\lambdab$$
Since $\left.\lambda\right|_{\partial_-{\cal B}}=\left.T\lambda\right|_{\partial_-{\cal B}}=0$, we have:
$$\left.T^m(p\lambda)\right|_{\partial_-{\cal B}}=
\sum_{k=2}^m\left(\begin{array}{c} 
m\\
k
\end{array}\right)\left.T^k\lambda\right|_{\partial_-{\cal B}}\left.T^{m-k}p\right|_{\partial_-{\cal B}}$$
so at most $m-2$ $T$-derivatives of $p$, thus the derived data of order only up to $m-1$ is involved 
here. Similarly since $\left.\lambdab\right|_{\partial_-{\cal B}}=0$, we have:
$$\left.T^m(q\lambdab)\right|_{\partial_-{\cal B}}=\left.q\right|_{\partial_-{\cal B}}
\left.T^m\lambdab\right|_{\partial_-{\cal B}}
+\sum_{k=1}^{m-1}\left(\begin{array}{c}
m\\
k
\end{array}\right)\left.T^k\lambdab\right|_{\partial_-{\cal B}}\left.T^{m-k}q\right|_{\partial_-{\cal B}}$$
and the sum is a known smooth function. It follows that:
\begin{equation}
\left.T^{m+1}\lambda\right|_{\partial_-{\cal B}}=\left.q\right|_{\partial_-{\cal B}}
\left.T^m\lambdab\right|_{\partial_-{\cal B}}+\xi
\label{5.143}
\end{equation}
Substituting this in \ref{5.141} and recalling \ref{5.97} we arrive at the relation: 
\begin{equation}
(m+2)(Tr)T^m\lambdab+(m+1)(T\lambdab)T^m r=\xi \ \ \mbox{: at $\partial_-{\cal B}$}
\label{5.144}
\end{equation}
Substituting finally \ref{5.140} and recalling \ref{5.97} and the 1st of \ref{4.216} we deduce:
\begin{equation}
\left.T^m\lambdab\right|_{\partial_-{\cal B}}=
-\frac{c_0 k}{6l}\frac{m(m+1)}{(m+2)}\left.\frac{\partial^{m-1}v}{\partial\tau^{m-1}}\right|_{\tau=0}
+\xi
\label{5.145}
\end{equation}

Next, we consider $T^{m+1}f$ at $\partial_-{\cal B}$. We apply $T^m$ to the equation 
(see \ref{4.119}):
\begin{equation}
cTf=\lambdab+\lambda \ \ \mbox{: along ${\cal K}$}
\label{5.146}
\end{equation}
and evaluate the result at $\partial_-{\cal B}$. Noting that 
$\left.T^m\lambda\right|_{\partial_-{\cal B}}$ and $\left.T^k c\right|_{\partial_-{\cal B}}$, 
$\left.T^k f\right|_{\partial_-{\cal B}}$ for $k=0,...,m$ are, by virtue of the inductive hypothesis, 
known smooth functions, we deduce:
\begin{equation}
\left.cT^{m+1}f\right|_{\partial_-{\cal B}}=\left.T^m\lambdab\right|_{\partial_-{\cal B}}+\xi
\label{5.147}
\end{equation}
Substituting \ref{5.145} we then obtain:
\begin{equation}
\left.T^{m+1}f\right|_{\partial_-{\cal B}}=-\frac{k}{6l}\frac{m(m+1)}{(m+2)}
\left.\frac{\partial^{m-1}v}{\partial\tau^{m-1}}\right|_{\tau=0}+\xi
\label{5.148}
\end{equation}
which in view of \ref{5.134} is equivalent to:
\begin{equation}
\left.\frac{\partial^{m-1}\hat{f}}{\partial\tau^{m-1}}\right|_{\tau=0}=
-\frac{k}{6l}\frac{1}{(m+2)}\left.\frac{\partial^{m-1}v}{\partial\tau^{m-1}}\right|_{\tau=0}+\xi
\label{5.149}
\end{equation}

Finally, we consider $T^{m+2}\delta^i$ at $\partial_-{\cal B}$. From the definition \ref{4.223},
\begin{equation}
\left.T^{m+2}\delta^i\right|_{\partial_-{\cal B}}=\left.T^{m+2}g^i\right|_{\partial_-{\cal B}}
-\left.N^i\right|_{\partial_-{\cal B}}\left.T^{m+2}f\right|_{\partial_-{\cal B}}
\label{5.150}
\end{equation}
We have (see \ref{4.119}):
\begin{equation}
cTg^i=\lambdab N^i+\lambda\Nb^i \ \ \mbox{: along ${\cal K}$}
\label{5.151}
\end{equation}
Multiplying by $c^{-1}$, applying $T^{m+1}$ and evaluating the result at $\partial_-{\cal B}$ we 
obtain, since $\left.\lambdab\right|_{\partial_-{\cal B}}=\left.\lambda\right|_{\partial_-{\cal B}}
=\left.T\lambda\right|_{\partial_-{\cal B}}=0$,
\begin{eqnarray}
&&\left.T^{m+2}g^i\right|_{\partial_-{\cal B}}=
\left. c^{-1}N^i T^{m+1}\lambdab\right|_{\partial_-{\cal B}}\label{5.152}\\
&&\hspace{18mm}+
\sum_{k=1}^m\left(\begin{array}{c}
m+1\\
k
\end{array}\right)\left.T^{m+1-k}(c^{-1}N^i)\right|_{\partial_-{\cal B}}\left.T^k\lambdab\right|_{\partial_-{\cal B}}\nonumber\\
&&\hspace{22mm}+\left. c^{-1}\Nb^i T^{m+1}\lambda\right|_{\partial_-{\cal B}}\nonumber\\
&&\hspace{18mm}+
\sum_{k=2}^m\left(\begin{array}{c}
m+1\\
k
\end{array}\right)\left. T^{m+1-k}(c^{-1}\Nb^i)\right|_{\partial_-{\cal B}}\left. T^k\lambda\right|_{\partial_-{\cal B}}\nonumber
\end{eqnarray}
Here, $\left.T^k\lambdab\right|_{\partial_-{\cal B}}$ for $k=1,...,m-1$ and 
$\left.T^k\lambda\right|_{\partial_-{\cal B}}$ for $k=2,...,m$ re already determined. 
Furthermore, in the first sum we have $\left.T^{m+1-k}(c^{-1}N^i)\right|_{\partial_-{\cal B}}$, 
$m+1-k\leq m$, which is already determined by the derived data of order up to $m$, and 
in the second sum we have $\left.T^{m+1-k}(c^{-1}\Nb^i)\right|_{\partial_-{\cal B}}$, 
$m+1-k\leq m-1$, which is a fortiori so determined. Consequently, \ref{5.152} is of the form: 
\begin{eqnarray}
&&\left.T^{m+2}g^i\right|_{\partial_-{\cal B}}=\left.c^{-1}N^i T^{m+1}\lambdab\right|_{\partial_-{\cal B}}
\label{5.153}\\
&&\hspace{24mm}+(m+1)\left.T(c^{-1}N^i)\right|_{\partial_-{\cal B}}\left.T^m\lambdab\right|_{\partial_-{\cal B}}\nonumber\\
&&\hspace{24mm}+\left.c^{-1}\Nb^i T^{m+1}\lambda\right|_{\partial_-{\cal B}}+\xi^i\nonumber
\end{eqnarray}
Similarly, we deduce from \ref{5.146}:
\begin{eqnarray}
&&\left.T^{m+2}f\right|_{\partial_-{\cal B}}=\left.c^{-1}T^{m+1}\lambdab\right|_{\partial_-{\cal B}}
\label{5.154}\\
&&\hspace{24mm}+(m+1)\left.T(c^{-1})\right|_{\partial_-{\cal B}}\left.T^m\lambdab\right|_{\partial_-{\cal B}}\nonumber\\
&&\hspace{24mm}+\left.c^{-1}T^{m+1}\lambda\right|_{\partial_-{\cal B}}+\xi\nonumber
\end{eqnarray}
Substituting \ref{5.153} and \ref{5.154} in \ref{5.150}, the terms involving 
$\left.T^{m+1}\lambdab\right|_{\partial_-{\cal B}}$ cancel and we obtain:
\begin{eqnarray}
&&\left.T^{m+2}\delta^i\right|_{\partial_-{\cal B}}=(m+1)\left.c^{-1}TN^i\right|_{\partial_-{\cal B}}\left.T^m\lambdab\right|_{\partial_-{\cal B}}\nonumber\\
&&\hspace{24mm}+\left.c^{-1}(\Nb^i-N^i)\right|_{\partial_-{\cal B}}\left.T^{m+1}\lambda\right|_{\partial_-{\cal B}}+\xi^i
\label{5.155}
\end{eqnarray}
Substituting for $\left.T^{m+1}\lambda\right|_{\partial_-{\cal B}}$ in terms of 
$\left.T^m\lambdab\right|_{\partial_-{\cal B}}$ from \ref{5.143} this becomes:
\begin{equation}
\left.T^{m+2}\delta^i\right|_{\partial_-{\cal B}}=
\left.c^{-1}\left((m+1)TN^i+q(\Nb^i-N^i)\right)\right|_{\partial_-{\cal B}}
\left.T^m\lambdab\right|_{\partial_-{\cal B}}+\xi^i
\label{5.156}
\end{equation}
From \ref{5.116}, \ref{5.97}, and \ref{4.199}:
\begin{equation}
\left.c^{-1}\left((m+1)TN^i+q(\Nb^i-N^i)\right)\right|_{\partial_-{\cal B}}
=c_0^{-1}l\left((m+2)S_0^i+2m^A\Omega^{\prime i}_{A 0}\right)
\label{5.157}
\end{equation}
Substituting this as well as \ref{5.145} in \ref{5.156} then yields:
\begin{equation}
\left.T^{m+2}\delta^i\right|_{\partial_-{\cal B}}=
-\frac{k}{6}\frac{m(m+1)}{(m+2)}\left((m+2)S_0^i+2m^A\Omega^{\prime i}_{A 0}\right)
\left.\frac{\partial^{m-1}v}{\partial\tau^{m-1}}\right|_{\tau=0}+\xi^i
\label{5.158}
\end{equation}
which in view of \ref{5.140} is equivalent to:
\begin{equation}
\left.\frac{\partial^{m-1}\hat{\delta}^i}{\partial\tau^{m-1}}\right|_{\tau=0}=
-\frac{k}{6}\frac{1}{(m+2)}\left(S_0^i+\frac{2m^A}{(m+2)}\Omega^{\prime i}_{A 0}\right)
\left.\frac{\partial^{m-1}v}{\partial\tau^{m-1}}\right|_{\tau=0}+\xi^i
\label{5.159}
\end{equation}

To determine $\left.(\partial^{m-1} v/\partial\tau^{m-1})\right|_{\tau=0}$ we substitute 
$$v=v(\tau,\vartheta), \ \ \ \gamma=\gamma(\tau,\vartheta)$$
in the regularized identification equations 
$$\hat{F}^i((\tau,\vartheta),(v,\gamma))=0,$$
differentiate implicitly with respect to $\tau$ $m-1$ times and evaluate the result at $\tau=0$. 
This gives:
\begin{equation}
\left.\frac{\partial\hat{F}^i}{\partial v}\frac{\partial^{m-1}v}{\partial\tau^{m-1}}\right|_{\tau=0}
+\left.\frac{\partial\hat{F}^i}{\partial\gamma^A}
\frac{\partial^{m-1}\gamma^A}{\partial\tau^{m-1}}\right|_{\tau=0}
+\left.\frac{\partial^{m-1}\hat{F}^i}{\partial\tau^{m-1}}\right|_{\tau=0}
+R^i_{m-1}=0
\label{5.160}
\end{equation}
The remainder, $R^i_{m-1}$, is a sum of terms of the form:
\begin{equation}
\left.\frac{\partial^{k+p+r}\hat{F}^i}
{\partial v^k\partial\gamma^{A_1}...\partial\gamma^{A_p}\partial\tau^r}
\frac{\partial^{l_1} v}{\partial\tau^{l_1}}...\frac{\partial^{l_k}v}{\partial\tau^{l_k}}
\frac{\partial^{q_1}\gamma^{A_1}}{\partial\tau^{q_1}}...
\frac{\partial^{q_p}\gamma^{A_p}}{\partial\tau^{q_p}}\right|_{\tau=0}
\label{5.161}
\end{equation}
where $k+p+r=2,...,m-1$, $k+p\geq 1$, and the $l_i:i=1,...,k$ and the $q_j:j=1,...,p$ are 
positive integers with: 
\begin{equation}
\sum_{i=1}^k l_i + \sum_{j=1}^p q_j + r = m-1
\label{5.162}
\end{equation}
Since $r\leq m-2$ the first factor in \ref{5.161} is a known smooth function on $\partial_-{\cal B}$. 
Either $k$ or $p$ may vanish but not both. Each of the sums in \ref{5.162} is at most $m-1-r$, 
so at most $m-2$ if $r\geq 1$. On the other hand if $r=0$ then $k+p\geq 2$ so either $p=0$ in 
which case $k\geq 2$ hence $\max_i l_i\leq m-2$, or $k=0$ in which case $p\geq 2$ hence 
$\max_j q_j\leq m-2$, or both $p$ and $k$ are positive in which case  again
\begin{equation}
\max_i l_i\leq m-2, \ \ \ \max_j q_j\leq m-2
\label{5.163}
\end{equation}
Thus in all cases \ref{5.163} holds, therefore the other factors in \ref{5.161} are all known smooth 
functions on $\partial_-{\cal B}$. Consequently, we can write:
\begin{equation}
R^i_{m-1}=\xi^i
\label{5.164}
\end{equation}

From the expression for $\hat{F}^i((\tau,\vartheta),(v,\gamma))$ of Proposition 4.5 we deduce:
\begin{equation}
\frac{\partial^{m-1}\hat{F}^i}{\partial\tau^{m-1}}((0,\vartheta),(v_0(\vartheta),\gamma_0(\vartheta)))
=-S_0^i(\vartheta)l(\vartheta)\frac{\partial^{m-1}\hat{f}}{\partial\tau^{m-1}}(0,\vartheta)
-\frac{\partial^{m-1}\hat{\delta}^i}{\partial\tau^{m-1}}(0,\vartheta)+\xi^i(\vartheta)
\label{5.165}
\end{equation}
Substituting \ref{5.149} and \ref{5.159} this becomes:
\begin{equation}
\left.\frac{\partial^{m-1}\hat{F}^i}{\partial\tau^{m-1}}\right|_{\tau=0}=
\frac{k}{3(m+2)}\left(S_0^i+\frac{m^A}{(m+2)}\Omega^{\prime i}_{A 0}\right)
\left.\frac{\partial^{m-1}v}{\partial\tau^{m-1}}\right|_{\tau=0}+\xi^i
\label{5.166}
\end{equation}
Substituting then \ref{5.121}, \ref{5.122}, \ref{5.166}, and \ref{5.164}, in \ref{5.160}, results 
in the spatial vector equation:
\begin{equation}
\frac{k}{3(m+2)}\left((m+3)S_0^i+\frac{m^A}{(m+2)}\Omega^{\prime i}_{A 0}\right)
\left.\frac{\partial^{m-1}v}{\partial\tau^{m-1}}\right|_{\tau=0}
+\Omega^{\prime i}_{A 0}\left.\frac{\partial^{m-1}\gamma^A}{\partial\tau^{m-1}}\right|_{\tau=0}
+\xi^i=0
\label{5.167}
\end{equation}
The $S_0$-component of this equation determines 
$\left.(\partial^{m-1}v/\partial\tau^{m-1})\right|_{\tau=0}$, 
after which the $\Omega_{A 0}$-components determine 
$\left.(\partial^{m-1}\gamma^A/\partial\tau^{m-1})\right|_{\tau=0}$, 
in terms of known smooth functions on $\partial_-{\cal B}$. Substituting the first 
in \ref{5.145} then so determines $\left.T^m\lambdab\right|_{\partial_-{\cal B}}$. 

To complete the derivation of the boundary conditions for the $m+1$ th derived data on 
$\partial_-{\cal B}$, what remains to be done is to determine 
$\left.T^m s_{NT}\right|_{\partial_-{\cal B}}$. Applying $T^m$ to \ref{5.69} along ${\cal K}$ 
we obtain:
\begin{eqnarray}
&&T^m s_{NT}=N^\mu T^{m+1}\beta^\prime_\mu+N^\mu T^{m+1}\triangle\beta_\mu\label{5.168}\\
&&\hspace{15mm}+\sum_{k=0}^{m-1}\left(\begin{array}{c}m\\k\end{array}\right)
(T^{m-k} N^\mu)T^{k+1}\beta^\prime_\mu\nonumber\\
&&\hspace{15mm}
+\sum_{k=0}^{m-1}\left(\begin{array}{c}m\\k\end{array}\right)(T^{m-k} N^\mu)T^{k+1}\triangle\beta_\mu
\nonumber
\end{eqnarray}
This is to be evaluated at $\partial_-{\cal B}$. Since at this point 
$$\left.\frac{\partial^k v}{\partial\tau^k}\right|_{\tau=0}, \ \ 
\left.\frac{\partial^k\gamma^A}{\partial\tau^k}\right|_{\tau=0} \ \ \mbox{: for $k=0,...,m-1$}$$
are already known smooth functions on $S^{n-1}$, by \ref{5.136} with $m+1$ in the role of $m$ we have:
\begin{equation}
\left. T^{m+1}\beta^\prime_\mu\right|_{\partial_-{\cal B}}=
(m+1)\left.\frac{\partial^m v}{\partial\tau^m}\right|_{\tau=0}
\left.\Lb\beta_\mu\right|_{\partial_-{\cal B}}+\xi_\mu
\label{5.169}
\end{equation}
Since $\left.N^\mu\Lb\beta_\mu\right|_{\partial_-{\cal B}}=\left.\lambda\mbox{tr}\sss\right|_{\partial_-{\cal B}}=0$ we then obtain simply:
\begin{equation}
\left.N^\mu T^{m+1}\beta^\prime_\mu\right|_{\partial_-{\cal B}}=\xi
\label{5.170}
\end{equation}
Applying $T^{m+1}$ to the expression \ref{4.4} for $\triangle\beta_\mu$ and evaluating the result at 
$\partial_-{\cal B}$, taking into account the facts that $\epb$, $\ep$ as well as $T\ep$ vanish 
at $\partial_-{\cal B}$ leads to the conclusion that:
\begin{equation}
\left.T^{m+1}\triangle\beta_\mu\right|_{\partial_-{\cal B}}=
-\left.\frac{h_{\mu\nu}}{2c}\right|_{\partial_-{\cal B}}
\left(\left.N^\nu T^{m+1}\epb\right|_{\partial_-{\cal B}}
+\left.\Nb^\nu T^{m+1}\ep\right|_{\partial_-{\cal B}}\right)+\xi_\mu
\label{5.171}
\end{equation}
By Proposition 4.2, 
$$\ep=-j\epb^2$$
Applying $T^{m+1}$ to this and evaluating the result at $\partial_-{\cal B}$ only the terms where 
each of the two $\epb$ factors receives at least one $T$ survive. This means that each of these 
factors receives at most $m$ $T$. Now at this point (see \ref{5.138}) the 
$\left.T^k\epb\right|_{\partial_-{\cal B}}$ for all $k\leq m$ have already been determined. 
It follows that the 2nd term in parenthesis in \ref{5.171} is a known smooth function on 
$\partial_-{\cal B}$. Therefore \ref{5.171} simplifies to:
\begin{equation}
\left.T^{m+1}\triangle\beta_\mu\right|_{\partial_-{\cal B}}=
\left.-\frac{h_{\mu\nu}}{2c}N^\nu T^{m+1}\epb\right|_{\partial_-{\cal B}}+\xi_\mu
\label{5.172}
\end{equation}
which implies:
\begin{equation}
\left.N^\mu T^{m+1}\triangle\beta_\mu\right|_{\partial_-{\cal B}}=\xi
\label{5.173}
\end{equation}
In view of \ref{5.170}, \ref{5.173} and the fact that the two sums in \ref{5.168} are known smooth 
functions on $\partial_-{\cal B}$, we conclude through \ref{5.168} that $T^m s_{NT}$ is determined 
as a smooth function on $\partial_-{\cal B}$. This completes the derivation of the boundary 
conditions for the $m+1$ th derived data on $\partial_-{\cal B}$. The inductive step is then complete 
and the inductive argument shows that the boundary conditions on $\partial_-{\cal B}$ for the 
derived data up to any desired finite order are determined in terms of known smooth functions on 
$S^{n-1}$. The linear system \ref{5.66} then determines the derived data of up to any desired finite 
order as known smooth functions on $\Cb_0$. 

\pagebreak

\chapter{The Variation Fields}

\section{The Bi-Variational Stress}

We have seen in Section 1.2 that a variation $\dot{\phi}$ through solutions of a solution $\phi$ 
of the nonlinear wave equation \ref{1.70} satisfies the equation of variation \ref{1.92}, which 
is simply the linear wave equation \ref{1.97} corresponding to the conformal acoustical metric 
$\tilde{h}$ of \ref{1.94}. The conformal factor $\Omega$ is given by \ref{1.96} in the case of 
3 spatial dimensions and in the general case of $n$ spatial dimensions by \ref{2.108}. 
In Section 1.2 we defined the variational stress $\dot{T}$ associated to a variation $\dot{\phi}$  
as a $T^1_1$ type tensorfield which is a quadratic form in $d\dot{\phi}$ given by \ref{1.98}. 
The divergence of $\dot{T}$ with respect to $\tilde{h}$ verifies the identity \ref{1.99}, hence 
by virtue of the equation of variation for any variation through solutions $\dot{\phi}$ the 
corresponding tensorfield $\dot{T}$ is divergence-free (\ref{1.100}). 

Now, ${\cal M}$ being the spacetime manifold,  at an arbitrary point $x\in{\cal M}$, 
$\dot{T}(x)$ is actually a quadratic form in $d\dot{\phi}(x)=\dot{v}\in T^*_x{\cal M}$, 
with values in the linear transformations of $T_x{\cal M}$, given according to \ref{1.98} by:
\begin{equation}
\dot{T}(\dot{v})=\dot{v}\otimes\dot{v}^{\sharp_{\tilde{h}}}-
\frac{1}{2}(\dot{v}\cdot\dot{v}^{\sharp_{\tilde{h}}})I_x
\label{6.1}
\end{equation}
$I_x$ being the identity transformation in $T_x{\cal M}$. Here we view the metric $\tilde{h}$ at 
$x$ as an isomorphism $\tilde{h}_x \ : \ T_x{\cal M}\rightarrow T^*_x{\cal M}$ and we denote by 
$\dot{v}^{\sharp_{\tilde{h}}}$ the inverse image of $\dot{v}$ under this isomorphism. Then the 
variational stress tensorfield associated to $\dot{\phi}$ is actually $\dot{T}\circ d\dot{\phi}$, 
$d\dot{\phi}$ being a section of $T^*{\cal M}$, that is, a 1-form on ${\cal M}$, and the identity 
\ref{1.99} reads:
\begin{equation}
\tilde{D}\cdot (\dot{T}\circ d\dot{\phi})=(\square_{\tilde{h}}\dot{\phi})d\dot{\phi}
\label{6.2}
\end{equation}

Now, $\dot{T}(x)$ being a quadratic form on $T^*_x{\cal M}$ there is a corresponding symmetric 
bilinear form on $T^*_x{\cal M}$ defined through the quadratic form by polarization. Denoting this 
symmetric bilinear form again by $\dot{T}(x)$, we have, for an arbitrary pair $\dot{v}^\prime, 
\dot{v}^{\prime\prime}\in T^*_x{\cal M}$:
\begin{equation}
\dot{T}(\dot{v}^\prime,\dot{v}^{\prime\prime})=
\frac{1}{2}(\dot{v}^\prime\otimes\dot{v}^{\prime\prime\sharp_{\tilde{h}}}
+\dot{v}^{\prime\prime}\otimes\dot{v}^{\prime\sharp_{\tilde{h}}}
-\tilde{h}^{-1}(\dot{v}^\prime,\dot{v}^{\prime\prime})I_x)
\label{6.3}
\end{equation}
noting that:
$$\dot{v}^\prime\cdot\dot{v}^{\prime\prime\sharp_{\tilde{h}}}=
\tilde{h}^{-1}(\dot{v}^\prime,\dot{v}^{\prime\prime})=
\dot{v}^{\prime\prime}\cdot\dot{v}^{\prime\sharp_{\tilde{h}}}$$
Then if $\dot{\phi}^\prime, \dot{\phi}^{\prime\prime}$ is a pair of variations we define the 
associated {\em bi-variational stress} to be the $T^1_1$-type tensorfield 
$\dot{T}\circ(d\dot{\phi}^\prime,d\dot{\phi}^{\prime\prime})$. In terms of arbitrary local coordinates 
on ${\cal M}$ this is given by: 
\begin{equation}
\dot{T}^\mu_\nu=\frac{1}{2}
\left(\partial_\nu\dot{\phi}^\prime(\tilde{h}^{-1})^{\mu\kappa}\partial_\kappa\dot{\phi}^{\prime\prime}
+\partial_\nu\dot{\phi}^{\prime\prime}(\tilde{h}^{-1})^{\mu\kappa}\partial_\kappa\dot{\phi}^\prime
-(\tilde{h}^{-1})^{\kappa\lambda}\partial_\kappa
\dot{\phi}^\prime\partial_\lambda\dot\phi^{\prime\prime}\delta^\mu_\nu\right)
\label{6.4}
\end{equation}
In analogy with \ref{6.2} we readily deduce the identity:
\begin{equation}
\tilde{D}\cdot (\dot{T}\circ (d\dot{\phi}^\prime,d\dot{\phi}^{\prime\prime})
=\frac{1}{2}\left((\square_{\tilde{h}}\dot{\phi}^\prime)d\dot{\phi}^{\prime\prime}
+(\square_{\tilde{h}}\dot{\phi}^{\prime\prime})d\dot{\phi}^\prime\right)
\label{6.5}
\end{equation}
that is, in terms of arbitrary local coordinates, 
\begin{equation}
\tilde{D}_\mu\dot{T}^\mu_\nu=\frac{1}{2}\left((\square_{\tilde{h}}\dot{\phi}^\prime)
\partial_\nu\dot{\phi}^{\prime\prime}+(\square_{\tilde{h}}\dot{\phi}^{\prime\prime})
\partial_\nu\dot{\phi}^\prime\right)
\label{6.6}
\end{equation}
In particular, if both $\dot{\phi}^\prime$, $\dot{\phi}^{\prime\prime}$ are variations through 
solutions, hence both satisfy the equation of variation \ref{1.92}, the associated bi-variational 
stress tensorfield $\dot{T}\circ(d\dot{\phi}^\prime,d\dot{\phi}^{\prime\prime})$ is divergence-free. 

Recall from Section 1.2 that the nonlinear wave equation \ref{1.70} is the Euler-Lagrange equation 
corresponding to the Lagrangian density \ref{1.74}. The Lagrangian is actually 
\begin{equation}
{\cal L}=Ld\mu_g
\label{6.7}
\end{equation}
$d\mu_g$ being the volume form of the background Minkowski metric $g$. The Lagrangian is actually 
defined on the bundle 
\begin{equation}
{\cal I}^{+*}=\bigcup_{x\in{\cal M}}{\cal I}^{*+}_x
\label{6.8}
\end{equation}
where ${\cal I}^{*+}_x$ is the interior of the positive cone in $T^*_x{\cal M}$ defined by $g_x$. 
The bundle ${\cal I}^{*+}$ is an open subbundle of $T^*{\cal M}$. Then $L(x)$ is a function on 
${\cal I}^{*+}_x$ while ${\cal L}(x)$ is an assignment to each 
$v\in {\cal I}^{*+}_x$ of a totally antisymmetric $m$-linear form in $T_x{\cal M}$, where 
$m=\mbox{dim}{\cal M}$. 
Given then a function $\phi$ on ${\cal M}$ such that at each $x\in{\cal M}$ $d\phi(x)$ belongs 
to ${\cal I}^{*+}_x$, the composition ${\cal L}\circ d\phi$ is a top degree form on ${\cal M}$, 
and for any domain ${\cal D}\subset{\cal M}$ the action corresponding to $\phi$ and to ${\cal D}$ 
is the integral 
\begin{equation}
\int_{{\cal D}}{\cal L}\circ d\phi
\label{6.9}
\end{equation}

The bi-variational stress is actually a concept of the general Lagrangian theory of mappings 
of a manifold ${\cal M}$ into another manifold ${\cal N}$ (see the monograph [Ch-A]). A Lagrangian 
in this theory is defined on the bundle 
\begin{equation}
{\cal V}=\bigcup_{(x,y)\in {\cal M}\times {\cal N}}{\cal V}_{(x,y)}
\label{6.10}
\end{equation}
where for each $(x,y)\in {\cal M}\times{\cal N}$, ${\cal V}_{(x,y)}$ is an open subset of 
${\cal L}(T_x{\cal M},T_y{\cal N})$, the space of linear maps of $T_x{\cal M}$ into $T_y{\cal N}$. 
The Lagrangian ${\cal L}$ at $(x,y)\in {\cal M}\times{\cal N}$ is an assignment to each 
$v\in {\cal V}_{(x,y)}$ of a totally antisymmetric $m$-linear form in $T_x{\cal M}$, 
$m=\mbox{dim}{\cal M}$. To represent things more simply we may choose a volume form $\ep$ on 
${\cal M}$ and write
\begin{equation}
{\cal L}(x,y)=L(x,y)\ep(x)
\label{6.11}
\end{equation}
where $L(x,y)$ is a function on ${\cal V}_{(x,y)}$. For each $(x,y)\in {\cal M}\times{\cal N}$ we 
define:
\begin{equation}
h(x,y)=\frac{\partial^2 L(x,y)}{\partial v^2}
\label{6.12}
\end{equation}
At each $v\in {\cal V}_{(x,y)}$ this is a symmetric bilinear form on 
${\cal L}(T_x{\cal M},T_y{\cal N})$. 
Given local coordinates $(x^\mu:\mu=1,...,m)$ in ${\cal M}$ and $(y^a:a=1,...,n)$ in ${\cal N}$, 
$n=\mbox{dim}{\cal N}$, we can expand an arbitrary $v\in {\cal V}_{(x,y)}$ as:
\begin{equation}
v=v_\mu^a\left.dx^\mu\right|_x\otimes\left.\frac{\partial}{\partial y^a}\right|_y
\label{6.13}
\end{equation}
The coefficients $(v_\mu^a:\mu=1,...,m; a=1,...,n)$ of the expansion constitute linear coordinates on 
${\cal V}_{(x,y)}$. We can similarly expand an arbitrary 
$\dot{v}\in {\cal L}(T_x{\cal M},T_y{\cal N})$. In these coordinates on ${\cal V}_{(x,y)}$, $h(x,y)$ 
is represented by:
\begin{equation}
h^{\mu\nu}_{ab}(x,y)=\frac{\partial^2 L(x,y)}{\partial v_\mu^a\partial v_\nu^b}
\label{6.14}
\end{equation}
Given a mapping 
$\phi \ : \ {\cal M}\rightarrow{\cal N}$ such that at each $x\in{\cal M}$, $d\phi(x)$ belongs 
to ${\cal V}_{(x,\phi(x))}$, the composition ${\cal L}\circ d\phi$ is a top degree form on ${\cal M}$, 
and for any domain ${\cal D}\subset{\cal M}$ the action corresponding to $\phi$ and to ${\cal D}$ 
is the integral 
\begin{equation}
\int_{{\cal D}}{\cal L}\circ d\phi
\label{6.15}
\end{equation}
Also, given such a mapping $\phi$, $h\circ d\phi$ is at each $x\in {\cal M}$ 
a symmetric bilinear form on 
${\cal L}(T_x{\cal M},T_{\phi(x)}{\cal N})$. In particular we can take $\phi$ to be a solution 
of the Euler-Lagrange equations deriving from the action \ref{6.15}. When such a background 
solution $\phi$ is understood fixed, we may write $h(x)$ for $h(d\phi(x))$. 

The {\em bi-variational stress} $\dot{T}(x)$ at an arbitrary point $x\in {\cal M}$ is then a 
symmetric bilinear form on ${\cal L}(T_x{\cal M},T_{\phi(x)}{\cal N})$ with values in the linear 
transformations on $T_x{\cal M}$. Thus for an arbitrary pair 
$(\xi,X)\in T^*_x{\cal M}\times T_x{\cal M}$, $\dot{T}(x)(\xi,X)$ is a symmetric bilinear form 
on ${\cal L}(T_x{\cal M},T_{\phi(x)}{\cal N})$ depending linearly on $\xi$ and on $X$. 
It is given by:
\begin{eqnarray}
&&\dot{T}(x)(\xi,X)=\frac{1}{2}\left(h(x)(\xi\otimes\dot{v}^\prime(X),\dot{v}^{\prime\prime})
+h(x)(\dot{v}^\prime,\xi\otimes\dot{v}^{\prime\prime})
-\xi(X)h(x)(\dot{v}^\prime,\dot{v}^{\prime\prime})\right)\nonumber\\
&&\hspace{30mm}: \ \forall \dot{v}^\prime, \dot{v}^{\prime\prime} \in 
{\cal L}(T_x{\cal M},T_{\phi(x)}{\cal N})
\label{6.16}
\end{eqnarray}
This has already been defined in [Ch-A] and called there {\em Noether transform} of the symmetric 
bilinear form $h(x)$. In the above presentation we chose a volume form $\ep$ on ${\cal M}$ 
(see \ref{6.11}) to represent things more simply. But then $\dot{T}$ depends on the choice of $\ep$. 
The {\em actual bi-variational stress} $\dot{{\cal T}}(x)$ at $x\in {\cal M}$ 
is a symmetric bilinear form on 
${\cal L}(T_x{\cal M},T_{\phi(x)}{\cal N})$ with values in the linear maps of $T_x{\cal M}$ into 
$\wedge_{m-1}(T_x{\cal M})$, the totally antisymmetric $m-1$-linear forms on $T_x{\cal M}$. This 
is related to $\dot{T}(x)$ as follows. If $X$ and $X_1,...,X_{m-1}$ are arbitrary vectors at $x$, 
then:
\begin{equation}
(\dot{{\cal T}}(x)\cdot X)(X_1,...,X_{m-1})=\ep(\dot{T}(x)\cdot X, X_1,...,X_{m-1})
\label{6.17}
\end{equation}
This is independent of the choice of volume form $\ep$ on ${\cal M}$. In terms of components in local 
coordinates, 
\begin{equation}
\dot{T}^\mu_\nu=\frac{1}{2}\left(h^{\mu\lambda}_{ab}
\dot{v}^{\prime a}_\nu\dot{v}^{\prime\prime b}_\lambda
+h^{\mu\lambda}_{ab}\dot{v}^{\prime\prime a}_\nu\dot{v}^{\prime b}_\lambda
-\delta^\mu_\nu h^{\kappa\lambda}_{ab}\dot{v}^{\prime a}_\kappa\dot{v}^{\prime\prime b}_\lambda\right)
\label{6.18}
\end{equation}
and:
\begin{equation}
\dot{{\cal T}}_{\nu,\lambda_1 . . . \lambda_{m-1}}=
\dot{T}^\mu_\nu\ep_{\mu \lambda_1 . . . \lambda_{m-1}}
\label{6.19}
\end{equation}
In the framework of the general Lagrangian theory of mappings of a manifold ${\cal M}$ into another 
manifold ${\cal N}$, a variation $\dot{\phi}$ of a background solution $\phi$ of the Euler-Lagrange 
equations is a section of the pullback bundle $\phi^*T{\cal N}$. To define the derivative of 
$\dot{\phi}$ at $x\in {\cal M}$ as a linear map of $T_x{\cal M}$ into $T_{\phi(x)}{\cal N}$ 
requires a connection in $T{\cal N}$ to define the horizontal part of $d\dot{\phi}$ and subtract it 
to obtain the vertical part which defines $D\dot{\phi}$, the covariant derivative of $\dot{\phi}$ 
with respect to the pullback connection on $\phi^*T{\cal N}$. To formulate the equations of variation 
satisfied by a variation through solutions $\dot{\phi}$ in fact requires the choice of such a 
connection. If $\dot{\phi}^\prime$, $\dot{\phi}^{\prime\prime}$ is a pair of variations of the 
background solution $\phi$ then $\dot{T}\circ(D\dot{\phi}^\prime,D\dot{\phi}^{\prime\prime})$ is 
a $T^1_1$-type tensorfield on ${\cal M}$. Given any vectorfield $X$ on ${\cal M}$, 
$\dot{T}\circ(D\dot{\phi}^\prime,D\dot{\phi}^{\prime\prime})\cdot X$ is a vectorfield on ${\cal M}$ 
and $\dot{{\cal T}}\circ(D\dot{\phi}^\prime,D\dot{\phi}^{\prime\prime})\cdot X$ is a $m-1$-form on 
${\cal M}$. The point is that if $\dot{\phi}^\prime$, $\dot{\phi}^{\prime\prime}$ are 
variations through solutions, then by virtue of the equations of variation 
\begin{equation}
d(\dot{{\cal T}}\circ(D\dot{\phi}^\prime,D\dot{\phi}^{\prime\prime})\cdot X)
\label{6.20}
\end{equation} 
does not involve the 2nd derivatives of $\dot{\phi}^\prime$, $\dot{\phi}^{\prime\prime}$. 
Equivalently, the divergence of the vectorfield 
$\dot{T}\circ(D\dot{\phi}^\prime,D\dot{\phi}^{\prime\prime})\cdot X$ does not contain these 2nd 
derivatives. Integrating then \ref{6.20} on a domain in ${\cal M}$ and applying Stokes' theorem, 
the integral of $\dot{{\cal T}}\circ(D\dot{\phi}^\prime,D\dot{\phi}^{\prime\prime})\cdot X$ on the boundary 
of the domain appears and a useful integral identity for  pairs of variations through solutions 
results. This topic belongs to the general theory of {\em compatible currents} expounded in [Ch-A]. 
In the present monograph shall not pursue the general theory any further. 

\vspace{5mm}

\section{The Variation Fields $V$ 
and the Associated 1-Forms $\s^{(V)}\theta^\mu$}

From this point on and throughout the remainder of this monograph indices from the middle of the 
Greek alphabet, that is $\mu, \nu, \kappa, \lambda$, shall refer to rectangular coordinates  
relative to an underlying Minkowskian structure on the spacetime manifold, while indices from 
the beginning of the Greek alphabet, that is $\alpha, \beta, \gamma, \delta$, shall refer to 
components relative to arbitrary local coordinates on the spacetime manifold. The indices 
from the middle of the Greek alphabet shall be raised and lowered using the Minkowski metric $g$, 
while the indices from the beginning of the Greek alphabet shall be raised and lowered using the 
conformal acoustical metric $\tilde{h}$. 

Let then $x^\mu:\mu=0,...,n$ be a system of rectangular coordinates and 
\begin{equation}
X_{(\mu)}=\frac{\partial}{\partial x^\mu} \ : \ \mu=0,...,n
\label{6.21}
\end{equation}
the corresponding vectorfields (see \ref{2.99}). These generate the spacetime translations, which are isometries 
of the underlying Minkowskian structure. As discussed in Section 2.4 these generate variations 
through solutions 
\begin{equation}
\dot{\phi}_\mu=X_{(\mu)}\phi=\beta_\mu \ : \ \mu=0,...,n
\label{6.22}
\end{equation}
which coincide with the rectangular components of the 1-form $\beta$ (see \ref{2.109}). Taking then 
in \ref{6.4} $\dot{\phi}^\prime=\dot{\phi}_\mu$, $\dot{\phi}^{\prime\prime}=\dot{\phi}_\nu$ we obtain 
the {\em bi-variational stress complex}, a symmetric matrix $\dot{T}_{\mu\nu}$ of $T^1_1$-type tensorfields on 
${\cal N}$ (the spacetime manifold is here ${\cal N}$, the domain of the new solution) given 
in terms of components in arbitrary local coordinates on ${\cal N}$ by:
\begin{equation}
\dot{T}^\alpha_{\beta,\mu\nu}=\frac{1}{2}\left(\partial_\beta\dot{\phi}_\mu(\tilde{h}^{-1})^{\alpha\gamma}
\partial_\gamma\dot{\phi}_\nu+\partial_\beta\dot{\phi}_\nu(\tilde{h}^{-1})^{\alpha\gamma}
\partial_\gamma\dot{\phi}_\mu-\delta^\alpha_\beta(\tilde{h}^{-1})^{\gamma\delta}
\partial_\gamma\dot{\phi}_\mu\partial_\delta\dot{\phi}_\nu\right)
\label{6.23}
\end{equation}
The $\dot{\phi}_\mu \ : \ \mu=0,...,n$ being variations through solutions satisfy the equation 
of variations \ref{1.97}, that is, we have:
\begin{equation}
\square_{\tilde{h}}\dot{\phi}_\mu=0 \ : \ \mu=0,...,n
\label{6.24}
\end{equation}
In Chapter 8 we shall consider higher order variations generated from these 1st order variations 
by applying strings of commutation fields. These higher order variations will satisfy an 
inhomogeneous version of \ref{1.97}. The bi-variational stress complex is of the same form 
\ref{6.23} for these higher order variations. With a view toward extending the treatment to the 
higher order variations to be considered in Chapter 8, we allow the $\dot{\phi}_\mu \ : \ \mu=0,...,n$ 
to satisfy the inhomogeneous equations:
\begin{equation}
\square_{\tilde{h}}\dot{\phi}_\mu=\rho_\mu \ : \ \mu=0,...,n
\label{6.25}
\end{equation}
Then by \ref{6.6} the divergence of the bi-variational stress complex is given by:
\begin{equation}
\tilde{D}_\alpha\dot{T}^\alpha_{\beta,\mu\nu}=\frac{1}{2}\left(\rho_\mu\partial_\beta\dot{\phi}_\nu
+\rho_\nu\partial_\beta\dot{\phi}_\mu\right)
\label{6.26}
\end{equation}

The usefulness of the concept of bi-variational stress in the context of a free boundary problem  
like ours, the free boundary here being the shock hypersurface ${\cal K}$, is in conjuction with 
the concept of variation fields to be presently introduced. A {\em variation field} is here simply a 
vectorfield $V$ on ${\cal N}$ which along ${\cal K}$ is tangential to ${\cal K}$. This can be 
expanded in terms of the translation fields $X_{(\mu)} \ : \ \mu=0,...,n$ (see \ref{6.21})) as:
\begin{equation}
V=V^\mu X_{(\mu)}
\label{6.27}
\end{equation}
The coefficients $V^\mu : \mu=0,...,n$ of the expansion are simply the components of $V$ in the 
rectangular coordinate system $(x^\mu:\mu=0,...,n)$. To the variation field $V$ we associate the 
column of 1-forms:
\begin{equation}
\s^{(V)}\theta^\mu=dV^\mu \ : \ \mu=0,...,n
\label{6.28}
\end{equation}
We call this column the {\em structure form} of $V$. To the variation field $V$ and to a solution 
$\phi$ of the nonlinear wave equation \ref{1.70} we associate the 1-form:
\begin{equation}
\s^{(V)}\xi=V^\mu d\dot{\phi}_\mu
\label{6.29}
\end{equation}
We say that this is associated to $\phi$ because the $\dot{\phi}_\mu : \mu=0,...,n$ are the 
variations of $\phi$ through translations. We remark that the 1-form $\s^{(V)}\xi$ is independent 
of the choice of rectangular coordinate system. For, if $(x^{\prime\mu}:\mu=0,...,n)$ is another 
rectangular coordinate system, there is a constant column $c$ and a constant matrix $\Lambda$ 
such that:
\begin{equation}
x^{\prime\mu}=c^\mu+\Lambda^\mu_\nu x^\nu
\label{6.30}
\end{equation}
and we have:
\begin{equation}
g_{\mu\nu}\Lambda^\mu_\kappa \Lambda^\nu_\lambda=g_{\kappa\lambda}
\label{6.31}
\end{equation}
that is, $\Lambda\in O(1,n)$. Therefore:
\begin{equation}
V^{\prime\mu}=Vx^{\prime\mu}=\Lambda^\mu_\nu Vx^\nu=\Lambda^\mu_\nu V^\nu
\label{6.32}
\end{equation}
while $X_{(\mu)}=\Lambda^\nu_\mu X^\prime_{(\nu)}$, 
that is:
\begin{equation}
X^\prime_{(\mu)}=(\Lambda^{-1})^\nu_\mu X_{(\nu)}
\label{6.33}
\end{equation}
hence:
\begin{equation}
\dot{\phi}^\prime_\mu=X^\prime_{(\mu)}\phi=(\Lambda^{-1})^\nu_\mu X_{(\nu)}\phi
=(\Lambda^{-1})^\nu_\mu\dot{\phi}_\nu
\label{6.34}
\end{equation}
Consequently,
\begin{equation}
V^{\prime\mu}d\dot{\phi}^\prime_\mu=\Lambda^\mu_\nu V^\nu(\Lambda^{-1})^\kappa_\mu d\dot{\phi}_\kappa
=\delta^\kappa_\nu V^\nu d\dot{\phi}_\kappa=V^\nu d\dot{\phi}_\nu
\label{6.35}
\end{equation}

Given now a variation field $V$, consider the $T^1_1$-type tensorfield $\s^{(V)}S$ given by:
\begin{equation}
\s^{(V)}S^\alpha_\beta=\dot{T}^\alpha_{\beta,\mu\nu}V^\mu V^\nu
\label{6.36}
\end{equation}
Substituting for $\dot{T}^\alpha_{\beta,\mu\nu}$ from \ref{6.23} we see that by virtue of the 
definition \ref{6.29} we have:
\begin{equation}
\s^{(V)}S^\alpha_\beta=\s^{(V)}\xi_\beta(\tilde{h}^{-1})^{\alpha\gamma}\s^{(V)}\xi_\gamma
-\frac{1}{2}\delta^\alpha_\beta(\tilde{h}^{-1})^{\gamma\delta}\s^{(V)}\xi_\gamma\s^{(V)}\xi_\delta
\label{6.37}
\end{equation}
In terms of the convention that indices from the beginning of the Greek alphabet are raised and 
lowered using $\tilde{h}$, this can be written more simply as:
\begin{equation}
\s^{(V)}S^\alpha_\beta=\s^{(V)}\xi_\beta\s^{(V)}\xi^\alpha
-\frac{1}{2}\delta^\alpha_\beta\s^{(V)}\xi_\gamma\s^{(V)}\xi^\gamma
\label{6.38}
\end{equation}
From the definition \ref{6.29} we have, in view of the definition \ref{6.28}, 
\begin{equation}
d\s^{(V)}\xi=\theta^\mu\wedge d\dot{\phi}_\mu 
\label{6.39}
\end{equation}
which can be written in the form:
\begin{equation}
\tilde{D}_\alpha\s^{(V)}\xi_\beta-\tilde{D}_\beta\s^{(V)}\xi_\alpha=
\s^{(V)}\theta^\mu_\alpha \ \partial_\beta\dot{\phi}_\mu-
\s^{(V)}\theta^\mu_\beta \ \partial_\alpha\dot{\phi}_\mu
\label{6.40}
\end{equation}
Also, by \ref{6.25} we have:
\begin{equation}
\tilde{D}_\alpha\s^{(V)}\xi^\alpha=V^\mu\rho_\mu+
\s^{(V)}\theta^\mu_\alpha \ \partial^\alpha\dot{\phi}_\mu 
\label{6.41}
\end{equation}
Using \ref{6.40} and \ref{6.41} we deduce from \ref{6.38}:
\begin{eqnarray}
&&\tilde{D}_\alpha\s^{(V)}S^\alpha_\beta=\s^{(V)}\xi^\alpha
\left(\s^{(V)}\theta^\mu_\alpha \ \partial_\beta\dot{\phi}_\mu-
\s^{(V)}\theta^\mu_\beta \ \partial_\alpha\dot{\phi}_\mu\right)
+\s^{(V)}\xi_\beta \ \s^{(V)}\theta^\mu_\alpha \ \partial^\alpha\dot{\phi}_\mu \nonumber\\
&&\hspace{20mm}+\s^{(V)}\xi_\beta \ V^\mu\rho_\mu
\label{6.42}
\end{eqnarray}
a formula which can also be deduced through \ref{6.36} using \ref{6.26} and \ref{6.23}. 

A basic requirement on the set of variation fields $V$ is that they span the tangent space to 
${\cal K}$ at each point. The simplest way to achieve this is to choose one of the variation 
fields, which we shall denote by $Y$, to be at each point of ${\cal N}$ in the linear span of 
$N$ and $\Nb$ and along ${\cal K}$ colinear to $T$, and to choose the 
other variation fields so that at each point of ${\cal N}$ they span the tangent space to 
the surface $S_{\ub,u}$ though that point. 

We thus set:
\begin{equation}
Y=\gamma N+\ogamma\Nb
\label{6.43}
\end{equation}
In view of the 2nd of \ref{4.5}, the requirement that $Y$ is along ${\cal K}$ colinear to $T$ 
reduces to:
\begin{equation}
\ogamma=r\gamma \ \mbox{: along ${\cal K}$}
\label{6.44}
\end{equation}
With the definition \ref{6.43}, the structure form of $Y$ is by \ref{3.29} - \ref{3.34} given by:
\begin{eqnarray}
&&\s^{(Y)}\stheta^\mu=(\gamma\sk+\ogamma\skb)\cdot\sd x^\mu+(\sd\gamma+\gamma k+\ogamma\kb)N^\mu 
+(\sd\ogamma+\gamma\ok+\ogamma\okb)\Nb^\mu\nonumber\\
&&\s^{(Y)}\theta^\mu_L=(\gamma\sm+\ogamma\snb)\cdot\sd x^\mu+(L\gamma+\gamma m+\ogamma\nb)N^\mu 
+(L\ogamma+\gamma\om+\ogamma\onb)\Nb^\mu\nonumber\\
&&\s^{(Y)}\theta^\mu_{\Lb}=(\gamma\sn+\ogamma\smb)\cdot\sd x^\mu+(\Lb\gamma+\gamma n+\ogamma\mb)
N^\mu+(\Lb\ogamma+\gamma\on+\ogamma\omb)\Nb^\mu\nonumber\\
&&\label{6.45}
\end{eqnarray}
Here, for any variation field $V$ we denote by $\s^{(V)}\stheta^\mu$ the $S$ 1-form which is the 
restriction of $\s^{(V)}\theta^\mu$ to the $S_{\ub,u}$. Also, we have taken account of the fact that
\begin{equation}
\Omega^\mu_A=\Omega_A x^\mu=\Omega_A\cdot\sd x^\mu 
\label{6.46}
\end{equation}
Thus \ref{3.29}, \ref{3.30} can be written in the form:
\begin{eqnarray}
&&\sd N^\mu=\sk\cdot\sd x^\mu+N^\mu k+\Nb^\mu\ok\nonumber\\
&&\sd\Nb^\mu=\skb\cdot\sd x^\mu+N^\mu\kb+\Nb^\mu\okb
\label{6.47}
\end{eqnarray}
\ref{3.31}, \ref{3.32} in the form:
\begin{eqnarray}
&&LN^\mu=\sm\cdot\sd x^\mu+m N^\mu+\om\Nb^\mu\nonumber\\
&&\Lb\Nb=\smb\cdot\sd x^\mu+\mb N^\mu+\omb\Nb^\mu 
\label{6.48}
\end{eqnarray}
and \ref{3.33}, \ref{3.34} in the form:
\begin{eqnarray}
&&\Lb N^\mu=\sn\cdot\sd x^\mu+nN^\mu+\on\Nb^\mu\nonumber\\
&&L\Nb^\mu=\snb\cdot\sd x^\mu+\nb N^\mu+\onb\Nb^\mu 
\label{6.49}
\end{eqnarray}
Optimal estimates are obtained if we choose:
\begin{equation}
\gamma=1
\label{6.50}
\end{equation}
and require $\ogamma$ to be constant along the integral curves of $L$:
\begin{equation}
L\ogamma=0 
\label{6.51}
\end{equation}
This, together with the boundary condition \ref{6.44} defines $\ogamma$ on ${\cal N}$. With these 
choices the expressions \ref{6.45} simplify accordingly. It is clear that $r$ being a smooth function 
on ${\cal K}$ (represented as a smooth function of $(\tau,\vartheta)$), $\ogamma$ defined in this 
way is a smooth function of the acoustical coordinates. 

We turn to the definition of the remaining variation fields. As noted above, at each point of 
${\cal N}$ these must span the tangent space to the surface $S_{\ub,u}$ through that point. 
Since the translation fields $X_{(\mu)} : \mu=0,...,n$ (see \ref{6.21}) span the tangent space to the 
spacetime manifold at each point, a suitable definition for these $S$-tangential variation fields is:
\begin{equation}
E_{(\mu)}=\Pi X_{(\mu)} \ : \ \mu=0,...,n
\label{6.52}
\end{equation}
where $\Pi$ is the $h$-orthogonal projection to the $S_{\ub,u}$. The rectangular components of 
the vectorfields $E_{(\mu)}$ are given by:
\begin{equation}
E_{(\mu)}^\nu=\Pi_\mu^\nu=\delta_\mu^\nu+
\frac{1}{2c}(N^\nu\Nb^\lambda+\Nb^\nu N^\lambda)h_{\lambda\mu}
\label{6.53}
\end{equation}
They are smooth functions of the acoustical coordinates. 
The structure forms $\s^{(E_{(\mu)})}\theta^\nu$ of the variation fields $E_{(\mu)}$ are deduced 
using \ref{3.84} - \ref{3.86}, and \ref{6.47} - \ref{6.49}. 

In the case $n=2$ there is a more natural choice for the remaining variation fields. What we have 
called surfaces $S_{\ub,u}$ reduce in this case to curves, and the unit tangent field $E$ to 
these oriented curves is the natural choice for the remaining variation field (see Section 3.5). 
The structure forms $\s^{(E)}\theta^\mu$ are then given by \ref{3.a31}, \ref{3.a33}, \ref{3.a38}:
\begin{eqnarray}
&&\s^{(E)}\stheta^\mu=\se E^\mu+eN^\mu+\oe\Nb^\mu\nonumber\\
&&\s^{(E)}\theta^\mu_L=\sf E^\mu+fN^\mu+\of\Nb^\mu\nonumber\\
&&\s^{(E)}\theta^\mu_{\Lb}=\sfb E^\mu+\fb N^\mu+\ofb\Nb^\mu 
\label{6.54}
\end{eqnarray} 

The following remark is intended to ensure that a proper non-relativistic limit 
can be extracted. Note that the variation fields $E_{(\mu)}$ or, in the case $n=2$, $E$, have the 
physical dimensions of $(\mbox{length})^{-1}$. On the other hand, by reason of the conditions 
$N^0=\Nb^0=1$ (see \ref{2.58}) the variation field $Y$ has in the non-relativistic theory the 
physical dimensions of $(\mbox{time})^{-1}$. However, adding the energies associated to the 
different variation fields is meaningful only when the different variation fields have the same 
physical dimensions. Choosing $\eta_0$ to be the sound speed in a reference constant state, we 
must then replace $Y$ as defined above by $\eta_0^{-1}Y$, which has the physical dimensions of 
$(\mbox{length})^{-1}$, in its role as a variation field. In the non-relativistic theory we can 
choose the relation of the unit of time to the unit of length so that $\eta_0=1$, so formally 
there is no change. In the relativistic theory however, the customary choice of relation of 
these units is so that the speed of light is set equal to 1, which would make $\eta_0$ a very 
small positive constant when the reference constant state is non-relativistic. Once we have  
made the point however it is straightforward to make the necessary changes and we shall ignore the 
factor $\eta_0^{-1}$ in the sequel. 

We shall presently show that, for the fundamental variations \ref{6.22}, control of all components 
of the associated 1-form $\s^{(V)}\xi$ 
for all variation fields $V$ provides control on all components of the symmetric 
2-covariant tensorfield $s$ (see \ref{3.36}). Consider first the case $n=2$, which is simpler. 
In this case, control of $\s^{(E)}\xi$ means control of 
\begin{eqnarray}
&&\s^{(E)}\xi(E)=\sss \label{6.a1}\\
&&\s^{(E)}\xi(L)=\rho\ss_N \label{6.a2}\\
&&\s^{(E)}\xi(\Lb)=\rhob\ss_{\Nb} \label{6.a3}
\end{eqnarray}
($\sss$, $\ss_N$, $\ss_{\Nb}$ being in this case the scalars \ref{3.a17}). Control of $\s^{(Y)}\xi$ 
means control of 
\begin{eqnarray}
&&\s^{(Y)}\xi(E)=\ss_N+\ogamma\ss_{\Nb} \label{6.a4}\\
&&\s^{(Y)}\xi(L)=s_{NL}+\ogamma s_{\Nb L} \label{6.a5}\\
&&\s^{(Y)}\xi(\Lb)=s_{N\Lb}+\ogamma s_{\Nb\Lb} \label{6.a6}
\end{eqnarray}
(see \ref{6.43}, \ref{6.50}). Now, by \ref{3.a18}:
\begin{equation}
s_{\Nb L}=\rho s_{\Nb N}=c\rho\sss, \ \ \ s_{N\Lb}=\rhob s_{N\Nb}=c\rhob\sss
\label{6.a7}
\end{equation}
Since $\sss$ is already controlled through \ref{6.a1}, we acquire control on $s_{NL}$ through 
\ref{6.a5} and on $s_{\Nb\Lb}$ through \ref{6.a6}. We thus achieve control on all components of 
$s$. For $n>2$ a similar argument applies with the variation fields $E_{(\mu)}:\mu=0,...,n$ in the 
role of $E$. That is, control of $\s^{(E_{(\mu)})}\xi$ for all $\mu=0,...,n$ means control of 
\begin{eqnarray}
&&\s^{(E_{(\mu)})}\sxi=E_{(\mu)}\cdot\sss \ \ \mbox{for all $\mu=0,...,n$} \label{6.a8}\\
&&\s^{(E_{(\mu)})}\xi(L)=\rho\ss_N(E_{(\mu)}) \ \  \mbox{for all $\mu=0,...,n$} \label{6.a9}\\
&&\s^{(E_{(\mu)})}\xi(\Lb)=\rhob\ss_{\Nb}(E_{(\mu)}) \ \ \mbox{for all $\mu=0,...,n$} \label{6.a10}
\end{eqnarray}
Since the $E_{(\mu)} \ : \ \mu=0,...,n$ span the tangent space to the $S_{\ub,u}$ at each point, 
we acquire full control on $\sss$ through \ref{6.a8}, on $\ss_N$ through \ref{6.a9} and on 
$\ss_{\Nb}$ through \ref{6.a10}. Also, control of $\s^{(Y)}\xi$ means control of 
\begin{eqnarray}
&&\s^{(Y)}\sxi=\ss_N+\ogamma\ss_{\Nb} \label{6.a11}\\
&&\s^{(Y)}\xi(L)=s_{NL}+\ogamma s_{\Nb L} \label{6.a12}\\
&&\s^{(Y)}\xi(\Lb)=s_{N\Lb}+\ogamma s_{\Nb\Lb} \label{6.a13}
\end{eqnarray}
By \ref{3.42}:
\begin{equation}
s_{\Nb L}=\rho s_{\Nb N}=c\rho\mbox{tr}\sss, \ \ \ s_{N\Lb}=\rhob s_{N\Nb}=s\rhob\mbox{tr}\sss
\label{6.a14}
\end{equation}
Since $\sss$ is already controlled, \ref{6.a12} then gives us control on $s_{NL}$ and \ref{6.a13} 
gives us control on $s_{\Nb\Lb}$. Thus again we achieve control on all components of 
$s$.

\vspace{5mm}

\section{The Fundamental Energy Identities}

Given now a vectorfield $X$, which we call {\em multiplier field}, we consider the vectorfield 
$\s^{(V)}P$ associated to $X$ and to a given variation field $V$ through $\s^{(V)}S$, defined by:
\begin{equation}
\s^{(V)}P^\alpha=-\s^{(V)}S^\alpha_\beta X^\beta
\label{6.55}
\end{equation}
We call $\s^{(V)}P$ the {\em energy current} associated to $X$ and to $V$. 
Let us denote by $\s^{(V)}Q$ the covariant divergence of $\s^{(V)}P$ (with respect to the conformal 
acoustical metric $\tilde{h}$):
\begin{equation}
\tilde{D}_\alpha\s^{(V)}P^\alpha=\s^{(V)}Q
\label{6.56}
\end{equation}
We have:
\begin{equation}
\s^{(V)}Q=-\s^{(V)}S^\alpha_\beta\tilde{D}_\alpha X^\beta-X^\beta\tilde{D}_\alpha\s^{(V)}S^\alpha_\beta
\label{6.57}
\end{equation}
Consider the $T^2_0$-type tensorfield corresponding to $\s^{(V)}S$ through the metric $\tilde{h}$:
\begin{equation}
\s^{(V)}S^{\alpha\beta}=\s^{(V)}S^\alpha_\gamma(\tilde{h}^{-1})^{\gamma\beta}
\label{6.58}
\end{equation}
Substituting from \ref{6.38} we obtain:
\begin{equation}
\s^{(V)}S^{\alpha\beta}=\s^{(V)}\xi^\alpha\s^{(V)}\xi^\beta-\frac{1}{2}(\tilde{h}^{-1})^{\alpha\beta}
(\tilde{h}^{-1})^{\gamma\delta}\s^{(V)}\xi_\gamma\s^{(V)}\xi_\delta
\label{6.59}
\end{equation}
which shows the $T^2_0$-type tensorfield to be {\em symmetric}. As a consequence, the first term 
on the right in \ref{6.57} can be written in the form:
\begin{equation}
-\frac{1}{2}\s^{(V)}S^{\alpha\beta}\left(\tilde{h}_{\alpha\gamma}\tilde{D}_\beta X^\gamma
+\tilde{h}_{\beta\gamma}\tilde{D}_\alpha X^\gamma\right)=-\frac{1}{2}\s^{(V)}S^{\alpha\beta}
\s^{(X)}\tilde{\pi}_{\alpha\beta}
\label{6.60}
\end{equation}
Here we denote by $\s^{(X)}\tilde{\pi}$ the {\em deformation tensor} of $X$ relative to $\tilde{h}$: 
\begin{equation}
\s^{(X)}\tilde{\pi}={\cal L}_X\tilde{h}
\label{6.61}
\end{equation}
that is, the rate of change of the metric $\tilde{h}$ under the flow generated by $X$ on ${\cal N}$. 
The last equality in \ref{6.62} is the identity:
\begin{equation}
({\cal L}_X\tilde{h})(Y,Z)=\tilde{h}(Y,\tilde{D}_Z X)+\tilde{h}(Z,\tilde{D}_Y X)
\label{6.62}
\end{equation}
$\tilde{D}$ being the covariant derivative operator associated to $\tilde{h}$. Substituting also for 
$\tilde{D}_\alpha\s^{(V)}S^\alpha_\beta$ from \ref{6.42} in the second term on the right in \ref{6.57} 
then yields:
\begin{equation}
\s^{(V)}Q=\s^{(V)}Q_1+\s^{(V)}Q_2+\s^{(V)}Q_3
\label{6.63}
\end{equation}
where:
\begin{eqnarray}
&&\s^{(V)}Q_1=-\frac{1}{2}\s^{(V)}S^{\alpha\beta}\s^{(X)}\tilde{\pi}_{\alpha\beta}
\label{6.64}\\
&&\s^{(V)}Q_2=-\s^{(V)}\xi^\alpha X^\beta
\left(\s^{(V)}\theta^\mu_\alpha \ \partial_\beta\dot{\phi}_\mu-
\s^{(V)}\theta^\mu_\beta \ \partial_\alpha\dot{\phi}_\mu\right)\nonumber\\
&&\hspace{15mm}-\s^{(V)}\xi_\beta X^\beta\ \s^{(V)}\theta^\mu_\alpha \ \partial^\alpha\dot{\phi}_\mu 
\label{6.65}\\
&&\s^{(V)}Q_3=-\s^{(V)}\xi_\beta X^\beta\ V^\mu\rho_\mu
\label{6.66}
\end{eqnarray}

Equation \ref{6.56} in arbitrary local coordinates reads:
\begin{equation}
\frac{1}{\sqrt{-\mbox{det}\tilde{h}}}\partial_\alpha
\left(\sqrt{-\mbox{det}\tilde h}\s^{(V)}P^\alpha\right)=
\s^{(V)}Q
\label{6.67}
\end{equation}
In particular in acoustical coordinates $(\ub,u,\vartheta^A:A=1,...,n-1)$, expanding $\s^{(V)}P$ in 
the associated coordinate frame field
$$\left(\frac{\partial}{\partial\ub},\frac{\partial}{\partial u},\frac{\partial}{\partial\vartheta^A}:
A=1,...,n-1\right)$$
we have:
\begin{equation}
\s^{(V)}P=\s^{(V)}P^{\ub}\frac{\partial}{\partial\ub}+\s^{(V)}P^u\frac{\partial}{\partial u}
+\s^{(V)}P^{\vartheta^A}\frac{\partial}{\partial\vartheta^A}
\label{6.68}
\end{equation}
Also, by \ref{1.94}, \ref{2.44}, 
\begin{equation}
\sqrt{-\mbox{det}\tilde{h}}=\Omega^{(n+1)/2}\sqrt{-\mbox{det}h}, \ \ \ 
\sqrt{-\mbox{det}h}=2a\sqrt{\mbox{det}\sh}
\label{6.69}
\end{equation}
Thus, equation \ref{6.67} expressed in acoustical coordinates reads:
\begin{eqnarray}
&&\frac{1}{\sqrt{\mbox{det}\sh}}\left\{
\frac{\partial}{\partial\ub}\left(2a\Omega^{(n+1)/2}\sqrt{\mbox{det}\sh}\s^{(V)}P^{\ub}\right)
+\frac{\partial}{\partial u}\left(2a\Omega^{(n+1)/2}\sqrt{\mbox{det}\sh}\s^{(V)}P^u\right)\right.\nonumber\\
&&\hspace{15mm}\left.+\frac{\partial}{\partial\vartheta^A}\left(2a\Omega^{(n+1)/2}\sqrt{\mbox{det}\sh}\s^{(V)}P^{\vartheta^A}\right)
\right\}=2a\Omega^{(n+1)/2}\s^{(V)}Q
\label{6.70}
\end{eqnarray}
Let $\s^{(V)}\sP$ be the $S$-vectorfield:
\begin{equation}
\s^{(V)}\sP=\s^{(V)}P^{\vartheta^A}\frac{\partial}{\partial\vartheta^A}
\label{6.71}
\end{equation}
We then have:
\begin{equation}
\frac{1}{\sqrt{\mbox{det}\sh}}\frac{\partial}{\partial\vartheta^A}
\left(2a\Omega^{(n+1)/2}\sqrt{\mbox{det}\sh}\s^{(V)}P^{\vartheta^A}\right)
=\sdiv(2a\Omega^{(n+1)/2}\s^{(V)}\sP)
\label{6.72}
\end{equation}
where we denote by $\sdiv$ the covariant divergence on $(S_{\ub,u},\sh)$. Integrating equation  
\ref{6.70} on $S_{\ub,u}$, represented as $S^{n-1}$, with the measure:
\begin{equation}
d\mu_{\sh}=\sqrt{\mbox{det}\sh}d\vartheta^1 . . . d\vartheta^{n-1}
\label{6.73}
\end{equation}
and noting that 
\begin{equation}
\int_{S_{\ub,u}}\sdiv(2a\Omega^{(n+1)/2}\s^{(V)}\sP)d\mu_{\sh}=0
\label{6.74}
\end{equation}
we obtain:
\begin{eqnarray}
&&\frac{\partial}{\partial\ub}\left(\int_{S_{\ub,u}}2a\Omega^{(n+1)/2}\s^{(V)}P^{\ub}d\mu_{\sh}\right)
+\frac{\partial}{\partial u}\left(\int_{S_{\ub,u}}2a\Omega^{(n+1)/2}\s^{(V)}P^u d\mu_{\sh}\right)
\nonumber\\
&&\hspace{55mm}=\int_{S_{\ub,u}}2a\Omega^{(n+1)/2}\s^{(V)}Q
\label{6.75}
\end{eqnarray}
Now, we can also expand the vectorfield $\s^{(V)}P$ in the frame field 
$(L,\Lb,\Omega_A:A=1,...,n-1)$:
\begin{equation}
\s^{(V)}P=\s^{(V)}P^L L+\s^{(V)}P^{\Lb}\Lb+\s^{(V)}P^A\Omega_A
\label{6.76}
\end{equation}
Since (see \ref{2.33}, \ref{2.36})
\begin{equation}
L=\frac{\partial}{\partial\ub}-b^A\frac{\partial}{\partial\vartheta^A}, \ \ \ 
\Lb=\frac{\partial}{\partial u}+b^A\frac{\partial}{\partial\vartheta^A}, \ \ \ 
\Omega_A=\frac{\partial}{\partial\vartheta^A} \ : \ A=1,...,n-1
\label{6.77}
\end{equation} 
comparing \ref{6.76} with \ref{6.68} we conclude that:
\begin{equation}
P^{\ub}=P^L, \ \ \ P^u=P^{\Lb}, \ \ \ 
P^{\vartheta^A}=P^A+b^A(P^{\Lb}-P^L) \ : \ A=1,...,n-1
\label{6.78}
\end{equation}
As a consequence, equation \ref{6.75} can be written in the form:
\begin{eqnarray}
&&\frac{\partial}{\partial\ub}\left(\int_{S_{\ub,u}}2a\Omega^{(n+1)/2}\s^{(V)}P^Ld\mu_{\sh}\right)
+\frac{\partial}{\partial u}\left(\int_{S_{\ub,u}}2a\Omega^{(n+1)/2}\s^{(V)}P^{\Lb} d\mu_{\sh}\right)
\nonumber\\
&&\hspace{55mm}=\int_{S_{\ub,u}}2a\Omega^{(n+1)/2}\s^{(V)}Q
\label{6.79}
\end{eqnarray}

Let now, with $u_1\geq\ub_1\geq 0$, $R_{\ub_1,u_1}$ be the trapezoid in the $(\ub,u)$ plane defined by:
\begin{equation}
R_{\ub_1,u_1}=\{(\ub,u) \ : \ u\in[\ub,u_1], \ \ub\in [0,\ub_1]\}
\label{6.80}
\end{equation}
We can also express $R_{\ub_1,u_1}$ as the union of the triangle: 
\begin{equation}
T_{\ub_1}=\{(\ub,u) \ : \ \ub\in[0,u], \ u\in[0,\ub_1]\}=R_{\ub_1,\ub_1}
\label{6.81}
\end{equation}
with the rectangle $[0,\ub_1]\times[\ub_1,u_1]$, the intersection of which is the segment 
$[0,\ub_1]\times\{\ub_1\}$:
\begin{equation}
R_{\ub_1,u_1}=T_{\ub_1}\bigcup\left([0,\ub_1]\times[\ub_1,u_1]\right)
\label{6.82}
\end{equation}
We denote by ${\cal R}_{\ub_1,u_1}$ the region in ${\cal N}$ corresponding to $R_{\ub_1,u_1}$:
\begin{equation}
{\cal R}_{\ub_1,u_1}=\bigcup_{(\ub,u)\in R_{\ub_1,u_1}}S_{\ub,u}
\label{6.83}
\end{equation}
This region is represented by $R_{\ub_1,u_1}\times S^{n-1}$. 

We integrate equation \ref{6.79} on the trapezoid $R_{\ub_1,u_1}$. Let us denote: 
\begin{eqnarray}
&&\Fb(\ub,u)=\int_{S_{\ub,u}}2a\Omega^{(n+1)/2}\s^{(V)}P^L d\mu_{\sh}\label{6.84}\\
&&F(\ub,u)=\int_{S_{\ub,u}}2a\Omega^{(n+1)/2}\s^{(V)}P^{\Lb} d\mu_{\sh}\label{6.85}
\end{eqnarray}
With the notation \ref{6.84} the first term on the left in \ref{6.79} is:
$$\frac{\partial\Fb}{\partial\ub}$$
To integrate this on $R_{\ub_1,u_1}$ we use the representation \ref{6.82}. We have:
$$\int_{T_{\ub_1}}\frac{\partial\Fb}{\partial\ub}d\ub du=\int_0^{\ub_1}\{\Fb(u,u)-\Fb(0,u)\}du,$$
$$\int_{[0,\ub_1]\times[\ub_1,u_1]}\frac{\partial\Fb}{\partial\ub}d\ub du=
\int_{\ub_1}^{u_1}\{\Fb(\ub_1,u)-\Fb(0,u)\}du$$
Therefore we obtain:
\begin{equation}
\int_{R_{\ub_1,u_1}}\frac{\partial\Fb}{\partial\ub}d\ub du=\int_{\ub_1}^{u_1}\Fb(\ub_1,u)du
+\int_0^{\ub_1}\Fb(u,u)du-\int_0^{u_1}\Fb(0,u)du
\label{6.86}
\end{equation}
With the notation \ref{6.85} the second term on the left in \ref{6.79} is:
$$\frac{\partial F}{\partial u}$$
To integrate this on $R_{\ub_1,u_1}$ we use the representation \ref{6.80}. We have:
\begin{equation}
\int_{R_{\ub_1,u_1}}\frac{\partial F}{\partial u}du d\ub=\int_0^{\ub_1}F(\ub,u_1)d\ub
-\int_0^{\ub_1}F(\ub,\ub)d\ub
\label{6.87}
\end{equation}
Denoting then also:
\begin{equation}
G(\ub,u)=\int_{S_{\ub,u}}2a\Omega^{(n+1)/2}\s^{(V)}Q d\mu_{\sh}
\label{6.88}
\end{equation}
equation \ref{6.79} integrated over $R_{\ub_1,u_1}$ is:
\begin{eqnarray}
&&\int_0^{\ub_1}F(\ub,u_1)d\ub+\int_{\ub_1}^{u_1}\Fb(\ub_1,u)du+\int_0^{\ub_1}(\Fb-F)(\tau,\tau)d\tau
\nonumber\\ 
&&\hspace{25mm}-\int_0^{u_1}\Fb(u,u)du=\int_{R_{\ub_1,u_1}}G(\ub,u)du d\ub
\label{6.89}
\end{eqnarray}
Recall that on $C_u$ the range of $\ub$ is $[0,u]$. We denote by $C_u^{\ub_1}$ the part of $C_u$ 
obtained by imposing the restriction $\ub\leq\ub_1$ on the range of $\ub$. So if $\ub_1\geq u$ this 
is no restriction and we have $C_u^{\ub_1}=C_u$. In view of \ref{6.85} first integral on the left in 
\ref{6.89} is an integral on $C_{u_1}^{\ub_1}$. Using $d\mu_{\sh} d\ub$ as the measure on the $C_u$, 
thus denoting, for an arbitrary function $f$ on ${\cal N}$, 
\begin{equation}
\int_{C_u^{\ub_1}}f=\int_0^{\ub_1}\left(\int_{S_{\ub,u}}fd\mu_{\sh}\right)d\ub
\label{6.90}
\end{equation}
the first integral on the left in \ref{6.89} is:
\begin{equation}
\s^{(V)}{\cal E}^{\ub_1}(u_1):=\int_{C_{u_1}^{\ub_1}}2a\Omega^{(n+1)/2}\s^{(V)}P^{\Lb}
\label{6.91}
\end{equation}
Recall that on $\Cb_{\ub}$ we have $u\geq\ub$. We denote by $\Cb_{\ub}^{u_1}$ the part of $\Cb_{\ub}$ 
obtained by imposing the restriction $u\leq u_1$. In view of \ref{6.84} the second integral on 
the left in \ref{6.89} is an integral on $\Cb_{\ub_1}^{u_1}$. Using $d\mu_{\sh}du$ as the measure on the 
$\Cb_{\ub}$, thus denoting, for an arbitrary function $f$ on ${\cal N}$,
\begin{equation}
\int_{\Cb_{\ub}^{u_1}}f=\int_{\ub}^{u_1}\left(\int_{S_{\ub,u}}fd\mu_{\sh}\right)du
\label{6.92}
\end{equation}
the second integral on the left in \ref{6.89} is:
\begin{equation}
\s^{(V)}\underline{{\cal E}}^{u_1}(\ub_1):=\int_{\Cb_{\ub_1}^{u_1}}2a\Omega^{(n+1)/2}\s^{(V)}P^L 
\label{6.93}
\end{equation}
Similarly, the fourth integral on the left in \ref{6.89} is:
\begin{equation}
\s^{(V)}\underline{{\cal E}}^{u_1}(0):=\int_{\Cb_0^{u_1}}2a\Omega^{(n+1)/2}\s^{(V)}P^L
\label{6.94}
\end{equation}
In regard to the third integral on the left in \ref{6.89}, in view of \ref{4.52} this is an integral 
on ${\cal K}^{\ub_1}$, where we denote by ${\cal K}^{\tau_1}$ the part of ${\cal K}$ obtained by 
imposing the restriction $\tau\leq\tau_1$. Using $d\mu_{\sh}d\tau$ as the measure on ${\cal K}$, 
thus denoting, for an arbitrary function $f$ on ${\cal N}$,
\begin{equation}
\int_{{\cal K}^{\tau_1}}=\int_0^{\tau_1}\left(\int_{S_{\tau,\tau}}fd\mu_{\sh}\right)d\tau
\label{6.95}
\end{equation}
in view of \ref{6.84}, \ref{6.85}, the third integral 
on the left in \ref{6.89} is:
\begin{equation}
\s^{(V)}{\cal F}^{\ub_1}:=
\int_{{\cal K}^{\ub_1}}2a\Omega^{(n+1)/2}\left(\s^{(V)}P^L-\s^{(V)}P^{\Lb}\right)
\label{6.96}
\end{equation}
Finally, using $d\mu_{\sh}dud\ub$ as the measure on ${\cal N}$, thus denoting for an arbitrary function 
$f$ on ${\cal N}$,
\begin{eqnarray}
&&\int_{{\cal R}_{\ub_1,u_1}}f=\int_0^{\ub_1}\left\{\int_{\ub}^{u_1}\left(\int_{S_{\ub,u}}d\mu_{\sh}
\right)du\right\}d\ub\label{6.97}\\
&&\hspace{10mm}=\int_0^{\ub_1}\left\{\int_0^u\left(\int_{S_{\ub,u}}d\mu_{\sh}\right)d\ub\right\}du
+\int_{\ub_1}^{u_1}\left\{\int_0^{\ub_1}\left(\int_{S_{\ub,u}}d\mu_{\sh}\right)d\ub\right\}du
\nonumber
\end{eqnarray}
according to the representations \ref{6.80} and \ref{6.82} respectively, the integral on the right 
in \ref{6.89} is:
\begin{equation}
\s^{(V)}{\cal G}^{\ub_1,u_1}:=\int_{{\cal R}_{\ub_1,u_1}}2a\Omega^{(n+1)/2}\s^{(V)}Q
\label{6.98}
\end{equation}
With the above definitions \ref{6.91}, \ref{6.93}, \ref{6.94}, \ref{6.96}, \ref{6.98}, equation 
\ref{6.89} reads:
\begin{equation}
\s^{(V)}{\cal E}^{\ub_1}(u_1)+\s^{(V)}\underline{{\cal E}}^{u_1}(\ub_1)+\s^{(V)}{\cal F}^{\ub_1}
-\s^{(V)}\underline{{\cal E}}^{u_1}(0)=\s^{(V)}{\cal G}^{\ub_1,u_1}
\label{6.99}
\end{equation}
This is the {\em fundamental energy identity} corresponding to the variation field $V$ and to the 
multiplier field $X$. We may call $\s^{(V)}{\cal E}^{\ub_1}(u)$ {\em energy} associated to 
$C_u^{\ub_1}$, $\s^{(V)}\underline{{\cal E}}^{u_1}(\ub)$ {\em energy} associated to $\Cb_{\ub}^{u_1}$, 
and we may call $\s^{(V)}{\cal F}^\tau_1$ {\em flux} associated to ${\cal K}^{\tau_1}$. 

Consider the integrants in the energies \ref{6.91} and \ref{6.93}. These integrants are the positive 
factor $2a\Omega^{(n+1)/2}$ times $\s^{(V)}P^{\Lb}$ and $\s^{(V)}P^L$ respectively. Taking the 
$\tilde{h}$-inner product of \ref{6.76} with $L$ and with $\Lb$ we obtain:
\begin{equation}
2a\Omega\s^{(V)}P^{\Lb}=-\tilde{h}(\s^{(V)}P,L), \ \ \ 
2a\Omega\s^{(V)}P^L=-\tilde{h}(\s^{(V)}P,\Lb)
\label{6.100}
\end{equation}
Now, from \ref{6.55} we have:
\begin{equation}
-\tilde{h}_{\alpha\gamma}\s^{(V)}P^\gamma=\s^{(V)}S_{\alpha\beta}X^\beta
\label{6.101}
\end{equation}
where
\begin{equation}
\s^{(V)}S_{\alpha\beta}=\tilde{h}_{\alpha\gamma}\s^{(V)}S^\gamma_\beta
=\s^{(V)}\xi_\alpha\s^{(V)}\xi_\beta-\frac{1}{2}\tilde{h}_{\alpha\beta}(\tilde{h}^{-1})^{\gamma\delta}
\s^{(V)}\xi_\gamma\s^{(V)}\xi_\delta
\label{6.102}
\end{equation}
This is the $T^0_2$-type tensorfield corresponding to $\s^{(V)}S$ through $\tilde{h}$. It is symmetric. 
Denoting by $\s^{(V)}S$ this tensorfield rather than the original $T^1_1$-type tensorfield, we have 
by \ref{6.101}, for an arbitrary vectorfield $Y$ on ${\cal N}$, 
\begin{equation}
-\tilde{h}(\s^{(V)}P,Y)=\s^{(V)}S(X,Y)
\label{6.103}
\end{equation}
Taking in particular $Y=L$ and $Y=\Lb$ we obtain by \ref{6.100}:
\begin{equation}
2a\Omega\s^{(V)}P^{\Lb}=\s^{(V)}S(X,L), \ \ \ 2a\Omega\s^{(V)}P^L=\s^{(V)}S(X,\Lb)
\label{6.104}
\end{equation}
As a consequence of the following proposition, the integrants in the energies \ref{6.91} and \ref{6.93} 
are positive-semidefinite if the multiplier field $X$ is chosen to be acoustically timelike 
future-directed. 

\vspace{2.5mm}

\noindent{\bf Proposition 6.1} \ \ \ Let $X,Y$ be a pair of future-directed non-spacelike vectors 
at a point of $({\cal N},\tilde{h})$. Then $\s^{(V)}S(X,Y)$ is a positive-semidefinite quadratic 
form in the covector $\s^{(V)}\xi$ at this point. Moreover, if $X,Y$ are timelike then this  
quadratic form is positive-definite. 

\vspace{2.5mm}

\noindent{\em Proof:} If $X$ and $Y$ are not colinear, their linear span is a timelike plane $P$. 
If on the other hand they are colinear, we can take $P$ to be any timelike plane containing $X,Y$. 
In both cases there are future-directed null vectors $E_+,E_-$ in $P$ such that 
$$\tilde{h}(E_+,E_-)=-2$$
and we have:
$$X=X^+E_++X^-E_-,\ \ \ Y=Y^+E_++Y^-E_-$$
where $X^+,X^-,Y^+,Y^-$ are all non-negative, all positive if $X,Y$ are timelike.  
Also, $P\bot$, the $\tilde{h}$-orthogonal complement of $P$, is a spacelike plane. We denote by 
$\tilde{\sh}$ the restriction of $\tilde{h}$ to $P\bot$, a positive-definite quadratic form on $P\bot$. 
Let us denote:
$$\s^{(V)}S_{++}=\s^{(V)}S(E_+,E_+), \ \ \ \s^{(V)}S_{+-}=\s^{(V)}S(E_+,E_-), \ \ \ 
\s^{(V)}S_{--}=\s^{(V)}S(E_-,E_-)$$
We then have:
$$\s^{(V)}S(X,Y)=X^+Y^+\s^{(V)}S_{++}+(X^+Y^-+X^-Y^+)\s^{(V)}S_{+-}+X^-Y^-\s^{(V)}S_{--}$$
Note that the coefficients of $\s^{(V)}S_{++},\s^{(V)}S_{+-},\s^{(V)}S_{--}$ here are all non-negative, 
all positive if $X,Y$ are timelike. Therefore the proposition follows if we show that $\s^{(V)}S_{++}$, 
$\s^{(V)}S_{+-}$, $\s^{(V)}S_{--}$ are all positive-semidefinite and their sum is positive-definite. 
In fact, denoting 
$$\s^{(V)}\xi_+=\s^{(V)}\xi(E_+), \ \ \ \s^{(V)}\xi_-=\s^{(V)}\xi(E_-)$$
from \ref{6.102} we have:
$$\s^{(V)}S_{++}=(\s^{(V)}\xi_+)^2, \ \ \ \s^{(V)}S_{--}=(\s^{(V)}\xi_-)^2$$
and, denoting by $\s^{(V)}\sxi$ the restriction of $\s^{(V)}\xi$ to $P\bot$, we have:
$$\tilde{h}^{-1}(\s^{(V)}\xi,\s^{(V)}\xi)=-\s^{(V)}\xi_+\s^{(V)}\xi_-
+\tilde{\sh}^{-1}(\s^{(V)}\sxi,\s^{(V)}\sxi)$$
hence:
$$\s^{(V)}S_{+-}=\tilde{\sh}^{-1}(\s^{(V)}\sxi,\s^{(V)}\sxi)$$
This being a positive-definite quadratic form on $(P\bot)^*$, the proposition follows. 

\vspace{2.5mm}

As for the integrant in the flux \ref{6.96}, by \ref{6.104} it is equal to the positive factor 
$\Omega^{(n-1)/2}$ times:
\begin{equation}
\s^{(V)}S(X,M), \ \ \mbox{where: $M=\Lb-L$}
\label{6.105}
\end{equation}
is a normal to ${\cal K}$ pointing to the interior of ${\cal N}$. In the next chapter we shall 
investigate the conditions on the multiplier field $X$ which make this flux integrant coercive 
when $\s^{(V)}\xi$ satisfies on ${\cal K}$ the boundary condition to be discussed in the next section. 

\vspace{5mm}

\section{The Boundary Condition on ${\cal K}$
for the 1-Forms $\s^{(V)}\xi$}

We recall from 
Sections 1.4 and 4.1 the covector field $\xi$ defined along ${\cal K}$ by the conditions that at each point on ${\cal K}$ 
the null space of $\xi$ is the tangent space to ${\cal K}$ at the point, $\xi(V)>0$ for any vector 
$V$ at the point which points to the interior of ${\cal N}$, and $\xi$ is normalized to be of unit magnitude with respect to the Minkowski metric $g$. Thus $\xi^\sharp$, the vectorfield along ${\cal K}$ corresponding to $\xi$ 
through $g$, is the unit normal to ${\cal K}$ with respect to $g$ pointing to the interior of 
${\cal N}$. To avoid confusion with the 1-form $\s^{(V)}\xi$,  
we denote here and in the following the covector field $\xi$ just referred to by $\zeta$. 
We also recall the function $\delta$ along ${\cal K}$ defined in the last part of Section 4.1. In the 
new notation, \ref{4.47}, which contains the definition of $\delta$, reads:
\begin{equation}
\triangle\beta=\delta\zeta
\label{6.106}
\end{equation} 

Consider now the jump condition, 1st of \ref{1.195}, which in the new notation reads:
\begin{equation}
\zeta_\mu\triangle I^\mu=0
\label{6.107}
\end{equation}
$I^\mu$ being the components of the particle current:
\begin{equation}
I^\mu=-G\beta^\mu 
\label{6.108}
\end{equation}
(see \ref{1.53}, \ref{1.72}). As we discussed in Section 1.6, \ref{6.107} is equivalent to the 
nonlinear jump condition \ref{1.329}. Let $V$ be an arbitrary variation field. Then 
$$V^\mu\frac{\partial}{\partial x^\mu}$$
is along ${\cal K}$ a differential operator interior to ${\cal K}$. Therefore \ref{6.107} implies:
\begin{equation}
V^\mu\frac{\partial}{\partial x^\mu}(\zeta_\nu\triangle I^\nu)=0 \ \ \mbox{: on ${\cal K}$}
\label{6.109}
\end{equation}
We have, from \ref{6.108}, 
\begin{eqnarray}
&&\frac{\partial I^\nu}{\partial x^\mu}=-(g^{-1})^{\nu\lambda}G\partial_\mu\beta_\lambda
+2\frac{dG}{d\sigma}\beta^\nu\beta^\lambda\partial_\mu\beta_\lambda\nonumber\\
&&\hspace{8mm}=-G(h^{-1})^{\nu\lambda}\partial_\mu\beta_\lambda
=-G(h^{-1})^{\nu\lambda}\partial_\lambda\beta_\mu 
\label{6.110}
\end{eqnarray}
by the 1st of \ref{4.20} and by \ref{1.d10}. Hence:
\begin{equation}
V^\mu\frac{\partial I^\nu}{\partial x^\mu}=-G(h^{-1})^{\nu\lambda}\s^{(V)}\xi_\lambda
\label{6.111}
\end{equation}
Thus, along ${\cal K}$ we have:
\begin{equation}
V^\mu\frac{\partial I_\pm^\nu}{\partial x^\mu}=-G_\pm(h_\pm^{-1})^{\mu\nu}\s^{(V)}\xi_{\pm\mu}
\label{6.112}
\end{equation}
Let us define along ${\cal K}$:
\begin{equation}
K_\pm^\mu=G_\pm(h_\pm^{-1})^{\mu\nu}\zeta_\nu
\label{6.113}
\end{equation}
These are normals to ${\cal K}$ relative to the acoustical metric corresponding to the future and the 
past solutions respectively, pointing to the interior of ${\cal N}$. Then by \ref{6.112} 
we have, along ${\cal K}$:
\begin{equation}
\zeta_\nu V^\mu\frac{\partial\triangle I^\nu}{\partial x^\mu}=-K_+^\mu\s^{(V)}\xi_{+\mu}
+K_-^\mu\s^{(V)}\xi_{-\mu}
\label{6.114}
\end{equation}
Consider next 
$$V^\mu\frac{\partial\zeta_\nu}{\partial x^\mu}\triangle I^\nu$$
which is the remainder of the left hand side of \ref{6.109}. By \ref{6.106},
\begin{eqnarray*}
&&V^\mu\frac{\partial\zeta_\nu}{\partial x^\mu}=V^\mu\frac{\partial}{\partial x^\mu}
\left(\frac{\triangle\beta_\nu}{\delta}\right)\\
&&\hspace{12mm}=\frac{V^\mu}{\delta}\frac{\partial\triangle\beta_\nu}{\partial x^\mu}
-\frac{V^\mu}{\delta^2}\triangle\beta_\nu\frac{\partial\delta}{\partial x^\mu}
\end{eqnarray*}
that is:
\begin{equation}
V^\mu\frac{\partial\zeta_\nu}{\partial x^\mu}=
\frac{V^\mu}{\delta}\frac{\partial\triangle\beta_\nu}{\partial x^\mu}
-\frac{V^\mu}{\delta}\frac{\partial\delta}{\partial x^\mu}\zeta_\nu
\label{6.115}
\end{equation}
If we multiply this by $\triangle I^\nu$ the 2nd term does not contribute by virtue of the condition 
\ref{6.107}. As for the 1st term, it is:
\begin{eqnarray}
&&\frac{V^\mu}{\delta}(\partial_\mu\beta_{+\nu}-\partial_\mu\beta_{-\nu})
=\frac{V^\mu}{\delta}\left((\partial_\nu\beta_\mu)_+-(\partial_\nu\beta_\mu)_-\right)\nonumber\\
&&\hspace{20mm}=\frac{1}{\delta}(\s^{(V)}\xi_{+\nu}-\s^{(V)}\xi_{-\nu}) \label{6.a15}
\end{eqnarray}
Therefore we obtain, along ${\cal K}$:
\begin{equation}
V^\mu\frac{\partial\zeta_\nu}{\partial x^\mu}\triangle I^\nu=
\frac{\triangle I^\nu}{\delta}(\s^{(V)}\xi_{+\nu}-\s^{(V)}\xi_{-\nu})
\label{6.116}
\end{equation}
Combining \ref{6.114} with \ref{6.116} we arrive through \ref{6.109} at the following proposition. 

\vspace{2.5mm}

\noindent{\bf Proposition 6.2} \ \ \ For any variation field $V$, the associated 1-form $\s^{(V)}\xi$ 
satisfies the following boundary condition on ${\cal K}$:
$$\s^{(V)}\xi_+(A_+)=\s^{(V)}\xi_-(A_-)$$
where $A_\pm$ are the following vectorfields along ${\cal K}$:
$$A_\pm=\frac{\triangle I}{\delta}-K_\pm$$
with $K_\pm$ being the normals to ${\cal K}$ defined by the future and past solution respectively, 
the rectangular components of which are given by \ref{6.113}. 

\vspace{2.5mm}

The next proposition plays a central role in the analysis of the coercivity of the integrant in the 
flux \ref{6.96}. 

\vspace{2.5mm}

\noindent{\bf Proposition 6.3} \ \ \ The vectorfield 
$$\frac{\triangle I}{\delta}$$
defined on ${\cal K}\setminus\partial_-{\cal B}$, is tangential to ${\cal K}$ and 
timelike future-directed with respect to the acoustical metric defined by the future solution. 

\vspace{2.5mm}

\noindent{\em Proof:} We begin by showing that at each point of ${\cal K}$, the vectors 
$u_+$, $u_-$, the fluid velocity of the future solution and the past solution respectively, and 
$\zeta^\sharp$ the Minkowski unit normal to ${\cal K}$ directed toward the domain of the 
future solution, all three lie in the same plane $P$, which is timelike relative to the Minkowski 
metric. This statement has already been established in Section 1.4 in the context of the full shock 
development problem and follows directly from the jump condition \ref{1.203}. However here we must 
establish the statement in the context of the restricted shock development problem. Let again 
$\beta^\sharp$ denote the vector with rectangular components $\beta^\mu=(g^{-1})^{\mu\nu}\beta_\nu$. 
Then the linear jump condition \ref{1.328} is:
\begin{equation}
g(\triangle\beta^\sharp,X)=0 \ \ : \ \forall X\in T{\cal K}
\label{6.117}
\end{equation}
This means that there is a function $\alpha$ on ${\cal K}$ such that:
\begin{equation}
\triangle\beta^\sharp=\alpha\zeta^\sharp
\label{6.118}
\end{equation}
Since according to \ref{1.53}:
\begin{equation}
\beta^\sharp=-h u
\label{6.119}
\end{equation}
equation \ref{6.118} reads:
\begin{equation}
h_+u_+=h_-u_--\alpha\zeta^\sharp
\label{6.120}
\end{equation}
hence ($h$ being positive) at any point of ${\cal K}$ the three vectors $u_+$, $u_-$, $\zeta^\sharp$ 
at the point do lie on the same plane.
This plane $P$ is timelike relative to the Minkowski metric as it contains the timelike vectors 
$u_+$, $u_-$. 

Let us denote by $H$ the tangent hyperplane to ${\cal K}$ at the point in question. Since 
${\cal K}$ is a timelike hypersurface relative to the Minkowski metric, $H$ is a 
timelike hyperplane relative to $g$. The plane $P$ is $g$-orthogonal to $H$ as it contains the 
$g$-normal to $H$ vector $\zeta^\sharp$. The plane $P$ being timelike, there is unique unit 
future-directed timelike relative to $g$ vector $U$ in $P$ which is orthogonal to $\zeta^\sharp$. 
The vector $U$ then belongs to $H$. 
We may expand $u_\pm$ as:
\begin{equation}
u_\pm=a_\pm U+b_\pm\zeta^\sharp
\label{6.121}
\end{equation}
Now the coefficients $a_\pm$, $b_\pm$ are all positive. That this is so for $a_\pm$ is obvious by 
the future-directed nature of $u_\pm$. Noting that according to \ref{1.199}:
\begin{equation}
b_\pm=u_{\bot\pm}
\label{6.122}
\end{equation}
the positivity of $b_\pm$ is \ref{1.200}. The jump condition 1st of \ref{1.195}, which coincides with 
\ref{6.107}, is \ref{1.201}, that is:
\begin{equation}
n_+b_+=n_-b_-
\label{6.123}
\end{equation}
Then in view of \ref{1.215} we have:
\begin{eqnarray}
&&b_+<b_- \ \mbox{: in the compressive case $\triangle p>0$}\nonumber\\ 
&&b_+>b_- \ \mbox{: in the depressive case $\triangle p<0$}
\label{6.124}
\end{eqnarray}
Now, $u_\pm$ are of unit magnitude relative to $g$, hence:
\begin{equation}
a_\pm=\sqrt{1+b_\pm^2}
\label{6.125}
\end{equation}

Since 
\begin{equation}
I_\pm=n_\pm u_\pm=n_\pm a_\pm U+n_\pm b_\pm \zeta^\sharp
\label{6.126}
\end{equation}
we have, by \ref{6.123}:
\begin{equation}
\triangle I:=I_+-I_-=(n_+a_+-n_-a_-)U
\label{6.127}
\end{equation}

There are two questions to be answered here. The first question is: $U$ is timelike future-directed 
relative to $g$. Is it also timelike future-directed relative to $h_+$? Here $h_+$ is the acoustical metric 
defined along ${\cal K}$ by the future solution. The second question is: what is the sign of the 
coefficient of $U$ in \ref{6.127}? 

To answer the first question, consider $h_{+\mu\nu}U^\mu U^\nu$. 
From \ref{1.51},
\begin{equation}
h_{+\mu\nu}=g_{\mu\nu}+(1-\eta_+^2)u_{+\mu}u_{+\nu}, \ \mbox{where $u_{+\mu}=g_{\mu\nu}u_+^\nu$}
\label{6.128}
\end{equation}
From \ref{6.121}, $u_{+\mu}U^\mu=-a_+$. Hence, substituting,  
\begin{equation}
h_{+\mu\nu}U^\mu U^\nu=-1+(1-\eta_+^2)a_+^2
\label{6.129}
\end{equation}
Here we bring in the fact that ${\cal K}$ is acoustically timelike relative to the future solution 
(called 2nd determinism condition in Section 1.4), expressed by \ref{1.237}. In view of \ref{1.238} 
and \ref{6.122} this reads:
\begin{equation}
(h_+^{-1})^{\mu\nu}\zeta_\mu\zeta_\nu=1-\left(\frac{1}{\eta_+^2}-1\right)b_+^2>0
\label{6.130}
\end{equation}
In view of \ref{6.125} this is equivalent to:
\begin{equation}
a_+^2<\frac{1}{1-\eta_+^2}
\label{6.131}
\end{equation}
which by \ref{6.129} is in turn equivalent to:
\begin{equation}
h_{+\mu\nu}U^\mu U^\nu<0
\label{6.132}
\end{equation}
that is to $U$ being timelike relative to $h_+$. So either $U$ is contained in the interior of the 
future null cone of $h_+$ or else it is contained in the interior of the past null cone of $h_+$. 
However, the null cone of $h_+$ being contained in the null cone of $g$, the second alternative 
implies that $U$ is past-directed timelike relative to $g$, contradicting what we already know. 
Thus only the first alternative remains and we have answered the first question in the affirmative.

To answer the second question, the sign of $n_+a_+-n_-a_-$, we multiply by $b_+$ to see that by virtue 
of \ref{6.123} the sign in question is that of $b_-a_+-b_+a_-$, or, diving by $a_+ a_-$, that of 
$$\frac{b_-}{a_-}-\frac{b_+}{a_+}=\frac{b_-}{\sqrt{1+b_-^2}}-\frac{b_+}{\sqrt{1+b_+^2}}$$
by \ref{6.125}. Therefore by \ref{6.124} the sign is positive in the compressive case and negative 
in the depressive case. 

We conclude that $\triangle I$ is relative to $h_+$ future-directed timelike in the compressive case, 
past-directed timelike in the depressive case. Recalling the conclusion reached at the end of 
Section 4.1 that $\delta$ is positive in the compressive case and negative in the depressive case, 
the proposition then follows. 

\vspace{2.5mm}

From this point, as in Section 4.1, for quantities defined on ${\cal K}$ by the future solution 
we omit the subscript + when there can be no misunderstanding. We thus write:
\begin{equation}
\triangle I^\mu=-\triangle(G\beta^\mu)=-G_-\triangle\beta^\mu-\beta^\mu\triangle G
\label{6.133}
\end{equation}
Since the vectorfield $\triangle I$ is tangential to ${\cal K}$ it can be expanded in terms of 
the frame field $T,\Omega_A:A=1,...,n-1$:
\begin{equation}
\triangle I=(\triangle I)^T+(\triangle I)^A\Omega_A
\label{6.134}
\end{equation}
and we have:
\begin{equation}
(\triangle I)^T=-\frac{1}{4a}h(\triangle I,T), \ \ \ 
(\triangle I)^A=(\sh^{-1})^{AB}h(\triangle I,\Omega_B)
\label{6.135}
\end{equation}
Now, we have:
\begin{equation}
h_{\mu\nu}\beta^\mu=(g_{\mu\nu}+H\beta_\mu\beta_\nu)\beta^\mu=(1-\sigma H)\beta_\nu=\eta^2\beta_\nu
\label{6.136}
\end{equation}
hence:
\begin{eqnarray}
&&h_{\mu\nu}\beta^\mu T^\nu=\eta^2\beta_\nu T^\nu=\eta^2\beta_T,\nonumber\\
&&h_{\mu\nu}\beta^\mu\Omega_B^\nu=\eta^2\beta_\nu\Omega_B^\nu=\eta^2\sbeta_B
\label{6.137}
\end{eqnarray}
where in the first we have denoted $\beta_T=\beta(T)$. 
Also, we have:
\begin{equation}
\triangle\beta^\mu=\left((h^{-1})^{\mu\nu}+F\beta^\mu\beta^\nu\right)\triangle\beta_\nu
=(h^{-1})^{\mu\nu}\triangle\beta_\nu+F\beta^\mu\nu
\label{6.138}
\end{equation}
recalling the 1st of the definitions \ref{4.8}. By virtue of the linear jump condition \ref{1.328} we 
then obtain:
\begin{eqnarray}
&&h_{\mu\nu}\triangle\beta^\mu T^\nu= F\nu h_{\mu\nu}\beta^\mu T^\nu=F\nu\eta^2\beta_T,\nonumber\\
&&h_{\mu\nu}\triangle\beta^\mu\Omega_B^\nu=F\nu h_{\mu\nu}\beta^\mu\Omega_B^\nu
=F\nu\eta^2\sbeta_B
\label{6.139}
\end{eqnarray}
using \ref{6.137}. Comparing \ref{6.133} - \ref{6.135} with \ref{6.137}, \ref{6.139} we see that 
$\triangle I$ is proportional to the quantity:
\begin{equation}
\Lambda=\triangle G +F\nu G_-
\label{6.140}
\end{equation}
In fact, by \ref{6.133} and the 1st of \ref{6.135}, \ref{6.137}, \ref{6.139}: 
\begin{equation}
(\triangle I)^T=\frac{\eta^2}{4a}\Lambda\beta_T
\label{6.141}
\end{equation}
while by \ref{6.133} and the 2nd of \ref{6.135}, \ref{6.137}, \ref{6.139}:
\begin{equation}
(\triangle I)^A=-\eta^2\Lambda(\sh^{-1})^{AB}\sbeta_B
\label{6.142}
\end{equation}
In conclusion, we have deduced the following formula for the vectorfield $\triangle I$:
\begin{equation}
\triangle I=\frac{\eta^2}{a}\Lambda\left(\frac{1}{4}\beta_T T-a\sbeta^\sharp\right)
\label{6.143}
\end{equation}
where $\sbeta^\sharp$ stands for the $S$ vectorfield corresponding through $\sh$ to the $S$ 1-form 
$\sbeta$:
$$\sbeta^\sharp=(\sh^{-1})^{AB}\sbeta_B\Omega_A$$

Consider now the vectorfield $K_+$ defined by \ref{6.113}. Omitting the subscript +, it being 
clear that we are considering a quantity defined along ${\cal K}$ by the future solution, the 
definition is:
\begin{equation}
K=G\hat{K}, \  \mbox{where $\hat{K}^\mu=(h^{-1})^{\mu\nu}\zeta_\nu$}
\label{6.144}
\end{equation}
Now, the vectorfield $\hat{K}$ is an interior normal to ${\cal K}$ relative to $({\cal N},h)$. 
So is the vectorfield:
\begin{equation}
M=\Lb-L
\label{6.145}
\end{equation}
Therefore there is a positive function $\kappa$ on ${\cal K}$ such that:
\begin{equation}
M=2\kappa\hat{K} \ \mbox{: on ${\cal K}$}
\label{6.146}
\end{equation} 
Let us denote:
\begin{equation}
w=(h^{-1})^{\mu\nu}\zeta_\mu\zeta_\nu
\label{6.147}
\end{equation}
By \ref{6.130} $w$ is a positive function on ${\cal K}$. Since 
\begin{equation}
h(M,M)=4a, \ \mbox{while} \ h(\hat{K},\hat{K})=w
\label{6.148}
\end{equation}
it follows that:
\begin{equation}
w=\frac{a}{\kappa^2}
\label{6.149}
\end{equation}

Consider next \ref{6.106}. By \ref{4.4}, 
\begin{equation}
a\triangle\beta_\mu=-\frac{1}{2}h_{\mu\nu}(\rhob\epb L^\nu+\rho\ep\Lb^\nu)
\label{6.150}
\end{equation}
According to \ref{4.5}, there is a function $q$ defined by:
\begin{equation}
q=\rhob\epb=-\rho\ep 
\label{6.151}
\end{equation}
therefore \ref{6.150} takes the form:
\begin{equation}
a\triangle\beta_\mu=\frac{q}{2}h_{\mu\nu}M^\nu
\label{6.152}
\end{equation}
or, substituting for $M$ from \ref{6.146} and the definition of $\hat{K}$ from \ref{6.144}, 
\begin{equation}
a\triangle\beta=q\kappa\zeta
\label{6.153}
\end{equation}
Comparing with \ref{6.106} we conclude that:
\begin{equation}
a\delta=\kappa q
\label{6.154}
\end{equation}

Going back to Proposition 6.2, denoting $A_+$ simply by $A$ and defining  
\begin{equation}
B=\kappa G^{-1}A
\label{6.155}
\end{equation}
we have, by \ref{6.144}, \ref{6.146} and \ref{6.154}, 
\begin{equation}
B=B_{||}+B_\bot
\label{6.156}
\end{equation}
where $B_\bot$, the part of the vectorfield $B$ which is $h$-orthogonal to ${\cal K}$ (relative to 
the future solution), is simply:
\begin{equation}
B_{\bot}=-\frac{1}{2}M
\label{6.157}
\end{equation}
while $B_{||}$, the part of the vectorfield $B$ which is tangential to ${\cal K}$, is:
\begin{equation}
B_{||}=\hat{\Lambda}\left(\frac{1}{4}\beta_T T-a\sbeta^\sharp\right)
\label{6.158}
\end{equation}
with
\begin{equation}
\hat{\Lambda}=\frac{\eta^2\Lambda}{Gq}
\label{6.159}
\end{equation}

Let us also define 
\begin{equation}
B_-=\kappa G^{-1}A_-
\label{6.160}
\end{equation}
Note that here while $A_-$ refers to the past solution, $\kappa$ and $G$ refer to the future solution. 
Decomposing:
\begin{equation}
B_-=B_{-||}+B_{-\bot}
\label{6.161}
\end{equation}
where $B_{-||}$ is the part of $B_-$ which is tangential to ${\cal K}$, while $B_{-\bot}$ is 
the part of $B_-$ which is $h_-$- orthogonal to ${\cal K}$ (that is, relative to the past solution), 
$B_{-||}$ coincides with $B_{||}$:
\begin{equation}
B_{-||}=B_{||}
\label{6.162}
\end{equation}
but $B_{-\bot}$ differs from $B_\bot$, being instead given by:
\begin{equation}
B_{-\bot}=-\frac{1}{2}M^\prime_-, \ \mbox{where} \ M^\prime_-=2\kappa G^{-1}K_-
\label{6.163}
\end{equation}
with $K_-$ given by \ref{6.113}. 

Then the boundary condition for the 1-form $\s^{(V)}\xi$ of Proposition 6.2 takes the form:
\begin{equation}
\s^{(V)}\xi(B)=\s^{(V)}b \ \ \mbox{: on ${\cal K}$}
\label{6.164}
\end{equation}
where: 
\begin{equation}
\s^{(V)}b=\s^{(V)}\xi_-(B_-)
\label{6.165}
\end{equation}

We shall presently analyze the behavior of the vectorfield $B_{||}$ for small $\tau$, that is 
in a neighborhood of $\partial_-{\cal B}$ in ${\cal K}$, using the results of Chapter 4. 
Since $G_-=G-\triangle G$, we can write the definition \ref{6.140} in the form:
\begin{equation}
\Lambda=(1-F\nu)\triangle G+F\nu G
\label{6.166}
\end{equation}
Now, $\triangle G$ is expressed by \ref{4.15} as $K(\sigma,\triangle\sigma)\triangle\sigma$ with 
$K(\sigma,\triangle\sigma)$ given by \ref{4.19}. Thus,
\begin{equation}
\triangle G=G^\prime\triangle\sigma-\frac{1}{2}G^{\prime\prime}(\triangle\sigma)^2
+O((\triangle\sigma)^3)
\label{6.167}
\end{equation}
By \ref{4.20} this is:
\begin{equation}
\triangle G=\frac{1}{2}GF\triangle\sigma-\frac{1}{4}G\left(F^\prime+\frac{1}{2}F^2\right)(\triangle\sigma)^2+O((\triangle\sigma)^3)
\label{6.168}
\end{equation}
Since by \ref{4.31} and Proposition 4.1 :
\begin{equation}
\mu=O(\nu^3)
\label{6.169}
\end{equation}
\ref{4.14} implies:
\begin{equation}
\triangle\sigma=-2\nu+F\nu^2+O(\nu^3)
\label{6.170}
\end{equation}
Substituting this in \ref{6.168} gives:
\begin{equation}
\triangle G=-GF\nu-GF^\prime\nu^2+O(\nu^3)
\label{6.171}
\end{equation}
Substituting this in \ref{6.166} the terms linear in $\nu$ cancel and we obtain:
\begin{equation}
\Lambda=G(F^2-F^\prime)\nu^2+O(\nu^3)
\label{6.172}
\end{equation}
or, by \ref{4.27}, 
\begin{equation}
\Lambda=-\eta^{-4}GH^\prime\nu^2+O(\nu^3)
\label{6.173}
\end{equation}

Next, from the 1st of \ref{4.37}, in view of \ref{4.40} and Proposition 4.2, 
\begin{equation}
\nu=-\frac{\eta^2}{2c}\beta_N\epb+O(\epb^2)
\label{6.174}
\end{equation}
Actually $\nu$ is expressed as a smooth function of $\epb$ and the quadruplet (see \ref{4.43}) 
$(\sigma,c,\beta_N,\beta_{\Nb})$ along ${\cal K}$ corresponding to the new solution. To find 
the form of $\nu$ for small $\tau$ however we need to assume a bound for $T^2\epb$ in a neighborhood 
of $\partial_-{\cal B}$ in ${\cal K}$. Then, since $\epb$ vanishes at $\partial_-{\cal B}$ we have:
\begin{equation}
\epb=\left.T\epb\right|_{\partial_-{\cal B}}\tau+O(\tau^2)
\label{6.175}
\end{equation}
From \ref{4.213} and \ref{4.154}:
\begin{equation}
\left.T\epb\right|_{\partial_-{\cal B}}=\left.2s_{\Nb\Lb}\right|_{\partial_-{\cal B}}
\label{6.176}
\end{equation}
hence:
\begin{equation}
\epb=\left.2s_{\Nb\Lb}\right|_{\partial_-{\cal B}}\tau+O(\tau^2)
\label{6.177}
\end{equation}
Substituting in \ref{6.174} we obtain:
\begin{equation}
\nu=\left.-\frac{\eta^2}{c}\beta_N s_{\Nb\Lb}\right|_{\partial_-{\cal B}}\tau+O(\tau^2)
\label{6.178}
\end{equation}
Now, by \ref{4.163}:
\begin{equation}
\left.\frac{\eta^2}{c}\beta_N s_{\Nb\Lb}\right|_{\partial_-{\cal B}}=
\left.\frac{\Lb H}{H^\prime}\right|_{\partial_-{\cal B}}
\label{6.179}
\end{equation}
and according to \ref{4.112}:
\begin{equation}
\left.\frac{1}{2}\beta_N^2\Lb H\right|_{\partial_-{\cal B}}=l
\label{6.180}
\end{equation}
Therefore \ref{6.178} takes the form:
\begin{equation}
\nu=-\frac{2l}{H^\prime_0}\left.\frac{1}{\beta_N^2}\right|_{\partial_-{\cal B}}\tau+O(\tau^2) \ \ 
\mbox{where} \ H^\prime_0=\left.H^\prime\right|_{\partial_-{\cal B}}
\label{6.181}
\end{equation}
Substituting this in \ref{6.173} we obtain:
\begin{equation}
\Lambda=-\frac{4G_0 l^2}{H^\prime_0}\left.\frac{1}{\eta^4\beta_N^4}\right|_{\partial_-{\cal B}}\tau^2
+O(\tau^3) \ \ \mbox{where} \ G_0=\left.G\right|_{\partial_-{\cal B}}
\label{6.182}
\end{equation}

Next, to find the form of $\lambdab$ along ${\cal K}$ for small $\tau$ we need to assume a bound for 
$T^2\lambdab$ in a neighborhood of $\partial_-{\cal B}$ in ${\cal K}$. Then since $\lambdab$ vanishes 
at $\partial_-{\cal B}$ we have:
\begin{equation}
\lambdab=\left.T\lambda\right|_{\partial_-{\cal B}}\tau+O(\tau^2) \ \mbox{: along ${\cal K}$}
\label{6.183}
\end{equation}
Substituting the 1st of \ref{4.216} then gives:
\begin{equation}
\lambdab=-\frac{c_0k}{3 l}\tau+O(\tau^2) \ \mbox{: along ${\cal K}$}
\label{6.184}
\end{equation}
To find the form of $\lambda$ along ${\cal K}$ for small $\tau$ we need to assume a bound for 
$T^3\lambda$ in a neighborhood of $\partial_-{\cal B}$ in ${\cal K}$. Then since $\lambda$, $T\lambda$ 
both vanish at $\partial_-{\cal B}$ we have:
\begin{equation}
\lambda=\left.\frac{1}{2}T^2\lambda\right|_{\partial_-{\cal B}}\tau^2+O(\tau^3) \ 
\mbox{: along ${\cal K}$}
\label{6.185}
\end{equation}
Substituting the 2nd of \ref{4.216} then gives:
\begin{equation}
\lambda=\frac{k}{6}\tau^2+O(\tau^3) \ \mbox{: along ${\cal K}$}
\label{6.186}
\end{equation}
The above imply that along ${\cal K}$:
\begin{eqnarray}
&&\rho=-\frac{k}{3l}\tau+O(\tau^2)\nonumber\\
&&\rhob=\frac{k}{6c_0}\tau^2+O(\tau^3)\nonumber\\
&&a=-\frac{k^2}{18l}\tau^3+O(\tau^4)
\label{6.187}
\end{eqnarray}
The 2nd of these together with the fact that from \ref{6.177}, in view of \ref{6.179}, \ref{6.180} 
we have:
\begin{equation}
\epb=\frac{4c_0 l}{H^\prime_0}\left.\frac{1}{\eta^2\beta_N^3}\right|_{\partial_-{\cal B}}\tau+O(\tau^2)
\label{6.188}
\end{equation}
implies that the function $q$ defined by \ref{6.151} has for small $\tau$ the form:
\begin{equation}
q=\frac{2kl}{3H^\prime_0}\left.\frac{1}{\eta^2\beta_N^3}\right|_{\partial_-{\cal B}}\tau^3+O(\tau^4)
\label{6.189}
\end{equation}

By \ref{6.182} and \ref{6.189} we obtain from \ref{6.159}:
\begin{equation}
\hat{\Lambda}=-\frac{6l}{k}\left.\frac{1}{\beta_N}\right|_{\partial_-{\cal B}}\tau^{-1}+O(1)
\label{6.190}
\end{equation}

Consider now the expression \ref{6.158} for $B_{||}$. It constitutes a decomposition of $B_{||}$, 
a vectorfield tangential to ${\cal K}$ into a term proportional to $T$, thus $h$-orthogonal to 
the $S_{\tau,\tau}$ sections of ${\cal K}$, and a $S$-vectorfield, thus tangential to these sections. 
The coefficient of $T$ is: 
\begin{equation}
\frac{1}{4}\hat{\Lambda}\beta_T=\frac{1}{2}+O(\tau)
\label{6.191}
\end{equation}
by \ref{6.190}, the fact that $\beta_T=\rho\beta_N+\rhob\beta_{\Nb}$, and the first two of \ref{6.187}. 
Since $h(T,T)=-4a$ \ref{6.191} should be compared with $(4a)^{-1/2}$ times the part of $B_{||}$ which 
is $S$ tangential, that is with:
\begin{equation}
-\frac{1}{2}\hat{\Lambda}\sqrt{a}\sbeta^\sharp
\label{6.192}
\end{equation}
Here, the coefficient of $\sbeta^\sharp$ is:
\begin{equation}
-\frac{1}{2}\hat{\Lambda}\sqrt{a}=-\frac{\sqrt{-2l}}{2}
\left.\frac{1}{\beta_N}\right|_{\partial_-{\cal B}}\tau^{1/2}+O(\tau^{3/2})
\label{6.193}
\end{equation}
by \ref{6.190} and the 3rd of \ref{6.187}. 

We shall now derive a suitable expression for the vectorfield $M_-^\prime$ which defines $B_{-\bot}$ 
through \ref{6.163}. We have:
\begin{equation}
M_-^\prime=2\kappa G^{-1}K_-
\label{6.194}
\end{equation}
The vectorfield $K_-$ is defined by \ref{6.113}. In analogy with \ref{6.144} we write:
\begin{equation}
K_-=G_-\hat{K}_-, \ \mbox{where $\hat{K}_-^\mu=(h_-^{-1})^{\mu\nu}\zeta_\nu$}
\label{6.195}
\end{equation}
Then in view of \ref{6.146} the difference $M-M_-^\prime$ is given by:
\begin{equation}
M-M_-^\prime=2\kappa\left(\hat{K}-\frac{G_-}{G}\hat{K}_-\right)
\label{6.196}
\end{equation}
Since 
$$h^{-1})^{\mu\nu}-(h_-^{-1})^{\mu\nu}=F\beta^\mu\beta^\nu-F_-\beta_-^\mu\beta_-^\nu$$
the rectangular components of the difference $M-M_-^\prime$ are:
\begin{equation}
M^\mu-M_-^{\prime\mu}=\frac{\triangle G}{G}M^\mu 
-2\kappa\frac{G_-}{G}\left(F\beta^\mu\beta^\nu-F_-\beta_-^\mu\beta_-^\nu\right)\zeta_\nu
\label{6.197}
\end{equation}
Consequently, denoting by $J$ the vectorfield with rectangular components
\begin{equation}
J^\mu=2\kappa\left(F\beta^\mu\beta^\nu-F_-\beta_-^\mu\beta_-^\nu\right)\zeta_\nu
\label{6.198}
\end{equation}
equation \ref{6.196} takes the form:
\begin{equation}
M-M_-^\prime=\frac{\triangle G}{G}M-\frac{G_-}{G}J
\label{6.199}
\end{equation}
We have:
\begin{equation}
\beta^\mu\beta^\nu-\beta_-^\mu\beta_-^\nu=\beta^\mu\triangle\beta^\nu+\beta^\nu\triangle\beta^\mu 
-\triangle\beta^\mu\triangle\beta^\nu
\label{6.200}
\end{equation}
Let us denote:
\begin{equation}
\beta^\mu\zeta_\mu=n
\label{6.201}
\end{equation}
Then in view of the fact that by \ref{6.106}
\begin{equation}
\triangle\beta^\mu \zeta_\mu=\delta
\label{6.202}
\end{equation}
we obtain from \ref{6.98}:
\begin{equation}
\left(\beta^\mu\beta^\nu-\beta_-^\mu\beta_-^\nu\right)\zeta_\nu=\delta\beta^\mu 
+(n-\delta)\triangle\beta^\mu 
\label{6.203}
\end{equation}
Then (in view of \ref{6.201}) \ref{6.198} becomes:
\begin{equation}
J^\mu=2\kappa\left((n\triangle F+F_-\delta)\beta^\mu+F_-(n-\delta)\triangle\beta^\mu\right)
\label{6.204}
\end{equation}

Consider now the function $n$ which is defined on ${\cal K}$. Recalling \ref{6.119} and the fact that 
$u_\bot$ is positive, $n$ is a negative function. We first relate $n$ to the positive function 
$w$ defined by \ref{6.147}. Substituting $(h^{-1})^{\mu\nu}=(g^{-1})^{\mu\nu}-F\beta^\mu\beta^\nu$, 
the definition \ref{6.147} takes in terms of $n$ the form:
\begin{equation}
w=1-Fn^2
\label{6.205}
\end{equation}
Now, evaluating \ref{6.106} on $\beta^\sharp$ we obtain $\beta^\mu\triangle\beta_\mu=\nu$ on the left 
(recall the definition 1st of \ref{4.8}) and $\delta n$ on the right, hence:
\begin{equation}
n=\frac{\nu}{\delta}
\label{6.206}
\end{equation}
Substituting this in \ref{6.205} gives:
\begin{equation}
w=1-\frac{F\nu^2}{\delta^2}
\label{6.207}
\end{equation}
Recalling \ref{4.13} this is:
\begin{equation}
w=\frac{\mu}{F\nu^2+\mu}
\label{6.208}
\end{equation}
We can use this expression to determine $w$. Then $n$ is determined through $w$ by \ref{6.205} 
and the fact that it is negative as:
\begin{equation}
n=-\frac{1}{\sqrt{F}}\sqrt{1-w}
\label{6.209}
\end{equation}

We now decompose the vectorfield $J$, given by \ref{6.204}, into its parts $J_{||}$ tangential 
to ${\cal K}$, and $J_\bot$ $h$-orthogonal, that is $h_+$-orthogonal, to ${\cal K}$:
\begin{equation}
J=J_{||}+J_\bot
\label{6.210}
\end{equation}
Since $J_\bot$ is colinear to $M$ and $h(M,M)=4a$ (see \ref{6.145}), we have: 
\begin{equation}
J_\bot=J^M M, \ \mbox{where $J^M=(4a)^{-1}h(J,M)$}
\label{6.211}
\end{equation}
Also $J_{||}$ can be decomposed into a part colinear to $T$ and a part $\sJ$ which is an $S$ 
vectorfield defined on ${\cal K}$:
\begin{equation}
J_{||}= J^T T+\sJ \ \mbox{where $J^T=-(4a)^{-1}h(J,T)$}
\label{6.212}
\end{equation}
in view of the fact that $h(T,T)=-4a$. Moreover, we can write:
\begin{equation}
\sJ=J^A\Omega_A \ \mbox{where $J^A=(\sh^{-1})^{AB}h(J,\Omega_B)$}
\label{6.213}
\end{equation}
We shall presently analyze the behavior of the vectorfield $J$ for small $\tau$, that is in a 
neighborhood of $\partial_-{\cal B}$ in ${\cal K}$. This reduces 
to the analysis of the behavior of the coefficients $J^M$, $J^T$, and $J^A$, 
for small $\tau$. 

By \ref{6.146} and \ref{6.144},
\begin{equation}
h(J,M)=2\kappa \zeta(J)
\label{6.214}
\end{equation}
From \ref{6.204} we obtain, in view of \ref{6.201}, \ref{6.202},
\begin{equation}
\zeta(J)=2\kappa\left(n^2\triangle F+(2n\delta-\delta^2)F_-\right)
\label{6.215}
\end{equation}
In view of \ref{6.149}, the coefficient $J^M$ in \ref{6.211} is:
\begin{equation}
J^M=\frac{1}{w}\left(n^2\triangle F+(2n\delta-\delta^2)F_-\right)
\label{6.216}
\end{equation}
Writing $F_-=F-\triangle F$, this reads:
\begin{equation}
J^M=\frac{1}{w}\left((n-\delta)^2\triangle F+(2n\delta-\delta^2)F\right)
\label{6.217}
\end{equation}

We first analyze the behavior of $w$, expressed by \ref{6.208} for small $\tau$. By \ref{4.12} and 
Proposition 4.2:
\begin{equation}
\mu=\frac{j}{c}\epb^3
\label{6.218}
\end{equation}
Then in view of \ref{4.39} and \ref{6.188} we have:
\begin{equation}
\mu=-\frac{8l^3}{H_0^{\prime 2}}\left.\frac{1}{\eta^4\beta_N^6}\right|_{\partial_-{\cal B}} \tau^3 
+O(\tau^4)
\label{6.219}
\end{equation}
Substituting \ref{6.181} and \ref{6.219} in \ref{6.208} yields:
\begin{equation}
w=-\frac{2l}{F_0}\left.\frac{1}{\eta^4\beta_N^2}\right|_{\partial_-{\cal B}}\tau+O(\tau^2) \ \ 
\mbox{where $F_0=\left.F\right|_{\partial_-{\cal B}}$}
\label{6.220}
\end{equation}
This implies, through \ref{6.209},
\begin{equation}
n=-\frac{1}{\sqrt{F_0}}+O(\tau)
\label{6.221}
\end{equation}
which in turn implies through \ref{6.206}, in view of \ref{6.181},
\begin{equation}
\delta=\frac{\nu}{n}=\frac{2\sqrt{F_0}l}{H_0^\prime}
\left.\frac{1}{\beta_N^2}\right|_{\partial_-{\cal B}}\tau+O(\tau^2)
\label{6.222}
\end{equation}
Since 
\begin{equation}
\triangle F=\frac{dF}{d\sigma}\triangle\sigma+O((\triangle\sigma)^2)
\label{6.223}
\end{equation}
by \ref{4.14} and \ref{6.181} we obtain:
\begin{equation}
\triangle F=\frac{4l F^\prime_0}{H^\prime_0}\left.\frac{1}{\beta_N^2}\right|_{\partial_-{\cal B}}\tau 
+O(\tau^2) \ \ \mbox{where $F^\prime_0=\left.F^\prime\right|_{\partial_-{\cal B}}$}
\label{6.224}
\end{equation}
The expressions \ref{6.221}, \ref{6.222}, and \ref{6.224}, give, recalling \ref{4.27},
\begin{equation}
(n-\delta)^2\triangle F+(2n\delta-\delta^2)F
=\frac{4l}{F_0}\left.\frac{1}{\eta^4\beta_N^2}\right|_{\partial_-{\cal B}}\tau+O(\tau^2)
\label{6.225}
\end{equation}
In view of \ref{6.220} and \ref{6.225} we conclude through \ref{6.217} that, simply:
\begin{equation}
J^M=-2+O(\tau)
\label{6.226}
\end{equation}

Consider next $J^T$ (see \ref{6.212}). Writing $F_-=F-\triangle F$, \ref{6.204} reads:
\begin{equation}
J^\mu=2\kappa\left\{(n-\delta)(\beta^\mu-\triangle\beta^\mu)\triangle F
+\left(\delta\beta^\mu+(n-\delta)\triangle\beta^\mu\right)F\right\}
\label{6.227}
\end{equation}
We have:
\begin{equation}
h_{\mu\nu}\beta^\mu T^\nu=(1-\sigma H)\beta_\nu T^\nu=\eta^2\beta_T
\label{6.228}
\end{equation}
and:
\begin{equation}
h_{\mu\nu}\triangle\beta^\mu T^\nu=(g_{\mu\nu}+H\beta_\mu\beta_\nu)\triangle\beta^\mu T^\nu
=\triangle\beta_\nu T^\nu+H\nu\beta_T=H\nu\beta_T
\label{6.229}
\end{equation}
by virtue of the linear jump condition \ref{1.328} with $X=T$. We then obtain:
\begin{equation}
J^T=-\frac{\kappa}{2a}\beta_T\left\{(n-\delta)(\eta^2-H\nu)\triangle F
+\left(\delta\eta^2+(n-\delta)H\nu\right)F\right\}
\label{6.230}
\end{equation}
By \ref{6.221}, \ref{6.222}, \ref{6.224}, and \ref{6.181}, taking also into account \ref{4.27}, 
we deduce:
\begin{equation}
(n-\delta)(\eta^2-H\nu)\triangle F
+\left(\delta\eta^2+(n-\delta)H\nu\right)F=-\frac{4l}{\sqrt{F_0}}
\left.\frac{1}{\eta^2\beta_N^2}\right|_{\partial_-{\cal B}}\tau+O(\tau^2)
\label{6.231}
\end{equation}
Using \ref{6.149} we express the factor $\kappa/a$ on the right in \ref{6.230} as:
\begin{equation}
\frac{\kappa}{a}=\frac{1}{\sqrt{aw}}
\label{6.232}
\end{equation}
Then by \ref{6.220} and the 3rd of \ref{6.187} we obtain:
\begin{equation}
\frac{\kappa}{a}=\frac{3\sqrt{F_0}}{k}
\left.\eta^2\beta_N\right|_{\partial_-{\cal B}}\tau^{-2}+O(\tau^{-1})
\label{6.233}
\end{equation}
Taking also into account that in view of the fact $\beta_T=\rho\beta_N+\rhob\beta_{\Nb}$ the 
first two of \ref{6.187} imply
\begin{equation}
\beta_T=-\frac{k}{3l}\left.\beta_N\right|_{\partial_-{\cal B}}\tau+O(\tau^2)
\label{6.234}
\end{equation}
we then conclude through \ref{6.230} that, simply:
\begin{equation}
J^T=-2+O(\tau)
\label{6.235}
\end{equation}

Consider finally $J^A$ (see \ref{6.213}). We have:
\begin{equation}
h_{\mu\nu}\beta^\mu\Omega_B^\nu=(1-\sigma H)\beta_\nu\Omega_B^\nu=\eta^2\sbeta_B
\label{6.236}
\end{equation}
and:
\begin{equation}
h_{\mu\nu}\triangle\beta^\mu\Omega_B^\nu=(g_{\mu\nu}+H\beta_\mu\beta_\nu)\triangle\beta^\mu\Omega_B^\nu
=\triangle\beta_\nu\Omega_B^\nu+H\nu\sbeta_B=H\nu\sbeta_B
\label{6.237}
\end{equation}
by virtue of the linear jump condition \ref{1.328} with $X=\Omega_B$. We then obtain:
\begin{equation}
J^A=2\kappa(\sh^{-1})^{AB}\sbeta_B\left\{(n-\delta)(\eta^2-H\nu)\triangle F
+\left(\delta\eta^2+(n-\delta)H\nu\right)F\right\}
\label{6.238}
\end{equation}
Now, by \ref{6.233} and the 3rd of \ref{6.187}:
\begin{equation}
\kappa=-\frac{k\sqrt{F_0}}{6l}\left.\eta^2\beta_N\right|_{\partial_-{\cal B}}\tau+O(\tau^2)
\label{6.239}
\end{equation}
This, together with \ref{6.231} yields the conclusion that:
\begin{equation}
J^A=\frac{4k}{3}(\sh^{-1})^{AB}\left.\frac{\sbeta_B}{\beta_N}\right|_{\partial_-{\cal B}}\tau^2
+O(\tau^3)
\label{6.240}
\end{equation}
or, in terms of $m^A$ defined by \ref{4.186}, 
\begin{equation}
J^A=\frac{4km^A}{3}\tau^2+O(\tau^3)
\label{6.241}
\end{equation}

We turn to the function $\s^{(V)}b$, defined by \ref{6.165}, associated to each of the variation fields 
$V$. We shall analyze the dependence of this function on the derivatives of the transformation 
functions $f, w, \psi$, or equivalently $\hat{f}, v, \gamma$, introduced in Chapter 4 (see \ref{4.54}, 
\ref{4.56}, \ref{4.65}, \ref{4.222}, \ref{4.240}, and Proposition 4.4). By \ref{6.162} and 
\ref{6.158},
\begin{equation}
B_{-||}=B_{||}=B_{||}^T T+B_{||}^A\Omega_A
\label{6.242}
\end{equation}
where
\begin{equation}
B_{||}^T=\frac{1}{4}\hat{\Lambda}\beta_T=\frac{1}{2}+O(\tau)
\label{6.243}
\end{equation}
by \ref{6.191}, while 
\begin{equation}
B_{||}^A=-a\hat{\Lambda}(\sh^{-1})^{AB}\sbeta_B=-\frac{km^A}{3}\tau^2+O(\tau^3)
\label{6.244}
\end{equation}
by \ref{6.190} and the 3rd of \ref{6.187}.  By \ref{6.163} and \ref{6.199},
\begin{equation}
B_{-\bot}=-\frac{1}{2}\frac{G_-}{G}(M+J)
\label{6.245}
\end{equation}
We then have:
\begin{equation}
B_{-\bot}=B_{-\bot}^M M+B_{-\bot}^T T+B_{-\bot}^A\Omega_A
\label{6.246}
\end{equation}
where 
\begin{equation}
B_{-\bot}^M=-\frac{1}{2}\frac{G_-}{G}(1+J^M)=\frac{1}{2}+O(\tau)
\label{6.247}
\end{equation}
by \ref{6.226}, 
\begin{equation}
B_{-\bot}^T=-\frac{1}{2}\frac{G_-}{G}J^T=1+O(\tau)
\label{6.248}
\end{equation}
by \ref{6.235}, and 
\begin{equation}
B_{-\bot}^A=-\frac{1}{2}\frac{G_-}{G}J^A=-\frac{2km^A}{3}\tau^2+O(\tau^3)
\label{6.249}
\end{equation}
Defining along ${\cal K}$ the vectorfield
\begin{equation}
\hat{B}_-=\rho^{-1}B_-
\label{6.250}
\end{equation}
the rectangular components of $\hat{B}_-$ are given by:
\begin{equation}
\hat{B}_-^\mu=\hat{B}_-^N N^\mu+\hat{B}_-^{\Nb}\Nb^\mu+\hat{B}_-^A\Omega_A^\mu
\label{6.251}
\end{equation}
and in view of the above and also of the first two of \ref{6.187} we have: 
\begin{eqnarray}
&&\hat{B}_-^N=B_{||}^T+B_{-\bot}^T-B_{-\bot}^M=1+O(\tau)\label{6.252}\\
&&\hat{B}_-^{\Nb}=\frac{\rhob}{\rho}(B_{||}^T+B_{-\bot}^T+B_{-\bot}^M)=-\frac{l}{c_0}\tau+O(\tau^2)
\label{6.253}\\
&&\hat{B}_-^A=\rho^{-1}(B_{||}^A+B_{-\bot}^A)=3lm^A\tau+O(\tau^2)\label{6.254}
\end{eqnarray}

Now, if $X$ is any vectorfield along ${\cal K}$, we have, for any variation field $V$, denoting 
as before by $\beta^\prime_\mu$ the rectangular components of the 1-form $\beta$ corresponding to the 
prior solution,
\begin{equation}
\s^{(V)}\xi_-(X)=\left.V^\mu X\beta^\prime_\mu\right|_{{\cal K}}=
\left.X^\mu V\beta^\prime_\mu\right|_{{\cal K}} 
\label{6.255}
\end{equation}
by \ref{2.88}. Hence, from the definition \ref{6.165} and the expansion \ref{6.251},
\begin{equation}
\s^{(V)}b=\hat{B}_-^N\left.N^\mu\rho V\beta^\prime_\mu\right|_{{\cal K}}
+\hat{B}_-^{\Nb}\left.\Nb^\mu\rho V\beta^\prime_\mu\right|_{{\cal K}}
+\hat{B}_-^A\left.\Omega_A^\mu\rho V\beta^\prime_\mu\right|_{{\cal K}}
\label{6.256}
\end{equation}

Consider first the case $V=Y$. From \ref{6.43}, \ref{6.44}, \ref{6.50}, since according to 
\ref{4.5} along ${\cal K}$ we have $\rhob=r\rho$, it holds:
\begin{equation}
\rho Y=T \ \mbox{: along ${\cal K}$}
\label{6.257}
\end{equation}
Hence in the case $V=Y$ we have:
\begin{equation}
\s^{(Y)}b=\hat{B}_-^N\left.N^\mu T\beta^\prime_\mu\right|_{{\cal K}} 
+\hat{B}_-^{\Nb}\left.\Nb^\mu T\beta^\prime_\mu\right|_{{\cal K}}
+\hat{B}_-^A\left.\Omega_A^\mu T\beta^\prime_\mu\right|_{{\cal K}}
\label{6.258}
\end{equation}
Since the functions $\beta^\prime_\mu$ are known smooth functions of $(t,u^\prime,\vartheta^\prime)$ 
and along ${\cal K}$ we have (see \ref{4.54}, 
\ref{4.56}, \ref{4.65}, \ref{4.222}, \ref{4.240}, and Proposition 4.4):
\begin{eqnarray}
&&t=f(0,\vartheta)+\tau^2\hat{f}(\tau,\vartheta)\nonumber\\
&&u^\prime=\tau v(\tau,\vartheta)\nonumber\\
&&\vartheta^{\prime A}=\vartheta^A+\tau^3\gamma^A(\tau,\vartheta) \ : \ A=1,...,n-1
\label{6.259}
\end{eqnarray}
while, along ${\cal K}$, $T=\partial/\partial\tau$, we have:
\begin{eqnarray}
&&\left.T\beta^\prime_\mu\right|_{{\cal K}}=
\left.\frac{\partial\beta^\prime_\mu}{\partial t}\right|_{{\cal K}}
\left(\tau^2\frac{\partial\hat{f}}{\partial\tau}+2\tau\hat{f}\right)\nonumber\\
&&\hspace{15mm}+\left.\frac{\partial\beta^\prime_\mu}{\partial u^\prime}\right|_{{\cal K}}
\left(\tau\frac{\partial v}{\partial\tau}+v\right)\nonumber\\
&&\hspace{15mm}+\left.\frac{\partial\beta^\prime_\mu}{\partial\vartheta^{\prime A}}\right|_{{\cal K}}
\left(\tau^3\frac{\partial\gamma^A}{\partial\tau}+3\tau^2\gamma^A\right)
\label{6.260}
\end{eqnarray}
This expression is to be substituted in \ref{6.258} to obtain $\s^{(Y)}b$. The 2nd term on the right 
in \ref{6.260} depends most sensitively on the 1st derivatives of the transformation functions 
$\hat{f}, v, \gamma$ in a neighborhood of $\partial_-{\cal B}$ in ${\cal K}$, as it involves the 
factor $\tau$ rather than the factor $\tau^2$ or $\tau^3$ as the 1st or 3rd term on the right 
respectively. However, in assessing the contribution of this term to the 1st term on the right in 
\ref{6.258} we must take account of the fact that by Proposition 4.3, \ref{4.102} and the 1st of 
\ref{4.103} we have:
\begin{equation}
\left.N^\mu\frac{\partial\beta^\prime_\mu}{\partial u^\prime}\right|_{\partial_-{\cal B}}
=\left.L^{\prime\mu} T^\prime\beta^\prime_\mu\right|_{\partial_-{\cal B}}
=\left.T^{\prime\mu}L^\prime\beta^\prime_\mu\right|_{\partial_-{\cal B}}=0
\label{6.261}
\end{equation}
consequently:
\begin{equation}
\left.N^\mu\frac{\partial\beta^\prime_\mu}{\partial u^\prime}\right|_{{\cal K}}=O(\tau)
\label{6.262}
\end{equation}
Similarly, in assessing the contribution of the 2nd term on the right in \ref{6.260} to the 3rd term 
on the right in \ref{6.258} we must take into account the fact that by Proposition 4.3, \ref{4.102} 
and the 1st of \ref{4.103} we have:
\begin{equation}
\left.\Omega_A^\mu\frac{\partial\beta^\prime_\mu}{\partial u^\prime}\right|_{\partial_-{\cal B}}
=\left.(\Omega^{\prime\mu}_A+\pi_A L^{\prime\mu})\frac{\partial\beta^\prime_\mu}{\partial u^\prime}\right|_{\partial_-{\cal B}}
=\left.T^{\prime\mu}(\Omega^\prime_A\beta^\prime_\mu+\pi_A L^\prime\beta^\prime_\mu)\right|_{\partial_-{\cal B}}=0
\label{6.263}
\end{equation}
consequently:
\begin{equation}
\left.\Omega_A^\mu\frac{\partial\beta^\prime_\mu}{\partial u^\prime}\right|_{{\cal K}}=O(\tau)
\label{6.264}
\end{equation}
 On the other hand, by \ref{4.154}: 
\begin{equation}
\left.\Nb^\mu\frac{\partial\beta^\prime_\mu}{\partial u^\prime}\right|_{\partial_-{\cal B}}
=\left.s_{\Nb\Lb}\right|_{\partial_-{\cal B}}
\label{6.265}
\end{equation}
does not vanish. In view of the results \ref{6.252} for the coefficients in \ref{6.256}, 
the dominant contributions to \ref{6.258} are three. First, the contribution of the 1st term in \ref{6.260} to the 1st term in \ref{6.258}. Second, the contribution of the 2nd term in \ref{6.260} 
to the 1st term in \ref{6.258}. And third, the similar contribution of the 2nd term in \ref{6.260} 
to the 2nd term in \ref{6.258}. 

Consider next, in the case $n>2$, the variation fields $E_{(\kappa)}$, defined by \ref{6.53}, 
which are $S$ vectorfields.  Using \ref{6.54} we can express these as:
\begin{equation}
E_{(\kappa)}=E_{(\kappa)}^A\Omega_A \ \ \mbox{where $E^A_{(\kappa)}=(\sh^{-1})^{AB}\Omega_B^\nu h_{\kappa\nu}$}
\label{6.266}
\end{equation}
Then in the case $V=E_{(\kappa)}$ from \ref{6.256} we have:
\begin{equation}
\s^{(E_{(\kappa)})}b=(\sh^{-1})^{AB}h_{\kappa\nu}\Omega_B^\nu\rho(\hat{B}_-^N N^\mu
\left.\Omega_A\beta^\prime_\mu\right|_{{\cal K}}
+\hat{B}_-^{\Nb}\Nb^\mu\left.\Omega_A\beta^\prime_\mu\right|_{{\cal K}}
+\hat{B}_-^C\Omega_C^\mu\left.\Omega_A\beta^\prime_\mu\right|_{{\cal K}})
\label{6.267}
\end{equation}
In view of \ref{6.259} and the fact that $\Omega_A=\partial/\partial\vartheta^A$ we have, 
denoting 
\begin{equation}
f_0(\vartheta)=f(0,\vartheta),
\label{6.268}
\end{equation}
\begin{eqnarray}
&&\left.\Omega_A\beta^\prime_\mu\right|_{{\cal K}}=
\left.\frac{\partial\beta^\prime_\mu}{\partial t}\right|_{{\cal K}}
\left(\frac{\partial f_0}{\partial\vartheta^A}+\tau^2\frac{\partial\hat{f}}{\partial\vartheta^A}\right)
\nonumber\\
&&\hspace{15mm}+\left.\frac{\partial\beta^\prime_\mu}{\partial u^\prime}\right|_{{\cal K}}
\tau\frac{\partial v}{\partial\vartheta^A}\nonumber\\
&&\hspace{15mm}+\left.\frac{\partial\beta^\prime_\mu}{\partial\vartheta^{\prime B}}\right|_{{\cal K}}
\left(\delta^B_A+\tau^3\frac{\partial\gamma^B}{\partial\vartheta^A}\right)\label{6.269}
\end{eqnarray}
This expression is to be substituted in \ref{6.267} to obtain $\s^{(E_{(\kappa)})}b$. 
The 2nd term on the right in \ref{6.269} depends most sensitively on the 1st derivatives of the transformation functions 
$\hat{f}, v, \gamma$ in a neighborhood of $\partial_-{\cal B}$ in ${\cal K}$, as it involves the 
factor $\tau$ rather than the factor $\tau^2$ or $\tau^3$ as the 1st or 3rd term on the right 
respectively. This term contributes to the 1st, 3rd, and 2nd terms in the parenthesis 
on the right in \ref{6.267} with the coefficients \ref{6.262}, \ref{6.264}, and 
$$\left.\Nb^\mu\frac{\partial\beta^\prime_\mu}{\partial u^\prime}\right|_{{\cal K}}$$
respectively, the last being $O(1)$ as \ref{6.265} does not vanish. 
In view of the results \ref{6.252} for the coefficients in \ref{6.256}, 
the dominant contributions to \ref{6.267} are three. First, the contribution of the 1st term in 
\ref{6.269} to the 1st term in parenthesis in \ref{6.267}. Second, the contribution 
of the 2nd term in \ref{6.269} to the 1st term in parenthesis in \ref{6.267}. And third, the similar 
contribution of the 2nd term in \ref{6.269} to the 2nd term in parenthesis in \ref{6.267}.
Note that here we have the partial derivatives of the transformation functions 
with respect to $\vartheta$ and there is an extra factor of $\rho=O(\tau)$ (see 1st of \ref{6.187}) 
on the right in \ref{6.267} in comparison with \ref{6.258} where the partial derivatives 
of the transformation functions with respect to $\tau$ are involved. 

In the case $n=2$ the variation field which is a $S$ vectorfield is $E$, given by \ref{3.a10}.
In place of \ref{6.267} we have in this case the simpler formula:
\begin{equation}
\s^{(E)}b=\sh^{-1/2}\rho(\hat{B}_-^N N^\mu\left.\Omega\beta^\prime_\mu\right|_{{\cal K}}
+\hat{B}_-^{\Nb}\Nb^\mu\left.\Omega\beta^\prime_\mu\right|_{{\cal K}}
+\hat{B}_-^{E}E^\mu\left.\Omega\beta^\prime_\mu\right|_{{\cal K}})
\label{6.271}
\end{equation}
Otherwise everything is the same as above. 

We remark that the transformation functions $\hat{f}, v, \gamma$ being connected to the coordinates 
$t, u^\prime, \vartheta^\prime$ through \ref{6.259}, their first derivatives should be considered 
of order 0 following the remark at the end of Section 2.2. However the way that the transformation 
functions enter the energy estimates is through the boundary condition \ref{6.164} via the functions 
$\s^{(V)}b$. Since the 1-forms $\s^{(V)}\xi$ are of order 1, we shall consider the 1st derivatives 
of the transformation functions to be also of order 1. 

\pagebreak

\chapter{The Multiplier Field}

\section{Coercivity at the Boundary. 
The Choice of Multiplier Field}

Recall from Section 6.3 that the flux integrant is $\Omega^{(n-1)/2}$ times \ref{6.105}. Recall also 
from Section 6.4 that the boundary condition for the 1-form $\s^{(V)}\xi$ is \ref{6.164}, with the 
vectorfield $B$ along ${\cal K}$ being given by \ref{6.156} - \ref{6.158}. By Proposition 6.3 
$B_{||}$ is timelike future-directed relative to $h$. Thus, there is a timelike future-directed 
vectorfield $E_0$ of unit magnitude relative to $h$, which is colinear to $B_{||}$:
\begin{equation}
B_{||}=\mu \sqrt{a}E_0
\label{7.1}
\end{equation}
where $\mu$ is a positive function on ${\cal K}$. Let 
\begin{equation}
E_n=-\frac{M}{2\sqrt{a}}
\label{7.2}
\end{equation}
Then $E_n$ is the unit, relative to $h$, outward, relative to ${\cal N}$, normal to ${\cal K}$. 
By \ref{6.156}, \ref{6.157} we have:
\begin{equation}
B=\sqrt{a}(\mu E_0+E_n)
\label{7.3}
\end{equation}
Consider a given point $p$ on ${\cal K}$. Let $\Sigma_p$ be the $h$-orthogonal complement of $E_0$ in 
$T_p{\cal K}$, a $n-1$ dimensional spacelike plane. Let $(E_A:A=1,...,n-1)$ be a basis for 
$\Sigma_p$ which is orthonormal relative to $h$. Then we can express, at $p$,
\begin{equation}
h^{-1}=-E_0\otimes E_0+\sum_{i=1}^n E_i\otimes E_i
\label{7.4}
\end{equation}
Let us denote: 
\begin{equation}
\s^{(V)}\xi_0=\s^{(V)}\xi(E_0), \ \ \ \s^{(V)}\xi_i=\s^{(V)}\xi(E_i) \ : \ i=1,...,n
\label{7.5}
\end{equation}
Setting then 
\begin{equation}
\s^{(V)}\hat{b}=\frac{\s^{(V)}b}{\sqrt{a}}
\label{7.6}
\end{equation}
in view of \ref{7.3} the boundary condition \ref{6.164} reads:
\begin{equation}
\s^{(V)}\xi_n=-\mu\s^{(V)}\xi_0+\s^{(V)}\hat{b}
\label{7.7}
\end{equation}

Recall from Section 6.3 that the multiplier field $X$ is required to be timelike future-directed 
relative to $h$. Expanding 
\begin{equation}
X=X^0 E_0+\sum_{i=1}^n X^i E_i
\label{7.8}
\end{equation}
the requirement is: 
\begin{equation}
X^0>\sqrt{\sum_{i=1}^n (X^i)^2}
\label{7.9}
\end{equation}
Since $\s^{(V)}S(X,M)$ is linear in $X$ as well as $M$, coercivity of $\s^{(V)}S(X,M)$ reduces, 
in view of \ref{7.2} to coercivity of $-\s^{(V)}S(\hat{X},E_n)$, where 
\begin{equation}
\hat{X}=\frac{X}{X^0}=E_0+\sum_{i=1}^n\hat{X}^i E_i, \ \ \ \hat{X}^i=\frac{X^i}{X^0} \ : \ i=1,...,n
\label{7.10}
\end{equation}
and the requirement \ref{7.9} becomes: 
\begin{equation}
\sum_{i=1}^n(\hat{X}^i)^2<1
\label{7.11}
\end{equation}
Denoting 
\begin{equation}
\s^{(V)}S_{0n}=\s^{(V)}S(E_0,E_n), \ \ \ \s^{(V)}S_{in}=\s^{(V)}S(E_i,E_n) \ : \ i=1,...,n
\label{7.12}
\end{equation}
since by \ref{7.4}:
\begin{equation}
h^{-1}(\s^{(V)}\xi,\s^{(V)}\xi)=-(\s^{(V)}\xi_0)^2+\sum_{i=1}^n(\s^{(V)}\xi_i)^2
\label{7.13}
\end{equation}
we have, from \ref{6.102}:
\begin{eqnarray}
&&\s^{(V)}S_{0n}=\s^{(V)}\xi_0 \s^{(V)}\xi_n\nonumber\\
&&\s^{(V)}S_{An}=\s^{(V)}\xi_A \s^{(V)}\xi_n \ : \ A=1,...,n-1\nonumber\\
&&\s^{(V)}S_{nn}=\frac{1}{2}((\s^{(V)}\xi_0)^2+(\s^{(V)}\xi_n)^2)
-\frac{1}{2}\sum_{A=1}^{n-1}(\s^{(V)}\xi_A)^2
\label{7.14}
\end{eqnarray}
By \ref{7.10}, 
\begin{equation}
-\s^{(V)}S(\hat{X},E_n)=-\s^{(V)}S_{0n}-\sum_{A=1}^{n-1}\hat{X}^A \s^{(V)}S_{An}
-\hat{X}^n \s^{(V)}S_{nn}
\label{7.15}
\end{equation}
Now, the coercivity of $-\s^{(V)}S(\hat{X},E_n)$ modulo the boundary condition \ref{7.7} reduces to 
the positive-definiteness of the quadratic form in 
$(\s^{(V)}\xi_0,\s^{(V)}\xi_1,...,\s^{(V)}\xi_{n-1})$ obtained by substituting for $\s^{(V)}\xi_n$ 
in terms of $\s^{(V)}\xi_0$ in $-\s^{(V)}S(\hat{X},E_n)$ from the homogeneous version of \ref{7.7}, 
that is, by substituting 
$$\s^{(V)}\xi_n=-\mu \s^{(V)}\xi_0$$
Denoting the resulting quadratic form $(1/2)q$, we have:
\begin{eqnarray}
&&\mu(\s^{(V)}\xi_0)^2+\sum_{A=1}^{n-1}\mu\s^{(V)}\xi_0 \s^{(V)}\xi_A \hat{X}^A\label{7.16}\\
&&\hspace{15mm}-\frac{1}{2}\left((1+\mu^2)(\s^{(V)}\xi_0)^2-\sum_{A=1}^{n-1}(\s^{(V)}\xi_A)^2\right)\hat{X}^n
=\frac{1}{2}q\nonumber
\end{eqnarray}
Let us denote:
\begin{equation}
\hat{X}^A=x^A \ : \ A=1,...,n-1 \ \ \ \hat{X}^n=y
\label{7.17}
\end{equation}
Remark that $(x^1,...,x^{n-1},y)$ are the components of the vector 
$$\sum_{i=1}^n \hat{X}^i E_i$$ 
in the orthonormal basis $(E_1,...,E_n)$ in the $n$ dimensional spacelike hyperplane $H_p$ which is 
the $h$-orthogonal complement of $E_0$ in $T_p{\cal N}$. 
The quadratic form $q$ is then given by:
\begin{equation}
q=\alpha(\s^{(V)}\xi_0)^2+2\s^{(V)}\xi_0\sum_{A=1}^{n-1}\beta^A\s^{(V)}\xi_A
+\gamma\sum_{A=1}^{n-1}(\s^{(V)}\xi_A)^2
\label{7.18}
\end{equation}
where:
\begin{eqnarray}
&&\alpha=2\mu-(1+\mu^2)y\nonumber\\
&&\beta^A=\mu x^A \ : \ A=1,...,n-1\nonumber\\
&&\gamma=y
\label{7.19}
\end{eqnarray}
We shall invoke the following elementary algebraic lemma. 

\vspace{2.5mm}

\noindent{\bf Lemma 7.1} \ \ \ The quadratic form $q$ is positive-definite if and only if: 
$$\gamma>0 \ \ \mbox{and} \ \ |\beta|^2-\alpha\gamma<0$$

\vspace{2.5mm}

\noindent{\em Proof:} {\bf Sufficiency:} Setting 
$$\beta=\sum_{A=1}^{n-1}\beta^A E_A, \ \ \ \s^{(V)}\zeta=\sum_{A=1}^{n-1}\s^{(V)}\xi_A E_A$$
we have 
$$\sum_{A=1}^{n-1}\beta^A\s^{(V)}\xi_A=\beta\cdot\s^{(V)}\zeta$$ 
and 
$$q=\alpha(\s^{(V)}\xi_0)^2+2\s^{(V)}\xi_0 \ \beta\cdot\s^{(V)}\zeta+\gamma|\s^{(V)}\zeta|^2$$
In view of the fact that 
$$|\beta\cdot\s^{(V)}\zeta|\leq |\beta||\s^{(V)}\zeta|\leq |\beta||\s^{(V)}\zeta|$$
if the conditions of the lemma are satisfied, then:
\begin{eqnarray*}
&&q\geq \alpha |\s^{(V)}\xi_0|^2-2|\beta||\s^{(V)}\xi_0||\s^{(V)}\zeta|
+\gamma|\s^{(V)}\zeta|^2\\
&&\hspace{3mm}\geq \alpha |\s^{(V)}\xi_0|^2-2\sqrt{\alpha\gamma}|\s^{(V)}\xi_0||\s^{(V)}\zeta|
+\gamma|\s^{(V)}\zeta|^2\\
&&\hspace{3mm}=\left(\sqrt{\alpha}|\s^{(V)}\xi_0|-\sqrt{\gamma}|\s^{(V)}\zeta|\right)^2\geq 0
\end{eqnarray*}
and the first inequality is strict unless one of $|\s^{(V)}\xi_0|$, $|\s^{(V)}\zeta|$ 
vanishes. Therefore $q>0$ unless  one of $|\s^{(V)}\xi_0|$, $|\s^{(V)}\zeta|$ 
vanishes. But if $q=0$ then $\sqrt{\alpha}|\s^{(V)}\xi_0|-\sqrt{\gamma}|\s^{(V)}\zeta|=0$, 
and one of $|\s^{(V)}\xi_0|$, $|\s^{(V)}\zeta|$ vanishes. 
Thus in this case $|\s^{(V)}\xi_0|$, $|\s^{(V)}\zeta|$ 
both vanish. We conclude that $q$ is positive-definite. 

{\bf Necessity:} That the positivity of $\gamma$ and $\alpha$ is necessary, is obvious by setting 
$\s^{(V)}\xi_0=0$ and $\s^{(V)}\zeta=0$ respectively. Suppose then that $|\beta|^2-\alpha\gamma\geq 0$. 
Let us set:
$$\s^{(V)}\zeta=-\frac{\beta}{|\beta|}|\s^{(V)}\zeta|, \ \ \ \s^{(V)}\xi_0>0$$
Then we have:
\begin{eqnarray*}
&&q=\alpha|\s^{(V)}\xi_0|^2-2|\beta||\s^{(V)}\xi_0||\s^{(V)}\zeta|+\gamma|\s^{(V)}\zeta|^2\\
&&\hspace{5mm}=\alpha(|\s^{(V)}\xi_0|-r_-|\s^{(V)}\zeta|)(|\s^{(V)}\xi_0|-r_+|\s^{(V)}\zeta|)
\end{eqnarray*}
where
$$r_\pm=\frac{|\beta|\pm\sqrt{|\beta|^2-\alpha\gamma}}{\alpha}>0$$
Therefore $q<0$ for $r_-|\s^{(V)}\zeta|<|\s^{(V)}\xi_0|<r_+|\s^{(V)}\zeta|$. 

\vspace{2.5mm}

In the notation \ref{7.17}, the condition \ref{7.11} reads:
\begin{equation}
|x|^2+y^2<1
\label{7.20}
\end{equation}
with $x=(x^1,...,x^{n-1})\in\mathbb{R}^{n-1}$, $y\in\mathbb{R}$. So this condition defines the 
open unit ball in $\mathbb{R}^{n-1}\times\mathbb{R}=\mathbb{R}^n$, the last representing the 
hyperplane $H_p$. The coefficients of the quadratic form $q$ being given here by \ref{7.19} the 
necessary and sufficient conditions of Lemma 7.1 here are:
\begin{equation}
y>0 \ \ \mbox{and} \ \ F=\mu^2 |x|^2-[2\mu-(1+\mu^2)y]y<0
\label{7.21}
\end{equation}
the function $F$ being $|\beta|^2-\alpha\gamma$. The 2nd of these conditions takes the form:
\begin{equation}
(1+\mu^2)|x|^2+\frac{(1+\mu^2)^2}{\mu^2}\left(y-\frac{\mu}{1+\mu^2}\right)^2<1
\label{7.22}
\end{equation}
This defines the interior of an oblate spheroid in $\mathbb{R}^n$, centered at 
\begin{equation}
y=\frac{\mu}{1+\mu^2} \ \ \mbox{on the $y$ axis}
\label{7.23}
\end{equation}
with minor axis the $y$ axis, minor radius 
\begin{equation}
r_m=\frac{\mu}{1+\mu^2}
\label{7.24}
\end{equation}
and major radius the radius 
\begin{equation}
r_M=\frac{1}{\sqrt{1+\mu^2}}
\label{7.25}
\end{equation}
of the circle of intersection of the spheroid with the plane 
$$y=\frac{\mu}{1+\mu^2}$$
The spheroid is tangent to the plane $y=0$. As a consequence, the 1st of conditions \ref{7.21} is 
redundant. Let us take as the pole on the unit sphere $S^{n-1}$ in $\mathbb{R}^n$ the point of 
intersection with the positive $y$ axis. Letting $\vartheta$ being the corresponding polar 
angle, on $S^{n-1}$ we have:
\begin{equation}
y=\cos\vartheta, \ \ |x|=\sin\vartheta, \ \ \vartheta\in[0,\pi]
\label{7.26}
\end{equation}
and:
\begin{equation}
F=\mu^2\sin^2\vartheta-[2\mu-(1+\mu^2)\cos\vartheta]\cos\vartheta=(\mu-\cos\vartheta)^2
\label{7.27}
\end{equation}
Since $F\geq 0$ with equality if and only if $\cos\vartheta=\mu$, this means that the unit sphere 
comprehends the spheroid, so the interior of the spheroid is contained in the open unit ball. 
If $\mu>1$ the spheroid is contained in the open unit ball. If $\mu=1$ the spheroid is tangent to the 
unit sphere at the pole. If $\mu<1$ the spheroid is tangent to the unit sphere at the circle 
$\cos\vartheta=\mu$. Note from \ref{7.24}, \ref{7.25} that:
$$\frac{r_m}{r_M}=\frac{\mu}{\sqrt{1+\mu^2}}$$
This increases from 0 to 1 as $\mu$ increases from 0 to $\infty$. Thus, as $\mu$ increases, the 
spheroid becomes increasingly more spherical, becoming round in the limit $\mu\rightarrow\infty$. 
And, as $\mu$ decreases, the spheroid becomes increasingly more oblate, flattening out in the limit 
$\mu\rightarrow 0$. Finally, from \ref{7.23} the height $y$ of the center of the spheroid above 
the $y=0$ plane increases with $\mu$ from 0 to 1 as $\mu$ increases from 0 to 1, then decreases with 
$\mu$ from 1 to 0 as $\mu$ increases further from 1 to $\infty$. 

From \ref{7.1} $\mu$ is given by:
\begin{equation}
\mu=\sqrt{-a^{-1}h(B_{||},B_{||})}
\label{7.28}
\end{equation}
Substituting for $B_{||}$ from \ref{6.158} we then obtain:
\begin{equation}
\mu=\hat{\mu}\sqrt{1-\varepsilon^2}
\label{7.29}
\end{equation}
where 
\begin{equation}
\hat{\mu}=\frac{1}{2}\beta_T\hat{\Lambda}, \ \ \ 
\varepsilon=\frac{|\sbeta|}{\beta_{\hat{T}}}
\label{7.30}
\end{equation}
with 
\begin{equation}
\hat{T}=\frac{T}{2\sqrt{a}}
\label{7.31}
\end{equation}
being the future-directed timelike vectorfield colinear to $T$ which is of unit magnitude 
relative to $h$. From \ref{6.191} we have:
\begin{equation}
\hat{\mu}=1+O(\tau)
\label{7.32}
\end{equation}
Also,  
\begin{equation}
\beta_{\hat{T}}=\frac{1}{2\sqrt{c}}\left(\frac{1}{\sqrt{r}}\beta_N+\sqrt{r}\beta_{\Nb}\right)
=\frac{1}{\sqrt{-2l}}\left.\beta_N\right|_{\partial_-{\cal B}}\tau^{-1/2}+O(\tau^{1/2})
\label{7.33}
\end{equation}
by \ref{4.215}, hence:
\begin{equation}
\varepsilon=O(\tau^{1/2})
\label{7.34}
\end{equation}
We then conclude through \ref{7.29} that:
\begin{equation}
\mu=1+O(\tau)
\label{7.35}
\end{equation}

Lemma 7.1 leads at once to the following proposition (see \ref{7.6}, \ref{7.7}).

\vspace{2.5mm}

\noindent{\bf Proposition 7.1} \ \ \ Suppose that the point 
$(x,y)\in\mathbb{R}^{n-1}\times\mathbb{R}$ defined by the coefficients of the expansion \ref{7.10} 
of the vectorfield $\hat{X}$ according to \ref{7.17} lies in the interior of the spheroid \ref{7.22} 
and bounded away from the boundary. Then there are positive constants $C$ and $C^\prime$ such 
that on ${\cal K}$: 
$$-\s^{(V)}S(\hat{X},E_n)\geq C^{-1}\left((\s^{(V)}\xi_0)^2+\sum_{A=1}^{n-1}(\s^{(V)}\xi_A)^2\right)
-C^\prime(\s^{(V)}\hat{b})^2$$

\vspace{2.5mm}

Defining on ${\cal K}$ the vectorfield: 
\begin{equation}
U=\beta_{\hat{T}}\hat{T}-\sbeta^\sharp
\label{7.36}
\end{equation}
equation \ref{6.158} takes the form:
\begin{equation}
B_{||}=a\hat{\Lambda}U
\label{7.37}
\end{equation}
Let $V$ be the vectorfield on ${\cal K}$:
\begin{equation}
V=\frac{|\sbeta|^2}{\beta_{\hat{T}}}\hat{T}-\sbeta^\sharp
\label{7.38}
\end{equation}
This vectorfield is $h$-orthogonal to $U$, hence spacelike, and we have:
$$h(V,V)=|\sbeta|^2\left(1-\frac{|\sbeta|^2}{\beta_{\hat{T}}^2}\right)$$
that is:
\begin{equation}
|V|=|\sbeta|\sqrt{1-\varepsilon^2}
\label{7.39}
\end{equation}
(see \ref{7.30}). Subtracting \ref{7.38} from \ref{7.36} yields:
\begin{equation}
\hat{T}=\frac{U-V}{\beta_{\hat{T}}(1-\varepsilon^2)}
\label{7.40}
\end{equation}
Comparing \ref{7.37} with \ref{7.1} we see that:
\begin{equation}
U=\frac{\mu}{\hat{\Lambda}\sqrt{a}}E_0
\label{7.41}
\end{equation}
Substituting for $\mu$ from \ref{7.29}, noting that the 1st of \ref{7.30} can be written in the form:
\begin{equation}
\hat{\mu}=\beta_{\hat{T}}\hat{\Lambda}\sqrt{a}
\label{7.42}
\end{equation}
we obtain:
\begin{equation}
U=\beta_{\hat{T}}\sqrt{1-\varepsilon^2}E_0
\label{7.43}
\end{equation}
Now, we can choose the unit vector $E_1$ of the orthonormal basis $(E_A : A=1,...,n-1)$ for 
$\Sigma_p$ at any given point $p\in{\cal K}$ to be colinear and in the same sense as $V$. Then 
by \ref{7.39} we have:
\begin{equation}
V=\beta_{\hat{T}}\varepsilon\sqrt{1-\varepsilon^2}E_1
\label{7.44}
\end{equation}
Substituting \ref{7.43} and \ref{7.44} in \ref{7.40} yields:
\begin{equation}
\hat{T}=\frac{E_0-\varepsilon E_1}{\sqrt{1-\varepsilon^2}}
\label{7.45}
\end{equation}
Let also $\hat{M}$ be the spacelike vectorfield which is colinear and in the same sense as $M$ and 
of unit magnitude relative to $h$:
\begin{equation}
\hat{M}=\frac{M}{2\sqrt{a}}
\label{7.46}
\end{equation}
So, \ref{7.2} reads: 
\begin{equation}
\hat{M}=-E_n
\label{7.47}
\end{equation}

We now choose the multiplier field $X$ along ${\cal K}$ to be:
\begin{equation}
X=\alpha_0 L+\Lb \ \ \mbox{: along ${\cal K}$}
\label{7.48}
\end{equation}
where $\alpha_0$ is a positive constant. Then $X$ is along ${\cal K}$ a future-directed timelike 
vectorfield. Since (see \ref{7.31}, \ref{7.46}):
\begin{equation}
L=\frac{1}{2}(T-M)=\sqrt{a}(\hat{T}-\hat{M}), \ \ \ 
\Lb=\frac{1}{2}(T+M)=\sqrt{a}(\hat{T}+\hat{M})
\label{7.49}
\end{equation}
the definition \ref{7.48} reads:
\begin{equation}
X=\sqrt{a}[(\alpha_0+1)\hat{T}-(\alpha_0-1)\hat{M}]
\label{7.50}
\end{equation}
or, substituting for $\hat{T}$, $\hat{M}$ from \ref{7.45}, \ref{7.47}, 
\begin{equation}
X=\sqrt{a}\left((\alpha_0+1)\frac{(E_0-\varepsilon E_1)}{\sqrt{1-\varepsilon^2}}
+(\alpha_0-1)E_n\right)
\label{7.51}
\end{equation}
Consequently, for our choice \ref{7.48} we have:
\begin{equation}
X^0=\sqrt{a}\frac{(\alpha_0+1)}{\sqrt{1-\varepsilon^2}}
\label{7.52}
\end{equation}
and:
\begin{equation}
x^1=-\varepsilon, \ \ x^2=...=x^{n-1}=0; \ \ \ \ 
y=\left(\frac{\alpha_0-1}{\alpha_0+1}\right)\sqrt{1-\varepsilon^2}
\label{7.53}
\end{equation}

In view of \ref{7.34}, \ref{7.35}, Proposition 7.1 applies if we fix 
\begin{equation}
\alpha_0>1
\label{7.54}
\end{equation}
and restrict $\tau$ to $\tau\leq\tau_0$ for some suitably small positive $\tau_0$. 
Since by \ref{7.46}, \ref{7.47}, \ref{7.52},
\begin{equation}
\s^{(V)}S(X,M)=-2\frac{(\alpha_0+1)}{\sqrt{1-\varepsilon^2}}a\s^{(V)}S(\hat{X},E_n)
\label{7.55}
\end{equation}
Proposition 7.1 yields:
\begin{equation}
\s^{(V)}S(X,M)\geq aC^{-1}\left((\s^{(V)}\xi_0)^2+\sum_{A=1}^{n-1}(\s^{(V)}\xi_A)^2\right)
-C^\prime(\s^{(V)}b)^2
\label{7.56}
\end{equation}
(see \ref{7.6}) for positive constants $C$ and $C^\prime$. 

By \ref{7.45} we have:
\begin{equation}
\s^{(V)}\xi_{\hat{T}}=\frac{\s^{(V)}\xi_0-\varepsilon\s^{(V)}\xi_1}{\sqrt{1-\varepsilon^2}}
\label{7.57}
\end{equation}
Let $e$ be a vectorfield along ${\cal K}$ which is tangential to the $S_{\tau,\tau}$ sections of 
${\cal K}$ and of unit magnitude with respect to $h$. We can expand: 
\begin{equation}
e=e^0 E_0 + \sum_{A=1}^{n-1} e^A E_A
\label{7.58}
\end{equation}
We have $h(\hat{T},e)=0$. Substituting for $\hat{T}$ from \ref{7.45} this gives:
\begin{equation}
e^0=-\varepsilon e^1
\label{7.59}
\end{equation}
We also have $h(e,e)=1$, which reads:
\begin{equation}
-(e^0)^2+\sum_{A=1}^{n-1}(e^A)^2=1
\label{7.60}
\end{equation}
Substituting for $e^0$ from \ref{7.59} this becomes:
\begin{equation}
(1-\varepsilon^2)(e^1)^2+\sum_{A=2}^{n-1}(e^A)^2=1
\label{7.61}
\end{equation}
From \ref{7.58},
\begin{eqnarray}
&&\s^{(V)}\sxi(e)=e^0\s^{(V)}\xi_0+\sum_{A=1}^{n-1}e^A\s^{(V)}\xi_A\nonumber\\
&&\hspace{10mm}=e^1(\s^{(V)}\xi_1-\varepsilon \s^{(V)}\xi_0)+\sum_{A=2}^{n-1}e^A\s^{(V)}\xi_A
\label{7.62}
\end{eqnarray}
by \ref{7.59}. By virtue of \ref{7.61} it then follows that:
\begin{equation}
|\s^{(V)}\sxi(e)|\leq\sqrt{\frac{(\s^{(V)}\xi_1-\varepsilon\s^{(V)}\xi_0)^2}{1-\varepsilon^2}
+\sum_{A=2}^{n-1}(\s^{(V)}\xi_A)^2}
\label{7.63}
\end{equation}
Setting  
$$e=\frac{\s^{(V)}\sxi^\sharp}{|\s^{(V)}\sxi|}, \ \ \mbox{we have $\s^{(V)}\sxi(e)=|\s^{(V)}\sxi|$}.$$
We thus obtain:
\begin{equation}
|\s^{(V)}\sxi|\leq\sqrt{\frac{(\s^{(V)}\xi_1-\varepsilon\s^{(V)}\xi_0)^2}{1-\varepsilon^2}
+\sum_{A=2}^{n-1}(\s^{(V)}\xi_A)^2}
\label{7.64}
\end{equation}

In view of \ref{7.57}, \ref{7.64}, the inequality \ref{7.56} implies:
\begin{equation}
\s^{(V)}S(X,M)\geq C^{-1}\left((\s^{(V)}\xi_T)^2+a|\s^{(V)}\sxi|^2\right)-C^\prime(\s^{(V)}b)^2
\label{7.65}
\end{equation}
for a different positive constant $C$. This is the desired coercivity inequality on ${\cal K}$.  

From this point we set, in accordance with \ref{7.54}, 
\begin{equation}
\alpha_0=3
\label{7.66}
\end{equation}
Then by \ref{7.53}, in view of \ref{7.34}, \ref{7.35}, and \ref{7.23}, the point $(x,y)$ in 
$\mathbb{R}^n$ defined by $X$ tends to the center of the spheroid as $\tau\rightarrow 0$. 
Moreover, we take \ref{7.48} to hold on all of ${\cal N}$, that is we define the multiplier 
field $X$ by:
\begin{equation}
X=3L+\Lb \ \ \mbox{: on ${\cal N}$}
\label{7.67}
\end{equation}

\vspace{5mm}

\section{The Deformation Tensor of the Multiplier Field.
The Error Integral Associated to $\s^{(V)}Q_1$}

We now consider the deformation tensor $\s^{(X)}\pi$ of the vectorfield $X$ relative to $h$:
\begin{equation}
\s^{(X)}\pi={\cal L}_X h
\label{7.68}
\end{equation}
Since $\tilde{H}=\Omega h$, then $\s^{(X)}\tilde{\pi}$, the deformation tensor of $X$ relative 
to $\tilde{h}$ (see \ref{6.61}), is expressed through $\s^{(X)}\pi$ by:
\begin{equation}
\s^{(X)}\tilde{\pi}=\Omega\s^{(X)}\pi+(X\Omega)h
\label{7.69}
\end{equation}
To calculate $\s^{(X)}\pi$ we use the identity:
\begin{equation}
({\cal L}_X h)(Y,Z)=h(Y,D_Z X)+h(Z,D_Y X)
\label{7.70}
\end{equation}
which applies to any triplet $X,Y,Z$ of vectorfields on ${\cal N}$. By \ref{7.67} we have:
\begin{equation}
\s^{(X)}\pi=3\s^{(L)}\pi+\s^{(\Lb)}\pi
\label{7.71}
\end{equation}
Using the table \ref{3.23} we then obtain:
$$\s^{(L)}\pi_{LL}=0, \ \ \s^{(L)}\pi_{\Lb\Lb}=0, \ \ \s^{(L)}\pi_{L\Lb}=-2La$$
Also, defining for any vectorfield $Y$ the $S$ 1-forms $\s^{(Y)}\spi_L$, $\s^{(Y)}\spi_{\Lb}$ by:
\begin{equation}
\s^{(Y)}\spi_L(\Omega_A)=\s^{(Y)}\pi(L,\Omega_A), \ \ \ 
\s^{(Y)}\spi_{\Lb}(\Omega_A)=\s^{(Y)}\pi(\Lb,\Omega_A)
\label{7.72}
\end{equation}
we have:
$$\s^{(L)}\spi_L=0, \ \ \ \s^{(L)}\spi_{\Lb}=2\eta-2\etab$$
Here in view of \ref{3.16} the last is $\sh\cdot Z$. 
Also, defining for any vectorfield $Y$ the symmetric 2-covariant $S$ tensorfield $\s^{(Y)}\sspi$ by:
\begin{equation}
\s^{(Y)}\sspi(\Omega_A,\Omega_B)=\s^{(Y)}\pi(\Omega_A,\Omega_B)
\label{7.73}
\end{equation}
we have:
$$\s^{(L)}\sspi=2\chi$$
We combine the above results into the following table:
\begin{eqnarray}
&&\s^{(L)}\pi_{LL}=\s^{(L)}\pi_{\Lb\Lb}=0, \ \ \s^{(L)}\pi_{L\Lb}=-2La\nonumber\\
&&\s^{(L)}\spi_L=0, \ \ \ \s^{(L)}\spi_{\Lb}=\sh\cdot Z\nonumber\\
&&\s^{(L)}\sspi=2\chi
\label{7.74}
\end{eqnarray}
By conjugation we obtain the table:
\begin{eqnarray}
&&\s^{(\Lb)}\pi_{\Lb\Lb}=\s^{(\Lb)}\pi_{LL}=0, \ \ \s^{(\Lb)}\pi_{\Lb L}=-2\Lb a\nonumber\\
&&\s^{(\Lb)}\spi_{\Lb}=0, \ \ \ \s^{(\Lb)}\spi_L=-\sh\cdot Z\nonumber\\
&&\s^{(\Lb)}\sspi=2\chib
\label{7.75}
\end{eqnarray}
From tables \ref{7.74} and \ref{7.75} we obtain through \ref{7.71} the table:
\begin{eqnarray}
&&\s^{(X)}\pi_{LL}=\s^{(X)}\pi_{\Lb\Lb}=0, \ \ \s^{(X)}\pi_{L\Lb}=-2Xa\nonumber\\
&&\s^{(X)}\spi_L=-\sh\cdot Z, \ \ \ \s^{(X)}\spi_{\Lb}=3\sh\cdot Z\nonumber\\
&&\s^{(X)}\sspi=2(3\chi+\chib)
\label{7.76}
\end{eqnarray}
Moreover, by \ref{7.69}:
\begin{eqnarray}
&&\s^{(X)}\tilde{\pi}_{LL}=\s^{(X)}\tilde{\pi}_{\Lb\Lb}=0, \s^{(X)}\tilde{\pi}_{L\Lb}=-2X(\Omega a)
\nonumber\\
&&\s^{(X)}\tilde{\spi}_L=-\Omega\sh\cdot Z, \ \ \ \s^{(X)}\tilde{\spi}_{\Lb}=3\Omega\sh\cdot Z
\nonumber\\
&&\s^{(X)}\tilde{\sspi}=2\Omega(3\chi+\chib)+(X\Omega)\sh
\label{7.77}
\end{eqnarray}
the $S$ 1-forms $\s^{(Y)}\tilde{\spi}_L$, $\s^{(Y)}\tilde{\spi}_{\Lb}$, and the symmetric 2-covariant 
$S$ tensorfield $\s^{(Y)}\tilde{\sspi}$ being defined for any vectorfield $Y$ in analogy with 
\ref{7.72}, \ref{7.73} with $\s^{(Y)}\tilde{\pi}$  in the role of $\s^{(Y)}\pi$. 

We turn to the error term $\s^{(V)}Q_1$, given by \ref{6.64}, which arises from $\s^{(X)}\tilde{\pi}$. 
By \ref{6.58} and \ref{6.102}:
\begin{equation}
\s^{(V)}S^{\alpha\beta}=\Omega^{-2}(h^{-1})^{\alpha\gamma}(h^{-1})^{\beta\delta}\s^{(V)}S_{\gamma\delta}
\label{7.78}
\end{equation}
and:
\begin{equation}
\s^{(V)}S_{\alpha\beta}=\s^{(V)}\xi_\alpha\s^{(V)}\xi_\beta
 -\frac{1}{2}h_{\alpha\beta}(h^{-1})^{\gamma\delta}\s^{(V)}\xi_\gamma\s^{(V)}\xi_\delta
\label{7.79}
\end{equation}
Using then \ref{7.77} we deduce the following expression for the error term $\s^{(V)}Q_1$:
\begin{eqnarray}
&&2a\Omega\s^{(V)}Q_1=(a^{-1}Xa+\Omega^{-1}X\Omega) a|\s^{(V)}\sxi|^2\nonumber\\
&&\hspace{15mm}-\frac{1}{2}\s^{(V)}\xi_{\Lb} \s^{(V)}\sxi(Z)
+\frac{3}{2}\s^{(V)}\xi_L \s^{(V)}\sxi(Z)\nonumber\\
&&\hspace{15mm}-2(3\chi+\chib)\cdot a(\s^{(V)}\sxi^\sharp\otimes\s^{(V)}\sxi^\sharp 
-\frac{1}{2}|\s^{(V)}\sxi|^2\sh^{-1})\nonumber\\
&&\hspace{15mm}-(3\mbox{tr}\chi+\mbox{tr}\chib)\s^{(V)}\xi_L\s^{(V)}\xi_{\Lb}
\label{7.80}
\end{eqnarray}

Now, by \ref{6.91}, \ref{6.100}, \ref{6.103}, with $X$ defined by \ref{7.67}, 
the energy $\s^{(V)}{\cal E}^{\ub_1}(u)$ is: 
\begin{eqnarray}
&&\s^{(V)}{\cal E}^{\ub_1}(u)=\int_{C_u^{\ub_1}}\Omega^{(n-1)/2}\s^{(V)}S(X,L)\nonumber\\
&&\hspace{17mm}=\int_{C_u^{\ub_1}}\Omega^{(n-1)/2}(3\s^{(V)}S_{LL}+\s^{(V)}S_{\Lb L})\nonumber\\
&&\hspace{17mm}=\int_{C_u^{\ub_1}}\Omega^{(n-1)/2}(3(\s^{(V)}\xi_L)^2+a|\s^{(V)}\sxi|^2)
\label{7.81}
\end{eqnarray}
Also, by \ref{6.93}, \ref{6.100}, \ref{6.103}, with $X$ defined by \ref{7.67}, the energy 
$\s^{(V)}\underline{{\cal E}}^{u_1}(\ub)$ is:
\begin{eqnarray}
&&\s^{(V)}\underline{{\cal E}}^{u_1}(\ub)=\int_{\Cb_{\ub}^{u_1}}\Omega^{(n-1)/2}\s^{(V)}S(X,\Lb)
\nonumber\\
&&\hspace{17mm}=\int_{\Cb_{\ub}^{u_1}}\Omega^{(n-1)/2}(3\s^{(V)}S_{L\Lb}+\s^{(V)}S_{\Lb\Lb})\nonumber\\
&&\hspace{17mm}=\int_{\Cb_{\ub}^{u_1}}\Omega^{(n-1)/2}(3a|\s^{(V)}\sxi|^2+(\s^{(V)}\xi_{\Lb})^2)
\label{7.82}
\end{eqnarray}

Recall from Proposition 5.1 that $\lambdab$ vanishes on $\Cb_0$. Integrating then from 
$\Cb_0$ along the integral curves of $T$ (see \ref{2.34}) we obtain (acoustical coordinates):
\begin{equation}
\lambdab(\ub,u,\vartheta)=\int_0^{\ub}(T\lambdab)(\tau,u-\ub+\tau,\vartheta)d\tau
\label{7.83}
\end{equation}
Recall also that according to the 1st of \ref{4.216}:
\begin{equation}
(T\lambdab)(0,0,\vartheta)=\left(-\frac{c_0 k}{3l}\right)(\vartheta)
\label{7.84}
\end{equation}
It follows that along $\Cb_0$: 
\begin{equation}
(T\lambdab)(0,u,\vartheta)=\left(-\frac{c_0 k}{3l}\right)(\vartheta)+O(u)
\label{7.85}
\end{equation}
Integrating $T^2\lambdab$ from $\Cb_0$ along the integral curves of $T$ gives, 
for all $\tau\in[0,\ub]$:
\begin{equation}
(T\lambdab)(\tau,u-\ub+\tau,\vartheta)=(T\lambdab)(0,u-\ub,\vartheta)
+\int_0^\tau(T^2\lambdab)(\tau^\prime,u-\ub+\tau^\prime,\vartheta)d\tau^\prime
\label{7.86}
\end{equation}
Assuming a bound on $T^2\lambdab$ in ${\cal R}_{\ub_1,u_1}$, the integral here is $O(\tau)$. 
Substituting for the 1st term on the right from \ref{7.85} then yields, for all $\tau\in[0,\ub]$:
\begin{equation}
(T\lambdab)(\tau,u-\ub+\tau,\vartheta)=\left(-\frac{c_0 k}{3l}\right)(\vartheta)+O(u)
\label{7.87}
\end{equation}
Substituting this on the right in \ref{7.83} we conclude that:
\begin{equation}
\lambdab(\ub,u,\vartheta)=\left(-\frac{c_0 k}{3l}\right)(\vartheta)\ub+O(u\ub)
\label{7.88}
\end{equation}

Recall next that $L\lambda$ vanishes on $\Cb_0$ (see \ref{5.38}). Integrating then from 
$\Cb_0$ along the integral curves of $T$ we obtain:
\begin{equation}
(L\lambda)(\ub,u,\vartheta)=\int_0^{\ub}(TL\lambda)(\tau,u-\ub+\tau,\vartheta)d\tau
\label{7.89}
\end{equation}
Since 
$$TL\lambda=T^2\lambda-T\Lb\lambda=T^2\lambda-\Lb L\lambda+Z\lambda-\Lb^2\lambda$$
in view of \ref{5.35} and \ref{5.38} we have:
\begin{equation}
TL\lambda=T^2\lambda-\Lb^2\lambda \ \ \mbox{: on $\Cb_0$}
\label{7.90}
\end{equation}
Recall that according to \ref{4.107} and the 2nd of \ref{4.216}:
\begin{equation}
(T^2\lambda)(0,0,\vartheta)=\frac{k(\vartheta)}{3}, \ \ \ 
(\Lb^2\lambda)(0,0,\vartheta)=\frac{k(\vartheta)}{2}
\label{7.91}
\end{equation}
It follows that along $\Cb_0$:
\begin{equation}
(TL\lambda)(0,u,\vartheta)=-\frac{k(\vartheta)}{6}+O(u)
\label{7.92}
\end{equation}
Integrating $T^2 L\lambda$ from $\Cb_0$ along the integral curves of $T$ gives, 
for all $\tau\in[0,\ub]$:
\begin{equation}
(TL\lambda)(\tau,u-\ub+\tau,\vartheta)=(TL\lambda)(0,u-\ub,\vartheta)
+\int_0^\tau(T^2 L\lambda)(\tau^\prime,u-\ub+\tau^\prime,\vartheta)d\tau^\prime
\label{7.93}
\end{equation}
Assuming a bound on $T^2 L\lambda$ in ${\cal R}_{\ub_1,u_1}$, the integral here is $O(\tau)$. 
Substituting for the 1st term on the right from \ref{7.92} then yields, for all $\tau\in[0,\ub]$:
\begin{equation}
(TL\lambda)(\tau,u-\ub+\tau,\vartheta)=-\frac{k(\vartheta)}{6}+O(u)
\label{7.94}
\end{equation}
Substituting this on the right in \ref{7.89} we conclude that:
\begin{equation}
(L\lambda)(\ub,u,\vartheta)=-\frac{k(\vartheta)}{6}\ub+O(u\ub)
\label{7.95}
\end{equation}

Integrating from $\Cb_0$ along the integral curves of $T$ we also obtain:
\begin{equation}
(\Lb\lambda)(\ub,u,\vartheta)=(\Lb\lambda)(0,u-\ub,\vartheta)
+\int_0^{\ub}(T\Lb\lambda)(\tau,u-\ub+\tau,\vartheta)d\tau
\label{7.96}
\end{equation}
By \ref{4.106}, \ref{4.107}, along $\Cb_0$ we have:
\begin{equation}
(\Lb\lambda)(0,u,\vartheta)=\frac{k(\vartheta)}{2}u+O(u^2)
\label{7.97}
\end{equation}
Since 
$$T\Lb\lambda=\Lb^2\lambda+L\Lb\lambda=\Lb^2\lambda+\Lb L\lambda-Z\lambda$$
in view of \ref{5.35} and \ref{5.38} we have:
\begin{equation}
T\Lb\lambda=\Lb^2\lambda \ \ \mbox{: on $\Cb_0$}
\label{7.98}
\end{equation}
It follows in view of \ref{4.107} that along $\Cb_0$:
\begin{equation}
(T\Lb\lambda)(0,u,\vartheta)=\frac{k(\vartheta)}{2}+O(u)
\label{7.99}
\end{equation}
Integrating $T^2\Lb\lambda$ from $\Cb_0$ along the integral curves of $T$ gives, 
for all $\tau\in[0,\ub]$:
\begin{equation}
(T\Lb\lambda)(\tau,u-\ub+\tau,\vartheta)=(T\Lb\lambda)(0,u-\ub,\vartheta)
+\int_0^\tau(T^2\Lb\lambda)(\tau^\prime,u-\ub+\tau^\prime,\vartheta)d\tau^\prime
\label{7.100}
\end{equation}
Assuming a bound on $T^2 \Lb\lambda$ in ${\cal R}_{\ub_1,u_1}$, the integral here is $O(\tau)$. 
Substituting for the 1st term on the right from \ref{7.99} then yields, for all $\tau\in[0,\ub]$:
\begin{equation}
(T\Lb\lambda)(\tau,u-\ub+\tau,\vartheta)=\frac{k(\vartheta)}{2}+O(u)
\label{7.101}
\end{equation}
Substituting this as well as \ref{7.97} on the right in \ref{7.96} we conclude that: 
\begin{equation}
(\Lb\lambda)(\ub,u,\vartheta)=\frac{k(\vartheta)}{2}u+O(u^2)
\label{7.102}
\end{equation}

Adding \ref{7.95} and \ref{7.102} we obtain:
\begin{equation}
(T\lambda)(\ub,u,\vartheta)=\frac{k(\vartheta)}{6}(3u-\ub)+O(u^2)
\label{7.103}
\end{equation}
Integrating from $\Cb_0$ along the integral curves of $T$ gives:
\begin{equation}
\lambda(\ub,u,\vartheta)=\lambda(0,u-\ub,\vartheta)
+\int_0^{\ub} (T\lambda)(\tau,u-\ub+\tau,\vartheta)d\tau
\label{7.104}
\end{equation}
Substituting \ref{7.103} and taking into account the fact that by \ref{4.106}, \ref{4.107}, we have,  
along $\Cb_0$:
$$\lambda(0,u,\vartheta)=\frac{k(\vartheta)}{4}u^2+O(u^3)$$
we conclude that:
\begin{equation}
\lambda(\ub,u,\vartheta)=\frac{k(\vartheta)}{12}(3u^2-\ub^2)+O(u^3)
\label{7.105}
\end{equation}

We collect the above results into the following proposition. 

\vspace{2.5mm}

\noindent{\bf Proposition 7.2} \ \ \ Assuming a bound on $T^2\lambdab$ in ${\cal R}_{\ub_1,u_1}$ 
we have: 
$$\lambdab(\ub,u,\vartheta)=\left(-\frac{c_0 k}{3l}\right)(\vartheta)\ub+O(u\ub) \ \ \mbox{: in ${\cal R}_{\ub_1,u_1}$}$$
Assuming a bound on $T^2 L\lambda$ in ${\cal R}_{\ub_1,u_1}$ we have: 
$$(L\lambda)(\ub,u,\vartheta)=-\frac{k(\vartheta)}{6}\ub+O(u\ub) \ \ 
\mbox{: in ${\cal R}_{\ub_1,u_1}$}$$
Assuming a bound on $T^2\Lb\lambda$ in ${\cal R}_{\ub_1,u_1}$ we have: 
$$(\Lb\lambda)(\ub,u,\vartheta)=\frac{k(\vartheta)}{2}u+O(u^2) \ \ \mbox{: in ${\cal R}_{\ub_1,u_1}$}$$
Moreover, under both of the last two assumptions we have: 
$$\lambda(\ub,u,\vartheta)=\frac{k(\vartheta)}{12}(3u^2-\ub^2)+O(u^3)
 \ \ \mbox{: in ${\cal R}_{\ub_1,u_1}$}$$
 
\vspace{2.5mm}

Let us introduce the following notation. Let $g$ be a positive function and $f$ another function 
on the same domain. We write $f\sim g$ if there are positive constants $C_m$ and $C_M$  such that: 
$C_m g\leq f\leq C_M g$. By Proposition 7.2 we then have:
\begin{equation}
\lambdab\sim \ub, \ \ \ \lambda\sim u^2
\label{7.106}
\end{equation}
Now, we may assume that: 
\begin{equation}
c\sim 1, \ \ \ \Omega\sim 1
\label{7.107}
\end{equation}
these being satisfied by the restrictions of these functions to $\Cb_0$. The above imply that:
\begin{equation}
\rho\sim \ub, \ \ \ \rhob\sim u^2, \ \ \ a\sim \ub u^2
\label{7.108}
\end{equation}

Since the $L\beta_\mu$ vanish on $\Cb_0$ (first corollary of Proposition 5.1), assuming a bound 
on $TL\beta_\mu$ in ${\cal R}_{\ub_1,u_1}$ we obtain:
\begin{equation}
L\beta_\mu=O(\ub), \ \ \mbox{hence} \ \ LH=O(\ub)
\label{7.109}
\end{equation}
Then in view of the 1st of \ref{7.106}, $\underline{q}/\lambdab$ is bounded, where the function 
$\underline{q}$ is defined in the statement of Proposition 3.3. Then from Proposition 3.3: 
\begin{equation}
\frac{\Lb\lambdab}{\lambdab} \ \ \mbox{is bounded}
\label{7.110}
\end{equation}

We now consider the error integral $\s^{(V)}{\cal G}_1^{\ub_1,u_1}$ corresponding to $\s^{(V)}Q_1$ 
(see \ref{6.98}):
\begin{equation}
\s^{(V)}{\cal G}_1^{\ub_1,u_1}=\int_{{\cal R}_{\ub_1,u_1}}2a\Omega^{(n+1)/2}\s^{(V)}Q_1
\label{7.111}
\end{equation}
We shall consider the contribution to this error integral of each of the five terms in the expression 
\ref{7.80}, starting from the last. The contribution of this 5th term is: 
\begin{equation}
-\int_{{\cal R}_{\ub_1,u_1}}(3\mbox{tr}\chi+\mbox{tr}\chib)\Omega^{(n-1/2}s^{(V)}\xi_L\s^{V)}\xi_{\Lb}
\label{7.112}
\end{equation}
Assuming bounds on $\mbox{tr}\chi$, $\mbox{tr}\chib$, this is bounded in absolute value by:
\begin{eqnarray*}
&&C\int_{{\cal R}_{\ub_1,u_1}}\sqrt{3}\Omega^{(n-1)/2}|\s^{(V)}\xi_L||\s^{(V)}\xi_{\Lb}|\\
&&\leq C\left(\int_{{\cal R}_{\ub_1,u_1}}3\Omega^{(n-1)/2}(\s^{(V)}\xi_L)^2\right)^{1/2} 
\left(\int_{{\cal R}_{\ub_1,u_1}}\Omega^{(n-1)/2}(\s^{(V)}\xi_{\Lb})^2\right)^{1/2}
\end{eqnarray*}
By \ref{7.81}, \ref{7.82}, 
\begin{eqnarray}
&&\int_{{\cal R}_{\ub_1,u_1}}3\Omega^{(n-1)/2}(\s^{(V)}\xi_L)^2
=\int_0^{u_1}\left\{\int_{C_u^{\ub_1}}3\Omega^{(n-1)/2}(\s^{(V)}\xi_L)^2\right\}du\nonumber\\
&&\hspace{25mm}\leq\int_0^{u_1}{\cal E}^{\ub_1}(u)du
\label{7.113}
\end{eqnarray}
\begin{eqnarray}
&&\int_{{\cal R}_{\ub_1,u_1}}\Omega^{(n-1)/2}(\s^{(V)}\xi_{\Lb})^2
=\int_0^{\ub_1}\left\{\int_{\Cb_{\ub}^{u_1}}\Omega^{(n-1)/2}(\s^{(V)}\xi_{\Lb})^2\right\}d\ub
\nonumber\\
&&\hspace{25mm}\leq\int_0^{\ub_1}\underline{{\cal E}}^{u_1}(\ub)d\ub
\label{7.114}
\end{eqnarray}
Substituting above we conclude that \ref{7.112} is bounded in absolute value by:
\begin{equation}
C\left(\int_0^{u_1}{\cal E}^{\ub_1}(u)du\right)^{1/2}
\left(\int_0^{\ub_1}\underline{{\cal E}}^{u_1}(\ub)d\ub\right)^{1/2}
\label{7.115}
\end{equation}

Assuming bounds on $|\chi|$ and $|\chib|$ the contribution of the 4th term in \ref{7.80} to the 
error integral \ref{7.111} is bounded in absolute value by:
\begin{equation}
C\int_{{\cal R}_{\ub_1,u_1}}3\Omega^{(n-1)/2}a|\s^{(V)}\sxi|^2
\leq C\int_0^{\ub_1}\underline{{\cal E}}^{u_1}(\ub)d\ub
\label{7.116}
\end{equation}
by \ref{7.82}. 

We turn to the contributions of the 2nd and 3rd terms in \ref{7.80}. In view of the fact that 
$Z$ vanishes on $\Cb_0$ (see \ref{5.35}), assuming a bound on $|\sL_T Z|$ in ${\cal R}_{\ub_1,u_1}$ 
we deduce:
\begin{equation}
|Z|=O(\ub) \ \ \mbox{in ${\cal R}_{\ub_1,u_1}$}
\label{7.117}
\end{equation}
Then in view of the 3rd of \ref{7.108} the contribution of the 3rd term in \ref{7.80} to the error 
integral \ref{7.111} is bounded in absolute value by:
\begin{eqnarray}
&&C\int_{{\cal R}_{\ub_1,u_1}}\ub\Omega^{(n-1)/2}|\s^{(V)}\xi_L||\s^{(V)}\sxi|\nonumber\\
&&\hspace{25mm}\leq C\int_{{\cal R}_{\ub_1,u_1}}2\sqrt{3}\Omega^{(n-1)/2}\frac{\sqrt{\ub}}{u}|\s^{(V)}\xi_L| \sqrt{a}|\s^{(V)}\sxi|
\nonumber\\
&&\hspace{25mm}\leq C\int_{{\cal R}_{\ub_1,u_1}}\frac{2\sqrt{3}}{\sqrt{u}}\Omega^{(n-1)/2}
|\s^{(V)}\xi_L| \sqrt{a}|\s^{(V)}\sxi| \ \ \mbox{: since $\ub\leq u$}\nonumber\\
&&\hspace{25mm}=\int_0^{u_1}\left\{\int_{C_u^{\ub_1}}2\sqrt{3}\Omega^{(n-1)/2}|\s^{(V)}\xi_L| 
\sqrt{a}|\s^{(V)}\sxi|\right\}\frac{du}{\sqrt{u}}
\nonumber\\
&&\hspace{35mm}\leq\int_0^{u_1}{\cal E}^{\ub_1}(u)\frac{du}{\sqrt{u}}
\label{7.118}
\end{eqnarray}
by \ref{7.81}. By \ref{7.117} and the 3rd of \ref{7.108} the contribution of the 2nd term in \ref{7.80} 
to the error integral \ref{7.111} is bounded in absolute value by:
\begin{eqnarray}
&&C\int_{{\cal R}_{\ub_1,u_1}}\ub\Omega^{(n-1)/2}|\s^{(V)}\xi_{\Lb}||\s^{(V)}\sxi|\nonumber\\
&&\hspace{25mm}\leq C\int_{{\cal R}_{\ub_1,u_1}}2\sqrt{3}\Omega^{(n-1)/2}\frac{\sqrt{\ub}}{u}
|\s^{(V)}\xi_{\Lb}| \sqrt{a}|\s^{(V)}\sxi|
\nonumber\\
&&\hspace{25mm}\leq C\int_{{\cal R}_{\ub_1,u_1}}\frac{2\sqrt{3}}{\sqrt{\ub}}\Omega^{(n-1)/2}
|\s^{(V)}\xi_{\Lb}| \sqrt{a}|\s^{(V)}\sxi| \ \ \mbox{: since $\ub\leq u$}\nonumber\\
&&\hspace{25mm}=\int_0^{\ub_1}\left\{\int_{\Cb_{\ub}^{u_1}}2\sqrt{3}\Omega^{(n-1)/2}|\s^{(V)}\xi_{\Lb}| 
\sqrt{a}|\s^{(V)}\sxi|\right\}\frac{d\ub}{\sqrt{\ub}}
\nonumber\\
&&\hspace{35mm}\leq\int_0^{\ub_1}\underline{{\cal E}}^{u_1}(\ub)\frac{d\ub}{\sqrt{\ub}}
\label{7.119}
\end{eqnarray}
by \ref{7.82}. 

Finally we consider the contribution of the 1st term in \ref{7.80} to the error integral \ref{7.111}. 
Consider first the coefficient:
\begin{equation}
a^{-1}Xa+\Omega^{-1}X\Omega=\lambdab^{-1}X\lambdab+\lambda^{-1}X\lambda-c^{-1}Xc+\Omega^{-1}X\Omega
\label{7.120}
\end{equation}
In view of \ref{7.107}, assuming bounds on $Xc$, $X\Omega$, the last two terms in \ref{7.120} are 
bounded. Therefore the partial contribution of these terms through the 1st term in \ref{7.80} is 
bounded as the contribution of the 4th term by \ref{7.116}. Since 
$X\lambda=3L\lambda+\Lb\lambda$ and by Proposition 7.2 $L\lambda$ is nonpositive for suitably small 
$u$, we only have to consider the partial contribution of 
\begin{equation}
\lambda^{-1}\Lb\lambda\sim u^{-1}
\label{7.121}
\end{equation}
(see Proposition 7.2). The partial contribution of this term, which is nonnegative, 
is then bounded by:
\begin{eqnarray}
&&C\int_{{\cal R}_{\ub_1,u_1}}u^{-1}\Omega^{(n-1)/2}a|\s^{(V)}\sxi|^2
=C\int_0^{u_1}\left\{\int_{C_u^{\ub_1}}\Omega^{(n-1)/2}a|\s^{(V)}\sxi|^2\right\}\frac{du}{u}
\nonumber\\
&&\hspace{35mm}\leq C\int_0^{u_1}{\cal E}^{\ub_1}(u)\frac{du}{u}
\label{7.122}
\end{eqnarray}
by \ref{7.81}.

Since $X\lambdab=3L\lambdab+\Lb\lambdab$, in view of \ref{7.110} we only have to consider the 
partial contribution of $3\lambdab^{-1}L\lambdab$, which since 
\begin{equation}
\lambdab^{-1}L\lambdab\sim \ub^{-1}
\label{7.123}
\end{equation}
(see \ref{7.87}, \ref{7.88}) is nonnegative and bounded by:
\begin{eqnarray}
&&C\int_{{\cal R}_{\ub_1,u_1}}3\ub^{-1}\Omega^{(n-1)/2}a|\s^{(V)}\sxi|^2
=C\int_0^{\ub_1}\left\{\int_{\Cb_{\ub}^{u_1}}3\Omega^{(n-1)/2}a|\s^{(V)}\sxi|^2\right\}\frac{d\ub}{\ub}
\nonumber\\
&&\hspace{35mm}\leq C\int_0^{\ub_1}\underline{{\cal E}}^{u_1}(\ub)\frac{d\ub}{\ub}
\label{7.124}
\end{eqnarray}
by \ref{7.82}. 

The integrals on the right in \ref{7.122} and \ref{7.124} are {\em borderline error integrals}. 
In Chapter 9 the method shall be introduced by which such integrals can be handled . 

\vspace{5mm}

\section{The Error Integral Associated to $\s^{(V)}Q_2$}

We turn to the error integral $\s^{(V)}{\cal G}_2^{\ub_1,u_1}$ corresponding to $\s^{(V)}Q_2$ 
(see \ref{6.98}):
\begin{equation}
\s^{(V)}{\cal G}_2^{\ub_1,u_1}=\int_{{\cal R}_{\ub_1,u_1}}2a\Omega^{(n+1)/2}\s^{(V)}Q_2
\label{7.125}
\end{equation}
where $\s^{(V)}Q_2$ is given by \ref{6.65}. Here we shall place the fundamental variations \ref{6.22} 
in the role of the $\dot{\phi}_\mu$. For these fundamental variations we have shown at the 
end of Section 6.2 that control of all components of the associated 1-form $\s^{(V)}\xi$ 
for all variation fields $V$ provides control on all components of the symmetric 
2-covariant tensorfield $s$. A similar statement will hold for the higher order variations defined 
in Chapter 8. The analogous error integral $\s^{(V)}{\cal G}_2^{\ub_1,u_1}$ corresponding to these 
higher order variations will then be handled in the same manner. 

We begin by quantifying in terms of the energies 
$\s^{(V)}{\cal E}$, $\s^{(V)}\underline{{\cal E}}$ the control on the components of $s$ discussed 
at the end of Section 6.2. Consider first the case $n=2$. By \ref{6.a1}, \ref{6.a2}, \ref{6.a3}, 
and \ref{7.81}, \ref{7.82}:
\begin{eqnarray}
&&\int_{C_u^{\ub_1}}\Omega^{(n-1)/2}\left\{3(\rho\ss_N)^2+a\sss^2\right\}=\s^{(E)}{\cal E}^{\ub_1}(u)\label{7.126}\\
&&\int_{\Cb_{\ub}^{u_1}}\Omega^{(n-1)/2}\left\{3a\sss^2+(\rhob\ss_{\Nb})^2\right\}=\s^{(E)}
\underline{{\cal E}}^{u_1}(\ub)\label{7.127}
\end{eqnarray}
By \ref{6.a7}, equations \ref{6.a5}, \ref{6.a6} take the form:
\begin{eqnarray}
&&s_{NL}+c\rho\ogamma\sss=\s^{(Y)}\xi_L \label{7.128}\\
&&\ogamma s_{\Nb\Lb}+c\rhob\sss=\s^{(Y)}\xi_{\Lb}\label{7.129}
\end{eqnarray}
Since, by \ref{7.81}, \ref{7.82}, 
\begin{equation}
\int_{C_u^{\ub_1}}3\Omega^{(n-1)/2}(\s^{(Y)}\xi_L)^2\leq\s^{(Y)}{\cal E}^{\ub_1}(u), \ \ \ 
\int_{\Cb_{\ub}^{u_1}}\Omega^{(n-1)/2}(\s^{(Y)}\xi_{\Lb})^2\leq\s^{(Y)}\underline{{\cal E}}^{u_1}(\ub)
\label{7.130}
\end{equation}
and $\sss$ is controlled by \ref{7.126}, \ref{7.127}, we acquire control on $s_{NL}$ 
and $s_{\Nb\Lb}$ through \ref{7.128}, \ref{7.129}. 

We note that since 
\begin{equation}
r\sim\tau \ \mbox{: on ${\cal K}$}
\label{7.131}
\end{equation}
\ref{6.44}, \ref{6.50}, and \ref{6.51} imply:
\begin{equation}
\ogamma\sim u\sim\sqrt{\rhob} \ \mbox{: on ${\cal N}$}
\label{7.132}
\end{equation}

For $n>2$ we shall make use of the following lemma. 

\vspace{2.5mm}

\noindent{\bf Lemma 7.2} \ \ \ {\bf a)} Let $\zeta$ be an $S$ 1-form. We then have:
$$\sum_{\mu}(\zeta(E_{(\mu)})^2\geq |\zeta|^2$$
{\bf b)} Let $\zeta$ be a 2-covariant $S$ tensorfield. We then have:
$$\sum_\mu|E_{(\mu)}\cdot\zeta|^2\geq |\zeta|^2$$

\vspace{2.5mm}

\noindent{\em Proof:} {\bf a)} Let us denote by $\zeta_\nu$ the rectangular components of $\zeta$. 
Then by \ref{6.53} we have:
\begin{equation}
\zeta(E_{(\mu)})=\zeta_\nu E_{(\mu)}^\nu=\zeta_\nu\left(\delta^\nu_\mu+\frac{1}{2c}(N^\nu\Nb^\lambda
+\Nb^\nu N^\lambda)h_{\lambda\mu}\right)=\zeta_\mu 
\label{7.133}
\end{equation}
since $\zeta_\nu N^\nu=\zeta_\nu\Nb^\nu=0$. 
Hence:
\begin{equation}
\sum_{\mu}(\zeta(E_{(\mu)})^2=\sum_\mu (\zeta_\mu)^2
\label{7.134}
\end{equation}
On the other hand, since 
\begin{equation}
(\sh^{-1})^{\mu\nu}=(h^{-1})^{\mu\nu}+\frac{1}{2c}(N^\mu\Nb^\nu+\Nb^\mu N^\nu)
\label{7.135}
\end{equation}
for the same reason,
\begin{equation}
|\zeta|^2:=\sh^{-1}(\zeta,\zeta)=
(\sh^{-1})^{\mu\nu}\zeta_\mu\zeta_\nu=(h^{-1})^{\mu\nu}\zeta_\mu\zeta_\nu
\label{7.136}
\end{equation}
But, denoting $\zeta_\mu\beta^\mu=\kappa$, we have:
\begin{equation}
(h^{-1})^{\mu\nu}\zeta_\mu\zeta_\nu=(g^{-1})^{\mu\nu}\zeta_\mu\zeta_\nu-F\kappa^2
\leq (g^{-1})^{\mu\nu}\zeta_\mu\zeta_\nu:=-(\zeta_0)^2+\sum_i(\zeta_i)^2\leq\sum_\mu(\zeta_\mu)^2
\label{7.137}
\end{equation}

{\bf b)} Let us denote by $\zeta_{\lambda\nu}$ the rectangular components of $\zeta$. 
Then by \ref{6.53} the rectangular components of $E_{(\mu)}\cdot\zeta$ are given by:
\begin{equation}
(E_{(\mu)}\cdot\zeta)_\nu=E^\lambda_{(\mu)}\zeta_{\lambda\nu}=\zeta_{\mu\nu}
\label{7.138}
\end{equation}
since $N^\lambda\zeta_{\lambda\nu}=\Nb^\lambda\zeta_{\lambda\nu}=0$. Hence:
\begin{equation}
\sum_{\mu}|E_{(\mu)}\cdot\zeta|^2=\sum_{\mu}(\sh^{-1})^{\nu\lambda}
(E_{(\mu)}\cdot\zeta)_{\nu}(E_{(\mu)}\cdot\zeta)_{\lambda}
=\sum_{\mu}(\sh^{-1})^{\nu\lambda}\zeta_{\mu\nu}\zeta_{\mu\lambda}
\label{7.139}
\end{equation}
On the other hand, denoting $\kappa_\nu=\beta^\mu\zeta_{\mu\nu}$, we have:
\begin{eqnarray}
&&|\zeta|^2=(\sh^{-1})^{\mu\kappa}(\sh^{-1})^{\nu\lambda}\zeta_{\mu\nu}\zeta_{\kappa\lambda}
=(h^{-1})^{\mu\kappa}(\sh^{-1})^{\nu\lambda}\zeta_{\mu\nu}\zeta_{\kappa\lambda}\nonumber\\
&&\hspace{15mm}=(g^{-1})^{\mu\kappa}(\sh^{-1})^{\nu\lambda}\zeta_{\mu\nu}\zeta_{\kappa\lambda}
-F(\sh^{-1})^{\nu\lambda}\kappa_\nu\kappa_\lambda\nonumber\\
&&\hspace{30mm}\leq\sum_{\mu}(\sh^{-1})^{\nu\lambda}\zeta_{\mu\nu}\zeta_{\mu\lambda}
\label{7.140}
\end{eqnarray}

\vspace{5mm}

By virtue of the above lemma and \ref{6.a8}, \ref{6.a9}, \ref{6.a10} we have:
\begin{eqnarray}
&&\sum_{\mu}|\s^{(E_{(\mu)})}\sxi|^2\geq |\sss|^2 \label{7.141}\\
&&\sum_{\mu}(\s^{(E_{(\mu)})}\xi(L))^2\geq \rho^2|\ss_N|^2 \label{7.142}\\
&&\sum_{\mu}(\s^{(E_{(\mu)})}\xi(\Lb))^2\geq \rhob^2|\ss_{\Nb}|^2 \label{7.143}
\end{eqnarray}
hence, comparing with \ref{7.81}, \ref{7.82}, 
\begin{eqnarray}
&&\int_{C_u^{\ub_1}}\Omega^{(n-1)/2}\left\{\rho^2 |\ss_N|^2+a|\sss|^2\right\}\leq 
\sum_{\mu}\s^{(E_{(\mu)})}{\cal E}^{\ub_1}(u) \label{7.144}\\
&&\int_{\Cb_{\ub}^{u_1}}\Omega^{(n-1)/2}\left\{\rhob^2 |\ss_{\Nb}|^2+a|\sss|^2\right\}\leq 
\sum_{\mu}\s^{(E_{(\mu)})}\underline{{\cal E}}^{u_1}(\ub) \label{7.145}
\end{eqnarray}
By \ref{6.a14}, equations \ref{6.a12}, \ref{6.a13} take the form:
\begin{eqnarray}
&&s_{NL}+c\rho\ogamma\mbox{tr}\sss=\s^{(Y)}\xi(L) \label{7.146}\\
&&\ogamma s_{\Nb\Lb}+c\rhob\mbox{tr}\sss=\s^{(Y)}\xi(\Lb) \label{7.147}
\end{eqnarray}
Since, by \ref{7.81}, \ref{7.82}, 
\begin{equation}
\int_{C_u^{\ub_1}}3\Omega^{(n-1)/2}(\s^{(Y)}\xi_L)^2\leq\s^{(Y)}{\cal E}^{\ub_1}(u), \ \ \ 
\int_{\Cb_{\ub}^{u_1}}\Omega^{(n-1)/2}(\s^{(Y)}\xi_{\Lb})^2\leq\s^{(Y)}\underline{{\cal E}}^{u_1}(\ub)
\label{7.148}
\end{equation}
and, in view of the elementary algebraic inequality
\begin{equation}
|\sss|^2\geq\frac{1}{n-1}(\mbox{tr}\sss)^2
\label{7.149}
\end{equation}
$\mbox{tr}\sss$ is controlled by \ref{7.144}, \ref{7.145}, we acquire control on $s_{NL}$, 
$s_{\Nb\Lb}$ through \ref{7.146}, \ref{7.147}.

The expression \ref{6.65} with the fundamental variations \ref{6.22} in the 
role of the $\dot{\phi}_\mu$, reads, recalling that in \ref{6.65} the indices 
from the beginning of the Greek alphabet are raised with respect to $\tilde{h}$, 
\begin{eqnarray}
&&\Omega\s^{(V)}Q_2=-(h^{-1})^{\alpha\gamma}\s^{(V)}\xi_\gamma X^\beta
\left(\s^{(V)}\theta^\mu_\alpha\partial_\beta\beta_\mu-
\s^{(V)}\theta^\mu_\beta\partial_\alpha\beta_\mu\right)\nonumber\\
&&\hspace{18mm}-\s^{(V)}\xi_\beta X^\beta(h^{-1})^{\alpha\gamma}
\s^{(V)}\theta^\mu_\alpha\partial_\gamma\beta_\mu 
\label{7.150}
\end{eqnarray}
Substituting the expansion \ref{2.38} of $h^{-1}$ and denoting, in accordance with the definition 
\ref{3.36} of $s$, 
\begin{equation}
L\beta_\mu=s_{\mu L}, \ \ \ \Lb\beta_\mu= s_{\mu\Lb}, \ \ \ \sd\beta_\mu=\ss_\mu 
\label{7.151}
\end{equation}
we obtain:
\begin{eqnarray}
&&a\Omega\s^{(V)}Q_2=3\s^{(V)}\xi_L\left(\s^{(V)}\theta^\mu_{\Lb} s_{\mu L}
-a\sh^{-1}(\s^{(V)}\stheta^\mu,\ss_\mu)\right) \label{7.152}\\
&&\hspace{18mm}+3a\left(-\sh^{-1}(\s^{(V)}\sxi,\s^{(V)}\stheta^\mu)s_{\mu L}
+\sh^{-1}(\s^{(V)}\sxi,\ss_\mu)\s^{(V)}\theta^\mu_L\right)\nonumber\\
&&\hspace{18mm}+\s^{(V)}\xi_{\Lb}\left(\s^{(V)}\theta^\mu_L s_{\mu\Lb}
-a\sh^{-1}(\s^{(V)}\stheta^\mu,\ss_\mu)\right)\nonumber\\
&&\hspace{18mm}+a\left(-\sh^{-1}(\s^{(V)}\sxi,\s^{(V)}\stheta^\mu)s_{\mu\Lb}
+\sh^{-1}(\s^{(V)}\sxi,\ss_\mu)\s^{(V)}\theta^\mu_{\Lb}\right)\nonumber
\end{eqnarray}
In the above we can expand, as in \ref{6.45}, 
\begin{eqnarray}
&&\s^{(V)}\theta^\mu_L=\s^{(V)}\stheta_L\cdot\sd x^\mu+\s^{(V)}\theta^{N}_L N^\mu 
+\s^{(V)}\theta^{\Nb}_L \Nb^\mu \nonumber\\
&&\s^{(V)}\theta^\mu_{\Lb}=\s^{(V)}\stheta_{\Lb}\cdot\sd x^\mu+\s^{(V)}\theta^N_{\Lb} N^\mu 
+\s^{(V)}\theta^{\Nb}_{\Lb} \Nb^\mu \nonumber\\
&&\s^{(V)}\stheta^\mu=\s^{(V)}\sstheta\cdot\sd x^\mu+\s^{(V)}\stheta^N N^\mu 
+\s^{(V)}\stheta^{\Nb} \Nb^\mu \label{7.153}
\end{eqnarray}
Here $\s^{(V)}\stheta_L$, $\s^{(V)}\stheta_{\Lb}$ are $S$ vectorfields, 
$\s^{(V)}\stheta^N$, $\s^{(V)}\stheta^{\Nb}$ are $S$ 1-forms, and $\s^{(V)}\sstheta$ is a 
$T^1_1$-type $S$ tensorfield. In terms of the above expansions we can express, in view of 
\ref{6.a14}, 
\begin{eqnarray}
&&\s^{(V)}\theta^\mu_L s_{\mu \Lb}=c\rhob\s^{(V)}\theta^N_L \mbox{tr}\sss
+\s^{(V)}\theta^{\Nb}_L s_{\Nb\Lb}+\rhob \s^{(V)}\stheta_L\cdot\ss_{\Nb}\nonumber\\
&&\s^{(V)}\theta^\mu_{\Lb} s_{\mu L}=\s^{(V)}\theta^N_{\Lb} s_{NL}
+c\rho\s^{(V)}\theta^{\Nb}_{\Lb}\mbox{tr}\sss+\rho \s^{(V)}\stheta_{\Lb}\cdot\ss_N 
\label{7.154}
\end{eqnarray}
and:
\begin{eqnarray}
&&\s^{(V)}\stheta^\mu s_{\mu L}=\s^{(V)}\stheta^N s_{NL}+c\rho\s^{(V)}\stheta^{\Nb}\mbox{tr}\sss
+\rho\s^{(V)}\sstheta\cdot\ss_N\nonumber\\
&&\s^{(V)}\stheta^\mu s_{\mu\Lb}=c\rhob\s^{(V)}\stheta^N\mbox{tr}\sss
+\s^{(V)}\stheta^{\Nb}s_{\Nb\Lb}+\rhob\s^{(V)}\sstheta\cdot\ss_{\Nb}
\label{7.155}
\end{eqnarray}
Moreover, we can express:
\begin{eqnarray}
&&\sh^{-1}(\s^{(V)}\stheta^\mu,\ss_\mu)=\sh^{-1}(\s^{(V)}\stheta^N,\ss_N)
+\sh^{-1}(\s^{(V)}\stheta^{\Nb},\ss_{\Nb})+\mbox{tr}(\s^{(V)}\sstheta\cdot\sss)\nonumber\\
&&\s^{(V)}\theta^\mu_L\ss_\mu=\s^{(V)}\stheta_L\cdot\sss+\s^{(V)}\theta^N_L\ss_N+
\s^{(V)}\theta^{\Nb}_L\ss_{\Nb}\nonumber\\
&&\s^{(V)}\theta^\mu_{\Lb}\ss_\mu=\s^{(V)}\stheta_{\Lb}\cdot\sss+\s^{(V)}\theta^N_{\Lb}\ss_N
+\s^{(V)}\theta^{\Nb}_{\Lb}\ss_{\Nb}
\label{7.156}
\end{eqnarray}

In the following we shall estimate the error integral \ref{7.125} in terms of the energies, 
using only the assumption that for each of the variation fields $V$: 
\begin{eqnarray}
&&|\s^{(V)}\sstheta|, \ |\s^{(V)}\stheta^N|, \ |\s^{(V)}\stheta^{\Nb}| \ \ : \ \mbox{are bounded}
\nonumber\\
&&|\s^{(V)}\stheta_{\Lb}|, \ \s^{(V)}\theta^N_{\Lb}, \ \s^{(V)}\theta^{\Nb}_{\Lb} \ \ : \ 
\mbox{are bounded}\nonumber\\
&&|\s^{(V)}\stheta_L|, \s^{(V)}\theta^N_L, \ \s^{(V)}\theta^{\Nb}_L \ \ : \ 
\mbox{are $O(\ub)$}
\label{7.157}
\end{eqnarray}
For each of the variation fields $V$ these assumptions follow from the fact that the rectangular 
components $V^\mu$ are smooth functions of the acoustical coordinates and that the $LV^\mu$ vanish 
on $\Cb_0$ by virtue of the vanishing corollaries of Proposition 5.1. See \ref{6.45} for $V=Y$, 
\ref{6.53} for $V=E_{(\mu)}$ (for $n>2$), \ref{6.54} for $V=E$ (case $n=2$). We shall show the 
estimates for $n>2$, the case $n=2$ being similar but simpler. 

By \ref{7.108} and the 3rd of \ref{7.157} we have, from the 1st of \ref{7.154},
\begin{equation}
|\s^{(V)}\theta^\mu_L s_{\mu\Lb}|\leq C(a|\sss|+\rho|s_{\Nb\Lb}|+a|\ss_{\Nb}|)
\label{7.158}
\end{equation}
By the 2nd of \ref{7.157} we have, from the 2nd of \ref{7.154},
\begin{equation}
|\s^{(V)}\theta^\mu_{\Lb} s_{\mu L}|\leq C(|s_{NL}|+\rho|\sss|+\rho|\ss_N|)
\label{7.159}
\end{equation}
By the 1st of \ref{7.157} we have, from \ref{7.155}, 
\begin{eqnarray}
&&|\s^{(V)}\stheta^\mu s_{\mu L}|\leq C(|s_{NL}|+\rho|\sss|+\rho|\ss_N|) \label{7.160}\\
&&|\s^{(V)}\stheta^\mu s_{\mu\Lb}|\leq C(\rhob|\sss|+|s_{\Nb\Lb}|+\rhob|\ss_{\Nb}|)\label{7.161}
\end{eqnarray}
Moreover, by the 1st of \ref{7.157} we have, from the 1st of \ref{7.156}, 
\begin{equation}
|\sh^{-1}(\s^{(V)}\stheta^\mu,\ss_\mu)|\leq C(|\ss_N|+|\ss_{\Nb}|+|\sss|)
\label{7.162}
\end{equation}
By \ref{7.108} and the 3rd of \ref{7.157} we have, from the 2nd of \ref{7.156}, 
\begin{equation}
|\s^{(V)}\theta^\mu_L\ss_\mu|\leq C\rho(|\sss|+|\ss_N|+|\ss_{\Nb}|)
\label{7.163}
\end{equation}
and by the 2nd of \ref{7.157} we have, from the 3rd of \ref{7.156}, 
\begin{equation}
|\s^{(V)}\theta^\mu_{\Lb}\ss_\mu|\leq C(|\sss|+|\ss_N|+|\ss_{\Nb}|)
\label{7.164}
\end{equation}

By \ref{7.159} and \ref{7.162}, the contribution of the 1st term in \ref{7.152} to the error 
integral \ref{7.125} is bounded by:
\begin{equation}
C\int_{{\cal R}_{\ub_1,u_1}}\Omega^{(n-1)/2}|\s^{(V)}\xi_L|(|s_{NL}|+\rho|\sss|+\rho|\ss_N|
+a|\ss_{\Nb}|)
\label{7.165}
\end{equation}
By \ref{7.160} and \ref{7.163} the contribution of the 2nd term in \ref{7.152} to the error 
integral \ref{7.125} is bounded by:
\begin{equation}
C\int_{{\cal R}_{\ub_1,u_1}}\Omega^{(n-1)/2}a|\s^{(V)}\sxi|(|s_{NL}|+\rho|\sss|+\rho|\ss_N|
+\rho|\ss_{\Nb}|)
\label{7.166}
\end{equation}
By \ref{7.158} and \ref{7.162}, the contribution of the 3rd term in \ref{7.152} to the error 
integral \ref{7.125} is bounded by:
\begin{equation}
C\int_{{\cal R}_{\ub_1,u_1}}\Omega^{(n-1)/2}|\s^{(V)}\xi_{\Lb}|(\rho|s_{\Nb\Lb}|+a|\sss|+a|\ss_N|
+a|\ss_{\Nb}|)
\label{7.167}
\end{equation}
Finally, by \ref{7.161} and \ref{7.164} the contribution of the 4th term in \ref{7.152} to the error 
integral \ref{7.125} is bounded by:
\begin{equation}
C\int_{{\cal R}_{\ub_1,u_1}}\Omega^{(n-1)/2}a|\s^{(V)}\sxi|(|s_{\Nb\Lb}|+|\sss|+|\ss_N|
+|\ss_{\Nb}|)
\label{7.168}
\end{equation}

To estimate \ref{7.165} and \ref{7.166} we first note that by \ref{7.146} and the fact that by 
\ref{7.108} and \ref{7.132} $\rho\ogamma\sim\rho\sqrt{\rhob}\sim\sqrt{\rho}\sqrt{a}$ and 
$\sqrt{\rho}\sim\sqrt{\ub}$ while $\ub\leq u$, we have:
\begin{equation}
|| s_{NL}||_{L^2(C_u^{\ub_1})}\leq ||\s^{(Y)}\xi_L||_{L^2(C_u^{\ub_1})}
+C\sqrt{u}||\sqrt{a}\sss||_{L^2(C_u^{\ub_1})}
\label{7.169}
\end{equation}
Hence, by \ref{7.144} and the 1st of \ref{7.148}:
\begin{equation}
\int_{C_u^{\ub_1}}\Omega^{(n-1)/2}s_{NL}^2\leq C\left(\s^{(Y)}{\cal E}^{\ub_1}(u)
+u\sum_\mu\s^{(E_{(\mu)})}{\cal E}^{\ub_1}(u)\right)
\label{7.170}
\end{equation}

In regard then to the partial contribution of the 1st term in parenthesis in \ref{7.165}, by \ref{7.170} we have:
\begin{eqnarray}
&&\int_{{\cal R}_{\ub_1,u_1}}\Omega^{(n-1)/2}|\s^{(V)}\xi_L||s_{NL}|
\leq \label{7.171}\\
&&\hspace{15mm}\int_0^{u_1}\left(\int_{C_u^{\ub_1}}\Omega^{(n-1)/2}(\s^{(V)}\xi_L)^2\right)^{1/2}
\left(\int_{C_u^{\ub_1}}\Omega^{(n-1)/2}s_{NL}^2\right)^{1/2}du\nonumber\\
&&\hspace{15mm}\leq C\int_0^{u_1}\left(\s^{(V)}{\cal E}^{\ub_1}(u)\right)^{1/2}
\left(\s^{(Y)}{\cal E}^{\ub_1}(u)+u\sum_\mu\s^{(E_{(\mu)})}{\cal E}^{\ub_1}(u)\right)^{1/2}du
\nonumber
\end{eqnarray}
To estimate the partial contribution of the 2nd term in parenthesis in \ref{7.165} we note that 
by \ref{7.108} 
$$\rho\sim\sqrt{\frac{\rho}{\rhob}}\sqrt{a}\sim\frac{\sqrt{\ub}}{u}\sqrt{a}$$
and $\ub\leq u$, hence:
\begin{eqnarray}
&&\int_{{\cal R}_{\ub_1,u_1}}\Omega^{(n-1)/2}|\s^{(V)}\xi_L|\rho|\sss|\leq \label{7.172}\\
&&\hspace{10mm} C\int_0^{u_1}\left(\int_{C_u^{\ub_1}}\Omega^{(n-1)/2}(\s^{(V)}\xi_L)^2\right)^{1/2}
\left(\int_{C_u^{\ub_1}}\Omega^{(n-1)/2}a|\sss|^2\right)^{1/2}\frac{du}{\sqrt{u}}\nonumber\\
&&\hspace{15mm}\leq C\int_0^{u_1}\left(\s^{(Y)}{\cal E}^{\ub_1}(u)\right)^{1/2}
\left(\sum_{\mu}\s^{(E_{(\mu)})}{\cal E}^{\ub_1}(u)\right)^{1/2}\frac{du}{\sqrt{u}}\nonumber
\end{eqnarray}
by \ref{7.144}. The partial contribution of the 3rd term in parenthesis in \ref{7.165} is 
estimated in a straightforward manner using \ref{7.144}:
\begin{equation}
\int_{{\cal R}_{\ub_1,u_1}}\Omega^{(n-1)/2}|\s^{(V)}\xi_L|\rho|\ss_N|
\leq  C\int_0^{u_1}\left(\s^{(Y)}{\cal E}^{\ub_1}(u)\right)^{1/2}
\left(\sum_{\mu}\s^{(E_{(\mu)})}{\cal E}^{\ub_1}(u)\right)^{1/2}du
\label{7.173}
\end{equation}
Finally, the partial contribution of the 4th term in parenthesis in \ref{7.165} is estimated 
using \ref{7.145}:
\begin{eqnarray}
&&\int_{{\cal R}_{\ub_1,u_1}}\Omega^{(n-1)/2}|\s^{(V)}\xi_L|a|\ss_{\Nb}|\leq \label{7.174}\\
&&\hspace{10mm}\left(\int_{{\cal R}_{\ub_1,u_1}}\Omega^{(n-1)/2}(\s^{(V)}\xi_L)^2\right)^{1/2}
\left(\int_{{\cal R}_{\ub_1,u_1}}\Omega^{(n-1)/2}a^2|\ss_{\Nb}|^2\right)^{1/2}\nonumber\\
&&\hspace{10mm}\leq C\left(\int_0^{u_1}\s^{(V)}{\cal E}^{\ub_1}(u)du\right)^{1/2}
\left(\int_0^{\ub_1}\ub^2\sum_{\mu}\s^{(E_{(\mu)})}\underline{{\cal E}}^{u_1}(\ub)d\ub\right)^{1/2}
\nonumber
\end{eqnarray}
since $a\sim\rho\rhob$ and $\rho\sim\ub$. 

In regard to \ref{7.166} we remark from \ref{7.81} that the $L^2(C_u^{\ub_1})$ bounds for 
$\sqrt{a}\s^{(V)}\sxi$ and $\s^{(V)}\xi_L$ in terms of $\sqrt{{\cal E}^{\ub_1}(u)}$ are similar. 
Thus the partial contributions of the first three terms in parenthesis in \ref{7.166} are 
estimated in the same way as the partial contributions of the corresponding terms in \ref{7.165}. 
In fact they are depressed relative to the former by the small factor $\sqrt{a}$. As for the partial 
contribution of the 4th term in \ref{7.166}, we estimate using \ref{7.82} and \ref{7.145} 
as follows:
\begin{eqnarray}
&&\int_{{\cal R}_{\ub_1,u_1}}\Omega^{(n-1)/2}a|\s^{(V)}\xi|\rho|\ss_{\Nb}|\leq \label{7.175}\\
&&\hspace{10mm}\int_0^{\ub_1}\left(\int_{\Cb_{\ub}^{u_1}}a|\s^{(V)}\sxi|^2\right)^{1/2}
\left(\int_{\Cb_{\ub}^{u_1}}a\rho^2|\ss_{\Nb}|^2\right)^{1/2}d\ub\nonumber\\
&&\hspace{10mm}\leq C\int_0^{\ub_1}\left(\s^{(V)}\underline{{\cal E}}^{u_1}(\ub)\right)^{1/2}
\left(\sum_{\mu}\s^{(E_{(\mu)})}\underline{{\cal E}}^{u_1}(\ub)\right)^{1/2}\sqrt{\ub}d\ub\nonumber
\end{eqnarray}
where we have used the fact that
$$\frac{a\rho^2}{\rhob^2}\sim\frac{\rho^3}{\rhob}\sim\frac{\ub^3}{u^2}$$
and $\ub^3/u^2\leq\ub$. 

To estimate \ref{7.167}  we first note that by \ref{7.147}, in view of the fact that by \ref{7.108} 
and \ref{7.132} 
$$\frac{\rho}{\ogamma}\sim\frac{\ub}{u}, \ \ \ 
\frac{\rho\rhob}{\ogamma}\sim\frac{a}{\ogamma}\sim\sqrt{\rho}\sqrt{a},$$
we have:
\begin{equation}
\rho|s_{\Nb\Lb}|\leq C(|\s^{(Y)}\xi_{\Lb}|+\sqrt{\rho}\sqrt{a}|\sss|)
\label{7.176}
\end{equation}
Hence by \ref{7.82} and \ref{7.145}:
\begin{equation}
\int_{\Cb_{\ub}^{u_1}}\Omega^{(n-1)/2}\rho^2 s_{\Nb\Lb}^2
\leq C\left(\s^{(Y)}\underline{{\cal E}}^{u_1}(\ub)
+\sum_{\mu}\s^{(E_{(\mu)})}\underline{{\cal E}}^{u_1}(\ub)\right)
\label{7.177}
\end{equation}
In regard then to the partial contribution of the 1st term in parenthesis in \ref{7.167}, by \ref{7.177} we have:
\begin{eqnarray}
&&\int_{{\cal R}_{\ub_1,u_1}}\Omega^{(n-1)/2}|\s^{(V)}\xi_{\Lb}|\rho|s_{\Nb\Lb}|\leq \label{7.178}\\
&&\hspace{15mm}\int_0^{\ub_1}\left(\int_{\Cb_{\ub}^{u_1}}\Omega^{(n-1)/2}(\s^{(V)}\xi_{\Lb})^2\right)^{1/2}\left(\int_{\Cb_{\ub}^{u_1}}\rho^2 s_{\Nb\Lb}^2\right)^{1/2}d\ub
\nonumber\\
&&\hspace{15mm}\leq C\int_0^{\ub_1}\left(\s^{(Y)}\underline{{\cal E}}^{u_1}(\ub)\right)^{1/2}
\left(\s^{(Y)}\underline{{\cal E}}^{u_1}(\ub)
+\sum_{\mu}\s^{(E_{(\mu)})}\underline{{\cal E}}^{u_1}(\ub)\right)^{1/2}d\ub\nonumber
\end{eqnarray}
The partial contribution of the 2nd term in parenthesis in \ref{7.167} is estimated using 
\ref{7.145} by:
\begin{eqnarray}
&&\int_{{\cal R}_{\ub_1,u_1}}\Omega^{(n-1)/2}|\s^{(V)}\xi_{\Lb}|a|\sss|\leq \label{7.179}\\
&&\hspace{15mm}\int_0^{\ub_1}\left(\int_{\Cb_{\ub}^{u_1}}\Omega^{(n-1)/2}(\s^{(V)}\xi_{\Lb})^2\right)^{1/2}\left(\int_{\Cb_{\ub}^{u_1}}\Omega^{(n-1)/2}
a^2|\sss|^2\right)^{1/2}d\ub\nonumber\\
&&\hspace{10mm}\leq Cu_1\int_0^{\ub_1}\left(\s^{(Y)}\underline{{\cal E}}^{u_1}(\ub)\right)^{1/2}
\left(\s^{(Y)}\underline{{\cal E}}^{u_1}(\ub)
+\sum_{\mu}\s^{(E_{(\mu)})}\underline{{\cal E}}^{u_1}(\ub)\right)^{1/2}\sqrt{\ub}d\ub\nonumber
\end{eqnarray}
(relative to \ref{7.145} there is an extra factor of $a\sim\ub u^2$ in the 2nd integral on $\Cb_{\ub}^{u_1}$). 
The partial contribution of the 3rd term in parenthesis in \ref{7.167} is estimated using \ref{7.144} 
as follows:
\begin{eqnarray}
&&\int_{{\cal R}_{\ub_1,u_1}}\Omega^{(n-1)/2}|\s^{(V)}\xi_{\Lb}|a|\ss_N|\leq \label{7.180}\\
&&\hspace{10mm}\left(\int_{{\cal R}_{\ub_1,u_1}}\Omega^{(n-1)/2}(\s^{(V)}\xi_{\Lb})^2\right)^{1/2}
\left(\int_{{\cal R}_{\ub_1,u_1}}\Omega^{(n-1)/2}a^2|\ss_N|^2\right)^{1/2}\nonumber\\
&&\hspace{15mm}\leq C\left(\int_0^{\ub_1}\s^{(V)}\underline{{\cal E}}^{u_1}(\ub)d\ub\right)^{1/2}
\left(\int_0^{u_1}u^4\sum_{\mu}\s^{(E_{(\mu)})}{\cal E}^{\ub_1}(u)du\right)^{1/2}\nonumber
\end{eqnarray}
for, in the 2nd integral, since $a^2\sim\rho^2\rhob^2$ there is an extra factor of $\rhob^2\sim u^4$ 
relative to \ref{7.144}. Finally, the partial contribution of the 4th term in parenthesis in \ref{7.167} is estimated using \ref{7.145}:
\begin{eqnarray}
&&\int_{{\cal R}_{\ub_1,u_1}}\Omega^{(n-1)/2}|\s^{(V)}\xi_{\Lb}|a|\ss_{\Nb}|\leq \label{7.181}\\
&&\hspace{10mm}\int_0^{\ub_1}\left(\int_{\Cb_{\ub}^{u_1}}\Omega^{(n-1)/2}(\s^{(V)}\xi_{\Lb})^2\right)^{1/2}\left(\int_{\Cb_{\ub}^{u_1}}a^2|\ss_{\Nb}|^2\right)^{1/2}\nonumber\\
&&\hspace{15mm}\leq C\int_0^{\ub_1}\left(\s^{(V)}\underline{{\cal E}}^{u_1}(\ub)\right)^{1/2}
\left(\sum_{\mu}\s^{(E_{(\mu)})}\underline{{\cal E}}^{u_1}(\ub)\right)^{1/2}\ub d\ub\nonumber
\end{eqnarray}
for, in the 2nd integral, since $a^2\sim\rho^2\rhob^2$ there is an extra factor of $\rho^2\sim\ub^2$ 
relative to \ref{7.145}.

To estimate \ref{7.168}  we first note that by \ref{7.147}, in view of the fact that by \ref{7.108} 
and \ref{7.132} 
$$\frac{\sqrt{a}}{\ogamma}\sim\sqrt{\rho}, \ \ \ 
\frac{\rhob}{\ogamma}\sim\sqrt{\rhob},$$
we have:
\begin{equation}
\sqrt{a}|s_{\Nb\Lb}|\leq C(\sqrt{\rho}|\s^{(Y)}\xi_{\Lb}|+\sqrt{\rhob}\sqrt{a}|\sss|)
\label{7.182}
\end{equation}
Hence by \ref{7.82} and \ref{7.145}:
\begin{equation}
\int_{\Cb_{\ub}^{u_1}}as_{\Nb\Lb}^2\leq C\left(\ub\s^{(Y)}\underline{{\cal E}}^{u_1}(\ub)
+u_1^2\sum_{\mu}\s^{(E_{(\mu)})}\underline{{\cal E}}^{u_1}(\ub)\right)
\label{7.183}
\end{equation}
In regard then to the partial contribution of the 1st term in parenthesis in \ref{7.168}, by \ref{7.183} we have:
\begin{eqnarray}
&&\int_{{\cal R}_{\ub_1,u_1}}\Omega^{(n-1)/2}a|\s^{(V)}\sxi||s_{\Nb\Lb}|\leq \label{7.184}\\
&&\hspace{10mm}\int_0^{\ub_1}\left(\int_{\Cb_{\ub}^{u_1}}a|\s^{(V)}\sxi|^2\right)^{1/2}
\left(\int_{\Cb_{\ub}^{u_1}}a s_{\Nb\Lb}^2\right)^{1/2} d\ub\nonumber\\
&&\hspace{10mm}\leq C\int_0^{\ub_1}\left(\s^{(V)}\underline{{\cal E}}^{u_1}(\ub)\right)^{1/2}
\left(\ub\s^{(Y)}\underline{{\cal E}}^{u_1}(\ub)
+u_1^2\sum_{\mu}\s^{(E_{(\mu)})}\underline{{\cal E}}^{u_1}(\ub)\right)^{1/2} d\ub \nonumber
\end{eqnarray}
The partial contribution of the 2nd term in parenthesis in \ref{7.168} is estimated using 
\ref{7.145} in a straightforward manner:
\begin{eqnarray}
&&\int_{{\cal R}_{\ub_1,u_1}}\Omega^{(n-1)/2}a|\s^{(V)}\sxi||\sss|\leq \label{7.185}\\
&&\hspace{15mm}\int_0^{\ub_1}\left(\int_{\Cb_{\ub}^{u_1}}\Omega^{(n-1)/2}a|\s^{(V)}\sxi|^2\right)^{1/2}
\left(\int_{\Cb_{\ub}^{u_1}}\Omega^{(n-1)/2}a|\sss|^2\right)^{1/2}d\ub\nonumber\\
&&\hspace{15mm}\leq C\int_0^{\ub_1}\left(\s^{(V)}\underline{{\cal E}}^{u_1}(\ub)\right)^{1/2}
\left(\sum_{\mu}\s^{(E_{(\mu)})}\underline{{\cal E}}^{u_1}(\ub)\right)^{1/2}d\ub\nonumber
\end{eqnarray}
We turn to the partial contribution of the 3rd term in parenthesis in \ref{7.168}. We estimate this 
using \ref{7.144} as follows:
\begin{eqnarray}
&&\int_{{\cal R}_{\ub_1,u_1}}\Omega^{(n-1)/2}a|\s^{(V)}\sxi||\ss_N|\leq \label{7.186}\\
&&\hspace{10mm}\left(\int_{{\cal R}_{\ub_1,u_1}}\Omega^{(n-1)/2}\frac{a}{\rho}
|\s^{(V)}\sxi|^2\right)^{1/2}\left(\int_{{\cal R}_{\ub_1,u_1}}\Omega^{(n-1)/2}a\rho
|\ss_N|^2\right)^{1/2}\nonumber\\
&&\hspace{10mm}\leq C\left(\int_0^{\ub_1}\s^{(V)}\underline{{\cal E}}^{u_1}(\ub)
\frac{d\ub}{\ub}\right)^{1/2}
\left(\int_0^{u_1}u^2\sum_{\mu}\s^{(E_{(\mu)})}{\cal E}^{\ub_1}(u)du\right)^{1/2}\nonumber
\end{eqnarray}
since $\rho\sim\ub$ while $a\rho\sim\rhob\rho^2\sim u^2\rho^2$. Here the 1st integral on the right 
is a {\em borderline integral}. Finally, the partial contribution of the 4th term in parenthesis in 
\ref{7.168} is estimated using \ref{7.145} as follows:
\begin{eqnarray}
&&\int_{{\cal R}_{\ub_1,u_1}}\Omega^{(n-1)/2}a|\s^{(V)}\sxi||\ss_{\Nb}|\leq \label{7.187}\\
&&\hspace{10mm}\int_0^{\ub_1}\left(\int_{\Cb_{\ub}^{u_1}}\Omega^{(n-1)/2}a|\s^{(V)}\sxi|^2\right)^{1/2}
\left(\int_{\Cb_{\ub}^{u_1}}\Omega^{(n-1)/2}a|\ss_{\Nb}|^2\right)^{1/2}d\ub\nonumber\\
&&\hspace{10mm}\leq C\int_0^{\ub_1}\left(\s^{(V)}\underline{{\cal E}}^{u_1}(\ub)\right)^{1/2}
\left(\sum_{\mu}\s^{(E_{(\mu)})}\underline{{\cal E}}^{u_1}(\ub)\right)^{1/2}\frac{d\ub}{\sqrt{\ub}}
\nonumber
\end{eqnarray}
since, comparing with \ref{7.145}, $a/\rhob^2\sim\rho/\rhob\sim\ub/u^2$ and the last is not greater than 
$1/\ub$. 

\pagebreak

\chapter{The Commutation Fields}

\section{The Commutation Fields and the Higher Order Variations}

We shall presently define the higher order variations. As noted in Chapter 6, the higher order 
variations are generated by applying to the fundamental variations \ref{6.22}, namely to the 
$\beta_\mu : \mu=0,...n$, a string of commutation fields $C$. The commutation fields $C$ are 
vectorfields on ${\cal N}$ which along ${\cal K}$ are tangential to ${\cal K}$. This is required 
so that the 1-forms $\s^{(V)}\xi$ associated to these higher order variations and to the variation 
fields $V$ satisfy to principal terms the boundary condition of Proposition 6.2, or equivalently 
\ref{6.164} - \ref{6.165}. Moreover, we require that at each point on ${\cal K}$ the set of 
commutation fields $C$ spans the tangent space to ${\cal K}$ at the point. Together with the analogous 
requirement on the set of variation fields (see Section 6.2), this ensures control on all derivatives 
of the $\beta_\mu$ when the wave system (see Section 2.3) is imposed. 

As first of the commutation fields we take the vectorfield $T$. The remaining commutation fields 
are required to span the tangent space to the $S_{\ub,u}$ at each point.  In the case $n=2$ only 
one remaining commutation field is needed and we take this to be the vectorfield $E$. For $n>2$ 
we take the remaining commutation fields to be the $E_{(\mu)} : \mu=0,...,n$. We note that 
$E$ in the case $n=2$ and the $E_{(\mu)} : \mu=0,...,n$ for $n>2$, play a dual role being 
both commutation fields as well as variation fields. These are however complemented by $T$ in 
the set of commutation fields, whereas they are complemented by $Y$ in the set of variation fields. 
As we have seen in Chapter 6, the error terms in the energy identities associated to a variation field 
$V$ are generated by the corresponding structure forms $\s^{(V)}\theta^\mu=dV^\mu$, 
$V^\mu : \mu=0,...,n$ being the rectangular components of $V$. On the other hand, as we shall presently 
see, the error terms in the energy identities associated to a commutation field $C$ are generated 
by the corresponding deformation tensors $\s^{(C)}\tilde{\pi}={\cal L}_C\tilde{h}$ through the 
source functions $\rho_\mu$. (See Section 6.3, in particular \ref{6.63} - \ref{6.66}.) 

In the case $n=2$ a higher order variation $\s^{(m,l)}\dot{\phi}_\mu$ corresponds to 
a pair of non-negative integers 
$(m,l)$, not both zero, where $m$ is the number of $T$ commutation fields and $l$ the number of 
$E$ commutation fields applied. It is defined by:
\begin{equation}
\s^{(m,l)}\dot{\phi}_\mu=E^l T^m \beta_\mu 
\label{8.1}
\end{equation}
$\s^{(0,0)}\dot{\phi}_\mu=\beta_\mu$ being the fundamental variations \ref{6.22}. For $n>2$ 
a higher order variation  $\s^{(m,\nu_1...\nu_l)}\dot{\phi}_\mu$ corresponds to a non-negative 
integer $m$ and to a string $\nu_1 . . . \nu_l$ of length $l$, where each index $\nu_i : i=1,...,l$ 
belongs to $\{0,...,n\}$. It is defined by:
\begin{equation}
\s^{(m,\nu_1...\nu_l)}\dot{\phi}_\mu=E_{(\nu_l)}...E_{(\nu_1)}T^m\beta_\mu 
\label{8.2}
\end{equation}

\vspace{5mm}

\section{The Recursion Formulas for the Source Functions}

The higher order variations satisfy inhomogeneous wave equations of the form \ref{6.25}. 
More precisely: 
\begin{eqnarray}
&&\square_{\tilde{h}}\s^{(m,l)}\dot{\phi}_\mu=\s^{(m,l)}\rho_\mu \ \ \mbox{: in the case $n=2$} 
\label{8.3}\\
&&\square_{\tilde{h}}\s^{(m,\nu_1...\nu_l)}\dot{\phi}_{\mu}=\s^{(m,\nu_1...\nu_l)}\rho_\mu 
\ \ \mbox{: for $n>2$} \label{8.4}
\end{eqnarray}
where the {\em source functions} $\s^{(m,l)}\rho_\mu$, $\s^{(m,\nu_1...\nu_l)}\rho_\mu$ 
shall be recursively determined using the following lemma. 

\vspace{2.5mm}

\noindent{\bf Lemma 8.1} \ \ \ Let $\psi$ be a solution of the inhomogeneous wave equation 
$$\square_{\tilde{h}}\psi=\rho$$
and let $C$ be an arbitrary vectorfield. Then 
$$\s^{(C)}\psi=C\psi$$
satisfies the inhomogeneous wave equation
$$\square_{\tilde{h}}\s^{(C)}\psi=\s^{(C)}\rho$$
where the source function $\s^{(C)}\rho$ is related to the source function $\rho$ by: 
$$\s^{(C)}\rho=C\rho+\tilde{\mbox{div}}\s^{(C)}J
+\frac{1}{2}\tilde{\mbox{tr}}\s^{(C)}\tilde\pi \rho$$
Here $\s^{(C)}J$ is the {\em commutation current} associated to $\psi$ and $C$, a vectorfield, 
given in arbitrary local coordinates by:
$$\s^{(C)}J^\alpha=\left(\s^{(C)}\tilde{\pi}^{\alpha\beta}
-\frac{1}{2}(\tilde{h}^{-1})^{\alpha\beta}\tilde{\mbox{tr}}\s^{(C)}\tilde{\pi}\right)
\partial_\beta\psi$$
with
$$\s^{(C)}\tilde{\pi}^{\alpha\beta}=(\tilde{h}^{-1})^{\alpha\gamma}(\tilde{h}^{-1})^{\beta\delta}
\s^{(C)}\tilde{\pi}_{\gamma\delta}$$
$\s^{(C)}\tilde{\pi}$ being ${\cal L}_C\tilde{h}$, the Lie derivative of the metric $\tilde{h}$ with 
respect to $C$. 

\vspace{2.5mm}

\noindent{\em Proof:} Let $f_s$ be the local 1-parameter group of diffeomorphisms generated by $C$. 
Denoting by $f_s^*$ the corresponding pullback, we have:
\begin{equation}
\square_{f_s^*\tilde{h}}f_s^*\psi=f_s^*\left(\square_{\tilde{h}}\psi\right)=f_s^*\rho
\label{8.5}
\end{equation}
Now, in arbitrary local coordinates we have:
\begin{equation}
\square_{\tilde{h}}\psi=\frac{1}{\sqrt{-\mbox{det}\tilde{h}}}\partial_\alpha
\left(\sqrt{-\mbox{det}\tilde{h}}(\tilde{h}^{-1})^{\alpha\beta}\partial_\beta\psi\right)
\label{8.6}
\end{equation}
and:
\begin{equation}
\square_{f_s^*\tilde{h}}f_s^*\psi=\frac{1}{\sqrt{-\mbox{det}f_s^*\tilde{h}}}\partial_\alpha
\left(\sqrt{-\mbox{det}f_s^*\tilde{h}}((f_s^*\tilde{h})^{-1})^{\alpha\beta}
\partial_\beta f_s^*\psi\right)
\label{8.7}
\end{equation}
Let us differentiate the expression \ref{8.7} with respect to $s$ at $s=0$. In view of the facts that:
$$\left.\frac{d}{ds}f_s^*\tilde{h}\right|_{s=0}={\cal L}_C\tilde{h}=\s^{(C)}\tilde{\pi},$$
$$\left.\frac{d}{ds}\sqrt{-\mbox{det}f_s^*\tilde{h}}\right|_{s=0}=
\frac{1}{2}\sqrt{-\mbox{det}\tilde{h}}\tilde{\mbox{tr}}\s^{(C)}\tilde{\pi},$$
$$\left.\frac{d}{ds}((f_s^*\tilde{h})^{-1})^{\alpha\beta}\right|_{s=0}=
-\s^{(C)}\tilde{\pi}^{\alpha\beta}$$
and: 
$$\left.\frac{d}{ds}f_s^*\psi\right|_{s=0}=\left.\frac{d}{ds}\psi\circ f_s\right|_{s=0}=C\psi$$
we obtain:
\begin{eqnarray}
&&\left.\frac{d}{ds}\square_{f_s^*\tilde{h}}f_s^*\psi\right|_{s=0}=
-\frac{1}{2}\tilde{\mbox{tr}}\s^{(C)}\tilde{\pi}\frac{1}{\sqrt{-\mbox{det}\tilde{h}}}
\partial_\alpha\left(\sqrt{-\mbox{det}\tilde{h}}(\tilde{h}^{-1})^{\alpha\beta}\partial_\beta\psi\right)
\nonumber\\
&&\hspace{25mm}+\frac{1}{\sqrt{-\mbox{det}\tilde{h}}}\partial_\alpha
\left\{\frac{1}{2}\tilde{\mbox{tr}}\s^{(C)}\tilde{\pi}\sqrt{-\mbox{det}\tilde{h}}
(\tilde{h}^{-1})^{\alpha\beta}\partial_\beta\psi\right. \label{8.8}\\
&&\hspace{25mm}\left.-\sqrt{-\mbox{det}\tilde{h}}\s^{(C)}\tilde{\pi}^{\alpha\beta}\partial_\beta\psi
+\sqrt{-\mbox{det}\tilde{h}}(\tilde{h}^{-1})^{\alpha\beta}\partial_\beta(C\psi)\right\}\nonumber
\end{eqnarray}
In view of \ref{8.6}, a similar expression with $C\psi$ in the role of $\psi$, and the fact that:
$$\frac{1}{\sqrt{-\mbox{det}\tilde{h}}}\partial_\alpha
\left(\sqrt{-\mbox{det}\tilde{h}}\s^{(C)}J^\alpha\right)=\tilde{\mbox{div}}\s^{(C)}J$$
\ref{8.8} is:
\begin{equation}
\left.\frac{d}{ds}\square_{f_s^*\tilde{h}}f_s^*\psi\right|_{s=0}=
-\frac{1}{2}\tilde{\mbox{tr}}\s^{(C)}\tilde{\pi}\square_{\tilde{h}}\psi
-\tilde{\mbox{div}}\s^{(C)}J+\square_{\tilde{h}}(C\psi)
\label{8.9}
\end{equation}
On the other hand by \ref{8.5} this is equal to:
\begin{equation}
\left.\frac{d}{ds}f_s^*\rho\right|_{s=0}=\left.\frac{d}{ds}\rho\circ f_s\right|_{s=0}=C\rho
\label{8.10}
\end{equation}
Comparing \ref{8.9} with \ref{8.10} the lemma follows. 

\vspace{2.5mm}

To derive a recursion formula for the sources $\s^{(m,l)}\rho_\mu$, or 
$\s^{(m,\nu_1...\nu_l)}\rho_\mu$, we define the {\em rescaled sources}:
\begin{eqnarray}
&&\s^{(m,l)}\tilde{\rho}_\mu=\Omega a\s^{(m,l)}\rho_\mu \ \ \mbox{: in the case $n=2$} \label{8.11}\\
&&\s^{(m,\nu_1...\nu_l)}\tilde{\rho}_\mu=\Omega a\s^{(m,\nu_1...\nu_l)}\rho_\mu \ \ 
\mbox{: for $n>2$} \label{8.12} 
\end{eqnarray}
The reason why we rescale the sources in this way is that the operator $\Omega a\square_{\tilde{h}}$ 
is regular in terms of the frame field $(L,\Lb,\Omega_A : A=1,...,n-1)$ even where $a$ vanishes, 
in contrast to the operator $\square_{\tilde{h}}$. To see this we first remark that since 
$$\sqrt{-\mbox{det}\tilde{h}}=\Omega^{(n+1)/2}\sqrt{-\mbox{det}h}, \ \ \ 
(\tilde{h}^{-1})^{\alpha\beta}=\Omega^{-1}(h^{-1})^{\alpha\beta}$$
the expression \ref{8.6} gives, for an arbitrary function $\psi$, 
\begin{equation}
\square_{\tilde{h}}\psi=\Omega^{-1}\square_h\psi+\frac{(n-1)}{2}\Omega^{-2}h^{-1}(d\Omega, d\psi)
\label{8.13}
\end{equation}
Now, $\square_h\psi$ is the trace relative to $h$ of $D^2\psi$, the covariant relative to $h$ 
Hessian of $\psi$. Thus by the expansion \ref{2.38} of $h^{-1}$ in terms of the frame field 
$(L,\Lb,\Omega_A : A=1,...,n-1)$ we have:
\begin{equation}
\square_h\psi=a^{-1}D^2\psi(L,\Lb)+(\sh^{-1})^{AB}D^2\psi(\Omega_A,\Omega_B)
\label{8.14}
\end{equation}
Moreover, 
\begin{equation}
D^2\psi(L,\Lb)=L(\Lb\psi)-(D_L\Lb)\psi=L(\Lb\psi)-2\etab^\sharp\cdot\sd\psi
\label{8.15}
\end{equation}
by \ref{3.15}. Alternatively, we can express 
\begin{equation}
D^2\psi(L,\Lb)=\Lb(L\psi)-(D_{\Lb}L)\psi=\Lb(L\psi)-2\eta^\sharp\cdot\sd\psi
\label{8.16}
\end{equation}
by \ref{8.14}. Furthermore, we have:
\begin{equation}
D^2\psi(\Omega_A,\Omega_B)=\Omega_A(\Omega_B\psi)-(D_{\Omega_A}\Omega_B)\psi
\label{8.17}
\end{equation}
On the other hand, $\sD^2\psi$ the covariant Hessian relative to $\sh$ of the restriction 
of $\psi$ to each $S_{\ub,u}$, a symmetric 2-covariant $S$ tensorfield, is given by:
\begin{equation}
\sD^2\psi(\Omega_A,\Omega_B)=\Omega_A(\Omega_B\psi)-(\sD_{\Omega_A}\Omega_B)\psi
\label{8.18}
\end{equation}
Subtracting \ref{8.18} from \ref{8.17} and using the fact that by \ref{3.26}:
\begin{equation}
D_{\Omega_A}\Omega_B=\sD_{\Omega_A}\Omega_B+\frac{1}{2a}(\chib_{AB}L+\chi_{AB}\Lb)
\label{8.19}
\end{equation}
we obtain:
\begin{equation}
D^2\psi(\Omega_A,\Omega_B)=\sD^2\psi(\Omega_A,\Omega_B)-\frac{1}{2a}(\chib_{AB}L\psi+\chi_{AB}\Lb\psi)
\label{8.20}
\end{equation}
Substituting \ref{8.15}, or \ref{8.16}, and \ref{8.20} yields:
\begin{equation}
a\square_h\psi=-L(\Lb\psi)+2\etab^\sharp\cdot\sd\psi+a\slap_{\sh}\psi
-\frac{1}{2}(\mbox{tr}\chib L\psi+\mbox{tr}\chi\Lb\psi) 
\label{8.21}
\end{equation}
or: 
\begin{equation}
a\square_h\psi=-\Lb(L\psi)+2\eta^\sharp\cdot\sd\psi+a\slap_{\sh}\psi
-\frac{1}{2}(\mbox{tr}\chib L\psi+\mbox{tr}\chi\Lb\psi)
\label{8.22}
\end{equation}
where $\slap_{\sh}\psi$ is the covariant Laplacian relative to $\sh$ of the restriction of $\psi$ 
to each $S_{\ub,u}$, the trace relative to $\sh$ of $\sD^2\psi$. These expressions show that 
$a\square_h$, in contrast to $\square_h$, is a regular operator in terms of the frame field 
$(L,\Lb,\Omega_A : A=1,...,n-1)$, thus a regular operator in acoustical coordinates, even where 
$a$ vanishes. By \ref{8.13} and the expansion \ref{2.38} we have:
\begin{equation}
\Omega a\square_{\tilde{h}}\psi=a\square_h\psi+\frac{(n-1)}{2}\Omega^{-1}
\left(-\frac{1}{2}((\Lb\Omega)L\psi+(L\Omega)\Lb\psi)+a\sh^{-1}(\sd\Omega,\sd\psi)\right)
\label{8.23}
\end{equation}
Consequently, $\Omega a\square_{\tilde{h}}$, in contrast to $\square_{\tilde{h}}$, is a regular 
operator even where $a$ vanishes. 

Going now back to Lemma 8.1 we define:
\begin{equation}
\tilde{\rho}=\Omega a\rho, \ \ \ \s^{(C)}\tilde{\rho}=\Omega a\s^{(C)}\rho
\label{8.24}
\end{equation}
Then according to the lemma $\s^{(C)}\tilde{\rho}$ is related to $\tilde{\rho}$ by:
\begin{equation}
\s^{(C)}\tilde{\rho}=C\tilde{\rho}+\s^{(C)}\sigma+\s^{(C)}\delta \tilde{\rho}
\label{8.25}
\end{equation}
where we denote:
\begin{equation}
\s^{(C)}\sigma=\Omega a\tilde{\mbox{div}}\s^{(C)}J
\label{8.26}
\end{equation}
and:
\begin{equation}
\s^{(C)}\delta=\frac{1}{2}\tilde{\mbox{tr}}\s^{(C)}\tilde{\pi}-a^{-1}Ca-\Omega^{-1}C\Omega
\label{8.27}
\end{equation}
Here $\s^{(C)}\sigma$ is associated to the commutation field $C$ and to the variation $\psi$. 

Consider first the case $n=2$. Then with $\s^{(m,l)}\dot{\phi}_\mu$, given by \ref{8.1} in the 
role of $\psi$, we have $\s^{(m,l)}\rho_\mu$, $\s^{(m,l)}\tilde{\rho}_\mu$ in the role of 
$\rho$, $\tilde{\rho}$ respectively. Let us denote by $\s^{(C,m,l)}\sigma_\mu$ the function 
$\s^{(C)}\sigma$ as above associated to the commutation field $C$ and to the variation 
$\s^{(m,l)}\dot{\phi}_\mu$. Taking $C=T$ we have $\s^{(m+1,l)}\dot{\phi}_\mu$ in the role of 
$\s^{(C)}\psi$ and $\s^{(m+1,l)}\tilde{\rho}_\mu$ in the role of $\s^{(C)}\tilde{\rho}$. Also, 
taking $C=E$ we have $\s^{(m,l+1)}\dot{\phi}_\mu$ in the role of $\s^{(C)}\psi$ and 
$\s^{(m,l+1)}\tilde{\rho}_\mu$ in the role of $\s^{(C)}\tilde{\rho}$. Then in the first case 
\ref{8.25} takes the form of the recursion formula:
\begin{equation}
\s^{(m+1,l)}\tilde{\rho}_\mu=T\s^{(m,l)}\tilde{\rho}_\mu+\s^{(T,m,l)}\sigma_\mu
+\s^{(T)}\delta\s^{(m,l)}\tilde{\rho}_\mu
\label{8.28}
\end{equation}
and in the second case \ref{8.25} takes the form of the recursion formula:
\begin{equation}
\s^{(m,l+1)}\tilde{\rho}_\mu=E\s^{(m,l)}\tilde{\rho}_\mu+\s^{(E,m,l)}\sigma_\mu
+\s^{(E)}\delta\s^{(m,l)}\tilde{\rho}_\mu 
\label{8.29}
\end{equation}
Given the fact that 
\begin{equation}
\s^{(0,0)}\tilde{\rho}_\mu=0
\label{8.30}
\end{equation}
as $\s^{(0,0)}\dot{\phi}_\mu=\beta_\mu$ satisfies the homogeneous wave equation 
relative to $\tilde{h}$, the recursion formulas 
\ref{8.28}, \ref{8.29} determine $\s^{(m,l)}\tilde{\rho}_\mu$ for all pairs of non-negative integers 
$(m,l)$. In fact \ref{8.28} with $l=0$ first determines $\s^{(m,0)}\tilde{\rho}_\mu$ for any $m$, 
after which \ref{8.29} determines $\s^{(m,l)}\tilde{\rho}_\mu$ for any $l$. 

For $n>2$, with $\s^{(m,\nu_1...\nu_l)}\dot{\phi}_\mu$, given by \ref{8.2} in the 
role of $\psi$, we have $\s^{(m,\nu_1...\nu_l)}\rho_\mu$, $\s^{(m,\nu_1...\nu_l)}\tilde{\rho}_\mu$ 
in the role of $\rho$, $\tilde{\rho}$ respectively. 
Let us denote by $\s^{(C,m,\nu_1...\nu_l)}\sigma_\mu$ the function 
$\s^{(C)}\sigma$ as above associated to the commutation field $C$ and to the variation 
$\s^{(m,\nu_1...\nu_l)}\dot{\phi}_\mu$. 
Taking $C=T$ we have $\s^{(m+1,\nu_1...\nu_l)}\dot{\phi}_\mu$ in the role of 
$\s^{(C)}\psi$ and $\s^{(m+1,\nu_1...\nu_l)}\tilde{\rho}_\mu$ in the role of 
$\s^{(C)}\tilde{\rho}$. Also, taking $C=E_{(\nu_{l+1})}$ we have 
$\s^{(m,\nu_1...\nu_{l+1})}\dot{\phi}_\mu$ in the role of $\s^{(C)}\psi$ and 
$\s^{(m,\nu_1...\nu_{l+1})}\tilde{\rho}_\mu$ in the role of $\s^{(C)}\tilde{\rho}$. 
Then in the first case \ref{8.27} takes the form of the recursion formula:
\begin{equation}
\s^{(m+1,\nu_1...\nu_l)}\tilde{\rho}_\mu=T\s^{(m,\nu_1...\nu_l)}\tilde{\rho}_\mu
+\s^{(T,m,\nu_1...\nu_l)}\sigma_\mu
+\s^{(T)}\delta\s^{(m,\nu_1...\nu_l)}\tilde{\rho}_\mu
\label{8.31}
\end{equation}
and in the second case \ref{8.27} takes the form of the recursion formula:
\begin{equation}
\s^{(m,\nu_1...\nu_{l+1})}\tilde{\rho}_\mu=E_{(\nu_{l+1})}\s^{(m,\nu_1...\nu_l)}\tilde{\rho}_\mu
+\s^{(E_{(\nu_{l+1})},m,\nu_1...\nu_l)}\sigma_\mu
+\s^{(E_{(\nu_{l+1})})}\delta\s^{(m,\nu_1...\nu_l)}\tilde{\rho}_\mu 
\label{8.32}
\end{equation}
Given again the fact that 
\begin{equation}
\s^{(0)}\tilde{\rho}_\mu=0
\label{8.33}
\end{equation}
as $\s^{(0)}\dot{\phi}_\mu=\beta_\mu$ satisfies the homogeneous wave equation 
relative to $\tilde{h}$, the recursion formulas 
\ref{8.31}, \ref{8.32} determine $\s^{(m,\nu_1...\nu_l)}\tilde{\rho}_\mu$ for all non-negative integers 
$m$ and all $\nu_i\in\{0,...,n\}$, $i=1,...,l$, $l$ being any positive integer. 
In fact \ref{8.31} with $l=0$ first determines $\s^{(m)}\tilde{\rho}_\mu$ for any $m$, 
after which \ref{8.32} determines $\s^{(m,\nu_1...\nu_l)}\tilde{\rho}_\mu$ for any $l$. 

Consider now  $\s^{(C)}\sigma$ (see \ref{8.26}). We express $\tilde{\mbox{div}}\s^{(C)}J$ in 
acoustical coordinates using the formulas \ref{6.56}, \ref{6.67}, \ref{6.68}, and \ref{6.70}, 
with the vectorfield $\s^{(C)}J$ in the role of the vectorfield $\s^{(V)}P$. Expanding then 
in terms of the acoustical coordinate frame field 
$(\partial/\partial\ub,\partial/\partial u,\partial/\partial\vartheta^A:A=1,...,n-1)$:
\begin{equation}
\s^{(C)}J=\s^{(C)}J^{\ub}\frac{\partial}{\partial\ub}+\s^{(C)}J^u\frac{\partial}{\partial u}
+\s^{(C)}J^{\vartheta^A}\frac{\partial}{\partial\vartheta^A}
\label{8.34}
\end{equation}
we have:
\begin{eqnarray}
&&\s^{(C)}\sigma=\frac{\Omega^{-(n-1)/2}}{\sqrt{\mbox{det}\sh}}
\left\{\frac{\partial}{\partial\ub}\left(a\Omega^{(n+1)/2}\sqrt{\mbox{det}\sh}\s^{(C)}J^{\ub}\right)
\right.\label{8.35}\\
&&\hspace{10mm}\left.
+\frac{\partial}{\partial u}\left(a\Omega^{(n+1)/2}\sqrt{\mbox{det}\sh}\s^{(C)}J^u\right)
+\frac{\partial}{\partial\vartheta^A}\left(a\Omega^{(n+1)/2}\sqrt{\mbox{det}\sh}
\s^{(C)}J^{\vartheta^A}\right)\right\}\nonumber
\end{eqnarray}
We can also expand the vectorfield $\s^{(C)}J$ as:
\begin{equation}
\s^{(C)}J=\s^{C)}J^L L+\s^{(C)}J^{\Lb}\Lb+\s^{(C)}\sJ
\label{8.36}
\end{equation}
where
\begin{equation}
\s^{(C)}\sJ=\s^{(C)}J^A\Omega_A
\label{8.37}
\end{equation}
is a $S$ vectorfield. So \ref{8.36} with \ref{8.37} substituted is the expansion of $\s^{(C)}J$ in 
terms of the frame field $(L,\Lb,\Omega_A:A=1,...,n-1)$. In view of the relations \ref{6.77}, 
comparing with \ref{8.34} we conclude that: 
\begin{equation}
\s^{(C)}J^{\ub}=\s^{(C)}J^L, \ \ \s^{(C)}J^u=\s^{(C)}J^{\Lb}, \ \ 
\s^{(C)}J^{\vartheta^A}=\s^{(C)}J^A+b^A(\s^{(C)}J^{\Lb}-\s^{(C)}J^L)
\label{8.38}
\end{equation}
in analogy with \ref{6.78}. Substituting these in \ref{8.35} and recalling again \ref{6.77} we obtain: 
\begin{eqnarray}
&&\s^{(C)}\sigma=\Omega^{-(n-1)/2}\left\{L(\Omega^{(n+1)/2}a\s^{(C)}J^L)
+\Lb(\Omega^{(n+1)/2}a\s^{(C)}J^{\Lb})\right.\nonumber\\
&&\hspace{35mm}\left.+\sdiv(\Omega^{(n+1)/2}a\s^{(C)}\sJ)\right\} \label{8.39}\\
&&\hspace{20mm}+\frac{\Omega}{\sqrt{\mbox{det}\sh}}
\left(\frac{\partial\sqrt{\mbox{det}\sh}}{\partial\ub}
-\frac{\partial}{\partial\vartheta^A}(\sqrt{\mbox{det}\sh}b^A)\right)a\s^{(C)}J^L\nonumber\\
&&\hspace{20mm}+\frac{\Omega}{\sqrt{\mbox{det}\sh}}
\left(\frac{\partial\sqrt{\mbox{det}\sh}}{\partial u}
+\frac{\partial}{\partial\vartheta^A}(\sqrt{\mbox{det}\sh}b^A)\right)a\s^{(C)}J^{\Lb}\nonumber
\end{eqnarray}
To express the last two terms we note that
\begin{eqnarray*}
&&\frac{1}{\sqrt{\mbox{det}\sh}}\left(\frac{\partial\sqrt{\mbox{det}\sh}}{\partial\ub}
-\frac{\partial}{\partial\vartheta^A}(\sqrt{\mbox{det}\sh}b^A)\right)\\
&&\hspace{15mm}=\frac{1}{2}(\sh^{-1})^{BC}\left(\frac{\partial\sh_{BC}}{\partial\ub}
-b^A\frac{\partial\sh_{BC}}{\partial\vartheta^A}\right)-\frac{\partial b^A}{\partial\vartheta^A}\\
&&\hspace{15mm}=\frac{1}{2}(\sh^{-1})^{BC}L\sh_{BC}-\frac{\partial b^A}{\partial\vartheta^A}
\end{eqnarray*}
On the other hand, since
$$[L,\Omega_B]=[-b^A\Omega_A,\Omega_B]=(\Omega_B b^A)\Omega_A,$$
we have:
$$\sh([L,\Omega_B],\Omega_C)+\sh(\Omega_B,[L,\Omega_C])=\sh_{AC}\Omega_B b^A+\sh_{AB}\Omega_C b^A$$
Hence:
$$\frac{1}{2}(\sh^{-1})^{BC}\left(\sh([L,\Omega_B],\Omega_C)+\sh(\Omega_B,[L,\Omega_C])\right)
=\Omega_A b^A=\frac{\partial b^A}{\partial\vartheta^A}$$
Since
$$(\sL_L\sh)_{BC}=L\sh_{BC}-\sh([L,\Omega_B],\Omega_C)-\sh(\Omega_B,[L,\Omega_C]),$$
we then conclude that:
$$\frac{1}{2}(\sh^{-1})^{BC}(L\sh_{BC})-\frac{\partial b^A}{\partial\vartheta^A}
=\frac{1}{2}(\sh^{-1})^{BC}(\sL_L\sh)_{BC}=(\sh^{-1})^{BC}\chi_{BC}=\mbox{tr}\chi$$
(see Proposition 3.1). Consequently:
\begin{equation}
\frac{1}{\sqrt{\mbox{det}\sh}}\left(\frac{\partial\sqrt{\mbox{det}\sh}}{\partial\ub}
-\frac{\partial}{\partial\vartheta^A}(\sqrt{\mbox{det}\sh}b^A)\right)=\mbox{tr}\chi
\label{8.40}
\end{equation}
Similarly we deduce:
\begin{equation}
\frac{1}{\sqrt{\mbox{det}\sh}}\left(\frac{\partial\sqrt{\mbox{det}\sh}}{\partial u}
+\frac{\partial}{\partial\vartheta^A}(\sqrt{\mbox{det}\sh}b^A)\right)=\mbox{tr}\chib
\label{8.41}
\end{equation}
In view of \ref{8.40} and \ref{8.41} the expression \ref{8.39} simplifies to:
\begin{eqnarray}
&&\s^{(C)}\sigma=\Omega^{-(n-1)/2}\left\{L(\Omega^{(n+1)/2}a\s^{(C)}J^L)
+\Lb(\Omega^{(n+1)/2}a\s^{(C)}J^{\Lb})\right.\nonumber\\
&&\hspace{35mm}\left.+\sdiv(\Omega^{(n+1)/2}a\s^{(C)}\sJ)\right\}\nonumber\\
&&\hspace{25mm}+\Omega\left(\mbox{tr}\chi a\s^{(C)}J^L+\mbox{tr}\chib a\s^{(C)}J^{\Lb}\right)
\label{8.42}
\end{eqnarray}

We shall now express the components of $a\s^{(C)}J$ in terms of the components of the deformation 
tensor $\s^{(C)}\tilde{\pi}$ in the $(L,\Lb,\Omega_A : A=1,...,n-1)$ frame. In analogy with \ref{7.72}, 
\ref{7.73} with the conformal acoustical metric $\tilde{h}$ in the role of the acoustical metric $h$ we define the $S$ 1-forms $\s^{(C)}\tilde{\spi}_L$, $\s^{(C)}\tilde{\spi}_{\Lb}$ by:
\begin{equation}
\s^{(C)}\tilde{\spi}_L(\Omega_A)=\s^{(C)}\tilde{\pi}(L,\Omega_A) , \ \ \ 
\s^{(C)}\tilde{\spi}_{\Lb}(\Omega_A)=\s^{(C)}\tilde{\pi}(\Lb,\Omega_A)
\label{8.43}
\end{equation}
and the symmetric 2-covariant $S$ tensorfield $\s^{(C)}\tilde{\sspi}$ by:
\begin{equation}
\s^{(C)}\tilde{\sspi}(\Omega_A,\Omega_B)=\s^{(C)}\tilde{\pi}(\Omega_A,\Omega_B)
\label{8.44}
\end{equation}
Noting then that in the definition of $\s^{(C)}J$ given in the statement of Lemma 8.1 the indices 
in $\s^{(C)}\tilde{\pi}^{\alpha\beta}$ are raised with respect to $\tilde{h}$, we deduce:
\begin{eqnarray}
&&\Omega^2 a\s^{(C)}J^L=\frac{1}{4a}\s^{(C)}\tilde{\pi}_{\Lb\Lb}L\psi
+\frac{1}{4}\mbox{tr}\s^{(C)}\tilde{\sspi}\Lb\psi
-\frac{1}{2}\sh^{-1}(\s^{(C)}\tilde{\spi}_{\Lb},\sd\psi) \label{8.45}\\
&&\Omega^2 a\s^{(C)}J^{\Lb}=\frac{1}{4a}\s^{(C)}\tilde{\pi}_{LL}\Lb\psi
+\frac{1}{4}\mbox{tr}\s^{(C)}\tilde{\sspi}L\psi
-\frac{1}{2}\sh^{-1}(\s^{(C)}\tilde{\spi}_L,\sd\psi) \label{8.46}\\
&&\Omega^2 a\s^{(C)}\sJ=-\frac{1}{2}\s^{(C)}\tilde{\spi}_{\Lb}^\sharp L\psi
-\frac{1}{2}\s^{(C)}\tilde{\spi}_L^\sharp\Lb\psi 
+\frac{1}{2}\s^{(C)}\tilde{\pi}_{L\Lb}(\sd\psi)^\sharp\nonumber\\
&&\hspace{20mm}+a\left(\s^{(C)}\tilde{\sspi}^\sharp-\frac{1}{2}\mbox{tr}\s^{(C)}\tilde{\sspi} I\right)
\cdot (\sd\psi)^\sharp \label{8.47}
\end{eqnarray}
In \ref{8.47} if $\theta$ is any 2-covariant $S$ tensorfield, we define $\theta^\sharp$ to be the 
$T^1_1$ type $S$ tensorfield defined by:
\begin{equation}
\sh(X,\theta^\sharp\cdot Y)=\theta(X,Y) \ \ \mbox{: for any pair of $S$ vectorfields $X,Y$}
\label{8.48}
\end{equation}
Also $I$ stands for the $T^1_1$ type $S$ tensorfield which at any $p\in S_{\ub,u}$ is the identity 
transformation in $T_p S_{\ub,u}$. Note finally that by \ref{7.69}, which holds for any vectorfield 
$X$, we have:
\begin{eqnarray}
&&\s^{(C)}\tilde{\pi}_{LL}=\Omega\s^{(C)}\pi_{LL}, \ \ 
\s^{(C)}\tilde{\pi}_{\Lb\Lb}=\Omega\s^{(C)}\pi_{\Lb\Lb}\nonumber\\
&&\hspace{8mm}\s^{(C)}\tilde{\pi}_{L\Lb}=\Omega(\s^{(C)}\pi_{L\Lb}-2a\Omega^{-1}C\Omega)\nonumber\\
&&\hspace{5mm}\s^{(C)}\tilde{\spi}_L=\Omega\s^{(C)}\spi_L, \ \ 
\s^{(C)}\tilde{\spi}_{\Lb}=\Omega\s^{(C)}\spi_{\Lb}\nonumber\\
&&\hspace{10mm}\s^{(C)}\tilde{\sspi}=\Omega(\s^{(C)}\sspi+\Omega^{-1}(C\Omega)\sh)
\label{8.49}
\end{eqnarray}

\vspace{5mm}

\section{The Deformation Tensors of the Commutation Fields}

Consider now the deformation tensors $\s^{(C)}\pi$ of the commutation fields $C$. For $C=T$ from 
tables \ref{7.74} and \ref{7.75} we obtain:
\begin{eqnarray}
&&\s^{(T)}\pi_{LL}=\s^{(T)}\pi_{\Lb\Lb}=0, \ \ \s^{(T)}\pi_{L\Lb}=-2Ta\nonumber\\
&&\s^{(T)}\spi_L=-\zeta, \ \ \s^{(T)}\spi_{\Lb}=\zeta \ \ 
\mbox{: where $\zeta E=Z$}\nonumber\\
&&\hspace{10mm}\s^{(T)}\sspi=2(\chi+\chib)
\label{8.50}
\end{eqnarray}

In the case $n=2$ the other commutation field is $E$. From \ref{3.a13} we deduce, 
in view of \ref{3.a12},
\begin{eqnarray}
&&\s^{(E)}\pi_{LL}=\s^{(E)}\pi_{\Lb\Lb}=0, \ \ \s^{(E)}\pi_{L\Lb}=-2Ea\nonumber\\
&&\hspace{5mm}\s^{(E)}\spi_L=-\chi, \ \ \s^{(E)}\spi_{\Lb}=-\chib\nonumber\\
&&\hspace{10mm}\s^{(E)}\sspi=0
\label{8.51}
\end{eqnarray}

For $n>2$ the remaining commutation fields are the $E_{(\mu)} : \mu=0,...,n$. We shall presently 
determine the deformation tensors of these. First, we apply the identity \ref{7.70} with $E_{(\mu)}$ in 
the role of $X$. Setting $(Y,Z)=(L,L)$, we obtain, since $E_{(\mu)}$ is a $S$ vectorfield, 
\begin{equation}
\s^{(E_{(\mu)})}\pi(L,L)=2h(D_L E_{(\mu)},L)=-2h(E_{(\mu)},D_L L)=0
\label{8.52}
\end{equation}
the integral curves of $L$ being geodesics. Similarly,
\begin{equation}
\s^{(E_{(\mu)})}\pi(\Lb,\Lb)=0
\label{8.53}
\end{equation}
Setting next $(Y,Z)=(L,\Lb)$ in \ref{7.70} we have, since $E_{(\mu)}$ is a $S$ vectorfield,
\begin{eqnarray}
&&\s^{(E_{(\mu)})}\pi(L,\Lb)=h(D_L E_{(\mu)},\Lb)+h(D_{\Lb}E_{(\mu)},L)\nonumber\\
&&\hspace{22mm}=-h(E_{(\mu)},D_L\Lb)-h(E_{(\mu)},D_{\Lb}L)\nonumber\\
&&\hspace{22mm}=-2(\eta+\etab)\cdot E_{(\mu)}=-2E_{(\mu)}a
\label{8.54}
\end{eqnarray}
by \ref{3.14}, \ref{3.15} and \ref{3.17}, \ref{3.18}. Note that the above depend only on the fact that 
$E_{(\mu)}$ is a $S$ vectorfield. To determine $\s^{(E_{(\mu)})}\spi_L$, $\s^{(E_{(\mu)})}\spi_L$ 
we proceed in a different manner. We have: 
\begin{eqnarray}
&&\s^{(E_{(\mu)})}\pi(\Omega_A,L)=({\cal L}_{E_{(\mu)}}h)(\Omega_A,L) \nonumber\\
&&\hspace{18mm}=E_{(\mu)}(h(\Omega_A,L))-h([E_{(\mu)},\Omega_A],L]-h(\Omega_A,[E_{(\mu)},L]) \nonumber\\
&&\hspace{23mm}=h(\Omega_A,[L,E_{(\mu)}])
\label{8.b1}
\end{eqnarray}
since $\Omega_A$ as well as $[E_{(\mu)},\Omega_A]$ are $S$ vectorfields, hence $h$-orthogonal to $L$. 
Thus, we can write:
\begin{equation}
\s^{(E_{(\mu)})}\spi_L^\sharp=[L,E_{(\mu)}]
\label{8.b2}
\end{equation}

Now, we can expand:
\begin{equation}
E_{(\mu)}=E_{(\mu)}^A\Omega_A
\label{8.71}
\end{equation}
Taking the $h$ - inner product with $\Omega_B$, in view of the definition \ref{6.52} 
(see \ref{6.21}) the left hand side becomes:
$$h(\partial/\partial x^\mu,\Omega_B)=h_{\mu\nu}\Omega_B^\nu$$
while the right hand side becomes $E_{(\mu)}^A\sh_{AB}$. It follows that the coefficients in \ref{8.71} 
are given by:
\begin{equation}
E_{(\mu)}^A=h_{\mu\nu}(\sh^{-1})^{AB}\Omega_B^\nu
\label{8.72}
\end{equation}
As a consequence, \ref{8.71} takes the form:
\begin{equation}
E_{(\mu)}=h_{\mu\nu}(\sd x^\nu)^\sharp
\label{8.a1}
\end{equation}
which constitutes an alternative expression for $E_{(\mu)}$. In view of this expression we have:
\begin{eqnarray}
&&[L,E_{(\mu)}]=h_{\mu\nu}(\sL_L(\sh^{-1}))\cdot\sd x^\nu+h_{\mu\nu}\sh^{-1}\cdot\sd Lx^\nu\nonumber\\
&&\hspace{20mm}+(Lh_{\mu\nu})(\sd x^\nu)^\sharp
\label{8.b3}
\end{eqnarray}
By Proposition 3.1, 
\begin{equation}
h_{\mu\nu}\sL_L(\sh^{-1})\cdot\sd x^\nu=-2\chi^\sharp\cdot E_{(\mu)}
\label{8.b4}
\end{equation}
By the 1st of \ref{6.47}, 
\begin{eqnarray}
&&\sd(Lx^\nu)=\sd(\rho N^\nu)=\rho\sd N^\nu+N^\nu\sd\rho \label{8.b5}\\
&&\hspace{12mm}=\rho(\sk\cdot\sd x^\nu+k N^\nu+\ok\Nb^\nu)+N^\nu c^{-1}
(\sd\lambdab-c^{-1}\lambdab \sd c) \nonumber
\end{eqnarray}
Also, 
\begin{equation}
(Lh_{\mu\nu})(\sd x^\nu)^\sharp=\beta_\mu(\sbeta^\sharp LH+\rho H\ss_N^\sharp)+\sbeta^\sharp L\beta_\mu 
\label{8.b6}
\end{equation}
We substitute in \ref{8.b5} for $\sd c$ from \ref{3.84} for $\sd\lambdab$ in terms of $\tetab$ from 
Proposition 3.2 and for $\sk$ in terms of $\tchi$ from \ref{3.76}. Substituting the result as well as \ref{8.b4}, \ref{8.b6} in \ref{8.b3}, and taking into account the fact that by \ref{6.53} we
can express 
\begin{equation}
\sbeta\cdot E_{(\mu)}=\beta_\mu+\frac{h_{\mu\nu}}{2c}(\beta_{\Nb}N^\nu+\beta_N\Nb^\nu)
\label{8.b7}
\end{equation}
we deduce (see \ref{8.b2}):
\begin{eqnarray}
&&\s^{(E_{(\mu)})}\spi_L=-\chi\cdot E_{(\mu)}+\frac{1}{c}h_{\mu\nu}N^\nu\tetab+H(L\beta_\mu)\sbeta 
\label{8.85}\\
&&\hspace{10mm}+\frac{1}{2}\left\{\left(\beta_\mu-\frac{h_{\mu\nu}\Nb^\nu}{2c}\beta_N\right)LH
+\rho\beta_N E_{(\mu)}H\right\}\sbeta-\frac{1}{2}\rho\beta_N\beta_\mu\sd H \nonumber
\end{eqnarray}
Similarly, 
\begin{eqnarray}
&&\s^{(E_{(\mu)})}\spi_{\Lb}=-\chib\cdot E_{(\mu)}+\frac{1}{c}h_{\mu\nu}\Nb^\nu\teta
+H(\Lb\beta_\mu)\sbeta \label{8.86}\\
&&\hspace{10mm}+\frac{1}{2}\left\{\left(\beta_\mu-
\frac{h_{\mu\nu}N^\nu}{2c}\beta_{\Nb}\right)\Lb H
+\rhob\beta_{\Nb} E_{(\mu)}H\right\}\sbeta-\frac{1}{2}\rhob\beta_{\Nb}\beta_\mu\sd H \nonumber 
\end{eqnarray}

Finally (setting $(Y,Z)=(\Omega_A,\Omega_B)$ in \ref{7.70}) we have:
\begin{eqnarray}
&&\s^{(E_{(\mu)})}\pi(\Omega_A,\Omega_B)=h(D_{\Omega_A} E_{(\mu)},\Omega_B)+
h(D_{\Omega_B} E_{(\mu)},\Omega_A)\nonumber\\
&&\hspace{25mm}=\sh(\sD_{\Omega_A}E_{(\mu)},\Omega_B)+\sh(\sD_{\Omega_B}E_{(\mu)},\Omega_A)
\label{8.b8}
\end{eqnarray}
In view of the expression \ref{8.a1} we have:
\begin{equation}
\sh(\sD_{\Omega_A}E_{(\mu)},\Omega_B)=(\sd h_{\mu\nu}\cdot\Omega_A)(\sd x^\nu\cdot\Omega_B)
+h_{\mu\nu}(\sD^2 x^\nu)(\Omega_A,\Omega_B)
\label{8.b9}
\end{equation}
According to \ref{3.105}:
\begin{equation}
\sD^2 x^\nu=\frac{1}{2c}(\tchib N^\nu+\tchi\Nb^\nu)-\tgamma^\nu
\label{8.b10}
\end{equation}
and according to \ref{3.101}:
\begin{equation}
\tgamma^\nu=\sgamma^C\Omega_C^\nu+\tgamma^N N^\nu+\tgamma^{\Nb}\Nb^\nu
\label{8.b11}
\end{equation}
From \ref{3.98} in view of \ref{8.a1} we obtain:
\begin{eqnarray}
&&h_{\mu\nu}\sgamma^C\Omega_C^\nu=-\frac{1}{2}\sbeta\otimes\sbeta E_{(\mu)}H \label{8.b12}\\
&&\hspace{20mm}+\frac{1}{2}(\sbeta\cdot E_{(\mu)})(\sbeta\otimes\sd H+\sd H\otimes\sbeta+2H\sss)
\nonumber
\end{eqnarray}
Also, according to \ref{3.100}:
\begin{eqnarray}
&&\tgamma^N=-\frac{1}{4c}\beta_{\Nb}(\sbeta\otimes\sd H+\sd H\otimes\sbeta+2 H\sss)\nonumber\\
&&\tgamma^{\Nb}=-\frac{1}{4c}\beta_N(\sbeta\otimes\sd H+\sd H\otimes\sbeta+2 H\sss) \label{8.b13}
\end{eqnarray}
Taking into account \ref{8.b7} we then conclude that:
\begin{eqnarray}
&&h_{\mu\nu}\sD^2 x^\nu=\frac{h_{\mu\nu}}{2c}(\tchib N^\nu+\tchi\Nb^\nu)
+\frac{1}{2}\sbeta\otimes\sbeta E_{(\mu)}H \nonumber\\
&&\hspace{15mm}-\frac{1}{2}\beta_\mu(\sbeta\otimes\sd H+\sd H\otimes\sbeta+2H\sss) 
\label{8.b14}
\end{eqnarray}
Also, we have:
$$\sd h_{\mu\nu}=\beta_\mu\beta_\nu\sd H+H(\beta_mu\sd\beta_\nu+\beta_\nu\sd\beta_\mu)$$
and (see \ref{2.88}):
$$\sd x^\nu\otimes\sd\beta_\nu=\sss=\sd\beta_\nu\otimes\sd\beta^\nu$$
Hence:
\begin{eqnarray}
&&\sd h_{\mu\nu}\otimes\sd x^\nu=\beta_\mu(\sd H\otimes\sbeta+H\sss)+H\sd\beta_\mu\otimes\sbeta 
\nonumber\\
&&\sd x^\nu\otimes\sd h_{\mu\nu}=\beta_\mu(\sbeta\otimes \sd H+H\sss)+H\sbeta\otimes\sd\beta_\mu 
\label{8.b15}
\end{eqnarray}
Combining \ref{8.b14}, \ref{8.b15} we conclude through \ref{8.b8}, \ref{8.b9} that:
\begin{eqnarray}
&&\s^{(E_{(\mu)})}\sspi=\frac{1}{c}h_{\mu\nu}\Nb^\nu\tchi+\frac{1}{c}h_{\mu\nu}N^\nu\tchib\nonumber\\
&&\hspace{15mm}+(E_{(\mu)}H)\sbeta\otimes\sbeta+H(\sbeta\otimes\sd\beta_\mu+\sd\beta_\mu\otimes\sbeta)
\label{8.89}
\end{eqnarray}
In conclusion, equations \ref{8.52}, \ref{8.53}, \ref{8.54}, \ref{8.85}, \ref{8.86}, \ref{8.89} 
give all components of $\s^{(E_{(\mu)})}\pi$. 

We note that the definition \ref{8.a1} implies the following commutation relation:
\begin{eqnarray}
&&[E_{(\mu)},E_{(\nu)}]f=h_{\mu\kappa}(\sd x^\kappa,\sd h_{\nu\lambda})_{\sh}
(\sd x^\lambda,\sd f)_{\sh}-h_{\nu\lambda}(\sd x^\lambda,\sd h_{\mu\kappa})_{\sh}
(\sd x^\kappa,\sd f)_{\sh} \nonumber\\
&&\hspace{25mm}+\left\{\sD^2 x^\lambda\cdot(\sd x^\kappa)^\sharp
-\sD^2 x^\kappa\cdot(\sd x^\lambda)^\sharp\right\}\cdot(\sd f)^\sharp 
\label{8.c1}
\end{eqnarray}
for an arbitrary function $f$ on ${\cal N}$. 

Now, from \ref{8.50} - \ref{8.53} we have:
\begin{equation}
\s^{(C)}\pi_{LL}=\s^{(C)}\pi_{\Lb\Lb}=0 \ \ \mbox{: for all commutation fields $C$}
\label{8.a2}
\end{equation}
We see then (in view of \ref{8.49}) from \ref{8.45}, \ref{8.46} that $a\s^{(C)}J^L$, $a\s^{(C)}J^{\Lb}$ 
are, like $a\s^{(C)}\sJ$, regular even where $a$ vanishes, and this holds for all commutation fields 
$C$. 

Since by \ref{8.49}:
\begin{equation}
\tilde{\mbox{tr}}\s^{(C)}\tilde{\pi}=-a^{-1}\s^{(C)}\pi_{L\Lb}+\mbox{tr}\s^{(C)}\sspi 
+(n+1)\Omega^{-1}C\Omega
\label{8.a3}
\end{equation}
the functions $\s^{(C)}\delta$, defined by \ref{8.27}, are expressed as:
\begin{equation}
\s^{(C)}\delta=-\frac{1}{2a}(\s^{(C)}\pi_{L\Lb}+2Ca)+\frac{1}{2}\mbox{tr}\s^{(C)}\sspi
+\frac{(n-1)}{2}\Omega^{-1}C\Omega
\label{8.a4}
\end{equation}
From \ref{8.50}, \ref{8.51}, \ref{8.54} we see that:
\begin{equation}
\s^{(C)}\pi_{L\Lb}+2Ca=0 \ \ \mbox{: for all commutation fields $C$}
\label{8.a5}
\end{equation}
Consequently the functions $\s^{(C)}\delta$ are regular even where $a$ vanishes, for all commutation 
fields $C$. It then follows that the recursion formulas \ref{8.28}, \ref{8.29} are likewise regular 
even where $a$ vanishes.

\vspace{5mm}

\section{The Principal Acoustical Error Terms}

We have seen in Chapter 3 that the irreducible 1st order acoustical quantities are  
$\tchi$, $\tchib$ and the first derivatives of $\lambda$, $\lambdab$. However of the last 
$L\lambda$, $\Lb\lambdab$ are not 1st order acoustical because they are expressed by Proposition 3.3 
as 1st order quantities not containing 1st order acoustical terms. We thus take the irreducible 
1st order acoustical quantities to be:
\begin{equation}
\tchi, \ \ \tchib, \ \ \sd\lambda, \ \ \sd\lambdab, \ \ T\lambda, \ \ T\lambdab
\label{8.90}
\end{equation}
All other 1st order acoustical quantities can be expressed in terms of these. In the following 
we shall denote by $[\s\s]_{P.A.}$ the principal acoustical part, that is for a quantity of order 
$k$ the part of the quantity containing the $k$th order acoustical terms. Consider 
the 1st order quantity $\sd c$. From \ref{3.84} we see that:
\begin{equation}
[\sd c]_{P.A.}=[c(k+\okb)]_{P.A.}
\label{8.91}
\end{equation}
In view of the fact that by \ref{3.53}, \ref{3.54} the principal acoustical parts of $\ok$, $\kb$ 
vanish, we obtain from \ref{3.46}:
\begin{equation}
[k]_{P.A.}=-[\sk\cdot\pi]_{P.A.}, \ \ \ [\okb]_{P.A.}=-[\skb\cdot\pi]_{P.A.}
\label{8.92}
\end{equation}
Now, by \ref{3.76}, \ref{3.77}:
\begin{equation}
[\sk]_{P.A.}=\tchi^\sharp, \ \ \ [\skb]_{P.A.}=\tchib^\sharp
\label{8.93}
\end{equation}
Hence:
\begin{equation}
[k]_{P.A.}=-\tchi\cdot\pi^\sharp, \ \ \ [\okb]_{P.A.}=-\tchib\cdot\pi^\sharp
\label{8.94}
\end{equation}
Substituting in \ref{8.91} we obtain:
\begin{equation}
[\sd c]_{P.A.}=-c(\tchi+\tchib)\cdot\pi^\sharp
\label{8.95}
\end{equation}
Consider the 1st order acoustical quantity $Tc$. From \ref{3.85}, \ref{3.86}, we see that, in view of 
the fact that by \ref{3.47} in conjunction with \ref{3.57}, \ref{3.58}, \ref{3.68}, \ref{3.69} the 
principal acoustical parts of $m$, $\omb$ vanish, 
\begin{equation}
[Tc]_{P.A.}=[c(n+\onb)]_{P.A.}
\label{8.96}
\end{equation}
Also, in view of the fact that by \ref{3.61}, \ref{3.62} the principal acoustical parts of $\on$, 
$\nb$ vanish, we obtain from \ref{3.48}:
\begin{equation}
[n]_{P.A.}=-[\sn\cdot\pi]_{P.A.}, \ \ \ [\onb]_{P.A.}=-[\snb\cdot\pi]_{P.A.}
\label{8.97}
\end{equation}
Now, by \ref{3.a6}, \ref{3.a7}, in view of \ref{8.94}:
\begin{equation}
[\sn]_{P.A.}=2(\sd\lambda+\lambda\tchi\cdot\pi^\sharp)^\sharp, \ \ \ [\snb]_{P.A.}=2(\sd\lambdab+\lambdab\tchib\cdot\pi^\sharp)^\sharp
\label{8.98}
\end{equation}
Hence:
\begin{equation}
[n]_{P.A.}=-2(\sd\lambda+\lambda\tchi\cdot\pi^\sharp)\cdot\pi^\sharp, \ \ \ [\onb]_{P.A.}=-2(\sd\lambdab+\lambdab\tchib\cdot\pi^\sharp)\cdot\pi^\sharp
\label{8.99}
\end{equation}
Substituting in \ref{8.96} we obtain:
\begin{equation}
[Tc]_{P.A.}=-2c\left(\sd\lambda+\sd\lambdab+(\lambda\tchi+\lambdab\tchib)\cdot\pi^\sharp\right)
\cdot\pi^\sharp
\label{8.100}
\end{equation} 
Consider next the 1st order quantities $\sd a$, $Ta$. Writing 
\begin{equation} 
a=c^{-1}\lambda\lambdab
\label{8.101}
\end{equation}
(see \ref{2.72}, \ref{2.74}) in view of \ref{8.95}, \ref{8.100} we obtain:
\begin{eqnarray}
&&[\sd a]_{P.A.}=\rho\sd\lambda+\rhob\sd\lambdab+a(\tchi+\tchib)\cdot\pi^\sharp\nonumber\\
&&[Ta]_{P.A.}=\rho T\lambda+\rhob T\lambdab+2a\left(\sd\lambda+\sd\lambdab+
(\lambda\tchi+\lambdab\tchib)\cdot\pi^\sharp\right)\cdot\pi^\sharp \label{8.102}
\end{eqnarray}
Consider next the 1st order quantities $\teta$, $\tetab$. By \ref{3.80}, \ref{3.81} together 
with \ref{8.98}:
\begin{equation}
[\teta]_{P.A.}=\sd\lambda+\lambda\tchi\cdot\pi^\sharp, \ \ \ [\tetab]_{P.A.}=\sd\lambdab+\lambdab\tchib\cdot\pi^\sharp
\label{8.103}
\end{equation}
Since by \ref{3.79}, \ref{3.82}:
\begin{equation}
[\eta]_{P.A.}=\rho[\teta]_{P.A.}, \ \ \ [\etab]_{P.A.}=\rhob[\tetab]_{P.A.}
\label{8.104}
\end{equation}
it then follows by \ref{3.16} that:
\begin{equation}
[Z]_{P.A.}=2\left(\rho\sd\lambda-\rhob\sd\lambdab+a(\tchi-\tchib)\cdot\pi^\sharp\right)^\sharp
\label{8.105}
\end{equation}
Note also that by \ref{3.75}, \ref{3.78}:
\begin{equation}
[\chi]_{P.A.}=\rho\tchi, \ \ \ [\chib]_{P.A.}=\rhob\tchib
\label{8.106}
\end{equation}

Now, the deformation tensors of the commutation fields are also 1st order quantities. Let us determine 
the principal acoustical parts. In view of \ref{8.a2}, \ref{8.a5} we only need to to consider the 
components $\s^{(C)}\spi_L$, $\s^{(C)}\spi_{\Lb}$, $\s^{(C)}\sspi$, for each commutation field $C$. 
From \ref{8.50}, in view of \ref{8.105}, \ref{8.106} we obtain:
\begin{eqnarray}
&&[\s^{(T)}\spi_L]_{P.A.}=-2\left(\rho\sd\lambda-\rhob\sd\lambdab
+a(\tchi-\tchib)\cdot\pi^\sharp\right)
\nonumber\\
&&[\s^{(T)}\spi_{\Lb}]_{P.A.}=2\left(\rho\sd\lambda-\rhob\sd\lambdab
+a(\tchi-\tchib)\cdot\pi^\sharp\right)
\nonumber\\
&&[\s^{(T)}\sspi]_{P.A.}=2(\rho\tchi+\rhob\tchib)
\label{8.107}
\end{eqnarray}
In the case $n=2$ from \ref{8.51} we obtain: 
\begin{equation}
[\s^{(E)}\spi_L]_{P.A.}=-\rho\tchi, \ \ \ [\s^{(E)}\spi_{\Lb}]_{P.A.}=-\rhob\tchib, \ \ \ 
[\s^{(E)}\sspi]_{P.A.}=0
\label{8.108}
\end{equation}
For $n>2$ from \ref{8.85}, \ref{8.86}, \ref{8.89}, in view of \ref{8.103}, \ref{8.106} we obtain:
\begin{eqnarray}
&&[\s^{(E_{(\mu)})}\spi_L]_{P.A.}=-\rho\tchi\cdot E_{(\mu)}+\frac{1}{c}h_{\mu\nu}N^\nu
(\sd\lambdab+\lambdab\tchib\cdot\pi^\sharp) \nonumber\\
&&[\s^{(E_{(\mu)})}\spi_{\Lb}]_{P.A.}=-\rhob\tchib\cdot E_{(\mu)}+\frac{1}{c}h_{\mu\nu}\Nb^\nu 
(\sd\lambda+\lambda\tchi\cdot\pi^\sharp) \nonumber\\
&&[\s^{(E_{(\mu)})}\sspi]_{P.A.}=\frac{1}{c}h_{\mu\nu}\Nb^\nu\tchi+\frac{1}{c}h_{\mu\nu}N^\nu\tchib
\label{8.109}
\end{eqnarray}

We now consider the rescaled sources \ref{8.11}, \ref{8.12} for $m+l=1$. In the case $n=2$, 
by \ref{8.30} these are given by \ref{8.28}, \ref{8.29} with $(m,l)=(0,0)$ as:
\begin{equation}
\s^{(1,0)}\tilde{\rho}_\mu=\s^{(T,0,0)}\sigma_\mu, \ \ \ 
\s^{(0,1)}\tilde{\rho}_\mu=\s^{(E,0,0)}\sigma_\mu
\label{8.110}
\end{equation}
For $n>2$, by \ref{8.33} these are similarly given by \ref{8.31}, \ref{8.32} with $(m,l)=(0,0)$ as:
\begin{equation}
\s^{(1)}\tilde{\rho}_\mu=\s^{(T,0)}\sigma_\mu, \ \ \ 
\s^{(0,\nu)}\tilde{\rho}_\mu=\s^{(E_{(\nu)},0)}\sigma_\mu 
\label{8.111}
\end{equation}
The right hand sides in \ref{8.110}, \ref{8.111} correspond to $\s^{(C)}\sigma$, defined by \ref{8.26} 
$C$ being $T,E$ in the case $n=2$, $T,E_{(\nu)}:\nu=0,...,n$ for $n>2$ and with the fundamental 
variations $\beta_\mu$ in the role of $\psi$. Now $\s^{(C)}\sigma$ is expressed by \ref{8.42}, 
the components of $\s^{(C)}J$ being expressed by \ref{8.45} - \ref{8.47}. So here $\s^{(C)}J$ is 
a quantity of order 1 and $\s^{(C)}\sigma$ is a quantity of order 2. Consequently the rescaled sources 
are for $m+l=1$ quantities of order 2. Let us determine their principal acoustical parts. 

From \ref{8.45} - \ref{8.47} and \ref{8.a5} with $C=T$, together with \ref{8.107} and the 2nd of 
\ref{8.102} we find:
\begin{eqnarray}
&&[\Omega a\s^{(T)}J^L]_{P.A.}=\frac{1}{2}(\rho\mbox{tr}\tchi+\rhob\mbox{tr}\tchib)\Lb\psi \label{8.112}\\
&&\hspace{20mm}-\left(\rho\sd\lambda-\rhob\sd\lambdab+a(\tchi-\tchib)\cdot\pi^\sharp\right)\cdot 
(\sd\psi)^\sharp\nonumber
\end{eqnarray}
\begin{eqnarray}
&&[\Omega a\s^{(T)}J^{\Lb}]_{P.A.}=\frac{1}{2}(\rho\mbox{tr}\tchi+\rhob\mbox{tr}\tchib)L\psi \label{8.113}\\
&&\hspace{20mm}+\left(\rho\sd\lambda-\rhob\sd\lambdab+a(\tchi-\tchib)\cdot\pi^\sharp\right)\cdot 
(\sd\psi)^\sharp\nonumber
\end{eqnarray}
\begin{eqnarray}
&&[\Omega a\s^{(T)}\sJ]_{P.A.}=
\left(\rho\sd\lambda-\rhob\sd\lambdab+a(\tchi-\tchib)\cdot\pi^\sharp\right)^\sharp(\Lb\psi-L\psi)
\label{8.114}\\
&&\hspace{20mm}+a\left(2(\rho\tchi+\rhob\tchib)^\sharp
-(\rho\mbox{tr}\tchi+\rhob\mbox{tr}\tchib)I\right)\cdot(\sd\psi)^\sharp\nonumber\\
&&\hspace{20mm}-\left\{(\rho T\lambda+\rhob T\lambdab)
+2a\left(\sd\lambda+\sd\lambdab+(\lambda\tchi+\lambdab\tchib)\cdot\pi^\sharp\right)\cdot\pi^\sharp
\right\}(\sd\psi)^\sharp\nonumber
\end{eqnarray}

In the case $n=2$, from \ref{8.45} - \ref{8.47} and \ref{8.a5} with $C=E$, together with \ref{8.108} 
and the 1st of \ref{8.102} we find: 
\begin{eqnarray}
&&[\Omega a\s^{(E)}J^L]_{P.A.}=\frac{1}{2}\rhob\tchib E\psi \label{8.115}\\
&&[\Omega a\s^{(E)}J^{\Lb}]_{P.A.}=\frac{1}{2}\rho\tchi E\psi \label{8.116}\\
&&[\Omega a\s^{(E)}\sJ]_{P.A.}=\frac{1}{2}\rhob\tchib L\psi+\frac{1}{2}\rho\tchi\Lb\psi \label{8.117}\\
&&\hspace{20mm}-\left(\rho\sd\lambda+\rhob\sd\lambdab+a(\tchi+\tchib)\pi\right)E\psi\nonumber
\end{eqnarray}
recalling from Section 3.5 that in the case $n=2$, for any function $f$ on ${\cal N}$, $\sd f$ 
is $Ef$ (in particular we have \ref{3.a19}). 

For $n>2$, from \ref{8.45} - \ref{8.47} and \ref{8.a5} with $C=E_{(\nu)}$, together with \ref{8.109} 
and the 1st of \ref{8.102} we find:
\begin{eqnarray}
&&[\Omega a\s^{(E_{(\nu)})}J^L]_{P.A.}=\frac{h_{\nu\kappa}}{4c}(\Nb^\kappa\mbox{tr}\tchi
+N^\kappa\mbox{tr}\tchib) \label{8.118}\\
&&\hspace{25mm}+\frac{1}{2}\left(\rhob\tchib\cdot E_{(\nu)}-\frac{h_{\nu\kappa}}{c}\Nb^\kappa
(\sd\lambda+\lambda\tchi\cdot\pi^\sharp)\right)\cdot(\sd\psi)^\sharp\nonumber
\end{eqnarray}
\begin{eqnarray}
&&[\Omega a\s^{(E_{(\nu)})}J^{\Lb}]_{P.A.}=\frac{h_{\nu\kappa}}{4c}(\Nb^\kappa\mbox{tr}\tchi
+N^\kappa\mbox{tr}\tchib) \label{8.119}\\
&&\hspace{25mm}+\frac{1}{2}\left(\rho\tchi\cdot E_{(\nu)}-\frac{h_{\nu\kappa}}{c}N^\kappa
(\sd\lambdab+\lambdab\tchib\cdot\pi^\sharp)\right)\cdot(\sd\psi)^\sharp\nonumber
\end{eqnarray}
\begin{eqnarray}
&&[\Omega a\s^{(E_{(\nu)})}\sJ]_{P.A.}=\frac{1}{2}
\left(\rhob\tchib\cdot E_{(\nu)}-\frac{h_{\nu\kappa}}{c}\Nb^\kappa
(\sd\lambda+\lambda\tchi\cdot\pi^\sharp)\right)^\sharp L\psi \label{8.120}\\
&&\hspace{25mm}+\frac{1}{2}\left(\rho\tchi\cdot E_{(\nu)}-\frac{h_{\nu\kappa}}{c}N^\kappa
(\sd\lambdab+\lambdab\tchib\cdot\pi^\sharp)\right)^\sharp\Lb\psi\nonumber\\
&&\hspace{25mm}-\left(\rho E_{(\nu)}\lambda+\rhob E_{(\nu)}\lambdab
+a(\tchi+\tchib)(E_{(\nu)},\pi^\sharp)\right)(\sd\psi)^\sharp\nonumber\\
&&\hspace{25mm}+\frac{ah_{\nu\kappa}}{c}\left(\Nb^\kappa\left(\tchi^\sharp-\frac{1}{2}\mbox{tr}\tchi I\right)+N^\kappa\left(\tchib^\sharp-\frac{1}{2}\mbox{tr}\tchib I\right)\right)\cdot(\sd\psi)^\sharp 
\nonumber
\end{eqnarray}

Consider now the 1st of \ref{8.111}. This reduces in the case $n=2$ to the 1st of \ref{8.110}, the 
difference being only in the notation. Here $\s^{(T,0)}\sigma_\mu$ is given by \ref{8.42} with 
$C=T$ and $\psi=\beta_\mu$. In general, we denote the commutation current corresponding to 
$C$ and to $\beta_\mu$ by $\s^{(C)}J_\mu$. By \ref{8.42} we then have:
\begin{equation}
[\s^{(1)}\tilde{\rho}_\mu]_{P.A.}=[L(\Omega a\s^{(T)}J_\mu^L)]_{P.A.}+
[\Lb(\Omega a\s^{(T)}J_\mu^{\Lb})]_{P.A.}+[\sdiv(\Omega a\s^{(T)}\sJ_\mu)]_{P.A.}
\label{8.121}
\end{equation}
Now, the 1st derivatives of $\s^{(C)}J_\mu$ are quantities of 2nd order. In determining their 
principal acoustical parts we must take into account the following. First that, 
as already noted, by Proposition 3.3 the principal acoustical parts of 
$L\lambda$, $\Lb\lambdab$ vanish, hence so do the principal acoustical parts of the 1st derivatives 
of $L\lambda$, $\Lb\lambdab$. It follows that we may take the irreducible 2nd derivatives of 
$\lambda$, $\lambdab$ to be:
\begin{equation}
\sD^2\lambda, \ \sD^2\lambdab, \ \ \sd T\lambda, \ \sd T\lambdab, \ \ T^2\lambda, T^2\lambdab
\label{8.122}
\end{equation}
Second that by Proposition 3.4 the principal acoustical parts of $\sL_L\tchi$, $\sL_{\Lb}\tchib$ 
vanish. Third, that by Proposition 3.5: 
\begin{eqnarray}
&&[\sL_{\Lb}\tchi]_{P.A.}=[\sD\sn\cdot\sh]_{P.A.}=2\sD^2\lambda+2\lambda[\sD\tchi]_{P.A}
\cdot\pi^\sharp \label{8.123}\\
&&[\sL_L\tchib]_{P.A.}=[\sD\snb\cdot\sh]_{P.A.}=2\sD^2\lambdab+2\lambdab[\sD\tchib]_{P.A.}
\cdot\pi^\sharp \label{8.124}
\end{eqnarray}
by \ref{8.98}. In terms of components in the $(\Omega_A:A=1,...,n-1)$ frame field for the $S_{\ub,u}$, 
the $AB$ components of the 2nd terms in \ref{8.123}, \ref{8.124} are:
$$2\lambda[\sD_A\tchi_{BC}]_{P.A.}\pi^C, \ \ \ 2\lambdab[\sD_A\tchib_{BC}]_{P.A.}\pi^C$$
Finally that by Proposition 3.6:
\begin{equation}
[\sD_A\tchi_{BC}-\sD_B\tchi_{AC}]_{P.A.}=0, \ \ \ [\sD_A\tchib_{BC}-\sD_B\tchib_{AC}]_{P.A.}=0
\label{8.125}
\end{equation}
As a consequence the remaining irreducible 2nd order acoustical quantities can be taken to be:
\begin{equation}
\sD\tchi, \ \ \ \sD\tchib 
\label{8.126}
\end{equation}
modulo however \ref{8.125}. Taking the trace of equations \ref{8.125} relative to $AC$ yields:
\begin{equation}
[\sdiv\tchi-\sd\mbox{tr}\tchi]_{P.A.}=0, \ \ \ [\sdiv\tchib-\sd\mbox{tr}\tchib]_{P.A.}=0
\label{8.127}
\end{equation}
Taking the trace of \ref{8.123}, \ref{8.124} and using \ref{8.127} we obtain:
\begin{eqnarray}
&&[\Lb\mbox{tr}\tchi]_{P.A.}=2\slap\lambda+2\lambda\sd\mbox{tr}\tchi\cdot\pi^\sharp \label{8.128}\\
&&[L\mbox{tr}\tchib]_{P.A.}=2\slap\lambdab+2\lambdab\sd\mbox{tr}\tchib\cdot\pi^\sharp \label{8.129}
\end{eqnarray}

Taking account of the above we deduce from \ref{8.112} - \ref{8.114} (noting that $\rho\lambda=\rhob\lambdab=a$):
\begin{eqnarray}
&&[L(\Omega a\s^{(T)}J_\mu^L)]_{P.A.}=(\rhob\slap\lambdab+a\sd\mbox{tr}\tchib\cdot\pi^\sharp)
\Lb\beta_\mu+\rhob\sd T\lambdab\cdot\sd\beta_\mu \nonumber\\
&&\hspace{27mm}+2a\left(\sD^2\lambdab+\lambdab[\sD\tchib]_{P.A.}\cdot\pi^\sharp\right)
(\pi^\sharp,(\sd\beta_\mu)^\sharp) \label{8.130}
\end{eqnarray}
\begin{eqnarray}
&&[\Lb(\Omega a\s^{(T)}J_\mu^{\Lb})]_{P.A.}=(\rho\slap\lambda+a\sd\mbox{tr}\tchi\cdot\pi^\sharp)
L\beta_\mu+\rho\sd T\lambda\cdot\sd\beta_\mu \nonumber\\
&&\hspace{27mm}+2a\left(\sD^2\lambda+\lambda[\sD\tchi]_{P.A.}\cdot\pi^\sharp\right)
(\pi^\sharp,(\sd\beta_\mu)^\sharp) \label{8.131}
\end{eqnarray}
\begin{eqnarray}
&&[\sdiv(\Omega a\s^{(T)}\sJ_\mu)]_{P.A.}=\left(\rho\slap\lambda-\rhob\slap\lambdab
+a(\sd\mbox{tr}\tchi-\sd\mbox{tr}\tchib)\cdot\pi^\sharp\right)(\Lb\beta_\mu-L\beta_\mu)\nonumber\\
&&\hspace{27mm}+a(\rho\sd\mbox{tr}\tchi+\rhob\sd\mbox{tr}\tchib)\cdot(\sd\beta_\mu)^\sharp
-(\rho\sd T\lambda+\rhob\sd T\lambdab)\cdot(\sd\beta_\mu)^\sharp \nonumber\\
&&\hspace{33mm}-2a(\sD^2\lambda+\sD^2\lambdab)(\pi^\sharp,(\sd\beta_\mu)^\sharp) \label{8.132}\\
&&\hspace{30mm}-2a\lambda[\sD_{(\sd\beta_\mu)^\sharp}\tchi]_{P.A.}(\pi^\sharp,\pi^\sharp)
-2a\lambdab[\sD_{(\sd\beta_\mu)^\sharp}\tchib]_{P.A.}(\pi^\sharp,\pi^\sharp)\nonumber
\end{eqnarray} 
In summing \ref{8.130}, \ref{8.131}, \ref{8.132} to obtain $[\s^{(1)}\tilde{\rho}_\mu]_{P.A.}$ 
according to \ref{8.121}, we see that the terms in $\sd T\lambda$, $\sd T\lambdab$ cancel. 
The terms in $\sD^2\lambda$, $\sD^2\lambdab$ likewise cancel. As for the terms in $\sD\tchi$, 
$\sD\tchib$, we note that 
$$[\sD_{(\sd\beta_\mu)^\sharp}\tchi]_{P.A.}(\pi^\sharp,\pi^\sharp)
=(\sd\beta_\mu)^A[\sD_A\tchi_{BC}]_{P.A.}\pi^B\pi^C$$
while
$$([\sD\tchi]_{P.A.}\cdot\pi^\sharp)(\pi^\sharp,(\sd\beta_\mu)^\sharp)
=[\sD_B\tchi_{AC}]_{P.A.}\pi^B\pi^C(\sd\beta_\mu)^A$$
hence by the 1st of \ref{8.125}: 
\begin{equation}
[\sD_{(\sd\beta_\mu)^\sharp}\tchi]_{P.A.}(\pi^\sharp,\pi^\sharp)=
([\sD\tchi]_{P.A.}\cdot\pi^\sharp)(\pi^\sharp,(\sd\beta_\mu)^\sharp)
\label{8.133}
\end{equation}
Similarly, by the 2nd of \ref{8.125}:
\begin{equation}
[\sD_{(\sd\beta_\mu)^\sharp}\tchib]_{P.A.}(\pi^\sharp,\pi^\sharp)=
([\sD\tchib]_{P.A.}\cdot\pi^\sharp)(\pi^\sharp,(\sd\beta_\mu)^\sharp)
\label{8.134}
\end{equation}
We see then that the terms in $[\sD\tchi]_{P.A.}$, $[\sD\tchib]_{P.A.}$ cancel as well. 
What remains is:
\begin{eqnarray}
&&[\s^{(1)}\tilde{\rho}_\mu]_{P.A.}=\rho\slap\lambda \ \Lb\beta_\mu+\rhob\slap\lambdab \ L\beta_\mu 
\label{8.135}\\
&&\hspace{10mm}+a\sd\mbox{tr}\tchi\cdot(\pi^\sharp\Lb\beta_\mu+\rho(\sd\beta_\mu)^\sharp)
+a\sd\mbox{tr}\tchib\cdot(\pi^\sharp L\beta_\mu+\rhob(\sd\beta_\mu)^\sharp)\nonumber
\end{eqnarray}
Note the absence of the 2nd order acoustical quantities $\sd T\lambda$, $\sd T\lambdab$ and 
$T^2\lambda$, $T^2\lambdab$.

Consider now the 2nd of \ref{8.110} (case $n=2$). Here $\s^{(E,0,1)}\sigma_\mu$ is given by \ref{8.42} with 
$C=E$ and $\psi=\beta_\mu$. By \ref{8.42} we then have:
\begin{equation}
[\s^{(0,1)}\tilde{\rho}_\mu]_{P.A.}=[L(\Omega a\s^{(E)}J_\mu^L)]_{P.A.}+
[\Lb(\Omega a\s^{(E)}J_\mu^{\Lb})]_{P.A.}+[\sdiv(\Omega a\s^{(E)}\sJ_\mu)]_{P.A.}
\label{8.136}
\end{equation}
Taking account of \ref{8.123}, \ref{8.124} which in the case $n=2$ reduce to:
\begin{eqnarray}
&&[\Lb\tchi]_{P.A.}=2E^2\lambda+2\lambda\pi E\tchi \label{8.137}\\
&&[L\tchib]_{P.A.}=2E^2\lambdab+2\lambdab\pi E\tchib \label{8.138}
\end{eqnarray}
we deduce from \ref{8.115} - \ref{8.117}:
\begin{eqnarray}
&&[L(\Omega a\s^{(E)}J_\mu^L]_{P.A.}=(\rhob E^2\lambdab+a\pi E\tchib)E\beta_\mu \label{8.139}\\
&&[\Lb(\Omega a\s^{(E)}J_\mu^{\Lb}]_{P.A.}=(\rho E^2\lambda+a\pi E\tchi)E\beta_\mu \label{8.140}\\
&&[\sdiv(\Omega a\s^{(E)}\sJ_\mu]_{P.A.}=\frac{1}{2}\rhob E\tchib \ L\beta_\mu 
+\frac{1}{2}\rho E\tchi \ \Lb\beta_\mu \label{8.141}\\
&&\hspace{30mm}-\left(\rho E^2\lambda+\rhob E^2\lambdab+a\pi(E\tchi+E\tchib)\right)E\beta_\mu \nonumber
\end{eqnarray}
Adding, the terms in $E^2\lambda$, $E^2\lambdab$ cancel, and we obtain:
\begin{equation}
[\s^{(0,1)}\tilde{\rho}_\mu]_{P.A.}=\frac{1}{2}\rhob E\tchib \ L\beta_\mu 
+\frac{1}{2}\rho E\tchi \ \Lb\beta_\mu 
\label{8.142}
\end{equation}

Consider finally the 2nd of \ref{8.111} (case $n>2$). Here $\s^{(E_{(\nu)},0)}\sigma_\mu$ is given by 
\ref{8.42} with $C=E_{(\nu)}$ and $\psi=\beta_\mu$. By \ref{8.42} we then have:
\begin{equation}
[\s^{(0,\nu)}\tilde{\rho}_\mu]_{P.A.}=[L(\Omega a\s^{(E_{(\nu)})}J_\mu^L)]_{P.A.}+
[\Lb(\Omega a\s^{(E_{(\nu)})}J_\mu^{\Lb})]_{P.A.}+[\sdiv(\Omega a\s^{(E_{(\nu)})}\sJ_\mu)]_{P.A.}
\label{8.143}
\end{equation}
Taking into account \ref{8.123}, \ref{8.124}, \ref{8.128}, \ref{8.129} 
we deduce from \ref{8.118} - \ref{8.119}:
\begin{eqnarray}
&&[L(\Omega a\s^{(E_{(\nu)})}J_\mu^L)]_{P.A.}=\frac{h_{\nu\kappa}N^\kappa}{2c}
(\slap\lambdab+\lambdab\sd\mbox{tr}\tchib\cdot\pi^\sharp)\Lb\beta_\mu \label{8.144}\\
&&\hspace{27mm}+\rhob(\sD^2\lambdab)(E_{(\nu)},(\sd\beta_\mu)^\sharp)
+a\left([\sD\tchi]_{P.A.}\cdot\pi^\sharp\right)(E_{(\nu)},(\sd\beta_\mu)^\sharp)\nonumber
\end{eqnarray}
\begin{eqnarray}
&&[\Lb(\Omega a\s^{(E_{(\nu)})}J_\mu^{\Lb})]_{P.A.}=\frac{h_{\nu\kappa}\Nb^\kappa}{2c}
(\slap\lambda+\lambda\sd\mbox{tr}\tchi\cdot\pi^\sharp)L\beta_\mu \label{8.145}\\
&&\hspace{27mm}+\rho(\sD^2\lambda)(E_{(\nu)},(\sd\beta_\mu)^\sharp)
+a\left([\sD\tchib]_{P.A.}\cdot\pi^\sharp\right)(E_{(\nu)},(\sd\beta_\mu)^\sharp)\nonumber
\end{eqnarray}
To calculate $[\sdiv(\Omega a\s^{(E_{(\nu)})}\sJ_\mu)]_{P.A.}$ from \ref{8.120}, we 
use the following. First, that by the 1st of \ref{8.127}, if $X$ is a $S$ vectorfield of order 0 
we have:
\begin{equation}
[\sdiv(\tchi\cdot X)^\sharp]_{P.A.}=[\sdiv\tchi]_{P.A.}\cdot X=X\mbox{tr}\tchi
\label{8.146}
\end{equation}
Similarly, by the 2nd of  \ref{8.127}:
\begin{equation}
[\sdiv(\tchib\cdot X)^\sharp]_{P.A.}=[\sdiv\tchib]_{P.A.}\cdot X
=X\mbox{tr}\tchib
\label{8.147}
\end{equation}
In particular the above apply taking $X=E_{(\nu)},\pi^\sharp$.
Second, that:
\begin{eqnarray}
&&[\sdiv\left((E_{(\nu)}\lambda)(\sd\psi)^\sharp\right)]_{P.A.}=
(\sD^2\lambda)(E_{(\nu)},(\sd\psi)^\sharp) \nonumber\\
&&[\sdiv\left((E_{(\nu)}\lambdab)(\sd\psi)^\sharp\right)]_{P.A.}=
(\sD^2\lambdab)(E_{(\nu)},(\sd\psi)^\sharp)
\label{8.148}
\end{eqnarray}
Third that:
\begin{eqnarray}
&&[\sdiv(\tchi(E_{(\nu)},\pi^\sharp)(\sd\psi)^\sharp)]_{P.A.}=
[\sD_{(\sd\psi)^\sharp}\tchi]_{P.A.}(E_{(\nu)},\pi^\sharp)\nonumber\\
&&[\sdiv(\tchib(E_{(\nu)},\pi^\sharp)(\sd\psi)^\sharp)]_{P.A.}=
[\sD_{(\sd\psi)^\sharp}\tchib]_{P.A.}(E_{(\nu)},\pi^\sharp)\label{8.149}
\end{eqnarray}
And fourth that by \ref{8.127}:
\begin{eqnarray}
&&[\sdiv((\tchi^\sharp-(1/2)\mbox{tr}\tchi I)\cdot(\sd\psi)^\sharp)]_{P.A.}
=[\sdiv\tchi-(1/2)\sd\mbox{tr}\tchi]_{P.A.}\cdot\sd\psi\nonumber\\
&&\hspace{40mm}=(1/2)\sd\mbox{tr}\tchi\cdot (\sd\psi)^\sharp\nonumber\\
&&[\sdiv((\tchib^\sharp-(1/2)\mbox{tr}\tchib I)\cdot(\sd\psi)^\sharp)]_{P.A.}
=[\sdiv\tchib-(1/2)\sd\mbox{tr}\tchib]_{P.A.}\cdot\sd\psi\nonumber\\
&&\hspace{40mm}=(1/2)\sd\mbox{tr}\tchib\cdot (\sd\psi)^\sharp
\label{8.150}
\end{eqnarray}
We obtain:
\begin{eqnarray}
&&\sdiv(\Omega a\s^{(E_{(\nu)})}\sJ_\mu)=
\frac{1}{2}\rhob (E_{(\nu)}\mbox{tr}\tchib)L\beta_\mu-\frac{h_{\nu\kappa}\Nb^\kappa}{2c} 
(\slap\lambda+\lambda\sd\mbox{tr}\tchi\cdot\pi^\sharp)L\beta_\mu\nonumber\\
&&\hspace{27mm}+\frac{1}{2}\rho (E_{(\nu)}\mbox{tr}\tchi)\Lb\beta_\mu
-\frac{h_{\nu\kappa}N^\kappa}{2c} 
(\slap\lambdab+\lambdab\sd\mbox{tr}\tchib\cdot\pi^\sharp)\Lb\beta_\mu\nonumber\\
&&\hspace{27mm}-\rho(\sD^2\lambda)(E_{(\nu)},(\sd\beta_\mu)^\sharp)
-\rhob(\sD^2\lambdab)(E_{(\nu)},(\sd\beta_\mu)^\sharp)\nonumber\\
&&\hspace{27mm}-a[\sD_{(\sd\beta_\mu)^\sharp}\tchi]_{P.A.}(E_{(\nu)},\pi^\sharp)
-a[\sD_{(\sd\beta_\mu)^\sharp}\tchib]_{P.A.}(E_{(\nu)},\pi^\sharp)\nonumber\\
&&\hspace{27mm}+\frac{ah_{\nu\kappa}}{2c}(\Nb^\kappa\sd\mbox{tr}\tchi+N^\kappa\sd\mbox{tr}\tchib)
\cdot(\sd\beta_\mu)^\sharp
\label{8.151}
\end{eqnarray}
In summing \ref{8.144}, \ref{8.145}, and \ref{8.151} to obtain $\s^{(0,\nu)}\tilde{\rho}_\mu$ 
according to \ref{8.143}, the terms in $\sD^2\lambda$, $\slap\lambda$, $\sD^2\lambdab$, $\slap\lambdab$ 
cancel. Also, noting that by \ref{8.125} we have:
$$\left([\sD\tchi]_{P.A.}\cdot\pi^\sharp\right)(E_{(\nu)},(\sd\beta_\mu)^\sharp)=
[\sD_{(\sd\beta_\mu)^\sharp}\tchi]_{P.A.}(E_{(\nu)},\pi^\sharp)$$
$$\left([\sD\tchib]_{P.A.}\cdot\pi^\sharp\right)(E_{(\nu)},(\sd\beta_\mu)^\sharp)=
[\sD_{(\sd\beta_\mu)^\sharp}\tchib]_{P.A.}(E_{(\nu)},\pi^\sharp)$$
the terms involving $[\sD\tchi]_{P.A.}$, $[\sD\tchib]_{P.A.}$ likewise cancel. Moreover, 
the terms involving $\sd\mbox{tr}\tchi\cdot\pi^\sharp$, $\sd\mbox{tr}\tchib\cdot\pi^\sharp$ 
cancel as well. What remains is:
\begin{eqnarray}
&&\s^{(0,\nu)}\tilde{\rho}_\mu=\frac{1}{2}\rhob(E_{(\nu)}\mbox{tr}\tchib)L\beta_\mu+
\frac{1}{2}\rho(E_{(\nu)}\mbox{tr}\tchi)\Lb\beta_\mu \nonumber\\
&&\hspace{13mm}+\frac{ah_{\nu\kappa}}{2c}(\Nb^\kappa\sd\mbox{tr}\tchi+N^\kappa\sd\mbox{tr}\tchib)
\cdot (\sd\beta_\mu)^\sharp
\label{8.152}
\end{eqnarray}

Now, the higher order variations $\s^{(m,l)}\dot{\phi}_\mu$ in the case $n=2$, 
$\s^{(m,\nu_1...\nu_l)}\dot{\phi}_\mu$ for $n>2$, given by \ref{8.1}, \ref{8.2} respectively, 
are quantities of order $m+l$, so the 1-forms $\s^{(V;m,l)}\xi$, $\s^{(V;m,\nu_1..\nu_l)}\xi$ 
associated to these and to the variation field $V$ are quantities of order $m+l+1$. So are the 
corresponding rescaled sources 
$\s^{(m,l)}\tilde{\rho}_\mu$, $\s^{(m,\nu_1...\nu_l)}\tilde{\rho}_\mu$. The principal acoustical error 
terms in the energy identities corresponding to $\s^{(m,l)}\dot{\phi}_\mu$ (case $n=2$), 
$\s^{(m,\nu_1...\nu_l)}\dot{\phi}_\mu$ (for $n>2$) are contained in $\s^{(V;m,l))}Q_3$, 
$\s^{(V;m,\nu_1...\nu_l)}Q_3$ given by \ref{6.66} with $\s^{(m,l)}\dot{\phi}_\mu$, 
$\s^{(m,\nu_1...\nu_l)}\dot{\phi_\mu}$ in the role of $\dot{\phi}$ and $\s^{(m,l)}\tilde{\rho}_\mu$, 
$\s^{(m,\nu_1...\nu_l)}\tilde{\rho}_\mu$ in the role of $\tilde{\rho}_\mu$. The multiplier field 
$X$ having been fixed by \ref{7.67}, we have:
\begin{eqnarray}
&&\Omega a\s^{(V;m,l)}Q_3=-(3\s^{(V;m,l)}\xi_L+\s^{(V;m,l)}\xi_{\Lb})V^\mu\s^{(m,l)}\tilde{\rho}_\mu 
 \ \ \mbox{: for $n=2$}\nonumber\\ 
&&\label{8.153}\\
&&\Omega a\s^{(V;m,\nu_1...\nu_l)}Q_3=-(3\s^{(V;m,\nu_1...\nu_l)}\xi_L
+\s^{(V;m,\nu_1...\nu_l)}\xi_{\Lb})V^\mu\s^{((m,\nu_1...\nu_l)}\tilde{\rho}_\mu\nonumber\\
&&\mbox{: for $n>2$} \label{8.154}
\end{eqnarray}
In view of the recursion formulas \ref{8.28}, \ref{8.29} in the case $n=2$, and \ref{8.31}, \ref{8.32} 
for $n>2$, we have, in the case $n=2$:
\begin{eqnarray}
&&[\s^{(0,l)}\tilde{\rho}_\mu]_{P.A.}=[E^{l-1}\s^{(0,1)}\tilde{\rho}_\mu]_{P.A.} \label{8.155}\\
&&[\s^{(m,l)}\tilde{\rho}_\mu]_{P.A.}=[E^l T^{m-1}\s^{(1,0)}\tilde{\rho}_\mu]_{P.A.} \ \ 
\mbox{: for $m\geq 1$} \nonumber
\end{eqnarray}
and for $n>2$:
\begin{eqnarray}
&&[\s^{(0,\nu_1...\nu_l)}\tilde{\rho}_\mu]_{P.A.}=
[E_{(\nu_l)}...E_{(\nu_2)}\s^{(0,\nu_1)}\tilde{\rho}_\mu]_{P.A.} \label{8.156}\\
&&[\s^{(m,\nu_1...\nu_l)}\tilde{\rho}_\mu]_{P.A.}=
[E_{(\nu_l)}...E_{(\nu_1)}T^{m-1}\s^{(1)}\tilde{\rho}]_{P.A.} \ \ \mbox{: for $m\geq 1$} \nonumber
\end{eqnarray}
From \ref{8.142} and \ref{8.135}, which in the case $n=2$ takes the form:
\begin{eqnarray}
&&[\s^{(1,0)}\tilde{\rho}_\mu]_{P.A.}=\rho E^2\lambda \ \Lb\beta_\mu+\rhob E^2\lambdab \ L\beta_\mu 
\label{8.157}\\
&&\hspace{10mm}+a E\tchi \ (\pi\Lb\beta_\mu+\rho E\beta_\mu)
+a E\tchib \ (\pi L\beta_\mu+\rhob E\beta_\mu) \nonumber
\end{eqnarray}
we then obtain, in the case $n=2$:
\begin{eqnarray}
&&[\s^{(0,l)}\tilde{\rho}_\mu]_{P.A.}=\frac{1}{2}\rhob E^l\tchib \ L\beta_\mu 
+\frac{1}{2}\rho E^l\tchi \ \Lb\beta_\mu \label{8.158}\\
&&[\s^{(m,l)}\tilde{\rho}_\mu]_{P.A.}=\rho E^{l+2}T^{m-1}\lambda \ \Lb\beta_\mu 
+\rhob E^{l+2} T^{m-1}\lambdab \ L\beta_\mu\nonumber\\
&&\hspace{15mm}+aE^{l+1}[T^{m-1}\tchi]_{P.A.} \ (\pi\Lb\beta_\mu+\rho E\beta_\mu) \nonumber\\
&&\hspace{10mm}+aE^{l+1}[T^{m-1}\tchib]_{P.A.} \ (\pi L\beta_\mu+\rhob E\beta_\mu) \ \ 
\mbox{: for $m\geq 1$} \label{8.159}
\end{eqnarray}
In the last, for $m\geq 2$ we must appeal to \ref{8.137}, \ref{8.138} to determine 
$[T^{m-1}\tchi]_{P.A.}$, $[T^{m-1}\tchi]_{P.A.}$. From \ref{8.152} and \ref{8.135} 
we obtain, through \ref{8.156}, for $n>2$:
\begin{eqnarray}
&&[\s^{(0,\nu_1...\nu_l)}\tilde{\rho}_\mu]_{P.A.}=\frac{1}{2}\rhob(E_{(\nu_l)}...E_{(\nu_1)}
\mbox{tr}\tchib)L\beta_\mu 
+\frac{1}{2}\rho(E_{(\nu_l)}...E_{(\nu_1)}\mbox{tr}\tchi)\Lb\beta_\mu\nonumber\\
&&\hspace{15mm}+\frac{ah_{\nu_1,\kappa}}{2c}
\left(\Nb^\kappa\sd(E_{(\nu_l)}...E_{(\nu_2)}\mbox{tr}\tchi)
+N^\kappa\sd(E_{(\nu_l)}...E_{(\nu_2)}\mbox{tr}\tchib)\right)\cdot(\sd\beta_\mu)^\sharp \nonumber\\ 
&&\label{8.160}\\
&&[\s^{(m,\nu_1...\nu_l)}\tilde{\rho}_\mu]_{P.A.}=
\rho (E_{(\nu_l)}...E_{(\nu_1)}T^{m-1}\slap\lambda) \ \Lb\beta_\mu \nonumber\\
&&\hspace{28mm}+\rhob (E_{(\nu_l)}...E_{(\nu_1)}T^{m-1}\slap\lambdab \ L\beta_\mu \nonumber\\
&&\hspace{15mm}+a \sd(E_{(\nu_l)}...E_{(\nu_1)}[T^{m-1}\mbox{tr}\tchi]_{P.A.})\cdot
(\pi^\sharp\Lb\beta_\mu+\rho(\sd\beta_\mu)^\sharp)\nonumber\\
&&\hspace{12mm}+a\sd(E_{(\nu_l)}...E_{(\nu-1}[T^{m-1}\mbox{tr}\tchib]_{P.A.})\cdot 
(\pi^\sharp L\beta_\mu+\rhob(\sd\beta_\mu)^\sharp) \ \ \mbox{: for $m\geq 1$} \nonumber\\
&&\label{8.161}
\end{eqnarray}
where in the last, for $m\geq 2$ we must appeal to \ref{8.128}, \ref{8.129} to determine 
$[T^{m-1}\tchi]_{P.A.}$, $[T^{m-1}\tchib]_{P.A.}$. 

\pagebreak

\chapter{The Power Series Approximation Method}

\section{Setup of the Truncated Power Series} 

From this point we shall focus on the case $n=2$, except in Chapter 11. All the arguments, 
with the exception of the top order estimates for the acoustical functions, readily 
extend to the case $n>2$. The top order estimates for the acoustical functions in the case $n>2$ 
is the topic of Chapter 11. 

Let $\tau$ be the parameter along the integral curves of $T$, measured from $\Cb_0$. Let us also 
introduce the parameter $\sigma$ which is constant along the integral curves of $T$ and equal to 
$u$ on $\Cb_0$. We then have:
\begin{equation}
\begin{array}{llll}
&\tau=\ub & \ \ \ &\ub=\tau\\
&\sigma=u-\ub & \ \ \ &u=\sigma+\tau
\end{array}
\label{9.1}
\end{equation}
Noting that the vectorfield $T$, being expressed in coordinates $(\ub,u,\vartheta)$ by  \ref{2.34}, is expressed in terms of $(\tau,\sigma,\vartheta)$ coordinates simply by:
\begin{equation}
T=\frac{\partial}{\partial\tau}
\label{9.2}
\end{equation}
we define the $N$th order approximants the $x^\mu$, $b$, and the $\beta_\mu$ by the power series 
expansions in $\tau$, truncated at order $N$:
\begin{eqnarray}
&&x^\mu_N(\tau,\sigma,\vartheta)=\sum_{n=0}^N\frac{1}{n!}(T^n x^\mu)(0,\sigma,\vartheta)\tau^n 
\nonumber\\
&&b_N(\tau,\sigma,\vartheta)=\sum_{n=0}^N\frac{1}{n!}(T^n b)(0,\sigma,\vartheta)\tau^n \nonumber\\
&&\beta_{\mu,N}(\tau,\sigma,\vartheta)=\sum_{n=0}^N\frac{1}{n!}(T^n\beta_\mu)(0,\sigma,\vartheta)
\tau^n \label{9.3}
\end{eqnarray}
expressed  in  $(\tau,\sigma,\vartheta)$). It is important to note in connection with the first of the above that the $(\tau,\sigma,\vartheta)$ or equivalently the $\ub,u,\vartheta)$ coordinates are fixed 
and what is being approximated is the functions expressing the rectangular coordinates in terms of these. The $\beta_{\mu,N}$ together with the $x_N^\mu$ define a 1-form on ${\cal N}$ by:
\begin{equation}
\beta_N=\beta_{\mu,N}dx_N^\mu
\label{9.a1}
\end{equation}
This is the $N$th approximant 1-form $\beta$. 

Defining, in analogy with \ref{2.129}, $\sigma_N$ the $N$th approximant $\sigma$ function, by:
\begin{equation}
\sigma_N=-(g^{-1})^{\mu\nu}\beta_{\mu,N}\beta_{\nu,N}
\label{9.4}
\end{equation}
We then set:
\begin{equation}
H_N=H(\sigma_N)
\label{9.5}
\end{equation}
and define the rectangular components $h_{\mu\nu,N}$ of$N$th approximant acoustical metric in analogy with \ref{1.d11} by:
\begin{equation}
h_{\mu\nu,N}=g_{\mu\nu}+H_N\beta_{\mu,N}\beta_{\nu,N}
\label{9.6}
\end{equation}
The $h_{\mu\nu,N}$ together with the $x_N^\mu$ define the $N$th approximant acoustical metric:
\begin{equation}
h_N=h_{\mu\nu,N}dx_N^\mu\otimes dx_N^\nu
\label{9.a2}
\end{equation}
a Lorentzian metric on ${\cal N}$. The inverse of \ref{9.6} is:
\begin{equation}
(h_N^{-1})^{\mu\nu}=(g^{-1})^{\mu\nu}-F_N\beta_N^\mu\beta_N^\nu
\label{9.b4}
\end{equation}
where: 
\begin{equation}
F_N=F(\sigma_N)
\label{9.b5}
\end{equation}

The vectorfield $\Omega$ is defined by the 3rd of \ref{3.a8}, so like the vectorfield $T$, it is 
independent of the approximation:
\begin{equation}
T=\frac{\partial}{\partial\ub}+\frac{\partial}{\partial u}, \ \ \ 
\Omega=\frac{\partial}{\partial\vartheta}
\label{9.7}
\end{equation}
However, the rectangular coordinates being given by the first of \ref{9.3}, $\Omega_N^\mu$, the 
rectangular components of $\Omega$, do depend on the approximation:
\begin{equation}
\Omega_N^\mu=\Omega x_N^\mu=\frac{\partial x_N^\mu}{\partial\vartheta}
\label{9.8}
\end{equation}
and so does the induced metric $\sh_N$, which is given in analogy with \ref{3.a9} by:
\begin{equation}
\sh_N=h_{\mu\nu,N}\Omega_N^\mu\Omega_N^\nu
\label{9.9}
\end{equation}
This being the squared magnitude of $\Omega$ relative to the $N$th approximant metric, the 
corresponding unit vectorfield is: 
\begin{equation}
E_N=\frac{1}{\sqrt{\sh_N}}\Omega
\label{9.10}
\end{equation}
This is then the $N$th approximant version of the vectorfield $E$ (see \ref{3.a10}). The $N$th 
approximant versions of the vectorfields $L$, $\Lb$ are given by:
\begin{equation}
L_N=\frac{\partial}{\partial\ub}-b_N\frac{\partial}{\partial\vartheta}, \ \ \ 
\Lb_N=\frac{\partial}{\partial u}+b_N\frac{\partial}{\partial\vartheta}
\label{9.11}
\end{equation}
in analogy with the first two of \ref{3.a8}. In analogy with the commutation relations \ref{3.a14} 
we have:
\begin{eqnarray}
&&[L_N,E_N]=-\chi_N E_N \ \ \ \ \ [\Lb_N,E_N]=-\chib_N E_N \nonumber\\
&&\hspace{25mm} [\Lb_N,L_N]=Z_N \label{9.b1}
\end{eqnarray}
Here, 
\begin{equation}
\chi_N=\frac{1}{2\sh_N}(L_N\sh_N-2\sh_N\Omega b_N) \ \ \ \ \ 
\chib_N=\frac{1}{2\sh_N}(\Lb_N\sh_N+2\sh_N\Omega b_N)
\label{9.b2}
\end{equation}
and $Z_N$ is the vectorfield:
\begin{equation}
Z_N=-(Tb_N)\Omega=\zeta_N E_N
\label{9.b3}
\end{equation}
(see \ref{8.50}).
The first variational formulas \ref{3.4}, \ref{3.5} reduce to \ref{9.b2} in the case $n=2$ when 
the $S$ vectorfield $b_N\partial/\partial\vartheta$ is placed in the role of the $S$ vectorfield $b$, 
and the symmetric 2-covariant $S$ tensorfields $\sh_N d\vartheta\otimes d\vartheta$, 
$\sh_N\chi_N d\vartheta\otimes d\vartheta$, $\sh_N\chib_N d\vartheta\otimes d\vartheta$ are placed 
in the roles of the symmetric 2-covariant $S$ tensorfields $\sh$, $\chi$, $\chib$ respectively. 
Also, \ref{2.35} reduces to \ref{9.b3} in the case $n=2$ when the $S$ vectorfield 
$b_N\partial/\partial\vartheta$ is placed in the role of the $S$ vectorfield $b$.

In analogy with the construction of Section 2.2, we define  
the components $N_N^\mu$ by the conditions:
\begin{equation}
h_{\mu\nu,N}N_N^\mu\Omega_N^\nu=0, \ \ \ h_{\mu\nu,N}N_N^\mu N_N^\nu=0, \ \ \ N_N^0=1
\label{9.12}
\end{equation} 
together with the condition that the vectorfield with rectangular components $N_N^\mu$ be outgoing. 
We also define the components $\Nb_N^\mu$ by the conditions:
\begin{equation}
h_{\mu\nu,N}\Nb_N^\mu\Omega_N^\nu=0, \ \ \ h_{\mu\nu,N}\Nb_N^\mu\Nb_N^\nu=0, \ \ \ \Nb_N^0=1
\label{9.13}
\end{equation}
together with the condition that the vectorfield with rectangular components $\Nb_N^\mu$ be incoming. 
(See \ref{2.49}, \ref{2.51}, \ref{2.52}, \ref{2.58}, \ref{2.59}.) Finally, we define the $N$th 
approximant functions $\rho_N$, $\rhob_N$ by:
\begin{equation}
\rho_N=L_N t_N, \ \ \ \rhob_N=\Lb_N t_N \ \ \ \mbox{: where $t_N=x^0_N$}
\label{9.14}
\end{equation}
(see \ref{2.61}).

We remark that from the results of Chapter 5 the coefficients of the expansions \ref{9.3} 
are known smooth functions of $(\sigma,\vartheta)$, hence the $N$th approximants $x_N^\mu$, 
$b_N$, $\beta_{\mu,N}$ are in $(\tau,\sigma,\vartheta)$ coordinates $N$th degree polynomials in $\tau$ 
with coefficients which are known smooth functions of $(\sigma,\vartheta)$. Moreover, $\sigma_N$ 
and $H_N$ are known positive smooth functions of $(\tau,\sigma,\vartheta)$ or equivalently of 
$(\ub,u,\vartheta)$. Then $h_{\mu\nu,N}$ are known smooth functions representing the rectangular 
components of a Lorentzian metric, and $F_N$ and the components $(h_N^{-1})^{\mu\nu}$ of the inverse 
are likewise known smooth functions. Also $N_N^\mu$, $\Nb_N^\mu$ defined according to the above 
through $h_{\mu\nu,N}$ by \ref{9.13}, \ref{9.14} are also known smooth functions. Moreover $\rho_N$ and 
$\rhob_N$ are known smooth functions, in fact polynomials in $\tau$ of degree $2N$ with coefficients 
which are known smooth functions of $(\sigma,\vartheta)$. Finally since $N_N^\mu$, $\Nb_N^\mu$ and 
\begin{equation}
E_N^\mu=\frac{\Omega_N^\mu}{\sqrt{\sh_N}}
\label{9.b6}
\end{equation}
constitute a null frame for $h_{\mu\nu,N}$, we have the expansion:
\begin{equation}
(h_N^{-1})^{\mu\nu}=-\frac{1}{2c_N}(N_N^\mu\Nb_N^\nu+\Nb_N^\mu N_N^\nu)+E_N^\mu E_N^\nu
\label{9.b7}
\end{equation}
where, recalling from \ref{2.71} the function $c$:
\begin{equation}
c=-\frac{1}{2}h_{\mu\nu}N^\mu\Nb^\nu
\label{9.30}
\end{equation}
$c_N$ is the corresponding $N$th approximant function:
\begin{equation}
c_N=-\frac{1}{2}h_{\mu\nu,N}N_N^\mu\Nb_N^\nu
\label{9.31}
\end{equation}

\vspace{5mm}

\section{Estimates for the Quantities by which the $N$th Approximants fail to satisfy 
the Characteristic and Wave Systems} 

Now the $N$th approximants \ref{9.3} do not satisfy with \ref{9.11}, \ref{9.14} and $N_N^\mu$, 
$\Nb_N^\mu$ defined as above the characteristic system \ref{2.62}. That is, the quantities 
\begin{eqnarray}
&&L_N x_N^i-\rho_N N_N^i:=\vep_N^i \nonumber\\
&&\Lb_N x_N^i-\rhob_N\Nb_N^i:=\vepb_N^i
\label{9.15}
\end{eqnarray}
do not vanish. By the above, the quantities $\vep_N^i$, $\vepb_N^i$ are known smooth functions of 
$(\tau,\sigma,\vartheta)$. We shall presently estimate them. 

By \ref{9.2} if $f$ is a arbitrary function on ${\cal N}$ in $(\tau,\sigma,\vartheta)$ coordinates 
we have
\begin{equation}
T^n f=\frac{\partial^n f}{\partial\tau^n}, \ \ \ 
\left.T^n f\right|_{\Cb_0}=\left.\frac{\partial^n f}{\partial\tau^n}\right|_{\tau=0}
\label{9.16}
\end{equation}
By \ref{9.3} we have:
\begin{eqnarray}
&&\left.\frac{\partial^n x_N^\mu}{\partial\tau^n}\right|_{\tau=0}
=\left.\frac{\partial^n x^\mu}{\partial\tau^n}\right|_{\tau=0}, \ \ \
\left.\frac{\partial^n b_N}{\partial\tau^n}\right|_{\tau=0}
=\left.\frac{\partial^n b}{\partial\tau^n}\right|_{\tau=0}, \nonumber\\
&&\left.\frac{\partial^n\beta_{\mu,N}}{\partial\tau^n}\right|_{\tau=0}
=\left.\frac{\partial^n\beta_\mu}{\partial\tau^n}\right|_{\tau=0}
\nonumber\\
&&\hspace{20mm} \mbox{: for $n=0,...,N$}
\label{9.17}
\end{eqnarray}
It follows that:
\begin{equation}
\left.\frac{\partial^n h_{\mu\nu,N}}{\partial\tau^n}\right|_{\tau=0}=
\left.\frac{\partial^n h_{\mu\nu}}{\partial\tau^n}\right|_{\tau=0} \ \ 
\mbox{: for $n=0,...,N$}
\label{9.18}
\end{equation}
and:
\begin{equation}
\left.\frac{\partial^n\Omega_N^\mu}{\partial\tau^n}\right|_{\tau=0}=
\left.\frac{\partial^n\Omega^\mu}{\partial\tau^n}\right|_{\tau=0} \ \ 
\mbox{: for $n=0,...,N$}
\label{9.19}
\end{equation}

\vspace{2.5mm} 

\noindent{\bf Lemma 9.1} \ \ \ We have:
\begin{eqnarray*}
&&\left.\frac{\partial^n N_N^\mu}{\partial\tau^n}\right|_{\tau=0}=
\left.\frac{\partial^n N^\mu}{\partial\tau^n}\right|_{\tau=0}, \ \ \ 
\left.\frac{\partial^n\Nb_N^\mu}{\partial\tau^n}\right|_{\tau=0}=
\left.\frac{\partial^n\Nb^\mu}{\partial\tau^n}\right|_{\tau=0}\\
&&\hspace{25mm}\mbox{: for $n=0,...,N$} 
\end{eqnarray*}

\vspace{2.5mm}

\noindent{\em Proof:} Taking $n=0$ in \ref{9.18}, \ref{9.19} and comparing \ref{9.12}, \ref{9.13} 
with \ref{2.49}, \ref{2.51}, \ref{2.52} we conclude that on $\Cb_0$ $N_N^\mu$ and $\Nb_N^\mu$ 
are subject to the same conditions as $N^\mu$ and $\Nb^\mu$ respectively. Since these conditions 
possess a unique solution, it follows that:
\begin{equation}
\left.N_N^\mu\right|_{\tau=0}=\left.N^\mu\right|_{\tau=0}, \ \ \ 
\left.\Nb_N^\mu\right|_{\tau=0}=\left.\Nb^\mu\right|_{\tau=0}
\label{9.20}
\end{equation}
We proceed by finite induction. Let 
\begin{eqnarray}
&&\left.\frac{\partial^n N_N^\mu}{\partial\tau^n}\right|_{\tau=0}=
\left.\frac{\partial^n N^\mu}{\partial\tau^n}\right|_{\tau=0}, \ \ \ 
\left.\frac{\partial^n\Nb_N^\mu}{\partial\tau^n}\right|_{\tau=0}=
\left.\frac{\partial^n\Nb^\mu}{\partial\tau^n}\right|_{\tau=0} \nonumber\\
&&\hspace{10mm}\mbox{: for $n=0,...,m$, where $0\leq m\leq N-1$} 
\label{9.21}
\end{eqnarray}
Differentiating \ref{2.49}, \ref{2.51}, \ref{2.52} with respect to $\tau$ at $\tau=0$ $m+1$ times, 
with $N^\mu$, $\Nb^\mu$ in the role of $M^\mu$, we obtain: 
\begin{eqnarray}
&&\left.h_{\mu\nu}\right|_{\tau=0}\left.\frac{\partial^{m+1}N^\mu}{\partial\tau^{m+1}}\right|_{\tau=0}
\left.\Omega^\nu\right|_{\tau=0}+A_m=0, \nonumber\\
&&\left.2h_{\mu\nu}\right|_{\tau=0}\left.\frac{\partial^{m+1}N^\mu}{\partial\tau^{m+1}}\right|_{\tau=0}
\left.N^\nu\right|_{\tau=0}+B_m=0, \nonumber\\
&&\left.\frac{\partial^{m+1}N^0}{\partial\tau^{m+1}}\right|_{\tau=0}=0 
\label{9.22}
\end{eqnarray}
and
\begin{eqnarray}
&&\left.h_{\mu\nu}\right|_{\tau=0}\left.\frac{\partial^{m+1}\Nb^\mu}{\partial\tau^{m+1}}\right|_{\tau=0}
\left.\Omega^\nu\right|_{\tau=0}+\underline{A}_m=0, \nonumber\\
&&\left.2h_{\mu\nu}\right|_{\tau=0}\left.\frac{\partial^{m+1}\Nb^\mu}{\partial\tau^{m+1}}\right|_{\tau=0}
\left.\Nb^\nu\right|_{\tau=0}+\underline{B}_m=0, \nonumber\\
&&\left.\frac{\partial^{m+1}\Nb^0}{\partial\tau^{m+1}}\right|_{\tau=0}=0 
\label{9.23}
\end{eqnarray}
Here $A_m$ depends on:
$$\left.\frac{\partial^n h_{\mu\nu}}{\partial\tau^n}\right|_{\tau=0}, \ \ 
\left.\frac{\partial^n\Omega^\mu}{\partial\tau^n}\right|_{\tau=0} \ \ \mbox{: for $n=0,...,m+1$}$$
and on:
$$\left.\frac{\partial^n N^\mu}{\partial\tau^n}\right|_{\tau=0} \ \ \mbox{: for $n=0,...,m$}$$
and $B_m$ depends on:
$$\left.\frac{\partial^n h_{\mu\nu}}{\partial\tau^n}\right|_{\tau=0} \ \ \mbox{: for $n=0,...,m+1$}$$
and on:
$$\left.\frac{\partial^n N^\mu}{\partial\tau^n}\right|_{\tau=0} \ \ \mbox{: for $n=0,...,m$}$$
Moreover, $\underline{A}_m$ and $\underline{B}_m$ are the same as $A_m$ and $B_m$ respectively 
but with $\Nb^\mu$ in the role of $N^\mu$. 

Differentiating \ref{9.12}, \ref{9.13} with respect to $\tau$ at $\tau=0$ $m+1$ times we obtain, in view of \ref{9.18}, \ref{9.19} for $n=0$ and \ref{9.20}, 
\begin{eqnarray}
&&\left.h_{\mu\nu}\right|_{\tau=0}\left.\frac{\partial^{m+1}N_N^\mu}{\partial\tau^{m+1}}\right|_{\tau=0}
\left.\Omega^\nu\right|_{\tau=0}+A_{m,N}=0, \nonumber\\
&&\left.2h_{\mu\nu}\right|_{\tau=0}\left.\frac{\partial^{m+1}N_N^\mu}{\partial\tau^{m+1}}\right|_{\tau=0}
\left.N^\nu\right|_{\tau=0}+B_{m,N}=0, \nonumber\\
&&\left.\frac{\partial^{m+1}N_N^0}{\partial\tau^{m+1}}\right|_{\tau=0}=0 
\label{9.24}
\end{eqnarray}
and
\begin{eqnarray}
&&\left.h_{\mu\nu}\right|_{\tau=0}\left.\frac{\partial^{m+1}\Nb_N^\mu}{\partial\tau^{m+1}}\right|_{\tau=0}
\left.\Omega^\nu\right|_{\tau=0}+\underline{A}_{m,N}=0, \nonumber\\
&&\left.2h_{\mu\nu}\right|_{\tau=0}\left.\frac{\partial^{m+1}\Nb_N^\mu}{\partial\tau^{m+1}}\right|_{\tau=0}
\left.\Nb^\nu\right|_{\tau=0}+\underline{B}_{m,N}=0, \nonumber\\
&&\left.\frac{\partial^{m+1}\Nb_N^0}{\partial\tau^{m+1}}\right|_{\tau=0}=0 
\label{9.25}
\end{eqnarray}
Here $A_{m,N}$ is the same as $A_m$ but with 
$$\left.\frac{\partial^n h_{\mu\nu,N}}{\partial\tau^n}\right|_{\tau=0}, \ \ 
\left.\frac{\partial^n\Omega_N^\mu}{\partial\tau^n}\right|_{\tau=0} \ \ : \ n=0,...,m+1$$
in the role of 
$$\left.\frac{\partial^n h_{\mu\nu}}{\partial\tau^n}\right|_{\tau=0}, \ \ 
\left.\frac{\partial^n\Omega^\mu}{\partial\tau^n}\right|_{\tau=0} \ \ : \ n=0,...,m+1$$
and 
$$\left.\frac{\partial^n N_N^\mu}{\partial\tau^n}\right|_{\tau=0} \ \ : \ n=0,...,m$$
in the role of 
$$\left.\frac{\partial^n N^\mu}{\partial\tau^n}\right|_{\tau=0} \ \ : \ n=0,...,m$$
Then by \ref{9.18}, \ref{9.19} and the inductive hypothesis \ref{9.21}:
\begin{equation}
A_{m,N}=A_m
\label{9.26}
\end{equation}
Also, $B_{m,N}$ is the same as $B_m$ but with 
$$\left.\frac{\partial^n h_{\mu\nu,N}}{\partial\tau^n}\right|_{\tau=0} \ \ : \ n=0,...,m+1$$
in the role of 
$$\left.\frac{\partial^n h_{\mu\nu}}{\partial\tau^n}\right|_{\tau=0} \ \ : \ n=0,...,m+1$$
and 
$$\left.\frac{\partial^n N_N^\mu}{\partial\tau^n}\right|_{\tau=0} \ \ : \ n=0,...,m$$
in the role of 
$$\left.\frac{\partial^n N^\mu}{\partial\tau^n}\right|_{\tau=0} \ \ : \ n=0,...,m$$
Then by \ref{9.18} and the inductive hypothesis \ref{9.21}:
\begin{equation}
B_{m,N}=B_m
\label{9.27}
\end{equation}
Similarly, placing $\Nb_N^\mu$ in the role of $N_N^\mu$ we deduce, using the inductive hypothesis 
\ref{9.21}, 
\begin{equation}
\underline{A}_{m,N}=\underline{A}_m, \ \ \ \underline{B}_{m,N}=\underline{B}_m
\label{9.28}
\end{equation}
In view of \ref{9.26}, \ref{9.27}, \ref{9.28}, we see, comparing \ref{9.24}, \ref{9.25} with 
\ref{9.22}, \ref{9.23} that $\left.\partial^{m+1}N_N^\mu/\partial\tau^{m+1}\right|_{\tau=0}$, 
$\left.\partial^{m+1}\Nb_N^\mu/\partial\tau^{m+1}\right|_{\tau=0}$ satisfy the same linear 
systems of equations as $\left.\partial^{m+1}N^\mu/\partial\tau^{m+1}\right|_{\tau=0}$, 
$\left.\partial^{m+1}\Nb^\mu/\partial\tau^{m+1}\right|_{\tau=0}$ respectively. Noting that the 
corresponding homogeneous systems admit only the trivial solution, it follows that:
\begin{equation}
\left.\frac{\partial^{m+1}N_N^\mu}{\partial\tau^{m+1}}\right|_{\tau=0}=
\left.\frac{\partial^{m+1}N^\mu}{\partial\tau^{m+1}}\right|_{\tau=0}, \ \ \ 
\left.\frac{\partial^{m+1}\Nb_N^\mu}{\partial\tau^{m+1}}\right|_{\tau=0}=
\left.\frac{\partial^{m+1}\Nb^\mu}{\partial\tau^{m+1}}\right|_{\tau=0}
\label{9.29}
\end{equation}
and the inductive step is complete. 

\vspace{2.5mm}

Lemma 9.1 implies that:
\begin{equation}
\left.\frac{\partial^n c_N}{\partial\tau^n}\right|_{\tau=0}=
\left.\frac{\partial^n c}{\partial\tau^n}\right|_{\tau=0} \ \ 
\mbox{: for $n=0,...,N$}
\label{9.32}
\end{equation}

Since from \ref{9.1} we have:
\begin{equation}
\frac{\partial}{\partial\ub}=\frac{\partial}{\partial\tau}-\frac{\partial}{\partial\sigma}, \ \ \ 
\frac{\partial}{\partial u}=\frac{\partial}{\partial\sigma}
\label{9.33}
\end{equation}
the vectorfields $L$, $\Lb$ are expressed in $(\tau,\sigma,\vartheta)$ coordinates by:
\begin{equation}
L=\frac{\partial}{\partial\tau}-\frac{\partial}{\partial\sigma}-b\frac{\partial}{\partial\vartheta}, 
\ \ \ \Lb=\frac{\partial}{\partial\sigma}+b\frac{\partial}{\partial\vartheta}
\label{9.34}
\end{equation}
and their $N$th approximant versions \ref{9.11} by:
\begin{equation}
L_N=\frac{\partial}{\partial\tau}-\frac{\partial}{\partial\sigma}
-b_N\frac{\partial}{\partial\vartheta}, \ \ \ 
\Lb_N=\frac{\partial}{\partial\sigma}+b_N\frac{\partial}{\partial\vartheta}
\label{9.35}
\end{equation}
It then follows from the first two of \ref{9.17}, that:
\begin{eqnarray}
&&\left.\frac{\partial^n(L_N x_N^i)}{\partial\tau^n}\right|_{\tau=0}=
\left.\frac{\partial^n(Lx^i)}{\partial\tau^n}\right|_{\tau=0} \ \ \mbox{: for $n=0,...,N-1$}\nonumber\\
&&\left.\frac{\partial^n(\Lb_N x_N^i)}{\partial\tau^n}\right|_{\tau=0}=
\left.\frac{\partial^n(\Lb x^i)}{\partial\tau^n}\right|_{\tau=0} \ \ 
\mbox{: for $n=0,...,N$}
\label{9.36}
\end{eqnarray}
and, through \ref{9.14}, that:
\begin{eqnarray}
&&\left.\frac{\partial^n\rho_N}{\partial\tau^n}\right|_{\tau=0}=
\left.\frac{\partial^n\rho}{\partial\tau^n}\right|_{\tau=0} \ \ \mbox{: for $n=0,...,N-1$} \nonumber\\
&&\left.\frac{\partial^n\rhob_N}{\partial\tau^n}\right|_{\tau=0}=
\left.\frac{\partial^n\rhob}{\partial\tau^n}\right|_{\tau=0} \ \ \mbox{: for $n=0,...,N$} 
\label{9.37}
\end{eqnarray}
The above together with Lemma 9.1 imply, in view of the definitions \ref{9.15}, 
\begin{eqnarray}
&&\left.\frac{\partial^n\vep_N^i}{\partial\tau^n}\right|_{\tau=0}=
\left.\frac{\partial^n(Lx^i-\rho N^i)}{\partial\tau^n}\right|_{\tau=0}=0 \ \ \mbox{: for $n=0,...,N-1$} 
\nonumber\\
&&\left.\frac{\partial^n\vepb_N^i}{\partial\tau^n}\right|_{\tau=0}=
\left.\frac{\partial^n(\Lb x^i-\rhob\Nb^i)}{\partial\tau^n}\right|_{\tau=0}=0 \ \ 
\mbox{: for $n=0,...,N$}
\label{9.38}
\end{eqnarray}
The quantities $\vep_N^i$, $\vepb_N^i$ being known smooth functions of $(\tau,\sigma,\vartheta)$ it follows 
that there are known smooth functions $\hat{\vep}_N^i$, $\hat{\vepb}_N^i$ of the same such that:
\begin{equation}
\vep_N^i=\tau^N\hat{\vep}_N^i, \ \ \ \vepb_N^i=\tau^{N+1}\hat{\vepb}_N^i
\label{9.39}
\end{equation}
To simplify the notation we introduce the following. 

\vspace{2.mm}

\noindent{\bf Definition 9.1:} \ \ Let $f$ be a known smooth function of $(\tau,\sigma,\vartheta)$ and 
let $p$ be a non-negative integer. We write: 
$$f=O(\tau^p)$$
if there is another known smooth function $\hat{f}$ of $(\tau,\sigma,\vartheta)$ such that:
$$f=\tau^p\hat{f}$$

\vspace{2.5mm}

Note that according to this definition, if $n=0,...,p$ and $m$, $l$ are non-negative integers we have:
$$\frac{\partial^{n+m+l}f}{\partial\tau^n\partial\sigma^m\partial\vartheta^l}=O(\tau^{p-n})$$ 

In the sense of the above definition we can state the results \ref{9.39} in the following form. 

\vspace{2.5mm}

\noindent{\bf Proposition 9.1} \ \ We have:
$$\vep_N^i=O(\tau^N), \ \ \ \vepb_N^i=O(\tau^{N+1})$$

\vspace{2.5mm} 

In the following we take, in accordance with \ref{9.14}, \ref{9.15}, $\vep_N^0=\vepb_N^0=0$. 

\vspace{2.5mm}

By \ref{9.b6}, \ref{9.8} we have:
\begin{equation}
E_N x_N^\mu=E_N^\mu
\label{9.b8}
\end{equation}
The question arises to find the vectorfields $N_N$, $\Nb_N$ such that:
\begin{equation}
N_N x_N^\mu=N_N^\mu, \ \ \ \Nb_N x_N^\mu=\Nb_N^\mu 
\label{9.b9}
\end{equation}
To answer this we consider the system \ref{9.15} with the 0-component included:
\begin{eqnarray}
&&L_N x_N^\mu=\rho_N N_N^\mu+\vep_N^\mu \nonumber\\
&&\Lb_N x_N^\mu=\rhob_N\Nb_N^\mu+\vepb_N^\mu 
\label{9.b10}
\end{eqnarray}
Now $E_N^\mu$, $N_N^\mu, \Nb_N^\mu$ constitute the rectangular components of a frame field.
Therefore, the rectangular components $V^\mu$ of an arbitrary vectorfield can be expanded as: 
\begin{equation}
V^\mu=V_N^E E_N^\mu+V_N^N N_N^\mu+V_N^{\Nb}\Nb_N^\mu
\label{9.81}
\end{equation}
and by \ref{9.9}, \ref{9.12}, \ref{9.13}, \ref{9.31} we have:
\begin{eqnarray}
&&V_N^E=h_{\mu\nu,N}V^\mu E_N^\nu \nonumber\\
&&V_N^N=-\frac{1}{2c_N}h_{\mu\nu,N}V^\mu \Nb_N^\nu \nonumber\\
&&V_N^{\Nb}=-\frac{1}{2c_N}h_{\mu\nu,N}V^\mu N_N^\nu \label{9.82}
\end{eqnarray}
Placing then $\vep_N^\mu$, $\vepb_N^\mu$ in the role of $V^\mu$ we expand:
\begin{eqnarray}
&&\vep_N^\mu=\vep_N^E E_N^\mu+\vep_N^N N_N^\mu+\vep_N^{\Nb}\Nb_N^\mu\nonumber\\
&&\vepb_N^\mu=\vepb_N^E E_N^\mu+\vepb_N^N N_N^\mu+\vepb_N^{\Nb}\Nb_N^\mu 
\label{9.b11}
\end{eqnarray}
Substituting these expansions in \ref{9.b10} and then substituting \ref{9.b8}, \ref{9.b9} yields 
the following linear system for the vectorfields $N_N$, $\Nb_N$:
\begin{eqnarray}
&&(\rho_N+\vep_N^N)N_N+\vep_N^{\Nb}\Nb_N=L_N-\vep_N^E E_N\nonumber\\
&&\vepb_N^N N_N+(\rhob_N+\vepb_N^{\Nb})\Nb_N=\Lb_N-\vepb_N^E E_N 
\label{9.b12}
\end{eqnarray}
Solving gives: 
\begin{eqnarray}
&&N_N=\frac{(\rhob_N+\vepb_N^{\Nb})(L_N-\vep_N^E E_N)-\vep_N^{\Nb}(\Lb_N-\vepb_N^E E_N)}
{(\rho_N+\vep_N^N)(\rhob_N+\vepb_N^{\Nb})-\vep_N^{\Nb}\vepb_N^N} \nonumber\\
&&\Nb_N=\frac{(\rho_N+\vep_N^N)(\Lb_N-\vepb_N^E E_N)-\vepb_N^N(L_N-\vep_N^E E_N)}
{(\rho_N+\vep_N^N)(\rhob_N+\vepb_N^{\Nb})-\vep_N^{\Nb}\vepb_N^N} \nonumber\\
&&\label{9.b13}
\end{eqnarray}

\vspace{2.5mm}

We turn to the wave system, expressed in general form by \ref{2.88}, \ref{2.86}:
\begin{eqnarray}
&&d\beta_\mu\wedge dx^\mu=0 \label{9.a3}\\
&&h^{-1}(d\beta_\mu, dx^\mu)=0 \label{9.a4}
\end{eqnarray}
Now, the $N$th approximants \ref{9.3} together with \ref{9.a4} fail to satisfy this system. 
That is, defining in regard to \ref{9.a3} the 2-form (see \ref{1.57}:
\begin{equation}
\omega_N=-d\beta_N=dx_N^\mu\wedge d\beta_{\mu,N}
\label{9.a5}
\end{equation}
this 2-form does not vanish. Also, defining first $a_N$, the $N$th approximant function $a$, 
in analogy with \ref{2.72}, by:
\begin{equation}
a_N=c_N\rho_N\rhob_N
\label{9.42}
\end{equation}
and then defining in regard to \ref{9.a4} the quantity:
\begin{equation}
\delta_N=a_N h_N^{-1}(dx_N^\mu,d\beta_{\mu,N})
\label{9.a6}
\end{equation}
this quantity does not vanish. The components of $\omega_N$ in the $(L_N,\Lb_N,\Omega)$ frame are:
\begin{eqnarray}
&&\omega_{L\Lb,N}=(L_N x_N^\mu)\Lb_N\beta_{\mu,N}-(\Lb_N x_N^\mu)L_N\beta_{\mu,N} \nonumber\\
&&\omega_{L\Omega,N}=(L_N x_N^\mu)\Omega\beta_{\mu,N}-(\Omega x_N^\mu)L_N\beta_{\mu,N} \nonumber\\
&&\omega_{\Lb\Omega,N}=(\Lb_N x_N^\mu)\Omega\beta_{\mu,N}-(\Omega x_N^\mu)\Lb_N\beta_{\mu,N} 
\label{9.44}
\end{eqnarray}
The quantities $\omega_{L\Lb,N},\omega_{L\Omega,N},\omega_{\Lb\Omega,N}$ and $\delta_N$ 
are known smooth functions of $(\tau,\sigma,\vartheta)$. We shall presently estimate them. 

We begin with \ref{9.44}. Noting that by \ref{9.35}
\begin{equation}
L_N x_N^\mu=\frac{\partial x_N^\mu}{\partial\tau}
-\frac{\partial x_N^\mu}{\partial\sigma}-b_N\frac{\partial x_N^\mu}{\partial\vartheta}, 
\ \ \ 
\Lb_N x_N^\mu=\frac{\partial x_N^\mu}{\partial\sigma}
+b_N\frac{\partial x_N^\mu}{\partial\vartheta}
\label{9.a7}
\end{equation}
and
\begin{equation}
L_N\beta_{\mu,N}=\frac{\partial\beta_{\mu,N}}{\partial\tau}
-\frac{\partial\beta_{\mu,N}}{\partial\sigma}-b_N\frac{\partial\beta_{\mu,N}}{\partial\vartheta}, 
\ \ \ 
\Lb_N\beta_{\mu,N}=\frac{\partial\beta_{\mu,N}}{\partial\sigma}
+b_N\frac{\partial\beta_{\mu,N}}{\partial\vartheta}
\label{9.45}
\end{equation}
\ref{9.17} imply:
\begin{eqnarray}
&&\left.\frac{\partial^n (L_N x_N^\mu)}{\partial\tau^n}\right|_{\tau=0}=
\left.\frac{\partial^n (L x^\mu)}{\partial\tau^n}\right|_{\tau=0} \ \ \mbox{: for $n=0,...,N-1$}
\nonumber\\
&&\left.\frac{\partial^n (\Lb_N x_N^\mu)}{\partial\tau^n}\right|_{\tau=0}=
\left.\frac{\partial^n (\Lb x^\mu)}{\partial\tau^n}\right|_{\tau=0} \ \ \mbox{: for $n=0,...,N$}
\label{9.a8}
\end{eqnarray}
and:
\begin{eqnarray}
&&\left.\frac{\partial^n (L_N\beta_{\mu,N})}{\partial\tau^n}\right|_{\tau=0}=
\left.\frac{\partial^n (L\beta_\mu)}{\partial\tau^n}\right|_{\tau=0} \ \ \mbox{: for $n=0,...,N-1$}
\nonumber\\
&&\left.\frac{\partial^n (\Lb_N\beta_{\mu,N})}{\partial\tau^n}\right|_{\tau=0}=
\left.\frac{\partial^n (\Lb\beta_\mu)}{\partial\tau^n}\right|_{\tau=0} \ \ \mbox{: for $n=0,...,N$}
\label{9.46}
\end{eqnarray}
Consequently:
\begin{eqnarray}
&&\left.\frac{\partial^n\omega_{L\Lb,N}}{\partial\tau^n}\right|_{\tau=0}
=\left.\frac{\partial^n((Lx^\mu)\Lb\beta_\mu-(\Lb x^\mu)L\beta_\mu)}{\partial\tau^n}\right|_{\tau=0}=0 \ \ \mbox{: for $n=0,...,N-1$}\nonumber\\
&&\left.\frac{\partial^n\omega_{L\Omega,N}}{\partial\tau^n}\right|_{\tau=0}
=\left.\frac{\partial^n((Lx^\mu)\Omega\beta_\mu-(\Omega x^\mu)L\beta_\mu)}{\partial\tau^n}\right|_{\tau=0}=0 \ \ \mbox{: for $n=0,...,N-1$}\nonumber\\
&&\left.\frac{\partial^n\omega_{\Lb\Omega,N}}{\partial\tau^n}\right|_{\tau=0}
=\left.\frac{\partial^n((\Lb x^\mu)\Omega\beta_\mu-(\Omega x^\mu)\Lb\beta_\mu)}{\partial\tau^n}\right|_{\tau=0}=0 \ \ \mbox{: for $n=0,...,N$}\nonumber\\
&&\label{9.47}
\end{eqnarray}
The quantities $\omega_{L\Lb,N},\omega_{L\Omega,N},\omega_{\Lb\Omega,N}$ being known 
smooth functions of $(\tau,\sigma,\vartheta)$ the results \ref{9.47} yield the following 
proposition. 

\vspace{2.5mm}

\noindent{\bf Proposition 9.2} \ \ We have:
$$\omega_{L\Lb,N}=O(\tau^N), \ \ \omega_{L\Omega,N}=O(\tau^N), 
\ \ \omega_{\Lb\Omega,N}=O(\tau^{N+1})$$

We turn to \ref{9.a6}. Actually we shall first estimate
\begin{equation}
\delta^\prime_N=a_N h_N^{\prime-1}(dx_N^\mu,d\beta_{\mu,N})
\label{9.a9}
\end{equation}
where $h^\prime_N$ is another $N$th approximant metric, defined as follows. 
Consider first the exact acoustical metric $h$. It can be viewed as being given by:  
\begin{equation}
h=h_{\mu\nu}dx^\mu\otimes dx^\nu
\label{9.134}
\end{equation}
where the components $h_{\mu\nu}$ are given by \ref{1.d11}. This way of viewing $h$ gives rise 
to the $N$th approximant $h_N$ given by \ref{9.a2}. However $h$ can also be viewed in a 
different way. Given $\sh$,
we define $E$ according to \ref{3.a10} and consider the co-frame field $(\zeta,\zetab,\szeta)$ 
(basis 1-forms) dual to 
the frame field $(L,\Lb,E)$:
\begin{equation}
\zeta=d\ub, \ \ \zetab=du, \ \ \szeta=\sqrt{\sh}(d\vartheta+b(d\ub-du))
\label{9.135}
\end{equation}
(see \ref{3.a8}). Given then also $a$, 
we can view $h$ as being given by:
\begin{equation}
h=-2a(\zeta\otimes\zetab+\zetab\otimes\zeta)+\szeta\otimes\szeta
\label{9.136}
\end{equation}
Similarly, in regard to $N$th approximants, we consider the co-frame field 
$(\zeta_N,\zetab_N,\szeta_N)$ dual to the frame field $(L_N,\Lb_N,E_N)$  (see \ref{9.10}, \ref{9.11}):
\begin{equation}
\zeta_N=d\ub, \ \ \zetab_N=du, \ \ \szeta_N=\sqrt{\sh_N}(d\vartheta+b_N(d\ub-du))
\label{9.137}
\end{equation}
We then define $h^\prime_N$, a different $N$th approximant acoustical metric, by:
\begin{equation}
h^\prime_N=-2a_N(\zeta_N\otimes\zetab_N+\zetab_N\otimes\zeta_N)+\szeta_N\otimes\szeta_N
\label{9.138}
\end{equation}
($a_N$ being given by \ref{9.42}). Note that both $N$th approximant metrics $h_N$ and $h^\prime_N$ 
are smooth metrics on the manifold ${\cal N}$, that is their components in the $(\ub,u,\vartheta)$ 
coordinates are smooth functions of these coordinates. The $N$th approximant metric $h^\prime_N$ 
being adapted to the frame field $(L_N,\Lb_N,E_N)$, we have:
\begin{eqnarray}
&&h^\prime_N(L_N,L_N)=h^\prime_N(\Lb_N,\Lb_N)=0, \nonumber\\
&&h^\prime_N(L_N,E_N)=h^\prime_N(\Lb_N,E_N)=0, \nonumber\\
&&h^\prime_N(L_N,\Lb_N)=-2a_N, \ \ h^\prime_N(E_N,E_N)=1\label{9.139}
\end{eqnarray}
The inverse $h_N^{\prime-1}$ is then given by:
\begin{eqnarray}
&&h_N^{\prime-1}=-\frac{1}{2a_N}(L_N\otimes\Lb_N+\Lb_N\otimes L_N)+E_N\otimes E_N\nonumber\\
&&\hspace{8mm}=-\frac{1}{2a_N}(L_N\otimes\Lb_N+\Lb_N\otimes L_N)+\sh_N^{-1}\Omega\otimes\Omega
\label{9.a10}
\end{eqnarray}
We then have:
\begin{equation}
\delta^\prime_N=-\frac{1}{2}((L_N x_N^\mu)\Lb_N\beta_{\mu,N}+(\Lb_N x_N^\mu)L_N\beta_{\mu,N})
+a_N\sh_N^{-1}(\Omega x_N^\mu)\Omega\beta_{\mu,N}
\label{9.a11}
\end{equation}

From \ref{9.32}, \ref{9.37} we deduce, through \ref{9.42},
\begin{equation}
\left.\frac{\partial^n a_N}{\partial\tau^n}\right|_{\tau=0}=
\left.\frac{\partial^n a}{\partial\tau^n}\right|_{\tau=0} \ \ \mbox{: for $n=0,...,N-1$}
\label{9.48}
\end{equation}
Also, from \ref{9.18}, \ref{9.19} we deduce, through \ref{9.9}, 
\begin{equation}
\left.\frac{\partial^n\sh_N}{\partial\tau^n}\right|_{\tau=0}=
\left.\frac{\partial^n\sh}{\partial\tau^n}\right|_{\tau=0} \ \ \mbox{: for $n=0,...,N$}
\label{9.49}
\end{equation}
Using \ref{9.a8}, \ref{9.46}, \ref{9.48}, \ref{9.49}, we deduce:
\begin{eqnarray}
&&\left.\frac{\partial^n\delta_N^\prime}{\partial\tau^n}\right|_{\tau=0}=
\left.\frac{\partial^n}{\partial\tau^n}\left\{-\frac{1}{2}((Lx^\mu)\Lb\beta_\mu+(\Lb x^\mu)L\beta_\mu)
-a\sh^{-1}(\Omega x^\mu)\Omega\beta_\mu\right\}\right|_{\tau=0}\nonumber\\
&&\hspace{18mm}=0 \ \ \ \mbox{: for $n=0,...,N-1$}
\label{9.50}
\end{eqnarray}

The quantity $\delta_N^\prime$ being a known 
smooth function of $(\tau,\sigma,\vartheta)$ the result \ref{9.50} yields the following 
lemma. 

\vspace{2.5mm}

\noindent{\bf Lemma 9.2} \ \ We have:
$$\delta_N^\prime=O(\tau^N)$$

\vspace{2.5mm}

To estimate then $\delta_N$ what remains to be done is to estimate the difference 
$h_N^{-1}-h_N^{\prime-1}$. Now, while the components of $h_N^\prime$ in the frame field 
$(L_N,\Lb_N,E_N)$ are given by \ref{9.139}, the components of $h_N$ in the same frame field 
are given by:
\begin{eqnarray}
&&h_N(L_N,L_N)=h_{\mu\nu,N}(L_N x_N^\mu)(L_N x_N^\nu)
=h_{\mu\nu,N}(\rho_N N_N^\mu+\vep_N^\mu)(\rho_N N_N^\nu+\vep_N^\nu)\nonumber\\
&&\hspace{20mm}=2h_{\mu\nu,N}\rho_N N_N^\mu\vep_N^\nu+h_{\mu\nu,N}\vep_N^\mu\vep_N^\nu
:=\delta_{LL,N} \label{9.140}
\end{eqnarray}
\begin{eqnarray}
&&h_N(\Lb_N,\Lb_N)=h_{\mu\nu,N}(\Lb_N x_N^\mu)(\Lb_N x_N^\nu)
=h_{\mu\nu,N}(\rhob_N\Nb_N^\mu+\vepb_N^\mu)(\rhob_N\Nb_N^\nu+\vepb_N^\nu)\nonumber\\
&&\hspace{20mm}=2h_{\mu\nu,N}\rhob_N\Nb_N^\mu\vepb_N^\nu+h_{\mu\nu,N}\vepb_N^\mu\vepb_N^\nu
:=\delta_{\Lb\Lb,N} \label{9.141}
\end{eqnarray}
\begin{eqnarray}
&&h_N(L_N,E_N)=h_{\mu\nu,N}(L_N x_N^\mu)(E_N x_N^\nu)
=h_{\mu\nu,N}(\rho_N N_N^\mu+\vep_N^\mu)E_N^\nu\nonumber\\
&&\hspace{20mm}=h_{\mu\nu,N}\vep_N^\mu E_N^\nu:=\delta_{LE,N} \label{9.142}
\end{eqnarray}
\begin{eqnarray}
&&h_N(\Lb_N,E_N)=h_{\mu\nu,N}(\Lb_N x_N^\mu)(E_N x_N^\nu)
=h_{\mu\nu,N}(\rhob_N\Nb_N^\mu+\vepb_N^\mu)E_N^\nu\nonumber\\
&&\hspace{20mm}=h_{\mu\nu,N}\vepb_N^\mu E_N^\nu:=\delta_{\Lb E,N} \label{9.143}
\end{eqnarray}
\begin{eqnarray}
&&h_N(L_N,\Lb_N)=h_{\mu\nu,N}(\Lb_N x_N^\mu)(\Lb_N x_N^\nu)
=h_{\mu\nu,N}(\rho_N N_N^\mu+\vep_N^\mu)(\rhob_N\Nb_N^\mu+\vepb_N^\nu)\nonumber\\
&&\hspace{20mm}=-2a_N+h_{\mu\nu,N}(\rho_N N_N^\mu\vepb_N^\nu+\rhob_N\Nb_N^\mu\vep_N^\nu)
+h_{\mu\nu,N}\vep_N^\mu\vepb_N^\nu\nonumber\\
&&\hspace{20mm}:=-2a_N+\delta_{L\Lb,N} \label{9.144}
\end{eqnarray}
\begin{equation}
h_N(E_N,E_N)=h_{\mu\nu,N}(E_N x_N^\mu)(E_N x_N^\nu)=
\frac{h_{\mu\nu,N}\Omega_N^\mu\Omega_N^\nu}{\sh_N}=1
\label{9.145}
\end{equation}
by \ref{9.15} and \ref{9.12}, \ref{9.13}. Defining then the smooth 2-covariant symmetric tensorfield 
$\delta_{h,N}$ on ${\cal N}$ by:
\begin{eqnarray}
&&\delta_{h,N}=\delta_{LL,N}\zeta_N\otimes\zeta_N+\delta_{\Lb\Lb}\zetab_N\otimes\zetab_N
+\delta_{L\Lb,N}(\zeta_N\otimes\zetab_N+\zetab_N\otimes\zeta_N)\nonumber\\
&&\hspace{10mm}+\delta_{LE,N}(\zeta_N\otimes\szeta_N+\szeta_N\otimes\zeta_N)
+\delta_{\Lb E,N}(\zetab_N\otimes\szeta_N+\szeta_N\otimes\zetab_N)\nonumber\\
&& \label{9.146}
\end{eqnarray}
the two $N$th approximant acoustical metrics are related by:
\begin{equation}
h_N=h^\prime_N+\delta_{h,N}
\label{9.147}
\end{equation} 
Proposition 9.1 yields the following estimates for the components of $\delta_{h,N}$:
\begin{eqnarray}
&&\delta_{LL,N}=O(\tau^{N+1}), \ \ \delta_{\Lb\Lb,N}=\rhob_N O(\tau^{N+1}), \ \ 
\delta_{L\Lb,N}=O(\tau^{N+2})+\rhob_N O(\tau^N) \nonumber\\
&&\delta_{LE,N}=O(\tau^N), \ \ \delta_{\Lb E,N}=O(\tau^{N+1}) \ \ \ \ \ \ \ \ (\delta_{EE,N}=0)
\label{9.a12}
\end{eqnarray}
We note that the quadratic error terms are absorbed because $N$ shall be chosen sufficiently large.
Here we have used the fact that since $\rho$ vanishes on $\Cb_0$ we have:
\begin{equation}
\rho_N=O(\tau)
\label{9.a13}
\end{equation}
In fact by the 1st of \ref{4.216}:
\begin{equation}
\left.\frac{\partial\rho}{\partial\tau}\right|_{(\tau,\sigma)=(0,0)}=-\frac{k}{3l}
\label{9.a14}
\end{equation}
which is positive. This implies that under a suitable restriction on $u=\tau+\sigma$ there is 
a positive constant $C$ such that:
\begin{equation}
\rho_N\geq C^{-1}\tau
\label{9.a15}
\end{equation}
In the estimates \ref{9.a12} we have kept the factor $\rhob_N$ where appropriate. We have:
\begin{equation}
\rhob_N=\left.\rhob\right|_{\tau=0}+\tau\left.\frac{\partial\rhob}{\partial\tau}\right|_{\tau=0}
+\frac{\tau^2}{2}\left.\frac{\partial^2\rhob}{\partial\tau^2}\right|_{\tau=0}+O(\tau^3)
\label{9.a16}
\end{equation}
By \ref{4.106}, \ref{4.107}, \ref{4.130} and the 2nd of \ref{4.216}:
\begin{eqnarray}
&&\left.\rhob\right|_{(\tau,\sigma)=(0,0)}=0, \ \ 
\left.\frac{\partial\rhob}{\partial\sigma}\right|_{(\tau,\sigma)=(0,0)}=0, \ \ 
\left.\frac{\partial^2\rhob}{\partial\sigma^2}\right|_{(\tau,\sigma)=(0,0)}=\frac{k}{2c_0} \nonumber\\
&&\left.\frac{\partial\rhob}{\partial\tau}\right|_{(\tau,\sigma)=(0,0)}=0, \ \ 
\left.\frac{\partial^2\rhob}{\partial\tau\partial\sigma}\right|_{(\tau,\sigma)=(0,0)}=\frac{k}{2c_0}
\nonumber\\
&&\left.\frac{\partial^2\rhob}{\partial\tau^2}\right|_{(\tau,\sigma)=(0,0)}=\frac{k}{3c_0}
\label{9.a17}
\end{eqnarray}
Hence:
\begin{equation}
\left.\rhob\right|_{\tau=0}+\tau\left.\frac{\partial\rhob}{\partial\tau}\right|_{\tau=0}
+\frac{\tau^2}{2}\left.\frac{\partial^2\rhob}{\partial\tau^2}\right|_{\tau=0}
=\frac{k}{12c_0}(3\sigma^2+6\sigma\tau+2\tau^2)+O(u^3)
\label{9.a18}
\end{equation}
Substituting in \ref{9.a16} we conclude that:
\begin{equation}
\rhob_N=\frac{k}{12c_0}(3u^2-\ub^2)+O(u^3)
\label{9.a19}
\end{equation}
which implies that under a suitable restriction on $u=\tau+\sigma$ there is 
a positive constant $C$ such that:
\begin{equation}
\rhob_N\geq C^{-1}u^2
\label{9.a20}
\end{equation}
Now \ref{9.a16} implies that if $m,n,l$ are non-negative integers there are constants $C_{m,n,l}$ 
such that:
\begin{equation}
\left|\frac{\partial^{m+n+l}\rhob_N}{\partial\tau^m\partial\sigma^n\partial\vartheta^l}\right|
\leq C_{m,n,l}\left\{\begin{array}{lll} u^2 \ &:& \ \mbox{if $m+n=0$}\\
u \ &:& \ \mbox{if $m+n=1$}\\ 
1 \ &:& \ \mbox{if $m+n\geq 2$} \end{array}\right.
\label{9.a21}
\end{equation}
This together with \ref{9.a20} implies that if $m,n,l$ are non-negative integers there are constants 
$C_{m,n,l}$ such that:
\begin{equation}
\left|\frac{\partial^{m+n+l}(\rhob_N^{-1})}{\partial\tau^m\partial\sigma^n\partial\vartheta^l}\right|
\leq C_{m,n,l} \ u^{-2-m-n}
\label{9.a22}
\end{equation}
On the other hand, \ref{9.a13} together with \ref{9.a15} implies that 
if $m,n,l$ are non-negative integers there are constants 
$C_{m,n,l}$ such that:
\begin{equation}
\left|\frac{\partial^{m+n+l}(\rho_N^{-1})}{\partial\tau^m\partial\sigma^n\partial\vartheta^l}\right|
\leq C_{m,n,l} \ \tau^{-1-m}
\label{9.a23}
\end{equation}
In fact for any positive integer $p$ we can write:
\begin{equation}
\rho_N^{-1}O(\tau^p)=O(\tau^{p-1})
\label{9.a24}
\end{equation}
in the sense of Definition 9.1.

Let us now define the symmetric 2-contravariant tensorfield  $\delta^\prime_{h,N}$ by:
\begin{equation}
h_N^{-1}=h_N^{\prime-1}+\delta^\prime_{h,N}
\label{9.a25}
\end{equation}
In terms of rectangular components \ref{9.a10} reads:
\begin{eqnarray}
&&(h_N^{\prime-1})^{\mu\nu}=-\frac{1}{2a_N}(L_N^\mu\Lb_N^\nu+\Lb_N^\mu L_N^\nu)+E_N^\mu E_N^\nu \nonumber\\
&&=-\frac{1}{2a_N}\left((\rho_N N_N^\mu+\vep_N^\mu)(\rhob_N\Nb_N^\mu+\vepb_N^\nu)
+(\rhob_N\Nb_N^\mu+\vepb_N^\mu)(\rho_N N_N^\nu+\vep_N^\nu)\right)\nonumber\\
&&\hspace{40mm}+E_N^\mu E_N^\nu \label{9.b14}
\end{eqnarray}
the last by \ref{9.b10}. Comparing with \ref{9.b7} we then obtain for the rectangular components 
of $\delta_{h,N}^{\prime \mu\nu}$:
\begin{eqnarray}
&&\delta_{h,N}^{\prime \mu\nu}=\frac{1}{2c_N}\left\{
\frac{1}{\rho_N}(\vep_N^\mu\Nb_N^\nu+\Nb_N^\mu\vep_N^\nu)+
\frac{1}{\rhob_N}(N_N^\mu\vepb_N^\nu+\vepb_N^\mu N_N^\nu)\right.\nonumber\\
&&\hspace{30mm}\left.+\frac{1}{\rho_N\rhob_N}(\vep_N^\mu\vepb_N^\nu+\vepb_N^\mu\vep_N^\nu)\right\}
\label{9.b15}
\end{eqnarray}
Substituting the expansions \ref{9.b11} gives the following expansion for $\delta^\prime_{h,N}$ 
in the frame field $N_N$, $\Nb_N$, $E_N$:
\begin{eqnarray}
&&\delta^\prime_{h,N}=\delta_N^{\prime NN}N_N\otimes N_N+\delta_N^{\prime\Nb\Nb}\Nb_N\otimes\Nb_N
+\delta_N^{\prime N\Nb}(N_N\otimes\Nb_N+\Nb_N\otimes N_N) \nonumber\\
&&\hspace{10mm}+\delta_N^{\prime NE}(N_N\otimes E_N+E_N\otimes N_N)
+\delta_N^{\prime \Nb E}(\Nb_N\otimes E_N+E_N\otimes\Nb_N) \nonumber\\
&&\hspace{10mm}+\delta_N^{\prime EE}E_N\otimes E_N
\label{9.b16}
\end{eqnarray}
with:
\begin{eqnarray}
&&\delta_N^{\prime NN}=\frac{\vepb_N^N}{c_N\rhob_N}\left(1+\frac{\vep_N^N}{\rho_N}\right)=
\rhob_N^{-1}O(\tau^{N+1})\nonumber\\
&&\delta_N^{\prime\Nb\Nb}=\frac{\vep_N^{\Nb}}{c_N\rho_N}\left(1+\frac{\vepb_N^{\Nb}}{\rhob_N}\right)=
O(\tau^{N-1})\nonumber\\
&&\delta_N^{\prime N\Nb}=\frac{1}{2c_N}\left(\frac{\vep_N^N}{\rho_N}+\frac{\vepb_N^{\Nb}}{\rhob_N}
+\frac{(\vep_N^N\vepb_N^{\Nb}+\vepb_N^N\vep_N^{\Nb})}{\rho_N\rhob_N}\right)=O(\tau^{N-1})+\rhob_N^{-1}O(\tau^{N+1})\nonumber\\
&&\delta_N^{\prime NE}=\frac{1}{2c_N\rhob_N}\left(\vepb_N^E+
\frac{(\vep_N^N\vepb_N^E+\vepb_N^N\vep_N^E)}{\rho_N}\right)=\rhob_N^{-1}O(\tau^{N+1})\nonumber\\
&&\delta_N^{\prime\Nb E}=\frac{1}{2c_N\rho_N}\left(\vep_N^E+
\frac{(\vep_N^{\Nb}\vep_N^E+\vepb_N^{\Nb}\vep_N^E)}{\rhob_N}\right)=O(\tau^{N-1})\nonumber\\
&&\delta_N^{\prime EE}=\frac{\vep_N^E\vepb_N^E}{c_N\rho_N\rhob_N}=\rhob_N^{-1}O(\tau^{2N}) 
\label{9.b17}
\end{eqnarray}
Substituting for $N_N$, $\Nb_N$ from \ref{9.b13} in \ref{9.b16} we obtain the expansion of
$\delta^\prime_{h,N}$ in the frame field $L_N$, $\Lb_N$, $E_N$:
\begin{eqnarray}
&&\delta^\prime_{h,N}=\delta_N^{\prime LL}L_N\otimes L_N+\delta_N^{\prime\Lb\Lb}\Lb_N\otimes\Lb_N
+\delta_N^{\prime L\Lb}(L_N\otimes\Lb_N+\Lb_N\otimes L_N) \nonumber\\
&&\hspace{10mm}+\delta_N^{\prime LE}(L_N\otimes E_N+E_N\otimes L_N)
+\delta_N^{\prime \Lb E}(\Lb_N\otimes E_N+E_N\otimes\Lb_N) \nonumber\\
&&\hspace{10mm}+\delta_N^{\prime EE}E_N\otimes E_N
\label{9.a26}
\end{eqnarray}
and \ref{9.b17} give: 
\begin{eqnarray}
&&\delta_N^{\prime LL}=\rhob_N^{-1}O(\tau^{N-1}), \ \ \
\delta_N^{\prime \Lb\Lb}=\rhob_N^{-2}O(\tau^{N-1})\nonumber\\
&&\hspace{8mm} \delta_N^{\prime L\Lb}=\rhob_N^{-1}O(\tau^{N-2})+\rhob_N^{-2}O(\tau^N)\nonumber\\
&&\delta_N^{\prime LE}=\rhob_N^{-1} O(\tau^N), \ \ \ 
\delta_N^{\prime\Lb E}=\rhob_N^{-1} O(\tau^{N-1}) \nonumber\\
&&\hspace{12mm}\delta_N^{\prime EE}=\rhob_N^{-1}O(\tau^{2N})
\label{9.a29}
\end{eqnarray}

We now consider:
$$a_N \delta^\prime_{h,N}(dx_N^\mu,d\beta_{\mu,N})$$ 
We write:
$$(L_N x_N^\mu)\Lb_N\beta_{\mu,N}+(\Lb_N x^\mu)L_N\beta_{\mu,N}=2(L_N x_N^\mu)\Lb_N\beta_{\mu,N}
-\omega_{L\Lb,N},$$ 
$$(L_N x_N^\mu)E_N\beta_{\mu,N}+(E_N x_N^\mu)L_N\beta_{\mu,N}=2(L_N x_N^\mu)E_N\beta_{\mu,N}
-\omega_{LE,N},$$
where $\omega_{LE,N}=\sh_N^{-1/2}\omega_{L\Omega,N}$,
$$(\Lb_N x_N^\mu)E_N\beta_{\mu,N}+(E_N x_N^\mu)\Lb_N\beta_{\mu,N}=2(\Lb_N x_N^\mu)E_N\beta_{\mu,N}
-\omega_{\Lb E,N},$$
where $\omega_{\Lb E,N}=\sh_N^{-1/2}\omega_{\Lb\Omega,N}$, 
and we express $L_N x_N^\mu$, $\Lb_N x_N^\mu$ using \ref{9.15} ($\vep_N^0=\vepb_N^0=0$). 
This gives:
\begin{eqnarray}
&&a_N \delta^\prime_{h,N}(dx_N^\mu,d\beta_{\mu,N})=\nonumber\\
&&\hspace{20mm} a_N\left\{\delta_N^{\prime LL}(\rho_N N_N^\mu+\vep_N^\mu)L_N\beta_{\mu,N}
+\delta_N^{\prime\Lb\Lb}(\rhob_N\Nb_N^\mu+\vepb_N^\mu)\Lb_N\beta_{\mu,N}\right.\nonumber\\
&&\hspace{27mm}+\delta_N^{\prime L\Lb}\left(2(\rho_N N_N^\mu+\vep_N^\mu)\Lb_N\beta_{\mu,N}-\omega_{L\Lb,N}\right)
\nonumber\\
&&\hspace{27mm}+\delta_N^{\prime LE}\left(2(\rho_N N_N^\mu+\vep_N^\mu)E_N\beta_{\mu,N}-\omega_{LE,N}\right)
\nonumber\\
&&\hspace{27mm}+\delta_N^{\prime\Lb E}\left(2(\rhob_N\Nb_N^\mu+\vepb_N^\mu)E_N\beta_{\mu,N}-\omega_{\Lb E,N}\right)\nonumber\\
&&\hspace{40mm}\left.+\delta_N^{\prime EE}(E_N x_N^\mu)E_N\beta_{\mu,N}\right\}
\label{9.a30}
\end{eqnarray}
Substituting the estimates \ref{9.a29}, by virtue of Propositions 9.1 and 9.2 the quadratic error terms 
are absorbed, and we obtain:
\begin{equation}
a_N \delta^\prime_{h,N}(dx_N^\mu,d\beta_{\mu,N})=O(\tau^N)+\rhob_N^{-1}O(\tau^{N+2})
\label{9.a31}
\end{equation}
Combining finally this result with Lemma 9.2 (see \ref{9.a6}, \ref{9.a25}) we deduce the following. 

\vspace{2.5mm}

\noindent{\bf Proposition 9.3} We have:
$$\delta_N=O(\tau^N)+\rhob_N^{-1}O(\tau^{N+2})$$

\vspace{5mm}

\section{Estimates for the Quantities by which the $N$th Approximants fail to satisfy the Boundary 
Conditions}

In connection with the $N$th approximants $x^\mu_N$ we define, in analogy with \ref{4.53}, in 
terms of $(\tau,\sigma,\vartheta)$ coordinates:
\begin{equation}
f_N(\tau,\vartheta)=x_N^0(\tau,0,\vartheta), \ \ \ g_N^i(\tau,\vartheta)=x_N^i(\tau,0,\vartheta) \ 
: \ i=1,2
\label{9.51}
\end{equation}
Having determined, in Section 5.3, 
\begin{equation}
\left.\frac{\partial^n v}{\partial\tau^n}\right|_{\tau=0}, \ \ 
\left.\frac{\partial^n\gamma}{\partial\tau^n}\right|_{\tau=0} \ \ \mbox{: for $n=0,...,N-2$}
\label{9.52}
\end{equation}
we define the $N$th order approximants of the transformation functions by:
\begin{eqnarray}
&&v_N(\tau,\vartheta)=\sum_{n=0}^{N-2}\frac{\tau^n}{n!}\frac{\partial^n v}{\partial\tau^n}(0,\vartheta)
\nonumber\\
&&\gamma_N(\tau,\vartheta)=\sum_{n=0}^{N-2}\frac{\tau^n}{n!}\frac{\partial^n\gamma}{\partial\tau^n}
(0,\vartheta) \label{9.53}
\end{eqnarray}
and (see \ref{4.53}, \ref{4.60} and Proposition 4.4):
\begin{eqnarray}
&&w_N(\tau,\vartheta)=\tau v_N(\tau,\vartheta) \nonumber\\
&&\psi_N(\tau,\vartheta)=\vartheta+\tau^3\gamma_N(\tau,\vartheta) 
\label{9.54}
\end{eqnarray}
We then set, in analogy with \ref{4.149}, 
\begin{equation}
\triangle_N\beta_\mu(\tau,\vartheta)=\beta_{\mu,N}(\tau,0,\vartheta)
-\beta^\prime_\mu(f_N(\tau,\vartheta),w_N(\tau,\vartheta), \psi_N(\tau,\vartheta))
\label{9.55}
\end{equation}
where the $\beta_{\mu,N}$  are in terms of $(\tau,\sigma,\vartheta)$ coordinates, while the 
$\beta^\prime_\mu$, which refer to the prior solution, are known smooth functions of 
$(t,u^\prime,\vartheta^\prime)$. 

Let us recall the linear boundary conditions \ref{4.1}, \ref{4.2}:
\begin{eqnarray}
&&\Omega^\mu\triangle\beta_\mu=0 \label{9.56}\\
&&T^\mu\triangle\beta_\mu=0 \label{9.57}
\end{eqnarray}
The $N$th approximants of $\Omega^\mu$ and $T^\mu$ are $\Omega_N^\mu$ given by \ref{9.8} and 
$T_N^\mu$ given similarly by:
\begin{equation}
T_N^\mu=Tx_N^\mu=\frac{\partial x_N^\mu}{\partial\tau}
\label{9.58}
\end{equation}
Now, the $N$th approximants \ref{9.55} together with $\Omega_N^\mu,T_N^\mu$ do not satisfy the boundary 
conditions \ref{9.56}, \ref{9.57}. That is, the quantities 
\begin{eqnarray}
&&\Omega_N^\mu\triangle_N\beta_\mu:=\iota_{\Omega,N} \label{9.59}\\
&&T_N^\mu\triangle_N\beta_\mu:=\iota_{T,N} \label{9.60}
\end{eqnarray}
do not vanish. By the above, the quantities $\iota_{\Omega,N},\iota_{T,N}$ are known smooth functions 
of $(\tau,\vartheta)$. We shall presently estimate them. 

By virtue of the 3rd of the definitions \ref{9.3} we have, in $(\tau,\sigma,\vartheta)$ coordinates,
\begin{equation}
\left.\frac{\partial^n(\beta_{\mu,N}(\tau,0,\vartheta))}{\partial\tau^n}\right|_{\tau=0}=
\left.\frac{\partial^n(\beta_\mu(\tau,0,\vartheta))}{\partial\tau^n}\right|_{\tau=0} 
\ \ \mbox{: for $n=0,...,N$}
\label{9.61}
\end{equation}
By virtue of \ref{9.51} and the 1st of the definitions \ref{9.3} we have: 
\begin{eqnarray}
&&\left.\frac{\partial^n f_N}{\partial\tau^n}\right|_{\tau=0}=
\left.\frac{\partial^n f}{\partial\tau^n}\right|_{\tau=0} \ \ \mbox{: for $n=0,...,N$} \nonumber\\
&&\left.\frac{\partial^n g_N^i}{\partial\tau^n}\right|_{\tau=0}=
\left.\frac{\partial^n g^i}{\partial\tau^n}\right|_{\tau=0} \ \ \mbox{: for $n=0,...,N$} \label{9.62}
\end{eqnarray}
By virtue of the definitions \ref{9.53} we have:
\begin{eqnarray}
&&\left.\frac{\partial^n v_N}{\partial\tau^n}\right|_{\tau=0}=
\left.\frac{\partial^n v}{\partial\tau^n}\right|_{\tau=0} \ \ \mbox{: for $n=0,...,N-2$} 
\nonumber\\
&&\left.\frac{\partial^n\gamma_N}{\partial\tau^n}\right|_{\tau=0}=
\left.\frac{\partial^n\gamma}{\partial\tau^n}\right|_{\tau=0} \ \ \mbox{: for $n=0,...,N-2$} 
\label{9.63}
\end{eqnarray}
It then follows from the definitions \ref{9.54} that:
\begin{eqnarray}
&&\left.\frac{\partial^n w_N}{\partial\tau^n}\right|_{\tau=0}=
\left.\frac{\partial^n w}{\partial\tau^n}\right|_{\tau=0} \ \ \mbox{: for $n=0,...,N-1$} \nonumber\\
&&\left.\frac{\partial^n\psi_N}{\partial\tau^n}\right|_{\tau=0}=
\left.\frac{\partial^n\psi}{\partial\tau^n}\right|_{\tau=0} \ \ \mbox{: for $n=0,...,N+1$} 
\label{9.64}
\end{eqnarray}
The 1st of \ref{9.64} holds also for $n=N-1$ because 
$\left.\partial^{N-1}v_N/\partial\tau^{N-1}\right|_{\tau=0}$, 
$\left.\partial^{N-1}v/\partial\tau^{N-1}\right|_{\tau=0}$ are not involved. Also, 
the 2nd of \ref{9.64} holds also for $n=N-1, \ N, N+1$ because neither
$\left.\partial^{N-1}\gamma_N/\partial\tau^{N-1}\right|_{\tau=0}$, 
$\left.\partial^{N-1}\gamma/\partial\tau^{N-1}\right|_{\tau=0}$, nor 
$\left.\partial^N\gamma_N/\partial\tau^N\right|_{\tau=0}$, 
$\left.\partial^N\gamma/\partial\tau^N\right|_{\tau=0}$, or 
$\left.\partial^{N+1}\gamma_N/\partial\tau^{N+1}\right|_{\tau=0}$, 
$\left.\partial^{N+1}\gamma/\partial\tau^{N+1}\right|_{\tau=0}$ are involved. 

Now, \ref{9.64} together with the 1st of \ref{9.62} imply:
\begin{eqnarray}
&&\left.\frac{\partial^n}{\partial\tau^n}\beta^\prime_\mu(f_N(\tau,\vartheta),w_N(\tau,\vartheta),\psi_N(\tau,\vartheta))\right|_{\tau=0} \label{9.65}\\
&&\hspace{10mm}=\left.\frac{\partial^n}{\partial\tau^n}(\beta^\prime_\mu(f(\tau,\vartheta),w(\tau,\vartheta),\psi(\tau,\vartheta))\right|_{\tau=0} \ \ \ \mbox{: for $n=0,...,N-1$} \nonumber
\end{eqnarray}
In the case $n=N$, since from \ref{9.53} $\left.\partial^{N-1}v_N/\partial\tau^{N-1}\right|_{\tau=0}$ 
vanishes, while $\left.\partial^{N-1}v/\partial\tau^{N-1}\right|_{\tau=0}$ may not vanish, we have:
\begin{equation}
\left.\frac{\partial^N w_N}{\partial\tau^N}\right|_{\tau=0}-
\left.\frac{\partial^N w}{\partial\tau^N}\right|_{\tau=0}=
-N\left.\frac{\partial^{N-1}v}{\partial\tau^{N-1}}\right|_{\tau=0}
\label{9.66}
\end{equation}
hence:
\begin{eqnarray}
&&\left.\frac{\partial^N}{\partial\tau^N}\beta^\prime_\mu(f_N(\tau,\vartheta),w_N(\tau,\vartheta),\psi_N(\tau,\vartheta))\right|_{\tau=0} \nonumber\\
&&\hspace{10mm}-\left.\frac{\partial^N}{\partial\tau^N}(\beta^\prime_\mu(f(\tau,\vartheta),w(\tau,\vartheta),\psi(\tau,\vartheta))\right|_{\tau=0}\nonumber\\
&&\hspace{20mm}=-N\frac{\partial\beta^\prime}{\partial u^\prime}(f(0,\vartheta),0,\vartheta)
\frac{\partial^{N-1}v}{\partial\tau^{N-1}}(0,\vartheta)
\label{9.67}
\end{eqnarray}
Equations \ref{9.61} together with \ref{9.65}, \ref{9.67} imply, through \ref{9.55}, 
\begin{equation}
\left.\frac{\partial^n\triangle_N\beta_\mu}{\partial\tau^n}\right|_{\tau=0}=
\left.\frac{\partial^n\triangle\beta_\mu}{\partial\tau^n}\right|_{\tau=0} \ \ 
\mbox{: for $n=0,...,N-1$}
\label{9.68}
\end{equation}
and: 
\begin{equation}
\left.\frac{\partial^N\triangle_N\beta_\mu}{\partial\tau^N}\right|_{\tau=0}=
\left.\frac{\partial^N\triangle\beta_\mu}{\partial\tau^N}\right|_{\tau=0}
+N\left.\frac{\partial\beta^\prime_\mu}{\partial u^\prime}\right|_{\partial_-{\cal B}}
\left.\frac{\partial^{N-1}v}{\partial\tau^{N-1}}\right|_{\tau=0}
\label{9.69}
\end{equation}
recalling that $\partial_-{\cal B}$ is represented in $(t,u^\prime,\vartheta^\prime)$ coordinates as 
the graph:
\begin{equation}
\partial_-{\cal B}=\{ (f(0,\vartheta),0,\vartheta) \ : \ \vartheta\in S^1\}
\label{9.70}
\end{equation}

Now, \ref{9.19} together with \ref{9.68} imply (see \ref{9.59}:
\begin{equation}
\left.\frac{\partial^n\iota_{\Omega,N}}{\partial\tau^n}\right|_{\tau=0}=
\left.\frac{\partial^n(\Omega^\mu\triangle\beta_\mu)}{\partial\tau^n}\right|_{\tau=0}=0 \ \ 
\mbox{: for $n=0,...,N-1$}
\label{9.71}
\end{equation}
Moreover, since by Proposition 4.3, \ref{4.102} and \ref{4.103}:
\begin{equation}
\left.\left(\Omega^\mu\frac{\partial\beta^\prime_\mu}{\partial u^\prime}\right)
\right|_{\partial_-{\cal B}}=
\left.\left((\Omega^{\prime\mu}+\pi L^{\prime\mu})T^\prime\beta^\prime_\mu\right)
\right|_{\partial_-{\cal B}}=\left.\left(T^{\prime\mu}(\Omega^\prime\beta^\prime_\mu
+\pi L^\prime\beta^\prime_\mu)\right)\right|_{\partial_-{\cal B}}=0
\label{9.72}
\end{equation}
in view of \ref{9.69} we also obtain:
\begin{equation}
\left.\frac{\partial^N\iota_{\Omega,N}}{\partial\tau^N}\right|_{\tau=0}=
\left.\frac{\partial^N(\Omega^\mu\triangle\beta_\mu)}{\partial\tau^N}\right|_{\tau=0}=0
\label{9.73}
\end{equation}

Next, by \ref{9.51}, \ref{9.58}, along ${\cal K}$ we have:
\begin{equation}
T_N^0=\frac{\partial f_N}{\partial\tau}, \ \ \ T_N^i=\frac{\partial g_N^i}{\partial\tau} \ : \ i=1,2
\label{9.74}
\end{equation}
hence along ${\cal K}$, by \ref{9.62}: 
\begin{equation}
\left.\frac{\partial^n T_N^\mu}{\partial\tau^n}\right|_{\tau=0}=
\left.\frac{\partial^n T^\mu}{\partial\tau^n}\right|_{\tau=0} \ \ \mbox{: for $n=0,...,N-1$}
\label{9.75}
\end{equation}
Noting that: 
\begin{equation}
\left.\triangle\beta_\mu\right|_{\tau=0}=0
\label{9.76}
\end{equation}
and that by \ref{4.120}:
\begin{equation}
\left.T^\mu\right|_{\tau=0}=0
\label{9.77}
\end{equation}
\ref{9.68} and \ref{9.75} imply:
\begin{equation}
\left.\frac{\partial^n\iota_{T,N}}{\partial\tau^n}\right|_{\tau=0}=
\left.\frac{\partial^n(T^\mu\triangle\beta_\mu)}{\partial\tau^n}\right|_{\tau=0}=0 \ \ 
\mbox{: for $n=0,...,N$}
\label{9.78}
\end{equation}

The quantities $\iota_{\Omega,N}$, $\iota_{T,N}$ being known smooth functions of 
$(\tau,\vartheta)$ equations \ref{9.71}, \ref{9.73}, and \ref{9.78} yield the following.

\vspace{2.5mm}

\noindent{\bf Proposition 9.4} \ \ We have:
$$\iota_{\Omega,N}=O(\tau^{N+1}), \ \ \ \iota_{T,N}=O(\tau^{N+1})$$

\vspace{2.5mm}

Let us define, in analogy with \ref{4.3}:
\begin{equation}
\ep_N=N_N^\mu\triangle_N\beta_\mu, \ \ \ \epb_N=\Nb_N^\mu\triangle_N\beta_\mu 
\label{9.79}
\end{equation}
Applying \ref{9.81}, \ref{9.82} taking 
$$V^\mu=(h_N^{-1})^{\mu\nu}\triangle_N\beta_\nu$$
we have, by \ref{9.59} and the definitions \ref{9.79},
\begin{eqnarray*}
&&V_N^E=\frac{\iota_{\Omega,N}}{\sqrt{\sh_N}}\\
&&V_N^N=-\frac{\epb_N}{2c_N}, \ \ V_N^{\Nb}=-\frac{\ep_N}{2c_N}
\end{eqnarray*}
we conclude that:
\begin{equation}
\triangle_N\beta_\mu=h_{\mu\nu,N}\left\{-\frac{1}{2c_N}(\epb_N N_N^\nu+\ep_N\Nb_N^\nu)
+\frac{\Omega_N^\nu}{\sh_N}\iota_{\Omega,N}\right\}
\label{9.83}
\end{equation}
Multiplying this with $T_N^\mu$ the left hand side becomes $\iota_{T,N}$ according to \ref{9.60} 
and we obtain:
\begin{equation}
\iota_{T,N}=h_{\mu\nu,N}T_N^\mu\left\{-\frac{1}{2c_N}(\epb_N N_N^\nu+\ep_N\Nb_N^\nu)
+\frac{\Omega_N^\nu}{\sh_N}\iota_{\Omega,N}\right\}
\label{9.84}
\end{equation}
Since from \ref{9.35}, \ref{9.58}, and \ref{9.15}:
\begin{eqnarray}
&&T_N^\mu=Tx_N^\mu=L_N x_N^\mu+\Lb_N x_N^\mu \nonumber\\
&&\hspace{5mm}=\rho_N N_N^\mu+\rhob_N\Nb_N^\mu+\vep_N^\mu+\vepb_N^\mu \label{9.85}
\end{eqnarray}
where (see \ref{9.14}):
\begin{equation}
\vep_N^0=\vepb_N^0=0
\label{9.86}
\end{equation}
we then deduce:
\begin{equation}
\iota_{T,N}=\nu_N+h_{\mu\nu,N}(\vep_N^\mu+\vepb_N^\mu)
\left\{-\frac{1}{2c_N}(\epb_N N_N^\nu+\ep_N\Nb_N^\nu)+\frac{\Omega_N^\nu}{\sh_N}\iota_{\Omega,N}\right\}
\label{9.87}
\end{equation}
where $\nu_N$ is the quantity:
\begin{equation}
\nu_N=\epb_N\rhob_N+\ep_N\rho_N
\label{9.88}
\end{equation}
In view of \ref{9.87}, Propositions 9.1 and 9.4 together with the fact that:
\begin{equation}
\epb_N=O(\tau), \ \ \ \ep_N=O(\tau^2)
\label{9.89}
\end{equation}
imply:
\begin{equation}
\nu_N=O(\tau^{N+1})
\label{9.90}
\end{equation}

However, a better estimate for $\nu_N$ can in fact be derived. From the definitions \ref{9.79}, 
by Lemma 9.1 and \ref{9.68} we obtain:
\begin{eqnarray}
&&\left.\frac{\partial^n\ep_N}{\partial\tau^n}\right|_{\tau=0}=
\left.\frac{\partial^n\ep}{\partial\tau^n}\right|_{\tau=0}, \ \ \ 
\left.\frac{\partial^n\epb_N}{\partial\tau^n}\right|_{\tau=0}=
\left.\frac{\partial^n\epb}{\partial\tau^n}\right|_{\tau=0} \nonumber\\
&&\hspace{20mm}\mbox{: for $n=0,...,N-1$} \label{9.91}
\end{eqnarray}
Moreover, since by Proposition 4.3, \ref{4.102}, \ref{4.103}:
\begin{equation}
\left.\left(N^\mu\frac{\partial\beta^\prime_\mu}{\partial u^\prime}\right)\right|_{\partial_-{\cal B}}
=\left.\left(L^{\prime\mu}T^\prime\beta^\prime_\mu\right)\right|_{\partial_-{\cal B}}=
\left.\left(T^{\prime\mu}L^\prime\beta^\prime_\mu\right)\right|_{\partial_-{\cal B}}=0
\label{9.92}
\end{equation}
in view of \ref{9.69} we also obtain:
\begin{equation}
\left.\frac{\partial^N\ep_N}{\partial\tau^N}\right|_{\tau=0}=
\left.\frac{\partial^N\ep}{\partial\tau^N}\right|_{\tau=0}
\label{9.93}
\end{equation}
Then, in view of the facts that:
\begin{equation}
\left.\rho\right|_{\tau=0}=0, \ \ \ \left.\ep\right|_{\tau=0}=\left.\frac{\partial\ep}{\partial\tau}\right|_{\tau=0}=0
\label{9.94}
\end{equation}
\ref{9.91}, \ref{9.93} together with \ref{9.37} imply:
\begin{equation}
\left.\frac{\partial^n(\rho_N\ep_N)}{\partial\tau^n}\right|_{\tau=0}=
\left.\frac{\partial^n(\rho\ep)}{\partial\tau^n}\right|_{\tau=0} \ \ \mbox{: for $n=0,...,N+1$}
\label{9.95}
\end{equation}
For, by \ref{9.94}, for a non-vanishing term at least two derivatives must fall on $\ep_N,\ep$ leaving at most $N-1$ 
derivatives to fall on $\rho_N,\rho$, and at least one derivative must fall on $\rho_N,\rho$ leaving 
at most $N$ derivatives to fall on $\ep_N,\ep$. Similarly, in view of the facts that:
\begin{equation}
\left.\rhob\right|_{\tau=0}=\left.\frac{\partial\rhob}{\partial\tau}\right|_{\tau=0}=0, \ \ \ 
\left.\epb\right|_{\tau=0}=0
\label{9.96}
\end{equation}
\ref{9.91} together with \ref{9.37} imply:
\begin{equation}
\left.\frac{\partial^n(\rhob_N\epb_N)}{\partial\tau^n}\right|_{\tau=0}=
\left.\frac{\partial^n(\rhob\epb)}{\partial\tau^n}\right|_{\tau=0} \ \ \mbox{: for $n=0,...,N+1$}
\label{9.97}
\end{equation}
For, by \ref{9.96}, for a non-vanishing term at least two derivatives must fall on $\rhob_N,\rhob$ 
leaving at most $N-1$ derivatives to fall on $\epb$, and at least one derivative must fall on 
$\epb_N,\epb$ leaving at most $N$ derivatives to fall on $\rhob_N,\rhob$. Adding \ref{9.95} and 
\ref{9.97} yields:
\begin{equation}
\left.\frac{\partial^n\nu_N}{\partial\tau^n}\right|_{\tau=0}=
\left.\frac{\partial^n(\rho\ep+\rhob\epb)}{\partial\tau^n}\right|_{\tau=0}=0 \ \ 
\mbox{: for $n=0,...,N+1$}
\label{9.98}
\end{equation}
(see \ref{4.5}). The quantity $\nu_N$ being a known smooth function of $(\tau,\vartheta)$ this 
yields the following proposition, which improves the estimate \ref{9.90}.

\vspace{2.5mm}

\noindent{\bf Proposition 9.5} \ \ We have:
$$\nu_N=O(\tau^{N+2})$$

\vspace{2.5mm}

Let us define, in analogy with the 2nd of \ref{4.5}, 
\begin{equation}
r_N=-\frac{\ep_N}{\epb_N}
\label{9.99}
\end{equation}
Let us also define:
\begin{equation}
\hat{\nu}_N=\frac{c_N\nu_N}{\epb_N}
\label{9.100}
\end{equation}
Also, in analogy with \ref{2.74},
\begin{equation}
\lambda_N=c_N\rhob_N, \ \ \ \lambdab_N=c_N\rho_N
\label{9.101}
\end{equation}
The definition \ref{9.88} can then we written in the form:
\begin{equation}
r_N\lambdab_N=\lambda_N-\hat{\nu}_N
\label{9.102}
\end{equation}
Comparing with \ref{4.113} we see that $\hat{\nu}_N$ is the quantity by which the boundary condition 
\ref{4.113} fails to be satisfied by the $N$th approximants. Now by \ref{4.155} and \ref{4.214}:
\begin{equation}
\left.\frac{\partial\epb}{\partial\tau}\right|_{\tau=0}=2\left.s_{\Nb\Lb}\right|_{\partial_-{\cal B}}
\label{9.103}
\end{equation}
and by \ref{4.163} and \ref{4.112}:
\begin{equation}
\left.\frac{\eta^2}{2c}\beta_N^3\frac{dH}{d\sigma}s_{\Nb\Lb}\right|_{\partial_-{\cal B}}=l
\label{9.104}
\end{equation}
is a smooth strictly negative function on $\partial_-{\cal B}$, hence $\left.s_{\Nb\Lb}\right|_{\partial_-{\cal B}}$ is 
a smooth nowhere vanishing function on $\partial_-{\cal B}$ of sign opposite to that of 
$\left.(dH/d\sigma)\right|_{\partial_-{\cal B}}$. Since 
$$\left.\epb_N\right|_{\tau=0}=\left.\epb\right|_{\tau=0}=0, \ \ \ 
\left.\frac{\partial\epb_N}{\partial\tau}\right|_{\tau=0}=
\left.\frac{\partial\epb}{\partial\tau}\right|_{\tau=0}$$
and $\epb_N$ is a known smooth function of $(\tau,\vartheta)$, it follows that:
\begin{equation}
\epb_N(\tau,\vartheta)=\tau\left(2\left.s_{\Nb\Lb}\right|_{\partial_-{\cal B}}(\vartheta)+O(\tau)\right)
\label{9.105}
\end{equation}
hence, with a suitable restriction on the size of $\tau$ we may conclude that 
$\tau/\epb_N$ is a known smooth function of $(\tau,\vartheta)$. Then Proposition 9.5 gives:
\begin{equation}
\hat{\nu}_N=O(\tau^{N+1})
\label{9.106}
\end{equation}

We turn to the nonlinear boundary condition \ref{1.329}. Actually, we have shown in Proposition 4.2 
that this condition is equivalent to:
\begin{equation}
\ep+j(\kappa,\epb)\epb^2=0
\label{9.107}
\end{equation}
where $\kappa$ stands for the quadruplet \ref{4.43}:
\begin{equation}
\kappa=(\sigma,c,\beta_N,\beta_{\Nb}) \ \ \mbox{: along ${\cal K}$}
\label{9.108}
\end{equation}
and $j$ is a given smooth function of its arguments. Defining then in connection with the $N$th 
approximants:
\begin{equation}
\kappa_N=(\sigma_N,c_N,\beta_{N,N},\beta_{\Nb,N}) \ \ \mbox{: along ${\cal K}$}
\label{9.109}
\end{equation}
where
\begin{equation}
\beta_{N,N}=\beta_{\mu,N}N_N^\mu, \ \ \ \beta_{\Nb,N}=\beta_{\mu,N}\Nb_N^\mu 
\label{9.110}
\end{equation}
the point is to estimate the quantity: 
\begin{equation}
\ep_N+j(\kappa_N,\epb_N)\epb_N^2:=d_N
\label{9.111}
\end{equation}
by which \ref{9.107} fails to be satisfied by the $N$th approximants. By the 3rd of \ref{9.17}, 
Lemma 9.1, and \ref{9.32}:
\begin{equation}
\left.\frac{\partial^n\kappa_N}{\partial\tau^n}\right|_{\tau=0}=
\left.\frac{\partial^n\kappa}{\partial\tau^n}\right|_{\tau=0} \ \ \mbox{: for $n=0,...,N$}
\label{9.112}
\end{equation}
Noting that in applying $\left.\partial^n/\partial\tau^n\right|_{\tau=0}$ to the 2nd term on 
the left hand side of \ref{9.111}, a non-vanishing term results only when each of the two $\epb_N$ 
factors receives at least one derivative leaving at most $n-1$ derivatives to fall on the other 
factor, using \ref{9.91}, \ref{9.93} we deduce:
\begin{equation}
\left.\frac{\partial^n d_N}{\partial\tau^n}\right|_{\tau=0}=
\left.\frac{\partial^n(\ep+j(\kappa,\epb)\epb^2)}{\partial\tau^n}\right|_{\tau=0}=0 \ \ 
\mbox{: for $n=0,...,N$}
\label{9.113}
\end{equation}
which, in view of the fact that $d_N$ is a known smooth function of $(\tau,\vartheta)$ yields 
the following proposition. 

\vspace{2.5mm}

\noindent{\bf Proposition 9.6} \ \ We have:
$$d_N=O(\tau^{N+1})$$

\vspace{5mm}

\section{Estimates for the Quantities by which the $N$th Approximants fail to satisfy the 
Identification Equations}

Let us define, in reference to \ref{9.51}, in analogy with \ref{4.222}, \ref{4.223}, \ref{4.240}, 
\ref{4.244}, 
\begin{equation}
\check{f}_N(\tau,\vartheta)=f_N(\tau,\vartheta)-f(0,\vartheta), \ \ \ 
\check{g}_N^i(\tau,\vartheta)=g_N^i(\tau,\vartheta)-g^i(0,\vartheta)
\label{9.114}
\end{equation}
\begin{equation}
\delta_N^i(\tau,\vartheta)=\check{g}_N^i(\tau,\vartheta)-N_0^i(\vartheta)\check{f}_N(\tau,\vartheta)
\label{9.115}
\end{equation}
\begin{equation}
\hat{f}_N(\tau,\vartheta)=\tau^{-2}\check{f}_N(\tau,\vartheta), \ \ \ 
\hat{\delta}_N^i(\tau,\vartheta)=\tau^{-3}\delta_N^i(\tau,\vartheta)
\label{9.116}
\end{equation}
We then define, in analogy with the expressions of Proposition 4.5 
(here we are considering the case of 2 spatial dimensions):
\begin{eqnarray}
&&\hat{F}_N^i((\tau,\vartheta),(v,\gamma))=S_0^i(\vartheta)v\left[l(\vartheta)\hat{f}_N(\tau,\vartheta)
+\frac{1}{6}k(\vartheta)v^2\right] \label{9.117}\\
&&\hspace{26mm}+\Omega^{\prime i}_{0}(\vartheta)\gamma-\hat{\delta}_N^i(\tau,\vartheta)
+\tau E_N^i((\tau,\vartheta),(v,\gamma)) \nonumber
\end{eqnarray}
the functions $E_N^i$ being given by:
\begin{eqnarray}
&&E_N^i((\tau,\vartheta),(v,\gamma))=\frac{1}{2}\left(a_0^i(\vartheta)\hat{f}_N(\tau,\vartheta)
+a_1^i(\vartheta)v^2\right)\hat{f}_N(\tau,\vartheta)\nonumber\\
&&\hspace{26mm}+\frac{1}{2}\tau a_2^i(\vartheta)v\hat{f}_N^2(\tau,\vartheta)
+\frac{1}{6}\tau^2 a_3^i(\vartheta)\hat{f}_N^3(\tau,\vartheta)\label{9.118}\\
&&\hspace{26mm}+v^4 A_0^i(\tau^2\hat{f}_N(\tau,\vartheta),\tau v,\vartheta)
+\tau v^3\hat{f}_N(\tau,\vartheta)A_1^i(\tau^2\hat{f}_N(\tau,\vartheta),\tau v,\vartheta)\nonumber\\
&&\hspace{13mm}+\tau^2 v^2\hat{f}_N^2(\tau,\vartheta)
A_2^i(\tau^2\hat{f}_N(\tau,\vartheta),\tau v,\vartheta)
+\tau^3 v\hat{f}_N^3(\tau,\vartheta) A_3^i(\tau^2\hat{f}_N(\tau,\vartheta),\tau v,\vartheta)\nonumber\\
&&\hspace{35mm}+\tau^4\hat{f}_N^4(\tau,\vartheta) A_4^i(\tau^2\hat{f}_N(\tau,\vartheta),\tau v,\vartheta)\nonumber\\
&&\hspace{13mm}+\gamma\left\{ v\Theta_{1}^i(\tau^2\hat{f}_N(\tau,\vartheta),\tau v,\tau^3\gamma,\vartheta)+\tau\hat{f}_N(\tau,\vartheta)\Theta_{0}^i(\tau^2\hat{f}_N(\tau,\vartheta),\tau v, 
\tau^3\gamma,\vartheta)\right.\nonumber\\
&&\hspace{30mm}+\left.\tau^2\gamma\Phi^i(\tau^2\hat{f}_N(\tau,\vartheta),\tau v,\tau^3\gamma,
\vartheta)\right\}\nonumber
\end{eqnarray}
If we insert, from \ref{9.53}, 
$$v=v_N(\tau,\vartheta), \ \ \ \gamma=\gamma_N(\tau,\vartheta)$$
into $\hat{F}^i_N$, the resulting quantities
\begin{equation}
\hat{F}_N^i((\tau,\vartheta),(v_N(\tau,\vartheta),\gamma_N(\tau,\vartheta))):=D^i_N(\tau,\vartheta)
\label{9.119}
\end{equation}
do not vanish, signifying that the $N$th approximants fail to satisfy the identification equations. 
We shall presently estimate these quantities, which are by their definition known smooth functions 
of $(\tau,\vartheta)$. 

According to the definitions of $f_N,g_N^i$ through \ref{9.51} and the 1st of \ref{9.3} these 
are polynomials of degree $N$ in $\tau$. In fact:
\begin{equation}
f_N=\sum_{n=0}^N\frac{\tau^n}{n!}\left.\frac{\partial^n f}{\partial\tau^n}\right|_{\tau=0}, \ \ \ 
g_N^i=\sum_{n=0}^N\frac{\tau^n}{n!}\left.\frac{\partial^n g^i}{\partial\tau^n}\right|_{\tau=0}
\label{9.120}
\end{equation}
In view of \ref{4.120} we then have:
\begin{equation}
\check{f}_N=\sum_{n=2}^N\frac{\tau^n}{n!}\left.\frac{\partial^n f}{\partial\tau^n}\right|_{\tau=0}, 
\ \ \ 
\check{g}_N^i=\sum_{n=2}^N\frac{\tau^n}{n!}\left.\frac{\partial^n g^i}{\partial\tau^n}\right|_{\tau=0}
\label{9.121}
\end{equation}
The first gives:
\begin{equation}
\hat{f}_N=\sum_{n=0}^{N-2}\frac{\tau^n}{(n+2)!}\left.\frac{\partial^{n+2}\check{f}}{\partial\tau^{n+2}}\right|_{\tau=0}
\label{9.122}
\end{equation}
On the other hand, since $\check{f}=\tau^2\hat{f}$, applying 
$\left.\partial^{n+2}/\partial\tau^{n+2}\right|_{\tau=0}$ only when exactly 2 derivatives fall on the 
factor $\tau^2$ do we obtain non-vanishing terms, hence:
\begin{equation}
\left.\frac{\partial^{n+2}\check{f}}{\partial\tau^{n+2}}\right|_{\tau=0}=
(n+2)(n+1)\left.\frac{\partial^n\hat{f}}{\partial\tau^n}\right|_{\tau=0}
\label{9.123}
\end{equation}
Substituting in \ref{9.122} then yields:
\begin{equation}
\hat{f}_N=\sum_{n=0}^{N-2}\frac{\tau^n}{n!}\left.\frac{\partial^n\hat{f}}{\partial\tau^n}\right|_{\tau=0}
\label{9.124}
\end{equation}
Next, substituting the expansions \ref{9.121} in \ref{9.115}, the $n=2$ terms cancel by \ref{4.131},  
hence we obtain:
\begin{equation}
\hat{\delta}_N^i=\sum_{n=0}^{N-3}\frac{\tau^n}{(n+3)!}\left.\frac{\partial^{n+3}\delta^i}{\partial\tau^{n+3}}\right|_{\tau=0}
\label{9.125}
\end{equation}
On the other hand, since $\delta^i=\tau^3\hat{\delta}^i$, applying 
$\left.\partial^{n+3}/\partial\tau^{n+3}\right|_{\tau=0}$ only when exactly 3 derivatives fall on the 
factor $\tau^3$ do we obtain non-vanishing terms, hence:
\begin{equation}
\left.\frac{\partial^{n+3}\delta^i}{\partial\tau^{n+3}}\right|_{\tau=0}=
(n+3)(n+2)(n+1)\left.\frac{\partial^n\hat{\delta}^i}{\partial\tau^n}\right|_{\tau=0}
\label{9.126}
\end{equation}
Substituting in \ref{9.125} then yields:
\begin{equation}
\hat{\delta}_N^i=\sum_{n=0}^{N-3}\frac{\tau^n}{n!}\left.\frac{\partial^n\hat{\delta}^i}{\partial\tau^n}
\right|_{\tau=0}
\label{9.127}
\end{equation}
Here we encounter a difficulty. For, while according to \ref{9.124} the partial derivatives of 
$\hat{f}_N$ with respect to $\tau$ of order up to $N-2$ agree at $\tau=0$ with those of $\hat{f}$, 
the partial derivatives of $\hat{\delta}_N^i$ with respect to $\tau$ of order only up 
to $N-3$ agree at $\tau=0$ with those of $\hat{\delta}^i$. In fact the partial derivatives of 
$\hat{\delta}_N^i$ with respect to $\tau$ of order $N-2$ vanish. 
On the other hand, as we have shown in Section 5.3, the pair 
$(\left.\partial^n v/\partial\tau^n\right|_{\tau=0},\left.\partial^n\gamma/\partial\tau^n\right|_{\tau=0}$ is derived from the linear equations obtained by setting to zero the $n$th derivative 
of $\hat{F}^i((\tau,\vartheta),(v(\tau,\vartheta),\gamma(\tau,\vartheta)))$ with respect to $\tau$ 
at $\tau=0$. Therefore for $n=N-2$, corresponding to the last terms in the expansions \ref{9.53}, 
$\left.\partial^{N-2}\hat{\delta}^i/\partial\tau^{N-2}\right|_{\tau=0}$ besides 
$\left.\partial^{N-2}\hat{f}/\partial\tau^{N-2}\right|_{\tau=0}$ is required. 

To overcome the difficulty we must redefine $\hat{\delta}_N^i$ to include the $n=N-2$ term. 
That is, we must redefine $\delta_N^i$ to include the term:
\begin{equation}
\frac{\tau^{N+1}}{(N+1)!}\left.\frac{\partial^{N+1}\delta^i}{\partial\tau^{N+1}}\right|_{\tau=0}
\label{9.128}
\end{equation}
Nevertheless we must show that this term depends only on the $N$th approximants. In fact we have, 
along ${\cal K}$:
\begin{eqnarray}
&&\frac{\partial\check{g}^i}{\partial\tau}=\left.\frac{\partial x^i}{\partial\tau}\right|_{\cal K}
=\left.\left(Lx^i+\Lb x^i\right)\right|_{\cal K}=\left.\left(\rho N^i+\rhob\Nb^i\right)\right|_{\cal K}
\nonumber\\
&&\frac{\partial\check{f}}{\partial\tau}=\left.\frac{\partial t}{\partial\tau}\right|_{\cal K}
=\left.\left(Lt+\Lb t\right)\right|_{\cal K}=\left.\left(\rho+\rhob\right)\right|_{\cal K}
\label{9.129}
\end{eqnarray}
Therefore we have:
\begin{eqnarray}
&&\left.\frac{\partial^{N+1}\delta^i}{\partial\tau^{N+1}}\right|_{\tau=0}=
\left.\frac{\partial^N(\rho N^i+\rhob\Nb^i)}{\partial\tau^N}\right|_{\tau=0}
-N_0^i\left.\frac{\partial^N(\rho+\rhob)}{\partial\tau^N}\right|_{\tau=0}\label{9.130}\\
&&\hspace{20mm}=\sum_{n=1}^N\left(\begin{array}{l} N\\n\end{array}\right)
\left.\frac{\partial^{N-n}\rho}{\partial\tau^{N-n}}\right|_{\tau=0}
\left.\frac{\partial^n N^i}{\partial\tau^n}\right|_{\tau=0}
+\left.\frac{\partial^N(\rhob(\Nb^i-N_0^i))}{\partial\tau^N}\right|_{\tau=0}\nonumber
\end{eqnarray}
Since $\left.\partial^N\rho/\partial\tau^N\right|_{\tau=0}$ cancels, by \ref{9.37} and Lemma 9.1 
this corresponds to what is obtained using the $N$th approximants, that is we have:
\begin{eqnarray}
&&\left.\frac{\partial^{N+1}\delta^i}{\partial\tau^{N+1}}\right|_{\tau=0}= \label{9.131}\\
&&\hspace{10mm}\sum_{n=1}^N\left(\begin{array}{l} N\\n\end{array}\right)
\left.\frac{\partial^{N-n}\rho_N}{\partial\tau^{N-n}}\right|_{\tau=0}
\left.\frac{\partial^n N_N^i}{\partial\tau^n}\right|_{\tau=0}
+\left.\frac{\partial^N(\rhob_N(\Nb_N^i-N_0^i))}{\partial\tau^N}\right|_{\tau=0}\nonumber
\end{eqnarray}
Redefining then $\delta_N^i$ by adding the term \ref{9.128} as given by \ref{9.131} we achieve:
\begin{equation}
\delta_N^i=\sum_{n=0}^{N+1}\frac{\tau^n}{n!}\left.\frac{\partial^n\delta^i}{\partial\tau^n}\right|_{\tau=0}, \ \ \ 
\hat{\delta}_N^i=\sum_{n=0}^{N-2}\frac{\tau^n}{n!}\left.\frac{\partial^n\hat{\delta}^i}{\partial\tau^n}\right|_{\tau=0}
\label{9.132}
\end{equation}
Then in view of \ref{9.53}, \ref{9.124} and the 2nd of \ref{9.132}, the partial derivatives of the left hand side of \ref{9.119} with respect to $\tau$ of order up to $N-2$ vanish at $\tau=0$ as these 
agree with the corresponding derivatives of 
$\hat{F}^i((\tau,\vartheta),(v(\tau,\vartheta),\gamma(\tau,\vartheta)))$ at $\tau=0$. This 
yields the following. 

\vspace{2.5mm}

\noindent{\bf Proposition 9.7} \ \ We have:
$$D_N^i=O(\tau^{N-1})$$

\vspace{5mm}

\section{Estimates for the Quantity by which the $\beta_{\mu,N}$ fail to satisfy the Wave Equation 
relative to $\tilde{h}_N$ and to $\tilde{h}^\prime_N$}

We shall apply the following extension of Proposition 2.1. 

\vspace{2.5mm}

\noindent{\bf Proposition 9.8} \ \ Let $({\cal M},g)$ be a $n+1$ dimensional Lorentzian manifold and 
$\beta$ a 1-form on ${\cal M}$ satisfying:
$$d\beta=-\omega, \ \ \ \mbox{div}_g(G\beta)=\rho$$
where with 
$$\sigma=-g^{-1}(\beta,\beta)>0,$$
$G>0$ is a given function of $\sigma$ such that 
$$\frac{2}{G}\frac{dG}{d\sigma}>0.$$
Here $\mbox{div}_g$ is the divergence operator associated to the metric $g$. 
Let $X$ be a vectorfield generating a 1-parameter group of isometries of $({\cal M},g)$. 
With $\eta>0$ defined by
$$\frac{1}{\eta^2}=1+\frac{2}{G}\frac{dG}{d\sigma},$$
so $\eta<1$, and $H>0$ by
$$\sigma H=1-\eta^2,$$
let $h,\tilde{h}$ be the Lorentzian metrics:
$$h=g+H\beta\otimes\beta, \ \ \tilde{h}=\Omega h, \ \ \Omega=\left(\frac{G}{\eta}\right)^{2/(n-1)}$$
Then the function $\beta(X)$ satisfies the equation:
$$\square_{\tilde{h}}\beta(X)=\kappa_X$$
where $\kappa_X$ is the function:
$$\kappa_X=\mbox{div}_{\tilde{h}}i_X\omega+\left(\frac{\eta^2}{G^{n+1}}\right)^{1/(n-1)}X\rho$$
Here $\mbox{div}_{\tilde{h}}$ is the divergence operator associated to the metric $\tilde{h}$. 

\vspace{2.5mm}

\noindent{\em Proof:} \ We go through the proof of Proposition 2.1 and make the necessary changes 
as required by the more general assumptions on $\beta$. First, \ref{2.135} changes to:
\begin{equation}
\nabla_\mu A^\mu=\rho
\label{9.148}
\end{equation}
As a consequence, \ref{2.138} changes to:
\begin{equation}
d\alpha=\rho\epsilon
\label{9.149}
\end{equation}
and, in view of the fact that ${\cal L}_X\epsilon=0$, \ref{2.139} changes to:
\begin{equation}
d{\cal L}_X\alpha=(X\rho)\epsilon
\label{9.150}
\end{equation}
Moreover, \ref{2.142} changes to:
\begin{equation}
{\cal L}_X\beta=d\beta(X)-i_X\omega
\label{9.151}
\end{equation}
As a consequence, \ref{2.145} changes to:
\begin{equation}
XG=-GFg^{-1}(\beta,d\beta(X)-i_X\omega)
\label{9.152}
\end{equation}
Then \ref{2.146} changes to:
\begin{equation}
{\cal L}_X(G\beta)=G\left(d\beta(X)-i_X\omega-Fg^{-1}(\beta,d\beta(X)-i_X\omega)\right)
\label{9.153}
\end{equation}
We then apply \ref{2.147} taking $\theta=d\beta(X)-i_X\omega$ to obtain \ref{2.148} with 
$B$ the vectorfield:
\begin{equation}
B=h^{-1}\cdot(d\beta(X)-i_X\omega)
\label{9.154}
\end{equation}
As a consequence, \ref{2.162} changes to:
\begin{equation}
{\cal L}_X\alpha=\s^\star(d\beta(X)-i_X\omega)
\label{9.155}
\end{equation}
Since \ref{2.139} has been changed to \ref{9.150}, we now obtain:
\begin{equation}
d\s^\star(d\beta(X)-i_X\omega)=(X\rho)\epsilon
\label{9.156}
\end{equation}
Since by \ref{2.156}, \ref{2.160} 
$$\epsilon=\eta^{-1}e=\eta^{-1}\Omega^{-(n+1)/2}\tilde{e},$$
\ref{9.156} is equivalent to:
\begin{equation}
\square_{\tilde{h}}\beta(X)=\mbox{div}_{\tilde{h}}i_X\omega+\eta^{-1}\Omega^{-(n+1)/2}X\rho
\label{9.157}
\end{equation}
which in view of the definition of $\Omega$ is the conclusion of the proposition. 

\vspace{2.5mm}

Taking in the above proposition $({\cal M},g)$ to be the Minkowski spacetime and taking the 
1-form $\beta_N$ (see \ref{9.a1}) in the role of the 1-form $\beta$, then we have $\sigma_N$, $G_N$, $H_N$ and 
$\eta_N$ in the roles of $\sigma$, $G$, $H$ and $\eta$, and the $N$th approximant acoustical metric 
$h_N$ (see \ref{9.a2}) in the role of the Lorentzian metric $h$. Also, we have 
\begin{equation}
\Omega_N=\left(\frac{G_N}{\eta_N}\right)^{2/(n-1)}, \ \ \ \tilde{h}_N=\Omega_N h_N
\label{9.158}
\end{equation}
in the roles of the conformal factor $\Omega$ and the Lorentzian metric $h_N$ respectively. 
Moreover, we have $\omega_N$ (see \ref{9.a5}) in the role of $\omega$ and $G_N a_N^{-1}\delta_N$ in 
the role of $\rho$. Taking the translation fields $X_{(\mu)}$ (see \ref{6.21}), we have 
the functions $\beta_{\mu,N}$ in the role of the function $\beta(X)$ and the conclusion of the 
proposition is:
\begin{equation}
\square_{\tilde{h}_N}\beta_{\mu,N}=\kappa_{\mu,N}
\label{9.159}
\end{equation}
where:
\begin{equation}
\kappa_{\mu,N}=\mbox{div}_{\tilde{h}_N}i_{X_{(\mu)}}\omega_N+
\left(\frac{\eta_N^2}{G_N^{n+1}}\right)^{1/(n-1)}X_{(\mu)}(G_N a_N^{-1}\delta_N)
\label{9.160}
\end{equation}
We proceed to estimate $\kappa_{\mu,N}$ in the case $n=2$. 

We first expand the translation fields $X_{(\mu)}$ in the $(L_N,\Lb_N,E_N)$ frame:
\begin{equation}
X_{(\mu)}=X^L_{\mu,N}L_N+X^{\Lb}_{\mu,N}\Lb_N+X^E_{\mu,N}E_N
\label{9.161}
\end{equation}
To determine the coefficients we first note that by \ref{9.139}:
\begin{eqnarray}
&&h^\prime_N(X_{(\mu)},L_N)=-2a_N X^{\Lb}_{\mu,N}\nonumber\\
&&h^\prime_N(X_{(\mu)},\Lb_N)=-2a_N X^L_{\mu,N}\nonumber\\
&&h^\prime_N(X_{(\mu)},E_N)=X^E_{\mu,N}
\label{9.162}
\end{eqnarray}
On the other hand we have (see \ref{9.15}):
\begin{eqnarray}
&&h_N(X_{(\mu)},L_N)=h_{\mu\nu,N}(L_N x_N^\nu)=h_{\mu\nu,N}(\rho_N N_N^\nu+\vep_N^\nu)\nonumber\\
&&h_N(X_{(\mu)},\Lb_N)=h_{\mu\nu,N}(\Lb_N x_N^\nu)=h_{\mu\nu,N}(\rhob_N\Nb_N^\nu+\vepb_N\Nb_N^\nu)
\nonumber\\
&&h_N(X_{(\mu)},E_N)=h_{\mu\nu,N}(E_N x_N^\mu)=h_{\mu\nu,N}E_N^\nu
\label{9.163}
\end{eqnarray}
Now, by \ref{9.147}, \ref{9.146} and the expansion \ref{9.161}, the differences of the left hand sides 
of {9.162} from the left hand sides of \ref{9.163} are given by:
\begin{eqnarray*}
&&\delta_{h,N}(X_{(\mu)},L_N)=\delta_{LL,N}X^L_{\mu,N}+\delta_{L\Lb,N}X^{\Lb}_{\mu,N}
+\delta_{LE,N}X^E_{\mu,N}\\
&&\delta_{h,N}(X_{(\mu)},\Lb_N)=\delta_{L\Lb,N}X^L_{\mu,N}+\delta_{\Lb\Lb,N}X^{\Lb}_{\mu,N}
+\delta_{\Lb E,N}X^E_{\mu,N}\\
&&\delta_{h,N}(X_{(\mu)},E_N)=\delta_{LE,N}X^L_{\mu,N}+\delta_{\Lb E,N}X^{\Lb}_{\mu,N}
+\delta_{EE,N}X^E_{\mu,N} 
\end{eqnarray*}
Equating these to the differences of the right hand sides we obtain for each $\mu=0,1,2$ the following 
system of 3 linear equations for the 3 coefficients $X^L_{\mu,N}$, $X^{\Lb}_{\mu,N}$, $X^E_{\mu,N}$:
\begin{eqnarray}
&&-2a_N X^{\Lb}_{\mu,N}=h_{\mu\nu,N}(\rho_N N_N^\nu+\vep_N^\nu)\label{9.164}\\
&&\hspace{20mm}-\delta_{LL,N}X^L_{\mu,N}-\delta_{L\Lb,N}X^{\Lb}_{\mu,N}-\delta_{LE,N}X^E_{\mu,N}
\nonumber
\end{eqnarray}
\begin{eqnarray}
&&-2a_N X^L_{\mu,N}=h_{\mu\nu}(\rhob_N\Nb_N^\nu+\vepb_N^\nu)\label{9.165}\\
&&\hspace{20mm}-\delta_{L\Lb,N}X^L_{\mu,N}-\delta_{\Lb\Lb,N}X^{\Lb}_{\mu,N}-\delta_{\Lb E,N}X^E_{\mu,N}
\nonumber
\end{eqnarray}
\begin{eqnarray}
&&X^E_{\mu,N}=h_{\mu\nu,N}E_N^\nu \label{9.166}\\
&&\hspace{20mm}-\delta_{LE,N}X^L_{\mu,N}-\delta_{\Lb E,N}X^{\Lb}_{\mu,N} \nonumber
\end{eqnarray}
It follows that, neglecting terms of the form $\rhob_N^{-1}O(\tau^{N+1})$, 
\begin{equation}
\rho_N X^L_{\mu,N}=O(1), \ \ \ \rhob_N X^{\Lb}_{\mu,N}=O(1), \ \ \ X^E_{\mu,N}=O(1)
\label{9.167}
\end{equation}

The estimates \ref{9.167} for the coefficients of the expansion \ref{9.161} of $X_{(\mu)}$, together 
with Proposition 9.3 allow us to estimate the 2nd term on the right in \ref{9.160}. To state the 
estimate in as simple form as possible we introduce the following definition. 

\vspace{2.5mm}

\noindent{\bf Definition 9.2} \ \ Let $f$ be a known function of $(\tau,\sigma,\vartheta)$ which is 
smooth except possibly at $\tau=\sigma=0$ and let $p,q$ be non-negative integers. We write: 
$$f={\bf O}(u^{-q}\tau^p)$$
if for each triplet of non-negative integers $m,n,l$ with $m\leq p$ there is a constant 
$C_{m,n,l}$ such that:
$$\left|\frac{\partial^{m+n+l}f}{\partial\tau^m\partial\sigma^n\partial\vartheta^l}\right|
\leq C_{m,n,l}u^{-q-n}\tau^{p-m}$$ 

\vspace{2.5mm}

Using then \ref{9.a22} we deduce, in regard to the 2nd term on the right in \ref{9.160},
\begin{equation}
\left(\frac{\eta_N^2}{G_N^{n+1}}\right)^{1/(n-1)}X_{(\mu)}(G_N a_N^{-1}\delta_N)
={\bf O}(u^{-2}\tau^{N-3})+{\bf O}(u^{-5}\tau^{N-1})
\label{9.168}
\end{equation}

We turn to the 1st term on the right in \ref{9.160}. We have:
\begin{eqnarray*}
&&(i_{X_{(\mu)}}\omega_N)(L_N)=X^{\Lb}_{\mu,N}\omega_{\Lb L,N}+X^E_{\mu,N}\omega_{EL,N}\\
&&(i_{X_{(\mu)}}\omega_N)(\Lb_N)=X^L_{\mu,N}\omega_{L\Lb,N}+X^E_{\mu,N}\omega_{E\Lb,N}\\
&&(i_{X_{(\mu)}}\omega_N)(E_N)=X^L_{\mu,N}\omega_{LE,N}+X^{\Lb}_{\mu,N}\omega_{\Lb E,N}
\end{eqnarray*}
From Proposition 9.2 and {9.167} we then deduce:
\begin{eqnarray}
&&(i_{X_{(\mu)}}\omega_N)(L_N)=\rhob_N^{-1}O(\tau^N) \nonumber\\
&&(i_{X_{(\mu)}}\omega_N)(\Lb_N)=O(\tau^{N-1}) \nonumber\\
&&(i_{X_{(\mu)}}\omega_N)(E_N)=O(\tau^{N-1})+\rhob_N^{-1}O(\tau^{N+1})
\label{9.169}
\end{eqnarray}

To determine the 1st term on the right in \ref{9.160}, we recall that if $({\cal M},g)$ is an  
arbitrary 3-dimensional Lorentzian manifold and $\xi$ is an arbitrary 1-form on ${\cal M}$, 
$\mbox{div}_g\xi$ is the function defined by:
\begin{equation}
d\s^{*_g}\xi=\epsilon\mbox{div}_g\xi
\label{9.170}
\end{equation}
where $\epsilon$ is the volume form of $g$ and $\s^{*_g}\xi$ is the Hodge dual of $\xi$ relative 
to $g$, the 2-form on ${\cal M}$ defined by:
\begin{equation}
\s^{*_g}\xi(X,Y)=\epsilon(g^{-1}\cdot\xi,X,Y) \ \ \mbox{: for any pair $X,Y$ of vectorfields on 
${\cal M}$}
\label{9.171}
\end{equation}
To determine then the 1st term on the right in \ref{9.160} we must place $({\cal N},\tilde{h}_N)$ 
in the role of the Lorentzian manifold $({\cal M}, g)$ and the 1-form $i_{X_{(\mu)}}\omega_N$ in 
the role of the 1-form $\xi$. Consider first the Lorentzian metrics $h^\prime_N$, $h_N$ on 
${\cal N}$. Let $e^\prime_N$, $e_N$ be the corresponding volume forms. With 
$$T_N=\frac{1}{2\sqrt{a_N}}(L_N+\Lb_N), \ \ \ M_N=\frac{1}{2\sqrt{a_N}}(L_N-\Lb_N),$$
$(T_N,M_N,E_N)$ is a positive orthonormal frame relative to $h^\prime_N$ (see \ref{9.139}). 
Therefore
$$e^\prime_N(T_N,M_N,E_N)=1$$
which implies:
\begin{equation}
e^\prime_N(L_N,\Lb_N,E_N)=-2a_N
\label{9.172}
\end{equation}
Similarly, $(N_N^\mu,\Nb_N^\mu,E_N^\mu)$ are the rectangular components of a null frame relative to 
$h_N$ with
$$h_{\mu\nu,N}N_N^\mu\Nb_N^\nu=-2c_N, \ \ \ h_{\mu\nu, N}E_N^\mu E_N^\nu=1$$
It follows that:
\begin{equation}
e_{\mu\nu\lambda,N}N_N^\mu\Nb_N^\nu E_N^\lambda=-2c_N
\label{9.173}
\end{equation}
or:
\begin{equation}
e_{\mu\nu\lambda,N}\rho_N N_N^\mu \rhob_N \Nb_N^\nu E_N^\lambda=-2a_N
\label{9.174}
\end{equation}
Now, by \ref{9.15},
$$\rho_N N_N^\mu=L_N^\mu-\vep_N^\mu, \ \ \ \rhob_N\Nb_N^\mu=\Lb_N^\mu-\vepb_N^\mu$$
Substituting in \ref{9.174} we then obtain:
\begin{equation}
e_N(L_N,\Lb_N,E_N)=-2a_N+v_N
\label{9.175}
\end{equation}
where:
\begin{eqnarray}
&&v_N=-\rhob_N e_{\mu\nu\lambda,N}\vep_N^\mu\Nb_N^\nu E_N^\lambda
-\rho_N e_{\mu\nu\lambda,N}N_N^\mu\vepb_N^\nu E_N^\lambda \nonumber\\
&&\hspace{10mm}-e_{\mu\nu\lambda,N}\vep_N^\mu\vepb_N^\nu E_N^\lambda 
\label{9.176}
\end{eqnarray}
In view of the fact that (see \ref{2.156}):
\begin{equation}
e_{\mu\nu\lambda,N}=\eta_N\epsilon_{\mu\nu\lambda}
\label{9.177}
\end{equation}
where $\epsilon_{\mu\nu\lambda}$ are the rectangular components of the volume form $\epsilon$ 
of the Minkowski metric $g$, namely the fully antisymmetric 3-dimensional symbol, Proposition 9.1 
implies that:
\begin{equation}
v_N=\rhob_N O(\tau^N)+O(\tau^{N+2})
\label{9.178}
\end{equation}

Let us set:
\begin{equation}
\xi=i_{X_{(\mu)}}\omega_N
\label{9.179}
\end{equation}
and consider $h_N^{\prime-1}\cdot\xi$. We can expand:
\begin{equation}
h_N^{\prime-1}\cdot\xi=(h_N^{\prime-1}\cdot\xi)^L L_N+(h_N^{\prime-1}\cdot\xi)^{\Lb}\Lb_N
+(h_N^{\prime-1}\cdot\xi)^E E_N
\label{9.180}
\end{equation}
By \ref{9.a10} and the estimates \ref{9.169} we obtain:
\begin{eqnarray}
&&(h_N^{\prime-1}\cdot\xi)^L=-\frac{\xi(\Lb_N)}{2a_N}=\rhob_N^{-1}O(\tau^{N-2}) \nonumber\\
&&(h_N^{\prime-1}\cdot\xi)^{\Lb}=-\frac{\xi(L_N)}{2a_N}=\rhob_N^{-2}O(\tau^{N-1}) \nonumber\\
&&(h_N^{\prime-1}\cdot\xi)^E=\xi(E_N)=O(\tau^{N-1})+\rhob_N^{-1}O(\tau^{N+1}) 
\label{9.181}
\end{eqnarray}
Then \ref{9.171} with $h^\prime_N$ in the role of $g$ gives:
\begin{eqnarray}
&&(\s^{*_{h^\prime_N}}\xi)(L_N,\Lb_N)=(h_N^{\prime-1}\cdot\xi)^E e^\prime_N(E_N,L_N,\Lb_N)
=-2a_N\xi(E_N)\nonumber\\
&&\hspace{7cm}=\rhob_N O(\tau^N)+O(\tau^{N+2})\nonumber\\
&&(\s^{*_{h^\prime_N}}\xi)(L_N,E_N)=(h_N^{\prime-1}\cdot\xi)^{\Lb} e^\prime_N(\Lb_N,L_N,E_N)
=-\xi(L_N)=\rhob_N^{-1}O(\tau^N)\nonumber\\
&&(\s^{*_{h^\prime_N}}\xi)(\Lb_N,E_N)=(h_N^{\prime-1}\cdot\xi)^L e^\prime_N(L_N,\Lb_N,E_N)
=\xi(\Lb_N)=O(\tau^{N-1})\nonumber\\
&&\label{9.182}
\end{eqnarray}
Next, we consider $h_N^{-1}\cdot\xi$. We similarly expand:
\begin{equation}
h_N^{-1}\cdot\xi=(h_N^{-1}\cdot\xi)^L L_N+(h_N^{-1}\cdot\xi)^{\Lb}\Lb_N
+(h_N^{-1}\cdot\xi)^E E_N
\label{9.183}
\end{equation}
Since by \ref{9.a25},
\begin{equation}
h_N^{-1}\cdot\xi=h_N^{\prime-1}\cdot\xi+\delta^\prime_{h,N}\cdot\xi
\label{9.184}
\end{equation}
we have, by \ref{9.a26}, \ref{9.182}, 
\begin{eqnarray}
&&(h_N^{-1}\cdot\xi)^L=-\frac{\xi(\Lb_N)}{2a_N}+\delta_N^{\prime LL}\xi(L_N)
+\delta_N^{\prime L\Lb}\xi(\Lb_N)+\delta_N^{\prime LE}\xi(E_N)\nonumber\\
&&(h_N^{-1}\cdot\xi)^{\Lb}=-\frac{\xi(L_N)}{2a_N}+\delta_N^{\prime L\Lb}\xi(L_N)
+\delta_N^{\prime \Lb\Lb}\xi(\Lb_N)+\delta_N^{\prime \Lb E}\xi(E_N)\nonumber\\
&&(h_N^{-1}\cdot\xi)^E=\xi(E_N)+\delta_N^{\prime LE}\xi(L_N)+\delta_N^{\prime \Lb E}\xi(\Lb_N)
+\delta_N^{\prime EE}\xi(E_N) \label{9.185}
\end{eqnarray}
hence, by the estimates \ref{9.a29} and \ref{9.169}:
\begin{eqnarray}
&&(h_N^{-1}\cdot\xi)^L=-\frac{\xi(\Lb_N)}{2a_N}+\rhob_N^{-1}O(\tau^{2N-3})+\rhob_N^{-2}O(\tau^{2N-1})
\nonumber\\
&&(h_N^{-1}\cdot\xi)^{\Lb}=-\frac{\xi(L_N)}{2a_N}+\rhob_N^{-2}O(\tau^{2N-2})+\rhob_N^{-3}O(\tau^{2N})
\nonumber\\
&&(h_N^{-1}\cdot\xi)^E=\xi(E_N)+\rhob_N^{-1}O(\tau^{2N-2})+\rhob_N^{-2}O(\tau^{2N}) \label{9.186}
\end{eqnarray}
Then \ref{9.171} with $h_N$ in the role of $g$ gives, in view of \ref{9.175} and the estimates 
\ref{9.169}, \ref{9.178},
\begin{eqnarray}
&&(\s^{*_{h_N}}\xi)(L_N,\Lb_N)=(h_N^{-1}\cdot\xi)^E e_N(E_N,L_N,\Lb_N)\nonumber\\
&&\hspace{25mm}=-2a_N\xi(E_N)+O(\tau^{2N-1})+\rhob_N^{-1}O(\tau^{2N+1})\nonumber\\
&&(\s^{*_{h_N}}\xi)(L_N,E_N)=(h_N^{-1}\cdot\xi)^{\Lb} e_N(\Lb_N,L_N,E_N)\nonumber\\
&&\hspace{25mm}=-\xi(L_N)+\rhob_N^{-1}O(\tau^{2N-1})+\rhob_N^{-2}O(\tau^{2N+1})\nonumber\\
&&(\s^{*_{h_N}}\xi)(\Lb_N,E_N)=(h_N^{-1}\cdot\xi)^L e_N(L_N,\Lb_N,E_N)\nonumber\\
&&\hspace{25mm}=\xi(\Lb_N)+O(\tau^{2N-2})+\rhob_N^{-1}O(\tau^{2N})\label{9.187}
\end{eqnarray}
or, again by the estimates \ref{9.169},
\begin{eqnarray}
&&(\s^{*_{h_N}}\xi)(L_N,\Lb_N)=\rhob_N O(\tau^N)+O(\tau^{N+2})+\rhob_N^{-1}O(\tau^{2N+1})\nonumber\\
&&(\s^{*_{h_N}}\xi)(L_N,E_N)=\rhob_N^{-1}O(\tau^N)+\rhob_N^{-2}O(\tau^{2N+1})\nonumber\\
&&(\s^{*_{h_N}}\xi)(\Lb_N,E_N)=O(\tau^{N-1})+\rhob_N^{-1}O(\tau^{2N}) 
\label{9.188}
\end{eqnarray}
Since (see \ref{9.158}):
\begin{equation}
\s^{*_{\tilde{h}_N}}\xi=\Omega_N^{1/2}\s^{*_{h_N}}\xi=\frac{G_N}{\eta_N}\s^{*{h_N}}\xi
\label{9.189}
\end{equation}
we also have:
\begin{eqnarray}
&&(\s^{*_{\tilde{h}_N}}\xi)(L_N,\Lb_N)=
\rhob_N O(\tau^N)+O(\tau^{N+2})+\rhob_N^{-1}O(\tau^{2N+1})\nonumber\\
&&(\s^{*_{\tilde{h}_N}}\xi)(L_N,E_N)=\rhob_N^{-1}O(\tau^N)+\rhob_N^{-2}O(\tau^{2N+1})\nonumber\\
&&(\s^{*_{\tilde{h}_N}}\xi)(\Lb_N,E_N)=O(\tau^{N-1})+\rhob_N^{-1}O(\tau^{2N}) 
\label{9.190}
\end{eqnarray}
According to \ref{9.170} with $\tilde{h}_N$ in the role of $g$:
\begin{equation}
\mbox{div}_{\tilde{h}_N}\xi=\frac{(d\s^{*_{\tilde{h}_N}}\xi)(L_N,\Lb_N,E_N)}
{\tilde{e}_N(L_N,\Lb_N,E_N)}
\label{9.191}
\end{equation}
and by \ref{9.175}:
\begin{equation}
\tilde{e}_N(L_N,\Lb_N,E_N)=\Omega_N^{3/2}e_N(L_N,\Lb_N,E_N)=\left(\frac{G_N}{\eta_N}\right)^3
(-2a_N+v_N)
\label{9.192}
\end{equation}
The 1-form $\xi$ being given by \ref{9.179}, \ref{9.191} is the 1st term on the right in \ref{9.160}. 
Now, if $\theta$ is an arbitrary 2-form on an arbitrary manifold ${\cal M}$, $d\theta$, the 
exterior derivative of $\theta$, is given by:
\begin{eqnarray}
&&(d\theta)(X,Y,Z)=X(\theta(Y,Z))+Y(\theta(Z,X))+Z(\theta(X,Y))\nonumber\\
&&\hspace{22.5mm}-\theta([X,Y],Z)-\theta([Y,Z],X)
-\theta([Z,X],Y)
\label{9.193}
\end{eqnarray}
where $X,Y,Z$ is an arbitrary triplet of vectorfields on ${\cal M}$. Taking the manifold ${\cal N}$ 
in the role of ${\cal M}$, the 2-form $\s^{*_{\tilde{h}_N}}\xi$ in the role of $\theta$ and the 
triplet of vectorfields $L_N,\Lb_N,E_N$ in the role of $X,Y,Z$, we obtain:
\begin{eqnarray*}
&&(d\s^{*_{\tilde{h}_N}}\xi)(L_N,\Lb_N,E_N)=L_N(\s^{*_{\tilde{h}_N}}\xi(\Lb_N,E_N))
+\Lb_N(\s^{*_{\tilde{h}_N}}\xi(E_N,L_N))\\
&&\hspace{50mm}+E_N(\s^{*_{\tilde{h}_N}}\xi(L_N,\Lb_N))\\
&&\hspace{35mm}-\s^{*_{\tilde{h}_N}}\xi([L_N,\Lb_N],E_N)-\s^{*_{\tilde{h}_N}}\xi([\Lb_N,E_N],L_N)\\
&&\hspace{50mm}-\s^{*_{\tilde{h}_N}}\xi([E_N,L_N],\Lb_N)
\end{eqnarray*}
Substituting from the commutation relations \ref{9.b1} this becomes:
\begin{eqnarray}
&&(d\s^{*_{\tilde{h}_N}}\xi)(L_N,\Lb_N,E_N)=L_N(\s^{*_{\tilde{h}_N}}\xi(\Lb_N,E_N))
-\Lb_N(\s^{*_{\tilde{h}_N}}\xi(L_N,E_N)) \nonumber\\
&&\hspace{50mm}+E_N(\s^{*_{\tilde{h}_N}}\xi(L_N,\Lb_N)) \nonumber\\
&&\hspace{35mm}+\chi_N\s^{*_{\tilde{h}_N}}\xi(\Lb_N,E_N)-\chib_N\s^{*_{\tilde{h}_N}}\xi(L_N,E_N)
\label{9.194}
\end{eqnarray}
Using the estimates \ref{9.190}, we then deduce through \ref{9.191}, \ref{9.192}, in view of the 
fact that by \ref{9.178}
\begin{equation}
1-\frac{v_N}{2a_N}=1+O(\tau^{N-1})+\rhob_N^{-1}O(\tau^{N+1})
\label{9.195}
\end{equation}
the following estimate in regard to the 1st term on the right in \ref{9.160}:
\begin{equation}
\mbox{div}_{\tilde{h}_N}(i_{X_{(\mu)}}\omega_N)={\bf O}(u^{-2}\tau^{N-3})+{\bf O}(u^{-5}\tau^{N-1})
\label{9.196}
\end{equation}
Combining finally this estimate with the estimate \ref{9.168} for the 2nd term on the right in 
\ref{9.160} yields the following proposition. 

\vspace{2.5mm}

\noindent{\bf Proposition 9.9} \ \ In the sense of Definition 9.2 we have:
$$\kappa_{\mu,N}={\bf O}(u^{-2}\tau^{N-3})+{\bf O}(u^{-5}\tau^{N-1})$$

\vspace{2.5mm}

Defining, in analogy with \ref{8.11} the rescaled functions
\begin{equation}
\tilde{\kappa}_{\mu,N}=\Omega_N a_N \kappa_{\mu,N}
\label{9.197}
\end{equation}
we have, according to \ref{9.159},
\begin{equation}
\Omega_N a_N\square_{\tilde{h}_N}\beta_{\mu,N}=\tilde{\kappa}_{\mu,N}
\label{9.198}
\end{equation}
and by Proposition 9.9:
\begin{equation}
\tilde{\kappa}_{\mu,N}={\bf O}(\tau^{N-2})+{\bf O}(u^{-3}\tau^N)
\label{9.199}
\end{equation}

Let now $f$ be a smooth function on ${\cal N}$, that is a $C^\infty$ function of the coordinates 
$(\ub,u,\vartheta)$ or equivalently of the coordinates $(\tau,\sigma,\vartheta)$. We want to 
express the difference:
\begin{equation}
a\Omega\square_{\tilde{h}}f-a_N\Omega_N\square_{\tilde{h}_N}f
\label{9.200}
\end{equation}
Consider first the difference:
\begin{equation}
a\square_h f-a_N\square_{h_N}f=a\square_h f-a_N\square_{h^\prime_N}f
-a_N(\square_{h_N}f-\square_{h^\prime_N}f)
\label{9.201}
\end{equation}
Now, if $({\cal M},g)$ is an arbitrary Lorentzian manifold an $f$ an arbitrary function on ${\cal M}$, 
then
\begin{equation}
\square_g f=\mbox{tr}(g^{-1}\cdot Ddf)
\label{9.202}
\end{equation}
$D$ being the covariant derivative operator associated to $g$, so $Ddf$ is the covariant Hessian of 
$f$ relative to $g$ a symmetric 2-covariant tensorfield on ${\cal M}$. If $X,Y$ is an arbitrary pair 
of vectorfields on ${\cal M}$ we have:
$$(Ddf)(X,Y)=(D_X df)(Y)=X(Yf)-(D_X Y)f=(Ddf)(Y,X)$$
Placing then $({\cal N},h)$ in the role of $({\cal M},g)$ and recalling that 
(see \ref{2.38}, \ref{3.a10}):
\begin{equation}
h^{-1}=-\frac{1}{2a}(L\otimes\Lb+\Lb\otimes L)+E\otimes E
\label{9.203}
\end{equation}
we deduce:
\begin{eqnarray}
&&a\square_h f=-\{\Lb(Lf)-(D_{\Lb}L)f\}+a\{ E(Ef)-(D_E E)f\}\nonumber\\ 
&&\hspace{10mm}=-\{L(\Lb f)-(D_L\Lb)f\}+a\{E(Ef)-(D_E E)f\}
\label{9.204}
\end{eqnarray}
Substituting for the connection coefficients of $D$ in the frame field $(L,\Lb,E)$ from  
\ref{3.a11}, \ref{3.a13} yields:
\begin{eqnarray}
&&a\square_h f=-\Lb(Lf)+aE(Ef)-\frac{1}{2}(\chib Lf+\chi\Lb f)+2\eta Ef\nonumber\\
&&\hspace{10mm}=-L(\Lb f)+aE(Ef)-\frac{1}{2}(\chib Lf+\chi\Lb f)+2\etab Ef
\label{9.205}
\end{eqnarray}
Now, the co-frame field $(\zeta_N,\zetab_N,\szeta_N)$ (see \ref{9.137}) being dual to the frame field 
$(L_N,\Lb_N,E_N)$ (see \ref{9.10}, \ref{9.11}), the metric $h^\prime_N$, defined by \ref{9.138}, bears 
the same relation to the frame field $(L_N,\Lb_N,E_N)$ and to the positive functions $a_N,\sh_N$ as 
the metric $h$ bears to the frame field $(L,\Lb,E)$ and to the positive functions $a,\sh$. 
It follows that, denoting by $D^\prime_N$ the covariant derivative operator associated to $h^\prime_N$, 
the connection coefficients of $D^\prime_N$ in the frame field $(L_N,\Lb_N,E_N)$ are given, in analogy 
to \ref{3.a11}, \ref{3.a13} by:
\begin{equation}
\begin{array}{llll}
&D^\prime_{N,L_N} L_N=a_N^{-1}(L_N a_N)L_N &\s &D^\prime_{N,\Lb}\Lb_N=a_N^{-1}(\Lb_N a_N)\Lb_N\\
&D^\prime_{N,\Lb_N} L_N=2\eta_N E_N &\s &D^\prime_{N,L_N}\Lb_N=2\etab_N E_N\\
&D^\prime_{N,E_N} L_N=\chi_N E_N+a_N^{-1}\etab_N L_N &\s &D^\prime_{N,E_N}\Lb_N=\chib_N E+a_N^{-1}\eta_N\Lb_N
\end{array}
\label{9.206}
\end{equation}
\begin{eqnarray}
&&D^\prime_{N,L_N} E_N=a_N^{-1}\etab_N L_N \ \ \ \ \ D^\prime_{N,\Lb_N} E_N=a_N^{-1}\eta_N\Lb_N \nonumber\\
&&\hspace{7mm} D^\prime_{N,E_N} E_N=(2a_N)^{-1}(\chib_N L_N+\chi_N\Lb_N)
\label{9.207}
\end{eqnarray}
Here $\chi_N,\chib_N$ are defined by \ref{9.b2}, so the first two of \ref{9.b1}, which are the 
analogues of the fist two of \ref{3.a14}, hold, and 
$\eta_N,\etab_N$ are defined by the analogue of \ref{3.a12}:
\begin{equation}
\eta_N+\etab_N=E_N a_N
\label{9.208}
\end{equation}
and by:
\begin{equation}
\eta_N-\etab_N=\frac{1}{2}\zeta_N
\label{9.209}
\end{equation}
with $\zeta_N$ defined by \ref{9.b3} and the last of \ref{9.b1}. That the connection coefficients of 
$D^\prime_N$ in the frame field $(L_N,\Lb_N,E_N)$ are given by \ref{9.206}, \ref{9.207} is a 
consequence of the fact that the connection coefficients of a torsion-free metric connection in 
a given frame field are determined by the metric coefficients in the frame field, here given by 
\ref{9.139}, and the commutation relations of the frame field, here given by \ref{9.b1}. 

Using \ref{9.206}, \ref{9.207} we deduce in analogy with \ref{9.205}:
\begin{eqnarray}
&&a_N\square_{h^\prime_N} f=-\Lb_N(L_Nf)+a_N E_N(E_N f)
-\frac{1}{2}(\chib_N L_N f+\chi_N\Lb_N f)+2\eta_N E_Nf\nonumber\\
&&\hspace{14mm}=-L_N(\Lb_N f)+a_N E_N(E_N f)-\frac{1}{2}(\chib_N L_N f+\chi_N\Lb_N f)+2\etab_N E_N f
\nonumber\\
&&\label{9.210}
\end{eqnarray}
The first difference on the right in \ref{9.201} is then obtained by subtracting \ref{9.210} from 
\ref{9.205}. 

We turn to the second difference on the right in \ref{9.201}, namely 
$$a_N(\square_{h_N}f-\square_{h^\prime_N}f)$$
To estimate this we recall the fact that if $f$ is an arbitrary function on an arbitrary 
Lorentzian manifold $({\cal M},g)$ we have:
\begin{equation}
\square_g f=\mbox{div}_g df
\label{9.211}
\end{equation}
and $\mbox{div}_g df$ is given by \ref{9.170} with the 1-form $df$ in the role of the 1-form $\xi$. 
Placing ${\cal N}$ in the role of ${\cal M}$ and $h_N$, $e_N$ or $h^\prime_N$, $e^\prime_N$ 
in the roles of $g$, $\epsilon$, we have:
\begin{equation}
d\s^{*h_N}df=e_N \square_{h_N}f, \ \ \ d\s^{*h^\prime_N}df=e^\prime_N\square_{h^\prime_N}f
\label{9.212}
\end{equation}
In view of \ref{9.172}, \ref{9.175} we then deduce:
\begin{eqnarray}
&&a_N(\square_{h_N}f-\square_{h^\prime_N}f)
=\nonumber\\
&&\hspace{18mm}-\frac{1}{2}\left(1-\frac{v_N}{2a_N}\right)^{-1}
\left\{(d(\s^{*h_N}df-\s^{*h^\prime_N}df))(L_N,\Lb_N,E_N)\right.\nonumber\\
&&\hspace{50mm}\left.-v_N\square_{h^\prime_N}f\right\}
\label{9.213}
\end{eqnarray}
Now from the first members of equations \ref{9.183}, \ref{9.187}, together with \ref{9.184} and 
\ref{9.172}, \ref{9.175}  we obtain:
\begin{eqnarray}
&&(\s^{*h_N}df-\s^{*h^\prime_N}df)(L_N,\Lb_N)=
-2a_N\left(1-\frac{v_N}{2a_N}\right)(\delta^\prime_{h,N}\cdot df)^E+v_N(h_N^{\prime-1}\cdot df)^E
\nonumber\\
&&(\s^{*h_N}df-\s^{*h^\prime_N}df)(L_N,E_N)=
2a_N\left(1-\frac{v_N}{2a_N}\right)(\delta^\prime_{h,N}\cdot df)^{\Lb}-v_N(h_N^{\prime-1}\cdot df)^{\Lb}
\nonumber\\
&&(\s^{*h_N}df-\s^{*h^\prime_N}df)(\Lb_N,E_N)=
-2a_N\left(1-\frac{v_N}{2a_N}\right)(\delta^\prime_{h,N}\cdot df)^L+v_N(h_N^{\prime-1}\cdot df)^L
\nonumber\\
&&\label{9.214}
\end{eqnarray}
We have:
\begin{eqnarray}
&&2a_N(\delta^\prime_{h,N}\cdot df)^E=2a_N(\delta_N^{\prime LE}L_N f+\delta_N^{\prime\Lb E}\Lb f
+\delta_N^{\prime EE}E_N f)\nonumber\\
&&2a_N(\delta^\prime_{h,N}\cdot df)^{\Lb}=2a_N(\delta_N^{\prime L\Lb} L_N f
+\delta_N^{\prime\Lb\Lb}\Lb_N f+\delta_N^{\prime\Lb E}E_N f)\nonumber\\
&&2a_N(\delta^\prime_{h,N}\cdot df)^L=2a_N(\delta_N^{\prime LL}L_N f+\delta_N^{\prime L\Lb}\Lb f
+\delta_N^{\prime LE}E_N f
\label{9.215}
\end{eqnarray}
Since $f$ is a $C^\infty$ function of $(\tau,\sigma,\vartheta)$ the estimates \ref{9.a29} then 
imply:
\begin{eqnarray}
&&2a_N(\delta^\prime_{h,N}\cdot df)^E=O(\tau^N)\nonumber\\
&&2a_N(\delta^\prime_{h,N}\cdot df)^{\Lb}=\rhob_N^{-1}O(\tau^N)+O(\tau^{N-1})\nonumber\\
&&2a_N(\delta^\prime_{h,N}\cdot df)^L=O(\tau^{N-1})+\rhob_N^{-1}O(\tau^{N+1})
\label{9.216}
\end{eqnarray}
In regard to the second terms in equations \ref{9.214} we have:
\begin{equation} 
2a_N(h_N^{\prime-1}\cdot df)^E=2a_N E_N f, \ \ \ 2a_N(h_N^{\prime-1}\cdot df)^{\Lb}=-L_N f, \ \ \ 
2a_N(h_N^{\prime-1}\cdot df)^L=-\Lb_N f
\label{9.217}
\end{equation}
and by \ref{9.178}:
\begin{equation}
\frac{v_N}{2a_N}=O(\tau^{N-1})+\rhob_N^{-1}O(\tau^{N+1})
\label{9.218}
\end{equation}
Defining the 2-form:
\begin{equation}
\theta=\s^{*h_N}df-\s^{*h^\prime_N}df
\label{9.219}
\end{equation}
the estimates \ref{9.216} together with the above yield, through \ref{9.214}, 
\begin{eqnarray}
&&\theta(L_N,\Lb_N)=O(\tau^N)\nonumber\\
&&\theta(L_N,E_N)=\rhob_N^{-1}O(\tau^N)+O(\tau^{N-1})\nonumber\\
&&\theta(\Lb_N,E_N)=O(\tau^{N-1})+\rhob_N^{-1}O(\tau^{N+1})
\label{9.220}
\end{eqnarray}
We now apply the formula \ref{9.193}. In view of the commutation relations \ref{9.b1} we obtain:
\begin{eqnarray}
&&(d\theta)(L_N,\Lb_N,E_N)=L_N(\theta(\Lb_N,E_N))+\Lb_N(\theta(E_N,L_N))+E_N(\theta(L_N,\Lb_N))
\nonumber\\
&&\hspace{32mm}+\chib_N\theta(E_N,L_N)-\chi_N\theta(E_N,\Lb_N)
\label{9.221}
\end{eqnarray}
Taking into account the fact that 
\begin{equation}
L_N\rhob_N=O(\tau) \ \ \mbox{as $\left.L_N\rhob_N\right|_{\tau=0}=\left.L\rhob\right|_{\tau=0}=0$}
\label{9.222}
\end{equation}
hence
\begin{equation}
L_N(\rhob_N^{-1})=\rhob_N^{-2}O(\tau)
\label{9.223}
\end{equation}
the estimates \ref{9.220} imply:
\begin{eqnarray}
&&L_N(\theta(\Lb_N,E_N))=O(\tau^{N-2})+\rhob_N^{-1}O(\tau^N)+\rhob_N^{-2}O(\tau^{N+2})\nonumber\\
&&\Lb_N(\theta(E_N,L_N))=\rhob_N^{-1}O(\tau^N)+O(\tau^{N-1})\nonumber\\
&&E_N(\theta(L_N,\Lb_N))=O(\tau^N) 
\label{9.224}
\end{eqnarray}
Also, taking into account the fact that
\begin{equation}
\chi_N=O(\tau) \ \ \mbox{as $\left.\chi_N\right|_{\tau=0}=\left.\chi\right|_{\tau=0}=0$}
\label{9.225}
\end{equation}
we have:
\begin{eqnarray}
&&\chib_N\theta(E_N,L_N)=\rhob_N^{-1}O(\tau^N)+O(\tau^{N-1})\nonumber\\
&&\chi_N\theta(E_N,\Lb_N)=O(\tau^N)+\rhob_N^{-1}O(\tau^{N+2})
\label{9.226}
\end{eqnarray}
We conclude that:
\begin{equation}
(d\theta)(L_N,\Lb_N,E_N)=O(\tau^{N-2})+\rhob_N^{-1}O(\tau^N)+\rhob_N^{-2}O(\tau^{N+2})
\label{9.227}
\end{equation}
This controls the first term in parenthesis on the right in \ref{9.213}. As for the second term, we note 
that by the formula \ref{9.210}, $a_N\square_{h_N^\prime}f$ is a $C^\infty$ function of 
$(\tau,\sigma,\vartheta)$. Then by \ref{9.218} this term is 
$$O(\tau^{N-1})+\rhob_N^{-1}O(\tau^{N+1})$$
In conclusion, we have established the following lemma.

\vspace{2.5mm}

\noindent{\bf Lemma 9.3} \ \ Let $f$ be a known $C^\infty$ function of the coordinates 
$(\tau,\sigma,\vartheta)$. We then have:
$$a_N(\square_{h_N}f-\square_{h_N^\prime}f)=
O(\tau^{N-2})+\rhob_N^{-1}O(\tau^N)+\rhob_N^{-2}O(\tau^{N+2})$$

\vspace{2.5mm}

We now turn to the difference \ref{9.200}. We express this in analogy with \ref{9.201} as:
\begin{equation}
a\Omega\square_{\tilde{h}}f-a_N\Omega_N\square_{\tilde{h}_N}f=
a\Omega\square_{\tilde{h}}f-a_N\Omega_N\square_{\tilde{h}^\prime_N}f 
-a_N\Omega_N(\square_{\tilde{h}_N}f-\square_{\tilde{h}^\prime_N}f)
\label{9.228}
\end{equation}
We shall apply the following. If $({\cal M},g)$ is an arbitrary $n+1$ dimensional Lorentzian manifold,  
and $\Omega$ an arbitrary positive function on ${\cal M}$, then with $\tilde{g}$ the conformal metric 
$\tilde{g}=\Omega g$ we have, for any function $f$ on ${\cal M}$,
\begin{equation}
\Omega\square_{\tilde{g}}f=\square_g f+\frac{(n-1)}{2}\Omega^{-1}g^{-1}(d\Omega,df)
\label{9.229}
\end{equation}
Applying this placing ${\cal N}$ in the role of ${\cal M}$, so here $n=2$, and placing $h$, $\Omega$ 
or $h^\prime_N$, $\Omega_N$ in the roles of $g$ and $\Omega$ respectively, we obtain:
\begin{eqnarray}
&&a\Omega\square_{\tilde{h}}f=a\square_h f+\frac{1}{2}\Omega^{-1}\left(-(1/2)((L\Omega)\Lb f
+(\Lb\Omega)Lf)+a(E\Omega)Ef\right)\nonumber\\
&&a_N\Omega_N\square_{\tilde{h}^\prime_N}f=a_N\square_{h^\prime_N}f+\frac{1}{2}\Omega_N^{-1}
\left(-(1/2)((L_N\Omega_N)\Lb_N f+(\Lb_N\Omega_N)L_N f)\right.\nonumber\\
&&\hspace{52mm}\left.+a_N(E_N\Omega_N)E_N f\right)\label{9.230}
\end{eqnarray}
The above together with \ref{9.205}, \ref{9.210} allow us to express the difference 
$$a\Omega\square_{\tilde{h}}f-a_N\Omega_N\square_{\tilde{h}^\prime_N}f$$
As for  
$$a_N\Omega_N(\square_{\tilde{h}_N}f-\square_{\tilde{h}^\prime_N}f)$$
we apply \ref{9.229} placing ${\cal N}$ in the role of ${\cal M}$, $h_N$ or $h^\prime_N$ in the roles 
of $g$, and $\Omega_N$ in the role of $\Omega$, to obtain:
\begin{eqnarray}
&&a_N\Omega_N\square_{\tilde{h}_N}f=a_N\square_{h_N}f
+\frac{1}{2}\Omega_N^{-1}a_N h_N^{-1}(d\Omega_N,df)\nonumber\\
&&a_N\Omega_N\square_{\tilde{h}^\prime_N}f=a_N\square_{h^\prime_N}f
+\frac{1}{2}\Omega_N^{-1}a_N h_N^{\prime-1}(d\Omega_N,df)
\label{9.231}
\end{eqnarray}
Taking the difference gives, in view of \ref{9.a25}, 
\begin{eqnarray}
&&a_N\Omega_N(\square_{\tilde{h}_N}f-\square_{\tilde{h}^\prime_N}f)
=a_N((\square_{h_N}f-\square_{h^\prime_N}f)\nonumber\\
&&\hspace{38mm}+\frac{1}{2}\Omega_N^{-1}a_N\delta^\prime_{h,N}(d\Omega_N,df)
\label{9.232}
\end{eqnarray}
The fact that $\Omega_N$ is a $C^\infty$ function of the coordinates $(\tau,\sigma,\vartheta)$ which is 
bounded from below by a positive constant together with the estimates \ref{9.a29} imply:
\begin{equation}
\frac{1}{2}\Omega_N^{-1}a_N\delta^\prime_{h,N}(d\Omega_N,df)=\rhob_N^{-1}O(\tau^N)+O(\tau^{N-1})
\label{9.233}
\end{equation}
Combining this with Lemma 9.3 yields the following. 

\vspace{2.5mm}

\noindent{\bf Lemma 9.4} \ \ Let $f$ be a known $C^\infty$ function of the coordinates 
$(\tau,\sigma,\vartheta)$. We then have:
$$a_N\Omega_N(\square_{\tilde{h}_N}f-\square_{\tilde{h}^\prime_N}f)=
O(\tau^{N-2})+\rhob_N^{-1}O(\tau^N)+\rhob_N^{-2}O(\tau^{N+2})$$

\vspace{2.5mm}

Placing in the above lemma the known smooth functions $\beta_{\mu,N}$ in the role of $f$ and recalling 
\ref{9.198}, \ref{9.199} yields the following proposition. 

\vspace{2.5mm}

\noindent{\bf Proposition 9.10} \ \ We have: 
$$\Omega_N a_N \square_{\tilde{h}^\prime_N}\beta_{\mu,N}=\tilde{\kappa}^\prime_{\mu,N}$$
where 
$$\tilde{\kappa}^\prime_{\mu,N}={\bf O}(\tau^{N-2})+{\bf O}(u^{-3}\tau^N)$$
in the sense of Definition 9.2.

\vspace{5mm}

\section{The Variation Differences $\s^{(m,l)}\check{\dot{\phi}}_\mu$ and the Rescaled Source Differences 
$\s^{(m,l)}\check{\tilde{\rho}}_\mu$}

Recalling the fundamental variations \ref{6.22}, we define the fundamental variation differences:
\begin{equation}
\check{\dot{\phi}}_\mu=\beta_\mu-\beta_{\mu,N}
\label{9.234}
\end{equation}
Since 
$$\Omega a\square_{\tilde{h}}\beta_\mu=0, \ \ \ \Omega_N a_N\square_{\tilde{h}^\prime_N}\beta_{\mu,N}=
\tilde{\kappa}^\prime_{\mu,N}$$
the fundamental variation differences satisfy:
\begin{equation}
\Omega a\square_{\tilde{h}}\check{\dot{\phi}}_\mu=
-(\Omega a\square_{\tilde{h}}-\Omega_N a_N\square_{\tilde{h}^\prime_N}) \beta_{\mu,N}
-\tilde{\kappa}^\prime_{\mu,N}
\label{9.235}
\end{equation}
Here on the right hand side the 1st term is estimated through \ref{9.230}, and \ref{9.205}, \ref{9.210} 
with $\beta_{\mu,N}$ in the role of $f$, while the 2nd term is estimated by Proposition 9.10. 

Similarly, recalling the higher order variations \ref{8.1}, we define the higher order variation 
differences:
\begin{equation}
\s^{(m,l)}\check{\dot{\phi}}_\mu=E^l T^m\beta_\mu-E_N^l T^m\beta_{\mu,N}
\label{9.236}
\end{equation}
By \ref{8.3} and \ref{8.11}:
\begin{equation}
\Omega a\square_{\tilde{h}}(E^l T^m\beta_\mu)=\s^{(m,l)}\tilde{\rho}_\mu 
\label{9.237}
\end{equation}
the rescaled sources $\s^{(m,l)}\tilde{\rho}_\mu$ being determined by the recursion \ref{8.28} with 
$l=0$, that is:
\begin{equation}
\s^{(m+1,0)}\tilde{\rho}_\mu=\hat{T}\s^{(m,0)}\tilde{\rho}_\mu+\s^{(T,m,0)}\sigma_\mu 
\label{9.238}
\end{equation}
and the homogeneous initial condition \ref{8.30}:
\begin{equation}
\s^{(0,0)}\tilde{\rho}_\mu=0
\label{9.239}
\end{equation}
followed by the recursion \ref{8.29}:
\begin{equation}
\s^{(m,l+1)}\tilde{\rho}_\mu=\hat{E}\s^{(m,l)}\tilde{\rho}_\mu+\s^{(E,m,l)}\sigma_\mu 
\label{9.240}
\end{equation}
In \ref{9.238} and \ref{9.240} we denote, for any commutation field $C$, here $T$ or $E$, by 
$\hat{C}$ the operator, acting on functions, 
\begin{equation}
\hat{C}=C+\s^{(C)}\delta
\label{9.241}
\end{equation}
the 2nd term being a multiplication operator. 

The $E_N^l T^m\beta_{\mu,N}$ similarly satisfy equations of the form:
\begin{equation}
\Omega_N a_N\square_{\tilde{h}^\prime_N}(E_N^l T^m\beta_{\mu,N})=\s^{(m,l)}\tilde{\rho}^\prime_{\mu,N}
\label{9.242}
\end{equation}
which follow from the equation of Proposition 9.10 by applying Lemma 8.1 with the metric 
$\tilde{h}^\prime_N$ in the role of the metric $\tilde{h}$. Then $\s^{(C)}\tilde{\pi}$ 
is replaced by: 
\begin{equation}
\s^{(C)}\tilde{\pi}_N={\cal L}_C\tilde{h}^\prime_N
\label{9.243}
\end{equation}
and $\s^{(C)}J$, the commutation current associated to $\psi$ and $C$, is replaced by $\s^{(C)}J_N$, 
a vectorfield given in arbitrary local coordinates by:
\begin{equation}
\s^{(C)}J_N=\left(\s^{(C)}\tilde{\pi}_N^{\alpha\beta}-\frac{1}{2}(\tilde{h}_N^{\prime-1})^{\alpha\beta}
\tilde{\mbox{tr}}\s^{(C)}\tilde{\pi}_N\right)\partial_\beta\psi
\label{9.244}
\end{equation}
where
\begin{equation}
\s^{(C)}\tilde{\pi}_N=
(\tilde{h}_N^{\prime-1})^{\alpha\gamma}(\tilde{h}_N^{\prime-1})^{\beta\delta}
\s^{(C)}\tilde{\pi}_{\gamma\delta,N}, \ \ 
\tilde{\mbox{tr}}\s^{(C)}\tilde{\pi}_N=
(\tilde{h}_N^{\prime-1})^{\alpha\beta}\tilde{\pi}_{\alpha\beta,N}
\label{9.245}
\end{equation}
With these replacements, Lemma 8.1 states that if $\psi$ is a solution of 
$$\square_{\tilde{h}^\prime_N}\psi=\rho$$ 
and $C$ is an arbitrary vectorfield, then 
$$\s^{(C)}\psi=C\psi$$
satisfies the equation
$$\square_{\tilde{h}^\prime_N}\s^{(C)}\psi=\s^{(C)}\rho$$
where 
$$\s^{(C)}\rho=C\rho+\tilde{\mbox{div}}_N\s^{(C)}J_N
+\frac{1}{2}\tilde{\mbox{tr}}\s^{(C)}\tilde{\pi}_N\rho$$
$\tilde{\mbox{div}}_N$ being the covariant divergence operator associated to the metric 
$\tilde{h}^\prime_N$. Defining then, in analogy with \ref{8.24}, 
\begin{equation}
\tilde{\rho}=\Omega_N a_N\rho, \ \ \ \s^{(C)}\tilde{\rho}=\Omega_N a_N\s^{(C)}\rho
\label{9.246}
\end{equation}
$\s^{(C)}\tilde{\rho}$ is related to $\tilde{\rho}$ according to the lemma by:
\begin{equation}
\s^{(C)}\tilde{\rho}=C\tilde{\rho}+\s^{(C)}\sigma_N+\s^{(C)}\delta_N\tilde{\rho}
\label{9.247}
\end{equation}
where
\begin{equation}
\s^{(C)}\sigma_N=\Omega_N a_N\tilde{\mbox{div}}_N\s^{(C)}J_N
\label{9.248}
\end{equation}
and
\begin{equation}
\s^{(C)}\delta_N=\frac{1}{2}\tilde{\mbox{tr}}\s^{(C)}\tilde{\pi}_N-a_N^{-1}Ca_N-\Omega_N^{-1}C\Omega_N
\label{9.249}
\end{equation}
The function $\s^{(C)}\sigma_N$ is associated to the commutation field $C$ and to the function $\psi$ 
through the metric $\tilde{h}^\prime_N$. 

From the above we deduce for the functions $\s^{(m,l)}\tilde{\rho}^\prime_{\mu,N}$ (see \ref{9.242}), 
in analogy with Section 8.2, the recursion formulas: 
\begin{equation}
\s^{(m+1,0)}\tilde{\rho}^\prime_{\mu,N}=\hat{T}_N\s^{(m,0)}\tilde{\rho}^\prime_{\mu,N}
+\s^{(T,m,0)}\sigma_{\mu,N}
\label{9.250}
\end{equation}
\begin{equation}
\s^{(m,l+1)}\tilde{\rho}^\prime_{\mu,N}=\hat{E}_N \s^{(m,l)}\tilde{\rho}^\prime_{\mu,N}
+\s^{(E_N,m,l)}\sigma_{\mu,N}
\label{9.251}
\end{equation}
Here $\hat{T}_N$ and $\hat{E}_N$ are the operators:
\begin{equation}
\hat{T}_N=T+\s^{(T)}\delta_N, \ \ \ \hat{E}_N=E_N+\s^{(E_N)}\delta_N
\label{9.252}
\end{equation}
Also, $\s^{(T,m,0)}\sigma_N$ is the function $\s^{(C)}\sigma_N$, as above, associated to the 
commutation field $C=T$ and to the function $T^m\beta_{\mu,N}$, and $\s^{(E_N,m,l)}\sigma_N$ is 
the function $\s^{(C)}\sigma_N$, as above, associated to the commutation field $C=E_N$ and to the 
function $E_N^l T^m\beta_{\mu,N}$. In contrast to the recursion \ref{9.238}, the initial condition 
for the recursion \ref{9.250} is inhomogeneous, as according to Proposition 9.10 we have:
\begin{equation}
\s^{(0,0)}\tilde{\rho}^\prime_{\mu,N}=\tilde{\kappa}^\prime_{\mu,N}
\label{9.253}
\end{equation}

To the recursions \ref{9.238}, \ref{9.240} and the recursions \ref{9.250}, \ref{9.251} we shall apply 
the following elementary lemma on linear recursions. 

\vspace{2.5mm}

\noindent{\bf Lemma 9.5} \ \ Let $(y_n:n=0,1,2,...)$ be a given sequence in a space $X$ and $A$ 
a linear operator in $X$. Suppose that $(x_n:n=0,1,2,...)$ is a sequence in $X$ satisfying the 
recursion:
$$x_{n+1}=Ax_n+y_n$$
Then for each $n=0,1,2,...$ we have:
$$x_n=A^n x_0+\sum_{m=0}^{n-1}A^{n-1-m}y_m$$

\vspace{2.5mm}

Applying this lemma to the recursion \ref{9.238}, in view of the homogeneous initial condition 
\ref{9.239} we deduce:
\begin{equation}
\s^{(m,0)}\tilde{\rho}_\mu=\sum_{n=0}^{m-1}\hat{T}^{m-1-n}\s^{(T,n,0)}\sigma_\mu 
\label{9.254}
\end{equation}
Applying then the lemma to the recursion \ref{9.240} we deduce:
\begin{equation}
\s^{(m,l)}\tilde{\rho}_\mu=\hat{E}^l\s^{(m,0)}\tilde{\rho}_\mu 
+\sum_{k=0}^{l-1}\hat{E}^{l-1-k}\s^{(E,m,k)}\sigma_\mu 
\label{9.255}
\end{equation}
Substituting for $\s^{(m,0)}\tilde{\rho}_\mu$ from \ref{9.254} then yields:
\begin{equation}
\s^{(m,l)}\tilde{\rho}_\mu=\sum_{n=0}^{m-1}\hat{E}^l\hat{T}^{m-1-n}\s^{(T,n,0)}\sigma_\mu
+\sum_{k=0}^{l-1}\hat{E}^{l-1-k}\s^{(E,m,k)}\sigma_\mu 
\label{9.256}
\end{equation}
Applying Lemma 9.5 to the recursion \ref{9.251}, in view of the inhomogeneous initial condition 
\ref{9.253} we deduce:
\begin{equation}
\s^{(m,0)}\tilde{\rho}^\prime_{\mu,N}=\hat{T}_N^m\tilde{\kappa}^\prime_{\mu,N}
+\sum_{n=0}^{m-1}\hat{T}_N^{m-1-n}\s^{(T,n,0)}\sigma_{\mu,N}
\label{9.257}
\end{equation}
Applying then Lemma 9.5 to the recursion \ref{9.252} we deduce:
\begin{equation}
\s^{(m,l)}\tilde{\rho}^\prime_{\mu,N}=\hat{E}_N^l\s^{(m,0)}\tilde{\rho}^\prime_{\mu,N}
+\sum_{k=0}^{l-1}\hat{E}_N^{l-1-k}\s^{(E_N,m,k)}\sigma_{\mu,N}
\label{9.258}
\end{equation}
Substituting for $\s^{(m,0)}\tilde{\rho}^\prime_{\mu,N}$ from \ref{9.257} then yields:
\begin{eqnarray}
&&\s^{(m,l)}\tilde{\rho}^\prime_{\mu,N}=\hat{E}_N^l\hat{T}_N^m\tilde{\kappa}^\prime_{\mu,N}
+\sum_{n=0}^{m-1}\hat{E}_N^l\hat{T}_N^{m-1-n}\s^{(T,n,0)}\sigma_{\mu,N}\nonumber\\
&&\hspace{30mm}+\sum_{k=0}^{l-1}\hat{E}_N^{l-1-k}\s^{(E_N,m,k)}\sigma_{\mu,N}
\label{9.259}
\end{eqnarray}
We can write this in the form:
\begin{equation}
\s^{(m,l)}\tilde{\rho}^\prime_{\mu,N}=\hat{E}_N^l\hat{T}_N^m\tilde{\kappa}^\prime_{\mu,N}
+\s^{(m,l)}\tilde{\rho}_{\mu,N}
\label{9.260}
\end{equation}
where $\s^{(m,l)}\tilde{\rho}_{\mu,N}$ is the solution of the recursions \ref{9.251}, \ref{9.252} 
with the inhomogeneous initial condition \ref{9.253} replaced by the corresponding homogeneous one:
\begin{equation}
\s^{(0,0)}\tilde{\rho}_{\mu,N}=0
\label{9.261}
\end{equation}
Thus, we have:
\begin{equation}
\s^{(m,l)}\tilde{\rho}_{\mu,N}=\sum_{n=0}^{m-1}\hat{E}_N^l\hat{T}_N^{m-1-n}\s^{(T,n,0)}\sigma_{\mu,N}
+\sum_{k=0}^{l-1}\hat{E}_N^{l-1-k}\s^{(E_N,m,k)}\sigma_{\mu,N}
\label{9.262}
\end{equation}
In conclusion, we have established the following proposition. 

\vspace{2.5mm}

\noindent{\bf Proposition 9.11} \ \ The variation differences $\s^{(m,l)}\check{\dot{\phi}}_\mu$ 
satisfy
$$\Omega a\square_{\tilde{h}}\s^{(m,l)}\check{\dot{\phi}}_\mu=\s^{(m,l)}\check{\tilde{\rho}}_\mu$$
where the rescaled source differences $\s^{(m,l)}\check{\tilde{\rho}}_\mu$ are given by:
$$\s^{(m,l)}\check{\tilde{\rho}}_\mu=-(\Omega a\square_{\tilde{h}}-\Omega_N\square_{\tilde{h}^\prime_N})
E_N^l T^m\beta_{\mu,N}+\s^{(m,l)}\tilde{\rho}_\mu-\s^{(m,l)}\tilde{\rho}_{\mu,N}
-\hat{E}_N^l\hat{T}_N^m\tilde{\kappa}^\prime_{\mu,N}$$
with $\s^{(m,l)}\tilde{\rho}_\mu$, $\s^{(m,l)}\tilde{\rho}_{\mu,N}$ given by \ref{9.256}, \ref{9.262} 
respectively. 

\vspace{2.5mm}

We remark that the first term in $\s^{(m,l)}\check{\tilde{\rho}}_\mu$ is estimated 
through \ref{9.230}, and \ref{9.205}, \ref{9.210} 
with the known $C^\infty$ function $E_N^l T^m\beta_{\mu,N}$ in the role of $f$, while the last term 
is estimated through Proposition 9.10. 

We denote
\begin{equation}
\s^{(C)}\pi_N={\cal L}_C h^\prime_N
\label{9.263}
\end{equation}
Recalling from Chapter 8 that in the case $n=2$ we denote  
\begin{eqnarray*}
&&\s^{(C)}\pi_{LL}=\s^{(C)}\pi(L,L), \ \ \s^{(C)}\pi_{\Lb\Lb}=\s^{(C)}\pi(\Lb,\Lb), 
\ \ \s^{(C)}\pi_{L\Lb}=\s^{(C)}\pi(L,\Lb)\\
&&\s^{(C)}\spi_L=\s^{(C)}\pi(L,E), \ \ \s^{(C)}\spi_{\Lb}=\s^{(C)}\pi(\Lb,E), 
\ \ \s^{(C)}\sspi=\s^{(C)}\pi(E,E)
\end{eqnarray*}
the components of $\s^{(C)}\pi$ in the $(L,\Lb,E)$ frame, we similarly denote 
\begin{eqnarray}
&&\s^{(C)}\pi_{LL,N}=\s^{(C)}\pi_N(L_N,L_N), \ \ \s^{(C)}\pi_{\Lb\Lb,N}=\s^{(C)}\pi_N(\Lb_N,\Lb_N) 
\nonumber\\ 
&&\hspace{16mm}\s^{(C)}\pi_{L\Lb,N}=\s^{(C)}\pi_N(L_N,\Lb_N)\nonumber\\
&&\s^{(C)}\spi_{L,N}=\s^{(C)}\pi_N(L_N,E_N), \ \ \s^{(C)}\spi_{\Lb,N}=\s^{(C)}\pi_N(\Lb_N,E_N)
\nonumber\\
&&\hspace{17mm}\s^{(C)}\sspi_N=\s^{(C)}\pi(E_N,E_N)
\label{9.264}
\end{eqnarray}
the components of $\s^{(C)}\pi_N$ in the $(L_N,\Lb_N,E_N)$ frame. Since with $X,Y$ an arbitrary 
pair of vectorfields we have:
\begin{equation}
\s^{(C)}\pi_N(X,Y)=C(h_N^\prime(X,Y))-h_N^\prime([C,X],Y)-h_N^\prime(X,[C,Y])
\label{9.265}
\end{equation}
taking $C$ equal to $T$ or $E_N$ and selecting the pair $X,Y$ out of the frame field $(L_N,\Lb_N,E_N)$, 
by reason of the commutation relations \ref{9.b1} and table \ref{9.139}, we deduce tables for the components of $\s^{(T)}\pi_N$ and $\s^{(E_N)}\pi_N$ analogous to the tables \ref{8.50} and \ref{8.51}:
\begin{eqnarray}
&&\s^{(T)}\pi_{LL,N}=\s^{(T)}\pi_{\Lb\Lb,N}=0, \ \ \s^{(T)}\pi_{L\Lb,N}=-2Ta_N\nonumber\\
&&\hspace{5mm}\s^{(T)}\spi_{L,N}=-\zeta_N, \ \ \s^{(T)}\spi_{\Lb,N}=\zeta_N \nonumber\\
&&\hspace{10mm}\s^{(T)}\sspi_N=2(\chi_N+\chib_N)
\label{9.266}
\end{eqnarray}
\begin{eqnarray}
&&\s^{(E_N)}\pi_{LL,N}=\s^{(E_N)}\pi_{\Lb\Lb,N}=0, \ \ \s^{(E_N)}\pi_{L\Lb,N}=-2E_N a_N\nonumber\\
&&\hspace{5mm}\s^{(E_N)}\spi_{L,N}=-\chi_N, \ \ \s^{(E_N)}\spi_{\Lb,N}=-\chib_N\nonumber\\
&&\hspace{10mm}\s^{(E_N)}\sspi_N=0
\label{9.267}
\end{eqnarray}
Since $\tilde{h}_N^\prime=\Omega_N h_N^\prime$, we have:
\begin{equation}
\s^{(C)}\tilde{\pi}_N=\Omega_N\s^{(C)}\pi_N+(C\Omega_N)h_N^\prime
\label{9.268}
\end{equation}
Hence, the components of $\s^{(C)}\tilde{\pi}_N$ in the frame field $(L_N,\Lb_N,E_N)$ are given 
in terms of the components of $\s^{(C)}\pi_N$ by:
\begin{eqnarray}
&&\s^{(C)}\tilde{\pi}_{LL,N}=\Omega_N\s^{(C)}\pi_{LL,N}, \ \ 
\s^{(C)}\tilde{\pi}_{\Lb\Lb,N}=\Omega_N\s^{(C)}\pi_{\Lb\Lb,N}\nonumber\\
&&\hspace{8mm}\s^{(C)}\tilde{\pi}_{L\Lb,N}=\Omega_N(\s^{(C)}\pi_{L\Lb,N}
-2a_N\Omega_N^{-1}C\Omega_N)\nonumber\\
&&\hspace{5mm}\s^{(C)}\tilde{\spi}_{L,N}=\Omega_N\s^{(C)}\spi_{L,N}, \ \ 
\s^{(C)}\tilde{\spi}_{\Lb,N}=\Omega_N\s^{(C)}\spi_{\Lb,N}\nonumber\\
&&\hspace{10mm}\s^{(C)}\tilde{\sspi}_N=\Omega_N(\s^{(C)}\sspi_N+\Omega_N^{-1}C\Omega_N)
\label{9.269}
\end{eqnarray}
(compare with \ref{8.49}). Since
\begin{equation}
\tilde{\mbox{tr}}\s^{(C)}\tilde{\pi}_N=-a_N^{-1}\s^{(C)}\pi_{L\Lb,N}+\s^{(C)}\sspi_N+3\Omega_N^{-1}C\Omega_N
\label{9.270}
\end{equation}
(compare with \ref{8.a3}, here $n=2$) the functions $\s^{(C)}\delta_N$, defined by \ref{9.249} are 
expressed as:
\begin{equation}
\s^{(C)}\delta_N=-\frac{1}{2a_N}(\s^{(C)}\pi_{L\Lb,N}+2Ca_N)+\frac{1}{2}\s^{(C)}\sspi_N
+\frac{1}{2}\Omega_N^{-1}C\Omega_N
\label{9.271}
\end{equation}
(compare with \ref{8.a4}, here $n=2$). 
From tables \ref{9.266}, \ref{9.267} we see that:
\begin{equation}
\s^{(C)}\pi_{L\Lb,N}+2Ca_N=0 \ : \ \mbox{for both commutation fields $C=T,E_N$}
\label{9.272}
\end{equation}
(compare with \ref{8.a5}). Therefore \ref{9.271} reduces to:
\begin{equation}
\s^{(C)}\delta_N=\frac{1}{2}\s^{(C)}\sspi_N
+\frac{1}{2}\Omega_N^{-1}C\Omega_N
\label{9.273}
\end{equation}
Consequently, the functions $\s^{(C)}\delta_N$ are smooth functions of the coordinates 
$(\tau,\sigma,\vartheta)$ for both commutation fields $C=T,E_N$. Proposition 9.10 then implies 
the following estimate in regard to the last term in $\s^{(m,l)}\check{\tilde{\rho}}_{\mu}$ 
as given by Proposition 9.11:
\begin{equation}
\hat{E}_N^l\hat{T}_N^m\tilde{\kappa}^\prime_{\mu,N}={\bf O}(\tau^{N-m-2})+{\bf O}(u^{-3}\tau^{N-m})
\label{9.274}
\end{equation}

\vspace{5mm}

\section{The Difference 1-Forms $\s^{(V;m,l)}\check{\xi}$. The Difference Energies and the 
Difference Energy Identities}

To the variation differences $\s^{(m,l)}\check{\dot{\phi}}_\mu : \mu=0,1,2$ and to each of 
the two variation fields $V=Y,E$ we associate the difference 1-form:
\begin{equation}
\s^{(V;m,l)}\check{\xi}=V^\mu d\s^{(m,l)}\check{\dot{\phi}}_\mu 
\label{9.275}
\end{equation}
(see Section 6.2). To the 1-form $\s^{(V;m,l)}\check{\xi}$ we associate, in analogy with 
\ref{6.102} (or \ref{7.79}), the symmetric 2-covariant tensorfield:
\begin{equation}
\s^{(V;m,l)}\check{S}=\s^{(V;m,l)}\check{\xi}\otimes\s^{(V;m,l)}\check{\xi}
-\frac{1}{2}\tilde{h}^{-1}(\s^{(V;m,l)}\check{\xi},\s^{(V;m,l)}\check{\xi})\tilde{h}
\label{9.276}
\end{equation}
With the multiplier field $X$ fixed according to \ref{7.67}, we define, in analogy with \ref{6.55}, 
the difference energy current, the vectorfield given in arbitrary local coordinates by:
\begin{equation}
\s^{(V;m,l)}\check{P}^\alpha=-\s^{(V;m,l)}\check{S}^\alpha_\beta X^\beta
\label{9.277}
\end{equation}
Here, as in Chapter 6, indices from the beginning of the Greek alphabet are raised with respect to 
$\tilde{h}$. In particular, 
\begin{equation}
\s^{(V;m,l)}\check{S}^\alpha_\beta=(\tilde{h}^{-1})^{\alpha\gamma}\s^{(V;m,l)}\check{S}_{\gamma\beta}
\label{9.278}
\end{equation}
We then have, in analogy with \ref{6.56}, 
\begin{equation}
\tilde{D}_\alpha\s^{(V;m,l)}\check{P}^\alpha=\s^{(V;m,l)}\check{Q}
\label{9.279}
\end{equation}
where, in view of Proposition 9.11, $\s^{(V;m,l)}\check{Q}$ is given, in analogy with \ref{6.63} - 
\ref{6.66}, by:
\begin{equation}
\s^{(V;m,l)}\check{Q}=\s^{(V;m,l)}\check{Q}_1+\s^{(V;m,l)}\check{Q}_2+\s^{(V;m,l)}\check{Q}_3
\label{9.280}
\end{equation}
where:
\begin{eqnarray}
&&\s^{(V;m,l)}\check{Q}_1=-\frac{1}{2}\s^{(V;m,l)}\check{S}^{\alpha\beta}\s^{(X)}\tilde{\pi}_{\alpha\beta}
\label{9.281}\\
&&\s^{(V;m,l)}\check{Q}_2=-\s^{(V;m,l)}\check{\xi}^\alpha X^\beta
\left(\s^{(V)}\theta^\mu_\alpha \ \partial_\beta\s^{(m,l)}\check{\dot{\phi}}_\mu-
\s^{(V)}\theta^\mu_\beta \ \partial_\alpha\s^{(m,l)}\check{\dot{\phi}}_\mu\right)\nonumber\\
&&\hspace{20mm}-\s^{(V;m,l)}\check{\xi}_\beta X^\beta\ \s^{(V)}\theta^\mu_\alpha \ \partial^\alpha\s^{(m,l)}\check{\dot{\phi}}_\mu 
\label{9.282}\\
&&\s^{(V;m,l)}\check{Q}_3=-\s^{(V;m,l)}\check{\xi}_\beta X^\beta\ V^\mu(\Omega a)^{-1}
\s^{(m,l)}\check{\tilde{\rho}}_\mu
\label{9.283}
\end{eqnarray}
indices from the beginning of the Greek alphabet being raised with respect to 
$\tilde{h}$. In particular, 
\begin{eqnarray}
&&\s^{(V;m,l)}\check{S}^{\alpha\beta}=(\tilde{h}^{-1})^{\alpha\gamma}(\tilde{h}^{-1})^{\beta\delta}
\s^{(V;m,l)}\check{S}_{\gamma\delta}\label{9.284}\\
&&\s^{(V;m,l)}\check{\xi}^\alpha=(\tilde{h}^{-1})^{\alpha\beta}\s^{(V;m,l)}\check{\xi}_\beta, 
\ \ \ \partial^\alpha\s^{(m,l)}\check{\dot{\phi}}_\mu=
(\tilde{h}^{-1})^{\alpha\beta}\partial_\beta\s^{(m,l)}\check{\dot{\phi}}_\mu\nonumber
\end{eqnarray}

From \ref{9.279} we deduce following the same path as that leading from \ref{6.56} to \ref{6.99}, 
the $(m,l)$ {\em difference energy identity} associated to the variation field $V$:
\begin{equation}
\s^{(V;m,l)}\check{{\cal E}}^{\ub_1}(u_1)+\s^{(V;m,l)}\underline{\check{{\cal E}}}^{u_1}(\ub_1)
+\s^{(V;m,l)}\check{{\cal F}}^{\ub_1}-\s^{(V;m,l)}\underline{\check{{\cal E}}}^{u_1}(0)=
\s^{(V;m,l)}\check{{\cal G}}^{\ub_1,u_1}
\label{9.285}
\end{equation}
Here $\s^{(V;m,l)}\check{{\cal E}}^{\ub_1}(u_1)$, 
$\s^{(V;m,l)}\underline{\check{{\cal E}}}^{u_1}(\ub_1)$ are the $(m,l)$ {\em difference energies}, 
associated to the variation field $V$, defined in analogy to \ref{6.91}, \ref{6.93} by (here $n=2$):
\begin{equation}
\s^{(V;m,l)}\check{{\cal E}}^{\ub_1}(u_1):=\int_{C_{u_1}^{\ub_1}}2a\Omega^{3/2}
\s^{(V;m,l)}\check{P}^{\Lb}
\label{9.286}
\end{equation}
\begin{equation}
\s^{(V;m,l)}\underline{\check{{\cal E}}}^{u_1}(\ub_1):=\int_{\Cb_{\ub_1}^{u_1}}2a\Omega^{3/2}\s^{(V;m,l)}\check{P}^L 
\label{9.287}
\end{equation}
Also, $\s^{(V;m,l)}\check{{\cal F}}^{\ub_1}$ is the $(m,l)$ {\em difference flux} associated to the 
variation field $V$, defined in analogy to \ref{6.96} by:
\begin{equation}
\s^{(V;m,l)}\check{{\cal F}}^{\ub_1}:=
\int_{{\cal K}^{\ub_1}}2a\Omega^{3/2}\left(\s^{(V;m,l)}\check{P}^L
-\s^{(V;m,l)}\check{P}^{\Lb}\right)
\label{9.288}
\end{equation}
Finally, $\s^{(V;m,l)}\check{{\cal G}}^{\ub_1,u_1}$ is the $(m,l)$ {\em difference error integral} 
associated to the variation field $V$, defined in analogy to \ref{6.98} by:
\begin{equation}
\s^{(V;m,l)}\check{{\cal G}}^{\ub_1,u_1}:=\int_{{\cal R}_{\ub_1,u_1}}2a\Omega^{3/2}
\s^{(V;m,l)}\check{Q}
\label{9.289}
\end{equation}
The multiplier field $X$ being fixed according to \ref{7.67}, we have, in analogy with \ref{7.81}, 
\ref{7.82}, 
\begin{equation}
\s^{(V;m,l)}\check{{\cal E}}^{\ub_1}(u)=\int_{C_u^{\ub_1}}\Omega^{1/2}(3(\s^{(V;m,l)}\check{\xi}_L)^2
+a(\s^{(V;m,l)}\check{\sxi})^2)
\label{9.290}
\end{equation}
\begin{equation}
\s^{(V;m,l)}\underline{\check{{\cal E}}}^{u_1}(\ub)=\int_{\Cb_{\ub}^{u_1}}\Omega^{1/2}
(3a(\s^{(V;m,l)}\check{\sxi})^2+(\s^{(V;m,l)}\check{\xi}_{\Lb})^2)
\label{9.291}
\end{equation}
Also, by the analogues of \ref{6.100}, \ref{6.103} we have:
\begin{equation}
\s^{(V;m,l)}\check{{\cal F}}^{\tau_1}=\int_{{\cal K}^{\tau_1}}\Omega^{1/2}\s^{(V;m,l)}\check{S}(X,M)
\label{9.292}
\end{equation}
$M$ being the interior normal \ref{6.145} to ${\cal K}$, so with $X$ fixed according to \ref{7.67}, 
the analogue of the coercivity inequality \ref{7.65} applies:
\begin{equation}
\s^{(V;m,l)}\check{S}(X,M)\geq C^{-1}\left((\s^{(V;m,l)}\check{\xi}_T)^2
+a(\s^{(V;m,l)}\check{\sxi})^2\right)-C^\prime(\s^{(V;m,l)}\check{b})^2
\label{9.a32}
\end{equation}
where, in analogy with \ref{6.164}, \ref{6.165}, $\s^{(V;m,l)}\check{b}$ represents the boundary values 
on ${\cal K}$ of $\s^{(V;m,l)}\check{\xi}(B)$:
\begin{equation}
\s^{(V;m,l)}\check{\xi}(B)=\s^{(V;m,l)}\check{b} \ : \ \mbox{on ${\cal K}$}
\label{9.a33}
\end{equation}
The functions $\s^{(V;m,l)}\check{b}$ shall be analyzed in Chapter 13. Defining then 
\begin{equation}
\s^{(V;m,l)}\check{{\cal F}}^{\prime\tau_1}=\s^{V;m,l)}\check{{\cal F}}^{\tau_1}
+2C^\prime\int_{{\cal K}^{\tau_1}}\Omega^{1/2}(\s^{(V;m,l)}\check{b})^2
\label{9.a34} 
\end{equation}
$C^\prime$ being the constant in the last term in \ref{9.a32}, we have:
\begin{eqnarray}
&&\s^{(V;m,l)}\check{{\cal F}}^{\prime\tau_1}\geq C^{-1}\int_{{\cal K}^{\tau_1}}\Omega^{1/2}\left((\s^{(V;m,l)}\check{\xi}_T)^2+a(\s^{(V;m,l)}\check{\sxi})^2\right)
\nonumber\\
&&\hspace{20mm}+C^\prime\int_{{\cal K}^{\tau_1}}\Omega^{1/2}(\s^{(V;m,l)}\check{b})^2
\label{9.a35}
\end{eqnarray}
Now, by \ref{6.156} - \ref{6.158}: 
\begin{equation}
L=B-\left(\frac{1}{4}\hat{\Lambda}\beta_T-\frac{1}{2}\right)T+a\hat{\Lambda}\sbeta E
\label{9.a36}
\end{equation}
hence by \ref{9.a33}:
\begin{equation}
\s^{(V;m,l)}\check{\xi}_L=\s^{(V;m,l)}\check{b}
-\left(\frac{1}{4}\hat{\Lambda}\beta_T-\frac{1}{2}\right)\s^{(V;m,l)}\check{\xi}_T
+a\hat{\Lambda}\sbeta\s^{(V;m,l)}\check{\sxi}
\label{9.a37}
\end{equation}
In view of \ref{6.191}, \ref{6.193} we can assume that:
\begin{equation}
\left|\frac{1}{2}\hat{\Lambda}\beta_T-1\right|\leq C\tau, \ \ 
\sqrt{a}|\hat{\Lambda}|\leq C\tau^{1/2} \ \ \mbox{: in ${\cal K}^\delta$}
\label{9.a38}
\end{equation}
It then follows from \ref{9.a33}, \ref{9.a35} that:
\begin{equation}
\int_{{\cal K}^{\tau_1}}
\left((\s^{(V;m,l)}\check{\xi}_L)^2+(\s^{(V;m,l)}\check{\xi}_{\Lb})^2
+a(\s^{(V;m,l)}\check{\sxi})^2\right)\leq C\s^{(V;m,l)}\check{{\cal F}}^{\prime\tau_1}
\label{9.a39}
\end{equation}
for a different constant $C$. 

Now, the definition of the $N$th order approximants \ref{9.3} implies, through \ref{9.236} that:
\begin{equation}
\left.\s^{(m,l)}\check{\dot{\phi}}_\mu\right|_{\Cb_0}=0 \ \ \mbox{: if $N\geq m$}
\label{9.293}
\end{equation}
hence for both $V=Y,E$:
\begin{equation}
\left.\s^{(V;m,l)}\check{\sxi}\right|_{\Cb_0}=\left.\s^{(V;m,l)}\check{\xi}_{\Lb}\right|_{\Cb_0}=0 
\ \ \mbox{: if $N\geq m$}
\label{9.294}
\end{equation}
Choosing then $N\geq m$, we have, from \ref{9.291},
\begin{equation}
\s^{(V;m,l)}\underline{\check{{\cal E}}}^{u_1}(0)=0
\label{9.295}
\end{equation}
Therefore the $(m,l)$ {\em difference energy identity} \ref{9.285} reduces to:
\begin{equation}
\s^{(V;m,l)}\check{{\cal E}}^{\ub_1}(u_1)+\s^{(V;m,l)}\underline{\check{{\cal E}}}^{u_1}(\ub_1)
+\s^{(V;m,l)}\check{{\cal F}}^{\ub_1}=
\s^{(V;m,l)}\check{{\cal G}}^{\ub_1,u_1}
\label{9.296}
\end{equation}

Finally, we define the weighted quantities: 
\begin{equation}
\s^{(V;m,l)}\cB(\ub_1,u_1)=\sup_{(\ub,u)\in R_{\ub_1,u_1}}\ub^{-2a_m}u^{-2b_m}
\s^{(V;m,l)}\cE^{\ub}(u)
\label{9.297}
\end{equation}
\begin{equation}
\s^{(V;m,l)}\cBb(\ub_1,u_1)=\sup_{(\ub,u)\in R_{\ub_1,u_1}}\ub^{-2a_m}u^{-2b_m}
\s^{(V;m,l)}\cEb^u(\ub)
\label{9.298}
\end{equation}
and:
\begin{equation}
\s^{(V;m,l)}\cA(\tau_1)=\sup_{\tau\in[0,\tau_1]}\tau^{-2(a_m+b_m)}\s^{(V;m,l)}\cF^{\prime\tau}
\label{9.299}
\end{equation}
the exponents $a_m$, $b_m$ being non-negative real numbers which are non-increasing with $m$. Remark that $(\ub_1,u_1)\in R_{\ub_2,u_2}$, 
that is $\ub_2\geq\ub_1$ and $u_2\geq u_1$, implies:
$$\s^{(V;m,l)}\cB(\ub_2,u_2)\geq\s^{(V;m,l)}\cB(\ub_1,u_1), \ \ \ 
\s^{(V;m,l)}\cBb(\ub_2,u_2)\geq\s^{(V;m,l)}\cBb(\ub_1,u_1)$$
and $\tau_2\geq\tau_1$ implies:
$$\s^{(V;m,l)}\cA(\tau_2)\geq\s^{(V;m,l)}\cA(\tau_1)$$ 
We also define:
\begin{equation}
\s^{(m,l)}\cB=\max_{V=Y,E}\s^{(V;m,l)}\cB, \ \ \ \s^{(m,l)}\cBb=\max_{V=Y,E}\s^{(V;m,l)}\cBb
\label{9.300}
\end{equation}
\begin{equation}
\s^{[m,l]}\cB=\max_{i=0,...,m}\s^{(i,l+m-i)}\cB, \ \ \ 
\s^{[m,l]}\cBb=\max_{i=0,...,m}\s^{(i,l+m-i)}\cBb
\label{9.301}
\end{equation}
and:
\begin{equation}
\s^{(m,l)}\cA=\max_{V=Y,E}\s^{(V;m,l)}\cA, \ \ \ 
\s^{[m,l]}\cA=\max_{i=0,...,m}\s^{(i,l+m-i)}\cA
\label{9.302}
\end{equation}

\pagebreak

\chapter{The Top Order Acoustical Estimates in the Case $d=2$}

\section{Regularization of the Propagation Equations for $\tchi$ and $\tchib$}

The variation differences $\s^{(m,l)}\check{\dot{\phi}}_\mu$, defined by \ref{9.236}, being of order 
$m+l$, the difference 1-forms $\s^{(V;m,l)}\check{\xi}$, defined by \ref{9.275}, are of order 
$m+l+1$, hence the $(m,l)$ difference energies $\s^{(V;m,l)}\cE^{\ub_1}(u)$, 
$\s^{(V;m,l)}\cEb^{u_1}(\ub)$, defined by \ref{9.290}, \ref{9.291}, are of order $m+l+1$. 
Let $n+1$ be the order of the top order estimates to be considered. Then for the top order 
difference energies we have $m+l=n$. The associated difference error integral, defined by \ref{9.289}, 
is, in view of \ref{9.280}, the sum:
\begin{equation}
\s^{(V;m,l)}G^{\ub_1,u_1}=\s^{(V;m,l)}G_1^{\ub_1,u_1}+\s^{(V;m,l)}G_2^{\ub_1,u_1}
+\s^{(V;m,l)}G_3^{\ub_1,u_1}
\label{10.1}
\end{equation}
where:
\begin{equation}
\s^{(V;m,l)}\cG_i^{\ub_1,u_1}=\int_{{\cal R}_{\ub_1,u_1}}2a\Omega^{3/2}\s^{(V;m,l)}\check{Q}_i 
\ \ : \ i=1,2,3
\label{10.2}
\end{equation}
the $\s^{(V;m,l)}\check{Q}_i:i=1,2,3$ being given by \ref{9.281}, \ref{9.282}, \ref{9.283}. 
The difference error integral containing the top order acoustical quantities is 
$\s^{(V;m,l)}\cG_3^{\ub_1,u_1}$, which according to \ref{9.283} is given by ($X$ has been fixed by 
\ref{7.67}):
\begin{equation}
\s^{(V;m,l)}\cG_3^{\ub_1,u_1}=-2\int_{{\cal R}_{\ub_1,u_1}}\Omega^{1/2}
(3\s^{(V;m,l)}\check{\xi}_L+\s^{(V;m,l)}\check{\xi}_{\Lb})V^\mu\s^{(m,l)}\check{\tilde{\rho}}_\mu 
\label{10.3}
\end{equation}
The principal acoustical part of $\s^{(m,l)}\check{\tilde{\rho}}_\mu$ is given by \ref{8.158}, 
\ref{8.159}. Since here $m+l=n$, the top order, $n+1$, acoustical quantities appearing are:
\begin{eqnarray}
&&\mbox{for $m=0$} \ \ : \ E^l\tchi, \ E^l\tchib \nonumber\\
&&\mbox{for $m\geq 1$} \ \ : \ E^{l+2}T^{m-1}\lambda, \ E^{l+2}T^{m-1}\lambdab 
\label{10.4}
\end{eqnarray}
The quantities $ET^{n-1}\lambda$, $ET^{n-1}\lambdab$, or $T^n\lambda$, $T^n\lambdab$, do not appear. 

We must estimate the top order acoustical quantities \ref{10.4}, or more precisely their 
difference from the corresponding 
$N$th approximants, in terms of the top order difference energies. In view of the 1st of the 
relations \ref{3.a21}, to control the departure of 
$\tchi$ from its initial values on $\Cb_0$ we must use the 2nd variation equation \ref{3.a44}, 
this equation being a propagation equation along the integral curves of $L$ 
(these intersect in the past $\Cb_0$).  
Similarly, in view of the 2nd of the relations \ref{3.a21}, to control the departure of $\tchib$ 
from its boundary values on ${\cal K}$ we must use the 2nd variation equation \ref{3.a45}, this 
equation being a propagation equation along the integral curves of $\Lb$  
(these intersect in the past ${\cal K}$). 
The principal term on the right hand side of equation \ref{3.a44} is the order 2 term 
$E\sm$. Similarly, the principal term on the right hand side of equation \ref{3.a45} is the 
order 2 term $E\smb$. To estimate $E^l\tchi$, $E^l\tchib$ requires that we commute equations 
\ref{3.a44}, \ref{3.a45} respectively with $E^l$. The principal terms on the right will then 
be $E^{l+1}\sm$, $E^{l+1}\smb$ respectively, which are of order $l+2$. At the top order 
we have $l=n$, so these principal terms, being then of order $n+2$, cannot be controlled in terms of 
the top order energies. The only way in which this difficulty can be overcome is if we can show that 
there are suitable 1st order quantities $v$ and $\underline{v}$ such that the differences 
\begin{equation}
E\sm-Lv \ \ \mbox{and} \ \ E\smb-\Lb\underline{v}
\label{10.5}
\end{equation}
are quantities of 1st order. The terms $Lv$ and $\Lb\underline{v}$ can then be brought to 
the left hand sides of equations \ref{3.a44} and \ref{3.a45} obtaining in this way propagation 
equations for the quantities $\sk-v$ and $\skb-\underline{v}$ along the integral curves of $L$ 
and $\Lb$ respectively, with right hand sides which are only of order 1. 

According to \ref{3.68} (see also \ref{3.a15}) we have:
\begin{equation}
\sm=-\beta_N\sbeta LH+\frac{1}{2}\rho\beta_N^2 EH-H\sbeta s_{NL}
\label{10.6}
\end{equation}
Note that the principal factors of the 1st and 3rd terms are both $L$ derivatives. Applying $E$ we 
then obtain, to principal terms, 
\begin{equation}
[E\sm]_{P.P}=\left[L(-\beta_N\sbeta EH-H\sbeta\ss_N)+\frac{1}{2}\rho\beta_N^2 E^2 H\right]_{P.P.}
\label{10.7}
\end{equation}
Now, suppose that: 
\begin{equation}
\left[\square_h H\right]_{P.P.}=0
\label{10.8}
\end{equation}
Then by \ref{9.205} with $H$ in the role of $f$ we can write:
\begin{equation}
[a E^2 H]_{P.P.}=[L\Lb H]_{P.P.}
\label{10.9}
\end{equation}
Since $a=\rho\lambda$ (see \ref{2.72}, \ref{2.74}) this implies:
\begin{equation}
\left[\frac{1}{2}\rho\beta_N^2 E^2 H\right]_{P.P.}=
\left[L\left(\frac{1}{2\lambda}\beta_N^2\Lb H\right)\right]_{P.P.}
\label{10.10}
\end{equation}
We can then set:
\begin{equation}
v=\frac{1}{2\lambda}\beta_N^2\Lb H-\beta_N\sbeta EH-H\sbeta\ss_N
\label{10.11}
\end{equation}
to achieve our objective in the case of $E\sm$. By the 1st of \ref{3.a21} we then have:
\begin{equation}
\sk-v=\tchi-\frac{1}{2}\lambda\beta_N^2\Lb H+\beta_N\sbeta EH
\label{10.12}
\end{equation}
However the function $v$ is singular at $\partial_-{\cal B}$ where $\lambda$ vanishes. To obtain a 
regular function we multiply by $\lambda$, defining instead:
\begin{equation}
\theta=\lambda(\sk-v)=\lambda\tchi+f
\label{10.13}
\end{equation}
where 
\begin{equation}
f=-\frac{1}{2}\beta_N^2\Lb H+\lambda\beta_N\sbeta EH
\label{10.14}
\end{equation}
As a consequence $\theta$ satisfies a propagation equation along the integral curves of $L$ of 
the form:
\begin{equation}
L\theta=R
\label{10.15}
\end{equation}
where $R$ is a quantity of order 1. 

In a similar manner we deduce using \ref{3.69}, that is 
\begin{equation}
\smb=-\beta_{\Nb}\sbeta\Lb H+\frac{1}{2}\rhob\beta_{\Nb}^2-H\sbeta s_{\Nb\Lb}
\label{10.16}
\end{equation}
that, provided that \ref{10.8} holds, the function 
\begin{equation}
\underline{v}=\frac{1}{2\lambdab}\beta_{\Nb}^2 LH-\beta_{\Nb}\sbeta EH-H\sbeta\ss_{\Nb}
\label{10.17}
\end{equation}
which is conjugate to $v$, achieves our objective in the case of $E\smb$. By the 2nd of \ref{3.a21} we then have:
\begin{equation}
\skb-\underline{v}=\tchib-\frac{1}{2}\lambdab\beta_{\Nb}^2 LH+\beta_{\Nb}\sbeta EH
\label{10.18}
\end{equation}
However the function $\underline{v}$ is singular at $\Cb_0$ where $\lambdab$ vanishes. To obtain a 
regular function we multiply by $\lambdab$, defining instead:
\begin{equation}
\thetab=\lambdab(\skb-\underline{v})=\lambdab\tchib+\fb
\label{10.19}
\end{equation}
where 
\begin{equation}
\fb=-\frac{1}{2}\beta_{\Nb}^2 LH+\lambdab\beta_{\Nb}\sbeta EH
\label{10.20}
\end{equation}
As a consequence $\thetab$ satisfies a propagation equation along the integral curves of $\Lb$ of 
the form:
\begin{equation}
\Lb\thetab=\Rb
\label{10.21}
\end{equation}
where $\Rb$ is a quantity of order 1. 

It remains for us to establish \ref{10.8}. In fact, more precisely, we have the following lemma. 

\vspace{2.5mm}

\noindent{\bf Lemma 10.1} \ \ We have: 
$$\square_h H=\left(H^{\prime\prime}-(1/2)\Omega^{-1}\Omega^\prime H^\prime\right)h^{-1}(d\sigma,d\sigma)
-2H^\prime(g^{-1})^{\mu\nu}h^{-1}(d\beta_\mu,d\beta_\nu)$$
where $H$ and $\Omega$ being a functions of $\sigma$ we denote
$$H^\prime=\frac{dH}{d\sigma}, \ H^{\prime\prime}=\frac{d^2 H}{d\sigma^2} \ \ \mbox{and} \ \ 
\Omega^\prime=\frac{d\Omega}{d\sigma}$$

\vspace{2.5mm}

\noindent{\em Proof:} Since $H$ is a function of $\sigma$ we have:
$$dH=H^\prime d\sigma, \ \ \ \square_{\tilde{h}}H=H^\prime\square_{\tilde{h}}\sigma
+H^{\prime\prime}\tilde{h}^{-1}(d\sigma,d\sigma)$$
Moreover, since 
$$\sigma=-(g^{-1})^{\mu\nu}\beta_\mu\beta_\nu$$
we have:
$$\square_{\tilde{h}}\sigma=-2(g^{-1})^{\mu\nu}\beta_\nu\square_{\tilde{h}}\beta_\mu 
-2(g^{-1})^{\mu\nu}\tilde{h}^{-1}(d\beta_\mu,d\beta_\nu)$$
Now the first term vanishes as $\square_{\tilde{h}}\beta_\mu=0$, therefore this simplifies to:
\begin{equation}
\square_{\tilde{h}}\sigma=-2(g^{-1})^{\mu\nu}\tilde{h}^{-1}(d\beta_\mu,d\beta_\nu)
\label{10.22}
\end{equation}
Substituting above we then obtain:
\begin{equation}
\square_{\tilde{h}}H=-2H^\prime(g^{-1})^{\mu\nu}\tilde{h}^{-1}(d\beta_\mu,d\beta_\nu)
+H^{\prime\prime}\tilde{h}^{-1}(d\sigma,d\sigma)
\label{10.23}
\end{equation}
On the other hand, by \ref{9.229}, taking $({\cal N},h)$ in the role of $({\cal M},g)$ and $H$ in the 
role of $f$, we have:
$$\Omega\square_{\tilde{h}}H=\square_h H+(1/2)\Omega^{-1}h^{-1}(d\Omega,dH)$$
Substituting \ref{10.23} and taking account of the fact that both $\Omega$ and $H$ being functions of 
$\sigma$ we have
$$h^{-1}(d\Omega,dH)=\Omega^\prime H^\prime h^{-1}(d\sigma,d\sigma)$$
yields the lemma. 

\vspace{2.5mm}

We shall presently discern the principal acoustical part of the order 1 quantity $R$ (see \ref{10.15}). 
We have (see \ref{10.13}): 
\begin{equation}
R=\lambda L\sk-L(\lambda v)+(L\lambda)\sk
\label{10.24}
\end{equation}
The principal term in $L\sk$ being, according to \ref{3.a44}, $E\sm$, let us first consider 
$\lambda E\sm-L(\lambda v)$. From \ref{10.6} and \ref{10.11}, using the first of the commutation 
relations \ref{3.a14} and the fact that by \ref{9.205}:
\begin{equation}
L\Lb H-aE^2 H=-a\square_h H-(1/2)(\chib LH+\chi\Lb H)+2\etab EH
\label{10.25}
\end{equation}
we deduce:
\begin{eqnarray}
&&\lambda E\sm-L(\lambda v)=-\lambda\chi\sbeta(\beta_N EH+H\ss_N) \label{10.26}\\
&&\hspace{20mm}+\frac{1}{2}\beta_N^2\left[a\square_h H+\frac{1}{2}(\chib LH+\chi\Lb H)-2\etab EH\right]
\nonumber\\
&&\hspace{20mm}+\lambda\left[-(LH)E(\beta_N\sbeta)+\frac{1}{2}(EH)E(\rho\beta_N^2)
-(L\beta_\mu)E(H\sbeta N^\mu)\right]\nonumber\\
&&\hspace{20mm}-\frac{1}{2}(\Lb H)L(\beta_N^2)+(EH)L(\lambda\beta_N\sbeta)
+(E\beta_\mu)L(\lambda H\sbeta N^\mu)\nonumber
\end{eqnarray}
This being of order 1, its principal acoustical (P.A.) part means the part containing order 1 
acoustical terms. According to Lemma 10.1:
\begin{eqnarray}
&&a\square_h H=(H^{\prime\prime}-(1/2)\Omega^{-1}\Omega^\prime H^\prime)(-(L\sigma)\Lb\sigma+a(E\sigma)E\sigma)\nonumber\\
&&\hspace{10mm}-2H^\prime(g^{-1})^{\mu\nu}(-(L\beta_\mu)\Lb\beta_\nu+a(E\beta_\mu)E\beta_\nu)
\label{10.27}
\end{eqnarray}
so this is of order 1 with vanishing P.A. part. By \ref{3.a15} and \ref{3.57}, \ref{10.6} 
and the 1st of \ref{3.47},
\begin{equation}
[LN^\mu]_{P.A.}=0 \ \ \mbox{hence} \ \ [L\beta_N]_{P.A.}=0
\label{10.28}
\end{equation}
By \ref{3.a33}, \ref{3.a35} - \ref{3.a37}, 
\begin{equation}
[LE^\mu]_{P.A.}=c^{-1}(E\lambdab+\pi\lambdab\tchib)N^\mu \ \ \mbox{hence} \ \ 
[L\sbeta]_{P.A.}=c^{-1}(E\lambdab+\pi\lambdab\tchib)\beta_N
\label{10.29}
\end{equation}
By \ref{3.46}, \ref{3.53}, \ref{3.54}, and \ref{3.a21}, we have:
\begin{equation}
[k]_{P.A.}=-\pi[\sk]_{P.A.}=-\pi\tchi, \ \ \ [\okb]_{P.A.}=-\pi[\skb]_{P.A.}=-\pi\tchib
\label{10.30}
\end{equation}
By \ref{3.a15} and the 1st of \ref{10.30},
\begin{equation}
[EN^\mu]_{P.A.}=\tchi(E^\mu-\pi N^\mu) \ \ \mbox{hence} \ \ [E\beta_N]_{P.A.}=\tchi(\sbeta-\pi\beta_N)
\label{10.31}
\end{equation}
Also, by \ref{3.a31}, \ref{3.a32},
\begin{equation}
[EE^\mu]_{P.A.}=\frac{1}{2c}(\tchib N^\mu+\tchi\Nb^\mu) \ \ \mbox{hence} \ \ 
[E\sbeta]_{P.A.}=\frac{1}{2c}(\tchib\beta_N+\tchi\beta_{\Nb})
\label{10.32}
\end{equation}
Since by the 2nd of \ref{3.a22} and of \ref{3.a25}
$$[\etab]_{P.A.}=\rhob[\tetab]_{P.A.}=\rhob(E\lambdab-\lambdab[\okb]_{P.A.})$$
substituting from the 2nd of \ref{10.30} we obtain:
\begin{equation}
[\etab]_{P.A.}=\rhob E\lambdab+a\pi\tchib
\label{10.33}
\end{equation}
Moreover, by \ref{3.a23} and \ref{10.30}:
\begin{equation}
[Ec]_{P.A.}=c[k+\okb]_{P.A.}=-c\pi(\tchi+\tchib)
\label{10.34}
\end{equation}
hence:
\begin{equation}
[E\rho]_{P.A.}=[E(c^{-1}\lambdab]_{P.A.}=c^{-1}E\lambdab+\pi\rho(\tchi+\tchib)
\label{10.35}
\end{equation}
Using the above results we deduce from \ref{10.26}:
\begin{eqnarray}
&&[\lambda E\sm-L(\lambda v)]_{P.A.}=\left(\frac{1}{2}\beta_N^2 EH+H\beta_N\ss_N\right)\rhob E\lambdab
\label{10.36}\\
&&\hspace{20mm}+\left\{\frac{1}{4}\rho\beta_N^2\Lb H-\rhob\left(\frac{1}{2}\beta_N\beta_{\Nb}
+c\sbeta(\sbeta-\pi\beta_N)\right)LH-\frac{1}{2}a\pi\beta_N^2 EH\right.\nonumber\\
&&\hspace{35mm}\left.+H\left[-2a\sbeta\ss_N+\rhob\left(c\pi\sbeta-\frac{1}{2}\beta_{\Nb}\right)s_{NL}
\right]\right\}\tchi\nonumber\\
&&\hspace{20mm}+\left\{-\frac{1}{4}\rhob\beta_N^2 LH+\frac{1}{2}a\pi\beta_N^2 EH
+H\beta_N\left(a\pi\ss_N-\frac{1}{2}\rhob s_{NL}\right)\right\}\tchib\nonumber
\end{eqnarray}
Going back to \ref{10.24}, we now consider the contributions of the remaining terms on the right in 
\ref{3.a44} to $[R]_{P.A.}$. In regard to the 2nd term we have, according to \ref{3.98} in the case 
of 2 spatial dimensions, 
\begin{equation}
\sgamma=\frac{1}{2}\sbeta^2 EH+H\sbeta\sss \ \ \mbox{hence} \ \ [\sgamma]_{P.A.}=0
\label{10.37}
\end{equation}
The results above allow us to determine the contributions of the remaining terms. Combining 
with \ref{10.36} and recalling Proposition 3.3 we arrive at, simply:
\begin{equation}
[R]_{P.A.}=-a\tchi^2+\left(2q\lambdab+(2p-\sbeta^2 LH)\lambda\right)\tchi
\label{10.38}
\end{equation}
Denoting by $\tilde{R}$ the remainder:
\begin{equation}
\tilde{R}=R-[R]_{P.A.}
\label{10.39}
\end{equation}
$\tilde{R}$ is of order 1 with vanishing principal acoustical part. A straightforward calculation shows 
that $\tilde{R}$ is of the form:
\begin{eqnarray}
&&\tilde{R}=c^{\mu\nu}_{L\Lb}(L\beta_\mu)\Lb\beta_\nu+\lambda c^{\mu\nu}_{LE}(L\beta_\mu)E\beta_\nu 
\nonumber\\
&&\hspace{7mm}+\lambdab c^{\mu\nu}_{\Lb E}(\Lb\beta_\mu)E\beta_\nu
+\lambda\lambdab c^{\mu\nu}_{EE}(E\beta_\mu)E\beta_\nu
\label{10.40}
\end{eqnarray}
where the coefficients are regular and of order 0. 

Similarly we arrive at analogous results for the conjugate function $\Rb$, the right hand side of 
\ref{10.21}:
\begin{equation}
\Rb=[\Rb]_{P.A.}+\tilde{\Rb}
\label{10.41}
\end{equation}
where:
\begin{equation}
[\Rb]_{P.A.}=-a\tchib^2+\left(2\qb\lambda+(2\pb-\sbeta^2 \Lb H)\lambdab\right)\tchib
\label{10.42}
\end{equation}
and $\tilde{\Rb}$ is of the form:
\begin{eqnarray}
&&\tilde{\Rb}=\underline{c}^{\mu\nu}_{\Lb L}(\Lb\beta_\mu)L\beta_\nu+\lambdab \underline{c}^{\mu\nu}_{\Lb E}(\Lb\beta_\mu)E\beta_\nu 
\nonumber\\
&&\hspace{7mm}+\lambda\underline{c}^{\mu\nu}_{LE}(L\beta_\mu)E\beta_\nu
+\lambdab\lambda\underline{c}^{\mu\nu}_{EE}(E\beta_\mu)E\beta_\nu
\label{10.43}
\end{eqnarray}
where the coefficients are regular and of order 0. 

We summarize the above in the following proposition. 

\vspace{2.5mm}

\noindent{\bf Proposition 10.1} \ \ The quantities: 
$$\theta=\lambda\tchi+f, \ \ \ \thetab=\lambdab\tchib+\fb$$
where 
$$f=-\frac{1}{2}\beta_N^2\Lb H+\lambda\beta_N\sbeta EH, \ \ \ 
\fb=-\frac{1}{2}\beta_{\Nb}^2 LH+\lambdab\beta_{\Nb}\sbeta EH$$
satisfy the following regularized propagation equations:
$$L\theta=R, \ \ \ \Lb\thetab=\Rb$$
where 
$$R=[R]_{P.A.}+\tilde{R}, \ \ \ \Rb=[\Rb]_{P.A.}+\tilde{\Rb}$$
are of order 1, their principal acoustical parts being given by
\begin{eqnarray*}
&&[R]_{P.A.}=-a\tchi^2+\left(2L\lambda-\sbeta^2(LH)\lambda\right)\tchi\\
&&[\Rb]_{P.A.}=-a\tchib^2+\left(2\Lb\lambdab-\sbeta^2(\Lb H)\lambdab\right)\tchib
\end{eqnarray*}
and the remainders $\tilde{R}$, $\tilde{\Rb}$ are of the form \ref{10.40}, \ref{10.43} respectively. 

\vspace{2.5mm}

Note that the only 1st order acoustical quantity appearing in $R$ is $\tchi$, that is $\theta$ itself, 
and the only 1st order acoustical quantity appearing in $\Rb$ is $\tchib$, that is $\thetab$ itself. 

\section{Regularization of the Propagation Equations for $E^2\lambda$ and $E^2\lambdab$}

The quantities $\lambda$ and $\lambdab$ satisfy the propagation equations of Proposition 3.3. 
The equation for $\lambda$ being a propagation equation along the integral curves of $L$ controls 
the departure of $\lambda$ from its initial values on $\Cb_0$. The equation for $\lambdab$ being 
a propagation equation along the integral curves of $\Lb$ controls the departure of $\lambdab$ from 
its boundary values on ${\cal K}$. While the quantities $\lambda$ and $\lambdab$ are of order 0, 
the right hand sides of these propagation equations are of order 1, and it is not possible to 
regularize these equations in the way that we regularized the propagation equations for $\tchi$ 
and $\tchib$. However, seeing that what must be estimated at the top order are the quantities 
\ref{10.4}, what is required is to regularize the propagation equations for $E^2\lambda$ and 
$E^2\lambdab$. These follow by applying $E^2$ to the propagation equations of Proposition 3.3 
and using the first two of the commutation relations \ref{3.a14}, which imply:
\begin{eqnarray*}
&&[L,E^2]=[L,E]E+E[L,E]=-\chi E^2-E\chi E \\
&&[\Lb,E^2]=[\Lb,E]E+E[\Lb E]=-\chib E^2-E\chib E
\end{eqnarray*}
We thus obtain the propagation equations:
\begin{eqnarray}
&&L(E^2\lambda)=-2\chi E^2\lambda-(E\chi)E\lambda+p E^2\lambda+q E^2\lambdab\nonumber\\
&&\hspace{15mm}+\lambda E^2 p+\lambdab E^2 q +2(Ep)E\lambda+2(Eq)E\lambdab
\label{10.44}
\end{eqnarray}
\begin{eqnarray}
&&\Lb(E^2\lambdab)=-2\chib E^2\lambdab-(E\chib)E\lambdab+\pb E^2\lambdab+\qb E^2\lambda\nonumber\\
&&\hspace{15mm}+\lambdab E^2\pb+\lambda E^2\qb +2(E\pb)E\lambdab+2(E\qb)E\lambda
\label{10.45}
\end{eqnarray}
The principal terms on the right in \ref{10.44} are the terms $\lambda E^2 p$, $\lambdab E^2 q$. 
The principal terms on the right in \ref{10.45} are the terms $\lambdab E^2\pb$, $\lambda E^2\qb$. 
These terms are all of order 3. 

Consider first equation \ref{10.44}. The function $p$ is given in the statement of Proposition 3.3:
\begin{equation}
p=m-\frac{1}{2c}(\beta_N\beta_{\Nb}LH+H\beta_{\Nb}s_{NL})
\label{10.46}
\end{equation}
$m$ being given in terms of $\sm$ and $\om$ by the 1st of \ref{3.47}:
\begin{equation}
m=-\pi\sm-\om 
\label{10.47}
\end{equation}
and $\om$ is given by \ref{3.57}:
\begin{equation}
\om=\frac{1}{4c}(\beta_N^2 LH+2H\beta_N s_{NL})
\label{10.48}
\end{equation}
while $\sm$ is given by \ref{10.6}. Thus $p$ is given by:
\begin{equation}
p=-\pi\sm-\frac{1}{4c}\beta_N(\beta_N+2\beta_{\Nb})LH-\frac{1}{2c}H(\beta_N+\beta_{\Nb})s_{NL}
\label{10.49}
\end{equation}
Applying $E^2$ we then obtain, to principal terms, 
\begin{eqnarray}
&&[E^2 p]_{P.P.}=-\pi[E^2\sm]_{P.P.} \label{10.50}\\
&&\hspace{10mm}-\left[L\left(\frac{1}{4c}\beta_N(\beta_N+2\beta_{\Nb})E^2 H
+\frac{1}{2c}H(\beta_N+\beta_{\Nb})N^\mu E^2\beta_\mu\right)\right]_{P.P.} \nonumber
\end{eqnarray}
Now, we have shown in the previous section that with $v$ defined by \ref{10.11} 
$\lambda E\sm-L(\lambda v)$ is a quantity of 1st order. It follows that $\lambda E^2\sm-LE(\lambda v)$ is a quantity of 2nd order. 
We can therefore write:
\begin{equation}
[\lambda E^2 p]_{P.P.}=[Lw_1]_{P.P.}
\label{10.51}
\end{equation}
Here $w_1$ is the function:
\begin{equation}
w_1=-\pi [E(\lambda v)]_{P.P}-\lambda\left(\frac{1}{4c}\beta_N(\beta_N+2\beta_{\Nb})E^2 H
+\frac{1}{2c}H(\beta_N+\beta_{\Nb})N^\mu E^2\beta_\mu\right)
\label{10.52}
\end{equation}
where by \ref{10.11}:
\begin{equation}
[E(\lambda v)]_{P.P.}=\frac{1}{2}\beta_N^2 E\Lb H-\lambda(\beta_N\sbeta E^2 H+H\sbeta N^\mu E^2\beta_\mu)
\label{10.53}
\end{equation}
Next, the function $q$ is also given in the statement of Proposition 3.3:
\begin{equation}
q=\frac{1}{4c}\beta_N^2\Lb H
\label{10.54}
\end{equation}
Applying $E^2$ we then obtain, to principal terms, 
\begin{equation}
[E^2 q]_{P.P.}=\left[\frac{1}{4c}\beta_N^2 E^2\Lb H\right]_{P.P.}
\label{10.55}
\end{equation}
Multiplying then by $\lambda\lambdab=ca$ we have:
\begin{equation}
[\lambda\lambdab E^2 q]_{P.P.}=\left[\frac{1}{4}\beta_N^2 aE^2\Lb H\right]_{P.P.}
=\left[\frac{1}{4}\beta_N^2 \Lb(aE^2 H)\right]_{P.P.}
\label{10.56}
\end{equation}
Now by \ref{10.25} and \ref{10.27} (Lemma 10.1) $aE^2 H-L\Lb H$ is of order 1. We can therefore 
write:
\begin{equation}
[\lambda\lambdab E^2 q]_{P.P.}=[Lw_2]_{P.P.}
\label{10.57}
\end{equation}
where $w_2$ is the function:
\begin{equation}
w_2=\frac{1}{4}\beta_N^2\Lb^2 H
\label{10.58}
\end{equation}
Let us define the function:
\begin{equation}
j=\lambda w_1+w_2
\label{10.59}
\end{equation}
Then according to the above we have:
\begin{equation}
[\lambda(\lambda E^2 p+\lambdab E^2 q)]_{P.P.}=[Lj]_{P.P.}
\label{10.60}
\end{equation}
In regard then to equation \ref{10.44}, the quantity:
\begin{equation}
\nu=\lambda E^2\lambda-j
\label{10.61}
\end{equation}
satisfies a propagation equation along the integral curves of $L$ of the form:
\begin{equation}
L\nu=K
\label{10.62}
\end{equation}
where $K$ is a quantity of order 2. Note that according to \ref{10.59} and \ref{10.52}, \ref{10.53}, 
\ref{10.58}, the function $j$ is given by:
\begin{eqnarray}
&&j=\frac{1}{4}\beta_N^2\Lb^2 H-\frac{1}{2}\lambda\pi\beta_N^2 E\Lb H\nonumber\\
&&\hspace{5mm}+\lambda^2\left\{\beta_N\left[\pi\sbeta-\frac{1}{4c}(\beta_N+2\beta_{\Nb})\right]E^2 H
\right.\nonumber\\
&&\hspace{15mm}\left.+H\left[\pi\sbeta-\frac{1}{2c}(\beta_N+\beta_{\Nb})\right]N^\mu E^2\beta_\mu
\right\}\label{10.63} 
\end{eqnarray}
In regard to equation \ref{10.45} we deduce in a similar manner that the quantity:
\begin{equation}
\nub=\lambdab E^2\lambdab-\jb
\label{10.64}
\end{equation}
with $\jb$ given by:
\begin{eqnarray}
&&\jb=\frac{1}{4}\beta_{\Nb}^2 L^2 H-\frac{1}{2}\lambdab\pi\beta_{\Nb}^2 ELH\nonumber\\
&&\hspace{5mm}+\lambdab^2\left\{\beta_{\Nb}\left[\pi\sbeta-\frac{1}{4c}(2\beta_N+\beta_{\Nb})\right]
E^2 H\right.\nonumber\\
&&\hspace{15mm}\left.+H\left[\pi\sbeta-\frac{1}{2c}(\beta_N+\beta_{\Nb})\right]\Nb^\mu E^2\beta_\mu
\right\}\label{10.65} 
\end{eqnarray}
satisfies a propagation equation along the integral curves of $\Lb$ of the form:
\begin{equation}
\Lb\nub=\Kb
\label{10.66}
\end{equation}
where $\Kb$ is a quantity of order 2. 

We shall presently discern the principal acoustical part of the order 2 quantity $K$, that is the part containing order 2 acoustical terms. Let us define, 
in connection with \ref{10.52}:
\begin{equation}
j_0=-\lambda^2\left(\frac{\beta_N(\beta_N+2\beta_{\Nb})}{4c}E^2 H
+\frac{(\beta_N+\beta_{\Nb})}{2c}N^\mu E^2\beta_\mu\right)
\label{10.67}
\end{equation}
Consider first 
\begin{equation}
P_1:=\lambda^2 E^2 p+\pi\lambda^2 E^2\sm-Lj_0
\label{10.68}
\end{equation}
Since $p$ is given by \ref{10.49}, we have:
\begin{eqnarray}
&&[P_1]_{P.A.}=-\lambda^2\sm[E^2\pi]_{P.A.}-\lambda^2\frac{(LH)}{4}
\left[E^2\left(\frac{\beta_N(\beta_N+2\beta_{\Nb})}{c}\right)\right]_{P.A.} \nonumber\\
&&\hspace{17mm}+\frac{\lambda^2}{4c}\beta_N(\beta_N+2\beta_{\Nb})\left[[L,E^2]H\right]_{P.A.} \nonumber\\
&&\hspace{17mm}-\frac{\lambda^2 H}{2}\left[E^2\left(\frac{(\beta_N+\beta_{\Nb})}{c}N^\mu\right)
\right]_{P.A.}L\beta_\mu \nonumber\\
&&\hspace{17mm}+\frac{\lambda^2 H}{2c}(\beta_N+\beta_{\Nb})N^\mu\left[[L,E^2]\beta_\mu\right]_{P.A.}
\label{10.69}
\end{eqnarray}
By \ref{3.105} and \ref{3.a27} with $x^0=t$ in the role of $f$ we have: 
\begin{equation}
[E^2\pi]_{P.A.}=\frac{1}{2c}(E\tchib+E\tchi)
\label{10.70}
\end{equation}
By \ref{10.31} and its conjugate we have:
\begin{equation}
[E^2\beta_N]_{P.A.}=(\sbeta-\pi\beta_N)E\tchi, \ \ \ 
[E^2\beta_{\Nb}]_{P.A.}=(\sbeta-\pi\beta_{\Nb})E\tchib
\label{10.71}
\end{equation}
and:
\begin{equation}
[E^2 N^\mu]_{P.A.}=(E^\mu-\pi N^\mu)E\tchi
\label{10.72}
\end{equation}
Also, by \ref{10.34}, 
\begin{equation}
[E^2 c]_{P.A.}=-c\pi(E\tchi+E\tchib)
\label{10.73}
\end{equation}
We then obtain:
\begin{eqnarray}
&&\left[E^2\left(\frac{\beta_N(\beta_N+2\beta_{\Nb})}{c}\right)\right]_{P.A.} \label{10.74}\\
&&\hspace{20mm}=\frac{1}{c}\left\{(2(\beta_N+\beta_{\Nb})\sbeta-\pi\beta_N^2)E\tchi
+(2\beta_N\sbeta+\pi\beta_N^2)E\tchib\right\}\nonumber
\end{eqnarray}
and:
\begin{eqnarray}
&&\left[E^2\left(\frac{(\beta_N+\beta_{\Nb})}{c}N^\mu\right)\right]_{P.A.}L\beta_\mu  
\nonumber\\
&&\hspace{20mm}=\frac{1}{c}\left\{\left[(\sbeta+\pi\beta_{\Nb})s_{NL}+(\beta_N+\beta_{\Nb})(\rho\ss_N-\pi s_{NL})
\right]E\tchi\right.\nonumber\\
&&\hspace{30mm}\left.+(\sbeta+\pi\beta_N)s_{NL}E\tchib\right\}
\label{10.75}
\end{eqnarray}
The first of the commutation relations \ref{3.a14} implies:
\begin{equation}
\left[[L,E^2]H\right]_{P.A.}=-\rho(EH)E\tchi
\label{10.76}
\end{equation}
and:
\begin{equation}
N^\mu\left[[L,E^2]\beta_\mu\right]_{P.A.}=-\rho\ss_N E\tchi
\label{10.77}
\end{equation}
Using the above results and recalling also the expression \ref{10.6} for $\sm$ we obtain 
the following:
\begin{eqnarray}
&&[P_1]_{P.A.}=\lambda^2\left\{\left[-\frac{\pi\beta_N^2}{4c}LH-\frac{\beta_N^2}{4c}\rho EH
-\frac{\pi\beta_N}{2c}Hs_{NL}\right]E\tchib\right.\nonumber\\
&&\hspace{23mm}+\left[\frac{1}{4c}(-2\sbeta\beta_{\Nb}+\pi\beta_N^2)LH
-\frac{\beta_N}{2c}(\beta_N+\beta_{\Nb})\rho EH\right.\nonumber\\
&&\hspace{28mm}\left.\left.+\frac{\pi\beta_N}{2c}H s_{NL}-\frac{(\beta_N+\beta_{\Nb})}{c} H\rho\ss_N
\right]E\tchi\right\}
\label{10.78}
\end{eqnarray}
Next, let us define:
\begin{equation}
j_1=\pi\lambda\left(-\frac{1}{2}\beta_N^2 E\Lb H+\lambda\sbeta(\beta_N E^2 H+H N^\mu E^2\beta_\mu)\right)
\label{10.79}
\end{equation}
We then have (see \ref{10.52}, \ref{10.53}):
\begin{equation}
j_0+j_1=\lambda w_1 
\label{10.80}
\end{equation}
Consider: 
\begin{equation}
P_2:=-\pi\lambda\left(\lambda^2 E^2\sm-LE(\lambda v)\right)
\label{10.81}
\end{equation}
In view of \ref{10.68} we then have:
\begin{eqnarray}
&&P_1+P_2=\lambda^2 E^2 p-L(\lambda w_1)\nonumber\\
&&\hspace{15mm}+L(\pi\lambda E(\lambda v)+j_1)-L(\pi\lambda)E(\lambda v)
\label{10.82}
\end{eqnarray}
Now, by \ref{10.11} and \ref{10.79}:
\begin{eqnarray}
&&\pi\lambda E(\lambda v)+j_1=\pi\lambda\beta_N(E\beta_N)\Lb H
-\pi\lambda\sbeta(\beta_N EH+H\ss_N)E\lambda\nonumber\\
&&\hspace{24mm}-\pi\lambda^2\left[(\beta_N E\sbeta+\sbeta E\beta_N)EH+\sbeta\ss_N EH\right.\nonumber\\
&&\hspace{34mm}\left.+H\ss_NE\sbeta
+H\sbeta(EN^\mu)E\beta_\mu\right]
\label{10.83}
\end{eqnarray}
This is a 1st order quantity. We must find the principal acoustical part of its $L$ derivative, 
a 2nd order quantity. By the propagation equation for $\lambda$ of Proposition 3.3 $LE\lambda$ has 
vanishing principal acoustical part. 
Also, by \ref{10.31} $EN^\mu$ and $E\beta_N$ only involve the 1st order acoustical 
quantity $\tchi$ and by the 2nd variation equation \ref{3.a44} $L\tchi$ has vanishing principal 
acoustical part. It then follows that:
\begin{equation}
\left[L(\pi\lambda E(\lambda v)+j_1)\right]_{P.A.}=-\pi\lambda^2(\beta_N EH+H\ss_N)[LE\sbeta]_{P.A.}
\label{10.84}
\end{equation}
By \ref{10.32} and the fact that $L\tchi$ has vanishing P.A. part, 
\begin{equation}
[LE\sbeta]_{P.A.}=\frac{\beta_N}{2c}[L\tchib]_{P.A.}=\frac{\beta_N}{c}(E^2\lambdab+\pi\lambdab E\tchib)
\label{10.85}
\end{equation}
the last by \ref{8.138}. We therefore have:
\begin{equation}
\left[L(\pi\lambda E(\lambda v)+j_1)\right]_{P.A.}=-\frac{\pi\lambda^2\beta_N}{c}(\beta_N EH+H\ss_N)
(E^2\lambdab+\pi\lambdab E\tchib)
\label{10.86}
\end{equation}
We then conclude through \ref{10.82} that:
\begin{eqnarray}
&&\left[\lambda^2 E^2 p-L(\lambda w_1)\right]_{P.A}=[P_1]_{P.A.}+[P_2]_{P.A.}\label{10.87}\\
&&\hspace{35mm}+\frac{\pi\lambda^2\beta_N}{c}(\beta_N EH+H\ss_N)
(E^2\lambdab+\pi\lambdab E\tchib)\nonumber
\end{eqnarray}
Since $E\sm$ has vanishing P.A. part, the P.A. part of $P_2$ is readily deduced from \ref{10.36}:
\begin{eqnarray}
&&[P_2]_{P.A}=-\pi\lambda\left\{\left(\frac{1}{2}\beta_N^2 EH+H\beta_N\ss_N\right)\rhob E^2\lambdab
\right. \label{10.88}\\
&&\hspace{15mm}+\left[-\frac{1}{4}\rhob\beta_N^2 LH+\frac{1}{2}a\pi\beta_N^2 EH+H\beta_N\left(a\pi\ss_N-\frac{1}{2}\rhob s_{NL}\right)\right]E\tchib\nonumber\\
&&\hspace{15mm}+\left[\frac{1}{4}\rho\beta_N^2\Lb H-\rhob\left(\frac{1}{2}\beta_N\beta_{\Nb}
+c\sbeta(\sbeta-\pi\beta_N)\right)LH\right.
\nonumber\\
&&\hspace{20mm}\left.\left.-\frac{1}{2}a\pi\beta_N^2 EH+H\left(-2a\sbeta\ss_N
+\rhob\left(c\pi\sbeta-\frac{1}{2}\beta_{\Nb}\right)s_{NL}\right)\right]E\tchi\right\} \nonumber
\end{eqnarray}
Consider finally (see \ref{10.58}):
\begin{equation}
Q:=\lambda\lambdab E^2 q-Lw_2
\label{10.89}
\end{equation}
By \ref{10.57} this is a quantity of order 2. We have:
\begin{eqnarray}
&&\left[\lambda\lambdab E^2q-\frac{1}{4}\beta_N^2 aE^2\Lb H\right]_{P.A.}=
\lambda\lambdab(\Lb H)\left[E^2\left(\frac{\beta_N^2}{4c}\right)\right]_{P.A.} \label{10.90}\\
&&\hspace{30mm}=\lambda\lambdab(\Lb H)\left\{\frac{\beta_N}{4c}(2\sbeta-\pi\beta_N)E\tchi
+\frac{\pi\beta_N^2}{4c}E\tchib\right\} \nonumber
\end{eqnarray}
the last by \ref{10.73} and the 1st of \ref{10.71}. Moreover, by the 2nd of the commutation relations 
\ref{3.a14}, 
\begin{equation}
\left[[aE^2,\Lb]H\right]_{P.A.}=a\rhob(EH)E\tchib
\label{10.91}
\end{equation}
On the other hand, applying $\Lb$ to \ref{10.25} we obtain:
\begin{equation}
\left[\Lb(aE^2 H)-L\Lb^2 H\right]_{P.A.}=\left[\frac{1}{2}((LH)\Lb\chib+(\Lb H)\Lb\chi)
-2(EH)\Lb\etab\right]_{P.A.}
\label{10.92}
\end{equation}
Now, by the 2nd variation equation \ref{3.a45} the P.A. part of $\Lb\tchib$ vanishes, while by 
\ref{8.137}:
\begin{equation}
[\Lb\tchi]_{P.A.}=2E^2\lambda+2\pi\lambda E\tchi
\label{10.93}
\end{equation} 
Also, by \ref{10.33} in view of the propagation equation for $\lambdab$ of Proposition 3.3 and 
the 2nd variation equation \ref{3.a45}, the P.A. part of $\Lb\etab$ vanishes. We then 
conclude through \ref{10.92} that:
\begin{equation}
\left[\Lb(aE^2 H)-L\Lb^2 H\right]_{P.A.}=\rho(\Lb H)(E^2\lambda+\pi\lambda E\tchi)
\label{10.94}
\end{equation}
Using the above results we deduce through \ref{10.89}:
\begin{eqnarray}
&&[Q]_{P.A.}=\frac{\lambdab\beta_N^2}{4c}(\Lb H)E^2\lambda+\left(\frac{a\pi\beta_N^2}{4}\Lb H+\frac{a\lambda\beta_N^2}{4c}EH\right)E\tchib
\nonumber\\
&&\hspace{15mm}+\frac{a\sbeta\beta_N}{2}(\Lb H)E\tchi
\label{10.95}
\end{eqnarray}
Now, by \ref{10.87} and \ref{10.89} we have, in view of \ref{10.59},
\begin{eqnarray}
&&\left[\lambda(\lambda E^2 p+\lambdab E^2 q)-Lj\right]_{P.A.}=[P_1]_{P.A.}+[P_2]_{P.A.}+[Q]_{P.A.}
\label{10.96}\\
&&\hspace{35mm}+\frac{\pi\lambda^2\beta_N}{c}(\beta_N EH+H\ss_N)
(E^2\lambdab+\pi\lambdab E\tchib)\nonumber
\end{eqnarray}
Substituting \ref{10.78}, \ref{10.88}, and \ref{10.95}, yields:
\begin{eqnarray}
&&\left[\lambda(\lambda E^2 p+\lambdab E^2 q)-Lj\right]_{P.A.}= \nonumber\\
&&\hspace{25mm}\frac{\lambda^2\pi}{2c}\beta_N^2(EH)E^2\lambdab
+\frac{\lambdab}{4c}\beta_N^2(\Lb H)E^2\lambda \nonumber\\
&&\hspace{25mm}+\frac{a\pi}{4}\beta_N^2(\Lb H+2\pi\lambda EH)E\tchib \label{10.97}\\
&&\hspace{25mm}+\lambda\left[\lambdab\frac{\beta_N}{2c}\left(\sbeta-\frac{\pi}{2}\beta_N\right)\Lb H 
\right.\nonumber\\
&&\hspace{30mm}+\lambda\left(\sbeta\left(\pi\sbeta-\pi^2\beta_N-\frac{\beta_{\Nb}}{2c}\right)
+\frac{\pi\beta_N}{4c}(\beta_N+2\beta_{\Nb})\right)LH\nonumber\\
&&\hspace{30mm}-\frac{a\beta_N}{2}\left(\frac{(\beta_N+\beta_{\Nb})}{c}-\pi^2\beta_N\right)EH\nonumber\\
&&\hspace{30mm}\left.+H\left(\frac{(\beta_N+\beta_{\Nb})}{2c}-\pi\sbeta\right)(-2a\ss_N+\pi\lambda s_{NL})\right]
E\tchi\nonumber
\end{eqnarray}
To this must be added the terms 
$$(L\lambda-2\chi\lambda)E^2\lambda-\lambda\rho(E\lambda)E\tchi+\lambda p E^2\lambda
+\lambda q E^2\lambdab$$
to obtain the principal acoustical part of $K$. This yields:
\begin{eqnarray}
&&[K]_{P.A.}=2(p\lambda+q\lambdab-\chi\lambda)E^2\lambda
+\frac{\lambda\beta_N^2}{4c}(\Lb H+2\pi\lambda EH)E^2\lambdab\nonumber\\
&&\hspace{17mm}+\lambda AE\tchi+a \oA E\tchib
\label{10.98}
\end{eqnarray}
where the coefficients $A$ and $\oA$ are given by:
\begin{eqnarray}
&&A=-\frac{\lambdab}{c}E\lambda+\frac{\lambdab\beta_N}{2c}\left(\sbeta-\frac{\pi}{2}\beta_N\right)\Lb H 
\nonumber\\
&&\hspace{10mm}+\lambda\left[\sbeta\left(\pi\sbeta-\pi^2\beta_N-\frac{\beta_{\Nb}}{2c}\right)
+\frac{\pi\beta_N}{4c}(\beta_N+2\beta_{\Nb})\right]LH\nonumber\\
&&\hspace{10mm}-\frac{a\beta_N}{2}\left(\frac{(\beta_N+\beta_{\Nb})}{c}-\pi^2\beta_N\right)EH\nonumber\\
&&\hspace{10mm}+H\left(\frac{(\beta_N+\beta_{\Nb})}{2c}-\pi\sbeta\right)(-2a\ss_N+\pi\lambda s_{NL})
\label{10.99}
\end{eqnarray}
and:
\begin{equation}
\oA=\frac{\pi\beta_N^2}{4}(\Lb H+2\pi\lambda EH)
\label{10.100}
\end{equation}
Denoting by $\tilde{K}$ the remainder:
\begin{equation}
\tilde{K}=K-[K]_{P.A.}
\label{10.101}
\end{equation}
$\tilde{K}$ is of order 2 with vanishing principal acoustical part. A straightforward calculation 
shows that the principal part of $\tilde{K}$ is of the form:
\begin{eqnarray}
&&[\tilde{K}]_{P.P.}=\left[\lambdab(\Lb\lambda)d^\mu_{1,EE}+\lambdab\lambda(E\lambda)d^\mu_{2,EE}
+\lambda^2(E\lambdab)d^\mu_{3,EE}\right.\nonumber\\
&&\hspace{25mm}+\lambda^2\lambdab\tchi d^\mu_{4,EE}+\lambda^2\lambdab\tchib d^\mu_{5,EE} \nonumber\\
&&\hspace{20mm}\left.+\lambda\lambdab d^{\mu\nu}_{6,EE}\Lb\beta_\nu
+\lambda^2 d^{\mu\nu}_{7,EE}L\beta_\nu+\lambda^2\lambdab d^{\mu\nu}_{8,EE}E\beta_\nu\right]E^2\beta_\mu \nonumber\\
&&\hspace{15mm}+\left[\lambda(E\lambda)d^\mu_{1,EL}+\lambda^2\tchi d^\mu_{2,EL}
+\lambda^2\tchib d^\mu_{3,EL}\right.\nonumber\\
&&\hspace{20mm}\left.+\lambda d^{\mu\nu}_{4,EL}\Lb\beta_\nu 
+\lambda^2 d^{\mu\nu}_{5,EL}E\beta_\nu\right]EL\beta_\mu \nonumber\\
&&\hspace{15mm}+\left[\lambda(E\lambdab)d^\mu_{1,E\Lb}+\lambdab(E\lambda)d^\mu_{2,E\Lb}
+\lambda\lambdab\tchi d^\mu_{3,E\Lb}+\lambda\lambdab\tchib d^\mu_{4,E\Lb}\right. \nonumber\\
&&\hspace{20mm}\left.+\lambdab d^{\mu\nu}_{5,E\Lb}\Lb\beta_\nu+\lambda d^{\mu\nu}_{6,E\Lb}L\beta_\nu
+\lambda\lambdab d^{\mu\nu}_{7,E\Lb}E\beta_\nu\right]E\Lb\beta_\mu \nonumber\\
&&\hspace{15mm}+\left[d^{\mu\nu}_{1,\Lb\Lb}L\beta_\nu
+\lambdab d^{\mu\nu}_{2,\Lb\Lb}E\beta_\nu\right]\Lb^2\beta_\mu 
\label{10.102}
\end{eqnarray}
where the coefficients are regular and of order 0. 

Similarly, we arrive at analogous results for the conjugate function $\Kb$, the right hand side of 
\ref{10.66}:
\begin{equation}
\Kb=[\Kb]_{P.A.}+\tilde{\Kb}
\label{10.103}
\end{equation}
Here:
\begin{eqnarray}
&&[\Kb]_{P.A.}=2(\pb\lambdab+\qb\lambda-\chib\lambdab)E^2\lambdab
+\frac{\lambdab\beta_{\Nb}^2}{4c}(LH+2\pi\lambdab EH)E^2\lambda\nonumber\\
&&\hspace{17mm}+\lambdab \oAb E\tchib+a \Ab E\tchi
\label{10.104}
\end{eqnarray}
where the coefficients $\oAb$ and $\Ab$ are given by:
\begin{eqnarray}
&&\oAb=-\frac{\lambda}{c}E\lambdab+\frac{\lambda\beta_{\Nb}}{2c}\left(\sbeta-\frac{\pi}{2}\beta_{\Nb}\right)LH 
\nonumber\\
&&\hspace{10mm}+\lambdab\left[\sbeta\left(\pi\sbeta-\pi^2\beta_{\Nb}-\frac{\beta_N}{2c}\right)
+\frac{\pi\beta_{\Nb}}{4c}(\beta_{\Nb}+2\beta_N)\right]\Lb H\nonumber\\
&&\hspace{10mm}-\frac{a\beta_{\Nb}}{2}\left(\frac{(\beta_N+\beta_{\Nb})}{c}-\pi^2\beta_{\Nb}\right)EH\nonumber\\
&&\hspace{10mm}+H\left(\frac{(\beta_N+\beta_{\Nb})}{2c}-\pi\sbeta\right)(-2a\ss_{\Nb}+\pi\lambdab  s_{\Nb\Lb})
\label{10.105}
\end{eqnarray}
and:
\begin{equation}
\Ab=\frac{\pi\beta_{\Nb}^2}{4}(LH+2\pi\lambdab EH)
\label{10.106}
\end{equation}
Also, the principal part of $\tilde{\Kb}$ is of the form: 
\begin{eqnarray}
&&[\tilde{\Kb}]_{P.P.}=\left[\lambda(L\lambdab)\db^\mu_{1,EE}+\lambda\lambdab(E\lambdab)\db^\mu_{2,EE}
+\lambdab^2(E\lambda)\db^\mu_{3,EE}\right.\nonumber\\
&&\hspace{25mm}+\lambdab^2\lambda\tchib\db^\mu_{4,EE}+\lambdab^2\lambda\tchi\db^\mu_{5,EE} \nonumber\\
&&\hspace{20mm}\left.+\lambdab\lambda \db^{\mu\nu}_{6,EE}L\beta_\nu
+\lambdab^2 \db^{\mu\nu}_{7,EE}\Lb\beta_\nu+\lambdab^2\lambda \db^{\mu\nu}_{8,EE}E\beta_\nu\right]E^2\beta_\mu \nonumber\\
&&\hspace{15mm}+\left[\lambdab(E\lambdab)\db^\mu_{1,E\Lb}
+\lambdab^2\tchib\db^\mu_{2,E\Lb}+\lambdab^2\tchi\db^\mu_{3,E\Lb}\right.\nonumber\\
&&\hspace{20mm}\left.+\lambdab \db^{\mu\nu}_{4,E\Lb}L\beta_\nu
+\lambdab^2 \db^{\mu\nu}_{5,E\Lb}E\beta_\nu\right]E\Lb\beta_\mu \nonumber\\
&&\hspace{15mm}+\left[\lambdab(E\lambda)\db^\mu_{1,EL}+\lambda(E\lambdab)\db^\mu_{2,EL}
+\lambdab\lambda\tchib\db^\mu_{3,EL}+\lambdab\lambda\tchi\db^\mu_{4,EL}\right. \nonumber\\
&&\hspace{20mm}\left.+\lambda \db^{\mu\nu}_{5,EL}L\beta_\nu+\lambdab \db^{\mu\nu}_{6,EL}\Lb\beta_\nu
+\lambdab\lambda \db^{\mu\nu}_{5,EL}E\beta_\nu\right]EL\beta_\mu \nonumber\\
&&\hspace{15mm}+\left[\db^{\mu\nu}_{1,LL}\Lb\beta_\nu+\lambda \db^{\mu\nu}_{2,LL}E\beta_\nu\right]L^2\beta_\mu 
\label{10.107}
\end{eqnarray}
where the coefficients are regular and of order 0. 

We summarize the above in the following proposition. 

\vspace{2.5mm}

\noindent{\bf Proposition 10.2} \ \ The quantities: 
$$\nu=\lambda E^2\lambda-j, \ \ \ \nub=\lambdab E^2\lambdab-\jb$$
where
\begin{eqnarray*}
&&j=\frac{1}{4}\beta_N^2\Lb^2 H-\frac{1}{2}\lambda\pi\beta_N^2 E\Lb H\\
&&\hspace{5mm}+\lambda^2\left\{\beta_N\left[\pi\sbeta-\frac{1}{4c}(\beta_N+2\beta_{\Nb})\right]E^2 H
\right.\\
&&\hspace{15mm}\left.+H\left[\pi\sbeta-\frac{1}{2c}(\beta_N+\beta_{\Nb})\right]N^\mu E^2\beta_\mu
\right\}
\end{eqnarray*}
and
\begin{eqnarray*}
&&\jb=\frac{1}{4}\beta_{\Nb}^2 L^2 H-\frac{1}{2}\lambdab\pi\beta_{\Nb}^2 ELH\\
&&\hspace{5mm}+\lambdab^2\left\{\beta_{\Nb}\left[\pi\sbeta-\frac{1}{4c}(2\beta_N+\beta_{\Nb})\right]
E^2 H\right.\\
&&\hspace{15mm}\left.+H\left[\pi\sbeta-\frac{1}{2c}(\beta_N+\beta_{\Nb})\right]\Nb^\mu E^2\beta_\mu
\right\}
\end{eqnarray*}
satisfy the following regularized propagation equations:
$$L\nu=K, \ \ \ \Lb\nub=\Kb$$
where
$$K=[K]_{P.A.}+\tilde{K}, \ \ \ \Kb=[\Kb]_{P.A.}+\tilde{\Kb}$$
are of order 2, and their principal acoustical parts are given by:
\begin{eqnarray*}
&&[K]_{P.A.}=2(L\lambda-\chi\lambda)E^2\lambda
+\frac{\lambda\beta_N^2}{4c}(\Lb H+2\pi\lambda EH)E^2\lambdab\\
&&\hspace{17mm}+\lambda AE\tchi+a \oA E\tchib
\end{eqnarray*}
and:
\begin{eqnarray*}
&&[\Kb]_{P.A.}=2(\Lb\lambdab-\chib\lambdab)E^2\lambdab
+\frac{\lambdab\beta_{\Nb}^2}{4c}(LH+2\pi\lambdab EH)E^2\lambda\nonumber\\
&&\hspace{17mm}+\lambdab \oAb E\tchib+a \Ab E\tchi
\end{eqnarray*}
with the coefficients $A, \oA$ and $\oAb, \Ab$ being given by \ref{10.99}, \ref{10.100} and 
\ref{10.105}, \ref{10.106} respectively. Also, the principal parts of the remainders $\tilde{K}$, 
$\tilde{\Kb}$ are of the form \ref{10.102}, \ref{10.107} respectively. 

\vspace{5mm}  

\section{The Structure Equations for the $N$th Approximants}

In the present section we shall derive analogues for the $N$th approximants of the acoustical 
structure equations of Chapter 3 for the case of 2 spatial dimensions and 
of the regularized propagation equations of the preceding two sections. 

We begin with the analogues 
of the propagation equations for $\lambda$, $\lambdab$ of Proposition 3.3, together with 
the analogue of the equation for $\zeta$. As we have seen in Section 3.3, these three equations 
are together deduced from the characteristic system \ref{2.62} by applying $\Lb$ to the equation 
for $Lx^\mu$, $L$ to the equation for $\Lb x^\mu$ and subtracting. For the $N$th approximants the 
characteristic system is replaced by the system \ref{9.15}:
\begin{eqnarray}
&&L_N x^\mu_N=\rho_N N^\mu_N+\vep^\mu_N \nonumber\\
&&\Lb_N x^\mu_N=\rhob_N\Nb^\mu_N+\vepb^\mu_N, \ \ \ \mbox{where} \ \ \vep^0_N=\vepb^0_N=0 
\label{10.108}
\end{eqnarray}
(see Proposition 9.1). Applying $\Lb_N$ to the first equation and $L_N$ to the second, we can 
express the right hand sides in terms of the $(E_N^\nu,N_N^\mu,\Nb_N^\mu)$ frame, defining 
analogues of the 5th and 6th of \ref{3.a15}:
\begin{eqnarray}
&&\Lb_N N^\mu_N=\sn_N E^\mu_N+n_N N^\mu_N+\on_N\Nb^\mu_N \nonumber\\
&&L_N\Nb^\mu_N=\snb_N E^\mu_N+\nb_N N^\mu_N+\onb_N\Nb^\mu_N 
\label{10.109}
\end{eqnarray}
and defining also, in regard to the error terms:
\begin{eqnarray}
&&\Lb_N\vep^\mu_N=\vep_{\Lb,N}^E E^\mu_N+\vep_{\Lb,N}^N N^\mu_N+\vep_{\Lb,N}^{\Nb}\Nb^\mu_N \nonumber\\
&&L_N\vepb^\mu_N=\vepb_{L,N}E^\mu_N+\vepb_{L,N}^N N^\mu_N+\vepb_{L,N}^{\Nb}\Nb^\mu_N 
\label{10.110}
\end{eqnarray}
We then obtain:
\begin{eqnarray}
&&\Lb_N L_N x^\mu_N=(\rho_N\sn_N+\vep_{\Lb,N}^E)E^\mu_N+(\Lb_N\rho_N+\rho_N n_N+\vep_{\Lb,N}^N)N^\mu_N
\nonumber\\
&&\hspace{20mm}+(\rho_N \on_N+\vep_{\Lb,N}^{\Nb})\Nb^\mu_N\nonumber\\
&&L_N\Lb_N x^\mu_N=(\rhob_N\snb_N+\vepb_{L,N}^E)E^\mu_N+(\rhob_N\nb_N+\vepb_{L,N}^N)N^\mu_N\nonumber\\
&&\hspace{20mm}+(L_N\rhob_N+\rhob_N\onb_N+\vepb_{L,N}^{\Nb})\Nb^\mu_N
\label{10.111}
\end{eqnarray}
Subtracting, results, in view of the 3rd of the commutation relations \ref{9.b1}, in the following:
\begin{eqnarray}
&&\left(-\zeta_N+\rho_N\sn_N-\rhob_N\snb_N+\vep_{\Lb,N}^E-\vepb_{L,N}^E\right)E^\mu_N\nonumber\\
&&+\left(\Lb_N\rho_N+\rho_N n_N-\rhob_N\nb_N+\vep_{\Lb,N}^N-\vepb_{L,N}^N\right)N^\mu_N\nonumber\\
&&+\left(-L_N\rhob_N-\rhob_N\onb_N+\rho_N\on_N -\vepb_{L,N}^{\Nb}+\vep_{\Lb,N}^{\Nb}\right)\Nb^\mu_N 
=0
\label{10.112}
\end{eqnarray}
This represents three equations, the vanishing of the coefficients of $E^\mu_N$, $N^\mu_N$, and 
$\Nb^\mu_N$.

The vanishing of the coefficient of $E^\mu_N$ is the following equation for $\zeta_N$:
\begin{equation}
\zeta_N=\rho_N\sn_N-\rhob_N\snb_N+\vep_{\Lb,N}^E-\vepb_{L,N}^E
\label{10.113}
\end{equation}
We shall presently determine $\sn_N$, $\snb_N$. According to \ref{10.109}:
\begin{equation}
\sn_N=h_{\mu\nu,N}E^\mu_N\Lb_N N^\nu_N, \ \ \ \snb_N=h_{\mu\nu,N}E^\mu_N L_N\Nb^\nu_N
\label{10.114}
\end{equation}
In view of the fact that $h_{\mu\nu,N}E^\mu_N\Nb^\nu_N=0$, we have:
\begin{equation}
\snb_N=-h_{\mu\nu,N}(L_N E^\mu_N)\Nb^\nu_N-(L_N h_{\mu\nu,N})E^\mu_N\Nb^\nu_N
\label{10.115}
\end{equation}
Now, by the 1st of the commutation relations \ref{9.b1}:
\begin{equation}
L_N E^\mu_N=E_N L^\mu_N-\chi_N E^\mu_N=E_N(\rho_N N^\mu_N+\vep^\mu_N)-\chi_N E^\mu_N
\label{10.116}
\end{equation}
We expand $E_N N^\mu_N$ and $E_N\Nb^\mu_N$ in the $(E^\mu_N,N^\mu_N,\Nb^\mu_N)$ frame in analogy with the first two of \ref{3.a15}:
\begin{eqnarray}
&&E_N N^\mu_N=\sk_N E^\mu_N+k_N N^\mu_N+\ok_N\Nb^\mu_N \nonumber\\
&&E_N\Nb^\mu_N=\skb_N E^\mu_N+\kb_N N^\mu_N+\okb_N\Nb^\mu_N 
\label{10.117}
\end{eqnarray}
We also expand the error term:
\begin{equation}
E_N\vep^\mu_N=\vep_{E,N}^E E^\mu_N+\vep_{E,N}^N N^\mu_N +\vep_{E,N}^{\Nb}\Nb^\mu_N
\label{10.118}
\end{equation}
We then obtain, through \ref{10.116},
\begin{equation}
h_{\mu\nu,N}(L_N E^\mu_N)\Nb^\nu_N=-2c_N(E_N\rho_N+\rho_N k_N+\vep_{E,N}^N)
\label{10.119}
\end{equation}
(see \ref{9.31}). By \ref{9.7} we have:
\begin{eqnarray}
&&(L_N h_{\mu\nu,N})E^\mu_N\Nb^\nu_N=(L_N H_N)\sbeta_N\beta_{\Nb,N} \label{10.120}\\
&&\hspace{20mm}+H_N\left\{E^\mu_N(L_N\beta_{\mu,N})\beta_{\Nb,N}+\sbeta_N\Nb^\nu_N L_N\beta_{\nu,N}
\right\}\nonumber
\end{eqnarray}
where we denote:
\begin{equation}
\sbeta_N=E_N^\mu\beta_{\mu,N}, \ \ \ \beta_{N,N}=N^\mu_N\beta_{\mu,N}, \ \ \ 
\beta_{\Nb,N}=\Nb^\mu_N\beta_{\mu,N}
\label{10.121}
\end{equation}
In view of the 2nd of \ref{9.44} we can write:
\begin{eqnarray}
&&E^\mu_NL_N\beta_{\mu,N}=L^\mu_N E_N\beta_{\mu,N}-\omega_{LE,N}\nonumber\\
&&\hspace{18mm}=(\rho_N N^\mu_N+\vep^\mu_N)E_N\beta_{\mu,N}-\omega_{LE,N}\nonumber\\
&&\hspace{18mm}=\rho_N\ss_{N,N}+\vep_N^\mu E_N\beta_{\mu,N}-\omega_{LE,N}
\label{10.122}
\end{eqnarray}
where we denote:
\begin{equation}
\ss_{N,N}=N_N^\mu E_N\beta_{\mu,N}, \ \ \ \ss_{\Nb,N}=\Nb^\mu_N E_N\beta_{\mu,N}
\label{10.123}
\end{equation}
Substituting the above in \ref{10.120} and the result, together with \ref{10.119}, in \ref{10.115} 
we obtain:
\begin{eqnarray}
&&\snb_N=2c_N(E_N\rho_N+\rho_N k_N+\vep_{E,N}^N)-\sbeta_N\beta_{\Nb,N}L_N H_N\nonumber\\
&&\hspace{8mm}-H_N\beta_{\Nb,N}\left(\rho_N\ss_{N,N}+\vep_N^\mu E_N\beta_{\mu,N}-\omega_{LE,N}\right)
\nonumber\\
&&\hspace{8mm}-H_N\sbeta_N\Nb^\nu_N L_N\beta_{\nu,N}
\label{10.124}
\end{eqnarray}
Similarly, we obtain the conjugate equation:
\begin{eqnarray}
&&\sn_N=2c_N(E_N\rhob_N+\rhob_N\okb_N+\vepb_{E,N}^{\Nb})-\sbeta_N\beta_{N,N}\Lb_N H_N\nonumber\\
&&\hspace{8mm}-H_N\beta_{N,N}\left(\rhob_N\ss_{\Nb,N}+\vepb_N^\mu E_N\beta_{\mu,N}
-\omega_{\Lb E,N}\right)
\nonumber\\
&&\hspace{8mm}-H_N\sbeta_N N^\nu_N\Lb_N\beta_{\nu,N}
\label{10.125}
\end{eqnarray}
Substituting \ref{10.125} and \ref{10.124} in \ref{10.113}, and writing 
\begin{eqnarray*}
&&\rhob_N\Nb^\nu_N L_N\beta_{\nu,N}-\rho_N N^\nu_N\Lb_N\beta_{\nu,N}=\\
&&\hspace{20mm}(\Lb_N x^\nu_N)L_N\beta_{\nu,N}-(L_N x^\nu_N)\Lb_N\beta_{\nu,N}
-\vepb_N^\nu L_N\beta_{\nu,N}+\vep_N^\nu\Lb_N\beta_{\nu,N} \\
&&\hspace{30mm}=-\omega_{L\Lb,N}-\vepb_N^\nu L_N\beta_{\nu,N}+\vep_N^\nu\Lb_N\beta_{\nu,N} 
\end{eqnarray*}
using the 1st of \ref{9.44} we obtain:
\begin{eqnarray}
&&\zeta_N=2c_N(\rho_N E_N\rhob_N-\rhob_N E_N\rho_N)+2c_N\rho_N\rhob_N(\okb_N-k_N)\nonumber\\
&&\hspace{10mm}+\sbeta_N(\beta_{\Nb,N}\rhob_N L_N H_N-\beta_{N,N}\rho_N\Lb_N H_N)\nonumber\\
&&\hspace{10mm}+\rho_N\rhob_N H_N(\beta_{\Nb,N}\ss_{N,N}-\beta_{N,N}\ss_{\Nb,N})+\vep_{\zeta,N}
\label{10.126}
\end{eqnarray}
where:
\begin{eqnarray}
&&\vep_{\zeta,N}=\vep_{\Lb,N}^E-\vepb_{L,N}^E-\omega_{L\Lb,N}
-\vepb_N^\nu L_N\beta_{\nu,N}+\vep_N^\nu\Lb_N\beta_{\nu,N} \nonumber\\
&&\hspace{10mm}+2c_N(\rho_N\vepb_{E,N}^{\Nb}-\rhob_N\vep_{E,N}^N)
+H_N(\beta_{\Nb,N}\rhob_N\vep^\mu_N-\beta_{N,N}\rho_N\vepb^\mu_N)E_N\beta_{\mu,N}\nonumber\\
&&\hspace{10mm}+H_N(\beta_{N,N}\rho_N\omega_{\Lb E,N}-\beta_{\Nb,N}\rhob_N\omega_{LE,N})
\label{10.127}
\end{eqnarray}
Using Propositions 9.1 and 9.2 we deduce:
\begin{equation}
\vep_{\zeta,N}=O(\tau^N)
\label{10.128}
\end{equation}

Next we shall determine $k_N$ and $\okb_N$. Since $N^0_N=\Nb^0_N=1$, the expansions \ref{10.117} 
imply, in analogy with \ref{3.46}, 
\begin{equation}
\sk_N\pi_N+k_N+\ok_N=0, \ \ \ \skb_N\pi_N+\kb_N+\okb_N=0 
\label{10.129}
\end{equation}
where:
\begin{equation}
\pi_N=E_N^0=E_N t_N
\label{10.130}
\end{equation}
Consider then $k_N$. This is determined through $\ok_N$ and $\sk_N$ by the 1st of \ref{10.129}. 
In regard to $\ok_N$, we have:
$$\ok_N=-\frac{1}{2c_N}h_{\mu\nu,N}(E_N N^\mu_N)N^\nu_N
=\frac{1}{4c_N}(E_N h_{\mu\nu,N})N_N^\mu N_N^\nu$$
By \ref{9.9} we have:
\begin{equation}
(E_N h_{\mu\nu,N})N_N^\mu N_N^\nu=\beta_{N,N}^2 E_N H_N+2H_N\beta_{N,N}\ss_{N,N}
\label{10.131}
\end{equation}
We thus obtain:
\begin{equation}
\ok_N=\frac{1}{4c_N}(\beta_{N,N}^2 E_N H_N+2H_N\beta_{N,N}\ss_{N,N})
\label{10.132}
\end{equation}
which is analogous to \ref{3.53}. In regard to $\sk$ we define, in analogy with the 1st of \ref{3.a21}, 
\begin{equation}
\tchi_N=\sk_N+H_N\sbeta_N\ss_{N,N}
\label{10.133}
\end{equation}

We proceed to relate $\tchi_N$ to $\chi_N$ as defined by the 1st of \ref{9.b2}. Since $\sh_N$ is 
defined by \ref{9.9} we have:
\begin{eqnarray}
&&L_N\sh_N-2\sh_N\Omega b_N=2h_{\mu\nu,N}\left(\frac{\partial\Omega_N^\mu}{\partial\ub}
-b_N\frac{\partial\Omega_N^\mu}{\partial\vartheta}-\frac{\partial b_N}{\partial\vartheta}\Omega_N^\mu 
\right)\Omega_N^\nu \nonumber\\
&&\hspace{32mm}+\left(\frac{\partial h_{\mu\nu,N}}{\partial\ub}-b_N\frac{\partial h_{\mu\nu,N}}{\partial\vartheta}\right)\Omega_N^\mu\Omega_N^\nu 
\label{10.134}
\end{eqnarray}
On the other hand, we have:
\begin{eqnarray}
&&\frac{\partial(L_N x_N^\mu)}{\partial\vartheta}=\frac{\partial}{\partial\vartheta}\left(
\frac{\partial x_N^\mu}{\partial\ub}-b_N\frac{\partial x_N^\mu}{\partial\vartheta}\right) \nonumber\\
&&\hspace{17mm}=\frac{\partial}{\partial\ub}\left(\frac{\partial x_N^\mu}{\partial\vartheta}\right)
-b_N\frac{\partial}{\partial\vartheta}\left(\frac{\partial x_N^\mu}{\partial\vartheta}\right)
-\frac{\partial b_N}{\partial\vartheta}\frac{\partial x_N^\mu}{\partial\vartheta} \nonumber\\
&&\hspace{17mm}=\frac{\partial\Omega_N^\mu}{\partial\ub}-b_N\frac{\partial\Omega_N^\mu}{\partial\vartheta}-\frac{\partial b_N}{\partial\vartheta}\Omega_N^\mu
\label{10.135}
\end{eqnarray}
Therefore the 1st term on the right in \ref{10.134} is:
\begin{eqnarray}
&&2h_{\mu\nu,N}\frac{\partial(L_N x_N^\mu)}{\partial\vartheta}\Omega_N^\nu=
2h_{\mu\nu,N}\frac{\partial}{\partial\vartheta}(\rho_N N_N^\mu+\vep_N^\mu) \Omega_N^\nu \nonumber\\
&&\hspace{25mm}=2h_{\mu\nu,N}\rho_N\frac{\partial N_N^\mu}{\partial\vartheta} \Omega_N^\nu 
+2h_{\mu\nu,N}\frac{\partial\vep_N^\mu}{\partial\vartheta} \Omega_N^\nu\nonumber\\
&&\hspace{25mm}=2\sh_N\rho_N h_{\mu\nu,N}(E_N N_N^\mu)E_N^\nu+2\sh_N h_{\mu\nu,N}(E_N\vep_N^\mu)E_N^\nu
\nonumber\\
&&\hspace{25mm}=2\sh_N(\rho_N \sk_N+\vep_{E,N}^E)
\label{10.136}
\end{eqnarray}
(by \ref{10.117}, \ref{10.118}). Also, the 2nd term on the right in \ref{10.134} is:
\begin{eqnarray}
&&(L_N h_{\mu\nu,N})\Omega_N^\mu\Omega_N^\nu=\sh_N\left\{(L_N H_N)\sbeta_N^2+
2H_N\sbeta_N E_N^\mu L_N\beta_{\mu,N}\right\}\nonumber\\
&&\hspace{15mm}=\sh_N\left\{(L_N H_N)\sbeta_N^2+2H_N\sbeta_N\left((\rho_N \ss_{N,N}+\vep_N^\mu E_N\beta_{\mu,N})-\omega_{LE,N}\right)\right\}\nonumber\\
&&\label{10.137}
\end{eqnarray}
by the 2nd of \ref{9.44}. Substituting \ref{10.136} and \ref{10.137} in \ref{10.134} we conclude 
through the 1st of \ref{9.b2}, in view of \ref{10.133}, that: 
\begin{equation}
\chi_N=\rho_N\tchi_N+\frac{1}{2}\sbeta_N^2 L_N H_N +\vep_{\chi,N}
\label{10.138}
\end{equation}
(compare with the 1st of \ref{3.a20}) where
\begin{equation}
\vep_{\chi,N}=\vep_{E,N}^E+H_N\sbeta_N(\vep_N^\mu E_N\beta_{\mu,N}-\omega_{LE,N})
\label{10.139}
\end{equation}
Using Propositions 9.1 and 9.2 we deduce:
\begin{equation}
\vep_{\chi,N}=O(\tau^N)
\label{10.140}
\end{equation}
In a similar manner, we derive through the 2nd of \ref{9.b2} the conjugate equation:
\begin{equation}
\chib_N=\rhob_N\tchib_N+\frac{1}{2}\sbeta_N^2\Lb_N H_N+\vep_{\chib,N}
\label{10.141}
\end{equation}
(compare with the 2nd of \ref{3.a20}) where $\tchib_N$ is defined in analogy with the 2nd of 
\ref{3.a21} by:
\begin{equation}
\tchib_N=\skb_N+H_N\sbeta_N\ss_{\Nb,N}
\label{10.142}
\end{equation}
and
\begin{equation}
\vep_{\chib,N}=\vepb_{E,N}^E+H_N\sbeta_N(\vepb_N^\mu E_N\beta_{\mu,N}-\omega_{\Lb E,N})
\label{10.143}
\end{equation}
Using Propositions 9.1 and 9.2 we deduce:
\begin{equation}
\vep_{\chib,N}=O(\tau^{N+1})
\label{10.144}
\end{equation}

We now go back to the 1st of \ref{10.129}. By \ref{10.132} and \ref{10.133} this gives:
\begin{eqnarray}
&&k_N=-\pi_N\tchi_N+\pi_N\sbeta_N H_N\ss_{N,N}\nonumber\\
&&\hspace{8mm}-\frac{\beta_{N,N}^2}{4c_N}E_N H_N-\frac{\beta_{N,N}}{2c_N}H_N\ss_{N,N}
\label{10.145}
\end{eqnarray}
Similarly, by the 2nd of \ref{10.129}, the conjugate of \ref{10.132}, which reads:
\begin{equation}
\kb_N=\frac{1}{4c_N}(\beta_{\Nb,N}^2 E_N H_N+2H_N\beta_{\Nb,N}\ss_{\Nb,N})
\label{10.146}
\end{equation}
and \ref{10.142}, we obtain:
\begin{eqnarray}
&&\okb_N=-\pi_N\tchib_N+\pi_N\sbeta_N H_N\ss_{\Nb,N}\nonumber\\
&&\hspace{8mm}-\frac{\beta_{\Nb,N}^2}{4c_N}E_N H_N-\frac{\beta_{\Nb,N}}{2c_N}H_N\ss_{\Nb,N}
\label{10.147}
\end{eqnarray}
(the conjugate of \ref{10.145}).

The difference $\okb_N-k_N$ substituted in \ref{10.126} from \ref{10.147}, \ref{10.145}, gives 
us an explicit expression for $\zeta_N$. Equations \ref{9.198}, \ref{9.199} then determine 
$\eta_N$, $\etab_N$. Also, \ref{10.147}, \ref{10.145} substituted in \ref{10.125}, 
\ref{10.124} give us explicit expressions for $\sn$, $\snb$. To express these quantities in 
terms of $E_N\lambda_N$, $E_N\lambdab_N$ (see \ref{9.101}) however requires that we also calculate 
$E_N c_N$. From \ref{9.31}, 
\begin{eqnarray}
&&E_N c_N=-\frac{1}{2}(E_N h_{\mu\nu,N})N_N^\mu\Nb_N^\nu-\frac{1}{2}h_{\mu\nu,N}(E_N N_N^\mu)\Nb_N^\nu
-\frac{1}{2}h_{\mu\nu,N}N_N^\mu(E_N N_N^\nu)\nonumber\\
&&\hspace{10mm}=-\frac{1}{2}\beta_{N,N}\beta_{\Nb,N}E_N H_N-\frac{1}{2}H_N(\beta_{\Nb,N}\ss_{N,N}
+\beta_{N,N}\ss_{\Nb,N})\nonumber\\
&&\hspace{15mm}+c_N(\okb_N+k_N)
\label{10.148}
\end{eqnarray}
and we can substitute for the sum $\okb_N+k_N$ from \ref{10.147}, \ref{10.145}. 

We now consider the other two equations represented by \ref{10.112}, namely the vanishing of the 
coefficients of $N_N^\mu$ and $\Nb_N^\mu$:
\begin{eqnarray}
&&\Lb_N\rho_N+\rho_N n_N-\rhob_N\nb_N+\vep_{\Lb,N}^N-\vepb_{L,N}^N=0 \label{10.149}\\
&&L_N\rhob_N+\rhob_N\onb_N-\rho_N\on_N +\vepb_{L,N}^{\Nb}-\vep_{\Lb,N}^{\Nb}=0 \label{10.150}
\end{eqnarray}
These will yield propagation equations for $\lambda_N$, $\lambdab_N$ analogous to those of 
Proposition 3.3. To derive these we must first calculate $L_N c_N$, $\Lb_N c_N$. We expand 
$L_N N_N^\mu$ and $\Lb_N\Nb_N^\mu$ in the $(E_N^\mu,N_N^\mu,\Nb_N^\mu)$ frame, in analogy with 
the 3rd and 4th of \ref{3.a15}:
\begin{eqnarray}
&&L_N N_N^\mu=\sm_N E_N^\mu+m_N N_N^\mu+\om_N\Nb_N^\mu \label{10.151}\\
&&\Lb_N\Nb_N^\mu=\smb_N E_N^\mu+\mb_N N_N^\mu+\omb_N\Nb_N^\mu \label{10.152} 
\end{eqnarray}
From \ref{9.31}, 
\begin{eqnarray}
&&L_N c_N=-\frac{1}{2}(L_N h_{\mu\nu,N})N_N^\mu\Nb_N^\nu-\frac{1}{2}h_{\mu\nu,N}(L_N N_N^\mu)\Nb_N^\nu
-\frac{1}{2}h_{\mu\nu,N}N_N^\mu(L_N N_N^\nu)\nonumber\\
&&\hspace{10mm}=-\frac{1}{2}\beta_{N,N}\beta_{\Nb,N}L_N H_N-\frac{1}{2}H_N(\beta_{\Nb,N}s_{NL,N}
+\beta_{N,N}\Nb_N^\nu L_N\beta_{\nu,N})\nonumber\\
&&\hspace{15mm}+c_N(m_N+\onb_N)
\label{10.153}
\end{eqnarray}
Here, we denote:
\begin{equation}
s_{NL,N}=N_N^\mu L_N\beta_{\mu,N}, \ \ \ s_{\Nb\Lb,N}=\Nb_N^\mu\Lb_N\beta_{\mu,N}
\label{10.154}
\end{equation}
Similarly, we obtain the conjugate equation: 
\begin{eqnarray}
&&\Lb_N c_N=-\frac{1}{2}\beta_{N,N}\beta_{\Nb,N}\Lb_N H_N-\frac{1}{2}H_N(\beta_{N,N}s_{\Nb\Lb,N}
+\beta_{\Nb,N}N_N^\nu \Lb_N\beta_{\nu,N})\nonumber\\
&&\hspace{15mm}+c_N(\omb_N+n_N)
\label{10.155}
\end{eqnarray}
Using equation \ref{10.150} together with equation \ref{10.153} to derive an equation for 
$L_N\lambda_N$, the terms in $\onb_N$ cancel and we obtain: 
\begin{equation}
L_N\lambda_N=p_N \lambda_N+\on_N\lambdab_N-\frac{1}{2}H_N\beta_{N,N}\rhob_N\Nb_N^\nu L_N\beta_{\nu,N} 
+\vep_{\Lb,N}^{\Nb}-\vepb_{L,N}^{\Nb}
\label{10.156}
\end{equation}
where:
\begin{equation}
p_N=m_N-\frac{1}{2c_N}(\beta_{N,N}\beta_{\Nb,N}L_N H_N+H_N\beta_{\Nb,N}s_{NL,N})
\label{10.157}
\end{equation}
Now, from the 1st of \ref{10.109}, 
\begin{eqnarray}
&&\on_N=-\frac{1}{2c_N}h_{\mu\nu,N}(\Lb_N N_N^\mu)N_N^\nu
=\frac{1}{4c_N}(\Lb_N h_{\mu\nu,N})N_N^\mu N_N^\nu \nonumber\\
&&\hspace{6mm}=\frac{1}{4c_N}\left\{\beta_{N,N}^2\Lb_N H_N+2H_N\beta_{N,N} N_N^\nu\Lb_N\beta_{\nu,N}\right\}
\label{10.158}
\end{eqnarray}
Setting then:
\begin{equation}
q_N=\frac{1}{4c_N}\beta_{N,N}^2\Lb_N H_N
\label{10.159}
\end{equation}
equation \ref{10.156} takes the form:
\begin{eqnarray}
&&L_N\lambda_N=p_N\lambda_N+q_N\lambdab_N \label{10.160}\\
&&\hspace{10mm}+\frac{1}{2}H_N\beta_{N,N}
\left\{\rho_N N_N^\nu\Lb_N\beta_{\nu,N}-\rhob_N\Nb_N^\nu L_N\beta_{\nu,N}\right\}
+\vep_{\Lb,N}^{\Nb}-\vepb_{L,N}^{\Nb}\nonumber
\end{eqnarray}
Here, the term in parenthesis is:
\begin{eqnarray*}
&&(L_N x_N^\nu)\Lb_N\beta_{\nu,N}-(\Lb_N x_N^\nu)L_N\beta_{\nu,N}-\vep_N^\nu\Lb_N\beta_{\nu,N}
+\vepb_N^\nu L_N\beta_{\nu,N}\\
&&\hspace{30mm}=\omega_{L\Lb,N}-\vep_N^\nu\Lb_N\beta_{\nu,N}
+\vepb_N^\nu L_N\beta_{\nu,N}
\end{eqnarray*}
by the first of \ref{9.44}. We conclude that $\lambda_N$ satisfies the propagation equation:
\begin{equation}
L_N\lambda_N=p_N\lambda_N+q_N\lambdab_N+\vep_{\lambda,N}
\label{10.161}
\end{equation}
where:
\begin{equation}
\vep_{\lambda,N}=\vep_{\Lb,N}^{\Nb}-\vepb_{L,N}^{\Nb}-\vep_N^\nu\Lb_N\beta_{\nu,N}
+\vepb_N^\nu L_N\beta_{\nu,N}+\omega_{L\Lb,N}
\label{10.162}
\end{equation}
Using Propositions 9.1 and 9.2 we deduce:
\begin{equation}
\vep_{\lambda,N}=O(\tau^N)
\label{10.163}
\end{equation}
Similarly, using equation \ref{10.149} together with equation \ref{10.155} we derive the conjugate 
equation, which is the following propagation equation for $\lambdab_N$:
\begin{equation}
\Lb_N\lambdab_N=\pb_N\lambdab_N+\qb_N\lambda_N+\vep_{\lambdab,N}
\label{10.164}
\end{equation}
Here: 
\begin{equation}
\pb_N=\omb_N-\frac{1}{2c_N}(\beta_{N,N}\beta_{\Nb,N}\Lb_N H_N+H_N\beta_{N,N}s_{\Nb\Lb,N}), \ \ \ 
\qb_N=\frac{1}{4c_N}\beta_{\Nb,N}^2 L_N H_N
\label{10.165}
\end{equation}
and:
\begin{equation}
\vep_{\lambdab,N}=\vepb_{L,N}^N-\vep_{\Lb,N}^N-\vepb_N^\nu L_N\beta_{\nu,N}
+\vep_N^\nu\Lb_N\beta_{\nu,N}-\omega_{L\Lb,N}
\label{10.166}
\end{equation}
Using Propositions 9.1 and 9.2 we deduce:
\begin{equation}
\vep_{\lambdab,N}=O(\tau^N)
\label{10.167}
\end{equation}
The coefficients $m_N$ and $\omb_N$ which enter \ref{10.157} and \ref{10.165} respectively, are 
determined through the coefficients $\sm_N$, $\om_N$ and $\smb_N$, $\mb_N$ respectively, by the relations 
\begin{equation}
\sm_N\pi_N+m_N+\om_N=0, \ \ \ \smb_N\pi_N+\mb_N+\omb_N=0
\label{10.168}
\end{equation}
which follow from \ref{10.151}, \ref{10.152} in view of $N_N^0=\Nb_N^0=1$. As for $\sm_N$ we have:
\begin{equation}
\sm_N=h_{\mu\nu,N}(L_N N_N^\mu)E_N^\nu=-(L_N h_{\mu\nu,N})N_N^\mu E_N^\nu
-h_{\mu\nu,N}N_N^\mu L_N E_N^\nu
\label{10.169}
\end{equation}
In view of the 1st of the commutation relations \ref{9.b1} we have:
\begin{eqnarray}
&&h_{\mu\nu,N}N_N^\mu L_N E_N^\nu=h_{\mu\nu,N}N_N^\mu E_N L_N^\nu
=h_{\mu\nu,N}N_N^\mu E_N (\rho_N N_N^\nu+\vep_N^\nu)\nonumber\\
&&\hspace{25mm}=\rho_N h_{\mu\nu,N}N_N^\mu E_N N_N^\nu-2c_N\vep_{E,N}^{\Nb}\nonumber\\
&&\hspace{25mm}=-\frac{1}{2}\rho_N (E_N h_{\mu\nu,N})N_N^\mu N_N^\nu-2c_N\vep_{E,N}^{\Nb}
\label{10.170}
\end{eqnarray}
We then obtain:
\begin{eqnarray}
&&\sm_N=-\beta_{N,N}\sbeta_N L_N H_N+\frac{1}{2}\rho_N\beta_{N,N}^2 E_N H_N
-H_N\sbeta_N s_{NL,N}\nonumber\\
&&\hspace{7mm}+H_N\beta_{N,N}(\omega_{LE,N}-\vep_N^\mu E_N\beta_{\mu,N})-2c_N\vep_{E,N}^{\Nb}
\label{10.171}
\end{eqnarray}
(compare with \ref{3.68}). Similarly, we have the conjugate equation:
\begin{eqnarray}
&&\smb_N=-\beta_{\Nb,N}\sbeta_N\Lb_N H_N+\frac{1}{2}\rhob_N\beta_{\Nb,N}^2 E_N H_N 
-H_N\sbeta_N s_{\Nb\Lb,N}\nonumber\\
&&\hspace{7mm}+H_N\beta_{\Nb,N}(\omega_{\Lb E,N}-\vepb_N^\mu E_N\beta_{\mu,N})
-2c_N\vepb_{E,N}^N
\label{10.172}
\end{eqnarray}
(compare with \ref{3.69}). For $\om_N$ we have:
$$\om_N=-\frac{1}{2c_N}h_{\mu\nu,N}(L_N N_N^\mu)N_N^\nu=\frac{1}{4c_N}(L_N h_{\mu\nu,N})N_N^\mu N_N^\nu$$
hence:
\begin{equation}
\om_N=\frac{1}{4c_N}(\beta_{N,N}^2 L_N H_N+2H_N\beta_{N,N} s_{NL,N})
\label{10.173}
\end{equation}
and by conjugation:
\begin{equation}
\mb_N=\frac{1}{4c_N}(\beta_{\Nb,N}^2\Lb_N H_N+2H_N\beta_{\Nb,N} s_{\Nb\Lb,N})
\label{10.174}
\end{equation}

We proceed to derive the analogues of the 2nd variation equations \ref{3.a44}, \ref{3.a45} for the 
$N$th approximants. We shall follow a somewhat different approach than that of Chapter 3, which 
relies more on the $(E,N,\Nb)$ frame and on the expansions \ref{3.a31}, \ref{3.a33}, \ref{3.a38}, 
besides the expansions \ref{3.a15}. This approach yields the following form of the 2nd variation 
equations:
\begin{eqnarray}
&&L\sk=E\sm-\chi\sk+\sm\se+m\sk+\om\skb+\sk\sf-2c(k\of+\ok f) \label{10.175}\\
&&\hspace{8mm}+(\sbeta\sk+\beta_N k+\beta_{\Nb}\ok)(\sbeta LH+\rho H\ss_N)
+H\sbeta(\rho\ss_N\sk+s_{NL}k+\rho s_{N\Nb}\ok) \nonumber
\end{eqnarray}
\begin{eqnarray}
&&\Lb\skb=E\smb-\chib\skb+\smb\se+\omb\skb+\mb\sk+\skb\sfb-2c(\okb\fb+\kb\ofb) \label{10.176}\\
&&\hspace{8mm}+(\sbeta\skb+\beta_{\Nb}\okb+\beta_N\kb)(\sbeta\Lb H+\rhob H\ss_{\Nb})
+H\sbeta(\rhob\ss_{\Nb}\skb+s_{\Nb\Lb}\okb+\rhob s_{N\Nb}\kb) \nonumber
\end{eqnarray}
From  the 1st of \ref{3.a30} and of \ref{3.a32} we see that:
\begin{equation}
\se+\sgamma=0
\label{10.177}
\end{equation}
Taking account also of the formulas \ref{3.a34} and their conjugates we see that equations 
\ref{10.175}, \ref{10.176} indeed coincide with equations \ref{3.a44}, \ref{3.a45}. 

According to the 1st of the expansions \ref{10.117} we have:
\begin{equation}
\sk_N=h_{\mu\nu,N}(E_N N_N^\mu)E_N^\nu
\label{10.178}
\end{equation}
Applying $L_N$ to this we obtain:
\begin{eqnarray}
&&L_N\sk_N=h_{\mu\nu,N}(L_N E_N N_N^\mu)E_N^\nu+h_{\mu\nu,N}(E_N N_N^\mu)(L_N E_N^\nu) \nonumber\\
&&\hspace{13mm}+(L_N h_{\mu\nu,N})(E_N N_N^\mu)E_N^\nu
\label{10.179}
\end{eqnarray}
By the 1st of the commutation relations \ref{9.b1} the 1st term on the right is:
\begin{equation}
h_{\mu\nu,N}(E_N L_N N_N^\mu)E_N^\nu-\chi_N\sk_N
\label{10.180}
\end{equation}
Applying $E_N$ to \ref{10.151} we obtain:
\begin{eqnarray}
&&E_N L_N N_N^\mu=(E_N\sm_N)E_N^\mu+(E_N m_N)N_N^\mu+(E_N\om_N)\Nb_N^\mu \nonumber\\
&&\hspace{20mm}+\sm_N(\se_N E_N^\mu+e_N N_N^\mu+\oe_N\Nb_N^\mu)\nonumber\\
&&\hspace{20mm}+m_N(\sk_N E_N^\mu+k_N N_N^\mu+\ok_N\Nb_N^\mu)\nonumber\\
&&\hspace{20mm}+\om_N(\skb_N E_N^\mu+\kb_N N_N^\mu+\okb_N\Nb_N^\mu)
\label{10.181}
\end{eqnarray}
Here, besides the expansions \ref{10.117} we have also substituted the expansion:
\begin{equation}
E_N E_N^\mu=\se_N E_N^\mu+e_N N_N^\mu+\oe_N\Nb_N^\mu 
\label{10.182}
\end{equation}
Consequently, \ref{10.180}, the 1st term on the right in \ref{10.179}, is:
\begin{equation}
E_N\sm_N+\sm_N\se_N+m_N\sk_N+\om_N\skb_N-\chi_N\sk_N
\label{10.183}
\end{equation}
Consider next the 2nd term on the right in \ref{10.179}. We expand $L_N E_N^\mu$ in the 
$(E_N^\mu,N_N^\mu,\Nb_N^\mu)$ frame:
\begin{equation}
L_N E_N^\mu=\sf_N E_N^\mu+f_N N_N^\mu+\of_N\Nb_N^\mu 
\label{10.184}
\end{equation}
Substituting the expansions 1st of \ref{10.117} and \ref{10.184} we find that the 2nd term on the right in 
\ref{10.179} is:
\begin{equation}
\sk_N\sf_N-2c_N(k_N\of_N+\ok_N f_N)
\label{10.185}
\end{equation}
In regard finally to the 3rd term on the right in \ref{10.179}, substituting the expansion 1st of \ref{10.117} 
we obtain that this term is:
\begin{eqnarray}
&&(L_N h_{\mu\nu,N})(\sk_N E_N^\mu+k_N N_N^\mu+\ok_N\Nb_N^\mu)E_N^\nu=\nonumber\\
&&\hspace{10mm}\sk_N\left(\sbeta_N^2 L_N H_N+2\sbeta_N H_N E_N^\mu L_N\beta_{\mu,N}\right)\nonumber\\
&&\hspace{7mm}+k_N\left(\sbeta_N\beta_{N,N}L_N H_N+H_N(\sbeta_N s_{NL,N}+\beta_{N,N}E_N^\mu L_N\beta_{\mu,N})\right)\nonumber\\
&&\hspace{7mm}+\ok_N\left(\sbeta_N\beta_{\Nb,N}L_N H_N+H_N(\sbeta_N\Nb_N^\mu L_N\beta_{\mu,N}
+\beta_{\Nb,N}E_N^\mu L_N\beta_{\mu,N})\right)\nonumber\\
&&\label{10.186}
\end{eqnarray}
Using also the 2nd of \ref{9.44} to express $E_N^\mu L_N\beta_{\mu,N}$ in terms of 
$\ss_{N,N}$, we conclude that $\sk_N$ satisfies the following $N$th approximant version of 
the 2nd variation equation \ref{10.175}:
\begin{eqnarray}
&&L_N\sk_N=E_N\sm_N-\chi_N\sk_N+\sm_N\se_N+m_N\sk_N+\om_N\skb_N\nonumber\\
&&\hspace{15mm}+\sk_N\sf_N-2c_N(k_N\of_N+\ok_N f_N)\nonumber\\
&&\hspace{15mm}+(\sbeta_N\sk_N+\beta_{N,N}k_N+\beta_{\Nb,N}\ok_N)(\sbeta_N L_N H_N+\rho_N\ss_{N,N})
\nonumber\\
&&\hspace{15mm}+H_N\sbeta_N(\rho_N\ss_{N,N}\sk_N+s_{NL,N}k_N+s_{\Nb L,N}\ok_N)+\vep_{L\sk,N} 
\nonumber\\
&&\label{10.187}
\end{eqnarray}
Here, 
\begin{equation}
s_{\Nb L,N}=\Nb_N^\mu L_N\beta_{\mu,N}
\label{10.188}
\end{equation}
and, 
\begin{equation}
\vep_{L\sk,N}=H_N(2\sbeta_N\sk_N+\beta_{N,N}k_N+\beta_{\Nb,N}\ok_N)
(\vep_N^\mu E_N\beta_{\mu,N}-\omega_{LE,N})=O(\tau^N)
\label{10.189}
\end{equation}
the estimate following from Propositions 9.1 and 9.2. In a similar manner we arrive at the conclusion 
that $\skb_N$ satisfies the following $N$th approximant version of the 2nd variation equation 
\ref{10.176}:
\begin{eqnarray}
&&\Lb_N\skb_N=E_N\smb_N-\chib_N\skb_N+\smb_N\se_N+\omb_N\skb_N+\mb_N\sk_N\nonumber\\
&&\hspace{15mm}+\skb_N\sfb_N-2c_N(\kb_N\fb_N+\kb_N\ofb_N)\nonumber\\
&&\hspace{15mm}+(\sbeta_N\skb_N+\beta_{\Nb,N}\okb_N+\beta_{N,N}\kb_N)(\sbeta_N\Lb_N H_N+\rhob_N\ss_{\Nb,N})
\nonumber\\
&&\hspace{15mm}+H_N\sbeta_N(\rhob_N\ss_{\Nb,N}\skb_N+s_{\Nb\Lb,N}\okb_N+s_{N\Lb,N}\kb)+\vep_{\Lb\skb,N}
\nonumber\\
&&\label{10.190}
\end{eqnarray}
Here, 
\begin{equation}
s_{N\Lb,N}=N_N^\mu\Lb_N\beta_{\mu,N}
\label{10.191}
\end{equation}
and, 
\begin{equation}
\vep_{\Lb\skb,N}=H_N(2\sbeta_N\skb_N+\beta_{\Nb,N}\okb_N+\beta_{N,N}\kb_N)
(\vepb_N^\mu E_N\beta_{\mu,N}-\omega_{\Lb E,N})=O(\tau^{N+1})
\label{10.192}
\end{equation}
the estimate following from Propositions 9.1 and 9.2. In deriving \ref{10.190} we have used the 
conjugate of the expansion \ref{10.184}: 
\begin{equation}
\Lb_N E_N^\mu=\sfb_N E_N^\mu+\fb_N N_N^\mu+\ofb_N \Nb_N^\mu 
\label{10.193}
\end{equation}
Note in regard to \ref{10.188}, \ref{10.191} that by \ref{9.a11} and \ref{9.44} we have:
\begin{eqnarray}
&&\frac{1}{2}(\rhob_N s_{\Nb L,N}+\rho_N s_{N\Lb,N})-a_N\sss_N=-\delta^\prime_N
-\frac{1}{2}(\vepb_N^\mu L_N\beta_{\mu,N}+\vep_N^\mu\Lb_N\beta_{\mu,N})\nonumber\\
&&\hspace{55mm}=O(\tau^N)\nonumber\\
&&\frac{1}{2}(\rhob_N s_{\Nb L,N}-\rho_N s_{N\Lb,N})=-\frac{1}{2}\omega_{L\Lb,N}
-\frac{1}{2}(\vepb_N^\mu L_N\beta_{\mu,N}-\vep_N^\mu\Lb_N\beta_{\mu,N})\nonumber\\
&&\hspace{55mm}=O(\tau^N)
\label{10.194}
\end{eqnarray}
the estimates following from Propositions 9.1, 9.2 and Lemma 9.2. Here,
\begin{equation}
\sss_N=E_N^\mu E_N\beta_{\mu,N}
\label{10.195}
\end{equation}
The coefficient $\se_N$, which appears in \ref{10.187} and \ref{10.190}, is given according to the 
expansion \ref{10.182} by:
$$\se_N=h_{\mu\nu,N}(E_N E_N^\mu)E_N^\nu=-\frac{1}{2}(E_N h_{\mu\nu,N})E_N^\mu E_N^\nu$$
$E_N$ being of unit magnitude with respect to $h_N$. We thus obtain:
\begin{equation}
\se_N=-\frac{1}{2}\sbeta_N^2 E_N H_N-\sbeta H_N \sss_N
\label{10.196}
\end{equation}
(compare with the 1st of \ref{3.a32}). As for the coefficients $\sf_N, f_N, \of_N$, which appear in \ref{10.187}, these are determined through the commutation relation, 1st of \ref{9.b1}, which implies:
$$L_N E_N^\mu-E_N L_N^\mu=-\chi_N E_N^\mu$$
substituting $L_N^\mu=\rho_N N_N^\mu+\vep_N^\mu$ and the expansions for 
$E_N N_N^\mu$ (1st of \ref{10.117}) and for $E_N\vep_N^\mu$ (\ref{10.118}) . We obtain:
\begin{equation}
\sf_N=\rho_N\sk_N-\chi_N+\vep_{E,N}^E, \ \ \ f_N=E_N\rho_N+\rho_N k_N+\vep_{E,N}^N, \ \ \ 
\of_N=\rho_N\ok_N+\vep_{E,N}^{\Nb}
\label{10.197}
\end{equation}
(compare with \ref{3.a34}). The coefficients $\sfb_N, \ofb_N, \fb_N$, which appear in \ref{10.190}, are 
given by the conjugates of these formulas:
\begin{equation}
\sfb_N=\rhob_N\skb_N-\chib_N+\vepb_{E,N}^E, \ \ \ \ofb_N=E_N\rhob_N+\rhob_N\okb_N+\vepb_{E,N}^{\Nb}, 
\ \ \ \fb_N=\rhob_N\kb_N+\vepb_{E,N}^N
\label{10.198}
\end{equation}

We proceed to derive analogues of the cross variation equations \ref{3.a46}, \ref{3.a47} for the 
$N$th approximants. Following the approach just used yields the following form of the cross variation 
equations:
\begin{eqnarray}
&&\Lb\sk=E\sn-\chib\sk+\sn\se+n\sk+\on\skb+\sk\sfb-2c(k\ofb+\ok\fb) \label{10.199}\\
&&\hspace{8mm}+(\sbeta\sk+\beta_N k+\beta_{\Nb}\ok)(\sbeta\Lb H+\rhob H\ss_{\Nb})
+H\sbeta(\rhob\ss_{\Nb}\sk+\rhob s_{N\Nb}k+s_{\Nb\Lb}\ok)\nonumber
\end{eqnarray}
\begin{eqnarray}
&&L\skb=E\snb-\chi\skb+\snb\se+\onb\skb+\nb\sk+\skb\sf-2c(\okb f+\kb\of) \label{10.200}\\
&&\hspace{8mm}+(\sbeta\skb+\beta_{\Nb}\okb+\beta_N\kb)(\sbeta LH+\rho H\ss_N)
+H\sbeta(\rho\ss_N\skb+\rho s_{N\Nb}\okb+s_{NL}\kb)\nonumber
\end{eqnarray}
In view of \ref{10.177} and also of formulas \ref{3.a34} and their conjugates, equations \ref{10.199}, 
\ref{10.200} are seen to coincide with equations \ref{3.a46}, \ref{3.a47}. 

Applying $\Lb_N$ to \ref{10.178} we obtain:
\begin{eqnarray}
&&\Lb_N\sk_N=h_{\mu\nu,N}(\Lb_N E_N N_N^\mu)E_N^\nu+h_{\mu\nu,N}(E_N N_N^\mu)(\Lb_N E_N^\nu) \nonumber\\
&&\hspace{13mm}+(\Lb_N h_{\mu\nu,N})(E_N N_N^\mu)E_N^\nu
\label{10.201}
\end{eqnarray}
By the 2nd of the commutation relations \ref{9.b1} the 1st term on the right is:
\begin{equation}
h_{\mu\nu,N}(E_N \Lb_N N_N^\mu)E_N^\nu-\chib_N\sk_N
\label{10.202}
\end{equation}
Applying $E_N$ to the 1st of \ref{10.109} we obtain:
\begin{eqnarray}
&&E_N \Lb_N N_N^\mu=(E_N\sn_N)E_N^\mu+(E_N n_N)N_N^\mu+(E_N\on_N)\Nb_N^\mu \nonumber\\
&&\hspace{20mm}+\sn_N(\se_N E_N^\mu+e_N N_N^\mu+\oe_N\Nb_N^\mu)\nonumber\\
&&\hspace{20mm}+n_N(\sk_N E_N^\mu+k_N N_N^\mu+\ok_N\Nb_N^\mu)\nonumber\\
&&\hspace{20mm}+\on_N(\skb_N E_N^\mu+\kb_N N_N^\mu+\okb_N\Nb_N^\mu)
\label{10.203}
\end{eqnarray}
Consequently, \ref{10.202}, the 1st term on the right in \ref{10.201}, is:
\begin{equation}
E_N\sn_N+\sn_N\se_N+n_N\sk_N+\on_N\skb_N-\chib_N\sk_N
\label{10.204}
\end{equation}
Consider next the 2nd term on the right in \ref{10.201}. The expansion of $\Lb_N E_N^\mu$ in the 
$(E_N^\mu,N_N^\mu,\Nb_N^\mu)$ frame is the conjugate of \ref{10.184}:
\begin{equation}
\Lb_N E_N^\mu=\sfb_N E_N^\mu+\fb_N N_N^\mu+\ofb_N \Nb_N^\mu 
\label{10.205}
\end{equation}
Substituting the expansions 1st of \ref{10.117} and \ref{10.205} we find that the 2nd term on the right in 
\ref{10.201} is:
\begin{equation}
\sk_N\sfb_N-2c_N(\ok_N\fb_N+k_N\ofb_N)
\label{10.206}
\end{equation}
In regard finally to the 3rd term on the right in \ref{10.201}, substituting the expansion \ref{10.151} 
we obtain that this term is:
\begin{eqnarray}
&&(\Lb_N h_{\mu\nu,N})(\sk_N E_N^\mu+k_N N_N^\mu+\ok_N\Nb_N^\mu)E_N^\nu=\nonumber\\
&&\hspace{10mm}\sk_N\left(\sbeta_N^2 \Lb_N H_N+2\sbeta_N H_N E_N^\mu \Lb_N\beta_{\mu,N}\right)\nonumber\\
&&\hspace{7mm}+k_N\left(\sbeta_N\beta_{N,N}\Lb_N H_N+H_N(\sbeta_N s_{N\Lb,N}+\beta_{N,N}E_N^\mu \Lb_N\beta_{\mu,N})\right)\nonumber\\
&&\hspace{7mm}+\ok_N\left(\sbeta_N\beta_{\Nb,N}\Lb_N H_N+H_N(\sbeta_N s_{\Nb\Lb,N}
+\beta_{\Nb,N}E_N^\mu \Lb_N\beta_{\mu,N})\right)\nonumber\\
&&\label{10.207}
\end{eqnarray}
Using also the 3rd of \ref{9.44} to express $E_N^\mu \Lb_N\beta_{\mu,N}$ in terms of 
$\ss_{\Nb,N}$, we conclude that $\sk_N$ satisfies the following $N$th approximant version of 
the cross variation equation \ref{10.199}:
\begin{eqnarray}
&&\Lb_N\sk_N=E_N\sn_N-\chib_N\sk_N+\sn_N\se_N+n_N\sk_N+\on_N\skb_N\nonumber\\
&&\hspace{15mm}+\sk_N\sfb_N-2c_N(\ok_N\fb_N+k_N\ofb_N)\nonumber\\
&&\hspace{15mm}+(\sbeta_N\sk_N+\beta_{N,N}k_N+\beta_{\Nb,N}\ok_N)(\sbeta_N\Lb_N H_N+\rhob_N\ss_{\Nb,N})
\nonumber\\
&&\hspace{15mm}+H_N\sbeta_N(\rhob_N\ss_{\Nb,N}\sk_N+s_{N\Lb,N}k_N+s_{\Nb\Lb,N}\ok_N)+\vep_{\Lb\sk,N}
\nonumber\\
&&\label{10.208}
\end{eqnarray}
Here:
\begin{equation}
\vep_{\Lb\sk,N}=H_N(2\sbeta_N\sk_N+\beta_{N,N}k_N+\beta_{\Nb,N}\ok_N)
(\vepb_N^\mu E_N\beta_{\mu,N}-\omega_{\Lb E,N})=O(\tau^{N+1})
\label{10.209}
\end{equation}
the estimate following from Propositions 9.1 and 9.2. In a similar manner we arrive at the conclusion 
that $\skb_N$ satisfies the following $N$th approximant version of the cross variation equation 
\ref{10.200}:
\begin{eqnarray}
&&L_N\skb_N=E_N\snb_N-\chi_N\skb_N+\snb_N\se_N+\onb_N\skb_N+\nb_N\sk_N\nonumber\\
&&\hspace{15mm}+\skb_N\sf_N-2c_N(\kb_N\of_N+\okb_N f_N)\nonumber\\
&&\hspace{15mm}+(\sbeta_N\skb_N+\beta_{\Nb,N}\okb_N+\beta_{N,N}\kb_N)(\sbeta_N L_N H_N+\rho_N\ss_{N,N})
\nonumber\\
&&\hspace{15mm}+H_N\sbeta_N(\rho_N\ss_{N,N}\skb_N+s_{\Nb L,N}\okb_N+s_{NL,N}\kb_N)+\vep_{L\skb,N}
\nonumber\\
&&\label{10.210}
\end{eqnarray}
Here:
\begin{equation}
\vep_{L\skb,N}=H_N(2\sbeta_N\skb_N+\beta_{\Nb,N}\okb_N+\beta_{N,N}\kb_N)
(\vep_N^\mu E_N\beta_{\mu,N}-\omega_{L E,N})=O(\tau^N)
\label{10.211}
\end{equation}
the estimate following from Propositions 9.1 and 9.2. 

We proceed to derive analogues for the $N$th approximant quantities:
\begin{eqnarray}
&&\theta_N=\lambda_N\tchi_N+f_N \nonumber\\
&&\thetab_N=\lambdab_N\tchib_N+\fb_N \label{10.212}
\end{eqnarray}
where: 
\begin{eqnarray}
&&f_N=-\frac{1}{2}\beta_{N,N}^2\Lb_N H_N+\lambda_N\beta_{N,N}\sbeta_N E_N H_N \nonumber\\
&&\fb_N=-\frac{1}{2}\beta_{\Nb,N}^2 L_N H_N+\lambdab_N\beta_{\Nb,N}\sbeta_N E_N H_N \label{10.213}
\end{eqnarray}
of the regularized propagation equations of Proposition 10.1.

Recall that the argument leading to the regularized propagation equation for $\theta$ starts from 
the observation that the principal term on the right hand side of the 2nd variation equation 
\ref{3.a44}, equivalently \ref{10.175}, is the term $E\sm$. The corresponding term on the right hand 
side of the analogous equation for the $N$th approximants, equation \ref{10.187}, is the term 
$E_N\sm_N$. The argument leading to the propagation equation for $\theta$ of Proposition 10.1 
relies on Lemma 10.1 to express the principal term $E\sm$ as $Lv$ up to a 1st order remainder. 
Now, while $\sm$ is given by \ref{10.6}, $\sm_N$ is given by \ref{10.171}, an analogous formula  
for the $N$th approximants with a $O(\tau^N)$ error term. A similar argument then applies to 
express $E_N\sm_N$ as $L_N v_N$ up to a suitable remainder, where in analogy with the definition 
\ref{10.11} of $v$, 
\begin{equation}
v_N=\frac{1}{2\lambda_N}\beta_{N,N}^2\Lb_N H_N-\beta_{N,N}\sbeta_N E_N H_N-H_N\sbeta_N\ss_{N,N}
\label{10.214}
\end{equation}
However the analogue of Lemma 10.1 concerns $\square_{\tilde{h}^\prime_N}H_N$ therefore we 
must consider:
$$\square_{\tilde{h}_N^\prime}\sigma_N=
-2(g^{-1})^{\mu\nu}\beta_{\nu,N}\square_{\tilde{h}_N^\prime}\beta_{\mu,N}
-2(g^{-1})^{\mu\nu}\tilde{h}_N^{\prime-1}(d\beta_{\mu,N},d\beta_{\nu,N})$$
Since, according to Proposition 9.10:
$$\square_{\tilde{h}_N^\prime}\beta_{\mu,N}=\Omega_N^{-1}a_N^{-1}\tilde{\kappa}^\prime_{\mu,N}$$
this is:
$$\square_{\tilde{h}_N^\prime}\sigma_N=
-2(g^{-1})^{\mu\nu}\left(\beta_{\nu,N}\Omega_N^{-1}a_N^{-1}\tilde{\kappa}^\prime_{\mu,N}
+\tilde{h}_N^{\prime-1}(d\beta_{\mu,N},d\beta_{\nu,N})\right)$$
As a consequence, the analogue of Lemma 10.1 for the $N$th approximants is:
\begin{eqnarray}
&&\square_{h_N^\prime}H_N=\left(H^{\prime\prime}_N-(1/2)\Omega_N^{-1}\Omega_N^\prime H^\prime_N\right)
h_N^{\prime-1}(d\sigma_N,d\sigma_N) \label{10.215}\\
&&\hspace{16mm}-2H^\prime_N(g^{-1})^{\mu\nu}
\left(\beta_{\nu,N}a_N^{-1}\tilde{\kappa}^\prime_{\mu,N}+h_N^{\prime-1}(d\beta_{\mu,N},d\beta_{\nu,N})\right) \nonumber
\end{eqnarray}
where (see \ref{9.6}):
\begin{equation}
H^\prime_N=\left(\frac{dH}{d\sigma}\right)(\sigma_N), \ \ \ 
H^{\prime\prime}_N=\left(\frac{d^2 H}{d\sigma^2}\right)(\sigma_N), \ \ \ 
\Omega^\prime_N=\left(\frac{d\Omega}{d\sigma}\right)(\sigma_N)
\label{10.216}
\end{equation}
Then in the role of \ref{10.27} for the $N$th approximants we have: 
\begin{eqnarray}
&&a_N\square_{h_N^\prime}H_N=\left(H^{\prime\prime}_N-(1/2)\Omega_N^{-1}\Omega^\prime_N H^\prime_N\right)\left(-(L_N\sigma_N)\Lb_N\sigma_N+a_N(E_N\sigma_N)E_N\sigma_N\right)\nonumber\\
&&\hspace{20mm}-2H^\prime_N(g^{-1})^{\mu\nu}\left(-(L_N\beta_{\mu,N})\Lb_N\beta_{\nu,N}
+a_N(E_N\beta_{\mu,N})E_N\beta_{\nu,N}\right)\nonumber\\
&&\hspace{20mm}-2H^\prime_N(g^{-1})^{\mu\nu}\beta_{\nu,N}\tilde{\kappa}^\prime_{\mu,N} 
\label{10.217}
\end{eqnarray}
and, by \ref{9.210}, we have:
\begin{equation}
L_N\Lb_N H_N-a_N E_N^2 H_N=-a_N\square_{h_N^\prime}H_N-(1/2)(\chib_N L_N H_N+\chi_N\Lb_N H_N)
+2\etab_N E_N H_N
\label{10.218}
\end{equation}
in place of \ref{10.25}. As a consequence of the additional last term on the right in \ref{10.217}, 
$\lambda_N E_N\sm_N-L_N(\lambda_N v_N)$
has an additional term of the form: 
$$-\beta_{N,N}^2H^\prime_N (g^{-1})^{\mu\nu}\beta_{\nu,N}\tilde{\kappa}^\prime_{\mu,N}$$ 
compared to the result of Section 10.1 for $\lambda E\sm-L(\lambda v)$. We are thus led to the 
conclusion that $\theta_N$ satisfies a propagation equation of the form:
\begin{equation}
L_N\theta_N=R_N+\lambda_N\vep_{L\sk,N}-\beta_{N,N}^2 H^\prime_N (g^{-1})^{\mu\nu}\beta_{\nu,N}\tilde{\kappa}^\prime_{\mu,N}
\label{10.219}
\end{equation}
where 
\begin{equation}
R_N=[R_N]_{P.A.}+\tilde{R}_N
\label{10.220}
\end{equation}
with
\begin{equation}
[R_N]_{P.A.}=-a_N\tchi_N^2+\left(2L_N\lambda_N-\sbeta_N^2(L_N H_N)\lambda_N\right)\tchi_N
\label{10.221}
\end{equation}
and the remainder $\tilde{R}_N$ given, up to $O(\tau^N)$ error terms resulting from 
$\vep_N^\mu$, $\vepb_N^\mu$, $\omega_{L\Lb,N}$, $\omega_{LE,N}$, $\omega_{\Lb E,N}$ and $\delta_N^\prime$ 
(see Propositions 9.1, 9.2 and Lemma 9.2), by what corresponds to the 
expression \ref{10.40} for the $N$th approximants. The 2nd term on the right in \ref{10.219} 
is contributed by the error term in \ref{10.187}. By Proposition 9.10 the 
3rd term on the right in \ref{10.219} is  
$${\bf O}(\tau^{N-2})+{\bf O}(u^{-3}\tau^N)$$
By conjugation, we conclude that $\thetab_N$ satisfies a propagation equation of the form:
\begin{equation}
\Lb_N\thetab_N=\Rb_N+\lambdab_N\vep_{\Lb\skb,N}-\beta_{\Nb,N}^2 H^\prime_N (g^{-1})^{\mu\nu}\beta_{\nu,N}\tilde{\kappa}^\prime_{\mu,N}
\label{10.222}
\end{equation}
where 
\begin{equation}
\Rb_N=[\Rb_N]_{P.A.}+\tilde{\Rb}_N
\label{10.223}
\end{equation}
with
\begin{equation}
[\Rb_N]_{P.A.}=-a_N\tchib_N^2+\left(2\Lb_N\lambdab_N-\sbeta_N^2(\Lb_N H_N)\lambdab_N\right)\tchib_N
\label{10.224}
\end{equation}
and the remainder $\tilde{\Rb}_N$ given, up to $O(\tau^N)$ error terms resulting from 
$\vep_N^\mu$, $\vepb_N^\mu$, $\omega_{L\Lb,N}$, $\omega_{LE,N}$, $\omega_{\Lb E,N}$ and 
$\delta_N^\prime$ 
(see Propositions 9.1, 9.2 and Lemma 9.2), by what corresponds to the 
expression \ref{10.43} for the $N$th approximants. 

 We finally derive analogues for the $N$th approximant quantities:
 \begin{eqnarray}
 &&\nu_N=\lambda_N E_N^2\lambda_N-j_N \nonumber\\
 &&\nub_N=\lambdab_N E_N^2\lambdab_N-\jb_N \label{10.225} 
 \end{eqnarray}
 where:
\begin{eqnarray}
&&j_N=\frac{1}{4}\beta_{N,N}^2\Lb_N^2 H_N-\frac{1}{2}\lambda_N\pi_N\beta_{N,N}^2 E_N\Lb_N H_N \label{10.226}\\
&&\hspace{5mm}+\lambda_N^2\left\{\beta_{N,N}\left[\pi_N\sbeta_N-\frac{1}{4c_N}(\beta_{N,N}+2\beta_{\Nb,N})\right]E_N^2 H_N
\right. \nonumber\\
&&\hspace{15mm}\left.+H_N\left[\pi_N\sbeta_N-\frac{1}{2c_N}(\beta_{N,N}+\beta_{\Nb,N})\right]N_N^\mu E_N^2\beta_{\mu,N}
\right\} \nonumber
\end{eqnarray}
and:
\begin{eqnarray}
&&\jb_N=\frac{1}{4}\beta_{\Nb,N}^2 L_N^2 H_N-\frac{1}{2}\lambdab_N\pi_N\beta_{\Nb,N}^2 E_N L_N H_N 
\label{10.227}\\
&&\hspace{5mm}+\lambdab_N^2\left\{\beta_{\Nb,N}\left[\pi_N\sbeta_N-\frac{1}{4c_N}(2\beta_{N,N}+\beta_{\Nb,N})\right]
E_N^2 H_N\right. \nonumber\\
&&\hspace{15mm}\left.+H_N\left[\pi_N\sbeta_N-\frac{1}{2c_N}(\beta_{N,N}+\beta_{\Nb,N})\right]\Nb_N^\mu E_N^2\beta_{\mu,N}
\right\} \nonumber
\end{eqnarray}

Recall that the argument leading to the regularized propagation equation for $\nu$ focusds on the 
terms $\lambda E^2 p+\lambdab E^2 q$, the principal terms on the right in the propagation equation 
\ref{10.44} for $E^2\lambda$. The quantity $p$ being given by \ref{10.49}, we expressed $\lambda E^2 p$ 
as $Lw_1$ to principal terms, with $w_1$ defined by \ref{10.51}, the point being that $\lambda E^2\sm$ 
is by the preceding argument equal to $LE(\lambda v)$ to principal terms. When a similar argument 
is applied to the $N$th approximants to express $\lambda_N E_N^2\sm_N$ in terms of 
$L_N E_N(\lambda_N v_N)$, an additional term of the form 
$$\pi_N\beta_{N,N}^2 E_N\left(H^\prime_N (g^{-1})^{\mu\nu}\beta_{\nu,N}\tilde{\kappa}^\prime_{\mu,N}\right)$$
results. Also, the quantity $q$ being given by \ref{10.54}, we expressed $\lambda\lambdab E^2 q$ 
as $Lw_2$ to principal terms, with $w_2$ defined by \ref{10.58}. When a similar argument is applied 
to the $N$th approximants to express $\lambda_N\lambdab_N E_N^2 q_N$, an additional term of the form 
$$-(1/2)\beta_{N,N}^2\Lb_N\left(H^\prime_N (g^{-1})^{\mu\nu}\beta_{\nu,N}\tilde{\kappa}^\prime_{\mu,N}\right)$$ 
results. We are thus led to the conclusion that $\nu_N$ satisfies a propagation equation of the form:
\begin{eqnarray}
&&L_N\nu_N=K_N+\lambda_N E_N^2\vep_{\lambda,N} \label{10.228}\\
&&\hspace{15mm}+\beta_{N,N}^2\left(\pi_N\lambda_N E_N-(1/2)\Lb_N\right)\
\left(H^\prime_N(g^{-1})^{\mu\nu}\beta_{\nu,N}\tilde{\kappa}^\prime_{\mu,N}\right) \nonumber
\end{eqnarray}
where 
\begin{equation}
K_N=[K_N]_{P.A.}+\tilde{K}_N
\label{10.229}
\end{equation}
with 
\begin{eqnarray}
&&[K_N]_{P.A.}=2(L_N\lambda_N-\chi_N\lambda_N)E_N^2\lambda_N \nonumber\\
&&\hspace{17mm}+\frac{\lambda_N\beta_{N,N}^2}{4c_N}(\Lb_N H_N+2\pi_N\lambda_N E_N H_N)E_N^2\lambdab_N
\nonumber\\
&&\hspace{17mm}+\lambda_N A_N E_N\tchi_N+a_N \oA_N E_N\tchib_N 
\label{10.230}
\end{eqnarray}
the coefficients $A_N$, $\oA_N$ being analogues for the $N$th approximants of the coefficients 
$A$, $\oA$ given by \ref{10.99}, \ref{10.100}. The principal part of the remainder $\tilde{K}_N$ 
is given by what corresponds to the expression \ref{10.102} for the $N$th approximants. 
The 2nd term on the right in \ref{10.228} is contributed by the error term in \ref{10.161}. By 
Proposition 9.10 the 3rd term on the right in \ref{10.228} is 
$${\bf O}(\tau^{N-2})+{\bf O}(u^{-3}\tau^N)$$
By conjugation, we conclude that $\nub$ satisfies a propagation equation of the form:
\begin{eqnarray}
&&\Lb_N\nub_N=\Kb_N+\lambdab_N E_N^2\vep_{\lambdab,N} \label{10.231}\\
&&\hspace{15mm}+\beta_{\Nb,N}^2\left(\pi_N\lambdab_N E_N-(1/2))L_N\right)\
\left(H^\prime_N(g^{-1})^{\mu\nu}\beta_{\nu,N}\tilde{\kappa}^\prime_{\mu,N}\right) \nonumber
\end{eqnarray}
where 
\begin{equation}
\Kb_N=[\Kb_N]_{P.A.}+\tilde{\Kb}_N
\label{10.232}
\end{equation}
with 
\begin{eqnarray}
&&[\Kb_N]_{P.A.}=2(\Lb_N\lambdab_N-\chib_N\lambdab_N)E_N^2\lambdab_N \nonumber\\
&&\hspace{17mm}+\frac{\lambdab_N\beta_{\Nb,N}^2}{4c_N}(L_N H_N+2\pi_N\lambdab_N E_N H_N)E_N^2\lambda_N
\nonumber\\
&&\hspace{17mm}+\lambdab_N\oAb_N E_N\tchib_N+a_N \Ab_N E_N\tchi_N 
\label{10.233}
\end{eqnarray}
the coefficients $\oAb_N$, $\Ab_N$ being analogues for the $N$th approximants of the coefficients 
$\oAb$, $\Ab$ given by \ref{10.105}, \ref{10.106}. The principal part of the remainder $\tilde{\Kb}_N$ 
is given by what corresponds to the expression \ref{10.107} for the $N$th approximants. 
The 2nd term on the right in \ref{10.231} is contributed by the error term in \ref{10.164}. 
By Proposition 9.10 the 3rd term in \ref{10.231} is:
$${\bf O}(\tau^{N-3})+{\bf O}(u^{-3}\tau^{N-1})$$

\vspace{5mm}

\section{The Propagation Equations for $\cth_l$, $\cthb_l$ and for $\cnu_{m,l}$, $\cnub_{m,l}$}

To control the higher order quantities $E^l\tchi$, $E^l\tchib$, we introduce 
the higher order versions 
\begin{equation}
\theta_l=\lambda E^l\tchi+E^l f, \ \ \ \thetab_l=\lambdab E^l\tchib+E^l\fb 
\label{10.234}
\end{equation}
of the quantities $\theta$, $\thetab$ of Proposition 10.1. Setting $l=n$ these will 
control the top order quantities in the 1st line of \ref{10.4} ($m=0$). To control the higher order 
quantities $E^{l+1}T^m\lambda$, $E^{l+1}T^m\lambdab$, where $l\geq 1$, we introduce the higher order 
versions 
\begin{equation}
\nu_{m,l}=\lambda E^{l-1}T^m E^2\lambda-E^{l-1}T^m j, \ \ \ 
\nub_{m,l}=\lambdab E^{l-1}T^m E^2\lambdab-E^{l-1}T^m\jb
\label{10.235}
\end{equation}
of the quantities $\nu$, $\nub$ of Proposition 10.2. Setting $m+l=n$, the quantities 
$\nu_{m-1,l+1}$, $\nub_{m-1,l+1}$ control the top order quantities in the 2nd line of \ref{10.4}\\ 
($m\geq 1$). The higher order quantities $\theta_l$, $\thetab_l$ will be shown to satisfy higher 
order versions of the regularized propagation equations of Proposition 10.1. Similarly, the 
higher order quantities $\nu_{m,l}$, $\nub_{m,l}$ will be shown to satisfy higher order versions 
of the regularized propagation equations of Proposition 10.2. We proceed to define 
the analogous higher order quantities:
\begin{eqnarray}
&&\theta_{l,N}=\lambda_N E_N^l\tchi_N+E_N^l f_N, \nonumber\\
&&\thetab_{l,N}=\lambdab_N E_N^l\tchib_N+E_N^l\fb_N
\label{10.236}
\end{eqnarray}
and:
\begin{eqnarray}
&&\nu_{m,l,N}=\lambda_N E_N^{l-1}T^m E_N^2\lambda_N-E_N^{l-1}T^m j_N, \nonumber\\
&&\nub_{m,l,N}=\lambdab_N E_N^{l-1}T^m E_N^2\lambdab_N-E_N^{l-1}T^m\jb_N
\label{10.237}
\end{eqnarray}
for the $N$th approximants. We shall show that $\theta_{l,N}$, $\thetab_{l,N}$ satisfy higher order 
versions of the propagation equations \ref{10.219}, \ref{10.222} respectively, 
while $\nu_{m,l,N}$, $\nub_{m,l,N}$ satisfy higher order versions propagation equations \ref{10.228}, 
\ref{10.231} respectively. Defining finally the acoustical difference quantities:
\begin{equation}
\cth_l=\theta_l-\theta_{l,N}, \ \ \ \cthb_l=\thetab_l-\thetab_{l,N}
\label{10.238}
\end{equation}
and:
\begin{equation}
\cnu_{m,l}=\nu_{m,l}-\nu_{m,l,N}, \ \ \ \cnub_{m,l}=\nub_{m,l}-\nub_{m,l,N}
\label{10.239}
\end{equation}
we shall deduce the corresponding propagation equations, on the basis of which the appropriate 
estimates shall be derived. 

To derive the propagation equation for $\theta_l$ we must apply $E^l$ to the propagation equation 
for $\theta$ of Proposition 10.1. In view of the expression for $[R]_{P.A.}$, $E^l[R]_{P.A}$ is 
the sum of: 
\begin{equation}
E^l\left((2L\lambda-\sbeta^2 \lambda LH)\tchi\right)=
\left(2L\lambda-\sbeta^2\lambda LH\right)E^l\tchi+A_l
\label{10.240}
\end{equation}
and: 
\begin{equation}
E^l\left(-a\tchi^2\right)=E^{l-1}\left(-2a\tchi E\tchi-(Ea)\tchi^2\right)
=-2a\tchi E^l\tchi+B_l
\label{10.241}
\end{equation}
In \ref{10.240}, the remainder $A_l$ is of order $l+1$ but has vanishing principal acoustical part. 
In fact, since $L\lambda=p\lambda+q\lambdab$, the principal part of $A_l$ is:
\begin{equation}
[A_l]_{P.P.}=\left[2\lambdab E^l q+\lambda(2 E^l p-\sbeta^2 E^l LH)\right]_{P.P.}\tchi
\label{10.242}
\end{equation}
In \ref{10.241}, the remainder $B_l$ is of order $l$. To derive a formula for $L(E^l\theta)$ we must 
also consider the commutator $[L,E^l]\theta$. We apply the following general formula for linear 
operators $A$, $B$:
\begin{equation}
[A,B^n]=\sum_{m=0}^{n-1}B^m[A,B]B^{n-1-m}
\label{10.243}
\end{equation}
In view of the 1st of the commutation relations \ref{3.a14} we then obtain:
\begin{eqnarray}
&&[L,E^l]\theta=\sum_{k=0}^{l-1}E^k[L,E]E^{l-1-k}\theta \nonumber\\
&&\hspace{13mm}=-\sum_{k=0}^{l-1}E^k(\chi E^{l-k}\theta)=
-l\chi E^l\theta+C_l \label{10.244}
\end{eqnarray}
where the remainder $C_l$ is of order $l$. Taking also into account that by the 1st of \ref{3.a20}:
$$2a\tchi+\sbeta^2\lambda H=2\lambda\chi$$ 
we deduce the equation:
\begin{eqnarray}
&&L(E^l\theta)=\left(2L\lambda-2\lambda\chi\right)E^l\tchi-l\chi E^l\theta \nonumber\\
&&\hspace{18mm}+A_l+B_l+C_l+E^l\tilde{R} 
\label{10.245}
\end{eqnarray}
Consider now the difference:
\begin{equation}
D_l=E^l\theta-\theta_l=\sum_{k=1}^l\left(\begin{array}{c} l\\k\end{array}\right)(E^k\lambda)E^{l-k}\tchi
\label{10.246}
\end{equation}
This difference is of order $l$. We have:
\begin{equation}
LD_l=\sum_{k=1}^l\left(\begin{array}{c} l\\k\end{array}\right)\left\{(LE^k\lambda)E^{l-k}\tchi
+(E^k\lambda)LE^{l-k}\tchi\right\}
\label{10.247}
\end{equation}
which is of order $l+1$ but with vanishing principal acoustical part. Only the extreme terms $k=1$, 
$k=l$ in the sum in \ref{10.247} contribute to the principal part. Since 
$$[LE^l\lambda]_{P.P.}=[\lambda E^l p+\lambdab E^l q]_{P.P.}$$
and by equation \ref{10.175}, the 1st of \ref{3.a21} and \ref{10.6}:
$$[LE^{l-1}\tchi]_{P.P.}=[-\sbeta\beta_N E^l LH+(1/2)\rho\beta_N^2 E^{l+1}H]_{P.P.}$$
we obtain:
\begin{eqnarray}
&&[LD_l]_{P.P.}=\left[l(E\lambda)\left(-\sbeta\beta_N E^l LH+(1/2)\rho\beta_N^2 E^{l+1}H\right)\right. \nonumber\\
&&\hspace{20mm}\left.+\tchi\left(\lambda E^l p+\lambdab E^l q\right)\right]_{P.P.} 
\label{10.248}
\end{eqnarray}
Substituting 
$$E^l\theta=\theta_l+D_l, \ \ \ E^l\tchi=\frac{\theta_l-E^l f}{\lambda}$$
in equation \ref{10.245} we arrive at the following propagation equation for $\theta_l$:
\begin{equation}
L\theta_l=\left(\frac{2L\lambda}{\lambda}-(l+2)\chi\right)\theta_l
-\frac{2L\lambda}{\lambda}E^l f
+\tilde{R}_l
\label{10.249}
\end{equation}
where:
\begin{equation}
\tilde{R}_l=E^l\tilde{R}+2\chi E^l f+A_l+B_l+C_l-LD_l-l\chi D_l
\label{10.250}
\end{equation}
Note that $\tilde{R}_l$ is of order $l+1$ but with vanishing principal acoustical part. In fact, from 
\ref{10.40}, \ref{10.242}, \ref{10.248}, together with the expressions for $p$, $q$ of Proposition 3.3 
and the expression for $f$ of Proposition 10.1, we find:
\begin{eqnarray}
&&[\tilde{R}_l]_{P.P.}=\left(\lambdab (E\lambda)c^\mu_{1,E}+\lambda\lambdab\tchi 
c_{2,E}^\mu+\lambda\lambdab c_{3,E}^{\mu\nu}(E\beta_\nu) \right.\nonumber\\
&&\hspace{35mm}\left.+\lambda c_{4,E}^{\mu\nu}(L\beta_\nu)+\lambdab c_{5,E}^{\mu\nu}(\Lb\beta_\nu)\right)E^{l+1}\beta_\mu 
\nonumber\\
&&\hspace{15mm}+\left((E\lambda)c_{1,L}^\mu+\lambda\tchi c_{2,L}^\mu 
+\lambda c_{3,L}^{\mu\nu}(E\beta_\nu)+c_{4,L}^{\mu\nu}(\Lb\beta_\nu)\right)
E^l L\beta_\mu \nonumber\\
&&\hspace{15mm}+\left(\lambdab\tchi c_{1,\Lb}^\mu+\lambdab c_{2,\Lb}^{\mu\nu}(E\beta_\nu)
+c_{3,\Lb}^{\mu\nu}L\beta_\nu\right)E^l\Lb\beta_\mu 
\label{10.251}
\end{eqnarray}
where the coefficients are regular and of order 0. Similarly, we deduce the following propagation 
equation for $\thetab_l$:
\begin{equation}
\Lb\thetab_l=\left(\frac{2\Lb\lambdab}{\lambdab}-(l+2)\chib\right)\thetab_l
-\frac{2\Lb\lambdab}{\lambdab}E^l\fb
+\tilde{\Rb}_l
\label{10.252}
\end{equation}
where $\tilde{\Rb}_l$ is of order $l+1$ but with vanishing principal acoustical part. In fact, 
\begin{eqnarray}
&&[\tilde{\Rb}_l]_{P.P.}=\left(\lambda (E\lambdab)\underline{c}^\mu_{1,E}+\lambda\lambdab\tchib 
\underline{c}_{2,E}^\mu+\lambda\lambdab\underline{c}_{3,E}^{\mu\nu}(E\beta_\nu) \right.\nonumber\\
&&\hspace{35mm}\left.+\lambdab\underline{c}_{4,E}^{\mu\nu}(\Lb\beta_\nu)+\lambda\underline{c}_{5,E}^{\mu\nu}(L\beta_\nu)\right)E^{l+1}\beta_\mu 
\nonumber\\
&&\hspace{15mm}+\left((E\lambdab)\underline{c}_{1,\Lb}^\mu+\lambdab\tchib\underline{c}_{2,\Lb}^\mu 
+\lambdab\underline{c}_{3,\Lb}^{\mu\nu}(E\beta_\nu)+\underline{c}_{4,\Lb}^{\mu\nu}(L\beta_\nu)\right)
E^l\Lb\beta_\mu \nonumber\\
&&\hspace{15mm}+\left(\lambda\tchib\underline{c}_{1,L}^\mu+\lambda\underline{c}_{2,L}^{\mu\nu}(E\beta_\nu)
+\underline{c}_{3,L}^{\mu\nu}\Lb\beta_\nu\right)E^lL\beta_\mu 
\label{10.253}
\end{eqnarray}
where the coefficients are regular and of order 0. 

From the propagation equation \ref{10.219} for the $N$th approximant quantity $\theta_N$, 
following the same argument as above which leads from the propagation equation for $\theta$ 
of Proposition 10.1 to equation \ref{10.249}, we deduce, using the 1st of the commutation relations 
\ref{9.b1} in the role of the 1st of the commutation relations \ref{3.a14}, the following propagation 
equation for the $N$th approximant quantity $\theta_{l,N}$: 
\begin{eqnarray}
&&L_N\theta_{l,N}=\left(\frac{2L_N\lambda_N}{\lambda_N}-(l+2)\chi_N\right)\theta_{l,N}
-\frac{2L_N\lambda_N}{\lambda}E_N^l f_N
+\tilde{R}_{l,N} \nonumber\\
&&\hspace{20mm}+E_N^l\left(\lambda_N\vep_{L\sk,N}-\beta_{N,N}^2 H^\prime_N (g^{-1})^{\mu\nu}\beta_{\nu,N}\tilde{\kappa}^\prime_{\mu,N}\right)
\label{10.254}
\end{eqnarray}
the last term resulting from the two error terms on the right in \ref{10.219}. Here $\tilde{R}_{l,N}$ 
is the analogue of $\tilde{R}_l$ for the $N$ approximants. Similarly, we deduce from the propagation 
equation \ref{10.222} for $\thetab_N$ the following propagation equation for $\thetab_{l,N}$: 
\begin{eqnarray}
&&\Lb_N\thetab_{l,N}=\left(\frac{2\Lb_N\lambdab_N}{\lambdab_N}-(l+2)\chib_N\right)\thetab_{l,N}
-\frac{2\Lb_N\lambdab_N}{\lambdab}E_N^l\fb_N
+\tilde{\Rb}_{l,N} \nonumber\\
&&\hspace{20mm}+E_N^l\left(\lambdab_N\vep_{\Lb\skb,N}-\beta_{\Nb,N}^2 H^\prime_N (g^{-1})^{\mu\nu}\beta_{\nu,N}\tilde{\kappa}^\prime_{\mu,N}\right)
\label{10.255}
\end{eqnarray}
the last term resulting from the two error terms on the right in \ref{10.222}. Here $\tilde{\Rb}_{l,N}$ 
is the analogue of $\tilde{\Rb}_l$ for the $N$ approximants. 

Finally, subtracting equation \ref{10.254} from equation \ref{10.249} we arrive at the following 
propagation equation for the acoustical difference quantity $\cth_l$ defined by the 1st of \ref{10.238}: 
\begin{eqnarray}
&&L\cth_l-\left(\frac{2L\lambda}{\lambda}-(l+2)\chi\right)\cth_l= \nonumber\\
&&\hspace{25mm}
-\left\{L-L_N-\frac{2L\lambda}{\lambda}+\frac{2L_N\lambda_N}{\lambda_N}+(l+2)(\chi-\chi_N)\right\}
\theta_{l,N}\nonumber\\
&&\hspace{25mm}-\frac{(2L\lambda)}{\lambda}E^l f+\frac{(2L_N\lambda_N)}{\lambda_N}E_N^l f_N
+\tilde{R}_l-\tilde{R}_{l,N}\nonumber\\
&&\hspace{25mm}-E_N^l\vep_{\theta,N}
\label{10.256}
\end{eqnarray}
where we denote by $\vep_{\theta,N}$ the error term:
\begin{equation}
\vep_{\theta,N}=\lambda_N\vep_{L\sk,N}-\beta_{N,N}^2 H^\prime_N (g^{-1})^{\mu\nu}\beta_{\nu,N}\tilde{\kappa}^\prime_{\mu,N}
\label{10.a1}
\end{equation}
Also, subtracting equation \ref{10.255} from equation \ref{10.252} we arrive at the following 
propagation equation for the acoustical difference quantity $\cthb_l$ defined by the 2nd of \ref{10.238}: 
\begin{eqnarray}
&&\Lb\cthb_l-\left(\frac{2\Lb\lambdab}{\lambdab}-(l+2)\chib\right)\cthb_l= \nonumber\\
&&\hspace{25mm}
-\left\{\Lb-\Lb_N-\frac{2\Lb\lambdab}{\lambdab}+\frac{2\Lb_N\lambdab_N}{\lambdab_N}+(l+2)(\chib-\chib_N)\right\}
\thetab_{l,N}\nonumber\\
&&\hspace{25mm}-\frac{(2\Lb\lambdab)}{\lambdab}E^l\fb+\frac{(2\Lb_N\lambdab_N)}{\lambdab_N}E_N^l\fb_N
+\tilde{\Rb}_l-\tilde{\Rb}_{l,N}\nonumber\\
&&\hspace{25mm}-E_N^l\vep_{\thetab,N}
\label{10.257}
\end{eqnarray}
where we denote by $\vep_{\thetab,N}$ the error term:
\begin{equation}
\vep_{\thetab,N}=\lambdab_N\vep_{\Lb\skb,N}-\beta_{\Nb,N}^2 H^\prime_N (g^{-1})^{\mu\nu}\beta_{\nu,N}\tilde{\kappa}^\prime_{\mu,N}
\label{10.a2}
\end{equation}
We remark that the 1st term on the right in each of \ref{10.256}, \ref{10.257} is only of order 1, 
the quantities $\theta_{l,N}$, $\thetab_{l,N}$ being known smooth functions of the coordinates 
$(\ub,u,\vartheta)$. 

We turn to the quantities $\nu_{m,l}$, $\nub_{m,l}$ defined by \ref{10.235}. However before we are 
ready to derive propagation equations for these higher order quantities, we must first re-express 
the 2nd term in $[K]_{P.A.}$ as given in the statement of Proposition 10.2, namely the term: 
\begin{equation}
\frac{\lambda\beta_N^2}{4c}(\Lb H+2\pi\lambda EH)E^2\lambdab
\label{10.258}
\end{equation}
For, this term as it stands will generate a term in the propagation equation for $\nu_{m,l}$ 
the contribution of which cannot properly be estimated. We must similarly re-express the 2nd term in 
$[\Kb]_{P.A.}$ as given in the statement of Proposition 10.2, namely the term:
\begin{equation}
\frac{\lambdab\beta_{\Nb}^2}{4c}(LH+2\pi\lambdab EH)E^2\lambda
\label{10.259}
\end{equation}
Let us set:
\begin{equation}
q^\prime=\frac{\beta_N^2}{4c}(\Lb H+2\pi\lambda EH)
\label{10.260}
\end{equation}
Then the term \ref{10.258} is:
\begin{equation}
\lambda q^\prime E^2\lambdab
\label{10.261}
\end{equation}
Now, since $\lambdab=c\rho=cLt$ we have:
$$E\lambdab=cELt+\rho Ec$$
and, by the 1st of the commutation relations \ref{3.a14}, 
$$ELt=LEt+\chi Et=L\pi+\chi\pi$$
Therefore: 
\begin{equation}
E\lambdab=cL\pi+c\chi\pi+\rho Ec
\label{10.262}
\end{equation}
Substituting for $Ec$ from \ref{3.a23} in which $k$ and $\okb$ are substituted from \ref{3.46}, 
\ref{3.53}, \ref{3.54} and \ref{3.a21}, that is:
\begin{eqnarray}
&&k=-\pi\tchi-\frac{\beta_N^2}{4c}EH+\left(\pi\sbeta-\frac{\beta_N}{2c}\right)H\ss_N\nonumber\\
&&\okb=-\pi\tchib-\frac{\beta_{\Nb}^2}{4c}EH+\left(\pi\sbeta-\frac{\beta_{\Nb}}{2c}\right)H\ss_{\Nb}
\label{10.263}
\end{eqnarray}
to obtain:
\begin{eqnarray}
&&Ec=-c\pi(\tchi+\tchib)-\frac{1}{4}(\beta_N+\beta_{\Nb})^2 EH \nonumber\\
&&\hspace{10mm}+H\left(c\pi\sbeta-\frac{1}{2}(\beta_N+\beta_{\Nb})\right)(\ss_N+\ss_{\Nb}) 
\label{10.264}
\end{eqnarray}
yields:
\begin{equation}
E\lambdab=cL\pi-\lambdab\pi\tchib+F
\label{10.265}
\end{equation}
where 
\begin{eqnarray}
&&F=\frac{1}{2}c\pi\sbeta^2 LH
-\frac{\lambdab}{4c}(\beta_N+\beta_{\Nb})^2 EH\nonumber\\
&&\hspace{10mm}+\lambdab H\left(\pi\sbeta-\frac{1}{2c}(\beta_N+\beta_{\Nb})\right)(\ss_N+\ss_{\Nb})
\label{10.266}
\end{eqnarray}
is of order 1 but has vanishing principal acoustical part. Note that by \ref{3.a27} with $t$ in the 
role of $f$ and \ref{3.105} for $\mu=0$:
\begin{equation}
E\pi=\frac{1}{2c}(\tchi+\tchib)-\tilde{\gamma}^0
\label{10.267}
\end{equation}
where, according to \ref{3.101},
\begin{equation}
\tilde{\gamma}^0=\sgamma\pi+\tilde{\gamma}^N+\tilde{\gamma}^{\Nb}
\label{10.268}
\end{equation}
hence, substituting from \ref{3.a30}, 
\begin{eqnarray}
&&\tilde{\gamma}^0=\frac{1}{2}\sbeta\left(\pi\sbeta-\frac{1}{c}(\beta_N+\beta_{\Nb})\right)EH 
\nonumber\\
&&\hspace{7mm}+H\left(\pi\sbeta-\frac{1}{2c}(\beta_N+\beta_{\Nb})\right)\sss
\label{10.269}
\end{eqnarray}
In particular, since 
\begin{equation}
[\tilde{\gamma}^0]_{P.A.}=0
\label{10.270}
\end{equation}
we have:
\begin{equation}
[E\pi]_{P.A.}=\frac{1}{2c}(\tchi+\tchib)
\label{10.271}
\end{equation}
Applying $E$ to \ref{10.265} we obtain:
\begin{equation}
E^2\lambdab=cLE\pi-\lambdab\pi E\tchib+F_1
\label{10.272}
\end{equation}
where
\begin{equation}
F_1=EF+(Ec)L\pi + (c\chi-\lambdab\tchib)E\pi-(E\lambdab)\pi\tchi
\label{10.273}
\end{equation}
is of order 2 with vanishing P.A. part. In fact $[F_1]_{P.P.}=[EF]_{P.P.}$ is of the form:
\begin{equation}
[F_1]_{P.P.}=\lambdab c_0^\mu E^2\beta_\mu+c_1^\mu EL\beta_\mu 
\label{10.274}
\end{equation}
where the coefficients are regular and of order 0. In view of the expression \ref{10.272} for 
$E^2\lambdab$ the term \ref{10.258}, that is \ref{10.261}, takes the form:
\begin{equation}
L(\lambda cq^\prime E\pi)-\lambda\lambdab q^\prime\pi E\tchib+F_2
\label{10.275}
\end{equation}
where 
\begin{equation}
F_2=\lambda q^\prime F_1-L(\lambda cq^\prime)E\pi
\label{10.276}
\end{equation}
is of order 2 with vanishing P.A. part. 
Replacing the 2nd term in $[K]_{P.A.}$ as given in the statement of Proposition 10.2, namely 
the term \ref{10.258}, by its equal \ref{10.275}, we see that in view of \ref{10.100} and 
\ref{10.260} the terms in $E\tchib$ cancel and we obtain: 
\begin{equation}
[K]_{P.A.}=[K^\prime]_{P.A.}+F_2
\label{10.277}
\end{equation}
where:
\begin{equation}
[K^\prime]_{P.A.}=2(L\lambda-\chi\lambda)E^2\lambda+L(\lambda cq^\prime E\pi) +\lambda AE\tchi 
\label{10.278}
\end{equation}
Defining then:
\begin{equation}
\tilde{K}^\prime=\tilde{K}+F_2
\label{10.279}
\end{equation}
we have:
\begin{equation}
K=[K^\prime]_{P.A.}+\tilde{K}^\prime
\label{10.280}
\end{equation}
Moreover it is evident that $[\tilde{K}^\prime]_{P.P.}$ is of the same form as $[\tilde{K}]_{P.P.}$, 
given by \ref{10.102}. An analogous treatment of \ref{10.259}, the 2nd term in $[\Kb]_{P.A.}$ as given 
in the statement of Proposition 10.2, yields the conclusion that:
\begin{equation}
\Kb=[\Kb^\prime]_{P.A.}+\tilde{\Kb}^\prime
\label{10.281}
\end{equation}
where:
\begin{equation}
[\Kb^\prime]_{P.A.}=2(\Lb\lambdab-\chib\lambdab)E^2\lambdab+\Lb(\lambdab c\qb^\prime E\pi) +\lambdab  \oAb E\tchib 
\label{10.282}
\end{equation}
where 
\begin{equation}
\qb^\prime=\frac{\beta_{\Nb}^2}{4c}(LH+2\pi\lambdab EH)
\label{10.283}
\end{equation}
and $[\tilde{\Kb}^\prime]_{P.P.}$ is of the same form as $[\tilde{\Kb}]_{P.P.}$, given by 
\ref{10.107}. In the following we shall consider the propagation equations for $\nu$ and $\nub$ 
of Proposition 10.2 in their fully regularized form, that is with $K$ expressed by \ref{10.280}, 
$[K^\prime]_{P.A.}$ by \ref{10.278}, and $\Kb$ expressed by \ref{10.281}, $[\Kb^\prime]_{P.A.}$ by 
\ref{10.282}. 

To derive the propagation equation for $\nu_{m,l}$, we first apply $T^m$ to the 
fully regularized propagation equation for $\nu=\nu_{0,1}$ to deduce the propagation equation for 
$\nu_{m,1}$. From \ref{10.278} we obtain:
\begin{eqnarray}
&&T^m[K^\prime]_{P.A.}=2(L\lambda-\chi\lambda)T^m E^2\lambda+LT^m(\lambda cq^\prime E\pi) \nonumber\\
&&\hspace{15mm}+\lambda A T^m E\tchi+M_m 
\label{10.284}
\end{eqnarray}
where the remainder $M_m$ is of order $m+1$. In fact, since by the 3rd of the 
commutation relations \ref{3.a14} and the general formula \ref{10.243} we have:
\begin{equation}
[T^m,L]=\sum_{i=0}^{m-1}T^i\zeta ET^{m-1-i}
\label{10.285}
\end{equation}
$M_m$ is given by:
\begin{eqnarray}
&&M_m=\sum_{i=0}^{m-1}T^i(\zeta E(T^{m-1-i}(\lambda cq^\prime E\pi)))
+2\sum_{i=1}^m\left(\begin{array}{c} m\\i \end{array}\right)(T^i(L\lambda-\chi\lambda))
T^{m-i}E^2\lambda \nonumber\\
&&\hspace{20mm}+\sum_{i=1}^m\left(\begin{array}{c} m\\i \end{array}\right)T^i(\lambda A) T^{m-i}E\tchi
\label{10.286}
\end{eqnarray}
To derive a formula for $L(T^m\nu)$ we must also consider the commutator $[L,T^m]\nu$. By \ref{10.285} 
this is:
\begin{equation}
[L,T^m]\nu=-\sum_{i=0}^{m-1}T^i(\zeta ET^{m-1-i}\nu)=-m\zeta ET^{m-1}\nu+N_m
\label{10.287}
\end{equation}
where the remainder $N_m$ is of order $m+1$. We thus deduce the equation:
\begin{eqnarray}
&&L(T^m\nu)=2(L\lambda-\chi\lambda)T^m E^2\lambda+LT^m(\lambda cq^\prime E\pi) \label{10.288}\\
&&\hspace{15mm}+\lambda A T^m E\tchi-m\zeta ET^{m-1}\nu+M_m+N_m+T^m\tilde{K}^\prime \nonumber
\end{eqnarray}
Consider the difference:
\begin{equation}
G_m=T^m\nu-\nu_{m,1}=\sum_{i=1}^m\left(\begin{array}{c} m\\i \end{array}\right)
(T^i\lambda)T^{m-i}E^2\lambda
\label{10.289}
\end{equation}
This difference is of order $m+1$. We have:
\begin{equation}
LG_m=\sum_{i=1}^m\left(\begin{array}{c} m\\i \end{array}\right)\left\{(LT^i\lambda)T^{m-i}E^2\lambda
+(T^i\lambda)LT^{m-i}E^2\lambda\right\}
\label{10.290}
\end{equation}
which is of order $m+2$ but with vanishing principal acoustical part. Only the term $i=1$ contributes 
to the principal part:
\begin{equation}
[LG_m]_{P.P.}=m(T\lambda)[LT^{m-1}E^2\lambda]_{P.P.}=m(T\lambda)\left[\lambda E^2 T^{m-1}p
+\lambdab E^2 T^{m-1}q\right]_{P.P.}
\label{10.291}
\end{equation}
Substituting 
$$L(T^m\nu)=L\nu_{m,1}+LG_m$$
in equation \ref{10.288} we arrive at the following propagation equation for $\nu_{m,1}$: 
\begin{eqnarray}
&&L\nu_{m,1}=2(L\lambda-\chi\lambda)T^m E^2\lambda+LT^m(\lambda cq^\prime E\pi) \label{10.292}\\
&&\hspace{10mm}+\lambda A T^m E\tchi-m\zeta ET^{m-1}\nu+M_m+N_m-LG_m+T^m\tilde{K}^\prime \nonumber
\end{eqnarray}
Next, we apply $E^{l-1}$ to this equation. Using 
\begin{equation}
[E^{l-1},L]=\sum_{k=0}^{l-2}E^k\chi E^{l-1-k}
\label{10.293}
\end{equation}
as follows from the 1st of the commutation relations \ref{3.a14} and the general formula \ref{10.243}, 
to re-express $E^{l-1}$ applied to the 2nd term on the right in \ref{10.292}, we obtain:
\begin{eqnarray}
&&E^{l-1}L\nu_{m,1}=2(L\lambda-\chi\lambda)E^{l-1}T^m E^2\lambda
+L(E^{l-1}T^m(\lambda cq^\prime E\pi)) \nonumber\\
&&\hspace{20mm}+\lambda A E^{l-1}T^m E\tchi-m\zeta E^l T^{m-1}\nu+U_{m,l} \nonumber\\
&&\hspace{20mm}+E^{l-1}(M_m+N_m-LG_m+T^m\tilde{K}^\prime)
\label{10.294}
\end{eqnarray}
where:
\begin{eqnarray}
&&U_{m,l}=2\sum_{k=1}^{l-1}\left(\begin{array}{c} l-1\\k \end{array}\right) E^k(L\lambda-\chi\lambda)
E^{l-1-k}T^m E^2\lambda \nonumber\\
&&\hspace{10mm}+\sum_{k=0}^{l-2}E^k(\chi E^{l-1-k}T^m(\lambda cq^\prime E\pi)) \nonumber\\
&&\hspace{10mm}+\sum_{k=1}^{l-1}\left(\begin{array}{c} l-1\\k \end{array}\right)E^k(\lambda A)
E^{l-1-k}T^m E\tchi \nonumber\\
&&\hspace{10mm}-m\sum_{k=1}^{l-1}\left(\begin{array}{c} l-1\\k \end{array}\right)(E^k\zeta)
E^{l-k}T^{m-1}\nu 
\label{10.295}
\end{eqnarray}
is of order $m+l$. The term $E^{l-1}(M_m+N_m)$ is also of order $m+l$, while the remaining terms on 
the right in \ref{10.294} are of top order $m+l+1$. Consider the difference:
\begin{equation}
F_{m,l}=E^{l-1}\nu_{m,1}-\nu_{m,l}=\sum_{k=1}^{l-1}\left(\begin{array}{c} l-1\\k \end{array}\right)
(E^k\lambda)E^{l-1-k}T^m E^2\lambda
\label{10.296}
\end{equation}
This difference is of order $m+l$. We have:
\begin{equation}
LF_{m,l}=\sum_{k=1}^{l-1}\left(\begin{array}{c} l-1\\k \end{array}\right)
\left\{(LE^k\lambda)E^{l-1-k}T^m E^2\lambda+(E^k\lambda)LE^{l-1-k}T^m E^2\lambda\right\}
\label{10.297}
\end{equation}
which is of order $m+l+1$ but with vanishing principal acoustical part. Only the term $k=1$ 
contributes to the principal part:
\begin{equation}
[LF_{m,l}]_{P.P.}=(l-1)(E\lambda)[\lambda E^l T^m p+\lambdab E^l T^m q]_{P.P.}
\label{10.298}
\end{equation}
Now, by \ref{10.293} we have:
\begin{eqnarray}
&&E^{l-1}L\nu_{m,1}=LE^{l-1}\nu_{m,1}+\sum_{k=0}^{l-2}E^k(\chi E^{l-1-k}\nu_{m,1}) \nonumber\\
&&\hspace{15mm}=LE^{l-1}\nu_{m,1}+(l-1)\chi E^{l-1}\nu_{m,1}+V_{m,l}
\label{10.299}
\end{eqnarray}
where $V_{m,l}$ is of order $m+l$. Substituting this as well as 
$$E^{l-1}\nu_{m,1}=\nu_{m,l}+F_{m,l}$$
on the left in \ref{10.294} we obtain the following propagation equation for $\nu_{m,l}$:
\begin{eqnarray}
&&L\nu_{m,l}=\left(2\frac{L\lambda}{\lambda}-(l+1)\chi\right)\nu_{m,l}
+2\frac{L\lambda}{\lambda}E^{l-1}T^m j+L(E^{l-1}T^m(\lambda cq^\prime E\pi))\nonumber\\
&&\hspace{15mm}+\lambda A E^{l-1}T^m E\tchi-m\zeta E^l T^{m-1}\nu+\tilde{K}^\prime_{m,l}
\label{10.300}
\end{eqnarray}
where:
\begin{eqnarray}
&&\tilde{K}^\prime_{m,l}=E^{l-1}T^m\tilde{K}^\prime-E^{l-1}LG_m-LF_{m,l}
-2\chi E^{l-1}T^m j\nonumber\\
&&\hspace{11mm}+E^{l-1}(M_m+N_m)+U_{m,l}-V_{m,l}-(l-1)\chi F_{m,l}
\label{10.301}
\end{eqnarray}
(the quantities in the 1st line on the right are of order $m+l+1$, those in the 2nd of order $m+l$).

We shall finally re-express the quantities $E^{l-1}T^m E\tchi$ and $E^l T^{m-1}\nu$ in the 4th and 5th 
terms on the right in \ref{10.300} in terms of $\nu_{m-i,l+i} \ : \ i=1,...,m$ and $\theta_{m,l}$. 
Consider first the quantity $E^l T^{m-1}\nu$. We have, 
$$T^{m-1}\nu=T^{m-1}(\lambda E^2\lambda-j)=\nu_{m-1,1}
+\sum_{i=1}^{m-1}\left(\begin{array}{c} m-1\\i \end{array}\right)(T^i\lambda)T^{m-1-i}E^2\lambda$$
Also, 
$$E^l\nu_{m-1,1}=\nu_{m-1,l+1}+\sum_{k=1}^l\left(\begin{array}{c} l\\k \end{array}\right) (E^k\lambda)
E^{l-k} T^{m-1} E^2 \lambda$$
Hence:
\begin{equation}
E^l T^{m-1}\nu=\nu_{m-1,l+1}+W_{m,l}
\label{10.302}
\end{equation}
where
\begin{eqnarray}
&&W_{m,l}=E^l\left(\sum_{i=1}^{m-1}
\left(\begin{array}{c} m-1\\i \end{array}\right)(T^i\lambda)T^{m-1-i}E^2\lambda\right) \nonumber\\
&&\hspace{10mm}+\sum_{k=1}^l\left(\begin{array}{c} l\\k \end{array}\right) 
(E^k\lambda) E^{l-k} T^{m-1} E^2 \lambda 
\label{10.303}
\end{eqnarray}
is of order $m+l$. 

Consider next the quantity $E^{l-1}T^m E\tchi$. We shall first express $T^n E\tchi$ for any 
positive integer $n$. Consider $T\tchi$. In view of the 1st of \ref{3.a21} 
we can express $T\tchi$ in terms of $T\sk=L\sk+\Lb\sk$, and $L\sk$ is expressed by the 2nd variation 
equation \ref{10.175}, $\Lb\sk$ by the cross variation equation \ref{10.199} as:
\begin{equation}
L\sk=E\sm+I, \ \ \ \Lb\sk=E\sn+J
\label{10.304}
\end{equation}
where $I$ and $J$ are quantities of order 1. Recall that $\sm$ is of order 1 with vanishing P.A. 
part, while from the 1st of \ref{3.a24}, substituting 
\begin{equation}
k=-\pi\tchi-\frac{\beta_N^2}{4c}EH+\left(\pi\sbeta-\frac{\beta_N}{2c}\right)H\ss_N
\label{10.305}
\end{equation}
(see \ref{10.245}), we obtain:
\begin{equation}
\sn=2E\lambda+2\pi\lambda\tchi+\tilde{\sn}
\label{10.306}
\end{equation}
where 
\begin{equation}
\tilde{\sn}=\lambda\left[\frac{\beta_N}{2c}(\beta_N+2\beta_{\Nb})EH
-2\left(\pi\sbeta-\frac{\beta_N}{2c}\right)H\ss_N\right]-\beta_N\sbeta\Lb H
\label{10.307}
\end{equation}
is of order 1 with vanishing P.A. part. Consequently:
\begin{eqnarray}
&&T\tchi=2E^2\lambda+2\pi\lambda E\tchi+E\sm+E\tilde{\sn}+T(\sbeta H\ss_N) \nonumber\\
&&\hspace{10mm}+\tchi E(2\pi\lambda)+I+J
\label{10.308}
\end{eqnarray}
Here the first two terms on the right constitute the principal acoustical part 
(compare with \ref{8.137}), 
the next three terms are of order 2 but with vanishing P.A. parts, while the last three terms 
are only of order 1. Since by the first two of the commutation relations \ref{3.a14} we have:
$$TE\tchi=ET\tchi-(\chi+\chib)E\tchi$$
\ref{10.308} implies:
\begin{equation}
TE\tchi=2E^3\lambda+2\pi\lambda E^2\tchi+S
\label{10.309}
\end{equation}
where $S$ is of order 3 but with vanishing P.A. part:
\begin{equation}
S=E^2\sm+E^2\tilde{\sn}+ET(\sbeta H\ss_N)+\left(2E(2\pi\lambda)-(\chi+\chib)\right)E\tchi
+\tchi E^2(2\pi\lambda)+EI+EJ
\label{10.310}
\end{equation}
The first three terms on the right contain the principal part of $S$. We apply $T^n$ to \ref{10.309} 
to obtain: 
\begin{equation}
T^{n+1}E\tchi=2\pi\lambda ET^n E\tchi+2ET^n E^2\lambda+S_n
\label{10.311}
\end{equation}
where:
\begin{eqnarray}
&&S_n=T^n S+\sum_{i=1}^n\left(\begin{array}{c} n\\i \end{array}\right)T^i(2\pi\lambda)T^{n-i}E^2\tchi
\nonumber\\
&&\hspace{10mm}-2\pi\lambda \sum_{i=0}^{n-1}T^i\left((\chi+\chib) ET^{n-1-i}E\tchi\right) \nonumber\\
&&\hspace{10mm}-2\sum_{i=0}^{n-1}T^i\left((\chi+\chib)ET^{n-1-i}E^2\lambda\right)
\label{10.312}
\end{eqnarray}
The last two sums represent $[E,T^n]E\tchi$ and $[E,T^n]E^2\lambda$, the commutator $[E,T^n]$ 
being expressed using the general formula \ref{10.243} and the first two of the commutation relations 
\ref{3.a14}. The first term on the right in \ref{10.312} is of order $n+3$ while each of the three 
sums is of order $n+2$. Now, equation \ref{10.311} represents a linear recursion to which Lemma 9.5 applies 
taking 
$$x_n=T^n E\tchi, \ \ \ A=2\pi\lambda E, \ \ \ y_n=2ET^n E^2\lambda+S_n$$
The conclusion is:
\begin{equation}
T^n E\tchi=(2\pi\lambda E)^n E\tchi+\sum_{i=0}^{n-1}(2\pi\lambda E)^{n-1-i}
\left(2ET^i E^2\lambda+S_i\right)
\label{10.313}
\end{equation}

Setting then $n=m$ in \ref{10.313} and applying $E^{l-1}$ gives:
\begin{eqnarray}
&&E^{l-1}T^m E\tchi=E^{l-1}(2\pi\lambda E)^m E\tchi+2\sum_{i=0}^{m-1}E^{l-1}(2\pi\lambda E)^{m-1-i}
ET^i E^2\lambda\nonumber\\
&&\hspace{25mm}+\sum_{i=0}^{m-1}E^{l-1}(2\pi\lambda E)^{m-1-i}S_i 
\label{10.314}
\end{eqnarray}
Now, with $n$ an arbitrary positive integer the operator
\begin{equation}
O_{n-1}=(2\pi\lambda E)^n-(2\pi\lambda)^n E^n
\label{10.315}
\end{equation}
is of order $n-1$ and vanishes for $n=1$. In fact, the principal part of $O_{n-1}$ is the operator:
\begin{equation}
[O_{n-1}]_{P.P.}=n(2\pi\lambda)^{n-1}(E(2\pi\lambda))E^{n-1}
\label{10.316}
\end{equation}
We can then write \ref{10.314} in the form:
\begin{equation}
E^{l-1}T^m E\tchi=(2\pi\lambda)^m E^{l+m}\tchi
+2\sum_{i=0}^{m-1}(2\pi\lambda)^{m-1-i}E^{l-1+m-i}T^i E^2\lambda+S_{m,l}
\label{10.317}
\end{equation}
where:
\begin{eqnarray}
&&S_{m,l}=\sum_{i=0}^{m-1}E^{l-1}(2\pi\lambda E)^{m-1-i}S_i \label{10.318}\\
&&\hspace{10mm}+E^{l-1}O_{m-1}E\tchi+2\sum_{i=0}^{m-2}E^{l-1}O_{m-2-i}ET^iE^2\lambda \nonumber\\
&&\hspace{10mm}+\sum_{k=1}^{l-1}\left(\begin{array}{c} l-1\\k \end{array}\right)
E^k\left((2\pi\lambda)^m\right)E^{l-k+m}\tchi \nonumber\\
&&\hspace{10mm}+2\sum_{i=0}^{m-1}\sum_{k=1}^{l-1}\left(\begin{array}{c} l-1\\k \end{array}\right)
E^k\left((2\pi\lambda)^{m-1-i}\right)E^{l-1-k+m-i}T^i E^2\lambda \nonumber
\end{eqnarray}
Here, each of the terms in the first sum is of top order, $m+l+1$, while the remaining sums are all 
of order $m+l$. Finally, in view of the definitions, 1st of \ref{10.234}, \ref{10.235}, we can express 
\ref{10.317} in the form:
\begin{eqnarray}
&&\lambda E^{l-1}T^m E\tchi=(2\pi\lambda)^m\theta_{l+m}
+2\sum_{i=0}^{m-1}(2\pi\lambda)^{m-1-i}\nu_{i,l+m-i} \label{10.319}\\
&&\hspace{24mm}+\lambda S_{m,l}-(2\pi\lambda)^m E^{l+m}f+2\sum_{i=0}^{m-1}
(2\pi\lambda)^{m-1-i}E^{l-1+m-i}T^i j
\nonumber
\end{eqnarray}

We substitute for the 4th and 5th terms on the right in \ref{10.300} from \ref{10.319} and \ref{10.302}
to bring the propagation equation for $\nu_{m,l}$ to its final form:
\begin{eqnarray}
&&L\nu_{m,l}=\left(2\frac{L\lambda}{\lambda}-(l+1)\chi\right)\nu_{m,l}
+2\frac{L\lambda}{\lambda}E^{l-1}T^m j+L(E^{l-1}T^m(\lambda cq^\prime E\pi))\nonumber\\
&&\hspace{15mm}+A\left\{(2\pi\lambda)^m\theta_{l+m}
+2\sum_{i=0}^{m-1}(2\pi\lambda)^{m-1-i}\nu_{i,l+m-i}\right\}-m\zeta\nu_{m-1,l+1}\nonumber\\
&&\hspace{15mm}+\tilde{K}^{\prime\prime}_{m,l}
\label{10.320}
\end{eqnarray}
where:
\begin{eqnarray}
&&\tilde{K}^{\prime\prime}_{m,l}=\tilde{K}^\prime_{m,l}-m\zeta W_{m,l} \label{10.321}\\
&&\hspace{10mm}+A\left\{\lambda S_{m,l}-(2\pi\lambda)^m E^{l+m}f
+2\sum_{i=0}^{m-1}(2\pi\lambda)^{m-1-i}E^{l-1+m-i}T^i j\right\}\nonumber
\end{eqnarray}

A straightforward calculation shows that the principal part of $\tilde{K}^{\prime\prime}_{m,l}$ is 
of the form:
\begin{eqnarray}
&&[\tilde{K}^{\prime\prime}_{m,l}]_{P.P.}=\left[\lambdab(\Lb\lambda)d^\mu_{1,EE}
+\lambdab\lambda(E\lambda)d^\mu_{2,EE}+\lambda^2(E\lambdab) d^\mu_{3,EE} \right. \nonumber\\
&&\hspace{27mm}+\lambdab\lambda^2\tchi d^\mu_{4,EE}+\lambdab\lambda^2\tchib d^\mu_{5,EE} \nonumber\\
&&\hspace{23mm}\left.+\lambdab\lambda d^{\mu\nu}_{6,EE}\Lb\beta_\nu+\lambda^2 d^{\mu\nu}_{7,EE}L\beta_\nu
+\lambdab\lambda^2 d^{\mu\nu}_{8,EE}E\beta_\nu\right]E^{l+1}T^m\beta_\mu \nonumber\\
&&\hspace{20mm}+\left[\lambda(E\lambda)d^\mu_{1,EL}+\lambda^2\tchi d^\mu_{2,EL}
+\lambda^2\tchib d^\mu_{3,EL}\right. \nonumber\\
&&\hspace{25mm}\left.+\lambda d^{\mu\nu}_{4,EL}\Lb\beta_\nu+\lambda^2 d^{\mu\nu}_{5,EL}E\beta_\nu
\right]E^l T^m L\beta_\mu\nonumber\\
&&\hspace{20mm}+\left[\lambda(E\lambdab)d^\mu_{1,E\Lb}+\lambdab(E\lambda)d^\mu_{2,E\Lb}
+\lambda\lambdab\tchi d^\mu_{3,E\Lb}+\lambda\lambdab\tchib d^\mu_{4,E,\Lb}\right.\nonumber\\
&&\hspace{23mm}\left.+\lambdab d^{\mu\nu}_{5,E\Lb}\Lb\beta_\nu+\lambda d^{\mu\nu}_{6,E\Lb}L\beta_\nu
+\lambda\lambdab d^{\mu\nu}_{7,E\Lb}E\beta_\nu\right]E^l T^m\Lb\beta_\mu\nonumber\\
&&\hspace{20mm}+\left[\lambdab\tchi d^\mu_{1,\Lb\Lb}+d^{\mu\nu}_{2,\Lb\Lb}L\beta_\nu
+\lambdab d^{\mu\nu}_{3,\Lb\Lb}\right]E^{l-1}T^m\Lb^2\beta_\mu\nonumber\\
&&+(T\lambda)\left[\lambda\lambdab g^\mu_{EE}E^{l+2}T^{m-1}\beta_\mu
+\lambda g^\mu_{EL}E^{l+1}T^{m-1}L\beta_\mu+\lambdab g^\mu_{E\Lb}E^{l+1}T^{m-1}\Lb\beta_\mu\right]
\nonumber\\
&&+A\sum_{i=0}^{m-1}(2\pi\lambda)^{m-1-i}\left\{
\left(\lambda^2 e^\mu_{1,EE}+\lambda\lambdab  e^\mu_{2,EE}\right)E^{l+1+m-i}T^i\beta_\mu
+\lambda e^\mu_{EL}E^{l+m-i}T^i L\beta_\mu \right.\nonumber\\
&&\hspace{40mm}\left.+\lambda e^\mu_{E\Lb}E^{l+m-i}T^i\Lb\beta_\mu 
+e^\mu_{\Lb\Lb}E^{l-1+m-i}T^i\Lb^2\beta_\mu\right\} \nonumber\\
&&\label{10.322}
\end{eqnarray}
Here the term proportional to $T\lambda$ represents $-[E^{l-1}LG_m]_{P.P.}$ (see \ref{10.301}), 
while the sum represents the principal part of the terms in parenthesis in \ref{10.321}. 
The coefficients $d$, $g$, and $e$, are all regular and of order 0. 

Similarly, we deduce the conjugate propagation equation for $\nub_{m,l}$:
\begin{eqnarray}
&&\Lb\nub_{m,l}=\left(2\frac{\Lb\lambdab}{\lambdab}-(l+1)\chib\right)\nub_{m,l}
+2\frac{\Lb\lambdab}{\lambdab}E^{l-1}T^m\jb+\Lb(E^{l-1}T^m(\lambdab c\qb^\prime E\pi))\nonumber\\
&&\hspace{15mm}+\oAb\left\{(2\pi\lambdab)^m\thetab_{l+m}
+2\sum_{i=0}^{m-1}(2\pi\lambdab)^{m-1-i}\nub_{i,l+m-i}\right\}+m\zeta\nub_{m-1,l+1}\nonumber\\
&&\hspace{15mm}+\tilde{\Kb}^{\prime\prime}_{m,l}
\label{10.323}
\end{eqnarray}
the principal part of $\tilde{\Kb}^{\prime\prime}_{m,l}$ being given by:
\begin{eqnarray}
&&[\tilde{\Kb}^{\prime\prime}_{m,l}]_{P.P.}=\left[\lambda(L\lambdab)\db^\mu_{1,EE}
+\lambda\lambdab(E\lambdab)\db^\mu_{2,EE}+\lambdab^2(E\lambda)\db^\mu_{3,EE} \right. \nonumber\\
&&\hspace{27mm}+\lambda\lambdab^2\tchib\db^\mu_{4,EE}+\lambda\lambdab^2\tchi\db^\mu_{5,EE} \nonumber\\
&&\hspace{23mm}\left.+\lambda\lambdab\db^{\mu\nu}_{6,EE}L\beta_\nu+\lambdab^2\db^{\mu\nu}_{7,EE}\Lb\beta_\nu
+\lambda\lambdab^2\db^{\mu\nu}_{8,EE}E\beta_\nu\right]E^{l+1}T^m\beta_\mu \nonumber\\
&&\hspace{20mm}+\left[\lambdab(E\lambdab)\db^\mu_{1,E\Lb}+\lambdab^2\tchib\db^\mu_{2,E\Lb}
+\lambdab^2\tchi\db^\mu_{3,E\Lb}\right. \nonumber\\
&&\hspace{25mm}\left.+\lambdab\db^{\mu\nu}_{4,E\Lb}L\beta_\nu+\lambdab^2\db^{\mu\nu}_{5,E\Lb}E\beta_\nu
\right]E^l T^m\Lb\beta_\mu\nonumber\\
&&\hspace{20mm}+\left[\lambdab(E\lambda)\db^\mu_{1,EL}+\lambda(E\lambdab)\db^\mu_{2,EL}
+\lambdab\lambda\tchib\db^\mu_{3,EL}+\lambdab\lambda\tchi\db^\mu_{4,E,L}\right.\nonumber\\
&&\hspace{23mm}\left.+\lambda\db^{\mu\nu}_{5,EL}L\beta_\nu+\lambdab\db^{\mu\nu}_{6,EL}\Lb\beta_\nu
+\lambdab\lambda\db^{\mu\nu}_{7,EL}E\beta_\nu\right]E^l T^m L\beta_\mu\nonumber\\
&&\hspace{20mm}+\left[\lambda\tchib\db^\mu_{1,LL}+\db^{\mu\nu}_{2,LL}\Lb\beta_\nu
+\lambda\db^{\mu\nu}_{3,LL}\right]E^{l-1}T^m L^2\beta_\mu\nonumber\\
&&+(T\lambdab)\left[\lambdab\lambda\gb^\mu_{EE}E^{l+2}T^{m-1}\beta_\mu
+\lambdab\gb^\mu_{E\Lb}E^{l+1}T^{m-1}\Lb\beta_\mu+\lambda\gb^\mu_{EL}E^{l+1}T^{m-1}L\beta_\mu\right]
\nonumber\\
&&+\oAb\sum_{i=0}^{m-1}(2\pi\lambdab)^{m-1-i}\left\{
\left(\lambdab^2\eb^\mu_{1,EE}+\lambdab\lambda\eb^\mu_{2,EE}\right)E^{l+1+m-i}T^i\beta_\mu
+\lambdab\eb^\mu_{E\Lb}E^{l+m-i}T^i\Lb\beta_\mu \right.\nonumber\\
&&\hspace{40mm}\left.+\lambdab\eb^\mu_{EL}E^{l+m-i}T^i L\beta_\mu 
+\eb^\mu_{LL}E^{l-1+m-i}T^i L^2\beta_\mu\right\} \nonumber\\
&&\label{10.324}
\end{eqnarray}
the coefficients $\db$, $\gb$, and $\eb$ being all regular and of order 0. 

To derive from the propagation equation \ref{10.228} for $\nu_N$ a propagation equation for the 
higher order $N$th approximant quantity $\nu_{m,l,N}$, we follow the same argument as above which 
leads from the propagation equation for $\nu$ of Proposition 10.2 to the propagation equation 
\ref{10.320} for $\nu_{m,l}$. Thus, we first re-express the 2nd term in $[K_n]_{P.A.}$ as given 
by \ref{10.230}, namely the term  
\begin{equation}
\lambda_N q_N^\prime E_N^2\lambdab_N, \ \ \mbox{where} \ \ 
q_N^\prime=\frac{\beta_{N,N}^2}{4c_N}(\Lb_N H_N+2\pi_N\lambda_N E_N H_N)
\label{10.325}
\end{equation}
using the analogue of \ref{10.265} for the $N$th approximants:
\begin{equation}
E_N\lambdab_N=c_N L_N\pi_N-\lambdab_N\pi_N\tchib_N+F_N
\label{10.326}
\end{equation}
$F_N$ being given by the analogue of \ref{10.266} for the $N$th approximants. To establish 
\ref{10.326} we substitute formulas \ref{10.145}, \ref{10.147} for $k_N$, $\okb_N$, which 
are analogous to formulas \ref{10.263} for $k$, $\okb$, into equation \ref{10.148} for $E_N c_N$, 
which is analogous to equation \ref{3.a23} for $Ec$. A decomposition for $K_N$ analogous to the 
decomposition \ref{10.280} for $K$ results:
\begin{equation}
K_N=[K^\prime_N]_{P.A.}+\tilde{K}^\prime_N
\label{10.327}
\end{equation}
with $[K^\prime_N]_{P.A.}$ analogous to $[K^\prime]_{P.A.}$:
\begin{equation}
[K^\prime_N]_{P.A.}=2(L_N-\chi_N\lambda_N)E_N^2\lambda_N+L_N(\lambda_N c_N q^\prime_N E_N\pi_N)
+\lambda_N A_N E_N\tchi_N
\label{10.328}
\end{equation}
We then follow the argument leading to the propagation equation \ref{10.300} with the commutation 
relations \ref{9.b1} in the role of the commutation relations \ref{3.a14}. The argument expressing 
in \ref{10.302} $E^l T^{m-1}\nu$ in terms of $\nu_{m-1,l+1}$ applies to 
express $E_N^l T^{m-1}\nu_N$ in terms of $\nu_{m-1,l+1,N}$. The next argument, which leads to the  
expression in \ref{10.319} of $\lambda E^{l-1}T^m E\tchi$ in terms of $\theta_{l+m}$ and 
$\nu_{i,l+m-i}$ for $i=0,...,m-1$ likewise applies to express $\lambda_N E_N^{l-1}T^m E_N\tchi_N$ 
in terms of $\theta_{l+m,N}$ and $\nu_{i,l+m-i,N}$ for $i=0,...,m-1$. However, in view of the fact 
that $L_N\sk_N$, $\Lb_N\sk_N$ are given by \ref{10.187}, \ref{10.209} respectively, and these 
equations contain the error terms $\vep_{L\sk,N}$, $\vep_{\Lb\sk,N}$, given by 
\ref{10.189}, \ref{10.209} respectively, in replacing equations \ref{10.304} by 
\begin{equation}
L_N\sk_N=E_N\sm_N+I_N, \ \ \ \Lb_N\sk_N=E_N\sn_N+J_N
\label{10.329}
\end{equation}
the remainders $I_N$, $J_N$ must be understood to include the error terms 
$\vep_{L\sk,N}$, $\vep_{\Lb\sk,N}$ respectively. In conclusion, $\nu_{m,l,N}$ satisfies the 
$N$th approximant analogue of the propagation equation \ref{10.320}:
\begin{eqnarray}
&&L_N\nu_{m,l,N}=\left(2\frac{L_N\lambda_N}{\lambda_N}-(l+1)\chi_N\right)\nu_{m,l,N}
+2\frac{L_N\lambda_N}{\lambda_N}E_N^{l-1}T^m j_N\nonumber\\
&&\hspace{20mm}+L_N(E_N^{l-1}T^m(\lambda_N c_N q_N^\prime E_N\pi_N))\nonumber\\
&&\hspace{15mm}+A_N\left\{(2\pi_N\lambda_N)^m\theta_{l+m,N}
+2\sum_{i=0}^{m-1}(2\pi_N\lambda_N)^{m-1-i}\nu_{i,l+m-i,N}\right\}\nonumber\\
&&\hspace{20mm}-m\zeta_N\nu_{m-1,l+1,N}+\tilde{K}^{\prime\prime}_{m,l,N}+E_N^{l-1}T^m\vep_{\nu,N}
\label{10.330}
\end{eqnarray}
where $\vep_{\nu,N}$ is the error term in \ref{10.228}:
\begin{equation}
\vep_{\nu,N}=\lambda_N E_N^2\vep_{\lambda,N}
+\beta_{N,N}^2\left(\pi_N\lambda_N E_N-(1/2)\Lb_N\right)
\left(H^\prime_N(g^{-1})^{\mu\nu}\beta_{\nu,N}\tilde{\kappa}^\prime_{\mu,N}\right) 
\label{10.331}
\end{equation}
In \ref{10.330}, $\tilde{K}^{\prime\prime}_{m,l,N}$ is the analogue of $\tilde{K}^{\prime\prime}_{m,l}$ 
for the $N$th approximants. The error term satisfies the estimate:
\begin{equation}
E_N^{l-1}T^m\vep_{\nu,N}={\bf O}(\tau^{N-2-m})+{\bf O}(u^{-3}\tau^{N-m})
\label{10.332}
\end{equation}

Similarly, we derive from the propagation equation \ref{10.231} for $\nub_N$ 
the following propagation equation for the higher order $N$th approximant quantity 
$\nub_{m,l,N}$:
\begin{eqnarray}
&&\Lb_N\nub_{m,l,N}=\left(2\frac{\Lb_N\lambdab_N}{\lambdab_N}-(l+1)\chib_N\right)\nub_{m,l,N}
+2\frac{\Lb_N\lambdab_N}{\lambdab_N}E_N^{l-1}T^m\jb_N\nonumber\\
&&\hspace{20mm}+\Lb_N(E_N^{l-1}T^m(\lambdab_N c_N\qb_N^\prime E_N\pi_N))\nonumber\\
&&\hspace{15mm}+\oAb_N\left\{(2\pi_N\lambdab_N)^m\thetab_{l+m,N}
+2\sum_{i=0}^{m-1}(2\pi_N\lambdab_N)^{m-1-i}\nub_{i,l+m-i,N}\right\}\nonumber\\
&&\hspace{20mm}+m\zeta_N\nub_{m-1,l+1,N}+\tilde{\Kb}^{\prime\prime}_{m,l,N}+E_N^{l-1}T^m\vep_{\nub,N}
\label{10.333}
\end{eqnarray}
where $\vep_{\nub,N}$ is the error term in \ref{10.231}:
\begin{equation}
\vep_{\nub,N}=\lambdab_N E_N^2\vep_{\lambdab,N}
+\beta_{\Nb,N}^2\left(\pi_N\lambdab_N E_N-(1/2)L_N\right)
\left(H^\prime_N(g^{-1})^{\mu\nu}\beta_{\nu,N}\tilde{\kappa}^\prime_{\mu,N}\right) 
\label{10.334}
\end{equation}
In \ref{10.333}, $\tilde{\Kb}^{\prime\prime}_{m,l,N}$ is the analogue of $\tilde{\Kb}^{\prime\prime}_{m,l}$ 
for the $N$th approximants. The error term satisfies the estimate:
\begin{equation}
E_N^{l-1}T^m\vep_{\nub,N}={\bf O}(\tau^{N-3-m})+{\bf O}(u^{-3}\tau^{N-1-m})
\label{10.335}
\end{equation}

Finally, subtracting equation \ref{10.330} from equation \ref{10.320} we arrive at the following 
propagation equation for the acoustical difference quantity $\cnu_{m,l}$ defined by the 1st of 
\ref{10.239}:
\begin{eqnarray}
&&L\cnu_{m,l}-\left(\frac{2L\lambda}{\lambda}-(l+1)\chi\right)\cnu_{m,l} \nonumber\\
&&\hspace{9mm}-L\left(E^{l-1}T^m(\lambda cq^\prime E\pi)-E_N^{l-1}T^m(\lambda_N c_N q^\prime_N E_N\pi_N)\right)
\nonumber\\
&&\hspace{20mm}=\frac{2L\lambda}{\lambda}E^{l-1}T^m j
-\frac{2L_N\lambda_N}{\lambda_N}E_N^{l-1}T^m j_N+\tilde{K}^{\prime\prime}_{m,l}
-\tilde{K}^{\prime\prime}_{m,l,N}\nonumber\\
&&\hspace{9mm}-m\zeta\cnu_{m-1,l+1}
+A\left\{(2\pi\lambda)^m\cth_{l+m}
+2\sum_{i=0}^{m-1}(2\pi\lambda)^{m-1-i}\cnu_{i,l+m-i}\right\} \nonumber\\
&&\hspace{9mm}-\left\{L-L_N-\frac{2L\lambda}{\lambda}+\frac{2L_N\lambda_N}{\lambda_N} 
+(l+1)(\chi-\chi_N)\right\}\nu_{m,l,N} \nonumber\\
&&\hspace{9mm}+(L-L_N)E_N^{l-1}T^m(\lambda_N c_N q^\prime_N E_N\pi_N)
-m(\zeta-\zeta_N)\nu_{m-1,l+1,N}\nonumber\\
&&\hspace{9mm}+2\sum_{i=0}^{m-1}\left\{A(2\pi\lambda)^{m-1-i}-A_N(2\pi_N\lambda_N)^{m-1-i}\right\}
\nu_{i,l+m-i,N}\nonumber\\
&&\hspace{9mm}+\left\{A(2\pi\lambda)^m-A_N(2\pi_N\lambda_N)^m\right\}\theta_{l+m,N}
-E_N^{l-1}T^m\vep_{\nu,N}
\label{10.336}
\end{eqnarray}
Also, subtracting equation \ref{10.333} from equation \ref{10.323} we arrive at the following 
propagation equation for the acoustical difference quantity $\cnub_{m,l}$ defined by the 2nd of 
\ref{10.239}:
\begin{eqnarray}
&&\Lb\cnub_{m,l}-\left(\frac{2\Lb\lambdab}{\lambdab}-(l+1)\chib\right)\cnub_{m,l} \nonumber\\
&&\hspace{9mm}-\Lb\left(E^{l-1}T^m(\lambdab c\qb^\prime E\pi)-E_N^{l-1}T^m(\lambdab_N c_N\qb^\prime_N E_N\pi_N)\right)
\nonumber\\
&&\hspace{20mm}=\frac{2\Lb\lambdab}{\lambdab}E^{l-1}T^m\jb
-\frac{2\Lb_N\lambdab_N}{\lambdab_N}E_N^{l-1}T^m\jb_N+\tilde{\Kb}^{\prime\prime}_{m,l}
-\tilde{\Kb}^{\prime\prime}_{m,l,N}\nonumber\\
&&\hspace{9mm}+m\zeta\cnub_{m-1,l+1}
+\oAb\left\{(2\pi\lambdab)^m\cthb_{l+m}
+2\sum_{i=0}^{m-1}(2\pi\lambdab)^{m-1-i}\cnub_{i,l+m-i}\right\} \nonumber\\
&&\hspace{9mm}-\left\{\Lb-\Lb_N-\frac{2\Lb\lambdab}{\lambdab}+\frac{2\Lb_N\lambdab_N}{\lambdab_N} 
+(l+1)(\chib-\chib_N)\right\}\nub_{m,l,N} \nonumber\\
&&\hspace{9mm}+(\Lb-\Lb_N)E_N^{l-1}T^m(\lambdab_N c_N\qb^\prime_N E_N\pi_N)
+m(\zeta-\zeta_N)\nub_{m-1,l+1,N}\nonumber\\
&&\hspace{9mm}+2\sum_{i=0}^{m-1}\left\{\oAb(2\pi\lambdab)^{m-1-i}
-\oAb_N(2\pi_N\lambdab_N)^{m-1-i}\right\}\nub_{i,l+m-i,N}\nonumber\\
&&\hspace{9mm}+\left\{\oAb(2\pi\lambdab)^m-\oAb_N(2\pi_N\lambdab_N)^m\right\}\thetab_{l+m,N}
-E_N^{l-1}T^m\vep_{\nub,N}
\label{10.337}
\end{eqnarray}
We remark that the 7th, 8th, 9th, 10th, and 11th terms on the right in each of \ref{10.336}, \ref{10.337} are only of order 1, 
the quantities $\nu_{m,l,N}$, $\nub_{m,l,N}$ and $\theta_{l,N}$, $\thetab_{l,N}$ being known smooth functions of the coordinates 
$(\ub,u,\vartheta)$. 

\vspace{5mm}

\section{Estimates for $\cth_l$}

Consider the propagation equation \ref{10.256} for $\cth_l$. Let us denote the right hand side 
by $\check{S}_l$, so the equation reads:
$$L\cth_l-\left(\frac{2L\lambda}{\lambda}-(l+2)\chi\right)\cth_l=\check{S}_l$$
or:
\begin{equation}
L(\lambda^{-2}\cth_l)+(l+2)\chi(\lambda^{-2}\cth_l)=\lambda^{-2}\check{S}_l
\label{10.338}
\end{equation}
We shall integrate this equation along the integral curves of $L$. Since (see 1st of \ref{2.33})
$$L=\frac{\partial}{\partial\ub}-b\frac{\partial}{\partial\vartheta}$$
in $(\ub,u,\vartheta)$ coordinates the integral curve of $L$ through $(0,u,\vartheta^\prime)\in\Cb_0$ is:
\begin{equation}
\ub \mapsto (\ub,u,\Phi_{\ub,u}(\vartheta^\prime))
\label{10.339}
\end{equation}
where $\Phi_{\ub,u}(\vartheta^\prime)$ is the solution of:
\begin{equation}
\frac{d\Phi_{\ub,u}(\vartheta^\prime)}{d\ub}=-b(\ub,u,\Phi_{\ub,u}(\vartheta^\prime)) , \ \ \ 
\Phi_{0,u}(\vartheta^\prime)=\vartheta^\prime
\label{10.340}
\end{equation}
The mapping 
\begin{equation}
\vartheta^\prime\mapsto\vartheta=\Phi_{\ub,u}(\vartheta^\prime)
\label{10.341}
\end{equation}
is an orientation preserving diffeomorphism of $S^1$, $\Phi_{0,u}$ being the identity on $S^1$. 
We can use $(\ub,u,\vartheta^\prime)$ as coordinates on ${\cal N}$ in place of the $(\ub,u,\vartheta)$.  
Recall Section 2.1, equations \ref{2.18} - \ref{2.22}. These coordinates are adapted to the flow of $L$, $\vartheta^\prime$ being constant along the integral curves of $L$. Hence in these coordinates 
$L$ is simply given by:
\begin{equation}
L=\frac{\partial}{\partial\ub}
\label{10.342}
\end{equation}
and $E$ is given by:
\begin{equation}
E=\frac{1}{\sqrt{\sh^\prime}}\frac{\partial}{\partial\vartheta^\prime}
\label{10.343}
\end{equation}
where $\sh^\prime d\vartheta^\prime\otimes d\vartheta^\prime$ is the induced metric on the 
$S_{\ub,u}$ in the same coordinates. The first of the commutation relations \ref{3.a14} 
in the $(\ub,u,\vartheta^\prime)$ coordinates reads:
\begin{equation}
\frac{\partial\sh^\prime}{\partial\ub}=2\chi\sh^\prime
\label{10.344}
\end{equation}
and we have:
\begin{equation}
\left.\sh^\prime\right|_{\ub=0}=\left.\sh\right|_{\ub=0} \ \ \mbox{as} \ \ \Phi_{0,u}=id
\label{10.345}
\end{equation}
Integrating then \ref{10.344} from $\Cb_0$ we obtain:
\begin{equation}
\sh^\prime(\ub,u,\vartheta^\prime)=\sh(0,u,\vartheta^\prime)
e^{2\int_0^{\ub}\chi(\ub^\prime,u,\vartheta^\prime)d\ub^\prime}
\label{10.346}
\end{equation}
In view of \ref{10.344}, equation \ref{10.338} can be written in $(\ub,u,\vartheta^\prime)$ 
coordinates in the form:
\begin{equation}
\frac{\partial}{\partial\ub}\left(\sh^{\prime(l+2)/2}\lambda^{-2}\cth_l\right)
=\sh^{\prime(l+2)/2}\lambda^{-2}\check{S}_l
\label{10.347}
\end{equation}
Integrating then from $\Cb_0$, noting that 
\begin{equation}
\left.\cth_l\right|_{\Cb_0}=0
\label{10.348}
\end{equation}
we obtain: 
\begin{equation}
\left(\sh^{\prime(l+2)/2}\lambda^{-2}\cth_l\right)(\ub_1,u,\vartheta^\prime)= 
\int_0^{\ub_1}\left(\sh^{\prime(l+2)/2}\lambda^{-2}\check{S}_l\right)(\ub,u,\vartheta^\prime)d\ub 
\label{10.349}
\end{equation}

Before we proceed to analyze \ref{10.349}, we make the following remark on the $L^2$ norm of the  
restriction to $S_{\ub,u}$ of a function on ${\cal N}$. 
If $f$ is a function on $S^1$ its $L^2$ norm is:
\begin{equation}
||f||_{L^2(S^1)}=\sqrt{\int_{\vartheta^\prime\in S^1}f^2(\vartheta^\prime)d\vartheta^\prime}
\label{10.350}
\end{equation}
On the other hand, if $f$ is a function on ${\cal N}$ and we represent $f$ in 
$(\ub,u,\vartheta^\prime)$ coordinates, the $L^2$ norm of $f$ on $S_{\ub,u}$ is:
\begin{equation}
||f||_{L^2(S_{\ub,u})}=\sqrt{\int_{\vartheta^\prime\in S^1}\left(f^2\sqrt{\sh^\prime}\right)
(\ub,u,\vartheta^\prime)d\vartheta^\prime}
\label{10.351}
\end{equation}
$\sqrt{\sh^\prime(\ub,u,\vartheta^\prime)}d\vartheta^\prime$ being the element of arc length on $S_{\ub,u}$. 
Comparing \ref{10.351} with \ref{10.350} we see that:
\begin{equation}
||f||_{L^2(S_{\ub,u})}=\left\|\left(f\sh^{\prime 1/4}\right)(\ub,u)\right\|_{L^2(S^1)}
\label{10.352}
\end{equation}
Here $g$ being an arbitrary function on ${\cal N}$ we denote by $g(\ub,u)$ the restriction of 
$g$ to $S_{\ub,u}$ represented in terms of the $\vartheta^\prime$ coordinate. 

With the above remark in mind, we rewrite \ref{10.349} in the form:
\begin{eqnarray}
&&\left(\cth_l\sh^{\prime 1/4}\right)(\ub_1,u,\vartheta^\prime)= \label{10.353}\\
&&\hspace{10mm}\int_0^{\ub_1}\left(\frac{\sh^\prime(\ub,u,\vartheta^\prime)}{\sh^\prime(\ub_1,u,\vartheta^\prime)}
\right)^{(2l+3)/4}\left(\frac{\lambda(\ub_1,u,\vartheta^\prime)}{\lambda(\ub,u,\vartheta^\prime)}
\right)^2\left(\check{S}_l\sh^{\prime 1/4}\right)(\ub,u,\vartheta^\prime)d\ub \nonumber
\end{eqnarray}
Now, by \ref{10.346}:
\begin{equation}
\frac{\sqrt{\sh^\prime(\ub,u,\vartheta^\prime)}}{\sqrt{\sh^\prime(\ub_1,u,\vartheta^\prime)}}
=e^{-\int_{\ub}^{\ub_1}\chi(\ub^\prime,u,\vartheta^\prime)d\ub^\prime}
\label{10.354}
\end{equation}
Assuming then a bound on $\chi$ of the form:
\begin{equation}
\sup_{S_{\ub,u}}|\chi|\leq C\ub \ : \ \forall \ub\in[0,\ub_1]
\label{10.355}
\end{equation}
we have:
\begin{equation}
e^{-\frac{1}{2}C(\ub_1^2-\ub^2)}\leq \frac{\sqrt{\sh^\prime(\ub,u,\vartheta^\prime)}}{\sqrt{\sh^\prime(\ub_1,u,\vartheta^\prime)}}\leq e^{\frac{1}{2}C(\ub_1^2-\ub^2)}
\label{10.356}
\end{equation}
Therefore:
\begin{equation}
\left(\frac{\sh^\prime(\ub,u,\vartheta^\prime)}{\sh^\prime(\ub_1,u,\vartheta^\prime)}
\right)^{(2l+3)/4}\leq k
\label{10.357}
\end{equation}
where $k$ is a constant greater than 1, but which can be chosen as close to 1 as we wish by 
suitably restricting $\ub_1$. 

By Proposition 7.2 $\lambda$ is non-increasing along the integral curves of $L$, hence:
\begin{equation}
\left(\frac{\lambda(\ub_1,u,\vartheta^\prime)}{\lambda(\ub,u,\vartheta^\prime)}\right)^2\leq 1
\label{10.358}
\end{equation}

In view of the above, \ref{10.353} implies:
\begin{equation}
\left|\left(\cth_l\sh^{\prime 1/4}\right)(\ub_1,u,\vartheta^\prime)\right|\leq k
\int_0^{\ub_1}\left|\left(\check{S}_l\sh^{\prime 1/4}\right)(\ub,u,\vartheta^\prime)\right|d\ub
\label{10.359}
\end{equation}
Taking $L^2$ norms with respect to $\vartheta^\prime\in S^1$ this implies:
\begin{equation}
\left\|\left(\cth_l\sh^{\prime 1/4}\right)(\ub_1,u)\right\|_{L^2(S^1)}\leq
 k\int_0^{\ub_1}\left\|\left(\check{S}_l\sh^{\prime 1/4}\right)(\ub,u)\right\|_{L^2(S^1)}d\ub
 \label{10.360}
 \end{equation}
or, in view of \ref{10.352}:
\begin{equation}
\|\cth_l\|_{L^2(S_{\ub_1,u})}\leq k\int_0^{\ub_1}\|\check{S}_l\|_{L^2(S_{\ub,u})}d\ub
\label{10.361}
\end{equation}

We shall presently estimate the contribution of the principal terms in $\check{S}_l$, the right hand 
side of \ref{10.256}, to the integral on the right in \ref{10.361}. These principal terms are the two 
difference terms:
\begin{equation}
-\frac{2L\lambda}{\lambda}\left(E^l f-E_N^l f_N\right)+\tilde{R}_l-\tilde{R}_{l,N}
\label{10.362}
\end{equation}
In fact, since by Proposition 7.2:
\begin{equation}
-\frac{2L\lambda}{\lambda}\sim\frac{\lambdab}{\lambda}\sim\frac{\ub}{u^2}
\label{10.363}
\end{equation}
the first of the two difference terms \ref{10.362} makes, as we shall see, the leading contribution. 
In fact the leading contribution comes from the 1st term in $f$ as given by Proposition 10.1, 
that is from:
\begin{equation}
-\frac{1}{2}\beta_N^2\left(E^l\Lb H-E_N^l\Lb_N H_N\right)
\label{10.364}
\end{equation}
More precisely, since 
$$\Lb H=H^\prime \Lb\sigma=-2H^\prime(g^{-1})^{\mu\nu}\beta_\nu\Lb\beta_\mu$$
hence
$$[E^l\Lb H]_{P.P.}=-2H^\prime(g^{-1})^{\mu\nu}\beta_\nu E^l\Lb\beta_\mu$$
and similarly for the $N$th approximants, the leading contribution comes from the principal part of \ref{10.364} which is:
\begin{equation}
\beta_N^2 H^\prime(g^{-1})^{\mu\nu}\beta_\nu
\Lb(E^l\beta_\mu-E_N^l\beta_{\mu,N})
\label{10.365}
\end{equation}

Now, if $X$ is an arbitrary vectorfield on ${\cal N}$ we can expand $X$ in the $(E,N,\Nb)$ frame, so 
if $X^\mu$ are the rectangular components of $X$ we have:
\begin{eqnarray}
&&X^\mu=X^E E^\mu+X^N N^\mu+X^{\Nb}\Nb^\mu \mbox{and:} \label{10.366}\\
&&X^E=h_{\mu\nu} X^\mu E^\nu, \ \ X^N=-(1/2c)h_{\mu\nu}X^\mu\Nb^\nu, \ \ 
X^{\Nb}=-(1/2c)h_{\mu\nu}X^\mu N^\nu \nonumber
\end{eqnarray}
and for arbitrary $\theta_\mu$ we have:
\begin{equation}
X^\mu\theta_\mu=X^E E^\mu\theta_\mu+X^N N^\mu\theta_\mu+X^{\Nb}\Nb^\mu\theta_\mu 
\label{10.367}
\end{equation}
In particular, for $X^\mu=(g^{-1})^{\mu\nu}\beta_\nu$ we have:
$$h_{\mu\nu}X^\mu=(1-\sigma H)\beta_\nu=\eta^2 \beta_\nu$$
($\eta$ being here the sound speed). Hence for arbitrary $\theta_\mu$ it holds:
\begin{equation}
(g^{-1})^{\mu\nu}\beta_\nu\theta_\mu=\eta^2\left(\sbeta E^\mu\theta_\mu
-(1/2c)\beta_{\Nb} N^\mu\theta_\mu-(1/2c)\beta_N\Nb^\mu\theta_\mu\right)
\label{10.368}
\end{equation}
It follows that \ref{10.365} is expressed as:
\begin{eqnarray}
&&\beta_N^2 H^\prime\eta^2\left\{\sbeta E^\mu\Lb(E^l\beta_\mu-E_N^l\beta_{\mu,N})\right. 
\label{10.369}\\
&&\hspace{5mm}\left.-(1/2c)\beta_{\Nb}N^\mu\Lb(E^l\beta_\mu-E_N^l\beta_{\mu,N})
-(1/2c)\beta_N\Nb^\mu\Lb(E^l\beta_\mu-E_N^l\beta_{\mu,N})\right\} \nonumber
\end{eqnarray}
As we shall see, it is the 3rd term in parenthesis which makes the leading partial contribution. 
In view of \ref{10.363} this partial contribution to the integral on the right in \ref{10.361} is 
bounded by:
\begin{equation}
Cu^{-2}\int_0^{\ub_1}\|\Nb^\mu\Lb(E^l\beta_\mu-E_N^l\beta_{\mu,N})\|_{L^2(S_{\ub,u})}\ub d\ub
\label{10.370}
\end{equation}

Now, according to the definition \ref{9.236} of the higher order variation differences we have: 
\begin{equation}
E^l\beta_\mu-E_N^l\beta_{\mu,N}=\s^{(0,l)}\check{\dot{\phi}}_\mu
\label{10.371}
\end{equation}
Thus, according to the definition \ref{9.275} for each of the two variation fields $V=Y,E$ we have:
\begin{equation}
V^\mu d(E^l\beta_\mu-E_N^l\beta_{\mu,N})=\s^{(V;0,l)}\check{\xi}
\label{10.372}
\end{equation}
In particular taking $V=Y$ we have:
\begin{equation}
N^\mu\Lb(E^l\beta_\mu-E_N^l\beta_{\mu,N})+\ogamma\Nb^\mu\Lb(E^l\beta_\mu-E_N^l\beta_{\mu,N})
=\s^{(Y;0,l)}\check{\xi}_{\Lb}
\label{10.373}
\end{equation}
Here, in view of the definition \ref{9.291} the right hand side is bounded in $L^2(\Cb_{\ub}^u)$ norm
by $C\sqrt{\s^{(Y;0,l)}\cEb^u(\ub)}$.  To estimate the first term on the left in $L^2(\Cb_{\ub}^u)$ 
norm we use the fact that:
\begin{equation}
N^\mu\Lb\beta_\mu=\lambda E^\mu E\beta_{\mu}
\label{10.374}
\end{equation}
This allows us to express:
\begin{eqnarray}
&&N^\mu\Lb E^l\beta_\mu=N^\mu E^l \Lb\beta_\mu+\mbox{terms of order $l$} \nonumber\\
&&\hspace{15mm}=E^l(N^\mu\Lb\beta_\mu)+\mbox{terms of order $l$} \nonumber\\
&&\hspace{15mm}=E^l(\lambda E^\mu E\beta_\mu)+\mbox{terms of order $l$} \nonumber\\
&&\hspace{15mm}=\lambda E^\mu E^{l+1}\beta_\mu+\mbox{terms of order $l$} 
\label{10.375}
\end{eqnarray}
Similarly for the $N$th approximants, but in this case there is a error term:
\begin{equation}
E_N^l\left(\rho_N^{-1}(-\delta_N^\prime-\omega_{L\Lb,N})\right)={\bf O}(\tau^{N-1})
\label{10.376}
\end{equation}
Writing also $N_N^\mu\Lb_N E_N^l\beta_{\mu,N}$ as $N^\mu \Lb E_N^l\beta_{\mu,N}$ up to the 0th 
order term \\$(N^\mu\Lb-N_N^\mu\Lb_N)E_N^l\beta_{\mu,N}$, we express in this way 
$$N^\mu\Lb(E^l\beta_\mu-E_N^l\beta_{\mu,N})$$ 
up to lower order terms as:
\begin{equation}
\lambda E^\mu E(E^l\beta_\mu-E_N^l\beta_{\mu,N})=\lambda\s^{(E;0,l)}\check{\sxi}
\label{10.377}
\end{equation}
by \ref{10.372} for $V=E$. In view of the definition \ref{9.291}, we have:
\begin{equation}
\|\sqrt{a}\s^{(E;0,l)}\check{\sxi}\|_{L^2(\Cb_{\ub}^u)}\leq C\sqrt{\s^{(E;0,l)}\cEb^u(\ub)}
\label{10.378}
\end{equation}
Recalling \ref{7.106} - \ref{7.108}, \ref{7.132}, the above imply that: 
\begin{equation}
\|\Nb^\mu\Lb(E^l\beta_\mu-E_N^l\beta_{\mu,N})\|_{L^2(S_{\ub,u})}\leq g_0(\ub,u)+g_1(\ub,u)
\label{10.379}
\end{equation}
where $g_0$ and $g_1$ are non-negative functions of $(\ub,u)$ such that, neglecting lower order 
terms,
\begin{equation}
\int_{\ub}^{u_1}u^2 g_0^2(\ub,u)du\leq C\s^{(Y;0,l)}\cEb^{u_1}(\ub)
\label{10.380}
\end{equation}
while:
\begin{equation}
\ub \int_{\ub}^{u_1}g_1^2(\ub,u)du\leq C\s^{(E;0,l)}\cEb^{u_1}(\ub)
\label{10.381}
\end{equation}
Then \ref{10.370} is bounded by:
\begin{equation}
C\left(G_0(\ub_1,u)+G_1(\ub_1,u)\right) 
\label{10.382}
\end{equation}
where: 
\begin{equation}
G_i(\ub_1,u)=u^{-2}\int_0^{\ub_1}g_i(\ub,u)\ub d\ub \ : \ i=0,1
\label{10.383}
\end{equation}
We shall bound 
\begin{equation}
\|G_i(\ub_1,\cdot)\|_{L^2([\ub_1,u_1])}:=\sqrt{\int_{\ub_1}^{u_1}G_i^2(\ub_1,u)du} \ 
\mbox{: for $i=0,1$}
\label{10.384}
\end{equation}
which, through \ref{10.361}, bound the contribution of \ref{10.370} to the $L^2(\Cb_{\ub_1}^{u_1})$ of 
$\cth_l$. We have (Schwartz inequality):
\begin{eqnarray}
&&G_0^2(\ub_1,u)\leq u^{-4}\int_0^{\ub_1}\ub^2 d\ub \cdot \int_0^{\ub_1} g_0^2(\ub,u)d\ub\nonumber\\
&&\hspace{15mm}=\frac{1}{3}\frac{\ub_1^3}{u^4}\int_0^{\ub_1}g_0^2(\ub,u)d\ub 
\label{10.385}
\end{eqnarray}
Hence:
\begin{eqnarray}
&&\int_{\ub_1}^{u_1}G_0^2(\ub_1,u)du\leq\frac{1}{3}\ub_1^3
\int_{\ub_1}^{u_1}\left\{\int_0^{\ub_1}g_0^2(\ub,u)d\ub\right\}\frac{du}{u^4}\nonumber\\
&&\hspace{25mm}=\frac{1}{3}\ub_1^3\int_0^{\ub_1}
\left\{\int_{\ub_1}^{u_1}u^{-4}g_0^2(\ub,u)du\right\}d\ub
\label{10.386}
\end{eqnarray}
(reversing the order of integration). In regard to the interior integral we have, for all 
$\ub\in[0,\ub_1]$:
\begin{eqnarray}
&&\int_{\ub_1}^{u_1}u^{-4}g_0^2(\ub,u)du
=\int_{\ub_1}^{u_1}\frac{\partial}{\partial u}\left(\int_{\ub_1}^u u^{\prime 2}g_0^2(\ub,u^\prime)
du^\prime\right)u^{-6}du\nonumber\\
&&\hspace{10mm}=u_1^{-6}\int_{\ub_1}^{u_1}u^{\prime 2}g_0^2(\ub,u^\prime)du^\prime
+6\int_{\ub_1}^{u_1}\left(\int_{\ub_1}^u u^{\prime 2}g_0^2(\ub,u^\prime)du^\prime\right)u^{-7}du 
\nonumber\\
&&\label{10.387}
\end{eqnarray}
Since $\ub_1\geq \ub$, by \ref{10.380} we have, for all $u\in[\ub_1,u_1]$:
\begin{eqnarray}
&&\int_{\ub_1}^u u^{\prime 2} g_0^2(\ub,u^\prime)du^\prime\leq 
\int_{\ub}^u u^{\prime 2}g_0^2(\ub,u^\prime)du^\prime \nonumber\\
&&\hspace{15mm}\leq C\s^{(Y;0,l)}\cEb^u(\ub)\leq C\s^{(Y;0,l)}\cBb(\ub_1,u_1)\cdot \ub^{2a_0} u^{2b_0}
\label{10.388}
\end{eqnarray}
in view of the definition \ref{9.298}. Then, assuming that 
\begin{equation}
b_0>3
\label{10.389}
\end{equation}
\ref{10.387} is bounded by:
\begin{equation}
C\s^{(Y;0,l)}\cBb(\ub_1,u_1)\ub^{2a_0}\left\{u_1^{2b_0-6}
+6\frac{(u_1^{2b_0-6}-\ub_1^{2b_0-6})}{(2b_0-6)}\right\}
\label{10.390}
\end{equation}
Taking in fact 
\begin{equation}
b_0\geq 4
\label{10.391}
\end{equation}
we have a bound of \ref{10.387} by:
\begin{equation}
C\s^{(Y,0,l)}\cBb(\ub_1,u_1)\ub^{2a_0} u_1^{2b_0-6}
\label{10.392}
\end{equation}
Then by \ref{10.386}:
\begin{eqnarray}
&&\int_{\ub_1}^{u_1}G_0^2(\ub_1,u)du\leq C\s^{(Y;0,l)}\cBb(\ub_1,u_1)u_1^{2b_0-6}\ub_1^3
\int_0^{\ub_1}\ub^{2a_0}d\ub \nonumber\\
&&\hspace{27mm}=\frac{C\s^{(Y;0,l)}\cBb(\ub_1,u_1)}{(2a_0+1)}
\cdot\ub_1^{2a_0} u_1^{2b_0}\cdot\frac{\ub_1^4}{u_1^6}
\label{10.393}
\end{eqnarray}
The last factor $\ub_1^4/u_1^6$ is homogeneous of degree -2. This signifies, as we shall see, that 
the contribution is {\em borderline}. To control it we shall have to choose $a_0$ appropriately 
large so that the coefficient $C/(2a_0+1)$ is suitably small. 

Consider next $G_1$. We have (Schwartz inequality):
\begin{eqnarray}
&&G_1^2(\ub_1,u)\leq u^{-4}\int_0^{\ub_1}\ub d\ub \cdot \int_0^{\ub_1}\ub g_1^2(\ub,u)d\ub\nonumber\\
&&\hspace{15mm}=\frac{1}{2}\frac{\ub_1^2}{u^4}\int_0^{\ub_1}\ub g_1^2(\ub,u)d\ub 
\label{10.394}
\end{eqnarray}
Hence:
\begin{eqnarray}
&&\int_{\ub_1}^{u_1}G_1^2(\ub_1,u)du\leq\frac{1}{2}\ub_1^2
\int_{\ub_1}^{u_1}\left\{\int_0^{\ub_1}\ub g_1^2(\ub,u)d\ub\right\}\frac{du}{u^4}\nonumber\\
&&\hspace{25mm}=\frac{1}{2}\ub_1^2\int_0^{\ub_1}
\left\{\int_{\ub_1}^{u_1}u^{-4}\ub g_1^2(\ub,u)du\right\}d\ub
\label{10.395}
\end{eqnarray}
(reversing the order of integration). In regard to the interior integral we have, for all 
$\ub\in[0,\ub_1]$:
\begin{eqnarray}
&&\int_{\ub_1}^{u_1}u^{-4}\ub g_1^2(\ub,u)du
=\int_{\ub_1}^{u_1}\frac{\partial}{\partial u}\left(\int_{\ub_1}^u \ub g_1^2(\ub,u^\prime)
du^\prime\right)u^{-4}du\nonumber\\
&&\hspace{10mm}=u_1^{-4}\int_{\ub_1}^{u_1}\ub g_1^2(\ub,u^\prime)du^\prime
+4\int_{\ub_1}^{u_1}\left(\int_{\ub_1}^u \ub g_1^2(\ub,u^\prime)du^\prime\right)u^{-5}du 
\nonumber\\
&&\label{10.396}
\end{eqnarray}
Since $\ub_1\geq \ub$, by \ref{10.381} we have, for all $u\in[\ub_1,u_1]$:
\begin{eqnarray}
&&\int_{\ub_1}^u \ub g_1^2(\ub,u^\prime)du^\prime\leq 
\int_{\ub}^u \ub g_1^2(\ub,u^\prime)du^\prime \nonumber\\
&&\hspace{15mm}\leq C\s^{(E;0,l)}\cEb^u(\ub)\leq C\s^{(E;0,l)}\cBb(\ub_1,u_1)\cdot \ub^{2a_0} u^{2b_0}
\label{10.397}
\end{eqnarray}
in view of the definition \ref{9.298}. Then \ref{10.396} is bounded by:
\begin{eqnarray}
&&C\s^{(E;0,l)}\cBb(\ub_1,u_1)\ub^{2a_0}\left\{u_1^{2b_0-4}
+4\frac{(u_1^{2b_0-4}-\ub_1^{2b_0-4})}{(2b_0-4)}\right\}\nonumber\\
&&\leq C\s^{(E;0,l)}\cBb(\ub_1,u_1)\ub^{2a_0}u_1^{2b_0-4}
\label{10.398}
\end{eqnarray}
(in view of \ref{10.391}), and by \ref{10.395}:
\begin{eqnarray}
&&\int_{\ub_1}^{u_1}G_1^2(\ub_1,u)du\leq C\s^{(E;0,l)}\cBb(\ub_1,u_1)u_1^{2b_0-4}\ub_1^2
\int_0^{\ub_1}\ub^{2a_0}d\ub \nonumber\\
&&\hspace{27mm}=\frac{C\s^{(E;0,l)}\cBb(\ub_1,u_1)}{(2a_0+1)}
\cdot\ub_1^{2a_0} u_1^{2b_0}\cdot\frac{\ub_1^3}{u_1^4}
\label{10.399}
\end{eqnarray}
The last factor $\ub_1^3/u_1^4$ is homogeneous of degree -1, which is one degree better than 
borderline. Restricting ourselves to ${\cal R}_{\delta,\delta}$ the contribution will come with an 
extra factor of $\delta$ and can be absorbed in the energy estimates by choosing $\delta$ 
appropriately small. 

We have now treated the partial contribution of the 3rd term in parenthesis in \ref{10.369}, 
which is the leading term. The partial contribution of the 2nd term is treated in the same way 
as the 1st term on the left in \ref{10.373}. Here we have an extra factor of $\ogamma\sim u$, 
therefore in place of \ref{10.379} we have:
\begin{equation}
\|N^\mu\Lb(E^l\beta_\mu-E_N^l\beta_{\mu,N})\|_{L^2(S_{\ub,u})}\leq  Cu g_1(\ub,u)
\label{10.400}
\end{equation}
Then the partial contribution of this term to $\|\cth_l\|^2_{L^2(\Cb_{\ub_1}^{u_1})}$ is bounded 
by \ref{10.399} with an extra $u_1^2$ factor, that is by:
\begin{equation}
\frac{C\s^{(E;0,l)}\cBb(\ub_1,u_1)}{(2a_0+1)}
\cdot\ub_1^{2a_0} u_1^{2b_0}\cdot\frac{\ub_1^3}{u_1^2}
\label{10.401}
\end{equation}
The last factor being homogeneous of degree 1, this is 3 degrees better than borderline, so the 
contribution will come with an extra factor of $\delta^3$. As for the partial contribution of the 
1st term in parenthesis in \ref{10.369}, since by \ref{10.372} for $V=E$ we have:
\begin{equation}
E^\mu\Lb(E^l\beta_\mu-E_N^l\beta_{\mu,N})=\s^{(E;0,l)}\check{\xi}_{\Lb}
\label{10.402}
\end{equation}
this partial contribution to the integral on the right in \ref{10.361} is bounded by $CG_2(\ub,u)$, 
where, as in \ref{10.383},
$$G_2(\ub_1,u)=u^{-2}\int_0^{\ub_1}g_2(\ub,u)\ub d\ub$$ 
and here:
\begin{equation}
g_2(\ub,u)=\|\s^{(E;0,l)}\check{\xi}_{\Lb}\|_{L^2(S_{\ub,u})}
\label{10.403}
\end{equation}
so in view of the definition \ref{9.291} we have:
\begin{equation}
\int_{\ub}^{u_1}g_2^2(\ub,u)du\leq C\s^{(E;0,l)}\cEb^{u_1}(\ub)
\label{10.404}
\end{equation}
(compare with \ref{10.380}). The corresponding partial contribution to 
$\|\cth_l\|^2_{L^2(\Cb_{\ub_1}^{u_1})}$ is then bounded by:
\begin{equation}
\frac{C\s^{(E;0,l)}\cBb(\ub_1,u_1)}{(2a_0+1)}
\cdot\ub_1^{2a_0} u_1^{2b_0}\cdot\frac{\ub_1^4}{u_1^4}
\label{10.405}
\end{equation}
The last factor being homogeneous of degree 0, this is 2 degrees better than borderline, so the 
contribution will come with an extra factor of $\delta^2$. This completes the treatment of the 
leading contribution to $\|\cth_l\|_{L^2(\Cb_{\ub_1}^{u_1})}$, which is that of the term \ref{10.364},  
through the expression \ref{10.369}. Comparing \ref{10.393}, \ref{10.401}, \ref{10.405}, 
we see that the dominant contribution is that of the $\Nb$-component, 
the least that of the $N$-component, while that of the $E$-component is intermediate, 
the extra factors being homogeneous of degree -2, 1, and 
0 respectively. 

We turn to the contribution coming from the 2nd term in $f$ as given by Proposition 10.1, that is 
from:
\begin{equation}
\lambda\beta_N\sbeta(E^{l+1}H-E_N^{l+1}H_N)
\label{10.406}
\end{equation}
Here, in comparison with \ref{10.364} we have an extra factor of $\lambda\sim u^2$. We consider the 
principal part of \ref{10.406}, which is:
\begin{equation}
-\lambda\beta_N\sbeta H^\prime(g^{-1})^{\mu\nu}\beta_\nu
E(E^l\beta_\mu-E_N^l\beta_{\mu,N})
\label{10.407}
\end{equation}
(compare with \ref{10.365}). By \ref{10.368} this is expressed as:
\begin{eqnarray}
&&-\lambda\beta_N\sbeta H^\prime\eta^2\left\{\sbeta E^\mu E(E^l\beta_\mu-E_N^l\beta_{\mu,N})\right. 
\label{10.408}\\
&&\hspace{5mm}\left.-(1/2c)\beta_{\Nb}N^\mu E(E^l\beta_\mu-E_N^l\beta_{\mu,N})
-(1/2c)\beta_N\Nb^\mu E(E^l\beta_\mu-E_N^l\beta_{\mu,N})\right\} \nonumber
\end{eqnarray}
Because of the extra factor of $\lambda$ the contribution of this to the integral on the right in 
\ref{10.361} is bounded by a constant multiple of:
\begin{eqnarray}
&&\int_0^{\ub_1}\|E^\mu E(E^l\beta_\mu-E_N^l\beta_{\mu,N})\|_{L^2(S_{\ub,u})}\ub d\ub \nonumber\\
&&+\int_0^{\ub_1}\|N^\mu E(E^l\beta_\mu-E_N^l\beta_{\mu,N})\|_{L^2(S_{\ub,u})}\ub d\ub \nonumber\\
&&+\int_0^{\ub_1}\|\Nb^\mu E(E^l\beta_\mu-E_N^l\beta_{\mu,N})\|_{L^2(S_{\ub,u})}\ub d\ub
\label{10.409}
\end{eqnarray}
Since 
$$E^\mu E(E^l\beta_\mu-E_N^l\beta_{\mu,N})=\s^{(E;0,l)}\check{\sxi}$$
the first integral is:
\begin{equation}
\tilde{G}_1(\ub_1,u)=\int_0^{\ub_1}\tilde{g}_1(\ub,u)\ub d\ub
\label{10.410}
\end{equation}
where:
\begin{equation}
\tilde{g}_1(\ub,u)=\|\s^{(E;0,l)}\check{\sxi}\|_{L^2(S_{\ub,u})}
\label{10.411}
\end{equation}
We have (Schwartz inequality):
\begin{eqnarray}
&&\tilde{G}_1^2(\ub_1,u)\leq \int_0^{\ub_1}\ub d\ub \cdot 
\int_0^{\ub_1}\tilde{g}_1^2(\ub,u)\ub d\ub\nonumber\\
&&\hspace{15mm}=\frac{1}{2}\ub_1^2\int_0^{\ub_1}\tilde{g}_1^2(\ub,u)\ub d\ub 
\label{10.412}
\end{eqnarray}
Hence:
\begin{eqnarray}
&&\int_{\ub_1}^{u_1}\tilde{G}_1^2(\ub_1,u)du\leq\frac{1}{2}\ub_1^2
\int_{\ub_1}^{u_1}\left\{\int_0^{\ub_1}\tilde{g}_1^2(\ub,u)\ub d\ub\right\}du\nonumber\\
&&\hspace{25mm}=\frac{1}{2}\ub_1^2\int_0^{\ub_1}
\left\{\int_{\ub_1}^{u_1}\tilde{g}_1^2(\ub,u)\ub du\right\}d\ub
\label{10.413}
\end{eqnarray}
(reversing the order of integration). In regard to the interior integral we have, for all 
$\ub\in[0,\ub_1]$:
\begin{eqnarray}
&&\int_{\ub_1}^{u_1}\ub\tilde{g}_1^2(\ub,u)du
=\int_{\ub_1}^{u_1}\frac{\partial}{\partial u}\left(\int_{\ub_1}^u 
u^{\prime 2}\ub\tilde{g}_1^2(\ub,u^\prime)du^\prime\right)u^{-2}du\nonumber\\
&&\hspace{10mm}=u_1^{-2}\int_{\ub_1}^{u_1}u^{\prime 2}\ub\tilde{g}_1^2(\ub,u^\prime)du^\prime
+2\int_{\ub_1}^{u_1}
\left(\int_{\ub_1}^u u^{\prime 2}\ub\tilde{g}_1^2(\ub,u^\prime)du^\prime\right)u^{-3}du 
\nonumber\\
&&\label{10.414}
\end{eqnarray}
By \ref{10.378}, since $a\sim u^2\ub$, for all $u\in[\ub,u_1$ and all $\ub\in[0,\ub_1]$ we have:
\begin{eqnarray}
&&\int_{\ub}^u u^{\prime 2}\ub\tilde{g}_1^2(\ub,u^\prime)du^\prime
\leq C\|\sqrt{a}\s^{(E;0,l)}\check{\sxi}\|^2_{L^2(\Cb_{\ub}^u)} \label{10.415}\\
&&\hspace{30mm}\leq C\s^{(E;0,l)}\cEb^u(\ub)\leq C\s^{(E;0,l)}\cBb(\ub_1,u_1)\ub^{2a_0} u^{2b_0}
\nonumber
\end{eqnarray}
Noting that since $\ub_1\geq \ub$ we have 
$$\int_{\ub_1}^u u^{\prime 2}\ub\tilde{g}_1^2(\ub,u^\prime)du^\prime\leq 
\int_{\ub}^u u^{\prime 2}\ub\tilde{g}_1^2(\ub,u^\prime)du^\prime,$$
we then conclude that \ref{10.414} is bounded by:
\begin{eqnarray}
&&C\s^{(E;0,l)}\cBb(\ub_1,u_1)\ub^{2a_0}\left\{u_1^{2b_0-2}
+2\frac{(u_1^{2b_0-2}-\ub_1^{2b_0-2})}{(2b_0-2)}\right\}\nonumber\\
&&\leq C\s^{(E;0,l)}\cBb(\ub_1,u_1)\ub^{2a_0}u_1^{2b_0-2}
\label{10.416}
\end{eqnarray}
(in view of \ref{10.391}), and by \ref{10.413}:
\begin{eqnarray}
&&\int_{\ub_1}^{u_1}\tilde{G}_1^2(\ub_1,u)du\leq C\s^{(E;0,l)}\cBb(\ub_1,u_1)u_1^{2b_0-2}\ub_1^2
\int_0^{\ub_1}\ub^{2a_0}d\ub \nonumber\\
&&\hspace{27mm}=\frac{C\s^{(E;0,l)}\cBb(\ub_1,u_1)}{(2a_0+1)}
\cdot\ub_1^{2a_0} u_1^{2b_0}\cdot\frac{\ub_1^3}{u_1^2}
\label{10.417}
\end{eqnarray}
This bounds the contribution of the 1st of the integrals \ref{10.409} to 
$\|\cth_l\|^2_{L^2(\Cb_{\ub_1}^{u_1})}$. The last factor in \ref{10.417} being homogeneous of degree 1, 
which is 3 degrees better than borderline, this contribution will come with an extra factor of 
$\delta^3$.  

Consider next the contribution of the 2nd of the integrals \ref{10.409}. Since $\rho\sim\ub$ this 
integral is bounded by a constant multiple of the integral
\begin{equation}
\int_0^{\ub_1}\|L^\mu E(E^l\beta_\mu-E_N^l\beta_{\mu,N})\|_{L^2(S_{\ub,u})}d\ub
\label{10.418}
\end{equation}
To estimate this we use the fact that:
\begin{equation}
L^\mu E\beta_\mu=E^\mu L\beta_\mu 
\label{10.419}
\end{equation}
This allows us to express: 
\begin{eqnarray}
&&L^\mu E^{l+1}\beta_\mu=E^l(L^\mu E\beta_\mu)+\mbox{terms of order $l$} \nonumber\\
&&\hspace{17mm}=E^l(E^\mu L\beta_\mu)+\mbox{terms of order $l$} \nonumber\\
&&\hspace{17mm}=E^\mu E^l L\beta_\mu+\mbox{terms of order $l$} \nonumber\\
&&\hspace{17mm}=E^\mu LE^l\beta_\mu+\mbox{terms of order $l$} 
\label{10.420}
\end{eqnarray}
Similarly for the $N$th approximants, but in this case there is a error term:
\begin{equation}
E_N^l\omega_{LE,N}={\bf O}(\tau^N)
\label{10.421}
\end{equation}
Writing also $L_N^\mu E_N^{l+1}\beta_{\mu,N}$ as $L^\mu E E_N^l\beta_{\mu,N}$ up to the 0th 
order term \\$(L^\mu E-L_N^\mu E_N)E_N^l\beta_{\mu,N}$, we express in this way 
$$L^\mu E(E^l\beta_\mu-E_N^l\beta_{\mu,N})$$ 
up to lower order terms as:
\begin{equation}
E^\mu L(E^l\beta_\mu-E_N^l\beta_{\mu,N})=\s^{(E;0,l)}\check{\xi}_L
\label{10.422}
\end{equation}
by \ref{10.372} for $V=E$. We can then bound the integral \ref{10.418} in terms of:
\begin{eqnarray}
&&\int_0^{\ub_1}\|\s^{(E;0,l)}\check{\xi}_L\|_{L^2(S_{\ub,u})}d\ub
\leq \sqrt{\ub_1}\|\s^{(E;0,l)}\check{\xi}_L\|_{L^2(C_u^{\ub_1})} \nonumber\\
&&\hspace{35mm}\leq C\sqrt{\ub_1}\sqrt{\s^{(E;0,l)}\cE^{\ub_1}(u)} \nonumber\\
&&\hspace{35mm}\leq C\sqrt{\ub_1}\sqrt{\s^{(E;0,l)}\cB(\ub_1,u_1)}\cdot \ub_1^{a_0} u^{b_0}
\label{10.423}
\end{eqnarray}
in view of \ref{9.290} and the definition \ref{9.297}. 
This bounds the contribution to $\|\cth_l\|_{L^2(S_{\ub_1,u})}$ for all $u\in[\ub_1,u_1]$. 
The corresponding contribution to $\|\cth_l\|^2_{L^2(\Cb_{\ub_1}^{u_1})}$ is then bounded by:
\begin{eqnarray}
&&C\s^{(E;0,l)}\cB(\ub_1,u_1)\cdot \ub_1^{2a_0+1}\int_{\ub_1}^{u_1}u^{2b_0}du \nonumber\\
&&\hspace{25mm}\leq \frac{C\s^{(E;0,l)}\cB(\ub_1,u_1)}{(2b_0+1)}\cdot \ub_1^{2a_0}u_1^{2b_0}
\cdot \ub_1 u_1
\label{10.424}
\end{eqnarray}
The last factor being homogeneous of degree 2, 
which is 4 degrees better than borderline, this contribution will come with an extra factor of 
$\delta^4$.  

Consider finally the contribution of the 3rd of the integrals \ref{10.409}. Since $\rhob\sim u^2$ this 
is bounded by a constant multiple of 
\begin{equation}
u^{-2}\int_0^{\ub_1}\|\Lb^\mu E(E^l\beta_\mu-E_N^l\beta_{\mu,N})\|_{L^2(S_{\ub,u})}\ub d\ub
\label{10.425}
\end{equation}
To estimate this we use the fact that:
\begin{equation}
\Lb^\mu E\beta_\mu=E^\mu \Lb\beta_\mu 
\label{10.426}
\end{equation}
(conjugate of \ref{10.419}). This allows us to express: 
\begin{eqnarray}
&&\Lb^\mu E^{l+1}\beta_\mu=E^l(\Lb^\mu E\beta_\mu)+\mbox{terms of order $l$} \nonumber\\
&&\hspace{17mm}=E^l(E^\mu\Lb\beta_\mu)+\mbox{terms of order $l$} \nonumber\\
&&\hspace{17mm}=E^\mu E^l\Lb\beta_\mu+\mbox{terms of order $l$} \nonumber\\
&&\hspace{17mm}=E^\mu\Lb E^l\beta_\mu+\mbox{terms of order $l$} 
\label{10.427}
\end{eqnarray}
Similarly for the $N$th approximants, but in this case there is a error term:
\begin{equation}
E_N^l\omega_{\Lb E,N}={\bf O}(\tau^{N+1})
\label{10.428}
\end{equation}
Writing also $\Lb_N^\mu E_N^{l+1}\beta_{\mu,N}$ as $\Lb^\mu E E_N^l\beta_{\mu,N}$ up to the 0th 
order term \\$(\Lb^\mu E-\Lb_N^\mu E_N)E_N^l\beta_{\mu,N}$, we express in this way 
$$\Lb^\mu E(E^l\beta_\mu-E_N^l\beta_{\mu,N})$$ 
up to lower order terms as:
\begin{equation}
E^\mu\Lb(E^l\beta_\mu-E_N^l\beta_{\mu,N})=\s^{(E;0,l)}\check{\xi}_{\Lb}
\label{10.429}
\end{equation}
by \ref{10.372} for $V=E$. Then in view of \ref{10.403} the integral \ref{10.425}, neglecting lower 
order terms, is bounded by $G_2(\ub_1,u)$. The contribution to $\|\cth_l\|^2_{L^2(\Cb_{\ub_1}^{u_1})}$ 
is then bounded as in \ref{10.405} by:
\begin{equation}
\frac{C\s^{(E;0,l)}\cBb(\ub_1,u_1)}{(2a_0+1)}
\cdot\ub_1^{2a_0} u_1^{2b_0}\cdot\frac{\ub_1^4}{u_1^4}
\label{10.436}
\end{equation}
The last factor being homogeneous of degree 0, which is 2 degrees better than borderline, this 
contribution will come with an extra factor of $\delta^2$. This completes the treatment of the 
contribution to $\|\cth_l\|_{L^2(\Cb_{\ub_1}^{u_1})}$ of the term \ref{10.406}, through the expression 
\ref{10.408}. Comparing \ref{10.417}, \ref{10.424}, \ref{10.436}, we see again that the 
dominant contribution is that of the $\Nb$-component, the least that of the $N$-component, while 
that of the $E$-component is intermediate, the extra factors being homogeneous of degree 0, 2, and 
1 respectively.

We now consider the remaining principal terms in $\check{S}_l$, the right hand side of \ref{10.256}, 
which are contained in 
\begin{equation}
\tilde{R}_l-\tilde{R}_{L,N}
\label{10.437}
\end{equation}
From \ref{10.251}, the principal part of this is of the form:
\begin{eqnarray}
&&\left(\lambdab (E\lambda)c^\mu_{1,E}+\lambda\lambdab\tchi 
c_{2,E}^\mu+\lambda\lambdab c_{3,E}^{\mu\nu}(E\beta_\nu) \right.\nonumber\\
&&\hspace{25mm}\left.+\lambda c_{4,E}^{\mu\nu}(L\beta_\nu)+\lambdab c_{5,E}^{\mu\nu}(\Lb\beta_\nu)\right)E(E^l\beta_\mu-E_N^l\beta_{\mu,N})
\nonumber\\
&&\hspace{5mm}+\left((E\lambda)c_{1,L}^\mu+\lambda\tchi c_{2,L}^\mu 
+\lambda c_{3,L}^{\mu\nu}(E\beta_\nu)+c_{4,L}^{\mu\nu}(\Lb\beta_\nu)\right)
L(E^l\beta_\mu-E_N^l\beta_{\mu,N}) \nonumber\\
&&\hspace{15mm}+\left(\lambdab\tchi c_{1,\Lb}^\mu+\lambdab c_{2,\Lb}^{\mu\nu}(E\beta_\nu)
+c_{3,\Lb}^{\mu\nu}L\beta_\nu\right)\Lb(E^l\beta_\mu-E_N^l\beta_{\mu,N}) \nonumber\\
&&\label{10.438} 
\end{eqnarray}
In reference to the formula \ref{10.367}, we take $X^\mu$ to be each of the coefficients
\begin{eqnarray*}
&&c^\mu_{1,E}, c^\mu_{2,E}, c_{3,E}^{\mu\nu}E\beta_\nu, c_{4,E}^{\mu\nu}L\beta_\nu, 
c_{5,E}^{\mu\nu}\Lb\beta_\nu, \\
&&c^\mu_{1,L}, c^\mu_{2,L}, c_{3,L}^{\mu\nu}E\beta_\nu, 
c_{4,L}^{\mu\nu}\Lb\beta_\nu, \\
&&c^\mu_{1,\Lb}, c_{2,\Lb}^{\mu\nu}E\beta_\nu, c_{3,\Lb}^{\mu\nu}L\beta_\nu
\end{eqnarray*}
Then \ref{10.438} is expressed as a linear combination of 
\begin{eqnarray*}
&&\Nb^\mu E(E^l\beta_\mu-E_N^l\beta_{\mu,N}), \ N^\mu E(E^l\beta_\mu-E_N^l\beta_{\mu,N}), 
\ E^\mu E(E^l\beta_\mu-E_N^l\beta_{\mu,N}) \\
&&\Nb^\mu L(E^l\beta_\mu-E_N^l\beta_{\mu,N}), \ N^\mu L(E^l\beta_\mu-E_N^l\beta_{\mu,N}), 
\ E^\mu L(E^l\beta_\mu-E_N^l\beta_{\mu,N}) \\
&&\Nb^\mu\Lb(E^l\beta_\mu-E_N^l\beta_{\mu,N}), \ N^\mu\Lb(E^l\beta_\mu-E_N^l\beta_{\mu,N}), 
\ E^\mu\Lb(E^l\beta_\mu-E_N^l\beta_{\mu,N})
\end{eqnarray*}
with appropriate coefficients. Assuming then that in ${\cal R}_{\delta,\delta}$ it holds:
\begin{equation}
|E\lambda|\leq C\lambda, \ \ |\tchi|\leq C
\label{10.439}
\end{equation}
and:
\begin{equation}
|E\beta_\mu|\leq C, \ \ |L\beta_\mu|\leq C\lambdab, \ \ |\Lb\beta_\mu|\leq C \ \ \mbox{: $\mu=0,1,2$}
\label{10.440}
\end{equation}
as is valid for the corresponding $N$th approximants, we conclude that \ref{10.438} is pointwise 
bounded by:
\begin{eqnarray}
&&C\lambdab\left\{|E^\mu E(E^l\beta_\mu-E_N^l\beta_{\mu,N})|+|N^\mu E(E^l\beta_\mu-E_N^l\beta_{\mu,N})|
\right.\nonumber\\
&&\hspace{40mm}\left.+|\Nb^\mu E(E^l\beta_\mu-E_N^l\beta_{\mu,N})|\right\} \nonumber\\
&&+C\left\{|E^\mu L(E^l\beta_\mu-E_N^l\beta_{\mu,N})|+|N^\mu L(E^l\beta_\mu-E_N^l\beta_{\mu,N})|
\right.\nonumber\\
&&\hspace{40mm}\left.+|\Nb^\mu L(E^l\beta_\mu-E_N^l\beta_{\mu,N})|\right\} \nonumber\\
&&+C\lambdab\left\{|E^\mu\Lb(E^l\beta_\mu-E_N^l\beta_{\mu,N})|
+|N^\mu\Lb(E^l\beta_\mu-E_N^l\beta_{\mu,N})|\right.\nonumber\\
&&\hspace{40mm}\left.+|\Nb^\mu\Lb(E^l\beta_\mu-E_N^l\beta_{\mu,N})|\right\} 
\label{10.441}
\end{eqnarray}
therefore the contribution of \ref{10.438}, through the integral on the right in \ref{10.361}, 
to $\|\cth_l\|_{L^2(S_{\ub_1,u})}$ is bounded by a constant multiple of:
\begin{eqnarray}
&&\int_0^{\ub_1}\left\{\|E^\mu E(E^l\beta_\mu-E_N^l\beta_{\mu,N})\|_{L^2(S_{\ub,u})}
+\|N^\mu E(E^l\beta_\mu-E_N^l\beta_{\mu,N})\|_{L^2(S_{\ub,u})}\right.\nonumber\\
&&\hspace{50mm}\left.+\|\Nb^\mu E(E^l\beta_\mu-E_N^l\beta_{\mu,N})\|_{L^2(S_{\ub,u})}
\right\}\ub d\ub\nonumber\\
&&+\int_0^{\ub_1}\left\{\|E^\mu L(E^l\beta_\mu-E_N^l\beta_{\mu,N})\|_{L^2(S_{\ub,u})}
+\|N^\mu L(E^l\beta_\mu-E_N^l\beta_{\mu,N})\|_{L^2(S_{\ub,u})}\right.\nonumber\\
&&\hspace{50mm}\left.+\|\Nb^\mu L(E^l\beta_\mu-E_N^l\beta_{\mu,N})\|_{L^2(S_{\ub,u})}
\right\}d\ub\nonumber\\
&&+\int_0^{\ub_1}\left\{\|E^\mu\Lb(E^l\beta_\mu-E_N^l\beta_{\mu,N})\|_{L^2(S_{\ub,u})}
+\|N^\mu\Lb(E^l\beta_\mu-E_N^l\beta_{\mu,N})\|_{L^2(S_{\ub,u})}\right.\nonumber\\
&&\hspace{50mm}\left.+\|\Nb^\mu\Lb(E^l\beta_\mu-E_N^l\beta_{\mu,N})\|_{L^2(S_{\ub,u})}
\right\}\ub d\ub\nonumber\\
&&\label{10.442}
\end{eqnarray}
Here the first integral coincides with the sum \ref{10.409}, which has already been treated. 
The third integral is similar to what is contributed by \ref{10.369}, but whereas \ref{10.369} 
contributes with an extra $\lambda^{-1}\sim u^{-2}$ factor, here this factor is missing (see  
\ref{10.370}). As a result we obtain a bound for the contribution of this integral to 
$\|\cth_l\|^2_{L^2(\Cb_{\ub_1}^{u_1})}$ by:
\begin{equation}
\frac{C\ub_1^{2a_0}u_1^{2b_0}}{(2a_0+1)} \left(\s^{(Y;0,l)}\cBb(\ub_1,u_1)\frac{\ub_1^4}{u_1^2}
+\s^{(E;0,l)}\cBb(\ub_1,u_1)\ub_1^3\right)
\label{10.443}
\end{equation}
which is 4 degrees better than borderline. 

What remains to be considered is the 2nd of the integrals \ref{10.442}. The contribution to this 
integral of the 1st term in the integrant is:
\begin{eqnarray}
&&\int_0^{\ub_1}\|\s^{(E;0,l)}\check{\xi}_L\|_{L^2(S_{\ub,u})}d\ub \leq 
\sqrt{\ub_1}\|\s^{(E;0,l)}\check{\xi}_L\|_{L^2(C_u^{\ub_1})} \nonumber\\
&&\hspace{30mm}\leq\sqrt{\ub_1}\sqrt{\s^{(E;0,l)}\cE^{\ub_1}(u)} \nonumber\\
&&\hspace{30mm}\leq \sqrt{\s^{(E;0,l)}\cB(\ub_1,u_1)}\ub_1^{a_0+\frac{1}{2}}u^{b_0} 
\label{10.444}
\end{eqnarray}
by \ref{9.290} and \ref{9.297}. This bounds the contribution to $\|\cth_l\|_{L^2(S_{\ub_1,u})}$. 
The corresponding contribution to $\|\cth_l\|^2_{L^2(\Cb_{\ub_1}^{u_1})}$ is then bounded by:
\begin{equation}
\frac{\s^{(E;0,l)}\cB(\ub_1,u_1)}{(2b_0+1)}\cdot\ub_1^{2a_0}u_1^{2b_0}\cdot\ub_1 u_1
\label{10.445}
\end{equation}
which is 4 degrees better than borderline. Next, we consider the contribution to the 2nd of the 
integrals \ref{10.442} of the 3rd term in the integrant, that is we consider the integral:
\begin{equation}
\int_0^{\ub_1}\|\Nb^\mu L(E^l\beta_\mu-E_N^l\beta_{\mu,N})\|_{L^2(S_{\ub,u})}d\ub
\label{10.446}
\end{equation}
To estimate this integral, we use the fact that:
\begin{equation}
\Nb^\mu L\beta_\mu=\lambdab E^\mu E\beta_{\mu}
\label{10.447}
\end{equation}
This allows us to express:
\begin{eqnarray}
&&\Nb^\mu LE^l\beta_\mu=\Nb^\mu E^l L\beta_\mu+\mbox{terms of order $l$} \nonumber\\
&&\hspace{15mm}=E^l(\Nb^\mu L\beta_\mu)+\mbox{terms of order $l$} \nonumber\\
&&\hspace{15mm}=E^l(\lambdab E^\mu E\beta_\mu)+\mbox{terms of order $l$} \nonumber\\
&&\hspace{15mm}=\lambdab E^\mu E^{l+1}\beta_\mu+\mbox{terms of order $l$} 
\label{10.448}
\end{eqnarray}
(\ref{10.447} and \ref{10.448} are the conjugates of \ref{10.374} and \ref{10.375}). 
Similarly for the $N$th approximants, but in this case there is an error term:
\begin{equation}
E_N^l\left(\rhob_N^{-1}(-\delta_N^\prime+\omega_{L\Lb,N})\right)={\bf O}(u^{-2}\tau^N)
\label{10.449}
\end{equation}
Writing also $\Nb_N^\mu L_N E_N^l\beta_{\mu,N}$ as $\Nb^\mu L E_N^l\beta_{\mu,N}$ up to the 0th 
order term \\$(\Nb^\mu L-\Nb_N^\mu L_N)E_N^l\beta_{\mu,N}$, we express in this way 
$$\Nb^\mu L(E^l\beta_\mu-E_N^l\beta_{\mu,N})$$ 
up to lower order terms as:
\begin{equation}
\lambdab E^\mu E(E^l\beta_\mu-E_N^l\beta_{\mu,N})=\lambdab\s^{(E;0,l)}\check{\sxi}
\label{10.450}
\end{equation}
Since $\lambdab\sim\ub$ we can then estimate the integral \ref{10.446}, up to lower order terms, 
in terms of the integral:
\begin{equation}
\int_0^{\ub_1}\|\s^{(E;0,l)}\check{\sxi}\|_{L^2(S_{\ub,u})}\ub d\ub=\tilde{G}_1(\ub,u)
\label{10.451}
\end{equation}
(see \ref{10.410}). Therefore the contribution to $\|\cth_l\|^2_{L^2(\Cb_{\ub_1}^{u_1})}$ is bounded 
by \ref{10.417}, which is 3 degrees better than borderline. Finally, we have the contribution to the 
2nd of the integrals \ref{10.442} of the 2nd term in the integrant:
\begin{equation}
\int_0^{\ub_1}\|N^\mu L(E^l\beta_\mu-E_N^l\beta_{\mu,N})\|_{L^2(S_{\ub,u})}d\ub
\label{10.452}
\end{equation}
To estimate this we write:
\begin{equation}
N^\mu L(E^l\beta_\mu-E_N^l\beta_{\mu,N})=\s^{(Y;0,l)}\check{\xi}_L
-\ogamma\Nb^\mu L(E^l\beta_\mu-E_N^l\beta_{\mu,N})
\label{10.453}
\end{equation}
The partial contribution to the integral \ref{10.452} of the 2nd term on the right is similar to the 
integral \ref{10.446}, but here, since $\ogamma\sim u$, there is an extra $u$ factor. So this 
contribution is bounded in terms of $u\tilde{G}_1(\ub_1,u)$, and the corresponding contribution 
to $\|\cth_l\|^2_{L^2(\Cb_{\ub_1}^{u_1})}$ by:
$$\frac{C\s^{(E;0,l)}\cBb(\ub_1,u_1)}{(2a_0+1)}
\cdot\ub_1^{2a_0} u_1^{2b_0}\cdot\ub_1^3$$
which is 5 degrees better than borderline. As for the partial contribution to the integral 
\ref{10.452} of the 1st term on the right in \ref{10.453}, this is:
\begin{eqnarray}
&&\int_0^{\ub_1}\|\s^{(Y;0,l)}\check{\xi}_L\|_{L^2(S_{\ub,u})}d\ub
\leq \sqrt{\ub_1}\|\s^{(Y;0,l)}\check{\xi}_L\|_{L^2(C_u^{\ub_1})} \nonumber\\
&&\hspace{30mm}\leq \sqrt{\ub_1}\sqrt{\s^{(Y;0,l)}\cE^{\ub_1}(u)} \nonumber\\
&&\hspace{30mm}\leq \sqrt{\s^{(Y;0,l)}\cB(\ub_1,u_1)}\ub_1^{a_0+\frac{1}{2}}u^{b_0}
\label{10.454}
\end{eqnarray}
This bounds the contribution to $\|\cth_l\|_{L^2(S_{\ub_1,u})}$. The corresponding contribution to 
$\|\cth_l\|^2_{L^2(\Cb_{\ub_1}^{u_1})}$ is then bounded by:
\begin{equation}
\frac{\s^{(Y;0,l)}\cB(\ub_1,u_1)}{(2b_0+1)}\cdot\ub_1^{2a_0}u_1^{2b_0}\cdot \ub_1 u_1
\label{10.455}
\end{equation}
which is 4 degrees better than borderline. 

We remark that \ref{10.424}, \ref{10.445} and \ref{10.455} while not the dominant of the contributions 
below the borderline, the last factor in these being $\ub_1 u_1$ while in \ref{10.399} the last factor is $\ub_1^3 u_1^{-4}$, the last factor in these three dominates near $\Cb_0$ 
having the least power of $\ub_1$. Now, while the only way to estimate \ref{10.452} is through 
\ref{10.453} leading to the estimate \ref{10.455}, \ref{10.418} can alternatively be estimated by 
writing: 
\begin{equation}
N^\mu E(E^l\beta_\mu-E_N^l\beta_{\mu,N})=\s^{(Y;0,l)}\check{\sxi}-\ogamma\Nb^\mu E(E^l\beta_\mu-
E_N^l\beta_{\mu,N})
\label{10.a3}
\end{equation}
and using the fact that 
\begin{eqnarray}
&&\int_{\ub}^u\|\s^{(Y;0,l)}\check{\sxi}\|^2_{L^2(S_{\ub,u^\prime})}\ub u^{\prime 2}du^\prime
\leq C\|\sqrt{a}\s^{(Y;0,l)}\check{\sxi}\|^2_{L^2(\Cb_{\ub}^u)}\nonumber\\
&&\hspace{18mm}\leq C\s^{(Y;0,l)}\cEb^u(\ub)\leq C\s^{(Y;0,l)}\cBb(\ub_1,u_1)\ub^{2a_0} u^{2b_0}
\label{10.a4}
\end{eqnarray}
The contribution of the 1st term on the right in \ref{10.a3} through \ref{10.418} to 
$\|\cth_l\|^2_{L^2(\Cb_{\ub_1}^{u_1})}$ is then shown to be bounded by:
\begin{equation}
\frac{C\s^{(Y;0,l)}\cBb(\ub_1,u_1)}{(2a_0+1)}\cdot\ub_1^{2a_0} u_1^{2b_0}\cdot\frac{\ub_1^3}{u_1^2}
\label{10.a5}
\end{equation}
Here the last factor is homogeneous of degree 1 which is worse than the last factor of \ref{10.424}, 
however the last factor here is better than the last factor in \ref{10.424} in a suitable 
neighborhood of $\Cb_0$. The contribution of the 2nd term on the right in \ref{10.a3} to 
$\|\cth_l\|_{L^2(S_{\ub_1,u})}$ is 
similar to \ref{10.425} but with a $u$ factor, since $\ogamma\sim u$. The corresponding contribution 
to $\|\cth_l\|^2_{L^2(\Cb_{\ub_1}^{u_1})}$ is then bounded as in \ref{10.436} but with a last factor 
of $\ub_1^4/u_1^2$. A similar estimate is obtained in place of \ref{10.445} for the contribution of 
$$\int_0^{\ub_1}\|E^\mu L(E^l\beta_\mu-E_N^l\beta_{\mu,N}\|_{L^2(S_{\ub,u})}d\ub$$
if we use  
$$E^\mu L\beta_\mu=L^\mu E\beta_\mu=\rho N^\mu E\beta_\mu$$
to express $E^\mu L(E^l\beta_\mu-E_N^l\beta_{\mu,N})$ up to lower order terms as 
$$\rho N^\mu E(E^l\beta_\mu-E_N^l\beta_{\mu,N})=\rho\left(\s^{(Y;0,l)}\check{\sxi}
-\ogamma\Nb^\mu E(E^l\beta_\mu-E_N^l\beta_{\mu,N})\right)$$

We summarize the above results in the following proposition. 

\vspace{2.5mm}

\noindent {\bf Proposition 10.3} \ \ Under the assumptions \ref{10.439}, \ref{10.440} the 
top order acoustical difference quantity $\cth_l$, $l=n$, satisfies to principal terms  the following estimate on the $\Cb_{\ub}$:
\begin{eqnarray*}
&&\|\cth_l\|_{L^2(\Cb_{\ub_1}^{u_1})}
\leq \ub_1^{a_0}u_1^{b_0}\left\{C\sqrt{\frac{\s^{(Y;0,l)}\cBb(\ub_1,u_1)}{2a_0+1}}
\cdot\frac{\ub_1^2}{u_1^3}\right.\\
&&\hspace{25mm}\left.+C\sqrt{\max\{\s^{(0,l)}\cB(\ub_1,u_1),\s^{(0,l)}\cBb(\ub_1,u_1)\}}
\cdot\frac{\ub_1^{1/2}}{u_1}\right\}
\end{eqnarray*}

\vspace{2.5mm}

Actually, the coefficient of the second term in parenthesis in the above estimate is of the form:
\begin{equation}
\left(\frac{1}{\sqrt{2a_0+1}}+\frac{1}{\sqrt{2b_0+1}}\right)C
\label{10.c1}
\end{equation}
but this refinement is not needed in the following. 

\vspace{2.5mm}

As we shall see in Section 10.8, the quantity $\cthb_l$ satisfies a boundary condition on ${\cal K}$ 
which involves the quantity $\cth_l$. To control appropriately the boundary values of $\cthb_l$ 
on ${\cal K}$ we need to derive an appropriate estimate for $\cth_l$ in $L^2({\cal K}^\tau)$. 
Setting $\ub_1=u$ in \ref{10.361} we have:
\begin{equation}
\|\cth_l\|_{L^2(S_{u,u})}\leq k\int_0^u\|\hat{S}_l\|_{L^2(S_{\ub,u})}d\ub
\label{10.456}
\end{equation}
and what we want to estimate is:
\begin{equation}
\|\cth_l\|^2_{L^2({\cal K}^\tau)}=\int_0^\tau \|\cth_l\|^2_{L^2(S_{u,u})}du
\label{10.457}
\end{equation}
We shall consider the leading contribution to the integral on the right in \ref{10.456}, which 
is that of the 3rd term in parenthesis in \ref{10.369} and comes from \ref{10.364} through  \ref{10.365}. This contribution is bounded by \ref{10.370} with $\ub_1=u$, that is by:
\begin{equation}
Cu^{-2}\int_0^u\|\Nb^\mu\Lb(E^l\beta_\mu-E_N^l\beta_{\mu,N})\|_{L^2(S_{\ub,u})}\ub d\ub
\label{10.458}
\end{equation}
Recalling \ref{10.379} - \ref{10.383}, \ref{10.458} is bounded by:
\begin{equation}
C(G_0(u,u)+G_1(u,u))
\label{10.459}
\end{equation}
Setting $\ub_1=u$ in \ref{10.385} we obtain:
\begin{equation}
G_0^2(u,u)\leq \frac{1}{3}u^{-1}\int_0^u g_0^2(\ub,u)d\ub
\label{10.460}
\end{equation}
Hence:
\begin{eqnarray}
&&\int_0^\tau G_0^2(u,u)du\leq\frac{1}{3}
\int_0^\tau\left\{\int_0^u g_0^2(\ub,u)d\ub\right\}\frac{du}{u}\nonumber\\
&&\hspace{25mm}=\frac{1}{3}\int_0^\tau
\left\{\int_{\ub}^\tau u^{-1}g_0^2(\ub,u)du\right\}d\ub
\label{10.461}
\end{eqnarray}
(reversing the order of integration). In regard to the interior integral we have, for all 
$\ub\in[0,\tau]$:
\begin{eqnarray}
&&\int_{\ub}^{\tau}u^{-1}g_0^2(\ub,u)du
=\int_{\ub}^{\tau}\frac{\partial}{\partial u}\left(\int_{\ub}^u u^{\prime 2}g_0^2(\ub,u^\prime)
du^\prime\right)u^{-3}du\nonumber\\
&&\hspace{10mm}=\tau^{-3}\int_{\ub}^{\tau}u^{\prime 2}g_0^2(\ub,u^\prime)du^\prime
+3\int_{\ub}^{\tau}\left(\int_{\ub}^u u^{\prime 2}g_0^2(\ub,u^\prime)du^\prime\right)u^{-4}du 
\nonumber\\
&&\label{10.462}
\end{eqnarray}
By \ref{10.380} with $(u,u^\prime)$ in the role of $(u_1,u)$ we have, for all $u\in[\ub,\tau]$:
\begin{equation}
\int_{\ub}^u u^{\prime 2} g_0^2(\ub,u^\prime)du^\prime\leq 
C\s^{(Y;0,l)}\cEb^u(\ub)\leq C\s^{(Y;0,l)}\cBb(\tau,\tau)\cdot \ub^{2a_0} u^{2b_0}
\label{10.463}
\end{equation}
in view of the definition \ref{9.298}. Then \ref{10.462} is bounded by:
\begin{equation}
C\s^{(Y;0,l)}\cBb(\tau,\tau)\ub^{2a_0}\left\{\tau^{2b_0-3}
+3\frac{(\tau^{2b_0-3}-\ub^{2b_0-3})}{(2b_0-3)}\right\}
\label{10.464}
\end{equation}
In view of \ref{10.391} we have a bound of \ref{10.462} by:
\begin{equation}
C\s^{(Y,0,l)}\cBb(\tau,\tau)\ub^{2a_0} \tau^{2b_0-3}
\label{10.465}
\end{equation}
Then by \ref{10.461}:
\begin{equation}
\int_0^{\tau}G_0^2(u,u)du\leq\frac{C\s^{(Y;0,l)}\cBb(\tau,\tau)}{(2a_0+1)}
\cdot\tau^{2(a_0+b_0)}\cdot\tau^{-2}
\label{10.466}
\end{equation}
The last factor $\tau^{-2}$ being of degree -2 signifies that 
the contribution is {\em borderline}. We shall also estimate the contribution of $G_1$. 
Setting $\ub_1=u$ in \ref{10.394} we obtain:
\begin{equation}
G_1^2(u,u)\leq\frac{1}{2}u^{-2}\int_0^u \ub g_1^2(\ub,u)d\ub 
\label{10.467}
\end{equation}
Hence:
\begin{eqnarray}
&&\int_0^{u_1}G_1^2(u,u)du\leq\frac{1}{2}
\int_0^\tau\left\{\int_0^u\ub g_1^2(\ub,u)d\ub\right\}\frac{du}{u^2}\nonumber\\
&&\hspace{25mm}=\frac{1}{2}\int_0^\tau
\left\{\int_{\ub}^\tau \ub u^{-2} g_1^2(\ub,u)du\right\}d\ub
\label{10.468}
\end{eqnarray}
(reversing the order of integration). In regard to the interior integral we have, for all 
$\ub\in[0,\tau]$:
\begin{eqnarray}
&&\int_{\ub}^{\tau}\ub u^{-2} g_1^2(\ub,u)du
=\int_{\ub}^{\tau}\frac{\partial}{\partial u}\left(\int_{\ub}^u \ub g_1^2(\ub,u^\prime)
du^\prime\right)u^{-2}du\nonumber\\
&&\hspace{10mm}=\tau^{-2}\int_{\ub}^{\tau}\ub g_1^2(\ub,u^\prime)du^\prime
+2\int_{\ub}^{\tau}\left(\int_{\ub}^u \ub g_1^2(\ub,u^\prime)du^\prime\right)u^{-3}du 
\nonumber\\
&&\label{10.469}
\end{eqnarray}
By \ref{10.381} with $(u,u^\prime)$ in the role of $(u_1,u)$ we have, for all $u\in[\ub,\tau]$: 
\begin{equation}
\int_{\ub}^u \ub g_1^2(\ub,u^\prime)du^\prime \leq C\s^{(E;0,l)}\cEb^u(\ub)
\leq C\s^{(E;0,l)}\cBb(\tau,\tau)\cdot \ub^{2a_0} u^{2b_0}
\label{10.470}
\end{equation}
in view of the definition \ref{9.298}. Then \ref{10.469} is bounded by:
\begin{eqnarray}
&&C\s^{(E;0,l)}\cBb(\tau,\tau)\ub^{2a_0}\left\{\tau^{2b_0-2}
+2\frac{(\tau^{2b_0-2}-\ub^{2b_0-2})}{(2b_0-2)}\right\}\nonumber\\
&&\leq C\s^{(E;0,l)}\cBb(\tau,\tau)\ub^{2a_0}\tau^{2b_0-2}
\label{10.471}
\end{eqnarray}
(in view of \ref{10.391}), and by \ref{10.468}:
\begin{equation}
\int_0^{\tau}G_1^2(u,u)du\leq \frac{C\s^{(E;0,l)}\cBb(\tau,\tau)}{(2a_0+1)}
\cdot\tau^{2(a_0+b_0)}\cdot\tau^{-1}
\label{10.472}
\end{equation}
The last factor $\tau^{-1}$ being of degree -1, which is one degree better than 
borderline, this contribution to $\|\cth\|^2_{L^2({\cal K}^\tau)}$ will come with an extra factor 
of $\delta$. 

Finally, all other contributions to $\|\cth\|^2_{L^2({\cal K}^\tau)}$ will come with extra factors 
of $\delta^2$ at least. We thus deduce the following proposition. 

\vspace{2.5mm}

\noindent {\bf Proposition 10.4} \ \ Under the assumptions \ref{10.439}, \ref{10.440} the 
top order acoustical difference quantity $\cth_l$, $l=n$, satisfies to principal terms the following estimate on ${\cal K}$:
\begin{eqnarray*}
&&\|\cth_l\|_{L^2({\cal K}^\tau)}
\leq \tau^{a_0+b_0}\left\{C\sqrt{\frac{\s^{(Y;0,l)}\cBb(\tau,\tau)}{2a_0+1}}
\cdot\frac{1}{\tau}\right.\\
&&\hspace{25mm}\left.+C\sqrt{\max\{\s^{(0,l)}\cB(\tau,\tau),\s^{(0,l)}\cBb(\tau,\tau)\}}
\cdot\frac{1}{\tau^{1/2}}\right\}
\end{eqnarray*}

\vspace{5mm}

\section{Estimates for $\cnu_{m-1,l+1}$}

Consider now the propagation equation \ref{10.336} for $\cnu_{m,l}$. Let us denote:
\begin{equation}
\check{\tau}_{m,l}=E^{l-1}T^m(\lambda cq^\prime E\pi)-E_N^{l-1}T^m(\lambda_N c_N q^\prime_N E_N\pi_N)
\label{10.473}
\end{equation}
This is a quantity of order $m+l$. 
Also, let us add to both sides of equation \ref{10.336} the term:
\begin{equation}
\left(\frac{2L\lambda}{\lambda}-(l+1)\chi\right)\check{\tau}_{m,l}
\label{10.474}
\end{equation}
Denoting the resulting right hand side with the term \ref{10.474} added by $\check{M}_{m,l}$, the  equation reads:
\begin{equation}
L(\cnu_{m,l}-\check{\tau}_{m,l})-\left(\frac{2L\lambda}{\lambda}-(l+1)\chi\right)
(\cnu_{m,l}-\check{\tau}_{m,l})
=\check{M}_{m,l}
\label{10.475}
\end{equation} 
Actually, comparing \ref{10.235} with \ref{10.4} we see that what is involved 
in the $(m,l)$ difference energy estimates for $m\geq 1$ is $\cnu_{m-1,l+1}$. 
Being at the top order, we have $m+l=n$, $m=1,...,n$. Thus we consider 
equation \ref{10.475} with $(m-1,l+1)$ in the role of $(m,l)$:
\begin{eqnarray}
&&L(\cnu_{m-1,l+1}-\check{\tau}_{m-1,l+1})-\left(\frac{2L\lambda}{\lambda}-(l+2)\chi\right)(\cnu_{m-1,l+1}-\check{\tau}_{m-1,l+1})\nonumber\\
&&\hspace{50mm}=\check{M}_{m-1,l+1} \nonumber\\
&&\hspace{30mm} : \ m+l=n, \ m=1,...,n
\label{10.476}
\end{eqnarray} 
This is exactly of the same form as the equation preceding \ref{10.338}, with the quantity 
$\cnu_{m-1,l+1}-\check{\tau}_{m-1,l+1}$ in the role of the quantity $\cth_l$ and $\check{M}_{m-1,l+1}$ 
in the role of $\check{S}_l$. Equation \ref{10.476} can then be written, in analogy with \ref{10.347},  
in the $(\ub,u,\vartheta^\prime)$ coordinates, which are adapted to the flow of $L$, in the form:
\begin{equation}
\frac{\partial}{\partial\ub}\left(\sh^{\prime(l+2)/2}\lambda^{-2}(\cnu_{m-1,l+1}-\check{\tau}_{m-1,l+1})\right)=\sh^{\prime(l+2)/2}\lambda^{-2}\check{M}_{m-1,l+1}
\label{10.477}
\end{equation}
Integrating from $\Cb_0$, noting that:
\begin{equation}
\left.\cnu_{m-1,l+1}\right|_{\Cb_0}=\left.\check{\tau}_{m-1,l+1}\right|_{\Cb_0}=0
\label{10.478}
\end{equation}
we obtain:
\begin{eqnarray}
&&\left(\sh^{\prime(l+2)/2}\lambda^{-2}(\cnu_{m-1,l+1}-\check{\tau}_{m-1,l+1})
\right)(\ub_1,u,\vartheta^\prime)
\nonumber\\
&&\hspace{40mm}=\int_0^{\ub_1}\left(\sh^{\prime(l+2)/2}\lambda^{-2}\check{M}_{m-1,l+1}\right)
(\ub,u,\vartheta^\prime)d\ub \nonumber\\
&&\label{10.479}
\end{eqnarray}
Following the argument leading from \ref{10.349} to the inequality \ref{10.361} we deduce from 
\ref{10.479}, in the same way, the inequality:
\begin{equation}
\|\cnu_{m-1,l+1}-\check{\tau}_{m-1,l+1}\|_{L^2(S_{\ub_1,u})}\leq 
k\int_0^{\ub_1}\|\check{M}_{m-1,l+1}\|_{L^2(S{\ub,u})}d\ub
\label{10.480}
\end{equation}

The quantity $\check{\tau}_{m-1,l+1}$ is of order $m+l=n$. We shall presently estimate the 
contribution of the principal terms in $\check{M}_{m-1,l+1}$ to the integral on the right in 
\ref{10.480}. These principal terms, which are of order $m+l+1=n+1$, are (see \ref{10.336}) 
the two difference terms: 
\begin{equation}
\frac{2L\lambda}{\lambda}\left(E^l T^{m-1}j-E_N^l T^{m-1}j_N\right)
+\tilde{K}^{\prime\prime}_{m-1,l+1}-\tilde{K}^{\prime\prime}_{m-1,l+1,N}
\label{10.481}
\end{equation}
as well as the terms:
\begin{equation}
-(m-1)\zeta\cnu_{m-2,l+2}+A\left\{(2\pi\lambda)^{m-1}\cth_{l+m}+2\sum_{i=0}^{m-2}(2\pi\lambda)^{m-2-i}
\cnu_{i,l+m-i}\right\}
\label{10.482}
\end{equation}
which involve the top order acoustical quantities entering the $(i,l+m-i)$ difference  energy estimates 
for $i=0,...,m-1$. In fact, in view of \ref{10.363}, the first of the two difference terms 
\ref{10.481} makes the leading contribution. This comes from the 1st term in $j$ as given by 
Proposition 10.2, that is from: 
\begin{equation}
\frac{1}{4}\beta_N^2(E^l T^{m-1}\Lb^2 H-E_N^l T^{m-1}\Lb_N^2 H_N)
\label{10.483}
\end{equation}
more precisely from the principal part of this which is (compare with \ref{10.365}):
\begin{equation}
-\frac{1}{2}\beta_N^2 H^\prime(g^{-1})^{\mu\nu}\beta_\nu
(E^l T^{m-1}\Lb^2\beta_\mu-E_N^l T^{m-1}\Lb^2_N\beta_{\mu,N})
\label{10.484}
\end{equation}
Here, recalling that
$$T=L+\Lb=L_N+\Lb_N$$
we write:
\begin{equation}
\Lb^2\beta_\mu=\Lb T\beta_\mu-\Lb L\beta_\mu, \ \ \ 
\Lb_N^2\beta_{\mu,N}=\Lb_N T\beta_{\mu,N}-\Lb_N L_N\beta_{\mu,N}
\label{10.485}
\end{equation}

To treat the 2nd terms in each of \ref{10.485} we recall from Section 9.5 that 
$$\Omega a\square_{\tilde{h}}\beta_\mu=0, \ \ \ 
\Omega_N a_N\square_{\tilde{h}^\prime_N}\beta_{\mu,N}=\tilde{\kappa}^\prime_{\mu,N}$$
where $\tilde{\kappa}^\prime_{\mu,N}$ satisfies the estimate of Proposition 9.10. Also 
(see \ref{9.230}), that $\Omega a\square_{\tilde{h}}\beta_\mu$ is equal to $a\square_h\beta_\mu$ 
up to 1st order terms, and $\Omega_N a_N\square_{\tilde{h}^\prime_N}\beta_{\mu,N}$ is equal to 
$a_N\square_{h^\prime_N}\beta_{\mu,N}$ up to the corresponding $N$th approximant terms. 
Moreover (see \ref{9.205}, \ref{9.210}), that $a\square_h\beta_\mu$ is equal to 
$-\Lb L\beta_\mu+aE^2\beta_\mu$ up to 1st order terms, and $a_N \square_{h^\prime_N}\beta_{\mu,N}$ 
is equal to $-\Lb_N L_N\beta_{\mu,N}+a_N E_N^2\beta_{\mu,N}$ up to the corresponding $N$th approximant 
terms. Consequently, we can express $\Lb L\beta_\mu$ as $a E^2\beta_\mu$ up to 1st order terms,  
and we can express $\Lb_N L_N\beta_{\mu,N}$ as $a_N E_N^2\beta_{\mu,N}$ up to the corresponding 
$N$th approximant terms minus the error term $\tilde{\kappa}^\prime_{\mu,N}$ which satisfies the 
estimate of Proposition 9.10. 

The contribution of the 1st terms in \ref{10.485} to \ref{10.484} is, to principal terms, 
\begin{equation}
-\frac{1}{2}\beta_N^2 H^\prime(g^{-1})^{\mu\nu}\beta_\nu\Lb(E^l T^m\beta_\mu-E_N^l T^m\beta_{\mu,N})
\label{10.486}
\end{equation}
while the contribution of the 2nd terms in \ref{10.485} to \ref{10.484} is, according to the above, 
given to principal terms by:
\begin{equation}
\frac{1}{2}\beta_N^2 H^\prime(g^{-1})^{\mu\nu}\beta_\nu a E(E^{l+1}T^{m-1}\beta_\mu
-E_N^{l+1}T^{m-1}\beta_{\mu,N})
\label{10.487}
\end{equation}
It is \ref{10.486} which makes the leading contribution to the integral on the right in \ref{10.480}. 
In view of \ref{10.368} this contribution is bounded by a constant multiple of 
\begin{eqnarray}
&&u^{-2}\int_0^{\ub_1}\left\{\|E^\mu\Lb(E^l T^m\beta_\mu-E_N^l T^m\beta_{\mu,N})\|_{L^2(S_{\ub,u})}
\right. \label{10.488}\\
&&\hspace{16mm}+\|N^\mu\Lb(E^l T^m\beta_\mu-E_N^l T^m\beta_{\mu,N})\|_{L^2(S_{\ub,u})} \nonumber\\
&&\hspace{16mm}\left.+\|\Nb^\mu\Lb(E^l T^m\beta_\mu-E_N^l T^m\beta_{\mu,N}\|_{L^2(S_{\ub,u})}\right\}
\ub d\ub \nonumber
\end{eqnarray}
As usual, the dominant contribution is that of the $\Nb$ component, 3rd term in the integrant, 
which as we shall see is borderline. 

We begin the estimates with the $E$ component, 1st term in the integrant. According to the definition \ref{9.275} for each of the variation fields $V=Y,E$ we have:
\begin{equation}
V^\mu d(E^l T^m\beta_\mu-E_N^l T^m\beta_{\mu,N})=\s^{(V;m,l)}\check{\xi}
\label{10.489}
\end{equation}
In particular taking $V=E$ we have:
\begin{equation}
E^\mu\Lb(E^l T^m\beta_\mu-E_N^l T^m\beta_{\mu,N})=\s^{(E;m,l)}\check{\xi}_{\Lb}
\label{10.490}
\end{equation}
Therefore the 1st term in the integrant in \ref{10.488} contributes $G(\ub_1,u)$, where, as in 
\ref{10.383}, 
$$G(\ub_1,u)=u^{-2}\int_0^{\ub_1}g(\ub,u)\ub d\ub$$ 
and here: 
\begin{equation}
g(\ub,u)=\|\s^{(E;m,l)}\check{\xi}_{\Lb}\|_{L^2(S_{\ub,u})}
\label{10.491}
\end{equation}
so in view of the definition \ref{9.291} we have:
\begin{equation}
\int_{\ub}^{u_1}g^2(\ub,u)du\leq C\s^{(E;m,l)}\cEb^{u_1}(\ub)
\label{10.492}
\end{equation}
As in \ref{10.386} it holds:
\begin{equation}
\int_{\ub_1}^{u_1}G^2(\ub_1,u)du\leq \frac{1}{3}\ub_1^3\int_0^{\ub_1}\left\{\int_{\ub_1}^{u_1}
u^{-4}g^2(\ub,u)du\right\}d\ub
\label{10.493}
\end{equation}
Here, in regard to the interior integral we have, for all $\ub\in[0,\ub_1]$: 
\begin{eqnarray}
&&\int_{\ub_1}^{u_1}u^{-4}g^2(\ub,u)du
=\int_{\ub_1}^{u_1}\frac{\partial}{\partial u}\left(\int_{\ub_1}^u g^2(\ub,u^\prime)
du^\prime\right)u^{-4}du\nonumber\\
&&\hspace{10mm}=u_1^{-4}\int_{\ub_1}^{u_1}g^2(\ub,u^\prime)du^\prime
+4\int_{\ub_1}^{u_1}\left(\int_{\ub_1}^u g^2(\ub,u^\prime)du^\prime\right)u^{-5}du 
\nonumber\\
&&\label{10.494}
\end{eqnarray}
Since $\ub_1\geq \ub$, by \ref{10.492} we have, for all $u\in[\ub_1,u_1]$:
\begin{eqnarray}
&&\int_{\ub_1}^u g^2(\ub,u^\prime)du^\prime\leq 
\int_{\ub}^u g^2(\ub,u^\prime)du^\prime \nonumber\\
&&\hspace{15mm}\leq C\s^{(E;m,l)}\cEb^u(\ub)\leq C\s^{(E;m,l)}\cBb(\ub_1,u_1)\cdot \ub^{2a_m} u^{2b_m}
\nonumber\\
&&\label{10.495}
\end{eqnarray}
in view of the definition \ref{9.298}. Then assuming that $b_m>2$, \ref{10.494} is bounded by:
$$C\s^{(E;m,l)}\cBb(\ub_1,u_1)\ub^{2a_m}\left\{u_1^{2b_m-4}+4\frac{(u_1^{2b_m-4}-\ub_1^{2b_m-4})}{(2b_m-4)}\right\}$$
or, taking $b_m\geq 3$, simply by:
$$C\s^{(E;m,l)}\cBb(\ub_1,u_1)\ub^{2a_m}u_1^{2b_m-4}$$
Then by \ref{10.493}:
\begin{equation}
\int_{\ub_1}^{u_1}G^2(\ub_1,u)du\leq \frac{C\s^{(E;m,l)}\cBb(\ub_1,u_1)}{(2a_m+1)}
\cdot\ub_1^{2a_m} u_1^{2b_m}\cdot\frac{\ub_1^4}{u_1^4}
\label{10.496}
\end{equation}
(compare with \ref{10.405}). 
This bounds the contribution to $\|\cnu_{m-1,l+1}-\check{\tau}_{m-1,l+1}\|^2_{L^2(\Cb_{\ub_1}^{u_1})}$. The last factor 
being homogeneous of degree 0 signifies that this contribution is 2 degrees better than borderline, 
hence will come with an extra factor of $\delta^2$. 

We turn to the contribution of the 2nd term in the integrant in \ref{10.488}, the $N$ component. 
Here we use \ref{10.374} which in analogy with \ref{10.375} allows us to express:
\begin{eqnarray}
&&N^\mu\Lb E^l T^m\beta_\mu=N^\mu E^l T^m\Lb\beta_\mu+\mbox{terms of order $l+m=n$} \nonumber\\
&&\hspace{22mm}=E^l T^m(N^\mu\Lb\beta_\mu)+\mbox{terms of order $n$} \nonumber\\
&&\hspace{22mm}=E^l T^m(\lambda E^\mu E\beta_\mu)+\mbox{terms of order $n$} \nonumber\\
&&\hspace{22mm}=\lambda E^\mu E^{l+1}T^m\beta_\mu+\mbox{terms of order $n$} 
\label{10.497}
\end{eqnarray}
and similarly for the $N$th approximants, but with error term:
\begin{equation}
E_N^l T^m\left(\rho_N^{-1}(-\delta^\prime_N-\omega_{L\Lb,N})\right)={\bf O}(\tau^{N-1-m})
\label{10.498}
\end{equation}
Writing also $N_N^\mu\Lb_N E_N^l T^m\beta_{\mu,N}$ as $N^\mu\Lb E_N^l T^m\beta_{\mu,N}$ up to the 
0th order term $(N^\mu\Lb-N_N^\mu\Lb_N)E_N^l T^m\beta_{\mu,N}$ we express in this way 
$$N^\mu\Lb(E^l T^m\beta_\mu-E_N^l T^m\beta_{\mu,N})$$
up to lower order terms as:
\begin{equation}
\lambda E^\mu E(E^l T^m\beta_\mu-E_N^l T^m\beta_{\mu,N})=\lambda \s^{(E;m,l)}\check{\sxi}
\label{10.499}
\end{equation}
by \ref{10.489} for $V=E$. 
In view of the definition \ref{9.291} we have:
\begin{equation}
\|\sqrt{a}\s^{(E;m,l)}\check{\sxi}\|_{L^2(C_{\ub}^u)}\leq C\sqrt{\s^{(E;m,l)}\cEb^u(\ub)}
\label{10.500}
\end{equation}
Since $\lambda\sim u^2$, $\lambdab\sim\ub$, $a\sim u^2\ub$, the 2nd term in the integrant in 
\ref{10.488} has a contribution bounded by $C\tilde{G}(\ub_1,u)$ where:
\begin{equation}
\tilde{G}(\ub_1,u)=u^{-1}\int_0^{\ub_1}\tilde{g}(\ub,u)\sqrt{\ub}d\ub, \ \ \ 
\tilde{g}(\ub,u)=\|\sqrt{a}\s^{(E;m,l)}\check{\sxi}\|_{L^2(S_{\ub,u})}
\label{10.501}
\end{equation}
thus by \ref{10.500}:
\begin{equation}
\int_{\ub}^u\tilde{g}^2(\ub,u^\prime)du^\prime\leq C\s^{(E;m,l)}\cEb^u(\ub)
\label{10.502}
\end{equation}
We have (Schwartz inequality):
\begin{eqnarray}
&&\tilde{G}^2(\ub_1,u)\leq u^{-2}\int_0^{\ub_1}\ub d\ub \cdot \int_0^{\ub_1} \tilde{g}^2(\ub,u)d\ub\nonumber\\
&&\hspace{15mm}=\frac{1}{2}\frac{\ub_1^2}{u^2}\int_0^{\ub_1}\tilde{g}^2(\ub,u)d\ub 
\label{10.503}
\end{eqnarray}
Hence:
\begin{eqnarray}
&&\int_{\ub_1}^{u_1}\tilde{G}^2(\ub_1,u)du\leq\frac{1}{2}\ub_1^2
\int_{\ub_1}^{u_1}\left\{\int_0^{\ub_1}\tilde{g}^2(\ub,u)d\ub\right\}\frac{du}{u^2}\nonumber\\
&&\hspace{25mm}=\frac{1}{2}\ub_1^2\int_0^{\ub_1}
\left\{\int_{\ub_1}^{u_1}u^{-2}\tilde{g}^2(\ub,u)du\right\}d\ub
\label{10.504}
\end{eqnarray}
(reversing the order of integration). In regard to the interior integral we have, for all 
$\ub\in[0,\ub_1]$:
\begin{eqnarray}
&&\int_{\ub_1}^{u_1}u^{-2}\tilde{g}^2(\ub,u)du
=\int_{\ub_1}^{u_1}\frac{\partial}{\partial u}\left(\int_{\ub_1}^u \tilde{g}^2(\ub,u^\prime)
du^\prime\right)u^{-2}du\nonumber\\
&&\hspace{10mm}=u_1^{-2}\int_{\ub_1}^{u_1}\tilde{g}^2(\ub,u^\prime)du^\prime
+2\int_{\ub_1}^{u_1}\left(\int_{\ub_1}^u \tilde{g}^2(\ub,u^\prime)du^\prime\right)u^{-3}du 
\nonumber\\
&&\label{10.505}
\end{eqnarray}
By \ref{10.502} in view of the definition \ref{9.298} we have:
\begin{equation}
\int_{\ub_1}^u \tilde{g}^2(\ub,u^\prime)du^\prime\leq C\s^{(E;m,l)}\cBb(\ub_1,u_1)\ub^{2a_m}u^{2b_m}
\label{10.506}
\end{equation}
It follows that \ref{10.505} is bounded by: 
\begin{eqnarray}
&&C\s^{(E;m,l)}\cBb(\ub_1,u_1)\ub^{2a_m}\left\{u_1^{2b_m-2}
+2\frac{(u_1^{2b_m-2}-\ub_1^{2b_m-2})}{(2b_m-2)}\right\}\nonumber\\
&&\leq C\s^{(E;m,l)}\cBb(\ub_1,u_1)\ub^{2a_m}u_1^{2b_m-2}
\label{10.507}
\end{eqnarray}
the last by virtue of the assumption above that $b_m>2$. Substituting in \ref{10.504} we then obtain:
\begin{equation}
\int_{\ub-1}^{u_1}\tilde{G}^2(\ub_1,u)du\leq\frac{C\s^{(E;m,l)}\cBb(\ub_1,u_1)}{(2a_m+1)}\cdot\ub_1^{2a_m}u_1^{2b_m}\cdot\frac{\ub_1^3}{u_1^2}
\label{10.508}
\end{equation}
(compare with \ref{10.401}). This bounds the contribution to $\|\cnu_{m-1,l+1}\|^2_{L^2(\Cb_{\ub_1}^{u_1})}$. The last factor 
being homogeneous of degree 1 signifies that this contribution is 3 degrees better than borderline, 
hence will come with an extra factor of $\delta^3$. 

We turn finally to the dominant contribution, that of the 3rd term in the integrant in \ref{10.488}, 
the $\Nb$ component. By \ref{10.489} with $V=Y$ we have:
\begin{equation}
N^\mu\Lb(E^l T^m\beta_\mu-E_N^l T^m\beta_{\mu,N})+\ogamma\Nb^\mu\Lb(E^l T^m-E_N^l T^m\beta_{\mu,N})
=\s^{(Y;m,l)}\check{\xi}_{\Lb}
\label{10.509}
\end{equation}
Since $\ogamma\sim u$, this implies:
\begin{eqnarray}
&&u\|\Nb^\mu\Lb(E^l T^m\beta_\mu-E_N^l T^m\beta_{\mu,N})\|_{L^2(S_{\ub,u})}\leq \nonumber\\
&&\hspace{10mm}C\left\{\|\s^{(Y;m,l)}\check{\xi}_{\Lb}\|_{L^2(S_{\ub,u})}
+\|N^\mu\Lb(E^l T^m\beta_\mu-E_N^l T^m \beta_{\mu,N})\|_{L^2(S_{\ub,u})}\right\}\nonumber\\
\label{10.510}
\end{eqnarray}
By the preceding, the 2nd term on the right is bounded, up to lower order terms, by:
$$\frac{Cu}{\sqrt{\ub}}\tilde{g}(\ub,u)$$
with $\tilde{g}$ as in \ref{10.501}. We are then led to a conclusion analogous to \ref{10.379} - 
\ref{10.383}, namely that:
\begin{equation}
\|\Nb^\mu\Lb(E^l T^m\beta_\mu-E_N^l T^m\beta_{\mu,N})\|_{L^2(S_{\ub,u})}\leq g_0(\ub,u)+g_1(\ub,u)
\label{10.511}
\end{equation}
where $g_0(\ub,u)$ is equal to $Cu^{-1}\|\s^{(Y;m,l)}\check{\xi}_{\Lb}\|_{L^2(S_{\ub,u})}$ hence 
by \ref{9.291} satisfies
\begin{equation}
\int_{\ub}^{u_1}u^2 g_0^2(\ub,u)du\leq C\s^{(Y;m,l)}\cEb^{u_1}(\ub)
\label{10.512}
\end{equation}
while $g_1(\ub,u)$ is equal, neglecting lower order terms, to  $(C/\sqrt{\ub})\tilde{g}(\ub,u)$, 
hence by \ref{10.502} satisfies
\begin{equation}
\ub\int_{\ub}^{u_1} g_1^2(\ub,u)du\leq C\s^{(E;m,l)}\cEb^{u_1}(\ub)
\label{10.513}
\end{equation}
The contribution of 3rd term in the integrant in \ref{10.488} is then bounded by:
\begin{equation}
C(G_0(\ub_1,u)+G_1(\ub,u))
\label{10.514}
\end{equation}
where:
\begin{equation}
G_i(\ub_1,u)=u^{-2}\int_0^{\ub_1}g_i(\ub,u)\ub d\ub \ : \ i=0,1
\label{10.515}
\end{equation}
We thus have a complete analogy with \ref{10.379} - \ref{10.383}. Therefore we deduce, assuming 
that 
\begin{equation}
b_m\geq 4
\label{10.516}
\end{equation}
(see \ref{10.391}), in analogy with \ref{10.393}, the estimate:
\begin{equation}
\int_{\ub_1}^{u_1} G_0^2(\ub_1,u)du\leq\frac{C\s^{(Y;m,l)}\cBb(\ub_1,u_1)}{(2a_m+1)}\cdot 
\ub_1^{2a_m}u_1^{2b_m}\cdot\frac{\ub_1^4}{u_1^6}
\label{10.517}
\end{equation}
This is the {\em borderline} contribution, the last factor being homogeneous of degree -2. 
To control it we shall have to choose $a_m$ appropriately large so that the coefficient 
$C/(2a_m+1)$ is suitably small. We also deduce, in analogy with \ref{10.399}, the estimate:
\begin{equation}
\int_{\ub_1}^{u_1} G_1^2(\ub_1,u)du\leq\frac{C\s^{(E;m,l)}\cBb(\ub_1,u_1)}{(2a_m+1)}\cdot 
\ub_1^{2a_m}u_1^{2b_m}\cdot\frac{\ub_1^3}{u_1^4}
\label{10.518}
\end{equation}
The last factor being homogeneous of degree -1, which is one degree better than borderline, 
this contribution will come with an extra factor of $\delta$. 

We turn to the contribution of \ref{10.487} to the integral on the right in \ref{10.480}. 
Here we shall only consider the dominant contribution, that of the $\Nb$ component, which 
is bounded by a constant multiple of
\begin{equation}
u^{-2}\int_0^{\ub_1}\|a\Nb^\mu E(E^{l+1}T^{m-1}\beta_\mu-E_N^{l+1}T^{m-1}\beta_{\mu,N})
\|_{L^2(S_{\ub,u})}\ub d\ub
\label{10.519}
\end{equation}
Here, writing $a\Nb^\mu=\lambdab \Lb^\mu$, we use \ref{10.426} which in analogy with \ref{10.427} 
allows us to express:
\begin{eqnarray}
&&\Lb^\mu E^{l+2}T^{m-1}\beta_\mu=\Lb^\mu E^{l+1}T^{m-1}E\beta_\mu+\mbox{terms of order $l+m=n$}
\nonumber\\
&&\hspace{25mm}=E^{l+1}T^{m-1}(\Lb^\mu E\beta_\mu)+\mbox{terms of order $n$}\nonumber\\
&&\hspace{25mm}=E^{l+1}T^{m-1}(E^\mu\Lb\beta_\mu)+\mbox{terms of order $n$}\nonumber\\
&&\hspace{25mm}=E^\mu(E^{l+1}T^{m-1}\Lb\beta_\mu+\mbox{terms of order $n$}\nonumber\\
&&\hspace{25mm}=E^\mu\Lb E^{l+1}T^{m-1}\beta_\mu+\mbox{terms of order $n$} 
\label{10.520}
\end{eqnarray}
and similarly for the $N$th approximants, but with an error therm:
\begin{equation}
E_N^{l+1}T^{m-1}\omega_{\Lb E,N}={\bf O}(\tau^{N+2-m})
\label{10.521}
\end{equation}
Writing also $\Lb_N^\mu E_N^{l+2}T^{m-1}\beta_{\mu,N}$ as $\Lb^\mu E E_N^{l+1}T^{m-1}\beta_{\mu,N}$ 
up to the 0th order term $(\Lb^\mu E-\Lb_N^\mu E_N)E_N^{l+1}T^{m-1}\beta_{\mu,N}$, we express in this 
way 
$$\Lb^\mu E(E^{l+1}T^{m-1}\beta_\mu-E_N^{l+1}T^{m-1}\beta_{\mu,N})$$ 
up to lower order terms as:
\begin{equation}
E^\mu\Lb(E^{l+1}T^{m-1}\beta_\mu-E_N^{l+1}T^{m-1}\beta_{\mu,N})=\s^{(E;m-1,l+1)}\check{\xi}_{\Lb}
\label{10.522}
\end{equation}
by \ref{10.490} with $(m-1,l+1)$ in the role of $(m,l)$. Then, neglecting lower order terms, since 
$\lambdab\sim\ub$, \ref{10.519} is bounded by $CF(\ub_1,u)$, where
\begin{equation}
F(\ub_1,u)=u^{-2}\int_0^{\ub_1}\|\s^{(E;m-1,l+1)}\check{\xi}_{\Lb}\|_{L^2(S_{\ub,u})}\ub^2 d\ub
\label{10.523}
\end{equation}
We have (Schwartz inequality):
\begin{equation}
F^2(\ub_1,u)\leq \frac{1}{5}\frac{\ub_1^5}{u^4}\int_0^{\ub_1}\|\s^{(E;m-1,l+1)}\check{\xi}_{\Lb}\|^2_{L^2(S_{\ub,u})}d\ub
\label{10.524}
\end{equation}
and (reversing the order of integration):
\begin{equation}
\int_{\ub_1}^{u_1}F^2(\ub_1,u)du\leq\frac{1}{5}\ub_1^5\int_0^{\ub_1}\left\{\int_{\ub_1}^{u_1}
\|\s^{(E;m-1,l+1)}\check{\xi}_{\Lb}\|^2_{L^2(S_{\ub,u})}u^{-4}du\right\}d\ub
\label{10.525}
\end{equation}
Here, the interior integral is equal to:
\begin{eqnarray}
&&\int_{\ub_1}^{u_1}\frac{\partial}{\partial u}\left(\int_{\ub_1}^u
\|\s^{(E;m-1,l+1)}\check{\xi}_{\Lb}\|^2_{L^2(S_{\ub,u^\prime})}du^\prime\right)u^{-4}du \label{10.526}\\
&&\hspace{10mm}=u_1^{-4}\int_{\ub_1}^{u_1}
\|\s^{(E;m-1,l+1)}\check{\xi}_{\Lb}\|^2_{L^2(S_{\ub,u^\prime})}du^\prime\nonumber\\
&&\hspace{20mm}+4\int_{\ub_1}^{u_1}\left(\int_{\ub_1}^u
\|\s^{(E;m-1,l+1)}\check{\xi}_{\Lb}\|^2_{L^2(S_{\ub,u^\prime})}du^\prime\right)u^{-5}du
\nonumber
\end{eqnarray}
Now, in view of the definitions \ref{9.291} and \ref{9.298} we have:
\begin{eqnarray}
&&\int_{\ub}^u\|\s^{(E;m-1,l+1)}\check{\xi}_{\Lb}\|^2_{L^2(S_{\ub,u^\prime})}du^\prime
=\|\s^{(E;m-1,l+1)}\check{\xi}_{\Lb}\|^2_{L^2(\Cb_{\ub}^u)}\nonumber\\
&&\hspace{10mm}\leq C\s^{(E;m-1,l+1)}\cEb^{u}(\ub)\leq 
\s^{(E;m-1,l+1)}\cBb(\ub_1,u_1)\ub^{2a_{m-1}}u^{2b_{m-1}}\nonumber\\
&&\label{10.527}
\end{eqnarray}
It follows that \ref{10.526} is bounded by:
\begin{eqnarray}
&&C\s^{(E;m-1,l+1)}\cBb(\ub_1,u_1)\ub^{2a_{m-1}}\left\{u_1^{2b_{m-1}-4}
+4\frac{(u_1^{2b_{m-1}-4}-\ub_1^{2b_{m-1}-4})}{(2b_{m-1}-4)}\right\}\nonumber\\
&&\hspace{20mm}\leq C\s^{(E;m-1,l+1)}\cBb(\ub_1,u_1)\ub^{2a_{m-1}}u_1^{2b_{m-1}-4}
\label{10.528}
\end{eqnarray}
since by \ref{10.516} and the non-increasing property of the $b_m$ we have $b_{m-1}\geq 4$. 
Substituting in \ref{10.525} the estimate \ref{10.528} for the interior integral yields:
\begin{equation}
\int_{\ub_1}^{u_1}F^2(\ub_1,u)du\leq\frac{C\s^{(E;m-1,l+1)}\cBb(\ub_1,u_1)}{(2a_{m-1}+1)}\cdot 
\ub_1^{2a_{m-1}}u_1^{2b_{m-1}}\cdot\frac{\ub_1^6}{u_1^4}
\label{10.529}
\end{equation}
The last factor being homogeneous of degree 2, in view of the non-increasing property of the 
exponents $a_m$, $b_m$, this contribution is at least 4 degrees better than borderline 
hence will come with an extra factor of at least $\delta^4$. 

We have now completed the investigation 
of the contribution to $\|\cnu_{m-1,l+1}-\check{\tau}_{m-1,l+1}\|_{L^2(\Cb_{\ub_1}^u)}$ of the term \ref{10.483}, which arises 
from the 1st term in $j$ as given by Proposition 10.2.

We now briefly consider the contributions of the terms which arise from the 2nd and 3rd terms in $j$ 
as given by Proposition 10.2. These have extra factors of $\lambda$ and $\lambda^2$ respectively 
therefore it is obvious that these contributions will be far from borderline. The principal part of the term arising from the 2nd term in $j$ is:
\begin{equation}
\frac{1}{2}\pi\lambda\beta_N^2 H^\prime (g^{-1})^{\mu\nu}\beta_\nu\Lb(E^{l+1}T^{m-1}\beta_\mu 
-E_N^{l+1}T^{m-1}\beta_{\mu,N})
\label{10.530}
\end{equation}
We consider the dominant contribution of this to the integral on the right in \ref{10.480}, which 
is that of the $\Nb$ component. Since $\lambda\sim u^2$, this contribution is bounded by a 
constant multiple of:
\begin{equation}
\int_0^{\ub_1}\|\Nb^\mu\Lb(E^{l+1}T^{m-1}\beta_\mu-E_N^{l+1}T^{m-1}\beta_{\mu,N})\|_{L^2(S_{\ub,u})}
\ub d\ub
\label{10.531}
\end{equation}
Now, by \ref{10.489} with $V=Y$ and $(m-1,l+1)$ in the role of $(m,l)$, 
\begin{eqnarray}
&&N^\mu\Lb(E^{l+1}T^{m-1}\beta_\mu-E_N^{l+1}T^{m-1}\beta_{\mu,N})\nonumber\\
&&\hspace{15mm}+\ogamma\Nb^\mu\Lb(E^{l+1}T^{m-1}\beta_\mu-E_N^{l+1}T^{m-1}\beta_{\mu,N})\nonumber\\
&&\hspace{50mm}=\s^{(Y;m-1,l+1)}\check{\xi}_{\Lb}
\label{10.532}
\end{eqnarray}
The contribution of the right hand side dominates that of the 1st term on the left (compare 
\ref{10.393} with \ref{10.399} in regard to \ref{10.373}). Since $\ogamma\sim u$, this 
contribution to \ref{10.531} is bounded by a constant multiple of
\begin{equation}
F_1(\ub_1,u)=u^{-1}\int_0^{\ub_1}\|\s^{(Y;m-1,l+1)}\check{\xi}_{\Lb}\|_{L^2(S_{\ub,u})}\ub d\ub
\label{10.533}
\end{equation}
We have (Schwartz inequality):
\begin{equation}
F_1^2(\ub_1,u)\leq \frac{1}{3}\frac{\ub_1^3}{u^2}\int_0^{\ub_1}
\|\s^{(Y;m-1,l+1)}\check{\xi}_{\Lb}\|^2_{L^2(S_{\ub,u})}d\ub
\label{10.534}
\end{equation}
and (reversing the order of integration):
\begin{equation}
\int_{\ub_1}^{u_1}F_1^2(\ub_1,u)du\leq\frac{1}{3}\ub_1^3\int_0^{\ub_1}\left\{\int_{\ub_1}^{u_1}
\|\s^{(Y;m-1,l+1)}\check{\xi}_{\Lb}\|^2_{L^2(S_{\ub,u})}u^{-2}du\right\}d\ub
\label{10.535}
\end{equation}
The interior integral is equal to:
\begin{eqnarray}
&&\int_{\ub_1}^{u_1}\frac{\partial}{\partial u}\left(\int_{\ub_1}^u
\|\s^{(Y;m-1,l+1)}\check{\xi}_{\Lb}\|^2_{L^2(S_{\ub,u^\prime})}du^\prime\right)u^{-2}du \label{10.536}\\
&&\hspace{10mm}=u_1^{-2}\int_{\ub_1}^{u_1}
\|\s^{(Y;m-1,l+1)}\check{\xi}_{\Lb}\|^2_{L^2(S_{\ub,u^\prime})}du^\prime\nonumber\\
&&\hspace{20mm}+2\int_{\ub_1}^{u_1}\left(\int_{\ub_1}^u
\|\s^{(Y;m-1,l+1)}\check{\xi}_{\Lb}\|^2_{L^2(S_{\ub,u^\prime})}du^\prime\right)u^{-3}du
\nonumber
\end{eqnarray}
In view of the definitions \ref{9.291} and \ref{9.298} we have:
\begin{eqnarray}
&&\int_{\ub}^u\|\s^{(Y;m-1,l+1)}\check{\xi}_{\Lb}\|^2_{L^2(S_{\ub,u^\prime})}du^\prime
=\|\s^{(Y;m-1,l+1)}\check{\xi}_{\Lb}\|^2_{L^2(\Cb_{\ub}^u)}\nonumber\\
&&\hspace{10mm}\leq C\s^{(Y;m-1,l+1)}\cEb^{u}(\ub)\leq 
\s^{(Y;m-1,l+1)}\cBb(\ub_1,u_1)\ub^{2a_{m-1}}u^{2b_{m-1}}\nonumber\\
&&\label{10.537}
\end{eqnarray}
hence \ref{10.536} is bounded by:
\begin{eqnarray}
&&C\s^{(Y;m-1,l+1)}\cBb(\ub_1,u_1)\ub^{2a_{m-1}}\left\{u_1^{2b_{m-1}-2}
+2\frac{(u_1^{2b_{m-1}-2}-\ub_1^{2b_{m-1}-2})}{(2b_{m-1}-2)}\right\}\nonumber\\
&&\hspace{20mm}\leq C\s^{(Y;m-1,l+1)}\cBb(\ub_1,u_1)\ub^{2a_{m-1}}u_1^{2b_{m-1}-2}
\label{10.538}
\end{eqnarray}
Substituting in \ref{10.535} the estimate \ref{10.538} for the interior integral yields:
\begin{equation}
\int_{\ub_1}^{u_1}F_1^2(\ub_1,u)du\leq\frac{C\s^{(Y;m-1,l+1)}\cBb(\ub_1,u_1)}{(2a_{m-1}+1)}\cdot 
\ub_1^{2a_{m-1}}u_1^{2b_{m-1}}\cdot\frac{\ub_1^4}{u_1^2}
\label{10.539}
\end{equation}
The last factor being homogeneous of degree 2, in view of the non-increasing property of the 
exponents $a_m$, $b_m$, this contribution is at least 4 degrees better than borderline 
hence will come with an extra factor of at least $\delta^4$.

The principal part of the term arising from the 3rd term in $j$ is:
\begin{eqnarray}
&&\lambda^2\left\{-\beta_N\left[\pi\sbeta-\frac{1}{4c}(\beta_N+2\beta_{\Nb})\right](g^{-1})^{\mu\nu}
\beta_\nu+H\left[\pi\sbeta-\frac{1}{2c}(\beta_N+\beta_{\Nb})\right]N^\mu\right\}\nonumber\\
&&\hspace{30mm}\cdot E(E^{l+1}T^{m-1}\beta_\mu-E_N^{l+1}T^{m-1}\beta_{\mu,N})
\label{10.540}
\end{eqnarray}
We shall consider the dominant partial contribution from this to the integral on the right in \ref{10.480}, that of the $\Nb$ component. Since $\lambda\Nb^\mu=c\Lb^\mu$ this contribution 
is bounded by a constant multiple of:
\begin{equation}
\int_0^{\ub_1}\|\Lb^\mu E(E^{l+1}T^{m-1}\beta_\mu-E_N^{l+1}T^{m-1}\beta_{\mu,N})\|_{L^2(S_{\ub,u})}
\ub d\ub
\label{10.541}
\end{equation}
which (see \ref{10.522}) up to lower order terms is:
\begin{equation}
\int_0^{\ub_1}\|\s^{(E;m-1,l+1)}\check{\xi}_{\Lb}\|_{L^2(S_{\ub,u})}\ub d\ub:=F_2(\ub_1,u)
\label{10.542}
\end{equation}
We have (Schwartz inequality):
\begin{equation}
F_2^2(\ub_1,u)\leq \frac{1}{3}\ub_1^3\int_0^{\ub_1}
\|\s^{(E;m-1,l+1)}\check{\xi}_{\Lb}\|^2_{L^2(S_{\ub,u})}d\ub
\label{10.543}
\end{equation}
and (reversing the order of integration):
\begin{eqnarray}
&&\int_{\ub_1}^{u_1}F_2^2(\ub_1,u)du\leq \frac{1}{3}\ub_1^3\int_0^{\ub_1}\left\{\int_{\ub_1}^{u_1}
\|\s^{(E;m-1,l+1)}\check{\xi}_{\Lb}\|^2_{L^2(S_{\ub,u})}du\right\}d\ub\nonumber\\
&&\hspace{25mm}\leq\frac{1}{3}\ub_1^3\int_0^{\ub_1}
\|\s^{(E;m-1,l+1)}\check{\xi}_{\Lb}\|^2_{L^2(\Cb_{\ub}^{u_1})}d\ub\nonumber\\
&&\hspace{25mm}\leq C\ub_1^3\int_0^{\ub_1}\s^{(E;m-1,l+1)}\cEb^{u_1}(\ub)d\ub\nonumber\\
&&\hspace{20mm}\leq C\ub_1^3\int_0^{\ub_1}\s^{(E;m-1,l+1)}\cBb(\ub_1,u_1)\ub^{2a_{m-1}}u_1^{2b_{m-1}}d\ub\nonumber\\
&&\hspace{20mm}=\frac{C\s^{(E;m-1,l+1)}\cBb(\ub_1,u_1)}{(2a_{m-1}+1)}\cdot\ub_1^{2a_{m-1}}u_1^{2b_{m-1}}\cdot\ub_1^4
\label{10.544}
\end{eqnarray}
by \ref{9.291} and \ref{9.298}. The last factor being of degree 4, 
in view of the non-increasing property of the exponents $a_m$, $b_m$, 
this contribution is at least 6 degrees better than borderline 
hence will come with an extra factor of at least $\delta^6$. 

We now consider the contribution of the principal part of the difference term 
\begin{equation}
\tilde{K}^{\prime\prime}_{m-1,l+1}-\tilde{K}^{\prime\prime}_{m-1,l+1,N}
\label{10.545}
\end{equation}
(see \ref{10.481}) to the integral on the right in \ref{10.480}. Augmenting the assumptions 
\ref{10.439}, \ref{10.440} with the assumption that in ${\cal R}_{\delta,\delta}$ it holds:
\begin{equation}
|\Lb\lambda|\leq Cu, \ \ \ |E\lambdab|\leq C\lambdab, \ \ \ |\tchib|\leq C
\label{10.546}
\end{equation}
we deduce from \ref{10.322} that the principal part of \ref{10.545} can be written in the form:
\begin{eqnarray}
&&[\tilde{K}^{\prime\prime}_{m-1,l+1}-\tilde{K}^{\prime\prime}_{m-1,l+1,N}]_{P.P.}=\nonumber\\
&&\hspace{17mm}\ub u k_1^\mu E(E^{l+1}T^{m-1}\beta_\mu-E_N^{l+1}T^{m-1}\beta_{\mu,N})\nonumber\\
&&\hspace{17mm}+u^2 k_2^\mu L(E^{l+1}T^{m-1}\beta_\mu-E_N^{l+1}T^{m-1}\beta_{\mu,N})\nonumber\\
&&\hspace{17mm}+\ub k_3^\mu \Lb(E^{l+1} T^{m-1}\beta_\mu-E_N^{l+1}T^{m-1}\beta_{\mu,N})\nonumber\\
&&\hspace{17mm}+\ub k_4^\mu(E^l T^{m-1}\Lb^2\beta_\mu-E_N^l T^{m-1}\Lb_N^2\beta_{\mu,N})\nonumber\\
&&\hspace{15mm}+u\left[\ub u^2 k^{\prime\mu}_1 E(E^{l+2}T^{m-2}\beta_\mu 
-E_N^{l+2} T^{m-2}\beta_{\mu,N})\right.\nonumber\\
&&\hspace{19mm}+u^2 k^{\prime\mu}_2 L(E^{l+2}T^{m-2}\beta_\mu
-E_N^{l+2} T^{m-2}\beta_{\mu,N})\nonumber\\
&&\hspace{19mm}\left.+\ub k^{\prime\mu}_3\Lb(E^{l+2}T^{m-2}\beta_\mu-E_N^{l+2}T^{m-2}\beta_{\mu,N})
\right]\nonumber\\
&&+\sum_{i=0}^{m-2}(hu^2)^{m-2-i}\ub\left\{u^3 k^{\prime\prime\mu}_1 E(E^{l+m-i}T^i\beta_\mu 
-E_N^{l+m-i} T^i\beta_{\mu,N})\right.\nonumber\\
&&\hspace{19mm}+u^2 k^{\prime\prime\mu}_2 L(E^{l+m-i} T^i\beta_\mu
-E_N^{l+m-i} T^i\beta_{\mu,N})\nonumber\\
&&\hspace{19mm}+u^2 k^{\prime\prime\mu}_3\Lb(E^{l+m-i}T^i\beta_\mu
-E_N^{l+m-i}T^i\beta_{\mu,N})\nonumber\\
&&\hspace{19mm}\left.+k^{\prime\prime\mu}_4(E^{l-1+m-i}T^i\Lb^2\beta_\mu
-E_N^{l-1+m-i}T^i\Lb_N^2\beta_{\mu,N})\right\} \nonumber\\
&&\label{10.547}
\end{eqnarray}
Here the coefficients $k$, $k^\prime$, $k^{\prime\prime}$ and the function $h$ in the sum are 
bounded by a constant. We have used the fact that by virtue of the assumptions \ref{10.439},  \ref{10.440} the coefficient $A$ defined by \ref{10.99} satisfies:
\begin{equation}
|A|\leq C\ub \ \ \mbox{: in ${\cal R}_{\delta,\delta}$}
\label{10.548}
\end{equation}

Consider first the first four terms on the right in \ref{10.547}. The 4th term is similar to 
the term \ref{10.484}. But whereas \ref{10.484} is multiplied in \ref{10.481} by $2L\lambda/\lambda$ 
and $-L\lambda\sim \ub$, in \ref{10.547} the coefficient of the 4th term is bounded by $C\ub$. 
Thus, since $\lambda\sim u^2$, while the dominant contribution of \ref{10.484} to 
$\|\cnu_{m-1,l+1}-\check{\tau}_{m-1,l+1}\|^2_{L^2(\Cb_{\ub_1}^{u_1})}$ is the borderline contribution \ref{10.517}, 
the dominant contribution of the 4th term in \ref{10.547} is bounded by:
\begin{equation}
\frac{C\s^{(Y;m,l)}\cBb(\ub_1,u_1)}{(2a_m+1)}\cdot 
\ub_1^{2a_m}u_1^{2b_m}\cdot\frac{\ub_1^4}{u_1^2}
\label{10.549}
\end{equation}
The last factor being homogeneous of degree 2 which is 4 degrees better than borderline, this 
will come with an extra factor of $\delta^4$. As for the 3rd term in \ref{10.547} this is similar 
to the term \ref{10.530}. Moreover since \ref{10.530} is multiplied by $2L\lambda/\lambda$ and 
we have an extra $\lambda$ factor and in \ref{10.547} the coefficient of the 3rd term is bounded 
by $C\ub$, the contributions are comparable. Thus (see \ref{10.539}) the dominant contribution of the 
3rd term is bounded by:
\begin{equation}
\frac{C\s^{(Y;m-1,l+1)}\cBb(\ub_1,u_1)}{(2a_{m-1}+1)}\cdot 
\ub_1^{2a_{m-1}}u_1^{2b_{m-1}}\cdot\frac{\ub_1^4}{u_1^2}
\label{10.550}
\end{equation}
which will come with an extra factor of at least $\delta^4$. Consider next the 2nd term in 
\ref{10.547}. The dominant partial contribution of this term is that of the $\Nb$ component which, 
up to a bounded coefficient, is: 
$$u^2\Nb^\mu L(E^{l+1}T^{m-1}\beta_\mu-E_N^{l+1}T^{m-1}\beta_{\mu,N})$$
so the corresponding contribution to the integral on the right in \ref{10.480} is bounded by a constant 
multiple of:
\begin{equation}
u^2\int_0^{\ub_1}\|\Nb^\mu L(E^{l+1}T^{m-1}\beta_\mu-E_N^{l+1}T^{m-1}\beta_{\mu,N})\|_{L^2(S_{\ub,u})}d\ub
\label{10.551}
\end{equation}
To estimate this we use \ref{10.447}, which allows us to express:
\begin{eqnarray}
&&\Nb^\mu LE^{l+1}T^{m-1}\beta_\mu=\Nb^\mu E^{l+1}T^{m-1} L\beta_\mu+\mbox{terms of order $l+m=n$} \nonumber\\
&&\hspace{28mm}=E^{l+1}T^{m-1}(\Nb^\mu L\beta_\mu)+\mbox{terms of order $n$} \nonumber\\
&&\hspace{28mm}=E^{l+1}T^{m-1}(\lambdab E^\mu E\beta_\mu)+\mbox{terms of order $n$} \nonumber\\
&&\hspace{28mm}=\lambdab E^\mu E^{l+2}T^{m-1}\beta_\mu+\mbox{terms of order $n$} 
\label{10.552}
\end{eqnarray}
Similarly for the $N$th approximants, but in this case there is an error term:
\begin{equation}
E_N^{l+1}T^{m-1}\left(\rhob_N^{-1}(-\delta^\prime_N+\omega_{L\Lb,N}0\right)={\bf O}(u^{-2}\tau^{N+1-m})
\label{10.553}
\end{equation}
We express in this way 
$$\Nb^\mu L(E^{l+1}T^{m-1}\beta_\mu-E_N^{l+1}T^{m-1}\beta_{\mu,N})$$
up to lower order terms as:
\begin{equation}
\lambdab E^\mu E(E^{l+1}T^{m-1}\beta_\mu-E_N^{l+1}T^{m-1}\beta_{\mu,N})=\lambdab \s^{(E;m-1,l+1)}\check{\sxi}
\label{10.554}
\end{equation}
Since $\lambdab\sim\ub$ we can then estimate \ref{10.551} in terms of 
\begin{equation}
u^2\int_0^{\ub_1}\|\s^{(E;m-1,l+1)}\check{\sxi}\|_{L^2(S_{\ub,u})}\ub d\ub
\label{10.555}
\end{equation}
This is analogous to \ref{10.451} but with an extra $u^2$ factor. The contribution 
to $\|\cnu_{m-1,l+1}-\check{\tau}_{m-1,l+1}\|^2_{L^2(\Cb_{\ub_1}^{u_1})}$ is then bounded analogously to \ref{10.417} 
but with an extra $u_1^4$ factor, that is by:
\begin{equation}
\frac{C\s^{(E;m-1,l+1)}\cBb(\ub_1,u_1)}{(2a_{m-1}+1)}
\cdot\ub_1^{2a_{m-1}} u_1^{2b_{m-1}}\cdot\ub_1^3 u_1^2
\label{10.556}
\end{equation}
which is at least 7 degrees better than borderline. However in the case of the 2nd term in 
\ref{10.547}, which is what we are presently considering, this being an $L$ derivative, we must also 
consider the $N$ component, that is:
\begin{equation}
u^2 N^\mu L(E^{l+1}T^{m-1}\beta_\mu-E_N^{l+1} T^{m-1}\beta_{\mu,N})
\label{10.a6}
\end{equation}
which can only be estimated by writing, in analogy with \ref{10.453}, 
\begin{eqnarray}
&&N^\mu L(E^{l+1}T^{m-1}\beta_\mu-E_N^{l+1} T^{m-1}\beta_{\mu,N})=\s^{(Y;m-1,l+1)}\check{\sxi} \label{10.a7}\\
&&\hspace{35mm}-\ogamma\Nb^\mu L(E^{l+1}T^{m-1}\beta_\mu-E_N^{l+1} T^{m-1}\beta_{\mu,N}) \nonumber
\end{eqnarray}
Here the 2nd term is the $\Nb$ component whose contribution we have already estimated but with 
an extra $\ogamma\sim u$ factor. As for the 1st term, its contribution to 
$\|\cnu_{m-1,l+1}-\check{\tau}_{m-1,l+1}\|_{L^2(S_{\ub,u})}$ is bounded by:
\begin{eqnarray}
&&u^2\int_0^{\ub_1}\|\s^{(Y;m-1,l+1)}\check{\xi}_L\|_{L^2(S_{\ub,u})}d\ub
\leq u^2\sqrt{\ub_1}\|\s^{(Y;m-1,l+1}\check{\xi}_L\|_{L^2(C_u^{\ub_1})}\nonumber\\
&&\leq Cu^2\sqrt{\ub_1}\sqrt{\s^{(Y;m-1,l+1)}\cE^{\ub_1}(u)}
\leq C\sqrt{\s^{(Y;m-1,l+1)}\cB(\ub_1,u_1)}\ub_1^{a_{m-1}+\frac{1}{2}}u^{b_{m-1}+2}\nonumber\\
&&\label{10.a8}
\end{eqnarray}
The corresponding contribution to $\|\cnu_{m-1,l+1}-\check{\tau}_{m-1,l+1}\|^2_{L^2(\Cb_{\ub_1}^{u_1})}$ is then bounded by:
\begin{equation}
\frac{C\s^{(Y;m-1,l+1)}\cB}{(2b_{m-1}+5)}\cdot \ub_1^{2a_{m-1}}u_1^{2b_{m-1}}\cdot\ub_1 u_1^5
\label{10.a9}
\end{equation}
Comparing with \ref{10.556} we see that while here the last factor is of degree 6 rather than 
5, it dominates the last factor in \ref{10.556} in a suitable neighborhood of $\Cb_0$.

The first term on the right in \ref{10.547} 
is similar to the term \ref{10.540}. But whereas \ref{10.540} is multiplied in \ref{10.481} by 
$2L\lambda/\lambda$ so the overall coefficient is $2\lambda L\lambda$ which is bounded by $C\ub u^2$, 
in \ref{10.547} the coefficient of the 1st term is bounded by $C\ub u$. It follows that the 
contribution of this term to the integral on the right in \ref{10.480} is bounded by:
\begin{equation}
Cu^{-1}\tilde{F}_2(\ub_1,u)
\label{10.557}
\end{equation}
and the corresponding contribution to $\|\cnu_{m-1,l+1}-\check{\tau}_{m-1,l+1}\|^2_{L^2(\Cb_{\ub_1}^{u_1})}$ is bounded 
by (see \ref{10.544}):
\begin{equation}
\frac{C\s^{(E;m-1,l+1)}\cBb(\ub_1,u_1)}{(2a_{m-1}+1)}\cdot\ub_1^{2a_{m-1}}u_1^{2b_{m-1}}
\cdot\frac{\ub_1^4}{u_1^2}
\label{10.558}
\end{equation}
The last factor being homogeneous of degree 2, which is 4 degrees better than borderline, this 
contribution will come with an extra factor of $\delta^4$. 

Consider next the 5th term on the right in \ref{10.547}. Here we have three terms in the square 
bracket which are analogous to the first three terms on the right in \ref{10.547} but with 
$(m-2,l+2)$ in the role of $(m-1,l+1)$. Taking account of the overall factor $u$, the 
coefficients of the three terms in the square bracket have extra factors of $u^2$, $u$, $u$, 
respectively, relative to the coefficients of the first three terms in \ref{10.547}. 
We then deduce (see \ref{10.558}, \ref{10.556} and \ref{10.a9}, \ref{10.550}) that the dominant contributions 
of these terms to $\|\cnu_{m-1,l+1}-\check{\tau}_{m-1,l+1}\|^2_{L^2(\Cb_{\ub_1}^{u_1})}$ are bounded by:
\begin{eqnarray}
&&\frac{C\s^{(E;m-2,l+2)}\cBb(\ub_1,u_1)}{(2a_{m-2}+1)}\cdot\ub_1^{2a_{m-2}}u_1^{2b_{m-2}}
\cdot\ub_1^4 u_1^2, \label{10.559}\\
&&\ub_1^{2a_{m-2}} u_1^{2b_{m-2}}\left(\frac{C\s^{(E;m-2,l+2)}\cBb(\ub_1,u_1)}{(2a_{m-2}+1)}
\cdot\cdot\ub_1^3 u_1^4\right.\nonumber\\
&&\hspace{28mm}\left.+\frac{C\s^{(Y;m-2,l+2)}\cB(\ub_1,u_1)}{(2b_{m-2}+5)}\cdot\ub_1 u_1^7\right), \label{10.560}\\
&&\frac{C\s^{(Y;m-2,l+2)}\cBb(\ub_1,u_1)}{(2a_{m-2}+1)}\cdot 
\ub_1^{2a_{m-2}}u_1^{2b_{m-2}}\cdot\ub_1^4 \label{10.561}
\end{eqnarray}
respectively. 

We finally turn to the sum in \ref{10.547}. Here we have four terms in the curly bracket which are 
analogous to the first four terms on the right in \ref{10.547} but with $(i,l+m-i)$ in the role of 
$(m-1,l+1)$. Taking account of the overall factor $\ub$, the coefficients of the four terms in 
the curly bracket have extra factors of $u^2$, $\ub$, $u^2$, $1$, respectively, relative to the 
coefficients of the first four terms in \ref{10.547}. Taking also into account the overall factor 
$(hu^2)^{m-2-i}$ we then deduce (see \ref{10.558}, \ref{10.556} and \ref{10.a9}, \ref{10.550}, \ref{10.549}) that 
the dominant partial contributions of these terms, through the $i$th term in the sum, to 
$\|\cnu_{m-1,l+1}-\check{\tau}_{m-1,l+1}\|^2_{L^2(\Cb_{\ub_1}^{u_1})}$ are bounded by:
\begin{eqnarray}
&&\frac{C\s^{(E;i,l+m-i)}\cBb(\ub_1,u_1)}{(2a_i+1)}\cdot\ub_1^{2a_i}u_1^{2b_i}
\cdot \ub_1^4 u_1^2 (M u_1^2)^{2(m-2-i)}, \label{10.562}\\
&&\ub_1^{2a_i} u_1^{2b_i}(M u_1^2)^{2(m-2-i)}\left(\frac{C\s^{(E;i,l+m-i)}\cBb(\ub_1,u_1)}{(2a_i+1)}
\cdot\ub_1^5 u_1^2 \right.\nonumber\\
&&\hspace{39mm}\left.+\frac{C\s^{(Y;i,l+m-i)}\cB(\ub_1,u_1)}{(2b_i+5)}\cdot\ub_1^3 u_1^5\right), \label{10.563}\\
&&\frac{C\s^{(Y;i,l+m-i)}\cBb(\ub_1,u_1)}{(2a_i+1)}\cdot 
\ub_1^{2a_i}u_1^{2b_i}\cdot\ub_1^4 u_1^2 (M u_1^2)^{2(m-2-i)}, \label{10.564}\\
&&\frac{C\s^{(Y;i+1,l+m-i-1)}\cBb(\ub_1,u_1)}{(2a_{i+1}+1)}\cdot 
\ub_1^{2a_{i+1}}u_1^{2b_{i+1}}\cdot\frac{\ub_1^4}{u_1^2} (M u_1^2)^{2(m-2-i)} \nonumber\\
&&\label{10.565}
\end{eqnarray}
where $M$ is the bound for $|h|$ in ${\cal R}_{\delta,\delta}$. 
Then, since taking $\delta$ small enough so that $M\delta^2\leq 1/2$ we have
$$\sum_{i=0}^{m-2}(M u_1^2)^{m-2-i}\leq \sum_{j=0}^\infty(M\delta^2)^j\leq 2,$$
the contribution of the sum in \ref{10.547} to $\|\cnu_{m-1,l+1}-\check{\tau}_{m-1,l+1}\|_{L^2(\Cb_{\ub_1}^{u_1})}$ is 
bounded by, recalling the definitions \ref{9.300}, \ref{9.301},
\begin{equation}
C\sqrt{\max\{\s^{[m-1,l+1]}\cB(\ub_1,u_1),\s^{[m-1,l+1]}\cBb(\ub_1,u_1)\}}\cdot\ub_1^{a_{m-1}}u_1^{b_{m-1}}\cdot \frac{\ub_1^{3/2} }{u_1^{1/2}}
\label{10.566}
\end{equation}

We finally consider the contribution of the terms \ref{10.482} to the integral on the right in 
\ref{10.480}. From the above results all other contributions to 
$\|\cnu_{m-1,l+1}-\check{\tau}_{m-1,l+1}\|_{L^2(\Cb_{\ub_1}^{u_1})}$ are bounded by:
\begin{eqnarray}
&&\ub_1^{a_m}u_1^{b_m}\left\{C\sqrt{\frac{\s^{(Y;m,l)}\cBb(\ub_1,u_1)}{(2a_m+1)}}
\cdot\frac{\ub_1^2}{u_1^3}\right. \label{10.567}\\
&&+C_0\sqrt{\max\{\s^{(m,l)}\cB(\ub_1,u_1),\s^{(m,l)}\cBb(\ub_1,u_1)\}}
\cdot\frac{\ub_1^{3/2}}{u_1^2}\nonumber\\
&&\left.+C_1\sqrt{\max\{\s^{[m-1,l+1]}\cB(\ub_1,u_1),\s^{[m-1,l+1]}\cBb(\ub_1,u_1)\}}
\cdot \ub_1^{1/2}u_1^{1/2}\right\}\nonumber
\end{eqnarray}
Here we have taken into account the non-increasing property of the exponents $a_m$, $b_m$ with $m$. 

Now, \ref{10.482} is equal up to lower order terms to $y_{m-1,l+1}$, where:
\begin{eqnarray}
&&y_{m-1,l+1}=-(m-1)\zeta(\cnu_{m-2,l+2}-\check{\tau}_{m-2,l+2})\label{10.568}\\
&&\hspace{10mm}+A\left\{(2\pi\lambda)^{m-1}\cth_{l+m}+2\sum_{i=0}^{m-2}(2\pi\lambda)^{m-2-i}
(\cnu_{i,l+m-i}-\check{\tau}_{i,l+m-i})\right\}\nonumber
\end{eqnarray}
Since $y_{m-1,l+1}$ involves the $\cnu_{i,l+m-i}-\check{\tau}_{i,l+m-i}$ for $i=0,...,m-2$, we must proceed by induction. In the case $m=1$ \ref{10.568} reduces to:
\begin{equation}
y_{0,l+1}=A\cth_{l+1}
\label{10.569}
\end{equation}
In view of \ref{10.548} the contribution of this to the integral on the right in \ref{10.480} is 
bounded by a constant multiple of:
\begin{equation}
\int_0^{\ub_1}\ub\|\cth_{l+1}\|_{L^2(S_{\ub,u})}d\ub
\leq\left(\frac{1}{3}\ub_1^3\int_0^{\ub_1}\|\cth_{l+1}\|^2_{L^2(S_{\ub,u})}d\ub\right)^{1/2}
\label{10.570}
\end{equation}
It follows that the corresponding contribution to $\|\cnu_{0,l+1}-\check{\tau}_{0,l+1}\|_{L^2(\Cb_{\ub_1}^{u_1})}$ is 
bounded by a constant multiple of:
\begin{equation}
\ub_1^{3/2}\left(\int_0^{\ub_1}\|\cth_{l+1}\|^2_{L^2(\Cb_{\ub}^{u_1})}d\ub\right)^{1/2}
\label{10.571}
\end{equation}
Substituting from Proposition 10.3, with $l+1$ in the role of $l$ (there $l=n$ while here $l+1=n$) and 
$\ub$ in the role of $\ub_1$, the estimate for $\|\cth_{l+1}\|_{L^2(\Cb_{\ub}^{u_1})}$, we conclude,  
in view of the monotonicity properties of the functions $\s^{(V;m,l)}\cB(\ub,u)$, 
$\s^{(V;m,l)}\cBb(\ub,u)$, that this is bounded by:
\begin{eqnarray}
&&\ub_1^{a_0}u_1^{b_0}\left\{C\sqrt{\frac{\s^{(Y;0,l+1)}\cBb(\ub_1,u_1)}{(2a_0+5)(2a_0+1)}}
\frac{\ub_1^4}{u_1^3}\right. \label{10.572}\\
&&\hspace{10mm}\left.+C\sqrt{\frac{\max\{\s^{(0,l+1)}\cB(\ub_1,u_1),\s^{(0,l+1)}\cBb(\ub_1,u_1)\}}{(2a_0+2)}}
\frac{\ub_1^{5/2}}{u_1}\right\}\nonumber\\
&&\leq\frac{C^\prime}{\sqrt{2a_0+2}}\sqrt{\max\{\s^{(0,l+1)}\cB(\ub_1,u_1),\s^{(0,l+1)}\cBb(\ub_1,u_1)\}}
\cdot\ub_1^{a_0}u_1^{b_0}\cdot\frac{\ub_1^{5/2}}{u_1^{3/2}} \nonumber
\end{eqnarray}
Adding then this bound for the contribution of \ref{10.569} to \ref{10.567} with $m=1$ we conclude that:
\begin{eqnarray}
&&\|\cnu_{0,l+1}-\check{\tau}_{0,l+1}\|_{L^2(\Cb_{\ub_1}^{u_1})}\leq \ub_1^{a_1}u_1^{b_1}\left\{
C\sqrt{\frac{\s^{(Y;1,l)}\cBb(\ub_1,u_1)}{(2a_1+1)}}\cdot\frac{\ub_1^2}{u_1^3}\right.\nonumber\\
&&\hspace{30mm}\left.+C\sqrt{\max\{\s^{[1,l]}\cB(\ub_1,u_1),\s^{[1,l]}\cBb(\ub_1,u_1)\}}\cdot\frac{\ub_1^{1/2}}{u_1}
\right\}\nonumber\\
&&\label{10.573}
\end{eqnarray}
(recall the definitions \ref{9.301}). Suppose then, as the inductive hypothesis, that for some 
$m\in\{2,...,n\}$ and all $j=1,...,m-1$ the following estimate holds for all 
$(\ub_1,u_1)\in R_{\delta,\delta}$:
\begin{eqnarray}
&&\|\cnu_{j-1,l+m-j+1}-\check{\tau}_{j-1,l+m-j+1}\|_{L^2(\Cb_{\ub_1}^{u_1})}\leq \ub_1^{a_j}u_1^{b_j}\left\{
C\sqrt{\frac{\s^{(Y;j,l+m-j)}\cBb(\ub_1,u_1)}{(2a_j+1)}}\cdot\frac{\ub_1^2}{u_1^3}\right.\nonumber\\
&&\hspace{15mm}\left.+C_{j,l+m-j}\sqrt{\max\{\s^{[j,l+m-j]}\cB(\ub_1,u_1),\s^{[j,l+m-j]}\cBb(\ub_1,u_1)\}}\cdot\frac{\ub_1^{1/2}}{u_1}
\right\}\nonumber\\
&&\label{10.574}
\end{eqnarray}
for certain constants $C_{j,l+m-j}$ to be determined recursively. 
Now, assumptions \ref{10.439}, \ref{10.440} and \ref{10.546} imply that:
\begin{equation}
|\zeta|\leq C\ub \ \ \mbox{: in ${\cal R}_{\delta,\delta}$}
\label{10.575}
\end{equation}
Then the contribution of the 1st term in \ref{10.568} to the integral on the right in \ref{10.480} 
is bounded by a constant multiple of:
\begin{eqnarray}
&&\int_0^{\ub_1}\ub\|\cnu_{m-2,l+2}-\check{\tau}_{m-2,l+2}\|_{L^2(S_{\ub,u})}d\ub\leq \label{10.576}\\
&&\hspace{35mm}\left(\frac{1}{3}\ub_1^3
\int_0^{\ub_1}\|\cnu_{m-2,l+2}-\check{\tau}_{m-2,l+2}\|^2_{L^2(S_{\ub,u})}d\ub\right)^{1/2}
\nonumber
\end{eqnarray}
It follows that the corresponding contribution to $\|\cnu_{m-1,l+1}-\check{\tau}_{m-1,l+1}\|_{L^2(\Cb_{\ub_1}^{u_1})}$ is 
bounded by a constant multiple of:
\begin{equation}
\ub_1^{3/2}\left(\int_0^{\ub_1}\|\cnu_{m-2,l+2}-\check{\tau}_{m-2,l+2}\|^2_{L^2(\Cb_{\ub}^{u_1})}d\ub\right)^{1/2}
\label{10.577}
\end{equation}
Substituting the estimate \ref{10.574} for $j=m-1$ with 
$\ub$ in the role of $\ub_1$, we conclude,  
in view of the monotonicity properties of the functions $\s^{(V;m,l)}\cB(\ub,u)$, 
$\s^{(V;m,l)}\cBb(\ub,u)$, that this is bounded by:
\begin{eqnarray}
&&\ub_1^{a_{m-1}}u_1^{b_{m-1}}\left\{C\sqrt{\frac{\s^{(Y;m-1,l+1)}\cBb(\ub_1,u_1)}
{(2a_{m-1}+5)(2a_{m-1}+1)}}
\cdot\frac{\ub_1^4}{u_1^3}\right.\nonumber\\
&&\hspace{10mm}\left.+C_{m-1,l+1}\sqrt{\frac{\max\{\s^{[m-1,l+1]}\cB(\ub_1,u_1),
\s^{[m-1,l+1]}\cBb(\ub_1,u_1)\}}{(2a_{m-1}+2)}}\cdot\frac{\ub_1^{5/2}}{u_1}\right\}\nonumber\\
&&\leq \frac{C_{m-1,l+1}}{\sqrt{2a_{m-1}+2}}\sqrt{\max\{\s^{[m-1,l+1]}\cB(\ub_1,u_1),
\s^{[m-1,l+1]}\cBb(\ub_1,u_1)\}}\cdot\frac{\ub_1^{5/2}}{u_1^{3/2}}\nonumber\\
&&\label{10.578}
\end{eqnarray}
as we can assume that:
\begin{equation}
C_{m-1,l+1}\geq\sqrt{\frac{(2a_{m-1}+2)}{(2a_{m-1}+5)(2a_{m-1}+1)}}
\label{10.c2}
\end{equation}
In view of \ref{10.578}, the contribution of the sum in \ref{10.568} to the integral on the right in 
\ref{10.480} is bounded by a constant multiple of:
\begin{equation}
\sum_{i=0}^{m-2}\left(\Pi u^2\right)^{m-2-i}\int_0^{\ub_1}\ub\|\cnu_{i,l+m-i}-\check{\tau}_{i,l+m-i}\|_{L^2(S_{\ub,u})}d\ub
\label{10.579}
\end{equation}
where $\Pi$ is a fixed positive bound for $|2\pi\lambda/u^2|$ in ${\cal R}_{\delta,\delta}$. It follows that the 
corresponding contribution to $\|\cnu_{m-1,l+1}-\check{\tau}_{m-1,l+1}\|_{L^2(\Cb_{\ub_1}^{u_1})}$ is 
bounded by a constant multiple of:
\begin{equation}
\ub_1^{3/2}\sum_{i=0}^{m-2}\left(\Pi u_1^2\right)^{m-2-i}
\left(\int_0^{\ub_1}\|\cnu_{i,l+m-i}-\check{\tau}_{i,l+m-i}\|^2_{L^2(\Cb_{\ub}^{u_1})}d\ub\right)^{1/2}
\label{10.580}
\end{equation}
Substituting the estimate \ref{10.574} in the simplified form:
\begin{eqnarray}
&&\|\cnu_{j-1,l+m-j+1}-\check{\tau}_{j-1,l+m-j+1}\|_{L^2(\Cb_{\ub_1}^{u_1})}\leq \nonumber\\
&&\left(\frac{C}{\sqrt{2a_j+1}}+C_{j,l+m-j}\right)
\sqrt{\max\{\s^{[j,l+m-j]}\cB(\ub_1,u_1),\s^{[j,l+m-j]}\cBb(\ub_1,u_1)\}}\nonumber\\
&&\hspace{38mm}\cdot \ub_1^{a_j}u_1^{b_j}
\cdot\frac{\ub_1^{1/2}}{u_1^{3/2}}
\label{10.581}
\end{eqnarray}
%where 
%\begin{equation}
%C^\prime_{j,l+m-j}=\max\left\{\frac{C}{\sqrt{2a_j+1}},C_{j,l+m-j}\right\}
%\label{10.582}
%\end{equation}
for $j=i+1$, $i=0,...,m-2$ with $\ub$ in the role of $\ub_1$, 
taking into account the non-decreasing property of the exponents $a_m$, $b_m$ with $m$, 
and taking $\Pi\delta^2\leq 1/2$ so that 
$$\sum_{i=0}^{m-2}(\Pi u_1^2)^{m-2-i}\leq \sum_{j=0}^\infty (\Pi\delta^2)^j\leq 2,$$
we deduce that \ref{10.580} is bounded by:
\begin{eqnarray}
&&\frac{2}{\sqrt{2a_{m-1}+2}}\left(\frac{C}{\sqrt{2a_{m-1}+1}}+C_{[m-1,l+1]}\right)\cdot\nonumber\\
&&\hspace{8mm}\sqrt{\max\{\s^{[m-1,l+1]}\cB(\ub_1,u_1),\s^{[m-1,l+1]}\cBb(\ub_1,u_1)\}}
\cdot \ub_1^{a_{m-1}}u_1^{b_{m-1}}\cdot\frac{\ub_1^{5/2}}{u_1^{3/2}}\nonumber\\
&&\label{10.583}
\end{eqnarray}
where
\begin{equation}
C_{[m-1,l+1]}=\max_{i=0,...,m-2}C_{i+1,l+m-i-1}
\label{10.584}
\end{equation}
Finally we have the contribution of the term 
$$A(2\pi\lambda)^{m-1}\cth_{l+m}$$
in \ref{10.568} to the integral on the right in \ref{10.480}. This is bounded by a constant multiple 
of:
\begin{equation}
(\Pi u^2)^{m-1}\int_0^{\ub_1}\ub\|\cth_{l+m}\|_{L^2(S_{\ub,u})}d\ub
\label{10.585}
\end{equation}
It follows that the corresponding contribution to $\|\cnu_{m-1,l+1}-\check{\tau}_{m-1,l+1}\|_{L^2(\Cb_{\ub_1}^{u_1})}$ is 
bounded by a constant multiple of:
\begin{equation}
\ub_1^{3/2}(\Pi u_1^2)^{m-1}\left(\int_0^{\ub_1}\|\cth_{l+m}\|^2_{L^2(\Cb_{\ub_1}^{u_1})}d\ub\right)^{1/2}
\label{10.586}
\end{equation}
Substituting from Proposition 10.3, with $l+m$ in the role of $l$ (there $l=n$ while here $l+m=n$) 
and $\ub$ in the role of $\ub_1$, the estimate for $\|\cth_{l+m}\|_{L^2(\Cb_{\ub}^{u_1})}$, we conclude 
that this is bounded by:
\begin{equation}
\frac{C^\prime}{\sqrt{2a_0+2}}\sqrt{\max\{\s^{(0,l)}\cB(\ub_1,u_1),\s^{(0,l)}\cBb(\ub_1,u_1)\}}\cdot\ub_1^{a_0}u_1^{b_0}\cdot(\Pi u_1^2)^{m-1}\frac{\ub_1^{5/2}}{u_1^{3/2}}
\label{10.587}
\end{equation}
The contribution of $y_{m-1,l+1}$ through the integral on the right in \ref{10.480} to 
$\|\cnu_{m-1,l+1}-\check{\tau}_{m-1,l+1}\|_{L^2(\Cb_{\ub_1}^{u_1})}$ is then bounded by a constant multiple of the sum of 
\ref{10.578}, \ref{10.583} and \ref{10.587}. Combining with the bound \ref{10.567} for all other 
contributions to $\|\cnu_{m-1,l+1}-\check{\tau}_{m-1,l+1}\|_{L^2(\Cb_{\ub_1}^{u_1})}$ yields \ref{10.574} for $j=m$,  
completing the inductive step, if we define $C_{m,l}$ according to: 
\begin{equation}
C_{m,l}=C_0+C_1+\frac{C\delta^{3/2}}{\sqrt{2a_{m-1}+1}} C_{[m-1,l+1]}
\label{10.c3}
\end{equation}
This is a recursion which determines the coefficients $C_{m,l}$. Choosing $\delta$ small enough or 
$a_{m-1}$ large enough so that the coefficient of $C_{[m-1,l+1]}$ on the right is not greater than 
$1/2$, this recursion together with the fact that from \ref{10.573} with $l+m-1$ in 
the role of $l$ we can take:
\begin{equation}
C_{1,l+m-1}=C
\label{10.c4}
\end{equation}
implies that $C_{m,l}$ is bounded by a constant $C$ which is independent of $m$ and $l$. We have thus proved the following proposition. 

\vspace{2.5mm}

\noindent{\bf Proposition 10.5} \ \ Under the assumptions \ref{10.439}, \ref{10.440}, and \ref{10.546}, 
the top order acoustical difference quantities $\cnu_{m-1,l+1}$, $m=1,...,n$, $m+l=n$, satisfy to 
principal terms the following estimates on the $\Cb_{\ub}$:
\begin{eqnarray*}
&&\|\cnu_{m-1,l+1}-\check{\tau}_{m-1,l+1}\|_{L^2(\Cb_{\ub_1}^{u_1})}\leq \ub_1^{a_m}u_1^{b_m}
\left\{C\sqrt{\frac{\s^{(Y;m,l)}\cBb(\ub_1,u_1)}{(2a_m+1)}}\cdot\frac{\ub_1^2}{u_1^3}\right.\\
&&\hspace{25mm}\left.+C\sqrt{\max\{\s^{[m,l]}\cB(\ub_1,u_1),
\s^{[m,l]}\cBb(\ub_1,u_1)\}}\cdot\frac{\ub_1^{1/2}}{u_1}
\right\}
\end{eqnarray*}

\vspace{2.5mm}

As we shall see in Section 10.8, the quantity $\cnub_{m-1,l+1}$ satisfies a boundary condition on 
${\cal K}$ which involves the quantity $\cnu_{m-1,l+1}$. To control appropriately the boundary values 
of $\cnub_{m-1,l+1}$ on ${\cal K}$ we need to derive an appropriate estimate for $\cnu_{m-1,l+1}$ 
in $L^2({\cal K}^\tau)$. Setting $\ub_1=u$ in \ref{10.480} we have:
\begin{equation}
\|\cnu_{m-1,l+1}-\check{\tau}_{m-1,l+1}\|_{L^2(S_{u,u})}\leq k\int_0^{u}\|\check{M}_{m-1,l+1}\|_{L^2(S{\ub,u})}d\ub
\label{10.588}
\end{equation}
and what we want to estimate is:
\begin{equation}
\|\cnu_{m-1,l+1}-\check{\tau}_{m-1,l+1}\|^2_{L^2({\cal K}^\tau)}=
\int_0^\tau\|\cnu_{m-1,l+1}-\check{\tau}_{m-1,l+1}\|^2_{L^2(S_{u,u})}du
\label{10.589}
\end{equation}
We shall consider the leading contribution to the integral on the right in \ref{10.588}, which is that 
of the $\Nb$ component of \ref{10.486} (see \ref{10.488}) and comes from \ref{10.483} through 
\ref{10.484}. This contribution is bounded by:
\begin{equation}
Cu^{-2}\int_0^u\|\Nb^\mu\Lb(E^l T^m\beta_\mu-E_N^lT^m\beta_{\mu,N}\|_{L^2(S_{\ub,u})}\ub d\ub
\label{10.590}
\end{equation}
Recalling \ref{10.511} - \ref{10.515}, \ref{10.590} is bounded by:
\begin{equation}
C(G_0(u,u)+G_1(u,u))
\label{10.591}
\end{equation}
This is completely analogous to \ref{10.456}. We therefore deduce, in analogy with \ref{10.463}, the 
estimate:
\begin{equation}
\int_0^\tau G_0^2(u,u)du\leq\frac{C\s^{(Y;m,l)}\cBb(\tau,\tau)}{(2a_m+1)}\cdot\tau^{2(a_m+b_m)}\cdot\tau^{-2}
\label{10.592}
\end{equation}
The last factor $\tau^{-2}$ being of degree -2 signifies that the contribution is {\em borderline}. 
We also deduce, in analogy with \ref{10.469} the estimate:
\begin{equation}
\int_0^\tau G_1^2(u,u)du\leq\frac{C\s^{(E;m,l)}\cBb(\tau,\tau)}{(2a_m+1)}\cdot\tau^{2(a_m+b_m)}\cdot\tau^{-1}
\label{10.593}
\end{equation}
The last factor $\tau^{-1}$ being of degree -1, which is one degree better than borderline, this 
contribution to $\|\cnu_{m-1,l+1}-\check{\tau}_{m-1,l+1}\|^2_{L^2({\cal K}^\tau)}$ will come with an extra factor of $\delta$. 

Finally, all other contributions to $\|\cnu_{m-1,l+1}-\check{\tau}_{m-1,l+1}\|^2_{L^2({\cal K}^\tau)}$ will come with extra 
factors of $\delta^2$ at least. We thus deduce the following proposition. 

\vspace{2.5mm} 

\noindent{\bf Proposition 10.6}\ \ Under the assumptions \ref{10.439}, \ref{10.440}, and \ref{10.546}, 
the top order acoustical difference quantities $\cnu_{m-1,l+1}$, $m=1,...,n$, $m+l=n$, satisfy 
to principal terms the following estimates on ${\cal K}^\tau$: 
\begin{eqnarray*}
&&\|\cnu_{m-1,l+1}\|_{L^2({\cal K}^\tau)}\leq \tau^{a_m+b_m}
\left\{C\sqrt{\frac{\s^{(Y;m,l)}\cBb(\tau,\tau)}{(2a_m+1)}}\cdot\frac{1}{\tau}\right.\\
&&\hspace{25mm}\left.+C\sqrt{\max\{\s^{[m,l]}\cB(\tau,\tau),
\s^{[m,l]}\cBb(\tau,\tau)\}}\cdot\frac{1}{\tau^{1/2}}
\right\}
\end{eqnarray*}

\vspace{5mm}

\section{Estimates for $\cthb_l$ and $\cnub_{m-1,l+1}$ in Terms of their Boundary Values on ${\cal K}$}

Consider now the propagation equation \ref{10.258} for $\cthb_l$. Let us denote the right hand side by 
$\check{\Sb}_l$, so that the equation reads: 
$$\Lb\cthb_l-\left(\frac{2\Lb\lambdab}{\lambdab}-(l+2)\chib\right)\cthb_l=\check{\Sb}_l$$
or:
\begin{equation}
\Lb(\lambdab^{-2}\cthb_l)+(l+2)\chib(\lambdab^{-2}\cthb_l)=\lambdab^{-2}\check{\Sb}_l
\label{10.594}
\end{equation}
We shall integrate this equation along the integral curves of $\Lb$. Since (see 2nd of \ref{2.33})
$$\Lb=\frac{\partial}{\partial u}+b\frac{\partial}{\partial\vartheta}$$
in $(\ub,u,\vartheta)$ coordinates the integral curve of $\Lb$ through $(\ub,\ub,\vartheta^{\prime\prime})\in {\cal K}$ is:
\begin{equation}
u\mapsto(\ub,u,\Phib_{\ub,u}(\vartheta^{\prime\prime}))
\label{10.595}
\end{equation}
where $\Phib_{\ub,u}(\vartheta^{\prime\prime})$ is the solution of:
\begin{equation}
\frac{d\Phib_{\ub,u}(\vartheta^{\prime\prime})}{du}=b(\ub,u,\Phib_{\ub,u}(\vartheta^{\prime\prime})), 
\ \ \ \Phib_{\ub,\ub}(\vartheta^{\prime\prime})=\vartheta^{\prime\prime}
\label{10.596}
\end{equation}
The mapping 
\begin{equation}
\vartheta^{\prime\prime}\mapsto\vartheta=\Phib_{\ub,u}(\vartheta^{\prime\prime})
\label{10.597}
\end{equation}
is an orientation preserving diffeomorphism of $S^1$. We can use $(\ub,u,\vartheta^{\prime\prime})$ 
as coordinates on ${\cal N}$ in place of the $(\ub,u,\vartheta)$. Recall Section 2.1, equations 
\ref{2.23} - \ref{2.27}. These coordinates are adapted to the flow of $\Lb$, $\vartheta^{\prime\prime}$ 
being constant along the integral  curves of $\Lb$. Hence in these coordinates $\Lb$ is simply given 
by:
\begin{equation}
\Lb=\frac{\partial}{\partial u}
\label{10.598}
\end{equation}
and $E$ is given by:
\begin{equation}
E=\frac{1}{\sqrt{\sh^{\prime\prime}}}\frac{\partial}{\partial\vartheta^{\prime\prime}}
\label{10.599}
\end{equation}
where $\sh^{\prime\prime}d\vartheta^{\prime\prime}\otimes d\vartheta^{\prime\prime}$ is the induced 
metric on the $S_{\ub,u}$ in the same coordinates. The second of the commutation relations \ref{3.a14} 
in the $(\ub,u,\vartheta^{\prime\prime})$ coordinates reads: 
\begin{equation}
\frac{\partial\sh^{\prime\prime}}{\partial u}=2\chib\sh^{\prime\prime}
\label{10.600}
\end{equation}
and we have:
\begin{equation}
\left.\sh^{\prime\prime}\right|_{u=\ub}=\left.\sh\right|_{u=\ub} \ \ \mbox{as} \ \ \Phib_{\ub,\ub}=id
\label{10.601}
\end{equation}
Integrating then \ref{10.600} from ${\cal K}$ we obtain:
\begin{equation}
\sh^{\prime\prime}(\ub,u,\vartheta^{\prime\prime})=\sh(\ub,\ub,\vartheta^{\prime\prime})
e^{2\int_{\ub}^u\chib(\ub,u^\prime,\vartheta^{\prime\prime})du^\prime}
\label{10.602}
\end{equation}
In view of \ref{10.600}, equation \ref{10.594} can be written in $(\ub,u,\vartheta^{\prime\prime})$ 
coordinates in the form:
\begin{equation}
\frac{\partial}{\partial u}\left(\sh^{\prime\prime(l+2)/2}\lambdab^{-2}\cthb_l\right)=
\sh^{\prime\prime(l+2)/2}\lambdab^{-2}\check{\Sb}_l
\label{10.603}
\end{equation}
Integrating from ${\cal K}$ we obtain:
\begin{eqnarray}
&&\left(\sh^{\prime\prime(l+2)/2}\lambdab^{-2}\cthb_l\right)(\ub,u_1,\vartheta^{\prime\prime})=
\left(\sh^{(l+2)/2}\lambdab^{-2}\cthb_l\right)(\ub,\ub,\vartheta^{\prime\prime})\nonumber\\
&&\hspace{30mm}+\int_{\ub}^{u_1}\left(\sh^{\prime\prime(l+2)/2}\lambdab^{-2}\check{\Sb}_l\right)
(\ub,u,\vartheta^{\prime\prime})du
\label{10.604}
\end{eqnarray}

Before we proceed to analyze \ref{10.604}, we make the following remark on the $L^2$ norm of the 
restriction to $S_{\ub,u}$ of a function on ${\cal N}$. If $f$ is a function on $S^1$ its $L^2$ norm 
is:
\begin{equation}
\|f\|_{L^2(S^1)}=\sqrt{\int_{\vartheta^{\prime\prime}\in S^1}f^2(\vartheta^{\prime\prime})
d\vartheta^{\prime\prime}}
\label{10.605}
\end{equation}
On the other hand, if $f$ is a function on ${\cal N}$ and we represent $f$ in 
$(\ub,u,\vartheta^{\prime\prime})$ coordinates, the $L^2$ norm of $f$ on $S_{\ub,u}$ is:
\begin{equation}
\|f\|_{L^2(S_{\ub,u})}=\sqrt{\int_{\vartheta^{\prime\prime}\in S^1}\left(f^2\sqrt{\sh^{\prime\prime}}
\right)(\ub,u,\vartheta^{\prime\prime})d\vartheta^{\prime\prime}}
\label{10.606}
\end{equation}
$\sqrt{\sh^{\prime\prime}(\ub,u,\vartheta^{\prime\prime})}d\vartheta^{\prime\prime}$ being the element 
of arc length on $S_{\ub,u}$. Comparing \ref{10.606} with \ref{10.605} we see that:
\begin{equation}
\|f\|_{L^2(S_{\ub,u})}=\left\|\left(f\sh^{\prime\prime 1/4}\right)(\ub,u)\right\|_{L^2(S^1)}
\label{10.607}
\end{equation}
Here $g$ being an arbitrary function on ${\cal N}$ we denote by $g(\ub,u)$ the restriction of $g$ to 
$S_{\ub,u}$ represented in terms of the $\vartheta^{\prime\prime}$ coordinate. 

With the above remark in mind, we rewrite \ref{10.604} in the form:
\begin{eqnarray}
&&\left(\cthb_l\sh^{\prime\prime 1/4}\right)(\ub,u_1,\vartheta^{\prime\prime})= \nonumber\\
&&\hspace{10mm}\left(\frac{\sh(\ub,\ub,\vartheta^{\prime\prime})}{\sh^{\prime\prime}(\ub,u_1,\vartheta^{\prime\prime})}\right)^{(2l+3)/4}\left(\frac{\lambdab(\ub,u_1,\vartheta^{\prime\prime})}
{\lambdab(\ub,\ub,\vartheta^{\prime\prime})}\right)^2
\left(\cthb\sh^{1/4}\right)(\ub,\ub,\vartheta^{\prime\prime}) \nonumber\\
&&\hspace{5mm}+\int_{\ub}^{u_1}\left(\frac{\sh^{\prime\prime}(\ub,u,\vartheta^{\prime\prime})}
{\sh^{\prime\prime}(\ub,u_1,\vartheta^{\prime\prime})}\right)^{(2l+3)/4}
\left(\frac{\lambdab(\ub,u_1,\vartheta^{\prime\prime})}{\lambdab(\ub,u,\vartheta^{\prime\prime})}\right)^2\left(\check{\Sb}_l\sh^{\prime\prime 1/4}\right)(\ub,u,\vartheta^{\prime\prime})du \nonumber\\
&&\label{10.608}
\end{eqnarray}
Now, by \ref{10.602}:
\begin{eqnarray}
&&\frac{\sqrt{\sh(\ub,\ub,\vartheta^{\prime\prime})}}
{\sqrt{\sh^{\prime\prime}(\ub,u_1,\vartheta^{\prime\prime})}}
=e^{-\int_{\ub}^{u_1}\chib(\ub,u,\vartheta^{\prime\prime})du}, \nonumber\\
&&\frac{\sqrt{\sh^{\prime\prime}(\ub,u,\vartheta^{\prime\prime})}}
{\sqrt{\sh^{\prime\prime}(\ub,u_1,\vartheta^{\prime\prime})}}
=e^{-\int_{u}^{u_1}\chib(\ub,u^\prime,\vartheta^{\prime\prime})du^\prime}
\label{10.609}
\end{eqnarray}
Assuming then a bound on $\chib$ of the form:
\begin{equation}
\sup_{S_{\ub,u}}|\chib|\leq C \ : \ \forall u\in[\ub,u_1]
\label{10.610}
\end{equation}
we have:
\begin{eqnarray}
&&e^{-C(u_1-\ub)}\leq\frac{\sqrt{\sh(\ub,\ub,\vartheta^{\prime\prime})}}
{\sqrt{\sh^{\prime\prime}(\ub,u_1,\vartheta^{\prime\prime})}}\leq e^{C(u_1-\ub)}, \nonumber\\
&&e^{-C(u_1-u)}\leq\frac{\sqrt{\sh^{\prime\prime}(\ub,u,\vartheta^{\prime\prime})}}
{\sqrt{\sh^{\prime\prime}(\ub,u_1,\vartheta^{\prime\prime})}}\leq e^{C(u_1-u)}
\label{10.611}
\end{eqnarray}

Recall the propagation equation for $\lambdab$ of Proposition 3.3. By the 2nd of the assumptions 
\ref{10.440} we have:
\begin{equation}
|\qb|\leq C\lambdab \ : \ \mbox{in} \ {\cal R}_{\delta,\delta}
\label{10.613}
\end{equation}
It then follows that:
\begin{equation}
|\Lb\lambdab|\leq C\lambdab \ : \ \mbox{in} \ {\cal R}_{\delta,\delta}
\label{10.614}
\end{equation}
Integrating along the integral curves of $\Lb$ we then obtain:
\begin{eqnarray}
&&e^{-C(u_1-\ub)}\leq \frac{\lambdab(\ub,u_1,\vartheta^{\prime\prime})}
{\lambdab(\ub,\ub,\vartheta^{\prime\prime})}\leq e^{C(u_1-\ub)} \nonumber\\
&&e^{-C(u_1-u)}\leq \frac{\lambdab(\ub,u_1,\vartheta^{\prime\prime})}
{\lambdab(\ub,u,\vartheta^{\prime\prime})}\leq e^{C(u_1-u)} 
\label{10.615}
\end{eqnarray}

In view of the above, \ref{10.608} implies:
\begin{eqnarray}
&&\left|\left(\cthb_l\sh^{\prime\prime 1/4}\right)(\ub,u_1,\vartheta^{\prime\prime})\right|\leq 
k\left|\left(\cthb_l\sh^{1/4}\right)(\ub,\ub,\vartheta^{\prime\prime})\right| \nonumber\\
&&\hspace{25mm}+k\int_{\ub}^{u_1}\left|\left(\check{\Sb}_l\sh^{\prime\prime 1/4}\right)
(\ub,u,\vartheta^{\prime\prime})\right|du 
\label{10.616}
\end{eqnarray}
where $k$ is a constant greater than 1, but which can be chosen as close to 1 as we wish by suitably 
restricting $u_1$. Taking $L^2$ norms with respect to $\vartheta^{\prime\prime}\in S^1$ this implies:
\begin{eqnarray}
&&\left\|\left(\cthb_l\sh^{\prime\prime 1/4}\right)(\ub,u_1)\right\|_{L^2(S^1)}\leq 
k\left\|\left(\cthb_l\sh^{1/4}\right)(\ub,\ub)\right\|_{L^2(S^1)} \nonumber\\
&&\hspace{25mm}+k\int_{\ub}^{u_1}\left\|\left(\check{\Sb}_l\sh^{\prime\prime 1/4}\right)(\ub,u)
\right\|_{L^2(S^1)}du
\label{10.617}
\end{eqnarray}
or, in view of \ref{10.607} (see also \ref{10.601}):
\begin{equation}
\|\cthb_l\|_{L^2(S_{\ub,u_1})}\leq k\|\cthb_l\|_{L^2(S_{\ub,\ub})}
+k\int_{\ub}^{u_1}\|\check{\Sb}_l\|_{L^2(S_{\ub,u})}du
\label{10.618}
\end{equation}
Taking then $L^2$ norms with respect to $\ub\in[0,\ub_1]$, since
\begin{eqnarray*}
&&\left(\int_0^{\ub_1}\|\cthb_l\|^2_{L^2(S_{\ub,u_1})}d\ub\right)^{1/2}=
\|\cthb_l\|_{L^2(C_{u_1}^{\ub_1})}, \\
&&\left(\int_0^{\ub_1}\|\cthb_l\|^2_{L^2(S_{\ub,\ub})}d\ub\right)^{1/2}=
\|\cthb_l\|_{L^2({\cal K}^{\ub_1})} 
\end{eqnarray*}
we obtain:
\begin{equation}
\|\cthb_l\|_{L^2(C_{u_1}^{\ub_1})}\leq k\|\cthb_l\|_{L^2({\cal K}^{\ub_1})}
+k\left\{\int_0^{\ub_1}\left(\int_{\ub}^{u_1}\|\check{\Sb}_l\|_{L^2(S_{\ub,u})}du\right)^2 d\ub\right\}^{1/2}
\label{10.619}
\end{equation}

We shall presently estimate the contribution of the principal terms in $\check{\Sb}_l$, the right 
hand side of \ref{10.258}, to 
\begin{equation}
\left\{\int_0^{\ub_1}\left(\int_{\ub}^{u_1}\|\check{\Sb}_l\|_{L^2(S_{\ub,u})}du\right)^2 d\ub\right\}^{1/2}
\label{10.a10}
\end{equation}
the 2nd term in \ref{10.619}. The 1st term, which 
represents the contribution of the boundary values of $\cthb_l$ on ${\cal K}$, shall be estimated in 
the next section. 

The principal terms in $\check{\Sb}_l$ are the two difference terms: 
\begin{equation}
-\frac{2\Lb\lambdab}{\lambdab}(E^l\fb-E_N^l\fb_N)+\tilde{\Rb}_l-\tilde{\Rb}_{L,N}
\label{10.620}
\end{equation}
In view of \ref{10.614}, here, in contrast to \ref{10.362}, the coefficient of the first difference 
term is bounded. 
Consider first the partial contribution of the 1st term in $\fb$ as given by Proposition 10.1:
\begin{equation}
-\frac{1}{2}\beta_{\Nb}^2(E^l LH-E_N^l L_N H_N)
\label{10.621}
\end{equation}
the principal part of which is:
\begin{equation}
\frac{1}{2}\beta_{\Nb}^2H^\prime(g^{-1})^{\mu\nu}\beta_\nu L(E^l\beta_\mu-E_N^l\beta_{\mu,N})
\label{10.622}
\end{equation}
As noted in the previous section, the dominant contribution is that of the $\Nb$ component. 
This contribution to the integral on the right in \ref{10.618} is bounded by a constant multiple of:
\begin{equation}
\int_{\ub}^{u_1}\|\Nb^\mu L(E^l\beta_\mu-E_N^l\beta_{\mu,N})\|_{L^2(S_{\ub,u})}du
\label{10.623}
\end{equation}
As in \ref{10.442} we express 
$$\Nb^\mu L(E^l\beta_\mu-E_N^l\beta_{\mu,N})$$ 
up to lower order terms as
$$\lambdab\s^{(E;0,l)}\check{\sxi}$$
Then the square of \ref{10.623} reduces up to lower order terms to:
\begin{equation}
\left(\int_{\ub}^{u_1}\|\lambdab\s^{(E;0,l)}\check{\sxi}\|_{L^2(S_{\ub,u})}du\right)^2
\label{10.624}
\end{equation}
Since $\lambdab\sim\ub$, $\lambda\sim u^2$, $a\sim \ub u^2$, this is bounded by (Schwartz inequality):
\begin{eqnarray}
&&C\ub\left(\int_{\ub}^{u_1}u^{-1}\|\sqrt{a}\s^{(E;0,l)}\check{\sxi}\|_{L^2(S_{\ub,u})}du\right)^2
\leq \label{10.625}\\
&&\hspace{10mm}C\ub\int_{\ub}^{u_1}u^{-2}du\cdot\int_{\ub}^{u_1}\|\sqrt{a}\s^{(E;0,l)}\check{\sxi}\|^2_{L^2(S_{\ub,u})}du\nonumber\\
&&\hspace{10mm}\leq C\|\sqrt{a}\s^{(E;0,l)}\check{\sxi}\|_{\Cb_{\ub}^{u_1}}\leq C\s^{(E;0,l)}\cEb^{u_1}(\ub)
\leq C\s^{(E;0,l)}\cBb(\ub_1,u_1)\ub^{2a_0} u_1^{2b_0}
\nonumber
\end{eqnarray}
by virtue of the definitions \ref{9.291} and \ref{9.298}. Then, integrating with respect to $\ub\in[0,\ub_1]$ we conclude that the corresponding contribution 
to \ref{10.a10} is bounded by:
\begin{equation}
C\sqrt{\frac{\s^{(E;0,l)}\cBb(\ub_1,u_1)}{2a_0+1}}\cdot\ub_1^{a_0}u_1^{b_0}\cdot\ub_1^{1/2}
\label{10.626}
\end{equation}
However, while the contribution of the $\Nb$ component is the dominant contribution overall, that is 
is terms of the degree of the  extra factor, as in \ref{10.455}, it is contribution of the $N$ 
component which has the least extra power of $\ub_1$, hence dominates in a suitable neighborhood of 
$\Cb_0$. This contribution to the integral on the right in \ref{10.618} is bounded by a constant 
multiple of:
\begin{equation}
\int_{\ub}^{u_1}\|N^\mu L(E^l\beta_\mu-E_N^l\beta_{\mu,N})\|_{L^2(S_{\ub,u})}du
\label{10.627}
\end{equation}
To estimate this we write:
$$N^\mu L(E^l\beta_\mu-E_N^l\beta_{\mu,N})=\s^{(Y;0,l)}\check{\xi}_L
-\ogamma\Nb^\mu L(E^l\beta_\mu-E_N^l\beta_{\mu,N})$$
Since $\gamma\sim u$ the contribution of the 2nd term is similar to \ref{10.623} but with an extra 
$u$ factor. As for the contribution of the 1st term we have (Schwartz inequality):
\begin{equation}
\left(\int_{\ub}^{u_1}\|\s^{(Y;0,l)}\check{\xi}_L\|_{L^2(S_{\ub,u})}du\right)^2
\leq u_1\int_{\ub}^{u_1}\|\s^{(Y;0,l)}\check{\xi}_L\|^2_{L^2(S_{\ub,u})}du
\label{10.628}
\end{equation}
Integrating with respect to $\ub\in[0,\ub_1]$ we have, reversing the order of integration,
\begin{eqnarray}
&&\int_0^{\ub_1}\left\{\int_{\ub}^{u_1}\|\s^{(Y;0,l)}\check{\xi}_L\|^2_{L^2(S_{\ub,u})}du\right\}d\ub
\nonumber\\
&&\hspace{15mm}=\int_0^{u_1}\left\{\int_0^{\min\{u,\ub_1\}}\|\s^{(Y;0,l)}\check{\xi}_L\|^2_{L^2(S_{\ub,u})}d\ub
\right\}du\nonumber\\
&&\hspace{15mm}=\int_0^{u_1}\|\s^{(Y;0,l)}\check{\xi}_L\|^2_{L^2(C_u^{\ub_1})}du\nonumber\\
&&\hspace{15mm}\leq C\int_0^{u_1}\s^{(Y;0,l)}\cE^{\ub_1}(u)du\leq C\s^{(Y;0,l)}\cB(\ub_1,u_1)
\ub_1^{2a_0}\cdot\frac{u_1^{2b_0+1}}{2b_0+1}\nonumber\\
&&\label{10.629}
\end{eqnarray}
by virtue of the definitions \ref{9.290} and \ref{9.297}. 
It follows that this contribution to \ref{10.a10} is bounded by:
\begin{equation}
C\sqrt{\frac{\s^{(Y;0,l)}\cB(\ub_1,u_1)}{2b_0+1}}\cdot\ub_1^{a_0}u_1^{b_0}\cdot u_1^{3/2}
\label{10.630}
\end{equation}

We turn to the partial contribution of the 2nd term in $\fb$ as given by Proposition 10.1:
\begin{equation}
\lambdab\beta_{\Nb}\sbeta(E^{l+1}H-E_N^{l+1}H_N)
\label{10.631}
\end{equation}
the principal part of which is:
\begin{equation}
-2\lambdab\beta_{\Nb}\sbeta H^\prime(g^{-1})^{\mu\nu}\beta_\nu E(E^l\beta_\mu-E_N^l\beta_{\mu,N})
\label{10.632}
\end{equation}
We shall estimate the partial contribution of the $\Nb$ component. Since $\lambdab\sim\ub$ this contribution 
to the integral on the right in \ref{10.618} is bounded by a constant multiple of:
\begin{equation}
\ub\int_{\ub}^{u_1}\|\Nb^\mu E(E^l\beta_\mu-E_N^l\beta_{\mu,N})\|_{L^2(S_{\ub,u})}du
\label{10.633}
\end{equation}
To estimate this we write $\Nb^\mu=\rhob^{-1}\Lb^\mu$ and express 
$$\Lb^\mu E(E^l\beta_\mu-E_N^l\beta_{\mu,N})$$
as in \ref{10.429} up to lower order terms as
$$\s^{(E;0,l)}\check{\xi}_{\Lb}$$
Then since $\rhob\sim u^2$ the square of \ref{10.633} is bounded up to lower order terms by 
a constant multiple of:
\begin{eqnarray}
&&\ub^2\left(\int_{\ub}^{u_1}u^{-2}\|\s^{(E;0,l)}\check{\xi}_{\Lb}\|_{L^2(S_{\ub,u})}du\right)^2 \label{10.634}\\
&&\hspace{20mm}\leq\ub^2\int_{\ub}^{u_1}u^{-4}du\cdot.\int_{\ub}^{u_1}
\|\s^{(E;0,l)}\check{\xi}_{\Lb}\|^2_{L^2(S_{\ub,u})}du\nonumber\\
&&\hspace{20mm}\leq \frac{1}{3}\ub^{-1}\|\s^{(E;0,l)}\check{\xi}_{\Lb}\|^2_{L^2(\Cb_{\ub}^{u_1})}
\nonumber\\
&&\hspace{20mm}\leq C\ub^{-1}\s^{(E;0,l)}\cEb^{u_1}(\ub)\leq C\ub^{2a_0-1}u_1^{2b_0}
\s^{(E;0,l)}\cBb(\ub_1,u_1)\nonumber
\end{eqnarray}
by the Schwartz inequality and the definitions \ref{9.291} and \ref{9.298}. Integrating with respect 
to $\ub\in[0,\ub_1]$ we conclude that the corresponding contribution to \ref{10.a10} is bounded by:
\begin{equation}
C\sqrt{\frac{\s^{(E;0,l)}\cBb(\ub_1,u_1)}{2a_0}}\cdot\ub_1^{a_0}u_1^{b_0}
\label{10.635}
\end{equation}
Here the extra factor is simply 1, so this dominates \ref{10.626}. This is in fact the dominant contribution to \ref{10.a10}. The partial contribution of the $E$ component of 
\ref{10.632} to the integral on the right in \ref{10.618} is bounded by a constant multiple of 
\begin{equation}
\ub\int_{\ub}^{u_1}\|E^\mu E(E^l\beta_\mu-E_N^l\beta_{\mu,N})\|_{L^2(S_{\ub,u})}du
\label{10.636}
\end{equation}
which since $a\sim\ub u^2$ is in turn bounded by a constant multiple of:
\begin{equation}
\sqrt{\ub}\int_{\ub}^{u_1}u^{-1}\|\sqrt{a}\s^{(E;0,l)}\check{\sxi}\|_{L^2(S_{\ub,u})}du
\label{10.637}
\end{equation} 
The square of this is bounded by:
\begin{eqnarray}
&&\ub\int_{\ub}^{u_1}u^{-2}du\cdot\int_{\ub}^{u_1}\|\sqrt{a}\s^{(E;0,l)}\check{\sxi}\|^2_{L^2(S_{\ub,u})}du \leq \|\sqrt{a}\s^{(E;0,l)}\check{\sxi}\|^2_{L^2(\Cb_{\ub}^{u_1})}\nonumber\\
&&\hspace{15mm}\leq C\s^{(E;0,l)}\cEb^{u_1}(\ub)\leq C\s^{(E;0,l)}\cBb(\ub_1,u_1)\ub^{2a_0}u_1^{2b_0}
\label{10.638}
\end{eqnarray}
by the Schwartz inequality and the definitions \ref{9.291} and \ref{9.298}. Integrating with respect 
to $\ub\in[0,\ub_1]$ we conclude that the corresponding contribution to \ref{10.a10} is bounded by:
\begin{equation}
C\sqrt{\frac{\s^{(E;0,l)}\cBb(\ub_1,u_1)}{2a_0+1}}\cdot\ub_1^{a_0}u_1^{b_0}\cdot\ub_1^{1/2}
\label{10.639}
\end{equation}
which is similar to \ref{10.626}. As for the partial contribution of the $N$ component of \ref{10.632}, 
this is bounded by a constant multiple of:
\begin{equation}
\int_{\ub}^{u_1}\|L^\mu E(E^l\beta_\mu-E_N^l\beta_{\mu,N})\|_{L^2(S_{\ub,u})}du
\label{10.640}
\end{equation}
Expressing $L^\mu E(E^l\beta_\mu-E_N^l\beta_{\mu,N})$ as in \ref{10.422} as $\s^{(E;0,l)}\check{\xi}_L$ 
up to lower order terms, we estimate this in terms of:
\begin{equation}
\int_{\ub}^{u_1}\|\s^{(E;0,l)}\check{\xi}_L\|_{L^2(S_{\ub,u})}du
\label{10.641}
\end{equation}
the square of which is bounded by:
\begin{equation}
u_1\int_{\ub}^{u_1}\|\s^{(E;0,l)}\check{\xi}_L\|^2_{L^2(S_{\ub,u})}du
\label{10.642}
\end{equation}
This is similar to \ref{10.628} with $E$ in the role of $Y$. Therefore this contribution to 
\ref{10.a10} is bounded by \ref{10.630} with $E$ in the role of $Y$, that is by:
\begin{equation}
C\sqrt{\frac{\s^{(E;0,l)}\cB(\ub_1,u_1)}{2b_0+1}}\cdot\ub_1^{a_0}u_1^{b_0}\cdot u_1^{3/2}
\label{10.643}
\end{equation}

We now turn to the other difference term in \ref{10.620}, namely:
\begin{equation}
\tilde{\Rb}_l-\tilde{\Rb}_{l,N}
\label{10.644}
\end{equation}
From \ref{10.253}, the principal part of this is of the form:
\begin{eqnarray}
&&\left(\lambda(E\lambdab)\underline{c}^\mu_{1,E}+\lambda\lambdab\tchib 
\underline{c}_{2,E}^\mu+\lambda\lambdab \underline{c}_{3,E}^{\mu\nu}(E\beta_\nu) \right.\nonumber\\
&&\hspace{25mm}\left.+\lambdab\underline{c}_{4,E}^{\mu\nu}(\Lb\beta_\nu)+
\lambda\underline{c}_{5,E}^{\mu\nu}(L\beta_\nu)\right)E(E^l\beta_\mu-E_N^l\beta_{\mu,N})
\nonumber\\
&&\hspace{5mm}+\left((E\lambdab)\underline{c}_{1,\Lb}^\mu+\lambdab\tchib\underline{c}_{2,\Lb}^\mu 
+\lambdab\underline{c}_{3,\Lb}^{\mu\nu}(E\beta_\nu)+\underline{c}_{4,\Lb}^{\mu\nu}(L\beta_\nu)\right)
\Lb(E^l\beta_\mu-E_N^l\beta_{\mu,N}) \nonumber\\
&&\hspace{15mm}+\left(\lambda\tchib\underline{c}_{1,L}^\mu+\lambda\underline{c}_{2,L}^{\mu\nu}(E\beta_\nu)
+\underline{c}_{3,L}^{\mu\nu}\Lb\beta_\nu\right)L(E^l\beta_\mu-E_N^l\beta_{\mu,N}) \nonumber\\
&&\label{10.645} 
\end{eqnarray}
The assumptions \ref{10.439}, \ref{10.440}, \ref{10.546} then imply that this is of the form:
\begin{equation}
\ub u^2 k_1^\mu E(E^l\beta_\mu-E_N^l\beta_{\mu,N})
+\ub k_2^\mu\Lb(E^l\beta_\mu-E_N^l\beta_{\mu,N})
+k_3^\mu L(E^l\beta_\mu-E_N^l\beta_{\mu,N}) 
\label{10.646}
\end{equation}
the coefficients $k$ being bounded by a constant. Here the 1st term is similar to the term \ref{10.632} 
but with an extra $u^2$ factor, while the 3rd term is similar to the term \ref{10.622}. Thus only the 
second term needs to be treated. We shall consider the dominant contribution, that of the $\Nb$ 
component. This contribution to the integral on the right in \ref{10.618} is bounded by a constant 
multiple of:
\begin{equation}
\ub\int_{\ub}^{u_1}\|\Nb^\mu\Lb(E^l\beta_\mu-E_N^l\beta_{\mu,N})\|_{L^2(S_{\ub,u})}du
\label{10.647}
\end{equation}
To estimate this we write:
\begin{equation}
\ogamma\Nb^\mu\Lb(E^l\beta_\mu-E_N^l\beta_{\mu,N})=\s^{(Y;0,l)}\check{\xi}_{\Lb}
-N^\mu\Lb(E^l\beta_\mu-E_N^l\beta_{\mu,N})
\label{10.648}
\end{equation}
Since $\ogamma\sim u$ it follows that \ref{10.647} is bounded by a constant multiple of:
\begin{eqnarray}
&&\ub\int_{\ub}^{u_1}u^{-1}\|\s^{(Y;0,l)}\check{\xi}_{\Lb}\|_{L^2(S_{\ub,u})}du 
\nonumber\\
&&+\ub\int_{\ub}^{u_1}u^{-1}\|N^\mu\Lb(\beta_\mu-E_N^l\beta_{\mu,N})\|_{L^2(S_{\ub,u})}du
\label{10.649}
\end{eqnarray}
The square of the first term here is bounded by:
\begin{eqnarray}
&&\ub^2\int_{\ub}^{u_1}u^{-2}du\cdot\int_{\ub}^{u_1}\|\s^{(Y;0,l)}\check{\xi}_{\Lb}\|^2_{L^2(S_{\ub,u})}du\leq \ub\|\s^{(Y;0,l)}\check{\xi}_{\Lb}\|^2_{L^2(\Cb_{\ub}^{u_1})}\nonumber\\
&&\hspace{15mm}\leq C\ub\s^{(Y;0,l)}\cEb^{u_1}(\ub)\leq C\ub^{2a_0+1}u_1^{2b_0}\cBb(\ub_1,u_1)
\label{10.650}
\end{eqnarray}
therefore the corresponding contribution to \ref{10.a10} is bounded by:
\begin{equation}
C\sqrt{\frac{\s^{(Y;0,l)}\cBb(\ub_1,u_1)}{2a_0+2}}\cdot\ub_1^{a_0}u_1^{b_0}\cdot\ub_1
\label{10.651}
\end{equation}
As for the second term in \ref{10.649} expressing $N^\mu\Lb(E^l\beta_\mu-E_N^l\beta_{\mu,N})$ 
as in \ref{10.377} as $\lambda \s^{(E;0,l)}\check{\sxi}$ up to lower order terms, this term 
is seen to be bounded by a constant multiple of:
\begin{equation}
\sqrt{\ub}\int_{\ub}^{u_1}\|\sqrt{a}\s^{(E;0,l)}\check{\sxi}\|_{L^2(S_{\ub,u})}du
\label{10.652}
\end{equation}
the square of which is bounded by:
\begin{eqnarray}
&&\ub u_1\int_{\ub}^{u_1}\|\sqrt{a}\s^{(E;0,l)}\check{\sxi}\|^2_{L^2(S_{\ub,u})}du=
\ub u_1\|\sqrt{a}\s^{(E;0,l)}\check{\sxi}\|^2_{L^2(\Cb_{\ub}^{u_1})}\nonumber\\
&&\hspace{15mm}\leq C\ub u_1\s^{(E;0,l)}\cEb^{u_1}(\ub)\leq C\ub^{2a_0+1} u_1^{2b_0+1}
\label{10.653}
\end{eqnarray}
therefore the corresponding contribution to \ref{10.a10} is bounded by:
\begin{equation}
C\sqrt{\frac{\s^{(Y;0,l)}\cBb(\ub_1,u_1)}{2a_0+2}}\cdot\ub_1^{a_0}u_1^{b_0}\cdot\ub_1 u_1^{1/2}
\label{10.654}
\end{equation}

We summarize the above results, in simplified form, in the following proposition. 

\vspace{2.5mm}

\noindent{\bf Proposition 10.7} \ \ The top order acoustical difference quantities $\cthb_l$ 
satisfy to principal terms the following estimate on the $C_u$:
$$\|\cthb_l\|_{L^2(C_{u_1}^{\ub_1})}\leq k\|\cthb_l\|_{L^2({\cal K}^{\ub_1})}
+\underline{C}\sqrt{\max\{\s^{(0,l)}\cB(\ub_1,u_1),\s^{(0,l)}\cBb(\ub_1,u_1)\}} \cdot 
\ub_1^{a_0}u_1^{b_0}$$

\vspace{2.5mm}

Consider now the propagation equation \ref{10.337} for $\cnub_{m,l}$. Let us denote:
\begin{equation}
\check{\taub}_{m,l}=E^{l-1}T^m(\lambdab c\qb^\prime E\pi)-E_N^{l-1}T^m(\lambdab_N c_N \qb^\prime_N E_N\pi_N)
\label{10.655}
\end{equation}
This is a quantity of order $m+l$. 
Also, let us add to both sides of equation \ref{10.337} the term:
\begin{equation}
\left(\frac{2\Lb\lambdab}{\lambdab}-(l+1)\chib\right)\check{\taub}_{m,l}
\label{10.656}
\end{equation}
Denoting the resulting right hand side with the term \ref{10.656} added by $\check{\Mb}_{m,l}$, the  equation reads:
\begin{equation}
\Lb(\cnub_{m,l}-\check{\taub}_{m,l})-\left(\frac{2\Lb\lambdab}{\lambdab}-(l+1)\chib\right)
(\cnub_{m,l}-\check{\taub}_{m,l})
=\check{\Mb}_{m,l}
\label{10.657}
\end{equation} 
Actually, comparing \ref{10.235} with \ref{10.4} we see that what is involved 
in the $(m,l)$ difference energy estimates for $m\geq 1$ is $\cnub_{m-1,l+1}$. 
Being at the top order, we have $m+l=n$, $m=1,...,n$. Thus we consider 
equation \ref{10.658} with $(m-1,l+1)$ in the role of $(m,l)$:
\begin{eqnarray}
&&\Lb(\cnub_{m-1,l+1}-\check{\taub}_{m-1,l+1})-\left(\frac{2\Lb\lambdab}{\lambdab}-(l+2)\chib\right)(\cnub_{m-1,l+1}-\check{\taub}_{m-1,l+1})\nonumber\\
&&\hspace{50mm}=\check{\Mb}_{m-1,l+1} \nonumber\\
&&\hspace{30mm} : \ m+l=n, \ m=1,...,n
\label{10.658}
\end{eqnarray} 
This is exactly of the same form as the equation preceding \ref{10.594}, with the quantity 
$\cnub_{m-1,l+1}-\check{\taub}_{m-1,l+1}$ in the role of the quantity $\cthb_l$ 
and $\check{\Mb}_{m-1,l+1}$ in the role of $\check{S}_l$. Equation \ref{10.658} can then be written, 
in analogy with \ref{10.603},  
in the $(\ub,u,\vartheta^{\prime\prime})$ coordinates, which are adapted to the flow of $\Lb$, in 
the form:
\begin{equation}
\frac{\partial}{\partial u}\left(\sh^{\prime\prime(l+2)/2}\lambdab^{-2}(\cnub_{m-1,l+1}-\check{\taub}_{m-1,l+1})\right)=\sh^{\prime\prime(l+2)/2}\lambdab^{-2}
\check{\Mb}_{m-1,l+1}
\label{10.659}
\end{equation}
(Equations \ref{10.655} - \ref{10.659} are the conjugates of equations \ref{10.473} - \ref{10.477}.) 
Integrating from ${\cal K}$ we obtain:
\begin{eqnarray}
&&\left(\sh^{\prime\prime(l+2)/2}\lambdab^{-2}(\cnub_{m-1,l+1}-\check{\taub}_{m-1,l+1})\right)
(\ub,u_1,\vartheta^{\prime\prime})\nonumber\\
&&\hspace{10mm}=\left(\sh^{(l+2)/2}\lambdab^{-2}(\cnub_{m-1,l+1}-\check{\taub}_{m-1,l+1})\right)
(\ub,\ub,\vartheta^{\prime\prime})\nonumber\\
&&\hspace{18mm}+\int_{\ub}^{u_1}\sh^{\prime\prime(l+2)/2}\lambdab^{-2}\check{\Mb}_{m-1,l+1}
(\ub,u,\vartheta^{\prime\prime})du
\label{10.670}
\end{eqnarray}
Following the argument leading from \ref{10.604} to the inequality \ref{10.618} we deduce from 
\ref{10.670}, in the same way, the inequality:
\begin{eqnarray}
&&\|\cnub_{m-1,l+1}-\check{\taub}_{m-1,l+1}\|_{L^2(S_{\ub,u_1})}\leq 
\|\cnub_{m-1,l+1}-\check{\taub}_{m-1,l+1}\|_{L^2(S_{\ub,\ub})}\nonumber\\
&&\hspace{45mm}+\int_{\ub}^{u_1}\|\check{\Mb}_{m-1,l+1}\|_{L^2(S_{\ub,u})}du
\label{10.671}
\end{eqnarray}
Taking then $L^2$ norms with respect to $\ub\in[0,\ub_1]$, since
\begin{eqnarray*}
&&\left(\int_0^{\ub_1}\|\cnub_{m-1,l+1}-\check{\taub}_{m-1,l+1}\|^2_{L^2(S_{\ub,u_1})}d\ub\right)^{1/2}=
\|\cnub_{m-1,l+1}-\check{\taub}_{m-1,l+1}\|_{L^2(C_{u_1}^{\ub_1})}, \\
&&\left(\int_0^{\ub_1}\|\cnub_{m-1,l+1}-\check{\taub}_{m-1,l+1}\|^2_{L^2(S_{\ub,\ub})}d\ub\right)^{1/2}=
\|\cnub_{m-1,l+1}-\check{\taub}_{m-1,l+1}\|_{L^2({\cal K}^{\ub_1})} 
\end{eqnarray*}
we obtain:
\begin{eqnarray}
&&\|\cnub_{m-1,l+1}-\check{\taub}_{m-1,l+1}\|_{L^2(C_{u_1}^{\ub_1})}
\leq k\|\cnub_{m-1,l+1}-\check{\taub}_{m-1,l+1}\|_{L^2({\cal K}^{\ub_1})}\nonumber\\
&&\hspace{38mm}+k\left\{\int_0^{\ub_1}\left(\int_{\ub}^{u_1}\|\check{\Mb}_{m-1,l+1}\|_{L^2(S_{\ub,u})}du\right)^2 d\ub\right\}^{1/2}\nonumber\\
&&\label{10.672}
\end{eqnarray}

The quantity $\check{\taub}_{m-1,l+1}$ is of order $m+l=n$. We shall presently estimate the 
contribution of the principal terms in $\check{\Mb}_{m-1,l+1}$ to
\begin{equation}
\left\{\int_0^{\ub_1}\left(\int_{\ub}^{u_1}\|\check{\Mb}_{m-1,l+1}\|_{L^2(S_{\ub,u})}du\right)^2 d\ub\right\}^{1/2}
\label{10.673}
\end{equation}
The first term on the right in \ref{10.672}, which represents the boundary values of 
$\cnub_{m-1,l+1}$ on ${\cal K}$, shall be estimated in the next section. 
The principal terms in $\check{\Mb}_{m-1,l+1}$ (see \ref{10.337}) are the two difference terms: 
\begin{equation}
\frac{2\Lb\lambdab}{\lambdab}\left(E^l T^{m-1}\jb-E_N^l T^{m-1}\jb_N\right)
+\tilde{\Kb}^{\prime\prime}_{m-1,l+1}-\tilde{\Kb}^{\prime\prime}_{m-1,l+1,N}
\label{10.674}
\end{equation}
as well as the terms:
\begin{equation}
(m-1)\zeta\cnub_{m-2,l+2}+\oAb\left\{(2\pi\lambdab)^{m-1}\cthb_{l+m}+2\sum_{i=0}^{m-2}(2\pi\lambdab)^{m-2-i}
\cnub_{i,l+m-i}\right\}
\label{10.675}
\end{equation}
which involve the top order acoustical quantities entering the $(i,l+m-i)$ difference  energy estimates 
for $i=0,...,m-1$. In view of \ref{10.614}, in \ref{10.674}, in contrast to \ref{10.481}, the 
coefficient of the 1st term is bounded. Consider first the partial contribution of the 1st term in 
$\jb$ as given by Proposition 10.2: 
\begin{equation}
\frac{1}{4}\beta_{\Nb}^2(E^l T^{m-1}L^2 H-E_N^l T^{m-1}L_N^2 H_N)
\label{10.676}
\end{equation}
the principal part of which is: 
\begin{equation}
-\frac{1}{2}\beta_{\Nb}^2 H^\prime(g^{-1})^{\mu\nu}\beta_\nu
(E^l T^{m-1}L^2\beta_\mu-E_N^l T^{m-1}L^2_N\beta_{\mu,N})
\label{10.677}
\end{equation}
Recalling that
$$T=L+\Lb=L_N+\Lb_N$$
we write:
\begin{equation}
L^2\beta_\mu=L T\beta_\mu-L\Lb\beta_\mu, \ \ \ 
L_N^2\beta_{\mu,N}=L_N T\beta_{\mu,N}-L_N\Lb_N\beta_{\mu,N}
\label{10.678}
\end{equation}
The contribution of the 1st terms in \ref{10.678} to \ref{10.677} is, to principal terms:
\begin{equation}
-\frac{1}{2}\beta_{\Nb}^2 H^\prime(g^{-1})^{\mu\nu}\beta_\nu L(E^l T^m\beta_\mu-E_N^l T^m\beta_{\mu,N})
\label{10.679}
\end{equation}
while the contribution of the 2nd terms in \ref{10.678} to \ref{10.677} is, following the argument 
of the paragraph after \ref{10.485}, seen to be  
given to principal terms by:
\begin{equation}
\frac{1}{2}\beta_{\Nb}^2 H^\prime(g^{-1})^{\mu\nu}\beta_\nu a E(E^{l+1}T^{m-1}\beta_\mu
-E_N^{l+1}T^{m-1}\beta_{\mu,N})
\label{10.680}
\end{equation}
The dominant partial contribution from \ref{10.679} to the integral on the right in \ref{10.671} 
is that of the $\Nb$ component, which is bounded by a constant multiple of:
\begin{equation}
\int_{\ub}^{u_1}\|\Nb^\mu L(E^l T^m\beta_\mu-E_N^l T^m\beta_{\mu,N})\|_{L^2(S_{\ub,u})}du
\label{10.681}
\end{equation}
To estimate this we use \ref{10.447} to express 
$$\Nb^\mu L(E^l T^m\beta_\mu-E_N^l T^m\beta_{\mu,N})$$
up to lower order terms as
$$\lambdab\s^{(E;m,l)}\check{\sxi}$$
Then the square of \ref{10.681} reduces up to lower order terms to:
\begin{equation}
\left(\int_{\ub}^{u_1}\|\lambdab\s^{(E;m,l)}\check{\sxi}\|_{L^2(S_{\ub,u})}du\right)^2
\label{10.682}
\end{equation}
This is similar to \ref{10.624} with $(m,l)$ in the role of $(0,l)$. Then the contribution to 
\ref{10.673} is bounded by (see \ref{10.626}):
\begin{equation}
C\sqrt{\frac{\s^{(E;m,l)}\cBb(\ub_1,u_1)}{2a_m+1}}\cdot\ub_1^{a_m}u_1^{b_m}\cdot\ub_1^{1/2}
\label{10.683}
\end{equation}
We also consider the partial contribution of the $N$ component of \ref{10.679}. This contribution to 
the integral on the right in \ref{10.671} is bounded by a constant multiple of:
\begin{equation}
\int_{\ub}^{u_1}\|N^\mu L(E^l T^m\beta_\mu-E_N^l T^m\beta_{\mu,N})\|_{L^2(S_{\ub,u})}du
\label{10.684}
\end{equation}
To estimate this we write:
$$N^\mu L(E^l T^m\beta_\mu-E_N^l T^m\beta_{\mu,N})=\s^{(Y;m,l)}\check{\xi}_L
-\ogamma\Nb^\mu L(E^l T^m\beta_\mu-E_N^l T^m\beta_{\mu,N})$$
Since $\gamma\sim u$ the contribution of the 2nd term is similar to \ref{10.681} but with an extra 
$u$ factor. As for the contribution of the 1st term, this is treated as in \ref{10.628} - \ref{10.630}, 
$(m,l)$ being placed here in the role of $(0,l)$. We conclude that this contribution to \ref{10.673} 
is bounded by:
\begin{equation}
C\sqrt{\frac{\s^{(Y;m,l)}\cB(\ub_1,u_1)}{2b_m+1}}\cdot\ub_1^{a_m}u_1^{b_m}\cdot u_1^{3/2}
\label{10.685}
\end{equation}
Similarly, the partial contribution of the $E$ component of \ref{10.679} to the integral on the right 
in \ref{10.671} being bounded by a constant multiple of:
$$\int_{\ub}^{u_1}\|\s^{(E;m,l)}\check{\xi}_L\|_{L^2(S_{\ub,u})}du$$
the corresponding contribution to \ref{10.673} is bounded by:
\begin{equation}
C\sqrt{\frac{\s^{(E;m,l)}\cB(\ub_1,u_1)}{2b_m+1}}\cdot\ub_1^{a_m}u_1^{b_m}\cdot u_1^{3/2}
\label{10.a11}
\end{equation}
Note that the extra factor in \ref{10.685}, \ref{10.a11}, although of higher degree, 
dominates the extra factor in \ref{10.683} 
in a suitable neighborhood of $\Cb_0$, as in \ref{10.630} compared to \ref{10.626}. The dominant 
partial contribution from \ref{10.680} to the integral on the right in \ref{10.671}, 
that of the $\Nb$ component, which is bounded by a constant multiple of:
\begin{equation}
\int_{\ub}^{u_1}\|a\Nb^\mu E(E^{l+1}T^{m-1}\beta_{\mu}-E_N^{l+1}T^{m-1}\beta_{\mu,N})\|_{L^2(S_{\ub,u})}du
\label{10.686}
\end{equation}
As in the treatment of \ref{10.519} we write $a\Nb^\mu=\lambdab \Lb^\mu$ and express
$$\Lb^\mu E(E^{l+1}T^{m-1}\beta_\mu-E_N^{l+1}T^{m-1}\beta_{\mu,N})$$ 
up to lower order terms as:
$$E^\mu\Lb(E^{l+1}T^{m-1}\beta_\mu-E_N^{l+1}T^{m-1}\beta_{\mu,N})=\s^{(E;m-1,l+1)}\check{\xi}_{\Lb}$$
Then since $\lambdab\sim\ub$, \ref{10.686} is seen to be bounded, up to lower order terms, by a 
constant multiple of:
\begin{equation}
\ub\int_{\ub}^{u_1}\|\s^{(E;m-1,l+1)}\check{\xi}_{\Lb}\|_{L^2(S_{\ub,u})}du 
\label{10.687}
\end{equation}
which is in turn bounded by:
\begin{eqnarray}
&&\ub\sqrt{u_1}\|\s^{(E;m-1,l+1)}\check{\xi}_{\Lb}\|_{L^2(\Cb_{\ub}^{u_1})}
\leq C\ub\sqrt{u_1}\sqrt{\s^{(E;m-1,l+1)}\cEb^{u_1}(\ub)} \nonumber\\
&&\hspace{30mm}\leq C\ub^{a_{m-1}+1}u_1^{b_{m-1}+\frac{1}{2}}\sqrt{\s^{(E;m-1,l+1)}\cBb(\ub_1,u_1)} 
\nonumber\\
&&\label{10.688}
\end{eqnarray}
We then conclude that the corresponding contribution to \ref{10.673} is bounded by:
\begin{equation}
C\sqrt{\frac{\s^{(E;m-1,l+1)}\cBb(\ub_1,u_1)}{2a_{m-1}+3}}\cdot\ub_1^{a_{m-1}}u_1^{b_{m-1}}\cdot\ub_1^{3/2}u_1^{1/2}
\label{10.689}
\end{equation}

We now briefly consider the contributions of the terms which arise from the 2nd and 3rd terms in $\jb$ 
as given by Proposition 10.2. These have extra factors of $\lambdab$ and $\lambdab^2$ respectively 
therefore it is obvious that these contributions will be dominated by the contributions already 
considered. The principal part of the term arising from the 2nd term in $\jb$ is:
\begin{equation}
\frac{1}{2}\pi\lambdab\beta_{\Nb}^2 H^\prime (g^{-1})^{\mu\nu}\beta_\nu L(E^{l+1}T^{m-1}\beta_\mu 
-E_N^{l+1}T^{m-1}\beta_{\mu,N})
\label{10.690}
\end{equation}
The dominant contribution of this, that of the $\Nb$ component, to the integral on the right in 
\ref{10.671} is bounded by a constant multiple of:
\begin{equation}
\int_{\ub}^{u_1}\|\lambdab \Nb^\mu L(E^{l+1}T^{m-1}\beta_\mu-E_N^{l+1}T^{m-1}\beta_{\mu,N})\|_{L^2(S_{\ub,u})}du
\label{10.691}
\end{equation} 
which is bounded in terms of (see \ref{10.552} - \ref{10.554}):
\begin{equation}
\int_{\ub}^{u_1}\|\lambdab^2 \s^{(E;m-1,l+1)}\check{\sxi}\|_{L^2(S_{\ub,u})}du
\label{10.692}
\end{equation}
This is bounded by a constant multiple of
\begin{eqnarray}
&&\ub^{3/2}\int_{\ub}^{u_1}u^{-1}\|\sqrt{a}\s^{(E;m-1,l+1)}\check{\sxi}\|_{L^2(S_{\ub,u})}du\nonumber\\
&&\leq \ub^{3/2}\left(\int_{\ub}^{u_1}u^{-2}du\right)^{1/2}\|\sqrt{a}\s^{(E;m-1,l+1)}\check{\sxi}
\|_{L^2(\Cb_{\ub}^{u_1})} \nonumber\\
&&\leq C\ub\sqrt{\cEb^{u_1}(\ub)}\leq C\ub^{a_{m-1}+1}u_1^{b_{m-1}}\s^{(E;m-1,l+1)}\cBb(\ub_1,u_1)
\label{10.693}
\end{eqnarray}
The corresponding contribution to \ref{10.673} is then bounded by:
\begin{equation}
C\sqrt{\frac{\s^{(E;m-1,l+1)}\cBb(\ub_1,u_1)}{2a_{m-1}+3}}\cdot \ub_1^{a_{m-1}}u_1^{b_{m-1}}\cdot 
\ub_1^{3/2}
\label{10.694}
\end{equation}

The principal part of the term arising from the 3rd term in $\jb$ is:
\begin{eqnarray}
&&\lambdab^2\left\{-\beta_{\Nb}\left[\pi\sbeta-\frac{1}{4c}(\beta_{\Nb}+2\beta_N)\right](g^{-1})^{\mu\nu}
\beta_\nu+H\left[\pi\sbeta-\frac{1}{2c}(\beta_N+\beta_{\Nb})\right]\Nb^\mu\right\}\nonumber\\
&&\hspace{30mm}\cdot E(E^{l+1}T^{m-1}\beta_\mu-E_N^{l+1}T^{m-1}\beta_{\mu,N})
\label{10.695}
\end{eqnarray}
The dominant contribution of this, that of the $\Nb$ component, to the integral on the right in 
\ref{10.671} is bounded by a constant multiple of:
\begin{equation}
\int_{\ub}^{u_1}\|\lambdab^2\Nb^\mu E(E^{l+1}T^{m-1}\beta_\mu-E_N^{l+1}T^{m-1}\beta_{\mu,N})
\|_{L^2(S_{\ub,u})}du
\label{10.696}
\end{equation}
which, since $\lambdab\sim\ub$, $\rhob\sim u^2$ and $\ub\leq u$, is in turn bounded by a constant 
multiple of:
\begin{equation}
\int_{\ub}^{u_1}\|\Lb^\mu E(E^{l+1}T^{m-1}\beta_\mu-E_N^{l+1}T^{m-1}\beta_{\mu,N})\|_{L^2(S_{\ub,u})}du
\label{10.697}
\end{equation}
which (see \ref{10.522}) up to lower order terms is:
\begin{eqnarray}
&&\int_{\ub_1}^{u_1}\|\s^{(E;m-1,l+1)}\check{\xi}_{\Lb}\|_{L^2(S_{\ub,u})}du
\leq u_1^{1/2}\|\s^{(E;m-1,l+1)}\check{\xi}_{\Lb}\|_{L^2(\Cb_{\ub}^{u_1})} \nonumber\\
&&\hspace{40mm}\leq  Cu_1^{1/2}\sqrt{\s^{(E;m-1,l+1)}\cEb^{u_1}(\ub)} \label{10.698}\\
&&\hspace{35mm}\leq C\ub^{a_{m-1}}u_1^{b_{m-1}+\frac{1}{2}}\sqrt{\s^{(E;m-1,l+1)}\cBb(\ub_1,u_1)} \nonumber
\end{eqnarray}
The corresponding contribution to \ref{10.673} is then bounded by:
\begin{equation}
C\sqrt{\frac{\s^{(E;m-1,l+1)}\cBb(\ub_1,u_1)}{2a_{m-1}+1}}\cdot\ub_1^{a_{m-1}}u_1^{b_{m-1}}\cdot 
\ub_1^{1/2}u_1^{1/2}
\label{10.699}
\end{equation}

We now turn to the other difference term in \ref{10.674}, namely:
\begin{equation}
\tilde{\Kb}^{\prime\prime}_{m-1,l+1}-\tilde{\Kb}^{\prime\prime}_{m-1,l+1,N}
\label{10.700}
\end{equation}
Augmenting the assumptions \ref{10.439}, \ref{10.440}, \ref{10.546} with the assumption that in 
${\cal R}_{\delta,\delta}$ it holds: 
\begin{equation}
|L\lambdab|\leq C
\label{10.701}
\end{equation}
we deduce from \ref{10.324} that the principal part of \ref{10.700} can be written in the form:
\begin{eqnarray}
&&[\tilde{\Kb}^{\prime\prime}_{m-1,l+1}-\tilde{\Kb}^{\prime\prime}_{m-1,l+1,N}]_{P.P.}=\nonumber\\
&&\hspace{17mm}u^2 \kb_1^\mu E(E^{l+1}T^{m-1}\beta_\mu-E_N^{l+1}T^{m-1}\beta_{\mu,N})\nonumber\\
&&\hspace{17mm}+\ub^2 \kb_2^\mu \Lb(E^{l+1}T^{m-1}\beta_\mu-E_N^{l+1}T^{m-1}\beta_{\mu,N})\nonumber\\
&&\hspace{17mm}+\ub \kb_3^\mu L(E^{l+1} T^{m-1}\beta_\mu-E_N^{l+1}T^{m-1}\beta_{\mu,N})\nonumber\\
&&\hspace{17mm}+\kb_4^\mu(E^l T^{m-1} L^2\beta_\mu-E_N^l T^{m-1}L_N^2\beta_{\mu,N})\nonumber\\
&&\hspace{15mm}+\left[u^2\ub \kb^{\prime\mu}_1 E(E^{l+2}T^{m-2}\beta_\mu 
-E_N^{l+2} T^{m-2}\beta_{\mu,N})\right.\nonumber\\
&&\hspace{19mm}+ \ub\kb^{\prime\mu}_2 \Lb(E^{l+2}T^{m-2}\beta_\mu
-E_N^{l+2} T^{m-2}\beta_{\mu,N})\nonumber\\
&&\hspace{19mm}\left.+u^2 \kb^{\prime\mu}_3 L(E^{l+2}T^{m-2}\beta_\mu-E_N^{l+2}T^{m-2}\beta_{\mu,N})
\right]\nonumber\\
&&+\sum_{i=0}^{m-2}(\hb\ub)^{m-2-i}\ub\left\{u\ub \kb^{\prime\prime\mu}_1 E(E^{l+m-i}T^i\beta_\mu 
-E_N^{l+m-i} T^i\beta_{\mu,N})\right.\nonumber\\
&&\hspace{19mm}+\ub \kb^{\prime\prime\mu}_2 \Lb(E^{l+m-i} T^i\beta_\mu
-E_N^{l+m-i} T^i\beta_{\mu,N})\nonumber\\
&&\hspace{19mm}+\ub \kb^{\prime\prime\mu}_3 L(E^{l+m-i}T^i\beta_\mu
-E_N^{l+m-i}T^i\beta_{\mu,N})\nonumber\\
&&\hspace{19mm}\left.+\kb^{\prime\prime\mu}_4(E^{l-1+m-i}T^i L^2\beta_\mu
-E_N^{l-1+m-i}T^i L_N^2\beta_{\mu,N})\right\} \nonumber\\
&&\label{10.702}
\end{eqnarray}
Here the coefficients $\kb$, $\kb^\prime$, $\kb^{\prime\prime}$ and the function $\hb$ in the sum are 
bounded by a constant. We have used the fact that by virtue of the assumptions \ref{10.439},  \ref{10.440} the coefficient $\oAb$ defined by \ref{10.105} satisfies:
\begin{equation}
|\oAb|\leq C\ub \ \ \mbox{: in ${\cal R}_{\delta,\delta}$}
\label{10.703}
\end{equation}

Consider first the first four terms on the right in \ref{10.702}. The 4th term is similar to the 
term \ref{10.677}. The 3rd term is similar to the term \ref{10.690}. The dominant partial 
contribution from the 2nd term, that of the $\Nb$ component, to the integral on the right in 
\ref{10.671} is bounded by a constant multiple of:
\begin{equation}
\ub^2\int_{\ub}^{u_1}\|\Nb^\mu\Lb(E^{l+1}T^{m-1}\beta_{\mu}-E_N^{l+1} T^{m-1}\beta_{\mu,N})\|_{L^2(S_{\ub,u})}du
\label{10.704}
\end{equation}
To estimate this we write:
\begin{eqnarray}
&&\ogamma\Nb^\mu\Lb(E^{l+1}T^{m-1}\beta_{\mu}-E_N^{l+1} T^{m-1}\beta_{\mu,N})
=\s^{(Y;m-1,l+1)}\check{\xi}_{\Lb} \label{10.705}\\
&&\hspace{40mm}-N^\mu\Lb(E^{l+1}T^{m-1}\beta_{\mu}-E_N^{l+1} T^{m-1}\beta_{\mu,N})\nonumber
\end{eqnarray}
The 2nd term on the right can be expressed up to lower order terms as $\lambda\s^{(E;m-1,l+1)}\check{\sxi}$. Its contribution to \ref{10.704} is then up to lower order 
terms bounded by a constant multiple of:
\begin{eqnarray}
&&\ub^{3/2}\int_{\ub}^{u_1}\|\sqrt{a}\s^{(E;m-1,l+1)}\check{\sxi}\|_{L^2(S_{\ub,u})}du \nonumber\\
&&\hspace{20mm}\leq \ub^{3/2} u_1^{1/2}\left(\int_{\ub}^{u_1}\|\sqrt{a}\s^{(E;m-1,l+1)}\check{\sxi}\|_{L^2(S_{\ub,u})}du\right)^{1/2} \nonumber\\
&&\hspace{20mm}=\ub^{3/2} u_1^{1/2}\|\sqrt{a}\s^{(E;m-1,l+1)}\check{\sxi}\|_{L^2(\Cb_{\ub}^{u_1})}
\nonumber\\
&&\hspace{20mm}\leq C\ub^{3/2}u_1^{1/2}\sqrt{\s^{(E;m-1,l+1)}\cEb^{u_1}(\ub)}\nonumber\\
&&\hspace{20mm}\leq C\ub^{a_{m-1}+\frac{3}{2}}u_1^{b_{m-1}+\frac{1}{2}}\sqrt{\s^{(E;m-1,l+1)}\cBb(\ub_1,u_1)}
\label{10.706}
\end{eqnarray}
This partial contribution to \ref{10.673} is then bounded by:
\begin{equation}
C\sqrt{\frac{\s^{(E;m-1,l+1)}\cBb(\ub_1,u_1)}{2a_{m-1}+4}}\cdot\ub_1^{a_{m-1}}u_1^{b_{m-1}}
\cdot \ub_1^2 u_1^{1/2}
\label{10.707}
\end{equation}
On the other hand, the contribution of the 1st term on the right in \ref{10.705} to \ref{10.704} 
is bounded by:
\begin{eqnarray}
&&\ub^2\int_{\ub}^{u_1}u^{-1}\|\s^{(Y;m-1,l+1)}\check{\xi}_{\Lb}\|_{L^2(S_{\ub,u})}du \nonumber\\
&&\hspace{20mm}\leq \ub^2\left(\int_{\ub}^{u-1}u^{-2}du\right)^{1/2}\left(\int_{\ub}^{u_1}
\|\s^{(Y;m-1,l+1)}\check{\xi}_{\Lb}\|^2_{L^2(S_{\ub,u})}du\right)^{1/2} \nonumber\\
&&\hspace{20mm}\leq\ub^{3/2}\|\s^{(Y;m-1,l+1)}\check{\xi}_{\Lb}\|_{L^2(\Cb_{\ub}^{u_1})} \nonumber\\
&&\hspace{20mm}\leq C\ub^{3/2}\sqrt{\s^{(Y;m-1,l+1)}\cEb^{u_1}(\ub)}\nonumber\\
&&\hspace{20mm}\leq C\ub^{a_{m-1}+\frac{3}{2}}u_1^{b_{m-1}}\sqrt{\s^{(Y;m-1,l+1)}\cBb(\ub_1,u_1)}
\label{10.708}
\end{eqnarray}
and the contribution to \ref{10.673} is bounded by:
\begin{equation}
C\sqrt{\frac{\s^{(Y;m-1,l+1)}\cBb(\ub_1,u_1)}{2a_{m-1}+4}}\cdot \ub_1^{a_{m-1}}u_1^{b_{m-1}}\cdot 
\ub_1^2
\label{10.709}
\end{equation}
We turn to the 1st term on the right in \ref{10.702}. Since $\rhob\sim u^2$ the contribution 
of the $\Nb$ component to the integral on the right in \ref{10.671} is bounded by a constant 
multiple of:
\begin{equation}
\int_{\ub}^{u_1}\|\Lb^\mu E(E^{l+1}T^{m-1}\beta_\mu-E_N^{l+1}T^{m-1}\beta_{\mu,N})\|_{L^2(S_{\ub,u})}du
\label{10.710}
\end{equation}
which (see \ref{10.522}) up to lower order terms is:
\begin{equation}
\int_{\ub}^{u_1}\|\s^{(E;m-1,l+1)}\check{\xi}_{\Lb}\|_{L^2(S_{\ub,u})}du
\label{10.711}
\end{equation}
This is bounded by:
\begin{eqnarray}
&&u_1^{1/2}\|\s^{(E;m-1,l+1)}\check{\xi}_{\Lb}\|_{L^2(\Cb_{\ub}^{u_1})}
\leq Cu_1^{1/2}\sqrt{\s^{(E;m-1,l+1)}\cEb^{u_1}(\ub)}\nonumber\\
&&\hspace{20mm}\leq C\ub^{a_{m-1}}u_1^{b_{m-1}+\frac{1}{2}}\sqrt{\s^{(E;m-1,l+1)}\cBb(\ub_1,u_1)}
\label{10.712}
\end{eqnarray}
The corresponding contribution to \ref{10.673} is then bounded by:
\begin{equation}
C\sqrt{\frac{\s^{(E;m-1,l+1)}\cBb(\ub_1,u_1)}{2a_{m-1}+1}}\cdot \ub_1^{a_{m-1}}u_1^{b_{m-1}}
\cdot \ub_1^{1/2}u_1^{1/2}
\label{10.713}
\end{equation}
Here, we also need to consider the contribution of the $E$ and $N$ components because 
these contributions will dominate that of the $\Nb$ component in a suitable neighborhood of $\Cb_0$. 
Since 
$$E^\mu E(E^{l+1}T^{m-1}\beta_\mu-E_N^{l+1}T^{m-1}\beta_{\mu,N})=\s^{(E;m-1,l+1)}\check{\sxi}$$
and $a\sim \ub u^2$, the contribution of the $E$ component to the integral on the right in \ref{10.671} 
is bounded by a constant multiple of:
\begin{eqnarray}
&&\ub^{-1/2}\int_{\ub}^{u_1}u\|\sqrt{a}\s^{(E;m-1,l+1)}\check{\sxi}\|_{L^2(S_{\ub,u})}du 
\nonumber\\
&&\hspace{10mm}\leq\ub_1^{-1/2}\left(\int_{\ub}^{u_1}u^2 du\right)^{1/2}\left(
\int_{\ub}^{u_1}u\|\sqrt{a}\s^{(E;m-1,l+1)}\check{\sxi}\|^2_{L^2(S_{\ub,u})}du\right)^{1/2}\nonumber\\
&&\hspace{10mm}=\frac{1}{\sqrt{3}}\ub_1^{-1/2}u_1^{3/2}\|\sqrt{a}\s^{(E;m-1,l+1)}\check{\sxi}\|_{L^2(\Cb_{\ub}^{u_1})}
\nonumber\\
&&\hspace{10mm}\leq C\ub_1^{-1/2}u_1^{3/2}\sqrt{\s^{(E;m-1,l+1)}\cEb^{u_1}(\ub)}\nonumber\\
&&\hspace{10mm}\leq C\ub_1^{a_{m-1}-\frac{1}{2}}u_1^{b_{m-1}+\frac{3}{2}}\sqrt{\s^{(E;m-1,l+1)}\cBb(\ub_1,u_1)} \label{10.714}
\end{eqnarray}
The corresponding contribution to \ref{10.673} is then bounded by:
\begin{equation}
C\sqrt{\frac{\s^{(E;m-1,l+1)}\cBb(\ub_1,u_1)}{2a_{m-1}}}\cdot \ub_1^{a_{m-1}}u_1^{b_{m-1}}
\cdot u_1^{3/2}
\label{10.715}
\end{equation}
The contribution of the $N$ component to the integral on the right in \ref{10.671} is bounded by 
a constant multiple of:
\begin{equation}
\int_{\ub}^{u_1}u^2\|N^\mu E(E^{l+1}T^{m-1}\beta_\mu-E_N^{l+1}T^{m-1}\beta_{\mu,N})\|_{L^2(S_{\ub,u})}
du
\label{10.716}
\end{equation}
To estimate this we write:
\begin{eqnarray}
&&N^\mu E(E^{l+1}T^{m-1}\beta_\mu-E_N^{l+1}T^{m-1}\beta_{\mu,N})=\s^{(Y;m-1,l+1)}\check{\sxi}
\label{10.717}\\
&&\hspace{40mm}-\ogamma\Nb^\mu E(E^{l+1}T^{m-1}\beta_\mu-E_N^{l+1}T^{m-1}\beta_{\mu,N})\nonumber
\end{eqnarray}
The partial contribution of the 2nd term on the right is similar to that of the $\Nb$ component 
but with an extra $u$ factor, while the partial contribution of the 1st term on the right is 
estimated in the same way as \ref{10.714} with $Y$ in the role of $E$, hence the corresponding 
contribution to \ref{10.673} is bounded by:
\begin{equation}
C\sqrt{\frac{\s^{(Y;m-1,l+1)}\cBb(\ub_1,u_1)}{2a_{m-1}}}\cdot \ub_1^{a_{m-1}}u_1^{b_{m-1}}
\cdot u_1^{3/2}
\label{10.718}
\end{equation}
Summarizing, the contribution to \ref{10.673} of the first four terms on the right in \ref{10.702} 
is bounded by:
\begin{eqnarray}
&&C\sqrt{\max\{\s^{(m,l)}\cB(\ub_1,u_1),\s^{(m,l)}\cBb(\ub_1,u_1)\}}\cdot\ub_1^{a_m}u_1^{b_m}\cdot 
u_1^{1/2} \label{10.719} \\
&&+C\sqrt{\max\{\s^{(m-1,l+1)}\cB(\ub_1,u_1),\s^{(m-1,l+1)}\cBb(\ub_1,u_1)\}}\cdot\ub_1^{a_{m-1}}u_1^{b_{m-1}}\cdot u_1 \nonumber
\end{eqnarray}

Consider next the 5th term on the right in \ref{10.702}, which consists of the three terms in 
the square bracket. The first of these terms is similar to the 1st term on the right in \ref{10.702}, 
but with $(m-2,l+2)$ in the role of $(m-1,l+1)$ and with an extra $\ub$ factor. Hence the partial 
contribution of this term to \ref{10.673} is bounded by:
\begin{equation}
C\sqrt{\max\{\s^{(m-2,l+2)}\cB(\ub_1,u_1),\s^{(m-2,l+2)}\cBb(\ub_1,u_1)\}}\cdot\ub_1^{a_{m-2}}u_1^{b_{m-2}}\cdot \ub_1 u_1
\label{10.720}
\end{equation} 
The second of the terms in the square bracket is similar to the 2nd term on the right in \ref{10.702}, 
but with $(m-2,l+2)$ in the role of $(m-1,l+1)$ and with one less $\ub$ factor. It follows, 
in view of \ref{10.707}, \ref{10.709} that the partial contribution of this term to \ref{10.673} is bounded by:
\begin{equation}
C\sqrt{\s^{(m-1,l+1)}\cBb(\ub_1,u_1)}\cdot \ub_1^{a_{m-2}}u_1^{b_{m-2}}\cdot 
\ub_1 \label{10.721}
\end{equation}
Finally the third term in the square bracket is similar to the 3rd term on the right in \ref{10.702}, 
but with $(m-2,l+2)$ in the role of $(m-1,l+1)$ and with the factor $u^2$ in place of the factor $\ub$. 
The partial contribution of the $\Nb$ component of this term to the integral on the right in 
\ref{10.671} is bounded by a constant multiple of:
\begin{equation}
\int_{\ub}^{u_1}\|\lambda \Nb^\mu L(E^{l+2}T^{m-2}\beta_\mu-E_N^{l+2}T^{m-2}\beta_{\mu,N})\|_{L^2(S_{\ub,u})}du
\label{10.722}
\end{equation}
Comparing with \ref{10.691} - \ref{10.697} we conclude that this is bounded up to lower order terms by:
\begin{eqnarray}
&&C\ub^{1/2}\int_{\ub}^{u_1}u\|\sqrt{a}\s^{(E;m-2,l+2)}\check{\sxi}\|_{L^2(S_{\ub,u})}du \label{10.723}\\
&&\leq C\ub^{1/2}\left(\int_{\ub}^{u_1}u^2du\right)^{1/2}\|\sqrt{a}\s^{(E;m-2,l+2)}\check{\sxi}
\|_{L^2(\Cb_{\ub}^{u_1})} \nonumber\\
&&\leq C\ub^{1/2}u_1^{3/2}\sqrt{\cEb^{u_1}(\ub)}\leq 
\ub^{a_{m-2}+\frac{1}{2}}u_1^{b_{m-2}+\frac{3}{2}}\s^{(E;m-2,l+2)}\cBb(\ub_1,u_1) \nonumber
\end{eqnarray}
The corresponding contribution to \ref{10.673} is then bounded by:
\begin{equation}
C\sqrt{\frac{\s^{(E;m-2,l+2)}\cBb(\ub_1,u_1)}{2a_{m-2}+2}}\cdot\ub_1^{a_{m-2}}u_1^{2b_{m-2}}\cdot 
\ub_1 u_1^{3/2}
\label{10.724}
\end{equation}
Here we must also consider the partial contributions of the $N$ and $E$ components. 
The contribution of the $N$ component to 
the integral on the right in \ref{10.671} is bounded by a constant multiple of: 
\begin{equation}
\int_{\ub}^{u_1}\|\lambda N^\mu L(E^{l+2}T^{m-2}\beta_\mu-E_N^{l+2}T^{m-2}\beta_{\mu,N})\|_{L^2(S_{\ub,u})}du
\label{10.725}
\end{equation}
To estimate this we write:
\begin{eqnarray}
&&N^\mu L(E^{l+2}T^{m-2}\beta_\mu-E_N^{l+2}T^{m-2}\beta_{\mu,N})=\s^{(Y;m-2,l+2)}\check{\xi}_L
\label{10.726}\\
&&\hspace{40mm}-\ogamma\Nb^\mu L(E^{l+2}T^{m-2}\beta_\mu-E_N^{l+2}T^{m-2}\beta_{\mu,N})\nonumber
\end{eqnarray}
The partial contribution to \ref{10.725} of the 2nd term on the right is similar to \ref{10.723} but with an extra $u_1$ factor, as $\ogamma\sim u$. On the other hand, the partial contribution to 
\ref{10.725} of the 1st term on the right is bounded in square by a constant multiple of:
\begin{equation}
\left(\int_{\ub}^{u_1}u^2\|\s^{(Y;m-2,l+2)}\check{\xi}_L\|_{L^2(S_{\ub,u})}du\right)^2
\leq \frac{1}{5}u_1^5\cdot\int_{\ub}^{u_1}
\|\s^{(Y;m-2,l+2)}\check{\xi}_L\|^2_{L^2(S_{\ub,u})}du
\label{10.727}
\end{equation}
Integrating with respect to $\ub\in[0,\ub_1]$ we have, reversing the order of integration, 
\begin{eqnarray}
&&\int_0^{\ub_1}\left\{\int_{\ub}^{u_1}\|\s^{(Y;m-2,l+2)}\check{\xi}_L\|^2_{L^2(S_{\ub,u})}du\right\}d\ub
\nonumber\\
&&\hspace{15mm}=\int_0^{u_1}\left\{\int_0^{\min\{u,\ub_1\}}\|\s^{(Y;m-2,l+2)}\check{\xi}_L\|^2_{L^2(S_{\ub,u})}d\ub
\right\}du\nonumber\\
&&\hspace{15mm}=\int_0^{u_1}\|\s^{(Y;m-2,l+2)}\check{\xi}_L\|^2_{L^2(C_u^{\ub_1})}du
\leq C\int_0^{u_1}\s^{(Y;m-2,l+2)}\cE^{\ub_1}(u)du\nonumber\\
&&\hspace{30mm}\leq C\s^{(Y;m-2,l+2)}\cB(\ub_1,u_1)
\ub_1^{2a_{m-2}}\cdot\frac{u_1^{2b_{m-2}+1}}{2b_{m-2}+1}\nonumber\\
&&\label{10.728}
\end{eqnarray}
(compare with \ref{10.629}). Hence the corresponding contribution to \ref{10.673} is bounded by:
\begin{equation}
C\sqrt{\frac{\s^{(Y;m-2,l+2)}\cB(u_1,\ub_1)}{2b_{m-2}+1}}\cdot \ub_1^{a_{m-2}}u_1^{b_{m-2}}\cdot 
u_1^3
\label{10.729}
\end{equation}
Similarly, the contribution to \ref{10.673} of the $E$ component is bounded by:
\begin{equation}
C\sqrt{\frac{\s^{(E;m-2,l+2)}\cB(u_1,\ub_1)}{2b_{m-2}+1}}\cdot \ub_1^{a_{m-2}}u_1^{b_{m-2}}\cdot 
u_1^3
\label{10.730}
\end{equation}
Summarizing, the contribution to \ref{10.673} of the 5th term on the right in \ref{10.702} is 
bounded by:
\begin{equation}
C\sqrt{\max\{\s^{(m-2,l+2)}\cB(\ub_1,u_1),\s^{(m-2,l+2)}\cBb(\ub_1,u_1)\}}
\cdot\ub_1^{a_{m-2}}u_1^{b_{m-2}}\cdot u_1
\label{10.731}
\end{equation}

We finally turn to the sum in \ref{10.702}. Here we have four terms in the curly bracket which are 
analogous to the first four terms in \ref{10.702} but with $(i,l+m-i)$ in the role of $(m-1,l+1)$. 
The coefficients of the four terms in the curly bracket have extra factors of $\ub^2/u\leq\ub$, 1, 
$\ub$, $\ub$, respectively, relative to the coefficients of the first four terms in \ref{10.702}. 
Taking also into account the overall factor $(\hb\ub)^{m-2-i}$ we then deduce that the contribution of 
the sum in \ref{10.702} to \ref{10.673} is bounded by:
\begin{eqnarray}
&&C\sum_{i=0}^{m-2}(\Mb\ub_1)^{m-2-i}\left\{\sqrt{\max\{\s^{(i,l+m-i)}\cB,
\s^{(i,l+m-i)}\cBb\}}\cdot\ub_1^{a_i}u_1^{b_i}\cdot\ub_1 u_1\right.\nonumber\\
&&\hspace{15mm}\left.+\sqrt{\max\{\s^{(i+1,l-1+m-i)}\cB, 
\s^{(i+1,l-1+m-i)}\cBb\}}\cdot\ub_1^{a_{i+1}}u_1^{b_{i+1}}\cdot\ub_1 u_1^{1/2}\right\}
\nonumber\\
&&\label{10.732}
\end{eqnarray}
$\Mb$ being the bound for $\hb$ in ${\cal R}_{\delta,\delta}$. In view of the non-increasing with 
$m$ property of the exponents $a_m$, $b_m$, the curly parenthesis here is for each $i=0,...,m-2$ 
bounded by:
$$\sqrt{\max\{\s^{[m-1,l+1]}\cB,\s^{[m-1,l+1]}\cBb\}}\cdot\ub_1^{a_{m-1}}u_1^{b_{m-1}}\cdot\ub_1 u_1^{1/2}$$
Since, taking $\delta$ small enough so that $\Mb\delta\leq 1/2$ we have
$$\sum_{i=0}^{m-2}(\Mb\ub_1)^{m-2-i}\leq\sum_{j=0}^\infty(\Mb\delta)^j\leq 2,$$
we conclude that \ref{10.732} is bounded by:
\begin{equation}
C\sqrt{\max\{\s^{[m-1,l+1]}\cB(\ub_1,u_1),\s^{[m-1,l+1]}\cBb(\ub_1,u_1)\}}\cdot\ub_1^{a_{m-1}}u_1^{b_{m-1}}\cdot\ub_1 u_1^{1/2}
\label{10.733}
\end{equation}

We finally come to the contribution to \ref{10.673} of the terms \ref{10.675}. Now, 
\ref{10.675} is equal up to lower order terms to $\yb_{m-1,l+1}$, where:
\begin{eqnarray}
&&\yb_{m-1,l+1}=(m-1)\zeta(\cnub_{m-2,l+2}-\check{\taub}_{m-2,l+2}) \label{10.734}\\
&&\hspace{15mm}+\oAb\left\{(2\pi\lambdab)^{m-1}\cthb_{l+m}+2\sum_{i=0}^{m-2}(2\pi\lambdab)^{m-2-i}
(\cnub_{i,l+m-i}-\check{\taub}_{i,l+m-i})\right\}\nonumber
\end{eqnarray}
The results derived above (\ref{10.683}, \ref{10.685}, \ref{10.a11}, \ref{10.689}, \ref{10.694}, 
\ref{10.699}, \ref{10.719}, \ref{10.731}, \ref{10.734}) can be summarized in the following 
estimate for \ref{10.673}:
\begin{eqnarray}
&&\left\{\int_0^{\ub_1}\left(\int_{\ub}^{u_1}\|\check{\Mb}_{m-1,l+1}\|_{L^2(S_{\ub,u})}du\right)^2 d\ub\right\}^{1/2}\nonumber\\
&&\hspace{20mm}\leq \left\{\int_0^{\ub_1}\left(\int_{\ub}^{u_1}\|\check{\yb}_{m-1,l+1}\|_{L^2(S_{\ub,u})}du\right)^2 d\ub\right\}^{1/2}\nonumber\\
&&\hspace{20mm}+\Cb\sqrt{\max\{\s^{[m,l]}\cB(\ub_1,u_1),\s^{[m,l]}\cBb(\ub_1,u_1)\}}\cdot 
\ub_1^{a_m}u_1^{b_m}\cdot u_1^{1/2}\nonumber\\
&&\label{10.735}
\end{eqnarray}
Substituting this in \ref{10.672} we obtain:
\begin{eqnarray}
&&\|\cnub_{m-1,l+1}-\check{\taub}_{m-1,l+1}\|_{L^2(C_{u_1}^{\ub_1})}\leq 
k\|\cnub_{m-1,l+1}-\check{\taub}_{m-1,l+1}\|_{L^2({\cal K}^{\ub_1})} \label{10.736}\\
&&\hspace{30mm}+k\left\{\int_0^{\ub_1}\left(\int_{\ub}^{u_1}\|\yb_{m-1,l+1}\|_{L^2(S_{\ub,u})}du\right)^2 d\ub\right\}^{1/2}\nonumber\\
&&\hspace{25mm}+\Cb\sqrt{\max\{\s^{[m,l]}\cB(\ub_1,u_1),\s^{[m,l]}\cBb(\ub_1,u_1)\}}\cdot 
\ub_1^{a_m}u_1^{b_m}\cdot u_1^{1/2}\nonumber
\end{eqnarray}

Since $\yb_{m-1,l+1}$ involves the $\cnu_{i,l+m-i}-\check{\taub}_{i,l+m-i}$ for $i=0,...,m-2$, we must proceed by induction. 
In the case $m=1$ \ref{10.734} reduces to:
\begin{equation}
\yb_{0,l+1}=\oAb\cthb_{l+1}
\label{10.737}
\end{equation}
In view of \ref{10.703}, with $C$ the bound for $|\oAb|/\ub$ in ${\cal R}_{\delta,\delta}$, the 2nd term on the right in \ref{10.736} is for $m=1$ bounded by $kC$ times:
\begin{eqnarray}
&&\left\{\int_0^{\ub_1}\ub^2\left(\int_{\ub}^{u_1}\|\cthb_{l+1}\|_{L^2(S_{\ub,u})}du\right)^2 d\ub\right\}^{1/2} \nonumber\\
&&\hspace{10mm}\leq u_1^{1/2}\left\{\int_0^{\ub_1}\left(\int_{\ub}^{u_1}\|\cthb_{l+1}\|_{L^2(S_{\ub,u})}du\right)\ub^2 d\ub\right\}^{1/2}\nonumber\\
&&\hspace{10mm}= u_1^{1/2}\left\{\int_0^{u_1}\left(\int_0^{\min\{u,\ub_1\}}\ub^2 
\|\cthb_{l+1}\|^2_{L^2(S_{\ub,u})}d\ub\right)du\right\}^{1/2}\nonumber\\
&&\hspace{20mm}=u_1^{1/2}\ub_1\left\{\int_0^{u_1}\|\cthb_{l+1}\|^2_{L^2(C_u^{\ub_1})}du
\right\}^{1/2}\label{10.738}
\end{eqnarray}
Substituting from Proposition 10.7, with $l+1$ in the role of $l$ (there $l=n$ while here $l+1=n$) 
and $u$ in the role of $u_1$, the estimate for $\|\cthb_{l+1}\|_{L^2(C_u^{\ub_1})}$, we conclude, in 
view of the monotonicity properties of the functions $\s^{(m,l)}\cB(\ub,u)$, $\s^{(m,l)}\cBb(\ub,u)$, 
that \ref{10.738} is bounded by:
\begin{eqnarray}
&&k\ub_1 u_1\|\cthb_{l+1}\|_{L^2({\cal K}^{\ub_1})} \label{10.739}\\
&&+\Cb\ub_1^{a_0}u_1^{b_0}
\sqrt{\frac{\max\{\s^{(0,l+1)}\cB(\ub_1,u_1),\s^{(0,l+1)}\cBb(\ub_1,u_1)\}}{2b_0+1}}\cdot\ub_1 u_1 \nonumber
\end{eqnarray}
Substituting the resulting bound for the 2nd term on the right in \ref{10.736} in the case $m=1$ we conclude that:
\begin{eqnarray}
&&\|\cnub_{0,l+1}-\check{\taub}_{0,l+1}\|_{L^2(C_{u_1}^{\ub_1})}\leq 
k\|\cnub_{0,l+1}-\check{\taub}_{0,l+1}\|_{L^2({\cal K}^{\ub_1})}
+C_0\ub_1 u_1 \|\cthb_{l+1}\|_{L^2({\cal K}^{\ub_1})} \nonumber\\
&&\hspace{20mm}+\Cb_1\ub_1^{a_1} u_1^{b_1}
\sqrt{\max\{\s^{[1,l]}\cB(\ub_1,u_1),\s^{[1,l]}\cBb(\ub_1,u_1)\}}\cdot u_1^{1/2}\nonumber\\
&& \label{10.740}
\end{eqnarray}
Here,
\begin{equation} 
C_0=k^2 C
\label{10.a12}
\end{equation}
 
Suppose then, as the inductive hypothesis, that for some $m\in\{2,...,n\}$ and all $j=1,...,m-1$ 
the following estimate holds for all $(\ub_1,u_1)\in R_{\delta,\delta}$:
\begin{eqnarray}
&&\|\cnub_{j-1,l+m-j+1}-\check{\taub}_{j-1,l+m-j+1}\|_{L^2(C_{u_1}^{\ub_1})}\leq k\|\cnub_{j-1,l+m-j+1}
-\check{\taub}_{j-1,l+m-j+1}\|_{L^2({\cal K}^{\ub_1})}\nonumber\\
&&\hspace{40mm}+C_{j-1}\ub_1 u_1(\Pib\ub_1)^{j-1}\|\cthb_{l+m}\|_{L^2({\cal K}^{\ub_1})}
\nonumber\\
&&\hspace{25mm}+C^\prime_{j-1}\ub_1 u_1\sum_{i=0}^{j-2}(\Pib\ub_1)^{j-2-i}
\|\cnub_{i,l+m-i}-\check{\taub}_{i,l+m-i}\|_{L^2({\cal K}^{\ub_1})}\nonumber\\
&&\hspace{10mm}+\Cb_{j,l+m-j}\ub_1^{a_j}u_1^{b_j}\sqrt{\max\{\s^{[j,l+m-j]}\cB(\ub_1,u_1),
\s^{[j,l+m-j]}\cBb(\ub_1,u_1)\}}\cdot u_1^{1/2} \nonumber\\
&&\label{10.741}
\end{eqnarray} 
for certain constants $\Cb_{j,l+m-j}$ to be determined recursively. 
Here $\Pib$ is a fixed positive bound for $|2\pi\lambdab/\ub|$ in ${\cal R}_{\delta,\delta}$. 
Now, by \ref{10.575}, denoting by $C^\prime$ the bound for $|\zeta|/\ub$ in 
${\cal R}_{\delta,\delta}$ the contribution of the 1st term in \ref{10.734} to the 2nd term 
on the right in \ref{10.736} is bounded by $(m-1)kC^\prime$ times:
\begin{eqnarray}
&&\left\{\int_0^{\ub_1}\ub^2\left(\int_{\ub}^{u_1}\|\cnub_{m-2,l+2}-\check{\taub}_{m-2,l+2}\|_{L^2(S_{\ub,u})}du\right)^2 d\ub\right\}^{1/2} \label{10.742}\\
&&\hspace{30mm}\leq u_1^{1/2}\ub_1\left\{\int_0^{u_1}\|\cnub_{m-2,l+2}-\check{\taub}_{m-2,l+2}\|^2_{L^2(C_u^{\ub_1})}du
\right\}^{1/2}\nonumber
\end{eqnarray}
(see \ref{10.738}). Substituting the estimate \ref{10.741} for $j=m-1$ with $u$ in the role of $u_1$, 
we deduce, taking account of the monotonicity properties of the functions $\s^{(m,l)}\cB(\ub,u)$, 
$\s^{(m,l)}\cBb(\ub,u)$, that \ref{10.742} is bounded by:
\begin{eqnarray}
&&k\ub_1 u_1\|\cnub_{m-2,l+2}-\check{\taub}_{m-2,l+2}\|_{L^2({\cal K}^{\ub_1})}
+C_{m-2}(\ub_1 u_1)^2 (\Pib\ub_1)^{m-2}\|\cthb_{l+m}\|_{L^2({\cal K}^{\ub_1})}\nonumber\\
&&\hspace{15mm}+C^\prime_{m-2}(\ub_1 u_1)^2\sum_{i=0}^{m-3}(\Pib\ub_1)^{m-3-i}\|\cnub_{i,l+m-i}-\check{\taub}_{i,l+m-i}\|
_{L^2({\cal K}^{\ub_1})}\nonumber\\
&&+\frac{\Cb_{m-1,l+1}}{\sqrt{2b_{m-1}+2}}\ub_1^{a_{m-1}}u_1^{b_{m-1}}
\sqrt{\max\{\s^{[m-1,l+1]}\cB(\ub_1,u_1),\s^{[m-1,l+1]}\cBb(\ub_1,u_1)\}} \cdot \ub_1 u_1^{3/2}
\nonumber\\
&&\label{10.743}
\end{eqnarray}
Next, the contribution of the sum in \ref{10.734} to the 2nd term on the right in \ref{10.736} is 
bounded by $2kC$ times:
\begin{eqnarray}
&&\sum_{i=0}^{m-2}(\Pib\ub_1)^{m-2-i}\left\{\int_0^{\ub_1}\ub^2\left(\int_{\ub}^{u_1}
\|\cnub_{i,l+m-i}-\check{\taub}_{i,l+m-i}\|_{L^2(S_{\ub,u})}du\right)^2 d\ub\right\}^{1/2} \nonumber\\
&&\hspace{20mm}\leq u_1^{1/2}\ub_1\sum_{i=0}^{m-2}(\Pib\ub_1)^{m-2-i}\left\{\int_0^{u_1}
\|\cnub_{i,l+m-i}-\check{\taub}_{i,l+m-i}\|^2_{L^2(C_u^{\ub_1})}du
\right\}^{1/2}\nonumber\\
&&\label{10.744}
\end{eqnarray}
(see \ref{10.738}). Substituting the estimate \ref{10.741} for $j=i+1=1,...,m-1$ with $u$ in the role 
of $u_1$,  we deduce, taking account of the monotonicity properties of the functions 
$\s^{(m,l)}\cB(\ub_1,u_1)$, $\s^{(m,l)}\cBb(\ub_1,u_1)$ and of the non-increasing with $m$ property 
of the exponents $a_m$, $b_m$, that \ref{10.744} is bounded by:
\begin{eqnarray}
&&k\ub_1 u_1\sum_{i=0}^{m-2}(\Pib\ub_1)^{m-2-i}
\|\cnub_{i,l+m-i}-\check{\taub}_{i,l+m-i}\|_{L^2({\cal K}^{\ub_1})} \nonumber\\
&&+\left(\sum_{i=0}^{m-2}C_i\right)\ub_1^2 u_1^2(\Pib\ub_1)^{m-2}
\|\cthb_{l+m}\|_{L^2({\cal K}^{\ub_1})}\nonumber\\
&&+\ub_1^2 u_1^2\sum_{i=1}^{m-2}\sum_{k=0}^{i-1}C^\prime_i(\Pib\ub_1)^{m-3-k}
\|\cnub_{k,l+m-k}-\check{\taub}_{k,l+m-k}\|_{L^2({\cal K}^{\ub_1})} \nonumber\\
&&+\ub_1 u_1\cdot\ub_1^{a_{m-1}} u_1^{b_{m-1}}\left(\sum_{i=0}^{m-2}(\Pib\ub_1)^{m-2-i}\Cb_{i+1,l+m-i-1}\right)\cdot \nonumber\\
&&\hspace{20mm}\sqrt{\frac{\max\{\s^{[m-1,l+1]}\cB(\ub_1,u_1),\s^{[m-1,l+1]}\cBb(\ub_1,u_1)\}}{2b_{m-1}+1}}\nonumber\\
&&\label{10.745}
\end{eqnarray}
Here the double sum is, reversing the order of summation, 
\begin{equation}
\sum_{k=0}^{m-3}\left(\sum_{i=k+1}^{m-2}C^\prime_i\right)(\Pib\ub_1)^{m-3-k}
\|\cnub_{k,l+m-k}-\check{\taub}_{k,l+m-k}\|_{L^2({\cal K}^{\ub_1})}
\label{10.746}
\end{equation}
Finally, we have the contribution of the term $\oAb(2\pi\lambdab)^{m-1}\cthb_{l+m}$ in \ref{10.734} to the 2nd term on the right in \ref{10.736}. This contribution is bounded by $kC$ times:
\begin{eqnarray}
&&(\Pib\ub_1)^{m-1}\left\{\int_0^{\ub_1}\ub^2\left(\int_{\ub}^{u_1}
\|\cthb_{l+m}\|_{L^2(S_{\ub,u})}du\right)^2 d\ub\right\}^{1/2} \nonumber\\
&&\hspace{20mm}\leq u_1^{1/2}\ub_1(\Pib\ub_1)^{m-1}\left\{\int_0^{u_1}
\|\cthb_{l+m}\|^2_{L^2(C_u^{\ub_1})}du\right\}^{1/2}\nonumber\\
&&\label{10.747}
\end{eqnarray}
Substituting from Proposition 10.7, with $l+m$ in the role of $l$ (there $l=n$ while here $l+m=n$) 
and $u$ in the role of $u_1$, the estimate for $\|\cthb_{l+m}\|_{L^2(C_u^{\ub_1})}$, we conclude, in 
view of the monotonicity properties of the functions $\s^{(m,l)}\cB(\ub,u)$, $\s^{(m,l)}\cBb(\ub,u)$, 
that \ref{10.747} is bounded by:
\begin{eqnarray}
&&k\ub_1 u_1(\Pib\ub_1)^{m-1}\|\cthb_{l+m}\|_{L^2({\cal K}^{\ub_1})} \label{10.748}\\
&&+\Cb\ub_1 u_1(\Pib\ub_1)^{m-1}\cdot\ub_1^{a_0}u_1^{b_0}
\sqrt{\frac{\max\{\s^{(0,l+1)}\cB(\ub_1,u_1),\s^{(0,l+1)}\cBb(\ub_1,u_1)\}}{2b_0+1}} \nonumber
\end{eqnarray}
The above, that is \ref{10.743}, \ref{10.745} (noting \ref{10.746}), \ref{10.748}, multiplied by 
the constants $(m-1)kC^\prime$, $2kC$, $kC$, respectively, and summed, bound the 2nd term on the right 
in \ref{10.736}. Substituting this bound in \ref{10.736} yields \ref{10.741} for $j=m$ 
completing the inductive step, if we define $C_{m-1}$, $C^\prime_{m-1}$, $C_{m,l}$ according to:
\begin{equation}
C_{m-1}=k^2 C+(m-1)kC^\prime\Pib^{-1}\delta C_{m-2}+2kC\Pib^{-1}\delta\left(\sum_{i=0}^{m-2}C_i\right)
\label{10.749}
\end{equation}
\begin{equation} 
C^\prime_{m-1}=2k^2 C+(m-1)kC^\prime\Pib^{-1}\delta C^\prime_{m-2}
+2kC\Pib^{-1}\delta\left(\sum_{i=1}^{m-2}C^\prime_i\right) \ \ \ (C^\prime_0=0)
\label{10.750}
\end{equation}
\begin{eqnarray}
&&\Cb_{m,l}=\Cb+\frac{(m-1)kC^\prime\delta^2}{\sqrt{2b_{m-1}+2}}\Cb_{m-1,l+1}
+kC\Cb\delta^{3/2}(\Pib\delta)^{m-1} \nonumber\\
&&\hspace{15mm}+\frac{2kC\delta^{3/2}}{\sqrt{2b_{m-1}+1}}\left(\sum_{i=0}^{m-2}
(\Pib\delta)^{m-2-i}\Cb_{i+1,l+m-i-1}\right)
\label{10.751}
\end{eqnarray}
These are recursions which determine the coefficients $C_{m-1}$, $C^\prime_{m-1}$, $C_{m,l}$. 
In regard to the recursion \ref{10.751}, taking $\Pib\delta\leq\frac{1}{2}$ and defining 
\begin{equation}
\Cb_{[m-1,l+1]}=\max_{i=0,...,m-2}\Cb_{i+1,l+m-i-1}
\label{10.c5}
\end{equation}
we have:
\begin{equation}
\sum_{i=0}^{m-2}(\Pib\delta)^{m-2-i}\Cb_{i+1,l+m-i-1}\leq 2\Cb_{[m-1,l+1]}
\label{10.c6}
\end{equation}
hence \ref{10.751} implies:
\begin{eqnarray}
&&\Cb_{m,l}\leq \Cb\left(1+\frac{k C}{2^{m-1}}\delta^{3/2}\right)\nonumber\\
&&\hspace{10mm}+\frac{k\delta^{3/2}}{\sqrt{2b_{m-1}+1}}\left[4C+(m-1)C^\prime\delta^{1/2}\right]
\Cb_{[m-1,l+1]}
\label{10.c7}
\end{eqnarray}
Choosing 
$b_{m-1}$ large enough so that the coefficient of $\Cb_{[m-1,l+1]}$ on the right is not greater than 
$1/2$, this recursive inequality together with the fact that from \ref{10.740} with $l+m-1$ in 
the role of $l$ we can take:
\begin{equation}
\Cb_{1,l+m-1}=\Cb_1
\label{10.c8}
\end{equation}
in turn implies that $\Cb_{m,l}$ is bounded by a constant $\Cb$ which is independent of $m$ and $l$. We have thus proved the following proposition. 

\vspace{2.5mm}

\noindent{\bf Proposition 10.8} \ \ The top order acoustical difference quantities 
$\cnub_{m-1,l+1}$, $m=1,...,n$, $m+l=n$, satisfy to principal terms 
the following estimates on the $C_u$:
\begin{eqnarray*}
&&\|\cnub_{m-1,l+1}-\check{\taub}_{m-1,l+1}\|_{L^2(C_{u_1}^{\ub_1})}\leq k\|\cnub_{m-1,l+1}
-\check{\taub}_{m-1,l+1}\|_{L^2({\cal K}^{\ub_1})}\\
&&\hspace{32mm}+C_{m-1}\ub_1 u_1(\Pib\ub_1)^{m-1}\|\cthb_{l+m}\|_{L^2({\cal K}^{\ub_1})}\\
&&\hspace{12mm}+C^\prime_{m-1}\ub_1 u_1\sum_{i=0}^{m-2}(\Pib\ub_1)^{m-2-i}
\|\cnub_{i,l+m-i}-\check{\taub}_{i,l+m-i}\|_{L^2({\cal K}^{\ub_1})}\\
&&\hspace{11mm}+\Cb\ub_1^{a_m}u_1^{b_m}\sqrt{\max\{\s^{[m,l]}\cB(\ub_1,u_1),
\s^{[m,l]}\cBb(\ub_1,u_1)\}}\cdot u_1^{1/2}
\end{eqnarray*} 

\vspace{5mm}

\section{Boundary Conditions on ${\cal K}$ and Preliminary Estimates for $\cthb_l$ and $\cnub_{m-1,l+1}$ on ${\cal K}$}

The boundary condition for $\lambdab$ on ${\cal K}$ has been established in Section 4.1. It is 
equation \ref{4.6}:
\begin{equation}
r\lambdab=\lambda \ \ \mbox{: on ${\cal K}$}
\label{10.752}
\end{equation}
Here $r$ is the jump ratio (recall \ref{4.3}):
\begin{equation}
r=-\frac{\ep}{\epb} \ \ \mbox{where} \ \ \ep=N^\mu\triangle\beta_\mu, \ \ \epb=\Nb^\mu\triangle\beta_\mu 
\label{10.753}
\end{equation}
This is a positive function on ${\cal K}$, vanishing at $\partial_-{\cal K}=\partial_-{\cal B}$, 
and expressed in terms of the jump $\epb$ through \ref{4.40} and Proposition 4.2. 

Up to this point we have not established an analogous boundary condition for $\tchib$. 
What we have done in Chapter 5 is to add the 2nd variation and cross variation equations for 
$\tchib$ (see \ref{10.176}, \ref{10.200}) to express $T\tchib$. This gives $T\tchib$ to principal acoustical terms as $2E^2\lambdab+2\pi\lambdab E\tchib$ (see \ref{8.138}). This equation restricted to 
${\cal K}$ determines $\tchib$ on ${\cal K}$, $\lambdab$ being given on ${\cal K}$ by the boundary 
condition \ref{10.752}, in view of the fact that $\tchib$ at 
$\partial_-{\cal K}=\partial_-{\cal B}$ is provided as part of the initial data. In Chapter 5 we 
actually applied $T^k$ to this equation and evaluated the result at $\partial_-{\cal B}$ obtaining  
in this way $\left.T^k\tchib\right|_{\partial_-{\cal B}}$ for $k=1,2,...$. This way of determining 
$\tchib$ on ${\cal K}$ involves however a loss of one derivative and for this reason is inapplicable 
at the top order. 

To derive the boundary condition for $\tchib$ on ${\cal K}$ which applies also at the top order, 
we consider the jump condition \ref{4.1}. In the case of 2 spatial dimensions this reduces to:
\begin{equation}
E^\mu\triangle\beta_\mu=0
\label{10.754}
\end{equation}
Applying $E$ to this we obtain:
$$(EE^\mu)\triangle\beta_\mu+E^\mu E\triangle\beta_\mu=0$$
Substituting for $EE^\mu$ from \ref{3.a31} this becomes, in view of \ref{10.754} and the definitions 
of $\ep$, $\epb$ from \ref{10.753},
\begin{equation}
e\ep+\oe\epb+E^\mu E\triangle\beta_\mu=0
\label{10.755}
\end{equation}
Substituting for $e$, $\oe$ from \ref{3.a32} and taking account of the definition of $r$ from 
\ref{10.753} we deduce:
\begin{equation}
r\left(\tchib+\beta_{\Nb}(\sbeta EH+H\sss)\right)=\tchi+\beta_N(\sbeta EH+H\sss)
+\frac{2c}{\epb}E^\mu E\triangle\beta_{\mu}
\label{10.756}
\end{equation}
This is the appropriate boundary condition for $\tchib$ on ${\cal K}$. 

The analogue of the boundary condition \ref{10.752} for the $N$th approximants is equation \ref{9.102}:
\begin{equation}
r_N\lambdab_N=\lambda_N-\hat{\nu}_N \ \ \mbox{where} \ \ \hat{\nu}_N=O(\tau^{N+1})
\label{10.757}
\end{equation}
the second being the estimate \ref{9.106}. To deduce the analogue of the boundary condition 
\ref{10.756} for the $N$th approximants we similarly consider equation \ref{9.59} in the form:
\begin{equation}
E_N^\mu\triangle_N\beta_{\mu,N}=\iota_{E,N}  \ \ \mbox{where} \ \ \iota_{E,N}=\frac{\iota_{\Omega,N}}{\sqrt{\sh_N}}=O(\tau^{N+1})
\label{10.758}
\end{equation}
the estimate according to Proposition 9.4. Applying $E_N$ to this we obtain:
$$(E_N E_N^\mu)\triangle_N\beta_{\mu,N}+E_N^\mu E_N\triangle_N\beta_{\mu,N}=0$$
Substituting for $E_N E_N^\mu$ the expansion \ref{10.182} this becomes, in view of \ref{10.758} and 
the definitions \ref{9.79},
\begin{equation}
e_N\ep_N+\oe_N\epb_N+E_N^\mu E_N\triangle_N\beta_{\mu,N}=-\se_N\iota_{E,N}
\label{10.759}
\end{equation}
The coefficient $\se_N$ is given by \ref{10.196}. In regard to the coefficient $e_N$  we have:
\begin{eqnarray*}
&&2c_N e_N=-h_{\mu\nu,N}\Nb_N^\nu E_N E_N^\mu\\
&&\hspace{11mm}=h_{\mu\nu,N}(E_N\Nb_N^\mu)E_N^\mu+(E_N h_{\mu\nu,N})E_N^\mu\Nb_N^\nu
\end{eqnarray*}
By the 1st of \ref{10.117} and by \ref{10.142},
$$h_{\mu\nu,N}(E_N \Nb_N^\nu)E_N^\mu=\skb_N=\tchib_N-H_N\sbeta_N\ss_{\Nb,N}$$
Also, 
$$(E_N h_{\mu\nu,N})E_N^\mu\Nb_N^\nu=(E_N H_N)\sbeta_N\beta_{\Nb,N}
+H_N(\sss_N\beta_{\Nb,N}+\ss_{\Nb,N}\sbeta_N)$$
Hence we obtain:
\begin{equation}
2c_N e_N=\tchib_N+\beta_{\Nb,N}(\sbeta_N E_N H_N+H_N\sss_N)
\label{10.760}
\end{equation}
Similarly,
\begin{equation}
2c_N\oe_N=\tchi_N+\beta_{N,N}(\sbeta_N E_N H_N+H_N\sss_N)
\label{10.761}
\end{equation}
Substituting \ref{10.760}, \ref{10.761} in \ref{10.759} and taking account of the definition 
\ref{9.99} of $r_N$ we deduce the boundary condition satisfied by $\tchib_N$:
\begin{eqnarray}
&&r_N\tchib_N=\tchi_N+(\beta_{N,N}-r_N\beta_{\Nb,N})(\sbeta_N E_N H_N+H_N\sss_N) \nonumber\\
&&\hspace{15mm}+2c_N\frac{E_N^\mu E_N\triangle_N\beta_{\mu,N}}{\epb_N}
+2c_N\frac{\se_N\iota_{E,N}}{\epb_N}
\label{10.762}
\end{eqnarray}

From the boundary conditions \ref{10.756} for $\tchib$ and \ref{10.752} for $\lambdab$ 
and the boundary conditions \ref{10.762} and \ref{10.757} satisfied by $\tchib_N$ and $\lambdab_N$
we shall deduce the boundary conditions for the top order acoustical difference quantities 
$\cthb_l$, $l=n$ and $\cnub_{m-1,l+1}$, $m=1,...,n$, $m+l=n$. Applying $E^l$, $l=n$ to the 
boundary condition \ref{10.756} we obtain:
\begin{eqnarray}
&&rE^l\tchib=E^l\tchi+(\beta_N-r\beta_{\Nb})(\sbeta E^{l+1}H+HE^\mu E^{l+1}\beta_\mu)\nonumber\\
&&\hspace{12mm}+2c\frac{E^\mu E^{l+1}\triangle\beta_\mu}{\epb}+{\cal L}_l \ \  \mbox{: on ${\cal K}$}
\label{10.763}
\end{eqnarray}
where ${\cal L}_l$ is of lower order, that is of order $l=n$. Similarly, applying $E_N^l$ to 
\ref{10.762} we obtain:
\begin{eqnarray}
&&r_N E_N^l\tchib_N=E_N^l\tchi_N+(\beta_{N,N}-r_N\beta_{\Nb,N})(\sbeta_N E_N^{l+1}H_N+H_N E_N^\mu E_N^{l+1}\beta_{\mu,N})\nonumber\\
&&\hspace{12mm}+2c_N\frac{E_N^\mu E_N^{l+1}\triangle_N\beta_{\mu,N}}{\epb_N}+{\cal L}_{l,N}
+E_N^l(2c_N\epb_N^{-1}\se_N\iota_{E,N})
\label{10.764}
\end{eqnarray}
where ${\cal L}_{l,N}$ corresponds to ${\cal L}_l$ for the $N$th approximants. In view of \ref{10.752}, 
equation \ref{10.763} takes in  terms of the definitions \ref{10.234} of $\theta_l$, $\thetab_l$ the 
form:
\begin{eqnarray}
&&\thetab_l-E^l\fb=r^{-2}(\theta_l-E^l f)
+\lambdab r^{-1}(\beta_N-r\beta_{\Nb})(\sbeta E^{l+1}H+HE^\mu E^{l+1}\beta_\mu)\nonumber\\
&&\hspace{20mm}+2c\lambdab\epb^{-1}r^{-1}E^\mu E^{l+1}\triangle\beta_\mu+\lambdab r^{-1}{\cal L}_l
\label{10.765}
\end{eqnarray}
Also, in view of \ref{10.757}, equation \ref{10.764} takes in terms of the definitions \ref{10.236} of 
$\theta_{l,N}$, $\thetab_{l,N}$ the form:
\begin{eqnarray}
&&\thetab_{l,N}-E_N^l\fb_N=r_N^{-2}[1-(\hat{\nu}_N/\lambda_N)](\theta_{l,N}-E_N^l f_N)\nonumber\\
&&\hspace{22mm}+\lambdab_N r_N^{-1}(\beta_{N,N}-r_N\beta_{\Nb,N})(\sbeta_N E_N^{l+1}H_N+H_N E_N^\mu E_N^{l+1}\beta_{\mu,N})\nonumber\\
&&\hspace{25mm}+2c_N\lambdab_N\epb_N^{-1}r_N^{-1}E_N^\mu E_N^{l+1}\triangle_N\beta_{\mu,N}
+\lambdab_N r_N^{-1}{\cal L}_{l,N}\nonumber\\
&&\hspace{25mm}+\lambdab_N r_N^{-1}E_N^l(2c_N\epb_N^{-1}\se_N\iota_{E,N})\label{10.766}
\end{eqnarray}
Subtracting then \ref{10.766} from \ref{10.765} we deduce the following boundary condition for 
$\cthb_l$ on ${\cal K}$:
\begin{eqnarray}
&&\cthb_l=r^{-2}(\cth_l-E^l f+E_N^l f_N)+E^l\fb-E_N^l\fb_N \nonumber\\
&&\hspace{5mm}+\lambdab r^{-1}(\beta_N-r\beta_{\Nb})\left[\sbeta(E^{l+1}H-E_N^{l+1}H_N)\right.
\nonumber\\
&&\hspace{45mm}\left.+H(E^\mu E^{l+1}\beta_\mu-E_N^\mu E_N^{l+1}\beta_{\mu,N})\right]\nonumber\\
&&\hspace{5mm}+2c\lambdab\epb^{-1} r^{-1}(E^\mu E^{l+1}\triangle\beta_\mu
-E_N^\mu E_N^{l+1}\triangle_N\beta_{\mu,N})+\check{{\cal L}}_l+\vep_{l,N} \nonumber\\
&& \label{10.767}
\end{eqnarray}
where $\check{{\cal L}}_l$ consists of lower order acoustical difference quantities and:
\begin{equation}
\vep_{l,N}=r_N^{-2}(\hat{\nu}_N/\lambda_N)(\theta_{l,N}-E_N^l f_N)-\lambdab_N r_N^{-1}
E_N^l(2c_N\epb_N^{-1}\se_N\iota_{E,N})=O(\tau^{N-3})
\label{10.768}
\end{equation}

Applying $E^l T^{m-1} E^2$, $m=1,...,n$, $m+l=n$, to \ref{10.752} we obtain:
\begin{equation}
rE^l T^{m-1}E^2\lambdab=E^l T^{m-1}E^2\lambda-\lambdab E^l T^{m-1}E^2 r+{\cal N}_{m,l}
\label{10.769}
\end{equation}
where ${\cal N}_{m,l}$ is of lower order, that is of order $m+l=n$. Similarly, applying 
$E_N^l T^{m-1}E_N^2$ to \ref{10.757} we obtain:
\begin{eqnarray}
&&r_N E_N^l T^{m-1}E_N^2\lambdab_N=E_N^l T^{m-1} E_N^2\lambda_N-\lambdab_N E_N^l T^{m-1} E_N^2 r_N \nonumber\\
&&\hspace{33mm}+{\cal N}_{m,l,N}-E_N^l T^{m-1}E_N^2\hat{\nu}_N
\label{10.770}
\end{eqnarray}
where ${\cal N}_{m,l,N}$ corresponds to ${\cal N}_{m,l}$ for the $N$th approximants. In view of 
\ref{10.752} equation \ref{10.769} takes in terms of the definitions \ref{10.235} of $\nu_{m,l}$, 
$\nub_{m,l}$ the form:
\begin{equation}
\nub_{m-1,l+1}+E^l T^{m-1}\jb=r^{-2}(\nu_{m-1,l+1}+E^l T^{m-1}j)-\lambdab^2 r^{-1} E^l T^{m-1} E^2 r
+\lambdab r^{-1}{\cal N}_{m,l}
\label{10.771}
\end{equation}
Also, in view of \ref{10.757}, equation \ref{10.770} takes in terms of the definitions \ref{10.237} 
of $\nu_{m,l,N}$, $\nub_{m,l,N}$ the form:
\begin{eqnarray}
&&\nub_{m-1,l+1,N}+E_N^l T^{m-1}\jb_N=r_N^{-2}[1-(\hat{\nu}_N/\lambda_N)](\nu_{m-1,l+1,N}
+E_N^l T^{m-1}j_N)\nonumber\\
&&\hspace{20mm}-\lambdab_N^2 r_N^{-1} E_N^l T^{m-1} E_N^2 r_N+\lambdab_N r_N^{-1}{\cal N}_{m,l,N}\nonumber\\
&&\hspace{20mm}-\lambdab_N r_N^{-1}E_N^l T^{m-1}E_N^2\hat{\nu}_N
\label{10.772}
\end{eqnarray}
Subtracting then \ref{10.772} from \ref{10.771} we deduce the following boundary condition for 
$\cnub_{m-1,l+1}$ on ${\cal K}$:
\begin{eqnarray}
&&\cnub_{m-1,l+1}=r^{-2}(\cnu_{m-1,l+1}+E^l T^{m-1}j-E_N^l T^{m-1}j_N)\nonumber\\
&&\hspace{20mm}-E^l T^{m-1}\jb+E_N^l T^{m-1}\jb_N\nonumber\\
&&\hspace{15mm}-\lambdab^2 r^{-1}(E^l T^{m-1} E^2 r-E_N^l T^{m-1} E_N^2 r_N)
+\check{{\cal N}}_{m,l}+\vep_{m,l,N} \nonumber\\
&&\label{10.773}
\end{eqnarray}
where ${\cal N}_{m,l}$ consists of lower order acoustical difference quantities and: 
\begin{eqnarray}
&&\vep_{m,l,N}=r_N^{-2}(\hat{\nu}_N/\lambda_N)(\nu_{m-1,l+1,N}+E_N^l T^{m-1}j_N) \nonumber\\
&&\hspace{15mm}-\lambdab_N r_N^{-1}E_N^l T^{m-1}E_N^2\hat{\nu}_N=O(\tau^{N+2-m})
\label{10.774}
\end{eqnarray}

We shall presently estimate the contribution of the principal terms on the right in \ref{10.767} 
and \ref{10.773} in $L^2({\cal K}^{\tau_1})$. Recalling \ref{4.215}, fixing a constant $k>1$ we 
introduce the fundamental assumption:
\begin{equation}
-\frac{l}{2c_0 k}\tau\leq r\leq -\frac{lk}{2c_0}\tau \ \ : \ \mbox{in ${\cal K}^\delta$}
\label{10.775}
\end{equation}
Consider the 1st term on the right in \ref{10.767}. By assumption \ref{10.775} we have:
\begin{eqnarray}
&&\|r^{-2}\cth_l\|^2_{L^2({\cal K}^{\tau_1})}\leq 
C\int_0^{\tau_1}\tau^{-4}\|\cth_l\|^2_{L^2(S_{\tau,\tau})}d\tau\nonumber\\
&&\hspace{30mm}=C\int_0^{\tau_1}\tau^{-4}\frac{d}{d\tau}\left(\int_0^\tau\|\cth_l\|^2_{L^2(S_{\tau^\prime,\tau^\prime})}\right)d\tau^\prime \label{10.776}\\
&&=C\left\{\tau_1^{-4}\int_0^{\tau_1}\|\cth_l\|^2_{L^2(S_{\tau,\tau})}d\tau
+4\int_0^{\tau_1}\tau^{-5}\left(\int_0^\tau\|\cth_l\|^2_{L^2(S_{\tau^\prime,\tau^\prime})}d\tau^\prime\right)d\tau\right\}\nonumber
\end{eqnarray}
Substituting the estimate for $\|\cth_l\|_{L^2({\cal K}^\tau)}$ of Proposition 10.4 we then obtain:
\begin{eqnarray}
&&\|r^{-2}\cth_l\|_{L^2({\cal K}^{\tau_1})}\leq\tau_1^{a_0+b_0}\left\{
C\sqrt{\frac{\s^{(Y;0,l)}\cBb(\tau_1,\tau_1)}{2a_0+1}}\cdot \frac{1}{\tau_1^3}\right.\nonumber\\
&&\hspace{27mm}\left.+C_{0,l}\sqrt{\max\{\s^{(0,l)}\cB(\tau_1,\tau_1),\s^{(0,l)}\cBb(\tau_1,\tau_1)\}}
\cdot\frac{1}{\tau_1^{5/2}}\right\} \nonumber\\
&&\label{10.777}
\end{eqnarray}
since $2a_0+2b_0-6>0$ in view of \ref{10.391}. In regard to 
$$\|r^{-2}(E^l f-E_N^lf_N)\|_{L^2({\cal K}^{\tau_1})}$$
the leading contribution comes from \ref{10.364}, in fact from the $\Nb$ component of \ref{10.367}. 
By assumption \ref{10.775} this contribution is bounded by:
\begin{equation}
C\left(\int_0^{\tau_1}\tau^{-4}\|\Nb^\mu\Lb(E^l\beta_\mu-E_N^l\beta_{\mu,N})
\|^2_{L^2(S_{\tau,\tau})}d\tau\right)^{1/2}
\label{10.778}
\end{equation}
To estimate this we use \ref{10.375}. Since $\ogamma\sim \tau$ on ${\cal K}$, \ref{10.778} 
is bounded by:
\begin{equation}
C\left(\int_0^{\tau_1}\tau^{-6}\|\ogamma \Nb^\mu\Lb(E^l\beta_\mu-E_N^l\beta_{\mu,N})
\|^2_{L^2(S_{\tau,\tau})}d\tau\right)^{1/2}
\label{10.779}
\end{equation}
By \ref{9.a39} and \ref{9.299} the contribution to this of the right hand side of \ref{10.375} is bounded in square by:
\begin{eqnarray}
&&C\int_0^{\tau_1}\tau^{-6}\|\s^{(Y;0,l)}\check{\xi}_{\Lb}\|^2_{L^2(S_{\tau,\tau})}d\tau\nonumber\\
&&\hspace{30mm}=C\int_0^{\tau_1}\tau^{-6}\frac{d}{d\tau}\left(\int_0^{\tau}
\|\s^{(Y;0,l)}\check{\xi}_{\Lb}\|^2_{S_{\tau^\prime,\tau^\prime})}d\tau^\prime\right)d\tau\nonumber\\
&&=C\left\{\tau_1^{-6}\int_0^{\tau_1}\|\s^{(Y;0,l)}\check{\xi}_{\Lb}\|^2_{L^2(S_{\tau,\tau})}d\tau
+6\int_0^{\tau_1}\tau^{-7}\left(\int_0^{\tau}\|\s^{(Y;0,l)}\check{\xi}_{\Lb}
\|^2_{S_{\tau^\prime,\tau^\prime})}d\tau^\prime\right)d\tau\right\}\nonumber\\
&&\hspace{30mm}\leq C\s^{(Y;0,l)}{\cal A}(\tau_1)\tau_1^{2a_0+2b_0-6}\label{10.780}
\end{eqnarray}
As for the contribution to \ref{10.779} of the 1st term on the left hand side of \ref{10.375}, 
recalling \ref{10.377} and noting that along ${\cal K}$ $\lambda\sim\tau^2$, $\lambdab\sim\tau$, 
$a\sim\tau^3$, this is up to lower order terms bounded in square by:
\begin{eqnarray}
&&C\int_0^{\tau_1}\tau^{-5}\|\sqrt{a}\s^{(E;0,l)}\check{\sxi}\|^2_{S_{\tau,\tau)})}d\tau\nonumber\\
&&\hspace{30mm}=C\int_0^{\tau_1}\tau^{-5}\frac{d}{d\tau}\left(\int_0^{\tau}
\|\sqrt{a}\s^{(E;0,l)}\check{\sxi}\|^2_{S_{\tau^\prime,\tau^\prime})}d\tau^\prime\right)d\tau\nonumber\\
&&=C\left\{\tau_1^{-5}\int_0^{\tau_1}\|\sqrt{a}\s^{(E;0,l)}\check{\sxi}\|^2_{S_{\tau,\tau)})}d\tau
+5\int_0^{\tau_1}\tau^{-6}\left(\int_0^{\tau}
\|\sqrt{a}\s^{(E;0,l)}\check{\sxi}\|^2_{S_{\tau^\prime,\tau^\prime})}d\tau^\prime\right)d\tau\right\}\nonumber\\
&&\hspace{30mm}\leq C\s^{(E;0,l)}{\cal A}(\tau_1)\tau_1^{2a_0+2b_0-5}\label{10.781}
\end{eqnarray}
by \ref{9.a39} and \ref{9.299}. We conclude that \ref{10.778} is bounded by:
\begin{equation}
C\tau_1^{a_0+b_0}\left\{\sqrt{\s^{(Y;0,l)}{\cal A}(\tau_1)}\cdot\frac{1}{\tau_1^3}
+\sqrt{\s^{(0,l)}{\cal A}(\tau_1)}\cdot\frac{1}{\tau_1^{5/2}}\right\}
\label{10.782}
\end{equation}
Of the remaining terms on the right hand side of \ref{10.767}, apart from the 4th, 
which involves the difference $E^\mu E^{l+1}\triangle\beta_{\mu}-E_N^\mu E_N^{l+1}\triangle_N\beta_{\mu,N}$, it is the 3rd term which contains the leading contribution. 
This contribution is bounded by, in view of \ref{10.770}, 
\begin{equation}
C\left(\int_0^{\tau_1}\|\Nb^\mu E(E^l\beta_\mu-E_N^l\beta_{\mu,N})\|^2_{L^2(S_{\tau,\tau})}d\tau
\right)^{1/2}
\label{10.783}
\end{equation}
Since along ${\cal K}$ $\rhob\sim\tau^2$, this is bounded up to lower order terms by:
\begin{equation}
C\left(\int_0^{\tau_1}\tau^{-4}\|\s^{(E;0,l)}\check{\xi}_{\Lb}\|^2_{L^2(S_{\tau,\tau})}d\tau
\right)^{1/2}
\label{10.784}
\end{equation}
which by \ref{9.a39} and \ref{9.299} is bounded in square by:
\begin{eqnarray}
&&C\int_0^{\tau_1}\tau^{-4}\frac{d}{d\tau}\left(\int_0^{\tau}
\|\s^{(E;0,l)}\check{\xi}_{\Lb}\|^2_{L^2(S_{\tau^\prime,\tau^\prime})}d\tau^\prime\right)d\tau
\nonumber\\
&&=C\left\{\tau_1^{-4}\int_0^{\tau_1}\|\s^{(E;0,l)}\check{\xi}_{\Lb}\|^2_{L^2(S_{\tau,\tau})}d\tau
\right.\nonumber\\
&&\hspace{25mm}\left.+4\int_0^{\tau_1}\tau^{-5}\left(\int_0^{\tau}
\|\s^{(E;0,l)}\check{\xi}_{\Lb}\|^2_{L^2(S_{\tau^\prime,\tau^\prime})}d\tau^\prime\right)d\tau\right\}
\nonumber\\
&&\hspace{30mm}\leq C\s^{(E;0,l)}{\cal A}(\tau_1)\tau_1^{2a_0+2b_0-4}\label{10.785}
\end{eqnarray}

Consider next the 1st term on the right in \ref{10.773}. By assumption \ref{10.775} we have:
\begin{eqnarray}
&&\|r^{-2}\cnu_{m-1,l+1}\|^2_{L^2({\cal K}^{\tau_1})}\leq 
C\int_0^{\tau_1}\tau^{-4}\|\cnu_{m-1,l+1}\|^2_{L^2(S_{\tau,\tau})}d\tau\nonumber\\
&&\hspace{20mm}=C\int_0^{\tau_1}\tau^{-4}\frac{d}{d\tau}\left(\int_0^\tau\|\cnu_{m-1,l+1}\|^2_{L^2(S_{\tau^\prime,\tau^\prime})}\right)d\tau^\prime \label{10.786}\\
&&=C\left\{\tau_1^{-4}\int_0^{\tau_1}\|\cnu_{m-1,l+1}\|^2_{L^2(S_{\tau,\tau})}d\tau\right.\nonumber\\
&&\hspace{30mm}\left.+4\int_0^{\tau_1}\tau^{-5}\left(\int_0^\tau\|\cnu_{m-1,l+1}\|^2_{L^2(S_{\tau^\prime,\tau^\prime})}d\tau^\prime\right)d\tau\right\}\nonumber
\end{eqnarray}
Substituting the estimate for $\|\cnu_{m-1,l+1}\|_{L^2({\cal K}^\tau)}$ of Proposition 10.6 we then obtain:
\begin{eqnarray}
&&\|r^{-2}\cnu_{m-1,l+1}\|_{L^2({\cal K}^{\tau_1})}\leq\tau_1^{a_m+b_m}\left\{
C\sqrt{\frac{\s^{(Y;m,l)}\cBb(\tau_1,\tau_1)}{2a_m+1}}\cdot \frac{1}{\tau_1^3}\right.\nonumber\\
&&\hspace{27mm}\left.+C_{m,l}\sqrt{\max\{\s^{[m,l]}\cB(\tau_1,\tau_1),\s^{[m,l]}\cBb(\tau_1,\tau_1)\}}
\cdot\frac{1}{\tau_1^{5/2}}\right\} \nonumber\\
&&\label{10.787}
\end{eqnarray}
since $2a_m+2b_m-6>0$ in view of \ref{10.516}. In regard to 
$$\|r^{-2}(E^l T^{m-1}j-E_N^l T^{m-1}j_N)\|_{L^2({\cal K}^{\tau_1})}$$
the leading contribution comes from \ref{10.486}, in fact from the $\Nb$ its component. 
This contribution is bounded by, in view of \ref{10.770}, 
\begin{equation}
C\left(\int_0^{\tau_1}\tau^{-4}\|\Nb^\mu\Lb(E^l T^m\beta_\mu-E_N^l T^m\beta_{\mu,N})\|
^2_{L^2(S_{\tau,\tau})}d\tau\right)^{1/2}
\label{10.788}
\end{equation}
This is analogous to \ref{10.778} with $(m,l)$ in the role of $(0,l)$. We conclude that \ref{10.788} 
is bounded by:
\begin{equation}
C\tau_1^{a_m+b_m}\left\{\sqrt{\s^{(Y;m,l)}{\cal A}(\tau_1)}\cdot\frac{1}{\tau_1^3}
+\sqrt{\s^{(m,l)}{\cal A}(\tau_1)}\cdot\frac{1}{\tau_1^{5/2}}\right\}
\label{10.789}
\end{equation}
In regard to 
$$\|E^l T^{m-1}\jb-E_N^l T^{m-1}\jb_N\|_{L^2({\cal K}^{\tau_1})}$$
the leading contribution comes from \ref{10.673}, in fact from its $\Nb$ component. 
This contribution is bounded by:
\begin{equation}
C\left(\int_0^{\tau_1}\|\Nb^\mu L(E^l T^m\beta_\mu-E_N^l T^m\beta_{\mu,N})\|
^2_{L^2(S_{\tau,\tau})}d\tau\right)^{1/2}
\label{10.790}
\end{equation}
Using \ref{10.447} to express 
$$\Nb^\mu L(E^l T^m\beta_\mu-E_N^l T^m\beta_{\mu,N})$$
up to lower order terms as 
$$\lambdab\s^{(E;m,l)}\check{\sxi}$$
and noting that along ${\cal K}$ $\lambdab\sim\tau$ while $\sqrt{a}\sim\tau^{3/2}$, \ref{10.790} 
is seen to be bounded up to lower order terms by:
\begin{equation}
C\left(\int_0^{\tau_1}\tau^{-1}\|\sqrt{a}\s^{(E;m,l)}\check{\sxi}\|^2_{L^2(S_{\tau,\tau})}d\tau\right)^{1/2}\leq C\tau_1^{a_m+b_m}\sqrt{\s^{(E;m,l)}{\cal A}(\tau_1)}\cdot \frac{1}{\tau_1^{1/2}}
\label{10.791}
\end{equation}
(compare with \ref{10.781}). 

We finally turn to the terms involving the jump $\triangle\beta_\mu$, that is the 4th term on the right 
in \ref{10.767} and the 3rd term on the right in \ref{10.773}. The last involves $\triangle\beta_\mu$ 
since by \ref{4.40} and Proposition 4.2:
\begin{equation}
r=j(\kappa,\epb)\epb, \ \ \ \epb=\Nb^\mu\triangle\beta_\mu 
\label{10.792}
\end{equation}
(recall that $\kappa$ stands for the quadruplet \ref{4.43}). We can assume that $j$ is bounded on 
${\cal K}^\delta$:
\begin{equation}
|j|\leq C \ : \ \mbox{on ${\cal K}^\delta$}
\label{10.793}
\end{equation}
In view of \ref{10.775} it then follows that the coefficient $2c\lambdab\epb^{-1}r^{-1}$ of the 4th 
term on the right in \ref{10.767} is bounded by $C\tau^{-1}$. Hence the contribution of this term to 
$\|\cthb_l\|_{L^2({\cal K}^{\tau_1})}$ is bounded by a constant multiple of:
\begin{equation}
\left(\int_0^{\tau_1}\tau^{-2}\|E^\mu E^{l+1}\triangle\beta_\mu-E_N^\mu E_N^{l+1}\triangle_N\beta_{\mu,N}\|^2_{L^2(S_{\tau,\tau})}d\tau\right)^{1/2}
\label{10.794}
\end{equation}
Now, from \ref{4.149} and \ref{9.55}:
\begin{eqnarray}
&&E^\mu E^{l+1}\triangle\beta_\mu-E_N^\mu E_N^{l+1}\triangle_N\beta_{\mu,N}=
\left.(E^\mu E^{l+1}\beta_\mu-E_N^\mu E_N^{l+1}\beta_{\mu,N})\right|_{{\cal K}}\nonumber\\
&&\hspace{20mm}-(E^\mu E^{l+1}\beta^\prime_\mu\circ(f,w,\psi)-E_N^\mu E_N^{l+1}\beta^\prime_\mu\circ(f_N,w_N,\psi_N))\nonumber\\
&&\label{10.795}
\end{eqnarray}
the transformation functions $f$, $w$, $\psi$ and their $N$th approximants being functions on 
${\cal K}$ (thus represented as functions of $(\tau,\vartheta)$). The principal part of the 1st term 
on the right in \ref{10.795} is $\left.\s^{(E;0,l)}\check{\sxi}\right|_{{\cal K}}$, and its 
contribution to \ref{10.794} is by \ref{9.a39} and \ref{9.299} bounded in square by a 
constant multiple of (recall that $a\sim\tau^3$ along ${\cal K}$):
\begin{equation}
\int_0^{\tau_1}\tau^{-5}\|\sqrt{a}\s^{(E;0,l)}\check{\sxi}\|^2_{L^2(S_{\tau,\tau})}d\tau
\leq C\s^{(E;0,l)}{\cal A}(\tau_1)\tau_1^{2a_0+2b_0-5}
\label{10.796}
\end{equation}
(see \ref{10.781}). In regard to the 2nd term on the right in \ref{10.795}, we recall from \ref{3.a8}, 
\ref{3.a10}, \ref{9.10} that:
\begin{equation}
E=\frac{1}{\sqrt{\sh}}\Omega, \ \ \ E_N=\frac{1}{\sqrt{\sh_N}}\Omega, \ \ \ 
\Omega=\frac{\partial}{\partial\vartheta}
\label{10.797}
\end{equation}
We also recall from Chapter 4 (see \ref{4.222}, \ref{4.240}, \ref{4.56}, \ref{4.65} and Proposition 4.4) that:
\begin{equation}
f(\tau,\vartheta)=f(0,\vartheta)+\tau^2\hat{f}(\tau,\vartheta), \ \ \ 
w(\tau,\vartheta)=\tau v(\tau,\vartheta), \ \ \ \psi(\tau,\vartheta)=\vartheta+\tau^3\gamma(\tau,\vartheta)
\label{10.798}
\end{equation}
We then have:
\begin{eqnarray}
&&\frac{\partial}{\partial\vartheta}(\beta^\prime_\mu\circ(f,w,\psi))=
\left(\frac{\partial\beta^\prime_\mu}{\partial t}\right)\circ(f,w,\psi)
\left(\frac{\partial f(0,\vartheta)}{\partial\vartheta}+\tau^2\frac{\partial\hat{f}(\tau,\vartheta)}{\partial\vartheta}\right)\nonumber\\
&&\hspace{30mm}+\left(\frac{\partial\beta^\prime_\mu}{\partial u^\prime}\right)\circ(f,w,\psi)
\tau\frac{\partial v(\tau,\vartheta)}{\partial\vartheta}\nonumber\\
&&\hspace{30mm}+\left(\frac{\partial\beta^\prime_\mu}{\partial\vartheta^\prime}\right)\circ(f,w,\psi)
\left(1+\tau^3\frac{\partial\gamma(\tau,\vartheta)}{\partial\vartheta}\right)
\label{10.799}
\end{eqnarray}
Also recall from Chapter 9 (see \ref{9.54} and 1st of \ref{9.114}, \ref{9.116}) that:
\begin{eqnarray}
&&\hspace{15mm}f_N(\tau,\vartheta)=f(0,\vartheta)+\tau^2\hat{f}_N(\tau,\vartheta), \nonumber\\
&&w_N(\tau,\vartheta)=\tau v_N(\tau,\vartheta), \ \ \ 
\psi_N(\tau,\vartheta)=\vartheta+\tau^3\gamma_N(\tau,\vartheta)
\label{10.800}
\end{eqnarray}
We then have:
\begin{eqnarray}
&&\frac{\partial}{\partial\vartheta}(\beta^\prime_\mu\circ(f_N,w_N,\psi_N))=
\left(\frac{\partial\beta^\prime_\mu}{\partial t}\right)\circ(f_N,w_N,\psi_N)
\left(\frac{\partial f(0,\vartheta)}{\partial\vartheta}+\tau^2\frac{\partial\hat{f}_N(\tau,\vartheta)}{\partial\vartheta}\right)\nonumber\\
&&\hspace{35mm}+\left(\frac{\partial\beta^\prime_\mu}{\partial u^\prime}\right)\circ(f_N,w_N,\psi_N)
\tau\frac{\partial v_N(\tau,\vartheta)}{\partial\vartheta}\nonumber\\
&&\hspace{35mm}+\left(\frac{\partial\beta^\prime_\mu}{\partial\vartheta^\prime}\right)
\circ(f_N,w_N,\psi_N)
\left(1+\tau^3\frac{\partial\gamma_N(\tau,\vartheta)}{\partial\vartheta}\right)
\label{10.801}
\end{eqnarray}
Defining then:
\begin{equation}
\chf=\hat{f}-\hat{f}_N, \ \ \ \cv=v-v_N, \ \ \ \cga=\gamma-\gamma_N
\label{10.802}
\end{equation}
the principal part of the 2nd term on the right in \ref{10.795} is:
\begin{eqnarray}
&&-\sh^{-(l+1)/2}E^\mu\left\{\left(\frac{\partial\beta^\prime_\mu}{\partial t}\right)\circ(f,w,\psi)
\tau^2\Omega^{l+1}\chf+\left(\frac{\partial\beta^\prime_\mu}{\partial u^\prime}\right)\circ(f,w,\psi)
\tau \Omega^{l+1}\cv\right.\nonumber\\
&&\hspace{50mm}\left.+\left(\frac{\partial\beta^\prime_\mu}{\partial\vartheta^\prime}\right)
\circ(f,w,\psi)\tau^3\Omega^{l+1}\cga\right\}
\label{10.803}
\end{eqnarray}
Now, by Proposition 4.3:
\begin{equation}
\left.\Omega^\mu\frac{\partial\beta^\prime_\mu}{\partial u^\prime}\right|_{\partial_-{\cal B}}
=\left.(\Omega^{\prime\mu}+\pi L^{\prime\mu})T^\prime\beta^\prime_\mu\right|_{\partial_-{\cal B}}
=\left. T^{\prime\mu}(\Omega^\prime\beta^\prime_\mu+\pi L^\prime\beta^\prime_\mu)\right|_{\partial_-{\cal B}}=0
\label{10.804}
\end{equation}
Therefore we can assume that:
\begin{equation}
\left|E^\mu\left(\frac{\partial\beta^\prime_\mu}{\partial u^\prime}\right)\circ(f,w,\psi)\right|\leq C\tau \ : \ \mbox{on ${\cal K}^\delta$}
\label{10.805}
\end{equation}
It then follows that the contribution of \ref{10.803} to \ref{10.794} is bounded by:
\begin{equation}
C\|\tau\Omega^{l+1}\chf\|_{L^2({\cal K}^{\tau_1})}
+C\|\tau\Omega^{l+1}\cv\|_{L^2({\cal K}^{\tau_1})}
+C\|\tau^2\Omega^{l+1}\cga\|_{L^2({\cal K}^{\tau_1})}
\label{10.806}
\end{equation}
We finally consider the 3rd term on the right in \ref{10.773}. In view of \ref{10.775} and the fact 
that along ${\cal K}$ $\lambdab\sim\tau$, the coefficient $\lambdab^2 r^{-1}$ of this term is 
bounded by $C\tau$. Hence the contribution of this term to 
$\|\cnub_{m-1,l+1}\|_{L^2({\cal K}^{\tau_1})}$ is, to principal terms, bounded by 
a constant multiple of:
\begin{equation}
\left(\int_0^{\tau_1}\tau^2\|E^{l+2}T^{m-1}r-E_N^{l+2}T^{m-1}r_N\|^2_{L^2(S_{\tau,\tau})}d\tau\right)^{1/2}
\label{10.807}
\end{equation}
Now, from \ref{9.99} and \ref{9.111} we have, in analogy with \ref{10.792}:
\begin{equation}
r_N=j(\kappa_N,\epb_N)\epb_N-\hat{d}_N 
\label{10.808}
\end{equation}
where $\kappa_N$ stands for the quadruplet of approximants \ref{9.109} while 
\begin{equation}
\hat{d}_N=\frac{d_N}{\epb_N}=O(\tau^N)
\label{10.809}
\end{equation}
by Proposition 9.6. The leading part of 
$$E^{l+2}T^{m-1}r-E_N^{l+2}T^{m-1}r_N$$
is:
$$j(\kappa,\epb)\left(E^{l+2}T^{m-1}\epb-E_N^{l+2}T^{m-1}\epb_N\right)$$
and the principal part of this is $j(\kappa,\epb)$ times:
\begin{eqnarray}
&&\Nb^\mu E^{l+2}T^{m-1}\triangle\beta_\mu
-\Nb_N^\mu E_N^{l+2}T^{m-1}\triangle_N\beta_{\mu,N}\nonumber\\
&&+\triangle\beta_\mu E^{l+2}T^{m-1}\Nb^\mu 
-\triangle_N\beta_{\mu,N} E_N^{l+2}T^{m-1}\Nb_N^\mu
\label{10.810}
\end{eqnarray}
Here, the first difference is:
\begin{eqnarray}
&&\Nb^\mu E^{l+2}T^{m-1}\triangle\beta_\mu-\Nb_N^\mu E_N^{l+2}T^{m-1}\triangle_N\beta_{\mu,N}=
\nonumber\\
&&\hspace{10mm}\left.(\Nb^\mu E^{l+2}T^{m-1}\beta_\mu-\Nb_N^\mu E_N^{l+2}T^{m-1}\beta_{\mu,N})\right|_{{\cal K}}\nonumber\\
&&-(\Nb^\mu E^{l+2}T^{m-1}\beta^\prime_\mu\circ(f,w,\psi)-\Nb_N^\mu E_N^{l+2}T^{m-1}\beta^\prime_\mu\circ(f_N,w_N,\psi_N))\nonumber\\
&&\label{10.811}
\end{eqnarray}
The 1st term on the right in \ref{10.811}, multiplied by $\left.\rhob\right|_{{\cal K}}\sim\tau^2$, 
is to principal terms equal to:
$$\left.E^\mu\Lb(E^{l+1}T^{m-1}\beta_\mu-E_N^{l+1}T^{m-1}\beta_{\mu,N})\right|_{{\cal K}}
=\left.\s^{(E;m-1,l+1)}\cxi_{\Lb}\right|_{{\cal K}}$$
and its contribution to \ref{10.807} is by \ref{9.a39} and \ref{9.299} bounded by:
\begin{eqnarray}
&&C\left(\int_0^{\tau_1}\tau^{-2}\|\s^{(E;m-1,l+1)}\cxi_{\Lb}\|^2_{L^2(S_{\tau,\tau})}d\tau\right)^{1/2}\nonumber\\
&&\hspace{25mm}\leq C\sqrt{\s^{(E;m-1,l+1)}{\cal A}(\tau_1)}\cdot\tau_1^{a_{m-1}+b_{m-1}-1}
\label{10.812}
\end{eqnarray}
In the case $m=1$ ($l=n-1$) the principal part of 2nd term on the right in \ref{10.811} is equal to:
\begin{eqnarray}
&&-\sh^{-(l+2)/2}\Nb^\mu\left\{\left(\frac{\partial\beta^\prime_\mu}{\partial t}\right)\circ(f,w,\psi)
\tau^2\Omega^{l+2}\chf+\left(\frac{\partial\beta^\prime_\mu}{\partial u^\prime}\right)\circ(f,w,\psi)
\tau \Omega^{l+2}\cv\right.\nonumber\\
&&\hspace{50mm}\left.+\left(\frac{\partial\beta^\prime_\mu}{\partial\vartheta^\prime}\right)
\circ(f,w,\psi)\tau^3\Omega^{l+2}\cga\right\}
\label{10.813}
\end{eqnarray}
(compare with \ref{10.803}), the contribution of which to \ref{10.807} is bounded by:
\begin{equation}
C\|\tau^3\Omega^{l+2}\chf\|_{L^2({\cal K}^{\tau_1})}
+C\|\tau^2\Omega^{l+2}\cv\|_{L^2({\cal K}^{\tau_1})}
+C\|\tau^4\Omega^{l+2}\cga\|_{L^2({\cal K}^{\tau_1})}
\label{10.814}
\end{equation}
In the case $m>1$ we use \ref{10.798} and \ref{10.800} to express:
\begin{eqnarray}
&&\frac{\partial}{\partial\tau}(\beta^\prime_\mu\circ(f,w,\psi))=
\left(\frac{\partial\beta^\prime_\mu}{\partial t}\right)\circ(f,w,\psi)
\left(\tau^2\frac{\partial\hat{f}(\tau,\vartheta)}{\partial\tau}
+2\tau\hat{f}(\tau,\vartheta)\right)\nonumber\\
&&\hspace{30mm}+\left(\frac{\partial\beta^\prime_\mu}{\partial u^\prime}\right)\circ(f,w,\psi)
\left(\tau\frac{\partial v(\tau,\vartheta)}{\partial\tau}+v(\tau,\vartheta)\right)\nonumber\\
&&\hspace{30mm}+\left(\frac{\partial\beta^\prime_\mu}{\partial\vartheta^\prime}\right)\circ(f,w,\psi)
\left(\tau^3\frac{\partial\gamma(\tau,\vartheta)}{\partial\tau}+3\tau^2\gamma(\tau,\vartheta)\right)
\nonumber\\
&&\label{10.815}
\end{eqnarray}
and:
\begin{eqnarray}
&&\frac{\partial}{\partial\tau}(\beta^\prime_\mu\circ(f_N,w_N,\psi_N))=
\left(\frac{\partial\beta^\prime_\mu}{\partial t}\right)\circ(f_N,w_N,\psi_N)
\left(\tau^2\frac{\partial\hat{f}_N(\tau,\vartheta)}{\partial\tau}
+2\tau\hat{f}_N(\tau,\vartheta)\right)\nonumber\\
&&\hspace{35mm}+\left(\frac{\partial\beta^\prime_\mu}{\partial u^\prime}\right)\circ(f_N,w_N,\psi_N)
\left(\tau\frac{\partial v_N(\tau,\vartheta)}{\partial\tau}+v_N(\tau,\vartheta)\right)\nonumber\\
&&\hspace{35mm}+\left(\frac{\partial\beta^\prime_\mu}{\partial\vartheta^\prime}\right)
\circ(f_N,w_N,\psi_N)
\left(\tau^3\frac{\partial\gamma_N(\tau,\vartheta)}{\partial\tau}
+3\tau^2\gamma_N(\tau,\vartheta)\right)\nonumber\\
&&\label{10.816}
\end{eqnarray}
Then in terms of the definitions \ref{10.802} the principal part of the 2nd term on the right in \ref{10.811} is in the case $m>1$ given by:
\begin{eqnarray}
&&-\Nb^\mu\sh^{-(l+2)/2}\left\{\left(\frac{\partial\beta^\prime_\mu}{\partial t}\right)\circ(f,w,\psi)
\cdot \Omega^{l+2}\left(\tau^2 T^{m-1}\chf\right.\right.\nonumber\\
&&\hspace{40mm}\left.+2(m-1)\tau T^{m-2}\chf+(m-1)(m-2)T^{m-3}\chf\right)\nonumber\\
&&\hspace{20mm}+\left(\frac{\partial\beta^\prime_\mu}{\partial u^\prime}\right)\circ(f,w,\psi) \cdot 
\Omega^{l+2}\left(\tau T^{m-1}\cv+(m-1)T^{m-2}\cv\right)\nonumber\\
&&\hspace{20mm}+\left(\frac{\partial\beta^\prime_\mu}{\partial\vartheta^\prime}\right)\circ(f,w,\psi) 
\cdot \Omega^{l+2}\left(\tau^3 T^{m-1}\cga+3(m-1)\tau^2 T^{m-2}\cga\right.\nonumber\\
&&\hspace{20mm}\left.\left.+3(m-1)(m-2)\tau T^{m-3}\cga
+(m-1)(m-2)(m-3)T^{m-4}\cga\right)\right\}\nonumber\\
&&\label{10.817}
\end{eqnarray}
Here we have kept the lower order $T$ derivatives of $\chf$, $\cv$, $\cga$, as these are multiplied 
by correspondingly lower powers of $\tau$. That is, the power of $\tau$ minus the power of $T$ is 
$3-m$ for all three terms involving $\chf$, $2-m$ for the two terms involving $\cv$, $4-m$ for 
all four terms involving $\cga$. It then follows that the contribution of \ref{10.817} to \ref{10.807} 
is bounded by:
\begin{eqnarray}
&&C\left(\|\tau^3\Omega^{l+2}T^{m-1}\chf\|_{L^2({\cal K}^{\tau_1})}
+2(m-1)\|\tau^2 \Omega^{l+2}T^{m-2}\chf\|_{L^2({\cal K}^{\tau_1})}\right.\nonumber\\
&&\hspace{35mm}\left.+(m-1)(m-2)\|\tau\Omega^{l+2}T^{m-3}\chf\|_{L^2({\cal K}^{\tau_1})}\right)\nonumber\\
&&\hspace{3mm}+C\left(\|\tau^2\Omega^{l+2}T^{m-1}\cv\|_{L^2({\cal K}^{\tau_1})}
+(m-1)\|\tau\Omega^{l+2}T^{m-2}\cv\|_{L^2({\cal K}^{\tau_1})}\right)\nonumber\\
&&+C\left(\|\tau^4\Omega^{l+2}T^{m-1}\cga\|_{L^2({\cal K}^{\tau_1})}
+3(m-1)\|\tau^3\Omega^{l+2}T^{m-2}\cga\|_{L^2({\cal K}^{\tau_1})}\right.\nonumber\\
&&\hspace{33mm}+3(m-1)(m-2)\|\tau^2\Omega^{l+2}T^{m-3}\cga\|_{L^2({\cal K}^{\tau_1})}\nonumber\\
&&\hspace{23mm}\left.+(m-1)(m-2)(m-3)\|\tau\Omega^{l+2}T^{m-4}\cga\|_{L^2({\cal K}^{\tau_1})}\right)\nonumber\\
&&\label{10.818}
\end{eqnarray}
We finally consider the second difference in \ref{10.810}: 
\begin{equation}
\triangle\beta_\mu E^{l+2}T^{m-1}\Nb^\mu-\triangle_N\beta_{\mu,N} E_N^{l+2}T^{m-1}\Nb_N^\mu 
\label{10.819}
\end{equation}
Here we shall only consider the contribution of the principal acoustical part. By \ref{3.a15}, 
\ref{3.a21}:
\begin{equation}
[E\Nb^\mu]_{P.A.}=\tchib(E^\mu-\pi\Nb^\mu)
\label{10.820}
\end{equation}
Also, by \ref{3.a15}, \ref{3.a24}:
\begin{equation}
[T\Nb^\mu]_{P.A.}=(2E\lambdab+2\pi\lambdab\tchib)(E^\mu-\pi\Nb^\mu)
\label{10.821}
\end{equation}
Moreover, we have (see \ref{8.138}): 
\begin{equation}
[T\tchib]_{P.A.}=2E^2\lambdab+2\pi\lambdab E\tchib
\label{10.822}
\end{equation}
It then follows that for $m>1$:
\begin{equation}
[T^{m-1}\Nb^\mu]_{P.A.}=\left\{2\sum_{i=0}^{m-2}(2\pi\lambdab)^i E^{i+1}T^{m-2-i}\lambdab
+(2\pi\lambdab)^{m-1}E^{m-2}\tchib\right\}(E^\mu-\pi\Nb^\mu)
\label{10.823}
\end{equation}
Similarly for the $N$th approximants (see Section 10.3). We then conclude that the principal 
acoustical part of \ref{10.819} can be put into the form:
\begin{eqnarray}
&&-\pi\epb\lambdab^{-1}\cthb_{l+1} \ \ \ \mbox{: for $m=1$} \nonumber\\
&&-\pi\epb\lambdab^{-1}\left\{2\sum_{i=0}^{m-2}(2\pi\lambdab)^i\cnub_{m-2-i,l+2+i}
+(2\pi\lambdab)^{m-1}\cthb_{l+m}\right\} \ \mbox{: for $m>1$} \nonumber\\
&&\label{10.824}
\end{eqnarray}
(recall that $m+l=n$). The corresponding contribution to \ref{10.807} is then bounded by:
\begin{eqnarray}
&&C\tau_1\|\cthb_{l+1}\|_{L^2({\cal K}^{\tau_1})} \ \ \ \mbox{: for $m=1$} \nonumber\\
&&C\tau_1\left\{2\sum_{i=0}^{m-2}(\Pib\tau_1)^i\|\cnub_{m-2-i,l+2+i}\|_{L^2({\cal K}^{\tau_1})}
+(\Pib\tau_1)^{m-1}\|\cthb_{l+m}\|_{L^2({\cal K}^{\tau_1})}\right\} \ \mbox{: for $m>1$} \nonumber\\
&&\label{10.825}
\end{eqnarray}
Recall the definition of the quantity $\Pib$ just after \ref{10.741}. Here we may take $\Pib$ to be 
a fixed positive bound for $|2\pi\lambdab/\tau|$ on ${\cal K}^\delta$. Note that the case $m=1$ 
can be included in the case $m>1$ (in the case $m=1$ the sum vanishes by definition).

We summarize the results of this section in the following proposition. 

\vspace{2.5mm}

\noindent{\bf Proposition 10.9} \ \ The boundary values on ${\cal K}$ of the top order acoustical 
difference quantities $\cthb_l$, $l=n$, and $\cnub_{m-1,l+1}$, $m=1,...,n$, $l=n-m$, satisfy 
to principal terms the following estimates: 
\begin{eqnarray*}
&&\|\cthb_l\|_{L^2({\cal K}^{\tau_1})}\leq\tau_1^{a_0+b_0}\left\{
C\sqrt{\frac{\s^{(Y;0,l)}\cBb(\tau_1,\tau_1)}{2a_0+1}}\cdot \frac{1}{\tau_1^3}\right.\\
&&\hspace{25mm}\left.+C\sqrt{\max\{\s^{(0,l)}\cB(\tau_1,\tau_1),\s^{(0,l)}\cBb(\tau_1,\tau_1)\}}
\cdot\frac{1}{\tau_1^{5/2}}\right\} \\
&&\hspace{20mm}+\tau_1^{a_0+b_0}\left\{C\sqrt{\s^{(Y;0,l)}{\cal A}(\tau_1)}\cdot\frac{1}{\tau_1^3}
+C\sqrt{\s^{(0,l)}{\cal A}(\tau_1)}\cdot\frac{1}{\tau_1^{5/2}}\right\}\\
&&\hspace{10mm}+C\|\tau\Omega^{l+1}\chf\|_{L^2({\cal K}^{\tau_1})}
+C\|\tau\Omega^{l+1}\cv\|_{L^2({\cal K}^{\tau_1})}
+C\|\tau^2\Omega^{l+1}\cga\|_{L^2({\cal K}^{\tau_1})}
\end{eqnarray*}
\begin{eqnarray*}
&&\|\cnub_{m-1,l+1}\|_{L^2({\cal K}^{\tau_1})}\leq C\tau_1 \left\{(\Pib\tau_1)^{m-1}\|\cthb_{l+m}\|_{L^2({\cal K}^{\tau_1})}\right.\\
&&\hspace{35mm}\left.+2\sum_{i=0}^{m-2}(\Pib\tau_1)^i\|\cnub_{m-2-i,l+2+i}\|_{L^2({\cal K}^{\tau_1})}
\right\}\\
&&\hspace{15mm}+\tau_1^{a_m+b_m}\left\{
C\sqrt{\frac{\s^{(Y;m,l)}\cBb(\tau_1,\tau_1)}{2a_m+1}}\cdot \frac{1}{\tau_1^3}\right.\\
&&\hspace{30mm}\left.+C\sqrt{\max\{\s^{[m,l]}\cB(\tau_1,\tau_1),\s^{[m,l]}\cBb(\tau_1,\tau_1)\}}
\cdot\frac{1}{\tau_1^{5/2}}\right\} \\
&&\hspace{20mm}+\tau_1^{a_m+b_m}\left\{C\sqrt{\s^{(Y;m,l)}{\cal A}(\tau_1)}\cdot\frac{1}{\tau_1^3}
+C\sqrt{\s^{[m,l]}{\cal A}(\tau_1)}\cdot\frac{1}{\tau_1^{5/2}}\right\}\\
&&\hspace{20mm}+C\left(\|\tau^3\Omega^{l+2}T^{m-1}\chf\|_{L^2({\cal K}^{\tau_1})}
+2(m-1)\|\tau^2 \Omega^{l+2}T^{m-2}\chf\|_{L^2({\cal K}^{\tau_1})}\right.\\
&&\hspace{55mm}\left.+(m-1)(m-2)\|\tau\Omega^{l+2}T^{m-3}\chf\|_{L^2({\cal K}^{\tau_1})}\right)\\
&&\hspace{23mm}+C\left(\|\tau^2\Omega^{l+2}T^{m-1}\cv\|_{L^2({\cal K}^{\tau_1})}
+(m-1)\|\tau\Omega^{l+2}T^{m-2}\cv\|_{L^2({\cal K}^{\tau_1})}\right)\\
&&\hspace{20mm}+C\left(\|\tau^4\Omega^{l+2}T^{m-1}\cga\|_{L^2({\cal K}^{\tau_1})}
+3(m-1)\|\tau^3\Omega^{l+2}T^{m-2}\cga\|_{L^2({\cal K}^{\tau_1})}\right.\\
&&\hspace{52mm}+3(m-1)(m-2)\|\tau^2\Omega^{l+2}T^{m-3}\cga\|_{L^2({\cal K}^{\tau_1})}\\
&&\hspace{42mm}\left.+(m-1)(m-2)(m-3)\|\tau\Omega^{l+2}T^{m-4}\cga\|_{L^2({\cal K}^{\tau_1})}\right)
\end{eqnarray*}

\pagebreak

\chapter{Outline of the Top Order Acoustical Estimates for more than 2 Spatial Dimensions.}

In treating the shock development problem in more than 2 spatial dimensions, there is only one 
essential change in our approach, everything else carrying over in a straightforward manner. 
This change is in regard to the top order acoustical estimates treated in Chapter 10 in the case 
of two spatial dimensions. In the present chapter we shall discuss this essential change without 
going into the details of the estimates, which, after this presentation, can be carried out in a 
straightforward manner following the argument of Chapter 10. 

We recall that in Chapters 2 - 8 we addressed the problem in any number of spatial 
dimensions. In Chapter 9 we restricted ourselves to the case of 2 spatial dimensions and the treatment 
in Chapter 10 was in this framework. The restriction to 2 spatial dimensions will be restored in 
the chapters succeeding the present one, namely Chapters 12 - 14, with the exception of the last 
part of the last section of Chapter 14, where we shall discuss the obvious change in the derivation 
of pointwise estimates from $L^2$ estimates on the $S_{\ub,u}$ as a result of the different form of the 
Sobolev inequality in higher dimensions and how the corresponding Sobolev constant is deduced 
using the bootstrap assumptions. In Chapters 2 - 8, where we addressed the general problem, we 
denoted by $n$ the number of spatial dimensions. However in Chapters 9 and 10, where we focused on 
the case of 2 spatial dimensions, we denoted by $n=m+l$ the length of the commutation string, that is 
of $E^l T^m$, at the top order. To avoid confusion we denote from now on the number of spatial 
dimensions by $d$. 

We shall first derive the higher dimensional analogue of Proposition 10.1. 
The point is that in higher dimensions only the propagation equation for $\mbox{tr}\tchi$ which 
corresponds to the trace of the the 2nd variation equation for $\tchi$ of Proposition 3.4 
can be regularized. We denote:
\begin{equation}
\sk_{\flat}=\sk\cdot\sh
\label{11.1}
\end{equation}
This is a 2-covariant $S$ tensorfield given by \ref{3.76}:
\begin{equation}
\sk_{\flat}=\tchi+\frac{1}{2}\beta_N\sbeta\wedge\sd H-H\ss_N\otimes\sbeta
\label{11.2}
\end{equation}
and satisfying the 2nd variation equation of Proposition 3.4. Taking the trace of \ref{11.2} 
we obtain:
\begin{equation}
\mbox{tr}\sk_{\flat}=\mbox{tr}\tchi-H(\sbeta,\ss_N)
\label{11.3}
\end{equation}
where for arbitrary p-covariant $S$ tensorfields $\alpha$, $\beta$ we denote by $(\alpha,\beta)$ 
their inner product relative to $\sh$:
\begin{equation}
(\alpha,\beta)=(\sh^{-1})^{A_1 B_1} \cdot \cdot \cdot (\sh^{-1})^{A_p B_p}
\alpha_{A_1 . . . A_p}\beta_{B_1 . . . B_p}
\label{11.4}
\end{equation}
Note that $\mbox{tr}\sk_{\flat}=\mbox{tr}\sk$ is the trace of $\sk$ as a linear transformation 
in $T_q S_{\ub,u}$ at each point $q$, while $\mbox{tr}\tchi=\mbox{tr}\tchi^\sharp$ is the 
trace of $\tchi^\sharp$ as such a linear transformation. The 1st variation equations of Proposition 
3.1 imply:
\begin{equation}
L(\mbox{tr}\sk_{\flat})=\mbox{tr}(\sL_L\sk_{\flat})-2(\chi,\sk_{\flat})
\label{11.5}
\end{equation}
As noted earlier, the principal term on the right in the 2nd variation equation for $\sk_{\flat}$ 
is the order 2 term $\sD\sm_{\flat}$ where:
\begin{equation}
\sm_{\flat}=-\beta_N(LH)\sbeta+\frac{1}{2}\rho\beta_N^2\sd H-H s_{NL}\sbeta
\label{11.6}
\end{equation}
by \ref{3.68}. Therefore taking the trace the principal term is $\sdiv\sm$. As in the proof of 
Proposition 10.1, the regularization of this propagation equation for $\mbox{tr}\tchi$ is by finding 
a 1st order quantity $v$ such that the difference
\begin{equation}
\sdiv\sm-Lv
\label{11.7}
\end{equation}
is a quantity of 1st order. Now, from \ref{11.6} we deduce:
\begin{equation}
\left[\sdiv\sm\right]_{P.P.}=\left[L(\sbeta,-\beta_N\sd H-H\ss_N)
+\frac{1}{2}\rho\beta_N^2\slap H\right]_{P.P.}
\label{11.8}
\end{equation}
To express the last term on the right as $L$ applied to a 1st order quantity up to a 1st order remainder 
we appeal to Lemma 10.1, which in $d$ spatial dimensions takes the form:

\vspace{2.5mm}

\noindent{\bf Lemma 11.1} \ \ We have: 
$$\square_h H=\left(H^{\prime\prime}-((d-1)/2)\Omega^{-1}\Omega^\prime H^\prime\right)h^{-1}(d\sigma,d\sigma)
-2H^\prime(g^{-1})^{\mu\nu}h^{-1}(d\beta_\mu,d\beta_\nu)$$
where $H$ and $\Omega$ being a functions of $\sigma$ we denote
$$H^\prime=\frac{dH}{d\sigma}, \ H^{\prime\prime}=\frac{d^2 H}{d\sigma^2} \ \ \mbox{and} \ \ 
\Omega^\prime=\frac{d\Omega}{d\sigma}$$

\vspace{2.5mm}

The proof is the same, the conformal transformation formula \ref{9.229} applying with 
$d$ in the role of $n$. 

Since for any function $f$ on ${\cal N}$, 
$$\square_h f=\mbox{tr}(h^{-1}\cdot Ddf),$$
substituting the expansion \ref{2.38} of $h^{-1}$ we obtain:
\begin{equation}
\square_h f=-a^{-1}(Ddf)(L,\Lb)+(\sh^{-1})^{AB}(Ddf)(\Omega_A,\Omega_B)
\label{11.9}
\end{equation}
Since for any pair of vectorfields $X$, $Y$ on ${\cal N}$,
$$(Ddf)(X,Y)=X(Yf)-(D_X Y)f,$$
taking $X=L$, $Y=\Lb$ we obtain, by \ref{3.23}:
\begin{equation}
(Ddf)(L,\Lb)=L(\Lb f)-2(\etab,\sd f)
\label{11.10}
\end{equation}
Also, taking $X=\Omega_A$, $Y=\Omega_B$, these being $S$ vectorfields \ref{3.26} applies and we obtain:
\begin{equation}
(Ddf)(\Omega_A,\Omega_B)=(\sD\sd f)(\Omega_A,\Omega_B)-\frac{1}{2a}(\chib_{AB}Lf+\chi_{AB}\Lb f)
\label{11.11}
\end{equation}
and, taking the trace:
\begin{equation}
(\sh^{-1})^{AB}(Ddf)(\Omega_A,\Omega_B)-\frac{1}{2a}(\mbox{tr}\chib Lf+\mbox{tr}\chi\Lb f)
\label{11.12}
\end{equation}
Substituting then \ref{11.10} and \ref{11.12} in \ref{11.9} yields the formula:
\begin{equation}
a\slap f-L(\Lb f)=a\square_h f+\frac{1}{2}(\mbox{tr}\chib Lf+\mbox{tr}\chi\Lb f)
-2(\etab,\sd f)
\label{11.13}
\end{equation}

The formula \ref{11.13} applies in particular taking $f=H$ and in this case shows that by virtue of 
Lemma 11.1 that 
$$a\slap H-L(\Lb H)$$
is actually of order 1. Let us in fact define the quantity:
\begin{equation}
M=\frac{1}{2}\beta_N^2(a\slap H-L(\Lb H))
\label{11.14}
\end{equation}
so this is actually a 1st order quantity. Going then back to \ref{11.8} the last term on the right 
is 
$$\left[\lambda^{-1} M+\frac{1}{2\lambda}\beta_N^2 L(\Lb H)\right]_{P.P.}=
\left[L\left(\frac{1}{2\lambda}\beta_N^2\Lb H\right)\right]_{P.P.}$$
as $\lambda^{-1} M$ does not contribute to the principal part. Setting then:
\begin{equation}
v=\frac{1}{2\lambda}\beta_N^2\Lb H-(\sbeta,\beta_N\sd H+H\ss_N)
\label{11.15}
\end{equation}
the difference \ref{11.7} is indeed a quantity of 1st order. By \ref{11.3}, 
\begin{equation}
\lambda(\mbox{tr}\sk_{\flat}-v)=\lambda\mbox{tr}\tchi+f:=\theta
\label{11.16}
\end{equation}
where
\begin{equation}
f:=-\frac{1}{2}\beta_N^2\Lb H+\lambda\beta_N(\sbeta,\sd H)
\label{11.17}
\end{equation}
(compare with \ref{10.13}, \ref{10.14}). The 1st order quantity $\theta$ then satisfies a propagation 
equation of the form:
\begin{equation}
L\theta=R
\label{11.18}
\end{equation}
where $R$ is again a 1st order quantity (compare with \ref{10.15}). 

We shall now discern the principal acoustical part of $R$. We have, by \ref{11.5}, 
\begin{equation}
R=\lambda\mbox{tr}(\sL_L\sk_{\flat})-L(\lambda v)-2\lambda(\chi,\sk_{\flat})
+(L\lambda)\mbox{tr}\sk_{\flat}
\label{11.19}
\end{equation}
We first discern the principal acoustical part of the 1st order quantity 
$$\lambda\sdiv\sm-L(\lambda v)$$
(see \ref{11.7}). Taking the divergence of \ref{11.6} we obtain:
\begin{equation}
\sdiv\sm=-\beta_N(\sbeta,\sd LH)-H(\sbeta, N^\mu \sd L\beta_\mu)+\frac{1}{2}\rho\beta_N^2\slap H+R_1
\label{11.20}
\end{equation}
where $R_1$ is a quantity of order 1 which up to a remainder with vanishing principal acoustical part 
is given by:
\begin{eqnarray}
&&R_1=-\left[\beta_N\sdiv\sbeta+(\sbeta,\sd\beta_N)\right]LH\nonumber\\
&&\hspace{8mm}+\frac{1}{2}\beta_N^2(\sd\rho,\sd H)+\rho\beta_N(\sd\beta_N,\sd H)\nonumber\\
&&\hspace{8mm}-Hs_{NL}\sdiv\sbeta-H(\sbeta,\sd N^\mu)L\beta_\mu 
\label{11.21}
\end{eqnarray}
On the other hand, multiplying \ref{1.15} by $\lambda$ and applying $L$ we obtain:
\begin{equation}
L(\lambda v)=\frac{1}{2}\beta_N^2 L(\Lb H)
-\lambda\left[\beta_N(\sbeta,\sd LH)+H(\sbeta,N^\mu\sd L\beta_\mu)\right]+\lambda R_2
\label{11.22}
\end{equation}
where $R_2$ is a quantity of order 1 which up to a remainder with vanishing principal acoustical part 
is given by:
\begin{equation}
R_2=2\rho\tchi^{\sharp\sharp}(\sbeta,\beta_N\sd H+H\ss_N)-(\sL_L\sbeta,\beta_N\sd H+H\ss_N)
\label{11.23}
\end{equation}
Subtracting \ref{11.22} from $\lambda$ times \ref{11.20}, yields, 
in view of the definition \ref{11.14},
\begin{equation}
\lambda\sdiv\sm-L(\lambda v)=M+\lambda(R_1-R_2)
\label{11.24}
\end{equation}
Now, by \ref{11.13} with $H$ in the role of $f$:
\begin{equation}
M=\frac{1}{2}\beta_N^2\left[a\square_h H+\frac{1}{2}(\mbox{tr}\chib LH+\mbox{tr}\chi\Lb H)
-2(\etab,\sd H)\right]
\label{11.25}
\end{equation}
and by Lemma 11.1 $\square_h H$ is an order 1 quantity with vanishing order 1 acoustical part. 
In view of the fact that by Proposition 3.2
\begin{equation}
[\etab]_{P.A.}=\rhob(\sd\lambdab+\lambdab\tchib\cdot\pi)
\label{11.26}
\end{equation}
we obtain:
\begin{equation}
[M]_{P.A.}=\frac{1}{4}(\rhob\mbox{tr}\tchib LH+\rho\mbox{tr}\tchi\Lb H)
-\rhob\beta_N^2(\sd\lambdab,\sd H)-a\beta_N^2\tchib^{\sharp\sharp}(\pi,\sd H)
\label{11.27}
\end{equation}
To obtain the principal acoustical parts of $R_1$ and $R_2$ we use the following formulas, which 
are readily deduced from the results of Chapter 3:
\begin{eqnarray}
&&[\sD\sbeta]_{P.A.}=\frac{1}{2c}(\beta_N\tchib+\beta_{\Nb}\tchi)\nonumber\\
&&[\sd\beta_N]_{P.A.}=\tchi^{\sharp}\cdot(\sbeta-\beta_N\pi)\nonumber\\
&&[\sL_L\sbeta]_{P.A.}=\beta_N\sd\rho+\rho\tchi^{\sharp}\cdot(\sbeta-\beta_N\pi)
\label{11.28}
\end{eqnarray}
In view of \ref{11.19}, \ref{11.24}, to arrive at the formula for $[R]_{P.A.}$ we must add to the 
above $\lambda$ times the trace of the remaining terms in the 2nd variation equation of 
Proposition 3.4 after the removal of the principal term $\sD\sm$, as well as the last two terms in 
\ref{11.19}. We arrive in this way at the following extension of Proposition 10.1 to higher 
dimensions. 

\vspace{2.5mm}

\noindent{\bf Proposition 11.1} \ \ The 1st order quantities: 
$$\theta=\lambda\mbox{tr}\tchi+f, \ \ \ \thetab=\lambdab\mbox{tr}\tchib+\fb$$
where 
$$f=-\frac{1}{2}\beta_N^2\Lb H+\lambda\beta_N(\sbeta,\sd H), \ \ \ 
\fb=-\frac{1}{2}\beta_{\Nb}^2 LH+\lambdab\beta_{\Nb}(\sbeta,\sd H)$$
satisfy the following regularized propagation equations:
$$L\theta=R, \ \ \ \Lb\thetab=\Rb$$
where $R$, $\Rb$ are again quantities of order 1, their principal acoustical parts being given by
\begin{eqnarray*}
&&[R]_{P.A.}=-a|\tchi|^2+2(L\lambda)\mbox{tr}\tchi-\lambda(LH)\tchi^{\sharp\sharp}(\sbeta,\sbeta)\\
&&[\Rb]_{P.A.}=-a|\tchib|^2+2(\Lb\lambdab)\mbox{tr}\tchib-\lambdab(\Lb H)\tchib^{\sharp\sharp}
(\sbeta,\sbeta)
\end{eqnarray*}

\vspace{2.5mm}

We note in connection with the boundary condition for $\thetab$ on ${\cal K}$, that it follows 
from the boundary condition for $\mbox{tr}\tchib$, which, as in the case $d=2$ (see \ref{10.754} - 
\ref{10.756}), in turn follows from the boundary condition \ref{4.1}. The last reads:
\begin{equation}
\triangle\beta_\mu \ \sd x^\mu=0
\label{11.b1}
\end{equation}
Applying $\sD$ to this gives:
\begin{equation}
\triangle\beta_\mu \ \sD^2 x^\mu+\sd\triangle\beta_\mu\otimes\sd x^\mu=0
\label{11.b2}
\end{equation}
Now $\sD^2 x^\mu$ is expressed in terms of $\tchi$, $\tchib$ by \ref{3.105}:
\begin{equation}
\sD^2 x^\mu=\frac{1}{2c}(\tchib N^\mu+\tchi\Nb^\mu)-\tgamma^\mu 
\label{11.b3}
\end{equation}
where by \ref{3.101}:
\begin{equation}
\tgamma^\mu=\sgamma\cdot\sd x^\mu+\tgamma^N N^\mu+\tgamma^{\Nb}\Nb^\mu 
\label{11.b4}
\end{equation}
In view of the boundary condition \ref{11.b1}, the 1st term in \ref{11.b4} does not contribute to the 
1st term in \ref{11.b2}. The other two terms in \ref{11.b4} are given by \ref{3.99}:
\begin{equation}
\tgamma^N=-\frac{1}{2c}\beta_{\Nb}\salpha, \ \ \tgamma^{\Nb}=-\frac{1}{2c}\beta_N\salpha \ \ 
\mbox{where} \ \ \salpha=\frac{1}{2}(\sbeta\otimes\sd H+\sd H\otimes\sbeta+2H\sss)
\label{11.b5}
\end{equation}
Substituting the above in \ref{11.b3} we obtain the boundary condition for $\tchib$:
\begin{equation}
r(\tchib+\beta_{\Nb}\salpha)=\tchi+\beta_N\salpha+\frac{2c}{\epb}\sd\triangle\beta_\mu\otimes\sd x^\mu 
\label{11.b6}
\end{equation}
(note that the last term is symmetric).

Next we shall derive the higher dimensional analogue of Proposition 10.2. 
The point again is that in higher dimensions only the propagation equation for $\slap\lambda$ which is 
derived from the propagation equation for $\lambda$ of Proposition 3.3 
can be regularized. To derive the propagation equation for $\slap\lambda$, we must derive a formula 
for the commutator $[L,\slap]$ applied to an arbitrary function $f$ on ${\cal N}$. The required 
formula is most readily derived using the representation, described in Section 2.1 
(see \ref{2.18} - \ref{2.21}), which is adapted to the flow of $L$. This representation was 
also used in Section 10.5 and 10.6 (see \ref{10.339} - \ref{10.345}) in the case $d=2$ the 
corresponding coordinate on $S^1$ being denoted there by $\vartheta^\prime$. Here the argument 
is intrinsic to a given $C_u$ consequently there is no confusion if, in the context of the 
representation adapted to the flow of $L$ on the given $C_u$, we denote by 
$(\vartheta^A : A=1,...,d-1)$ the corresponding local coordinates on $S^{d-1}$. 
We then have:
\begin{equation}
L=\frac{\partial}{\partial\ub}, \ \ \sh=\sh_{AB}\sd\vartheta^A\otimes\sd\vartheta^B, 
\ \ \sL_L\sh=\frac{\partial\sh_{AB}}{\partial\ub}\sd\vartheta^A\otimes\sd\vartheta^B
\label{11.29}
\end{equation}
hence according to the 1st variation equations of Proposition 3.1 :
\begin{equation}
\frac{\partial\sh_{AB}}{\partial\ub}=2\chi_{AB}
\label{11.30}
\end{equation}
and, if $\sGamma^C_{AB}$ are the connection coefficients on $(S_{\ub,u},\sh)$ in these coordinates, 
\begin{equation}
\frac{\partial\sGamma^C_{AB}}{\partial\ub}=\sD_A\chi_B^C+\sD_B\chi_A^C-\sD^C\chi_{AB}
\label{11.31}
\end{equation}
The right hand side of \ref{11.31} represents a type $S^1_2$ $S$ tensorfield, that is, at each 
$S_{\ub,u}$ a section of the tensor bundle
$$\bigcup_{q\in S_{\ub,u}} S^1_2(T_q S_{\ub,u})$$
over $S_{\ub,u}$, where $S^1_2(T_q S_{\ub,u})$ denotes the space of symmetric bilinear forms on 
$T_q S_{\ub,u}$ with values in $T_q S_{\ub,u}$. Since in the above local coordinates we have, 
for an arbitrary function $f$ on ${\cal N}$, 
\begin{equation}
\slap f=(\sh^{-1})^{AB}\left(\frac{\partial^2 f}{\partial\vartheta^A\partial\vartheta^B}
-\sGamma^C_{AB}\frac{\partial f}{\partial\vartheta^C}\right)
\label{11.32}
\end{equation}
we then obtain:
\begin{equation}
\frac{\partial}{\partial\ub}\slap f=-2\chi^{AB}(\sD^2 f)_{AB}-X^A(\sd f)_A+\slap\frac{\partial f}{\partial\ub}
\label{11.33}
\end{equation}
where:
\begin{equation}
X^A=2\sD^B\chi_B^A-(\sd\mbox{tr}\chi)^A
\label{11.34}
\end{equation}
The above and a similar argument for the conjugate quantities establishes the following lemma. 

\vspace{2.5mm}

\noindent{\bf Lemma 11.2 :} \ \ For an arbitrary function $f$ on ${\cal N}$ we have:
\begin{eqnarray*}
&&[L,\slap]f=-2\chi^{\sharp\sharp}\cdot\sD^2 f-Xf\\
&&[\Lb,\slap]f=-2\chib^{\sharp\sharp}\cdot\sD^2 f-\underline{X}f
\end{eqnarray*}
where $X$, $\underline{X}$ are the $S$ vectorfields given by:
\begin{eqnarray*}
&&X^A=2\sD^B\chi_B^A-(\sd\mbox{tr}\chi)^A\\
&&\underline{X}^A=2\sD^B\chib_B^A-(\sd\mbox{tr}\chib)^A
\end{eqnarray*}

\vspace{2.5mm}

Note that the $S$ vectorfields $X$, $\underline{X}$ are quantities of order 2 and that by the 
Codazzi equations of Proposition 3.6 their principal acoustical parts are given by (see \ref{8.127}):
\begin{equation}
[X]_{P.A}=\rho(\sd\mbox{tr}\tchi)^\sharp, \ \ [\underline{X}]_{P.A.}=\rhob(\sd\mbox{tr}\tchib)^\sharp
\label{11.35}
\end{equation}
Applying $\slap$ to the propagation equation for $\lambda$ of Proposition 3.3 and using Lemma 11.2 
we deduce the following propagation equation for $\slap\lambda$ :
\begin{eqnarray}
&&L(\slap\lambda)=-2\chi^{\sharp\sharp}\cdot\sD^2\lambda-X\lambda+p\slap\lambda+q\slap\lambdab
\nonumber\\
&&\hspace{15mm}+\lambda\slap p+\lambdab\slap q+2(\sd p,\sd\lambda)+2(\sd q,\sd\lambdab) 
\label{11.36}
\end{eqnarray}
The principal part of the right hand side is contained in the order 3 terms:
$$\lambda\slap p+\lambdab\slap q$$
As in the proof of Proposition 10.2, the regularization of this propagation equation for 
$\slap\lambda$ is by finding a 2nd order quantity $w_1$ such that the difference
\begin{equation}
\hat{P}:=\lambda\slap p-Lw_1
\label{11.37}
\end{equation}
is a quantity of 2nd order, and by finding a 2nd order quantity $w_2$ such that the difference
\begin{equation}
Q:=\lambda\lambdab\slap q-Lw_2
\label{11.38}
\end{equation}
Setting then, as in \ref{10.59}, 
\begin{equation}
j=\lambda w_1+w_2
\label{11.39}
\end{equation}
a quantity of order 2, the difference
\begin{equation}
\lambda(\lambda\slap p+\lambdab\slap q)-Lj
\label{11.40}
\end{equation}
is a quantity of 2nd order, hence multiplying \ref{11.36} by $\lambda$ and subtracting $Lj$ from 
both sides we obtain a regularized propagation equation for the quantity
\begin{equation}
\nu=\lambda\slap\lambda-j
\label{11.41}
\end{equation}
In fact $w_1$ and $w_2$ will have vanishing principal acoustical parts, hence so will $j$. 

Now, $p$ is given in accordance with Proposition 3.3 by \ref{10.46}, $m$ being given in terms of 
$\sm$ and $\om$ by the 1st of \ref{3.47}:
\begin{equation}
m=-\pi\cdot\sm-\om 
\label{11.42}
\end{equation}
Recall that in higher dimensions $\pi=\sd t$ is a $S$ 1-form while $\sm$ is a $S$ vectorfield, 
however $\om$, a scalar, is still given in accordance with \ref{3.57} by \ref{10.48}. 
The higher dimensional analogue of the formula \ref{10.49} for $p$ is then:
\begin{equation}
p=-\pi\cdot\sm-\frac{1}{4c}\beta_N(\beta_N+2\beta_{\Nb})LH-\frac{1}{2c}H(\beta_N+\beta_{\Nb})s_{NL}
\label{11.43}
\end{equation}
This formula gives:
\begin{eqnarray}
&&\left[\slap p\right]_{P.P.}=-\left[\pi^\sharp\cdot\slap\sm_\flat\right]_{P.P.} \label{11.44}\\
&&\hspace{13mm}-\left[ L\left(\frac{1}{4c}\beta_N(\beta_N+2\beta_{\Nb})\slap H
+\frac{1}{2c}H(\beta_N+\beta_{\Nb})N^\mu\slap\beta_\mu\right)\right]_{P.P.}\nonumber
\end{eqnarray}
In deducing a formula for $\slap\sm_\flat$ from \ref{11.6}, to assess the contribution of the 2nd 
term  we use the fact that for any function $f$ on ${\cal N}$ we have:
\begin{equation}
\slap(\sd f)-\sd(\slap f)=\sRic\cdot(\sd f)^\sharp
\label{11.45}
\end{equation}
a general fact in Riemannian geometry, realized for the Riemannian manifold $(S_{\ub,u},\sh)$. 
Here, $\sRic$ is the Ricci curvature of $(S_{\ub,u},\sh)$, which by the Gauss equation of 
Proposition 3.7 is a quantity of order 2 with vanishing principal acoustical part. We apply \ref{11.45} 
with $H$ in the role of $f$. This yields:
\begin{equation}
\slap\sm_\flat=-\left[\sL_L\left(\sbeta(\beta_N\slap H+HN^\mu\slap\beta_\mu)\right)\right]_{P.P.}
+\left[\frac{1}{2}\rho\beta_N^2\sd\slap H\right]_{P.P.}
\label{11.46}
\end{equation}
Multiplying then by $\lambda$ we substitute 
\begin{equation}
\frac{1}{2}\beta_N^2 a\sd\slap H=\frac{1}{2}\beta_N^2\sd L(\Lb H)+\sd M
+(L\Lb H)\sd\left(\frac{1}{2}\beta_N^2\right)
-(\slap H)\sd\left(\frac{1}{2}\beta_N^2 a\right)
\label{11.47}
\end{equation}
from \ref{11.14}. Here the last two terms on the right, like the 2nd term, are of order 2. We then 
conclude that:
\begin{equation}
\left[\lambda\slap\sm_\flat-\sL_L\sxi\right]_{P.P.}=0
\label{11.48}
\end{equation}
where $\xi$ is the $S$ 1-form
\begin{equation}
\xi=\frac{1}{2}\beta_N^2\sd(\Lb H)-\lambda\sbeta\left(\beta_N\slap H+H N^\mu\slap\beta_\mu\right)
\label{11.49}
\end{equation}
Therefore by \ref{11.44} setting:
\begin{equation}
w_1=-\pi^\sharp\cdot\xi-\lambda\left(\frac{1}{4c}\beta_N(\beta_N+2\beta_{\Nb})\slap H
+\frac{1}{2c}H(\beta_N+\beta_{\Nb})N^\mu\slap\beta_\mu\right)
\label{11.50}
\end{equation}
that is:
\begin{eqnarray}
&&w_1=-\frac{1}{2}\beta_N^2\pi^\sharp\cdot\sd H \nonumber\\
&&\hspace{8mm}+\lambda\beta_N\left[(\pi,\sbeta)-\frac{1}{4c}(\beta_N+2\beta_{\Nb})\right]\slap H
\nonumber\\
&&\hspace{8mm}+\lambda\left[(\pi,\sbeta)-\frac{1}{2c}(\beta_N+\beta_{\Nb})\right]H N^\mu\slap\beta_\mu 
\label{11.51}
\end{eqnarray}
the difference \ref{11.37} is indeed a quantity of order 2. 

Consider next $q$, given by \ref{10.54}: 
\begin{equation}
q=\frac{1}{4c}\beta_N^2\Lb H
\label{11.52}
\end{equation}
Applying $\slap$ and multiplying by $\lambda\lambdab=ca$ the principal part is, in view of Lemma 11.2, 
\begin{equation}
[\lambda\lambdab\slap q]_{P.P.}=\left[\frac{1}{4}\beta_N^2 a\Lb\slap H\right]_{P.P.}
=\frac{1}{2}\left[\Lb\left(\frac{1}{2}\beta_N^2 a\slap H\right)\right]_{P.P.}
\label{11.53}
\end{equation}
Here we substitute
\begin{equation}
\Lb\left(\frac{1}{2}\beta_N^2 a\slap H\right)=\frac{1}{2}\beta_N^2 L(\Lb^2 H)+\frac{1}{2}\Lb M
+\frac{1}{2}\beta_N^2 Z\Lb H+(L\Lb H)\sd\left(\frac{1}{2}\beta_N^2\right)
\label{11.54}
\end{equation}
from \ref{11.14}. Here the last two terms on the right, like the 2nd term, are of order 2. Setting 
then: 
\begin{equation}
w_2=\frac{1}{4}\beta_N^2\Lb^2 H
\label{11.55}
\end{equation}
the difference \ref{11.38} is indeed a quantity of order 2. Then with  $j$ defined according 
to \ref{11.39}, thus by \ref{11.51} and \ref{11.55} being a 2nd order quantity with vanishing 
principal acoustical part, the 2nd order quantity $\nu$ defined by \ref{11.41} satisfies 
a propagation equation of the form:
\begin{equation}
L\nu=K
\label{11.56}
\end{equation}
where $K$ is again a 2nd order quantity. 

We shall now discern the principal acoustical part of $K$. From \ref{11.36} - \ref{11.38} we have:
\begin{eqnarray}
&&K=-2\lambda\chi^{\sharp\sharp}\cdot\sD^2\lambda+(L\lambda+p\lambda)\slap\lambda
-\lambda X\lambda+\lambda q\slap\lambdab \nonumber\\
&&\hspace{8mm}+\lambda\hat{P}+Q-(L\lambda)w_1+2\lambda(\sd p,\sd\lambda)+2\lambda(\sd q,\sd\lambdab)
\label{11.57}
\end{eqnarray}
The principal acoustical part of the first four terms is obvious, while the last three terms have 
vanishing principal acoustical part. Thus, what remains to be determined is the principal acoustical 
part of $\hat{P}$ and of $Q$. 

Revisiting the argument showing that with $w_1$ defined by 
\ref{11.51} $\hat{P}$ is actually a quantity of order 2, we find:
\begin{equation}
\hat{P}=-\pi^\sharp\cdot\sd M+\lambda\left[-\sm\cdot\slap\pi
-\frac{(\pi,\sd H)}{2}\slap\left(\frac{\beta_N^2\lambdab}{c}\right)+S_1-S_2\right]
\label{11.58}
\end{equation}
where $S_1$ and $S_2$ are quantities of order 2 which up to a remainder with vanishing principal 
acoustical part are given by:
\begin{eqnarray}
&&S_1=\frac{\beta_N}{4c}(\beta_N+2\beta_{\Nb})[L,\slap]H
+\frac{H}{2c}(\beta_N+\beta_{\Nb})N^\mu[L,\slap]\beta_\mu\nonumber\\
&&\hspace{8mm}-\frac{(LH)}{4}\slap\left(\frac{\beta_N(\beta_N+2\beta_{\Nb})}{c}\right)
-\frac{H(L\beta_\mu)}{2}\slap\left(\frac{(\beta_N+\beta_{\Nb})N^\mu}{c}\right)\nonumber\\
&&\label{11.59}
\end{eqnarray}
\begin{eqnarray}
&&S_2=(\pi,\sbeta)\left(\beta_N[L,\slap]H+H N^\mu [L,\slap]\beta_\mu\right)\nonumber\\
&&\hspace{8mm}-\pi^\sharp\cdot\left((LH)\slap(\beta_N\sbeta)+H(L\beta_\mu)\slap(\sbeta N^\mu)\right)
\label{11.60}
\end{eqnarray}
To obtain the principal acoustical part of $\slap\pi$, which appears in \ref{11.58}, 
we recall from Chapter 3 that $\pi=\sd x^0$ and by \ref{3.105}:
\begin{equation}
\slap x^0=\frac{1}{2c}(\mbox{tr}\tchib+\mbox{tr}\tchi)-\mbox{tr}\tilde{\gamma}^0 
\label{11.61}
\end{equation}
where by \ref{3.101}, in view of \ref{3.98}, \ref{3.100}, $\tilde{\gamma}^{\mu}_{AB}$ 
is of order 1 with vanishing principal acoustical part, hence:
\begin{equation}
[\slap x^0]_{P.A.}=\frac{1}{2c}(\mbox{tr}\tchib+\mbox{tr}\tchi)
\label{11.62}
\end{equation}
The formula \ref{11.45} then implies that:
\begin{equation}
\left[\slap\pi\right]_{P.A.}=\frac{1}{2c}(\sd\mbox{tr}\tchib+\sd\mbox{tr}\tchi)
\label{11.63}
\end{equation}
To obtain the principal acoustical parts of $S_1$ and $S_2$ we use, besides Lemma 11.2, 
the following formulas, which are readily deduced from the results of Chapter 3:
\begin{eqnarray}
&&\left[\slap N^\mu\right]_{P.A.}=(\sd\mbox{tr}\tchi,\sd x^\mu-N^\mu\pi)\nonumber\\
&&\left[c^{-1}\slap c\right]_{P.A.}=-(\pi,\sd\mbox{tr}\tchi+\sd\mbox{tr}\tchib)\nonumber\\
&&\left[\slap\sbeta\right]_{P.A.}=\frac{1}{2c}
(\beta_N\sd\mbox{tr}\tchib+\beta_{\Nb}\sd\mbox{tr}\tchi)\nonumber\\
&&\left[\slap\beta_N\right]_{P.A.}=(\sbeta-\beta_N\pi,\sd\mbox{tr}\tchi)\nonumber\\
&&\left[\slap\sbeta_{\Nb}\right]_{P.A.}=(\sbeta-\beta_{\Nb}\pi,\sd\mbox{tr}\tchib)
\label{11.64}
\end{eqnarray}
The principal acoustical part of $\sd M$, which appears in \ref{11.58}, is obtained directly from 
\ref{11.27}. We arrive in this way at the conclusion that $[\hat{P}]_{P.A.}$ is of the form:
\begin{eqnarray}
&&\left[\hat{P}\right]_{P.A.}=\oB_1\slap\lambdab+\oB_2\cdot\sD^2\lambdab\nonumber\\
&&\hspace{15mm}+\oA_{1,1}\cdot\sd\mbox{tr}\tchib+\oA_{1,2}\cdot\sD\tchib+A_1\cdot\sd\mbox{tr}\tchi
\label{11.65}
\end{eqnarray}
Here:
\begin{eqnarray}
&&\oB_1=-\frac{1}{2}\rhob\beta_N^2(\pi,\sd H)\nonumber\\
&&\oB_2=\frac{1}{2}\rhob\beta_N^2
\left(\pi^\sharp\otimes(\sd H)^\sharp+(\sd H)^\sharp\otimes\pi^\sharp\right) 
\label{11.66}
\end{eqnarray}
and:
\begin{eqnarray}
&&\oA_{1,1}=-\frac{1}{2}a\beta_N^2\left(\frac{1}{2c}(\sd H)^\sharp+(\pi,\sd H)\pi^\sharp\right)
\nonumber\\
&&\oA_{1,2}=\frac{1}{2}a\beta_N^2\pi^\sharp\otimes\left(\pi^\sharp\otimes(\sd H)^\sharp
+(\sd H)^\sharp\otimes\pi^\sharp\right)
\label{11.67}
\end{eqnarray}
(note that $\oA_{1,2}$ and $\sD\tchib$ are sections of dual tensor bundles over the $S_{\ub,u}$).  
Moreover, the coefficient $A_1$ is a $S$ vectorfield of the form:
\begin{equation}
A_1=-\frac{\rho}{4}\beta_N^2(\Lb H)\pi^\sharp+aA_{1,1}+\lambda A_{1,2}^\mu L\beta_\mu 
\label{11.68}
\end{equation}
where $A_{1,1}$ and $A_{1,2}$ are of order 1 with vanishing order 1 acoustical part.  

Revisiting the argument showing that with $w_2$ defined by 
\ref{11.55} $Q$ is actually a quantity of order 2, we find:
\begin{equation}
Q=\frac{1}{2}\Lb M+S_3
\label{11.69}
\end{equation}
where $S_3$ is a  quantity of order 2 which up to a remainder with vanishing principal 
acoustical part is given by:
\begin{equation}
S_3=\frac{ac}{4}(\Lb H)\slap\left(\frac{\beta_N^2}{c}\right)-\frac{a}{4}\beta_N^2[\Lb,\slap]H
\label{11.70}
\end{equation}
The principal acoustical part of $\Lb M$ is obtained directly from \ref{11.27}. To obtain the 
principal acoustical part of $S_3$ we use Lemma 11.2 and the formulas \ref{11.64}. We arrive in 
this way at the conclusion that $[Q]_{P.A.}$ is of the form:
\begin{equation}
\left[Q\right]_{P.A.}=q\lambdab\slap\lambda+A_2\cdot\sd\mbox{tr}\tchi+\oA_2\cdot\sd\mbox{tr}\tchib
\label{11.71}
\end{equation}
where:
\begin{eqnarray}
&&A_2=\frac{a}{2}\beta_N(\Lb H)\sbeta^\sharp\nonumber\\
&&\oA_2=\frac{a}{4}\beta_N^2\left((\Lb H)\pi^\sharp+\rhob(\sd H)^\sharp\right)
\label{11.72}
\end{eqnarray}

Taking into account the above results \ref{11.65} - \ref{11.71} we arrive through \ref{11.57} at 
the following expression for $[K]_{P.A.}$ which extends \ref{10.98} - \ref{10.100} to higher 
dimensions:
\begin{eqnarray}
&&\left[K\right]_{P.A.}=2(L\lambda)\slap\lambda-2\lambda\chi^{\sharp\sharp}\cdot\sD^2\lambda
+\left(\lambda A-a(\sd\lambda)^\sharp\right)\cdot\sd\mbox{tr}\tchi\nonumber\\
&&\hspace{15mm}+\lambda q\slap\lambdab+\frac{\lambda^2\beta_N^2}{2c}
\left\{2(\sD^2\lambdab)(\pi^\sharp,(\sd H)^\sharp)-(\pi,\sd H)\slap\lambdab\right\}
\nonumber\\
&&\hspace{15mm}+\lambda\lambdab q(\pi,\sd\mbox{tr}\tchib)
-\frac{1}{2}a\lambda\beta_N^2(\pi,\sd H)(\pi,\sd\mbox{tr}\tchib)\nonumber\\
&&\hspace{15mm}+a\lambda\beta_N^2(\sD_{\pi^\sharp}\tchib)(\pi^\sharp,(\sd H)^\sharp)
\label{11.73}
\end{eqnarray}
where:
\begin{eqnarray}
&&A=A_1+\lambda^{-1}A_2\nonumber\\
&&\hspace{4mm}=\frac{\lambdab}{2c}\beta_N(\Lb H)\left(\sbeta^\sharp-\frac{1}{2}\beta_N\pi^\sharp\right)
+a A_{1,1}+\lambda A_{1,2}^\mu L\beta_\mu \label{11.74}
\end{eqnarray}

However, the estimates involving $\nu$ in Chapter 10 where not deduced using the form \ref{10.98} 
of $[K]_{P.A.}$, which appears in the statement of Proposition 10.2, but rather the form 
\ref{10.278} (the function $F_2$ is of order 2 but with vanishing P.A. part) as the sum of $L\tau$ 
where 
\begin{equation}
\tau=\lambda cq^\prime E\pi
\label{11.75}
\end{equation}
($q^\prime$ being defined by \ref{10.260}) and a remainder which depends only on $E^2\lambda$ and $E\tchi$, the dependence on $E^2\lambdab$ and 
$E\tchib$ having been absorbed in the term $L\tau$. Transferring this term to the left hand side 
a propagation equation for the quantity $\nu-\tau$ results, 
and it is on on the basis of this equation that the associated estimates of Chapter 10 were derived 
(see Section 10.6, \ref{10.473} - \ref{10.476}). Although $\tau$ is a 1st order quantity  
so this is a lower order correction to $\nu$, it results in making the derivation of the estimates 
considerably easier. We shall presently show that this correction extends to higher dimensions. 

We begin with the higher dimensional analogue of \ref{10.262}. Recalling that $\lambdab=c\rho=cLt$ 
we have:
$$\sd\lambdab=c\sd Lt+(Lt)\sd c$$
Since $\sd Lt=\sL_L\sd t=\sL_L\pi$, this is:
\begin{equation}
\sd\lambdab=c\sL_L\pi+\rho\sd c
\label{11.76}
\end{equation}
Applying $\sD$ to \ref{11.76} and taking account of the fact that 
$[\sd c]_{P.A.}=-c(\tchi^\sharp+\tchib^\sharp)\cdot\pi$
we obtain:
\begin{equation}
\sD^2\lambdab=c\sD(\sL_L\pi)-\lambdab\sD\left((\tchi^\sharp+\tchib^\sharp)\cdot\pi\right)
\label{11.77}
\end{equation}
Now the formula \ref{11.31} implies:
\begin{equation}
(\sD(\sL_L\pi))_{AB}-(\sL_L(\sD\pi))_{AB}=(\sD_A\chi_{BC}+\sD_B\chi_{AC}-\sD_C\chi_{AB})\pi^C
\label{11.78}
\end{equation}
hence:
$$\left[(\sD(\sL_L\pi))_{AB}-(\sL_L(\sD\pi))_{AB}\right]_{P.A.}=
\left[\rho(\sD_A\tchi_{BC}+\sD_B\tchi_{AC}-\sD^C\tchi_{AB})\pi^C\right]_{P.A.}$$
In view of the Codazzi equations of Proposition 3.6 this simplifies to: 
$$\left[(\sD(\sL_L\pi))_{AB}-(\sL_L(\sD\pi))_{AB}\right]_{P.A.}=
\left[\rho\pi^C\sD_C\tchi_{AB}\right]_{P.A.}$$
that is:
\begin{equation}
\sD(\sL_L\pi)-\sL_L(\sD\pi)=\rho\sD_{\pi^\sharp}\tchi+F_1
\label{11.79}
\end{equation}
where $F_1$ is of order 2 but with vanishing P.A. part. In regard to the 2nd term on the right in 
\ref{11.77} we have:
\begin{eqnarray*}
&&\left[\sD\left((\tchi^\sharp+\tchib^\sharp)\cdot\pi\right)_{AB}\right]_{P.A.}
=\left[(\sD_A\tchi_{BC}+\sD_A\tchib_{BC})\pi^C\right]_{P.A.}\\
&&\hspace{40mm}=\left[\pi^C(\sD_C\tchi_{AB}+\sD_C\tchib_{AB}\right]_{P.A.}
\end{eqnarray*}
the last again by virtue of the Codazzi equations of Proposition 3.6. That is:
\begin{equation}
\sD\left((\tchi^\sharp+\tchib^\sharp)\cdot\pi\right)=\sD_{\pi^\sharp}\tchi+\sD_{\pi^\sharp}\tchib
+F_2
\label{11.80}
\end{equation}
where $F_2$ is of order 2 but with vanishing P.A. part. Substituting \ref{11.79} and \ref{11.80} in 
\ref{11.77} we obtain the extension of the formula \ref{10.272} to higher dimensions:
\begin{equation}
\sD^2\lambdab=c\sL_L(\sD\pi)-\lambdab\sD_{\pi^\sharp}\tchib+F
\label{11.81}
\end{equation}
where 
$$F=cF_1-\lambdab F_2$$
is again of order 2 with vanishing P.A. part. 

Going now back to \ref{11.73}, recalling that $a\lambda=\lambda^2\lambdab/c$, the last term combines 
with the first part of the 5th term to:
\begin{eqnarray}
&&\frac{\lambda^2\beta_N^2}{c}(\sD^2\lambdab+\lambdab\sD_{\pi^\sharp}\tchib)(\pi^\sharp,(\sd H)^\sharp)
\nonumber\\
&&=\lambda^2\beta_N^2 (\sL_L(\sD\pi))(\pi^\sharp,(\sd H)^\sharp) + F_3
\label{11.82}
\end{eqnarray}
where we have substituted for $\sD^2\lambdab$ from \ref{11.81}. Here
$$F_3=\frac{\lambda^2\beta_N^2}{c}F(\pi^\sharp,(\sd H)^\sharp)$$
is again of order 2 with vanishing P.A. part.  To evaluate 
$(\sL_L(\sD\pi))(\pi^\sharp,(\sd H)^\sharp)$ we recall from Section 3.1 that 
if $X$, $Y$ is a pair of $S$ vectorfields 
and $\theta$ a 2-covariant $S$ tensorfield we have:
\begin{eqnarray*}
&&(\sL_L\theta)(X,Y)=({\cal L}_L\theta)(X,Y)\\
&&\hspace{20mm}=L(\theta(X,Y))-\theta([L,X],Y)-\theta(X,[L,Y])
\end{eqnarray*}
Here we have $\pi^\sharp$, $(\sd H)^\sharp$ in the role of $X$, $Y$ and $\sL_L(\sD\pi)$ in the role 
of $\theta$. Then 
$$[L,\pi^\sharp]=-2\chi^{\sharp\sharp}\cdot\pi+(\sL_L\pi)^\sharp$$
is of order 1, while
$$[L,(\sd H)^\sharp]=-2\chi^{\sharp\sharp}\cdot\sd H+(\sd LH)^\sharp$$
is of order 2 but with vanishing P.A. part. We conclude that:
\begin{equation}
(\sL_L(\sD\pi))(\pi^\sharp,(\sd H)^\sharp)=L\left((\sD_{\pi^\sharp}\pi^\sharp)\cdot\sd H\right)
+F_3
\label{11.83}
\end{equation}
where $F_3$ is of order 2 with vanishing P.A. part. Here we recall the following fact in 
Riemannian geometry: if $f$ is a function on a Riemannian manifold $({\cal M},g)$ and
\begin{equation} 
X=(df)^\sharp
\label{11.84}
\end{equation}
is the corresponding gradient vectorfield, then $\nabla_X X$, the acceleration of 
the integral curves of $X$, is given by:
\begin{equation}
\nabla_X X=\frac{1}{2}\left(d(|df|^2)\right)^\sharp
\label{11.85}
\end{equation}
Taking then $(S_{\ub,u},\sh)$ in the role of $({\cal M},g)$ and $t$ in the role of $f$ hence 
$\pi^\sharp$ in the role of $X$, we obtain:
\begin{equation}
\sD_{\pi^\sharp}\pi^\sharp=\frac{1}{2}\left(\sd(|\pi|^2)\right)^\sharp
\label{11.86}
\end{equation}
Therefore \ref{11.83} takes the form:
\begin{equation}
(\sL_L(\sD\pi))(\pi^\sharp,(\sd H)^\sharp)=
\frac{1}{2}L\left(\sd H,\sd(|\pi|^2)\right)+F_3
\label{11.87}
\end{equation}
Consequently \ref{11.82} is equal to:
\begin{equation}
L\left\{\frac{1}{2}\lambda^2\beta_N^2\left(\sd H,\sd(|\pi|^2)\right)\right\}+F_4
\label{11.88}
\end{equation}
where $F_4$ is of order 2 but with vanishing P.A. part. 

Going going again back to \ref{11.73}, the next to last term combines 
with the second part of the 5th term to:
\begin{eqnarray}
&&-\frac{1}{2}\frac{\lambda^2\beta_N^2}{c}(\pi,\sd H)\left(\slap\lambda
+\lambdab\pi^\sharp\mbox{tr}\tchib\right)\nonumber\\
&&=-\frac{1}{2}\lambda^2\beta_N^2(\pi,\sd H)\left\{L(\sdiv\pi)+\chi^{\sharp\sharp}\cdot\sD\pi
+\frac{\mbox{tr}F}{c}\right\}\nonumber\\
&&=-\frac{1}{2}L\left\{\lambda^2\beta_N^2(\pi,\sd H)\sdiv\pi\right\}+F_5
\label{11.89}
\end{eqnarray}
where we have substituted for $\slap\lambdab$ from the trace of \ref{11.81}. Here again 
$F_5$ is of order 2 but with vanishing P.A. part. Also, the 4th and 6th terms in \ref{11.73} combine to: 
\begin{eqnarray}
&&q\lambda\left(\slap\lambda
+\lambdab\pi^\sharp\mbox{tr}\tchib\right)\nonumber\\
&&=cq\lambda\left\{L(\sdiv\pi)+\chi^{\sharp\sharp}\cdot\sD\pi
+\frac{\mbox{tr}F}{c}\right\}\nonumber\\
&&=L\left(cq\lambda\sdiv\pi\right)+F_6
\label{11.90}
\end{eqnarray}
by the same substitution, the function $F_6$ being of order 2 with vanishing P.A. part. 

In conclusion, we have succeeded to express the sum of all the terms in \ref{11.73} which involve 
$\sD^2\lambdab$, $\sD\tchib$ in the form $L\tau$ up to a remainder which is of order 2 but 
with vanishing principal acoustical part. We arrive in this way at the following extension of 
Proposition 10.2, as reformulated by \ref{10.277} - \ref{10.280} and their conjugates, to higher 
dimensions.

\vspace{2.5mm}

\noindent{\bf Proposition 11.2} \ \ The 2nd order quantities: 
$$\nu=\lambda\slap\lambda-j, \ \ \ \nub=\lambdab\slap\lambdab-\jb$$
where 
\begin{eqnarray*}
&&j=\frac{1}{4}\beta_N^2\Lb^2 H-\frac{1}{2}\lambda\beta_N^2\pi^\sharp\cdot\sd H \\
&&\hspace{6mm}+\lambda^2\left\{\beta_N\left[(\pi,\sbeta)-\frac{1}{4c}(\beta_N+2\beta_{\Nb})\right]
\slap H\right. \\
&&\hspace{6mm}\left.+H\left[(\pi,\sbeta)-\frac{1}{2c}(\beta_N+\beta_{\Nb})\right]N^\mu\slap\beta_\mu 
\right\}
\end{eqnarray*}
and
\begin{eqnarray*}
&&\jb=\frac{1}{4}\beta_{\Nb}^2 L^2 H-\frac{1}{2}\lambdab\beta_{\Nb}^2\pi^\sharp\cdot\sd H \\
&&\hspace{6mm}+\lambdab^2\left\{\beta_{\Nb}\left[(\pi,\sbeta)-\frac{1}{4c}(2\beta_N+\beta_{\Nb})\right]
\slap H\right. \\
&&\hspace{6mm}\left.+H\left[(\pi,\sbeta)-\frac{1}{2c}(\beta_N+\beta_{\Nb})\right]\Nb^\mu\slap\beta_\mu 
\right\}
\end{eqnarray*}
together with the 1st order quantities:
\begin{eqnarray*}
&&\tau=\lambda cq\sdiv\pi+\frac{1}{2}\lambda^2\beta_N^2\left\{\left(\sd H,\sd(|\pi|^2)\right)
-(\pi,\sd H)\sdiv\pi\right\}\\
&&\taub=\lambdab c\qb\sdiv\pi+\frac{1}{2}\lambdab^2\beta_{\Nb}^2\left\{\left(\sd H,\sd(|\pi|^2)\right)
-(\pi,\sd H)\sdiv\pi\right\}
\end{eqnarray*}
satisfy the following regularized propagation equations:
$$L(\nu-\tau)=I, \ \ \ \Lb(\nub-\taub)=\Ib$$
where $I$, $\Ib$ are quantities of order 2, their principal acoustical parts being given by
\begin{eqnarray*}
&&[I]_{P.A.}=2(L\lambda)\slap\lambda-2\lambda\chi^{\sharp\sharp}\cdot\sD^2\lambda
+(\lambda A-a(\sd\lambda)^\sharp)\cdot\sd\mbox{tr}\tchi\\
&&[\Ib]_{P.A.}=2(\Lb\lambdab)\slap\lambdab-2\lambdab\chib^{\sharp\sharp}\cdot\sD^2\lambdab
+(\lambdab\oAb-a(\sd\lambdab)^\sharp)\cdot\sd\mbox{tr}\tchib
\end{eqnarray*}
(Here $A$ is given by \ref{11.74} and $\oAb$ by its conjugate.)

\vspace{2.5mm}

Going back to Proposition 11.1 we introduce the $S$ 1-forms:
\begin{equation}
\theta_1=\lambda\sd\mbox{tr}\tchi+\sd f, \ \ \ \thetab_1=\lambdab\sd\mbox{tr}\tchib+\sd\fb
\label{11.91}
\end{equation}
These are quantities of order 2. It then follows from Proposition 11.1 that these quantities satisfy 
propagation equations of the form:
\begin{equation}
\sL_L\theta_1=R_1, \ \ \ \sL_L\thetab_1=\Rb_1
\label{11.92}
\end{equation}
where $R_1$, $\Rb_1$ are likewise $S$ 1-forms which are quantities of order 2 and we have:
\begin{eqnarray}
&&[R_1]_{P.A.}=[\sd R]_{P.A.}=2(L\lambda)\sd\mbox{tr}\tchi-a\left[\sd\left(|\tchi|^2\right)\right]_{P.A.}
-\lambda(LH)\left[\sd\left(\tchi^{\sharp\sharp}(\sbeta,\sbeta)\right)\right]_{P.A.}\nonumber\\
&&[\Rb_1]_{P.A.}=[\sd\Rb]_{P.A.}=2(\Lb\lambdab)\sd\mbox{tr}\tchib-a\left[\sd\left(|\tchib|^2\right)\right]_{P.A.}
-\lambda(\Lb H)\left[\sd\left(\tchib^{\sharp\sharp}(\sbeta,\sbeta)\right)\right]_{P.A.}\nonumber\\
\label{11.93}
\end{eqnarray}
By \ref{8.142} to control the P.A. part of the $(0,l)$ sources as in the 1st of \ref{8.156} 
we take strings of commutation fields $E_{(\nu)}$ of length $l-1$, defining the $S$ 1-forms:
\begin{eqnarray}
&&\s^{(\nu_1...\nu_{l-1})}\theta_l=\lambda \sd\left(E_{(\nu_{l-1})} . . . E_{(\nu_1)}
\mbox{tr}\tchi\right)+\sd\left(E_{(\nu_{l-1})} . . . E_{(\nu_1)}f\right)\nonumber\\
&&\s^{(\nu_1...\nu_{l-1})}\thetab_l=\lambdab \sd\left(E_{(\nu_{l-1})} . . . E_{(\nu_1)}
\mbox{tr}\tchib\right)+\sd\left(E_{(\nu_{l-1})} . . . E_{(\nu_1)}\fb\right)\nonumber\\
&&\label{11.94}
\end{eqnarray}
These are quantities of order $l+1$. It then follows that these quantities satisfy propagation 
equations of the form:
\begin{eqnarray}
&&\sL_L\s^{(\nu_1...\nu_{l-1})}\theta_l=\s^{(\nu_1...\nu_{l-1})}R_l\nonumber\\
&&\sL_L\s^{(\nu_1...\nu_{l-1})}\thetab_l=\s^{(\nu_1...\nu_{l-1})}\Rb_l \label{11.95}
\end{eqnarray}
where $\s^{(\nu_1...\nu_{l-1})}R_l$, $\s^{(\nu_1...\nu_{l-1})}\Rb_l$ 
are likewise $S$ 1-forms which are quantities of order $l+1$ and we have:
\begin{eqnarray}
&&\left[\s^{(\nu_1...\nu_{l-1})}R_l\right]_{P.A.}=
\left[\sd\left(E_{(\nu_{l-1})} . . . E_{(\nu_1)}R\right)\right]_{P.A.}\nonumber\\
&&\hspace{20mm}=2(L\lambda)\sd\left(E_{(\nu_{l-1})} . . . E_{(\nu_1)}
\mbox{tr}\tchi\right)-a\left[\sd\left(E_{(\nu_{l-1})} . . . E_{(\nu_1)}\left(|\tchi|^2\right)\right)\right]_{P.A.}\nonumber\\
&&\hspace{30mm}
-\lambda(LH)\left[\sd\left(E_{(\nu_{l-1})} . . . E_{(\nu_1)}\left(\tchi^{\sharp\sharp}(\sbeta,\sbeta)\right)\right)\right]_{P.A.}\nonumber\\
&&\label{11.96}
\end{eqnarray}
\begin{eqnarray}
&&\left[\s^{(\nu_1...\nu_{l-1})}\Rb_l\right]_{P.A.}=
\left[\sd\left(E_{(\nu_{l-1})} . . . E_{(\nu_1)}\Rb\right)\right]_{P.A.}\nonumber\\
&&\hspace{20mm}=2(\Lb\lambdab)\sd\left(E_{(\nu_{l-1})} . . . E_{(\nu_1)}
\mbox{tr}\tchib\right)
-a\left[\sd\left(E_{(\nu_{l-1})} . . . E_{(\nu_1)}\left(|\tchib|^2\right)\right)\right]_{P.A.}
\nonumber\\
&&\hspace{30mm}
-\lambdab(\Lb H)\left[\sd\left(E_{(\nu_{l-1})} . . . E_{(\nu_1)}\left(\tchib^{\sharp\sharp}(\sbeta,\sbeta)\right)\right)\right]_{P.A.}\nonumber\\
&&\label{11.97}
\end{eqnarray}

By \ref{8.135} to control the P.A. part of the $(m,l)$ sources as in the 2nd of \ref{8.156} 
we take strings of commutation fields of length $m+l-1$ of the form 
$E_{(\nu_l)}... E_{(\nu_1)}T^{m-1}$, defining the functions:
\begin{eqnarray}
&&\s^{(\nu_1 ... \nu_l)}\nu_{m-1,l+1}=\lambda E_{(\nu_l)} ... E_{(\nu_1)}T^{m-1}\slap\lambda
-E_{(\nu_l)} ... E_{(\nu_1)}T^{m-1}j\nonumber\\
&&\s^{(\nu_1 ... \nu_l)}\nub_{m-1,l+1}=\lambdab E_{(\nu_l)} ... E_{(\nu_1)}T^{m-1}\slap\lambdab
-E_{(\nu_l)} ... E_{(\nu_1)}T^{m-1}\jb \nonumber\\
&&\label{11.98}
\end{eqnarray}
These are quantities of order $m+l+1$. We also define the functions:
\begin{eqnarray}
&&\s^{(\nu_1 ... \nu_l)}\tau_{m-1,l+1}=E_{(\nu_l)} ... E_{(\nu_1)}T^{m-1}\tau\nonumber\\
&&\s^{(\nu_1 ... \nu_l)}\taub_{m-1,l+1}=E_{(\nu_l)} ... E_{(\nu_1)}T^{m-1}\taub \label{11.99}
\end{eqnarray}
which are quantities of order $m+l$. It then follows from Proposition 11.2 that these 
quantities satisfy propagation equations of the form:
\begin{eqnarray}
&&L\left(\s^{(\nu_1 ... \nu_l)}\nu_{m-1,l+1}-\s^{(\nu_1 ... \nu_l)}\tau_{m-1,l+1}\right)
=\s^{(\nu_1 ... \nu_l)}I_{m-1,l+1}\nonumber\\
&&\Lb\left(\s^{(\nu_1 ... \nu_l)}\nub_{m-1,l+1}-\s^{(\nu_1 ... \nu_l)}\taub_{m-1,l+1}\right)
=\s^{(\nu_1 ... \nu_l)}\Ib_{m-1,l+1}\nonumber\\
\label{11.100}
\end{eqnarray}
where $\s^{(\nu_1 ... \nu_l)}I_{m-1,l+1}$, $\s^{(\nu_1 ... \nu_l)}\Ib_{m-1,l+1}$ are  quantities 
of order $m+l+1$ and we have:
\begin{eqnarray}
&&\left[\s^{(\nu_1 ... \nu_l)}I_{m-1,l+1}\right]_{P.A.}=
\left[E_{(\nu_l)} ... E_{(\nu_1)}T^{m-1}I\right]_{P.A.}-(m-1)Z\s^{(\nu_1 ... \nu_l)}\nu_{m-2,l+1}
\nonumber\\
&&\hspace{14mm}=2(L\lambda)E_{(\nu_l)} ... E_{(\nu_1)}T^{m-1}\slap\lambda
-2\lambda\left[E_{(\nu_l)} ... E_{(\nu_1)}T^{m-1}\left(\chi^{\sharp\sharp}\cdot\sD^2\lambda\right)
\right]_{P.A.}\nonumber\\
&&\hspace{8mm}+(\lambda A-a(\sd\lambda)^\sharp)\cdot \left[\sd(E_{(\nu_l)} ... E_{(\nu_1)}T^{m-1}\mbox{tr}\tchi)\right]_{P.A.}-(m-1)Z\s^{(\nu_1 ... \nu_l)}\nu_{m-2,l+1}\nonumber\\
&&\label{11.101}
\end{eqnarray}
\begin{eqnarray}
&&\left[\s^{(\nu_1 ... \nu_l)}\Ib_{m-1,l+1}\right]_{P.A.}=
\left[E_{(\nu_l)} ... E_{(\nu_1)}T^{m-1}\Ib\right]_{P.A.}+(m-1)Z\s^{(\nu_1 ... \nu_l)}\nub_{m-2,l+1}
\nonumber\\
&&\hspace{14mm}=2(\Lb\lambdab)E_{(\nu_l)} ... E_{(\nu_1)}T^{m-1}\slap\lambdab
-2\lambdab\left[E_{(\nu_l)} ... E_{(\nu_1)}T^{m-1}\left(\chib^{\sharp\sharp}\cdot\sD^2\lambdab\right)
\right]_{P.A.}\nonumber\\
&&\hspace{8mm}+(\lambdab\oAb-a(\sd\lambda)^\sharp)\cdot \left[\sd(E_{(\nu_l)} ... E_{(\nu_1)}T^{m-1}\mbox{tr}\tchib)\right]_{P.A.}+(m-1)Z\s^{(\nu_1 ... \nu_l)}\nub_{m-2,l+1}\nonumber\\
&&\label{11.102}
\end{eqnarray}

Proceeding as in Chapter 10 we then consider the analogous propagation equations for the 
corresponding $N$th approximants. These contain error terms similar to those in \ref{10.254}, 
\ref{10.255}, \ref{10.330}, \ref{10.333}. 
Subtracting, we then obtain propagation equations for the acoustical difference quantities:
\begin{eqnarray}
&&\s^{(\nu_1...\nu_{l-1})}\cth_l=\s^{(\nu_1...\nu_{l-1})}\theta_l-\s^{(\nu_1...\nu_{l-1})}\theta_{l,N}\nonumber\\
&&\s^{(\nu_1...\nu_{l-1})}\cthb_l=\s^{(\nu_1...\nu_{l-1})}\thetab_l-\s^{(\nu_1...\nu_{l-1})}\thetab_{l,N} 
\label{11.103}
\end{eqnarray}
where at the top order $l=n$, and:
\begin{eqnarray}
&&\s^{(\nu_1 ... \nu_l)}\cnu_{m-1,l+1}=\s^{(\nu_1 ... \nu_l)}\nu_{m-1,l+1}-\s^{(\nu_1 ... \nu_l)}\nu_{m-1,l+1,N}\nonumber\\
&&\s^{(\nu_1 ... \cnub_l)}\cnu_{m-1,l+1}=\s^{(\nu_1 ... \nu_l)}\nub_{m-1,l+1}-\s^{(\nu_1 ... \nu_l)}\nub_{m-1,l+1,N}
\label{11.104}
\end{eqnarray}
where at the top order $m=1,...,n$, $l=n-m$. These propagation equations are treated in a similar 
manner as the corresponding equations in Chapter 10, after the essential replacements to be presently 
discussed are made. 

Consider the last two terms on the right in \ref{11.96}, \ref{11.97}. 
In regard to the next to last term in \ref{11.96}, we have:
\begin{equation}
\left[E_{(\nu_{(l-1)}} ... E_{(\nu_1)}\left(|\tchi|^2\right)\right]_{P.A.}
=2(\tchi,\sL_{E_{(\nu_{l-1})}} ... \sL_{E_{(\nu_1)}}\tchi)
\label{11.105}
\end{equation}
hence this term is:
\begin{equation}
-2a\left(\tchi,\sD\sL_{E_{(\nu_{l-1})}} ... \sL_{E_{(\nu_1)}}\tchi\right)
\label{11.106}
\end{equation}
In regard to the last term in \ref{11.96}, we have:
\begin{equation}
\left[E_{(\nu_{(l-1)}} ... E_{(\nu_1)}\left(\tchi^{\sharp\sharp}(\sbeta,\sbeta)\right)\right]_{P.A.}
=\left(\sL_{E_{(\nu_{l-1})}} ... \sL_{E_{(\nu_1)}}\tchi\right)^{\sharp\sharp}(\sbeta,\sbeta)
\label{11.107}
\end{equation}
hence this term is:
\begin{equation}
-\lambda(LH)\left(\sD\sL_{E_{(\nu_{l-1})}} ... \sL_{E_{(\nu_1)}}\tchi\right)^{\sharp\sharp}(\sbeta,\sbeta)
\label{11.108}
\end{equation}
Then, differences from the corresponding $N$th approximants being taken, the last two terms in 
\ref{11.96} become, up to lower order terms,
\begin{equation}
-2a\left(\tchi,\sD\s^{(\nu_1 ... \nu_{l-1})}\ctchi_{l-1}\right)
-\lambda(LH)\left(\sD\s^{(\nu_1 ... \nu_{l-1})}\ctchi_{l-1}\right)^{\sharp\sharp}
(\sbeta,\sbeta)
\label{11.109}
\end{equation}
Similarly, differences from the corresponding $N$th approximants being taken, the last two terms in 
\ref{11.97} become, up to lower order terms,
\begin{equation}
-2a\left(\tchib,\sD\s^{(\nu_1 ... \nu_{l-1})}\ctchib_{l-1}\right)
-\lambdab(\Lb H)\left(\sD\s^{(\nu_1 ... \nu_{l-1})}\ctchib_{l-1}\right)^{\sharp\sharp}
(\sbeta,\sbeta)
\label{11.110}
\end{equation}
In the above we denote:
\begin{eqnarray}
&&\s^{(\nu_1 ... \nu_{l-1})}\ctchi_{l-1}=\sL_{E_{(\nu_{l-1})}} ... \sL_{E_{(\nu_1)}}\tchi-
\sL_{E_{(\nu_{l-1}),N}} ... \sL_{E_{(\nu_1),N}}\tchi_N\nonumber\\
&&\s^{(\nu_1 ... \nu_{l-1})}\ctchib_{l-1}=\sL_{E_{(\nu_{l-1})}} ... \sL_{E_{(\nu_1)}}\tchib-
\sL_{E_{(\nu_{l-1}),N}} ... \sL_{E_{(\nu_1),N}}\tchib_N\nonumber\\
&&\label{11.a1}
\end{eqnarray}
the vectorfields $E_{(\mu),N}$ being the $N$th approximant versions of the vectorfields $E_{(\mu)}$. 
To estimate the contribution of \ref{11.109} to the integral on the right in \ref{10.361} 
we must estimate:
\begin{equation}
\|\lambda\sD\s^{(\nu_1 ... \nu_{l-1})}\ctchi_{l-1}\|_{L^2(S_{\ub,u})}
\label{11.111}
\end{equation}
Similarly, to estimate the contribution of \ref{11.110} to the integral on the right in \ref{10.618} 
we must estimate:
\begin{equation}
\|\lambdab\sD\s^{(\nu_1 ... \nu_{l-1})}\ctchib_{l-1}\|_{L^2(S_{\ub,u})}
\label{11.112}
\end{equation}
We shall show below how to estimate \ref{11.111} in terms of:
\begin{equation}
\|\s^{(\nu_1...\nu_{l-1})}\cth_l\|_{L^2(S_{\ub,u})}
\label{11.113}
\end{equation}
and how to estimate \ref{11.112} in terms of:
\begin{equation}
\|\s^{(\nu_1...\nu_{l-1})}\cthb_l\|_{L^2(S_{\ub,u})}
\label{11.114}
\end{equation}
using elliptic theory on $S_{\ub,u}$ in connection with the Codazzi equations (Proposition 3.6). 

Consider the second term on the right in each of \ref{11.101}, \ref{11.102}. 
In regard to the second term in \ref{11.101}, we have:
\begin{equation}
\left[E_{(\nu_{(l)}} ... E_{(\nu_1)}T^{m-1}\left(\tchi^{\sharp\sharp}\cdot\sD^2\lambda\right)\right]_{P.A.}
=\tchi^{\sharp\sharp}\cdot\sD^2(E_{(\nu_l)} ... E_{(\nu_1)}T^{m-1}\lambda)
\label{11.115}
\end{equation}
hence this term is:
\begin{equation}
-2\lambda\tchi^{\sharp\sharp}\cdot\sD^2(E_{(\nu_l)} ... E_{(\nu_1)}T^{m-1}\lambda)
\label{11.116}
\end{equation}
Then, differences from the corresponding $N$th approximants being taken, the second term in 
\ref{11.101} becomes, up to lower order terms, 
\begin{equation}
-2\lambda\tchi^{\sharp\sharp}\cdot\sD^2\s^{(\nu_1 ... \nu_l)}\cla_{m-1,l}
\label{11.117}
\end{equation}
Similarly, differences from the corresponding $N$th approximants being taken, the second term in 
\ref{11.102} becomes, up to lower order terms, 
\begin{equation}
-2\lambdab\tchib^{\sharp\sharp}\cdot\sD^2\s^{(\nu_1 ... \nu_l)}\clab_{m-1,l}
\label{11.118}
\end{equation}
In the above we denote:
\begin{eqnarray}
&&\s^{(\nu_1 ... \nu_l)}\cla_{m-1,l}=E_{(\nu_l)} ... E_{(\nu_1)}T^{m-1}\lambda
-E_{(\nu_l),N} ... E_{(\nu_1),N}T^{m-1}\lambda_N\nonumber\\
&&\s^{(\nu_1 ... \nu_l)}\clab_{m-1,l}=E_{(\nu_l)} ... E_{(\nu_1)}T^{m-1}\lambdab
-E_{(\nu_l),N} ... E_{(\nu_1),N}T^{m-1}\lambdab_N\nonumber\\
&&\label{11.a2}
\end{eqnarray}
To estimate the contribution of \ref{11.117} to the integral on the right in \ref{10.480} we must 
estimate:
\begin{equation}
\|\lambda\sD^2\s^{(\nu_1 ... \nu_l)}\cla_{m-1,l}\|_{L^2(S_{\ub,u})}
\label{11.119}
\end{equation}
To estimate the contribution of \ref{11.118} to the integral on the right in \ref{10.671} we must 
estimate:
\begin{equation}
\|\lambdab\sD^2\s^{(\nu_1 ... \nu_l)}\clab_{m-1,l}\|_{L^2(S_{\ub,u})}
\label{11.120}
\end{equation}
We shall show below how to estimate \ref{11.119} in terms of:
\begin{equation}
\|\s^{(\nu_1 ... \nu_l)}\cnu_{m-1,l+1}\|_{L^2(S_{\ub,u})}
\label{11.121}
\end{equation}
and how to estimate \ref{11.120} in terms of:
\begin{equation}
\|\s^{(\nu_1 ... \nu_l)}\cnub_{m-1,l+1}\|_{L^2(S_{\ub,u})}
\label{11.122}
\end{equation}
using elliptic theory on $S_{\ub,u}$. 

The essential replacements then consist of the replacement of \ref{11.111}, \ref{11.112} by their 
bounds in terms of \ref{11.113}, \ref{11.114} in assessing their contributions to the integrals  
in \ref{10.361}, \ref{10.618}, and of the replacement of \ref{11.119}, \ref{11.120} 
by their bounds in terms of \ref{11.121}, \ref{11.122} in assessing their contributions to the 
integrals in \ref{10.480}, \ref{10.671}. 

The Codazzi equations of Proposition 3.6 are of the form, in terms of components, writing 
$\sk_{AB}$ for $(\sk_\flat)_{AB}$ (see \ref{11.2}):
\begin{eqnarray}
&&\sD_A\sk_{BC}-\sD_B\sk_{AC}=l_{ABC}\nonumber\\
&&\sD_A\skb_{BC}-\sD_B\skb_{AC}=\lb_{ABC}
\label{11.123}
\end{eqnarray}
where $l_{ABC}$, $\lb_{ABC}$ are the components, antisymmetric in the first two indices, of type 
$\wedge_{2,1}$ $S$ tensorfields $l$, $\lb$, sections at each $S_{\ub,u}$ of the tensor bundle
$$\bigcup_{q\in S_{\ub,u}} \wedge_{2,1}(T_q S_{\ub,u})$$
over $S_{\ub,u}$, where $\wedge_{2,1}(T_q S_{\ub,u})$ denotes the space of antisymmetric bilinear 
forms on $T_q S_{\ub,u}$ with values in $T^*_q S_{\ub,u}$. The $S$ tensorfields $l$, $\lb$ are only 
of order 1. Substituting the expression \ref{11.2} of $\sk_\flat$, $\skb_\flat$ in terms of 
$\tchi$, $\tchib$, equations \ref{11.123} become:
\begin{eqnarray}
&&\sD_A\tchi_{BC}-\sD_B\tchi_{AC}=\psi_{ABC}\nonumber\\
&&\sD_A\tchib_{BC}-\sD_B\tchib_{AC}=\psib_{ABC}
\label{11.124}
\end{eqnarray}
where $\psi_{ABC}$, $\psib_{ABC}$ are the components of $S$ tensorfields $\psi$, $\psib$ of the 
same type which are of order 2 but with vanishing P.A. part. We apply the string 
$\sL_{E_{(\nu_{l-1})}} ... \sL_{E_{(\nu_1)}}$ to equations \ref{11.124} and note the fact that the 
commutator $[\sD,\sL_{E_{(\mu)}}]$ applied to a $S$ tensorfield of any type can be expressed in terms 
of: 
\begin{equation}
\s^{(E_{(\mu)})}\spi_{1,AB}^C=\frac{1}{2}\left(\sD_A\s^{(E_{(\mu)})}\sspi_B^C
+\sD_B\s^{(E_{(\mu)})}\sspi_A^C-\sD^C\s^{(E_{(\mu)})}\sspi_{AB}\right)
\label{11.125}
\end{equation}
the Lie derivative of the $S_{\ub,u}$ connection $\sGamma$ (compare with \ref{11.31}) with respect 
to the $S$ vectorfield $E_{(\mu)}$. Here (see \ref{8.b8}):
\begin{equation}
\s^{(E_{(\mu)})}\sspi=\sL_{E_{(\mu)}}\sh 
\label{11.126}
\end{equation}
is a quantity of order 1 expressed by \ref{8.89}. Hence $\s^{(E_{(\mu)})}\spi_1$, a type 
$S^1_2$ $S$ tensorfield, is of order 2. We obtain:
\begin{eqnarray}
&&\sD_A\left(\sL_{E_{(\nu_{l-1})}} ... \sL_{E_{(\nu_1)}}\tchi\right)_{BC}-
\sD_B\left(\sL_{E_{(\nu_{l-1})}} ... \sL_{E_{(\nu_1)}}\tchi\right)_{AC}
=\left(\s^{(\nu_1 ... \nu_{l-1})}\psi_{l-1}\right)_{ABC}\nonumber\\
&&\sD_A\left(\sL_{E_{(\nu_{l-1})}} ... \sL_{E_{(\nu_1)}}\tchib\right)_{BC}-
\sD_B\left(\sL_{E_{(\nu_{l-1})}} ... \sL_{E_{(\nu_1)}}\tchib\right)_{AC}
=\left(\s^{(\nu_1 ... \nu_{l-1})}\psib_{l-1}\right)_{ABC}\nonumber\\
\label{11.127}
\end{eqnarray}
where $\s^{(\nu_1 ... \nu_{l-1})}\psi_{l-1}$, $\s^{(\nu_1 ... \nu_{l-1})}\psib_{l-1}$ are 
type $\wedge_{2,1}$ $S$ tensorfields of order $l+1$, with
\begin{eqnarray}
&&\left[\s^{(\nu_1 ... \nu_{l-1})}\psi_{l-1}\right]_{P.P.}=
\left[\sL_{E_{(\nu_{l-1})}} ... \sL_{E_{(\nu_1)}}\psi\right]_{P.P.}\nonumber\\
&&\left[\s^{(\nu_1 ... \nu_{l-1})}\psib_{l-1}\right]_{P.P.}=
\left[\sL_{E_{(\nu_{l-1})}} ... \sL_{E_{(\nu_1)}}\psib\right]_{P.P.}
\label{11.128}
\end{eqnarray}
Analogous equations hold for the $N$th approximants. Taking differences then yields, in view of the 
definitions \ref{11.a1}, 
\begin{eqnarray}
&&\sD_A\left(\s^{(\nu_1 ... \nu_{l-1})}\ctchi_{l-1}\right)_{BC}
-\sD_B\left(\s^{(\nu_1 ... \nu_{l-1})}\ctchi_{l-1}\right)_{AC}\nonumber\\
&&\hspace{37mm}=
\left(\s^{(\nu_1 ... \nu_{l-1})}\check{\psi}_{l-1}\right)_{ABC}+\svep_{ABC,N}\nonumber\\
&&\sD_A\left(\s^{(\nu_1 ... \nu_{l-1})}\ctchib_{l-1}\right)_{BC}
-\sD_B\left(\s^{(\nu_1 ... \nu_{l-1})}\ctchib_{l-1}\right)_{AC}\nonumber\\
&&\hspace{37mm}=
\left(\s^{(\nu_1 ... \nu_{l-1})}\check{\psib}_{l-1}\right)_{ABC}+\svepb_{ABC,N}\nonumber\\
&&\label{11.129}
\end{eqnarray}
where $\s^{(\nu_1 ... \nu_{l-1})}\check{\psi}_{l-1}$, $\s^{(\nu_1 ... \nu_{l-1})}\check{\psib}_{l-1}$ 
are difference quantities of order $l+1$ with vanishing principal acoustical parts, and
\begin{equation}
\svep_N, \svepb_N = O(\tau^{N+1})
\label{11.130}
\end{equation}
correspond to the errors by which the $N$th aproximants fail to satisfy the Codazzi equations 
(see Lemma 9.1). Each of \ref{11.129} constitutes an elliptic system on $S_{\ub,u}$, the 
first for $\s^{(\nu_1 ... \nu_{l-1})}\ctchi_{l-1}$ given its trace and the second for 
$\s^{(\nu_1 ... \nu_{l-1})}\ctchib_{l-1}$ given its trace. We note that by virtue of the 
definitions \ref{11.a1} these traces are up to lower order terms given by:
\begin{eqnarray}
&&E_{(\nu_{l-1})} ... E_{(\nu_1)}\mbox{tr}\tchi
-E_{(\nu_{l-1}),N} ... E_{(\nu_1),N}\mbox{tr}_N\tchi_N\nonumber\\
&&E_{(\nu_{l-1})} ... E_{(\nu_1)}\mbox{tr}\tchib
-E_{(\nu_{l-1}),N} ... E_{(\nu_1),N}\mbox{tr}_N\tchib_N
\label{11.131}
\end{eqnarray}
respectively. (Here $\mbox{tr}_N$ signifies the trace of the 
corresponding through $\sh_N$ $T^1_1$ type $S$ tensorfields.) In view of the definitions 
\ref{11.94} and \ref{11.103} applying $\lambda\sd$ to the 1st of \ref{11.131} and $\lambdab\sd$ 
to the 2nd, we obtain, up to a 0th order term,
\begin{eqnarray}
&&\s^{(\nu_1...\nu_{l-1})}\cth_l-\sd\left(E_{(\nu_{l-1})} . . . E_{(\nu_1)}f
-E_{(\nu_{l-1}),N} . . . E_{(\nu_1),N}f_N\right)\nonumber\\
&&\s^{(\nu_1...\nu_{l-1})}\cthb_l-\sd\left(E_{(\nu_{l-1})} . . . E_{(\nu_1)}\fb
-E_{(\nu_{l-1}),N} . . . E_{(\nu_1),N}\fb_N\right) \label{11.132}
\end{eqnarray}
respectively. We are now in a position to apply to each of \ref{11.129} the following lemma. 

\vspace{2.5mm}

\noindent{\bf Lemma 11.3} \ \ Let $({\cal M},g)$ be a compact Riemannian manifold. We denote by 
$d\mu_g$ and by $\nabla$ the associated volume element and covariant derivative operator. Let 
$\theta$ be a symmetric 2-covariant tensorfield on ${\cal M}$ satisfying the system, in terms of 
components in an arbitrary local frame field,
$$\nabla_i\theta_{jk}-\nabla_j\theta_{ik}=f_{ijk}$$
where $f$ is a type $\wedge_{2,1}$ tensorfield on ${\cal M}$. Let also $\mu$ be an arbitrary 
non-negative weight function on ${\cal M}$. Then the following integral inequality holds:
\begin{eqnarray*}
&&\int_{\cal M}\mu^2\left\{\frac{1}{2}|\nabla\theta|^2+\mbox{Ric}_{ij}\theta^i_k\theta^{jk}
-R_{ijkl}\theta^{ik}\theta^{jl}\right\}d\mu_g\leq\\
&&\hspace{10mm}\int_{\cal M}
\left\{\mu^2\left(|d\mbox{tr}\theta+\mbox{tr}f|^2
+\frac{1}{2}|f|^2+|\mbox{tr}f|^2\right)+3|d\mu|^2|\theta|^2\right\}d\mu_g
\end{eqnarray*}
Here $R_{ijkl}$ are the components of the curvature tensor of $({\cal M},g)$, 
$\mbox{Ric}_{ij}=(g^{-1})^{kl}R_{ikjl}$ those of the Ricci curvature, and $\mbox{f}$ is the 1-form on 
${\cal M}$ with components $(\mbox{tr}f)_i=(g^{-1})^{jk}f_{jik}$. 

\vspace{2.5mm}

\noindent{\em Proof:} We have (indices are lowered and raised with respect to $g_{ij}$ and 
$(g^{-1})^{ij}$):
\begin{eqnarray*}
&&|f|^2=f_{ijk}f^{ijk}=(\nabla_i\theta_{jk})(\nabla^i\theta^{jk})
+(\nabla_j\theta_{ik})(\nabla^j\theta^{ik})\\
&&\hspace{25mm}-(\nabla_j\theta_{ik})(\nabla^i\theta^{jk})-(\nabla_i\theta_{jk})(\nabla^j\theta^{ik})
\end{eqnarray*}
that is:
\begin{equation}
|f|^2=2|\nabla\theta|^2-2(\nabla_j\theta_{ik})(\nabla^i\theta^{jk})
\label{11.133}
\end{equation}
In regard to the 2nd term on the right we write:
\begin{equation}
(\nabla_j\theta_{ik})(\nabla^i\theta^{jk})=\nabla_j(\theta_{ik}\nabla^i\theta^{jk})-
\theta_{ik}\nabla_j\nabla^i\theta^{jk}
\label{11.134}
\end{equation}
and we have:
\begin{equation}
\nabla_j\nabla_i\theta^{jk}-\nabla_i\nabla_j\theta^{jk}=\mbox{Ric}_{li}\theta^{lk}+R^k_{lji}\theta^{jl}
\label{11.135}
\end{equation}
Hence, substituting in \ref{11.134}, 
\begin{eqnarray}
&&(\nabla_j\theta_{ik})(\nabla^i\theta^{jk})=\nabla_j(\theta_{ik}\nabla^i\theta^{jk})-
\theta^i_k\nabla_i\nabla_j\theta^{jk}\nonumber\\
&&\hspace{27mm}-\mbox{Ric}_{li}\theta^i_k\theta^{lk}-R^k_{lji}\theta^i_k\theta^{jl}
\label{11.136}
\end{eqnarray}
In regard to the 2nd term on the right we write:
\begin{equation}
\theta^i_k\nabla_i\nabla_j\theta^{jk}=\nabla_i(\theta^i_k\nabla_j\theta^{jk})
-(\nabla_i\theta^i_k)(\nabla_j\theta^{jk})
\label{11.137}
\end{equation}
Defining then the vectorfield $J$ with components
\begin{equation}
J^i=\theta^j_k\nabla_j\theta^{ik}-\theta^i_k\nabla_j\theta^{jk}
\label{11.138}
\end{equation}
we arrive at the identity:
\begin{eqnarray}
&&(\nabla_j\theta_{ik})(\nabla^i\theta^{jk})=\nabla_i J^i+(\nabla_i\theta^i_k)(\nabla_j\theta^{jk})
\nonumber\\
&&\hspace{27mm}-\mbox{Ric}_{li}\theta^i_k\theta^{lk}-R^k_{lji}\theta^i_k\theta^{jl}
\label{11.139}
\end{eqnarray}
Next we contract the equations with $(g^{-1})^{ik}$ to obtain:
\begin{equation}
\nabla^i\theta_{ij}-\nabla_j\mbox{tr}\theta=(\mbox{tr}f)_j
\label{11.140}
\end{equation}
As a consequence the 2nd term on the right in \ref{11.139} is:
$$|d\mbox{tr}\theta+\mbox{tr}f|^2$$
Substituting, \ref{11.139} takes the form:
\begin{eqnarray}
&&(\nabla_j\theta_{ik})(\nabla^i\theta^{jk})=\mbox{div}J+|d\mbox{tr}\theta+\mbox{tr}f|^2
\nonumber\\
&&\hspace{27mm}-\mbox{Ric}_{li}\theta^i_k\theta^{lk}-R^k_{lji}\theta^i_k\theta^{jl}
\label{11.141}
\end{eqnarray}
On the other hand from \ref{11.133} we have:
\begin{equation}
(\nabla_j\theta_{ik})(\nabla^i\theta^{jk})=|\nabla\theta|^2-\frac{1}{2}|f|^2
\label{11.142}
\end{equation}
Equating the two expressions yields the identity:
\begin{eqnarray}
&&|\nabla\theta|^2=\mbox{div}J+|d\mbox{tr}\theta+\mbox{tr}f|^2+\frac{1}{2}|f|^2\nonumber\\
&&\hspace{12mm}-\mbox{Ric}_{ij}\theta^i_k\theta^{jk}+R_{ijkl}\theta^{ik}\theta^{jl}
\label{11.143}
\end{eqnarray}
We multiply this identity by $\mu^2$ and integrate over ${\cal M}$. 
Writing 
$$\mu^2\mbox{div}J=\mbox{div}(\mu^2 J)-2\mu J\cdot d\mu$$
The integral over ${\cal M}$ of $\mbox{div}(\mu^2 J)$ vanishes ${\cal M}$ being a compact 
manifold without boundary. Hence:
\begin{equation}
\int_{{\cal M}}\mu^2\mbox{div}J d\mu_g=-2\int_{\cal M}\mu(J\cdot d\mu)d\mu_g
\label{11.144}
\end{equation}
Moreover, in view of \ref{11.138} we can bound:
\begin{equation}
|J|\leq |\theta|\left(|\nabla\theta|+|\mbox{tr}f|\right)
\label{11.145}
\end{equation}
therefore:
\begin{eqnarray}
&&-2\int_{\cal M}\mu(J\cdot d\mu)d\mu_g\leq\int_{{\cal M}}\mu|J||d\mu|d\mu_g\nonumber\\
&&\hspace{30mm}\leq 2\int_{{\cal M}}|d\mu||\theta|\left(\mu|\nabla\theta|+\mu|\mbox{tr}f|\right)d\mu_g
\label{11.146}
\end{eqnarray}
Applying the inequalities
\begin{eqnarray*}
&&2|d\mu||\theta|\mu|\nabla\theta|\leq\frac{1}{2}\mu^2|\nabla\theta|^2+2|d\mu|^2|\theta|^2\\
&&2|d\mu||\theta|\mu|\mbox{tr}f|\leq\mu^2|\mbox{tr}f|^2+|d\mu|^2|\theta|^2
\end{eqnarray*}
we then obtain:
\begin{eqnarray}
&&-2\int_{\cal M}\mu(J\cdot d\mu)d\mu_g\leq \label{11.147}\\
&&\hspace{15mm}\int_{{\cal M}}\left\{\mu^2\left(\frac{1}{2}|\nabla\theta|^2
+|\mbox{tr}f|^2\right)+3|d\mu|^2|\theta|^2\right\}d\mu_g \nonumber
\end{eqnarray}
In view of \ref{11.144}, \ref{11.146}, integrating $\mu^2$ times \ref{11.143} over ${\cal M}$ 
yields the conclusion of the lemma. 

\vspace{2.5mm}

We apply Lemma 11.3 with the Riemannian manifold $(S_{\ub,u},\sh)$ in the role of $({\cal M},g)$. 
The curvature terms are controlled through the Gauss equations of Proposition 3.7. In connection 
with the 1st of \ref{11.129} we take $\s^{(\nu_1 ... \nu_{l-1})}\ctchi_{l-1}$ in the role of $\theta$ 
and $\lambda$ in the role of $\mu$. In connection with the 2nd of \ref{11.129} we take 
$\s^{(\nu_1 ... \nu_{l-1})}\ctchib_{l-1}$ in the role of $\theta$ and $\lambdab$ in the role of $\mu$. 
In view of \ref{11.131}, \ref{11.132}, we then obtain the desired estimates of 
\ref{11.111}, \ref{11.112} in terms of \ref{11.113}, \ref{11.114} respectively.

To estimate \ref{11.119}, \ref{11.120} in terms of \ref{11.121}, \ref{11.122} we apply the 
string $E_{(\nu_l)} ... E_{(\nu_1)}T^{m-1}$ to $\slap\lambda$, $\slap\lambdab$ and recall the 
definitions \ref{11.98}. Lemma 11.2 implies, in view of the last of \ref{8.50}, comparing with 
\ref{11.31} and its conjugate and with \ref{11.125}, that for an arbitrary function $f$ on ${\cal N}$:
\begin{equation}
[T,\slap]f=-\s^{(T)}\sspi^{\sharp\sharp}\cdot\sD^2 f-\mbox{tr}\s^{(T)}\spi_1\cdot\sd f
\label{11.148}
\end{equation}
where:
\begin{eqnarray}
&&\hspace{3mm}\s^{(T)}\spi^C_{1,AB}=\frac{1}{2}\left(\sD_A\s^{(T)}\sspi^C_B+\sD_B\s^{(T)}\sspi^C_A
-\sD^C\s^{(T)}\sspi_{AB}\right) \nonumber\\ 
&&(\mbox{tr}\s^{(T)}\spi_1)^A=\sD^B\s^{(T)}\sspi_{AB}-\frac{1}{2}\sD^A\mbox{tr}\s^{(T)}\sspi
\label{11.149}
\end{eqnarray}
Similarly we deduce, in view of \ref{11.125}, that for an arbitrary function $f$ on ${\cal N}$:
\begin{equation}
[E_{(\mu)},\slap]f=-\s^{(E_{(\mu)})}\sspi^{\sharp\sharp}\cdot\sD^2 f-\mbox{tr}\s^{(E_{(\mu)})}\spi_1\cdot\sd f
\label{11.150}
\end{equation}
where:
\begin{equation}
(\mbox{tr}\s^{(E_{(\mu)})}\spi_1)^A=\sD^B\s^{(E_{(\mu)})}\sspi_{AB}-\frac{1}{2}\sD^A\mbox{tr}\s^{(E_{(\mu)})}\sspi
\label{11.151}
\end{equation}
Using the above we obtain:
\begin{eqnarray}
&&\slap E_{(\nu_l)} ... E_{(\nu_1)}T^{m-1}\lambda= E_{(\nu_l)} ... E_{(\nu_1)}T^{m-1}\slap\lambda
+\s^{(\nu_1 ... \nu_l)}c_{m-1,l}\nonumber\\
&&\slap E_{(\nu_l)} ... E_{(\nu_1)}T^{m-1}\lambdab= E_{(\nu_l)} ... E_{(\nu_1)}T^{m-1}\slap\lambdab
+\s^{(\nu_1 ... \nu_l)}\underline{c}_{m-1,l} \nonumber\\
&&\label{11.152}
\end{eqnarray}
where the remainders $\s^{(\nu_1 ... \nu_l)}c_{m-1,l}$, $\s^{(\nu_1 ... \nu_l)}\underline{c}_{m-1,l}$ 
are of order $m+l$ only. Similar formulas hold for the corresponding $N$th approximants.  Taking 
differences we obtain, in view of the definitions \ref{11.a2},
\begin{eqnarray}
&&\left[\slap\s^{(\nu_1 ... \nu_l)}\cla_{m-1,l}\right]_{P.P.}
=\nonumber\\
&&\hspace{15mm}\left[E_{(\nu_l)} ... E_{(\nu_1)}T^{m-1}\slap\lambda
-E_{(\nu_l,N)} ... E_{(\nu_1,N)}T^{m-1}\slap_N\lambda_N\right]_{P.P.}
\nonumber\\
&&\left[\slap\s^{(\nu_1 ... \nu_l)}\clab_{m-1,l}\right]_{P.P.}
=\nonumber\\
&&\hspace{15mm}\left[E_{(\nu_l)} ... E_{(\nu_1)}T^{m-1}\slap\lambdab
-E_{(\nu_l,N)} ... E_{(\nu_1,N)}T^{m-1}\slap_N\lambdab_N\right]_{P.P.}
\nonumber\\
&&\label{11.153}
\end{eqnarray}
therefore, in view of the definitions \ref{11.98} and \ref{11.104}:
\begin{eqnarray}
&&\lambda\slap\s^{(\nu_1 ... \nu_l)}\cla_{m-1,l}
=\nonumber\\
&&\hspace{15mm}\s^{(\nu_1 ... \nu_l)}\cnu_{m-1,l+1}+E_{(\nu_l)} ... E_{(\nu_1)}T^{m-1}j
-E_{(\nu_l,N)} ... E_{(\nu_1,N)}T^{m-1}j_N\nonumber\\
&&\hspace{40mm}+\lambda\s^{(\nu_1 ... \nu_l)}\check{n}_{m-1,l+1}\nonumber\\
&&\lambdab\slap\s^{(\nu_1 ... \nu_l)}\clab_{m-1,l}
=\nonumber\\
&&\hspace{15mm}\s^{(\nu_1 ... \nu_l)}\cnub_{m-1,l+1}+E_{(\nu_l)} ... E_{(\nu_1)}T^{m-1}\jb
-E_{(\nu_l,N)} ... E_{(\nu_1,N)}T^{m-1}\jb_N\nonumber\\
&&\hspace{40mm}+\lambdab\s^{(\nu_1 ... \nu_l)}\check{\nb}_{m-1,l+1}\nonumber\\
&&\label{11.154}
\end{eqnarray}
where the difference quantities $\s^{(\nu_1 ... \nu_l)}\check{n}_{m-1,l+1}$, 
$\s^{(\nu_1 ... \nu_l)}\check{\nb}_{m-1,l+1}$ are of order $m+l$ only. We are now in a position to 
apply to each of \ref{11.154} the following lemma. 

\vspace{2.5mm}

\noindent{\bf Lemma 11.4} \ \ Let $({\cal M},g)$ be a compact Riemannian manifold and 
let $\phi$ be a function on ${\cal M}$ satisfying the equation:
$$\triangle_g\phi=\rho$$
Let also $\mu$ be an arbitrary 
non-negative weight function on ${\cal M}$. Then the following integral inequality holds:
$$\int_{\cal M}\mu^2\left\{\frac{1}{2}|\nabla^2 \phi|^2+\mbox{Ric}((d\phi)^\sharp,(d\phi)^\sharp)
\right\}d\mu_g
\leq\int_{\cal M}\left\{2\mu^2\rho^2+3|d\mu|^2|\theta|^2\right\}d\mu_g$$

\vspace{2.5mm}

\noindent{\em Proof:} We have:
\begin{equation}
|\nabla^2\phi|^2=(\nabla_i\nabla_j\phi)(\nabla^i\nabla^j\phi)
=\nabla_i(\nabla_j\phi\nabla^i\nabla^j\phi)-\nabla^j\phi\triangle_g\nabla_j\phi
\label{11.155}
\end{equation}
In regard to the 2nd term on the right we have (see \ref{11.45}):
$$\triangle_g d\phi-d\triangle_g\phi=\mbox{Ric}\cdot(d\phi)^\sharp$$
Therefore, by virtue of the equation satisfied by $\phi$, the 2nd term on the right in \ref{11.155} is:
\begin{eqnarray}
&&-\nabla^i\phi\nabla_i\rho-\mbox{Ric}_{ij}\nabla^i\phi\nabla^j\phi\nonumber\\
&&=-\nabla_i(\rho\nabla_i\phi)+\rho\triangle_g\phi-\mbox{Ric}_{ij}\nabla^i\phi\nabla^j\phi 
\label{11.156}
\end{eqnarray}
and the 2nd term in the last expression is $\rho^2$. Substituting then in \ref{11.155} we arrive at 
the identity:
\begin{equation}
|\nabla^2\phi|^2=\mbox{div}I+\rho^2-\mbox{Ric}((d\phi)^\sharp,(d\phi)^\sharp)
\label{11.157}
\end{equation}
where $I$ is the vectorfield with components:
\begin{equation}
I^i=\nabla_j\phi\nabla^i\nabla^j\phi-\rho\nabla^i\phi
\label{11.158}
\end{equation}
We multiply the identity \ref{11.157} by $\mu^2$ and integrate over ${\cal M}$. 
Writing 
$$\mu^2\mbox{div}I=\mbox{div}(\mu^2 I)-2\mu I\cdot d\mu$$
The integral over ${\cal M}$ of $\mbox{div}(\mu^2 I)$ vanishes ${\cal M}$ being a compact 
manifold without boundary. Hence:
\begin{equation}
\int_{{\cal M}}\mu^2\mbox{div}I d\mu_g=-2\int_{\cal M}\mu(I\cdot d\mu)d\mu_g
\label{11.159}
\end{equation}
Moreover, in view of \ref{11.158} we can bound:
\begin{equation}
|I|\leq |d\phi|\left(|\nabla^2\phi|+|\rho|\right)
\label{11.160}
\end{equation}
therefore:
\begin{eqnarray}
&&-2\int_{\cal M}\mu(I\cdot d\mu)d\mu_g\leq\int_{{\cal M}}\mu|I||d\mu|d\mu_g\nonumber\\
&&\hspace{30mm}\leq 2\int_{{\cal M}}|d\mu||d\phi|\left(\mu|\nabla^2\phi|+\mu|\rho|\right)d\mu_g
\label{11.161}
\end{eqnarray}
Applying the inequalities
\begin{eqnarray*}
&&2|d\mu||d\phi|\mu|\nabla^2\phi|\leq\frac{1}{2}\mu^2|\nabla^2\phi|^2+2|d\mu|^2|d\phi|^2\\
&&2|d\mu||d\phi|\mu|\rho|\leq\mu^2|\rho|^2+|d\mu|^2|\rho|^2
\end{eqnarray*}
we then obtain:
\begin{eqnarray}
&&-2\int_{\cal M}\mu(I\cdot d\mu)d\mu_g\leq \label{11.162}\\
&&\hspace{15mm}\int_{{\cal M}}\left\{\mu^2\left(\frac{1}{2}|\nabla^2\phi|^2
+|\rho|^2\right)+3|d\mu|^2|d\phi|^2\right\}d\mu_g \nonumber
\end{eqnarray}
In view of \ref{11.159}, \ref{11.162}, integrating $\mu^2$ times \ref{11.157} over ${\cal M}$ 
yields the conclusion of the lemma.

\vspace{2.5mm}

We apply Lemma 11.4 with the Riemannian manifold $(S_{\ub,u},\sh)$ in the role of $({\cal M},g)$. 
The curvature term is controlled through the Gauss equations of Proposition 3.7. In connection 
with the 1st of \ref{11.154} we take $\s^{(\nu_1 ... \nu_l)}\cla_{m-1,l}$ in the role of $\phi$ 
and $\lambda$ in the role of $\mu$. In connection with the 2nd of \ref{11.154} we take 
$\s^{(\nu_1 ... \nu_l)}\clab_{m-1,l}$ in the role of $\phi$ and $\lambdab$ in the role of $\mu$. 
This yields the desired estimates of 
\ref{11.119}, \ref{11.120} in terms of \ref{11.121}, \ref{11.122} respectively. 

\pagebreak

\chapter{The Top Order Estimates for the Transformation Functions and the Next to Top Order 
Acoustical Estimates}

\section{The Propagation Equations for the Next to Top Order Acoustical Difference Quantities 
$(\s^{(n-1)}\ctchi, \s^{(n-1)}\ctchib)$ and $(\s^{(m,n-m)}\cla, \s^{(m,n-m)}\clab)$ : $m=0,...,n$}

Recall the transformation functions $f$, $w$, $\psi$ introduced in Chapter 4, Section 4.2. These are 
functions on ${\cal K}$ which to each point $(\tau,\vartheta)$ on ${\cal K}$ associate a point 
$(t,u^\prime,\vartheta^\prime)$ in the domain of the prior solution by:
$$t=f(\tau,\vartheta), \ \ \ u^\prime=w(\tau,\vartheta), \ \ \ \vartheta^\prime=\psi(\tau,\vartheta)$$
(see \ref{4.53} - \ref{4.57}). The transformation functions are subject to the identification 
equations \ref{4.57} which express the condition that the rectangular coordinates of identified 
points coincide. The regularized form of the identification equations is given in Section 4.3, 
Proposition 4.5 in terms of the re-scaled transformation functions $\hat{f}$, $v$, $\gamma$, 
related to the original functions by:
$$f(\tau,\vartheta)=f(0,\vartheta)+\tau^2\hat{f}(\tau,\vartheta), \ \ \ 
w(\tau,\vartheta)=\tau v(\tau,\vartheta), \ \ \ \psi(\tau,\vartheta)=\vartheta+\tau^3\gamma(\tau,\vartheta)$$
The corresponding $N$th approximants $\hat{f}_N$, $v_N$, $\gamma_N$ where introduced in Chapter 9, 
Sections 9.3 (see \ref{9.51} - \ref{9.54}) and 9.4 (see \ref{9.114} - \ref{9.116}). In Section 9.3 
estimates where derived for the quantities by which the jumps \ref{9.55} defined through these $N$th approximants fail to satisfy the boundary conditions (see Propositions 9.4, 9.5, 9.6), and in Section 9.4 estimates where derived for the quantities by which the $N$th approximants fail to satisfy the 
regularized identification equations (see Proposition 9.7). 

Let us define the differences:
\begin{equation}
\chf=\hat{f}-\hat{f}_N, \ \ \ \cv=v-v_N, \ \ \ \cga=\gamma-\gamma_N
\label{12.1}
\end{equation}
The main objective of the present chapter is to derive estimates for the top order derivatives of 
these transformation function differences. Since (see \ref{9.2}, \ref{9.7})
\begin{equation}
T=\frac{\partial}{\partial\tau}, \ \ \ \Omega=\frac{\partial}{\partial\vartheta}
\label{12.2}
\end{equation}
this means estimates for: 
$$(T^{m+1}\Omega^{n-m}\chf, T^{m+1}\Omega^{n-m}\cv, T^{m+1}\Omega^{n-m}\cga) \ : \ m=0,...,n$$
as well as for:
$$(\Omega^{n+1}\chf, \Omega^{n+1}\cv, \Omega^{n+1}\cga)$$
As we shall see, the derivation of these estimates necessarily involves the derivation of approriate 
estimates for the next to top order acoustical difference quantities, that is for:
\begin{equation}
\s^{(l-1)}\ctchi=E^{l-1}\tchi-E_N^{l-1}\tchi_N, \ \ \ \s^{(l-1)}\ctchib=E^{l-1}\tchib-E_N^{l-1}\tchib_N
 \ \ \ : \ \mbox{for $l=n$}
\label{12.3}
\end{equation}
and for:
\begin{eqnarray}
&&\s^{(m,l)}\cla=E^l T^m\lambda-E_N^l T^m\lambda_N, \ \ \ 
\s^{(m,l)}\clab=E^l T^m\lambdab-E_N^l T^m\lambdab_N \nonumber\\
&& \hspace{20mm} : \ \mbox{for $m=0,...,n$; \ $m+l=n$}
\label{12.4}
\end{eqnarray}
The derivation of these next to top order acoustical estimates is the other objective of the present 
chapter. 

In this first section we shall derive the propagation equations satisfied by the quantities 
\ref{12.3} and \ref{12.4}. In regard to the quantities \ref{12.3} we begin with the 2nd variation 
equations \ref{10.175}, \ref{10.176} and the corresponding equations \ref{10.187}, \ref{10.190} 
satisfied by the $N$th approximants. 

We shall first bring the 2nd variation equation \ref{10.175} to a suitable form. In the following 
we shall denote by $R$ a generic term of order 1 with vanishing order 1 acoustical part. 
By the 1st of \ref{3.a21} we have:
$$L\tchi=L\sk+L(H\sbeta\ss_N)$$
Using the facts, from Chapter 3, that $[LN^\mu]_{P.A.}=0$, while:
\begin{equation}
[LE^\mu]_{P.A.}=[f]_{P.A.}N^\mu=(c^{-1} E\lambdab+\pi\rho\tchib)N^\mu 
\label{12.5}
\end{equation}
(see \ref{3.a33}, \ref{3.a35} - \ref{3.a37}), we deduce: 
\begin{equation}
L(H\sbeta\ss_N)=H\sbeta N^\mu LE\beta_\mu+H(c^{-1}E\lambdab+\pi\rho\tchib)\beta_N\ss_N+R
\label{12.6}
\end{equation}
The principal term on the right hand side of \ref{10.175} is the term $E\sm$, and $\sm$ is given by 
\ref{10.6}. Using the facts, from Chapter 3, that:
\begin{eqnarray}
&&[EN^\mu]_{P.A.}=\tchi(E^\mu-\pi N^\mu) \label{12.7}\\
&&[EE^\mu]_{P.A.}=\frac{1}{2c}(\tchib N^\mu+\tchi\Nb^\mu) \label{12.8}
\end{eqnarray}
(see \ref{3.a31}, \ref{3.a32}) and, since $\rho=c^{-1}\lambdab$ and by \ref{3.a23}
\begin{equation}
c^{-1}[Ec]_{P.A.}=[k+\okb]_{P.A.}=-\pi(\tchi+\tchib),
\label{12.9}
\end{equation}
that:
\begin{equation}
[E\rho]_{P.A.}=c^{-1}E\lambdab+\pi\rho(\tchi+\tchib)
\label{12.10}
\end{equation}
we deduce:
\begin{eqnarray}
&&E\sm=-\beta_N\sbeta ELH+\frac{1}{2}\rho\beta_N^2 E^2 H-H\sbeta N^\mu EL\beta_\mu \nonumber\\
&&\hspace{10mm}+\frac{1}{2c}\beta_N^2(EH)E\lambdab+\frac{1}{2}\beta_N\left[-\frac{1}{c}
(\beta_N LH+Hs_{NL})+\pi\rho\beta_N EH\right]\tchib \nonumber\\
&&\hspace{10mm}+\left\{\left(-\sbeta^2+\pi\sbeta\beta_N-\frac{1}{2c}\beta_N\beta_{\Nb}\right)LH
+\rho\beta_N\left(\sbeta-\frac{1}{2}\pi\beta_N\right)EH\right.\nonumber\\
&&\hspace{20mm}\left.+H\left(\pi\sbeta-\frac{1}{2c}\beta_{\Nb}\right)s_{NL}
-H\rho\sbeta\ss_N\right\}\tchi+R
\label{12.11}
\end{eqnarray}
Using the results of Chapter 3 we also find that the remainder of the right hand side of \ref{10.175}, 
after the principal term $E\sm$ has been subtracted, is given by:
\begin{eqnarray}
&&-\frac{\beta_N}{2c}(\beta_N EH+2H\ss_N)E\lambdab\nonumber\\
&&+\frac{\beta_N}{2}\left[\frac{1}{2c}(\beta_N LH+2Hs_{NL})-\pi\rho(\beta_N EH+2H\ss_N)\right]\tchib
\nonumber\\
&&+\left\{m-\pi\sbeta\beta_N LH+\frac{1}{2}\pi\rho\beta_N^2 EH-\pi H\sbeta s_{NL}
+2\rho H\sbeta\ss_N\right\}\tchi-\rho\tchi^2\nonumber\\
&&+R
\label{12.12}
\end{eqnarray}
Adding \ref{12.6}, \ref{12.11}, \ref{12.12} and noting that by the first of the commutation 
relations \ref{3.a14} we have:
$$N^\mu (LE\beta_\mu-EL\beta_\mu)=-\chi\ss_N,$$ 
and also recalling the definition of the function $p$ from the statement of Proposition 3.3, 
we bring the 2nd variation equation \ref{10.175} to the desired form:
\begin{equation}
L\tchi=A\tchi-\rho\tchi^2+\oA\tchib-\beta_N\sbeta ELH+\frac{1}{2}\rho\beta_N^2 E^2 H+R
\label{12.13}
\end{equation}
where:
\begin{eqnarray}
&&A=p-\sbeta^2 LH+\rho\beta_N\sbeta EH \label{12.14}\\
&&\oA=-\frac{\beta_N^2}{4c}LH \label{12.15}
\end{eqnarray}
Note that the terms in $E\lambdab$ cancel. 

Similarly, using \ref{10.142}, \ref{10.184}, the 1st of \ref{10.117}, \ref{10.182}, \ref{10.148} 
together with \ref{10.145}, \ref{10.146}, and noting also that by the 1st of the commutation relations 
\ref{9.b1} we have 
$$N_N^\mu(L_N E_N\beta_{\mu,N}-E_N L_N\beta_{\mu,N})=-\chi_N\ss_{N,N}$$
while $\chi_N$ is related to $\tchi_N$ by \ref{10.138}, we bring the $N$th approximant 
2nd variation equation \ref{10.187} to the analogous form: 
\begin{eqnarray}
&&L_N\tchi_N=A_N\tchi_N-\rho_N\tchi_N^2+\oA_N\tchib_N-\beta_{N,N}\sbeta_N E_N L_N H_N
+\frac{1}{2}\rho_N\beta_N^2 E_N^2 H_N\nonumber\\
&&\hspace{14mm}+R_N+\varepsilon_{L\tchi,N}\label{12.16}
\end{eqnarray}
where:
\begin{eqnarray}
&&A_N=p_N-\sbeta_N^2 L_N H_N+\rho_N\beta_{N,N}\sbeta_N E_N H_N \label{12.17}\\
&&\oA_N=-\frac{\beta_{N,N}^2}{4c_N}L_N H_N\label{12.18}
\end{eqnarray}
$R_N$ is the $N$th approximant version of $R$ in \ref{12.13}, and:
\begin{equation}
\varepsilon_{L\tchi,N}=\varepsilon_{L\sk,N}-H_N\sbeta_N\varepsilon_{\chi,N}=O(\tau^N) 
\label{12.19}
\end{equation}
by \ref{10.140} and \ref{10.189}. 

In a similar manner we bring the conjugate 2nd variation equation, \ref{10.176} to the form:
\begin{equation}
\Lb\tchib=\oAb\tchib-\rhob\tchib^2+\Ab\tchi-\beta_{\Nb}\sbeta E\Lb H
+\frac{1}{2}\rhob\beta_{\Nb}^2 E^2 H+\Rb
\label{12.20}
\end{equation}
where:
\begin{eqnarray}
&&\oAb=\pb-\sbeta^2\Lb H+\rhob\beta_{\Nb}\sbeta EH \label{12.21}\\
&&\Ab=-\frac{\beta_{\Nb}^2}{4c}\Lb H \label{12.22}
\end{eqnarray} 
and $\Rb$ is a quantity of order 1 with vanishing order 1 acoustical part. 
Equation \ref{12.20} is the conjugate of equation \ref{12.13}. We also similarly bring 
the $N$th approximant conjugate 2nd variation equation \ref{10.190} to the form conjugate 
to \ref{12.16}: 
\begin{eqnarray}
&&\Lb_N\tchib_N=\oAb_N\tchib_N-\rhob_N\tchib_N^2+\Ab_N\tchi_N-\beta_{\Nb,N}\sbeta_N E_N\Lb_N H_N
+\frac{1}{2}\rhob_N\beta_{\Nb}^2 E_N^2 H_N\nonumber\\
&&\hspace{14mm}+\Rb_N+\varepsilon_{\Lb\tchib,N}\label{12.23}
\end{eqnarray}
where:
\begin{eqnarray}
&&\oAb_N=\pb_N-\sbeta_N^2\Lb_N H_N+\rhob_N\beta_{\Nb,N}\sbeta_N E_N H_N \label{12.24}\\
&&\Ab_N=-\frac{\beta_{\Nb,N}^2}{4c_N}\Lb_N H_N\label{12.25}
\end{eqnarray}
$\Rb_N$ is the $N$th approximant version of $\Rb$ in \ref{12.20}, and:
\begin{equation}
\varepsilon_{\Lb\tchib,N}=\varepsilon_{\Lb\skb,N}-H_N\sbeta_N\varepsilon_{\chib,N}=O(\tau^{N+1}) 
\label{12.26}
\end{equation}
by \ref{10.144} and \ref{10.192}. 

Applying $E^{n-1}$ to equation \ref{12.13} and using the 1st of the commutation relations \ref{3.a14}, 
which implies that: 
$$[L,E^{n-1}]\tchi=-\sum_{i=0}^{n-2}E^{n-2-i}\left(\chi E^{i+1}\tchi\right)$$
we deduce the following propagation equation for $E^{n-1}\tchi$:
\begin{eqnarray}
&&LE^{n-1}\tchi=\left(A-(n+1)\chi\right)E^{n-1}\tchi+\oA E^{n-1}\tchib\nonumber\\
&&\hspace{17mm}-\beta_N\sbeta E^n LH+\frac{1}{2}\rho\beta_N^2 E^{n+1}H+R_{n-1} \label{12.27}
\end{eqnarray}
where $R_{n-1}$ is of order $n$, but the highest order acoustical quantity occurring in $R_{n-1}$ is 
only of order $n-1$. Similarly, applying $E_N^{n-1}$ to equation \ref{12.16} and using the 1st of the 
commutation relations \ref{9.b1}, which implies that:
$$[L_N,E_N^{n-1}]\tchi_N=-\sum_{i=0}^{n-2}E_N^{n-2-i}\left(\chi_N E_N^{i+1}\tchi_N\right)$$
we deduce:
\begin{eqnarray}
&&L_N E_N^{n-1}\tchi_N=\left(A_N-(n+1)\chi_N\right)E_N^{n-1}\tchi_N+\oA_N E_N^{n-1}\tchib_N\nonumber\\
&&\hspace{21mm}-\beta_{N,N}\sbeta_N E_N^n L_N H_N+\frac{1}{2}\rho_N\beta_{N,N}^2 E_N^{n+1}H_N 
+R_{n-1,N}\nonumber\\
&&\hspace{21mm}+E_N^{n-1}\vep_{L\tchi,N}
\label{12.28}
\end{eqnarray}
where $R_{n-1,N}$ is the $N$th approximant version of $R_{n-1}$. Subtracting equation \ref{12.28} from 
\ref{12.27} we arrive at the following propagation equation for the acoustical difference quantity 
$\s^{(n-1)}\ctchi$ defined by the 1st of \ref{12.3}:
\begin{eqnarray}
&&L\s^{(n-1)}\ctchi+(n+1)\chi\s^{(n-1)}\ctchi=A\s^{(n-1)}\ctchi+\oA\s^{(n-1)}\ctchib\nonumber\\
&&\hspace{5mm}-\left(L-L_N-A+A_N+(n+1)(\chi-\chi_N)\right)E_N^{n-1}\tchi_N 
+(\oA-\oA_N)E_N^{n-1}\tchib_N\nonumber\\
&&\hspace{5mm}-\beta_N\sbeta E^n LH+\beta_{N,N}\sbeta_N E_N^n L_N H_N 
+\frac{1}{2}\rho\beta_N^2 E^{n+1}H-\frac{1}{2}\rho_N\beta_{N,N}^2 E_N^{n+1}H_N\nonumber\\
&&\hspace{35mm}+R_{n-1}-R_{n-1,N}-E_N^{n-1}\vep_{L\tchi,N}\label{12.29}
\end{eqnarray}

Applying likewise $E^{n-1}$ to equation \ref{12.20} and using the 2nd of the commutation relations \ref{3.a14}, 
which implies that: 
$$[\Lb,E^{n-1}]\tchib=-\sum_{i=0}^{n-2}E^{n-2-i}\left(\chib E^{i+1}\tchib\right)$$
we deduce the following propagation equation for $E^{n-1}\tchib$ (the conjugate of \ref{12.27}):
\begin{eqnarray}
&&\Lb E^{n-1}\tchib=\left(\oAb-(n+1)\chib\right)E^{n-1}\tchib+\Ab E^{n-1}\tchi\nonumber\\
&&\hspace{17mm}-\beta_{\Nb}\sbeta E^n\Lb H+\frac{1}{2}\rhob\beta_{\Nb}^2 E^{n+1}H+\Rb_{n-1} \label{12.30}
\end{eqnarray}
where $\Rb_{n-1}$ is of order $n$, but the highest order acoustical quantity 
occurring in $\Rb_{n-1}$ is only of order $n-1$. 
Similarly, applying $E_N^{n-1}$ to equation \ref{12.23} and using the 2nd of the 
commutation relations \ref{9.b1}, which implies that:
$$[\Lb_N,E_N^{n-1}]\tchib_N=-\sum_{i=0}^{n-2}E_N^{n-2-i}\left(\chib_N E_N^{i+1}\tchib_N\right)$$
we deduce:
\begin{eqnarray}
&&\Lb_N E_N^{n-1}\tchib_N=\left(\oAb_N-(n+1)\chib_N\right)E_N^{n-1}\tchib_N+\Ab_N E_N^{n-1}\tchi_N\nonumber\\
&&\hspace{21mm}-\beta_{\Nb,N}\sbeta_N E_N^n\Lb_N H_N+\frac{1}{2}\rhob_N\beta_{\Nb,N}^2 E_N^{n+1}H_N 
+\Rb_{n-1,N}\nonumber\\
&&\hspace{21mm}+E_N^{n-1}\vep_{\Lb\tchib,N}
\label{12.31}
\end{eqnarray}
where $\Rb_{n-1,N}$ is the $N$th approximant version of $\Rb_{n-1}$. 
Subtracting equation \ref{12.31} from 
\ref{12.30} we arrive at the following propagation equation for the acoustical difference quantity 
$\s^{(n-1)}\ctchib$ defined by the 2nd of \ref{12.3}:
\begin{eqnarray}
&&\Lb\s^{(n-1)}\ctchib+(n+1)\chib\s^{(n-1)}\ctchib=\oAb\s^{(n-1)}\ctchib
+\Ab\s^{(n-1)}\ctchi\nonumber\\
&&\hspace{5mm}-\left(\Lb-\Lb_N-\oAb+\oAb_N+(n+1)(\chib-\chib_N)\right)E_N^{n-1}\tchib_N 
+(\Ab-\Ab_N)E_N^{n-1}\tchi_N\nonumber\\
&&\hspace{5mm}-\beta_{\Nb}\sbeta E^n\Lb H+\beta_{\Nb,N}\sbeta_N E_N^n\Lb_N H_N 
+\frac{1}{2}\rhob\beta_{\Nb}^2 E^{n+1}H-\frac{1}{2}\rhob_N\beta_{\Nb,N}^2 E_N^{n+1}H_N\nonumber\\
&&\hspace{35mm}+\Rb_{n-1}-\Rb_{n-1,N}-E_N^{n-1}\vep_{\Lb\tchib,N}\label{12.32}
\end{eqnarray}

We now turn to the quantities \ref{12.4}. In regard to these quantities we must first deduce from 
the propagation equations of Proposition 3.3 propagation equations for $E\lambda$, $E\lambdab$ and 
for $T\lambda$, $T\lambdab$. 

Applying $E$ to the propagation equation for $\lambda$ of Proposition 3.3 and using the first of the 
commutation relations \ref{3.a14} we obtain:
\begin{equation}
LE\lambda+\chi E\lambda=p E\lambda+q E\lambdab+\lambda Ep+\lambdab Eq
\label{12.33}
\end{equation}
The principal terms on the right are contained in $Ep$, $Eq$. The quantities $p$ and $q$ being of 
1st order with vanishing 1st order acoustical part, $Ep$ and $Eq$ are of 2nd order with vanishing 
2nd order acoustical part. The principal part of $Eq$ is:
\begin{equation}
[Eq]_{P.P.}=\frac{1}{4c}\beta_N^2 E\Lb H
\label{12.34}
\end{equation}
In regard to $Ep$ we first consider $Em$. By \ref{3.47}, \ref{3.57}, \ref{3.68} :
\begin{eqnarray}
&&m=-\pi\sm-\om\nonumber\\
&&\hspace{4mm}=\pi (\beta_N\sbeta LH-\frac{1}{2}\rho\beta_N^2 EH +H\sbeta s_{NL})\nonumber\\
&&\hspace{7mm}-\frac{1}{4c}(\beta_N^2 LH+2H\beta_N s_{NL}) \label{12.35}
\end{eqnarray}
The principal part of $Em$ is:
\begin{eqnarray}
&&[Em]_{P.P.}=\beta_N\left(\pi\sbeta-\frac{\beta_N}{4c}\right)ELH-\frac{1}{2}\pi\rho\beta_N^2 E^2 H 
\nonumber\\
&&\hspace{17mm}+H\left(\pi\sbeta-\frac{\beta_N}{2c}\right)N^\mu EL\beta_\mu \label{12.36}
\end{eqnarray}
The difference $p-m$ being given by:
\begin{equation}
p-m=-\frac{1}{2c}(\beta_N\beta_{\Nb}LH+H\beta_{\Nb}s_{NL})
\label{12.37}
\end{equation}
the principal part of $Ep-Em$ is:
\begin{equation}
[Ep-Em]_{P.P.}=-\frac{1}{2c}(\beta_N\beta_{\Nb}ELH+\beta_{\Nb}N^\mu EL\beta_\mu)
\label{12.38}
\end{equation}
The principal part of $Ep$ is then:
\begin{eqnarray}
&&[Ep]_{P.P.}=\beta_N\left(\pi\sbeta-\frac{\beta_N}{4c}-\frac{\beta_{\Nb}}{2c}\right)ELH
-\frac{1}{2}\pi\rho\beta_N^2 E^2 H \nonumber\\
&&\hspace{16mm}+H\left(\pi\sbeta-\frac{(\beta_N+\beta_{\Nb})}{2c}\right)N^\mu EL\beta_\mu \label{12.39}
\end{eqnarray}
We must also determine the principal acoustical parts of the remainders 
$$Eq-[Eq]_{P.P.}, \ \ \ Ep-[Ep]_{P.P.}$$
that is the terms in these remainders containing acoustical quantities of order 1. Using \ref{12.7}, 
which implies:
\begin{equation}
[E\beta_N]_{P.A.}=(\sbeta-\pi\beta_N)\tchi
\label{12.40}
\end{equation}
and \ref{12.9}, we find that the principal acoustical part of the remainder $Eq-[Eq]_{P.P.}$ is:
\begin{equation}
\frac{\beta_N}{2c}(\Lb H)\left[\left(\sbeta-\frac{\pi}{2}\beta_N\right)\tchi
+\frac{\pi}{2}\beta_N\tchib\right]
\label{12.41}
\end{equation}
To determine the principal acoustical part of the remainder $Ep-[Ep]_{P.P.}$ we first determine 
the principal acoustical part of the remainder $Em-[Em]_{P.P.}$. Using \ref{10.273}, together with 
\ref{12.7}, which besides \ref{12.40} also implies:
\begin{equation}
[EN^\mu]_{P.A.}L\beta_\mu=(\rho\ss_N-\pi s_{NL})\tchi,
\label{12.42}
\end{equation}
\ref{12.8}, which implies:
\begin{equation}
[E\sbeta]_{P.A.}=\frac{1}{2c}(\beta_N\tchib+\beta_{\Nb}\tchi),
\label{12.43}
\end{equation}
\ref{12.9} and \ref{12.10}, we find that the principal acoustical part of the 
remainder $Em-[Em]_{P.P.}$ is given by: 
\begin{eqnarray}
&&-\frac{\pi}{2c}\beta_N^2(EH)E\lambdab\nonumber\\
&&+\left[\frac{\beta_N}{2c}\left(\sbeta+\frac{\pi}{2}\beta_N\right)LH
-\frac{1}{2}\rho\beta_N^2\left(\frac{1}{2c}+\pi^2\right)EH+\frac{H\sbeta}{2c}s_{NL}\right]\tchib\nonumber\\
&&+\left\{\pi\left(\sbeta^2-\pi\sbeta\beta_N+\frac{\beta_N^2}{4c}+\frac{\beta_N\beta_{\Nb}}{2c}\right)
LH\right.\nonumber\\
&&\hspace{5mm}+\rho\beta_N\left[-\pi\sbeta+\frac{1}{2}\left(\pi^2-\frac{1}{2c}\right)\beta_N\right]EH
\nonumber\\
&&\hspace{5mm}\left.+\pi H\left(\frac{(\beta_N+\beta_{\Nb})}{2c}-\pi\sbeta\right)s_{NL}
+H\left(\pi\sbeta-\frac{\beta_N}{2c}\right)\rho\ss_N\right\}\tchi \label{12.44}
\end{eqnarray}
Using also the fact, from Chapter 3, that:
\begin{equation}
[E\Nb^\mu]_{P.A.}=\tchib(E^\mu-\pi\Nb^\mu)
\label{12.45}
\end{equation}
(conjugate of \ref{12.7}) which implies:
\begin{equation}
[E\beta_{\Nb}]_{P.A.}=(\sbeta-\pi\beta_{\Nb})\tchib
\label{12.46}
\end{equation}
(conjugate of \ref{12.40}), we find from \ref{12.37}, \ref{12.38} that the principal acoustical part of the remainder 
$Ep-Em-[Ep-Em]_{P.P.}$ is given by:
\begin{equation}
-\frac{\sbeta}{2c}(\beta_N LH+H s_{NL})\tchib
-\frac{\beta_{\Nb}}{2c}(\sbeta LH+H\rho\ss_N)\tchi
\label{12.47}
\end{equation}
Adding \ref{12.44} and \ref{12.47} we obtain the principal acoustical part of the remainder 
$Ep-[Ep]_{P.P.}$ as:
\begin{eqnarray}
&&-\frac{\pi}{2c}\beta_N^2(EH)E\lambdab\nonumber\\
&&+\left[\frac{\pi\beta_N^2}{4c}LH
-\frac{1}{2}\rho\beta_N^2\left(\frac{1}{2c}+\pi^2\right)EH\right]\tchib\nonumber\\
&&+\left\{\left(\pi\sbeta^2-\pi^2\sbeta\beta_N-\frac{\sbeta\beta_{\Nb}}{2c}
+\frac{\pi\beta_N^2}{4c}+\frac{\pi\beta_N\beta_{\Nb}}{2c}\right)LH\right.\nonumber\\
&&\hspace{5mm}+\rho\beta_N\left[-\pi\sbeta+\frac{1}{2}\left(\pi^2-\frac{1}{2c}\right)\beta_N\right]EH
\nonumber\\
&&\hspace{5mm}\left.+H\left(\frac{(\beta_N+\beta_{\Nb})}{2c}-\pi\sbeta\right)
(\pi s_{NL}-\rho\ss_N)\right\}\tchi \label{12.48}
\end{eqnarray}
Defining then 
\begin{equation}
\tilde{q}=\frac{\beta_N^2}{4c}(\Lb H-2\pi\lambda EH)
\label{12.49}
\end{equation}
(compare with \ref{10.262}) we conclude from the above that equation \ref{12.33} takes the form:
\begin{eqnarray}
&&LE\lambda=(p-\chi)E\lambda+\tilde{q}E\lambdab+B\tchi+\oB\tchib\nonumber\\
&&\hspace{10mm}+\frac{\rho}{4}\beta_N^2 E\Lb H
+\rhob\beta_N\left(c\pi\sbeta-\frac{\beta_N}{4}-\frac{\beta_{\Nb}}{2}\right)ELH\nonumber\\
&&\hspace{10mm}-\frac{\pi}{2}a\beta_N^2 E^2 H+\rhob H\left(c\pi\sbeta-\frac{(\beta_N+\beta_{\Nb})}{2}\right)N^\mu EL\beta_\mu\nonumber\\
&&\hspace{10mm}+S
\label{12.50}
\end{eqnarray}
Here:
\begin{eqnarray}
&&\oB=\beta_N^2\left[\frac{\pi\rho}{4}\Lb H+\frac{\pi\rhob}{4}LH
-\frac{1}{2}a\left(\pi^2+\frac{1}{2c}\right)EH\right]\label{12.51}\\
&&B=\frac{\rho}{2}\beta_N\left(\sbeta-\frac{\pi}{2}\beta_N\right)\Lb H\nonumber\\
&&\hspace{7mm}+\rhob\left(c\pi\sbeta^2-c\pi^2\sbeta\beta_N-\frac{1}{2}\sbeta\beta_{\Nb}
+\frac{\pi}{4}\beta_N^2+\frac{\pi}{2}\beta_N\beta_{\Nb}\right)LH\nonumber\\
&&\hspace{7mm}+a\beta_N\left[-\pi\sbeta+\frac{1}{2}\left(\pi^2-\frac{1}{2c}\right)\beta_N\right]EH
\nonumber\\
&&\hspace{7mm}+H\left(\frac{(\beta_N+\beta_{\Nb})}{2}-c\pi\sbeta\right)\left(\pi\rhob s_{NL}
-\frac{a}{c}\ss_N\right) \label{12.52}
\end{eqnarray}
and $S$ is of order 1, with no order 1 acoustical quantities occurring in $S$. 

Applying $E$ to the propagation equation for $\lambdab$ of Proposition 3.3 and using the second of the 
commutation relations \ref{3.a14} we derive in a similar manner the conjugate equation:
\begin{eqnarray}
&&\Lb E\lambdab=(\pb-\chib)E\lambdab+\tilde{\qb}E\lambda+\oBb\tchib+\Bb\tchi\nonumber\\
&&\hspace{10mm}+\frac{\rhob}{4}\beta_{\Nb}^2 ELH
+\rho\beta_{\Nb}\left(c\pi\sbeta-\frac{\beta_{\Nb}}{4}-\frac{\beta_N}{2}\right)E\Lb H\nonumber\\
&&\hspace{10mm}-\frac{\pi}{2}a\beta_{\Nb}^2 E^2 H+\rho H\left(c\pi\sbeta-\frac{(\beta_N+\beta_{\Nb})}{2}\right)\Nb^\mu E\Lb\beta_\mu\nonumber\\
&&\hspace{10mm}+\Sb
\label{12.53}
\end{eqnarray}
Here:
\begin{equation}
\tilde{\qb}=\frac{\beta_{\Nb}^2}{4c}(LH-2\pi\lambdab EH)
\label{12.54}
\end{equation}
(compare with \ref{10.285}) and: 
\begin{eqnarray}
&&\Bb=\beta_{\Nb}^2\left[\frac{\pi\rhob}{4}LH+\frac{\pi\rho}{4}\Lb H
-\frac{1}{2}a\left(\pi^2+\frac{1}{2c}\right)EH\right]\label{12.55}\\
&&\oBb=\frac{\rhob}{2}\beta_{\Nb}\left(\sbeta-\frac{\pi}{2}\beta_{\Nb}\right)LH\nonumber\\
&&\hspace{7mm}+\rho\left(c\pi\sbeta^2-c\pi^2\sbeta\beta_{\Nb}-\frac{1}{2}\sbeta\beta_N
+\frac{\pi}{4}\beta_{\Nb}^2+\frac{\pi}{2}\beta_{\Nb}\beta_N\right)\Lb H\nonumber\\
&&\hspace{7mm}+a\beta_{\Nb}\left[-\pi\sbeta+\frac{1}{2}\left(\pi^2-\frac{1}{2c}\right)\beta_{\Nb}\right]EH
\nonumber\\
&&\hspace{7mm}+H\left(\frac{(\beta_N+\beta_{\Nb})}{2}-c\pi\sbeta\right)\left(\pi\rho s_{\Nb\Lb}
-\frac{a}{c}\ss_{\Nb}\right) \label{12.56}
\end{eqnarray}
and $\Sb$ is of order 1, with no order 1 acoustical quantities occurring in $S$. 

Applying $T$ to the propagation equation for $\lambda$ of Proposition 3.3 and using the fact 
that by the third of the commutation relations \ref{3.a14} 
$$[L,T]=-\zeta E$$
we obtain:
\begin{equation}
LT\lambda+\zeta E\lambda=pT\lambda+qT\lambdab+\lambda Tp+\lambdab Tq
\label{12.57}
\end{equation}
The principal terms on the right are contained in $Tp$, $Tq$. The quantities $p$ and $q$ being of 
1st order with vanishing 1st order acoustical part, $Tp$ and $Tq$ are of 2nd order with vanishing 
2nd order acoustical part. The principal part of $Tq$ is:
\begin{equation}
[Tq]_{P.P.}=\frac{1}{4c}\beta_N^2 T\Lb H
\label{12.58}
\end{equation}
In regard to $Tp$ we first consider $Tm$. From \ref{12.35} the principal part of $Tm$ is:
\begin{eqnarray}
&&[Tm]_{P.P.}=\beta_N\left(\pi\sbeta-\frac{\beta_N}{4c}\right)TLH-\frac{1}{2}\pi\rho\beta_N^2 TEH 
\nonumber\\
&&\hspace{17mm}+H\left(\pi\sbeta-\frac{\beta_N}{2c}\right)N^\mu TL\beta_\mu \label{12.59}
\end{eqnarray}
The difference $p-m$ being given by \ref{12.37}, the principal part of $Tp-Tm$ is:
\begin{equation}
[Tp-Tm]_{P.P.}=-\frac{1}{2c}(\beta_N\beta_{\Nb}TLH+\beta_{\Nb}N^\mu TL\beta_\mu)
\label{12.60}
\end{equation}
The principal part of $Tp$ is then:
\begin{eqnarray}
&&[Tp]_{P.P.}=\beta_N\left(\pi\sbeta-\frac{\beta_N}{4c}-\frac{\beta_{\Nb}}{2c}\right)TLH
-\frac{1}{2}\pi\rho\beta_N^2 TEH \nonumber\\
&&\hspace{16mm}+H\left(\pi\sbeta-\frac{(\beta_N+\beta_{\Nb})}{2c}\right)N^\mu TL\beta_\mu \label{12.61}
\end{eqnarray}
We must also determine the principal acoustical parts of the remainders 
$$Tq-[Tq]_{P.P.}, \ \ \ Tp-[Tp]_{P.P.}$$
that is the terms in these remainders containing acoustical quantities of order 1. From the 3rd and 
5th of \ref{3.a15} we have:
\begin{equation}
[TN^\mu]_{P.A.}=[\Lb N^\mu]_{P.A}=[\sn]_{P.A.}(E^\mu-\pi N^\mu)
\label{12.62}
\end{equation}
(taking account of \ref{3.61} and the 1st of \ref{3.48}). Similarly, from the 4th and 6th of 
\ref{3.a15} we have:
\begin{equation}
[T\Nb^\mu]_{P.A.}=[L\Nb^\mu]_{P.A.}=[\snb]_{P.A.}(E^\mu-\pi\Nb^\mu)
\label{12.63}
\end{equation}
(taking account of \ref{3.62} and the 2nd of \ref{3.48}). 
Introducing the quantities:
\begin{equation}
\mu=E\lambda+\pi\lambda\tchi, \ \ \ \mub=E\lambdab+\pi\lambdab\tchib
\label{12.64}
\end{equation}
from \ref{3.a24} we have:
\begin{equation}
[\sn]_{P.A.}=2\mu, \ \ \ [\snb]_{P.A.}=2\mub
\label{12.65}
\end{equation}
Therefore \ref{12.62}, \ref{12.63} take the form:
\begin{equation}
[TN^\mu]_{P.A.}=2\mu(E^\mu-\pi N^\mu), \ \ \ [T\Nb^\mu]_{P.A.}=2\mub(E^\mu-\pi\Nb^\mu)
\label{12.66}
\end{equation}
These imply:
\begin{equation}
[T\beta_N]_{P.A.}=2\mu(\sbeta-\pi\beta_N), \ \ \ [T\beta_{\Nb}]_{P.A.}=2\mub(\sbeta-\pi\beta_{\Nb})
\label{12.67}
\end{equation}
and:
\begin{equation}
[TN^\mu]_{P.A.}L\beta_\mu=2\mu(\rho\ss_N-\pi s_{NL})
\label{12.68}
\end{equation}
Moreover by \ref{3.85}, \ref{3.86}, 
\begin{equation}
c^{-1}[Tc]_{P.A.}=[n+\onb]_{P.A.}=-\pi[\sn+\snb]_{P.A.}=-2\pi(\mu+\mub)
\label{12.69}
\end{equation}
hence, since $\rho=c^{-1}\lambdab$,
\begin{equation}
[T\rho]_{P.A.}=c^{-1}T\lambdab+2\pi\rho(\mu+\mub)
\label{12.70}
\end{equation}
Using \ref{12.67} and \ref{12.69} we find that the principal acoustical part of the remainder 
$Tq-[Tq]_{P.P.}$ is:
\begin{equation}
\frac{\beta_N}{c}(\Lb H)\left[\left(\sbeta-\frac{\pi}{2}\beta_N\right)\mu 
+\frac{\pi}{2}\beta_N\mub\right]
\label{12.71}
\end{equation}
To determine the principal acoustical part of the remainder $Tp-[Tp]_{P.P.}$ we first determine the 
principal acoustical part of the remainder $Tm-[Tm]_{P.P.}$. By \ref{3.a33}, \ref{3.a38}:
\begin{equation}
TE^\mu=(\sf+\sfb)E^\mu+(f+\fb)N^\mu+(\of+\ofb)\Nb^\mu 
\label{12.72}
\end{equation}
and by \ref{3.a35} - \ref{3.a37}, \ref{3.a39} - \ref{3.a41}:
\begin{eqnarray}
&&[\sf]_{P.A.}=0, \ \ \ [f]_{P.A.}=c^{-1}\mub, \ \ \ [\of]_{P.A.}=0 \nonumber\\
&&[\sfb]_{P.A.}=0, \ \ \ [\fb]_{P.A.}=0, \ \ \ [\ofb]_{P.A.}=c^{-1}\mu 
\label{12.73}
\end{eqnarray}
Hence:
\begin{equation}
[TE^\mu]_{P.A.}=c^{-1}\mub N^\mu+c^{-1}\mu\Nb^\mu 
\label{12.74}
\end{equation}
The 0-component of this is (see \ref{3.a19}):
\begin{equation}
[T\pi]_{P.A.}=c^{-1}(\mub+\mu)
\label{12.75}
\end{equation}
Also \ref{12.74} implies:
\begin{equation}
[T\sbeta]_{P.A.}=c^{-1}(\beta_N\mub+\beta_{\Nb}\mu)
\label{12.76}
\end{equation}
Using the 1st of \ref{12.67} and \ref{12.68}, \ref{12.69}, \ref{12.70} and \ref{12.75}, \ref{12.76} 
we find that the principal acoustical part of $Tm-[Tm]_{P.P.}$ is given by:
\begin{eqnarray}
&&-\frac{\pi}{2c}\beta_N^2(EH)T\lambdab \nonumber\\
&&+\left[\frac{\beta_N}{c}\left(\sbeta+\frac{\pi}{2}\beta_N\right)LH
-\rho\beta_N^2\left(\pi^2+\frac{1}{2c}\right)EH+\frac{H\sbeta}{c}s_{NL}\right]\mub\nonumber\\
&&+\left\{\pi\left(2\sbeta^2-2\pi\sbeta\beta_N+\frac{\beta_N^2}{2c}+\frac{\beta_N\beta_{\Nb}}{c}\right)
LH\right.\nonumber\\
&&\hspace{5mm}+\rho\beta_N\left[-2\pi\sbeta+\left(\pi^2-\frac{1}{2c}\right)\beta_N\right]EH \nonumber\\
&&\hspace{5mm}\left.+\pi H\left(\frac{(\beta_N+\beta_{\Nb})}{c}-2\pi\sbeta\right)s_{NL}
+H\left(2\pi\sbeta-\frac{\beta_N}{c}\right)\rho\ss_N\right\}\mu \nonumber\\
&&\label{12.77}
\end{eqnarray}
Using \ref{12.67}, \ref{12.68}, \ref{12.69} we find from \ref{12.37}, \ref{12.60} that the principal 
part of the remainder $Tp-Tm-[Tp-Tm]_{P.P.}$ is given by:
\begin{equation}
-\frac{\sbeta}{c}(\beta_N LH+Hs_{NL})\mub-\frac{\beta_{\Nb}}{c}(\sbeta LH+H\rho\ss_N)
\label{12.78}
\end{equation}
Adding \ref{12.77} and \ref{12.78} we obtain the principal part of the remainder $Tp-[Tp]_{P.P.}$ as:
\begin{eqnarray}
&&-\frac{\pi}{2c}\beta_N^2(EH)T\lambdab\nonumber\\
&&+\left[\frac{\pi\beta_N^2}{2c}LH
-\rho\beta_N^2\left(\frac{1}{2c}+\pi^2\right)EH\right]\mub\nonumber\\
&&+\left\{\left(2\pi\sbeta^2-2\pi^2\sbeta\beta_N-\frac{\sbeta\beta_{\Nb}}{c}
+\frac{\pi\beta_N^2}{2c}+\frac{\pi\beta_N\beta_{\Nb}}{c}\right)LH\right.\nonumber\\
&&\hspace{5mm}+\rho\beta_N\left[-2\pi\sbeta+\left(\pi^2-\frac{1}{2c}\right)\beta_N\right]EH
\nonumber\\
&&\hspace{5mm}\left.+H\left(\frac{(\beta_N+\beta_{\Nb})}{c}-2\pi\sbeta\right)
(\pi s_{NL}-\rho\ss_N)\right\}\mu \label{12.79}
\end{eqnarray}
We conclude from the above, in view of the definitions \ref{12.49}, \ref{12.51}, \ref{12.52}, 
that equation \ref{12.57} takes the form:
\begin{eqnarray}
&&LT\lambda=pT\lambda+\tilde{q}T\lambdab-\zeta E\lambda+2B\mu+2\oB\mub\nonumber\\
&&\hspace{10mm}+\frac{\rho}{4}\beta_N^2 T\Lb H
+\rhob\beta_N\left(c\pi\sbeta-\frac{\beta_N}{4}-\frac{\beta_{\Nb}}{2}\right)TLH\nonumber\\
&&\hspace{10mm}-\frac{\pi}{2}a\beta_N^2 TEH+\rhob H\left(c\pi\sbeta-\frac{(\beta_N+\beta_{\Nb})}{2}\right)N^\mu TL\beta_\mu\nonumber\\
&&\hspace{10mm}+Q
\label{12.80}
\end{eqnarray}
Here $Q$ is of order 1, with no order 1 acoustical quantities occurring in $Q$. 

Applying $T$ to the propagation equation for $\lambdab$ of Proposition 3.3  and using the fact that by the third 
of the commutation relations \ref{3.a14} 
$$[\Lb,T]=\zeta E$$ 
we derive in a similar manner the conjugate equation:
\begin{eqnarray}
&&\Lb T\lambdab=\pb T\lambdab+\tilde{\qb}T\lambda+\zeta E\lambdab+2\oBb\mub+2\Bb\mu\nonumber\\
&&\hspace{10mm}+\frac{\rhob}{4}\beta_{\Nb}^2 TLH
+\rho\beta_{\Nb}\left(c\pi\sbeta-\frac{\beta_{\Nb}}{4}-\frac{\beta_N}{2}\right)T\Lb H\nonumber\\
&&\hspace{10mm}-\frac{\pi}{2}a\beta_{\Nb}^2 TEH+\rho H\left(c\pi\sbeta-\frac{(\beta_N+\beta_{\Nb})}{2}\right)\Nb^\mu T\Lb\beta_\mu\nonumber\\
&&\hspace{10mm}+\Qb
\label{12.81}
\end{eqnarray}
where $\Qb$ is of order 1, with no order 1 acoustical quantities occurring in $\Qb$.
 
Applying $E_N$ to the propagation equation \ref{10.161} for the $N$th approximant quantity $\lambda_N$ and using 
the first of the commutation relations \ref{9.b1} we obtain the $N$th approximant version of equation \ref{12.33}:
\begin{equation}
L_N E_N\lambda_N+\chi_N E_N\lambda_N=p_N E_N\lambda_N+q_N E_N\lambdab_N +\lambda_N E_N p_N 
+\lambdab_N E_N q_N +E_N\vep_{\lambda,N}
\label{12.82}
\end{equation}
Proceeding from this equation following the path leading from equation \ref{12.33} to equation \ref{12.50}, using the results of Section 10.3 in place 
of the results of Chapter 3, we deduce the $N$th approximant version of equation \ref{12.50}: 
\begin{eqnarray}
&&L_N E_N\lambda_N=(p_N-\chi_N)E_N\lambda_N+\tilde{q}_N E_N\lambdab_N+B_N\tchi_N+\oB_N\tchib_N\nonumber\\
&&\hspace{20mm}+\frac{\rho_N}{4}\beta_{N,N}^2 E_N\Lb_N H_N \nonumber\\
&&\hspace{20mm}+\rhob_N\beta_{N,N}\left(c_N\pi_N\sbeta_N-\frac{\beta_{N,N}}{4}-\frac{\beta_{\Nb,N}}{2}\right)E_N L_N H_N\nonumber\\
&&\hspace{20mm}-\frac{\pi_N}{2}a_N\beta_{N,N}^2 E_N^2 H_N \nonumber\\
&&\hspace{20mm}+\rhob_N H_N\left(c_N\pi_N\sbeta_N-\frac{(\beta_{N,N}+\beta_{\Nb,N})}{2}\right)
N_N^\mu E_N L_N\beta_{\mu,N}\nonumber\\
&&\hspace{20mm}+S_N +\vep_{E\lambda,N}
\label{12.83}
\end{eqnarray}
Here: 
\begin{eqnarray}
&&\tilde{q}_N=\frac{\beta_{N,N}^2}{4c_N}(\Lb_N H_N-2\pi_N\lambda_N E_N H_N)\label{12.84}\\
&&\oB_N=\beta_{N,N}^2\left[\frac{\pi_N\rho_N}{4}\Lb_N H_N+\frac{\pi_N\rhob_N}{4}L_N H_N
-\frac{1}{2}a_N\left(\pi_N^2+\frac{1}{2c_N}\right)E_N H_N\right]\label{12.85}\\
&&B_N=\frac{\rho_N}{2}\beta_{N,N}\left(\sbeta_N-\frac{\pi_N}{2}\beta_{N,N}\right)\Lb_N H_N\nonumber\\
&&\hspace{7mm}+\rhob_N\left(c_N\pi_N\sbeta_N^2-c_N\pi_N^2\sbeta_N\beta_{N,N}-\frac{1}{2}\sbeta_N\beta_{\Nb,N}
+\frac{\pi_N}{4}\beta_{N,N}^2+\frac{\pi_N}{2}\beta_{N,N}\beta_{\Nb,N}\right)L_N H_N\nonumber\\
&&\hspace{7mm}+a_N\beta_{N,N}\left[-\pi_N\sbeta_N+\frac{1}{2}\left(\pi_N^2-\frac{1}{2c_N}\right)\beta_{N,N}\right]E_N H_N
\nonumber\\
&&\hspace{7mm}+H_N\left(\frac{(\beta_{N,N}+\beta_{\Nb,N})}{2}-c_N\pi_N\sbeta_N\right)\left(\pi_N\rhob_N s_{NL,N}
-\frac{a_N}{c_N}\ss_{N,N}\right) \label{12.86}
\end{eqnarray}
Also, $S_N$ is the $N$th approximant version of $S$, and $\vep_{E\lambda,N}$ is an error term satisfying the estimate:
\begin{equation}
\vep_{E\lambda,N}=O(\tau^N)
\label{12.87}
\end{equation}
Similarly, we deduce the $N$th approximant version of the conjugate equation \ref{12.53}:
\begin{eqnarray}
&&\Lb_N E_N\lambdab_N=(\pb_N-\chib_N)E_N\lambdab_N+\tilde{\qb}_N E_N\lambda_N+\oBb_N\tchib_N+
\Bb_N\tchi_N\nonumber\\
&&\hspace{10mm}+\frac{\rhob_N}{4}\beta_{\Nb,N}^2 E_N L_N H_N
+\rho_N\beta_{\Nb,N}\left(c_N\pi_N\sbeta_N-\frac{\beta_{\Nb,N}}{4}-\frac{\beta_{N,N}}{2}\right)
E\_N Lb_N H_N\nonumber\\
&&\hspace{10mm}-\frac{\pi_N}{2}a_N\beta_{\Nb,N}^2 E_N^2 H_N+\rho_N H_N\left(c_N\pi_N\sbeta_N-\frac{(\beta_{N,N}+\beta_{\Nb,N})}{2}\right)\Nb_N^\mu E_N\Lb_N\beta_{\mu,N}\nonumber\\
&&\hspace{10mm}+\Sb_N+\vep_{E\lambdab,N}
\label{12.88}
\end{eqnarray}
Here:
\begin{eqnarray}
&&\tilde{\qb}_N=\frac{\beta_{\Nb,N}^2}{4c_N}(L_N H_N-2\pi_N\lambdab_N E_N H_N)\label{12.89}\\
&&\Bb_N=\beta_{\Nb,N}^2\left[\frac{\pi_N\rhob_N}{4}L_N H_N+\frac{\pi_N\rho_N}{4}\Lb_N H_N
-\frac{1}{2}a_N\left(\pi_N^2+\frac{1}{2c_N}\right)E_N H_N\right]\label{12.90}\\
&&\oBb_N=\frac{\rhob_N}{2}\beta_{\Nb,N}\left(\sbeta_N-\frac{\pi_N}{2}\beta_{\Nb,N}\right)L_N H_N\nonumber\\
&&\hspace{7mm}+\rho_N\left(c_N\pi_N\sbeta_N^2-c_N\pi_N^2\sbeta_N\beta_{\Nb,N}
-\frac{1}{2}\sbeta_N\beta_{N,N}
+\frac{\pi_N}{4}\beta_{\Nb,N}^2+\frac{\pi_N}{2}\beta_{\Nb,N}\beta_{N,N}\right)\Lb_N H_N\nonumber\\
&&\hspace{7mm}+a_N\beta_{\Nb,N}\left[-\pi_N\sbeta_N+\frac{1}{2}\left(\pi_N^2-\frac{1}{2c_N}\right)\beta_{\Nb,N}\right]E_N H_N
\nonumber\\
&&\hspace{7mm}+H_N\left(\frac{(\beta_{N,N}+\beta_{\Nb,N})}{2}-c_N\pi_N\sbeta_N\right)
\left(\pi_N\rho_N s_{\Nb\Lb,N}-\frac{a_N}{c_N}\ss_{\Nb,N}\right) \label{12.91}
\end{eqnarray}
Also, $\Sb_N$ is the $N$th approximant version of $\Sb$, and $\vep_{E\lambdab,N}$ is an error term satisfying the estimate:
\begin{equation}
\vep_{E\lambdab,N}=O(\tau^N)
\label{12.92}
\end{equation}

Applying $T$ to the propagation equation \ref{10.161} for the $N$th approximant quantity $\lambda_N$ 
and using the fact that by the third of the commutation relations \ref{9.b1} (see \ref{9.b3}):
$$[L_N,T]=-\zeta_N E_N,$$
we obtain the $N$th approximant version of equation \ref{12.57}:
\begin{equation}
L_N T\lambda_N+\zeta_N E_N\lambda_N=p_N T\lambda_N+q_N T\lambdab_N+\lambda_N Tp_N+\lambdab_N Tq_N
\label{12.93}
\end{equation}
Proceeding from this equation following the path leading from equation \ref{12.57} to equation 
\ref{12.80} , using the results of Section 10.3 in place of the results of Chapter 3, and  
defining the $N$th approximant version of the quantities $\mu, \mub$ (see \ref{12.64})
\begin{equation}
\mu_N=E_N\lambda+\pi_N\lambda_N\tchi_N, \ \ \ 
\mub_N=E_N\lambdab_N+\pi_N\lambdab_N\tchib_N, 
\label{12.94}
\end{equation}
we deduce the $N$th approximant version of equation \ref{12.80}:
\begin{eqnarray}
&&L_N T\lambda_N=p_N T\lambda_N+\tilde{q}_N T\lambdab_N-\zeta_N E_N\lambda_N
+2B_N\mu_N+2\oB_N\mub_N\nonumber\\
&&\hspace{17mm}+\frac{\rho_N}{4}\beta_{N,N}^2 T\Lb_N H_N-\frac{\pi_N}{2}a_N\beta_{N,N}^2 TE_N H_N
\nonumber\\
&&\hspace{17mm}+\rhob_N\beta_{N,N}\left(c_N\pi_N\sbeta_N-\frac{\beta_{N,N}}{4}-\frac{\beta_{\Nb,N}}{2}\right)
TL_N H_N\nonumber\\
&&\hspace{17mm}+\rhob_N H_N
\left(c_N\pi_N\sbeta_N-\frac{(\beta_{N,N}+\beta_{\Nb,N})}{2}\right)N_N^\mu TL_N\beta_{\mu,N}
\nonumber\\
&&\hspace{17mm}+Q_N+\vep_{T\lambda,N}
\label{12.95}
\end{eqnarray}
Here $Q_N$ is the $N$th approximant version of $Q$ and $\vep_{T\lambda,N}$ is an error term satisfying 
the estimate:
\begin{equation}
\vep_{T\lambda,N}=O(\tau^{N-1})
\label{12.96}
\end{equation}
Similarly, we deduce the $N$th approximant version of the conjugate equation \ref{12.81}:
\begin{eqnarray}
&&\Lb_N T\lambdab_N=\pb_N T\lambdab_N+\tilde{\qb}_N T\lambda_N+\zeta_N E_N\lambdab_N
+2\oBb_N\mub_N+2\Bb_N\mu_N\nonumber\\
&&\hspace{17mm}+\frac{\rhob_N}{4}\beta_{\Nb,N}^2 TL_N H_N-\frac{\pi_N}{2}a_N\beta_{\Nb,N}^2 TE_N H_N
\nonumber\\
&&\hspace{17mm}+\rho_N\beta_{\Nb,N}\left(c_N\pi_N\sbeta_N-\frac{\beta_{\Nb,N}}{4}-\frac{\beta_{N,N}}{2}\right)
T\Lb_N H_N\nonumber\\
&&\hspace{17mm}+\rho_N H_N
\left(c_N\pi_N\sbeta_N-\frac{(\beta_{N,N}+\beta_{\Nb,N})}{2}\right)\Nb_N^\mu T\Lb_N\beta_{\mu,N}
\nonumber\\
&&\hspace{17mm}+\Qb_N+\vep_{T\lambdab,N}
\label{12.97}
\end{eqnarray}
Here $\Qb_N$ is the $N$th approximant version of $\Qb$ and $\vep_{T\lambdab,N}$ is an error term satisfying 
the estimate:
\begin{equation}
\vep_{T\lambdab,N}=O(\tau^{N-1})
\label{12.98}
\end{equation} 

Applying $E^{n-1}$ to equation \ref{12.50} and using the 1st of the commutation relations  
\ref{3.a14}, which implies that:
$$[L,E^{n-1}]E\lambda=-\sum_{i=0}^{n-2}E^{n-2-i}(\chi E^{i+2}\lambda)$$
we deduce the following propagation equation for $E^n\lambda$:
\begin{eqnarray}
&&LE^n\lambda=(p-n\chi)E^n\lambda+\tilde{q}E^n\lambdab+(B-\rho E\lambda)E^{n-1}\tchi+\oB E^{n-1}\tchib\nonumber\\
&&\hspace{10mm}+\frac{\rho}{4}\beta_N^2 E^n\Lb H
+\rhob\beta_N\left(c\pi\sbeta-\frac{\beta_N}{4}-\frac{\beta_{\Nb}}{2}\right)E^n LH\nonumber\\
&&\hspace{10mm}-\frac{\pi}{2}a\beta_N^2 E^{n+1} H
+\rhob H\left(c\pi\sbeta-\frac{(\beta_N+\beta_{\Nb})}{2}\right)N^\mu E^n L\beta_\mu\nonumber\\
&&\hspace{10mm}+S_{n-1}
\label{12.99}
\end{eqnarray}
where $S_{n-1}$ is of order $n$, but the highest order acoustical quantity occurring in $S_{n-1}$ 
is only of order $n-1$. Similarly, applying $E_N^{n-1}$ to equation \ref{12.83} and using the 1st of the 
commutation relations \ref{9.b1}, which implies that:
$$[L_N, E_N^{n-1}]E_N\lambda_N=-\sum_{i=0}^{n-2}E_N^{n-2-i}(\chi_N E_N^{i+2}\lambda_N)$$
we deduce:
\begin{eqnarray}
&&L_N E_N^n\lambda_N=(p_N-n\chi_N)E_N^n\lambda_N+\tilde{q}_N E_N^n\lambdab_N\nonumber\\
&&\hspace{20mm}+(B_N-\rho_N E_N\lambda_N)E_N^{n-1}\tchi_N+\oB_N E_N^{n-1}\tchib_N\nonumber\\
&&\hspace{20mm}+\frac{\rho_N}{4}\beta_{N,N}^2 E_N^n\Lb_N H_N \nonumber\\
&&\hspace{20mm}+\rhob_N\beta_{N,N}\left(c_N\pi_N\sbeta_N-\frac{\beta_{N,N}}{4}-\frac{\beta_{\Nb,N}}{2}\right)E_N^n L_N H_N\nonumber\\
&&\hspace{20mm}-\frac{\pi_N}{2}a_N\beta_{N,N}^2 E_N^{n+1} H_N \nonumber\\
&&\hspace{20mm}+\rhob_N H_N\left(c_N\pi_N\sbeta_N-\frac{(\beta_{N,N}+\beta_{\Nb,N})}{2}\right)
N_N^\mu E_N^n L_N\beta_{\mu,N}\nonumber\\
&&\hspace{20mm}+S_{n-1,N} +E_N^{n-1}\vep_{E\lambda,N}
\label{12.100}
\end{eqnarray}
where $S_{n-1,N}$ is the $N$th approximant analogue of $S_{n-1}$. 
Subtracting equation \ref{12.100} from \ref{12.99} we arrive at the following propagation equation 
for the acoustical difference quantity $\s^{(0,n)}\cla$ defined by the 1st of \ref{12.4} for 
$m=0, l=n$:
\begin{eqnarray}
&&L\s^{(0,n)}\cla=(p-n\chi)\s^{(0,n)}\cla+\tilde{q}\s^{(0,n)}\clab+(B-\rho E\lambda)s^{(n-1)}\ctchi
+\oB\s^{(n-1)}\ctchib\nonumber\\
&&\hspace{15mm}-(L-L_N-p+p_N+n(\chi-\chi_N))E_N^n\lambda_N+(\tilde{q}-\tilde{q}_N)E_N^n\lambdab_N\nonumber\\
&&\hspace{15mm}+(B-B_N-\rho E\lambda+\rho_N E_N\lambda_N)E_N^{n-1}\tchi_N+\oB E_N^{n-1}\ctchib_N\nonumber\\
&&\hspace{15mm}+\frac{\rho}{4}\beta_N^2 E^n\Lb H
-\frac{\rho_N}{4}\beta_{N,N}^2 E_N^n\Lb_N H_N \nonumber\\
&&\hspace{15mm}+\rhob\beta_N\left(c\pi\sbeta-\frac{\beta_N}{4}-\frac{\beta_{\Nb}}{2}\right)E^n LH
\nonumber\\
&&\hspace{30mm}-\rhob_N\beta_{N,N}\left(c_N\pi_N\sbeta_N-\frac{\beta_{N,N}}{4}-\frac{\beta_{\Nb,N}}{2}\right)E_N^n L_N H_N\nonumber\\
&&\hspace{15mm}-\frac{\pi}{2}a\beta_N^2 E^{n+1} H 
+\frac{\pi_N}{2}a_N\beta_{N,N}^2 E_N^{n+1} H_N \nonumber\\
&&\hspace{15mm}+\rhob H\left(c\pi\sbeta-\frac{(\beta_N+\beta_{\Nb})}{2}\right)N^\mu E^n L\beta_\mu
\nonumber\\
&&\hspace{30mm}-\rhob_N H_N\left(c_N\pi_N\sbeta_N-\frac{(\beta_{N,N}+\beta_{\Nb,N})}{2}\right)
N_N^\mu E_N^n L_N\beta_{\mu,N}\nonumber\\
&&\hspace{15mm}+S_{n-1}-S_{n-1,N}-E_N^{n-1}\vep_{E\lambda,N}
\label{12.101}
\end{eqnarray}

Applying likewise $E^{n-1}$ to equation \ref{12.53} and using the 2nd of the commutation relations  
\ref{3.a14}, which implies that:
$$[\Lb,E^{n-1}]E\lambdab=-\sum_{i=0}^{n-2}E^{n-2-i}(\chib E^{i+2}\lambdab)$$
we deduce the following propagation equation for $E^n\lambdab$ (the conjugate of \ref{12.99}):
\begin{eqnarray}
&&\Lb E^n\lambdab=(\pb-n\chib)E^n\lambdab+\tilde{\qb}E^n\lambda+(\oBb-\rhob E\lambdab)E^{n-1}\tchib+\Bb E^{n-1}\tchi\nonumber\\
&&\hspace{10mm}+\frac{\rhob}{4}\beta_{\Nb}^2 E^n LH
+\rho\beta_{\Nb}\left(c\pi\sbeta-\frac{\beta_{\Nb}}{4}-\frac{\beta_N}{2}\right)E^n\Lb H\nonumber\\
&&\hspace{10mm}-\frac{\pi}{2}a\beta_{\Nb}^2 E^{n+1} H
+\rho H\left(c\pi\sbeta-\frac{(\beta_N+\beta_{\Nb})}{2}\right)\Nb^\mu E^n\Lb\beta_\mu\nonumber\\
&&\hspace{10mm}+\Sb_{n-1}
\label{12.102}
\end{eqnarray}
where $\Sb_{n-1}$ is of order $n$, but the highest order acoustical quantity occurring in $\Sb_{n-1}$ 
is only of order $n-1$. Similarly, applying $E_N^{n-1}$ to equation \ref{12.88} and using the 2nd of the 
commutation relations \ref{9.b1}, which implies that:
$$[\Lb_N, E_N^{n-1}]E_N\lambdab_N=-\sum_{i=0}^{n-2}E_N^{n-2-i}(\chib_N E_N^{i+2}\lambdab_N)$$
we deduce:
\begin{eqnarray}
&&\Lb_N E_N^n\lambdab_N=(\pb_N-n\chib_N)E_N^n\lambdab_N+\tilde{\qb}_N E_N^n\lambda_N\nonumber\\
&&\hspace{20mm}+(\oBb_N-\rhob_N E_N\lambdab_N)E_N^{n-1}\tchib_N+\Bb_N E_N^{n-1}\tchi_N\nonumber\\
&&\hspace{20mm}+\frac{\rhob_N}{4}\beta_{\Nb,N}^2 E_N^n L_N H_N \nonumber\\
&&\hspace{20mm}+\rho_N\beta_{\Nb,N}\left(c_N\pi_N\sbeta_N-\frac{\beta_{\Nb,N}}{4}-\frac{\beta_{N,N}}{2}\right)E_N^n\Lb_N H_N\nonumber\\
&&\hspace{20mm}-\frac{\pi_N}{2}a_N\beta_{\Nb,N}^2 E_N^{n+1} H_N \nonumber\\
&&\hspace{20mm}+\rho_N H_N\left(c_N\pi_N\sbeta_N-\frac{(\beta_{N,N}+\beta_{\Nb,N})}{2}\right)
\Nb_N^\mu E_N^n\Lb_N\beta_{\mu,N}\nonumber\\
&&\hspace{20mm}+\Sb_{n-1,N} +E_N^{n-1}\vep_{E\lambdab,N}
\label{12.103}
\end{eqnarray}
where $\Sb_{n-1,N}$ is the $N$th approximant analogue of $\Sb_{n-1}$. 
Subtracting equation \ref{12.103} from \ref{12.102} we arrive at the following propagation equation 
for the acoustical difference quantity $\s^{(0,n)}\clab$ defined by the 2nd of \ref{12.4} for 
$m=0, l=n$:
\begin{eqnarray}
&&\Lb\s^{(0,n)}\clab=(\pb-n\chib)\s^{(0,n)}\clab+\tilde{\qb}\s^{(0,n)}\cla+(\oBb-\rhob E\lambdab)\s^{(n-1)}\ctchib
+\Bb\s^{(n-1)}\ctchi\nonumber\\
&&\hspace{15mm}-(\Lb-\Lb_N-\pb+\pb_N+n(\chib-\chib_N))E_N^n\lambdab_N+(\tilde{\qb}-\tilde{\qb}_N)E_N^n\lambda_N\nonumber\\
&&\hspace{15mm}+(\oBb-\oBb_N-\rhob E\lambdab+\rhob_N E_N\lambdab_N)E_N^{n-1}\tchib_N+\Bb E_N^{n-1}\ctchi_N\nonumber\\
&&\hspace{15mm}+\frac{\rhob}{4}\beta_{\Nb}^2 E^n LH
-\frac{\rhob_N}{4}\beta_{\Nb,N}^2 E_N^n L_N H_N \nonumber\\
&&\hspace{15mm}+\rho\beta_{\Nb}\left(c\pi\sbeta-\frac{\beta_{\Nb}}{4}-\frac{\beta_N}{2}\right)E^n\Lb H
\nonumber\\
&&\hspace{30mm}-\rho_N\beta_{\Nb,N}\left(c_N\pi_N\sbeta_N-\frac{\beta_{\Nb,N}}{4}-\frac{\beta_{N,N}}{2}\right)E_N^n\Lb_N H_N\nonumber\\
&&\hspace{15mm}-\frac{\pi}{2}a\beta_{\Nb}^2 E^{n+1} H 
+\frac{\pi_N}{2}a_N\beta_{\Nb,N}^2 E_N^{n+1} H_N \nonumber\\
&&\hspace{15mm}+\rho H\left(c\pi\sbeta-\frac{(\beta_N+\beta_{\Nb})}{2}\right)\Nb^\mu E^n\Lb\beta_\mu
\nonumber\\
&&\hspace{30mm}-\rho_N H_N\left(c_N\pi_N\sbeta_N-\frac{(\beta_{N,N}+\beta_{\Nb,N})}{2}\right)
\Nb_N^\mu E_N^n\Lb_N\beta_{\mu,N}\nonumber\\
&&\hspace{15mm}+\Sb_{n-1}-\Sb_{n-1,N}-E_N^{n-1}\vep_{E\lambdab,N}
\label{12.104}
\end{eqnarray}

For $m=1,...,n$, we apply $T^{m-1}$ to \ref{12.80} and use the fact that by the 3rd of the commutation relations 
\ref{3.154}:
$$[L,T^{m-1}]T\lambda=-\sum_{j=0}^{m-2}T^{m-2-j}(\zeta ET^{j+1}\lambda)$$
Also, in view of \ref{8.137}, \ref{8.138} we have:
\begin{eqnarray}
&&[T\tchi]_{P.A.}=2E^2\lambda+2\pi\lambda E\tchi \label{12.105}\\
&&[T\tchib]_{P.A.}=2E^2\lambdab+2\pi\lambdab E\tchib \label{12.106}
\end{eqnarray}
Using these we establish by induction that for $m=1,...,n$: 
\begin{eqnarray}
&&[T^{m-1}\mu]_{P.A.}=\sum_{j=0}^{m-1}(2\pi\lambda)^j E^{j+1}T^{m-1-j}\lambda +
\frac{1}{2}(2\pi\lambda)^m E^{m-1}\tchi \nonumber\\
&&\hspace{19mm}:=\mu_{m-1} \label{12.107}\\
&&[T^{m-1}\mub]_{P.A.}=\sum_{j=0}^{m-1}(2\pi\lambdab)^j E^{j+1}T^{m-1-j}\lambdab + 
\frac{1}{2}(2\pi\lambdab)^m E^{m-1}\tchib \nonumber\\
&&\hspace{19mm}:=\mub_{m-1} \label{12.108} 
\end{eqnarray}
We then deduce the following propagation equation for $T^m\lambda$, $m=1,...,n$:
\begin{eqnarray}
&&LT^m\lambda=pT^m\lambda+\tilde{q}T^m\lambdab-m\zeta ET^{m-1}\lambda+2B\mu_{m-1}+2\oB\mub_{m-1}\nonumber\\
&&\hspace{10mm}+\frac{\rho}{4}\beta_N^2 T^m\Lb H
+\rhob\beta_N\left(c\pi\sbeta-\frac{\beta_N}{4}-\frac{\beta_{\Nb}}{2}\right)T^m LH\nonumber\\
&&\hspace{10mm}-\frac{\pi}{2}a\beta_N^2 T^m EH+\rhob H\left(c\pi\sbeta-\frac{(\beta_N+\beta_{\Nb})}{2}\right)N^\mu T^m L\beta_\mu\nonumber\\
&&\hspace{10mm}+Q_{m-1}
\label{12.109}
\end{eqnarray}
where $Q_{m-1}$ is of order $m$, but the highest order acoustical quantity occurring in $Q_{m-1}$ is 
only of order $m-1$. We then apply $E^{n-m}$ to \ref{12.109} and use the fact that by the 1st of 
the commutation relations \ref{3.154}:
$$[L,E^{n-m}]T^m\lambda=-\sum_{i=0}^{n-m-1}E^{n-m-1-i}(\chi E^{i+1}T^m\lambda)$$
Noting that:
\begin{eqnarray}
&&[E^{n-m}\mu_{m-1}]_{P.A.}=\sum_{j=0}^{m-1}(2\pi\lambda)^j E^{n-m+j+1}T^{m-1-j}\lambda +
\frac{1}{2}(2\pi\lambda)^m E^{n-1}\tchi \nonumber\\
&&\hspace{19mm}:=\mu_{m-1,n-m} \label{12.110}\\
&&[E^{n-m}\mub_{m-1}]_{P.A.}=\sum_{j=0}^{m-1}(2\pi\lambdab)^j E^{n-m+j+1}T^{m-1-j}\lambdab + 
\frac{1}{2}(2\pi\lambdab)^m E^{n-1}\tchib \nonumber\\
&&\hspace{19mm}:=\mub_{m-1,n-m} \label{12.111} 
\end{eqnarray}
we deduce the following propagation equation for $E^{n-m}T^m\lambda$:
\begin{eqnarray}
&&LE^{n-m}T^m\lambda=(p-(n-m)\chi)E^{n-m}T^m\lambda+\tilde{q}E^{n-m}T^m\lambdab
\nonumber\\
&&\hspace{20mm}-m\zeta E^{n-m+1}T^{m-1}\lambda+2B\mu_{m-1,n-m}+2\oB\mub_{m-1,n-m}\nonumber\\
&&\hspace{10mm}+\frac{\rho}{4}\beta_N^2 E^{n-m}T^m\Lb H
+\rhob\beta_N\left(c\pi\sbeta-\frac{\beta_N}{4}-\frac{\beta_{\Nb}}{2}\right)E^{n-m}T^m LH\nonumber\\
&&\hspace{10mm}-\frac{\pi}{2}a\beta_N^2 E^{n-m}T^m EH+\rhob H\left(c\pi\sbeta-\frac{(\beta_N+\beta_{\Nb})}{2}\right)N^\mu E^{n-m}T^m L\beta_\mu\nonumber\\
&&\hspace{20mm}+Q_{m-1,n-m}
\label{12.112}
\end{eqnarray}
where $Q_{m-1,n-m}$ is of order $n$, but the highest order acoustical quantity occurring in 
$Q_{m-1,n-m}$ is only of order $n-1$. Similarly applying first $T^{m-1}$ to equation \ref{12.95} 
and then applying $E_N^{n-m}$ to the resulting equation, we deduce the analogous equation for the 
$N$th approximants:
\begin{eqnarray}
&&L_N E_N^{n-m}T^m\lambda_N=(p_N-(n-m)\chi_N)E_N^{n-m}T^m\lambda_N+\tilde{q}_N E_N^{n-m}T^m\lambdab_N
\nonumber\\
&&\hspace{18mm}-m\zeta_N E_N^{n-m+1}T^{m-1}\lambda_N+2B_N\mu_{m-1,n-m,N}+2\oB_N\mub_{m-1,n-m,N}\nonumber\\
&&\hspace{28mm}+\frac{\rho_N}{4}\beta_{N,N}^2 E_N^{n-m}T^m\Lb_N H_N \nonumber\\
&&\hspace{18mm}+\rhob_N\beta_{N,N}\left(c_N\pi_N\sbeta_N-\frac{\beta_{N,N}}{4}
-\frac{\beta_{\Nb,N}}{2}\right)E_N^{n-m}T^m L_N H_N\nonumber\\
&&\hspace{28mm}-\frac{\pi_N}{2}a_N\beta_{N,N}^2 E_N^{n-m}T^m E_N H_N\nonumber\\
&&\hspace{18mm}+\rhob_N H_N\left(c_N\pi_N\sbeta_N-\frac{(\beta_{N,N}+\beta_{\Nb,N})}{2}\right)N_N^\mu E_N^{n-m}T^m 
L_N\beta_{\mu,N}\nonumber\\
&&\hspace{28mm}+Q_{m-1,n-m,N}+E_N^{n-m}T^{m-1}\vep_{T\lambda,N}
\label{12.113}
\end{eqnarray}
Here, the quantities $\mu_{m-1,n-m,N}$, $\mub_{m-1,n-m,N}$ are the $N$th approximant versions of the 
quantities $\mu_{m-1,n-m}$, $\mub_{m-1,n-m}$:
\begin{eqnarray}
&&\mu_{m-1,n-m,N}=\sum_{j=0}^{m-1}(2\pi_N\lambda_N)^j E_N^{n-m+j+1}T^{m-1-j}\lambda_N +
\frac{1}{2}(2\pi_N\lambda_N)^m E_N^{n-1}\tchi_N \nonumber\\
&&\label{12.114}\\
&&\mub_{m-1,n-m,N}=\sum_{j=0}^{m-1}(2\pi_N\lambdab_N)^j E_N^{n-m+j+1}T^{m-1-j}\lambdab_N + 
\frac{1}{2}(2\pi_N\lambdab_N)^m E_N^{n-1}\tchib_N \nonumber\\
&&\label{12.115}
\end{eqnarray} 
and $Q_{m-1,n-m,N}$ is the $N$th approximant version of $Q_{m-1,n-m}$. 
Subtracting equation \ref{12.113} from \ref{12.112} we arrive at the following propagation equation 
for the acoustical difference quantity $\s^{(m,n-m)}\cla$ defined by the 1st of \ref{12.4} for 
$m=1,...,n$, $l=n-m$:
\begin{eqnarray}
&&L\s^{(m,n-m)}\cla=(p-(n-m)\chi)\s^{(m,n-m)}\cla+\tilde{q}\s^{(m,n-m)}\clab-m\zeta\s^{(m-1,n-m+1)}\cla
\nonumber\\
&&\hspace{30mm}+2B\cmu_{m-1,n-m}+2\oB\cmub_{m-1,n-m}\nonumber\\
&&\hspace{20mm}+(L-L_N+p-p_N-(n-m)(\chi-\chi_N))E_N^{n-m}T^m\lambda_N\nonumber\\
&&\hspace{20mm}+(\tilde{q}-\tilde{q}_N)E_N^{n-m}T^m\lambdab_N
-m(\zeta-\zeta_N)E_N^{n-m+1}T^{m-1}\lambda_N\nonumber\\
&&\hspace{30mm}+2(B-B_N)\mu_{m-1,n-m,N}+2(\oB-\oB_N)\mub_{m-1,n-m,N}\nonumber\\
&&\hspace{20mm}+\frac{\rho}{4}\beta_N^2 E^{n-m}T^m\Lb H-\frac{\rho_N}{4}\beta_{N,N}^2 E_N^{n-m}T^m\Lb_N H_N
\nonumber\\
&&\hspace{20mm}+\rhob\beta_N\left(c\pi\sbeta-\frac{\beta_N}{4}-\frac{\beta_{\Nb}}{2}\right)E^{n-m}T^m LH
\nonumber\\
&&\hspace{30mm}-\rhob_N\beta_{N,N}\left(c_N\pi_N\sbeta_N-\frac{\beta_{N,N}}{4}
-\frac{\beta_{\Nb,N}}{2}\right)E_N^{n-m}T^m L_N H_N \nonumber\\
&&\hspace{20mm}-\frac{\pi}{2}a\beta_N^2 E^{n-m}T^m EH+\frac{\pi_N}{2}a_N\beta_{N,N}^2 E_N^{n-m}T^m 
E_N H_N\nonumber\\
&&\hspace{20mm}+\rhob H\left(c\pi\sbeta-\frac{(\beta_N+\beta_{\Nb})}{2}\right)
N^\mu E^{n-m}T^m L\beta_\mu\nonumber\\
&&\hspace{30mm}-\rhob_N H_N\left(c_N\pi_N\sbeta_N-\frac{(\beta_{N,N}+\beta_{\Nb,N})}{2}\right)
N_N^\mu E_N^{n-m}T^m L_N\beta_{\mu,N}\nonumber\\
&&\hspace{20mm}+Q_{m-1,n-m}-Q_{m-1,n-m,N}-E_N^{n-m}T^{m-1}\vep_{T\lambda,N}
\label{12.116}
\end{eqnarray}
Here we have introduced the acoustical difference quantities:
\begin{equation}
\cmu_{m-1,n-m}=\mu_{m-1,n-m}-\mu_{m-1,n-m,N}, \ \ \ 
\cmub_{m-1,n-m}=\mub_{m-1,n-m}-\mub_{m-1,n-m,N}
\label{12.117}
\end{equation}
By \ref{12.110}, \ref{12.111} and \ref{12.114}, \ref{12.115}:
\begin{eqnarray}
&&\cmu_{m-1,n-m}=\sum_{j=0}^{m-1}(2\pi\lambda)^j \s^{(m-1-j,n-m+j+1)}\cla +
\frac{1}{2}(2\pi\lambda)^m \s^{(n-1)}\ctchi\nonumber\\
&&\hspace{20mm}+\sum_{j=0}^{m-1}\left((2\pi\lambda)^j-(2\pi_N\lambda_N)^j\right)
E_N^{n-m+j+1}T^{m-1-j}\lambda_N\nonumber\\
&&\hspace{20mm}+\frac{1}{2}\left((2\pi\lambda)^m-(2\pi_N\lambda_N)^m\right)
E_N^{n-1}\tchi_N \label{12.118}\\
&&\cmub_{m-1,n-m}=\sum_{j=0}^{m-1}(2\pi\lambdab)^j \s^{(m-1-j,n-m+j+1)}\clab +
\frac{1}{2}(2\pi\lambdab)^m \s^{(n-1)}\ctchib\nonumber\\
&&\hspace{20mm}+\sum_{j=0}^{m-1}\left((2\pi\lambdab)^j-(2\pi_N\lambdab_N)^j\right)
E_N^{n-m+j+1}T^{m-1-j}\lambdab_N\nonumber\\
&&\hspace{20mm}+\frac{1}{2}\left((2\pi\lambdab)^m-(2\pi_N\lambdab_N)^m\right)
E_N^{n-1}\tchib_N \label{12.119}
\end{eqnarray}
In a similar manner we deduce the propagation equation 
for the acoustical difference quantity $\s^{(m,n-m)}\clab$ defined by the 2nd of \ref{12.4} for 
$m=1,...,n$, $l=n-m$, the conjugate of \ref{12.116}:
\begin{eqnarray}
&&\Lb\s^{(m,n-m)}\clab=(\pb-(n-m)\chib)\s^{(m,n-m)}\clab+\tilde{\qb}\s^{(m,n-m)}\cla
+m\zeta\s^{(m-1,n-m+1)}\clab
\nonumber\\
&&\hspace{30mm}+2\oBb\cmub_{m-1,n-m}+2\Bb\cmu_{m-1,n-m}\nonumber\\
&&\hspace{20mm}+(\Lb-\Lb_N+\pb-\pb_N-(n-m)(\chib-\chib_N))E_N^{n-m}T^m\lambdab_N\nonumber\\
&&\hspace{20mm}+(\tilde{\qb}-\tilde{\qb}_N)E_N^{n-m}T^m\lambda_N
+m(\zeta-\zeta_N)E_N^{n-m+1}T^{m-1}\lambdab_N\nonumber\\
&&\hspace{30mm}+2(\oBb-\oBb_N)\mub_{m-1,n-m,N}+2(\Bb-\Bb_N)\mu_{m-1,n-m,N}\nonumber\\
&&\hspace{20mm}+\frac{\rhob}{4}\beta_{\Nb}^2 E^{n-m}T^m LH-\frac{\rhob_N}{4}\beta_{\Nb,N}^2 
E_N^{n-m}T^m L_N H_N
\nonumber\\
&&\hspace{20mm}+\rho\beta_{\Nb}\left(c\pi\sbeta-\frac{\beta_{\Nb}}{4}-\frac{\beta_N}{2}\right)
E^{n-m}T^m\Lb H
\nonumber\\
&&\hspace{30mm}-\rho_N\beta_{\Nb,N}\left(c_N\pi_N\sbeta_N-\frac{\beta_{\Nb,N}}{4}
-\frac{\beta_{N,N}}{2}\right)E_N^{n-m}T^m\Lb_N H_N \nonumber\\
&&\hspace{20mm}-\frac{\pi}{2}a\beta_{\Nb}^2 E^{n-m}T^m EH+\frac{\pi_N}{2}a_N\beta_{\Nb,N}^2 E_N^{n-m}T^m 
E_N H_N\nonumber\\
&&\hspace{20mm}+\rho H\left(c\pi\sbeta-\frac{(\beta_N+\beta_{\Nb})}{2}\right)
\Nb^\mu E^{n-m}T^m\Lb\beta_\mu\nonumber\\
&&\hspace{30mm}-\rho_N H_N\left(c_N\pi_N\sbeta_N-\frac{(\beta_{N,N}+\beta_{\Nb,N})}{2}\right)
\Nb_N^\mu E_N^{n-m}T^m\Lb_N\beta_{\mu,N}\nonumber\\
&&\hspace{20mm}+\Qb_{m-1,n-m}-\Qb_{m-1,n-m,N}-E_N^{n-m}T^{m-1}\vep_{T\lambdab,N}
\label{12.120}
\end{eqnarray}
Here $\Qb_{m-1,n-m}$ is of order $n$, but the highest order acoustical quantity occurring in 
$\Qb_{m-1,n-m}$ is only of order $n-1$, and $\Qb_{m-1,n-m,N}$ is the $N$th approximant version 
of $\Qb_{m-1,n-m}$. 

\vspace{5mm}

\section{Estimates for $(\s^{(n-1)}\ctchi, \s^{(n-1)}\ctchib)$ and $(\s^{(0,n)}\cla, \s^{(0,n)}\clab)$}

In the next section we shall derive estimates for $(\Omega^n T\chf, \Omega^n T\cv, T\Omega^n T\cga)$. 
As we shall see, the derivation of these estimates necessarily involves the derivation of appropriate 
estimates on ${\cal K}$ for $\s^{(n-1)}\ctchib$ and $\s^{(0,n)}\clab$, which in turn involve the 
derivation of appropriate estimates for $\s^{(n-1)}\ctchi$ and for $\s^{(0,n)}\cla$ respectively. 

We begin with the derivation of the estimates for $\s^{(n-1)}\ctchi$. 
We consider the propagation equation for $\s^{(n-1)}\ctchi$, equation \ref{12.29}. Let us denote 
the right hand side by $\s^{(n-1)}\chF$, so this equation reads:
\begin{equation}
L\s^{(n-1)}\ctchi+(n+1)\chi\s^{(n-1)}\ctchi=\s^{(n-1)}\chF
\label{12.121}
\end{equation}
Here as in Section 10.3 we use the $(\ub,u,\vartheta^\prime)$ coordinates on ${\cal N}$ which are 
adapted to the flow of $L$. In view of \ref{10.344}, equation \ref{12.121} takes in these coordinates 
the form:
\begin{equation}
\frac{\partial}{\partial\ub}\left(\sh^{\prime(n+1)/2}\s^{(n-1)}\ctchi\right)=
\sh^{\prime(n+1)/2}\s^{(n-1)}\chF
\label{12.122}
\end{equation}
Integrating then from $\Cb_0$, noting that 
\begin{equation}
\left.\s^{(n-1)}\ctchi\right|_{\Cb_0}=0
\label{12.123}
\end{equation}
we obtain:
\begin{equation}
\left(\sh^{\prime(n+1)/2}\s^{(n-1)}\ctchi\right)(\ub_1,u,\vartheta^\prime)=
\int_0^{\ub_1}\left(\sh^{\prime(n+1)/2}\s^{(n-1)}\chF\right)(\ub,u,\vartheta^\prime)d\ub
\label{12.124}
\end{equation}
In view of the remark \ref{10.352} we rewrite this in the form:
\begin{eqnarray*}
&&\left(\s^{(n-1)}\ctchi\sh^{\prime 1/4}\right)(\ub_1,u,\vartheta^\prime)=\\
&&\hspace{20mm}\int_0^{\ub_1}\left(\frac{\sh^\prime(\ub,u,\vartheta^\prime)}
{\sh^\prime(\ub_1,u,\vartheta^\prime)}\right)^{(2n+1)/4}\left(\s^{(n-1)}\chF\sh^{\prime 1/4}\right)
(\ub,u,\vartheta^\prime)d\ub
\end{eqnarray*}
By \ref{10.357} we have:
\begin{equation}
\left(\frac{\sh^\prime(\ub,u,\vartheta^\prime)}
{\sh^\prime(\ub_1,u,\vartheta^\prime)}\right)^{(2n+1)/4}\leq k
\label{12.125}
\end{equation}
where $k$ is a constant greater than 1, but which can be chosen as close to 1 as we wish by suitably 
restricting $\ub_1$. Therefore \ref{12.124} implies:
\begin{equation}
\left|\left(\s^{(n-1)}\ctchi\sh^{\prime 1/4}\right)(\ub_1,u,\vartheta^\prime)\right|\leq k
\int_0^{\ub_1}\left|\left(\s^{(n-1)}\chF\sh^{\prime 1/4}\right)(\ub,u,\vartheta^\prime)\right|d\ub
\label{12.126}
\end{equation}
Taking $L^2$ norms with respect to $\vartheta^\prime\in S^1$ this implies:
\begin{equation}
\left\|\left(\s^{(n-1)}\ctchi\sh^{\prime 1/4}\right)(\ub_1,u)\right\|_{L^2(S^1)}\leq k
\int_0^{\ub_1}\left\|\left(\s^{(n-1)}\chF\sh^{\prime 1/4}\right)(\ub,u)\right\|_{L^2(S^1)}d\ub
\label{12.127}
\end{equation}
or, in view of \ref{10.352}:
\begin{equation}
\|\s^{(n-1)}\ctchi\|_{L^2(S_{\ub_1,u})}\leq k\int_0^{\ub_1}\|\s^{(n-1)}\chF\|_{L^2(S_{\ub,u})}d\ub
\label{12.128}
\end{equation}

We shall presently estimate the leading contributions to the integral on the right in \ref{12.128}. 
Consider first the $n$th order acoustical difference terms: 
\begin{equation}
A\s^{(n-1)}\ctchi+\oA\s^{(n-1)}\ctchib
\label{12.129}
\end{equation}
The coefficients $A$ and $\oA$ being given by \ref{12.14} and \ref{12.15} respectively, the assumptions 
\ref{10.440} imply: 
\begin{equation}
|A|, |\oA| \leq C\ub \ \mbox{: in ${\cal R}_{\delta,\delta}$}
\label{12.130}
\end{equation}
It follows that the contribution of the terms \ref{12.129} is bounded by:
\begin{equation}
C\int_0^{\ub_1}\ub\|\s^{(n-1)}\ctchi\|_{L^2(S_{\ub,u})}d\ub + C\int_0^{\ub_1}\ub\|\s^{(n-1)}\ctchib\|_{L^2(S_{\ub,u})}d\ub
\label{12.131}
\end{equation}
Here the first can obviously be absorbed in the integral inequality \ref{12.128}. We shall show in 
the sequel that the second can also be absorbed. 

The other leading contributions are those of the principal difference terms:
\begin{equation}
-\beta_N\sbeta E^n LH+\beta_{N,N}\sbeta_N E_N^n L_N H_N 
+\frac{1}{2}\rho\beta_N^2 E^{n+1}H-\frac{1}{2}\rho_N\beta_{N,N}^2 E_N^{n+1} H_N
\label{12.132}
\end{equation}
Here we must keep track not only of the actually principal terms, namely those containing the 
derivatives of the $\beta_\mu$ of order $n+1$, but also of the terms containing acoustical quantities 
of order $n$. 
As we shall see the dominant contribution comes from the 2nd difference term. In fact, since 
$$EH=-2H^\prime (g^{-1})^{\mu\nu}\beta_\nu E\beta_\mu, \ \ \mbox{hence} \ \ 
[E^{n+1}H]_{P.P.}=-2H^\prime (g^{-1})^{\mu\nu}\beta_\nu E^{n+1}\beta_\mu$$
and similarly for the $N$th approximants, the actual principal difference term in question is:
\begin{equation}
-\rho\beta_N^2 H^\prime (g^{-1})^{\mu\nu}\beta_\nu E^{n+1}\beta_\mu 
+\rho_N\beta_{N,N}^2H_N^\prime (g^{-1})^{\mu\nu}\beta_{\nu,N} E_N^{n+1}\beta_{\mu,N}
\label{12.133}
\end{equation}
By the formula \ref{10.368} and an analogous formula for the $N$th approximants, \ref{12.133} is 
expressed as the sum of differences of $E$, $N$, and $\Nb$ components:
\begin{eqnarray}
&&-\beta_N^2\sbeta\eta^2 H^\prime\rho E^\mu E^{n+1}\beta_\mu 
+\beta_{N,N}^2\sbeta_N\eta_N^2 H^\prime_N \rho_N E_N^\mu E_N^{n+1}\beta_{\mu,N}\nonumber\\
&&+\beta_N^2\frac{\beta_{\Nb}}{2c}\eta^2 H^\prime\rho N^\mu E^{n+1}\beta_\mu 
-\beta_{N,N}^2\frac{\beta_{\Nb,N}}{2c_N}\eta_N^2 H_N^\prime\rho_N N_N^\mu E_N^{n+1}\beta_{\mu,N}
\nonumber\\
&&+\beta_N^2\frac{\beta_N}{2c}\eta^2 H^\prime\rho\Nb^\mu E^{n+1}\beta_\mu 
-\beta_{N,N}^2\frac{\beta_{N,N}}{2c_N}\eta_N^2 H^\prime_N\rho_N\Nb_N^\mu E_N^{n+1}\beta_{\mu,N}
\nonumber\\&&\label{12.134}
\end{eqnarray}
The dominant contribution is that of the $\Nb$ component difference. To estimate this we write: 
\begin{eqnarray}
&&\rhob \Nb^\mu E^{n+1}\beta_\mu=\Lb^\mu E^{n+1}\beta_\mu\label{12.135}\\
&&=E^n(\Lb^\mu E\beta_\mu)-(E^n\Lb^\mu)E\beta_\mu-\sum_{i=1}^{n-1}\left(\begin{array}{c} n\\i \end{array}\right)(E^i\Lb^\mu)E^{n+1-i}\beta_\mu\nonumber
\end{eqnarray}
Here the sum is of order $n$ with vanishing principal acoustical part. Such terms can be ignored. 
As for the 2nd term on the right, we have, since $\Lb^\mu=\rhob\Nb^\mu$, by the conjugates of \ref{12.7} and \ref{12.10}:
\begin{equation}
[(E^n\Lb^\mu)E\beta_\mu]_{P.A.}=c^{-1}\ss_{\Nb} E^n\lambda+\pi\rhob\ss_{\Nb}E^{n-1}\tchi
+\rhob\sss E^{n-1}\tchib
\label{12.136}
\end{equation}
Since $\Lb^\mu E\beta_\mu=E^\mu\Lb\beta_\mu$, the 1st term on the right in \ref{12.135} is:
\begin{eqnarray}
&&E^n(E^\mu\Lb\beta_\mu)=E^\mu E^n\Lb\beta_\mu \label{12.137}\\
&&\hspace{23mm}+(E^n E^\mu)\Lb\beta_\mu +\sum_{i=1}^{n-1}\left(\begin{array}{c} n\\i \end{array}\right)
(E^i E^\mu)E^{n-i}\Lb\beta_\mu \nonumber
\end{eqnarray}
Here the sum is of order $n$ with vanishing principal acoustical part, hence can be ignored. As for the 
2nd term on the right, by \ref{12.8} we have:
\begin{equation}
[(E^n E^\mu)\Lb\beta_\mu]_{P.A.}=\frac{1}{2}\left(c^{-1}s_{\Nb\Lb}E^{n-1}\tchi
+\rhob\sss E^{n-1}\tchib\right) 
\label{12.138}
\end{equation}
Finally, by the 2nd of the commutation relations \ref{3.a14} the 1st term on the right in \ref{12.137} 
is:
\begin{eqnarray}
&&E^\mu\Lb E^n\beta_\mu+E^\mu\sum_{i=0}^{n-1}E^{n-1-i}(\chib E^{i+1}\beta_\mu) \nonumber\\
&&=E^\mu\Lb E^n\beta_\mu+\sss E^{n-1}\chib+E^\mu\sum_{j=1}^{n-1}\left(\begin{array}{c} n-1\\j\end{array}
\right)(E^{n-1-j}\chib)E^{j+1}\beta_\mu\nonumber\\
&&\hspace{20mm}+E^\mu\sum_{i=1}^{n-1}E^{n-1-i}(\chib E^{i+1}\beta_\mu)\label{12.139}
\end{eqnarray}
Here the two sums on the right are of order $n$ with vanishing principal acoustical part, hence can be 
ignored. Since 
$$[E^{n-1}\chib]_{P.A.}=\rhob E^{n-1}\tchib$$
we conclude from the above that 
$$\rhob\Nb^\mu E^{n+1}\beta_\mu$$
is, up to terms which can be ignored, given by:
\begin{equation}
E^\mu\Lb E^n\beta_\mu+\frac{1}{2}\rhob\sss E^{n-1}\tchib
+\left(\frac{1}{2c}s_{\Nb\Lb}-\pi\rhob\ss_{\Nb}\right)E^{n-1}\tchi-\frac{1}{c}\ss_{\Nb}E^n\lambda
\label{12.140}
\end{equation}
Similar formulas hold for the $N$th approximants. It follows that up to terms which can be ignored 
the last the differences \ref{12.134} is given by:
\begin{equation}
\frac{\beta_N^3}{2c}\eta^2 H^\prime\left\{\frac{\rho}{\rhob}\s^{(E;0,n)}\cxi_{\Lb}
+\frac{1}{2}\rho\sss\s^{(n-1)}\ctchib+\rho\left(\frac{s_{\Nb\Lb}}{2\lambda}-\pi\ss_N\right)\s^{(n-1)}\ctchi 
-\frac{\rho}{\lambda}\ss_N\s^{(0,n)}\cla\right\}
\label{12.141}
\end{equation}
The contribution of this to the integral on the right in \ref{12.128} is bounded by:
\begin{eqnarray}
&&\frac{C}{u^2}\int_0^{\ub_1}\ub\|\s^{(E;0,n)}\cxi_{\Lb}\|_{L^2(S_{\ub,u})}d\ub \nonumber\\
&&+C\int_0^{\ub_1}\ub\|\s^{(n-1)}\ctchib\|_{L^2(S_{\ub,u})}d\ub
+\frac{C}{u^2}\int_0^{\ub_1}\ub\|\s^{(n-1)}\ctchi\|_{L^2(S_{\ub,u})}d\ub \nonumber\\
&&+\frac{C}{u^2}\int_0^{\ub_1}\ub\|\s^{(0,n)}\cla\|_{L^2(S_{\ub,u})}d\ub \label{12.142}
\end{eqnarray}
Here the 2nd integral is the same as the 2nd integral in \ref{12.131}. The 3rd integral is the same 
as the 1st integral in \ref{12.131}, but here this comes with the factor $u^{-2}$ which makes the 
contribution {\em borderline}. As for the contribution of the 4th integral in \ref{12.142}, we shall 
show in the sequel that this contribution, like that of the 2nd integral, can be absorbed. 

In regard to \ref{12.128}, we shall estimate $\|\s^{(n-1)}\ctchi\|_{L^2(\Cb_{\ub_1}^{u_1})}$ and 
$\|\s^{(n-1)}\ctchib\|_{L^2({\cal K}^\tau)}$. We shall make use of the following lemma. 

\vspace{2.5mm}

\noindent{\bf Lemma 12.1} \ \ Let $f$ be a non-negative function on $R_{\delta,\delta}$. 

\noindent{\bf a)} For all $(\ub_1,u_1)\in R_{\delta,\delta}$ we have:
$$\left\|\int_0^{\ub_1}f(\ub,\cdot)d\ub\right\|_{L^2(\ub_1,u_1)}
\leq\int_0^{\ub_1}\|f(\ub,\cdot)\|_{L^2(\ub_1,u_1)}d\ub$$

\noindent{\bf b)} Let us define on $[0,\delta]$ the function: 
$$g(u)=\int_0^uf(\ub,u)d\ub$$ 
Then for all $\tau\in [0,\delta]$ we have:
$$\|g\|_{L^2(0,\tau)}\leq\int_0^\tau\|f(\ub,\cdot)\|_{L^2(\ub,\tau)}d\ub$$

\noindent{\em Proof:} Part {\bf a} is standard. To establish part {\bf b} we extend for each 
$u\in [0,\tau]$ the function $f(\cdot,u)$ to $[0,\tau]$ by defining:
$$\of(\ub,u)=\left\{\begin{array}{cccc} f(\ub,u)&:&\mbox{if}& \ub\in [0,u]\\
0&:&\mbox{if}& \ub\in(u,\tau] \end{array}\right\}$$
We then have:
$$g(u)=\int_0^\tau\of(\ub,u)d\ub \ \ : \ u\in[0,\tau]$$
It follows that:
$$\|g\|_{L^2(0,\tau)}\leq\int_0^\tau\|\of(\ub,\cdot)\|_{L^2(0,\tau)}d\ub$$
Since 
$$\|\of(\ub,\cdot)\|^2_{L^2(0,\tau)}=\int_0^\tau \of^2(\ub,u)du=\int_{\ub}^\tau f^2(\ub,u)du
=\|f(\ub,\cdot)\|^2_{L^2(\ub,\tau)}$$
part {\bf b} then follows.

\vspace{2.5mm}

In regard to the contribution to $\|\s^{(n-1)}\ctchi\|_{L^2(\Cb_{\ub_1}^{u_1})}$ of the first of 
\ref{12.142}, we set:
\begin{equation}
f(\ub,u)=\ub u^{-2}\|\s^{(E;0,n)}\cxi_{\Lb}\|_{L^2(S_{\ub,u})}
\label{12.143}
\end{equation}
We have:
\begin{eqnarray}
&&\|f(\ub,\cdot)\|^2_{L^2(\ub_1,u_1)}=\int_{\ub_1}^{u_1}f^2(\ub,u)du\leq\int_{\ub}^{u_1}f^2(\ub,u)du 
\nonumber\\
&&=\ub^2\int_{\ub}^{u_1}u^ {-4}\|\s^{(E;0,n)}\cxi_{\Lb}\|^2_{L^2(S_{\ub,u})}du\nonumber\\
&&=\ub^2\int_{\ub}^{u_1}u^{-4}\frac{\partial}{\partial u}
\left(\int_{\ub}^u\|\s^{(E;0,n)}\cxi_{\Lb}\|^2_{L^2(S_{\ub,u^\prime})}du^\prime\right)du\nonumber\\
&&=\ub^2\left\{u_1^{-4}\int_{\ub}^{u_1}\|\s^{(E;0,n)}\cxi_{\Lb}\|^2_{L^2(S_{\ub,u})}du\right.\nonumber\\
&&\hspace{20mm}\left.+4\int_{\ub}^{u_1}u^{-5}\left(\int_{\ub}^u\|\s^{(E;0,n)}\cxi_{\Lb}
\|^2_{L^2(S_{\ub,u^\prime})}du^\prime\right)du\right\}\nonumber\\
&&\leq C\ub^2\left\{u_1^{-4}\s^{(E;0,n)}\cEb^{u_1}(\ub)+4\int_{\ub}^{u_1}u^{-5}
\s^{(E;0,n)}\cEb^u(\ub)du\right\}
\nonumber\\
&&\leq C\ub^{2a_0+2}\left\{u_1^{2b_0-4}+4\frac{(u_1^{2b_0-4}-\ub^{2b_0-4})}{(2b_0-4)}\right\}
\s^{(E;0,n)}\cBb(\ub_1,u_1)\nonumber\\
&&\leq C\ub^{2a_0+2}u_1^{2b_0-4}\s^{(E;0,n)}\cBb(\ub_1,u_1)
\label{12.144}
\end{eqnarray}
by \ref{9.291}, \ref{9.298} and the condition \ref{10.391}. Then by Lemma 12.1a the contribution to 
$\|\s^{(n-1)}\ctchi\|_{L^2(\Cb_{\ub_1}^{u_1})}$ of the first of \ref{12.142} is bounded by:
\begin{equation}
C\sqrt{\s^{(E;0,n)}\cBb(\ub_1,u_1)}u_1^{b_0-2}\int_0^{\ub_1}\ub^{a_0+1}d\ub
\leq\frac{C\sqrt{\s^{(E;0,n)}\cBb(\ub_1,u_1)}}{(a_0+2)}\cdot\ub_1^{a_0}u_1^{b_0}\cdot\frac{\ub_1^2}{u_1^2}
\label{12.145}
\end{equation}
On the other hand, setting $\ub_1=u$ in \ref{12.128} and taking the $L^2$ norm with respect to $u$ on 
$[0,\tau]$ yields: 
\begin{equation}
\|\s^{(n-1)}\ctchi\|_{L^2({\cal K}^\tau)}\leq 
k\left\{\int_0^\tau\left(\int_0^u\|\s^{(n-1)}\chF\|_{L^2(S_{\ub,u})}d\ub\right)^2 du\right\}^{1/2}
\label{12.146}
\end{equation}
Since by \ref{12.144} with $\tau$ in the role of $u_1$ we have:
\begin{equation}
\|f(\ub,\cdot)\|^2_{L^2(\ub,\tau)}=\int_{\ub}^\tau f^2(\ub,u)du
\leq C\ub^{2a_0+2}\tau^{2b_0-4}\s^{(E;0,n)}\cBb(\tau,\tau)
\label{12.147}
\end{equation}
by Lemma 12.1b the contribution to $\|\s^{(n-1)}\ctchi\|_{L^2({\cal K}^\tau)}$ of the first of 
\ref{12.142} is bounded by:
\begin{equation}
C\sqrt{\s^{(E;0,n)}\cBb(\tau,\tau)}\tau^{b_0-2}\int_0^{\tau}\ub^{a_0+1}d\ub
\leq\frac{C\sqrt{\s^{(E;0,n)}\cBb(\tau,\tau)}}{(a_0+2)}\cdot\tau^{a_0+b_0}
\label{12.148}
\end{equation}
In regard to the contribution to $\|\s^{(n-1)}\ctchi\|_{L^2(\Cb_{\ub_1}^{u_1})}$ of the third of 
\ref{12.142} we set:
\begin{equation}
f(\ub,u)=\ub u^{-2}\|\s^{(n-1)}\ctchi\|_{L^2(S_{\ub,u})}
\label{12.149}
\end{equation}
In analogy with \ref{12.144} we obtain:
\begin{equation}
\int_{\ub}^{u_1}f^2(\ub,u)du=\ub^2\left\{u_1^{-4}\|\s^{(n-1)}\ctchi\|^2_{L^2(\Cb_{\ub}^{u_1})}
+4\int_{\ub}^{u_1}u^{-5}\|\s^{(n-1)}\ctchi\|^2_{L^2(\Cb_{\ub}^u)}du\right\}
\label{12.150}
\end{equation}
hence applying Lemma 12.1a the contribution in question is seen to be bounded by:
\begin{equation}
\int_0^{\ub_1}\ub\left\{u_1^{-4}\|\s^{(n-1)}\ctchi\|^2_{L^2(\Cb_{\ub}^{u_1})}
+4\int_{\ub}^{u_1}u^{-5}\|\s^{(n-1)}\ctchi\|^2_{L^2(\Cb_{\ub}^u)}du\right\}^{1/2}d\ub
\label{12.151}
\end{equation}
Also, replacing $u_1$ by $\tau$ and applying Lemma 12.1b the corresponding contribution to 
$\|\s^{(n-1)}\ctchi\|_{L^2({\cal K}^\tau)}$ is seen to be bounded by: 
\begin{equation}
\int_0^{\tau}\ub\left\{\tau^{-4}\|\s^{(n-1)}\ctchi\|^2_{L^2(\Cb_{\ub}^{\tau})}
+4\int_{\ub}^{\tau}u^{-5}\|\s^{(n-1)}\ctchi\|^2_{L^2(\Cb_{\ub}^u)}du\right\}^{1/2}d\ub
\label{12.152}
\end{equation}

Consider now the first of the differences \ref{12.134}, the $E$ component difference. Up to terms 
which can be ignored this difference is:
\begin{equation}
-\beta_N^2\sbeta\eta^2 H^\prime \sqrt{\frac{\rho}{\lambda}}\sqrt{a}\s^{(E;0,n)}\csxi
\label{12.153}
\end{equation}
The contribution of this to the integral on the right in \ref{12.128} is bounded by:
\begin{equation}
\frac{C}{u}\int_0^{\ub_1}\sqrt{\ub}\|\sqrt{a}\s^{(E;0,n)}\csxi\|_{L^2(S_{\ub,u})}d\ub
\label{12.154}
\end{equation}
Setting then:
\begin{equation}
f(\ub,u)=\sqrt{\ub}u^{-1}\|\sqrt{a}\s^{(E;0,n)}\csxi\|_{L^2(S_{\ub,u})}
\label{12.155}
\end{equation}
we have:
\begin{eqnarray}
&&\int_{\ub}^{u_1}f^2(\ub,u)du=\ub\int_{\ub}^{u_1}u^{-2}
\|\sqrt{a}\s^{(E;0,n)}\csxi\|^2_{L^2(S_{\ub,u})}du \nonumber\\
&&=\ub\int_{\ub}^{u_1}u^{-2}\frac{\partial}{\partial u}
\left(\int_{\ub}^u\|\sqrt{a}\s^{(E;0,n)}\csxi\|^2_{L^2(S_{\ub,u^\prime})}du^\prime\right)du\nonumber\\
&&=\ub\left\{u_1^{-2}\int_{\ub}^{u_1}\|\sqrt{a}\s^{(E;0,n)}\csxi\|^2_{L^2(S_{\ub,u})}du\right.
\nonumber\\
&&\hspace{20mm}\left.+2\int_{\ub}^{u_1}u^{-3}\left(\int_{\ub}^u
\|\sqrt{a}\s^{(E;0,n)}\csxi\|^2_{L^2(S_{\ub,u^\prime})}du^\prime\right)du\right\}\nonumber\\
&&\leq C\ub\left\{u_1^{-2}\s^{(E;0,n)}\cEb^{u_1}(\ub)+2\int_{\ub}^{u_1}u^{-3}\cEb^u(\ub)du\right\}
\nonumber\\
&&\leq C\ub^{2a_0+1}\left\{u_1^{2b_0-2}+2\frac{(u_1^{2b_0-2}-\ub^{2b_0-2})}{(2b_0-2)}\right\}
\s^{(E;0,n)}\cBb(\ub_1,u_1)\nonumber\\
&&\leq C\ub^{2a_0+1}u_1^{2b_0-2}\s^{(E;0,n)}\cBb(\ub_1,u_1)
\label{12.156}
\end{eqnarray}
by \ref{9.291} and \ref{9.298}. Then by Lemma 12.1a the contribution to 
$\|\s^{(n-1)}\ctchi\|_{L^2(\Cb_{\ub_1}^{u_1})}$ of \ref{12.154} is bounded by:
\begin{equation}
C\sqrt{\s^{(E;0,n)}\cBb(\ub_1,u_1)}u_1^{b_0-1}\int_0^{\ub_1}\ub^{a_0+\frac{1}{2}}d\ub
\leq\frac{C\sqrt{\s^{(E;0,n)}\cBb(\ub_1,u_1)}}{\left(a_0+\frac{3}{2}\right)}\cdot\ub_1^{a_0}u_1^{b_0}\cdot\frac{\ub_1^{3/2}}{u_1}
\label{12.157}
\end{equation}
Here the last factor, $\ub_1^{3/2}/u_1$, is homogeneous of degree $1/2$. This is to be compared with 
\ref{12.145} where the last factor $\ub_1^2/u_1^2$ is homogeneous of degree 0. 

Consider finally the second of the differences \ref{12.134}, the $N$ component difference. 
Following a sequence of formulas which are the conjugates of the formulas in the sequence 
\ref{12.135} - \ref{12.140} we conclude that 
$$\rho N^\mu E^{n+1}\beta_\mu$$
is, up to terms which can be ignored, given by:
\begin{equation}
E^\mu L E^n\beta_\mu-\frac{1}{2}\rho\sss E^{n-1}\tchi
-\left(\frac{1}{2c}s_{NL}+\pi\rho\ss_N\right)E^{n-1}\tchib-\frac{1}{c}\ss_N E^n\lambdab
\label{12.158}
\end{equation}
Similar formulas hold for the $N$th approximants. It follows that up to terms which can be ignored 
the second the differences \ref{12.134} is given by:
\begin{equation}
\beta_N^2\frac{\beta_{\Nb}}{2c}\eta^2 H^\prime\left\{\s^{(E;0,n)}\cxi_L
-\frac{1}{2}\rho\sss\s^{(n-1)}\ctchi-\left(\frac{1}{2c}s_{NL}+\pi\rho\ss_N\right)\s^{(n-1)}\ctchib-\frac{1}{c}\ss_N \s^{(0,n)}\clab\right\}
\label{12.159}
\end{equation}
The contribution of this to the integral on the right in \ref{12.128} is bounded by:
\begin{eqnarray}
&&C\int_0^{\ub_1}\|\s^{(E;0,n)}\cxi_L\|_{L^2(S_{\ub,u})}d\ub
+C\int_0^{\ub_1}\ub\|\s^{(n-1)}\ctchi\|_{L^2(S_{\ub,u})}d\ub\nonumber\\
&&+C\int_0^{\ub_1}\ub\|\s^{(n-1)}\tchib\|_{L^2(S_{\ub,u})}d\ub
+C\int_0^{\ub}\|\s^{(0,n)}\clab\|_{L^2(S_{\ub,u})}d\ub
\label{12.160}
\end{eqnarray}
Here the 2nd integral is the same as the 3rd integral in \ref{12.142} but with an extra factor of $u^2$, hence its contribution can be absorbed. The 3rd integral is the same 
as the 2nd integral in \ref{12.142}. As for the contribution of 
the 4th integral in \ref{12.160}, we shall 
show in the sequel that this contribution, like that of the 3rd integral, can be absorbed. Finally, 
the 1st integral in \ref{12.160} is bounded by:
\begin{eqnarray}
&&\ub_1^{1/2}\|\s^{(E;0,n)}\cxi_L\|_{L^2(C_u^{\ub_1})}\leq C\ub_1^{1/2}\sqrt{\s^{(E;0,n)}\cE^{\ub_1}(u)}
\nonumber\\
&&\leq C\ub_1^{a_0+\frac{1}{2}}u^{b_0}\sqrt{\s^{(E;0,n)}\cB(\ub_1,u_1)}
\label{12.161}
\end{eqnarray}
hence its contribution to $\|\s^{(n-1)}\ctchi\|_{L^2(\Cb_{\ub_1}^{u_1})}$ is bounded by:
\begin{equation}
\frac{C\sqrt{\s^{(E;0,n)}\cB(\ub_1,u_1)}}{\sqrt{2b_0+1}}\cdot \ub_1^{a_0} u_1^{b_0}\cdot 
\ub_1^{1/2}u_1^{1/2}
\label{12.162}
\end{equation}
Note that here the last factor, $\ub_1^{1/2}u_1^{1/2}$, is homogeneous of degree 1. 

Consider next the 1st of the principal difference terms \ref{12.132}. Since 
$$LH=-2H^\prime (g^{-1})^{\mu\nu}\beta_\nu L\beta_\mu, \ \ \mbox{hence} \ \ 
[E^n LH]_{P.P.}=-2H^\prime (g^{-1})^{\mu\nu}\beta_\nu E^n L\beta_\mu$$
and similarly for the $N$th approximants, the actual principal difference term in question is:
\begin{equation}
2\beta_N\sbeta H^\prime (g^{-1})^{\mu\nu}\beta_\nu E^n L\beta_\mu 
-2\beta_{N,N}\sbeta_N H^\prime_N (g^{-1})^{\mu\nu}\beta_{\nu,N}E_N^n L_N\beta_{\mu,N}
\label{12.163}
\end{equation}
By the formula \ref{10.368} and the analogous formula for the $N$th approximants, \ref{12.163} is 
expressed as the sum of differences of $E$, $N$, and $\Nb$ components:
\begin{eqnarray}
&&2\beta_N\sbeta^2\eta^2 H^\prime E^\mu E^n L\beta_\mu 
-2\beta_{N,N}\sbeta_N^2\eta_N^2 H^\prime_N E_N^\mu E_N^n L_N\beta_{\mu,N}\nonumber\\
&&-\beta_N\sbeta\frac{\beta_{\Nb}}{c}\eta^2 H^\prime N^\mu E^n L\beta_\mu 
+\beta_{N,N}\sbeta_N\frac{\beta_{\Nb,N}}{c_N}\eta_N^2 H_N^\prime N_N^\mu E_N^n L_N\beta_{\mu,N}
\nonumber\\
&&-\beta_N\sbeta\frac{\beta_N}{c}\eta^2 H^\prime\Nb^\mu E^n L\beta_\mu 
-\beta_{N,N}\sbeta_N\frac{\beta_{N,N}}{c_N}\eta_N^2 H^\prime_N\Nb_N^\mu E_N^n L_N\beta_{\mu,N}
\nonumber\\&&\label{12.164}
\end{eqnarray}
To estimate the contribution of the $\Nb$ component difference we write: 
\begin{equation}
\Nb^\mu E^n L\beta_\mu=E^n(\Nb^\mu L\beta_\mu)-(E^n\Nb^\mu)L\beta_\mu-\sum_{i=1}^{n-1}\left(\begin{array}{c} n\\i \end{array}\right)(E^i\Nb^\mu)
E^{n-i}L\beta_\mu \label{12.165}
\end{equation}
Here the sum is of order $n$ with vanishing principal acoustical part, hence can be ignored. 
As for the 2nd term on the right, by the conjugate of \ref{12.7} we have:
\begin{equation}
[(E^n\Nb^\mu)L\beta_\mu]_{P.A.}=\rho(\ss_N-\pi c\sss)E^{n-1}\tchib
\label{12.166}
\end{equation}
Since $\Nb^\mu L\beta_\mu=\lambdab E^\mu E\beta_\mu$, the 1st term on the right in \ref{12.165} is:
\begin{eqnarray}
&&E^n(\lambdab E^\mu E\beta_\mu)=\lambdab E^\mu E^{n+1}\beta_\mu+(E^n\lambdab)\sss
+\lambdab(E^n E^\mu)E\beta_\mu \label{12.167}\\
&&+\sum_{i=1}^{n-1}\left(\begin{array}{c} n\\i \end{array}\right)
(E^i\lambdab)(E^{n-i}E^\mu)E\beta_\mu 
+\sum_{i=1}^{n-1}\left(\begin{array}{c} n\\i \end{array}\right)
(E^i (\lambdab E^\mu))E^{n+1-i}\beta_\mu \nonumber
\end{eqnarray}
Here the two sums are of order $n$ with vanishing principal acoustical part, hence can be ignored. 
As for the 3rd term on the right, by \ref{12.8} we have:
\begin{equation}
[(E^n E^\mu)E\beta_\mu]_{P.A.}=\frac{1}{2c}(\ss_{\Nb}E^{n-1}\tchi+\ss_N E^{n-1}\tchib)
\label{12.168}
\end{equation}
We conclude from the above that
$$\Nb^\mu E^n L\beta_\mu$$
is, up to terms which can be ignored, given by:
\begin{equation}
\lambdab E^\mu E^{n+1}\beta_\mu+\sss E^n\lambdab-\rho\left(\frac{\ss_N}{2}-\pi c\sss\right)E^{n-1}\tchib
+\rho\frac{\ss_{\Nb}}{2}E^{n-1}\tchi
\label{12.169}
\end{equation}
Similar formulas hold for the $N$th approximants. It follows that up to terms which can be ignored 
the last the differences \ref{12.164} is given by:
\begin{equation}
-\beta_N\sbeta\frac{\beta_N}{c}\eta^2 H^\prime\left\{\lambdab \s^{(E;0,n)}\csxi
+\sss \s^{(0,n)}\clab-\rho\left(\frac{\ss_N}{2}-\pi c\sss\right)\s^{(n-1)}\ctchib 
+\rho\frac{\ss_{\Nb}}{2}\s^{(n-1)}\ctchi\right\}
\label{12.170}
\end{equation}
The contribution of this to the integral on the right in \ref{12.128} is bounded by:
\begin{eqnarray}
&&\frac{C}{u}\int_0^{\ub_1}\sqrt{\ub}\|\sqrt{a}\s^{(E;0,n)}\csxi\|_{L^2(S_{\ub,u})}d\ub
+C\int_0^{\ub_1}\|\s^{(0,n)}\clab\|_{L^2(S_{\ub,u})}d\ub\nonumber\\
&&+C\int_0^{\ub_1}\ub\|\s^{(n-1)}\ctchib\|_{L^2(S_{\ub,u})}d\ub
+C\int_0^{\ub_1}\ub\|\s^{(n-1)}\ctchi\|_{L^2(S_{\ub,u})}d\ub\nonumber\\
&&\label{12.171}
\end{eqnarray}
Here the 1st term coincides with \ref{12.154}. The 2nd term coincides with the 4th term 
in \ref{12.160}. The 3rd term coincides with the 2nd term in \ref{12.142}. 
Comparing the 4th term with the 3rd term in \ref{12.142} we see that the contribution of this term here 
can be absorbed. 

Consider next the first of the differences \ref{12.164}, the $E$ component difference. To estimate this 
we use the 1st of the commutation relations \ref{3.a14} to express $E^\mu E^n L\beta_\mu$ as:
\begin{eqnarray}
&&E^\mu LE^n\beta_\mu+\sum_{i=0}^{n-1}E^\mu E^{n-1-i}(\chi E^{i+1}\beta_\mu) \nonumber\\
&&=E^\mu LE^n\beta_\mu+\sss E^{n-1}\chi+\sum_{j=1}^{n-1}\left(\begin{array}{c} n-1\\j\end{array}
\right)(E^{n-1-j}\chi)E^\mu E^{j+1}\beta_\mu\nonumber\\
&&\hspace{20mm}+\sum_{i=1}^{n-1}E^\mu E^{n-1-i}(\chi E^{i+1}\beta_\mu)\label{12.172}
\end{eqnarray}
Here the two sums on the right are of order $n$ with vanishing principal acoustical part, hence can be 
ignored. Since 
$$[E^{n-1}\chi]_{P.A.}=\rho E^{n-1}\tchi$$
we conclude from that, up to terms which can be ignored, $E^\mu E^n L\beta_\mu$ is:
\begin{equation}
E^\mu LE^n\beta_\mu+\rho\sss E^{n-1}\tchi
\label{12.173}
\end{equation}
A similar formula holds for the $N$th approximants. It follows that up to terms which can be ignored 
the first the differences \ref{12.164} is given by:
\begin{equation}
2\beta_N\sbeta^2\eta^2 H^\prime\left(\s^{(E;0,n)}\cxi_L+\rho\sss\s^{(n-1)}\ctchi\right)
\label{12.174}
\end{equation}
The contribution of this to the integral on the right in \ref{12.128} is bounded by:
\begin{equation}
C\int_0^{\ub_1}\|\s^{(E;0,n)}\cxi_L\|_{L^2(S_{\ub,u})}d\ub
+C\int_0^{\ub_1}\ub\|\s^{(n-1)}\ctchi\|_{L^2(S_{\ub,u})}d\ub
\label{12.175}
\end{equation}
Here the 1st term coincides with the 1st term in \ref{12.160} while the contribution of the 2nd term 
can be absorbed. 

Consider finally the second of the differences \ref{12.164}, the $N$ component difference. To estimate 
this we use the 1st of the commutation relations \ref{3.a14} to express $N^\mu E^n L\beta_\mu$ as:
\begin{eqnarray}
&&N^\mu LE^n\beta_\mu+\sum_{i=0}^{n-1}N^\mu E^{n-1-i}(\chi E^{i+1}\beta_\mu) \nonumber\\
&&=N^\mu LE^n\beta_\mu+\ss_N E^{n-1}\chi+\sum_{j=1}^{n-1}\left(\begin{array}{c} n-1\\j\end{array}
\right)(E^{n-1-j}\chi)N^\mu E^{j+1}\beta_\mu\nonumber\\
&&\hspace{20mm}+\sum_{i=1}^{n-1}N^\mu E^{n-1-i}(\chi E^{i+1}\beta_\mu)\label{12.176}
\end{eqnarray}
Again the two sums are of order $n$ with vanishing principal acoustical part and can be ignored. 
Then $N^\mu E^n L\beta_\mu$ is, up to terms which can be ignored, 
\begin{equation}
N^\mu LE^n\beta_\mu+\rho\ss_N E^{n-1}\tchi
\label{12.177}
\end{equation}
A similar formula holds for the $N$th approximants. It follows that up to terms which can be ignored 
the second of the differences \ref{12.164} is given by:
\begin{equation}
-\beta_N\sbeta\frac{\beta_{\Nb}}{c}\eta^2 H^\prime
\left\{N^\mu L(E^n\beta_\mu-E_N^n\beta_{\mu,N})+\rho\ss_N\s^{(n-1)}\ctchi\right\}
\label{12.178}
\end{equation}
The contribution of the 2nd term in parenthesis to the integral on the right in \ref{12.128} is 
bounded by:
\begin{equation}
C\int_0^{\ub_1}\ub\|\s^{(n-1)}\ctchi\|_{L^2(S_{\ub,u})}d\ub
\label{12.179}
\end{equation}
which coincides with the 2nd term in \ref{12.175} and can be absorbed. 
To estimate the contribution of the 1st term in parenthesis in \ref{12.178} we use the fact that 
\begin{equation}
N^\mu L(E^n\beta_\mu-E_N^n\beta_{\mu,N})+\ogamma\Nb^\mu L(E^n\beta_\mu-E_N^n\beta_{\mu,N})
=\s^{(Y;0,n)}\cxi_L
\label{12.180}
\end{equation}
In regard to the 2nd term on the left, we can express $\Nb^\mu L(E^n\beta_\mu-E_N^n\beta_{\mu,N})$ as
$$\Nb^\mu E^n L\beta_\mu-\Nb_N^\mu E_N^n L_N\beta_{\mu,N} -\rho\ss_{\Nb}E^{n-1}\tchi
+\rho_N\ss_{\Nb,N}E_N^{n-1}\tchi_N$$
up to terms which can be ignored. Since $\ogamma\sim u$, the contribution of the 1st of these
differences is bounded  by $u$ times \ref{12.171} while the contribution of the 2nd difference is 
bounded by $u$ times \ref{12.179}. Therefore the contribution of the 2nd term on the left in 
\ref{12.180} can be absorbed. We are left with the contribution of the right hand side of \ref{12.180}, 
through the 1st term in parenthesis in \ref{12.178}, to the integral on the right in \ref{12.128}. 
This contribution is bounded by a constant multiple of:
\begin{eqnarray}
&&\int_0^{\ub_1}\|\s^{(Y;0,n)}\cxi_L\|_{L^2(S_{\ub,u})}d\ub\leq 
\ub_1^{1/2}\|\s^{(Y;0,n)}\cxi_L\|_{L^2(C_u^{\ub_1})}\nonumber\\
&&\leq C\ub_1^{1/2}\sqrt{\s^{(Y;0,n)}\cE^{\ub_1}(u)}
\leq C\ub_1^{a_0+\frac{1}{2}}u^{b_0}\sqrt{\s^{(Y;0,n)}\cB(\ub_1,u_1)} \label{12.181}
\end{eqnarray}
hence its contribution to $\|\s^{(n-1)}\ctchi\|_{L^2(\Cb_{\ub_1}^{u_1})}$ is bounded by:
\begin{equation}
\frac{C\sqrt{\s^{(Y;0,n)}\cB(\ub_1,u_1)}}{\sqrt{2b_0+1}}\cdot \ub_1^{a_0}u_1^{b_0}\cdot 
\ub_1^{1/2}u_1^{1/2}
\label{12.182}
\end{equation}
(compare with \ref{12.162}). 

Now, in regard to the contribution of the integral 
$$\int_0^{\ub_1}\ub\|\s^{(n-1)}\ctchib\|_{L^2(S_{\ub,u})}d\ub$$ 
through the 2nd of \ref{12.131}, \ref{12.142} and the 3rd of \ref{12.160}, \ref{12.171} to 
$\|\s^{(n-1)}\ctchi\|_{L^2(\Cb_{\ub_1}^{u_1})}$, we remark that this is bounded by a constant multiple 
of:
\begin{equation}
\ub_1^{1/2}u_1^{1/2}\sup_{u\in[\ub_1,u_1]}\|\ub\s^{(n-1)}\ctchib\|_{L^2(C_u^{\ub_1})}
\label{12.183}
\end{equation}
and the corresponding contribution to $\|\s^{(n-1)}\ctchi\|_{L^2({\cal K}^\tau)}$ by a constant 
multiple of: 
\begin{equation}
\tau\sup_{u\in[0,\tau]}\|\ub\s^{(n-1)}\ctchib\|_{L^2(C_u^u)}
\label{12.184}
\end{equation}
Similarly, in regard to the contribution of the integral 
$$\int_0^{\ub_1}\|\s^{(0,n)}\clab\|_{L^2(S{\ub,u})}d\ub$$
through the 4th of \ref{12.160} and the 2nd \ref{12.171} to 
$\|\s^{(n-1)}\ctchi\|_{L^2(\Cb_{\ub_1}^{u_1})}$, we remark that this is bounded by a constant multiple 
of:
\begin{equation}
\ub_1^{1/2}u_1^{1/2}\sup_{u\in[\ub_1,u_1]}\|\s^{(n-1)}\clab\|_{L^2(C_u^{\ub_1})}
\label{12.185}
\end{equation}
and the corresponding contribution to $\|\s^{(n-1)}\ctchi\|_{L^2({\cal K}^\tau)}$ by a constant 
multiple of: 
\begin{equation}
\tau\sup_{u\in[0,\tau]}\|\s^{(n-1)}\clab\|_{L^2(C_u^u)}
\label{12.186}
\end{equation}
Also, in regard to the contribution of 
$$\frac{1}{u^2}\int_0^{\ub_1}\ub\|\s^{(0,n)}\cla\|_{L^2(S_{\ub,u})}d\ub$$
through the 4th of \ref{12.142} to 
$\|\s^{(n-1)}\ctchi\|_{L^2(\Cb_{\ub_1}^{u_1})}$, setting 
$$f(\ub,u)=\ub u^{-2}\|\s^{(0,n)}\cla\|_{L^2(S_{\ub,u})}$$
we deduce in analogy with \ref{12.149} - \ref{12.152} that this contribution is bounded by a constant 
multiple of:
\begin{equation}
\int_0^{\ub_1}\ub\left\{u_1^{-4}\|\s^{(0,n)}\cla\|^2_{L^2(\Cb_{\ub}^{u_1})}
+4\int_{\ub}^{u_1}u^{-5}\|\s^{(0,n)}\cla\|^2_{L^2(\Cb_{\ub}^u)}du\right\}^{1/2}d\ub
\label{12.187}
\end{equation}
and the corresponding contribution to $\|\s^{(n-1)}\ctchi\|_{L^2({\cal K}^\tau)}$ by a constant 
multiple of:
\begin{equation}
\int_0^{\tau}\ub\left\{\tau^{-4}\|\s^{(0,n)}\cla\|^2_{L^2(\Cb_{\ub}^{\tau})}
+4\int_{\ub}^{\tau}u^{-5}\|\s^{(0,n)}\cla\|^2_{L^2(\Cb_{\ub}^u)}du\right\}^{1/2}d\ub
\label{12.188}
\end{equation}

We combine the above results in the form in which they will be used in the sequel in the following 
lemma. 

\vspace{2.5mm}

\noindent{\bf Lemma 12.2} \ \ The next to the top order acoustical quantity 
$\|\s^{(n-1)}\ctchi\|_{L^2(\Cb_{\ub}^u)}$ satisfies to principal terms the following inequality:
\begin{eqnarray*}
&&\|\s^{(n-1)}\ctchi\|_{L^2(\Cb_{\ub_1}^{u_1})}\leq \\
&&\hspace{10mm}C\int_0^{\ub_1}\ub\left\{u_1^{-4}\|\s^{(n-1)}\ctchi\|^2_{L^2(\Cb_{\ub}^{u_1})}
+4\int_{\ub}^{u_1}u^{-5}\|\s^{(n-1)}\ctchi\|^2_{L^2(\Cb_{\ub}^u)}du\right\}^{1/2}d\ub\\
&&\hspace{10mm}+C\int_0^{\ub_1}\ub\left\{u_1^{-4}\|\s^{(0,n)}\cla\|^2_{L^2(\Cb_{\ub}^{u_1})}
+4\int_{\ub}^{u_1}u^{-5}\|\s^{(0,n)}\cla\|^2_{L^2(\Cb_{\ub}^u)}du\right\}^{1/2}d\ub\\
&&\hspace{28mm}+C\ub_1^{1/2}u_1^{1/2}\sup_{u\in[\ub_1,u_1]}\|\ub\s^{(n-1)}\ctchib\|_{L^2(C_u^{\ub_1})}\\
&&\hspace{28mm}+C\ub_1^{1/2}u_1^{1/2}\sup_{u\in[\ub_1,u_1]}\|\s^{(n-1)}\clab\|_{L^2(C_u^{\ub_1})}\\
&&\hspace{28mm}+C\frac{\sqrt{\s^{(E;0,n)}\cBb(\ub_1,u_1)}}{(a_0+2)}\cdot\ub_1^{a_0}u_1^{b_0}
\cdot\frac{\ub_1^2}{u_1^2}\\
&&\hspace{28mm}+C\sqrt{\max\{\s^{(0,n)}\cB(\ub_1,u_1),\s^{(0,n)}\cBb(\ub_1,u_1)\}}
\cdot\ub_1^{a_0}u_1^{b_0}\cdot \ub_1^{1/2}
\end{eqnarray*}
Moreover, $\|\s^{(n-1)}\ctchi\|_{L^2({\cal K}^\tau)}$ satisfies:
\begin{eqnarray*}
&&\|\s^{(n-1)}\ctchi\|_{L^2({\cal K}^\tau)}\leq \\
&&\hspace{10mm}C\int_0^{\tau}\ub\left\{\tau^{-4}\|\s^{(n-1)}\ctchi\|^2_{L^2(\Cb_{\ub}^{\tau})}
+4\int_{\ub}^{\tau}u^{-5}\|\s^{(n-1)}\ctchi\|^2_{L^2(\Cb_{\ub}^u)}du\right\}^{1/2}d\ub\\
&&\hspace{10mm}+C\int_0^{\tau}\ub\left\{\tau^{-4}\|\s^{(0,n)}\cla\|^2_{L^2(\Cb_{\ub}^{\tau})}
+4\int_{\ub}^{\tau}u^{-5}\|\s^{(0,n)}\cla\|^2_{L^2(\Cb_{\ub}^u)}du\right\}^{1/2}d\ub\\
&&\hspace{28mm}+C\tau\sup_{u\in[0,\tau]}\|\ub\s^{(n-1)}\ctchib\|_{L^2(C_u^u)}\\
&&\hspace{28mm}+C\tau\sup_{u\in[0,\tau]}\|\s^{(n-1)}\clab\|_{L^2(C_u^u)}\\
&&\hspace{28mm}+C\frac{\sqrt{\s^{(E;0,n)}\cBb(\tau,\tau)}}{(a_0+2)}\cdot\tau^{a_0+b_0}\\
&&\hspace{28mm}+C\sqrt{\max\{\s^{(0,n)}\cB(\tau,\tau),\s^{(0,n)}\cBb(\tau,\tau)\}}
\cdot\tau^{a_0+b_0}\cdot \tau^{1/2}
\end{eqnarray*}

\vspace{2.5mm}

We turn to the derivation of the estimates for $\s^{(0,n)}\cla$. We consider the propagation equation 
for $\s^{(0,n)}\cla$, equation \ref{12.101}. This equation takes the form:
\begin{equation}
L\s^{(0,n)}\cla+n\chi\s^{(0,n)}\s^{(0,n)}\cla=\s^{(0,n)}\chG
\label{12.189}
\end{equation}
In terms of the $(\ub,u,\vartheta^\prime)$ coordinates on ${\cal N}$, which are adapted to the flow of 
$L$, this reads:
\begin{equation}
\frac{\partial}{\partial\ub}\left(\sh^{\prime n/2}\s^{(0,n)}\cla\right)=\sh^{\prime n/2}\s^{(0,n)}\chG
\label{12.190}
\end{equation}
Integrating from $\Cb_0$, noting that
\begin{equation}
\left.\s^{(0,n)}\cla\right|_{\Cb_0}=0
\label{12.191}
\end{equation}
we obtain:
\begin{equation}
\left(\sh^{\prime n/2}\s^{(0,n)}\cla\right)(\ub_1,u,\vartheta^\prime)=
\int_0^{\ub_1}\left(\sh^{\prime n/2}\s^{(0,n)}\chG\right)(\ub,u,\vartheta^\prime)d\ub
\label{12.192}
\end{equation}
In view of the remark \ref{10.352} we rewrite this in the form:
\begin{eqnarray*}
&&\left(\s^{(0,n)}\cla\sh^{\prime 1/4}\right)(\ub_1,u,\vartheta^\prime)=\\
&&\hspace{20mm}\int_0^{\ub_1}\left(\frac{\sh^\prime(\ub,u,\vartheta^\prime)}
{\sh^\prime(\ub_1,u,\vartheta^\prime)}\right)^{(2n-1)/4}\left(\s^{(0,n)}\chG\sh^{\prime 1/4}\right)
(\ub,u,\vartheta^\prime)d\ub
\end{eqnarray*}
By \ref{10.356} we have:
\begin{equation}
\left(\frac{\sh^\prime(\ub,u,\vartheta^\prime)}
{\sh^\prime(\ub_1,u,\vartheta^\prime)}\right)^{(2n-1)/4}\leq k
\label{12.193}
\end{equation}
where $k$ is a constant greater than 1, but which can be chosen as close to 1 as we wish by suitably 
restricting $\ub_1$. Therefore \ref{12.192} implies:
\begin{equation}
\left|\left(\s^{(0,n)}\cla\sh^{\prime 1/4}\right)(\ub_1,u,\vartheta^\prime)\right|\leq k
\int_0^{\ub_1}\left|\left(\s^{(0,n)}\chG\sh^{\prime 1/4}\right)(\ub,u,\vartheta^\prime)\right|d\ub
\label{12.194}
\end{equation}
Taking $L^2$ norms with respect to $\vartheta^\prime\in S^1$ this implies:
\begin{equation}
\left\|\left(\s^{(0,n)}\cla\sh^{\prime 1/4}\right)(\ub_1,u)\right\|_{L^2(S^1)}\leq k
\int_0^{\ub_1}\left\|\left(\s^{(0,n)}\chG\sh^{\prime 1/4}\right)(\ub,u)\right\|_{L^2(S^1)}d\ub
\label{12.195}
\end{equation}
or, in view of \ref{10.352}:
\begin{equation}
\|\s^{(0,n)}\cla\|_{L^2(S_{\ub_1,u})}\leq k\int_0^{\ub_1}\|\s^{(0,n)}\chG\|_{L^2(S_{\ub,u})}d\ub
\label{12.196}
\end{equation}

We shall presently estimate the leading contributions to the integral on the right in \ref{12.196}. 
Consider first the $n$th order acoustical difference terms:
\begin{equation}
p\s^{(0,n)}\cla+\tilde{q}\s^{(0,n)}\clab+(B-\rho E\lambda)\s^{(n-1)}\ctchi+\oB\s^{(n-1)}\ctchib
\label{12.197}
\end{equation}
The coefficients $B$ and $\oB$ being given by \ref{12.52} and \ref{12.51} respectively, the 
assumptions \ref{10.439}, \ref{10.440} imply:
\begin{equation}
|B|, \ |\oB|  \leq C\ub \ \mbox{: in ${\cal R}_{\delta,\delta}$}
\label{12.198}
\end{equation}
It follows that the contribution of the terms \ref{12.197} is bounded by:
\begin{eqnarray}
&&C\int_0^{\ub_1}\ub\|\s^{(0,n)}\cla\|_{L^2(S_{\ub,u})}d\ub
+C\int_0^{\ub_1}\|\s^{(0,n)}\clab\|_{L^2(S_{\ub,u})}d\ub \nonumber\\
&&+C\int_0^{\ub_1}\ub\|\s^{(n-1)}\ctchi\|_{L^2(S_{\ub,u})}d\ub
+C\int_0^{\ub_1}\ub\|\s^{(n-1)}\ctchib\|_{L^2(S_{\ub,u})}d\ub \nonumber\\
&&\label{12.199}
\end{eqnarray}
Here the first can obviously be absorbed in the integral inequality \ref{12.196}. The third can 
likewise be absorbed. We shall see that the second gives a borderline contribution and the same 
holds for the fourth. 

The other leading contributions are those of the principal difference terms:
\begin{eqnarray}
&&\frac{\rho}{4}\beta_N^2 E^n\Lb H
-\frac{\rho_N}{4}\beta_{N,N}^2 E_N^n\Lb_N H_N \nonumber\\
&&+\rhob\beta_N\left(c\pi\sbeta-\frac{\beta_N}{4}-\frac{\beta_{\Nb}}{2}\right)E^n LH
\nonumber\\
&&\hspace{15mm}-\rhob_N\beta_{N,N}\left(c_N\pi_N\sbeta_N-\frac{\beta_{N,N}}{4}-\frac{\beta_{\Nb,N}}{2}\right)E_N^n L_N H_N\nonumber\\
&&+\rhob H\left(c\pi\sbeta-\frac{(\beta_N+\beta_{\Nb})}{2}\right)N^\mu E^n L\beta_\mu
\nonumber\\
&&\hspace{15mm}-\rhob_N H_N\left(c_N\pi_N\sbeta_N-\frac{(\beta_{N,N}+\beta_{\Nb,N})}{2}\right)
N_N^\mu E_N^n L_N\beta_{\mu,N}\nonumber\\
&&-\frac{\pi}{2}a\beta_N^2 E^{n+1} H 
+\frac{\pi_N}{2}a_N\beta_{N,N}^2 E_N^{n+1} H_N \label{12.200}
\end{eqnarray}
Here we must keep track not only of the actually principal terms, namely those containing derivatives 
of the $\beta_\mu$ of order $n+1$, but also of the terms containing acoustical quantities of order $n$. 
The terms of order $n$ with vanishing principal acoustical part can be ignored. As we shall see
the dominant contribution comes from the 1st of the four difference terms in \ref{12.200}. In fact, 
since 
$$\Lb H=-2H^\prime (g^{-1})^{\mu\nu}\beta_\nu \Lb\beta_\mu, \ \ \mbox{hence} \ \ 
[E^n\Lb H]_{P.P.}=-2H^\prime (g^{-1})^{\mu\nu}\beta_\nu E^n\Lb\beta_\mu$$
and similarly for the $N$th approximants, the actual principal difference term in question is:
\begin{equation}
-\frac{\rho}{2}\beta_N^2 H^\prime(g^{-1})^{\mu\nu}\beta_\nu E^n\Lb\beta_\mu 
+\frac{\rho_N}{2}\beta_{N,N}^2 H_N^\prime(g^{-1})^{\mu\nu}\beta_{\nu,N}E_N^n\Lb_N\beta_{\mu,N}
\label{12.201}
\end{equation}
By the formula \ref{10.368} and the analogous formula for the $N$th approximants, \ref{12.201} is 
expressed as the sum of differences of $E$, $N$, and $\Nb$ components:
\begin{eqnarray}
&&-\frac{1}{2}\beta_N^2\sbeta\eta^2 H^\prime\rho E^\mu E^n\Lb\beta_\mu 
+\frac{1}{2}\beta_{N,N}^2\sbeta_N\eta_N^2 H_N^\prime\rho_N E_N^\mu E_N^n\Lb_N\beta_{\mu,N}\nonumber\\
&&+\frac{1}{2}\beta_N^2\frac{\beta_{\Nb}}{2c}\eta^2 H^\prime\rho N^\mu E^n\Lb\beta_\mu 
-\frac{1}{2}\beta_{N,N}^2\frac{\beta_{\Nb,N}}{2c_N}\eta_N^2 H_N^\prime\rho_N N_N^\mu E_N^n\Lb_N\beta_{\mu,N}\nonumber\\
&&+\frac{1}{2}\beta_N^2\frac{\beta_N}{2c}\eta^2 H^\prime\rho\Nb^\mu E^n\Lb\beta_\mu 
-\frac{1}{2}\beta_{N,N}^2\frac{\beta_{N,N}}{2c_N}\eta_N^2 H_N^\prime\rho_N\Nb_N^\mu E_N^n\Lb_N\beta_{\mu,N} \nonumber\\
&&\label{12.202}
\end{eqnarray}
The dominant contribution is that of the $\Nb$ component difference. To estimate it we must use the 
fact that 
\begin{equation}
N^\mu\Lb(E^n\beta_\mu-E_N^n\beta_{\mu,N})+\ogamma\Nb^\mu\Lb(E^n\beta_\mu-E_N^n\beta_{\mu,N})
=\s^{(Y;0,n)}\cxi_{\Lb}
\label{12.203}
\end{equation}
Since the 1st term on the left is related to the $N$ component difference, the second of the 
differences \ref{12.202}, we begin with that difference. To estimate its contribution we follow 
the argument leading from \ref{12.165} to \ref{12.169} but with the original quantities replaced by 
their conjugates to conclude that 
$$N^\mu E^n\Lb\beta_\mu$$
is, up to terms which can be ignored, given by:
\begin{equation}
\lambda E^\mu E^{n+1}\beta_\mu+\sss E^n\lambda-\rhob\left(\frac{\ss_{\Nb}}{2}-\pi c\sss\right)
E^{n-1}\tchi+\rhob\frac{\ss_{\Nb}}{2}E^{n-1}\tchib
\label{12.204}
\end{equation}
Similarly for the $N$th approximants. It follows that up to terms which can be ignored the second of 
the differences \ref{12.202} is given by:
\begin{equation}
\frac{1}{2}\beta_N^2\frac{\beta_{\Nb}}{2c}\eta^2 H^\prime\rho\left\{\lambda\s^{(E;0,n)}\csxi
+\sss \s^{(0,n)}\cla-\rhob\left(\frac{\ss_{\Nb}}{2}-\pi c\sss\right)
\s^{(n-1)}\ctchi+\rhob\frac{\ss_{\Nb}}{2}\s^{(n-1)}\ctchib\right\}
\label{12.205}
\end{equation}
The contribution of this to the integral on the right in \ref{12.196} is bounded by:
\begin{eqnarray}
&&Cu\int_0^{\ub_1}\sqrt{\ub}\|\sqrt{a}\s^{(E;0,n)}\csxi\|_{L^2(S_{\ub,u})}d\ub
+C\int_0^{\ub_1}\ub\|\s^{(0,n)}\cla\|_{L^2(S_{\ub,u})}d\ub\nonumber\\
&&+Cu^2\int_0^{\ub_1}\ub\|\s^{(n-1)}\ctchi\|_{L^2(S_{\ub,u})}d\ub
+Cu^2\int_0^{\ub_1}\ub\|\s^{(n-1)}\ctchib\|_{L^2(S_{\ub,u})}d\ub \nonumber\\
&&\label{12.206}
\end{eqnarray}
Here the 2nd term coincides with the 1st term in \ref{12.199} and can be absorbed. Comparing the 
3rd and 4th terms with the corresponding terms in \ref{12.199} we see that the contributions of these 
terms here can be absorbed. As for the 1st term, setting 
\begin{equation}
f(\ub,u)=\sqrt{\ub} u\|\sqrt{a}\s^{(E;0,n)}\csxi\|_{L^2(S_{\ub,u})}
\label{12.207}
\end{equation}
we have:
\begin{eqnarray}
&&\|f(\ub,\cdot)\|^2_{L^2(\ub_1,u_1)}=\int_{\ub_1}^{u_1}f^2(\ub,u)du\leq \int_{\ub}^{u_1}f^2(\ub,u)du
\nonumber\\
&&\leq\ub u_1^2\int_{\ub}^{u_1}\|\sqrt{a}\s^{(E;0,n)}\csxi\|^2_{L^2(S_{\ub,u})}du
=\ub u_1^2\|\sqrt{a}\s^{(E;0,n)}\sxi\|^2_{L^2(\Cb_{\ub}^{u_1})}\nonumber\\
&&\leq C\ub u_1^2\s^{(E;0,n)}\cEb^{u_1}(\ub)\leq C\ub^{2a_0+1}u_1^{2b_0+2}\s^{(E;0,n)}\cBb(\ub_1,u_1)
\label{12.208}
\end{eqnarray}
Then by Lemma 12.1a the contribution to $\|\s^{(0,n)}\cla\|_{L^2(\Cb_{\ub_1}^{u_1})}$ of the 1st of 
\ref{12.206} is bounded by:
\begin{equation}
C\sqrt{\s^{(E;0,n)}\cBb(\ub_1,u_1)} u_1^{b_0+1}\int_0^{\ub_1}\ub^{a_0+\frac{1}{2}}d\ub
\leq \frac{C\sqrt{\s^{(E;0,n)}\cBb(\ub_1,u_1)}}{\left(a_0+\frac{3}{2}\right)}\cdot 
\ub_1^{a_0+\frac{3}{2}}u_1^{b_0+1}
\label{12.209}
\end{equation}
On the other hand, setting $u_1=u$ in \ref{12.196} and taking the $L^2$ norm with respect to $u$ on 
$[0,\tau]$ yields:
\begin{equation}
\|\s^{(0,n)}\cla\|_{L^2({\cal K}^\tau)}\leq k\left\{\int_0^\tau
\left(\int_0^u\|\s^{(0,n)}\chG\|_{L^2(S_{\ub,u})}d\ub\right)^2 du\right\}^{1/2}
\label{12.210}
\end{equation}
Since by \ref{12.208} with $\tau$ in the role of $u_1$ we have:
\begin{equation}
\|f(\ub,\cdot)\|^2_{L^2(\ub,\tau)}=\int_{\ub}^{\tau}f^2(\ub,u)du\leq C\ub^{2a_0+1}\tau^{2b_0+2}
\s^{(E;0,n)}\cBb(\tau,\tau)
\label{12.211}
\end{equation}
by Lemma 12.1b the contribution to $\|\s^{(0,n)}\cla\|_{L^2({\cal K}^{\tau})}$ of 1st of \ref{12.206} 
is bounded by:
\begin{equation}
C\sqrt{\s^{(E;0,n)}\cBb(\tau,\tau)} \tau^{b_0+1}\int_0^{\tau}\ub^{a_0+\frac{1}{2}}d\ub
\leq\frac{C\sqrt{\s^{(E;0,n)}\cBb(\tau,\tau)}}{\left(a_0+\frac{3}{2}\right)}\cdot 
\tau^{a_0+b_0+\frac{5}{2}}
\label{12.212}
\end{equation}

We turn to the last of the differences \ref{12.202}, the $\Nb$ component difference. By the 2nd of the 
commutation relations \ref{3.a14} we have:
\begin{eqnarray}
&&E^n\Lb\beta_\mu=\Lb E^n\beta_\mu+\sum_{i=0}^{n-1}E^{n-1-i}(\chib E^{i+1}\beta_\mu)\nonumber\\
&&\hspace{12mm}=\Lb E^n\beta_\mu+(E\beta_\mu)E^{n-1}\chib
+\sum_{j=1}^{n-1}\left(\begin{array}{c} n-1\\j\end{array}
\right)(E^{n-1-j}\chib)E^{j+1}\beta_\mu\nonumber\\
&&\hspace{28mm}+\sum_{i=1}^{n-1}(\chib E^{i+1}\beta_\mu)
\label{12.213}
\end{eqnarray}
The two sums on the right can be ignored being of order $n$ with vanishing principal acoustical part. 
Hence 
$$\Nb^\mu E^n\Lb\beta_\mu$$
is equal to:
\begin{equation}
\Nb^\mu\Lb E^n\beta_\mu+\ss_{\Nb}\rhob E^{n-1}\tchib
\label{12.214}
\end{equation}
up to terms which can be ignored. Similarly for the $N$th approximants. It follows that up to terms 
which can be ignored the last of the differences \ref{12.202} is given by:
\begin{equation}
\frac{\beta_N^3}{4c}\eta^2 H^\prime\rho\left\{\Nb^\mu\Lb(E^n\beta_\mu-E_N^n\beta_{\mu,N})
+\rhob\ss_{\Nb}\s^{(n-1)}\ctchib\right\}
\label{12.215}
\end{equation}
The contribution of this to the integral on the right in \ref{12.196} is bounded by:
\begin{eqnarray}
&&C\int_0^{\ub_1}\ub\|\Nb^\mu\Lb(E^n\beta_\mu-E_N^n\beta_{\mu,N})\|_{L^2(S_{\ub,u})}d\ub\nonumber\\
&&+Cu^2\int_0^{\ub_1}\ub\|\s^{(n-1)}\ctchib\|_{L^2(S_{\ub,u})}d\ub
\label{12.216}
\end{eqnarray}
Here the 2nd term coincides with the 4th term in \ref{12.206} and can be absorbed. To estimate the 
1st term we appeal to \ref{12.203}. Since $\ogamma\sim u$, this term is then seen to be bounded by:
\begin{eqnarray}
&&Cu^{-1}\int_0^{\ub_1}\ub\|\s^{(Y;0,n)}\cxi_{\Lb}\|_{L^2(S_{\ub,u})}d\ub \nonumber\\
&&+Cu^{-1}\int_0^{\ub_1}\ub\|N^\mu\Lb(E^n\beta_\mu-E_N^n\beta_{\mu,N})\|_{L^2(S_{\ub,u})}d\ub
\label{12.217}
\end{eqnarray}
In regard to the 2nd term here we can express $N^\mu\Lb(E^n\beta_\mu-E_N^n\beta_{\mu,N})$ as 
$$N^\mu E^n\Lb\beta_\mu-N_N^\mu E_N^n\Lb_N\beta_{\mu,N}-\rhob\ss_N E^{n-1}\tchib
+\rhob_N\ss_{N,N}E_N^{n-1}\tchib_N$$
up to terms which can be ignored. Hence the contribution of the 2nd term in \ref{12.217} is bounded 
in terms of:
\begin{eqnarray}
&&Cu^{-1}\int_0^{\ub_1}\ub\|N^\mu E^n\Lb\beta_\mu-N_N^\mu E_N^n\Lb_N\beta_{\mu,N}\|_{L^2(S_{\ub,u})}d\ub
\nonumber\\
&&+Cu\int_0^{u_1}\ub\|\s^{(n-1)}\ctchib\|_{L^2(S_{\ub,u})}d\ub
\label{12.218}
\end{eqnarray}
Here the 2nd term can be absorbed (compare with the 4th of \ref{12.199}) while the 1st is similar 
to the contribution of the $N$ component difference but with an extra $u^{-1}$ factor, hence is 
bounded in terms of \ref{12.206} but with an extra $u^{-1}$ factor. Then, up to terms which can be 
absorbed, \ref{12.218} is bounded by:
$$C\int_0^{\ub_1}\sqrt{\ub}\|\sqrt{a}\s^{(E;0,n)}\csxi\|_{L^2(S_{\ub,u})}d\ub$$
the contribution of which to $\|\s^{(0,n)}\cla\|_{L^2(\Cb_{\ub_1}^{u_1})}$ is bounded by:
\begin{equation}
\frac{C\sqrt{\s^{(E;0,n)}\cBb(\ub_1,u_1)}}{\left(a_0+\frac{3}{2}\right)}\cdot 
\ub_1^{a_0+\frac{3}{2}}u_1^{b_0}
\label{12.219}
\end{equation}
(compare with \ref{12.209}) and to $\|\s^{(0,n)}\cla\|_{L^2({\cal K}^{\tau})}$ by: 
\begin{equation}
\frac{C\sqrt{\s^{(E;0,n)}\cBb(\tau,\tau)}}{\left(a_0+\frac{3}{2}\right)}\cdot 
\tau^{a_0+b_0+\frac{3}{2}}
\label{12.220}
\end{equation}
(compare with \ref{12.212}). What remains to be considered is then the first of \ref{12.217}. 
Setting 
\begin{equation}
f(\ub,u)=\ub u^{-1}\|\s^{(Y;0,n)}\cxi_{\Lb}\|_{L^2(S_{\ub,u})}
\label{12.221}
\end{equation}
we have: 
\begin{eqnarray}
&&\|f(\ub,\cdot)\|^2_{L^2(\ub_1,u_1)}=\int_{\ub_1}^{u_1}f^2(\ub,u)du\leq\int_{\ub}^{u_1}f^2(\ub,u)du 
\nonumber\\
&&=\ub^2\int_{\ub}^{u_1}u^ {-2}\|\s^{(Y;0,n)}\cxi_{\Lb}\|^2_{L^2(S_{\ub,u})}du\nonumber\\
&&=\ub^2\int_{\ub}^{u_1}u^{-2}\frac{\partial}{\partial u}
\left(\int_{\ub}^u\|\s^{(Y;0,n)}\cxi_{\Lb}\|^2_{L^2(S_{\ub,u^\prime})}du^\prime\right)du\nonumber\\
&&=\ub^2\left\{u_1^{-2}\int_{\ub}^{u_1}\|\s^{(Y;0,n)}\cxi_{\Lb}\|^2_{L^2(S_{\ub,u})}du\right.\nonumber\\
&&\hspace{20mm}\left.+2\int_{\ub}^{u_1}u^{-3}\left(\int_{\ub}^u\|\s^{(Y;0,n)}\cxi_{\Lb}
\|^2_{L^2(S_{\ub,u^\prime})}du^\prime\right)du\right\}\nonumber\\
&&\leq C\ub^2\left\{u_1^{-2}\s^{(Y;0,n)}\cEb^{u_1}(\ub)+2\int_{\ub}^{u_1}u^{-3}
\s^{(Y;0,n)}\cEb^u(\ub)du\right\}
\nonumber\\
&&\leq C\ub^{2a_0+2}\left\{u_1^{2b_0-2}+2\frac{(u_1^{2b_0-2}-\ub^{2b_0-2})}{(2b_0-2)}\right\}
\s^{(Y;0,n)}\cBb(\ub_1,u_1)\nonumber\\
&&\leq C\ub^{2a_0+2}u_1^{2b_0-2}\s^{(Y;0,n)}\cBb(\ub_1,u_1)
\label{12.222}
\end{eqnarray}
by \ref{9.291}, \ref{9.298} (and \ref{10.391}). Then by Lemma 12.1a the contribution to 
$\|\s^{(0,n)}\cla\|_{L^2(\Cb_{\ub_1}^{u_1})}$ of the first of \ref{12.217} is bounded by:
\begin{equation}
C\sqrt{\s^{(Y;0,n)}\cBb(\ub_1,u_1)}u_1^{b_0-1}\int_0^{\ub_1}\ub^{a_0+1}d\ub
\leq\frac{C\sqrt{\s^{(Y;0,n)}\cBb(\ub_1,u_1)}}{(a_0+2)}\cdot\ub_1^{a_0}u_1^{b_0}\cdot\frac{\ub_1^2}{u_1}
\label{12.223}
\end{equation}
On the other hand, in view of \ref{12.210}, since by \ref{12.222} with $\tau$ in the role of $u_1$ we have:
\begin{equation}
\|f(\ub,\cdot)\|^2_{L^2(\ub,\tau)}=\int_{\ub}^\tau f^2(\ub,u)du
\leq C\ub^{2a_0+2}\tau^{2b_0-2}\s^{(Y;0,n)}\cBb(\tau,\tau)
\label{12.224}
\end{equation}
by Lemma 12.1b the contribution to $\|\s^{(0,n)}\cla\|_{L^2({\cal K}^\tau)}$ of the first of 
\ref{12.217} is bounded by:
\begin{equation}
C\sqrt{\s^{(Y;0,n)}\cBb(\tau,\tau)}\tau^{b_0-1}\int_0^{\tau}\ub^{a_0+1}d\ub
\leq\frac{C\sqrt{\s^{(Y;0,n)}\cBb(\tau,\tau)}}{(a_0+2)}\cdot\tau^{a_0+b_0+1}
\label{12.225}
\end{equation}

Consider finally the first of the differences \ref{12.202}, the $E$ component difference. Using 
\ref{12.213} we find that, up to terms which can be ignored, this is expressed as:
\begin{equation}
-\frac{1}{2}\beta_N^2\sbeta\eta^2 H^\prime\rho\left(\s^{(E;0,n)}\cxi_{\Lb}
+\rhob\sss\s^{(n-1)}\ctchib\right)
\label{12.226}
\end{equation}
The contribution of this to the integral \ref{12.196} is bounded by:
\begin{equation}
C\int_0^{\ub_1}\ub\|\s^{(E;0,n)}\cxi_{\Lb}\|_{L^2(S_{\ub,u})}d\ub
+Cu^2\int_0^{\ub_1}\ub\|\s^{(n-1)}\ctchib\|_{L^2(S_{\ub,u})}d\ub
\label{12.227}
\end{equation}
Here the contribution of the 2nd term can be absorbed, while the contribution of the 1st term 
to $\|\s^{(0,n)}\cla\|_{L^2(\Cb_{\ub_1}^{u_1})}$ is shown to be bounded by: 
\begin{equation}
\frac{C\sqrt{\s^{(E;0,n)}\cBb(\ub_1,u_1)}}{(a_0+2)}\cdot\ub_1^{a_0}u_1^{b_0}\cdot\ub_1^2
\label{12.228}
\end{equation}
and to $\|\s^{(0,n)}\cla\|_{L^2({\cal K}^\tau)}$ by:
\begin{equation}
\frac{C\sqrt{\s^{(E;0,n)}\cBb(\tau,\tau)}}{(a_0+2)}\cdot\tau^{a_0+b_0+2}
\label{12.229}
\end{equation}

We turn to the second of the principal difference terms \ref{12.200}. In fact the actual principal 
difference term in question is:
\begin{eqnarray}
&&-2\rhob\beta_N\left(c\pi\sbeta-\frac{\beta_N}{4}-\frac{\beta_{\Nb}}{2}\right)
H^\prime(g^{-1})^{\mu\nu}\beta_\nu E^n L\beta_\mu
\nonumber\\
&&\hspace{7mm}+2\rhob_N\beta_{N,N}\left(c_N\pi_N\sbeta_N-\frac{\beta_{N,N}}{4}-\frac{\beta_{\Nb,N}}{2}\right)H_N^\prime(g^{-1})^{\mu\nu}\beta_{\nu,N}E_N^n L_N\beta_{\mu,N}\nonumber\\
&&\label{12.230}
\end{eqnarray}
This is similar to \ref{12.163} but with extra $\rhob$, $\rhob_N$ factors.
It is expressed as the sum of differences of $E$, $N$, and $\Nb$ components:
\begin{eqnarray}
&&-2\rhob\sbeta\eta^2\beta_N\left(c\pi\sbeta-\frac{\beta_N}{4}-\frac{\beta_{\Nb}}{2}\right)
H^\prime E^\mu E^n L\beta_\mu\nonumber\\
&&\hspace{7mm}+2\rhob_N\sbeta_N\eta_N^2\beta_{N,N}\left(c_N\pi_N\sbeta_N-\frac{\beta_{N,N}}{4}-\frac{\beta_{\Nb,N}}{2}\right)H_N^\prime E_N^\mu E_N^n L_N\beta_{\mu,N}\nonumber\\
&&+\rhob\frac{\beta_{\Nb}}{c}\eta^2\beta_N
\left(c\pi\sbeta-\frac{\beta_N}{4}-\frac{\beta_{\Nb}}{2}\right)
H^\prime N^\mu E^n L\beta_\mu\nonumber\\
&&\hspace{7mm}-\rhob_N\frac{\beta_{\Nb,N}}{c_N}\eta_N^2\beta_{N,N}\left(c_N\pi_N\sbeta_N-\frac{\beta_{N,N}}{4}-\frac{\beta_{\Nb,N}}{2}\right)H_N^\prime N_N^\mu E_N^n L_N\beta_{\mu,N}\nonumber\\
&&+\rhob\frac{\beta_N}{c}\eta^2\beta_N
\left(c\pi\sbeta-\frac{\beta_N}{4}-\frac{\beta_{\Nb}}{2}\right)
H^\prime \Nb^\mu E^n L\beta_\mu\nonumber\\
&&\hspace{7mm}-\rhob_N\frac{\beta_{N,N}}{c_N}\eta_N^2\beta_{N,N}\left(c_N\pi_N\sbeta_N-\frac{\beta_{N,N}}{4}-\frac{\beta_{\Nb,N}}{2}\right)
H_N^\prime\Nb_N^\mu E_N^n L_N\beta_{\mu,N}\nonumber\\
&&\label{12.231}
\end{eqnarray}
which is similar to the differences \ref{12.164} but with extra $\rhob$, $\rhob_N$ factors. 
The contribution of the $\Nb$ component difference to the integral on the right in \ref{12.196} 
is then bounded by $Cu^2$ times \ref{12.171}, that is by:
\begin{eqnarray}
&&Cu\int_0^{\ub_1}\sqrt{\ub}\|\sqrt{a}\s^{(E;0,n)}\csxi\|_{L^2(S_{\ub,u})}d\ub
+Cu^2\int_0^{\ub_1}\|\s^{(0,n)}\lambdab\|_{L^2(S_{\ub,u})}d\ub\nonumber\\
&&+Cu^2\int_0^{\ub_1}\ub\|\s^{(n-1)}\ctchib\|_{L^2(S_{\ub,u})}d\ub
+Cu^2\int_0^{\ub_1}\ub\|\s^{(n-1)}\ctchi\|_{L^2(S_{\ub,u})}d\ub\nonumber\\
&&\label{12.232}
\end{eqnarray}
Here the last three can be absorbed while the contribution of the first to 
$\|\s^{(0,n)}\cla\|_{L^2(\Cb_{\ub_1}^{u_1})}$ can be shown to be bounded by 
(see \ref{12.154} - \ref{12.157}): 
\begin{equation}
\frac{C\sqrt{\s^{(E;0,n)}\cBb(\ub_1,u_1)}}{\left(a_0+\frac{3}{2}\right)}\cdot \ub_1^{a_0}u_1^{b_0}
\cdot \ub_1^{3/2} u_1
\label{12.233}
\end{equation}
and to $\|\s^{(0,n)}\cla\|_{L^2({\cal K}^{\tau})}$ by:
\begin{equation}
\frac{C\sqrt{\s^{(E;0,n)}\cBb(\tau,\tau)}}{\left(a_0+\frac{3}{2}\right)}\cdot \tau^{a_0+b_0+\frac{5}{2}}
\label{12.234}
\end{equation}
The contribution of the $E$ component difference to the integral on the right in \ref{12.196} is  
bounded by $Cu^2$ times \ref{12.175}, that is by:
\begin{equation}
Cu^2\int_0^{\ub_1}\|\s^{(E;0,n)}\cxi_L\|_{L^2(S_{\ub,u})}d\ub
+Cu^2\int_0^{\ub_1}\ub\|\s^{(n-1)}\ctchi\|_{L^2(S_{\ub,u})}d\ub
\label{12.235}
\end{equation}
Here the 2nd term can be absorbed while the 1st term is bounded according to \ref{12.161} by:
\begin{equation}
C\ub_1^{a_0+\frac{1}{2}}u^{b_0+2}\sqrt{\s^{(E;0,n)}\cB(\ub_1,u_1)}
\label{12.236}
\end{equation}
It follows that the corresponding contribution to $\|\s^{(0,n)}\cla\|_{L^2(\Cb_{\ub_1}^{u_1})}$ is 
bounded by:
\begin{equation}
\frac{C\sqrt{\s^{(E;0,n)}\cB(\ub_1,u_1)}}{\sqrt{2b_0+5}}\cdot\ub_1^{a_0}u_1^{b_0}\cdot\ub_1^{1/2}u_1^{5/2}
\label{12.237}
\end{equation}
and to $\|\s^{(0,n)}\cla\|_{L^2({\cal K}^{\tau})}$ by:
\begin{equation}
\frac{C\sqrt{\s^{(E;0,n)}\cB(\tau,\tau)}}{\sqrt{2b_0+5}}\cdot\tau^{a_0+b_0+3}
\label{12.238}
\end{equation}
Finally the contribution of the $N$ component difference to the integral on the right in \ref{12.196} 
is, up to terms which can be absorbed, bounded by:
\begin{equation}
Cu^2\int_0^{\ub_1}\|\s^{(Y;0,n)}\cxi_L\|_{L^2(S_{\ub,u})}d\ub\leq C\ub_1^{a_0+\frac{1}{2}}u_1^{b_0+2}
\sqrt{\s^{(Y;0,n)}\cB(\ub_1,u_1)}
\label{12.239}
\end{equation}
by \ref{12.181}. It follows that the corresponding contribution to 
$\|\s^{(0,n)}\cla\|_{L^2(\Cb_{\ub_1}^{u_1})}$ is 
bounded by:
\begin{equation}
\frac{C\sqrt{\s^{(Y;0,n)}\cB(\ub_1,u_1)}}{\sqrt{2b_0+5}}\cdot\ub_1^{a_0}u_1^{b_0}\cdot\ub_1^{1/2}u_1^{5/2}
\label{12.240}
\end{equation}
and to $\|\s^{(0,n)}\cla\|_{L^2({\cal K}^{\tau})}$ by:
\begin{equation}
\frac{C\sqrt{\s^{(Y;0,n)}\cB(\tau,\tau)}}{\sqrt{2b_0+5}}\cdot\tau^{a_0+b_0+3}
\label{12.241}
\end{equation}

The third of the principal difference terms \ref{12.200} is similar to the $N$ component difference 
in \ref{12.231}. 

We finally consider the fourth of the principal difference terms \ref{12.200}. In fact the actual 
principal difference term in question is:
\begin{equation}
\pi a\beta_N^2 H^\prime (g^{-1})^{\mu\nu}\beta_\nu E^{n+1}\beta_\mu
-\pi_N a_N\beta_{N,N}^2 H_N^\prime (g^{-1})^{\mu\nu}\beta_{\nu,N}E_N^{n+1}\beta_{\mu,N}
\label{12.242}
\end{equation}
This is similar to \ref{12.133} but, since $a=c\rho\rhob$, $a_N=c_N\rho_N\rhob_N$, 
with extra $\rhob$, $\rhob_N$ factors.
It is expressed as the sum of differences of $E$, $N$, and $\Nb$ components: 
\begin{eqnarray}
&&a\sbeta\eta^2\pi\beta_N^2 H^\prime E^\mu E^{n+1}\beta_\mu 
-a_N\sbeta_N\eta_N^2\pi_N\beta_{N,N}^2 H_N^\prime E_N^\mu E_N^{n+1}\beta_{\mu,N}\nonumber\\
&&-a\frac{\beta_{\Nb}}{2c}\eta^2\pi\beta_N^2 H^\prime N^\mu E^{n+1}\beta_\mu 
+a_N\frac{\beta_{\Nb,N}}{2c_N}\eta_N^2\pi_N\beta_{N,N}^2 H^\prime_N N_N^\mu E_N^{n+1}\beta_{\mu,N}
\nonumber\\
&&-a\frac{\beta_N}{2c}\eta^2\pi\beta_N^2 H^\prime \Nb^\mu E^{n+1}\beta_\mu 
+a_N\frac{\beta_{N,N}}{2c_N}\eta_N^2\pi_N\beta_{N,N}^2 H^\prime_N \Nb_N^\mu E_N^{n+1}\beta_{\mu,N}
\nonumber\\
&&\label{12.243}
\end{eqnarray}
which is similar to the differences \ref{12.134} but with extra $\rhob$, $\rhob_N$ factors.
The contribution of the $\Nb$ component difference to the integral on the right in \ref{12.196} 
is then bounded by $Cu^2$ times \ref{12.142}, that is by:
\begin{eqnarray}
&&C\int_0^{\ub_1}\ub\|\s^{(E;0,n)}\cxi_{\Lb}\|_{L^2(S_{\ub,u})}d\ub \nonumber\\
&&+Cu^2\int_0^{\ub_1}\ub\|\s^{(n-1)}\ctchib\|_{L^2(S_{\ub,u})}d\ub
+C\int_0^{\ub_1}\ub\|\s^{(n-1)}\ctchi\|_{L^2(S_{\ub,u})}d\ub \nonumber\\
&&+C\int_0^{\ub_1}\ub\|\s^{(0,n)}\cla\|_{L^2(S_{\ub,u})}d\ub \label{12.244}
\end{eqnarray}
Here the last three can be absorbed while the contribution of the first to 
$\|\s^{(0,n)}\cla\|_{L^2(\Cb_{\ub_1}^{u_1})}$ can be shown to be bounded by (see \ref{12.143} - 
\ref{12.145}): 
\begin{equation}
\frac{C\sqrt{\s^{(E;0,n)}\cBb(\ub_1,u_1)}}{(a_0+2)}\cdot \ub_1^{a_0}u_1^{b_0}\cdot\ub_1^2
\label{12.245}
\end{equation}
and to $\|\s^{(0,n)}\cla\|_{L^2({\cal K}^{\tau})}$ by:
\begin{equation}
\frac{C\sqrt{\s^{(E;0,n)}\cBb(\tau,\tau)}}{(a_0+2)}\cdot\tau^{a_0+b_0+2}
\label{12.246}
\end{equation}
The contribution of the $E$ component difference to the integral on the right in \ref{12.196} is  
bounded by $Cu^2$ times \ref{12.154}, that is by:
\begin{equation}
Cu\int_0^{\ub_1}\sqrt{\ub}\|\sqrt{a}\s^{(E;0,n)}\csxi\|_{L^2(S_{\ub,u})}d\ub
\label{12.247}
\end{equation}
the contribution of which to $\|\s^{(0,n)}\cla\|_{L^2(\Cb_{\ub_1}^{u_1})}$ is bounded by 
(see \ref{12.155} - \ref{12.157}):
\begin{equation}
\frac{C\sqrt{\s^{(E;0,n)}\cBb(\ub_1,u_1)}}{\left(a_0+\frac{3}{2}\right)}\cdot \ub_1^{a_0}u_1^{b_0}
\cdot\ub_1^{3/2}u_1
\label{12.248}
\end{equation}
and to $\|\s^{(0,n)}\cla\|_{L^2({\cal K}^{\tau})}$ by:
\begin{equation}
\frac{C\sqrt{\s^{(E;0,n)}\cBb(\tau,\tau)}}{\left(a_0+\frac{3}{2}\right)}\cdot \tau^{a_0+b_0+\frac{5}{2}}
\label{12.249}
\end{equation}
Finally, the contribution of the $N$ component difference to the integral on the right in \ref{12.196} 
is bounded by $Cu^2$ times \ref{12.160}, that is by:
\begin{eqnarray}
&&Cu^2\int_0^{\ub_1}\|\s^{(E;0,n)}\cxi_L\|_{L^2(S_{\ub,u})}d\ub
+Cu^2\int_0^{\ub_1}\ub\|\s^{(n-1)}\ctchi\|_{L^2(S_{\ub,u})}d\ub\nonumber\\
&&+Cu^2\int_0^{\ub_1}\ub\|\s^{(n-1)}\tchib\|_{L^2(S_{\ub,u})}d\ub
+Cu^2\int_0^{\ub}\|\s^{(0,n)}\clab\|_{L^2(S_{\ub,u})}d\ub\nonumber\\
&&\label{12.250}
\end{eqnarray}
Here the last three can be absorbed while the contribution of the first to 
$\|\s^{(0,n)}\cla\|_{L^2(\Cb_{\ub_1}^{u_1})}$ is bounded by (see \ref{12.161}, \ref{12.162}):
\begin{equation}
\frac{C\sqrt{\s^{(E;0,n)}\cB(\ub_1,u_1)}}{\sqrt{2b_0+5}}\cdot\ub_1^{a_0}u_1^{b_0}\cdot
\ub_1^{1/2}u_1^{5/2}
\label{12.251}
\end{equation}
and to $\|\s^{(0,n)}\cla\|_{L^2({\cal K}^{\tau})}$ by:
\begin{equation}
\frac{C\sqrt{\s^{(E;0,n)}\cB(\tau,\tau)}}{\sqrt{2b_0+5}}\cdot \tau^{a_0+b_0+3}
\label{12.252}
\end{equation}

Returning to \ref{12.199}, in regard to the 1st and 3rd terms we set 
\begin{equation}
f(\ub,u)=\ub\|\s^{(0,n)}\cla\|_{L^2(S_{\ub,u})} \ \ \mbox{and} \ \ 
f(\ub,u)=\ub\|\s^{(n-1)}\ctchi\|_{L^2(S_{\ub,u})}
\label{12.253}
\end{equation}
respectively. Then since 
\begin{equation}
\|f(\ub,\dot)\|_{L^2(\ub,u_1)}\leq\ub\|\s^{(0,n)}\cla\|_{L^2(\Cb_{\ub}^{u_1})} \ \ \mbox{and} \ \ 
\|f(\ub,\cdot)\|_{L^2(\ub,u_1)}\leq\ub\|\s^{(n-1)}\ctchi\|_{L^2(\Cb_{\ub}^{u_1})}
\label{12.254}
\end{equation}
respectively, applying Lemma 12.1a we obtain that the contributions of these terms to 
$\|\s^{(0,n)}\cla\|_{L^2(\Cb_{\ub_1}^{u_1})}$ are bounded by 
\begin{equation}
C\int_0^{\ub_1}\ub\|\s^{(0,n)}\cla\|_{L^2(\Cb_{\ub}^{u_1})}d\ub \ \ \mbox{and} \ \ 
C\int_0^{\ub_1}\ub\|\s^{(n-1)}\ctchi\|_{L^2(\Cb_{\ub}^{u_1})}d\ub 
\label{12.255}
\end{equation}
respectively. Also, setting $u_1=\tau$ in \ref{12.254} and applying Lemma 12.1b we obtain that 
the corresponding contributions to $\|\s^{(0,n)}\cla\|_{L^2({\cal K}^{\tau})}$ are bounded by 
\begin{equation}
C\int_0^{\tau} \ub\|\s^{(0,n)}\cla\|_{L^2(\Cb_{\ub}^{\tau})}d\ub \ \ \mbox{and} \ \ 
C\int_0^{\tau} \ub\|\s^{(n-1)}\ctchi\|_{L^2(\Cb_{\ub}^{\tau})}d\ub 
\label{12.256}
\end{equation}
respectively. 

We combine the above results in the form in which they will be used in the sequel in the following 
lemma. 

\vspace{2.5mm}

\noindent{\bf Lemma 12.3} \ \ The next to the top order acoustical quantity 
$\|\s^{(0,n)}\cla\|_{L^2(\Cb_{\ub}^u)}$ satisfies to principal terms the following inequality:
\begin{eqnarray*}
&&\|\s^{(0,n)}\cla\|_{L^2(\Cb_{\ub_1}^{u_1})}\leq 
C\left\{\int_{\ub_1}^{u_1}\left(\int_0^{\ub_1}
\|\s^{(0,n)}\clab\|_{L^2(S_{\ub,u})}d\ub\right)^2 du\right\}^{1/2} \\
&&\hspace{27mm}+C\left\{\int_{\ub_1}^{u_1}\left(\int_0^{\ub_1}\ub
\|\s^{(n-1)}\ctchib\|_{L^2(S_{\ub,u})}d\ub\right)^2 du\right\}^{1/2} \\
&&\hspace{25mm}+C\int_0^{\ub_1}\ub\|\s^{(0,n)}\cla\|_{L^2(\Cb_{\ub}^{u_1})}d\ub 
+C\int_0^{\ub_1}\ub\|\s^{(n-1)}\ctchi\|_{L^2(\Cb_{\ub}^{u_1})}d\ub \\
&&\hspace{27mm}+C\frac{\sqrt{\s^{(Y;0,n)}\cBb(\ub_1,u_1)}}{(a_0+2)}\cdot \ub_1^{a_0}u_1^{b_0}
\cdot\frac{\ub_1^2}{u_1}\\
&&\hspace{25mm}+C\sqrt{\max\{\s^{(0,n)}\cB(\ub_1,u_1),\s^{(0,n)}\cBb(\ub_1,u_1)\}}\cdot
\ub_1^{a_0}u_1^{b_0}\cdot \ub_1^{1/2}u_1
\end{eqnarray*}
Moreover $\|\s^{(0,n)}\cla\|_{L^2({\cal K}^{\tau})}$ satisfies: 
\begin{eqnarray*}
&&\|\s^{(0,n)}\cla\|_{L^2({\cal K}^{\tau})}\leq 
C\left\{\int_0^{\tau}\left(\int_0^u 
\|\s^{(0,n)}\clab\|_{L^2(S_{\ub,u})}d\ub\right)^2 du\right\}^{1/2} \\
&&\hspace{27mm}+C\left\{\int_0^{\tau}\left(\int_0^u \ub
\|\s^{(n-1)}\ctchib\|_{L^2(S_{\ub,u})}d\ub\right)^2 du\right\}^{1/2} \\
&&\hspace{25mm}+C\int_0^{\tau}\ub\|\s^{(0,n)}\cla\|_{L^2(\Cb_{\ub}^{\tau})}d\ub 
+C\int_0^{\tau}\ub\|\s^{(n-1)}\ctchi\|_{L^2(\Cb_{\ub}^{\tau})}d\ub \\
&&\hspace{27mm}+C\frac{\sqrt{\s^{(Y;0,n)}\cBb(\tau,\tau)}}{(a_0+2)}\cdot \tau^{a_0+b_0}\cdot\tau\\
&&\hspace{25mm}+C\sqrt{\max\{\s^{(0,n)}\cB(\tau,\tau),\s^{(0,n)}\cBb(\tau,\tau)\}}\cdot
\tau^{a_0+b_0}\cdot \tau^{3/2}
\end{eqnarray*}

\vspace{2.5mm}

We turn to the derivation of estimates for $\s^{(n-1)}\ctchib$ and for $\s^{(n-1)}\clab$ in ${\cal N}$ 
in terms of their boundary values on ${\cal K}$. We begin with the derivation of the estimate for 
$\s^{(n-1)}\ctchib$. We consider the propagation equation for $\s^{(n-1)}\ctchib$, equation \ref{12.32}. 
Let us denote the right hand side by $\s^{(n-1)}\chFb$, so this equation reads:
\begin{equation}
\Lb\s^{(n-1)}\ctchib+(n+1)\chib\s^{(n-1)}\ctchib=\s^{(n-1)}\chFb 
\label{12.257}
\end{equation}
Here as in Section 10.7 we use the $(\ub,u,\vartheta^{\prime\prime})$ coordinates on ${\cal N}$ which 
are adapted to the flow of $\Lb$. In view of \ref{10.600}, equation \ref{12.257} takes in these 
coordinates the form:
\begin{equation}
\frac{\partial}{\partial u}\left(\sh^{\prime\prime(n+1)/2}\s^{(n-1)}\ctchib\right)=
\sh^{\prime\prime(n+1)/2}\s^{(n-1)}\chFb
\label{12.258}
\end{equation}
Integrating from ${\cal K}$ we obtain:
\begin{eqnarray}
&&\left(\sh^{\prime\prime(n+1)/2}\s^{(n-1)}\ctchib\right)(\ub,u_1,\vartheta^{\prime\prime})
=\left(\sh^{\prime\prime(n+1)/2}\s^{(n-1)}\ctchib\right)(\ub,\ub,\vartheta^{\prime\prime})\nonumber\\
&&\hspace{30mm}+\int_{\ub}^{u_1}\left(\sh^{\prime\prime(n+1)/2}\s^{(n-1)}\chFb\right)
(\ub,u,\vartheta^{\prime\prime})du 
\label{12.259}
\end{eqnarray}
In view of the remark \ref{10.607} we rewrite this in the form:
\begin{eqnarray}
&&\left(\s^{(n-1)}\ctchib\sh^{\prime\prime 1/4}\right)(\ub,u_1,\vartheta^{\prime\prime})=
\left(\frac{\sh(\ub,\ub,\vartheta^{\prime\prime})}
{\sh^{\prime\prime}(\ub,u_1,\vartheta^{\prime\prime})}\right)^{(2n+1)/4}
\s^{(n-1)}\ctchib(\ub,\ub,\vartheta^{\prime\prime})\nonumber\\
&&\hspace{20mm}+\int_{\ub}^{u_1}\left(\frac{\sh^{\prime\prime}(\ub,u,\vartheta^{\prime\prime})}
{\sh^{\prime\prime}(\ub,u_1,\vartheta^{\prime\prime})}\right)^{(2n+1)/4}
\left(\s^{(n-1)}\chFb\sh^{\prime\prime 1/4}\right)(\ub,u,\vartheta^{\prime\prime})du\nonumber\\
&&\label{12.260}
\end{eqnarray}
In view of \ref{10.611} this implies:
\begin{eqnarray}
&&\left|\left(\s^{(n-1)}\ctchib\sh^{\prime\prime 1/4}\right)(\ub,u_1,\vartheta^{\prime\prime})\right|
\leq k\left|\left(\s^{(n-1)}\ctchib\sh^{1/4}\right)(\ub,\ub,\vartheta^{\prime\prime})\right|\nonumber\\
&&\hspace{30mm}+k\int_{\ub}^{u_1}\left|\left(\s^{(n-1)}\chFb\sh^{\prime\prime 1/4}\right)
(\ub,u,\vartheta^{\prime\prime})\right|du
\label{12.261}
\end{eqnarray}
where $k$ is a constant greater than 1, but which can be chosen as close to 1 as we wish by suitably 
restricting $u_1$. Taking $L^2$ norms with respect to $\vartheta^{\prime\prime}\in S^1$ we obtain:
\begin{eqnarray}
&&\left\|\left(\s^{(n-1)}\ctchib\sh^{\prime\prime 1/4}\right)(\ub,u_1)\right\|_{L^2(S^1)}
\leq k\left\|\left(\s^{(n-1)}\ctchib\sh^{1/4}\right)(\ub,\ub)\right\|_{L^2(S^1)}\nonumber\\
&&\hspace{25mm}+k\int_{\ub}^{u_1}\left\|\left(\s^{(n-1)}\chFb\sh^{\prime\prime 1/4}\right)(\ub,u)
\right\|_{L^2(S^1)}du
\label{12.262}
\end{eqnarray}
or, in view of \ref{10.607} (see also \ref{10.601}):
\begin{equation}
\|\s^{(n-1)}\ctchib\|_{L^2(S_{\ub,u_1})}\leq k\|\s^{(n-1)}\ctchib\|_{L^2(S_{\ub,\ub})}
+k\int_{\ub}^{u_1}\|\s^{(n-1)}\chFb\|_{L^2(S_{\ub,u})}du
\label{12.263}
\end{equation}
Taking then $L^2$ norms with respect to $\ub\in[0,\ub_1]$ we obtain:
\begin{eqnarray}
&&\|\s^{(n-1)}\ctchib\|_{L^2(C_{u_1}^{\ub_1})}\leq k\|\s^{(n-1)}\ctchib\|_{L^2({\cal K}^{\ub_1})} 
\label{12.264}\\
&&\hspace{27mm}+k\left\{\int_0^{\ub_1}\left(\int_{\ub}^{u_1}\|\s^{(n-1)}\chFb\|_{L^2(S_{\ub,u})}du\right)^2 d\ub\right\}^{1/2} \nonumber
\end{eqnarray}

We shall presently estimate the leading contributions to 
\begin{equation}
\left\{\int_0^{\ub_1}\left(\int_{\ub}^{u_1}\|\s^{(n-1)}\chFb\|_{L^2(S_{\ub,u})}du\right)^2 d\ub
\right\}^{1/2}
\label{12.265}
\end{equation}
Consider first the $n$th order acoustical acoustical difference terms:
\begin{equation}
\oAb\s^{(n-1)}\ctchib+\Ab\s^{(n-1)}\ctchi
\label{12.266}
\end{equation}
The coefficients $\oAb$ and $\Ab$ being given by \ref{12.21} ad \ref{12.22} respectively, the assumptions 
\ref{10.440} imply: 
\begin{equation}
|\oAb|, |\Ab| \leq C \ \mbox{: in ${\cal R}_{\delta,\delta}$}
\label{12.267}
\end{equation}
It follows that the contribution of the terms \ref{12.266} to the integral on the right in \ref{12.263} 
is bounded by:
\begin{equation}
C\int_{\ub}^{u_1}\|\s^{(n-1)}\ctchib\|_{L^2(S_{\ub,u})}du
+C\int_{\ub}^{u_1}\|\s^{(n-1)}\ctchi\|_{L^2(S_{\ub,u})}du
\label{12.268}
\end{equation}
Here the first can obviously be absorbed in the integral inequality \ref{12.263} while the second is bounded by $C\sqrt{u_1}\|\s^{(n-1)}\ctchi\|_{L^2(\Cb_{\ub}^{u_1})}$ hence 
its contribution to \ref{12.265} is bounded by: 
\begin{equation}
Cu_1^{1/2}\left\{\int_0^{\ub_1}\|\s^{(n-1)}\ctchi\|^2_{L^2(\Cb_{\ub}^{u_1})}d\ub\right\}^{1/2}
\label{12.269}
\end{equation}
We shall see in the sequel that this can be absorbed. 

The other leading contributions are those of the principal difference terms:
\begin{equation}
-\beta_{\Nb}\sbeta E^n \Lb H+\beta_{\Nb,N}\sbeta_N E_N^n \Lb_N H_N 
+\frac{1}{2}\rhob\beta_{\Nb}^2 E^{n+1}H-\frac{1}{2}\rhob_N\beta_{\Nb,N}^2 E_N^{n+1} H_N
\label{12.270}
\end{equation}
Here again we must keep track not only of the actually principal terms, 
namely those containing derivatives of the $\beta_\mu$ of order $n+1$, 
but also of the terms containing acoustical quantities of order $n$. 
The terms of order $n$ with vanishing principal acoustical part can be ignored. 
As we shall see the dominant contribution comes from the 1st of the difference terms \ref{12.270}. 
The actual principal difference term in question is:
\begin{equation}
2\beta_{\Nb}\sbeta\eta^2 H^\prime(g^{-1})^{\mu\nu}\beta_\nu E^n \Lb\beta_\mu 
-2\beta_{\Nb,N}\sbeta_N\eta_N^2 H^\prime_N(g^{-1})^{\mu\nu}\beta_{\nu,N}E_N^n\Lb_N\beta_{\mu,N}
\label{12.271}
\end{equation}
By \ref{10.368} and the analogous formula for the $N$th approximants this is expressed as the sum of 
differences of $E$, $N$, and $\Nb$ components:
\begin{eqnarray}
&&2\beta_{\Nb}\sbeta^2\eta^2 H^\prime E^\mu E^n\Lb\beta_\mu 
-2\beta_{\Nb,N}\sbeta_N^2\eta_N^2 H^\prime_N E_N^\mu E_N^n\Lb_N\beta_{\mu,N}\nonumber\\
&&-\beta_{\Nb}\sbeta\frac{\beta_{\Nb}}{c}\eta^2 H^\prime N^\mu E^n\Lb\beta_\mu 
+\beta_{\Nb,N}\sbeta_N\frac{\beta_{\Nb,N}}{c_N}\eta_N^2 H^\prime_N N_N^\mu E_N^n\Lb_N\beta_{\mu,N}
\nonumber\\
&&-\beta_{\Nb}\sbeta\frac{\beta_N}{c}\eta^2 H^\prime \Nb^\mu E^n\Lb\beta_\mu 
+\beta_{\Nb,N}\sbeta_N\frac{\beta_{N,N}}{c_N}\eta_N^2 H^\prime_N \Nb_N^\mu E_N^n\Lb_N\beta_{\mu,N}
\nonumber\\&&\label{12.272}
\end{eqnarray}
The dominant contribution is that of the $\Nb$ component difference. To estimate it we must use 
\ref{12.203}, which suggests that we consider first the $N$ component difference. By \ref{12.204} 
and a similar expression for the $N$th approximants, the $N$ component difference is seen to be 
given by 
\begin{equation}
-\beta_{\Nb}\sbeta\frac{\beta_{\Nb,N}}{c}\eta^2 H^\prime\left\{\lambda\s^{(E;0,n)}\csxi
+\sss \s^{(0,n)}\cla-\rhob\left(\frac{\ss_{\Nb}}{2}-\pi c\sss\right)
\s^{(n-1)}\ctchi+\rhob\frac{\ss_{\Nb}}{2}\s^{(n-1)}\ctchib\right\}
\label{12.273}
\end{equation}
The contribution of this to the integral on the right in \ref{12.263} is bounded by:
\begin{eqnarray}
&&\frac{C}{\sqrt{\ub}}\int_{\ub}^{u_1}u\|\sqrt{a}\s^{(E;0,n)}\csxi\|_{L^2(S_{\ub,u})}du
+C\int_{\ub}^{u_1}\|\s^{(0,n)}\cla\|_{L^2(S_{\ub,u})}du\nonumber\\
&&+C\int_{\ub}^{u_1}u^2\|\s^{(n-1)}\ctchi\|_{L^2(S_{\ub,u})}du
+C\int_{\ub}^{u_1}u^2\|\s^{(n-1)}\ctchib\|_{L^2(S_{\ub,u})}du \nonumber\\
&&\label{12.274}
\end{eqnarray}
Here the last two terms are dominated by the corresponding terms in \ref{12.268}, while the second  
is bounded by $C\sqrt{u_1}\|\s^{(0,n)}\cla\|_{L^2(\Cb_{\ub}^{u_1})}$ hence its contribution to 
\ref{12.265} is bounded by:
\begin{equation}
u_1^{1/2}\left\{\int_0^{\ub_1}\|\s^{(0,n)}\cla\|^2_{L^2(\Cb_{\ub}^{u_1})}d\ub\right\}^{1/2}
\label{12.275}
\end{equation}
On the other hand the first of \ref{12.274} is bounded by:
\begin{eqnarray}
&&\frac{Cu_1^{3/2}}{\ub^{1/2}}\|\sqrt{a}\s^{(E;0,n)}\csxi\|_{L^2(\Cb_{\ub}^{u_1})}\nonumber\\
&&\leq\frac{Cu_1^{3/2}}{\ub^{1/2}}\sqrt{\s^{(E;0,n)}\cEb^{u_1}(\ub)}
\leq \frac{Cu_1^{3/2}}{\ub^{1/2}}\ub^{a_0}u_1^{b_0}\sqrt{\s^{(E;0,n)}\cBb(\ub_1,u_1)}\nonumber\\
&&\label{12.276}
\end{eqnarray}
hence its contribution to \ref{12.265} is bounded by:
\begin{equation}
\frac{C\sqrt{\s^{(E;0,n)}\cBb(\ub_1,u_1)}}{\sqrt{2a_0}}\cdot\ub_1^{a_0}u_1^{b_0}\cdot u_1^{3/2}
\label{12.277}
\end{equation}
We turn to the $\Nb$ component difference, last of \ref{12.272}. By \ref{12.213} - \ref{12.215} 
this is given by, up to terms which can be ignored, 
\begin{equation}
-\beta_{\Nb}\sbeta\frac{\beta_N}{c}\eta^2 H^\prime 
\left\{\Nb^\mu\Lb(E^n\beta_\mu-E_N^n\beta_{\mu,N})+\rhob\ss_{\Nb}\s^{(n-1)}\ctchib\right\}
\label{12.278}
\end{equation}
hence its contribution to \ref{12.265} is bounded by:
\begin{eqnarray}
&&C\int_{\ub}^{u_1}\|\Nb^\mu\Lb(E^n\beta_\mu-E_N^n\beta_{\mu,N})\|_{L^2(S_{\ub,u})}du \nonumber\\
&&+C\int_{\ub}^{u_1}u^2\|\s^{(n-1)}\ctchib\|_{L^2(S_{\ub,u})}du
\label{12.279}
\end{eqnarray}
The second is dominated by the corresponding term in \ref{12.268}. As for the first, by \ref{12.203}, 
in view of the fact that $\ogamma\sim u$, this term is bounded by:
\begin{eqnarray}
&&C\int_{\ub}^{u_1}u^{-1}\|\s^{(Y;0,n)}\cxi_{\Lb}\|_{L^2(S_{\ub,u})}du\nonumber\\
&&+C\int_{\ub}^{u_1}u^{-1}\|N^\mu\Lb(E^n\beta_\mu-E_N^n\beta_{\mu,N})\|_{L^2(S_{\ub,u})}du
\label{12.280}
\end{eqnarray}
The second of \ref{12.280} is bounded by
\begin{eqnarray}
&&C\int_{\ub}^{u_1}u^{-1}\|N^\mu E^n\Lb\beta_\mu-N_N^\mu E_N^n\Lb_N\beta_{\mu,N}\|_{L^2(S_{\ub,u})}du
\nonumber\\
&&+C\int_{\ub}^{u_1}u\|\s^{(n-1)}\ctchib\|_{L^2(S_{\ub,u})}d\ub
\label{12.281}
\end{eqnarray}
(compare with \ref{12.218}). Here the 2nd term is dominated by the corresponding term in \ref{12.268}, 
while (see \ref{12.273}) the 1st is bounded by: 
\begin{eqnarray}
&&\frac{C}{\sqrt{\ub}}\int_{\ub}^{u_1}\|\sqrt{a}\s^{(E;0,n)}\csxi\|_{L^2(S_{\ub,u})}du
+C\int_{\ub}^{u_1}u^{-1}\|\s^{(0,n)}\cla\|_{L^2(S_{\ub,u})}du\nonumber\\
&&+C\int_{\ub}^{u_1}u\|\s^{(n-1)}\ctchi\|_{L^2(S_{\ub,u})}du
+C\int_{\ub}^{u_1}u\|\s^{(n-1)}\ctchib\|_{L^2(S_{\ub,u})}du \nonumber\\
&&\label{12.282}
\end{eqnarray}
In the above the last two terms are dominated by the corresponding terms in \ref{12.268}. The 
2nd term is bounded by (Schwartz inequality)
$$C\left(\int_{\ub}^{u_1}u^{-2}du\right)^{1/2}\|\s^{(0,n)}\cla\|_{L^2(\Cb_{\ub}^{u_1})}
\leq C\ub^{-1/2}\|\s^{(0,n)}\cla\|_{L^2(\Cb_{\ub}^{u_1})},$$ 
hence its contribution to \ref{12.265} is bounded by
\begin{equation}
C\left\{\int_0^{\ub_1}\ub^{-1}\|\s^{(0,n)}\cla\|^2_{L^2(\Cb_{\ub}^{u_1})}d\ub\right\}^{1/2}
\label{12.283}
\end{equation}
The 1st term is bounded by 
\begin{eqnarray}
&&\frac{Cu_1^{1/2}}{\ub^{1/2}}\|\sqrt{a}\s^{(E;0,n)}\csxi\|_{L^2(\Cb_{\ub}^{u_1})}\nonumber\\
&&\leq\frac{Cu_1^{1/2}}{\ub^{1/2}}\sqrt{\s^{(E;0,n)}\cEb^{u_1}(\ub)}
\leq \frac{Cu_1^{1/2}}{\ub^{1/2}}\ub^{a_0}u_1^{b_0}\sqrt{\s^{(E;0,n)}\cBb(\ub_1,u_1)}\nonumber\\
&&\label{12.284}
\end{eqnarray}
hence its contribution to \ref{12.265} is bounded by 
\begin{equation}
\frac{C\sqrt{\s^{(E;0,n)}\cBb(\ub_1,u_1)}}{\sqrt{2a_0}}\cdot\ub_1^{a_0}u_1^{b_0}\cdot u_1^{1/2}
\label{12.285}
\end{equation}
What remains to be considered in regard to the $\Nb$ component difference is the first of \ref{12.280}. 
This is bounded by (Schwartz inequality) 
\begin{eqnarray*}
&&C\left(\int_{\ub}^{u_1}u^{-2}du\right)^{1/2}\|\s^{(Y;0,n)}\cxi_{\Lb}\|_{L^2(\Cb_{\ub}^{u_1})}\\
&&\leq C\ub^{-1/2}\sqrt{\s^{(Y;0,n)}\cEb^{u_1}(\ub)}\leq C\ub^{a_0-\frac{1}{2}}u_1^{b_0}
\sqrt{\s^{(Y;0,n)}\cBb(\ub_1,u_1)}
\end{eqnarray*}
hence its contribution to \ref{12.265} is bounded by
\begin{equation}
\frac{C\sqrt{\s^{(Y;0,n)}\cBb(\ub_1,u_1)}}{\sqrt{2a_0}}\cdot\ub_1^{a_0}u_1^{b_0}
\label{12.286}
\end{equation}
Finally, up to terms which can be ignored, the $E$ component difference, first of \ref{12.272}, 
is given by (see \ref{12.226}):
\begin{equation}
2\beta_{\Nb}\sbeta^2\eta^2 H^\prime\left\{\s^{(E;0,n)}\cxi_{\Lb}+\rhob\sss \s^{(n-1)}\ctchib\right\}
\label{12.287}
\end{equation}
The contribution of this to the integral on the right in \ref{12.263} is bounded by:
\begin{equation}
C\int_{\ub}^{u_1}\|\s^{(E;0,n)}\cxi_{\Lb}\|_{L^2(S_{\ub,u})}du
+C\int_{\ub}^{u_1}u^2\|\s^{(n-1)}\ctchib\|_{L^2(S_{\ub,u})}du
\label{12.288}
\end{equation}
Here the 2nd term is dominated by the corresponding term in \ref{12.268} while the 1st term is 
bounded by 
\begin{eqnarray*}
&&Cu_1^{1/2}\|\s^{(E;0,n)}\cxi_{\Lb}\|_{L^2(\Cb_{\ub}^{u_1})}\\
&&\leq Cu_1^{1/2}\sqrt{\s^{(E;0,n)}\cEb^{u_1}(\ub)}\leq C\ub^{a_0} u_1^{b_0+\frac{1}{2}}
\sqrt{\s^{(E;0,n)}\cBb(\ub_1,u_1)}
\end{eqnarray*}
hence its contribution to \ref{12.265} is bounded by
\begin{equation}
\frac{C\sqrt{\s^{(E;0,n)}\cBb(\ub_1,u_1)}}{\sqrt{2a_0+1}}\cdot\ub_1^{a_0}u_1^{b_0}
\cdot\ub_1^{1/2}u_1^{1/2}
\label{12.289}
\end{equation}

Consider next the 2nd of the principal difference terms \ref{12.270}. The actual principal difference 
term in question is: 
\begin{equation}
-\rhob\beta_{\Nb}^2\eta^2 H^\prime (g^{-1})^{\mu\nu}\beta_\nu E^{n+1}\beta_\mu 
+\rhob_N\beta_{\Nb,N}^2\eta_N^2 H^\prime_N (g^{-1})^{\mu\nu}\beta_{\nu,N} E_N^{n+1}\beta_{\mu,N}
\label{12.290}
\end{equation}
By the formula \ref{10.368} and the analogous formula for the $N$th approximants this is expressed as a 
sum of differences of $E$, $N$ and $\Nb$ componets:
\begin{eqnarray}
&&-\beta_{\Nb}^2\sbeta\eta^2 H^\prime\rhob E^\mu E^{n+1}\beta_\mu 
+\beta_{\Nb,N}^2\sbeta_N\eta_N^2 H^\prime_N\rhob_N E_N^\mu E_N^{n+1}\beta_{\mu,N} \nonumber\\
&&+\beta_{\Nb}^2\frac{\beta_{\Nb}}{2c}\eta^2 H^\prime\rhob N^\mu E^{n+1}\beta_\mu 
-\beta_{\Nb,N}^2\frac{\beta_{\Nb,N}}{2c_N}\eta_N^2 H^\prime_N\rhob_N N_N^\mu E_N^{n+1}\beta_{\mu,N}
\nonumber\\
&&+\beta_{\Nb}^2\frac{\beta_N}{2c}\eta^2 H^\prime\rhob\Nb^\mu E^{n+1}\beta_\mu 
-\beta_{\Nb,N}^2\frac{\beta_{N,N}}{2c_N}\eta_N^2 H^\prime_N\rhob_N\Nb_N^\mu E_N^{n+1}\beta_{\mu,N}
\nonumber\\&&\label{12.291}
\end{eqnarray}
We first consider the $\Nb$ component difference. In view of \ref{12.140} we can express this as:
\begin{equation}
\frac{\beta_{\Nb}^2\beta_N}{2c}\eta^2 H^\prime\left\{\s^{(E;0,n)}\cxi_{\Lb}
+\frac{1}{2}\rhob\sss \s^{(n-1)}\ctchib+\left(\frac{s_{\Nb\Lb}}{2c}-\pi\rhob\ss_{\Nb}\right)
\s^{(n-1)}\ctchi-\frac{\ss_{\Nb}}{c}\s^{(0,n)}\cla\right\}
\label{12.292}
\end{equation}
The contribution of this to the integral on the right in \ref{12.263} is bounded by:
\begin{eqnarray}
&&C\int_{\ub}^{u_1}\|\s^{(E;0,n)}\cxi_{\Lb}\|_{L^2(S_{\ub,u})}du\nonumber\\
&&+C\int_{\ub}^{u_1}u^2\|\s^{(n-1)}\ctchib\|_{L^2(S_{\ub,u})}du 
+C\int_{\ub}^{u_1}\|\s^{(n-1)}\ctchi\|_{L^2(S_{\ub,u})}du \nonumber\\
&&+C\int_{\ub}^{u_1}\|\s^{(0,n)}\cla\|_{L^2(S_{\ub,u})}du \label{12.293}
\end{eqnarray}
Here the 2nd is dominated by while the 2nd coincides with the corresponding terms in \ref{12.268}. 
The 4th coincides with the 2nd term in \ref{12.274}. The 1st term coincides with the 1st term in 
\ref{12.288}. 

We turn to the $N$ component difference, second of \ref{12.291}. This is given by
\begin{equation}
\frac{\beta_{\Nb}^3}{2c}\eta^2 H^\prime\rhob N^\mu E(E^n\beta_\mu-E_N^n\beta_{\mu,N})
\label{12.294}
\end{equation}
up to terms which can be ignored. We have:
\begin{equation}
N^\mu E(E^n\beta_\mu-E_N^n\beta_{\mu,N})=\s^{(Y;0,n)}\csxi-\ogamma\Nb^\mu E(E^n\beta_\mu-E_N^n\beta_{\mu,N})
\label{12.295}
\end{equation}
Since $\ogamma\sim u$, the contribution of the 2nd term on the right in \ref{12.295}, through 
\ref{12.294}, to the integral on the right in \ref{12.263} is bounded by $Cu_1$ times the 
corresponding contribution of the $\Nb$ component difference estimated above. The contribution 
of the 1st term on the right is bounded by:
\begin{eqnarray}
&&\frac{C}{\sqrt{\ub}}\int_{\ub}^{u_1}u\|\sqrt{a}\s^{(Y;0,n)}\csxi\|_{L^2(S_{\ub,u})}du 
\leq \frac{C}{\sqrt{\ub}}\left(\int_{\ub}^{u_1}u^2 du\right)^{1/2}
\|\sqrt{a}\s^{(Y;0,n)}\csxi\|_{L^2(\Cb_{\ub}^{u_1})}\nonumber\\
&&\hspace{40mm}\leq\frac{Cu_1^{3/2}}{\ub^{1/2}}\sqrt{\s^{(Y;0,n)}\cEb^{u_1}(\ub)}
\label{12.296}
\end{eqnarray}
hence its contribution to \ref{12.265} is bounded by:
\begin{equation}
Cu_1^{3/2}\left(\int_0^{\ub_1}\frac{\s^{(Y;0,n)}\cEb^{u_1}(\ub)}{\ub}d\ub\right)^{1/2}
\leq \frac{C\sqrt{\s^{(Y;0,n)}\cBb(\ub_1,u_1)}}{\sqrt{2a_0}}\cdot\ub_1^{a_0}u_1^{b_0}\cdot u_1^{3/2}
\label{12.297}
\end{equation}
Finally the $E$ component difference, first of \ref{12.291}, is given by
\begin{equation}
-\beta_N^2\sbeta\eta^2 H^\prime\rhob\s^{(E;0,n)}\csxi
\label{12.298}
\end{equation}
up to terms which can be ignored. Therefore its contribution to \ref{12.265} is similarly bounded by:
\begin{equation}
\frac{C\sqrt{\s^{(E;0,n)}\cBb(\ub_1,u_1)}}{\sqrt{2a_0}}\cdot\ub_1^{a_0}u_1^{b_0}\cdot u_1^{3/2}
\label{12.299}
\end{equation}

Absorbing the first of \ref{12.268} in the integral inequality \ref{12.263}, we summarize the 
above results in the form in which they will be used in the sequel in the following lemma. 

\vspace{2.5mm}

\noindent{\bf Lemma 12.4} \ \ The next to the top order acoustical quantity 
$\|\s^{(n-1)}\ctchib\|_{L^2(C_u^{\ub})}$ satisfies to principal terms the following inequality:
\begin{eqnarray*}
&&\|\s^{(n-1)}\ctchib\|_{L^2(C_{u_1}^{\ub_1})}\leq k\|\s^{(n-1)}\ctchib\|_{L^2({\cal K}^{\ub_1})}\\
&&\hspace{27mm}+Cu_1^{1/2}\left\{\int_0^{\ub_1}\|\s^{(n-1)}\ctchi\|^2_{L^2(\Cb_{\ub}^{u_1})}d\ub
\right\}^{1/2}\\
&&\hspace{27mm}+C\left\{\int_0^{\ub_1}\ub^{-1}\|\s^{(0,n)}\cla\|^2_{L^2(\Cb_{\ub}^{u_1})}d\ub
\right\}^{1/2}\\
&&\hspace{27mm}+\Cb\sqrt{\max\{\s^{(0,n)}\cB(\ub_1,u_1),\s^{(0,n)}\cBb(\ub_1,u_1)\}}
\cdot\ub_1^{a_0}u_1^{b_0}
\end{eqnarray*}

Here $k$ is again a constant greater than 1 but which can be chosen as close to 1 as we wish by 
suitably restricting $\delta$.

\vspace{2.5mm}

We turn to the derivation of the estimate for $\s^{(0,n)}\clab$. We consider the propagation equation 
for $\s^{(0,n)}\clab$, equation \ref{12.104}. This equation takes the form:
\begin{equation}
\Lb\s^{(0,n)}\clab+n\chib\s^{(0,n)}\s^{(0,n)}\clab=\s^{(0,n)}\chGb
\label{12.300}
\end{equation}
In terms of the $(\ub,u,\vartheta^{\prime\prime})$ coordinates on ${\cal N}$, which are adapted to the 
flow of $\Lb$, this reads:
\begin{equation}
\frac{\partial}{\partial u}\left(\sh^{\prime\prime n/2}\s^{(0,n)}\clab\right)
=\sh^{\prime\prime n/2}\s^{(0,n)}\chGb
\label{12.301}
\end{equation}
Integrating from ${\cal K}$ we obtain:
\begin{eqnarray}
&&\left(\sh^{\prime\prime n/2}\s^{(0,n)}\clab\right)(\ub,u_1,\vartheta^{\prime\prime})
=\left(\sh^{\prime\prime n/2}\s^{(0,n)}\clab\right)(\ub,\ub,\vartheta^{\prime\prime})\nonumber\\
&&\hspace{30mm}+\int_{\ub}^{u_1}\left(\sh^{\prime\prime n/2}\s^{(0,n)}\chGb\right)
(\ub,u,\vartheta^{\prime\prime})du 
\label{12.302}
\end{eqnarray}
In view of the remark \ref{10.607} we rewrite this in the form:
\begin{eqnarray}
&&\left(\s^{(0,n)}\clab\sh^{\prime\prime 1/4}\right)(\ub,u_1,\vartheta^{\prime\prime})=
\left(\frac{\sh(\ub,\ub,\vartheta^{\prime\prime})}
{\sh^{\prime\prime}(\ub,u_1,\vartheta^{\prime\prime})}\right)^{(2n-1)/4}
\s^{(0,n)}\clab(\ub,\ub,\vartheta^{\prime\prime})\nonumber\\
&&\hspace{20mm}+\int_{\ub}^{u_1}\left(\frac{\sh^{\prime\prime}(\ub,u,\vartheta^{\prime\prime})}
{\sh^{\prime\prime}(\ub,u_1,\vartheta^{\prime\prime})}\right)^{(2n-1)/4}
\left(\s^{(0,n)}\chGb\sh^{\prime\prime 1/4}\right)(\ub,u,\vartheta^{\prime\prime})du\nonumber\\
&&\label{12.303}
\end{eqnarray}
In view of \ref{10.611} this implies:
\begin{eqnarray}
&&\left|\left(\s^{(0,n)}\clab\sh^{\prime\prime 1/4}\right)(\ub,u_1,\vartheta^{\prime\prime})\right|
\leq k\left|\left(\s^{(0,n)}\clab\sh^{1/4}\right)(\ub,\ub,\vartheta^{\prime\prime})\right|\nonumber\\
&&\hspace{30mm}+k\int_{\ub}^{u_1}\left|\left(\s^{(0,n)}\chGb\sh^{\prime\prime 1/4}\right)
(\ub,u,\vartheta^{\prime\prime})\right|du
\label{12.304}
\end{eqnarray}
where $k$ is a constant greater than 1, but which can be chosen as close to 1 as we wish by suitably 
restricting $u_1$. Taking $L^2$ norms with respect to $\vartheta^{\prime\prime}\in S^1$ we obtain:
\begin{eqnarray}
&&\left\|\left(\s^{(0,n)}\clab\sh^{\prime\prime 1/4}\right)(\ub,u_1)\right\|_{L^2(S^1)}
\leq k\left\|\left(\s^{(0,n)}\clab\sh^{1/4}\right)(\ub,\ub)\right\|_{L^2(S^1)}\nonumber\\
&&\hspace{25mm}+k\int_{\ub}^{u_1}\left\|\left(\s^{(0,n)}\chGb\sh^{\prime\prime 1/4}\right)(\ub,u)
\right\|_{L^2(S^1)}du
\label{12.305}
\end{eqnarray}
or, in view of \ref{10.607} (see also \ref{10.601}):
\begin{equation}
\|\s^{(0,n)}\clab\|_{L^2(S_{\ub,u_1})}\leq k\|\s^{(0,n)}\clab\|_{L^2(S_{\ub,\ub})}
+k\int_{\ub}^{u_1}\|\s^{(0,n)}\chGb\|_{L^2(S_{\ub,u})}du
\label{12.306}
\end{equation}
Taking then $L^2$ norms with respect to $\ub\in[0,\ub_1]$ we obtain:
\begin{eqnarray}
&&\|\s^{(0,n)}\clab\|_{L^2(C_{u_1}^{\ub_1})}\leq k\|\s^{(0,n)}\clab\|_{L^2({\cal K}^{\ub_1})} 
\label{12.307}\\
&&\hspace{27mm}+k\left\{\int_0^{\ub_1}\left(\int_{\ub}^{u_1}\|\s^{(0,n)}\chGb\|_{L^2(S_{\ub,u})}du\right)^2 d\ub\right\}^{1/2} \nonumber
\end{eqnarray}

We shall presently estimate the leading contributions to 
\begin{equation}
\left\{\int_0^{\ub_1}\left(\int_{\ub}^{u_1}\|\s^{(0,n)}\chGb\|_{L^2(S_{\ub,u})}du\right)^2 d\ub
\right\}^{1/2}
\label{12.308}
\end{equation}
Consider first the $n$th order acoustical difference terms:
\begin{equation}
\pb\s^{(0,n)}\clab+\tilde{\qb}\s^{(0,n)}\cla+(\oBb-\rhob E\lambdab)\s^{(n-1)}\ctchib
+\Bb\s^{(n-1)}\ctchi
\label{12.309}
\end{equation}
The coefficients $\oBb$ and $\Bb$ being given by \ref{12.56} and \ref{12.55} respectively, the 
assumptions \ref{10.439}, \ref{10.440} imply:
\begin{equation}
|\oBb|, \ |\Bb|  \leq C\ub \ \mbox{: in ${\cal R}_{\delta,\delta}$}
\label{12.310}
\end{equation}
Also, 
\begin{equation}
|\pb|\leq C, \ |\tilde{\qb}|\leq C\ub \ \ \mbox{: in ${\cal R}_{\delta,\delta}$}
\label{12.311}
\end{equation}
It follows that the contribution of the terms \ref{12.309} to the integral on the right in \ref{12.306} 
is bounded by:
\begin{eqnarray}
&&C\int_{\ub}^{u_1}\|\s^{(0,n)}\clab\|_{L^2(S_{\ub,u})}du
+C\ub\int_{\ub}^{u_1}\|\s^{(0,n)}\cla\|_{L^2(S_{\ub,u})}du\nonumber\\
&&+C\ub\int_{\ub}^{u_1}\|\s^{(n-1)}\ctchib\|_{L^2(S_{\ub,u})}du
+C\ub\int_{\ub}^{u_1}\|\s^{(n-1)}\ctchi\|_{L^2(S_{\ub,u})}du \nonumber\\
&&\label{12.312}
\end{eqnarray}
Here the first can obviously be absorbed in the integral inequality \ref{12.306} while the second is 
bounded by $C\ub\sqrt{u_1}\|\s^{(0,n)}\cla\|_{L^2(\Cb_{\ub}^{u_1})}$ hence its contribution to 
\ref{12.308} is bounded by:
\begin{equation}
Cu_1^{1/2}\left\{\int_0^{\ub_1}\ub^2\|\s^{(0,n)}\cla\|^2_{L^2(\Cb_{\ub}^{u_1})}d\ub\right\}^{1/2}
\label{12.313}
\end{equation}
Similarly the fourth of \ref{12.312} is bounded by 
$C\ub\sqrt{u_1}\|\s^{(n-1)}\ctchi\|_{L^2(\Cb_{\ub}^{u_1})}$ hence its contribution to \ref{12.308} is 
bounded by:
\begin{equation}
Cu_1^{1/2}\left\{\int_0^{\ub_1}\ub^2\|\s^{(n-1)}\ctchi\|^2_{L^2(\Cb_{\ub}^{u_1})}d\ub\right\}^{1/2}
\label{12.314}
\end{equation}
To estimate the contribution of the third of \ref{12.312} to \ref{12.308} we shall use the following 
lemma. 

\vspace{2.5mm}

\noindent{\bf Lemma 12.5} \ \ Let $f$ be a non-negative function on $R_{\delta,\delta}$. For 
$(\ub_1,u_1)\in R_{\delta,\delta}$ let us define on $[0,\ub_1]$ the function:
$$g(\ub)=\int_{\ub}^{u_1}f(\ub,u)du$$
We then have:
\begin{eqnarray*}
&&\|g\|_{L^2(0,\ub_1)}\leq \int_0^{\ub_1}\|f(\cdot,u)\|_{L^2(0,u)}du\\
&&\hspace{18mm}+\int_{\ub_1}^{u_1}\|f(\cdot,u)\|_{L^2(0,\ub_1)}du
\end{eqnarray*}

\noindent{\em Proof:} We extend for each $\ub\in[0,\ub_1$ the function $f(\ub,\cdot)$ to $[0,u_1]$ 
by defining:
$$\of(\ub,u)=\left\{\begin{array}{cccc} f(\ub,u)&:&\mbox{if}& u\in[\ub,u_1]\\
0&:&\mbox{if}& u\in[0,\ub) \end{array}\right\}$$
We then have:
$$g(\ub)=\int_0^{u_1}\of(\ub,u)du \ \ : \ \ub\in[0,\ub_1]$$
It follows that:
$$\|g\|_{L^2(0,\ub_1)}\leq\int_0^{u_1}\|\of(\cdot,u)\|_{L^2(0,\ub_1)}du$$
Since 
$$\|\of(\cdot,u)\|^2_{L^2(0,\ub_1)}=\int_0^{\ub_1}\of^2(\ub,u)d\ub
=\int_0^{\min\{u,\ub_1\}}f^2(\ub,u)d\ub$$
the lemma follows. 

\vspace{2.5mm}

We apply Lemma 12.5 taking 
\begin{equation}
f(\ub,u)=\ub\|\s^{(n-1)}\ctchib\|_{L^2(S_{\ub,u})}
\label{12.315}
\end{equation}
Then the third of \ref{12.312} is $Cg(\ub)$ and its contribution to \ref{12.308} is 
$C\|g\|_{L^2(0,\ub_1)}$. The lemma then yields that this is bounded by:
\begin{equation}
C\int_0^{\ub_1}\|\ub\s^{(n-1)}\ctchib\|_{L^2(C_u^u)}du
+C\int_{\ub_1}^{u_1}\|\ub\s^{(n-1)}\ctchib\|_{L^2(C_u^{\ub_1})}du
\label{12.316}
\end{equation}

The other leading contributions are those of the principal difference terms:
\begin{eqnarray}
&&\hspace{15mm}\frac{\rhob}{4}\beta_{\Nb}^2 E^n LH
-\frac{\rhob_N}{4}\beta_{\Nb,N}^2 E_N^n L_N H_N \nonumber\\
&&\hspace{15mm}+\rho\beta_{\Nb}\left(c\pi\sbeta-\frac{\beta_{\Nb}}{4}-\frac{\beta_N}{2}\right)E^n\Lb H
\nonumber\\
&&\hspace{30mm}-\rho_N\beta_{\Nb,N}\left(c_N\pi_N\sbeta_N-\frac{\beta_{\Nb,N}}{4}-\frac{\beta_{N,N}}{2}\right)E_N^n\Lb_N H_N\nonumber\\
&&\hspace{15mm}+\rho H\left(c\pi\sbeta-\frac{(\beta_N+\beta_{\Nb})}{2}\right)\Nb^\mu E^n\Lb\beta_\mu
\nonumber\\
&&\hspace{30mm}-\rho_N H_N\left(c_N\pi_N\sbeta_N-\frac{(\beta_{N,N}+\beta_{\Nb,N})}{2}\right)
\Nb_N^\mu E_N^n\Lb_N\beta_{\mu,N}\nonumber\\
&&\hspace{15mm}-\frac{\pi}{2}a\beta_{\Nb}^2 E^{n+1} H 
+\frac{\pi_N}{2}a_N\beta_{\Nb,N}^2 E_N^{n+1} H_N \label{12.317}
\end{eqnarray}
(conjugate of \ref{12.200}). Here again we must keep track not only of the actually principal terms, 
namely those containing derivatives of the $\beta_\mu$ of order $n+1$, 
but also of the terms containing acoustical quantities of order $n$. 
The terms of order $n$ with vanishing principal acoustical part can be ignored. Consider the 
first of the above difference terms. The actual principal difference term in question is:
\begin{equation}
-\frac{\rhob}{2}\beta_{\Nb}^2 H^\prime(g^{-1})^{\mu\nu}\beta_\nu E^n L\beta_\mu 
+\frac{\rhob_N}{2}\beta_{\Nb,N}^2 H^\prime_N(g^{-1})^{\mu\nu}\beta_{\nu,N} E_N^n L_N\beta_{\mu,N}
\label{12.318}
\end{equation}
By the formula \ref{10.368} and the analogous formula for the $N$th approximants, this is expressed 
as differences of $E$, $N$, and $\Nb$ components:
\begin{eqnarray}
&&-\frac{\rhob}{2}\beta_{\Nb}^2\sbeta\eta^2 H^\prime E^\mu E^n L\beta_\mu 
+\frac{\rhob_N}{2}\beta_{\Nb,N}^2\sbeta_N\eta_N^2 H^\prime_N E_N^\mu E_N^n L_N\beta_{\mu,N} \nonumber\\
&&+\frac{\rhob}{2}\beta_{\Nb}^2\frac{\beta_{\Nb}}{2c}\eta^2 H^\prime N^\mu E^n L\beta_\mu 
-\frac{\rhob_N}{2}\beta_{\Nb,N}^2\frac{\beta_{\Nb,N}}{2c_N}\eta_N^2 H^\prime_N N_N^\mu E_N^n L_N\beta_{\mu,N} \nonumber\\
&&+\frac{\rhob}{2}\beta_{\Nb}^2\frac{\beta_N}{2c}\eta^2 H^\prime\Nb^\mu E^n L\beta_\mu 
-\frac{\rhob_N}{2}\beta_{\Nb,N}^2\frac{\beta_{N,N}}{2c_N}\eta_N^2 H^\prime_N\Nb_N^\mu E_N^n L_N\beta_{\mu,N} \nonumber\\
&&\label{12.319}
\end{eqnarray}
This is similar to the differences \ref{12.164} but with extra $\rhob$, $\rhob_N$ factors. 
Thus in view of \ref{12.170} the contribution of the $\Nb$ component difference to the integral on the 
right in \ref{12.306} is bounded by:
\begin{eqnarray}
&&C\sqrt{\ub}\int_{\ub}^{u_1}u\|\sqrt{a}\s^{(E;0,n)}\csxi\|_{L^2(S_{\ub,u})}du
+C\int_{\ub}^{u_1}u^2\|\s^{(0,n)}\clab\|_{L^2(S_{\ub,u})}du\nonumber\\
&&+C\ub\int_{\ub}^{u_1}u^2\|\s^{(n-1)}\ctchib\|_{L^2(S_{\ub,u})}du
+C\ub\int_{\ub}^{u_1}u^2\|\s^{(n-1)}\ctchi\|_{L^2(S_{\ub,u})}du\nonumber\\
&&\label{12.320}
\end{eqnarray}
Here the last three terms are dominated by the corresponding terms in \ref{12.312} while the 1st term 
is bounded by 
\begin{equation}
C\ub^{1/2}u_1^{3/2}\sqrt{\s^{(E;0,n)}\cEb^{u_1}(\ub)}
\label{12.321}
\end{equation}
hence its contribution to \ref{12.308} is bounded by:
\begin{equation}
\frac{C\sqrt{\s^{(E;0,n)}\cBb(\ub_1,u_1)}}{\sqrt{2a_0+2}}\cdot \ub_1^{a_0}u_1^{b_0}\cdot\ub_1 u_1^{3/2}
\label{12.322}
\end{equation}
In view of \ref{12.178} the contribution of the $N$ component difference to the integral on the right 
in \ref{12.306} is bounded by:
\begin{equation}
C\int_{\ub}^{u_1}u^2\|N^\mu L(E^n\beta_\mu-E_N^n\beta_{\mu,N})\|{L^2(S_{\ub,u})}du
+C\ub\int_{\ub}^{u_1}u^2\|\s^{(n-1)}\ctchi\|_{L^2(S_{\ub,u})}du
\label{12.323}
\end{equation}
The 2nd coincides with the last of \ref{12.320}. To estimate the contribution of the 1st we use 
\ref{12.180}. Since $\ogamma\sim u$ the partial contribution of the 2nd term on the left in 
\ref{12.180} to the 1st of \ref{12.323} is bounded by $Cu$ times the contribution of the 
$\Nb$ component difference to the integral on the right in \ref{12.306}. On the other hand the 
partial contribution of the right hand side of \ref{12.180} to the 1st of \ref{12.323} is:
\begin{equation}
C\int_{\ub}^{u_1}u^2\|\s^{(Y;0,n)}\cxi_L\|_{L^2(S_{\ub,u})}du
\label{12.324}
\end{equation}
To estimate the contribution of \ref{12.324} to \ref{12.308} we apply Lemma 12.5 taking 
\begin{equation}
f(\ub,u)=u^2\|\s^{(Y;0,n)}\cxi_L\|_{L^2(S_{\ub,u})}
\label{12.325}
\end{equation}
Then \ref{12.324} is $Cg(\ub)$ and its contribution to \ref{12.308} is $C\|g\|_{L^2(0,\ub_1)}$. 
The lemma then yields that this is bounded by:
\begin{eqnarray}
&&C\int_0^{\ub_1}u^2\|\s^{(Y;0,n)}\cxi_L\|_{L^2(C_u^u)}du
+C\int_{\ub_1}^{u_1} u^2\|\s^{(Y;0,n)}\cxi_L\|_{L^2(C_u^{\ub_1})}du\nonumber\\
&&\leq C\int_0^{\ub_1}u^2\sqrt{\s^{(Y;0,n)}\cE^u(u)}du
+C\int_{\ub_1}^{u_1}u^2\sqrt{\s^{(Y;0,n)}\cE^{\ub_1}(u)}du\nonumber\\
&&\leq C\sqrt{\s^{(Y;0,n)}\cB(\ub_1,\ub_1)}\int_0^{\ub_1}u^{a_0+b_0+2}du\nonumber\\
&&\hspace{30mm}+C\sqrt{\s^{(Y;0,n)}\cB(\ub_1,u_1)}\cdot\ub_1^{a_0}\cdot\int_{\ub_1}^{u_1}u^{b_0+2}du\nonumber\\
&&\leq C\sqrt{\s^{(Y;0,n)}\cB(\ub_1,u_1)}\cdot\left\{\frac{\ub_1^{a_0+b_0+3}}{a_0+b_0+3}
+\ub_1^{a_0}\cdot\frac{(u_1^{b_0+3}-\ub_1^{b_0+3})}{b_0+3}\right\}\nonumber\\
&&\leq C\frac{\sqrt{\s^{(Y;0,n)}\cB(\ub_1,u_1)}}{b_0+3}\cdot\ub_1^{a_0}u_1^{b_0}\cdot u_1^3
\label{12.326}
\end{eqnarray}
In view of \ref{12.174} the contribution of the $E$ component difference to the integral on the right 
in \ref{12.306} is bounded by:
\begin{equation}
C\int_{\ub}^{u_1}u^2\|\s^{(E;0,n)}\cxi_L\|_{L^2(S_{\ub,u})}du
+C\ub\int_{\ub}^{u_1}u^2\|\s^{(n-1)}\ctchi\|_{L^2(S_{\ub,u})}du
\label{12.327}
\end{equation}
Here the 2nd term is dominated by the corresponding term in \ref{12.312} while the 1st term is 
similar to \ref{12.324} and its contribution to \ref{12.308} is bounded by:
\begin{equation}
C\frac{\sqrt{\s^{(E;0,n)}\cB(\ub_1,u_1)}}{b_0+3}\cdot\ub_1^{a_0}u_1^{b_0}\cdot u_1^3
\label{12.328}
\end{equation}

Consider now the second of the principal difference terms \ref{12.317}. The actual principal difference 
term in question is:
\begin{eqnarray}
&&-2\rho\beta_{\Nb}\left(c\pi\sbeta-\frac{\beta_{\Nb}}{4}-\frac{\beta_{N}}{2}\right)
H^\prime(g^{-1})^{\mu\nu}\beta_\nu E^n\Lb\beta_\mu
\nonumber\\
&&\hspace{7mm}+2\rho_N\beta_{\Nb,N}\left(c_N\pi_N\sbeta_N-\frac{\beta_{\Nb,N}}{4}-\frac{\beta_{N,N}}{2}\right)H_N^\prime(g^{-1})^{\mu\nu}\beta_{\nu,N}E_N^n\Lb_N\beta_{\mu,N}\nonumber\\
&&\label{12.329}
\end{eqnarray}
This is similar to \ref{12.271} but with extra $\rho$, $\rho_N$ factors.
It is similarly expressed as the sum of differences of $E$, $N$, and $\Nb$ components. 
It follows that up to terms which either coincide or are dominated by the contributions of 
the last three of \ref{12.312}, the contribution to \ref{12.308} of the $N$ component difference 
is bounded by:
\begin{equation}
\frac{C\sqrt{\s^{(E;0,n)}\cBb(\ub_1,u_1)}}{\sqrt{2a_0+2}}\cdot\ub_1^{a_0}u_1^{b_0}\cdot\ub_1 u_1^{3/2}
\label{12.330}
\end{equation}
Up to terms which are dominated by the contributions of the last two of \ref{12.312}, 
the contribution to \ref{12.308} of the $\Nb$ component difference is bounded by: 
\begin{eqnarray}
&&C\left\{\int_0^{\ub_1}\ub\|\s^{(0,n)}\cla\|^2_{L^2(\Cb_{\ub}^{u_1})}d\ub\right\}^{1/2} \nonumber\\
&&+\frac{C\sqrt{\s^{(E;0,n)}\cBb(\ub_1,u_1)}}{\sqrt{2a_0+2}}\cdot\ub_1^{a_0}u_1^{b_0}\cdot 
\ub_1 u_1^{1/2} \nonumber\\
&&+\frac{C\sqrt{\s^{(Y;0,n)}\cBb(\ub_1,u_1)}}{\sqrt{2a_0+2}}\cdot\ub_1^{a_0}u_1^{b_0}\cdot\ub_1 
\label{12.331}
\end{eqnarray}
Also, up to terms which are dominated by the contributions of the last two of \ref{12.312}, 
the contribution to \ref{12.308} of the $E$ component difference is bounded by: 
\begin{equation}
\frac{C\sqrt{\s^{(E;0,n)}\cBb(\ub_1,u_1)}}{\sqrt{2a_0+3}}\cdot\ub_1^{a_0}u_1^{b_0}
\cdot\ub_1^{3/2}u_1^{1/2}
\label{12.332}
\end{equation}
Moreover, the third of the principal difference terms \ref{12.317} is similar to the $\Nb$ component difference. 

Consider finally the last of the principal difference terms \ref{12.317}. The actual principal 
difference term in question is:
\begin{equation}
\pi a\beta_{\Nb}^2 H^\prime (g^{-1})^{\mu\nu}\beta_\nu E^{n+1}\beta_\mu
-\pi_N a_N\beta_{\Nb,N}^2 H^\prime_N (g^{-1})^{\mu\nu}\beta_{\nu,N} E_N^{n+1}\beta_{\mu,N}
\label{12.333}
\end{equation}
This is similar to \ref{12.290} but, since $a=c\rho\rhob$, $a_N=c_N\rho_N\rhob_N$, with extra 
$\rho$, $\rho_N$ factors. It is similarly expressed as the sum of differences 
of $E$, $N$, and $\Nb$ components. 
It follows that up to terms which either coincide or are dominated by the contributions of 
the last three of \ref{12.312}, the contribution to \ref{12.308} of the $\Nb$ component difference 
is bounded by:
\begin{equation}
\frac{C\sqrt{\s^{(E;0,n)}\cBb(\ub_1,u_1)}}{\sqrt{2a_0+3}}\cdot\ub_1^{a_0}u_1^{b_0}\cdot
\ub_1^{3/2}u_1^{1/2}
\label{12.334}
\end{equation}
Similarly, the contribution to \ref{12.308} of the $N$ component difference is bounded by:
\begin{eqnarray}
&&\frac{C\sqrt{\s^{(Y;0,n)}\cBb(\ub_1,u_1)}}{\sqrt{2a_0+2}}\cdot\ub_1^{a_0}u_1^{b_0}\cdot\ub_1 u_1^{3/2}
\nonumber\\
&&+\frac{C\sqrt{\s^{(E;0,n)}\cBb(\ub_1,u_1)}}{\sqrt{2a_0+3}}\cdot\ub_1^{a_0}u_1^{b_0}\cdot
\ub_1^{3/2}u_1^{3/2}
\label{12.335}
\end{eqnarray}
and the contribution to \ref{12.308} of the $E$ component difference is bounded by:
\begin{equation}
\frac{C\sqrt{\s^{(E;0,n)}\cBb(\ub_1,u_1)}}{\sqrt{2a_0+2}}\cdot\ub_1^{a_0}u_1^{b_0}\cdot\ub_1u_1^{3/2}
\label{12.336}
\end{equation}

Absorbing the first of \ref{12.312} in the integral inequality \ref{12.306}, we summarize the above 
results in the form in which they will be used in the sequel in the following lemma. 

\vspace{2.5mm}

\noindent{\bf Lemma 12.6} \ \ The next to the top order acoustical quantity 
$\|\s^{(0,n)}\clab\|_{L^2(C_u^{\ub})}$ satisfies to principal terms the following inequality:
\begin{eqnarray*}
&&\|\s^{(0,n)}\clab\|_{L^2(C_{u_1}^{\ub_1})}\leq k\|\s^{(0,n)}\clab\|_{L^2({\cal K}^{\ub_1})}\\
&&\hspace{27mm}+Cu_1^{1/2}\left\{\int_0^{\ub_1}\ub^2\|\s^{(n-1)}\ctchi\|^2_{L^2(\Cb_{\ub}^{u_1})}d\ub\right\}^{1/2}\\
&&\hspace{20mm}+C\int_0^{\ub_1}\|\ub\s^{(n-1)}\ctchib\|_{L^2(C_u^u)}du
+C\int_{\ub_1}^{u_1}\|\ub\s^{(n-1)}\ctchib\|_{L^2(C_u^{\ub_1})}du\\
&&\hspace{27mm}+C\left\{\int_0^{\ub_1}\ub\|\s^{(0,n)}\cla\|^2_{L^2(\Cb_{\ub}^{u_1})}d\ub\right\}^{1/2}\\
&&\hspace{27mm}+\Cb\sqrt{\max\{\s^{(0,n)}\cB(\ub_1,u_1),\s^{(0,n)}\cBb(\ub_1,u_1)\}}\cdot\ub_1^{a_0}u_1^{b_0}
\cdot u_1
\end{eqnarray*}

Here $k$ is again a constant greater than 1 but which can be chosen as close to 1 as we wish by 
suitably restricting $\delta$. 

\vspace{2.5mm}

To proceed we must bring in the boundary conditions for $\tchib$ and $\lambdab$ on ${\cal K}$, 
equations \ref{10.756} and \ref{10.752}. These express on ${\cal K}$ $\tchib$ in terms of $\tchi$ 
and $\lambdab$ in terms of $\lambda$ respectively. These conditions will allow us to estimate 
the boundary terms $\|\s^{(n-1)}\ctchib\|_{L^2({\cal K}^{\ub_1})}$ and 
$\|\s^{(0,n)}\clab\|_{L^2({\cal K}^{\ub_1})}$ in Lemmas 12.4 and 12.6 in terms of 
$\|\s^{(n-1)}\ctchi\|_{L^2{\cal K}^{\ub_1})}$ and $\|\s^{(0,n)}\cla\|_{L^2({\cal K}^{\ub_1})}$ 
respectively. 

Equations \ref{10.763} and \ref{10.764} with $l=n-1$ read:
\begin{eqnarray}
&&rE^{n-1}\tchib=E^{n-1}\tchi+2c\frac{E^\mu E^n\triangle\beta_\mu}{\epb}\nonumber\\
&&\hspace{12mm}+(\beta_N-r\beta_{\Nb})(\sbeta E^n H+HE^\mu E^n\beta_\mu)+{\cal L}_{n-1} \ \  \mbox{: on ${\cal K}$}\nonumber\\
&&\label{12.337}
\end{eqnarray}
and:
\begin{eqnarray}
&&r_N E_N^{n-1}\tchib_N=E_N^{n-1}\tchi_N+2c_N\frac{E_N^\mu E_N^n\triangle_N\beta_{\mu,N}}{\epb_N}\nonumber\\
&&\hspace{12mm}+(\beta_{N,N}-r_N\beta_{\Nb,N})(\sbeta_N E_N^n H_N+H_N E_N^\mu E_N^n\beta_{\mu,N})
+{\cal L}_{n-1,N}\nonumber\\
&&\hspace{12mm}+E_N^{n-1}(2c_N\epb_N^{-1}\se_N\iota_{E,N}) \ \  \mbox{: on ${\cal K}$}
\label{12.338}
\end{eqnarray}
where ${\cal L}_{n-1}$ is of order $n-1$ and ${\cal L}_{n-1,N}$ is the analogous quantity 
for the $N$th approximants. Subtracting \ref{12.338} from \ref{12.337} we obtain:
\begin{eqnarray}
&&r\s^{(n-1)}\ctchib=\s^{(n-1)}\ctchi+\frac{2c}{\epb}(E^\mu E^n\triangle\beta_\mu
-E_N^\mu E_N^n\triangle_N\beta_{\mu,N})\nonumber\\
&&\hspace{15mm}-(r-r_N)E_N^{n-1}\tchib_N+\left(\frac{2c}{\epb}-\frac{2c_N}{\epb_N}\right)
E_N^\mu E_N^n\triangle_N\beta_{\mu,N}\nonumber\\
&&\hspace{15mm}+(\beta_N-r\beta_{\Nb})(\sbeta E^n H+HE^\mu E^n\beta_\mu)\nonumber\\
&&\hspace{30mm}-(\beta_{N,N}-r_N\beta_{\Nb,N})(\sbeta_N E_N^n H_N+H_N E_N^\mu E_N^n\beta_{\mu,N})\nonumber\\
&&\hspace{15mm}+{\cal L}_{n-1}-{\cal L}_{n-1,N}+\tilde{\vep}_{n-1,N}
\label{12.339}
\end{eqnarray}
where $\tilde{\vep}_{n-1,N}$ is the error term:
\begin{equation}
\tilde{\vep}_{n-1,N}=-E_N^{n-1}(2c_N\epb^{-1}\se_N\iota_{E,N})=O(\tau^N)
\label{12.340}
\end{equation}
(see \ref{10.758}). In \ref{12.339} the difference term 
\begin{eqnarray*}
&&(\beta_N-r\beta_{\Nb})(\sbeta E^n H+HE^\mu E^n\beta_\mu)\nonumber\\
&&\hspace{20mm}-(\beta_{N,N}-r_N\beta_{\Nb,N})(\sbeta_N E_N^n H_N+H_N E_N^\mu E_N^n\beta_{\mu,N})
\end{eqnarray*}
is of order $n$ but does not contain acoustical quantities of order $n$. Therefore this difference 
term can be ignored when estimating the principal contributions. The other terms on the 
right hand side of \ref{12.339} can likewise be ignored, with the exception of the first two. 
The 2nd term, which involves the jump $\triangle\beta_\mu$, is also of order $n$ and does not contain acoustical terms of order $n$, however 
in view of the fact that the factor $\epb^{-1}$ is only bounded by $C\tau^{-1}$ (see \ref{10.775}, 
\ref{10.792}, \ref{10.793}) we must appeal to estimates for the $T$ derivative of 
$E^\mu E^n\triangle\beta_\mu-E_N^\mu E_N^n\triangle_N\beta_{\mu,N}$, a quantity of order $n+1$, 
to deduce an appropriate estimate. 

We shall estimate $\|r\s^{(n-1)}\ctchib\|_{L^2({\cal K}^{\tau_1})}$. Since 
$|\epb^{-1}|\leq C\tau^{-1}$ the contribution to this of the 2nd term on the right in \ref{12.339} 
is bounded by a constant multiple of:
\begin{equation}
\left(\int_0^{\tau_1}\tau^{-2}\|E^\mu E^n\triangle\beta_\mu-E_N^\mu E_N^n\triangle_N\beta_{\mu,N}\|^2_
{L^2(S_{\tau,\tau})}d\tau\right)^{1/2}
\label{12.341}
\end{equation}
According to \ref{10.795} with $l=n-1$ we have:
\begin{eqnarray}
&&E^\mu E^n\triangle\beta_\mu-E_N^\mu E_N^n\triangle_N\beta_{\mu,N}=
\left.(E^\mu E^n\beta_\mu-E_N^\mu E_N^n\beta_{\mu,N})\right|_{{\cal K}}\nonumber\\
&&\hspace{20mm}-(E^\mu E^n\beta^\prime_\mu\circ(f,w,\psi)-E_N^\mu E_N^n\beta^\prime_\mu\circ(f_N,w_N,\psi_N))\nonumber\\
&&\label{12.342}
\end{eqnarray}
To derive an appropriate estimate for the first term on the right in \ref{12.342}, we appeal to the 
fact that: 
\begin{equation}
\left.(E^\mu E^n\beta_\mu-E_N^\mu E_N^n\beta_{\mu,N})\right|_{\partial_-{\cal K}}=0 , \ \ 
\partial_-{\cal K}=S_{0,0}
\label{12.343}
\end{equation}
to express, in $(\tau,\vartheta)$ coordinates on ${\cal K}$, since $T=\partial/\partial\tau$ in these 
coordinates,
\begin{equation}
(E^\mu E^n\beta_\mu-E_N^\mu E_N^n\beta_{\mu,N})(\tau,\vartheta)
=\int_0^{\tau}(T(E^\mu E^n\beta_\mu-E_N^\mu E_N^n\beta_{\mu,N}))(\tau^\prime,\vartheta)
d\tau^\prime
\label{12.344}
\end{equation}

For functions on ${\cal K}$ we have the following analogue of \ref{10.350} - \ref{10.352} and of 
\ref{10.605} - \ref{10.607}. If $f$ is a function on $S^1$ its $L^2$ norm is:
\begin{equation}
\|f\|_{L^2(S^1)}=\sqrt{\int_{\vartheta\in S^1}f^2(\vartheta)d\vartheta}
\label{12.345}
\end{equation}
On the other hand, if $f$ is a function on ${\cal K}$ and we represent $f$ in $(\tau,\vartheta)$ 
coordinates, the $L^2$ norm of $f$ on $S_{\tau,\tau}$ is:
\begin{equation}
\|f\|_{L^2(S_{\tau,\tau})}=\sqrt{\int_{\vartheta\in S^1}\left(f^2\sqrt{\sh}\right)(\tau,\vartheta)
d\vartheta}
\label{12.346}
\end{equation}
$\sqrt{\sh(\tau,\vartheta)}d\vartheta$ being the element of arc length on $S_{\tau,\tau}$. 
Comparing \ref{12.345} with \ref{12.346} we see that:
\begin{equation}
\|f\|_{L^2(S_{\tau,\tau})}=\left\|\left(f\sh^{1/4}\right)(\tau)\right\|_{L^2(S^1)}
\label{12.347}
\end{equation}
Here $g$ being an arbitrary function on ${\cal K}$ we denote by $g(\tau)$ the restriction of $g$ to 
$S_{\tau,\tau}$ represented in terms of the $\vartheta$ coordinate. 

With the above remark in mind, we rewrite \ref{12.344} in the form:
\begin{eqnarray}
&&\left((E^\mu E^n\beta_\mu-E_N^\mu E_N^n\beta_{\mu,N})\sh^{1/4}\right)(\tau,\vartheta)= 
\label{12.348}\\
&&\hspace{15mm}\int_0^{\tau}\left(\frac{\sh(\tau,\vartheta)}{\sh(\tau^\prime,\vartheta)}\right)^{1/4}
\left(T(E^\mu E^n\beta_\mu-E_N^\mu E_N^n\beta_{\mu,N})\sh^{1/4}\right)(\tau^\prime,\vartheta)
d\tau^\prime \nonumber
\end{eqnarray}
Since in $(\tau,\vartheta)$ coordinates on ${\cal K}$ the vectorfields $T$ and $E$ are given by:
\begin{equation}
T=\frac{\partial}{\partial\tau}, \ \ \ E=\frac{1}{\sqrt{\sh}}\frac{\partial}{\partial\vartheta}
\label{12.349}
\end{equation}
Since $T=L+\Lb$ in ${\cal N}$, the first two of the commutation relations \ref{3.a14} give:
\begin{equation}
[T,E]=-2(\chi+\chib)E \ \ \mbox{: in ${\cal N}$}
\label{12.350}
\end{equation}
On ${\cal K}$, by \ref{12.349} this reads:
\begin{equation}
\frac{\partial\sh}{\partial\tau}=2(\chi+\chib)\sh
\label{12.351}
\end{equation}
In view of the fact that 
\begin{equation}
\left.\sh\right|_{S_{0,0}}=1
\label{12.352}
\end{equation}
$\vartheta$ being arc length along $\partial_-{\cal K}$, integrating \ref{12.351} we obtain:
\begin{equation}
\sqrt{\sh(\tau,\vartheta)}=e^{\int_0^\tau(\chi+\chib)(\tau^\prime,\vartheta)d\tau^\prime}
\label{12.353}
\end{equation}
Now by \ref{10.355} and \ref{10.607}:
\begin{equation}
\sup_{S_{\tau,\tau}}|\chi+\chib|\leq C \  : \forall \tau\in[0,\delta]
\label{12.354}
\end{equation}
It follows that:
\begin{equation}
e^{-C(\tau-\tau^\prime)}\leq \frac{\sqrt{\sh(\tau,\vartheta)}}{\sqrt{\sh(\tau^\prime,\vartheta)}}
\leq e^{C(\tau-\tau^\prime)}
\label{12.355}
\end{equation}
therefore, in regard to \ref{12.348}, 
\begin{equation}
\left(\frac{\sh(\tau,\vartheta)}{\sh(\tau^\prime,\vartheta)}\right)^{1/4}\leq k
\label{12.356}
\end{equation}
where $k$ is a constant greater than 1, but which can be chosen as close to 1 as we wish by suitably 
restricting $\delta$. Consequently \ref{12.348} implies:
\begin{eqnarray}
&&\left|\left((E^\mu E^n\beta_\mu-E_N^\mu E_N^n\beta_{\mu,N})\sh^{1/4}\right)(\tau,\vartheta)\right|
\leq \label{12.357}\\
&&\hspace{25mm}
k\int_0^{\tau}\left|\left(T(E^\mu E^n\beta_\mu-E_N^\mu E_N^n\beta_{\mu,N})\sh^{1/4}\right)
(\tau^\prime,\vartheta)\right|d\tau^\prime \nonumber
\end{eqnarray}
Taking $L^2$ norms with respect to $\vartheta\in S^1$ this implies:
\begin{eqnarray}
&&\left\|\left((E^\mu E^n\beta_\mu-E_N^\mu E_N^n\beta_{\mu,N})\sh^{1/4}\right)(\tau)\right\|_{L^2(S^1)}
\leq \label{12.358}\\
&&\hspace{25mm}k\int_0^{\tau}\left\|\left(T(E^\mu E^n\beta_\mu
-E_N^\mu E_N^n\beta_{\mu,N})\sh^{1/4}\right)(\tau^\prime)\right\|_{L^2(S^1)}d\tau^\prime \nonumber
\end{eqnarray}
or, in view of \ref{12.347}:
\begin{equation}
\|E^\mu E^n\beta_\mu-E_N^\mu E_N^n\beta_{\mu,N}\|_{L^2(S_{\tau,\tau})}\leq 
k\int_0^{\tau}\|T(E^\mu E^n\beta_\mu-E_N^\mu E_N^n\beta_{\mu,N})\|_{L^2(S_{\tau^\prime,\tau^\prime})}
d\tau^\prime
\label{12.359}
\end{equation}
Up to terms of order $n$ with vanishing $n$th order acoustical part 
$T(E^\mu E^n\beta_\mu-E_N^\mu E_N^n\beta_{\mu,N})$ is equal to $\s^{(E;0,n)}\cxi_T$ . Hence the 
integral on the right in \ref{12.359} is bounded up to terms which can be ignored by:
\begin{eqnarray}
&&\int_0^{\tau}\|\s^{(E;0,n)}\cxi_T\|_{L^2(S_{\tau^\prime,\tau^\prime})}d\tau^{\prime}
\leq \tau^{1/2}\|\s^{(E;0,n)}\cxi_T\|_{L^2({\cal K}^{\tau})} \nonumber\\
&&\leq C\tau^{1/2}\sqrt{\s^{(E;0,n)}\cF^{\prime\tau}}\leq C\tau^{a_0+b_0+\frac{1}{2}}\sqrt{\s^{(E;0,n)}\cA(\tau_1)} 
\label{12.360}
\end{eqnarray}
for all $\tau\in [0,\tau_1]$, by \ref{9.a35} and \ref{9.299}. Therefore, denoting, for any 
$m=0,...,n$, 
\begin{equation}
c_m=a_m+b_m
\label{12.361}
\end{equation}
the contribution to \ref{12.341} of the first term on the right in \ref{12.342} is bounded, 
to leading terms, by:
\begin{equation}
C\sqrt{\s^{(E;0,n)}\cA(\tau_1)}\left(\int_0^{\tau_1}\tau^{2c_0-1}d\tau\right)^{1/2}
=\frac{C\sqrt{\s^{(E;0,n)}\cA(\tau_1)}}{\sqrt{2c_0}}\cdot\tau_1^{c_0}
\label{12.362}
\end{equation}

We turn to the second term on the right in \ref{12.342}. By \ref{10.798} - \ref{10.802} and 
\ref{10.803} with $l=n-1$ the principal part of the second term on the right in \ref{12.342} is:
\begin{eqnarray}
&&-\sh^{-n/2}E^\mu\left\{\left(\frac{\partial\beta^\prime_\mu}{\partial t}\right)\circ(f,w,\psi)
\tau^2\Omega^n\chf+\left(\frac{\partial\beta^\prime_\mu}{\partial u^\prime}\right)\circ(f,w,\psi)
\tau \Omega^n\cv\right.\nonumber\\
&&\hspace{55mm}\left.+\left(\frac{\partial\beta^\prime_\mu}{\partial\vartheta^\prime}\right)
\circ(f,w,\psi)\tau^3\Omega^n\cga\right\}\nonumber\\
&&\label{12.363}
\end{eqnarray}
In view of \ref{10.805}, it follows that the principal part of the 2nd term on the right in 
\ref{12.342} is pointwise bounded by:
\begin{equation}
C|\tau^2\Omega^n\chf|+C|\tau^2\Omega^n\cv|+C|\tau^3\Omega^n\cga|
\label{12.364}
\end{equation}
Now the argument leading from \ref{12.344} to \ref{12.359} applies with an arbitrary function 
$f$ on ${\cal K}$ which vanishes on $\partial_-{\cal K}$ in the role of 
$E^\mu E^n\beta_\mu-E_N^\mu E_N^n\beta_{\mu,N}$. That is, for any function $f$ on ${\cal K}$ 
vanishing on $\partial_-{\cal K}=S_{0,0}$, so that
\begin{equation}
f(\tau,\vartheta)=\int_0^{\tau}(Tf)(\tau^\prime,\vartheta)d\tau^\prime
\label{12.365}
\end{equation}
we have:
\begin{eqnarray}
&&\|f\|_{L^2(S_{\tau,\tau})}\leq k\int_0^{\tau}\|Tf\|_{L^2(S_{\tau^\prime,\tau^\prime})}d\tau^\prime
\nonumber\\
&&\hspace{18mm}\leq k\tau^{1/2}\|Tf\|_{L^2({\cal K}^{\tau})}
\label{12.366}
\end{eqnarray}
Taking then successively $\Omega^n\chf$, $\Omega^n\cv$, $\Omega^n\cga$ in the role of the function $f$ 
we conclude that the principal part of the 2nd tern on the right in \ref{12.342} is bounded in 
$L^2(S_{\tau,\tau})$ by:
\begin{equation}
C\tau^{5/2}\|\Omega^n T\chf\|_{L^2({\cal K}^{\tau})}
+C\tau^{5/2}\|\Omega^n T\cv\|_{L^2({\cal K}^{\tau})}
+C\tau^{7/2}\|\Omega^n T\cga\|_{L^2({\cal K}^{\tau})}
\label{12.367}
\end{equation}
(recall that $[T,\Omega]=0$). Therefore the contribution to \ref{12.341} of the second term on the 
right in \ref{12.342} is bounded, to leading terms, by:
\begin{eqnarray}
&&C\left(\int_0^{\tau_1}\tau^3\|\Omega^n T\chf\|^2_{L^2({\cal K}^{\tau})}d\tau\right)^{1/2}
+C\left(\int_0^{\tau_1}\tau^3\|\Omega^n T\cv\|^2_{L^2({\cal K}^{\tau})}d\tau\right)^{1/2} \nonumber\\
&&+C\left(\int_0^{\tau_1}\tau^5\|\Omega^n T\cga\|^2_{L^2({\cal K}^{\tau})}d\tau\right)^{1/2} 
\label{12.368}
\end{eqnarray}
The above yield through \ref{12.339} the following lemma. 

\vspace{2.5mm}

\noindent{\bf Lemma 12.7} \ \ The boundary values on ${\cal K}$ of the next to the top order 
acoustical difference quantity $\s^{(n-1)}\ctchib$ satisfy, to leading terms, the inequality:
\begin{eqnarray*}
&&\|r\s^{(n-1)}\ctchib\|_{L^2({\cal K}^{\tau_1})}\leq k\|\s^{(n-1)}\ctchi\|_{L^2({\cal K}^{\tau_1})}\\
&&\hspace{30mm}+\frac{C\sqrt{\s^{(E;0,n)}\cA(\tau_1)}}{\sqrt{2c_0}}\cdot\tau_1^{c_0}\\
&&\hspace{30mm}+C\left(\int_0^{\tau_1}\tau^3\|\Omega^n T\chf\|^2_{L^2({\cal K}^{\tau})}d\tau\right)^{1/2}
\\
&&\hspace{30mm}+C\left(\int_0^{\tau_1}\tau^3\|\Omega^n T\cv\|^2_{L^2({\cal K}^{\tau})}d\tau\right)^{1/2}
\\
&&\hspace{30mm}+C\left(\int_0^{\tau_1}\tau^5\|\Omega^n T\cga\|^2_{L^2({\cal K}^{\tau})}d\tau\right)^{1/2} 
\end{eqnarray*}

\vspace{2.5mm}

Equations \ref{10.769} and \ref{10.770} with $m=1$, $l=n-2$ read:
\begin{equation}
rE^n\lambdab=E^n\lambda-\lambdab E^n r+{\cal N}_{1,n-2} \ \ \mbox{: on ${\cal K}$}
\label{12.369}
\end{equation}
and:
\begin{eqnarray}
&&r_N E_N^n\lambdab_N=E_N^n\lambda_N-\lambdab_N E_N^n r_N+{\cal N}_{1,n-2,N} \nonumber\\
&&\hspace{18mm}-E_N^n\hat{\nu}_N \ \ \mbox{: on ${\cal K}$}
\label{12.370}
\end{eqnarray}
where ${\cal N}_{1,n-2}$ is of order $n-1$ and ${\cal N}_{1,n-2,N}$ is the analogous quantity for the 
$N$th approximants. Subtracting \ref{12.370} from \ref{12.369} we obtain:
\begin{eqnarray}
&&r\s^{(0,n)}\clab=\s^{(0,n)}\cla-\lambdab(E^n r-E_N^n r_N) \nonumber\\
&&\hspace{17mm}-(r-r_N)E_N^n\lambda_N-(\lambdab-\lambdab_N)E_N^n r_N \nonumber\\
&&\hspace{17mm}+{\cal N}_{1,n-2}-{\cal N}_{1,n-2,N}+E_N^n\hat{\nu}_N
\label{12.371}
\end{eqnarray}
Here by \ref{10.757} the last term is $O(\tau^{N+1})$. In estimating 
$\|\s^{(0,n)}\clab\|_{L^2({\cal K}^{\tau_1})}$ only the first two terms on the right in \ref{12.371} 
need be considered. Since along ${\cal K}$ $\lambdab\sim\tau$ while according to \ref{10.775} 
$r\sim\tau$, the contribution of the 2nd term on the right in \ref{12.371} to 
$\|\s^{(0,n)}\clab\|_{L^2({\cal K}^{\tau_1})}$ is bounded by a constant multiple of:
\begin{equation}
\left(\int_0^{\tau_1}\|E^n r-E_N^n r_N\|^2_{L^2(S_{\tau,\tau})}d\tau\right)^{1/2}
\label{12.372}
\end{equation}
By \ref{10.792} and \ref{10.808} the leading part of $E^n r-E_N^n r_N$ is:
$$j(\kappa,\epb)(E^n\epb-E_N^N\epb_N)$$
and the principal part of this is $j(\kappa,\epb)$ times:
\begin{equation}
\Nb^\mu E^n\triangle\beta_\mu-\Nb_N^\mu E_N^n\triangle_N\beta_{\mu,N}
+\triangle\beta_\mu E^n\Nb^\mu-\triangle_N\beta_{\mu,N} E_N^n\Nb_N^\mu
\label{12.373}
\end{equation}
(compare with \ref{10.810}). Here, the first difference is: 
\begin{eqnarray}
&&\Nb^\mu E^n\triangle\beta_\mu-\Nb_N^\mu E_N^n\triangle_N\beta_{\mu,N}=
\left.(\Nb^\mu E^n\beta_\mu-\Nb_N^\mu E_N^n\beta_{\mu,N})\right|_{{\cal K}}\nonumber\\
&&\hspace{30mm}-(\Nb^\mu E^n\beta^\prime_\mu\circ(f,w,\psi)-\Nb_N^\mu E_N^n\beta^\prime_\mu\circ(f_N,w_N,\psi_N))\nonumber\\
&&\label{12.374}
\end{eqnarray}
(compare with \ref{10.811}). Applying \ref{12.366} to the function 
$\Nb^\mu E^n\beta_\mu-\Nb_N^\mu E_N^n\beta_{\mu,N}$ on ${\cal K}$ we obtain:
\begin{equation}
\|\Nb^\mu E^n\beta_\mu-\Nb_N^\mu E_N^n\beta_{\mu,N}\|_{L^2(S_{\tau,\tau})}\leq 
k\tau^{1/2}\|T(\Nb^\mu E^n\beta_\mu-\Nb_N^\mu E_N^n\beta_{\mu,N})
\|_{L^2({\cal K}^{\tau})}
\label{12.375}
\end{equation}
Up to terms of order $n$ with vanishing $n$th order acoustical part $T(\Nb^\mu E^n\beta_\mu-\Nb_N^\mu E_N^n\beta_{\mu,N})$ is equal to:
\begin{equation}
\Nb^\mu T(E^n\beta_\mu-E_N^n\beta_{\mu,N})=\Nb^\mu L(E^n\beta_\mu-E_N^n\beta_{\mu,N})
+\Nb^\mu\Lb(E^n\beta_\mu-E_N^n\beta_{\mu,N})
\label{12.376}
\end{equation}
In regard to the first term on the right, this is equal to 
$$\Nb^\mu E^n L\beta_\mu-\Nb_N^\mu E_N^n L_N\beta_{\mu,N}-\ss_{\Nb}\rho\s^{(n-1)}\ctchi$$
up to terms which can be ignored. Since $\Nb^\mu L\beta_\mu=\lambdab E^\mu E\beta_\mu$ and similarly 
for the $N$th approximants, this is in turn equal to 
\begin{equation}
\lambdab\s^{(E;0,n)}\csxi+\sss\s^{(0,n)}\clab
+\frac{1}{2}\rho(\ss_N\s^{(n-1)}\ctchib-\ss_{\Nb}\s^{(n-1)}\ctchi)
\label{12.377}
\end{equation}
Since along ${\cal K}$ $a\sim\tau^3$, in regard to the 1st of \ref{12.377} we have:
\begin{eqnarray}
&&\|\lambdab\s^{(E;0,n)}\csxi\|^2_{L^2({\cal K}^{\tau})}\leq 
C\int_0^{\tau}\tau^{\prime-1}\|\sqrt{a}\s^{(E;0,n)}\csxi\|^2_{L^2(S_{\tau^\prime,\tau^\prime})}
d\tau^\prime\nonumber\\
&&\hspace{30mm}=C\int_0^{\tau}\tau^{\prime-1}\frac{\partial}{\partial\tau^\prime}
\left(\|\sqrt{a}\s^{(E;0,n)}\csxi\|^2_{L^2({\cal K}^{\tau^\prime})}\right)d\tau^\prime\nonumber\\
&&\hspace{10mm}=C\left\{\tau^{-1}\|\sqrt{a}\s^{(E;0,n)}\csxi\|^2_{L^2({\cal K}^{\tau})}
+\int_0^{\tau}\tau^{\prime-2}\|\sqrt{a}\s^{(E;0,n)}\csxi\|^2_{L^2({\cal K}^{\tau^\prime})}d\tau^\prime
\right\}\nonumber\\
&&\hspace{30mm}\leq C\s^{(E;0,n)}\cA(\tau_1)\tau^{2c_0-1}
\label{12.378}
\end{eqnarray}
for all $\tau\in [0,\tau_1]$, by \ref{9.a35} and \ref{9.299}. Therefore the contribution of this term 
through \ref{12.375} to \ref{12.372} is bounded by:
\begin{equation}
C\sqrt{\s^{(E;0,n)}\cA(\tau_1)}\left(\int_0^{\tau_1}\tau^{2c_0}d\tau\right)^{1/2}
\leq\frac{C\sqrt{\s^{(E;0,n)}\cA(\tau_1)}}{\sqrt{2c_0+1}}\cdot\tau_1^{c_0}\cdot\tau_1^{1/2}
\label{12.379}
\end{equation}
The contribution of the 2nd of \ref{12.377} through \ref{12.375} to \ref{12.372} is bounded by:
\begin{equation}
C\left(\int_0^{\tau_1}\tau\|\s^{(0,n}\clab\|^2_{L^2({\cal K}^{\tau})}d\tau\right)^{1/2}
\label{12.380}
\end{equation}
Since along ${\cal K}$ $\rho\sim r$, the contribution of the 3rd of \ref{12.377} through \ref{12.375} 
to \ref{12.372} is bounded by:
\begin{equation}
C\left(\int_0^{\tau_1}\tau\|r\s^{(n-1)}\ctchib\|^2_{L^2({\cal K}^{\tau})}d\tau\right)^{1/2}
+C\left(\int_0^{\tau_1}\tau^3\|\s^{(n-1)}\ctchi\|^2_{L^2({\cal K}^{\tau})}d\tau\right)^{1/2}
\label{12.381}
\end{equation}
In regard to the second term on the right in \ref{12.376} we write:
\begin{equation}
N^\mu\Lb(E^n\beta_\mu-E_N^n\beta_{\mu,N})+\ogamma\Nb^\mu\Lb(E^n\beta_\mu-E_N^n\beta_{\mu,N})
=\s^{(Y;0,n)}\cxi_{\Lb}
\label{12.382}
\end{equation}
Since along ${\cal K}$ we have $\ogamma\sim\tau$, the contribution of the right hand side of 
\ref{12.382} to $\|T(\Nb^\mu E^n\beta_\mu-\Nb_N^\mu E_N^n\beta_{\mu,N})\|^2_{L^2({\cal K}^{\tau})}$ 
is bounded by a constant multiple of:
\begin{eqnarray}
&&\int_0^{\tau}\tau^{\prime-2}\|\s^{(Y;0,n)}\cxi_{\Lb}\|^2_{L^2(S_{\tau^\prime,\tau^\prime})}d\tau^\prime\nonumber\\
&&=\int_0^{\tau}\tau^{\prime-2}\frac{\partial}{\partial\tau^\prime}
\left(\|\s^{(Y;0,n)}\cxi_{\Lb}\|^2_{L^2({\cal K}^{\tau^\prime})}\right)d\tau^\prime\nonumber\\
&&=\tau^{-2}\|\s^{(Y;0,n)}\cxi_{\Lb}\|^2_{L^2({\cal K}^{\tau})}
+2\int_0^{\tau}\tau^{\prime-3}\|\s^{(Y;0,n)}\cxi_{\Lb}\|^2_{L^2({\cal K}^{\tau^\prime})}d\tau^\prime
\nonumber\\
&&\leq C\s^{(Y;0,n)}\cA(\tau_1)\tau^{2c_0-2}
\label{12.383}
\end{eqnarray}
for all $\tau\in [0,\tau_1]$, by \ref{9.a39} and \ref{9.299}. The corresponding contribution to 
\ref{12.372} is then bounded by:
\begin{equation}
C\sqrt{\s^{(Y;0,n)}\cA(\tau_1)}\left(\int_0^{\tau_1}\tau^{2c_0-1}d\tau\right)^{1/2}
\leq\frac{C\sqrt{\s^{(Y;0,n)}\cA(\tau_1)}}{\sqrt{2c_0}}\cdot\tau_1^{c_0}
\label{12.384}
\end{equation}
In regard to the 1st term on the left in \ref{12.382}, this term is given, up to terms which can be 
ignored, by the conjugate of \ref{12.377}:
\begin{equation}
\lambda\s^{(E;0,n)}\csxi+\sss\s^{(0,n)}\cla+\frac{1}{2}\rhob(\ss_{\Nb}\s^{(n-1)}\ctchi
-\ss_N\s^{(n-1)}\ctchib)
\label{12.385}
\end{equation}
Since along ${\cal K}$ $\lambda/\ogamma=\lambdab$, the contribution of the 1st of \ref{12.385} 
coincides with that of the 1st of \ref{12.377}. Since along ${\cal K}$ $\ogamma\sim\tau$, the 
contribution of the 2nd of \ref{12.385} to 
$\|T(\Nb^\mu E^n\beta_\mu-\Nb_N^\mu E_N^n\beta_{\mu,N})\|^2_{L^2({\cal K}^{\tau})}$ 
is bounded by a constant multiple of:
\begin{eqnarray}
&&\int_0^{\tau}\tau^{\prime-2}\frac{\partial}{\partial\tau^\prime}
\left(\|\s^{(0,n)}\cla\|^2_{L^2({\cal K}^{\tau^\prime})}\right)d\tau^\prime\nonumber\\
&&=\tau^{-2}\|\s^{(0,n)}\cla\|^2_{L^2({\cal K}^{\tau})}
+2\int_0^{\tau}\tau^{\prime-3}\|\s^{(0,n)}\cla\|^2_{L^2({\cal K}^{\tau^\prime})}d\tau^\prime \nonumber\\
&&\label{12.386}
\end{eqnarray}
The corresponding contribution to \ref{12.372} is then bounded by:
\begin{equation}
\left\{\int_0^{\tau_1}\left(\tau^{-1}\|\s^{(0,n)}\cla\|^2_{L^2({\cal K}^{\tau})}
+2\tau\int_0^{\tau}\tau^{\prime-3}\|\s^{(0,n)}\cla\|^2_{L^2({\cal K}^{\tau^\prime})}d\tau^\prime\right)
d\tau\right\}^{1/2}
\label{12.387}
\end{equation}
Finally, since along ${\cal K}$ $\rhob/\ogamma=\rho$, the contribution of the 3rd of \ref{12.385} 
is bounded by \ref{12.381}. 

We turn to the second term on the right in \ref{12.374}. By \ref{10.798} - \ref{10.802} and 
\ref{10.803} with $l=n-1$ the principal part of the second term on the right in \ref{12.374} is:
\begin{eqnarray}
&&-\sh^{-n/2}\Nb^\mu\left\{\left(\frac{\partial\beta^\prime_\mu}{\partial t}\right)\circ(f,w,\psi)
\tau^2\Omega^n\chf+\left(\frac{\partial\beta^\prime_\mu}{\partial u^\prime}\right)\circ(f,w,\psi)
\tau \Omega^n\cv\right.\nonumber\\
&&\hspace{55mm}\left.+\left(\frac{\partial\beta^\prime_\mu}{\partial\vartheta^\prime}\right)
\circ(f,w,\psi)\tau^3\Omega^n\cga\right\}\nonumber\\
&&\label{12.388}
\end{eqnarray}
It follows that the principal part of the 2nd term on the right in 
\ref{12.374} is pointwise bounded by:
\begin{equation}
C|\tau^2\Omega^n\chf|+C|\tau\Omega^n\cv|+C|\tau^3\Omega^n\cga|
\label{12.389}
\end{equation}
Following the argument leading from \ref{12.364} to \ref{12.367} we conclude that the principal 
part of the 2nd tern on the right in \ref{12.374} is bounded in $L^2(S_{\tau,\tau})$ by:
\begin{equation}
C\tau^{5/2}\|\Omega^n T\chf\|_{L^2({\cal K}^{\tau})}
+C\tau^{3/2}\|\Omega^n T\cv\|_{L^2({\cal K}^{\tau})}
+C\tau^{7/2}\|\Omega^n T\cga\|_{L^2({\cal K}^{\tau})}
\label{12.390}
\end{equation}
The corresponding contribution to \ref{12.372} is then bounded by:
\begin{eqnarray}
&&C\left(\int_0^{\tau_1}\tau^5\|\Omega^n T\chf\|^2_{L^2({\cal K}^{\tau})}d\tau\right)^{1/2}
+C\left(\int_0^{\tau_1}\tau^3\|\Omega^n T\cv\|^2_{L^2({\cal K}^{\tau})}d\tau\right)^{1/2} \nonumber\\
&&+C\left(\int_0^{\tau_1}\tau^7\|\Omega^n T\cga\|^2_{L^2({\cal K}^{\tau})}d\tau\right)^{1/2} 
\label{12.391}
\end{eqnarray}

Finally, we consider the second of the differences \ref{12.373}. In view of \ref{12.45} and the 
boundary condition $E^\mu \triangle\beta_\mu=0$, the principal acoustical part of this difference is:
\begin{equation}
-\pi\epb\s^{(n-1)}\ctchib
\label{12.392}
\end{equation}
Since $|\epb|\leq Cr$, the contribution of this to \ref{12.372} is bounded by:
\begin{equation}
C\left(\int_0^{\tau_1}\|r\s^{(n-1)}\ctchib\|^2_{L^2(S_{\tau,\tau})}d\tau\right)^{1/2}
=C\|r\s^{(n-1)}\ctchib\|_{L^2({\cal K}^{\tau_1})}
\label{12.393}
\end{equation}

We also note that, since $r\sim\tau$, the contribution to 
$\|\s^{(0,n)}\clab\|_{L^2({\cal K}^{\tau_1})}$ of the first term on the right in \ref{12.371} 
is bounded by a constant multiple of \ref{12.386} with $\tau_1$ in the role of $\tau$. 
The above yield through \ref{12.371} the following lemma. 

\vspace{2.5mm}

\noindent{\bf Lemma 12.8} \ \ The boundary values on ${\cal K}$ of the next to the top order 
acoustical difference quantity $\s^{(0,n)}\clab$ satisfy, to leading terms, the inequality:
\begin{eqnarray*}
&&\|\s^{(0,n)}\clab\|_{L^2({\cal K}^{\tau_1})}\leq 
C\left(\tau_1^{-2}\|\s^{(0,n)}\cla\|^2_{L^2({\cal K}^{\tau_1})}+2\int_0^{\tau_1}\tau^{-3}
\|\s^{(0,n)}\cla\|^2_{L^2({\cal K}^{\tau})}d\tau\right)^{1/2}\\
&&\hspace{27mm}+C\|r\s^{(n-1)}\ctchib\|_{L^2({\cal K}^{\tau_1})}\\
&&\hspace{27mm}+\frac{C\sqrt{\s^{(Y;0,n)}\cA(\tau_1)}}{\sqrt{2c_0}}\cdot\tau_1^{c_0}\\
&&\hspace{27mm}+C\left(\int_0^{\tau_1}\tau^5\|\Omega^n T\chf\|^2_{L^2({\cal K}^{\tau})}d\tau
\right)^{1/2}\\
&&\hspace{27mm}+C\left(\int_0^{\tau_1}\tau^3\|\Omega^n T\cv\|^2_{L^2({\cal K}^{\tau})}d\tau\right)^{1/2}
\\
&&\hspace{27mm}+C\left(\int_0^{\tau_1}\tau^7\|\Omega^n T\cga\|^2_{L^2({\cal K}^{\tau})}d\tau\right)^{1/2}
\end{eqnarray*}

\vspace{2.5mm}

We now define, for arbitrary $\tau_1\in(0,\delta]$, the following quantities in regard to the transformation 
functions:
\begin{eqnarray}
&&\s^{(1,n)}P(\tau_1)=\sup_{\tau\in[0,\tau_1]}\left\{\tau^{-c_0+2}
\|\Omega^n T\chf\|_{L^2({\cal K}^{\tau})}\right\} \label{12.394}\\
&&\s^{(1,n)}V(\tau_1)=\sup_{\tau\in[0,\tau_1]}\left\{\tau^{-c_0+2}
\|\Omega^n T\cv\|_{L^2({\cal K}^{\tau})}\right\} \label{12.395}\\
&&\s^{(1,n)}\Ga(\tau_1)=\sup_{\tau\in[0,\tau_1]}\left\{\tau^{-c_0+2}
\|\Omega^n T\cga\|_{L^2({\cal K}^{\tau})}\right\} \label{12.396}
\end{eqnarray}
We also define, for arbitrary $\tau_1\in(0,\delta]$, the following quantities 
in regard to the boundary values on ${\cal K}$ of $\s^{(n-1)}\ctchi$, $\s^{(n-1)}\ctchib$, 
$\s^{(0,n)}\cla$, $\s^{(0,n)}\clab$:
\begin{eqnarray}
&&\s^{(n-1)}\oX(\tau_1)=\sup_{\tau\in[0,\tau_1]}\left\{\tau^{-c_0}
\|\s^{(n-1)}\ctchi\|_{L^2({\cal K}^{\tau})}\right\} \label{12.397}\\
&&\s^{(n-1)}\oXb(\tau_1)=\sup_{\tau\in[0,\tau_1]}\left\{\tau^{-c_0}
\|r\s^{(n-1)}\ctchib\|_{L^2({\cal K}^{\tau})}\right\} \label{12.398}\\
&&\s^{(0,n)}\oLa(\tau_1)=\sup_{\tau\in[0,\tau_1]}\left\{\tau^{-c_0-1}
\|\s^{(0,n)}\cla\|_{L^2({\cal K}^{\tau})}\right\} \label{12.399}\\
&&\s^{(0,n)}\oLab(\tau_1)=\sup_{\tau\in[0,\tau_1]}\left\{\tau^{-c_0}
\|\s^{(0,n)}\clab\|_{L^2({\cal K}^{\tau})}\right\} \label{12.400}
\end{eqnarray}
Then by Lemma 12.7 with $\tau\in (0,\tau_1]$ in the role of $\tau$, substituting 
\begin{equation}
\|\s^{(n-1)}\ctchi\|_{L^2({\cal K}^{\tau})}\leq \s^{(n-1)}\oX(\tau_1)\cdot \tau^{c_0}
\label{12.401}
\end{equation}
\begin{eqnarray}
&&\|\Omega^n T\chf\|_{L^2({\cal K}^{\tau})}\leq \s^{(1,n)}P(\tau_1)\cdot \tau^{c_0-2} \nonumber\\
&&\|\Omega^n T\cv\|_{L^2({\cal K}^{\tau})}\leq \s^{(1,n)}V(\tau_1)\cdot \tau^{c_0-2} \nonumber\\
&&\|\Omega^n T\cga\|_{L^2({\cal K}^{\tau})}\leq \s^{(1,n)}\Ga(\tau_1)\cdot \tau^{c_0-2} \label{12.402}
\end{eqnarray}
and noting that $\s^{(E;0,n)}\cA(\tau)\leq\s^{(E;0,n)}\cA(\tau_1)$, we obtain: 
\begin{eqnarray}
&&\|r\s^{(n-1)}\ctchib\|_{L^2({\cal K}^{\tau})}\leq k\s^{(n-1)}\oX(\tau_1)\cdot\tau^{c_0}\nonumber\\
&&\hspace{30mm}+\frac{C\sqrt{\s^{(E;0,n)}\cA(\tau_1)}}{\sqrt{2c_0}}\cdot\tau^{c_0}\nonumber\\
&&\hspace{30mm}+\frac{C\left(\s^{(1,n)}P(\tau_1)+\s^{(1,n)}V(\tau_1)\right)}{\sqrt{2c_0}}\cdot\tau^{c_0}
\nonumber\\
&&\hspace{30mm}+\frac{C\s^{(1,n)}\Ga(\tau_1)}{\sqrt{2c_0+2}}\cdot\tau^{c_0}\cdot\tau \label{12.403}
\end{eqnarray}
Multiplying both sides by $\tau^{-c_0}$ and taking the supremum over $\tau\in(0,\tau_1]$ then yields:
\begin{eqnarray}
&&\s^{(n-1)}\oXb(\tau_1)\leq k\s^{(n-1)}\oX(\tau_1)+\frac{C\sqrt{\s^{(E;0,n)}\cA(\tau_1)}}{\sqrt{2c_0}}
\nonumber\\
&&\hspace{22mm}+\frac{C\left(\s^{(1,n)}P(\tau_1)+\s^{(1,n)}V(\tau_1)\right)}{\sqrt{2c_0}}
+\frac{C\s^{(1,n)}\Ga(\tau_1)}{\sqrt{2c_0+2}}\cdot\tau_1 \nonumber\\
&&\label{12.404}
\end{eqnarray}
Also, by Lemma 12.8 with $\tau\in (0,\tau_1]$ in the role of $\tau$, substituting 
\begin{equation}
\|r\s^{(n-1)}\ctchib\|_{L^2({\cal K}^{\tau})}\leq\s^{(n-1)}\oXb(\tau_1)\cdot\tau^{c_0}
\label{12.405}
\end{equation}
and 
\begin{equation}
\|\s^{(0,n)}\cla\|_{L^2({\cal K}^{\tau})}\leq\s^{(0,n)}\oLa(\tau_1)\cdot\tau^{c_0+1}
\label{12.406}
\end{equation}
which implies 
\begin{equation}
\left(\tau^{-2}\|\s^{(0,n)}\cla\|^2_{L^2({\cal K}^{\tau})}
+2\int_0^{\tau}\tau^{\prime-3}\|\s^{(0,n)}\cla\|^2_{L^2({\cal K}^{\tau^\prime})}d\tau^\prime\right)^{1/2}
\leq C\s^{(0,n)}\oLa(\tau_1)\cdot\tau^{c_0}
\label{12.407}
\end{equation}
and noting that $\s^{(Y;0,n)}\cA(\tau)\leq\s^{(Y;0,n)}\cA(\tau_1)$, we obtain, substituting also 
\ref{12.402}:
\begin{eqnarray}
&&\|\s^{(0,n)}\clab\|_{L^2({\cal K}^{\tau})}\leq C\s^{(0,n)}\oLa(\tau_1)\cdot\tau^{c_0}
+C\s^{(n-1)}\oXb(\tau_1)\cdot\tau^{c_0}\nonumber\\
&&\hspace{27mm}+\frac{C\sqrt{\s^{(Y;0,n)}\cA(\tau_1)}}{\sqrt{2c_0}}\cdot\tau^{c_0}
+\frac{C\s^{(1,n)}V(\tau_1)}{\sqrt{2c_0}}\cdot\tau^{c_0}\nonumber\\
&&\hspace{27mm}+\frac{C\s^{(1,n)}P(\tau_1)}{\sqrt{2c_0+2}}\cdot\tau^{c_0+1}
+\frac{C\s^{(1,n)}\Ga(\tau_1)}{\sqrt{2c_0+4}}\cdot\tau^{c_0+2}\nonumber\\
&&\label{12.408}
\end{eqnarray}
Multiplying both sides by $\tau^{-c_0}$ and taking the supremum over $\tau\in(0,\tau_1]$ then yields:
\begin{eqnarray}
&&\s^{(0,n)}\oLab(\tau_1)\leq C\s^{(0,n)}\oLa(\tau_1)+C\s^{(n-1)}\oXb(\tau_1)\nonumber\\
&&\hspace{20mm}+\frac{C\sqrt{\s^{(Y;0,n)}\cA(\tau_1)}}{\sqrt{2c_0}}
+\frac{C\s^{(1,n)}V(\tau_1)}{\sqrt{2c_0}}\nonumber\\
&&\hspace{20mm}+\frac{C\s^{(1,n)}P(\tau_1)}{\sqrt{2c_0+2}}\cdot\tau_1
+\frac{C\s^{(1,n)}\Ga(\tau_1)}{\sqrt{2c_0+4}}\cdot\tau_1^2 
\label{12.409}
\end{eqnarray}

Note that since $r\sim\tau$ we have:
\begin{eqnarray*}
&&\|\s^{(n-1)}\ctchib\|^2_{L^2({\cal K}^{\tau_1})}=
\int_0^{\tau_1}\|\s^{(n-1)}\tchib\|^2_{L^2(S_{\tau,\tau})}d\tau\\
&&\hspace{26mm}\leq C\int_0^{\tau_1}\tau^{-2}\|r\s^{(n-1)}\ctchib\|^2_{L^2(S_{\tau,\tau})}d\tau\\
&&\hspace{20mm}=C\int_0^{\tau_1}\tau^{-2}\frac{\partial}{\partial\tau}\left(\|r\s^{(n-1)}\ctchib\|
^2_{L^2({\cal K}^{\tau})}\right)d\tau\\
&&\hspace{20mm}=C\left\{\tau_1^{-2}\|r^{(n-1)}\ctchib\|^2_{L^2({\cal K}^{\tau_1})}
+2\int_0^{\tau_1}\tau^{-3}\|r^{(n-1)}\ctchib\|^2_{L^2({\cal K}^{\tau})}\right\}\\
&&\hspace{32mm}\leq C\tau_1^{2c_0-2}\oXb^2(\tau_1)
\end{eqnarray*}
that is:
\begin{equation}
\|\s^{(n-1)}\ctchib\|_{L^2({\cal K}^{\tau_1})}\leq C\tau_1^{c_0-1}\oXb(\tau_1)
\label{12.410}
\end{equation}

In regard to the next to the top order acoustical difference quantities $\s^{(n-1)}\ctchi$, 
and $\s^{(0,n)}\cla$ in ${\cal N}$, we define for arbitrary 
$(\ub_1,u_1)\in R_{\delta,\delta}$, the quantities: 
\begin{equation}
\s^{(n-1)}X(\ub_1,u_1)=\sup_{(\ub,u)\in R_{\ub_1,u_1}}\left\{\ub^{-a_0}u^{-b_0}
\|\s^{(n-1)}\ctchi\|_{L^2(\Cb_{\ub}^u)}\right\} \label{12.411}
\end{equation}
\begin{equation}
\s^{(0,n)}\La(\ub_1,u_1)=\sup_{(\ub,u)\in R_{\ub_1,u_1}}
\left\{\ub^{-a_0-\frac{1}{2}}u^{-b_0-\frac{1}{2}}
\|\s^{(0,n)}\cla\|_{L^2(\Cb_{\ub}^u)}\right\} \label{12.412}
\end{equation}

Consider Lemma 12.4. Substituting 
$$\|\s^{(n-1)}\ctchi\|_{L^2(\Cb_{\ub}^{u_1})}\leq \s^{(n-1)}X(\ub_1,u_1)\ub^{a_0}u_1^{b_0}$$
we obtain, in regard to the 2nd term on the right hand side of the inequality,
\begin{equation}
u_1^{1/2}\left\{\int_0^{\ub_1}\|\s^{(n-1)}\ctchi\|^2_{L^2(\Cb_{\ub}^{u_1})}d\ub\right\}^{1/2}
\leq \frac{\s^{(n-1)}X(\ub_1,u_1)}{\sqrt{2a_0+1}}\ub_1^{a_0}u_1^{b_0}\cdot\ub_1^{1/2}u_1^{1/2}
\label{12.413}
\end{equation}
Also, substituting 
$$\|\s^{(0,n)}\cla\|_{L^2(\Cb_{\ub}^{u_1})}\leq \s^{(0,n)}\La(\ub_1,u_1)
\ub^{a_0+\frac{1}{2}}u_1^{b_0+\frac{1}{2}}$$
we obtain, in regard to the 3rd term on the right hand side of the inequality,
\begin{equation}
\left\{\int_0^{\ub_1}\ub^{-1}\|\s^{(0,n)}\cla\|^2_{L^2(\Cb_{\ub}^{u_1})}d\ub\right\}^{1/2}
\leq \frac{\s^{(0,n)}\La(\ub_1,u_1)}{\sqrt{2a_0+1}}\ub_1^{a_0}u_1^{b_0}\cdot\ub_1^{1/2}u_1^{1/2}
\label{12.414}
\end{equation}
In regard to the 1st term on the right hand side of the inequality we substitute \ref{12.410} with 
$\ub_1$ in the role of $\tau_1$. In view of the fact that since $\rho\sim\ub$ we have
$$\|\rho\s^{(n-1)}\ctchib\|_{L^2(C_{u_1}^{\ub_1})}\leq C\ub_1\|\s^{(n-1)}\ctchib\|
_{L^2(C_{u_1}^{\ub_1})}$$
we then conclude that:
\begin{eqnarray}
&&\|\rho\s^{(n-1)}\ctchib\|_{L^2(C_{u_1}^{\ub_1})}\leq C\s^{(n-1)}\oXb(\ub_1)\ub_1^{c_0}\nonumber\\
&&\hspace{29mm}+\frac{C\s^{(n-1)}X(\ub_1,u_1)}{\sqrt{2a_0+1}}\ub_1^{a_0}u_1^{b_0}\
\cdot\ub_1^{3/2}u_1^{1/2}\nonumber\\
&&\hspace{29mm}+\frac{C\s^{(0,n)}\La(\ub_1,u_1)}{\sqrt{2a_0+1}}\ub_1^{a_0}u_1^{b_0}\cdot
\ub_1^{3/2}u_1^{1/2}
\nonumber\\
&&\hspace{29mm}+\Cb\sqrt{\max\{\s^{(0,n)}\cB(\ub_1,u_1),\s^{(0,n)}\cBb(\ub_1,u_1)\}}\cdot
\ub_1^{a_0}u_1^{b_0}\cdot\ub_1 \nonumber\\
&&\label{12.415}
\end{eqnarray}

Consider Lemma 12.6. The 1st term on the right hand side of the inequality is bounded by:
$$k\s^{(0,n)}\oLab(\ub_1)\ub_1^{c_0}$$
The 2nd term is bounded by:
$$C\s^{(n-1)}X(\ub_1,u_1)\frac{\ub_1^{a_0+\frac{3}{2}}u_1^{b_0+\frac{1}{2}}}{\sqrt{2a_0+3}}$$
and the 5th term is bounded by:
$$C\s^{(0,n)}\La(\ub_1,u_1)\frac{\ub_1^{a_0+\frac{3}{2}}u_1^{b_0+\frac{1}{2}}}{\sqrt{2a_0+3}}$$
Since $\rho\sim\ub$, the 3rd and 4th terms are bounded by:
$$C\int_0^{\ub_1}\|\rho\s^{(n-1)}\ctchib\|_{L^2(C_u^u)}du
+C\int_{\ub_1}^{u_1}\|\rho\s^{(n-1)}\ctchib\|_{L^2(C_u^{\ub_1})}du$$
Substituting the estimate \ref{12.415} we find that this is in turn bounded by:
\begin{eqnarray*}
&&C\oXb(\ub_1)\ub_1^{c_0}u_1+\frac{C\s^{(n-1)}X(\ub_1,u_1)}{\sqrt{2a_0+1}}
\frac{\ub_1^{a_0+\frac{3}{2}}u_1^{b_0+\frac{3}{2}}}{\left(b_0+\frac{3}{2}\right)}\\
&&+\frac{C\s^{(0,n)}\La(\ub_1,u_1)}{\sqrt{2a_0+1}}\frac{\ub_1^{a_0+\frac{3}{2}}u_1^{b_0+\frac{3}{2}}}{(b_0+\frac{3}{2})}\\
&&+\Cb\sqrt{\max\{\s^{(0,n)}\cB(\ub_1,u_1),\s^{(0,n)}\cBb(\ub_1,u_1)\}}
\frac{\ub_1^{a_0+1}u_1^{b_0+1}}{(b_0+1)}
\end{eqnarray*}
We then conclude that:
\begin{eqnarray}
&&\|\s^{(0,n)}\clab\|_{L^2(C_{u_1}^{\ub_1})}\leq k\s^{(0,n)}\oLab(\ub_1)\ub_1^{c_0}\nonumber\\
&&\hspace{28mm}+\frac{C\s^{(n-1)}X(\ub_1,u_1)}{\sqrt{2a_0+1}}\ub_1^{a_0}u_1^{b_0}
\cdot\ub_1^{3/2}u_1^{1/2}\nonumber\\
&&\hspace{28mm}+\frac{C\s^{(0,n)}\La(\ub_1,u_1)}{\sqrt{2a_0+1}}\ub_1^{a_0}u_1^{b_0}\cdot
\ub_1^{3/2}u_1^{1/2}
\nonumber\\
&&\hspace{28mm}+\Cb\sqrt{\max\{\s^{(0,n)}\cB(\ub_1,u_1),\s^{(0,n)}\cBb(\ub_1,u_1)\}}\cdot 
\ub_1^{a_0}u_1^{b_0}\cdot u_1\nonumber\\
&&\label{12.416}
\end{eqnarray}

We turn to Lemma 12.2. We first consider the inequality satisfied by 
$\|\s^{(n-1)}\ctchi\|_{L^2(\Cb_{\ub_1}^{u_1})}$. By definitions \ref{12.411}, \ref{12.412}, for all 
$(\ub,u)\in R_{\ub_1,u_1}$ we have:
\begin{eqnarray}
&&\|\s^{(n-1)}\ctchi\|_{L^2(\Cb_{\ub}^u)}\leq \s^{(n-1)}X(\ub_1,u_1)\ub^{a_0}u^{b_0} \label{12.417}\\
&&\|\s^{(0,n)}\cla\|_{L^2(\Cb_{\ub}^u)}\leq \s^{(0,n)}\La(\ub_1,u_1)
\ub^{a_0+\frac{1}{2}}u^{b_0+\frac{1}{2}} \label{12.418}
\end{eqnarray}
In regard to the 1st term on the right hand side of the inequality in question, substituting 
\ref{12.417} we obtain:
\begin{eqnarray*}
&&u_1^{-4}\|\s^{(n-1)}\ctchi\|^2_{L^2(\Cb_{\ub}^{u_1})}+
4\int_{\ub}^{u_1}u^{-5}\|\s^{(n-1)}\ctchi\|^2_{L^2(\Cb_{\ub}^u)}du\\
&&\hspace{30mm}\leq C\left(\s^{(n-1)}X(\ub_1,u_1)\right)^2\ub^{2a_0}u_1^{2b_0-4}
\end{eqnarray*}
hence the 1st term is bounded by:
\begin{equation}
\frac{C\s^{(n-1)}X(\ub_1,u_1)}{a_0+2}\cdot\ub_1^{a_0+2}u_1^{b_0-2}
\label{12.419}
\end{equation}
In regard to the 2nd term, substituting \ref{12.418} we obtain:
\begin{eqnarray*}
&&u_1^{-4}\|\s^{(0,n)}\cla\|^2_{L^2(\Cb_{\ub}^{u_1})}+
4\int_{\ub}^{u_1}u^{-5}\|\s^{(0,n)}\cla\|^2_{L^2(\Cb_{\ub}^u)}du\\
&&\hspace{30mm}\leq C\left(\s^{(0,n)}\La(\ub_1,u_1)\right)^2\ub^{2a_0+1}u_1^{2b_0-3}
\end{eqnarray*}
hence the 2nd term is bounded by:
\begin{equation}
\frac{C\s^{(0,n)}\La(\ub_1,u_1)}{a_0+\frac{5}{2}}\cdot\ub_1^{a_0+\frac{5}{2}}u_1^{b_0-\frac{3}{2}}
\label{12.420}
\end{equation}
By the estimate \ref{12.415} the 3rd term on the right hand side of the same inequality is bounded by:
\begin{eqnarray}
&&C\s^{(n-1)}\oXb(\ub_1)\ub_1^{c_0}\cdot\ub_1^{1/2}u_1^{1/2}\nonumber\\
&&+\frac{C\s^{(n-1)}X(\ub_1,u_1)}{\sqrt{2a_0+1}}\ub_1^{a_0}u_1^{b_0}\cdot\ub_1^2 u_1
+\frac{C\s^{(0,n)}\La(\ub_1,u_1)}{\sqrt{2a_0+1}}\ub_1^{a_0}u_1^{b_0}\cdot\ub_1^2 u_1
\nonumber\\
&&+\Cb\sqrt{\max\{\s^{(0,n)}\cB(\ub_1,u_1),\s^{(0,n)}\cBb(\ub_1,u_1)\}}\cdot
\ub_1^{a_0}u_1^{b_0}\cdot\ub_1^{3/2}u_1^{1/2}\label{12.421}
\end{eqnarray}
By the estimate \ref{12.416} the 4th term on the right hand side of the same inequality is bounded by:
\begin{eqnarray}
&&C\s^{(0,n)}\oLab(\ub_1)\ub_1^{c_0}\cdot\ub_1^{1/2}u_1^{1/2}\nonumber\\
&&+\frac{C\s^{(n-1)}X(\ub_1,u_1)}{\sqrt{2a_0+1}}\ub_1^{a_0}u_1^{b_0}\cdot\ub_1^2 u_1
+\frac{C\s^{(0,n)}\La(\ub_1,u_1)}{\sqrt{2a_0+1}}\ub_1^{a_0}u_1^{b_0}\cdot\ub_1^2 u_1
\nonumber\\
&&+\Cb\sqrt{\max\{\s^{(0,n)}\cB(\ub_1,u_1),\cBb(\ub_1,u_1)\}}\cdot\ub_1^{a_0}u_1^{b_0}\cdot 
\ub_1^{1/2}u_1^{3/2}\label{12.422}
\end{eqnarray}
Replacing in the inequality which results by substituting the above bounds $(\ub_1,u_1)$ by any 
$(\ub,u)\in R_{\ub_1,u_1}$ and noting that by the definitions {12.411}, \ref{12.412} we have:
\begin{eqnarray*}
&&\s^{(n-1)}X(\ub,u)\leq \s^{(n-1)}X(\ub_1,u_1)\\
&&\s^{(0,n)}\La(\ub,u)\leq \s^{(0,n)}\La(\ub_1,u_1)
\end{eqnarray*}
we obtain:
\begin{eqnarray}
&&\|\s^{(n-1)}\ctchi\|_{L^2(\Cb_{\ub}^u)}\leq \ub^{a_0}u^{b_0}\left\{ 
C\left(\s^{(n-1)}\oXb(\ub_1)+\s^{(0,n)}\oLab(\ub_1)\right)\left(\frac{\ub}{u}\right)^{b_0}
\ub^{1/2}u^{1/2}
\right.\nonumber\\
&&\hspace{32mm}+C\left(\frac{1}{a_0+2}\left(\frac{\ub}{u}\right)^2+\frac{1}{\sqrt{2a_0+1}}\ub^2 u\right)
\s^{(n-1)}X(\ub_1,u_1)\nonumber\\
&&\hspace{32mm}+C\left(\frac{1}{a_0+\frac{5}{2}}\left(\frac{\ub}{u}\right)^{3/2}+\frac{1}{\sqrt{2a_0+1}}
\ub u\right)\ub\s^{(0,n)}\La(\ub_1,u_1)\nonumber\\
&&\hspace{32mm}+C\frac{\sqrt{\s^{(E;0,n)}\cBb(\ub_1,u_1)}}{(a_0+2)}\left(\frac{\ub}{u}\right)^2
\nonumber\\
&&\hspace{32mm}\left. +C\sqrt{\max\{\s^{(0,n)}\cB(\ub_1,u_1),\s^{(0,n)}\cBb(\ub_1,u_1)\}}\ub^{1/2}
\right\}\nonumber\\
&&\label{12.423}
\end{eqnarray}
Since $\ub/u\leq 1$, multiplying both sides by $\ub^{-a_0}u^{-b_0}$ and taking the supremum 
over $(\ub,u)\in R_{\ub_1,u_1}$ we obtain:
\begin{eqnarray}
&&\s^{(n-1)}X(\ub_1,u_1)\leq C\left(\s^{(n-1)}\oXb(\ub_1)+\s^{(0,n)}\oLab(\ub_1)\right)
\ub_1^{1/2}u_1^{1/2}\nonumber\\
&&\hspace{25mm}+C\left(\frac{1}{a_0+2}+\frac{1}{\sqrt{2a_0+1}}\ub_1^2 u_1\right)
\s^{(n-1)}X(\ub_1,u_1)\nonumber\\
&&\hspace{25mm}+C\left(\frac{1}{a_0+\frac{5}{2}}+\frac{1}{\sqrt{2a_0+1}}\ub_1 u_1\right)\ub_1
\s^{(0,n)}\La(\ub_1,u_1)\nonumber\\
&&\hspace{25mm}+C\frac{\sqrt{\s^{(E;0,n)}\cBb(\ub_1,u_1)}}{(a_0+2)}\nonumber\\
&&\hspace{20mm}+C\sqrt{\max\{\s^{(0,n)}\cB(\ub_1,u_1),\s^{(0,n)}\cBb(\ub_1,u_1)\}}\ub_1^{1/2}
\label{12.424}
\end{eqnarray}
Choosing $a_0$ large enough so that, in regard to the coefficient of 
$\s^{(n-1)}X(\ub_1,u_1)$, 
\begin{equation}
C\left(\frac{1}{a_0+2}+\frac{1}{\sqrt{2a_0+1}}\delta^3\right)\leq \frac{1}{2}
\label{12.425}
\end{equation}
this implies:
\begin{eqnarray}
&&\s^{(n-1)}X(\ub_1,u_1)\leq C\left(\s^{(n-1)}\oXb(\ub_1)+\s^{(0,n)}\oLab(\ub_1)\right)
\ub_1^{1/2}u_1^{1/2}\nonumber\\
&&\hspace{25mm}+C\left(\frac{1}{a_0+\frac{5}{2}}+\frac{1}{\sqrt{2a_0+1}}\ub_1 u_1\right)\ub_1
\s^{(0,n)}\La(\ub_1,u_1)\nonumber\\
&&\hspace{25mm}+C\frac{\sqrt{\s^{(E;0,n)}\cBb(\ub_1,u_1)}}{(a_0+2)}\nonumber\\
&&\hspace{20mm}+C\sqrt{\max\{\s^{(0,n)}\cB(\ub_1,u_1),\s^{(0,n)}\cBb(\ub_1,u_1)\}}\ub_1^{1/2}
\label{12.426}
\end{eqnarray}
(for new constants $C$). 

We turn to the inequality of Lemma 12.2 satisfied by $\|\s^{(n-1)}\ctchi\|_{L^2({\cal K}^{\tau})}$. 
Now the right hand side of this inequality coincides with the right hand side of the inequality 
satisfied by $\|\s^{(n-1)}\ctchi\|_{L^2(\Cb_{\ub_1}^{u_1})}$ if we set $\ub_1=u_1=\tau$. Therefore 
$\|\s^{(n-1)}\ctchi\|_{L^2({\cal K}^{\tau})}$ is bounded by the right hand side of \ref{12.423} in 
which we have set $\ub=\ub_1=\tau$, $u=u_1=\tau$:
\begin{eqnarray}
&&\|\s^{(n-1)}\ctchi\|_{L^2({\cal K}^{\tau})}\leq \tau^{c_0}\left\{ 
C\left(\s^{(n-1)}\oXb(\tau)+\s^{(0,n)}\oLab(\tau)\right)\tau
\right.\nonumber\\
&&\hspace{32mm}+C\left(\frac{1}{a_0+2}+\frac{1}{\sqrt{2a_0+1}}\tau^3\right)
\s^{(n-1)}X(\tau,\tau)\nonumber\\
&&\hspace{32mm}+C\left(\frac{1}{a_0+\frac{5}{2}}+\frac{1}{\sqrt{2a_0+1}}
\tau^2\right)\tau\s^{(0,n)}\La(\tau,\tau)\nonumber\\
&&\hspace{32mm}+C\frac{\sqrt{\s^{(E;0,n)}\cBb(\tau,\tau)}}{(a_0+2)}
\nonumber\\
&&\hspace{32mm}\left. +C\sqrt{\max\{\s^{(0,n)}\cB(\tau,\tau),\s^{(0,n)}\cBb(\tau,\tau)\}}\tau^{1/2}
\right\}\nonumber\\
&&\label{12.427}
\end{eqnarray}
Let $\tau_1\in(0,\delta]$. For any $\tau\in[0,\tau_1]$ the above inequality holds if we replace 
$\s^{(n-1)}\oXb(\tau)$, $\s^{(0,n)}\oLab(\tau)$, $\s^{(n-1)}X(\tau,\tau)$, $\s^{(0,n)}\La(\tau,\tau)$ , 
$\s^{(E;0,n)}\cBb(\tau,\tau)$, $\s^{(0,n)}\cB(\tau,\tau)$, $\s^{(0,n)}\cBb(\tau,\tau)$ on the right 
by $\s^{(n-1)}\oXb(\tau_1)$, $\s^{(0,n)}\oLab(\tau_1)$, $\s^{(n-1)}X(\tau_1,\tau_1)$, 
$\s^{(0,n)}\La(\tau_1,\tau_1)$, $\s^{(E;0,n)}\cBb(\tau_1,\tau_1)$, $\s^{(0,n)}\cB(\tau_1,\tau_1)$, 
$\s^{(0,n)}\cBb(\tau_1,\tau_1)$, respectively, these being non-decreasing functions of $\tau$. 
Multiplying then both sides by $\tau^{-c_0}$ and taking the supremum over $\tau\in[0,\tau_1]$ 
we obtain (see definition \ref{12.397}):
\begin{eqnarray}
&&\s^{(n-1)}\oX(\tau_1)\leq C\left(\s^{(n-1)}\oXb(\tau_1)+\s^{(0,n)}\oLab(\tau_1)\right)\tau_1
\nonumber\\
&&\hspace{15mm}+C\left(\frac{1}{a_0+2}+\frac{1}{\sqrt{2a_0+1}}\tau_1^3\right)
\s^{(n-1)}X(\tau_1,\tau_1)\nonumber\\
&&\hspace{15mm}+C\left(\frac{1}{a_0+\frac{5}{2}}+\frac{1}{\sqrt{2a_0+1}}
\tau_1^2\right)\tau_1\s^{(0,n)}\La(\tau_1,\tau_1)\nonumber\\
&&\hspace{15mm}+C\frac{\sqrt{\s^{(E;0,n)}\cBb(\tau_1,\tau_1)}}{(a_0+2)}
\nonumber\\
&&\hspace{15mm}+C\sqrt{\max\{\s^{(0,n)}\cB(\tau_1,\tau_1),\s^{(0,n)}\cBb(\tau,_1\tau_1)\}}\tau_1^{1/2}\label{12.428}
\end{eqnarray}

We turn to Lemma 12.3. We first consider the inequality satisfied by 
$\|\s^{(0,n)}\cla\|_{L^2(\Cb_{\ub_1}^{u_1})}$. The first term on the right hand side of this inequality 
is bounded by (Schwartz inequality):
\begin{equation}
C\left\{\int_{\ub_1}^{u_1}\ub_1\|\s^{(0,n)}\clab\|^2_{L^2(C_u^{\ub_1})}du\right\}^{1/2}
\label{12.429}
\end{equation}
By the estimate \ref{12.416} this is bounded by:
\begin{eqnarray}
&&C\s^{(0,n)}\oLab(\ub_1)\ub_1^{c_0}\cdot \ub_1^{1/2}u_1^{1/2}\nonumber\\
&&+C\frac{\s^{(n-1)}X(\ub_1,u_1)}{\sqrt{2a_0+1}\sqrt{2b_0+2}}\ub_1^{a_0}u_1^{b_0}\cdot\ub_1^2 u_1
+C\frac{\s^{(0,n)}\La(\ub_1,u_1)}{\sqrt{2a_0+1}\sqrt{2b_0+2}}\ub_1^{a_0}u_1^{b_0}\cdot
\ub_1^2 u_1\nonumber\\
&&+\Cb\frac{\sqrt{\max\{\s^{(0,n)}\cB(\ub_1,u_1),\s^{(0,n)}\cBb(\ub_1,u_1)\}}}{\sqrt{2b_0+3}}
\cdot\ub_1^{a_0}u_1^{b_0}
\cdot\ub_1^{1/2} u_1^{3/2}\nonumber\\
&&\label{12.430}
\end{eqnarray}
Similarly, the second term on the right hand side of the inequality in question is bounded by:
\begin{equation}
C\left\{\int_{\ub_1}^{u_1}\ub_1^{1/2}\|\rho\s^{(n-1)}\ctchib\|_{L^2(C_u^{\ub_1})}du\right\}^{1/2}
\label{12.431}
\end{equation}
which by the estimate \ref{12.415} is bounded by:
\begin{eqnarray}
&&C\s^{(n-1)}\oXb(\ub_1)\ub_1^{c_0}\cdot\ub_1^{1/2}u_1^{1/2}\nonumber\\
&&+C\frac{\s^{(n-1)}X(\ub_1,u_1)}{\sqrt{2a_0+1}\sqrt{2b_0+2}}\ub_1^{a_0}u_1^{b_0}\cdot\ub_1^2 u_1 
+ C\frac{\s^{(0,n)}\La(\ub_1,u_1)}{\sqrt{2a_0+1}\sqrt{2b_0+2}}\ub_1^{a_0}u_1^{b_0}\cdot 
\ub_1^2 u_1\nonumber\\
&&+\Cb\frac{\sqrt{\max\{\s^{(0,n)}\cB(\ub_1,u_1),\s^{(0,n)}\cBb(\ub_1,u_1)\}}}{\sqrt{2b_0+1}}
\cdot \ub_1^{a_0}u_1^{b_0}\cdot \ub_1^{3/2}u_1^{1/2}\nonumber\\
&&\label{12.432}
\end{eqnarray}
In view of the definitions \ref{12.411}, \ref{12.412} the 3rd and 4th terms in the same inequality 
are bounded by:
\begin{equation}
C\frac{\s^{(0,n)}\La(\ub_1,u_1)}{a_0+\frac{5}{2}}\ub_1^{a_0}u_1^{b_0}\cdot\ub_1^{5/2}u_1^{1/2}
+C\frac{\s^{(n-1)}X(\ub_1,u_1)}{a_0+2}\ub_1^{a_0}u_1^{b_0}\cdot\ub_1^2
\label{12.433}
\end{equation}
Replacing in the inequality which results by substituting the above bounds $(\ub_1,u_1)$ by any 
$(\ub,u)\in R_{\ub_1,u_1}$ we obtain:
\begin{eqnarray}
&&\|\s^{(0,n)}\cla\|_{L^2(\Cb_{\ub}^u)}\leq \ub^{a_0+\frac{1}{2}}u^{b_0+\frac{1}{2}}\left\{
C\left(\s^{(0,n)}\oLab(\ub_1)+\s^{(n-1)}\oXb(\ub_1)\right)\left(\frac{\ub}{u}\right)^{b_0}
\right.\nonumber\\
&&\hspace{14mm}+C\left(\frac{1}{a_0+\frac{5}{2}}\left(\frac{\ub}{u}\right)^{1/2}
+\frac{1}{\sqrt{2a_0+1}\sqrt{2b_0+2}}\right)\ub^{3/2}u^{1/2}\s^{(0,n)}\La(\ub_1,u_1)\nonumber\\
&&\hspace{14mm}+C\left(\frac{1}{a_0+2}\left(\frac{\ub}{u}\right)^{1/2}+\frac{\ub^{1/2}u^{1/2}}{\sqrt{2a_0+1}\sqrt{2b_0+2}}\right)\ub\s^{(n-1)}X(\ub_1,u_1)\nonumber\\
&&\hspace{28mm}+C\frac{\sqrt{\s^{(Y;0,n)}\cBb(\ub_1,u_1)}}{(a_0+2)}\left(\frac{\ub}{u}\right)^{3/2}
\nonumber\\
&&\hspace{28mm}\left.+C\sqrt{\max\{\s^{(0,n)}\cB(\ub_1,u_1),\s^{(0,n)}\cBb(\ub_1,u_1)\}}\cdot u^{1/2}\right\}
\nonumber\\
&&\label{12.434}
\end{eqnarray}
Since $\ub/u\leq 1$, multiplying both sides by $\ub^{-a_0-\frac{1}{2}}u^{-b_0-\frac{1}{2}}$ 
and taking the supremum over $(\ub,u)\in R_{\ub_1,u_1}$ we obtain:
\begin{eqnarray}
&&\s^{(0,n)}\La(\ub_1,u_1)\leq C\s^{(0,n)}\oLab(\ub_1)+C\s^{(n-1)}\oXb(\ub_1)\nonumber\\
&&\hspace{17mm}+C\left(\frac{1}{a_0+\frac{5}{2}}
+\frac{1}{\sqrt{2a_0+1}\sqrt{2b_0+2}}\right)\ub_1^{3/2}u_1^{1/2}\s^{(0,n)}\La(\ub_1,u_1)\nonumber\\
&&\hspace{17mm}+C\left(\frac{1}{a_0+2}+\frac{\ub_1^{1/2}u_1^{1/2}}{\sqrt{2a_0+1}\sqrt{2b_0+2}}\right)
\ub_1\s^{(n-1)}X(\ub_1,u_1)\nonumber\\
&&\hspace{25mm}+C\frac{\sqrt{\s^{(Y;0,n)}\cBb(\ub_1,u_1)}}{(a_0+2)}\nonumber\\
&&\hspace{25mm}+C\sqrt{\max\{\s^{(0,n)}\cB(\ub_1,u_1),\s^{(0,n)}\cBb(\ub_1,u_1)\}}\cdot u_1^{1/2}
\nonumber\\
&&\label{12.435}
\end{eqnarray}
Subjecting $\delta$ to an appropriate smallness condition so that the coefficient of 
$\s^{(0,n)}\La(\ub_1,u_1)$ on the right does not exceed $1/2$ this implies: 
\begin{eqnarray}
&&\s^{(0,n)}\La(\ub_1,u_1)\leq C\s^{(0,n)}\oLab(\ub_1)+C\s^{(n-1)}\oXb(\ub_1)\nonumber\\
&&\hspace{17mm}+C\left(\frac{1}{a_0+2}+\frac{\ub_1^{1/2}u_1^{1/2}}{\sqrt{2a_0+1}\sqrt{2b_0+2}}\right)
\ub_1\s^{(n-1)}X(\ub_1,u_1)\nonumber\\
&&\hspace{25mm}+C\frac{\sqrt{\s^{(Y;0,n)}\cBb(\ub_1,u_1)}}{(a_0+2)}\nonumber\\
&&\hspace{25mm}+C\sqrt{\max\{\s^{(0,n)}\cB(\ub_1,u_1),\s^{(0,n)}\cBb(\ub_1,u_1)\}}\cdot u_1^{1/2}
\nonumber\\
&&\label{12.436}
\end{eqnarray}
(for new constants $C$). 

Substituting this in the inequality \ref{12.426} we obtain:
\begin{eqnarray}
&&\s^{(n-1)}X(\ub_1,u_1)\leq C\left(\s^{(n-1)}\oXb(\ub_1)+\s^{(0,n)}\oLab(\ub_1)\right)
\ub_1^{1/2}u_1^{1/2}\nonumber\\
&&+C\left(\frac{1}{a_0+\frac{5}{2}}+\frac{1}{\sqrt{2a_0+1}}\ub_1 u_1\right)
\left(\frac{1}{a_0+2}+\frac{\ub_1^{1/2}u_1^{1/2}}{\sqrt{2a_0+1}\sqrt{2b_0+2}}\right)
\ub_1^2\s^{(n-1)}X(\ub_1,u_1)\nonumber\\
&&\hspace{25mm}+C\frac{\sqrt{\s^{(E;0,n)}\cBb(\ub_1,u_1)}}{(a_0+2)}\nonumber\\
&&\hspace{20mm}+C\sqrt{\max\{\s^{(0,n)}\cB(\ub_1,u_1),\s^{(0,n)}\cBb(\ub_1,u_1)\}}\ub_1^{1/2}
\label{12.437}
\end{eqnarray}
Subjecting $\delta$ to an appropriate smallness condition so that the coefficient of 
$\s^{(n-1)}X(\ub_1,u_1)$ on the right does not exceed $1/2$ this implies: 
\begin{eqnarray}
&&\s^{(n-1)}X(\ub_1,u_1)\leq C\left(\s^{(n-1)}\oXb(\ub_1)+\s^{(0,n)}\oLab(\ub_1)\right)
\ub_1^{1/2}u_1^{1/2}\nonumber\\
&&\hspace{25mm}+C\frac{\sqrt{\s^{(E;0,n)}\cBb(\ub_1,u_1)}}{(a_0+2)}\nonumber\\
&&\hspace{20mm}+C\sqrt{\max\{\s^{(0,n)}\cB(\ub_1,u_1),\s^{(0,n)}\cBb(\ub_1,u_1)\}}\ub_1^{1/2}
\nonumber\\
&&\label{12.438}
\end{eqnarray}
(for new constants $C$). Substituting \ref{12.438} in \ref{12.436} then yields:
\begin{eqnarray}
&&\s^{(0,n)}\La(\ub_1,u_1)\leq C\s^{(0,n)}\oLab(\ub_1)+C\s^{(n-1)}\oXb(\ub_1)\nonumber\\
&&\hspace{25mm}+C\frac{\sqrt{\s^{(Y;0,n)}\cBb(\ub_1,u_1)}}{(a_0+2)}\nonumber\\
&&\hspace{25mm}+C\sqrt{\max\{\s^{(0,n)}\cB(\ub_1,u_1),\s^{(0,n)}\cBb(\ub_1,u_1)\}}\cdot u_1^{1/2}
\nonumber\\
&&\label{12.439}
\end{eqnarray}
(for new constants $C$). 

We turn to the inequality of Lemma 12.3 satisfied by $\|\s^{(0,n)}\cla\|_{L^2({\cal K}^{\tau})}$. 
The first term on the right hand side of this inequality is bounded by (Schwartz inequality):
\begin{equation}
C\left\{\int_0^{\tau}u\|\s^{(0,n)}\clab\|^2_{L^2(C_u^u)}du\right\}^{1/2} 
\label{12.440}
\end{equation}
According to the estimate \ref{12.416} with $\ub_1=u_1=u$ we have:
\begin{eqnarray}
&&\|\s^{(0,n)}\clab\|_{L^2(C_u^u)}\leq k\s^{(0,n)}\oLab(u)u^{c_0}\nonumber\\
&&\hspace{20mm}+\frac{C}{\sqrt{2a_0+1}}\left(\s^{(n-1)}X(u,u)+\s^{(0,n)}\La(u,u)\right)u^{c_0+2}
\nonumber\\
&&\hspace{20mm}+\Cb\sqrt{\max\{\s^{(0,n)}\cB(u,u),\s^{(0,n)}\cBb(u,u)\}}\cdot u^{c_0+1}\nonumber\\
&&\label{12.441}
\end{eqnarray}
The contribution of the 1st term on the right in \ref{12.441} to \ref{12.440} is bounded by:
\begin{equation}
C\s^{(0,n)}\oLab(\tau)\left(\int_0^{\tau}u^{2c_0+1}du\right)^{1/2}=\frac{C}{\sqrt{2c_0+2}}\s^{(0,n)}\oLab(\tau)\tau^{c_0+1}
\label{12.442}
\end{equation}
Since $\s^{(n-1)}X(u,u)$ and $\s^{(0,n)}\La(u,u)$ are non-decreasing functions of $u$, the 
contribution of the 2nd term on the right in \ref{12.441} to \ref{12.440} is bounded by:
\begin{eqnarray*}
&&\frac{C}{\sqrt{2a_0+1}}\left(\s^{(n-1)}X(\tau,\tau)+\s^{(0,n)}\La(\tau,\tau)\right)
\left(\int_0^{\tau} u^{2c_0+5}\right)^{1/2}\\
&&\hspace{20mm}=\frac{C}{\sqrt{2a_0+1}}\left(\s^{(n-1)}X(\tau,\tau)+\s^{(0,n)}\La(\tau,\tau)\right)
\frac{\tau^{c_0+3}}{\sqrt{2c_0+6}}
\end{eqnarray*}
Substituting the estimates \ref{12.438}, \ref{12.439} with $\ub_1=u_1=\tau$, this is in turn bounded by:
\begin{eqnarray}
&&\frac{1}{\sqrt{2a_0+1}\sqrt{2c_0+6}}
\left[C\left(\s^{(n-1)}\oXb(\tau)+\s^{(0,n)}\oLab(\tau)\right)\right.\nonumber\\
&&\hspace{30mm}\left.+C\sqrt{\max\{\s^{(0,n)}\cB(\tau,\tau),\s^{(0,n)}\cBb(\tau,\tau)\}}\right]
\tau^{c_0+3}\nonumber\\
&&\label{12.443}
\end{eqnarray}
Finally, the contribution of the 3rd term on the right in \ref{12.441} to \ref{12.440} is bounded by:
\begin{equation}
\frac{\Cb}{\sqrt{2c_0+4}}\sqrt{\max\{\s^{(0,n)}\cB(\tau,\tau),\s^{(0,n)}\cBb(\tau,\tau)\}}\cdot\tau^{c_0+2}
\label{12.444}
\end{equation}

The second term on the right hand side of the inequality of Lemma 12.3 satisfied by 
$\|\s^{(0,n)}\cla\|_{L^2({\cal K}^{\tau})}$ is bounded by (Schwartz inequality):
\begin{equation}
C\left\{\int_0^{\tau} u\|\rho\s^{(n-1)}\ctchib\|^2_{L^2(C_u^u)}du\right\}^{1/2}
\label{12.445}
\end{equation}
According to the estimate \ref{12.415} with $\ub_1=u_1=u$ we have:
\begin{eqnarray}
&&\|\rho\s^{(n-1)}\ctchib\|_{L^2(C_u^u)}\leq C\s^{(0,n)}\oXb(u)u^{c_0}\nonumber\\
&&\hspace{20mm}+\frac{C}{\sqrt{2a_0+1}}\left(\s^{(n-1)}X(u,u)+\s^{(0,n)}\La(u,u)\right)u^{c_0+2}
\nonumber\\
&&\hspace{20mm}+\Cb\sqrt{\max\{\s^{(0,n)}\cB(u,u),\s^{(0,n)}\cBb(u,u)\}}\cdot u^{c_0+1}\nonumber\\
&&\label{12.446}
\end{eqnarray}
The contribution of the 1st term on the right in \ref{12.446} to \ref{12.445} is bounded by:
\begin{equation}
C\s^{(n-1)}\oXb(\tau)\left(\int_0^{\tau}u^{2c_0+1}du\right)^{1/2}=\frac{C}{\sqrt{2c_0+2}}\s^{(n-1)}\oXb(\tau)\tau^{c_0+1}
\label{12.447}
\end{equation}
The 2nd and 3rd terms on the right in \ref{12.446} are similar to the corresponding terms in \ref{12.441}. 

The 3rd and 4th terms on the right hand side of the inequality of Lemma 12.3 satisfied by 
$\|\s^{(0,n)}\cla\|_{L^2({\cal K}^{\tau})}$ coincide with the corresponding terms on the right hand 
side of the inequality of the same lemma satisfied by $\|\s^{(0,n)}\cla\|_{L^2(\Cb_{\ub_1}^{u_1})}$ 
if we set $\ub_1=u_1=\tau$. Thus by \ref{12.433} these terms are bounded by:
\begin{eqnarray}
&&\frac{C}{a_0+2}\left(\s^{(0,n)}\La(\tau,\tau)+\s^{(n-1)}X(\tau,\tau)\right)\tau^{c_0+2}\nonumber\\
&&\leq \frac{1}{a_0+2}\left[C\left(\s^{(n-1)}\oXb(\tau)+\s^{(0,n)}\oLab(\tau)\right)\right.\nonumber\\
&&\hspace{20mm}\left.+C\sqrt{\max\{\s^{(0,n)}\cB(\tau,\tau),\s^{(0,n)}\cBb(\tau,\tau)\}}\right]
\tau^{c_0+2}\label{12.448}
\end{eqnarray}
where we have substituted from \ref{12.438}, \ref{12.439}. Combining the above results we conclude 
that $\|\s^{(0,n)}\cla\|_{L^2({\cal K}^{\tau})}$ satisfies the inequality:
\begin{eqnarray}
&&\|\s^{(0,n)}\cla\|_{L^2({\cal K}^{\tau})}\leq \tau^{c_0+1}\left\{
K(\tau)\left(\s^{(0,n)}\oLab(\tau)+\s^{(n-1)}\oXb(\tau)\right)\right.\nonumber\\
&&\hspace{26mm}+C\frac{\sqrt{\s^{(Y;0,n)}\cBb(\tau,\tau)}}{(a_0+2)}\nonumber\\
&&\hspace{26mm}\left.+C\sqrt{\max\{\s^{(0,n)}\cB(\tau,\tau),\s^{(0,n)}\cBb(\tau,\tau)\}}\tau^{1/2}
\right\}\nonumber\\
&&\label{12.449}
\end{eqnarray}
(for a new constant $C$). Here:
\begin{equation}
K(\tau)=C\left(\frac{1}{\sqrt{2c_0+2}}+\frac{\tau}{a_0+2}+\frac{\tau^2}{\sqrt{2a_0+1}\sqrt{2c_0+6}}
\right)
\label{12.450}
\end{equation}
Given now $\tau_1\in(0,\delta]$ and letting $\tau$ range in $[0,\tau_1]$, we may replace 
$K(\tau)$, $\oLab(\tau)$, $\oXb(\tau)$, $\s^{(Y;0,n)}\cBb(\tau,\tau)$, $\s^{(0,n)}\cB(\tau,\tau)$, 
$\s^{(0,n)}\cBb(\tau,\tau)$ on the right in \ref{12.449} by $K(\tau_1)$, $\oLab(\tau_1)$, 
$\oXb(\tau_1)$, $\s^{(Y;0,n)}\cBb(\tau_1,\tau_1)$, $\s^{(0,n)}\cB(\tau_1,\tau_1)$, 
$\s^{(0,n)}\cBb(\tau_1,\tau_1)$, respectively, these being all non-decreasing functions of $\tau$. 
Multiplying then by $\tau^{-c_0-1}$ and taking the supremum over $\tau\in[0,\tau_1]$ we obtain 
the inequality (see definition \ref{12.399}):
\begin{eqnarray}
&&\s^{(0,n)}\oLa(\tau_1)\leq K(\tau_1)\left(\s^{(0,n)}\oLab(\tau_1)+\s^{(n-1)}\oXb(\tau_1)\right)\nonumber\\
&&\hspace{15mm}+C\frac{\sqrt{\s^{(Y;0,n)}\cBb(\tau_1,\tau_1)}}{(a_0+2)}\nonumber\\
&&\hspace{15mm}+C\sqrt{\max\{\s^{(0,n)}\cB(\tau_1,\tau_1),\s^{(0,n)}\cBb(\tau,\tau)\}}\tau_1^{1/2}
\label{12.451}
\end{eqnarray}
Also, substituting the estimates \ref{12.438}, \ref{12.439} in \ref{12.428} we obtain the inequality:
\begin{eqnarray}
&&\s^{(0,n)}\oX(\tau_1)\leq C\tau_1\left(\s^{(n-1)}\oXb(\tau_1)+\s^{(0,n)}\oLab(\tau_1)\right)
\nonumber\\
&&\hspace{15mm}+C\frac{\sqrt{\s^{(E;0,n)}\cBb(\tau_1,\tau_1)}}{(a_0+2)}\nonumber\\
&&\hspace{15mm}+C\sqrt{\max\{\s^{(0,n)}\cB(\tau_1,\tau_1),\s^{(0,n)}\cBb(\tau_1,\tau_1)\}}\tau_1^{1/2}
\label{12.452}
\end{eqnarray}
(for new constants $C$). 

We have now deduced four inequalities connecting the four boundary next to top order acoustical 
difference quantities $\s^{(n-1)}\oX$, $\s^{(n-1)}\oXb$, $\s^{(0,n)}\oLa$, $\s^{(0,n)}\oLab$, 
the inequalities \ref{12.404}, \ref{12.409}, \ref{12.451}, \ref{12.452}. Substituting \ref{12.452} 
in \ref{12.404} we obtain:
\begin{eqnarray*}
&&\s^{(n-1)}\oXb(\tau_1)\leq C\tau_1\s^{(n-1)}\oXb(\tau_1)+C\tau_1\s^{(0,n)}\oLab(\tau_1)\\
&&\hspace{15mm}+\frac{C\left(\s^{(1,n)}P(\tau_1)+\s^{(1,n)}V(\tau_1)\right)}{\sqrt{2c_0}}
+\frac{C\s^{(1,n)}\Ga(\tau_1)}{\sqrt{2c_0+2}}\cdot\tau_1\\
&&\hspace{15mm}+\frac{C\sqrt{\s^{(E;0,n)}\cA(\tau_1)}}{\sqrt{2c_0}}
+\frac{C\sqrt{\s^{(E;0,n)}\cBb(\tau_1,\tau_1)}}{(a_0+2)}\\
&&\hspace{15mm}+C\sqrt{\max\{\s^{(0,n)}\cB(\tau_1,\tau_1),\s^{(0,n)}\cBb(\tau_1,\tau_1)\}}
\cdot\tau_1^{1/2}
\end{eqnarray*}
Subjecting $\delta$ to the smallness condition, in regard to the coefficient of 
$\s^{(n-1)}\oXb(\tau_1)$ on the right, 
\begin{equation}
C\delta\leq \frac{1}{2}
\label{12.453}
\end{equation}
this implies: 
\begin{eqnarray}
&&\s^{(n-1)}\oXb(\tau_1)\leq C\tau_1\s^{(0,n)}\oLab(\tau_1) \label{12.454}\\
&&\hspace{15mm}+\frac{C\left(\s^{(1,n)}P(\tau_1)+\s^{(1,n)}V(\tau_1)\right)}{\sqrt{2c_0}}
+\frac{C\s^{(1,n)}\Ga(\tau_1)}{\sqrt{2c_0+2}}\cdot\tau_1\nonumber\\
&&\hspace{15mm}+\frac{C\sqrt{\s^{(E;0,n)}\cA(\tau_1)}}{\sqrt{2c_0}}
+\frac{C\sqrt{\s^{(E;0,n)}\cBb(\tau_1,\tau_1)}}{(a_0+2)}\nonumber\\
&&\hspace{15mm}+C\sqrt{\max\{\s^{(0,n)}\cB(\tau_1,\tau_1),\s^{(0,n)}\cBb(\tau_1,\tau_1)\}}
\cdot\tau_1^{1/2}\nonumber
\end{eqnarray}
Substituting \ref{12.451} and \ref{12.454} in \ref{12.409} we obtain:
\begin{eqnarray*}
&&\s^{(0,n)}\oLab(\tau_1)\leq \left[CK(\tau_1)+C\tau_1\right]\s^{(0,n)}\oLab\\
&&\hspace{15mm}+\frac{C\left(\s^{(1,n)}P(\tau_1)+\s^{(1,n)}V(\tau_1)\right)}{\sqrt{2c_0}}
+\frac{C\s^{(1,n)}\Ga(\tau_1)}{\sqrt{2c_0+2}}\cdot\tau_1\\
&&\hspace{15mm}+\frac{C\sqrt{\s^{(Y;0,n)}\cA(\tau_1)}}{\sqrt{2c_0}}
+\frac{C\sqrt{\s^{(E;0,n)}\cA(\tau_1)}}{\sqrt{2c_0}}\\
&&\hspace{15mm}+\frac{C\sqrt{\s^{(Y;0,n)}\cBb(\tau_1,\tau_1)}}{(a_0+2)}
+\frac{C\sqrt{\s^{(E;0,n)}\cBb(\tau_1,\tau_1)}}{(a_0+2)}\\
&&\hspace{15mm}+C\sqrt{\max\{\s^{(0,n)}\cB(\tau_1,\tau_1),\s^{(0,n)}\cBb(\tau_1,\tau_1)\}}
\cdot\tau_1^{1/2}
\end{eqnarray*}
The coefficient of $\s^{(0,n)}\oLab$ on the right does not exceed $CK(\delta)+C\delta$. 
Recalling \ref{12.450} this is not greater than 
\begin{equation}
C\left(\frac{1}{\sqrt{2c_0+2}}+\delta\right)
\label{12.455}
\end{equation}
(for a new constant $C$). Choosing then $c_0$ large enough so that:
\begin{equation}
\frac{C}{\sqrt{2c_0+2}}\leq\frac{1}{4}
\label{12.456}
\end{equation}
and subjecting $\delta$ to the smallness condition:
\begin{equation}
C\delta\leq \frac{1}{4}
\label{12.457}
\end{equation}
\ref{12.455} is not greater than 1/2. The above inequality then implies:
\begin{eqnarray}
&&\s^{(0,n)}\oLab(\tau_1)\leq \frac{C\left(\s^{(1,n)}P(\tau_1)+\s^{(1,n)}V(\tau_1)\right)}{\sqrt{2c_0}}
+\frac{C\s^{(1,n)}\Ga(\tau_1)}{\sqrt{2c_0+2}}\cdot\tau_1\nonumber\\
&&\hspace{18mm}+\frac{C\sqrt{\s^{(Y;0,n)}\cA(\tau_1)}}{\sqrt{2c_0}}
+\frac{C\sqrt{\s^{(E;0,n)}\cA(\tau_1)}}{\sqrt{2c_0}}\nonumber\\
&&\hspace{18mm}+\frac{C\sqrt{\s^{(Y;0,n)}\cBb(\tau_1,\tau_1)}}{(a_0+2)}
+\frac{C\sqrt{\s^{(E;0,n)}\cBb(\tau_1,\tau_1)}}{(a_0+2)}\nonumber\\
&&\hspace{18mm}+C\sqrt{\max\{\s^{(0,n)}\cB(\tau_1,\tau_1),\s^{(0,n)}\cBb(\tau_1,\tau_1)\}}
\cdot\tau_1^{1/2}\nonumber\\
&&\label{12.458}
\end{eqnarray}
(for new constants $C$). Substituting finally \ref{12.458} in \ref{12.454} yields:
\begin{eqnarray}
&&\s^{(n-1)}\oXb(\tau_1)\leq \frac{C\left(\s^{(1,n)}P(\tau_1)+\s^{(1,n)}V(\tau_1)\right)}{\sqrt{2c_0}}
+\frac{C\s^{(1,n)}\Ga(\tau_1)}{\sqrt{2c_0+2}}\cdot\tau_1\nonumber\\
&&\hspace{20mm}+\frac{C\sqrt{\s^{(E;0,n)}\cA(\tau_1)}}{\sqrt{2c_0}}
+\frac{C\sqrt{\s^{(Y;0,n)}\cA(\tau_1)}}{\sqrt{2c_0}}\cdot\tau_1\nonumber\\
&&\hspace{20mm}+\frac{C\sqrt{\s^{(E;0,n)}\cBb(\tau_1,\tau_1)}}{(a_0+2)}
+\frac{C\sqrt{\s^{(Y;0,n)}\cBb(\tau_1,\tau_1)}}{(a_0+2)}\cdot\tau_1\nonumber\\
&&\hspace{20mm}+C\sqrt{\max\{\s^{(0,n)}\cB(\tau_1,\tau_1),\s^{(0,n)}\cBb(\tau_1,\tau_1)\}}
\cdot\tau_1^{1/2}\nonumber\\
&&\label{12.459}
\end{eqnarray}
(for new constants $C$). We summarize the above results in regard to the next to 
top order boundary acoustical difference quantities in the following proposition. 

\vspace{2.5mm}

\noindent{\bf Proposition 12.1} \ \ Choosing the exponents $a_0$ and $c_0$ suitably large and subjecting 
$\delta$ to a suitable smallness condition, the next to top order boundary acoustical difference 
quantities $\s^{(n-1)}\oXb$, $\s^{(0,n)}\oLab$, $\s^{(n-1)}\oX$, $\s^{(0,n)}\oLa$ satisfy 
to principal terms the estimates:
\begin{eqnarray*}
&&\s^{(n-1)}\oXb(\tau_1)\leq \frac{C\left(\s^{(1,n)}P(\tau_1)+\s^{(1,n)}V(\tau_1)\right)}{\sqrt{2c_0}}
+\frac{C\s^{(1,n)}\Ga(\tau_1)}{\sqrt{2c_0+2}}\cdot\tau_1\\
&&\hspace{20mm}+\frac{C\left(\sqrt{\s^{(E;0,n)}\cA(\tau_1)}+\tau_1\sqrt{\s^{(Y;0,n)}\cA(\tau_1)}\right)}
{\sqrt{2c_0}}\\
&&\hspace{20mm}+\frac{C\sqrt{\s^{(E;0,n)}\cBb(\tau_1,\tau_1)}}{a_0+2}\\
&&\hspace{20mm}+C\sqrt{\max\{\s^{(0,n)}\cB(\tau_1,\tau_1),\s^{(0,n)}\cBb(\tau_1,\tau_1)\}}
\cdot\tau_1^{1/2}
\end{eqnarray*}
\begin{eqnarray*}
&&\s^{(0,n)}\oLab(\tau_1)\leq \frac{C\left(\s^{(1,n)}P(\tau_1)+\s^{(1,n)}V(\tau_1)\right)}{\sqrt{2c_0}}
+\frac{C\s^{(1,n)}\Ga(\tau_1)}{\sqrt{2c_0+2}}\cdot\tau_1\\
&&\hspace{18mm}+\frac{C\left(\sqrt{\s^{(Y;0,n)}\cA(\tau_1)}+\sqrt{\s^{(E;0,n)}\cA(\tau_1)}
\right)}{\sqrt{2c_0}}\\
&&\hspace{18mm}+\frac{C\left(\sqrt{\s^{(Y;0,n)}\cBb(\tau_1,\tau_1)}
+\sqrt{\s^{(E;0,n)}\cBb(\tau_1,\tau_1)}\right)}{a_0+2}\\
&&\hspace{18mm}+C\sqrt{\max\{\s^{(0,n)}\cB(\tau_1,\tau_1),\s^{(0,n)}\cBb(\tau_1,\tau_1)\}}
\cdot\tau_1^{1/2}
\end{eqnarray*}
\begin{eqnarray*}
&&\s^{(n-1)}\oX(\tau_1)\leq 
\frac{C\left(\s^{(1,n)}P(\tau_1)+\s^{(1,n)}V(\tau_1)\right)}{\sqrt{2c_0}}\cdot\tau_1
+\frac{C\s^{(1,n)}\Ga(\tau_1)}{\sqrt{2c_0+2}}\cdot\tau_1^2\\
&&\hspace{20mm}+\frac{C\left(\sqrt{\s^{(Y;0,n)}\cA(\tau_1)}+\sqrt{\s^{(E;0,n)}\cA(\tau_1)}
\right)}{\sqrt{2c_0}}\cdot\tau_1\\
&&\hspace{20mm}+\frac{C\sqrt{\s^{(E;0,n)}\cBb(\tau_1,\tau_1)}}{a_0+2}\\
&&\hspace{20mm}+C\sqrt{\max\{\s^{(0,n)}\cB(\tau_1,\tau_1),\s^{(0,n)}\cBb(\tau_1,\tau_1)\}}
\cdot\tau_1^{1/2}
\end{eqnarray*}
\begin{eqnarray*}
&&\s^{(0,n)}\oLa(\tau_1)\leq K(\tau_1)
\left[\frac{C\left(\s^{(1,n)}P(\tau_1)+\s^{(1,n)}V(\tau_1)\right)}{\sqrt{2c_0}}
+\frac{C\s^{(1,n)}\Ga(\tau_1)}{\sqrt{2c_0+2}}\cdot\tau_1\right.\\
&&\hspace{37mm}+\frac{C\left(\sqrt{\s^{(Y;0,n)}\cA(\tau_1)}+\sqrt{\s^{(E;0,n)}\cA(\tau_1)}
\right)}{\sqrt{2c_0}}\\
&&\hspace{65mm}\left.+\frac{C\sqrt{\s^{(E;0,n)}\cBb(\tau_1,\tau_1)}}{a_0+2}\right]\\
&&\hspace{20mm}+\frac{C\sqrt{\s^{(Y;0,n)}\cBb(\tau_1,\tau_1)}}{a_0+2}\\
&&\hspace{20mm}+C\sqrt{\max\{\s^{(0,n)}\cB(\tau_1,\tau_1),\s^{(0,n)}\cBb(\tau_1,\tau_1)\}}
\cdot\tau_1^{1/2}
\end{eqnarray*}

\vspace{2.5mm}

To complete the next to top order acoustical difference quantities we define, in addition to 
the quantities $\s^{(n-1)}X(\ub_1,u_1)$, $\s^{(0,n)}\La(\ub_1,u_1)$ defined by 
\ref{12.411}, \ref{12.412}, the quantities:
\begin{equation}
\s^{(n-1)}\Xb(\ub_1,u_1)=\sup_{(\ub,u)\in R_{\ub_1,u_1}}\left\{\ub^{-a_0-1}u^{-b_0+1}
\|\rho\s^{(n-1)}\ctchib\|_{L^2(C_u^{\ub})}\right\}
\label{12.460}
\end{equation}
\begin{equation}
\s^{(0,n)}\Lab(\ub_1,u_1)=\sup_{(\ub,u)\in R_{\ub_1,u_1}}\left\{\ub^{-a_0}u^{-b_0}
\|\s^{(0,n)}\clab\|_{L^2(C_u^{\ub})}\right\}
\label{12.461}
\end{equation}
Replacing $(\ub_1,u_1)$ in inequality \ref{12.415} by any $(\ub,u)\in R_{\ub_1,u_1}$ but keeping 
the quantities $\oXb(\ub_1)$, $\s^{(n-1)}X(\ub_1,u_1)$, $\s^{(0,n)}\La(\ub_1,u_1)$, 
$\s^{(0,n)}\cB(\ub_1,u_1)$, $\s^{(0,n)}\cBb(\ub_1,u_1)$ on the right as they stand, then multiplying 
by $\ub^{-a_0-1}u^{-b_0+1}$ and taking the supremum over $(\ub,u)\in R_{\ub_1,u_1}$, we obtain:
\begin{eqnarray}
&&\s^{(n-1)}\Xb(\ub_1,u_1)\leq C\s^{(n-1)}\oXb(\ub_1)\nonumber\\
&&\hspace{25mm}+\frac{C\left(\s^{(n-1)}X(\ub_1,u_1)+\s^{(0,n)}\La(\ub_1,u_1)\right)}{\sqrt{2a_0+1}}
\cdot\ub_1^{1/2}u_1^{3/2}\nonumber\\
&&\hspace{25mm}+\Cb\sqrt{\max\{\s^{(0,n)}\cB(\ub_1,u_1),\s^{(0,n)}\cBb(\ub_1,u_1)\}}\cdot u_1
\nonumber\\
&&\label{12.462}
\end{eqnarray}
Replacing $(\ub_1,u_1)$ in inequality \ref{12.416} by any $(\ub,u)\in R_{\ub_1,u_1}$ but keeping 
the quantities $\oLab(\ub_1)$, $\s^{(n-1)}X(\ub_1,u_1)$, $\s^{(0,n)}\La(\ub_1,u_1)$, 
$\s^{(0,n)}\cB(\ub_1,u_1)$, $\s^{(0,n)}\cBb(\ub_1,u_1)$ on the right as they stand, then multiplying 
by $\ub^{-a_0}u^{-b_0}$ and taking the supremum over $(\ub,u)\in R_{\ub_1,u_1}$, we obtain:
\begin{eqnarray}
&&\s^{(0,n)}\Lab(\ub_1,u_1)\leq k\s^{(0,n)}\oLab(\ub_1)\nonumber\\
&&\hspace{23mm}+\frac{C\left(\s^{(n-1)}X(\ub_1,u_1)+\s^{(0,n)}\La(\ub_1,u_1)\right)}{\sqrt{2a_0+1}}
\cdot\ub_1^{3/2}u_1^{1/2}\nonumber\\
&&\hspace{23mm}+\Cb\sqrt{\max\{\s^{(0,n)}\cB(\ub_1,u_1),\s^{(0,n)}\cBb(\ub_1,u_1)\}}\cdot u_1
\nonumber\\
&&\label{12.463}
\end{eqnarray}

In view of the estimates \ref{12.438}, \ref{12.439} for $\s^{(n-1)}X(\ub_1,u_1)$, 
$\s^{(0,n)}\La(\ub_1,u_1)$ and the estimates of Proposition 12.1 for $\s^{(n-1)}\oXb(\ub_1)$, 
$\s^{(0,n)}\oLab(\ub_1)$, we arrive at the following proposition in regard to the next to 
top order interior acoustical difference quantities. 

\vspace{2.5mm}

\noindent{\bf Proposition 12.2} \ \ Choosing the exponents $a_0$ and $c_0$ suitably large and subjecting 
$\delta$ to a suitable smallness condition, the next to top order interior acoustical difference 
quantities $\s^{(n-1)}X$, $\s^{(0,n)}\La$, $\s^{(n-1)}\Xb$, $\s^{(0,n)}\Lab$ satisfy 
to principal terms the estimates:
\begin{eqnarray*}
&&\s^{(n-1)}X(\ub_1,u_1)\leq C\ub_1^{1/2}u_1^{1/2}\left\{\frac{
\left(\s^{(1,n)}P(\ub_1)+\s^{(1,n)}V(\ub_1)\right)}{\sqrt{2c_0}}
+\frac{\s^{(1,n)}\Ga(\ub_1)}{\sqrt{2c_0+2}}\ub_1\right.\\
&&\hspace{55mm}\left.+\frac{\sqrt{\s^{(Y;0,n)}\cA(\ub_1)}+\sqrt{\s^{(E;0,n)}\cA(\ub_1)}}{\sqrt{2c_0}}
\right\}\\
&&\hspace{27mm}+\frac{C\sqrt{\s^{(E;0,n)}\cBb(\ub_1,u_1)}}{a_0+2}\\
&&\hspace{27mm}+C\sqrt{\max\{\s^{(0,n)}\cB(\ub_1,u_1),\s^{(0,n)}\cBb(\ub_1,u_1)\}}\cdot\ub_1^{1/2}
\end{eqnarray*}
\begin{eqnarray*}
&&\s^{(0,n)}\La(\ub_1,u_1)\leq \frac{C\left(\s^{(1,n)}P(\ub_1)+\s^{(1,n)}V(\ub_1)\right)}{\sqrt{2c_0}}
+\frac{C\s^{(1,n)}\Ga(\ub_1)}{\sqrt{2c_0+2}}\ub_1\\
&&\hspace{25mm}+\frac{C\left(\sqrt{\s^{(Y;0,n)}\cA(\ub_1)}+\sqrt{\s^{(E;0,n)}\cA(\ub_1)}\right)}{\sqrt{2c_0}}\\
&&\hspace{25mm}+\frac{C\left(\sqrt{\s^{(Y;0,n)}\cBb(\ub_1,u_1)}
+\sqrt{\s^{(E;0,n)}\cBb(\ub_1,\ub_1)}\right)}{a_0+2}\\
&&\hspace{25mm}+C\sqrt{\max\{\s^{(0,n)}\cB(\ub_1,u_1),\s^{(0,n)}\cBb(\ub_1,u_1)\}}\cdot u_1^{1/2}
\end{eqnarray*}
\begin{eqnarray*}
&&\s^{(n-1)}\Xb(\ub_1,u_1)\leq \frac{C\left(\s^{(1,n)}P(\ub_1)+\s^{(1,n)}V(\ub_1)\right)}{\sqrt{2c_0}}
+\frac{C\s^{(1,n)}\Ga(\ub_1)}{\sqrt{2c_0+2}}\cdot\ub_1\\
&&\hspace{27mm}+\frac{C\sqrt{\s^{(E;0,n)}\cA(\ub_1)}}{\sqrt{2c_0}}
+\frac{C\sqrt{\s^{(Y;0,n)}\cA(\ub_1)}}{\sqrt{2c_0}}\cdot\ub_1^{1/2}u_1^{1/2}\\
&&\hspace{27mm}+\frac{C\sqrt{\s^{(E;0,n)}\cBb(\ub_1,\ub_1)}}{a_0+2}\\
&&\hspace{27mm}+\Cb\sqrt{\max\{\s^{(0,n)}\cB(\ub_1,u_1),\s^{(0,n)}\cBb(\ub_1,u_1)\}}\cdot u_1^{1/2}
\end{eqnarray*}
\begin{eqnarray*}
&&\s^{(0,n)}\Lab(\ub_1,u_1)\leq \frac{C\left(\s^{(1,n)}P(\ub_1)+\s^{(1,n)}V(\ub_1)\right)}{\sqrt{2c_0}}
+\frac{C\s^{(1,n)}\Ga(\ub_1)}{\sqrt{2c_0+2}}\cdot\ub_1\\
&&\hspace{25mm}+\frac{C\left(\sqrt{\s^{(Y;0,n)}\cA(\ub_1)}+\sqrt{\s^{(E;0,n)}\cA(\ub_1)}\right)}
{\sqrt{2c_0}}\\
&&\hspace{25mm}+\frac{C\left(\sqrt{\s^{(Y;0,n)}\cBb(\ub_1,\ub_1)}+
\sqrt{\s^{(E;0,n)}\cBb(\ub_1,\ub_1)}\right)}{a_0+2}\\
&&\hspace{25mm}+\Cb\sqrt{\max\{\s^{(0,n)}\cB(\ub_1,u_1),\s^{(0,n)}\cBb(\ub_1,u_1)\}}\cdot u_1^{1/2}
\end{eqnarray*}

\vspace{5mm}

\section{Estimates for $(T\Omega^n\chf, T\Omega^n\cv, T\Omega^n\cga)$}

Let us recall that the transformation function $f$ satisfies the first of \ref{4.119}, which, 
since $\rho=c^{-1}\lambdab$, $\rhob=c^{-1}\lambda$, reads:
\begin{equation}
Tf=\rho+\rhob \ : \ \mbox{along ${\cal K}$}
\label{12.464}
\end{equation}
In view of the fact that $\rho=c^{-1}\lambdab$, $\rhob=c^{-1}\lambda$ and 
the boundary condition \ref{4.6} (or \ref{10.752}) this takes the form:
\begin{equation}
Tf=\frac{(1+r)}{c}\lambdab \ : \ \mbox{along ${\cal K}$}
\label{12.465}
\end{equation}
By \ref{9.51} and \ref{9.14} the $N$th approximate transformation function $f_N$ satisfies:
\begin{equation}
Tf_N=\rho_N+\rhob_N \ : \ \mbox{along ${\cal K}$}
\label{12.466}
\end{equation}
In view of \ref{9.101} and \ref{9.102} this can be written in the form:
\begin{equation}
Tf_N=\frac{(1+r_N)}{c_N}\lambdab_N+\vep_N
\label{12.467}
\end{equation}
where: 
\begin{equation}
\vep_N=\frac{\hat{\nu}_N}{c_N}=O(\tau^{N+1})
\label{12.468}
\end{equation}
by \ref{9.106}. Let us presently denote by $\cf$ the difference:
\begin{equation}
\cf=f-f_N
\label{12.469}
\end{equation}
Note that this is does not agree with \ref{4.222}. In the present section $\cf$ stands for 
the difference \ref{12.469}. Subtracting \ref{12.467} from \ref{12.465} we obtain:
\begin{equation}
T\cf=\frac{(1+r)}{c}\lambdab-\frac{(1+r_N)}{c_N}\lambdab_N-\vep_N \ : \ \mbox{along ${\cal K}$}
\label{12.470}
\end{equation}

We apply $\Omega^n$ to \ref{12.470} to obtain:
\begin{equation}
\Omega^n T\cf=\Omega^n\left[\frac{(1+r)}{c}\lambdab\right]
-\Omega^n\left[\frac{(1+r_N)}{c_N}\lambdab_N\right]-\Omega^n\vep_N
\label{12.471}
\end{equation}
and from \ref{12.468} we have:
\begin{equation}
\Omega^n\vep_N=O(\tau^{N+1})
\label{12.472}
\end{equation}
The principal part of the right hand side of \ref{12.471} is contained in: 
\begin{equation}
\frac{(1+r)}{c}\Omega^n(\lambdab-\lambdab_N)
+\lambdab\Omega^n\left[\frac{(1+r)}{c}-\frac{(1+r_N)}{c_N}\right]
\label{12.473}
\end{equation}
In fact, the principal part of the first of \ref{12.473} is:
\begin{equation} 
\frac{(1+r)}{c}\sh^{n/2}\s^{(0,n)}\clab
\label{12.474}
\end{equation}
while the principal part of the second of \ref{12.473} is contained in:
\begin{eqnarray}
&&\frac{\lambdab}{c}\sh^{n/2}(E^n r-E_N^n r_N)\nonumber\\
&&-\lambdab\frac{(1+r)}{c}\sh^{n/2}\left(E^{n-1}(c^{-1}Ec)-E_N^{n-1}(c_N^{-1}E_N c_N)\right)
\label{12.475}
\end{eqnarray}
Given now $\tau_1\in (0,\delta]$, \ref{12.474} is bounded in $L^2({\cal K}^{\tau})$ for 
$\tau\in[0,\tau_1]$ by:
\begin{equation}
C\s^{(0,n)}\oLab(\tau_1)\tau^{c_0}
\label{12.476}
\end{equation}
On the other hand, since $\lambdab\sim\tau$ along ${\cal K}$, the 1st of \ref{12.475} is bounded in 
$L^2({\cal K}^{\tau})$ by:
\begin{equation}
C\tau\|E^n r-E_N^n r_N\|_{L^2({\cal K}^{\tau})}
\label{12.477}
\end{equation}
which is $C\tau$ times \ref{12.372} with $\tau$ in the role of $\tau_1$. By the estimates following 
\ref{12.372} and leading to Lemma 12.8 we have, to leading terms:
\begin{eqnarray}
&&\|E^n r-E_N^n r_N\|_{L^2({\cal K}^{\tau_1})}\leq  C\|r\s^{(n-1)}\ctchib\|_{L^2({\cal K}^{\tau_1})}
+\frac{C\s^{(Y;0,n)}\cA(\tau_1)}{\sqrt{2c_0}}\cdot\tau_1^{c_0}\nonumber\\
&&\hspace{37mm}
+C\left(\int_0^{\tau_1}\tau^5\|\Omega^n T\chf\|^2_{L^2({\cal K}^{\tau})}d\tau\right)^{1/2}\nonumber\\
&&\hspace{37mm}
+C\left(\int_0^{\tau_1}\tau^3\|\Omega^n T\cv\|^2_{L^2({\cal K}^{\tau})}d\tau\right)^{1/2}\nonumber\\
&&\hspace{37mm}
+C\left(\int_0^{\tau_1}\tau^7\|\Omega^n T\cga\|^2_{L^2({\cal K}^{\tau})}d\tau\right)^{1/2}
\label{12.478}
\end{eqnarray}
Substituting the definitions \ref{12.398} and \ref{12.394}-\ref{12.396} we then obtain, to leading 
terms:
\begin{eqnarray}
&&\|E^n r-E_N^n r_N\|_{L^2({\cal K}^{\tau_1})}\leq \tau_1^{c_0}\left\{\s^{(n-1)}\oXb(\tau_1)
+\frac{C\s^{(Y;0,n)}\cA(\tau_1)}{\sqrt{2c_0}}\right.\nonumber\\
&&\hspace{60mm}+\frac{C\s^{(1,n)}V(\tau_1)}{\sqrt{2c_0}}\nonumber\\
&&\hspace{45mm}\left.+\frac{C\s^{(1,n)}P(\tau_1)}{\sqrt{2c_0+2}}\cdot\tau_1
+\frac{C\s^{(1,n)}\Ga(\tau_1)}{\sqrt{2c_0+4}}\cdot\tau_1^2\right\}\nonumber\\
&&\label{12.479}
\end{eqnarray}
In regard to the 2nd of \ref{12.475}, by \ref{3.a23} and \ref{10.148} we have:
$$[c^{-1}Ec-c_N^{-1}E_N c_N]_{P.A.}=-\pi(\tchi+\tchib)+\pi_N(\tchi_N+\tchib_N)$$
hence the principal acoustical part of the 2nd of \ref{12.475} is:
\begin{equation}
\lambdab\frac{(1+r)}{c}\sh^{n/2}\pi\left(\s^{(n-1)}\ctchi+\s^{(n-1)}\ctchib\right)
\label{12.480}
\end{equation}
In view of the definitions \ref{12.397}, \ref{12.398}, this is bounded in $L^2({\cal K}^{\tau})$ 
for $\tau\in[0,\tau_1]$ by: 
\begin{equation}
C\tau^{c_0+1}\s^{(n-1)}\oX(\tau_1)+C\tau^{c_0}\s^{(n-1)}\oXb(\tau_1)
\label{12.481}
\end{equation}

The above lead to the conclusion that $\Omega^n T\cf$ satisfies the estimate:
\begin{eqnarray}
&&\sup_{\tau\in[0,\tau_1]}\left\{\tau^{-c_0}\|\Omega^n T\cf\|_{L^2({\cal K}^{\tau})}\right\}\leq 
C\s^{(0,n)}\oLab(\tau_1)+C\s^{(n-1)}\oXb(\tau_1)+C\tau_1\s^{(n-1)}\oX(\tau_1)\nonumber\\
&&\hspace{40mm}+C\tau_1\left(\frac{\s^{(1,n)}V(\tau_1)}{\sqrt{2c_0}}+\frac{\s^{(1,n)}P(\tau_1)}{\sqrt{2c_0+2}}\tau_1
+\frac{\s^{(1,n)}\Ga(\tau_1)}{\sqrt{2c_0+4}}\tau_1^2\right)\nonumber\\
&&\hspace{40mm}+C\tau_1\frac{\s^{(Y;0,n)}\cA(\tau_1)}{\sqrt{2c_0}}
\label{12.482}
\end{eqnarray}
for all $\tau_1\in(0,\delta]$. 
Substituting the estimates of Proposition 12.1 for $\s^{(0,n)}\oLab(\tau_1)$, 
$\s^{(n-1)}\oXb(\tau_1)$, $\s^{(n-1)}\oX(\tau_1)$ then yields:
\begin{eqnarray}
&&\sup_{\tau\in[0,\tau_1]}\left\{\tau^{-c_0}\|\Omega^n T\cf\|_{L^2({\cal K}^{\tau})}\right\}
\leq \frac{C\left(\s^{(1,n)}P(\tau_1)+\s^{(1,n)}V(\tau_1)\right)}{\sqrt{2c_0}}
+\frac{C\s^{(1,n)}\Ga(\tau_1)}{\sqrt{2c_0+2}}\tau_1 \nonumber\\
&&\hspace{45mm}+\frac{C\left(\sqrt{\s^{(Y;0,n)}\cA(\tau_1)}+\sqrt{\s^{(E;0,n)}\cA(\tau_1)}
\right)}{\sqrt{2c_0}}\nonumber\\
&&\hspace{45mm}+\frac{C\left(\sqrt{\s^{(Y;0,n)}\cBb(\tau_1,\tau_1)}
+\sqrt{\s^{(E;0,n)}\cBb(\tau_1,\tau_1)}\right)}{a_0+2}\nonumber\\
&&\hspace{45mm}+C\sqrt{\max\{\s^{(0,n)}\cB(\tau_1,\tau_1),\s^{(0,n)}\cBb(\tau_1,\tau_1)\}}
\cdot\tau_1^{1/2}\nonumber\\
&&\label{12.483}
\end{eqnarray}

By the first of \ref{10.798}, \ref{10.800}, \ref{10.802} and by \ref{12.469} we have:
\begin{equation}
\chf=\tau^{-2}\cf 
\label{12.484}
\end{equation}
hence:
\begin{equation}
T\chf=\tau^{-2}T\cf-2\tau^{-3}\cf 
\label{12.485}
\end{equation}
and:
\begin{equation}
\Omega^n T\chf=\tau^{-2}\Omega^n T\cf-2\tau^{-3}\Omega^n\cf
\label{12.486}
\end{equation}
To estimate the 2nd term on the right we apply \ref{12.366} with $\Omega^n\cf$ in the role of $f$. This 
gives:
\begin{eqnarray}
&&\|\Omega^n\cf\|_{L^2(S_{\tau,\tau})}\leq k\tau^{1/2}\|\Omega^n T\cf\|_{L^2({\cal K}^{\tau})}
\nonumber\\
&&\hspace{20mm}\leq k\tau^{c_0+\frac{1}{2}}
\sup_{\tau\in[0,\tau_1]}\left\{\tau^{-c_0}\|\Omega^n T\cf\|_{L^2({\cal K}^{\tau})}\right\}
\label{12.487}
\end{eqnarray}
Therefore, replacing $\tau$ by $\tau^{\prime}\in [0,\tau]$, 
\begin{eqnarray}
&&\|\tau^{\prime-3}\Omega^n\cf\|_{L^2({\cal K}^{\tau})}
=\left(\int_0^{\tau}\tau^{\prime-6}\|\Omega^n\cf\|^2_{L^2(S_{\tau^\prime,\tau^\prime})}d\tau^\prime
\right)^{1/2}\nonumber\\
&&\hspace{20mm}\leq k\sup_{\tau\in[0,\tau_1]}\left\{\tau^{-c_0}\|\Omega^n T\cf\|_{L^2({\cal K}^{\tau})}
\right\}\cdot\left(\int_0^{\tau}\tau^{\prime 2c_0-5}d\tau^\prime\right)^{1/2}\nonumber\\
&&\hspace{20mm}=\frac{k\tau^{c_0-2}}{\sqrt{2c_0-4}}\cdot
\sup_{\tau\in[0,\tau_1]}\left\{\tau^{-c_0}\|\Omega^n T\cf\|_{L^2({\cal K}^{\tau})}
\right\}\nonumber\\
&&\label{12.488}
\end{eqnarray}
As for the 1st term on the right in \ref{12.486}, we have:
\begin{eqnarray}
&&\|\tau^{\prime-2}\Omega^n T\cf\|^2_{L^2({\cal K}^{\tau})}
=\int_0^{\tau}\tau^{\prime-4}\|\Omega^n T\cf\|^2_{L^2(S_{\tau^\prime,\tau^\prime})}d\tau^\prime
\nonumber\\
&&\hspace{30mm}=\int_0^{\tau}\tau^{\prime-4}\frac{\partial}{\partial\tau^{\prime}}\left(
\|\Omega^n T\cf\|^2_{L^2({\cal K}^{\tau^\prime})}\right)d\tau^\prime\nonumber\\
&&\hspace{20mm}=\tau^{-4}\|\Omega^n T\cf\|^2_{L^2({\cal K}^{\tau})}
+4\int_0^{\tau}\tau^{\prime-5}\|\Omega^n T\cf\|^2_{L^2({\cal K}^{\tau^\prime})}d\tau^\prime\nonumber\\
&&\hspace{20mm}\leq\left(\sup_{\tau\in[0,\tau_1]}\left\{\tau^{-c_0}\|\Omega^n T\cf\|_{L^2({\cal K}^{\tau})}
\right\}\right)^2\cdot\left(1+\frac{4}{2c_0-4}\right)\tau^{2c_0-4}\nonumber\\
&&\label{12.489}
\end{eqnarray}
It then follows through \ref{12.486} that for all $\tau\in[0,\tau_1]$:
\begin{equation}
\|\Omega^n T\chf\|_{L^2({\cal K}^{\tau})}\leq C\tau^{c_0-2}\sup_{\tau\in[0,\tau_1]}\left\{\tau^{-c_0}\|\Omega^n T\cf\|_{L^2({\cal K}^{\tau})}
\right\}
\label{12.490}
\end{equation}
Multiplying then by $\tau^{-c_0+2}$ and taking the supremum over $\tau\in[0,\tau_1]$ we obtain, 
in view of the definition \ref{12.394}, 
\begin{equation}
\s^{(1,n)}P(\tau_1)\leq 
C\sup_{\tau\in[0,\tau_1]}\left\{\tau^{-c_0}\|\Omega^n T\cf\|_{L^2({\cal K}^{\tau})}\right\}
\label{12.491}
\end{equation}
Substituting on the right the estimate \ref{12.483} and choosing $c_0$ large enough so that in 
reference to the resulting coefficient of $\s^{(1,n)}P(\tau_1)$ on the right it holds:
\begin{equation}
\frac{C}{\sqrt{2c_0}}\leq \frac{1}{2}
\label{12.492}
\end{equation}
we conclude that:
\begin{eqnarray}
&&\s^{(1,n)}P(\tau_1)
\leq \frac{C\s^{(1,n)}V(\tau_1)}{\sqrt{2c_0}}
+\frac{C\s^{(1,n)}\Ga(\tau_1)}{\sqrt{2c_0+2}}\tau_1 \nonumber\\
&&\hspace{18mm}+\frac{C\left(\sqrt{\s^{(Y;0,n)}\cA(\tau_1)}+\sqrt{\s^{(E;0,n)}\cA(\tau_1)}
\right)}{\sqrt{2c_0}}\nonumber\\
&&\hspace{18mm}+\frac{C\left(\sqrt{\s^{(Y;0,n)}\cBb(\tau_1,\tau_1)}
+\sqrt{\s^{(E;0,n)}\cBb(\tau_1,\tau_1)}\right)}{a_0+2}\nonumber\\
&&\hspace{18mm}+C\sqrt{\max\{\s^{(0,n)}\cB(\tau_1,\tau_1),\s^{(0,n)}\cBb(\tau_1,\tau_1)\}}
\cdot\tau_1^{1/2}\nonumber\\
&&\label{12.493}
\end{eqnarray}
(for new constants $C$). 

Consider next the second of \ref{4.119}:
\begin{equation}
Tg^i=\rho N^i+\rhob\Nb^i \ : \ \mbox{along ${\cal K}$}
\label{12.494}
\end{equation}
In view of the fact that $\rho=c^{-1}\lambdab$, $\rhob=c^{-1}\lambda$ and 
the boundary condition \ref{4.6} (or \ref{10.752}) this takes the form:
\begin{equation}
Tg^i=\frac{(N^i+r\Nb^i)}{c}\lambdab \ : \ \mbox{along ${\cal K}$}
\label{12.495}
\end{equation}
By \ref{9.51} and \ref{9.15} $g_N^i$ the $N$th approximate analogues of the $g^i$ satisfy:
\begin{equation}
Tg_N^i=\rho_N N_N^i+\rhob_N\Nb_N^i=\vep_N^i+\vepb_N^i
\label{12.496}
\end{equation}
In view of \ref{9.101} and \ref{9.147} this can be written in the form: 
\begin{equation}
Tg_N^i=\frac{(N_N^i+r_N\Nb_N^i)}{c_N}\lambdab_N + \vep_N^{\prime i}
\label{12.497}
\end{equation}
where:
\begin{equation}
\vep_N^{\prime i}=\vep_N^i+\vepb_N^i+\frac{\hat{\nu}_N}{c_N}\Nb_N^i=O(\tau^N)
\label{12.498}
\end{equation}
by Proposition 9.1 and \ref{9.106}. Let us presently denote by $\check{g}^i$ the difference:
\begin{equation}
\check{g}^i=g^i-g_N^i
\label{12.499}
\end{equation}
Note that this does not agree with \ref{4.222}. In the present section $\check{g}^i$ stands for the 
difference \ref{12.499}. Subtracting \ref{12.497} from \ref{12.495} we obtain:
\begin{equation}
T\check{g}^i=\frac{(N^i+r\Nb^i)}{c}\lambdab-\frac{(N_N^i+r_N\Nb_N^i)}{c_N}\lambdab_N-\vep_N^{\prime i} 
\ : \ \mbox{along ${\cal K}$}
\label{12.500}
\end{equation}
Recalling the definition \ref{4.223} and the corresponding definition for the $N$th approximants,  
\ref{9.115}, let us define:
\begin{equation}
\check{\delta}^i=\delta^i-\delta_N^i
\label{12.501}
\end{equation}
By \ref{4.222} and \ref{9.114} we have, in terms of the present notations \ref{12.469}, 
\ref{12.499}, 
\begin{equation}
\check{\delta}^i=\check{g}^i-N_0^i\cf
\label{12.502}
\end{equation}
Then by \ref{12.470} and \ref{12.500} the $\check{\delta}^i$ satisfy: 
\begin{equation}
T\check{\delta}^i=U^i\lambdab-U_N^i\lambdab_N-\tilde{\vep}_N^i
\label{12.503}
\end{equation}
where:
\begin{eqnarray}
&&U^i=\frac{1}{c}\left[N^i-N_0^i+r(\Nb^i-N_0^i)\right] \nonumber\\
&&U_N^i=\frac{1}{c_N}\left[N_N^i-N_0^i+r_N(\Nb_N^i-N_0^i)\right] 
\label{12.504}
\end{eqnarray}
and:
\begin{equation}
\tilde{\vep}_N^i=\vep_N^{\prime i}-\vep_N N_0^i=O(\tau^N)
\label{12.505}
\end{equation}

We apply $\Omega^n$ to \ref{12.503} to obtain:
\begin{equation}
\Omega^n T\check{\delta}^i=\Omega^n(U^i\lambdab-U_N^i\lambdab_N)-\Omega^n\tilde{\vep}_N^i
\label{12.506}
\end{equation}
and from \ref{12.505} we have:
\begin{equation}
\Omega^n\tilde{\vep}_N^i=O(\tau^N)
\label{12.507}
\end{equation}
The principal part of the right hand side of \ref{12.506} is contained in:
\begin{equation}
U^i\Omega^n(\lambdab-\lambdab_N)+\lambdab\Omega^n(U^i-U_N^i)
\label{12.508}
\end{equation}
In fact, the principal part of the first of \ref{12.506} is:
\begin{equation}
\sh^{n/2}U^i\s^{(0,n)}\clab
\label{12.509}
\end{equation}
while the principal part of the second of \ref{12.506} is contained in:
\begin{equation}
\lambdab\sh^{n/2}(E^n U^i-E_N^n U_N^i)
\label{12.510}
\end{equation}
or, more precisely, in:
\begin{eqnarray}
&&\frac{\lambdab}{c}\sh^{n/2}\left[(E^n N^i-E_N^n N_N^i)+r(E^n\Nb^i-E_N^n\Nb_N^i)\right.\nonumber\\
&&\hspace{33mm}+\left. (\Nb^i-N_0^i)(E^n r-E_N^n r_N)\right]\nonumber\\
&&\hspace{5mm}-\lambdab U^i\sh^{n/2}\left[E^{n-1}(c^{-1}Ec)-E_N^{n-1}(c_N^{-1}E_N c_N)\right] 
\label{12.511}
\end{eqnarray}
Now,  
$$(N^i-N_0^i)(\tau,\vartheta)=\int_0^{\tau}(TN^i)(\tau^\prime,\vartheta)d\tau^\prime$$
From \ref{3.a15}:
\begin{equation}
TN^i=(\sm+\sn)E^i+(m+n)N^i+(\om+\on)\Nb^i
\label{12.512}
\end{equation}
and (see 1st of \ref{3.a24} and of \ref{3.48}, and 1st of \ref{3.a21} and of \ref{3.46}) we have:
\begin{equation}
[TN^i]_{P.A.}=[\sn]_{P.A.}(E^i-\pi N^i), \ \ \ [\sn]_{P.A.}=2E\lambda+2\pi\lambda\tchi
\label{12.513}
\end{equation}
The assumptions \ref{10.439}, \ref{10.440} then imply:
\begin{equation}
|TN^i|\leq C \ : \ \mbox{in ${\cal R}_{\delta,\delta}$}
\label{12.514}
\end{equation}
Hence along ${\cal K}$, for $\tau\in[0,\delta]$ we have: 
\begin{equation}
|N^i-N_0^i|\leq C\tau 
\label{12.515}
\end{equation}
which together with \ref{10.775} implies (see 1st of \ref{12.504}): 
\begin{equation}
|U^i|\leq C\tau \ : \ \forall\tau\in[0,\delta]
\label{12.516}
\end{equation}
Given then $\tau_1\in(0,\delta]$, \ref{12.509} is bounded in $L^2({\cal K}^{\tau})$ for 
$\tau\in[0,\tau_1]$ by:
\begin{equation}
C\s^{(0,n)}\oLab(\tau_1)\tau^{c_0+1}
\label{12.517}
\end{equation}
Since 
\begin{eqnarray}
&&[E^n N^i-E_N^n N_N^i]_{P.A.}=\s^{(n-1)}\ctchi (E^i-\pi N^i) \nonumber\\
&&[E^n\Nb^i-E_N^n\Nb_N^i]_{P.A.}=\s^{(n-1)}\ctchib (E^i-\pi\Nb^i) \label{12.518}
\end{eqnarray}
the contribution of the first two terms in parenthesis in the 1st of \ref{12.511} to the 
$L^2({\cal K}^{\tau})$ norm of \ref{12.511} is bounded for $\tau\in[0,\tau_1]$, to leading terms, by (see \ref{12.397}, 
\ref{12.398}): 
\begin{equation}
C\left(\s^{(n-1)}\oX(\tau_1)+\s^{(n-1)}\oXb(\tau_1)\right)\tau^{c_0+1}
\label{12.519}
\end{equation}
Moreover, by \ref{12.479} with $\tau$ in the role of $\tau_1$ the contribution of the third term in 
parenthesis in the 1st of \ref{12.511} to the 
$L^2({\cal K}^{\tau})$ norm of \ref{12.511} is bounded for $\tau\in[0,\tau_1]$, to leading terms, by:
\begin{eqnarray}
&&C\left\{\s^{(n-1)}\oXb(\tau_1)+\frac{C\s^{(Y;0,n)}\cA(\tau_1)}{\sqrt{2c_0}}\right.\label{12.520}\\
&&\hspace{10mm}\left.+\frac{C\s^{(1,n)}V(\tau_1)}{\sqrt{2c_0}}
+\frac{C\s^{(1,n)}P(\tau_1)}{\sqrt{2c_0+2}}\tau+\frac{C\s^{(1,n)}\Ga(\tau_1)}{\sqrt{2c_0+4}}\tau^2
\right\}\tau^{c_0+1}\nonumber
\end{eqnarray}
Also, by \ref{12.480}, \ref{12.481} and \ref{12.516}, the 2nd of \ref{12.511} is bounded in 
$L^2({\cal K}^{\tau})$ for $\tau\in[0,\tau_1]$ by:
\begin{equation}
C\left(\s^{(n-1)}\oXb(\tau_1)+\s^{(n-1)}\oX(\tau_1)\tau\right)\tau^{c_0+1}
\label{12.521}
\end{equation}
The above lead to the conclusion that $\Omega^n T\check{\delta}^i$ satisfies the estimate:
\begin{eqnarray}
&&\sup_{\tau\in[0,\tau_1]}\left\{\tau^{-c_0-1}\|\Omega^n T\check{\delta}^i\|_{L^2({\cal K}^{\tau})}\right\}\leq 
C\s^{(0,n)}\oLab(\tau_1)+C\s^{(n-1)}\oXb(\tau_1)\nonumber\\
&&\hspace{45mm}+C\s^{(n-1)}\oX(\tau_1)
+\frac{C\s^{(1,n)}V(\tau_1)}{\sqrt{2c_0}}\nonumber\\
&&\hspace{45mm}+\frac{C\s^{(1,n)}P(\tau_1)}{\sqrt{2c_0+2}}\tau_1
+\frac{C\s^{(1,n)}\Ga(\tau_1)}{\sqrt{2c_0+4}}\tau_1^2\nonumber\\
&&\hspace{45mm}+\frac{C\s^{(Y;0,n)}\cA(\tau_1)}{\sqrt{2c_0}} \label{12.522}
\end{eqnarray}
Substituting the estimates of Proposition 12.1 for $\s^{(0,n)}\oLab(\tau_1)$, $\s^{(n-1)}\oXb(\tau_1)$, 
$\s^{(n-1)}\oX(\tau_1)$ as well as the estimate \ref{12.493} for $\s^{(1,n)}P(\tau_1)$ then yields:
\begin{eqnarray}
&&\sup_{\tau\in[0,\tau_1]}\left\{\tau^{-c_0-1}\|\Omega^n T\check{\delta}^i\|_{L^2({\cal K}^{\tau})}\right\}\leq \frac{C\s^{(1,n)}V(\tau_1)}{\sqrt{2c_0}}
+\frac{C\s^{(1,n)}\Ga(\tau_1)}{\sqrt{2c_0+2}}\tau_1 \nonumber\\
&&\hspace{45mm}+\frac{C\left(\sqrt{\s^{(Y;0,n)}\cA(\tau_1)}+\sqrt{\s^{(E;0,n)}\cA(\tau_1)}
\right)}{\sqrt{2c_0}}\nonumber\\
&&\hspace{40mm}+\frac{C\left(\sqrt{\s^{(Y;0,n)}\cBb(\tau_1,\tau_1)}
+\sqrt{\s^{(E;0,n)}\cBb(\tau_1,\tau_1)}\right)}{a_0+2}\nonumber\\
&&\hspace{40mm}+C\sqrt{\max\{\s^{(0,n)}\cB(\tau_1,\tau_1),\s^{(0,n)}\cBb(\tau_1,\tau_1)\}}
\cdot\tau_1^{1/2}\nonumber\\
\label{12.523}
\end{eqnarray}
(for new constants $C$). 

Let us define the difference:
\begin{equation}
\chdl^i=\hdl^i-\hdl_N^i=\tau^{-3}\cdl^i
\label{12.524}
\end{equation}
(see \ref{12.469}, \ref{12.501} and \ref{4.223}, \ref{9.116}). Then we have:
\begin{equation}
T\chdl^i=\tau^{-3}T\cdl^i-3\tau^{-4}\cdl^i
\label{12.525}
\end{equation}
and:
\begin{equation}
\Omega^n T\chdl^i=\tau^{-3}\Omega^n T\cdl^i-3\tau^{-4}\Omega^n\cdl^i
\label{12.526}
\end{equation}
To estimate the 2nd term on the right we apply \ref{12.366} with $\Omega^n\cdl^i$ in the role of $f$. 
This gives:
\begin{eqnarray}
&&\|\Omega^n\cdl^i\|_{L^2(S_{\tau,\tau})}\leq k\tau^{1/2}\|\Omega^n T\cdl^i\|_{L^2({\cal K}^{\tau})}
\nonumber\\
&&\hspace{20mm}\leq k\tau^{c_0+\frac{3}{2}}
\sup_{\tau\in[0,\tau_1]}\left\{\tau^{-c_0-1}\|\Omega^n T\cdl^i\|_{L^2({\cal K}^{\tau})}\right\}
\label{12.527}
\end{eqnarray}
Therefore, replacing $\tau$ by $\tau^{\prime}\in [0,\tau]$, 
\begin{eqnarray}
&&\|\tau^{\prime-4}\Omega^n\cdl^i\|_{L^2({\cal K}^{\tau})}
=\left(\int_0^{\tau}\tau^{\prime-8}\|\Omega^n\cdl^i\|^2_{L^2(S_{\tau^\prime,\tau^\prime})}d\tau^\prime
\right)^{1/2}\nonumber\\
&&\hspace{20mm}\leq k\sup_{\tau\in[0,\tau_1]}\left\{\tau^{-c_0-1}\|\Omega^n T\cdl^i\|_{L^2({\cal K}^{\tau})}
\right\}\cdot\left(\int_0^{\tau}\tau^{\prime 2c_0-5}d\tau^\prime\right)^{1/2}\nonumber\\
&&\hspace{20mm}=\frac{k\tau^{c_0-2}}{\sqrt{2c_0-4}}\cdot
\sup_{\tau\in[0,\tau_1]}\left\{\tau^{-c_0-1}\|\Omega^n T\cdl^i\|_{L^2({\cal K}^{\tau})}
\right\}\nonumber\\
&&\label{12.528}
\end{eqnarray}
As for the 1st term on the right in \ref{12.526}, we have:
\begin{eqnarray}
&&\|\tau^{\prime-3}\Omega^n T\cdl^i\|^2_{L^2({\cal K}^{\tau})}
=\int_0^{\tau}\tau^{\prime-6}\|\Omega^n T\cdl^i\|^2_{L^2(S_{\tau^\prime,\tau^\prime})}d\tau^\prime
\nonumber\\
&&\hspace{30mm}=\int_0^{\tau}\tau^{\prime-6}\frac{\partial}{\partial\tau^{\prime}}\left(
\|\Omega^n T\cdl^i\|^2_{L^2({\cal K}^{\tau^\prime})}\right)d\tau^\prime\nonumber\\
&&\hspace{20mm}=\tau^{-6}\|\Omega^n T\cdl^i\|^2_{L^2({\cal K}^{\tau})}
+6\int_0^{\tau}\tau^{\prime-7}\|\Omega^n T\cdl^i\|^2_{L^2({\cal K}^{\tau^\prime})}d\tau^\prime\nonumber\\
&&\hspace{20mm}\leq\left(\sup_{\tau\in[0,\tau_1]}\left\{\tau^{-c_0-1}\|\Omega^n T\cdl^i\|_{L^2({\cal K}^{\tau})}
\right\}\right)^2\cdot\left(1+\frac{6}{2c_0-4}\right)\tau^{2c_0-4}\nonumber\\
&&\label{12.529}
\end{eqnarray}
It then follows through \ref{12.526} that for all $\tau[0,\tau_1]$:
\begin{equation}
\|\Omega^n T\chdl^i\|_{L^2({\cal K}^{\tau})}\leq C\tau^{c_0-2}\sup_{\tau\in[0,\tau_1]}
\left\{\tau^{-c_0-1}\|\Omega^n T\cdl^i\|_{L^2({\cal K}^{\tau})}\right\}
\label{12.530}
\end{equation}
Multiplying then by $\tau^{-c_0+2}$ and taking the supremum over $\tau\in[0,\tau_1]$ we obtain:
\begin{equation}
\sup_{\tau\in[0,\tau_1]}
\left\{\tau^{-c_0+2}\|\Omega^n T\chdl^i\|_{L^2({\cal K}^{\tau})}\right\}\leq 
C\sup_{\tau\in[0,\tau_1]}
\left\{\tau^{-c_0-1}\|\Omega^n T\cdl^i\|_{L^2({\cal K}^{\tau})}\right\}
\label{12.531}
\end{equation}
Hence, defining:
\begin{equation}
\s^{(1,n)}Q(\tau_1)=\sup_{\tau\in[0,\tau_1]}
\left\{\tau^{-c_0+2}\sum_i\|\Omega^n T\chdl^i\|_{L^2({\cal K}^{\tau})}\right\} 
\label{12.532}
\end{equation}
we have, by \ref{12.523}:
\begin{eqnarray}
&&\s^{(1,n)}Q(\tau_1)\leq \frac{C\s^{(1,n)}V(\tau_1)}{\sqrt{2c_0}}
+\frac{C\s^{(1,n)}\Ga(\tau_1)}{\sqrt{2c_0+2}}\tau_1 \nonumber\\
&&\hspace{20mm}+\frac{C\left(\sqrt{\s^{(Y;0,n)}\cA(\tau_1)}+\sqrt{\s^{(E;0,n)}\cA(\tau_1)}
\right)}{\sqrt{2c_0}}\nonumber\\
&&\hspace{20mm}+\frac{C\left(\sqrt{\s^{(Y;0,n)}\cBb(\tau_1,\tau_1)}
+\sqrt{\s^{(E;0,n)}\cBb(\tau_1,\tau_1)}\right)}{a_0+2}\nonumber\\
&&\hspace{20mm}+C\sqrt{\max\{\s^{(0,n)}\cB(\tau_1,\tau_1),\s^{(0,n)}\cBb(\tau_1,\tau_1)\}}
\cdot\tau_1^{1/2}\nonumber\\
\label{12.533}
\end{eqnarray}
(for new constants $C$). 

We now bring in from Proposition 4.5 the identification equations:
\begin{equation}
\hat{F}^i((\tau,\vartheta), (v(\tau,\vartheta),\gamma(\tau,\vartheta)))=0
\label{12.534}
\end{equation}
and from Section 9.4 the corresponding equations, \ref{9.119}, for the $N$th approximants:
\begin{equation}
\hat{F}^i_N((\tau,\vartheta),(v_N(\tau,\vartheta),\gamma_N(\tau,\vartheta)))=D^i_N(\tau,\vartheta)
\label{12.535}
\end{equation}
Recall that according to Proposition 9.7 we have:
\begin{equation}
D^i_N(\tau,\vartheta)=O(\tau^{N-1})
\label{12.536}
\end{equation}
Subtracting \ref{12.535} from \ref{12.534} we obtain:
\begin{equation}
\hat{F}^i((\tau,\vartheta), (v(\tau,\vartheta),\gamma(\tau,\vartheta)))
-\hat{F}^i_N((\tau,\vartheta),(v_N(\tau,\vartheta),\gamma_N(\tau,\vartheta)))=-D^i_N(\tau,\vartheta)
\label{12.537}
\end{equation}
Differentiating these equations implicitly with respect to $\tau$ yields:
$$\frac{\partial\hat{F}^i}{\partial v}\frac{\partial v}{\partial\tau}
+\frac{\partial\hat{F}^i}{\partial\gamma}\frac{\partial\gamma}{\partial\tau}
+\frac{\partial\hat{F}^i}{\partial\tau}
-\frac{\partial\hat{F}^i_N}{\partial v}\frac{\partial v_N}{\partial\tau}
-\frac{\partial\hat{F}^i_N}{\partial\gamma}\frac{\partial\gamma_N}{\partial\tau}
-\frac{\partial\hat{F}^i_N}{\partial\tau}=-\frac{\partial D^i_N}{\partial\tau}$$
Here, the argument of the partial derivatives of $\hat{F}^i$ is 
$((\tau,\vartheta),(v(\tau,\vartheta),\gamma(\tau,\vartheta)))$, while the argument of the partial 
derivatives of $\hat{F}^i_N$ is $((\tau,\vartheta),(v_N(\tau,\vartheta),\gamma_N(\tau,\vartheta)))$. 
In view of the definitions \ref{12.1} these equations take the form:
\begin{eqnarray}
&&\frac{\partial\hat{F}^i}{\partial v}\frac{\partial\cv}{\partial\tau}
+\frac{\partial\hat{F}^i}{\partial\gamma}\frac{\partial\cga}{\partial\tau}=
-\left(\frac{\partial\hat{F}^i}{\partial v}-\frac{\partial\hat{F}^i_N}{\partial v}\right)
\frac{\partial v_N}{\partial\tau}
-\left(\frac{\partial\hat{F}^i}{\partial\gamma}-\frac{\partial\hat{F}^i_N}{\partial\gamma}\right)
\frac{\partial\gamma_N}{\partial\tau}\nonumber\\
&&\hspace{45mm}-\left(\frac{\partial\hat{F}^i}{\partial\tau}-\frac{\partial\hat{F}^i_N}{\partial\tau}\right)-\frac{\partial D^i_N}{\partial\tau}
\label{12.538}
\end{eqnarray}
We remark that, the transformation functions $\hat{f}$, $v$, $\gamma$ being of order 0, the left 
hand side here is of order 1, the first two terms on the right hand side are of order 0, while 
the third term on the right hand side is of order 1. From Proposition 4.5:
\begin{eqnarray}
&&\frac{\partial\hat{F}^i}{\partial\tau}((\tau,\vartheta),(v,\gamma))=
S_0^i(\vartheta)l(\vartheta)v\frac{\partial\hat{f}}{\partial\tau}(\tau,\vartheta)
-\frac{\partial\hat{\delta}^i}{\partial\tau}(\tau,\vartheta)\nonumber\\
&&\hspace{23mm}+\tau\frac{\partial E^i}{\partial\tau}((\tau,\vartheta),(v,\gamma))
+E^i((\tau,\vartheta),(v,\gamma))
\label{12.539}
\end{eqnarray}
and from \ref{9.117}:
\begin{eqnarray}
&&\frac{\partial\hat{F}^i_N}{\partial\tau}((\tau,\vartheta),(v,\gamma))=
S_0^i(\vartheta)l(\vartheta)v\frac{\partial\hat{f}_N}{\partial\tau}(\tau,\vartheta)
-\frac{\partial\hat{\delta}_N^i}{\partial\tau}(\tau,\vartheta)\nonumber\\
&&\hspace{23mm}+\tau\frac{\partial E_N^i}{\partial\tau}((\tau,\vartheta),(v,\gamma))
+E_N^i((\tau,\vartheta),(v,\gamma))
\label{12.540}
\end{eqnarray}
Therefore the third term on the right in \ref{12.538} is:
\begin{eqnarray}
&&-\left(\frac{\partial\hat{F}^i}{\partial\tau}-\frac{\partial\hat{F}^i_N}{\partial\tau}\right)=
-S_0^i l v\frac{\partial\chf}{\partial\tau}+\frac{\partial\chdl^i}{\partial\tau}
-S_0^i l\cv\frac{\partial\hat{f}_N}{\partial\tau}\nonumber\\
&&\hspace{32mm}-\tau\left(\frac{\partial E^i}{\partial\tau}-\frac{\partial E^i_N}{\partial\tau}\right)
-\left(E^i-E^i_N\right) 
\label{12.541}
\end{eqnarray}
Here, the 3rd and 5th terms on the right are of order 0 while, in view of the factor $\tau$, 
the 4th term can be absorbed imposing a smallness condition on $\delta$. 

We conclude from the above that 
equations \ref{12.538} to leading terms read, recalling that $T=\partial/\partial\tau$, 
\begin{equation}
\frac{\partial\hat{F}^i}{\partial v}T\cv+\frac{\partial\hat{F}^i}{\partial\gamma}T\cga=
-S_0^i l v T\chf+T\chdl^i
\label{12.542}
\end{equation}
Differentiating equations \ref{12.538} implicitly with respect to $\vartheta$ $n$ times we obtain, 
recalling that $\Omega=\partial/\partial\vartheta$, equations of the form:
\begin{equation}
\frac{\partial\hat{F}^i}{\partial v}\Omega^n T\cv
+\frac{\partial\hat{F}^i}{\partial\gamma}\Omega^n T\cga=R^i_{0,n}
\label{12.543}
\end{equation}
where the leading part of $R^i_{0,n}$ is:
\begin{equation}
-S_0^i l v \Omega^n T\chf+\Omega^n T\chdl^i
\label{12.544}
\end{equation}
By Proposition 4.5 and \ref{4.241}:
\begin{equation}
\frac{\partial\hat{F}^i}{\partial v}((0,\vartheta),(-1,\gamma_0(\vartheta)))
=\frac{1}{3}k(\vartheta)S_0^i(\vartheta), \ \ \ 
\frac{\partial\hat{F}^i}{\partial\gamma}((0,\vartheta),(-1,\gamma_0(\vartheta)))
=\Omega_0^{\prime i}(\vartheta)
\label{12.545}
\end{equation}
Now, we can assume that for all $\tau\in[0,\delta]$ it holds:
\begin{equation}
|\hat{f}(\tau,\vartheta)|\leq C, \ \ |v(\tau,\vartheta)+1|\leq C\tau, \ \ 
|\gamma(\tau,\vartheta)-\gamma_0(\vartheta)|\leq C\tau \ \ : \ \forall\vartheta\in S^1
\label{12.546}
\end{equation}
It then follows that 
$$\frac{\partial E^i}{\partial v}((\tau,\vartheta),(v(\tau,\vartheta),\gamma(\tau,\vartheta))), \ 
\frac{\partial E^i}{\partial\gamma}((\tau,\vartheta),(v(\tau,\vartheta),\gamma(\tau,\vartheta)))$$
are accordingly bounded, and defining: 
\begin{eqnarray}
&&e^i(\tau,\vartheta)=\frac{\partial\hat{F}^i}{\partial v}((\tau,\vartheta),
(v(\tau,\vartheta),\gamma(\tau,\vartheta)))-\frac{1}{3}k(\vartheta)S_0^i(\vartheta)\nonumber\\
&&h^i(\tau,\vartheta)=\frac{\partial\hat{F}^i}{\partial\gamma}((\tau,\vartheta)),
(v(\tau,\vartheta),\gamma(\tau,\vartheta)))-\Omega_0^{\prime i}(\vartheta) 
\label{12.547}
\end{eqnarray}
we have, for all $\tau\in[0,\delta]$:
\begin{equation}
\max_i|e^i(\tau,\vartheta)|\leq C\tau, \ \ \max_i|h^i(\tau,\vartheta)|\leq C\tau \ \ : \ 
\forall\vartheta\in S^1
\label{12.548}
\end{equation}
As already noted in the course of the proof of Proposition 4.4, for each $\vartheta\in S^1$ the vectors 
\begin{equation}
S_0(\vartheta)=S_0^i(\vartheta)\frac{\partial}{\partial x^i}, \ \ 
\Omega_0^\prime(\vartheta)=\Omega_0^{\prime i}\frac{\partial}{\partial x^i}
\label{12.549}
\end{equation}
are linearly independent and span $\mathbb{R}^2$. Consequently any $V^i(\tau,\vartheta)$ 
can be expanded as:
\begin{equation}
V^i(\tau,\vartheta)=V_{\bot}(\tau,\vartheta)S_0^i(\vartheta)
+V_{||}(\tau,\vartheta)\Omega_0^{\prime i}(\vartheta)
\label{12.550}
\end{equation}
$V_{\bot}(\tau,\vartheta)$ and $V_{||}(\tau,\vartheta)$ being respectively the $S_0(\vartheta)$ 
and $\Omega_0^\prime(\vartheta)$ components of the vector 
$$V(\tau,\vartheta)=V^i(\tau,\vartheta)\frac{\partial}{\partial x^i}$$
In particular, we can expand:
\begin{eqnarray}
&&e^i(\tau,\vartheta)=e_{\bot}(\tau,\vartheta)S_0^i(\vartheta)
+e_{||}(\tau,\vartheta)\Omega_0^{\prime i}(\vartheta) \nonumber\\
&&h^i(\tau,\vartheta)=h_{\bot}(\tau,\vartheta)S_0^i(\vartheta)
+h_{||}(\tau,\vartheta)\Omega_0^{\prime i}(\vartheta)
\label{12.551}
\end{eqnarray}
and by \ref{12.548} we have:
\begin{equation}
|e_{\bot}(\tau,\vartheta)|, |e_{||}(\tau,\vartheta)|, 
|h_{\bot}(\tau,\vartheta)|, |h_{||}(\tau,\vartheta)| \ \leq \ C\tau \ \ : \ \forall\vartheta\in S^1, 
\ \forall\tau\in[0,\delta]
\label{12.552}
\end{equation}
We can also expand, in reference to the right hand side of \ref{12.543}, 
\begin{equation}
R^i_{0,n}=R_{0,n\bot}S_0^i+R_{0,n||}\Omega_0^{\prime i}
\label{12.553}
\end{equation}
Then in view of \ref{12.547} and \ref{12.551}, equations \ref{12.543} decompose into:
\begin{eqnarray}
&&\left(\frac{k}{3}+e_{\bot}\right)\Omega^n T\cv+h_{\bot}\Omega^n T\cga=R_{0,n\bot} \nonumber\\
&&e_{||}\Omega^n T\cv+(1+h_{||})\Omega^n T\cga=R_{0,n||} 
\label{12.554}
\end{eqnarray}
their $S_0$ and $\Omega_0^\prime$ components respectively. Solving for 
$(\Omega^n T\cv,\Omega^n T\cga)$ we obtain:
\begin{eqnarray}
&&\Omega^n T\cv=\frac{R_{0,n\bot}-(h_{\bot}/(1+h_{||}))R_{0,n||}}{(k/3)+e_{\bot}
-((e_{||}h_{\bot})/(1+h_{||}))} \nonumber\\
&&\Omega^n T\cga=\frac{R_{0,n||}-(e_{||}/((k/3)+e_{\bot}))R_{0,n\bot}}{1+h_{||}
-((e_{||}h_{\bot})/((k/3)+e_{\bot}))} 
\label{12.555}
\end{eqnarray}
In view of \ref{12.552} these imply that, pointwise: 
\begin{eqnarray}
&&|\Omega^n T\cv|\leq C\sum_i|R^i_{0,n}| \nonumber\\
&&|\Omega^n T\cga|\leq C\sum_i|R^i_{0,n}| 
\label{12.556}
\end{eqnarray}
Hence, taking $L^2$ norms on ${\cal K}^{\tau}$:
\begin{eqnarray}
&&\|\Omega^n T\cv\|_{L^2({\cal K}^{\tau})}\leq C\sum_i\|R^i_{0,n}\|_{L^2({\cal K}^{\tau})}\nonumber\\
&&\|\Omega^n T\cga\|_{L^2({\cal K}^{\tau})}\leq C\sum_i\|R^i_{0,n}\|_{L^2({\cal K}^{\tau})}
\label{12.557}
\end{eqnarray}
Now, by \ref{12.544} and the definitions \ref{12.394} and \ref{12.532} we have, to leading terms, 
\begin{eqnarray}
&&\sum_i\|R^i_{0,n}\|_{L^2({\cal K}^{\tau})}\leq C\left(\|\Omega^n T\chf\|_{L^2({\cal K}^{\tau})}
+\sum_i\|\Omega^n T\chdl^i\|_{L^2({\cal K}^{\tau})}\right) \nonumber\\
&&\hspace{27mm}\leq C\tau^{c_0-2}\left(\s^{(1,n)}P(\tau_1)+\s^{(1,n)}Q(\tau_1)\right)
\label{12.558}
\end{eqnarray}
for all $\tau\in[0,\tau_1]$, $\tau_1\in(0,\delta]$. Substituting in \ref{12.557}, multiplying by 
$\tau^{-c_0+2}$ and taking the supremum over $\tau\in[0,\tau_1]$ then yields, in view of the 
definitions \ref{12.395}, \ref{12.396}, 
\begin{eqnarray}
&&\s^{(1,n)}V(\tau_1)\leq C\left(\s^{(1,n)}P(\tau_1)+\s^{(1,n)}Q(\tau_1)\right)\nonumber\\
&&\s^{(1,n)}\Ga(\tau_1)\leq C\left(\s^{(1,n)}P(\tau_1)+\s^{(1,n)}Q(\tau_1)\right)
\label{12.559}
\end{eqnarray}
Adding these inequalities and substituting on the right the estimates \ref{12.493} and \ref{12.533} for $\s^{(1,n)}P(\tau_1)$ 
and $\s^{(1,n)}Q(\tau_1)$ respectively, we obtain:
\begin{eqnarray}
&&\s^{(1,n)}V(\tau_1)+\s^{(1,n)}\Ga(\tau_1)\leq \frac{C\s^{(1,n)}V(\tau_1)}{\sqrt{2c_0}}
+\frac{C\s^{(1,n)}\Ga(\tau_1)}{\sqrt{2c_0+2}}\tau_1 \nonumber\\
&&\hspace{35mm}+\frac{C\left(\sqrt{\s^{(Y;0,n)}\cA(\tau_1)}+\sqrt{\s^{(E;0,n)}\cA(\tau_1)}
\right)}{\sqrt{2c_0}}\nonumber\\
&&\hspace{30mm}+\frac{C\left(\sqrt{\s^{(Y;0,n)}\cBb(\tau_1,\tau_1)}
+\sqrt{\s^{(E;0,n)}\cBb(\tau_1,\tau_1)}\right)}{a_0+2}\nonumber\\
&&\hspace{30mm}+C\sqrt{\max\{\s^{(0,n)}\cB(\tau_1,\tau_1),\s^{(0,n)}\cBb(\tau_1,\tau_1)\}}
\cdot\tau_1^{1/2}\nonumber\\
\label{12.560}
\end{eqnarray}
Choosing $c_0$ large enough so that in regard to the coefficient of $\s^{(1,n)}V(\tau_1)$ on 
the right we have:
\begin{equation}
\frac{C}{\sqrt{2c_0}}\leq\frac{1}{2}
\label{12.561}
\end{equation}
the coefficient of $\s^{(1,n)}\Ga(\tau_1)$ on the right is a fortiori not greater than $1/2$ and 
the inequality \ref{12.560} then implies:
\begin{eqnarray}
&&\s^{(1,n)}V(\tau_1)+\s^{(1,n)}\Ga(\tau_1)\leq \frac{C\left(\sqrt{\s^{(Y;0,n)}\cA(\tau_1)}+\sqrt{\s^{(E;0,n)}\cA(\tau_1)}
\right)}{\sqrt{2c_0}}\nonumber\\
&&\hspace{40mm}+\frac{C\left(\sqrt{\s^{(Y;0,n)}\cBb(\tau_1,\tau_1)}
+\sqrt{\s^{(E;0,n)}\cBb(\tau_1,\tau_1)}\right)}{a_0+2}\nonumber\\
&&\hspace{40mm}+C\sqrt{\max\{\s^{(0,n)}\cB(\tau_1,\tau_1),\s^{(0,n)}\cBb(\tau_1,\tau_1)\}}
\cdot\tau_1^{1/2}\nonumber\\
\label{12.562}
\end{eqnarray}
(for new constants $C$). Substituting this in \ref{12.493} we arrive the 
following proposition.

\vspace{2.5mm}

\noindent{\bf Proposition 12.3} \ \ Recalling the definitions \ref{12.394} - \ref{12.396}, the 
transformation function differences satisfy to principal terms the following estimates: 
\begin{eqnarray*}
&&\s^{(1,n)}P(\tau_1),\s^{(1,n)}V(\tau_1),\s^{(1,n)}\Ga(\tau_1) \leq 
\frac{C\left(\sqrt{\s^{(Y;0,n)}\cA(\tau_1)}+\sqrt{\s^{(E;0,n)}\cA(\tau_1)}
\right)}{\sqrt{2c_0}}\\
&&\hspace{45mm}+\frac{C\left(\sqrt{\s^{(Y;0,n)}\cBb(\tau_1,\tau_1)}
+\sqrt{\s^{(E;0,n)}\cBb(\tau_1,\tau_1)}\right)}{a_0+2}\\
&&\hspace{40mm}+C\sqrt{\max\{\s^{(0,n)}\cB(\tau_1,\tau_1),\s^{(0,n)}\cBb(\tau_1,\tau_1)\}}
\cdot\tau_1^{1/2}
\end{eqnarray*}
for all $\tau_1\in(0,\delta]$. 

\vspace{5mm}

\section{Estimates for $(\s^{(m,n-m)}\cla,\s^{(m,n-m)}\clab)$ for $m=1,...,n$}

In the next section we shall derive estimates for 
$(\Omega^{n-m}T^m\chf,\Omega^{n-m}T^m\cv,\Omega^{n-m}T^m\cga)$ for $m=1,...,n$. As we shall see, 
the derivation of these estimates necessarily involves the derivation of appropriate estimates on 
${\cal K}$ for $\s^{(m,n-m)}\clab$, $: m=1,...,n$, which in turn involve the derivation of appropriate 
estimates for $\s^{(m,n-m)}\cla$, $: m=1,...,n$. 

We begin with the derivation of the estimates for $\s^{(m,n-m)}\cla$. We consider the propagation 
equation for $\s^{(m,n-m)}\cla$, equation \ref{12.116}. This equation takes the form:
\begin{equation}
L\s^{(m,n-m)}\cla+(n-m)\chi\s^{(m,n-m)}\cla=\s^{(m,n-m)}\chG
\label{12.563}
\end{equation}
Proceeding as in the derivation of the inequality \ref{12.196} from equation \ref{12.189}, we deduce 
from equation \ref{12.563} the inequality:
\begin{equation}
\|\s^{(m,n-m)}\cla\|_{L^2(S_{\ub_1,u})}\leq k\int_0^{\ub_1}\|\s^{(m,n-m)}\chG\|_{L^2(S_{\ub,u})}d\ub
\label{12.564}
\end{equation}
where $k$ is a constant greater than 1, but which can be chosen as close to 1 as we wish by suitably 
restricting $\ub_1$. 

We shall presently estimate the leading contributions to the integral on the right in \ref{12.564}. 
Consider first the $n$th order acoustical difference terms:
\begin{eqnarray}
&&p\s^{(m,n-m)}\cla+\tilde{q}\s^{(m,n-m)}\clab-m\zeta\s^{(m-1,n-m+1)}\cla\nonumber\\
&&+2B\cmu_{m-1,n-m}+2\oB\cmub_{m-1,n-m} \label{12.565}
\end{eqnarray}
In view of \ref{10.575} and \ref{12.198} the contribution of the terms \ref{12.565} is bounded by:
\begin{eqnarray}
&&C\int_0^{\ub_1}\ub\|\s^{(m,n-m)}\cla\|_{L^2(S_{\ub,u})}d\ub
+C\int_0^{\ub_1}\|\s^{(m,n-m)}\clab\|_{L^2(S_{\ub,u})}d\ub\nonumber\\
&&+Cm\int_0^{\ub_1}\ub\|\s^{(m-1,n-m+1)}\cla\|_{L^2(S_{\ub,u})}d\ub\nonumber\\
&&+C\int_0^{\ub_1}\ub\|\cmu_{m-1,n-m}\|_{L^2(S_{\ub,u})}d\ub
+C\int_0^{\ub_1}\ub\|\cmub_{m-1,n-m}\|_{L^2(S_{\ub,u})}d\ub\nonumber\\ 
&&\label{12.566}
\end{eqnarray}
Here the first and third can obviously be absorbed in the integral inequality \ref{12.564}. 
We shall see that the second gives a borderline contribution. As for the last two of \ref{12.566}, 
the principal part of $\cmu_{m-1,n-m}$ and of $\cmub_{m-1,n-m}$ consists of the first two terms on the 
right in each of \ref{12.118} and \ref{12.119}. Thus to principal terms 
$\|\cmu_{m-1,n-m}\|_{L^2(S_{\ub,u})}$ is bounded by:
\begin{equation}
\sum_{j=0}^{m-1}\|(2\pi\lambda)^j \s^{(m-1-j,n-m+j+1)}\cla\|_{L^2(S_{\ub,u})}
+\frac{1}{2}\|(2\pi\lambda)^m \s^{(n-1)}\ctchi\|_{L^2(S_{\ub,u})}
\label{12.567}
\end{equation}
and $\|\cmub_{m-1,n-m}\|_{L^2(S_{\ub,u})}$ is bounded by:
\begin{equation}
\sum_{j=0}^{m-1}\|(2\pi\lambdab)^j \s^{(m-1-j,n-m+j+1)}\clab\|_{L^2(S_{\ub,u})}
+\frac{1}{2}\|(2\pi\lambdab)^m \s^{(n-1)}\ctchib\|_{L^2(S_{\ub,u})}
\label{12.568}
\end{equation}
Thus, with $\Pi$ being as in \ref{10.579} a fixed positive bound for $|2\pi\lambda/u^2|$ in 
${\cal R}_{\delta,\delta}$, the fourth of \ref{12.566} is to principal terms bounded by:
\begin{eqnarray}
&&C\sum_{j=0}^{m-1}(\Pi u^2)^j \int_0^{\ub_1}\ub\|\s^{(m-1-j,n-m+j+1)}\cla\|_{L^2(S_{\ub,u})}d\ub
\nonumber\\
&&+\frac{C}{2}(\Pi u^2)^m\int_0^{\ub_1}\ub\|\s^{(n-1)}\ctchi\|_{L^2(S_{\ub,u})}d\ub
\label{12.569}
\end{eqnarray}
Also, with $\Pib$ being as in \ref{10.741} a fixed positive bound for $|2\pi\lambdab/\ub|$ in 
${\cal R}_{\delta,\delta}$, the fifth of \ref{12.566} is to principal terms bounded by:
\begin{eqnarray}
&&C\sum_{j=0}^{m-1}\Pib^j \int_0^{\ub_1}\ub^{j+1}\|\s^{(m-1-j,n-m+j+1)}\clab\|_{L^2(S_{\ub,u})}d\ub
\nonumber\\
&&+\frac{C}{2}\Pib^m\int_0^{\ub_1}\ub^{m+1}\|\s^{(n-1)}\ctchib\|_{L^2(S_{\ub,u})}d\ub
\label{12.570}
\end{eqnarray}
The contributions of \ref{12.569} and \ref{12.570} will be seen to be depressed by a factor of 
at least $\ub_1$ relative to the leading contributions. However since \ref{12.569} and \ref{12.570} 
depend on $(\s^{(i,n-i)}\cla,\s^{(i,n-i)}\clab)$ for $i=2,...,m-1$, their estimates must proceed 
by induction. 

Besides the contribution of the terms \ref{12.565}, the other leading contributions to the integral 
on the right in \ref{12.564} are those of the principal difference terms:
\begin{eqnarray}
&&\frac{\rho}{4}\beta_N^2 E^{n-m}T^m\Lb H-\frac{\rho_N}{4}\beta_{N,N}^2 E_N^{n-m}T^m\Lb_N H_N
\nonumber\\
&&+\rhob\beta_N\left(c\pi\sbeta-\frac{\beta_N}{4}-\frac{\beta_{\Nb}}{2}\right)E^{n-m}T^m LH
\nonumber\\
&&\hspace{10mm}-\rhob_N\beta_{N,N}\left(c_N\pi_N\sbeta_N-\frac{\beta_{N,N}}{4}
-\frac{\beta_{\Nb,N}}{2}\right)E_N^{n-m}T^m L_N H_N \nonumber\\
&&+\rhob H\left(c\pi\sbeta-\frac{(\beta_N+\beta_{\Nb})}{2}\right)
N^\mu E^{n-m}T^m L\beta_\mu\nonumber\\
&&\hspace{10mm}-\rhob_N H_N\left(c_N\pi_N\sbeta_N-\frac{(\beta_{N,N}+\beta_{\Nb,N})}{2}\right)
N_N^\mu E_N^{n-m}T^m L_N\beta_{\mu,N}\nonumber\\
&&-\frac{\pi}{2}a\beta_N^2 E^{n-m}T^m EH+\frac{\pi_N}{2}a_N\beta_{N,N}^2 E_N^{n-m}T^m 
E_N H_N 
\label{12.571}
\end{eqnarray}
Here we must keep track not only of the actually principal terms, namely those containing derivatives 
of the $\beta_\mu$ of order $n+1$, but also of the terms containing acoustical quantities of order $n$. 
The terms of order $n$ with vanishing principal acoustical part can be ignored. As we shall see
the dominant contribution comes from the 1st of the four difference terms in \ref{12.571}. In fact, 
the actual principal difference term in question is:
\begin{equation}
-\frac{\rho}{2}\beta_N^2 H^\prime (g^{-1})^{\mu\nu}\beta_\nu E^{n-m}T^m\Lb\beta_\mu 
+\frac{\rho_N}{2}\beta_{N,N}^2 H^\prime_N (g^{-1})^{\mu\nu}\beta_{\nu,N}E_N^{n-m} T^m\Lb_N\beta_{\mu,N}
\label{12.572}
\end{equation}
By the formula \ref{10.368} and the analogous formula for the $N$th approximants, \ref{12.572} is 
expressed as the sum of differences of $E$, $N$, and $\Nb$ components:
\begin{eqnarray}
&&-\frac{1}{2}\beta_N^2\sbeta\eta^2 H^\prime\rho E^\mu E^{n-m}T^m\Lb\beta_\mu 
+\frac{1}{2}\beta_{N,N}^2\sbeta_N\eta_N^2 H_N^\prime\rho_N E_N^\mu E_N^{n-m}T^m\Lb_N\beta_{\mu,N}\nonumber\\
&&+\frac{1}{2}\beta_N^2\frac{\beta_{\Nb}}{2c}\eta^2 H^\prime\rho N^\mu E^{n-m}T^m\Lb\beta_\mu 
-\frac{1}{2}\beta_{N,N}^2\frac{\beta_{\Nb,N}}{2c_N}\eta_N^2 H_N^\prime\rho_N N_N^\mu E_N^{n-m}T^m\Lb_N\beta_{\mu,N}\nonumber\\
&&+\frac{1}{2}\beta_N^2\frac{\beta_N}{2c}\eta^2 H^\prime\rho\Nb^\mu E^{n-m}T^m\Lb\beta_\mu 
-\frac{1}{2}\beta_{N,N}^2\frac{\beta_{N,N}}{2c_N}\eta_N^2 H_N^\prime\rho_N\Nb_N^\mu E_N^{n-m}T^m\Lb_N\beta_{\mu,N} \nonumber\\
&&\label{12.573}
\end{eqnarray}
The dominant contribution is that of the $\Nb$ component difference. To estimate it we must use the 
fact that 
\begin{eqnarray}
&&N^\mu\Lb(E^{n-m}T^m\beta_\mu-E_N^{n-m}T^m\beta_{\mu,N})
+\ogamma\Nb^\mu\Lb(E^{n-m}T^m\beta_\mu-E_N^{n-m}T^m\beta_{\mu,N})\nonumber\\
&&\hspace{40mm}=\s^{(Y;m,n-m)}\cxi_{\Lb}
\label{12.574}
\end{eqnarray}
Since the 1st term on the left is related to the $N$ component difference, the second of the 
differences \ref{12.573}, we begin with that difference. To estimate its contribution we remark 
that in view of the fact that $N^\mu\Lb\beta_\mu=\lambda E^\mu E\beta_\mu$, 
$$N^\mu E^{n-m}T^m\Lb\beta_\mu$$
is, up to terms which can be ignored, equal to:
\begin{equation}
\lambda E^\mu E^{n-m}T^m E\beta_\mu+(E\beta_\mu)E^{n-m}T^m(\lambda E^\mu)
-(\Lb\beta_\mu)E^{n-m}T^m N^\mu
\label{12.575}
\end{equation}
In regard to the 1st of \ref{12.575}, by the first two of the commutation relations \ref{3.a14} we 
have:
\begin{equation}
[T^m, E]\beta_\mu=-\sum_{j=0}^{m-1}T^{m-1-j}\left((\chi+\chib)ET^j\beta_\mu\right)
\label{12.576}
\end{equation}
This gives:
\begin{equation}
\left[[T^m,E]\beta_\mu\right]_{P.A.}=-(E\beta_\mu)
\left(\rho[T^{m-1}\tchi]_{P.A.}+\rhob[T^{m-1}\tchib]_{P.A.}\right)
\label{12.577}
\end{equation}
and:
\begin{equation}
\left[E^{n-m}[T^m,E]\beta_\mu\right]_{P.A.}=-(E\beta_\mu)
\left(\rho[E^{n-m}T^{m-1}\tchi]_{P.A.}+\rhob[E^{n-m}T^{m-1}\tchib]_{P.A.}\right)
\label{12.578}
\end{equation}
Now from \ref{12.64}, \ref{12.105} - \ref{12.108}, \ref{12.110} - \ref{12.111}:
\begin{equation}
[E^{n-m}T^{m-1}\tchi]_{P.A.}=\left\{ \begin{array}{lll}
E^{n-1}\tchi &:& \mbox{if $m=1$}\\
2\mu_{m-2,n-m+1} &:& \mbox{if $m\geq 2$}
\end{array} \right.
\label{12.579}
\end{equation}
\begin{equation}
[E^{n-m}T^{m-1}\tchib]_{P.A.}=\left\{ \begin{array}{lll} 
E^{n-1}\tchib &:& \mbox{if $m=1$}\\
2\mub_{m-2,n-m+1} &:& \mbox{if $m\geq 2$}
\end{array} \right.
\label{12.580}
\end{equation}
It follows that up to terms which can be ignored the 1st of \ref{12.575} is equal to:
\begin{eqnarray}
&&\lambda E^\mu E^{n+1-m}T^m\beta_\mu \nonumber\\
&&-\lambda\sss\rho\cdot\left\{ \begin{array}{lll}
E^{n-1}\tchi &:& \mbox{if $m=1$}\\
2\mu_{m-2,n-m+1} &:& \mbox{if $m\geq 2$}
\end{array} \right.\nonumber\\
&&-\lambda\sss\rhob\cdot\left\{ \begin{array}{lll}
E^{n-1}\tchib &:& \mbox{if $m=1$}\\
2\mub_{m-2,n-m+1} &:& \mbox{if $m\geq 2$}
\end{array} \right. \label{12.581}
\end{eqnarray}
In regard to the 2nd of \ref{12.575}, from \ref{12.74}, \ref{12.107} - \ref{12.108}, 
\ref{12.110} - \ref{12.111}: 
\begin{equation}
[E^{n-m}T^m E^\mu]_{P.A.}=c^{-1}N^\mu\mub_{m-1,n-m}+c^{-1}\Nb^\mu\mu_{m-1,n-m} 
\label{12.582}
\end{equation}
Hence up to terms which can be ignored the 2nd of \ref{12.575} is equal to:
\begin{equation}
\sss E^{n-m}T^m\lambda+c^{-1}\lambda(\ss_N\mub_{m-1,n-m}+\ss_{\Nb}\mu_{m-1,n-m})
\label{12.583}
\end{equation}
In regard to the 3rd of \ref{12.575}, from the 1st of \ref{12.66} and from \ref{12.107}, \ref{12.110}:
\begin{equation}
[E^{n-m}T^m N^\mu]_{P.A.}=2(E^\mu-\pi N^\mu)\mu_{m-1,n-m}
\label{12.584}
\end{equation}
Hence up to terms which can be ignored the 3rd of \ref{12.575} is equal to:
\begin{equation}
2\lambda(\pi\sss-c^{-1}\ss_{\Nb})\mu_{m-1,n-m}
\label{12.585}
\end{equation}
Similar formulas hold for the $N$th approximants. It follows that up to terms which can be ignored 
the second of the differences \ref{12.573} (the $N$ component difference) is given by:
\begin{eqnarray}
&&\frac{1}{2}\beta_N^2\frac{\beta_{\Nb}}{2c}\eta^2 H^\prime\rho\left[ \lambda \s^{(E;m,n-m)}\csxi 
+\sss\s^{(m,n-m)}\cla\right.\nonumber\\
&&\hspace{27mm}+\rhob\left(\ss_N\cmub_{m-1,n-m}+(c\pi\sss-\ss_{\Nb})\cmu_{m-1,n-m}\right)
\nonumber\\
&&\hspace{27mm}-\lambda\sss\rho\cdot\left\{ \begin{array}{lll}
\s^{(n-1)}\ctchi &:& \mbox{if $m=1$}\\
2\cmu_{m-2,n-m+1} &:& \mbox{if $m\geq 2$}
\end{array} \right.\nonumber\\
&&\hspace{27mm}\left.-\lambda\sss\rhob\cdot\left\{ \begin{array}{lll}
\s^{(n-1)}\ctchib &:& \mbox{if $m=1$}\\
2\cmub_{m-2,n-m+1} &:& \mbox{if $m\geq 2$}
\end{array} \right. \hspace{2mm}\right]
\label{12.586}
\end{eqnarray}
The contribution of this to the integral on the right in \ref{12.564} is bounded by:
\begin{eqnarray}
&&Cu\int_0^{\ub_1}\sqrt{\ub}\|\sqrt{a}\s^{(E;m,n-m)}\csxi\|_{L^2(S_{\ub,u})}d\ub 
+C\int_0^{\ub_1}\ub\|\s^{(m,n-m)}\cla\|_{L^2(S_{\ub,u})}d\ub\nonumber\\
&&+Cu^2\int_0^{\ub_1}\ub\|\cmub_{m-1,n-m}\|_{L^2(S_{\ub,u})}d\ub
+Cu^2\int_0^{\ub_1}\ub\|\cmu_{m-1,n-m}\|_{L^2(S_{\ub,u})}d\ub\nonumber\\
&&+C\s^{(m,n-m)}I(\ub_1,u)\label{12.587}
\end{eqnarray}
where:
\begin{eqnarray}
&&\mbox{for $m=1$:} \hspace{3mm} \s^{(1,n-1)}I(\ub_1,u)=u^2
\int_0^{\ub_1}\ub^2\|\s^{(n-1)}\ctchi\|_{L^2(S_{\ub,u})}d\ub\nonumber\\
&&\hspace{48mm}+u^4\int_0^{\ub_1}\ub\|\s^{(n-1)}\ctchib\|_{L^2(S_{\ub,u})}d\ub\nonumber\\
&&\mbox{for $m\geq 2$:} \hspace{3mm} \s^{(m,n-m)}I(\ub_1,u)=u^2
\int_0^{\ub_1}\ub^2\|\cmu_{m-2,n-m+1}\|_{L^2(S_{\ub,u})}d\ub\nonumber\\
&&\hspace{48mm}+u^4\int_0^{\ub_1}\ub\|\cmub_{m-2,n-m+1}\|_{L^2(S_{\ub,u})}d\ub\nonumber\\
&&\label{12.588}
\end{eqnarray}
Now the 1st of \ref{12.587} is analogous to the 1st of \ref{12.206} with $\s^{(E;m,n-m)}\csxi$ in 
the role of $\s^{(E;0,n)}\csxi$. In analogy with \ref{12.209} its contribution to 
$\|\s^{(m,n-m)}\cla\|_{L^2(\Cb_{\ub_1}^{u_1})}$ is bounded by:
\begin{equation}
\frac{C\sqrt{\s^{(E;m,n-m)}\cBb(\ub_1,u_1)}}{\left(a_m+\frac{3}{2}\right)}\cdot
\ub_1^{a_m+\frac{3}{2}}u_1^{b_m+1}
\label{12.589}
\end{equation}
while in analogy  with \ref{12.212} its contribution to $\|\s^{(m,n-m)}\cla\|_{L^2({\cal K}^{\tau})}$ 
is bounded by:
\begin{equation}
\frac{C\sqrt{\s^{(E;m,n-m)}\cBb(\tau,\tau)}}{\left(a_m+\frac{3}{2}\right)}\cdot
\tau^{a_m+b_m+\frac{5}{2}}
\label{12.590}
\end{equation}
The 2nd of \ref{12.587} coincides with the 1st of \ref{12.566} and can be absorbed. The 3rd of 
\ref{12.587} is the same as the 5th of \ref{12.566} but with an extra $u^2$ factor. The 4th of 
\ref{12.587} is likewise the same as the 4th of \ref{12.566} but with an extra $u^2$ factor. 
As for $\s^{(m,n-m)}I(\ub_1,u)$, given by \ref{12.588}, consider first the case $m=1$. 
In regard to the 1st of the integrals in $\s^{(1,n-1)}I(\ub_1,u)$, setting 
\begin{equation}
f(\ub,u)=\ub^2\|\s^{(n-1)}\ctchi\|_{L^2(S_{\ub,u})}
\label{12.591}
\end{equation}
this integral is:
\begin{equation}
\int_0^{\ub_1}f(\ub,u)d\ub
\label{12.592}
\end{equation}
Since by the definition \ref{12.411}:
\begin{equation}
\|f(\ub,\cdot)\|_{L^2(\ub,u_1)}=\ub^2\|\s^{(n-1)}\ctchi\|_{L^2(\Cb_{\ub}^{u_1})}
\leq \ub^{a_0+2}u_1^{b_0}\s^{(n-1)}X(\ub_1,u_1)
\label{12.593}
\end{equation}
we have:
\begin{equation}
\int_0^{\ub_1}\|f(\ub,\cdot)\|_{L^2(\ub,u_1)}d\ub
\leq \frac{\ub_1^{a_0+3}u_1^{b_0}}{(a_0+3)}\s^{(n-1)}X(\ub_1,u_1)
\label{12.594}
\end{equation}
By Lemma 12.1a this bounds 
$$\left\|\int_0^{\ub_1}f(\ub,\cdot)d\ub\right\|_{L^2(\ub_1,u_1)}$$
It follows that the contribution of the 1st term in $\s^{(1,n-1)}I(\ub_1,u)$ to 
$\|\s^{(1,n-1)}\cla\|_{L^2(\Cb_{\ub_1}^{u_1})}$ is bounded by:
\begin{equation}
\frac{\ub_1^{a_0+3}u_1^{b_0+2}}{(a_0+3)}\s^{(n-1)}X(\ub_1,u_1)
\label{12.595}
\end{equation}
Also, setting $\ub_1=u$, $u_1=\tau$,
\begin{equation}
\int_0^{\tau}\|f(\ub,\cdot)\|_{L^2(\ub,\tau)}d\ub
\leq \frac{\tau^{a_0+b_0+3}}{(a_0+3)}\s^{(n-1)}X(\tau,\tau)
\label{12.596}
\end{equation}
By Lemma 12.1b this bounds $\|g\|_{L^2(0,\tau)}$ where
$$g(u)=\int_0^u f(\ub,u)d\ub$$
It follows that the contribution of the 1st term in $\s^{(1,n-1)}I(u,u)$ to 
$\|\s^{(1,n-1)}\cla\|_{L^2({\cal K}^{\tau})}$ is bounded by:
\begin{equation}
\frac{\tau^{a_0+b_0+5}}{(a_0+3)}\s^{(n-1)}X(\tau,\tau)
\label{12.597}
\end{equation}
In regard to the 2nd of the integrals in $\s^{(1,n-1)}I(\ub_1,u)$, in view of the definition 
\ref{12.460} we have:
\begin{eqnarray}
&&\int_0^{\ub_1}\ub\|\s^{(n-1)}\ctchib\|_{L^2(S_{\ub,u})}d\ub
\leq C\ub_1^{1/2}\|\rho\s^{(n-1)}\ctchib\|_{L^2(C_u^{\ub_1})}\nonumber\\
&&\hspace{40mm}\leq C\ub_1^{a_0+\frac{3}{2}}u^{b_0-1}\s^{(n-1)}\Xb(\ub_1,u_1) \label{12.598}
\end{eqnarray}
It follows that the contribution of the 2nd term in $\s^{(1,n-1)}I(\ub_1,u_1)$ to 
$\|\s^{(1,n-1)}\cla\|_{L^2(\Cb_{\ub_1}^{u_1})}$ is bounded by:
\begin{equation}
C\frac{\ub_1^{a_0+\frac{3}{2}}u_1^{b_0+\frac{7}{2}}}{\sqrt{2b_0+7}}\s^{(n-1)}\Xb(\ub_1,u_1)
\label{12.599}
\end{equation}
Also, setting $\ub_1=u$ in \ref{12.598} we deduce that the corresponding contribution to 
$\|\s^{(1,n-1)}\cla\|_{L^2({\cal K}^{\tau})}$ is bounded by:
\begin{equation}
C\frac{\tau^{c_0+5}}{\sqrt{2c_0+10}}\s^{(n-1)}\Xb(\tau,\tau)
\label{12.600}
\end{equation}
Consider now the case $m\geq 2$. The 1st integral in $\s^{(m,n-m)}I(\ub_1,u)$ is bounded by 
$\ub_1$ times the integral in the 4th term in \ref{12.566} with $m-1$ in the role of $m$. Therefore 
the 1st term in $\s^{(m,n-m)}I(\ub_1,u)$ is bounded by $u^2\ub_1$ times \ref{12.569} with $m-1$ 
in the role of $m$. Also, the 2nd integral in $\s^{(m,n-m)}I(\ub_1,u)$ is similar to the integral 
in the 5th term in \ref{12.566} with $m-1$ in the role of $m$. Therefore the 2nd term in 
$\s^{(m,n-m)}I(\ub_1,u)$ is bounded by $u^4$ times \ref{12.570} with $m-1$ in the role of $m$. 

We turn to the last of the differences \ref{12.573}, the $\Nb$ component difference. By the third 
of the commutation relations \ref{3.a14} we have:
\begin{equation}
[T^m,\Lb]\beta_\mu=-\sum_{j=0}^{m-1}T^{m-1-j}(\zeta ET^j\beta_\mu)
\label{12.601}
\end{equation}
This gives:
\begin{equation}
\left[[T^m,\Lb]\beta_\mu\right]_{P.A.}=-(E\beta_\mu)[T^{m-1}\zeta]_{P.A.}
\label{12.602}
\end{equation}
and:
\begin{equation}
E^{n-m}\left[[T^m,\Lb]\beta_\mu\right]_{P.A.}=-(E\beta_\mu)[E^{n-m}T^{m-1}\zeta]_{P.A.}
\label{12.603}
\end{equation}
Since $\zeta=2(\eta-\etab)$, we have, by \ref{3.a22},
$$[\zeta]_{P.A.}=2[\rho\teta-\rhob\tetab]_{P.A.}$$
hence by \ref{3.a25} and the definitions \ref{12.64}:
\begin{equation}
[\zeta]_{P.A.}=2(\rho\mu-\rhob\mub)
\label{12.604}
\end{equation}
It follows that:
\begin{equation}
E^{n-m}\left[[T^m,\Lb]\beta_\mu\right]_{P.A.}=-2(E\beta_\mu)\left\{\rho\mu_{m-1,n-m}
-\rhob\mub_{m-1,n-m}\right\}
\label{12.605}
\end{equation}
Hence 
$$\Nb^\mu E^{n-m}T^m\Lb\beta_\mu$$
is equal to:
\begin{equation}
\Nb^\mu\Lb E^{n-m}T^m\beta_\mu-2\ss_{\Nb}\left\{\rho\mu_{m-1,n-m}-\rhob\mub_{m-1,n-m}\right\}
\label{12.606}
\end{equation}
up to terms which can be ignored. Similarly for the $N$th approximants. It follows that up to terms 
which can be ignored the last of the differences \ref{12.573} is given by:
\begin{eqnarray}
&&\frac{\beta_N^3}{4c}\eta^2 H^\prime\rho\left\{\Nb^\mu\Lb(E^{n-m}T^m\beta_\mu-E_N^{n-m}T^m\beta_{\mu,N})\right.\nonumber\\
&&\hspace{25mm}\left.-2\ss_{\Nb}\left(\rho\cmu_{m-1,n-m}-\rhob\cmub_{m-1,n-m}\right)\right\}
\label{12.607}
\end{eqnarray}
The contribution of this to the integral on the right in \ref{12.564} is bounded by:
\begin{eqnarray}
&&C\int_0^{\ub_1}\ub\|\Nb^\mu\Lb(E^{n-m}T^m\beta_\mu-E_N^{n-m}T^m\beta_{\mu,N})\|_{L^2(S_{\ub,u})}d\ub
\nonumber\\
&&+C\int_0^{\ub_1}\ub^2\|\cmu_{m-1,n-m}\|_{L^2(S_{\ub,u})}d\ub
+Cu^2\int_0^{\ub_1}\ub\|\cmub_{m-1,n-m}\|_{L^2(S_{\ub,u})}d\ub\nonumber\\
&&\label{12.608}
\end{eqnarray}
Here the integrant of the 2nd integral is $\ub$ times the integrant of the 4th integral in \ref{12.566}, 
while the 3rd term is $u^2$ times the 5th term in \ref{12.566}. To estimate the 1st term we appeal to 
\ref{12.574}. Since $\ogamma\sim u$, this term is seen to be bounded by:
\begin{eqnarray}
&&Cu^{-1}\int_0^{\ub_1}\ub\|\s^{(Y;m,n-m)}\cxi_{\Lb}\|_{L^2(S_{\ub,u})}d\ub\nonumber\\
&&+Cu^{-1}\int_0^{\ub_1}\ub\|N^\mu\Lb(E^{n-m}T^m\beta_\mu-E_N^{n-m}T^m\beta_{\mu,N})\|_{L^2(S_{\ub,u})}
d\ub\nonumber\\
&&\label{12.609}
\end{eqnarray}
In regard to the 2nd term here we can express $N^\mu\Lb(E^{n-m}T^m\beta_\mu-E_N^{n-m}T^m\beta_{\mu,N})$ 
as
\begin{eqnarray*}
&&N^\mu E^{n-m}T^m\Lb\beta_\mu-N_N^\mu E_N^{n-m}T^m\beta_{\mu,N}\\
&&+2\ss_N(\rho\mu_{m-1,n-m}-\rhob\mu_{m-1,n-m})-2\ss_{N,N}(\rho_N\mu_{m-1,n-m,N}-\rhob_N
\mub_{m-1,n-m,N})
\end{eqnarray*}
up to terms which can be ignored. Hence the contribution of the 2nd term in \ref{12.609} is bounded in 
terms of:
\begin{eqnarray}
&&Cu^{-1}\int_0^{\ub_1}\ub\|N^\mu E^{n-m}T^m\Lb\beta_\mu
-N_N^\mu E_N^{n-m}T^m\Lb_N\beta_{\mu,N}\|_{L^2(S_{\ub,u})}d\ub\nonumber\\
&&+Cu^{-1}\int_0^{\ub_1}\ub^2\|\cmu_{m-1,n-m}\|_{L^2(S_{\ub,u})}d\ub
+Cu\int_0^{\ub_1}\ub\|\cmub_{m-1,n-m}\|_{L^2(S_{\ub,u})}d\ub\nonumber\\
&&\label{12.610}
\end{eqnarray}
Here the 3rd term is $u$ times the 5th term in \ref{12.566}, the 2nd term is dominated by the 4th term 
in \ref{12.566}, while the 1st is similar to the contribution of the $N$ component difference 
estimated above but with an extra $u^{-1}$ factor. It follows that, up to terms which can be absorbed, 
\ref{12.610} is bounded by:
$$C\int_0^{\ub_1}\sqrt{\ub}\|\sqrt{a}\s^{(E;m,n-m)}\csxi\|_{L^2(S_{\ub,u})}d\ub$$
the contribution of which to $\|\s^{(m,n-m)}\cla\|_{L^2(\Cb_{\ub_1}^{u_1})}$ is bounded by:
\begin{equation}
\frac{C\sqrt{\s^{(E;m,n-m)}\cBb(\ub_1,u_1)}}{\left(a_m+\frac{3}{2}\right)}\cdot
\ub_1^{a_m+\frac{3}{2}}u_1^{b_m}
\label{12.611}
\end{equation}
(compare with \ref{12.589}) and to $\|\s^{(m,n-m)}\cla\|_{L^2({\cal K}^{\tau})}$ by:
\begin{equation}
\frac{C\sqrt{\s^{(E;m,n-m)}\cBb(\tau,\tau)}}{\left(a_m+\frac{3}{2}\right)}\cdot
\tau^{a_m+b_m+\frac{3}{2}}
\label{12.612}
\end{equation}
(compare with \ref{12.590}). What is left is then the first of \ref{12.605} which is analogous to the 
first of \ref{12.217} with $\s^{(Y;m,n-m)}\cxi_{\Lb}$ in the role of $\s^{(Y;0,n)}\cxi_{\Lb}$. In 
analogy with \ref{12.223} its contribution to $\|\s^{(m,n-m)}\cla\|_{L^2(\Cb_{\ub_1}^{u_1})}$ is 
bounded by:
\begin{equation}
\frac{C\sqrt{\s^{(Y;m,n-m)}\cBb(\ub_1,u_1)}}{(a_m+2)}\cdot\ub_1^{a_m}u_1^{b_m}\cdot\frac{\ub_1^2}{u_1}
\label{12.613}
\end{equation}
while in analogy with \ref{12.225} its contribution to $\|\s^{(m,n-m)}\cla\|_{L^2({\cal K}^{\tau})}$ is 
bounded by:
\begin{equation}
\frac{C\sqrt{\s^{(Y;m,n-m)}\cBb(\tau,\tau)}}{(a_m+2)}\cdot\tau^{a_m+b_m+1}
\label{12.614}
\end{equation}

Consider finally the first of the differences \ref{12.573}, the $E$ component difference. Using 
\ref{12.605} we find that, up to terms which can be ignored, this is expressed as:
\begin{equation}
-\frac{1}{2}\beta_N^2\sbeta\eta^2 H^\prime\rho\left\{\s^{(E;m,n-m)}\cxi_{\Lb}
-2\sss\left(\rho\cmu_{m-1,n-m}-\rhob\cmub_{m-1,n-m}\right)\right\}
\label{12.615}
\end{equation}
The contribution of this to the integral on the right in \ref{12.564} is bounded by:
\begin{eqnarray}
&&C\int_0^{\ub_1}\ub\|\s^{(E;m,n-m)}\cxi_{\Lb}\|_{L^2(S_{\ub,u})}d\ub \label{12.616}\\
&&+C\int_0^{\ub_1}\ub^2\|\cmu_{m-1,n-m}\|_{L^2(S_{\ub,u})}d\ub
+Cu^2\int_0^{\ub_1}\ub\|\cmub_{m-1,n-m}\|_{L^2(S_{\ub,u})}d\ub \nonumber
\end{eqnarray}
Here the contribution of the last two terms can be absorbed, while the contribution of the first term 
to $\|\s^{(m,n-m)}\cla\|_{L^2(\Cb_{\ub_1}^{u_1})}$ is bounded by:
\begin{equation}
\frac{C\sqrt{\s^{(E;m,n-m)}\cBb(\ub_1,u_1)}}{(a_m+2)}\cdot\ub_1^{a_m}u_1^{b_m}\cdot\ub_1^2
\label{12.617}
\end{equation}
and to $\|\s^{(m,n-m)}\cla\|_{L^2({\cal K}^{\tau})}$ by:
\begin{equation}
\frac{C\sqrt{\s^{(E;m,n-m)}\cBb(\tau,\tau)}}{(a_m+2)}\cdot\tau^{a_m+b_m+2}
\label{12.618}
\end{equation}

We turn to the second of the principal difference terms \ref{12.571}. In fact the actual principal 
difference term in question is:
\begin{eqnarray}
&&-2\rhob\beta_N\left(c\pi\sbeta-\frac{\beta_N}{4}-\frac{\beta_{\Nb}}{2}\right)
H^\prime(g^{-1})^{\mu\nu}\beta_\nu E^{n-m}T^m L\beta_\mu
\nonumber\\
&&+2\rhob_N\beta_{N,N}\left(c_N\pi_N\sbeta_N-\frac{\beta_{N,N}}{4}-\frac{\beta_{\Nb,N}}{2}\right)H_N^\prime(g^{-1})^{\mu\nu}\beta_{\nu,N}E_N^{n-m}T^m L_N\beta_{\mu,N}\nonumber\\
&&\label{12.619}
\end{eqnarray}
This is expressed as the sum of differences of $E$, $N$, and $\Nb$ components:
\begin{eqnarray}
&&-2\rhob\sbeta\eta^2\beta_N\left(c\pi\sbeta-\frac{\beta_N}{4}-\frac{\beta_{\Nb}}{2}\right)
H^\prime E^\mu E^{n-m}T^m L\beta_\mu\nonumber\\
&&+2\rhob_N\sbeta_N\eta_N^2\beta_{N,N}\left(c_N\pi_N\sbeta_N-\frac{\beta_{N,N}}{4}-\frac{\beta_{\Nb,N}}{2}\right)H_N^\prime E_N^\mu E_N^{n-m}T^m L_N\beta_{\mu,N}\nonumber\\
&&+\rhob\frac{\beta_{\Nb}}{c}\eta^2\beta_N
\left(c\pi\sbeta-\frac{\beta_N}{4}-\frac{\beta_{\Nb}}{2}\right)
H^\prime N^\mu E^{n-m}T^m L\beta_\mu\nonumber\\
&&-\rhob_N\frac{\beta_{\Nb,N}}{c_N}\eta_N^2\beta_{N,N}\left(c_N\pi_N\sbeta_N-\frac{\beta_{N,N}}{4}-\frac{\beta_{\Nb,N}}{2}\right)H_N^\prime N_N^\mu E_N^{n-m}T^m L_N\beta_{\mu,N}\nonumber\\
&&+\rhob\frac{\beta_N}{c}\eta^2\beta_N
\left(c\pi\sbeta-\frac{\beta_N}{4}-\frac{\beta_{\Nb}}{2}\right)
H^\prime \Nb^\mu E^{n-m}T^m L\beta_\mu\nonumber\\
&&-\rhob_N\frac{\beta_{N,N}}{c_N}\eta_N^2\beta_{N,N}\left(c_N\pi_N\sbeta_N-\frac{\beta_{N,N}}{4}-\frac{\beta_{\Nb,N}}{2}\right)
H_N^\prime\Nb_N^\mu E_N^{n-m}T^m L_N\beta_{\mu,N}\nonumber\\
&&\label{12.620}
\end{eqnarray}
Consider first the $\Nb$ component difference. Using the conjugates of \ref{12.575} - \ref{12.585} 
we find that up to terms which can be ignored this is given by:
\begin{eqnarray}
&&\frac{\beta_N^2}{c}\eta^2\left(c\pi\sbeta-\frac{\beta_N}{4}-\frac{\beta_{\Nb}}{2}\right)
H^\prime\rhob\left[ \lambdab \s^{(E;m,n-m)}\csxi 
+\sss\s^{(m,n-m)}\clab\right.\nonumber\\
&&\hspace{27mm}+\rho\left(\ss_{\Nb}\cmu_{m-1,n-m}+(c\pi\sss-\ss_N)\cmub_{m-1,n-m}\right)
\nonumber\\
&&\hspace{27mm}-\lambdab\sss\rhob\cdot\left\{ \begin{array}{lll}
\s^{(n-1)}\ctchib &:& \mbox{if $m=1$}\\
2\cmub_{m-2,n-m+1} &:& \mbox{if $m\geq 2$}
\end{array} \right.\nonumber\\
&&\hspace{27mm}\left.-\lambdab\sss\rho\cdot\left\{ \begin{array}{lll}
\s^{(n-1)}\ctchi &:& \mbox{if $m=1$}\\
2\cmu_{m-2,n-m+1} &:& \mbox{if $m\geq 2$}
\end{array} \right. \hspace{2mm}\right]
\label{12.621}
\end{eqnarray}
(compare with \ref{12.586}). The contribution of the $\Nb$ component difference to the integral 
on the right in \ref{12.564} is then bounded by:
\begin{eqnarray}
&&Cu\int_0^{\ub_1}\sqrt{\ub}\|\sqrt{a}\s^{(E;m,n-m)}\csxi\|_{L^2(S_{\ub,u})}d\ub
+Cu^2\int_0^{\ub_1}\|\s^{(m,n-m)}\clab\|_{L^2(S_{\ub,u})}d\ub\nonumber\\
&&+Cu^2\int_0^{\ub_1}\ub\|\cmub_{m-1,n-m}\|_{L^2(S_{\ub,u})}d\ub
+Cu^2\int_0^{\ub_1}\ub\|\cmu_{m-1,n-m}\|_{L^2(S_{\ub,u})}d\ub\nonumber\\
&&+C\s^{(m,n-m)}I(\ub_1,u)\label{12.622}
\end{eqnarray}
Here the first as well as the last three terms coincide with the corresponding terms in \ref{12.587},  
while the second term can be absorbed into the second term in \ref{12.566}. 

In regard to the $E$ component difference (1st of the differences \ref{12.620}), by the 
conjugates of \ref{12.601} - \ref{12.605} 
$$E^\mu E^{n-m}T^m L\beta_\mu$$
is expressed as 
\begin{equation}
E^\mu LE^{n-m}T^m\beta_\mu+2\sss\left\{\rho\mu_{m-1,n-m}-\rhob\mub_{m-1,n-m}\right\}
\label{12.623}
\end{equation}
up to terms which can be ignored. Similarly for the $N$th approximants. It follows that up to terms 
which can be ignored the $E$ component difference is:
\begin{eqnarray}
&&-2\beta_N\sbeta\eta^2\left(c\pi\sbeta-\frac{\beta_N}{4}-\frac{\beta_{\Nb}}{2}\right)H^\prime\rhob
\left\{\s^{(E;m,n-m)}\cxi_L\right.\nonumber\\
&&\hspace{50mm}\left.+2\sss\left(\rho\cmu_{m-1,n-m}-\rhob\cmub_{m-1,n-m}\right)\right\}\nonumber\\
&&\label{12.624}
\end{eqnarray}
The contribution of this to the integral on the right in \ref{12.564} is bounded by:
\begin{eqnarray}
&&Cu^2\int_0^{\ub_1}\|\s^{(E;m,n-m)}\cxi_L\|_{L^2(S_{\ub,u})}d\ub \label{12.625}\\
&&+Cu^2\int_0^{\ub_1}\ub\|\cmu_{m-1,n-m}\|_{L^2(S_{\ub,u})}d\ub
+Cu^4\int_0^{\ub_1}\|\cmub_{m-1,n-m}\|_{L^2(S_{\ub,u})}d\ub \nonumber
\end{eqnarray}
Here the 1st term is bounded by:
\begin{equation}
C\ub_1^{a_m+\frac{1}{2}}u^{b_m+2}\sqrt{\s^{(E;m,n-m)}\cB(\ub_1,u_1)}
\label{12.626}
\end{equation}
hence its contribution to $\|\s^{(m,n-m)}\cla\|_{L^2(\Cb_{\ub_1}^{u_1})}$ is bounded by:
\begin{equation}
\frac{C\sqrt{\s^{(E;m,n-m)}\cB(\ub_1,u_1)}}{\sqrt{2b_m+5}}\cdot\ub_1^{a_m}u_1^{b_m}\cdot
\ub_1^{1/2}u_1^{5/2}
\label{12.627}
\end{equation}
and to $\|\s^{(m,n-m)}\cla\|_{L^2({\cal K}^{\tau})}$ by:
\begin{equation}
\frac{C\sqrt{\s^{(E;m,n-m)}\cB(\tau,\tau)}}{\sqrt{2b_m+5}}\cdot\tau^{a_m+b_m+3}
\label{12.628}
\end{equation}
The 2nd term in \ref{12.625} is $u^2$ times the 4th term in \ref{12.566}. As for the 3rd term, 
the integral is similar to that of the 5th term in \ref{12.566} but the $\ub$ factor in the integrant 
is here absent. Then in view of \ref{12.570} the integral in the 3rd term of \ref{12.625} 
is to principal terms bounded by:
\begin{eqnarray}
&&C\sum_{j=0}^{m-1}\Pib^j \int_0^{\ub_1}\ub^{j}\|\s^{(m-1-j,n-m+j+1)}\clab\|_{L^2(S_{\ub,u})}d\ub
\nonumber\\
&&+\frac{C}{2}\Pib^m\int_0^{\ub_1}\ub^{m}\|\s^{(n-1)}\ctchib\|_{L^2(S_{\ub,u})}d\ub
\label{12.629}
\end{eqnarray}

Finally, in regard to the $N$ component difference (2nd of the differences \ref{12.620}), by the 
conjugates of \ref{12.601} - \ref{12.605} 
$$N^\mu E^{n-m}T^m L\beta_\mu$$
is expressed as 
\begin{equation}
N^\mu LE^{n-m}T^m\beta_\mu+2\ss_N\left\{\rho\mu_{m-1,n-m}-\rhob\mub_{m-1,n-m}\right\}
\label{12.630}
\end{equation}
up to terms which can be ignored. Similarly for the $N$th approximants. It follows that up to terms 
which can be ignored the $N$ component difference is:
\begin{eqnarray}
&&\frac{\beta_N\beta_{\Nb}}{c}\eta^2\left(c\pi\sbeta-\frac{\beta_N}{4}-\frac{\beta_{\Nb}}{2}\right)
H^\prime\rhob
\left\{N^\mu L(E^{n-m}T^m\beta_\mu-E_N^{n-m}T^m\beta_{\mu,N})\right.\nonumber\\
&&\hspace{50mm}\left.+2\ss_N\left(\rho\cmu_{m-1,n-m}-\rhob\cmub_{m-1,n-m}\right)\right\}\nonumber\\
&&\label{12.631}
\end{eqnarray}
The contribution of this to the integral on the right in \ref{12.564} is bounded by:
\begin{eqnarray}
&&Cu^2\int_0^{\ub_1}\|N^\mu L(E^{n-m}T^m\beta_\mu-E_N^{n-m}T^m\beta_{\mu,N})\|_{L^2(S_{\ub,u})}d\ub\nonumber\\
&&+Cu^2\int_0^{\ub_1}\ub\|\cmu_{m-1,n-m}\|_{L^2(S_{\ub,u})}d\ub
+Cu^4\int_0^{\ub_1}\|\cmub_{m-1,n-m}\|_{L^2(S_{\ub,u})}d\ub \nonumber\\
&&\label{12.632}
\end{eqnarray}
Here the last two terms are the same as the last two terms in \ref{12.625}. To estimate the 
first term we use the fact that 
\begin{eqnarray}
&&N^\mu L(E^{n-m}T^m\beta_\mu-E_N^{n-m}T^m\beta_{\mu,N})
+\ogamma\Nb^\mu L(E^{n-m}T^m\beta_\mu-E_N^{n-m}T^m\beta_{\mu,N})\nonumber\\
&&\hspace{45mm}=\s^{(Y;m,n-m)}\cxi_L
\label{12.633}
\end{eqnarray}
Since 
\begin{equation}
\int_0^{\ub_1}\|\s^{(Y;m,n-m)}\cxi_L\|_{L^2(S_{\ub,u})}d\ub\leq C\ub_1^{a_m+\frac{1}{2}}u_1^{b_m}
\sqrt{\s^{(Y;m,n-m)}\cB(\ub_1,u_1)}
\label{12.634}
\end{equation}
the contribution of the right hand side of \ref{12.633} through the first term in \ref{12.632} to 
$\|\s^{(m,n-m)}\cla\|_{L^2(\Cb_{\ub_1}^{u_1})}$ is bounded by:
\begin{equation}
\frac{C\sqrt{\s^{(Y;m,n-m)}\cB(\ub_1,u_1)}}{\sqrt{2b_m+5}}\cdot\ub_1^{a_m}u_1^{b_m}\cdot
\ub_1^{1/2}u_1^{5/2}
\label{12.635}
\end{equation}
and the corresponding contribution to $\|\s^{(m,n-m)}\cla\|_{L^2({\cal K}^{\tau})}$ by:
\begin{equation}
\frac{C\sqrt{\s^{(Y;m,n-m)}\cB(\tau,\tau)}}{\sqrt{2b_m+5}}\cdot\tau^{a_m+b_m+3}
\label{12.636}
\end{equation}
In regard to the 2nd term on the left hand side of \ref{12.633}, using the conjugates of 
\ref{12.601}- \ref{12.605} we can express this term as
\begin{equation}
\Nb^\mu E^{n-m}T^m L\beta_\mu-\Nb_N^\mu E_N^{n-m}T^m L_N\beta_{\mu,N}
-2\ss_{\Nb}\left\{\rho\cmu_{m-1,n-m}-\rhob\cmub_{m-1,n-m}\right\}
\label{12.637}
\end{equation}
up to terms which can be ignored. The contribution of the principal difference to the first term 
in \ref{12.632} is similar to the contribution of the $\Nb$ component difference to the integral 
on the right in \ref{12.564}, but with an extra $u$ factor as $\ogamma\sim u$. As for the 
contribution of the remainder, this is similar to the last two terms in \ref{12.632} but with an 
extra $u$ factor. 

The third of the principal difference terms \ref{12.571} is similar to the $N$ component difference 
in \ref{12.620}.  

We finally consider the fourth of the principal difference terms \ref{12.571}. In fact the actual 
principal difference term in question is:
\begin{equation}
\pi a\beta_N^2 H^\prime(g^{-1})^{\mu\nu}\beta_\nu E^{n-m}T^m E\beta_\mu
-\pi_N a_N\beta_{N,N}^2 H^\prime_N (g^{-1})^{\mu\nu}\beta_{\nu,N}E_N^{n-m}T^m E_N\beta_{\mu,N}
\label{12.638}
\end{equation}
This is expressed as a sum of differences of $E$, $N$, and $\Nb$ components:
\begin{eqnarray}
&&a\pi\beta_N^2\sbeta\eta^2 H^\prime E^\mu E^{n-m}T^m E\beta_\mu \nonumber\\
&&\hspace{20mm}-a_N\pi_N\beta_{N,N}^2\beta_N\eta_N^2 H^\prime_N E_N^\mu E_N^{n-m}T^m E_N\beta_{\mu,N}\nonumber\\
&&-a\pi\beta_N^2\frac{\beta_{\Nb}}{2c}\eta^2 H^\prime N^\mu E^{n-m}T^m E\beta_\mu \nonumber\\
&&\hspace{20mm}+a_N\pi_N\beta_{N,N}^2\frac{\beta_{\Nb,N}}{2c_N}\eta_N^2 H^\prime_N N_N^\mu E_N^{n-m}T^m E_N\beta_{\mu,N}\nonumber\\
&&-a\pi\beta_N^2\frac{\beta_N}{2c}\eta^2 H^\prime\Nb^\mu E^{n-m}T^m E\beta_\mu \nonumber\\
&&\hspace{20mm}+a_N\pi_N\beta_{N,N}^2\frac{\beta_{N,N}}{2c_N}\eta_N^2 H^\prime_N\Nb_N^\mu E_N^{n-m}T^m E_N\beta_{\mu,N} \nonumber\\
&&\label{12.639}
\end{eqnarray}

Consider first the $\Nb$ component difference. Writing $\rhob\Nb^\mu=\Lb^\mu$ and noting that 
$\Lb^\mu E\beta_\mu=E^\mu\Lb\beta_\mu$, we can express 
$$\rhob\Nb^\mu E^{n-m}T^m E\beta_\mu$$
as 
\begin{equation}
E^{n-m}T^m(E^\mu\Lb\beta_\mu)-(E\beta_\mu)E^{n-m}T^m\Lb^\mu
\label{12.640}
\end{equation}
up to terms which can be ignored. Moreover, by \ref{12.66} and \ref{12.69}, 
\begin{equation}
[T\Lb^\mu]_{P.A.}=2\rhob\mub E^\mu+(c^{-1}T\lambda+2\pi\rhob\mu)\Nb^\mu 
\label{12.641}
\end{equation}
hence:
\begin{equation}
[E^{n-m}T^m\Lb^\mu]_{P.A.}=2\rhob\mub_{m-1,n-m}E^\mu+(c^{-1}E^{n-m}T^m\lambda
+2\pi\rhob\mu_{m-1,n-m})\Nb^\mu 
\label{12.642}
\end{equation}
Therefore up to terms which can be ignored the second term in \ref{12.640} is:
\begin{equation}
-2\sss\rhob\mub_{m-1,n-m}-\ss_{\Nb}(c^{-1}E^{n-m}T^m\lambda+2\pi\rhob\mu_{m-1,n-m})
\label{12.643}
\end{equation}
On the other hand, the first term in \ref{12.640} is:
\begin{equation}
E^\mu E^{n-m}T^m\Lb\beta_\mu+(\Lb\beta_\mu)[E^{n-m}T^m E^\mu]_{P.A.}
\label{12.644}
\end{equation}
up to terms which can be ignored. Here, by \ref{12.605} the 1st term is
\begin{equation}
E^\mu\Lb E^{n-m}T^m\beta_\mu-2\sss\left\{\rho\mu_{m-1,n-m}-\rhob\mub_{m-1,n-m}\right\}
\label{12.645}
\end{equation}
and by \ref{12.582} the 2nd term is
\begin{equation}
\sss\rhob\mub_{m-1,n-m}+c^{-1}s_{\Nb\Lb}\mu_{m-1,n-m}
\label{12.646}
\end{equation}
Similar formulas hold for the $N$th approximants. It follows that up to terms which can be ignored the 
$\Nb$ component difference is:
\begin{eqnarray}
&&-\frac{\pi\beta_N^3}{2}\eta^2 H^\prime\rho\left\{\s^{(E;m,n-m)}\cxi_{\Lb}
-c^{-1}\ss_{\Nb}\s^{(m,n-m)}\cla\right.\nonumber\\
&&\hspace{20mm}\left.+(c^{-1}s_{\Nb\Lb}-2\rho\sss-2\pi\rhob\ss_{\Nb})\cmu_{m-1,n-m}
+\rhob\sss\cmub_{m-1,n-m}\right\}\nonumber\\
&&\label{12.647}
\end{eqnarray}
The contribution of this to the integral on the right in \ref{12.564} is bounded by:
\begin{eqnarray}
&&C\int_0^{\ub_1}\ub\|\s^{(E;m,n-m)}\cxi_{\Lb}\|_{L^2(S_{\ub,u})}d\ub 
+C\int_0^{\ub_1}\ub\|\s^{(m,n-m)}\cla\|_{L^2(S_{\ub,u})}d\ub\nonumber\\
&&+C\int_0^{\ub_1}\ub\|\cmu_{m-1,n-m}\|_{L^2(S_{\ub,u})}d\ub
+Cu^2\int_0^{\ub_1}\ub\|\cmub_{m-1,n-m}\|_{L^2(S_{\ub,u})}d\ub\nonumber\\
&&\label{12.648}
\end{eqnarray}
Here the first is analogous to the first of \ref{12.244} with $\s^{(E;m,n-m)}\cxi_{\Lb}$ in the role of 
$\s^{(E;0,n)}\cxi_{\Lb}$. In analogy with \ref{12.245} its contribution to 
$\|\s^{(m,n-m)}\cla\|_{\Cb_{\ub_1}^{u_1})}$ is bounded by:
\begin{equation}
\frac{C\sqrt{\s^{(E;m,n-m)}\cBb(\ub_1,u_1)}}{(a_m+2)}\cdot\ub_1^{a_m}u_1^{b_m}\cdot\ub_1^2
\label{12.649}
\end{equation}
and to $\|\s^{(m,n-m)}\cla\|_{L^2({\cal K}^{\tau})}$ by:
\begin{equation}
\frac{C\sqrt{\s^{(E;m,n-m)}}\cBb(\tau,\tau)}{(a_m+2)}\cdot\tau^{a_m+b_m+2}
\label{12.650}
\end{equation}
The 3rd term in \ref{12.648} coincides with the 4th term in \ref{12.566} while the remaining terms 
can be absorbed.

Consider next the $E$ component difference (1st of \ref{12.639}). Using \ref{12.578} - \ref{12.580} 
we find that up to terms which can be ignored this difference is:
\begin{eqnarray}
&&\pi\beta_N^2\eta^2 H^\prime a\left[\s^{(E;m,n-m)}\csxi\right.
\nonumber\\
&&\hspace{24mm}-\sss\rho\cdot\left\{ \begin{array}{lll}
\s^{(n-1)}\ctchi &:& \mbox{if $m=1$}\\
2\cmu_{m-2,n-m+1} &:& \mbox{if $m\geq 2$}
\end{array} \right.\nonumber\\
&&\hspace{24mm}\left.-\sss\rhob\cdot\left\{ \begin{array}{lll}
\s^{(n-1)}\ctchib &:& \mbox{if $m=1$}\\
2\cmub_{m-2,n-m+1} &:& \mbox{if $m\geq 2$}
\end{array} \right. \hspace{2mm}\right]
\label{12.651}
\end{eqnarray}
The contribution of this to the integral on the right in \ref{12.564} is bounded by:
\begin{equation}
Cu\int_0^{\ub_1}\sqrt{\ub}\|\sqrt{a}\s^{(E;m,n-m)}\csxi\|_{L^2(S_{\ub,u})}d\ub 
+C\s^{(m,n-m)}I(\ub_1,u)
\label{12.652}
\end{equation}
(see \ref{12.588}). The two terms here coincide with the first and last terms in \ref{12.587}. 

Finally, apart from the factors $\beta_N^2$, $\beta_{N,N}^2$, the $N$ component difference (2nd of 
\ref{12.639}) is conjugate to the $\Nb$ component difference, hence comparing with \ref{12.647} the 
$N$ component difference is given up to terms which can be ignored by:
\begin{eqnarray}
&&-\frac{\pi\beta_N^2\beta_{\Nb}}{2}\eta^2 H^\prime\rhob\left\{\s^{(E;m,n-m)}\cxi_L
-c^{-1}\ss_N\s^{(m,n-m)}\clab\right.\nonumber\\
&&\hspace{20mm}\left.+(c^{-1}s_{NL}-2\rhob\sss-2\pi\rho\ss_N)\cmub_{m-1,n-m}
+\rho\sss\cmu_{m-1,n-m}\right\}\nonumber\\
&&\label{12.653}
\end{eqnarray}
The contribution of this to the integral on the right in \ref{12.564} is bounded by:
\begin{eqnarray}
&&Cu^2\int_0^{\ub_1}\|\s^{(E;m,n-m)}\cxi_L\|_{L^2(S_{\ub,u})}d\ub
+Cu^2\int_0^{\ub_1}\|\s^{(m,n-m)}\clab\|_{L^2(S_{\ub,u})}d\ub\nonumber\\
&&+Cu^3\int_0^{\ub_1}\|\cmub_{m-1,n-m}\|_{L^2(S_{\ub,u})}d\ub
+Cu^2\int_0^{\ub_1}\ub\|\cmu_{m-1,n-m}\|_{L^2(S_{\ub,u})}d\ub\nonumber\\
&&\label{12.654}
\end{eqnarray}
(note that by assumptions \ref{10.440} $|s_{NL}|\leq C\ub$ in ${\cal R}_{\delta,\delta}$). Here the 
1st term coincides with the 1st term of \ref{12.625}, the 2nd term coincides with the corresponding 
term of \ref{12.622}, the integral in the 3rd term coincides with the integral in the 3rd term of 
\ref{12.625} which is bounded by \ref{12.629}, and the 4th term coincides with the 2nd term of 
\ref{12.625}.

Returning to \ref{12.566}, in regard to the 1st term we apply Lemma 12.1 setting
$$f(\ub,u)=\ub\|\s^{(m,n-m)}\cla\|_{L^2(S_{\ub,u})}$$
Then, since 
$$\|f(\ub, \cdot)\|_{L^2(\ub,u_1)}=\ub\|\s^{(m,n-m)}\cla\|_{L^2(\Cb_{\ub}^{u_1})},$$
by Lemma 12.1a the contribution of this term to $\|\s^{(m,n-m)}\cla\|_{L^2(\Cb_{\ub_1}^{u_1})}$ 
is bounded by:
\begin{equation}
C\int_0^{\ub_1}\ub\|\s^{(m,n-m)}\cla\|_{L^2(\Cb_{\ub}^{u_1})}
\label{12.655}
\end{equation}
and by Lemma 12.1b the corresponding contribution to $\|\s^{(m,n-m)}\cla\|_{L^2({\cal K}^{\tau})}$ 
is bounded by:
\begin{equation}
C\int_0^{\tau}\ub\|\s^{(m,n-m)}\cla\|_{L^2(\Cb_{\ub}^{\tau})}d\ub
\label{12.656}
\end{equation}
Similarly in regard to the 3rd term in \ref{12.566} with $m-1$ in the role of $m$. In regard to 
the 4th term, we consider \ref{12.569} and we apply Lemma 12.1 setting 
$$f(\ub,u)=\left\{ \begin{array}{l} u^{2j}\ub\|\s^{(m-1-j,n-m+j+1)}\cla\|_{L^2(S_{\ub,u})} \ : \ 
j=0,...,m-1\\
u^{2m}\ub\|\s^{(n-1)}\ctchi\|_{L^2(S_{\ub,u})}\end{array}\right.$$
Then by Lemma 12.1a the contribution of the 4th term to $\|\s^{(m,n-m)}\cla\|_{L^2(\Cb_{\ub_1}^{u_1})}$ 
is bounded by:
\begin{eqnarray}
&&C\sum_{j=0}^{m-1}(\Pi u_1^2)^j\int_0^{\ub_1}\ub\|\s^{(m-1-j,n-m+j+1)}\cla\|_{L^2(\Cb_{\ub}^{u_1})}d\ub
\nonumber\\
&&+C(\Pi u_1^2)^m\int_0^{\ub_1}\ub\|\s^{(n-1)}\ctchi\|_{L^2(\Cb_{\ub}^{u_1})}d\ub
\label{12.657}
\end{eqnarray}
and by Lemma 12.1b the corresponding contribution to $\|\s^{(m,n-m)}\cla\|_{L^2({\cal K}^{\tau})}$ 
is bounded by:
\begin{eqnarray}
&&C\sum_{j=0}^{m-1}(\Pi\tau^2)^j\int_0^{\tau}\ub\|\s^{(m-1-j,n-m+j+1)}\cla\|_{L^2(\Cb_{\ub}^{\tau})}d\ub
\nonumber\\
&&+C(\Pi\tau^2)^m\int_0^{\tau}\ub\|\s^{(n-1)}\ctchi\|_{L^2(\Cb_{\ub}^{\tau})}d\ub
\label{12.658}
\end{eqnarray}
The 2nd term in \ref{12.566} is bounded by (Schwartz inequality):
\begin{equation}
C\ub_1^{1/2}\|\s^{(m,n-m)}\clab\|_{L^2(C_u^{\ub_1})}
\label{12.659}
\end{equation}
hence its contribution to $\|\s^{(m,n-m)}\cla\|_{L^2(\Cb_{\ub_1}^{u_1})}$ 
is bounded by:
\begin{equation}
C\ub_1^{1/2}\left(\int_{\ub_1}^{u_1}\|\s^{(m,n-m)}\clab\|^2_{L^2(C_u^{\ub_1})}du\right)^{1/2}
\label{12.660}
\end{equation}
and to $\|\s^{(m,n-m)}\cla\|_{L^2({\cal K}^{\tau})}$  by:
\begin{equation}
C\left(\int_0^{\tau}u\|\s^{(m,n-m)}\clab\|^2_{L^2(C_u^u)}du\right)^{1/2}
\label{12.661}
\end{equation}
In regard to the 5th term in \ref{12.566} and the 3rd term in \ref{12.625}, both of these are 
dominated by $u$ times \ref{12.629}, therefore by virtue of the Schwarz inequality by:
\begin{eqnarray}
&&Cu\left\{\sum_{j=0}^{m-1}\Pib^j\left(\frac{\ub_1^{2j+1}}{2j+1}\right)^{1/2}
\|\s^{(m-1-j,n-m+j+1)}\clab\|_{L^2(C_u^{\ub_1})}\right.\nonumber\\
&&\hspace{10mm}\left.+\Pib^m\left(\frac{\ub_1^{2m-1}}{2m-1}\right)^{1/2}
\|\rho\s^{(n-1)}\ctchib\|_{L^2(C_u^{\ub_1})}\right\}
\label{12.662}
\end{eqnarray}
Consequently the contribution of both terms to $\|\s^{(m,n-m)}\cla\|_{L^2(\Cb_{\ub_1}^{u_1})}$ is 
bounded by:
\begin{eqnarray}
&&C\left\{\sum_{j=0}^{m-1}\Pib^j\left(\frac{\ub_1^{2j+1}}{2j+1}\right)^{1/2}
\left(\int_{\ub_1}^{u_1}u^2\|\s^{(m-1-j,n-m+j+1)}\clab\|^2_{L^2(C_u^{\ub_1})}du\right)^{1/2}\right.\nonumber\\
&&\hspace{10mm}\left.+\Pib^m\left(\frac{\ub_1^{2m-1}}{2m-1}\right)^{1/2}
\left(\int_{\ub_1}^{u_1}u^2\|\rho\s^{(n-1)}\ctchib\|^2_{L^2(C_u^{\ub_1})}du\right)^{1/2}\right\}
\label{12.663}
\end{eqnarray}
and to $\|\s^{(m,n-m)}\cla\|_{L^2({\cal K}^{\tau})}$ by:
\begin{eqnarray}
&&C\left\{\sum_{j=0}^{m-1}\frac{\Pib^j}{\sqrt{2j+1}}\left(\int_0^{\tau}u^{2j+3}
\|\s^{(m-1-j,n-m+j)}\clab\|^2_{L^2(C_u^u)}du\right)^{1/2}\right.\nonumber\\
&&\hspace{10mm}\left.+\frac{\Pib^m}{\sqrt{2m-1}}\left(\int_0^{\tau}u^{2m+1}
\|\rho\s^{(n-1)}\ctchib\|^2_{L^2(C_u^u)}du\right)^{1/2}\right\}
\label{12.664}
\end{eqnarray}
Let us define, in analogy with \ref{12.412}, \ref{12.461}, 
\begin{eqnarray}
&&\s^{(m,n-m)}\La(\ub_1,u_1)=\sup_{(\ub,u)\in R_{\ub_1,u_1}}
\left\{\ub^{-a_m-\frac{1}{2}}u^{-b_m-\frac{1}{2}}\|\s^{(m,n-m)}\cla\|_{L^2(\Cb_{\ub}^u)}\right\}
\nonumber\\
&&\label{12.665}\\
&&\s^{(m,n-m)}\Lab(\ub_1,u_1)=\sup_{(\ub,u)\in R_{\ub_1,u_1}}
\left\{\ub^{-a_m}u^{-b_m}\|\s^{(m,n-m)}\clab\|_{L^2(C_u^{\ub})}\right\}\nonumber\\
&&\label{12.666}
\end{eqnarray}
We then have, in regard to \ref{12.657} (or \ref{12.658}), 
\begin{eqnarray}
&&\int_0^{\ub_1}\ub\|\s^{(m-1-j,n-m+j+1)}\cla\|_{L^2(\Cb_{\ub}^{u_1})}d\ub\nonumber\\
&&\hspace{10mm}\leq \frac{\ub_1^{a_{m-1-j}+\frac{5}{2}}u_1^{b_{m-1-j}+\frac{1}{2}}}
{\left(a_{m-1-j}+\frac{5}{2}\right)}\s^{(m-1-j,n-m+j+1)}\La(\ub_1,u_1)
\label{12.667}
\end{eqnarray}
Also, in view of the definition \ref{12.411}, 
\begin{equation}
\int_0^{\ub_1}\ub\|\s^{(n-1)}\ctchi\|_{L^2(\Cb_{\ub}^{u_1})}d\ub
\leq\frac{\ub_1^{a_0+2}u_1^{b_0}}{(a_0+2)}\s^{(n-1)}X(\ub_1,u_1)
\label{12.668}
\end{equation}
On the other hand, in regard to \ref{12.663}, 
\begin{eqnarray}
&&\left(\int_{\ub_1}^{u_1} u^2\|\s^{(m-1-j,n-m+j+1)}\clab\|^2_{L^2(C_u^{\ub_1})}du\right)^{1/2}
\nonumber\\
&&\hspace{10mm}\leq\frac{\ub_1^{a_{m-1-j}}u_1^{b_{m-1-j}+\frac{3}{2}}}{\sqrt{2b_{m-1-j}+3}}
\s^{(m-1-j,n-m+j+1)}\Lab(\ub_1,u_1)
\label{12.669}
\end{eqnarray}
and, in view of the definition \ref{12.460},
\begin{equation}
\left(\int_{\ub_1}^{u_1} u^2\|\rho\s^{(n-1)}\ctchib\|^2_{L^2(C_u^{\ub_1})}du\right)^{1/2}\leq\frac{\ub_1^{a_0+1}u^{b_0+\frac{1}{2}}}{\sqrt{2b_0+1}}\s^{(n-1)}\Xb(\ub_1,u_1)
\label{12.670}
\end{equation}
Also, in regard to \ref{12.664},
\begin{eqnarray}
&&\left(\int_0^{\tau}u^{2j+3}\|\s^{(m-1-j,n-m+j+1)}\clab\|^2_{L^2(C_u^u)}du\right)^{1/2}\nonumber\\
&&\hspace{10mm}\leq \frac{\tau^{c_{m-1-j}+j+2}}{\sqrt{2c_{m-1-j}+2j+4}}
\s^{(m-1-j,n-m+j+1)}\Lab(\tau,\tau)
\label{12.671}
\end{eqnarray}
and:
\begin{equation}
\left(\int_0^{\tau}u^{2m+1}\|\rho\s^{(n-1)}\ctchib\|^2_{L^2(C_u^u)}du\right)^{1/2}
\leq\frac{\tau^{c_0+m+1}}{\sqrt{2c_0+2m+2}}\s^{(n-1)}\Xb(\tau,\tau)
\label{12.672}
\end{equation}

Recalling that the exponents $a_m$, $b_m$ are non-increasing in $m$ we combine the above results 
in the form in which they will be used in the sequel in the following lemma. 

\vspace{2.5mm}

\noindent{\bf Lemma 12.9} \ \ The next to the top order acoustical quantities 
$\|\s^{(m,n-m)}\cla\|_{L^2(\Cb_{\ub}^u)} \ : \ m=1,...,n$ satisfy to principal terms the following inequalities:
\begin{eqnarray*}
&&\|\s^{(m,n-m)}\cla\|_{L^2(\Cb_{\ub_1}^{u_1})}\leq 
C\ub_1^{1/2}\left(\int_{\ub_1}^{u_1}\|\s^{(m,n-m)}\clab\|^2_{L^2(C_u^{\ub_1})}du\right)^{1/2}\\
&&\hspace{27mm}+C\int_0^{\ub_1}\ub\|\s^{(m,n-m)}\cla\|_{L^2(\Cb_{\ub}^{u_1})}\\
&&\hspace{27mm}+Cm\frac{\ub_1^{a_m+\frac{5}{2}}u_1^{b_m+\frac{1}{2}}}{\left(a_m+\frac{5}{2}\right)}
\s^{(m-1,n-m+1)}\La(\ub_1,u_1)\\
&&\hspace{24mm}+C\sum_{j=0}^{m-1}(\Pi u_1^2)^j \ \frac{\ub_1^{a_m+\frac{5}{2}}u_1^{b_m+\frac{1}{2}}}
{\left(a_m+\frac{5}{2}\right)}\s^{(m-1-j,n-m+j+1)}\La(\ub_1,u_1)\\
&&\hspace{27mm}+C(\Pi u_1^2)^m \ \frac{\ub_1^{a_m+2}u_1^{b_m}}{(a_m+2)}\s^{(n-1)}X(\ub_1,u_1)\\
&&\hspace{20mm}+C\sum_{j=0}^{m-1}\Pib^j\left(\frac{\ub_1^{2j+1}}{2j+1}\right)^{1/2}
\frac{\ub_1^{a_m}u_1^{b_m+\frac{3}{2}}}{\sqrt{2b_m+3}}\s^{(m-1-j,n-m+j+1)}\Lab(\ub_1,u_1)\\
&&\hspace{24mm}+C\Pib^m\left(\frac{\ub_1^{2m-1}}{2m-1}\right)^{1/2}
\frac{\ub_1^{a_m+1}u_1^{b_m+\frac{1}{2}}}{\sqrt{2b_m+1}}\s^{(n-1)}\Xb(\ub_1,u_1)\\
&&\hspace{27mm}+C\frac{\sqrt{\s^{(Y;m,n-m)}\cBb(\ub_1,u_1)}}{(a_m+2)}\cdot \ub_1^{a_m} u_1^{b_m}\cdot 
\frac{\ub_1^2}{u_1}\\
&&\hspace{20mm}+C\sqrt{\max\{\s^{(m,n-m)}\cB(\ub_1,u_1),\s^{(m,n-m)}\cBb(\ub_1,u_1)\}}\cdot 
\ub_1^{a_m}u_1^{b_m}\cdot\ub_1^{1/2} u_1
\end{eqnarray*}
Moreover, the $\|\s^{(m,n-m)}\cla\|_{L^2({\cal K}^{\tau})} \ : \ m=1,...,n$ satisfy:
\begin{eqnarray*}
&&\|\s^{(m,n-m)}\cla\|_{L^2({\cal K}^{\tau})}\leq 
C\left(\int_0^{\tau}u\|\s^{(m,n-m)}\clab\|^2_{L^2(C_u^u)}du\right)^{1/2}\\
&&\hspace{27mm}+C\int_0^{\tau}\ub\|\s^{(m,n-m)}\cla\|_{L^2(\Cb_{\ub}^{\tau})}d\ub\\
&&\hspace{27mm}+Cm\frac{\tau^{c_m+3}}{\left(a_m+\frac{5}{2}\right)}
\s^{(m-1,n-m+1)}\La(\tau,\tau)\\
&&\hspace{24mm}+C\sum_{j=0}^{m-1}(\Pi\tau^2)^j \ \frac{\tau^{c_m+3}}
{\left(a_m+\frac{5}{2}\right)}\s^{(m-1-j,n-m+j+1)}\La(\tau,\tau)\\
&&\hspace{27mm}+C(\Pi\tau^2)^m \ \frac{\tau^{c_m+2}}{(a_m+2)}\s^{(n-1)}X(\tau,\tau)\\
&&\hspace{20mm}+C\sum_{j=0}^{m-1}\frac{\Pib^j}{\sqrt{2j+1}}
\frac{\tau^{c_m+j+2}}{\sqrt{2c_m+2j+4}}\s^{(m-1-j,n-m+j+1)}\Lab(\tau,\tau)\\
&&\hspace{24mm}+C\frac{\Pib^m}{\sqrt{2m-1}}
\frac{\tau^{c_m+m+1}}{\sqrt{2c_m+2m+2}}\s^{(n-1)}\Xb(\tau,\tau)\\
&&\hspace{27mm}+C\frac{\sqrt{\s^{(Y;m,n-m)}\cBb(\tau,\tau)}}{(a_m+2)}\cdot \tau^{c_m+1}\\
&&\hspace{20mm}+C\sqrt{\max\{\s^{(m,n-m)}\cB(\tau,\tau),\s^{(m,n-m)}\cBb(\tau,\tau)\}}\cdot 
\tau^{c_m+\frac{3}{2}}
\end{eqnarray*}

\vspace{2.5mm}

We turn to the derivation of the estimates for $\s^{(m,n-m)}\clab \ : \ m=1,...,n$. We consider the 
propagation equation for $\s^{(m,n-m)}\clab$, equation \ref{12.120}. This equation takes the form:
\begin{equation}
\Lb\s^{(m,n-m)}\clab+(n-m)\chib\s^{(m,n-m)}\clab=\s^{(m,n-m)}\chGb
\label{12.673}
\end{equation}
where $\s^{(m,n-m)}\chGb$ is the conjugate of $\s^{(m,n-m)}\chG$. 
Proceeding as in the derivation of the inequalities \ref{12.306} and \ref{12.307} from equation 
\ref{12.300}, we deduce from equation \ref{12.673} the inequalities:
\begin{equation}
\|\s^{(m,n-m)}\clab\|_{L^2(S_{\ub,u_1})}\leq k\|\s^{(m,n-m)}\clab\|_{L^2(S_{\ub,\ub})}
+k\int_{\ub}^{u_1}\|\s^{(m,n-m)}\chGb\|_{L^2(S_{\ub,u})}du
\label{12.674}
\end{equation}
and:
\begin{eqnarray}
&&\|\s^{(m,n-m)}\clab\|_{L^2(C_{u_1}^{\ub_1})}\leq k\|\s^{(m,n-m)}\clab\|_{L^2({\cal K}^{\ub_1})} 
\label{12.675}\\
&&\hspace{27mm}+k\left\{\int_0^{\ub_1}\left(\int_{\ub}^{u_1}\|\s^{(m,n-m)}\chGb\|_{L^2(S_{\ub,u})}du\right)^2 d\ub\right\}^{1/2} \nonumber
\end{eqnarray}
where $k$ is a constant greater than 1, but which can be chosen as close to 1 as we wish by suitably 
restricting $u_1$. 

We shall presently estimate the leading contributions to 
\begin{equation}
\left\{\int_0^{\ub_1}\left(\int_{\ub}^{u_1}\|\s^{(m,n-m)}\chGb\|_{L^2(S_{\ub,u})}du\right)^2 
d\ub\right\}^{1/2}
\label{12.676}
\end{equation}
Consider first the $n$th order acoustical difference terms:
\begin{eqnarray}
&&\pb\s^{(m,n-m)}\clab+\tilde{\qb}\s^{(m,n-m)}\cla+m\zeta\s^{(m-1,n-m+1)}\clab\nonumber\\
&&+2\oBb\cmub_{m-1,n-m}+2\Bb\cmu_{m-1,n-m} \label{12.677}
\end{eqnarray}
In view of \ref{12.310}, \ref{12.311} the contribution of the terms \ref{12.677} to the integral on the 
right in \ref{12.674} is bounded by:
\begin{eqnarray}
&&C\int_{\ub}^{u_1}\|\s^{(m,n-m)}\clab\|_{L^2(S_{\ub,u})}du
+C\ub\int_{\ub}^{u_1}\|\s^{(m,n-m)}\cla\|_{L^2(S_{\ub,u})}du\nonumber\\
&&+Cm\ub\int_{\ub}^{u_1}\|\s^{(m-1,n-m+1)}\clab\|_{L^2(S_{\ub,u})}du\nonumber\\
&&+C\ub\int_{\ub}^{u_1}\|\cmub_{m-1,n-m}\|_{L^2(S_{\ub,u})}du
+C\ub\int_{\ub}^{u_1}\|\cmu_{m-1,n-m}\|_{L^2(S_{\ub,u})}du\nonumber\\
&&\label{12.678}
\end{eqnarray}
Here the first can obviously be absorbed in the integral inequality\ref{12.674} while the second is 
bounded by $C\ub\sqrt{u_1}\|\s^{(m,n-m)}\cla\|_{L^2(\Cb_{\ub}^{u_1})}$ hence its contribution to 
\ref{12.676} is bounded by:
\begin{equation}
Cu_1^{1/2}\left\{\int_0^{\ub_1}\ub^2\|\s^{(m,n-m)}\cla\|^2_{L^2(\Cb_{\ub}^{u_1})}d\ub\right\}^{1/2}
\label{12.679}
\end{equation}
Similarly the last of \ref{12.678} is bounded by 
$C\ub\sqrt{u_1}\|\cmu_{m-1,n-m}\|_{L^2(\Cb_{\ub}^{u_1})}$ hence its contribution to \ref{12.676} is 
bounded by:
\begin{equation}
Cu_1^{1/2}\left\{\int_0^{\ub_1}\ub^2\|\cmu_{m-1,n-m}\|^2_{L^2(\Cb_{\ub}^{u_1})}d\ub\right\}^{1/2}
\label{12.680}
\end{equation}
By \ref{12.567}, taking $L^2$ norms with respect to $u$ on $[\ub,u_1]$, 
$\|\cmu_{m-1,n-m}\|_{L^2(\Cb_{\ub}^{u_1})}$ is bounded to principal terms by:
\begin{equation}
\sum_{j=0}^{m-1}(\Pi u_1^2)^j \ \|\s^{(m-1-j,n-m+j+1)}\cla\|_{L^2(\Cb_{\ub}^{u_1})}
+\frac{1}{2}(\Pi u_1^2)^m \ \|\s^{(n-1)}\ctchi\|_{L^2(\Cb_{\ub}^{u_1})}
\label{12.681}
\end{equation}
which in view of the definitions \ref{12.665} and \ref{12.411} is bounded by:
\begin{eqnarray}
&&\sum_{j=0}^{m-1}(\Pi u_1^2)^j \ \ub^{a_{m-1-j}+\frac{1}{2}}u_1^{b_{m-1-j}+\frac{1}{2}}
\s^{(m-1-j,n-m+j+1)}\La(\ub_1,u_1)\nonumber\\
&&+\frac{1}{2}(\Pi u_1^2)^m \ \ub^{a_0} u_1^{b_0}\s^{(n-1)}X(\ub_1,u_1)
\label{12.682}
\end{eqnarray}
Substituting in \ref{12.680} yields a bound of the last by $C$ times:
\begin{eqnarray}
&&\sum_{j=0}^{m-1}(\Pi u_1^2)^j \ \frac{\ub_1^{a_{m-1-j}+2}u_1^{b_{m-j-1}+1}}{\sqrt{2a_{m-1-j}+4}}
\s^{(m-1-j,n-m+j+1)}\La(\ub_1,u_1)\nonumber\\
&&+\frac{1}{2}(\Pi u_1^2)^m \ \frac{\ub_1^{a_0+\frac{3}{2}}u_1^{b_0+\frac{1}{2}}}{\sqrt{2a_0+3}}
\s^{(n-1)}X(\ub_1,u_1)
\label{12.683}
\end{eqnarray}
To estimate the contribution of the third of \ref{12.678} to \ref{12.676} we apply Lemma 12.5 taking 
\begin{equation}
f(\ub,u)=\|\s^{(m-1,n-m+1)}\clab\|_{L^2(S_{\ub,u})}
\label{12.684}
\end{equation}
Then the third of \ref{12.678} is $Cm\ub g(\ub)$ and its contribution to \ref{12.676} is bounded by 
$Cm\ub_1\|g\|_{L^2(0,\ub_1)}$. The lemma then yields that this is bounded by $Cm\ub_1$ times:
\begin{equation}
\int_0^{\ub_1}\|\s^{(m-1,n-m+1)}\clab\|_{L^2(C_u^u)}du
+\int_{\ub_1}^{u_1}\|\s^{(m-1,n-m+1)}\clab\|_{L^2(C_u^{\ub_1})}du
\label{12.685}
\end{equation}
which in view of the definition \ref{12.666} is bounded by:
\begin{equation}
\frac{\ub_1^{a_{m-1}}u_1^{b_{m-1}+1}}{b_{m-1}+1}\s^{(m-1,n-m+1)}\Lab(\ub_1,u_1)
\label{12.686}
\end{equation}
Similarly, to estimate the contribution of the fourth of \ref{12.678} to \ref{12.676} we apply 
Lemma 12.5 taking 
\begin{equation}
f(\ub,u)=\|\cmub_{m-1,n-m}\|_{L^2(S_{\ub,u})}
\label{12.687}
\end{equation}
Then the fourth of \ref{12.678} is $C\ub g(\ub)$ and its contribution to \ref{12.676} is bounded by 
$C\ub_1\|g\|_{L^2(0,\ub_1)}$. The lemma then yields that this is bounded by $C\ub_1$ times:
\begin{equation}
\int_0^{\ub_1}\|\cmub_{m-1,n-m}\|_{L^2(C_u^u)}du
+\int_{\ub_1}^{u_1}\|\cmub_{m-1,n-m}\|_{L^2(C_u^{\ub_1})}du
\label{12.688}
\end{equation}
By \ref{12.678}, taking $L^2$ norms with respect to $\ub$ on $[0,\ub_1]$, 
$\|\cmub_{m-1,n-m}\|_{L^2(C_u^{\ub_1})}$ is bounded to principal terms by:
\begin{equation}
\sum_{j=1}^{m-1}\Pib^j \ub_1^j \|\s^{(m-1-j,n-m+j+1)}\clab\|_{L^2(C_u^{\ub_1})}
+C\Pib^m \ub_1^{m-1}\|\rho\s^{(n-1)}\ctchib\|_{L^2(C_u^{\ub_1})}
\label{12.689}
\end{equation}
which in view of the definitions \ref{12.666} and \ref{12.460} is bounded by:
\begin{eqnarray}
&&\sum_{j=0}^{m-1}\Pib^j \ub_1^{j+a_{m-1-j}}u^{b_{m-1-j}} \s^{(m-1-j,n-m+j+1)}\Lab(\ub_1,u_1)\nonumber\\
&&+C\Pib^m \ub_1^{m+a_0}u^{b_0-1} \s^{(n-1)}\Xb(\ub_1,u_1) 
\label{12.690}
\end{eqnarray}
It then follows that \ref{12.688} is bounded by:
\begin{eqnarray}
&&\sum_{j=0}^{m-1}\Pib^j \frac{\ub_1^{j+a_{m-1-j}}u_1^{b_{m-1-j}+1}}{b_{m-1-j}+1} 
\s^{(m-1-j,n-m+j+1)}\Lab(\ub_1,u_1)\nonumber\\
&&+C\Pib^m \frac{\ub_1^{m+a_0}u_1^{b_0}}{b_0} \s^{(n-1)}\Xb(\ub_1,u_1)
\label{12.691}
\end{eqnarray}

Besides the contribution of the terms \ref{12.678}, the other leading contributions are those of the 
principal difference terms: 
\begin{eqnarray}
&&\frac{\rhob}{4}\beta_{\Nb}^2 E^{n-m}T^m LH-\frac{\rhob_N}{4}\beta_{\Nb,N}^2 E_N^{n-m}T^m L_N H_N
\nonumber\\
&&+\rho\beta_{\Nb}\left(c\pi\sbeta-\frac{\beta_{\Nb}}{4}-\frac{\beta_N}{2}\right)E^{n-m}T^m\Lb H
\nonumber\\
&&\hspace{10mm}-\rho_N\beta_{\Nb,N}\left(c_N\pi_N\sbeta_N-\frac{\beta_{\Nb,N}}{4}
-\frac{\beta_{N,N}}{2}\right)E_N^{n-m}T^m\Lb_N H_N \nonumber\\
&&+\rho H\left(c\pi\sbeta-\frac{(\beta_N+\beta_{\Nb})}{2}\right)
\Nb^\mu E^{n-m}T^m\Lb\beta_\mu\nonumber\\
&&\hspace{10mm}-\rho_N H_N\left(c_N\pi_N\sbeta_N-\frac{(\beta_{N,N}+\beta_{\Nb,N})}{2}\right)
\Nb_N^\mu E_N^{n-m}T^m\Lb_N\beta_{\mu,N}\nonumber\\
&&-\frac{\pi}{2}a\beta_{\Nb}^2 E^{n-m}T^m EH+\frac{\pi_N}{2}a_N\beta_{\Nb,N}^2 E_N^{n-m}T^m 
E_N H_N 
\label{12.692}
\end{eqnarray}
(conjugate of \ref{12.571}). Here again we must keep track not only of the actually principal terms, namely those containing derivatives 
of the $\beta_\mu$ of order $n+1$, but also of the terms containing acoustical quantities of order $n$. 
The terms of order $n$ with vanishing principal acoustical part can be ignored. Consider the first of 
the above difference terms. The actual principal difference term in question is:
\begin{equation}
-\frac{\rhob}{2}\beta_{\Nb}^2 H^\prime(g^{-1})^{\mu\nu}\beta_\nu E^{n-m}T^m L\beta_\mu 
+\frac{\rhob_N}{2}\beta_{\Nb,N}^2 H^\prime_N(g^{-1})^{\mu\nu}\beta_{\nu,N} E_N^{n-m}T^m L_N\beta_{\mu,N}
\label{12.693}
\end{equation}
This is similar to the principal difference \ref{12.619}. 
It is similarly expressed as differences of $E$, $N$, and $\Nb$ components:
\begin{eqnarray}
&&-\frac{\rhob}{2}\beta_{\Nb}^2\sbeta\eta^2 H^\prime E^\mu E^{n-m}T^m L\beta_\mu 
+\frac{\rhob_N}{2}\beta_{\Nb,N}^2\sbeta_N\eta_N^2 H^\prime_N E_N^\mu E_N^{n-m}T^m L_N\beta_{\mu,N} \nonumber\\
&&+\frac{\rhob}{2}\beta_{\Nb}^2\frac{\beta_{\Nb}}{2c}\eta^2 H^\prime N^\mu E^{n-m}T^m L\beta_\mu 
-\frac{\rhob_N}{2}\beta_{\Nb,N}^2\frac{\beta_{\Nb,N}}{2c_N}\eta_N^2 H^\prime_N N_N^\mu E_N^{n-m}T^m L_N\beta_{\mu,N} \nonumber\\
&&+\frac{\rhob}{2}\beta_{\Nb}^2\frac{\beta_N}{2c}\eta^2 H^\prime\Nb^\mu E^{n-m}T^m L\beta_\mu 
-\frac{\rhob_N}{2}\beta_{\Nb,N}^2\frac{\beta_{N,N}}{2c_N}\eta_N^2 H^\prime_N\Nb_N^\mu E_N^{n-m}T^m L_N\beta_{\mu,N} \nonumber\\
&&\label{12.694}
\end{eqnarray}
Comparing with \ref{12.624}, the $E$ component difference here is:
\begin{equation}
-\frac{1}{2}\beta_{\Nb}^2\sbeta\eta^2 H^\prime\rhob
\left\{\s^{(E;m,n-m)}\cxi_L + 2\sss\left(\rho\cmu_{m-1,n-m}-\rhob\cmub_{m-1,n-m}\right)\right\}
\label{12.695}
\end{equation}
up to terms which can be ignored. 
The contribution of this to the integral on the right in \ref{12.674} is bounded by:
\begin{eqnarray}
&&C\int_{\ub}^{u_1}u^2\|\s^{(E;m,n-m)}\cxi_L\|_{L^2(S_{\ub,u})}du \label{12.696}\\
&&+C\ub\int_{\ub}^{u_1}u^2 \|\cmu_{m-1,n-m}\|_{L^2(S_{\ub,u})}du
+C\int_{\ub}^{u_1}u^4\|\cmub_{m-1,n-m}\|_{L^2(S_{\ub,u})}du \nonumber
\end{eqnarray}
To estimate the contribution of the 1st term to \ref{12.676} we apply Lemma 12.5 taking
\begin{equation}
f(\ub,u)=u^2\|\s^{(E;m,n-m)}\cxi_L\|_{L^2(S_{\ub,u})}
\label{12.697}
\end{equation}
Then this term is $Cg(\ub)$ and its contribution to \ref{12.676} is $C\|g\|_{L^2(0,\ub_1)}$. 
The lemma then yields that this is bounded by:
\begin{eqnarray}
&&C\int_0^{\ub_1}u^2\|\s^{(E;m,n-m)}\cxi_L\|_{L^2(C_u^u)}du
+C\int_{\ub_1}^{u_1} u^2\|\s^{(E;m,n-m)}\cxi_L\|_{L^2(C_u^{\ub_1})}du\nonumber\\
&&\leq C\int_0^{\ub_1}u^2\sqrt{\s^{(E;m,n-m)}\cE^u(u)}du
+C\int_{\ub_1}^{u_1}u^2\sqrt{\s^{(E;m,n-m)}\cE^{\ub_1}(u)}du\nonumber\\
&&\leq C\sqrt{\s^{(E;m,n-m)}\cB(\ub_1,\ub_1)}\int_0^{\ub_1}u^{a_m+b_m+2}du\nonumber\\
&&\hspace{30mm}+C\sqrt{\s^{(E;m,n-m)}\cB(\ub_1,u_1)}\cdot\ub_1^{a_m}\cdot\int_{\ub_1}^{u_1}u^{b_m+2}du\nonumber\\
&&\leq C\sqrt{\s^{(E;m,n-m)}\cB(\ub_1,u_1)}\cdot\left\{\frac{\ub_1^{a_m+b_m+3}}{a_m+b_m+3}
+\ub_1^{a_m}\cdot\frac{(u_1^{b_m+3}-\ub_1^{b_m+3})}{b_m+3}\right\}\nonumber\\
&&\leq C\frac{\sqrt{\s^{(E;m,n-m)}\cB(\ub_1,u_1)}}{b_m+3}\cdot\ub_1^{a_m}u_1^{b_m}\cdot u_1^3
\label{12.698}
\end{eqnarray}
(this is analogous to \ref{12.326}). The 2nd term in \ref{12.696} is dominated by $u_1^2$ times 
the last term in \ref{12.678}, while the contribution of the 3rd term in \ref{12.696} to \ref{12.676} 
is dominated by $u_1^4$ times \ref{12.688}. 

Comparing with \ref{12.621} the $\Nb$ component difference in \ref{12.694} is: 
\begin{eqnarray}
&&\frac{\beta_{\Nb}^2\beta_N}{4c}\eta^2
H^\prime\rhob\left[ \lambdab \s^{(E;m,n-m)}\csxi 
+\sss\s^{(m,n-m)}\clab\right.\nonumber\\
&&\hspace{27mm}+\rho\left(\ss_{\Nb}\cmu_{m-1,n-m}+(c\pi\sss-\ss_N)\cmub_{m-1,n-m}\right)
\nonumber\\
&&\hspace{27mm}-\lambdab\sss\rhob\cdot\left\{ \begin{array}{lll}
\s^{(n-1)}\ctchib &:& \mbox{if $m=1$}\\
2\cmub_{m-2,n-m+1} &:& \mbox{if $m\geq 2$}
\end{array} \right.\nonumber\\
&&\hspace{27mm}\left.-\lambdab\sss\rho\cdot\left\{ \begin{array}{lll}
\s^{(n-1)}\ctchi &:& \mbox{if $m=1$}\\
2\cmu_{m-2,n-m+1} &:& \mbox{if $m\geq 2$}
\end{array} \right. \hspace{2mm}\right]
\label{12.699}
\end{eqnarray}
The contribution of this to the integral on the right in \ref{12.674} is bounded by: 
\begin{eqnarray}
&&C\sqrt{\ub}\int_{\ub}^{u_1}u\|\sqrt{a}\s^{(E;m,n-m)}\csxi\|_{L^2(S_{\ub,u})}du
+C\int_{\ub}^{u_1}u^2\|\s^{(m,n-m)}\clab\|_{L^2(S_{\ub,u})}du\nonumber\\
&&+C\ub\int_{\ub}^{u_1}u^2\|\cmub_{m-1,n-m}\|_{L^2(S_{\ub,u})}du
+C\ub\int_{\ub}^{u_1}u^2\|\cmu_{m-1,n-m}\|_{L^2(S_{\ub,u})}du\nonumber\\
&&+C\s^{(m,n-m)}J(\ub,u_1)
\label{12.700}
\end{eqnarray}
where:
\begin{eqnarray}
&&\mbox{for $m=1$:} \hspace{3mm} \s^{(1,n-1)}J(\ub,u_1)=\ub^2
\int_{\ub}^{u_1}u^2\|\s^{(n-1)}\ctchi\|_{L^2(S_{\ub,u})}du\nonumber\\
&&\hspace{48mm}+\ub\int_{\ub}^{u_1}u^4\|\s^{(n-1)}\ctchib\|_{L^2(S_{\ub,u})}du\nonumber\\
&&\mbox{for $m\geq 2$:} \hspace{3mm} \s^{(m,n-m)}J(\ub,u_1)=\ub^2
\int_{\ub}^{u_1}u^2\|\cmu_{m-2,n-m+1}\|_{L^2(S_{\ub,u})}du\nonumber\\
&&\hspace{48mm}+\ub\int_{\ub}^{u_1}u^4\|\cmub_{m-2,n-m+1}\|_{L^2(S_{\ub,u})}du\nonumber\\
&&\label{12.701}
\end{eqnarray}
The first term in \ref{12.700} is bounded by a constant multiple of:
\begin{eqnarray}
&&\ub^{1/2}u_1^{3/2}\|\sqrt{a}\s^{(E;m,n-m)}\csxi\|_{L^2(\Cb_{\ub}^{u_1})}\leq 
C\ub^{1/2}u_1^{3/2}\sqrt{\s^{(E;m,n-m)}\cEb^{u_1}(\ub)}\nonumber\\
&&\hspace{20mm}\leq C\ub^{a_m+\frac{1}{2}}u_1^{b_m+\frac{3}{2}}\sqrt{\s^{(E;m,n-m)}\cBb(\ub_1,u_1)}
\label{12.702}
\end{eqnarray}
hence its contribution to \ref{12.676} is bounded by:
\begin{equation}
C\frac{\ub_1^{a_m+1}}{\sqrt{2a_m+2}}u_1^{b_m+\frac{3}{2}}\sqrt{\s^{(E;m,n-m)}\cBb(\ub_1,u_1)}
\label{12.703}
\end{equation}
The second term in \ref{12.700} can be absorbed into the first term in \ref{12.678}. The third 
and fourth terms can likewise be absorbed into the corresponding terms in \ref{12.678}. As for 
$\s^{(m,n-m)}J(\ub,u_1)$, given by \ref{12.701}, consider first the case $m=1$. The 
1st term in $\s^{(1,n-1)}J(\ub,u_1)$ is bounded by
\begin{equation}
\ub^2 u_1^{5/2}\|\s^{(n-1)}\ctchi\|_{L^2(\Cb_{\ub}^{u_1})}
\leq \ub^{a_0+2}u_1^{b_0+\frac{5}{2}}\s^{(n-1)}X(\ub_1,u_1)
\label{12.704}
\end{equation}
in view of the definition \ref{12.411}. Hence its contribution to \ref{12.676} is bounded by:
\begin{equation}
\frac{\ub_1^{a_0+\frac{5}{2}}}{\sqrt{2a_0+5}}u_1^{b_0+\frac{5}{2}}\s^{(n-1)}X(\ub_1,u_1)
\label{12.705}
\end{equation}
To estimate the contribution to \ref{12.676} of the 2nd term in $\s^{(1,n-1)}J(\ub,u_1)$ we 
apply Lemma 12.5 taking 
\begin{equation}
f(\ub,u)=\ub u^4\|\s^{(n-1)}\ctchib\|_{L^2(S_{\ub,u})}
\label{12.706}
\end{equation}
Then the term in question is $g(\ub)$. Since 
\begin{eqnarray}
&&\|f(\cdot, u)\|_{L^2(0,\ub_1)}\leq Cu^4\|\rho\s^{(n-1)}\ctchib\|_{L^2(C_u^{\ub_1})}\nonumber\\
&&\hspace{20mm}\leq C\ub_1^{a_0+1}u^{b_0+3}\s^{(n-1)}\Xb(\ub_1,u_1) 
\label{12.707}
\end{eqnarray}
in view of the definition \ref{12.460}, using Lemma 12.5 we conclude that the contribution 
to \ref{12.676} is bounded by:
\begin{equation}
C\ub_1^{a_0+1}\frac{u_1^{b_0+4}}{b_0+4}\s^{(n-1)}\Xb(\ub_1,u_1)
\label{12.708}
\end{equation}
In the case $m\geq 2$, the 1st term in $\s^{(m,n-m)}J(\ub,u_1)$ is bounded by $\ub u_1^2$ times 
the 5th term in \ref{12.678} with $m-1$ in the role of $m$, while the 2nd term in 
$\s^{(m,n-m)}J(\ub,u_1)$ is bounded by $u_1^4$ times the 4th term in \ref{12.678} with $m-1$ in the 
role of $m$. 

Comparing with \ref{12.631} the $N$ component difference in \ref{12.694} is:
\begin{eqnarray}
&&\frac{\beta_{\Nb}^3}{4c}\eta^2 H^\prime\rhob
\left\{N^\mu L(E^{n-m}T^m\beta_\mu-E_N^{n-m}T^m\beta_{\mu,N})\right.\nonumber\\
&&\hspace{30mm}\left.+2\ss_N\left(\rho\cmu_{m-1,n-m}-\rhob\cmub_{m-1,n-m}\right)\right\}\nonumber\\
&&\label{12.709}
\end{eqnarray}
up to terms which can be ignored. 
The contribution of this to the integral on the right in \ref{12.674} is bounded by:
\begin{eqnarray}
&&C\int_{\ub}^{u_1}u^2\|N^\mu L(E^{n-m}T^m\beta_\mu-E_N^{n-m}T^m\beta_{\mu,N})\|_{L^2(S_{\ub,u})}du
\nonumber\\
&&+C\ub\int_{\ub}^{u_1}u^2\|\cmu_{m-1,n-m}\|_{L^2(S_{\ub,u})}du
+C\int_{\ub}^{u_1}u^4\|\cmub_{m-1,n-m}\|_{L^2(S_{\ub,u})}du\nonumber\\
&&\label{12.710}
\end{eqnarray}
The last two terms here are the same as the last two terms in \ref{12.696}. To estimate the first term 
we appeal to \ref{12.633}. The contribution of the right hand side of \ref{12.633} through the 
first of \ref{12.710} is analogous to the first of \ref{12.696} with the variation field $Y$ in the 
role of the variation field $E$. Hence this contribution to \ref{12.676} is bounded by 
(see \ref{12.698}):
\begin{equation}
C\frac{\sqrt{\s^{(Y;m,n-m)}\cB(\ub_1,u_1)}}{b_m+3}\cdot\ub_1^{a_m}u_1^{b_m}\cdot u_1^3
\label{12.711}
\end{equation}
As for the 2nd term on the left had side of \ref{12.633}, this is expressed in the form \ref{12.637} 
up to terms which can be ignored. The contribution of the principal difference in \ref{12.637} is 
similar to the contribution of the $\Nb$ component difference but with an extra $u_1$ factor 
(as $\ogamma\sim u$). As for the contribution of the remainder, this is similar to the last two 
terms in \ref{12.710} but with an extra $u_1$ factor. 

Consider next the second of the principal difference terms \ref{12.692}. In fact the actual principal 
difference term in question is:
\begin{eqnarray}
&&-2\rho\beta_{\Nb}\left(c\pi\sbeta-\frac{\beta_{\Nb}}{4}-\frac{\beta_N}{2}\right)H^\prime
(g^{-1})^{\mu\nu}\beta_\nu E^{n-m}T^m\Lb\beta_\mu\nonumber\\
&&+2\rho_N\beta_{\Nb,N}\left(c_N\pi_N\sbeta_N-\frac{\beta_{\Nb,N}}{4}-\frac{\beta_{N,N}}{2}\right)
H_N^\prime(g^{-1})^{\mu\nu}\beta_{\nu,N} E_N^{n-m}T^m\Lb_N\beta_{\mu,N}\nonumber\\
&&\label{12.712}
\end{eqnarray}
This is similar to the principal difference \ref{12.572}. It is similarly expressed as differences 
of $E$, $N$ and $\Nb$ components:
\begin{eqnarray}
&&-2\rho\sbeta\eta^2\beta_{\Nb}\left(c\pi\sbeta-\frac{\beta_{\Nb}}{4}-\frac{\beta_N}{2}\right)
H^\prime E^\mu E^{n-m}T^m\Lb\beta_\mu\nonumber\\
&&+2\rho_N\sbeta_N\eta_N^2\beta_{\Nb,N}\left(c_N\pi_N\sbeta_N-\frac{\beta_{\Nb,N}}{4}
-\frac{\beta_{N,N}}{2}\right)H_N^\prime E_N^\mu E_N^{n-m}T^m\Lb_N\beta_{\mu,N}\nonumber\\
&&+\rho\frac{\beta_{\Nb}}{c}\eta^2\beta_{\Nb}
\left(c\pi\sbeta-\frac{\beta_{\Nb}}{4}-\frac{\beta_N}{2}\right)
H^\prime N^\mu E^{n-m}T^m\Lb\beta_\mu\nonumber\\
&&-\rho_N\frac{\beta_{\Nb,N}}{c_N}\eta_N^2\beta_{\Nb,N}\left(c_N\pi_N\sbeta_N-\frac{\beta_{\Nb,N}}{4}-\frac{\beta_{N,N}}{2}\right)
H_N^\prime N_N^\mu E_N^{n-m}T^m\Lb_N\beta_{\mu,N}\nonumber\\
&&+\rho\frac{\beta_N}{c}\eta^2\beta_{\Nb}
\left(c\pi\sbeta-\frac{\beta_{\Nb}}{4}-\frac{\beta_N}{2}\right)
H^\prime \Nb^\mu E^{n-m}T^m\Lb\beta_\mu\nonumber\\
&&-\rho_N\frac{\beta_{N,N}}{c_N}\eta_N^2\beta_{\Nb,N}\left(c_N\pi_N\sbeta_N-\frac{\beta_{\Nb,N}}{4}-\frac{\beta_{N,N}}{2}\right)H_N^\prime\Nb_N^\mu E_N^{n-m}T^m\Lb_N\beta_{\mu,N}\nonumber\\
&&\label{12.713}
\end{eqnarray}
Comparing with \ref{12.615} the $E$ component difference here is:
\begin{eqnarray}
&&-2\beta_{\Nb}\left(c\pi\sbeta-\frac{\beta_{\Nb}}{4}-\frac{\beta_N}{2}\right)\sbeta\eta^2 H^\prime\rho\left\{\s^{(E;m,n-m)}\cxi_{\Lb}\right.\nonumber\\
&&\hspace{50mm}\left.-2\sss\left(\rho\cmu_{m-1,n-m}-\rhob\cmub_{m-1,n-m}\right)\right\}\nonumber\\
&&\label{12.714}
\end{eqnarray}
up to terms which can be ignored. The contribution of this to the integral on the right in \ref{12.674} 
is bounded by:
\begin{eqnarray}
&&C\ub\int_{\ub}^{u_1}\|\s^{(E;m,n-m)}\cxi_{\Lb}\|_{L^2(S_{\ub,u})}du\label{12.715}\\
&&+C\ub^2\int_{\ub}^{u_1}\|\cmu_{m-1,n-m}\|_{L^2(S_{\ub,u})}du
+C\ub\int_{\ub}^{u_1}u^2\|\cmub_{m-1,n-m}\|_{L^2(S_{\ub,u})}du\nonumber
\end{eqnarray}
Here the 2nd term is dominated by $\ub$ times the 5th term in \ref{12.678} while the 3rd term is 
dominated by $u_1^2$ times the 4th term in \ref{12.678}. The first term in \ref{12.715} is bounded by a constant multiple of:
\begin{eqnarray}
&&\ub u_1^{1/2}\|\s^{(E;m,n-m)}\cxi_{\Lb}\|_{L^2(\Cb_{\ub}^{u_1})}\leq 
C\ub u_1^{1/2}\sqrt{\s^{(E;m,n-m)}\cEb^{u_1}(\ub)}\nonumber\\
&&\hspace{20mm}\leq C\ub^{a_m+1}u_1^{b_m+\frac{1}{2}}\sqrt{\s^{(E;m,n-m)}\cBb(\ub_1,u_1)}
\label{12.716}
\end{eqnarray}
hence its contribution to \ref{12.676} is bounded by:
\begin{equation}
C\frac{\ub_1^{a_m+\frac{3}{2}}}{\sqrt{2a_m+3}}u_1^{b_m+\frac{1}{2}}\sqrt{\s^{(E;m,n-m)}\cBb(\ub_1,u_1)}
\label{12.717}
\end{equation}

Comparing with \ref{12.586}, the $N$ component difference in \ref{12.713} is, up to terms which can be 
ignored, 
\begin{eqnarray}
&&\frac{\beta_{\Nb}^2}{c}\left(c\pi\sbeta-\frac{\beta_{\Nb}}{4}-\frac{\beta_N}{2}\right)\eta^2 H^\prime\rho\left[ \lambda \s^{(E;m,n-m)}\csxi 
+\sss\s^{(m,n-m)}\cla\right.\nonumber\\
&&\hspace{27mm}+\rhob\left(\ss_N\cmub_{m-1,n-m}+(c\pi\sss-\ss_{\Nb})\cmu_{m-1,n-m}\right)
\nonumber\\
&&\hspace{27mm}-\lambda\sss\rho\cdot\left\{ \begin{array}{lll}
\s^{(n-1)}\ctchi &:& \mbox{if $m=1$}\\
2\cmu_{m-2,n-m+1} &:& \mbox{if $m\geq 2$}
\end{array} \right.\nonumber\\
&&\hspace{27mm}\left.-\lambda\sss\rhob\cdot\left\{ \begin{array}{lll}
\s^{(n-1)}\ctchib &:& \mbox{if $m=1$}\\
2\cmub_{m-2,n-m+1} &:& \mbox{if $m\geq 2$}
\end{array} \right. \hspace{2mm}\right]
\label{12.718}
\end{eqnarray}
The contribution of this to the integral on the right in \ref{12.674} is bounded by:
\begin{eqnarray}
&&C\sqrt{\ub}\int_{\ub}^{u_1}u\|\sqrt{a}\s^{(E;m,n-m)}\csxi\|_{L^2(S_{\ub,u})}du
+C\ub\int_{\ub}^{u_1}\|\s^{(m,n-m)}\cla\|_{L^2(S_{\ub,u})}du\nonumber\\
&&+C\ub\int_{\ub}^{u_1}u^2\|\cmub_{m-1,n-m}\|_{L^2(S_{\ub,u})}du
+C\ub\int_{\ub}^{u_1}u^2\|\cmu_{m-1,n-m}\|_{L^2(S_{\ub,u})}du\nonumber\\
&&+C\s^{(m,n-m)}J(\ub,u_1) \label{12.719}
\end{eqnarray}
The 1st, 3rd, 4th, and 5th terms here coincide with the corresponding term in \ref{12.700}, while 
the 2nd term here coincides with the 2nd term in \ref{12.678}. 

We turn to the last of the differences \ref{12.713}, the $\Nb$ component difference. 
Comparing with \ref{12.607} this is given by:
\begin{eqnarray}
&&\frac{\beta_{\Nb}\beta_N}{c}\left(c\pi\sbeta-\frac{\beta_{\Nb}}{4}-\frac{\beta_N}{2}\right)\eta^2 H^\prime\rho\left\{\Nb^\mu\Lb(E^{n-m}T^m\beta_\mu-E_N^{n-m}T^m\beta_{\mu,N})\right.\nonumber\\
&&\hspace{55mm}\left.-2\ss_{\Nb}\left(\rho\cmu_{m-1,n-m}-\rhob\cmub_{m-1,n-m}\right)\right\}\nonumber\\
&&\label{12.720}
\end{eqnarray}
The contribution of this to the integral on the right in \ref{12.674} is bounded by:
\begin{eqnarray}
&&C\ub\int_{\ub}^{u_1}\|\Nb^\mu\Lb(E^{n-m}T^m\beta_\mu-E_N^{n-m}T^m\beta_{\mu,N})\|_{L^2(S_{\ub,u})}du
\nonumber\\
&&+C\ub^2\int_{\ub}^{u_1}\|\cmu_{m-1,n-m}\|_{L^2(S_{\ub,u})}du
+C\ub\int_{\ub}^{u_1}u^2\|\cmub_{m-1,n-m}\|_{L^2(S_{\ub,u})}du\nonumber\\
&&\label{12.721}
\end{eqnarray}
Here the last two terms are dominated by the corresponding terms in \ref{12.678}. To estimate the first 
term we appeal to \ref{12.574}. Since $\ogamma\sim u$, this term is seen to be bounded by:
\begin{eqnarray}
&&C\ub\int_{\ub}^{u_1}u^{-1}\|\s^{(Y;m,n-m)}\cxi_{\Lb}\|_{L^2(S_{\ub,u})}du\nonumber\\
&&+C\ub\int_{\ub}^{u_1}u^{-1}\|N^\mu\Lb(E^{n-m}T^m\beta_\mu
-E_N^{n-m}T^m\beta_{\mu,N})\|_{L^2(S_{\ub,u})}du\nonumber\\
&&\label{12.722}
\end{eqnarray}
Following the treatment of the 2nd term in \ref{12.609}, the 2nd term here is shown to be bounded 
in terms of:
\begin{eqnarray}
&&C\ub\int_{\ub}^{u_1}u^{-1}\|N^\mu E^{n-m}T^m\Lb\beta_\mu
-N_N^\mu E_N^{n-m}T^m\Lb_N\beta_{\mu,N}\|_{L^2(S_{\ub,u})}du\nonumber\\
&&+C\ub^2\int_{\ub}^{u_1}u^{-1}\|\cmu_{m-1,n-m}\|_{L^2(S_{\ub,u})}du 
+C\ub\int_{\ub}^{u_1}u\|\cmub_{m-1,n-m}\|_{L^2(S_{\ub,u})}du\nonumber\\
&&\label{12.723}
\end{eqnarray}
Here the last two terms are dominated by the corresponding terms in \ref{12.678}, while the 1st is 
similar to the contribution of the $N$ component difference estimated above but with an 
extra $u^{-1}$ factor in the integrant. It follows that up to terms which can be absorbed the 1st 
term in \ref{12.723} is bounded by a constant multiple of:
\begin{eqnarray}
&&\sqrt{\ub}\int_{\ub}^{u_1}\|\sqrt{a}\s^{(E;m,n-m)}\csxi\|_{L^2(S_{\ub,u})}du
\leq C\ub^{1/2}u_1^{1/2}\sqrt{\s^{(E;m,n-m)}\cEb^{u_1}(\ub)}\nonumber\\
&&\hspace{20mm}\leq C\ub^{a_m+\frac{1}{2}}u_1^{b_m+\frac{1}{2}}\sqrt{\s^{(E;m,n-m)}\cBb(\ub_1,u_1)}
\label{12.724}
\end{eqnarray}
hence the corresponding contribution to \ref{12.676} is bounded by:
\begin{equation}
C\frac{\ub_1^{a_m+1}u_1^{b_m+\frac{1}{2}}}{\sqrt{2a_m+2}}\sqrt{\s^{(E;m,n-m)}\cBb(\ub_1,u_1)}
\label{12.725}
\end{equation}
What remains is then the first of \ref{12.722}. This is bounded by (Schwartz inequality):
\begin{eqnarray}
&&C\ub\left(\int_{\ub}^{u_1}u^{-2}du\right)^{1/2}\|\s^{(Y;m,n-m)}\cxi_{\Lb}\|_{L^2(\Cb_{\ub}^{u_1})}
\nonumber\\
&&\leq C\ub^{1/2}\sqrt{\s^{(Y;m,n-m)}\cEb^{u_1}(\ub)}\leq C\ub^{a_m+\frac{1}{2}}u_1^{b_m}
\sqrt{\s^{(Y;m,n-m)}\cBb(\ub_1,u_1)}\nonumber\\
&&\label{12.726}
\end{eqnarray}
hence the corresponding contribution to \ref{12.676} is bounded by:
\begin{equation}
C\frac{\ub_1^{a_m+1}u_1^{b_m}}{\sqrt{2a_m+2}}\sqrt{\s^{(Y;m,n-m)}\cBb(\ub_1,u_1)}
\label{12.727}
\end{equation}

The third of the principal difference terms \ref{12.692} is similar to the $\Nb$ component difference 
in \ref{12.713}. 

We finally consider the fourth of the principal difference terms \ref{12.692}. In fact the actual 
principal difference term in question is:
\begin{equation}
\pi a\beta_{\Nb}^2 H^\prime(g^{-1})^{\mu\nu}\beta_\nu E^{n-m}T^m E\beta_\mu
-\pi_N a_N\beta_{\Nb,N}^2 H^\prime_N (g^{-1})^{\mu\nu}\beta_{\nu,N}E_N^{n-m}T^m E_N\beta_{\mu,N}
\label{12.728}
\end{equation}
This is similar to the principal difference \ref{12.638}. It is similarly expressed as differences of 
$E$, $N$, and $\Nb$ components:
\begin{eqnarray}
&&a\pi\beta_{\Nb}^2\sbeta\eta^2 H^\prime E^\mu E^{n-m}T^m E\beta_\mu \nonumber\\
&&\hspace{20mm}-a_N\pi_N\beta_{\Nb,N}^2\beta_N\eta_N^2 H^\prime_N E_N^\mu E_N^{n-m}T^m E_N\beta_{\mu,N}\nonumber\\
&&-a\pi\beta_{\Nb}^2\frac{\beta_{\Nb}}{2c}\eta^2 H^\prime N^\mu E^{n-m}T^m E\beta_\mu \nonumber\\
&&\hspace{20mm}+a_N\pi_N\beta_{\Nb,N}^2\frac{\beta_{\Nb,N}}{2c_N}\eta_N^2 H^\prime_N N_N^\mu E_N^{n-m}T^m E_N\beta_{\mu,N}\nonumber\\
&&-a\pi\beta_{\Nb}^2\frac{\beta_N}{2c}\eta^2 H^\prime\Nb^\mu E^{n-m}T^m E\beta_\mu \nonumber\\
&&\hspace{20mm}+a_N\pi_N\beta_{\Nb,N}^2\frac{\beta_{N,N}}{2c_N}\eta_N^2 H^\prime_N\Nb_N^\mu E_N^{n-m}T^m E_N\beta_{\mu,N} \nonumber\\
&&\label{12.729}
\end{eqnarray}
Comparing with \ref{12.651} the $E$ component difference is, up to terms which can be ignored, 
\begin{eqnarray}
&&\pi\beta_{\Nb}^2\eta^2 H^\prime a\left[\s^{(E;m,n-m)}\csxi\right.
\nonumber\\
&&\hspace{24mm}-\sss\rho\cdot\left\{ \begin{array}{lll}
\s^{(n-1)}\ctchi &:& \mbox{if $m=1$}\\
2\cmu_{m-2,n-m+1} &:& \mbox{if $m\geq 2$}
\end{array} \right.\nonumber\\
&&\hspace{24mm}\left.-\sss\rhob\cdot\left\{ \begin{array}{lll}
\s^{(n-1)}\ctchib &:& \mbox{if $m=1$}\\
2\cmub_{m-2,n-m+1} &:& \mbox{if $m\geq 2$}
\end{array} \right. \hspace{2mm}\right]
\label{12.730}
\end{eqnarray}
The contribution of this to the integral on the right in \ref{12.674} is bounded by:
\begin{equation}
C\sqrt{\ub}\int_{\ub}^{u_1}u\|\sqrt{a}\s^{(E;m,n-m)}\csxi\|_{L^2(S_{\ub,u})}du+C\s^{(m,n-m)}J(\ub,u_1)
\label{12.731}
\end{equation}
The two terms here coincide with the first and last terms in \ref{12.719}. 

We turn to the $\Nb$ component difference (last of \ref{12.729}). Comparing with \ref{12.647} this 
difference is, up to terms which can be ignored, 
\begin{eqnarray}
&&-\frac{\pi\beta_{\Nb}^2\beta_N}{2}\eta^2 H^\prime\rho\left\{\s^{(E;m,n-m)}\cxi_{\Lb}
-c^{-1}\ss_{\Nb}\s^{(m,n-m)}\cla\right.\nonumber\\
&&\hspace{20mm}\left.+(c^{-1}s_{\Nb\Lb}-2\rho\sss-2\pi\rhob\ss_{\Nb})\cmu_{m-1,n-m}
+\rhob\sss\cmub_{m-1,n-m}\right\}\nonumber\\
&&\label{12.732}
\end{eqnarray}
The contribution of this to the integral on the right in \ref{12.674} is bounded by:
\begin{eqnarray}
&&C\ub\int_{\ub}^{u_1}\|\s^{(E;m,n-m)}\cxi_{\Lb}\|_{L^2(S_{\ub,u})}du
+C\ub\int_{\ub}^{u_1}\|\s^{(m,n-m)}\cla\|_{L^2(S_{\ub,u})}du\nonumber\\
&&+C\ub\int_{\ub}^{u_1}\|\cmu_{m-1,n-m}\|_{L^2(S_{\ub,u})}du
+C\ub\int_{\ub}^{u_1}u^2\|\cmub_{m-1,n-m}\|_{L^2(S_{\ub,u})}du\nonumber\\
&&\label{12.733}
\end{eqnarray}
Here the 2nd and 3rd terms coincide with the 2nd and 5th terms in \ref{12.678} respectively, 
while the 4th term here is dominated by $u_1^2$ times the 4th term in \ref{12.678}. The 1st term 
in \ref{12.733} is bounded by a constant multiple of:
\begin{eqnarray}
&&\ub u_1^{1/2}\|\s^{(E;m,n-m)}\cxi_{\Lb}\|_{L^2(\Cb_{\ub}^{u_1})}
\leq C\ub u_1^{1/2}\sqrt{\s^{(E;m,n-m)}\cEb^{u_1}(\ub)}\nonumber\\
&&\hspace{20mm}\leq C\ub^{a_m+1}u_1^{b_m+\frac{1}{2}}\sqrt{\s^{(E;m,n-m)}\cBb(\ub_1,u_1)}
\label{12.734}
\end{eqnarray}
hence its contribution to \ref{12.676} is bounded by:
\begin{equation}
C\frac{\ub_1^{a_m+\frac{3}{2}}u_1^{b_m+\frac{1}{2}}}{\sqrt{2a_m+3}}\sqrt{\s^{(E;m,n-m)}\cBb(\ub_1,u_1)}
\label{12.735}
\end{equation}

Finally, in regard to the $N$ component difference (second of \ref{12.729}), using \ref{12.578} - 
\ref{12.580} we find that up to terms which can be ignored this difference is expressed as:
\begin{eqnarray}
&&-\frac{\pi\beta_{\Nb}^3}{2c}\eta^2 H^\prime a
\left[N^\mu E(E^{n-m}T^m\beta_\mu-E_N^{n-m}T^m\beta_{\mu,N})\right.
\nonumber\\
&&\hspace{24mm}-\ss_N\rho\cdot\left\{ \begin{array}{lll}
\s^{(n-1)}\ctchi &:& \mbox{if $m=1$}\\
2\cmu_{m-2,n-m+1} &:& \mbox{if $m\geq 2$}
\end{array} \right.\nonumber\\
&&\hspace{24mm}\left.-\ss_N\rhob\cdot\left\{ \begin{array}{lll}
\s^{(n-1)}\ctchib &:& \mbox{if $m=1$}\\
2\cmub_{m-2,n-m+1} &:& \mbox{if $m\geq 2$}
\end{array} \right. \hspace{2mm}\right]
\label{12.736}
\end{eqnarray}
The contribution of this to the integral on the right in \ref{12.674} is bounded by:
\begin{eqnarray}
&&C\int_{\ub}^{u_1}\|a N^\mu E(E^{n-m}T^m\beta_\mu-E_N^{n-m}T^m\beta_{\mu,N})\|_{L^2(S_{\ub,u})}du
\nonumber\\
&&+C\s^{(m,n-m)}J(\ub,u_1)
\label{12.737}
\end{eqnarray}
To estimate the contribution of the first term here we appeal to the fact that
\begin{eqnarray}
&&N^\mu E(E^{n-m}T^m\beta_\mu-E_N^{n-m}T^m\beta_{\mu,N})
+\ogamma\Nb^\mu E(E^{n-m}T^m\beta_\mu-E_N^{n-m}T^m\beta_{\mu,N})\nonumber\\
&&\hspace{50mm}=\s^{(Y;m,n-m)}\csxi
\label{12.738}
\end{eqnarray}
The contribution 2nd term on the left, through the first term in \ref{12.737}, to \ref{12.676} 
is bounded, up to terms which can be absorbed, by $u_1$ times the contribution of the $\Nb$ 
component difference estimated above, that is by:
\begin{equation}
C\frac{\ub_1^{a_m+\frac{3}{2}}u_1^{b_m+\frac{3}{2}}}{\sqrt{2a_m+3}}\sqrt{\s^{(E;m,n-m)}\cBb(\ub_1,u_1)}
\label{12.739}
\end{equation}
As for the contribution of the right hand side, since
\begin{eqnarray}
&&\int_{\ub}^{u_1}\|a\s^{(Y;m,n-m)}\csxi\|_{L^2(S_{\ub,u})}du
\leq C\ub^{1/2}u_1^{3/2}\|\sqrt{a}\s^{(Y;m,n-m)}\csxi\|_{L^2(\Cb_{\ub}^{u_1})}\nonumber\\
&&\leq C\ub^{1/2}u_1^{3/2}\sqrt{\s^{(Y;m,n-m)}\cEb^{u_1}(\ub)}
\leq C\ub^{a_m+\frac{1}{2}}u_1^{b_m+\frac{3}{2}}\sqrt{\s^{(Y;m,n-m)}\cBb(\ub_1,u_1)}\nonumber\\
&&\label{12.740}
\end{eqnarray}
the contribution to \ref{12.676} is bounded by:
\begin{equation}
C\frac{\ub_1^{a_m+1}u_1^{b_m+\frac{3}{2}}}{\sqrt{2a_m+2}}\sqrt{\s^{(Y;m,n-m)}\cBb(\ub_1,u_1)}
\label{12.741}
\end{equation}

Absorbing the first of \ref{12.678} in the integral inequality \ref{12.674}, and recalling that the 
exponents $a_m$, $b_m$ are non-increasing in $m$, 
we summarize the above results in the form in which they will be used in the sequel in the 
following lemma.

\vspace{2.5mm}

\noindent{\bf Lemma 12.10} \ \ The next to the top order acoustical quantities 
$\|\s^{(m,n-m)}\clab\|_{L^2(C_u^{\ub})}$ satisfy to principal terms the following inequalities:
\begin{eqnarray*}
&&\|\s^{(m,n-m)}\clab\|_{L^2(C_{u_1}^{\ub_1})}\leq k\|\s^{(m,n-m)}\clab\|_{L^2({\cal K}^{\ub_1})}\\
&&\hspace{27mm}+Cu_1^{1/2}\left\{\int_0^{\ub_1}\ub^2\|\s^{(m,n-m)}\cla\|^2_{L^2(\Cb_{\ub}^{u_1})}d\ub
\right\}^{1/2}\\
&&\hspace{27mm}+Cm\frac{\ub_1^{a_m+1}u_1^{b_m+1}}{(b_m+1)}\s^{(m-1,n-m+1)}\Lab(\ub_1,u_1)\\
&&\hspace{24mm}+C\sum_{j=0}^{m-1}(\Pi u_1^2)^j \ \frac{\ub_1^{a_m+2}u_1^{b_m+1}}{\sqrt{2a_m+4}}
\s^{(m-1-j,n-m+j+1)}\La(\ub_1,u_1)\\
&&\hspace{27mm}+C(\Pi u_1^2)^m \ \frac{\ub_1^{a_m+\frac{3}{2}}u_1^{b_m+\frac{1}{2}}}{\sqrt{2a_m+3}}
\s^{(n-1)}X(\ub_1,u_1)\\
&&\hspace{24mm}+Cu_1\sum_{j=0}^{m-1}\Pib^j \ \frac{\ub_1^{a_m+j}u_1^{b_m+1}}{(b_m+1)}
\s^{(m-1-j,n-m+j+1)}\Lab(\ub_1,u_1)\\
&&\hspace{27mm}+Cu_1\Pib^m \ \frac{\ub_1^{a_m+m}u_1^{b_m}}{b_m}\s^{(n-1)}\Xb(\ub_1,u_1)\\
&&\hspace{20mm}+\Cb\sqrt{\max\{\s^{(m,n-m)}\cB(\ub_1,u_1),\s^{(m,n-m)}\cBb(\ub_1,u_1)\}}
\cdot \ub_1^{a_m}u_1^{b_m}\cdot u_1
\end{eqnarray*}

Here $k$ is again a constant greater than 1 but which can be chosen as close to 1 as we wish by 
suitably restricting $\delta$. 

\vspace{2.5mm}

To proceed we must bring in the boundary condition \ref{10.752} for $\lambda$ on ${\cal K}$ and the 
corresponding boundary condition \ref{10.757} for $\lambdab_N$. Applying $T$ to \ref{10.752} and to 
\ref{10.757} and subtracting gives:
\begin{equation}
T\cla=rT\clab+\lambdab T\check{r}+(Tr)\clab+(T\lambdab)\check{r}+{\cal M}_{-1,0}
\label{12.742}
\end{equation}
where:
\begin{equation}
{\cal M}_{-1,0}=\check{r} T\clab+\clab T\check{r}-T\hat{\nu}_N
\label{12.743}
\end{equation}
Here and in the following we denote:
\begin{equation}
\check{r}=r-r_N
\label{12.744}
\end{equation}
The first two terms in \ref{12.742} are of order 1 with coefficients, $r$ and $\lambdab$, which are 
$\sim\tau$. The first two terms in remainder ${\cal M}_{-1,0}$ are also of order 1, but with 
coefficients, $\check{r}$ and $\clab$, for which we can assume that:
\begin{equation}
|\check{r}|, |\clab| \leq C\tau^2 \ : \ \mbox{on ${\cal K}^{\delta}$}
\label{12.745}
\end{equation}
Applying $T^2$ to \ref{10.752} and to \ref{10.757} and subtracting gives:
\begin{equation}
T^2\cla=rT^2\clab+\lambdab T^2\check{r}+2(Tr)T\clab+2(T\lambdab)T\check{r}+{\cal M}_{0,0}
\label{12.746}
\end{equation}
where: 
\begin{equation}
{\cal M}_{0,0}=-2(T\check{r})(T\clab)+\check{r}T^2\lambdab_N+\clab T^2 r_N-T^2\hat{\nu}_N
\label{12.747}
\end{equation}
The first two terms in \ref{12.746} are of order 2, with coefficients $r$ and $\lambdab$, which are 
$\sim\tau$. The next two terms in \ref{12.746} are of order 1, but, considering $2(Tr)$ and 
$2(T\lambdab)$ to be the coefficients, these coefficients are $\sim 1$. The remainder ${\cal M}_{0,0}$ 
is of order 1. Moreover, the order 1 part of ${\cal M}_{0,0}$, $-2(T\check{r})(T\clab)$, can be 
absorbed in the 3rd and 4th terms in \ref{12.746} as we can assume that:
\begin{equation}
|T\check{r}|, |T\clab| \leq C\tau \ : \ \mbox{on ${\cal K}^{\delta}$}
\label{12.748}
\end{equation}
Applying finally $T^m$ for $m\geq 3$ to \ref{10.752} and to \ref{10.757} and subtracting 
gives:
\begin{equation}
T^m\cla=rT^m\clab+\lambdab T^m\check{r}+m(Tr)T^{m-1}\clab+m(T\lambdab)T^{m-1}\check{r}
+{\cal M}_{m-2,0}
\label{12.749}
\end{equation}
where the remainder ${\cal M}_{m-2,0}$ is of order $m-2$. The first two terms in \ref{12.749} are of order m, with coefficients $r$ and $\lambdab$, which are $\sim\tau$. The next two terms are of order 
$m-1$, but their coefficients $m(Tr)$ and $m(T\lambdab)$ are $\sim 1$.
We may incorporate the formula \ref{12.749} the formulas \ref{12.746} and \ref{12.742} 
as the cases $m=2$ and $m=1$ respectively, however the statement 
in regard to the order of the remainder does not hold in these cases. 

We now apply $\Omega^{n-m}$ to \ref{12.750} to obtain:
\begin{eqnarray}
&&\Omega^{n-m}T^m\cla=r\Omega^{n-m}T^m\clab+\lambdab\Omega^{n-m}T^m\check{r}\nonumber\\
&&\hspace{20mm}+m(Tr)\Omega^{n-m}T^{m-1}\clab+m(T\lambdab)\Omega^{n-m}T^{m-1}\check{r}\nonumber\\
&&\hspace{20mm}+{\cal M}_{m-2,n-m} \label{12.750}
\end{eqnarray}
where:
\begin{eqnarray}
&&{\cal M}_{m-2,n-m}=\Omega^{n-m}{\cal M}_{m-2,0}\nonumber\\
&&\hspace{20mm}+\sum_{i=1}^{n-m}\left(\begin{array}{c} n-m\\i\end{array}\right)(\Omega^i r)
\Omega^{n-m-i}T^m\clab\nonumber\\
&&\hspace{20mm}+\sum_{i=1}^{n-m}\left(\begin{array}{c} n-m\\i\end{array}\right)(\Omega^i\lambdab)
\Omega^{n-m-i}T^m\check{r}\nonumber\\
&&\hspace{20mm}+m\sum_{i=1}^{n-m}\left(\begin{array}{c} n-m\\i\end{array}\right)(\Omega^i Tr)
\Omega^{n-m-i}T^{m-1}\clab\nonumber\\
&&\hspace{20mm}+m\sum_{i=1}^{n-m}(\Omega^i T\lambdab)\Omega^{n-m-i}T^{m-1}\check{r} 
\label{12.751}
\end{eqnarray}
For $m\geq 2$ the highest order term in the 1st sum in \ref{12.751} is the $i=1$ term:
$$(n-m)(\Omega r)\Omega^{n-m-1}T^m\clab$$
This is of order $n-1$ but with a coefficient $\Omega r$ which is bounded in absolute value by 
$C\tau$ on ${\cal K}^{\delta}$. Similarly, for $m\geq 2$ the highest order term in the 2nd sum 
in \ref{12.751} is the $i=1$ term:
$$(n-m)(\Omega\lambdab)\Omega^{n-m-i}T^m\check{r}$$
This is of order $n-1$ but with a coefficient $\Omega\lambdab$ which is bounded in absolute value by 
$C\tau$ on ${\cal K}^{\delta}$. It follows that these terms, like the remainders of the 
1st and 2nd sums, can be absorbed in the estimates.  
For $m\geq 3$ the 3rd and 4th sums in \ref{12.751} are only of order $n-2$ and can likewise 
be absorbed in the estimates. On the other hand, for $m=2$ the highest order term in the 3rd sum is 
the $i=n-2$ term:
$$2(\Omega^{n-2}Tr)T\clab$$
This is of order $n-1$, the principal part being:
$$2(T\clab)\Omega^{n-2}T\check{r}$$
however the coefficient $T\clab$ is, according to \ref{12.748}, bounded in absolute value by $C\tau$ 
on ${\cal K}^{\delta}$. Similarly, for $m=2$ the highest order term in the 4th sum is the $i=n-2$ term: 
$$2(\Omega^{n-2}T\lambdab)T\check{r}$$
This is of order $n-1$, the principal part being:
$$2(T\check{r})\Omega^{n-2}T\clab$$
however the coefficient $T\check{r}$ is, according to \ref{12.748}, bounded in absolute value by 
$C\tau$ on ${\cal K}^{\delta}$. It follows that these terms, like the remainders of the 3rd and 4th 
sums, can be absorbed in the estimates. Consider finally the case $m=1$. In this case the highest 
order terms in the 1st sum are the terms $i=1$ and $i=n-1$:
$$(n-1)(\Omega r)\Omega^{n-2}T\clab \ \mbox{and} \ (\Omega^{n-1}r)T\clab$$
The first is of order $n-1$ but with a coefficient bounded in absolute value by $C\tau$, while the 
principal part of the second is:
$$(T\clab)\Omega^{n-1}\check{r}$$
which is also of order $n-1$ but with a coefficient bounded in absolute value by $C\tau$ (\ref{12.748}). 
In the case $m=1$ the highest order terms in the 2nd sum are the terms $i=1$ and $i=n-1$:
$$(n-1)(\Omega\lambdab)\Omega^{n-2}T\check{r} \ \mbox{and} \ (\Omega^{n-1}\lambdab T\check{r}$$
The first is of order $n-1$ but with a coefficient bounded in absolute value by $C\tau$, while 
the principal part of the second is:
$$(T\check{r})\Omega^{n-1}\clab$$
which is also of order $n-1$ but with a coefficient bounded in absolute value by $C\tau$ (\ref{12.748}). 
It follows that also in the case $m=1$ the 1st and 2nd sums in \ref{12.751} can be absorbed in the 
estimates. Furthermore, in the case $m=1$ the highest order terms in the 3rd sum are the terms $i=n-1$ 
and $i=n-2$, 
$$\clab\Omega^{n-1}T\check{r} \ \mbox{and} \ (n-1)(\Omega\clab)\Omega^{n-2}T\check{r}$$
being the corresponding principal parts. The first is of order $n$ but its coefficient is, 
according to \ref{12.745}, bounded in absolute value by $C\tau^2$ in ${\cal K}^{\delta}$. The 
second is of order $n-1$ and, assuming that:
\begin{equation}
|\Omega\check{r}|, |\Omega\clab| \leq C\tau \ : \ \mbox{on ${\cal K}^{\delta}$}
\label{12.752}
\end{equation}
its coefficient is bounded in absolute value by $C\tau$ on ${\cal K}^{\delta}$. Similarly, in the case 
$m=1$ the highest order terms in the 4th sum are the terms $i=n-1$ and $i=n-2$, 
$$\check{r}\Omega^{n-1}T\clab \ \mbox{and} \ (n-1)(\Omega\check{r})\Omega^{n-2}T\clab$$
being the corresponding principal parts. The first is of order $n$ but its coefficient is, according 
to \ref{12.745}, bounded in absolute value by $C\tau^2$ in ${\cal K}^{\delta}$, while the second is 
of order $n-1$ with a coefficient bounded in absolute value by $C\tau$. It follows that these terms, 
like the remainders of the 3rd and 4th sums, can be absorbed in the estimates. 

Consider now $\Omega^{n-m}{\cal M}_{m-2,0}$, the first term on the right in \ref{12.751}. Since 
for $m\geq 3$ ${\cal M}_{m-2,0}$ is of order $m-2$, this term is for $m\geq 3$ of order $n-2$. 
In the following we can ignore terms of order $n-2$ or less.  In the case $m=2$, from \ref{12.747} 
the first term on the right in \ref{12.751} is:
\begin{equation}
-2(T\check{r})\Omega^{n-2}T\clab-2(T\clab)\Omega^{n-2}T\check{r}
\label{12.753}
\end{equation}
up to terms which can be ignored. Both terms here are of order $n-1$ but with coefficients which 
are, according to \ref{12.748} bounded in absolute value by $C\tau$ in ${\cal K}^{\delta}$. So 
these terms can be absorbed in the estimates. In the case $m=1$, from \ref{12.743} the first term 
on the right in \ref{12.751} is:
\begin{eqnarray}
&&\check{r}\Omega^{n-1}T\clab+\clab\Omega^{n-1}T\check{r}\nonumber\\
&&+(n-1)(\Omega\check{r})\Omega^{n-2}T\clab+(n-1)(\Omega\clab)\Omega^{n-2}T\check{r}\nonumber\\
&&+(T\clab)\Omega^{n-1}\check{r}+(T\check{r})\Omega^{n-1}\clab \label{12.754}
\end{eqnarray}
up to terms which can be ignored. Here the first two terms are of order $n$ but, according to 
\ref{12.745}, with coefficients bounded in absolute vale by $C\tau^2$ in ${\cal K}^{\delta}$. 
the next two terms are of order $n-1$ but, according to \ref{12.752}, with coefficients bounded 
in absolute value by $C\tau$ in ${\cal K}^{\delta}$. The last two terms are also of order $n-1$, 
but, according to \ref{12.748}, with coefficients bounded in absolute value by $C\tau$ in 
${\cal K}^{\delta}$. Thus all these terms can be absorbed in the estimates. 

We conclude from the above that for any $m=1,...,n$ the remainder term ${\cal M}_{m-2,n-m}$ in 
\ref{12.750} can be absorbed in the estimates. We can then express \ref{12.750}, multiplying by 
$\sh^{-(n-m)/2}$, in the form:
\begin{eqnarray}
&&\s^{(m,n-m)}\cla=r\s^{(m,n-m)}\clab+\clab\sh^{-(n-m)/2}\Omega^{n-m}T^m\check{r} \nonumber\\
&&\hspace{18mm}+m(Tr)\s^{(m-1,n-m)}\clab+m(T\lambdab)\sh^{-(m-n)/2}\Omega^{n-m}T^{m-1}\check{r} \nonumber\\
&&\label{12.755}
\end{eqnarray}
up to terms which can be absorbed. In estimating $\s^{(m,n-m)}\clab$ in $L^2({\cal K}^{\tau_1})$ 
we are to divide through by $r$. Since $r\sim\tau$ the contribution of the left hand side 
of \ref{12.755} to $\|\s^{(m,n-m)}\clab\|_{L^2({\cal K}^{\tau_1})}$ is bounded by: 
\begin{eqnarray}
&&C\left\{\int_0^{\tau_1}\tau^{-2}\|\s^{(m,n-m)}\cla\|^2_{L^2(S_{\tau,\tau})}d\tau\right\}^{1/2}
\nonumber\\
&&=C\left\{\int_0^{\tau_1}\tau^{-2}\frac{\partial}{\partial\tau}
\left(\|\s^{(m,n-m)}\cla\|^2_{L^2({\cal K}^{\tau})}\right)d\tau\right\}^{1/2}\nonumber\\
&&=C\left\{\tau_1^{-2}\|\s^{(m,n-m)}\cla\|^2_{L^2({\cal K}^{\tau_1})}
+2\int_0^{\tau_1}\tau^{-3}\|\s^{(m,n-m)}\cla\|^2_{L^2({\cal K}^{\tau})}d\tau\right\}^{1/2}\nonumber\\
&&\label{12.756}
\end{eqnarray}
Since $\lambdab/r\sim 1$ along ${\cal K}$, the contribution of the 2nd term on the right hand side of \ref{12.755} to $\|\s^{(m,n-m)}\clab\|_{L^2({\cal K}^{\tau_1})}$ is bounded by: 
\begin{equation}
C\|\Omega^{n-m}T^m\check{r}\|_{L^2({\cal K}^{\tau_1})}
\label{12.757}
\end{equation}
By \ref{12.366}:
\begin{equation}
\|\Omega^{n-m}T^{m-1}\clab\|_{L^2(S_{\tau,\tau})}\leq k\tau^{1/2}
\|\Omega^{n-m}T^m\clab\|_{L^2({\cal K}^{\tau})}
\label{12.758}
\end{equation}
It follows that the contribution of the 3rd term on the right in \ref{12.755} to 
$\|\s^{(m,n-m)}\clab\|_{L^2({\cal K}^{\tau_1})}$ is bounded by:
\begin{equation}
Cm\left\{\int_0^{\tau_1}\tau^{-1}\|\s^{(m,n-m)}\clab\|^2_{L^2({\cal K}^{\tau})}d\tau\right\}^{1/2}
\label{12.759}
\end{equation}
up to terms which can be absorbed. Again by \ref{12.366}:
\begin{equation}
\|\Omega^{n-m}T^{m-1}\check{r}\|_{L^2(S_{\tau,\tau})}\leq k\tau^{1/2}
\|\Omega^{n-m}T^m\check{r}\|_{L^2({\cal K}^{\tau})}
\label{12.760}
\end{equation}
It follows that the contribution of the 4th term on the right in \ref{12.755} to 
$\|\s^{(m,n-m)}\clab\|_{L^2({\cal K}^{\tau_1})}$ is bounded by:
\begin{equation}
Cm\left\{\int_0^{\tau_1}\tau^{-1}\|\Omega^{n-m}T^m\check{r}\|^2_{L^2({\cal K}^{\tau})}
d\tau\right\}^{1/2}
\label{12.761}
\end{equation}

We must now derive an appropriate estimate for $\|\Omega^{n-m}T^m\check{r}\|_{L^2({\cal K}^{\tau_1})}$, 
or equivalently for $\|E^{n-m}T^m r-E_N^{n-m}T^m r_N\|_{L^2({\cal K}^{\tau_1})}$, for 
any $\tau_1\in(0,\delta]$. By \ref{10.792} and 
\ref{10.708} the leading part of $E^{n-m}T^m r-E_N^{n-m}T^m r_N$ is:
$$j(\kappa,\epb)(E^{n-m}T^m\epb-E_N^{n-m}T^m\epb_N)$$
and the leading part of this is $j(\kappa,\epb)$ times:
\begin{eqnarray}
&&\Nb^\mu E^{n-m}T^m\triangle\beta_\mu-\Nb_N^\mu E_N^{n-m}T^m\triangle_N\beta_{\mu,N}\nonumber\\
&&+\triangle\beta_\mu E^{n-m}T^m\Nb^\mu-\triangle_N\beta_{\mu,N} E_N^{n-m}T^m\Nb_N^\mu\nonumber\\
&&+m(T\triangle\beta_\mu)E^{n-m}T^{m-1}\Nb^\mu-m(T\triangle_N\beta_{\mu,N})E_N^{n-m}T^{m-1}\Nb_N^\mu 
\nonumber\\
&&\label{12.762}
\end{eqnarray}
Here we have three difference terms. The third difference term is of one order lower than the 
second difference term, however whereas the coefficients $\triangle\beta_\mu$, 
$\triangle_N\beta_{\mu,N}$ in the second difference are bounded in absolute value by $C\tau$, 
those in the third difference $T\triangle\beta_\mu$, $T\triangle_N\beta_{\mu,N}$ are only 
bounded in absolute value by $C$. For this reason the third difference term must be included 
in the leading part. 

The first of the difference terms \ref{12.762} is:
\begin{eqnarray}
&&\Nb^\mu E^{n-m}T^m\triangle\beta_\mu-\Nb_N^\mu E_N^{n-m}T^m\triangle_N\beta_{\mu,N}\nonumber\\
&&\hspace{20mm}=\left.(\Nb^\mu E^{n-m}T^m\beta_\mu-\Nb_N^\mu E_N^{n-m}T^m\beta_{\mu,N})
\right|_{{\cal K}}\nonumber\\
&&\hspace{20mm}-\left(\Nb^\mu E^{n-m}T^m\beta^\prime_\mu\circ(f,w,\psi)
-\Nb_N^\mu E_N^{n-m}T^m\beta^\prime_\mu\circ(f_N,w_N,\psi_N)\right)\nonumber\\
&&\label{12.763}
\end{eqnarray}
Applying \ref{12.366} to the function $\Nb^\mu E^{n-m}T^m\beta_\mu-\Nb_N^\mu E_N^{n-m}T^m\beta_{\mu,N}$ 
on ${\cal K}$ we obtain:
\begin{eqnarray}
&&\|\Nb^\mu E^{n-m}T^m\beta_\mu-\Nb_N^\mu E_N^{n-m}T^m\beta_{\mu,N}\|_{L^2(S_{\tau,\tau})}
\nonumber\\
&&\hspace{20mm}\leq k\tau^{1/2}\|T(\Nb^\mu E^{n-m}T^m\beta_\mu-\Nb_N^\mu E_N^{n-m}T^m\beta_{\mu,N})\|
_{L^2({\cal K}^{\tau})}\nonumber\\
&&\label{12.764}
\end{eqnarray}
Up to terms of order $n$ with vanishing $n$th order acoustical part 
$T(\Nb^\mu E^{n-m}T^m\beta_\mu-\Nb_N^\mu E_N^{n-m}T^m\beta_{\mu,N})$ is equal to:
\begin{eqnarray}
&&\Nb^\mu T(E^{n-m}T^m\beta_\mu-E_N^{n-m}T^m\beta_{\mu,N})=
\Nb^\mu L(E^{n-m}T^m\beta_\mu-E_N^{n-m}T^m\beta_{\mu,N})\nonumber\\
&&\hspace{42mm}+\Nb^\mu\Lb(E^{n-m}T^m\beta_\mu-E_N^{n-m}T^m\beta_{\mu,N})
\label{12.765}
\end{eqnarray}
In regard to the first term on the right, using the conjugate of \ref{12.605} and its analogue for 
the $N$th approximants we find that this term is equal to:
\begin{equation}
\Nb^\mu E^{n-m}T^m L\beta_\mu-\Nb_N^\mu E_N^{n-m}T^m L_N\beta_{\mu,N}
+2\ss_{\Nb}(\rhob\cmub_{m-1,n-m}-\rho\cmu_{m-1,n-m})
\label{12.766}
\end{equation}
up to terms which can be ignored. In view of the fact that 
$\Nb^\mu L\beta_\mu=\lambdab E^mu E\beta_\mu$, while by the 2nd of \ref{12.66} and by \ref{12.108}, 
\ref{12.111}
\begin{equation}
[E^{n-m}T^m\Nb^\mu]_{P.A.}=2\mub_{m-1,n-m}(E^\mu-\pi\Nb^\mu)
\label{12.767}
\end{equation}
we can express:
\begin{equation}
\Nb^\mu E^{n-m}T^m L\beta_\mu=E^{n-m}T^m(\lambdab E^\mu E\beta_\mu)
-2\rho(\ss_N-c\pi\sss)\mub_{m-1,n-m}
\label{12.768}
\end{equation}
up to terms which can be ignored. Moreover by \ref{12.582} the 1st term on the right in \ref{12.768} 
is:
\begin{equation}
\lambdab E^\mu E^{n-m}T^m E\beta_\mu+\sss E^{n-m}T^m\lambdab
+\rho(\ss_N\mub_{m-1,n-m}+\ss_{\Nb}\mu_{m-1,n-m})
\label{12.769}
\end{equation}
up to terms which can be ignored. As for the 1st term here, it is given by \ref{12.581} with 
the overall factor $\lambda$ replaced by $\lambdab$. Similar formulas hold for the $N$th approximants. 
It follows that up to terms which can be ignored the principal difference in \ref{12.766} is 
given by:
\begin{eqnarray}
&&\lambdab\s^{(E;m,n-m)}\csxi+\sss \s^{(m,n-m)}\clab\nonumber\\
&&-\lambdab\sss\rho\cdot\left\{ \begin{array}{lll}
\s^{(n-1)}\ctchi &:& \mbox{if $m=1$}\\
2\cmu_{m-2,n-m+1} &:& \mbox{if $m\geq 2$}
\end{array} \right.\nonumber\\
&&-\lambdab\sss\rhob\cdot\left\{ \begin{array}{lll}
\s^{(n-1)}\ctchib &:& \mbox{if $m=1$}\\
2\cmub_{m-2,n-m+1} &:& \mbox{if $m\geq 2$}
\end{array} \right. \nonumber\\
&&+\rho(\ss_N\cmub_{m-1,n-m}+\ss_{\Nb}\cmu_{m-1,n-m})\nonumber\\
&&-2\rho(\ss_N-c\pi\sss)\cmub_{m-1,n-m}
\label{12.770}
\end{eqnarray}
Consequently, the $L^2$ norm on ${\cal K}^{\tau}$ of the first term on the right in \ref{12.765} 
is bounded by:
\begin{eqnarray}
&&\|\lambdab\s^{(E;m,n-m)}\csxi||_{L^2({\cal K}^{\tau})}+C\|\s^{(m,n-m)}\clab\|_{L^2({\cal K}^{\tau})}
\nonumber\\
&&+C\tau^2\left\{ \begin{array}{lll} 
\|\s^{(n-1)}\ctchi\|_{L^2({\cal K}^{\tau})} &:& \mbox{if $m=1$}\\
2\|\cmu_{m-2,n-m+1}\|_{L^2({\cal K}^{\tau})} &:& \mbox{if $m\geq 2$} \end{array}\right.\nonumber\\
&&+C\tau^2\left\{ \begin{array}{lll} 
\|r\s^{(n-1)}\ctchib\|_{L^2({\cal K}^{\tau})} &:& \mbox{if $m=1$}\\
2\|r\cmub_{m-2,n-m+1}\|_{L^2({\cal K}^{\tau})} &:& \mbox{if $m\geq 2$} \end{array}\right.\nonumber\\
&&+C\|r\cmub_{m-1,n-m}\|_{L^2({\cal K}^{\tau})}+C\tau\|\cmu_{m-1,n-m}\|_{L^2({\cal K}^{\tau})}
\label{12.771}
\end{eqnarray}
Since along ${\cal K}$ $a\sim\tau^3$, in regard to the 1st of \ref{12.771} we have:
\begin{eqnarray}
&&\|\lambdab\s^{(E;m,n-m)}\csxi\|^2_{L^2({\cal K}^{\tau})}\leq 
C\int_0^{\tau}\tau^{\prime-1}\|\sqrt{a}\s^{(E;m,n-m)}\csxi\|^2_{L^2(S_{\tau^\prime,\tau^\prime})}
d\tau^\prime\nonumber\\
&&\hspace{30mm}=C\int_0^{\tau}\tau^{\prime-1}\frac{\partial}{\partial\tau^\prime}
\left(\|\sqrt{a}\s^{(E;m,n-m)}\csxi\|^2_{L^2({\cal K}^{\tau^\prime})}\right)d\tau^\prime\nonumber\\
&&\hspace{10mm}=C\left\{\tau^{-1}\|\sqrt{a}\s^{(E;m,n-m)}\csxi\|^2_{L^2({\cal K}^{\tau})}
+\int_0^{\tau}\tau^{\prime-2}\|\sqrt{a}\s^{(E;m,n-m)}\csxi\|^2_{L^2({\cal K}^{\tau^\prime})}d\tau^\prime
\right\}\nonumber\\
&&\hspace{30mm}\leq C\s^{(E;m,n-m)}\cA(\tau_1)\tau^{2c_m-1}
\label{12.772}
\end{eqnarray}
for all $\tau\in (0,\tau_1]$, by \ref{9.a35} and \ref{9.299}. Therefore the contribution of this term 
through \ref{12.764} to $\|E^{n-m}T^m r-E_N^{n-m}T^m r_N\|_{L^2({\cal K}^{\tau_1})}$ is bounded by:
\begin{equation}
C\sqrt{\s^{(E;m,n-m)}\cA(\tau_1)}\left(\int_0^{\tau_1}\tau^{2c_m}d\tau\right)^{1/2}
\leq\frac{C\sqrt{\s^{(E;m,n-m)}\cA(\tau_1)}}{\sqrt{2c_m+1}}\cdot\tau_1^{c_m}\cdot\tau_1^{1/2}
\label{12.773}
\end{equation}
Also, the contribution of the 2nd of \ref{12.771} through \ref{12.764} to 
$\|E^{n-m}T^m r-E_N^{n-m}T^m r_N\|_{L^2({\cal K}^{\tau_1})}$ is bounded by: 
\begin{equation}
C\left(\int_0^{\tau_1}\tau\|\s^{(m,n-m)}\clab\|^2_{L^2({\cal K}^{\tau})}d\tau\right)^{1/2}
\label{12.774}
\end{equation}
In regard to the second term on the right in \ref{12.765} we write:
\begin{eqnarray}
&&N^\mu\Lb(E^{n-m}T^m\beta_\mu-E_N^{n-m}T^m\beta_{\mu,N})
+\ogamma\Nb^\mu\Lb(E^{n-m}T^m\beta_\mu-E_N^{n-m}T^m\beta_{\mu,N})\nonumber\\
&&\hspace{56mm}=\s^{(Y;m,n-m)}\cxi_{\Lb}
\label{12.775}
\end{eqnarray}
Since along ${\cal K}$ we have $\ogamma\sim\tau$, the contribution of the right hand side of 
\ref{12.775} to 
$\|T(\Nb^\mu E^{n-m}T^m\beta_\mu-\Nb_N^\mu E_N^{n-m}T^m\beta_{\mu,N})\|^2_{L^2({\cal K}^{\tau})}$ 
is bounded by a constant multiple of:
\begin{eqnarray}
&&\int_0^{\tau}\tau^{\prime-2}\|\s^{(Y;m,n-m)}\cxi_{\Lb}\|^2_{L^2(S_{\tau^\prime,\tau^\prime})}d\tau^\prime\nonumber\\
&&=\int_0^{\tau}\tau^{\prime-2}\frac{\partial}{\partial\tau^\prime}
\left(\|\s^{(Y;m,n-m)}\cxi_{\Lb}\|^2_{L^2({\cal K}^{\tau^\prime})}\right)d\tau^\prime\nonumber\\
&&=\tau^{-2}\|\s^{(Y;m,n-m)}\cxi_{\Lb}\|^2_{L^2({\cal K}^{\tau})}
+2\int_0^{\tau}\tau^{\prime-3}\|\s^{(Y;m,n-m)}\cxi_{\Lb}\|^2_{L^2({\cal K}^{\tau^\prime})}d\tau^\prime
\nonumber\\
&&\leq C\s^{(Y;0,n)}\cA(\tau_1)\tau^{2c_m-2}
\label{12.776}
\end{eqnarray}
for all $\tau\in [0,\tau_1]$, by \ref{9.a39} and \ref{9.299}. The corresponding contribution to 
$\|E^{n-m}T^m r-E_N^{n-m}T^m r_N\|_{L^2({\cal K}^{\tau_1})}$ is then bounded by:
\begin{equation}
C\sqrt{\s^{(Y;m,n-m)}\cA(\tau_1)}\left(\int_0^{\tau_1}\tau^{2c_m-1}d\tau\right)^{1/2}
\leq\frac{C\sqrt{\s^{(Y;m,n-m)}\cA(\tau_1)}}{\sqrt{2c_m}}\cdot\tau_1^{c_m}
\label{12.777}
\end{equation}
In regard to the 1st term on the left in \ref{12.775}, this term is given, up to terms which can be 
ignored, by the conjugates of \ref{12.766}, \ref{12.770}, that is by:
\begin{eqnarray}
&&\lambda\s^{(E;m,n-m)}\csxi+\sss \s^{(m,n-m)}\cla\nonumber\\
&&-\lambda\sss\rhob\cdot\left\{ \begin{array}{lll}
\s^{(n-1)}\ctchib &:& \mbox{if $m=1$}\\
2\cmub_{m-2,n-m+1} &:& \mbox{if $m\geq 2$}
\end{array} \right.\nonumber\\
&&-\lambda\sss\rho\cdot\left\{ \begin{array}{lll}
\s^{(n-1)}\ctchi &:& \mbox{if $m=1$}\\
2\cmu_{m-2,n-m+1} &:& \mbox{if $m\geq 2$}
\end{array} \right. \nonumber\\
&&+\rhob(\ss_{\Nb}\cmu_{m-1,n-m}+\ss_N\cmub_{m-1,n-m})-2\rhob(\ss_{\Nb}-c\pi\sss)\cmu_{m-1,n-m}
\nonumber\\
&&+2\ss_N(\rho\cmu_{m-1,n-m}-\rhob\cmub_{m-1,n-m})\label{12.778}
\end{eqnarray}
Since along ${\cal K}$ $\lambda/\ogamma=\lambdab$, the contribution of the 1st of \ref{12.778} 
coincides with that of the 1st of \ref{12.770}. Since along ${\cal K}$ $\ogamma\sim\tau$, the 
contribution of the 2nd of \ref{12.778} to $\|T(\Nb^\mu E^{n-m}T^m\beta_\mu-\Nb_N^\mu E_N^{n-m}T^m
\beta_{\mu,N})\|^2_{L^2({\cal K}^{\tau})}$ is bounded by a constant multiple of:
\begin{eqnarray}
&&\int_0^{\tau}\tau^{\prime-2}\frac{\partial}{\partial\tau^\prime}\left(\|\s^{(m,n-m)}\cla\|^2
_{L^2({\cal K}^{\tau^\prime})}\right)d\tau^\prime\nonumber\\
&&=\tau^{-2}\|\s^{(m,n-m)}\cla\|^2_{L^2({\cal K}^{\tau})}+2\int_0^{\tau}\tau^{\prime-3}
\|\s^{(m,n-m)}\cla\|^2_{L^2({\cal K}^{\tau^\prime})}d\tau^\prime\nonumber\\
&&\label{12.779}
\end{eqnarray}
The corresponding contribution to $\|E^{n-m}T^m r-E_N^{n-m}T^m r_N\|_{L^2({\cal K}^{\tau_1})}$ is then bounded by:
\begin{equation}
\left\{\int_0^{\tau_1}\left(\tau^{-1}\|\s^{(m,n-m)}\cla\|^2_{L^2({\cal K}^{\tau})}
+2\tau\int_0^{\tau}\tau^{\prime-3}\|\s^{(m,n-m)}\cla\|^2_{L^2({\cal K}^{\tau^\prime})}d\tau^\prime\right)
d\tau\right\}^{1/2}
\label{12.780}
\end{equation}
Finally, the contribution of the remaining terms in \ref{12.778} to \\
$\|T(\Nb^\mu E^{n-m}T^m\beta_\mu-\Nb_N^\mu E_N^{n-m}T^m\beta_{\mu,N})\|_{L^2({\cal K}^{\tau})}$ 
is bounded by:
\begin{eqnarray}
&&C\tau^2\left\{ \begin{array}{lll} 
\|\s^{(n-1)}\ctchi\|_{L^2({\cal K}^{\tau})} &:& \mbox{if $m=1$}\\
2\|\cmu_{m-2,n-m+1}\|_{L^2({\cal K}^{\tau})} &:& \mbox{if $m\geq 2$} \end{array}\right.\nonumber\\
&&+C\tau^2\left\{ \begin{array}{lll} 
\|r\s^{(n-1)}\ctchib\|_{L^2({\cal K}^{\tau})} &:& \mbox{if $m=1$}\\
2\|r\cmub_{m-2,n-m+1}\|_{L^2({\cal K}^{\tau})} &:& \mbox{if $m\geq 2$} \end{array}\right.\nonumber\\
&&+C\|\cmu_{m-1,n-m}\|_{L^2({\cal K}^{\tau})}+C\|r\cmub_{m-1,n-m}\|_{L^2({\cal K}^{\tau})}
\label{12.781}
\end{eqnarray}

We turn to the second term on the right in \ref{12.763} (compare with \ref{10.811}). By \ref{10.817} 
with $(m,n-m)$ in the role of $(m-1,l+2)$ the leading part of the second term on the right in 
\ref{12.763} is:
\begin{eqnarray}
&&-\Nb^\mu\sh^{-(n-m)/2}\left\{\left(\frac{\partial\beta^\prime_\mu}{\partial t}\right)\circ(f,w,\psi)
\cdot \Omega^{n-m}\left(\tau^2 T^m\chf\right.\right.\nonumber\\
&&\hspace{40mm}\left.+2m\tau T^{m-1}\chf+m(m-1)T^{m-2}\chf\right)\nonumber\\
&&\hspace{20mm}+\left(\frac{\partial\beta^\prime_\mu}{\partial u^\prime}\right)\circ(f,w,\psi) \cdot 
\Omega^{n-m}\left(\tau T^m\cv+m T^{m-1}\cv\right)\nonumber\\
&&\hspace{20mm}+\left(\frac{\partial\beta^\prime_\mu}{\partial\vartheta^\prime}\right)\circ(f,w,\psi) 
\cdot \Omega^{n-m}\left(\tau^3 T^m\cga+3m\tau^2 T^{m-1}\cga\right.\nonumber\\
&&\hspace{20mm}\left.\left.+3m(m-1)\tau T^{m-2}\cga
+m(m-1)(m-2)T^{m-3}\cga\right)\right\}\nonumber\\
&&\label{12.782}
\end{eqnarray}
It follows that the leading part of the second term on the right in \ref{12.763} is pointwise 
bounded by:
\begin{eqnarray}
&&C\left(\tau^2|\Omega^{n-m}T^m\chf|+2m\tau|\Omega^{n-m}T^{m-1}\chf|+m(m-1)|\Omega^{n-m}T^{m-2}\chf|
\right)\nonumber\\
&&+C\left(\tau|\Omega^{n-m}T^m\cv|+m|\Omega^{n-m}T^{m-2}\cv|\right)\nonumber\\
&&+C\left(\tau^3|\Omega^{n-m}T^m\cga|+3m\tau^2|\Omega^{n-m}T^{m-1}\cga|\right.\nonumber\\
&&\hspace{7mm}\left.
+3m(m-1)\tau|\Omega^{n-m}T^{m-2}\cga|+m(m-1)(m-2)|\Omega^{n-m}T^{m-3}\cga|\right) \nonumber\\
&&\label{12.783}
\end{eqnarray}
Let us revisit the inequality \ref{12.366}, which is satisfied for any function $f$ on ${\cal K}$ 
which vanishes on $\partial_-{\cal K}=S_{0,0}$. This inequality implies:
\begin{equation}
\|f\|_{L^2({\cal K}^{\tau})}\leq k\left(\int_0^{\tau}\tau^\prime\|Tf\|^2_{L^2({\cal K}^{\tau^\prime})}
d\tau^\prime\right)^{1/2}
\label{12.784}
\end{equation}
If also $Tf$ vanishes on $S_{0,0}$, the inequality \ref{12.784} holds with $Tf$ in the role of $f$. 
Substituting the resulting bound for $\|Tf\|_{L^2({\cal K}^{\tau})}$ in \ref{12.366} yields:
\begin{equation}
\|f\|_{L^2(S_{\tau,\tau})}\leq k^2\tau^{1/2}\left(\int_0^{\tau}\tau^\prime
\|T^2 f\|^2_{L^2({\cal K}^{\tau^\prime})}d\tau^\prime\right)^{1/2}
\label{12.785}
\end{equation}
Applying inequality \ref{12.366} with $\Omega^{n-m}T^m\cv$ in the role of $f$ and inequality 
\ref{12.785} with $\Omega^{n-m}T^{m-1}\cv$ in the role of $f$ we conclude that the 2nd of 
\ref{12.783} is bounded in $L^2(S_{\tau,\tau})$ by a constant multiple of:
\begin{eqnarray}
&&\tau^{3/2}\|\Omega^{n-m}T^{m+1}\cv\|_{L^2({\cal K}^{\tau})}\nonumber\\
&&+m\tau^{1/2}\left(\int_0^{\tau}\tau^\prime\|\Omega^{n-m}T^{m+1}\cv\|^2_{L^2({\cal K}^{\tau^\prime})}
d\tau^\prime\right)^{1/2}\label{12.786}
\end{eqnarray}
If also $T^2 f$ vanishes on $S_{0,0}$, the inequality \ref{12.784} holds with $T^2 f$ in the role of 
$f$. Substituting the resulting bound for $\|T^2 f\|_{L^2({\cal K}^{\tau})}$ in \ref{12.785} yields:
\begin{equation}
\|f\|_{L^2(S_{\tau,\tau})}\leq k^3\tau^{1/2}\left[\int_0^{\tau}\tau^\prime\left(\int_0^{\tau^\prime}
\tau^{\prime\prime}\|T^3 f\|^2_{L^2({\cal K}^{\tau^{\prime\prime}})}d\tau^{\prime\prime}\right)d\tau^\prime\right]^{1/2}
\label{12.787}
\end{equation}
Applying inequality \ref{12.366} with $\Omega^{n-m}T^m\chf$ in the role of $f$, inequality \ref{12.785} 
with $\Omega^{n-m}T^{m-1}\chf$ in the role of $f$ and inequality \ref{12.787} with 
$\Omega^{n-m}T^{m-2}\chf$ in the role of $f$, we conclude that the 1st of \ref{12.783} is bounded in 
$L^2(S_{\tau,\tau})$ by a constant multiple of:
\begin{eqnarray}
&&\tau^{5/2}\|\Omega^{n-m}T^{m+1}\chf\|_{L^2({\cal K}^{\tau})}\nonumber\\
&&+2m\tau^{3/2}\left(\int_0^{\tau}\tau^\prime\|\Omega^{n-m}T^{m+1}\chf\|^2_{L^2({\cal K}^{\tau^\prime})}
d\tau^\prime\right)^{1/2}\nonumber\\
&&+m(m-1)\tau^{1/2}\left[\int_0^{\tau}\tau^\prime\left(\int_0^{\tau^\prime}\tau^{\prime\prime}
\|\Omega^{n-m}T^{m+1}\chf\|^2_{L^2({\cal K}^{\tau^{\prime\prime}})}d\tau^{\prime\prime}\right)
d\tau^\prime\right]^{1/2}\nonumber\\
&&\label{12.788}
\end{eqnarray}
If finally $T^3 f$ vanishes on $S_{0,0}$ as well, the inequality \ref{12.784} holds with $T^3 f$ in the 
role of $f$. Substituting the resulting bound for $\|T^3 f\|_{L^2({\cal K}^{\tau})}$ in \ref{12.787} 
yields:
\begin{equation}
\|f\|_{L^2(S_{\tau,\tau})}\leq k^4\tau^{1/2}\left\{\int_0^{\tau}\tau^\prime\left[\int_0^{\tau^\prime}
\tau^{\prime\prime}\left(\int_0^{\tau^{\prime\prime}}\tau^{\prime\prime\prime}
\|T^4 f\|^2_{L^2({\cal K}^{\tau^{\prime\prime\prime}})}d\tau^{\prime\prime\prime}\right)
d\tau^{\prime\prime}\right]d\tau^\prime\right\}^{1/2}
\label{12.789}
\end{equation}
Applying then inequality \ref{12.366} with $\Omega^{n-m}T^m\cga$ in the role of $f$, inequality 
\ref{12.785} with $\Omega^{n-m}T^{m-1}\cga$ in the role of $f$, inequality \ref{12.787} with 
$\Omega^{n-m}T^{m-2}\cga$ in the role of $f$, and inequality \ref{12.789} with 
$\Omega^{n-m}T^{m-3}\cga$ in the role of $f$, we conclude that the 3rd of \ref{12.783} is bounded 
in $L^2(S_{\tau,\tau})$ by a constant multiple of:
\begin{eqnarray}
&&\tau^{7/2}\|\Omega^{n-m}T^{m+1}\cga\|_{L^2({\cal K}^{\tau})}\nonumber\\
&&+3m\tau^{5/2}\left(\int_0^{\tau}\tau^\prime\|\Omega^{n-m}T^{m+1}\cga\|^2_{L^2({\cal K}^{\tau^\prime})}
d\tau^\prime\right)^{1/2}\nonumber\\
&&+3m(m-1)\tau^{3/2}\left[\int_0^{\tau}\tau^\prime\left(\int_0^{\tau^\prime}\tau^{\prime\prime}
\|\Omega^{n-m}T^{m+1}\cga\|^2_{L^2({\cal K}^{\tau^{\prime\prime}})}d\tau^{\prime\prime}\right)
d\tau^\prime\right]^{1/2}\nonumber\\
&&+m(m-1)(m-2)\tau^{1/2}\cdot\nonumber\\
&&\cdot\left\{\int_0^{\tau}\tau^\prime\left[\int_0^{\tau^\prime}
\tau^{\prime\prime}\left(\int_0^{\tau^{\prime\prime}}\tau^{\prime\prime\prime}
\|\Omega^{n-m}T^{m+1}\cga\|^2_{L^2({\cal K}^{\tau^{\prime\prime\prime}})}d\tau^{\prime\prime\prime}\right)
d\tau^{\prime\prime}\right]d\tau^\prime\right\}^{1/2}\nonumber\\
&&\label{12.790}
\end{eqnarray}

Consider finally the last two of the differences \ref{12.762}. In view of \ref{12.767} and the boundary 
condition $E^\mu\triangle\beta_\mu=0$, the principal acoustical part of the second difference is:
\begin{equation}
-2\pi\epb\cmub_{m-1,n-m}
\label{12.791}
\end{equation}
Since $|\epb|\leq C\tau$, this is bounded in $L^2({\cal K}^{\tau_1})$ by:
\begin{equation}
C\tau_1\|\cmub_{m-1,n-m}\|_{L^2({\cal K}^{\tau_1})}
\label{12.792}
\end{equation}
Also, by \ref{12.366} the third of the differences \ref{12.762} is bounded in $L^2(S_{\tau,\tau})$ 
by:
\begin{equation}
Cm\tau^{1/2}\|\cmub_{m-1,n-m}\|_{L^2({\cal K}^{\tau})}
\label{12.793}
\end{equation}
therefore in $L^2({\cal K}^{\tau})$ by:
\begin{equation}
Cm\left(\int_0^{\tau_1}\tau\|\cmub_{m-1,n-m}\|^2_{L^2({\cal K}^{\tau})}d\tau\right)^{1/2}
\label{12.794}
\end{equation}

We shall first summarize the above results in regard to the estimate for 
$\|\Omega^{n-m}T^m\check{r}\|_{L^2({\cal K}^{\tau_1})}$ which they imply. Let us define, in analogy 
with \ref{12.399}, \ref{12.400}, for arbitrary $\tau_1\in(0,\delta]$, the following quantities in regard to the boundary values on 
${\cal K}$ of $\s^{(m,n-m)}\cla$ and $\s^{(m,n-m)}\clab$:
\begin{equation}
\s^{(m,n-m)}\oLa(\tau_1)=\sup_{\tau\in(0,\tau_1]}\left\{\tau^{-c_m-1}
\|\s^{(m,n-m)}\cla\|_{L^2({\cal K}^{\tau})}\right\}
\label{12.795}
\end{equation}
\begin{equation}
\s^{(m,n-m)}\oLab(\tau_1)=\sup_{\tau\in(0,\tau_1]}\left\{\tau^{-c_m}
\|\s^{(m,n-m)}\clab\|_{L^2({\cal K}^{\tau})}\right\}
\label{12.796}
\end{equation}
In view of the definition \ref{12.796}, \ref{12.774} is bounded by:
\begin{equation}
C\frac{\tau_1^{c_m+1}}{\sqrt{2c_m+2}}\s^{(m,n-m)}\oLab(\tau_1)
\label{12.797}
\end{equation}
and in view of the definition \ref{12.795}, \ref{12.780} is bounded by:
\begin{equation}
C\frac{\tau_1^{c_m+1}}{\sqrt{2c_m+2}}\s^{(m,n-m)}\oLa(\tau_1)
\label{12.798}
\end{equation}
According to the definitions \ref{12.397}, \ref{12.398} we have, for all $\tau\in(0,\tau_1]$:
\begin{eqnarray}
&&\|\s^{(n-1)}\ctchi\|_{L^2({\cal K}^{\tau})}\leq \tau^{c_0}\s^{(n-1)}\oX(\tau_1)\label{12.799}\\
&&\|r\s^{(n-1)}\ctchib\|_{L^2({\cal K}^{\tau})}\leq \tau^{c_0}\s^{(n-1)}\oXb(\tau_1)\label{12.800}
\end{eqnarray}
Moreover, \ref{12.567} implies:
\begin{eqnarray}
&&\|\cmu_{m-1,n-m}\|_{L^2({\cal K}^{\tau})}\leq \sum_{j=0}^{m-1}(\Pi\tau^2)^j \ 
\|\s^{(m-1-j,n-m+j+1)}\cla\|_{L^2({\cal K}^{\tau})}\nonumber\\
&&\hspace{32mm}+\frac{1}{2}(\Pi\tau^2)^m \ \|\s^{(n-1)}\ctchi\|_{L^2({\cal K}^{\tau})} 
\label{12.801}
\end{eqnarray}
hence in view of the definitions \ref{12.397}, \ref{12.399}, \ref{12.795} we have, for all 
$\tau\in(0,\tau_1]$:
\begin{eqnarray}
&&\|\cmu_{m-1,n-m}\|_{L^2({\cal K}^{\tau})}\leq \sum_{j=0}^{m-1}(\Pi\tau^2)^j \ \tau^{c_{m-1-j}+1}
\s^{(m-1-j,n-m+j+1)}\oLa(\tau_1)\nonumber\\
&&\hspace{32mm}+\frac{1}{2}(\Pi\tau^2)^m \ \tau^{c_0}\s^{(n-1)}\oX(\tau_1)
\label{12.802}
\end{eqnarray}
Also, \ref{12.568} implies:
\begin{eqnarray}
&&\|\cmub_{m-1,n-m}\|_{L^2({\cal K}^{\tau})}\leq \sum_{j=0}^{m-1}(\Pib\tau)^j \ 
\|\s^{(m-1-j,n-m+j+1)}\clab\|_{L^2({\cal K}^{\tau})}\nonumber\\
&&\hspace{32mm}+C(\Pib\tau)^{m-1} \ \|r\s^{(n-1)}\ctchib\|_{L^2({\cal K}^{\tau})}
\label{12.803}
\end{eqnarray}
hence in view of the definitions \ref{12.398}, \ref{12.400}, \ref{12.796} we have, for all 
$\tau\in(0,\tau_1]$:
\begin{eqnarray}
&&\|\cmub_{m-1,n-m}\|_{L^2({\cal K}^{\tau})}\leq \sum_{j=0}^{m-1}(\Pib\tau)^j \ 
\tau^{c_{m-1-j}}\s^{(m-1-j,n-m+j+1)}\oLab(\tau_1)\nonumber\\
&&\hspace{32mm}+C(\Pib\tau)^{m-1} \ \tau^{c_0}\s^{(n-1)}\oXb(\tau_1) 
\label{12.804}
\end{eqnarray}

Next we define, in analogy with \ref{12.394} - \ref{12.396}, for arbitrary $\tau_1\in(0,\delta]$, 
the following quantities in regard to the transformation functions:
\begin{eqnarray}
&&\s^{(m+1,n-m)}P(\tau_1)=\sup_{\tau\in(0,\tau_1]}\left\{\tau^{-c_m+2}\|\Omega^{n-m}T^{m+1}\chf\|
_{L^2({\cal K}^{\tau})}\right\}\label{12.805}\\
&&\s^{(m+1,n-m)}V(\tau_1)=\sup_{\tau\in(0,\tau_1]}\left\{\tau^{-c_m+2}\|\Omega^{n-m}T^{m+1}\cv\|
_{L^2({\cal K}^{\tau})}\right\}\label{12.806}\\
&&\s^{(m+1,n-m)}\Ga(\tau_1)=\sup_{\tau\in(0,\tau_1]}\left\{\tau^{-c_m+2}\|\Omega^{n-m}T^{m+1}\cga\|
_{L^2({\cal K}^{\tau})}\right\}\label{12.807}
\end{eqnarray}
Then for all $\tau\in(0,\tau_1]$ \ref{12.788} is bounded by:
\begin{equation}
\tau^{c_m+\frac{1}{2}}\left\{1+\frac{2m}{\sqrt{2c_m-2}}+\frac{m(m-1)}{\sqrt{2c_m(2c_m-2)}}\right\}
\s^{(m+1,n-m)}P(\tau_1)
\label{12.808}
\end{equation}
\ref{12.786} is bounded by:
\begin{equation}
\tau^{c_m-\frac{1}{2}}\left\{1+\frac{m}{\sqrt{2c_m-2}}\right\}\s^{(m+1,n-m)}V(\tau_1)
\label{12.809}
\end{equation}
and \ref{12.790} is bounded by:
\begin{eqnarray}
&&\tau^{c_m+\frac{3}{2}}\left\{1+\frac{3m}{\sqrt{2c_m-2}}+\frac{3m(m-1)}{\sqrt{2c_m(2c_m-2)}}
\right.\nonumber\\
&&\hspace{15mm}\left.+\frac{m(m-1)(m-2)}{\sqrt{(2c_m+2)2c_m(2c_m-2)}}\right\}\s^{(m+1,n-m)}\Ga(\tau_1)
\label{12.810}
\end{eqnarray}
Choosing:
\begin{equation}
c_n-2\geq n^2
\label{12.811}
\end{equation}
the coefficients in \ref{12.808} - \ref{12.810} are bounded by a constant for all $m=1,...,n$. 
Recalling that the exponents $c_m$ are non-increasing with $m$ , we deduce from the above the following 
estimate for $\|\Omega^{n-m}T^m\check{r}\|_{L^2({\cal K}^{\tau_1})}$, for all $\tau_1\in(0,\delta]$: 
\begin{equation}
\|\Omega^{n-m}T^m\check{r}\|_{L^2({\cal K}^{\tau_1})}\leq \tau_1^{c_m}\s^{(m,n-m)}R(\tau_1)
\label{12.812}
\end{equation}
where:
\begin{eqnarray}
&&\s^{(m,n-m)}R(\tau_1)=C\frac{\tau_1}{\sqrt{2c_m+2}}\left(\s^{(m,n-m)}\oLab(\tau_1)+\s^{(m,n-m)}\oLa(\tau_1)\right)\nonumber\\
&&\hspace{10mm}+C\tau_1\left\{\sum_{j=0}^{m-1}(\Pib\tau_1)^j \ 
\s^{(m-1-j,n-m+j+1)}\oLab(\tau_1)+(\Pib\tau_1)^{m-1} \ \s^{(n-1)}\oXb(\tau_1)\right\}\nonumber\\
&&+C\tau_1^2\left\{\sum_{j=0}^{m-1}\frac{(\Pi\tau_1^2)^j}{\sqrt{2c_m+4j+4}}
\s^{(m-1-j,n-m+j+1)}\oLa(\tau_1)
+\frac{(\Pi\tau_1^2)^{m-1}\tau_1}{\sqrt{2c_m+4m+2}}\s^{(n-1)}\oX(\tau_1)\right\}\nonumber\\
&&\hspace{20mm}+C\left(\frac{\tau_1\s^{(m+1,n-m)}P(\tau_1)}{\sqrt{2c_m+2}}
+\frac{\s^{(m+1,n-m)}V(\tau_1)}{\sqrt{2c_m}}+\frac{\tau_1^2\s^{(m+1,n-m)}\Ga(\tau_1)}{\sqrt{2c_m+4}}
\right)\nonumber\\
&&\hspace{20mm}+C\frac{\sqrt{\s^{(Y;m,n-m)}\cA(\tau_1)}}{\sqrt{2c_m}}
+C\tau_1^{1/2}\frac{\sqrt{\s^{(E;m,n-m)}\cA(\tau_1)}}{\sqrt{2c_m+1}}
\nonumber\\&&\label{12.813}
\end{eqnarray}

Going now back to \ref{12.755} - \ref{12.761}, in view of the definition \ref{12.795}, \ref{12.756} 
is bounded by:
\begin{equation}
C\tau_1^{c_m}\s^{(m,n-m)}\oLa(\tau_1)
\label{12.814}
\end{equation}
and, in view of the definition \ref{12.796}, \ref{12.759} is bounded by:
\begin{equation}
Cm\frac{\tau^{c_m}}{\sqrt{2c_m}}\s^{(m,n-m)}\oLab(\tau_1)
\label{12.815}
\end{equation}
Moreover, \ref{12.761} is bounded by:
\begin{equation}
\frac{Cm}{\sqrt{2c_m}}\tau_1^{c_m}\sup_{\tau\in(0,\tau_1]}\left\{\tau^{-c_m}\|\Omega^{n-m}T^m\check{r}
\|_{L^2({\cal K}^{\tau})}\right\}\leq \frac{Cm}{\sqrt{2c_m}}\tau_1^{c_m}\s^{(m,n-m)}R(\tau_1)
\label{12.816}
\end{equation}
For, by \ref{12.812} with $\tau\in(0,\tau_1]$ in the role of $\tau_1$ and the fact that 
$\s^{(m,n-m)}R(\tau)$ we have:
$$\|\Omega^{n-m}T^m\check{r}\|_{L^2({\cal K}^{\tau})}\leq \tau^{c_m}\s^{(m,n-m)}R(\tau_1) 
\ : \ \forall \tau\in(0,\tau_1]$$
By \ref{12.811} the factor $m/\sqrt{c_m}$ in \ref{12.816} is less than 1. The above together with 
\ref{12.812} imply through \ref{12.755} the following estimate for $\s^{(m,n-m)}\clab$ in 
$L^2({\cal K}^{\tau})$:
\begin{eqnarray}
&&\|\s^{(m,n-m)}\clab\|_{L^2({\cal K}^{\tau_1})}\leq \tau_1^{c_m}\left\{C\s^{(m,n-m)}\oLa(\tau_1)
+\frac{Cm}{\sqrt{2c_m}}\s^{(m,n-m)}\oLab(\tau_1)\right.\nonumber\\
&&\left.\hspace{45mm}+C\s^{(m,n-m)}R(\tau_1)\right\}
\label{12.817}
\end{eqnarray}
Replacing $\tau_1$ by any $\tau\in(0,\tau_1]$ the inequality holds a fortiori if we keep the 
argument of the quantity in parenthesis fixed to $\tau_1$, this quantity being non-decreasing 
in $\tau$. Multiplying then by $\tau^{-c_m}$ and taking the supremum over $\tau\in(0,\tau_1]$ 
yields, in view of the definition \ref{12.796}, 
\begin{eqnarray}
&&\s^{(m,n-m)}\oLab(\tau_1)\leq C\s^{(m,n-m)}\oLa(\tau_1)
+\frac{Cm}{\sqrt{2c_m}}\s^{(m,n-m)}\oLab(\tau_1)\nonumber\\
&&\hspace{40mm}+C\s^{(m,n-m)}R(\tau_1)
\label{12.818}
\end{eqnarray}
Substituting the expression \ref{12.813} for $\s^{(m,n-m)}R(\tau_1)$ the coefficient of 
$\s^{(m,n-m)}\oLab(\tau_1)$ in the right hand side becomes:
\begin{equation}
\frac{Cm}{\sqrt{2c_m}}+\frac{C\tau_1}{\sqrt{2c_m+2}}
\label{12.819}
\end{equation}
Choosing then $c_n$ large enough so that:
\begin{equation}
\frac{Cn}{\sqrt{2c_n}}+\frac{C\delta}{\sqrt{2c_n+2}}\leq \frac{1}{2}
\label{12.820}
\end{equation}
in view of the non-increasing with $m$ property of the exponents $c_m$ \ref{12.819} does not 
exceed 1/2 for all $m=1,...,n$. Hence the quantity $\s^{(m,n-m)}\oLab(\tau_1)$ can be eliminated 
from the right hand side of \ref{12.818} and this yields the following lemma. 

\vspace{2.5mm}

\noindent{\bf Lemma 12.11} \ \ On ${\cal K}$, the next to top order acoustical difference quantities 
$\|\s^{(m,n-m)}\clab\|_{L^2({\cal K}^{\tau})} \ : \ m=1,...,n$ satisfy to principal terms 
the following inequalities:
\begin{eqnarray*}
&&\s^{(m,n-m)}\oLab(\tau_1)\leq C\s^{(m,n-m)}\oLa(\tau_1)\\
&&\hspace{10mm}+C\tau_1\left\{\sum_{j=0}^{m-1}(\Pib\tau_1)^j \ 
\s^{(m-1-j,n-m+j+1)}\oLab(\tau_1)+(\Pib\tau_1)^{m-1} \ \s^{(n-1)}\oXb(\tau_1)\right\}\\
&&+C\tau_1^2\left\{\sum_{j=0}^{m-1}\frac{(\Pi\tau_1^2)^j}{\sqrt{2c_m+4j+4}}
\s^{(m-1-j,n-m+j+1)}\oLa(\tau_1)
+\frac{(\Pi\tau_1^2)^{m-1}\tau_1}{\sqrt{2c_m+4m+2}}\s^{(n-1)}\oX(\tau_1)\right\}\\
&&\hspace{20mm}+C\left(\frac{\tau_1\s^{(m+1,n-m)}P(\tau_1)}{\sqrt{2c_m+2}}
+\frac{\s^{(m+1,n-m)}V(\tau_1)}{\sqrt{2c_m}}+\frac{\tau_1^2\s^{(m+1,n-m)}\Ga(\tau_1)}{\sqrt{2c_m+4}}
\right)\\
&&\hspace{20mm}+C\frac{\sqrt{\s^{(Y;m,n-m)}\cA(\tau_1)}}{\sqrt{2c_m}}
+C\tau_1^{1/2}\frac{\sqrt{\s^{(E;m,n-m)}\cA(\tau_1)}}{\sqrt{2c_m+1}}
\end{eqnarray*}

\vspace{2.5mm}

Consider Lemma 12.10. The 1st term on the right hand side of the inequality satisfied by 
$\|\s^{(m,n-m)}\clab\|_{L^2(C_{u_1}^{\ub_1})}$ is bounded by:
\begin{equation}
k\s^{(m,n-m)}\oLab(\ub_1)\ub_1^{c_m}
\label{12.821}
\end{equation}
In view of the definition \ref{12.665} the 2nd term on the right hand side of the same inequality 
is bounded by:
\begin{equation}
C\frac{\ub_1^{a_m+2}u_1^{b_m+1}}{\sqrt{2a_m+4}}\s^{(m,n-m)}\La(\ub_1,u_1)
\label{12.822}
\end{equation}
Substituting the above we obtain:
\begin{eqnarray}
&&\|\s^{(m,n-m)}\clab\|_{L^2(C_{u_1}^{\ub_1})}\leq k\s^{(m,n-m)}\oLab(\ub_1)\ub_1^{c_m}\nonumber\\
&&\hspace{27mm}+C\frac{\ub_1^{a_m+2}u_1^{b_m+1}}{\sqrt{2a_m+4}}\s^{(m,n-m)}\La(\ub_1,u_1)\nonumber\\
&&\hspace{27mm}+Cm\frac{\ub_1^{a_m+1}u_1^{b_m+1}}{(b_m+1)}\s^{(m-1,n-m+1)}\Lab(\ub_1,u_1)\nonumber\\
&&\hspace{24mm}+C\sum_{j=0}^{m-1}(\Pi u_1^2)^j \ \frac{\ub_1^{a_m+2}u_1^{b_m+1}}{\sqrt{2a_m+4}}
\s^{(m-1-j,n-m+j+1)}\La(\ub_1,u_1)\nonumber\\
&&\hspace{27mm}+C(\Pi u_1^2)^m \ \frac{\ub_1^{a_m+\frac{3}{2}}u_1^{b_m+\frac{1}{2}}}{\sqrt{2a_m+3}}
\s^{(n-1)}X(\ub_1,u_1)\nonumber\\
&&\hspace{24mm}+Cu_1\sum_{j=0}^{m-1}\Pib^j \ \frac{\ub_1^{a_m+j}u_1^{b_m+1}}{(b_m+1)}
\s^{(m-1-j,n-m+j+1)}\Lab(\ub_1,u_1)\nonumber\\
&&\hspace{27mm}+Cu_1\Pib^m \ \frac{\ub_1^{a_m+m}u_1^{b_m}}{b_m}\s^{(n-1)}\Xb(\ub_1,u_1)\nonumber\\
&&\hspace{20mm}+\Cb\sqrt{\max\{\s^{(m,n-m)}\cB(\ub_1,u_1),\s^{(m,n-m)}\cBb(\ub_1,u_1)\}}
\cdot \ub_1^{a_m}u_1^{b_m}\cdot u_1\nonumber\\
&&\label{12.823}
\end{eqnarray}

We turn to Lemma 12.9. We first consider the inequality satisfied by 
$\|\s^{(m,n-m)}\cla\|_{L^2(\Cb_{\ub_1}^{u_1})}$. The 1st term on the right is by virtue of the estimate 
\ref{12.823} bounded by:
\begin{eqnarray}
&&k\s^{(m,n-m)}\oLab(\ub_1)\ub_1^{c_m+\frac{1}{2}}u_1^{\frac{1}{2}}\nonumber\\
&&+C\frac{\ub_1^{a_m+\frac{5}{2}}u_1^{b_m+\frac{3}{2}}}{\sqrt{2a_m+4}\sqrt{2b_m+3}}\s^{(m,n-m)}\La(\ub_1,u_1)\nonumber\\
&&+Cm\frac{\ub_1^{a_m+\frac{3}{2}}u_1^{b_m+\frac{3}{2}}}{(b_m+1)\sqrt{2b_m+3}}\s^{(m-1,n-m+1)}\Lab(\ub_1,u_1)\nonumber\\
&&+C\sum_{j=0}^{m-1}(\Pi u_1^2)^j \ \frac{\ub_1^{a_m+\frac{5}{2}}u_1^{b_m+\frac{3}{2}}}{\sqrt{2a_m+4}\sqrt{2b_m+4j+3}}
\s^{(m-1-j,n-m+j+1)}\La(\ub_1,u_1)\nonumber\\
&&+C(\Pi u_1^2)^m \ \frac{\ub_1^{a_m+2}u_1^{b_m+1}}{\sqrt{2a_m+3}\sqrt{2b_m+4m+2}}
\s^{(n-1)}X(\ub_1,u_1)\nonumber\\
&&+Cu_1\sum_{j=0}^{m-1}\Pib^j \ \frac{\ub_1^{a_m+j+\frac{1}{2}}u_1^{b_m+\frac{3}{2}}}
{(b_m+1)\sqrt{2b_m+5}}
\s^{(m-1-j,n-m+j+1)}\Lab(\ub_1,u_1)\nonumber\\
&&+Cu_1\Pib^m \ \frac{\ub_1^{a_m+m+\frac{1}{2}}u_1^{b_m}+\frac{1}{2}}{b_m\sqrt{2b_m+3}}\s^{(n-1)}\Xb(\ub_1,u_1)\nonumber\\
&&+\Cb\sqrt{\max\{\s^{(m,n-m)}\cB(\ub_1,u_1),\s^{(m,n-m)}\cBb(\ub_1,u_1)\}}
\cdot\frac{\ub_1^{a_m+\frac{1}{2}}u_1^{b_m+\frac{3}{2}}}{\sqrt{2b_m+3}}\nonumber\\
&&\label{12.824}
\end{eqnarray}
The 2nd term on the right hand side of the same inequality is bounded by:
\begin{equation}
C\frac{\ub_1^{a_m+\frac{5}{2}}u_1^{b_m+\frac{1}{2}}}{\left(a_m+\frac{5}{2}\right)}
\s^{(m,n-m)}\La(\ub_1,u_1)
\label{12.825}
\end{equation}
Replacing in the inequality which results by substituting the above $(\ub_1,u_1)$ by any 
$(\ub,u)\in R_{\ub_1,u_1}$, multiplying both sides by $\ub^{-a_m-\frac{1}{2}}u^{-b_m-\frac{1}{2}}$ 
and taking the supremum over $(\ub,u)\in R_{\ub_1,u_1}$, taking account of the fact that $\ub/u\leq 1$, 
we obtain: 
\begin{eqnarray}
&&\s^{(m,n-m)}\La(\ub_1,u_1)\leq C\s^{(m,n-m)}\oLab(\ub_1)\nonumber\\
&&\hspace{20mm}+C\ub_1^2\left(\frac{1}{\left(a_m+\frac{5}{2}\right)}+\frac{u_1}
{\sqrt{2a_m+4}\sqrt{2b_m+3}}\right)\s^{(m,n-m)}\La(\ub_1,u_1)\nonumber\\
&&\hspace{25mm}+C\frac{m\ub_1^2}{\sqrt{2a_m+4}}\s^{(m-1,n-m+1)}\La(\ub_1,u_1)\nonumber\\
&&\hspace{25mm}+C\frac{m\ub_1 u_1}{(b_m+1)\sqrt{2b_m+3}}\s^{(m-1,n-m+1)}\Lab(\ub_1,u_1)\nonumber\\
&&\hspace{20mm}+C\sum_{j=0}^{m-1}(\Pi u_1^2)^j \ \frac{\ub_1^2}{\sqrt{2a_m+4}}\s^{(m-1-j,n-m+j+1)}\Lambda(\ub_1,u_1)\nonumber\\
&&\hspace{35mm}+C(\Pi u_1^2)^m  \ \frac{\ub_1}{\sqrt{2a_m+3}}\s^{(n-1)}X(\ub_1,u_1)\nonumber\\
&&\hspace{20mm}+C\sum_{j=0}^{m-1}(\Pib\ub_1)^j \ \frac{u_1}{\sqrt{2b_m+3}}\s^{m-1-j,n-m+j+1)}\Lab
(\ub_1,u_1)\nonumber\\
&&\hspace{35mm}+C(\Pib\ub_1)^m \ \frac{1}{\sqrt{2b_m+1}}\s^{(n-1)}\Xb(\ub_1,u_1)\nonumber\\
&&\hspace{35mm}+C\frac{\sqrt{\s^{(Y;m,n-m)}\cBb(\ub_1,u_1)}}{(a_m+2)}\nonumber\\
&&\hspace{20mm}+Cu_1^{1/2}\sqrt{\max\{\s^{(m,n-m)}\cB(\ub_1,u_1),\s^{(m,n-m)}\cBb(\ub_1,u_1)\}}
\nonumber\\
&&\label{12.826}
\end{eqnarray}
Subjecting $\delta$ to an appropriate smallness condition or $a_m$ to an appropriate largeness 
condition so that the coefficient of 
$\s^{(m,n-m)}\La(\ub_1,u_1)$ on the right does not exceed 1/2 this implies:
\begin{eqnarray}
&&\s^{(m,n-m)}\La(\ub_1,u_1)\leq C\s^{(m,n-m)}\oLab(\ub_1)\nonumber\\
&&\hspace{25mm}+C\frac{m\ub_1^2}{\sqrt{2a_m+4}}\s^{(m-1,n-m+1)}\La(\ub_1,u_1)\nonumber\\
&&\hspace{25mm}+C\frac{m\ub_1 u_1}{(b_m+1)\sqrt{2b_m+3}}\s^{(m-1,n-m+1)}\Lab(\ub_1,u_1)\nonumber\\
&&\hspace{20mm}+C\sum_{j=0}^{m-1}(\Pi u_1^2)^j \ \frac{\ub_1^2}{\sqrt{2a_m+4}}\s^{(m-1-j,n-m+j+1)}\Lambda(\ub_1,u_1)\nonumber\\
&&\hspace{35mm}+C(\Pi u_1^2)^m  \ \frac{\ub_1}{\sqrt{2a_m+3}}\s^{(n-1)}X(\ub_1,u_1)\nonumber\\
&&\hspace{20mm}+C\sum_{j=0}^{m-1}(\Pib\ub_1)^j \ \frac{u_1}{\sqrt{2b_m+3}}\s^{(m-1-j,n-m+j+1)}\Lab
(\ub_1,u_1)\nonumber\\
&&\hspace{35mm}+C(\Pib\ub_1)^m \ \frac{1}{\sqrt{2b_m+1}}\s^{(n-1)}\Xb(\ub_1,u_1)\nonumber\\
&&\hspace{35mm}+C\frac{\sqrt{\s^{(Y;m,n-m)}\cBb(\ub_1,u_1)}}{(a_m+2)}\nonumber\\
&&\hspace{20mm}+Cu_1^{1/2}\sqrt{\max\{\s^{(m,n-m)}\cB(\ub_1,u_1),\s^{(m,n-m)}\cBb(\ub_1,u_1)\}}
\nonumber\\
&&\label{12.827}
\end{eqnarray}
(for new constants $C$). We now consider the inequality of Lemma 12.9 satisfied by 
$\|\s^{(m,n-m)}\cla\|_{L^2({\cal K}^{\tau})}$. The most important term on the right hand side is 
the first term:
\begin{equation}
C\left(\int_0^{\tau} u\|\s^{(m,n-m)}\clab\|^2_{L^2(C_u^u)}du\right)^{1/2}
\label{12.828}
\end{equation}
We substitute the estimate \ref{12.823} with $\ub_1=u_1=u$. The leading contribution to \ref{12.828} 
comes from the boundary term:
$$k\s^{(m,n-m)}\oLab(u) u^{c_m}$$
This contribution is:
\begin{equation}
C\left(\int_0^{\tau}u^{2c_m+1}du\right)^{1/2}\s^{(m,n-m)}\oLab(\tau)
=C\frac{\tau^{c_m+1}}{\sqrt{2c_m+2}}\s^{(m,n-m)}\oLab(\tau)
\label{12.829}
\end{equation}
The second term on the right hand side of the inequality of Lemma 12.9 satisfied by 
$\|\s^{(m,n-m)}\cla\|_{L^2({\cal K}^{\tau})}$ coincides with the corresponding term in the 
inequality satisfied by $\|\s^{(m,n-m)}\cla\|_{L^2(\Cb_{\ub_1}^{u_1})}$ if we set $\ub_1=u_1=\tau$. 
So this term is bounded by (see \ref{12.825}):
\begin{equation}
C\frac{\tau^{c_m+3}}{\left(a_m+\frac{5}{2}\right)}\s^{(m,n-m)}\La(\tau,\tau)
\label{12.830}
\end{equation}
Here we substitute the estimate \ref{12.827} with $\ub_1=u_1=\tau$. We multiply the inequality resulting 
from the above substitutions by $\tau^{-c_m-1}$ and take the supremum over $\tau\in(0,\tau_1]$ to 
obtain:
\begin{eqnarray}
&&\s^{(m,n-m)}\oLa(\tau_1)\leq C\left(\frac{1}{\sqrt{2c_m+2}}+\frac{\tau_1^2}{\sqrt{2a_m+4}}\right)
\s^{(m,n-m)}\oLab(\tau_1)\nonumber\\
&&\hspace{15mm}+Cm\tau_1^2\left(\frac{\s^{(m-1,n-m+1)}\Lab(\tau_1,\tau_1)}{\sqrt{2c_m+6} \ (b_m+1)}
+\frac{\s^{(m-1,n-m+1)}\La(\tau_1,\tau_1)}{\sqrt{2a_m+4}}\right)\nonumber\\
&&\hspace{15mm}+C\sum_{j=0}^{m-1}(\Pi\tau_1^2)^j \ \frac{\tau_1^2}{\sqrt{2a_m+4}}
\s^{(m-1-j,n-m+j+1)}\La(\tau_1,\tau_1)\nonumber\\
&&\hspace{35mm}+C(\Pi\tau_1^2)^m \ \frac{\tau_1}{\sqrt{2a_m+3}}\s^{(n-1)}X(\tau_1,\tau_1)\nonumber\\
&&\hspace{15mm}+C\sum_{j=0}^{m-1}(\Pib\tau_1)^j \ \frac{\tau_1}{\sqrt{2c_m+4}}
\s^{(m-1-j,n-m+j+1)}\Lab(\tau_1,\tau_1)\nonumber\\
&&\hspace{35mm}+C(\Pib\tau_1)^m \ \frac{1}{\sqrt{2c_m+4}}\s^{(n-1)}\Xb(\tau_1,\tau_1)\nonumber\\
&&\hspace{25mm}+C\frac{\sqrt{\s^{(Y;m,n-m)}\cBb(\tau_1,\tau_1)}}{(a_m+2)}\nonumber\\
&&\hspace{15mm}+C\tau_1^{1/2}\sqrt{\max\{\s^{(m,n-m)}\cB(\tau_1,\tau_1),
\s^{(m,n-m)}\cBb(\tau_1,\tau_1)\}}\nonumber\\
&&\label{12.831}
\end{eqnarray}

Substituting this bound for $\s^{(m,n-m)}\oLa(\tau_1)$ in the inequality for $\s^{(m,n-m)}\oLab(\tau_1)$ 
of Lemma 12.11, the coefficient of $\s^{(m,n-m)}\oLab(\tau_1)$ is:
\begin{equation}
C\left(\frac{1}{\sqrt{2c_m+2}}+\frac{\delta^2}{\sqrt{2a_m+4}}\right)
\label{12.832}
\end{equation}
Choosing then $c_m$ large enough so that:
\begin{equation}
\frac{C}{\sqrt{2c_m+2}}\leq\frac{1}{4}
\label{12.833}
\end{equation}
and $a_m$ large enough so that:
\begin{equation}
\frac{C\delta^2}{\sqrt{2a_m+4}}\leq\frac{1}{4}
\label{12.834}
\end{equation}
\ref{12.832} is not greater than 1/2. The inequality then implies the estimate for 
$\s^{(m,n-m)}\oLab(\tau_1)$ of the following proposition. Substituting this estimate in \ref{12.831} 
yields the estimate for $\s^{(m,n-m)}\La(\tau_1)$ of the proposition. 

\vspace{2.5mm}

\noindent{\bf Proposition 12.4} \ \ Choosing the exponents $a_m$ and $c_m$ suitably large, 
the next to top order boundary acoustical 
difference quantities $\s^{(m,n-m)}\oLab$, $\s^{(m,n-m)}\oLa$ \ : \ $m=1,...,n$ satisfy 
to principal terms the estimates:
\begin{eqnarray*}
&&\s^{(m,n-m)}\oLab(\tau_1)\leq C\tau_1\left\{\sum_{j=0}^{m-1}(\Pib\tau_1)^j \ 
\s^{(m-1-j,n-m+j+1)}\oLab(\tau_1)+(\Pib\tau_1)^{m-1} \ \s^{(n-1)}\oXb(\tau_1)\right\}\\
&&\hspace{10mm}+C\tau_1^2\left\{\sum_{j=0}^{m-1}\frac{(\Pi\tau_1^2)^j}{\sqrt{2c_m+4j+4}}
\s^{(m-1-j,n-m+j+1)}\oLa(\tau_1)
+\frac{(\Pi\tau_1^2)^{m-1}\tau_1}{\sqrt{2c_m+4m+2}}\s^{(n-1)}\oX(\tau_1)\right\}\\
&&\hspace{25mm}+Cm\tau_1^2\left(\frac{\s^{(m-1,n-m+1)}\Lab(\tau_1,\tau_1)}{\sqrt{2c_m+6} \ (b_m+1)}
+\frac{\s^{(m-1,n-m+1)}\La(\tau_1,\tau_1)}{\sqrt{2a_m+4}}\right)\\
&&\hspace{25mm}+\frac{C\tau_1^2}{\sqrt{2a_m+4}}\sum_{j=0}^{m-1}(\Pi\tau_1^2)^j \ 
\s^{(m-1-j,n-m+j+1)}\La(\tau_1,\tau_1)\\
&&\hspace{40mm}+\frac{C\tau_1}{\sqrt{2a_m+3}}(\Pi\tau_1^2)^m \ \s^{(n-1)}X(\tau_1,\tau_1)\\
&&\hspace{25mm}+\frac{C\tau_1}{\sqrt{2c_m+4}}\sum_{j=0}^{m-1}(\Pib\tau_1)^j \ 
\s^{(m-1-j,n-m+j+1)}\Lab(\tau_1,\tau_1)\\
&&\hspace{40mm}+\frac{C}{\sqrt{2c_m+4}}(\Pib\tau_1)^m \ \s^{(n-1)}\Xb(\tau_1,\tau_1)\\
&&\hspace{25mm}+C\left(\frac{\tau_1\s^{(m+1,n-m)}P(\tau_1)}{\sqrt{2c_m+2}}
+\frac{\s^{(m+1,n-m)}V(\tau_1)}{\sqrt{2c_m}}+\frac{\tau_1^2\s^{(m+1,n-m)}\Ga(\tau_1)}{\sqrt{2c_m+4}}
\right)\\
&&\hspace{25mm}+C\frac{\sqrt{\s^{(Y;m,n-m)}\cA(\tau_1)}}{\sqrt{2c_m}}
+C\tau_1^{1/2}\frac{\sqrt{\s^{(E;m,n-m)}\cA(\tau_1)}}{\sqrt{2c_m+1}}\\
&&\hspace{40mm}+C\frac{\sqrt{\s^{(Y;m,n-m)}\cBb(\tau_1,\tau_1)}}{(a_m+2)}\\
&&\hspace{25mm}+C\tau_1^{1/2}\sqrt{\max\{\s^{(m,n-m)}\cB(\tau_1,\tau_1),
\s^{(m,n-m)}\cBb(\tau_1,\tau_1)\}}
\end{eqnarray*}
\begin{eqnarray*}
&&\s^{(m,n-m)}\oLa(\tau_1)\leq \frac{C\tau_1}{\sqrt{2a_m+4}}\left\{\sum_{j=0}^{m-1}(\Pib\tau_1)^j \ 
\s^{(m-1-j,n-m+j+1)}\oLab(\tau_1)+(\Pib\tau_1)^{m-1} \ \s^{(n-1)}\oXb(\tau_1)\right\}\\
&&\hspace{10mm}+\frac{C\tau_1^2}{\sqrt{2a_m+4}}\left\{\sum_{j=0}^{m-1}
\frac{(\Pi\tau_1^2)^j}{\sqrt{2c_m+4j+4}}\s^{(m-1-j,n-m+j+1)}\oLa(\tau_1)
+\frac{(\Pi\tau_1^2)^{m-1}\tau_1}{\sqrt{2c_m+4m+2}}\s^{(n-1)}\oX(\tau_1)\right\}\\
&&\hspace{25mm}+Cm\tau_1^2\left(\frac{\s^{(m-1,n-m+1)}\Lab(\tau_1,\tau_1)}{\sqrt{2c_m+6} \ (b_m+1)}
+\frac{\s^{(m-1,n-m+1)}\La(\tau_1,\tau_1)}{\sqrt{2a_m+4}}\right)\\
&&\hspace{25mm}+\frac{C\tau_1^2}{\sqrt{2a_m+4}}\sum_{j=0}^{m-1}(\Pi\tau_1^2)^j \ 
\s^{(m-1-j,n-m+j+1)}\La(\tau_1,\tau_1)\\
&&\hspace{40mm}+\frac{C\tau_1}{\sqrt{2a_m+3}}(\Pi\tau_1^2)^m \ \s^{(n-1)}X(\tau_1,\tau_1)\\
&&\hspace{25mm}+\frac{C\tau_1}{\sqrt{2c_m+4}}\sum_{j=0}^{m-1}(\Pib\tau_1)^j \ 
\s^{(m-1-j,n-m+j+1)}\Lab(\tau_1,\tau_1)\\
&&\hspace{40mm}+\frac{C}{\sqrt{2c_m+4}}(\Pib\tau_1)^m \ \s^{(n-1)}\Xb(\tau_1,\tau_1)\\
&&\hspace{25mm}+\frac{C}{\sqrt{2a_m+4}}\left(\frac{\tau_1\s^{(m+1,n-m)}P(\tau_1)}{\sqrt{2c_m+2}}
+\frac{\s^{(m+1,n-m)}V(\tau_1)}{\sqrt{2c_m}}+\frac{\tau_1^2\s^{(m+1,n-m)}\Ga(\tau_1)}{\sqrt{2c_m+4}}
\right)\\
&&\hspace{25mm}+\frac{C}{\sqrt{2a_m+4}}\left(\frac{\sqrt{\s^{(Y;m,n-m)}\cA(\tau_1)}}{\sqrt{2c_m}}
+\tau_1^{1/2}\frac{\sqrt{\s^{(E;m,n-m)}\cA(\tau_1)}}{\sqrt{2c_m+1}}\right)\\
&&\hspace{40mm}+C\frac{\sqrt{\s^{(Y;m,n-m)}\cBb(\tau_1,\tau_1)}}{(a_m+2)}\\
&&\hspace{25mm}+C\tau_1^{1/2}\sqrt{\max\{\s^{(m,n-m)}\cB(\tau_1,\tau_1),
\s^{(m,n-m)}\cBb(\tau_1,\tau_1)\}}
\end{eqnarray*}

\vspace{2.5mm}

Substituting the estimates of the above proposition first in \ref{12.827} and then in \ref{12.823} 
we deduce the following proposition. 

\vspace{2.5mm}

\noindent{\bf Proposition 12.5} \ \ Choosing the exponents $a_m$ and $c_m$ suitably large, the next to top order interior acoustical 
difference quantities $\s^{(m,n-m)}\La$, $\s^{(m,n-m)}\Lab$ \ : \ $m=1,...,n$ satisfy 
to principal terms the estimates:
\begin{eqnarray*}
&&\s^{(m,n-m)}\La(\ub_1,u_1)\leq C\ub_1\left\{\sum_{j=0}^{m-1}(\Pib\ub_1)^j \ 
\s^{(m-1-j,n-m+j+1)}\oLab(\ub_1)+(\Pib\ub_1)^{m-1} \ \s^{(n-1)}\oXb(\ub_1)\right\}\\
&&\hspace{10mm}+C\ub_1^2\left\{\sum_{j=0}^{m-1}\frac{(\Pi\ub_1^2)^j}{\sqrt{2c_m+4j+4}}
\s^{(m-1-j,n-m+j+1)}\oLa(\ub_1)
+\frac{(\Pi\ub_1^2)^{m-1}\ub_1}{\sqrt{2c_m+4m+2}}\s^{(n-1)}\oX(\ub_1)\right\}\\
&&\hspace{15mm}+\frac{Cm\ub_1 u_1}{(b_m+1)\sqrt{2b_m+3}}\s^{(m-1,n-m+1)}\Lab(\ub_1,u_1)
+\frac{Cm\ub_1^2}{\sqrt{2a_m+4}}\s^{(m-1,n-m)}\La(\ub_1,u_1)\\
&&\hspace{30mm}+\frac{C\ub_1^2}{\sqrt{2a_m+4}}\sum_{j=0}^{m-1}(\Pi u_1^2)^j \ 
\s^{(m-1-j,n-m+j+1)}\La(\ub_1,u_1)\\
&&\hspace{50mm}+\frac{C\ub_1}{\sqrt{2a_m+3}}(\Pi u_1^2)^m \ \s^{(n-1)}X(\ub_1,u_1)\\
&&\hspace{30mm}+\frac{Cu_1}{\sqrt{2b_m+3}}\sum_{j=0}^{m-1}(\Pib\ub_1)^j \ 
\s^{(m-1-j,n-m+j+1)}\Lab(\ub_1,u_1)\\
&&\hspace{50mm}+\frac{C}{\sqrt{2b_m+1}}(\Pib\ub_1)^m \ \s^{(n-1)}\Xb(\ub_1,u_1)\\
&&\hspace{25mm}+C\left(\frac{\ub_1\s^{(m+1,n-m)}P(\ub_1)}{\sqrt{2c_m+2}}
+\frac{\s^{(m+1,n-m)}V(\ub_1)}{\sqrt{2c_m}}+\frac{\ub_1^2\s^{(m+1,n-m)}\Ga(\ub_1)}{\sqrt{2c_m+4}}
\right)\\
&&\hspace{25mm}+C\frac{\sqrt{\s^{(Y;m,n-m)}\cA(\ub_1)}}{\sqrt{2c_m}}
+C\ub_1^{1/2}\frac{\sqrt{\s^{(E;m,n-m)}\cA(\ub_1)}}{\sqrt{2c_m+1}}\\
&&\hspace{50mm}+C\frac{\sqrt{\s^{(Y;m,n-m)}\cBb(\ub_1,u_1)}}{(a_m+2)}\\
&&\hspace{30mm}+C u_1^{1/2}\sqrt{\max\{\s^{(m,n-m)}\cB(\ub_1,u_1),\s^{(m,n-m)}\cBb(\ub_1,u_1)\}}
\end{eqnarray*}
\begin{eqnarray*}
&&\s^{(m,n-m)}\Lab(\ub_1,u_1)\leq C\ub_1\left\{\sum_{j=0}^{m-1}(\Pib\ub_1)^j \ 
\s^{(m-1-j,n-m+j+1)}\oLab(\ub_1)+(\Pib\ub_1)^{m-1} \ \s^{(n-1)}\oXb(\ub_1)\right\}\\
&&\hspace{10mm}+C\ub_1^2\left\{\sum_{j=0}^{m-1}\frac{(\Pi\ub_1^2)^j}{\sqrt{2c_m+4j+4}}
\s^{(m-1-j,n-m+j+1)}\oLa(\ub_1)
+\frac{(\Pi\ub_1^2)^{m-1}\ub_1}{\sqrt{2c_m+4m+2}}\s^{(n-1)}\oX(\ub_1)\right\}\\
&&\hspace{20mm}+\frac{Cm\ub_1 u_1}{(b_m+1)}\s^{(m-1,n-m+1)}\Lab(\ub_1,u_1)
+\frac{Cm\ub_1^4 u_1}{(2a_m+4)}\s^{(m-1,n-m+1)}\La(\ub_1,u_1)\\
&&\hspace{30mm}+\frac{C\ub_1^2 u_1}{\sqrt{2a_m+4}}\sum_{j=0}^{m-1}(\Pi u_1^2)^j \ 
\s^{(m-1-j,n-m+j+1)}\La(\ub_1,u_1)\\
&&\hspace{50mm}+\frac{\ub_1^{3/2}u_1^{1/2}}{\sqrt{2a_m+3}}(\Pi u_1^2)^m \ \s^{(n-1)}X(\ub_1,u_1)\\
&&\hspace{30mm}+\frac{Cu_1^2}{(b_m+1)}\sum_{j=0}^{m-1}(\Pib\ub_1)^j \ \s^{(m-1,n-m+1)}\Lab(\ub_1,u_1)\\
&&\hspace{50mm}+\frac{Cu_1}{b_m}(\Pib\ub_1)^m \ \s^{(n-1)}\Xb(\ub_1,u_1)\\
&&\hspace{25mm}+\Cb u_1\sqrt{\max\{\s^{(m,n-m)}\cB(\ub_1,u_1),\s^{(m,n-m)}\cBb(\ub_1,u_1)\}}
\end{eqnarray*}

\vspace{2.5mm}

\section{Estimates for $(T^{m+1}\Omega^{n-m}\chf, T^{m+1}\Omega^{n-m}\cv, T^{m+1}\Omega^{n-m}\cga)$, 
for $m=1,...,n$}

We now apply $\Omega^{n-m}T^m$ to \ref{12.470} to deduce that $\Omega^{n-m}T^{m+1}\cf$ is to 
principal terms given by:
\begin{equation}
\frac{(1+r)}{c}\Omega^{n-m}T^m\clab+\frac{\lambdab}{c}\Omega^{n-m}T^m\check{r}
-\frac{\lambdab (1+r)}{c}\Omega^{n-m}T^{m-1}(c^{-1}Tc-c_N^{-1}Tc_N)
\label{12.835}
\end{equation}
In fact, the principal part of the first of \ref{12.835} is:
\begin{equation}
\frac{(1+r)}{c}\sh^{(n-m)/2}\s^{(m,n-m)}\clab
\label{12.836}
\end{equation}
In view of the definition \ref{12.796} this is bounded in $L^2({\cal K}^{\tau})$ by:
\begin{equation}
C\tau^{c_m}\s^{(m,n-m)}\oLab(\tau)
\label{12.837}
\end{equation}
By \ref{12.812} the second of \ref{12.835} is bounded in $L^2({\cal K}^{\tau})$ by:
\begin{equation}
C\tau^{c_m+1}\s^{(m,n-m)}R(\tau)
\label{12.838}
\end{equation}
$\s^{(m,n-m)}R(\tau)$ been given by \ref{12.813}. In regard to the third of \ref{12.835} we 
recall that according to \ref{12.69}
$$c^{-1}[Tc]_{P.A.}=-2\pi(\mu\mub)$$
hence by \ref{12.107}, \ref{12.108}:
$$[T^{m-1}(c^{-1}Tc)]_{P.A.}=-2\pi[T^{m-1}\mu+T^{m-1}\mub]_{P.A.}=-2\pi(\mu_{m-1}+\mub_{m-1})$$
In view of \ref{12.110}, \ref{12.111} it then follows that the principal acoustical part of 
the third of \ref{12.835} is:
\begin{equation}
\frac{2\pi\lambdab(1+r)}{c}\sh^{(n-m)/2}(\cmu_{m-1,n-m}+\cmub_{m-1,n-m})
\label{12.839}
\end{equation}
By \ref{12.802}, \ref{12.804} this is bounded in $L^2({\cal K}^{\tau})$ by:
\begin{eqnarray}
&&C\tau^{c_m+1}\left\{\tau\sum_{j=0}^{m-1}(\Pi\tau^2)^j \ \s^{(m-1-j,n-m+j+1)}\oLa(\tau)
+(\Pi\tau^2)^m \ \s^{(n-1)}\oX(\tau)\right.\nonumber\\
&&\hspace{15mm}+\left.\sum_{j=0}^{m-1}(\Pib\tau)^j \ \s^{(m-1-j,n-m+j+1)}\oLab(\tau)
+(\Pib\tau)^{m-1} \ \s^{(n-1)}\oXb(\tau)\right\}\nonumber\\
\label{12.840}
\end{eqnarray}
Substituting the estimates  of Proposition 12.4 we conclude from the 
above that:
\begin{eqnarray}
&&\sup_{\tau\in(0,\tau_1]}\left\{\tau^{-c_m}\|\Omega^{n-m}T^{m+1}\cf\|_{L^2({\cal K}^{\tau})}\right\}
\leq C\frac{\s^{(m+1,n-m)}V(\tau_1)}{\sqrt{2c_m}}\nonumber\\
&&\hspace{40mm}+C\frac{\tau_1\s^{(m+1,n-m)}P(\tau_1)}{\sqrt{2c_m+2}}
+C\frac{\tau_1^2\s^{(m+1,n-m)}\Ga(\tau_1)}{\sqrt{2c_m+4}}\nonumber\\
&&\hspace{35mm}+C\frac{\sqrt{\s^{(Y;m,n-m)}\cA(\tau_1)}}{\sqrt{2c_m}}
+C\tau_1^{1/2}\frac{\sqrt{\s^{(E;m,n-m)}\cA(\tau_1)}}{\sqrt{2c_m+1}}\nonumber\\
&&\hspace{50mm}+C\frac{\sqrt{\s^{(Y;m,n-m)}\cBb(\tau_1,\tau_1)}}{(a_m+2)}\nonumber\\
&&\hspace{25mm}+C\tau_1^{1/2}\sqrt{\max\{\s^{(m,n-m)}\cB(\tau_1,\tau_1),
\s^{(m,n-m)}\cBb(\tau_1,\tau_1)\}}\nonumber\\
&&\hspace{25mm}+C\tau_1^2\s^{(m,n-m)}L(\tau_1)+C\tau_1\s^{(m,n-m)}\underline{L}(\tau_1)
+C\tau_1\s^{(m,n-m)}M(\tau_1)\nonumber\\
&&\label{12.841}
\end{eqnarray}
where $\s^{(m,n-m)}L(\tau_1)$ refers to the quantities 
$\s^{(n-1)}\oX(\tau_1)$ and $\s^{(i,n-i)}\oLa(\tau_1)$ for $i=0,...,m-1$, 
$\s^{(m,n-m)}\underline{L}(\tau_1)$ refers to the quantities $\s^{(n-1)}\oXb(\tau_1)$ and 
$\s^{(i,n-i)}\oLab(\tau_1)$ for $i=0,...,m-1$:
\begin{eqnarray}
&&\s^{(m,n-m)}L(\tau_1)=\sum_{j=0}^{m-1}(\Pi\tau_1^2)^j \s^{(m-1-j,n-m+j+1)}\oLa(\tau_1)\nonumber\\
&&\hspace{42.5mm}+\tau_1(\Pi \tau_1^2)^{m-1} \ \s^{(n-1)}\oX(\tau_1),\nonumber\\
&&\s^{(m,n-m)}\underline{L}(\tau_1)=\sum_{j=0}^{m-1}(\Pib\tau_1)^j \ \s^{(m-1-j,n-m+j+1)}\oLab(\tau_1)\nonumber\\
&&\hspace{42.5mm}+(\Pib\tau_1)^{m-1}\s^{(n-1)}\oXb(\tau_1) \label{12.842}
\end{eqnarray}
and $\s^{(m,n-m)}M(\tau_1)$ refers to the quantities 
$\s^{(n-1)}X(\tau_1,\tau_1)$, $\s^{(n-1)}\Xb(\tau_1,\tau_1)$ and $\s^{(i,n-i)}\La(\tau_1,\tau_1)$, 
$\s^{(i,n-1)}\Lab(\tau_1,\tau_1)$ for $i=0,...,m-1$:
\begin{eqnarray}
&&\s^{(m,n-m)}M(\tau_1)=m\tau_1\left(\frac{\s^{(m-1,n-m+1)}\La(\tau_1,\tau_1)}{\sqrt{2a_m+4}}
+\frac{\s^{(m-1,n-m+1)}\Lab(\tau_1,\tau_1)}{(b_m+1)\sqrt{2c_m+6}}\right)\nonumber\\
&&\hspace{25mm}+\frac{1}{\sqrt{2a_m+4}}\left\{\tau_1\sum_{j=0}^{m-1}(\Pi\tau_1^2)^j \ 
\s^{(m-1-j,n-m+j+1)}\La(\tau_1,\tau_1)\right.\nonumber\\
&&\hspace{60mm}\left.+\tau_1^2(\Pi\tau_1^2)^{m-1} \ \s^{(n-1)}X(\tau_1,\tau_1)\right\}\nonumber\\
&&\hspace{25mm}+\frac{1}{\sqrt{2c_m+4}}\left\{\sum_{j=0}^{m-1}(\Pib\tau_1)^j \ \s^{(m-1-j,n-m+j+1)}
\Lab(\tau_1,\tau_1)\right.\nonumber\\
&&\hspace{62.5mm}\left.+(\Pib\tau_1)^{m-1}\s^{(n-1)}\Xb(\tau_1,\tau_1)\right\}\nonumber\\
&&\label{12.843}
\end{eqnarray}

Next, recalling \ref{12.484} and noting that, for every non-negative integer $i$ we have
$$\frac{\partial^i (\tau^{-2})}{\partial\tau^i}=(-1)^i (i+1)! \tau^{-i-2}$$
we obtain:
\begin{equation}
\Omega^{n-m}T^{m+1}\chf=\sum_{i=0}^{m+1}\left(\begin{array}{c}m+1\\i\end{array}\right)(-1)^i(i+1)!
\tau^{-i-2}\Omega^{n-m}T^{m+1-i}\cf
\label{12.844}
\end{equation}
It follows that:
\begin{equation}
\|\Omega^{n-m}T^{m+1}\chf\|_{L^2({\cal K}^{\tau_1})}\leq \sum_{i=0}^{m+1}
\left(\begin{array}{c}m+1\\i\end{array}\right)(i+1)!
\|\tau^{-i-2}\Omega^{n-m}T^{m+1-i}\cf\|_{L^2({\cal K}^{\tau_1})}
\label{12.845}
\end{equation}
Consider first the $i=0$ term in the sum in \ref{12.845}. We have:
\begin{eqnarray}
&&\|\tau^{-2}\Omega^{n-m}T^{m+1}\cf\|^2_{L^2({\cal K}^{\tau_1})}=
\int_0^{\tau_1}\tau^{-4}\|\Omega^{n-m}T^{m+1}\cf\|^2_{L^2(S_{\tau,\tau})}d\tau\nonumber\\
&&\hspace{30mm}=\int_0^{\tau_1}\tau^{-4}\frac{\partial}{\partial\tau}
\left(\|\Omega^{n-m}T^{m+1}\cf\|^2_{L^2({\cal K}^{\tau})}\right)d\tau\nonumber\\
&&\hspace{15mm}=\tau_1^{-4}\|\Omega^{n-m}T^{m+1}\cf\|^2_{L^2({\cal K}^{\tau_1})}
+4\int_0^{\tau_1}\tau^{-5}\|\Omega^{n-m}T^{m+1}\cf\|^2_{L^2({\cal K}^{\tau})}d\tau\nonumber\\
&&\hspace{15mm}\leq\left(\sup_{\tau\in(0,\tau_1]}\left\{\tau^{-c_m}\|\Omega^{n-m}T^{m+1}\cf\|_{L^2({\cal K}^{\tau})}
\right\}\right)^2\cdot\left(1+\frac{2}{c_m-2}\right)\tau_1^{2c_m-4}\nonumber\\
&&\label{12.846}
\end{eqnarray}
Consider next the $i=1$ term in the sum in \ref{12.845}. Applying \ref{12.366} with 
$\Omega^{n-m}T^m\cf$ in the role of $f$ gives, for all $\tau\in(0,\tau_1]$,
\begin{eqnarray}
&&\|\Omega^{n-m}T^m\cf\|_{L^2(S_{\tau,\tau})}\leq k\tau^{1/2}
\|\Omega^{n-m}T^{m+1}\cf\|_{L^2({\cal K}^{\tau})}\nonumber\\
&&\hspace{25mm}\leq k\tau^{c_m+\frac{1}{2}}\sup_{\tau\in(0,\tau_1]}\left\{\tau^{-c_m}
\|\Omega^{n-m}T^{m+1}\cf\|_{L^2({\cal K}^{\tau})}\right\}\nonumber\\
&&\label{12.847}
\end{eqnarray}
hence:
\begin{eqnarray}
&&\|\tau^{-3}\Omega^{n-m}T^m\cf\|_{L^2({\cal K}^{\tau_1})}
=\left(\int_0^{\tau_1}\tau^{-6}\|\Omega^{n-m}T^m\cf\|^2_{L^2(S_{\tau,\tau})}d\tau\right)^{1/2}
\nonumber\\
&&\hspace{20mm}\leq k\sup_{\tau\in(0,\tau_1]}\left\{\tau^{-c_m}
\|\Omega^{n-m}T^{m+1}\cf\|_{L^2({\cal K}^{\tau})}\right\}\left(\int_0^{\tau_1}\tau^{2c_m-5}d\tau
\right)^{1/2}\nonumber\\
&&\hspace{20mm}=\frac{k\tau_1^{c_m-2}}{\sqrt{2(c_m-2)}}\sup_{\tau\in(0,\tau_1]}\left\{\tau^{-c_m}
\|\Omega^{n-m}T^{m+1}\cf\|_{L^2({\cal K}^{\tau})}\right\}\nonumber\\
&&\label{12.848}
\end{eqnarray}
From \ref{12.847} we proceed by induction on $i$ to show that for $i=1,...,m+1$ there are positive 
constants $k_i$ such that for all $\tau\in(0,\tau_1]$ we have:
\begin{equation}
\|\Omega^{n-m}T^{m+1-i}\cf\|_{L^2(S_{\tau,\tau})}\leq k_i\tau^{c_m+i-\frac{1}{2}}
\sup_{\tau\in(0,\tau_1]}\left\{\tau^{-c_m}\|\Omega^{n-m}T^{m+1}\cf\|_{L^2({\cal K}^{\tau})}\right\}
\label{12.849}
\end{equation}
According to \ref{12.847} this is true for $i=1$ with 
\begin{equation}
k_1=k
\label{12.850}
\end{equation}
Let it be true for $i=j-1$, $j\geq 2$. Applying then the first inequality in \ref{12.366} with $\Omega^{n-m}T^{m+1-j}\cf$ 
in the role of $f$ gives, for all $\tau\in(0,\tau_1]$:
\begin{eqnarray}
&&\|\Omega^{n-m}T^{m+1-j}\cf\|_{L^2(S_{\tau,\tau})}\leq k\tau^{1/2}
\|\Omega^{n-m}T^{m+2-j}\cf\|_{L^2({\cal K}^{\tau})}\nonumber\\
&&\hspace{30mm}=k\int_0^{\tau}
\|\Omega^{n-m}T^{m+2-j}\cf\|_{L^2(S_{\tau^\prime,\tau^\prime})}d\tau^\prime\nonumber\\
&&\hspace{20mm}\leq kk_{j-1}\sup_{\tau\in(0,\tau_1]}\left\{\tau^{-c_m}
\|\Omega^{n-m}T^{m+1}\cf\|_{L^2({\cal K}^{\tau})}\right\}\int_0^{\tau}\tau^{\prime c_m+j-\frac{3}{2}}
d\tau^\prime\nonumber\\
&&\hspace{20mm}=\frac{k k_{j-1}}{c_m+j-\frac{1}{2}}\tau^{c_m+j-\frac{1}{2}}
\sup_{\tau\in(0,\tau_1]}\left\{\tau^{-c_m}\|\Omega^{n-m}T^{m+1}\cf\|_{L^2({\cal K}^{\tau})}\right\}
\nonumber\\
&&\label{12.851}
\end{eqnarray}
therefore \ref{12.849} is true for $i=j$ if we set:
\begin{equation}
k_j=\frac{kk_{j-1}}{c_m+j-\frac{1}{2}}
\label{12.852}
\end{equation}
This is a recursion with starting point \ref{12.850}. It follows that for each $i=1,...,m+1$:
\begin{equation}
k_i= \frac{k^i}{\left(c_m+\frac{3}{2}\right)\cdot\cdot\cdot\left(c_m+i-\frac{1}{2}\right)}
\label{12.853}
\end{equation}
By virtue of \ref{12.849} we can estimate:
\begin{eqnarray}
&&\|\tau^{-i-2}\Omega^{n-m}T^{m+1-i}\cf\|_{L^2({\cal K}^{\tau_1})}=
\left(\int_0^{\tau_1}\tau^{-2i-4}\|\Omega^{n-m}T^{m+1-i}\cf\|^2_{L^2(S_{\tau,\tau})}d\tau\right)^{1/2}
\nonumber\\
&&\hspace{30mm}\leq k_i\sup_{\tau\in(0,\tau_1]}\left\{\tau^{-c_m}\|\Omega^{n-m}T^{m+1}\cf\|_{L^2({\cal K}^{\tau})}\right\}\left(\int_0^{\tau_1}\tau^{2c_m-5}d\tau\right)^{1/2}\nonumber\\
&&\hspace{30mm}=\frac{k_i\tau_1^{c_m-2}}{\sqrt{2(c_m-2)}}\sup_{\tau\in(0,\tau_1]}\left\{\tau^{-c_m}\|\Omega^{n-m}T^{m+1}\cf\|_{L^2({\cal K}^{\tau})}\right\}\nonumber\\
&&\label{12.854}
\end{eqnarray}
Together with \ref{12.846} this implies through \ref{12.845} that:
\begin{equation}
\|\Omega^{n-m}T^{m+1}\chf\|_{L^2({\cal K}^{\tau_1})}\leq K_m \tau_1^{c_m-2}
\sup_{\tau\in(0,\tau_1]}\left\{\tau^{-c_m}\|\Omega^{n-m}T^{m+1}\cf\|_{L^2({\cal K}^{\tau})}\right\}
\label{12.855}
\end{equation}
where:
\begin{equation}
K_m=\sqrt{1+\frac{2}{c_m-2}}+\frac{1}{\sqrt{2(c_m-2)}}\sum_{i=1}^{m+1}\left(\begin{array}{c}m+1\\i
\end{array}\right)(i+1)! \  k_i
\label{12.856}
\end{equation}
Note that for $c_m\geq 3$ we have:
\begin{equation}
K_m\leq C_m
\label{12.857}
\end{equation}
where $C_m$ is a constant which depends only on $m$ and is independent of the exponent $c_m$. 
In fact we have:
\begin{equation}
K_m\leq C \ : \ \mbox{independent of $m$}
\label{12.a1}
\end{equation}
if $c_m$ is chosen large enough, depending on $m$.
The estimate \ref{12.855} then implies, replacing $\tau_1$ by any $\tau\in(0,\tau_1]$, in view of 
the definition \ref{12.805}, 
\begin{equation}
\s^{(m+1,n-m)}P(\tau_1)\leq C\sup_{\tau\in(0,\tau_1]}\left\{\tau^{-c_m}\|\Omega^{n-m}T^{m+1}\cf\|_{L^2({\cal K}^{\tau})}\right\}
\label{12.858}
\end{equation}
Substituting on the right the estimate \ref{12.841} and choosing $c_m$ large enough so that in 
reference to the resulting coefficient of $\s^{(m+1,n-m)}P(\tau_1)$ on the right it holds:
\begin{equation}
\frac{C\delta}{\sqrt{2c_m+2}}\leq\frac{1}{2}
\label{12.859}
\end{equation}
we conclude that:
\begin{eqnarray}
&&\s^{(m+1,n-m)}P(\tau_1)
\leq C\left\{\frac{\s^{(m+1,n-m)}V(\tau_1)}{\sqrt{2c_m}}
+\frac{\tau_1^2\s^{(m+1,n-m)}\Ga(\tau_1)}{\sqrt{2c_m+4}}\right.\nonumber\\
&&\hspace{35mm}+\frac{\sqrt{\s^{(Y;m,n-m)}\cA(\tau_1)}}{\sqrt{2c_m}}
+\tau_1^{1/2}\frac{\sqrt{\s^{(E;m,n-m)}\cA(\tau_1)}}{\sqrt{2c_m+1}}\nonumber\\
&&\hspace{70mm}\left.+\frac{\sqrt{\s^{(Y;m,n-m)}\cBb(\tau_1,\tau_1)}}{(a_m+2)}\right\}\nonumber\\
&&\hspace{25mm}+C\tau_1^{1/2}\sqrt{\max\{\s^{(m,n-m)}\cB(\tau_1,\tau_1),
\s^{(m,n-m)}\cBb(\tau_1,\tau_1)\}}\nonumber\\
&&\hspace{25mm}+C\left\{\tau_1^2\s^{(m,n-m)}L(\tau_1)+\tau_1\s^{(m,n-m)}\underline{L}(\tau_1)
+\tau_1\s^{(m,n-m)}M(\tau_1)
\right\}\nonumber\\
&&\label{12.860}
\end{eqnarray}

We turn to the derivation of an appropriate estimate for $\Omega^{n-m}T^{m+1}\cdl^i$ in 
$L^2({\cal K}^{\tau})$. Applying $\Omega^{n-m}T^m$ to equation \ref{12.503} we obtain, 
to leading terms, 
\begin{eqnarray}
&&\Omega^{n-m}T^{m+1}\cdl^i=U^i\Omega^{n-m}T^m\clab+\lambdab\Omega^{n-m}T^m\check{U}^i\nonumber\\
&&\hspace{25mm}+m(TU^i)\Omega^{n-m}T^{m-1}\clab+m(T\lambdab)\Omega^{n-m}T^{m-1}\check{U}^i\nonumber\\
&&\label{12.861}
\end{eqnarray}
where we denote:
\begin{equation}
\check{U}^i=U^i-U_N^i
\label{12.862}
\end{equation}
In the right hand side of \ref{12.861}, while the first two terms are the principal terms being of 
order $n$, their coefficients are bounded in absolute value by $C\tau$, 
while the next two terms even though of order $n-1$ make a comparable contribution because 
their coefficients are bounded in absolute value only by $C$. 

The principal part of the first term on the right in \ref{12.861} is:
\begin{equation}
\sh^{(n-m)/2}U^i\s^{(m,n-m)}\clab
\label{12.863}
\end{equation}
In view of \ref{12.516} this is bounded in $L^2({\cal K}^{\tau})$ by:
\begin{equation}
C\tau^{c_m+1}\oLab(\tau)
\label{12.864}
\end{equation}
In regard to the third term on the right in \ref{12.861} we apply inequality \ref{12.366} 
with $\Omega^{n-m}T^{m-1}\clab$ in the role of $f$ to deduce:
\begin{equation}
\|\Omega^{n-m}T^{m-1}\clab\|_{L^2({\cal K}^{\tau})}\leq k
\left(\int_0^{\tau}\tau^\prime\|\Omega^{n-m}T^m\clab\|^2_{L^2({\cal K}^{\tau^\prime})}d\tau^\prime
\right)^{1/2}
\label{12.865}
\end{equation}
It follows that the third term on the right in \ref{12.861} is bounded in $L^2({\cal K}^{\tau})$ by:
\begin{equation}
\frac{Cm\tau^{c_m+1}}{\sqrt{2c_m+2}}\s^{(m,n-m)}\oLab(\tau)
\label{12.866}
\end{equation}
In regard to the second term on the right in \ref{12.861} applying $\Omega^{n-m}T^m$ to \ref{12.504}, 
we obtain, to leading terms, 
\begin{eqnarray}
&&\Omega^{n-m}T^m\check{U}^i=c^{-1}\Omega^{n-m}T^m\check{N}^i+c^{-1}r\Omega^{n-m}T^m\check{\Nb}^i
\nonumber\\
&&\hspace{23mm}+c^{-1}(\Nb^i-N_0^i)\Omega^{n-m}T^m\check{r}\nonumber\\
&&\hspace{23mm}-U^i c^{-1}\Omega^{n-m}T^m\check{c}
-m(TU^i)c^{-1}\Omega^{n-m}T^{m-1}\check{c}\nonumber\\
&&\label{12.867}
\end{eqnarray}
where we denote:
\begin{equation}
\check{N}^i=N^i-N_N^i, \ \ \ \check{\Nb}^i=\Nb^i-\Nb_N^i, \ \ \ \check{c}=c-c_N
\label{12.868}
\end{equation}
In regard to the first two terms on the right in \ref{12.867}, the principal parts of 
$\Omega^{n-m}T^m\check{N}^i$, $\Omega^{n-m}T^m\check{\Nb}^i$ are 
\begin{equation}
\sh^{(n-m)/2}(E^{n-m}T^m N^i-E_N^{n-m}T^m N_N^i), \ \ \ 
\sh^{(n-m)/2}(E^{n-m}T^m\Nb^i-E_N^{n-m}T^m\Nb_N^i)
\label{12.869}
\end{equation}
respectively, and by \ref{12.66}, \ref{12.107}, \ref{12.108}, \ref{12.110}, \ref{12.111} the 
principal acoustical parts of these are
\begin{equation}
2\sh^{(n-m)/2}(E^i-\pi N^i)\cmu_{m-1,n-m}, \ \ \ 
2\sh^{(n-m)/2}(E^i-\pi\Nb^i)\cmub_{m-1,n-m}
\label{12.870}
\end{equation}
respectively. Hence by \ref{12.802} the 1st term on the right in \ref{12.867} is bounded in 
$L^2({\cal K}^{\tau})$ by:
\begin{equation}
C\tau^{c_m+1}\s^{(m,n-m)}L(\tau)
\label{12.871}
\end{equation}
and by \ref{12.804} the 2nd term on the right in \ref{12.867} is bounded in $L^2({\cal K}^{\tau})$ by:
\begin{equation}
C\tau^{c_m+1}\s^{(m,n-m}\underline{L}(\tau)
\label{12.872}
\end{equation}
On the other hand, by \ref{12.812} the 3rd term on the right in \ref{12.867} is bounded in 
$L^2({\cal K}^{\tau})$ by:
\begin{equation}
C\tau^{c_m}\s^{(m,n-m)}R(\tau)
\label{12.873}
\end{equation}
As for the 4th term on the right in \ref{12.867}, by \ref{12.69} together with \ref{12.107}, 
\ref{12.108}, \ref{12.110}, \ref{12.111} the principal acoustical part of 
$c^{-1}\Omega^{n-m}T^m\check{c}$ is:
\begin{equation}
-2\pi \sh^{(n-m)/2}(\cmu_{m-1,n-m}+\cmub_{m-1,n-m})
\label{12.874}
\end{equation}
hence by \ref{12.802}, \ref{12.804} and in view of \ref{12.516} the 4th term on the right in 
\ref{12.867} is bounded in $L^2({\cal K}^{\tau})$ by:
\begin{equation}
C\tau^{c_m+2}\s^{(m,n-m)}L(\tau)
+C\tau^{c_m+1}\s^{(m,n-m)}\underline{L}(\tau)
\label{12.875}
\end{equation}
Moreover, using the inequality \ref{12.366} we can show that the 5th term on the right in \ref{12.867}
is bounded in $L^2({\cal K}^{\tau})$ by $Cm/\sqrt{2c_m+2}$ times \ref{12.875}. Note that by condition 
\ref{12.820} $m/\sqrt{2c_m}$ does not exceed a fixed constant. In view of the above results, 
recalling also \ref{12.813} and Proposition 12.4 we conclude that:
\begin{eqnarray}
&&\|\Omega^{n-m}T^m\check{U}^i\|_{L^2({\cal K}^{\tau})}\leq C\tau^{c_m}\left\{
\frac{\s^{(m+1,n-m)}V(\tau)}{\sqrt{2c_m}}\right.\nonumber\\
&&\hspace{47mm}+\frac{\tau\s^{(m+1,n-m)}P(\tau)}{\sqrt{2c_m+2}}
+\frac{\tau^2\s^{(m+1,n-m)}\Ga(\tau)}{\sqrt{2c_m+4}}\nonumber\\
&&\hspace{40mm}+\frac{\sqrt{\s^{(Y;m,n-m)}\cA(\tau)}}{\sqrt{2c_m}}
+\frac{\tau^{1/2}\sqrt{\s^{(E;m,n-m)}\cA(\tau)}}{\sqrt{2c_m+1}}\nonumber\\
&&\hspace{27mm}+\tau\s^{(m,n-m)}L(\tau)+\frac{\tau^2}{\sqrt{2c_m+2}}\s^{(m,n-m)}L(\tau)
\nonumber\\
&&\hspace{47mm}+\frac{\tau}{\sqrt{2c_m+2}}\frac{\sqrt{\s^{(Y;m,n-m)}\cBb(\tau,\tau)}}{(a_m+2)}
\nonumber\\
&&\hspace{27mm}\left.+\frac{C\tau^{3/2}}{\sqrt{2c_m+2}}\sqrt{\max\{\s^{(m,n-m)}\cB(\tau,\tau),
\s^{(m,n-m)}\cBb(\tau,\tau)\}}\right\}\nonumber\\
&&\label{12.876}
\end{eqnarray}
The second term on the right in \ref{12.861} is then bounded in $L^2({\cal K}^{\tau})$ by $C\tau$ 
times the right hand side of \ref{12.876}, and using inequality \ref{12.366} we deduce that the 
fourth term on the right in \ref{12.861} is bounded in $L^2({\cal K}^{\tau})$ by $Cm\tau/\sqrt{2c_m+2}$ 
times the right hand side of \ref{12.876}. Combining with the previous results \ref{12.864} and 
\ref{12.866} and recalling the estimate for $\s^{(m,n-m)}\oLab(\tau)$ of Proposition 12.4 and the 
fact that $m/\sqrt{2c_m}$ does not exceed a fixed constant we conclude that $\Omega^{n-m}T^{m+1}\cdl^i$ 
satisfies the estimate:
\begin{eqnarray}
&&\sup_{\tau\in(0,\tau_1]}\left\{\tau^{-c_m-1}\|\Omega^{n-m}T^{m+1}\cdl^i\|_{L^2({\cal K}^{\tau})}
\right\}\leq C\frac{\s^{(m+1,n-m)}V(\tau_1)}{\sqrt{2c_m}}\nonumber\\
&&\hspace{40mm}+C\frac{\tau_1\s^{(m+1,n-m)}P(\tau_1)}{\sqrt{2c_m+2}}
+C\frac{\tau_1^2\s^{(m+1,n-m)}\Ga(\tau_1)}{\sqrt{2c_m+4}}\nonumber\\
&&\hspace{35mm}+C\frac{\sqrt{\s^{(Y;m,n-m)}\cA(\tau_1)}}{\sqrt{2c_m}}
+C\tau_1^{1/2}\frac{\sqrt{\s^{(E;m,n-m)}\cA(\tau_1)}}{\sqrt{2c_m+1}}\nonumber\\
&&\hspace{50mm}+C\frac{\sqrt{\s^{(Y;m,n-m)}\cBb(\tau_1,\tau_1)}}{(a_m+2)}\nonumber\\
&&\hspace{25mm}+C\tau_1^{1/2}\sqrt{\max\{\s^{(m,n-m)}\cB(\tau_1,\tau_1),
\s^{(m,n-m)}\cBb(\tau_1,\tau_1)\}}\nonumber\\
&&\hspace{25mm}+C\tau_1\s^{(m,n-m)}L(\tau_1)+C\tau_1\s^{(m,n-m)}\underline{L}(\tau_1)
+C\tau_1\s^{(m,n-m)}M(\tau_1)\nonumber\\
&&\label{12.877}
\end{eqnarray}

Next, recalling \ref{12.524} and noting that, for every non-negative integer $i$ we have
$$\frac{\partial^i (\tau^{-3})}{\partial\tau^i}=(-1)^i \frac{(i+2)!}{2} \tau^{-i-3}$$
we obtain:
\begin{equation}
\Omega^{n-m}T^{m+1}\chdl^i=\sum_{i=0}^{m+1}\left(\begin{array}{c}m+1\\i\end{array}\right)(-1)^i
\frac{(i+2)!}{2}\tau^{-i-3}\Omega^{n-m}T^{m+1-i}\cdl^i
\label{12.878}
\end{equation}
It follows that:
\begin{equation}
\|\Omega^{n-m}T^{m+1}\chdl^i\|_{L^2({\cal K}^{\tau_1})}\leq \sum_{i=0}^{m+1}
\left(\begin{array}{c}m+1\\i\end{array}\right)\frac{(i+2)!}{2}
\|\tau^{-i-3}\Omega^{n-m}T^{m+1-i}\cdl^i\|_{L^2({\cal K}^{\tau_1})}
\label{12.879}
\end{equation}
Consider first the $i=0$ term in the sum in \ref{12.879}. We have:
\begin{eqnarray}
&&\|\tau^{-3}\Omega^{n-m}T^{m+1}\cdl^i\|^2_{L^2({\cal K}^{\tau_1})}=
\int_0^{\tau_1}\tau^{-6}\|\Omega^{n-m}T^{m+1}\cdl^i\|^2_{L^2(S_{\tau,\tau})}d\tau\nonumber\\
&&\hspace{30mm}=\int_0^{\tau_1}\tau^{-6}\frac{\partial}{\partial\tau}
\left(\|\Omega^{n-m}T^{m+1}\cdl^i\|^2_{L^2({\cal K}^{\tau})}\right)d\tau\nonumber\\
&&\hspace{15mm}=\tau_1^{-6}\|\Omega^{n-m}T^{m+1}\cdl^i\|^2_{L^2({\cal K}^{\tau_1})}
+6\int_0^{\tau_1}\tau^{-7}\|\Omega^{n-m}T^{m+1}\cdl^i\|^2_{L^2({\cal K}^{\tau})}d\tau\nonumber\\
&&\hspace{15mm}\leq\left(\sup_{\tau\in(0,\tau_1]}\left\{\tau^{-c_m-1}\|\Omega^{n-m}T^{m+1}\cdl^i\|_{L^2({\cal K}^{\tau})}
\right\}\right)^2\cdot\left(1+\frac{3}{c_m-2}\right)\tau_1^{2c_m-4}\nonumber\\
&&\label{12.880}
\end{eqnarray}
Consider next the $i=1$ term in the sum in \ref{12.879}. Applying \ref{12.366} with 
$\Omega^{n-m}T^m\cdl^i$ in the role of $f$ gives, for all $\tau\in(0,\tau_1]$,
\begin{eqnarray}
&&\|\Omega^{n-m}T^m\cdl^i\|_{L^2(S_{\tau,\tau})}\leq k\tau^{1/2}
\|\Omega^{n-m}T^{m+1}\cdl^i\|_{L^2({\cal K}^{\tau})}\nonumber\\
&&\hspace{25mm}\leq k\tau^{c_m+\frac{3}{2}}\sup_{\tau\in(0,\tau_1]}\left\{\tau^{-c_m-1}
\|\Omega^{n-m}T^{m+1}\cdl^i\|_{L^2({\cal K}^{\tau})}\right\}\nonumber\\
&&\label{12.881}
\end{eqnarray}
hence:
\begin{eqnarray}
&&\|\tau^{-4}\Omega^{n-m}T^m\cdl^i\|_{L^2({\cal K}^{\tau_1})}
=\left(\int_0^{\tau_1}\tau^{-8}\|\Omega^{n-m}T^m\cdl^i\|^2_{L^2(S_{\tau,\tau})}d\tau\right)^{1/2}
\nonumber\\
&&\hspace{20mm}\leq k\sup_{\tau\in(0,\tau_1]}\left\{\tau^{-c_m-1}
\|\Omega^{n-m}T^{m+1}\cdl^i\|_{L^2({\cal K}^{\tau})}\right\}\left(\int_0^{\tau_1}\tau^{2c_m-5}d\tau
\right)^{1/2}\nonumber\\
&&\hspace{20mm}=\frac{k\tau_1^{c_m-2}}{\sqrt{2(c_m-2)}}\sup_{\tau\in(0,\tau_1]}\left\{\tau^{-c_m-1}
\|\Omega^{n-m}T^{m+1}\cdl^i\|_{L^2({\cal K}^{\tau})}\right\}\nonumber\\
&&\label{12.882}
\end{eqnarray}
From \ref{12.881} we proceed by induction on $i$ to show that for $i=1,...,m+1$ there are positive 
constants $l_i$ such that for all $\tau\in(0,\tau_1]$ we have:
\begin{equation}
\|\Omega^{n-m}T^{m+1-i}\cdl^i\|_{L^2(S_{\tau,\tau})}\leq l_i\tau^{c_m+i+\frac{1}{2}}
\sup_{\tau\in(0,\tau_1]}\left\{\tau^{-c_m}-1\|\Omega^{n-m}T^{m+1}\cdl^i\|_{L^2({\cal K}^{\tau})}\right\}
\label{12.883}
\end{equation}
According to \ref{12.881} this is true for $i=1$ with 
\begin{equation}
l_1=k
\label{12.884}
\end{equation}
Let it be true for $i=j-1$, $j\geq 2$. Applying then the first inequality in \ref{12.366} with $\Omega^{n-m}T^{m+1-j}\cdl^i$ 
in the role of $f$ gives, for all $\tau\in(0,\tau_1]$:
\begin{eqnarray}
&&\|\Omega^{n-m}T^{m+1-j}\cdl^i\|_{L^2(S_{\tau,\tau})}\leq k\tau^{1/2}
\|\Omega^{n-m}T^{m+2-j}\cdl^i\|_{L^2({\cal K}^{\tau})}\nonumber\\
&&\hspace{30mm}=k\int_0^{\tau}
\|\Omega^{n-m}T^{m+2-j}\cdl^i\|_{L^2(S_{\tau^\prime,\tau^\prime})}d\tau^\prime\nonumber\\
&&\hspace{20mm}\leq kl_{j-1}\sup_{\tau\in(0,\tau_1]}\left\{\tau^{-c_m-1}
\|\Omega^{n-m}T^{m+1}\cdl^i\|_{L^2({\cal K}^{\tau})}\right\}\int_0^{\tau}\tau^{\prime c_m+j-\frac{1}{2}}
d\tau^\prime\nonumber\\
&&\hspace{20mm}=\frac{k l_{j-1}}{c_m+j+\frac{1}{2}}\tau^{c_m+j+\frac{1}{2}}
\sup_{\tau\in(0,\tau_1]}\left\{\tau^{-c_m-1}\|\Omega^{n-m}T^{m+1}\cdl^i\|_{L^2({\cal K}^{\tau})}\right\}
\nonumber\\
&&\label{12.885}
\end{eqnarray}
therefore \ref{12.883} is true for $i=j$ if we set:
\begin{equation}
l_j=\frac{kl_{j-1}}{c_m+j+\frac{1}{2}}
\label{12.886}
\end{equation}
This is a recursion with starting point \ref{12.884}. It follows that for each $i=1,...,m+1$:
\begin{equation}
l_i= \frac{k^i}{\left(c_m+\frac{5}{2}\right)\cdot\cdot\cdot\left(c_m+i+\frac{1}{2}\right)}
\label{12.887}
\end{equation}
By virtue of \ref{12.883} we can estimate:
\begin{eqnarray}
&&\|\tau^{-i-3}\Omega^{n-m}T^{m+1-i}\cdl^i\|_{L^2({\cal K}^{\tau_1})}=
\left(\int_0^{\tau_1}\tau^{-2i-6}\|\Omega^{n-m}T^{m+1-i}\cdl^i\|^2_{L^2(S_{\tau,\tau})}d\tau\right)^{1/2}
\nonumber\\
&&\hspace{30mm}\leq l_i\sup_{\tau\in(0,\tau_1]}\left\{\tau^{-c_m-1}\|\Omega^{n-m}T^{m+1}\cdl^i\|_{L^2({\cal K}^{\tau})}\right\}\left(\int_0^{\tau_1}\tau^{2c_m-5}d\tau\right)^{1/2}\nonumber\\
&&\hspace{30mm}=\frac{l_i\tau_1^{c_m-2}}{\sqrt{2(c_m-2)}}\sup_{\tau\in(0,\tau_1]}\left\{\tau^{-c_m-1}\|\Omega^{n-m}T^{m+1}\cdl^i\|_{L^2({\cal K}^{\tau})}\right\}\nonumber\\
&&\label{12.888}
\end{eqnarray}
Together with \ref{12.880} this implies through \ref{12.879} that:
\begin{equation}
\|\Omega^{n-m}T^{m+1}\chdl^i\|_{L^2({\cal K}^{\tau_1})}\leq L_m \tau_1^{c_m-2}
\sup_{\tau\in(0,\tau_1]}\left\{\tau^{-c_m-1}\|\Omega^{n-m}T^{m+1}\cdl^i\|_{L^2({\cal K}^{\tau})}\right\}
\label{12.889}
\end{equation}
where:
\begin{equation}
L_m=\sqrt{1+\frac{3}{c_m-2}}+\frac{1}{\sqrt{2(c_m-2)}}\sum_{i=1}^{m+1}\left(\begin{array}{c}m+1\\i
\end{array}\right)\frac{(i+2)!}{2} \  l_i
\label{12.890}
\end{equation}
Note that for $c_m\geq 3$ we have:
\begin{equation}
L_m\leq C_m
\label{12.891}
\end{equation}
where $C_m$ is a constant which depends only on $m$ and is independent of the exponent $c_m$. 
In fact we have:
\begin{equation}
L_m\leq C \ : \ \mbox{independent of $m$}
\label{12.a2}
\end{equation}
if $c_m$ is chosen large enough, depending on $m$.
Defining then:
\begin{equation}
\s^{(m+1,n-m)}Q(\tau_1)=\sup_{\tau\in(0,\tau_1]}
\left\{\tau^{-c_m+2}\sum_i\|\Omega^{n-m}T^{m+1}\chdl^i\|_{L^2({\cal K}^{\tau})}\right\} 
\label{12.892}
\end{equation}
and substituting the estimate \ref{12.877} in which the estimate \ref{12.860} for 
$\s^{(m+1,n-m)}P(\tau_1)$ is substituted we conclude that:
\begin{eqnarray}
&&\s^{(m+1,n-m)}Q(\tau_1)
\leq C\left\{\frac{\s^{(m+1,n-m)}V(\tau_1)}{\sqrt{2c_m}}
+\frac{\tau_1^2\s^{(m+1,n-m)}\Ga(\tau_1)}{\sqrt{2c_m+4}}\right.\nonumber\\
&&\hspace{35mm}+\frac{\sqrt{\s^{(Y;m,n-m)}\cA(\tau_1)}}{\sqrt{2c_m}}
+\tau_1^{1/2}\frac{\sqrt{\s^{(E;m,n-m)}\cA(\tau_1)}}{\sqrt{2c_m+1}}\nonumber\\
&&\hspace{70mm}\left.+\frac{\sqrt{\s^{(Y;m,n-m)}\cBb(\tau_1,\tau_1)}}{(a_m+2)}\right\}\nonumber\\
&&\hspace{25mm}+C\tau_1^{1/2}\sqrt{\max\{\s^{(m,n-m)}\cB(\tau_1,\tau_1),
\s^{(m,n-m)}\cBb(\tau_1,\tau_1)\}}\nonumber\\
&&\hspace{25mm}+C\left\{\tau_1\s^{(m,n-m)}L(\tau_1)+\tau_1\s^{(m,n-m)}\underline{L}(\tau_1)
+\tau_1\s^{(m,n-m)}M(\tau_1)
\right\}\nonumber\\
&&\label{12.893}
\end{eqnarray}

We now bring in the identification equations \ref{12.537}. Differentiating equations \ref{12.538} 
implicitly, first with respect to $\tau$ $m$ times and then with respect to $\vartheta$ $n-m$ times, 
we obtain, recalling that $T=\partial/\partial\tau$, $\Omega=\partial/\partial\vartheta$, equations 
of the form:
\begin{equation}
\frac{\partial\hat{F}^i}{\partial v}\Omega^{n-m}T^{m+1}\cv+
\frac{\partial\hat{F}^i}{\partial\gamma}\Omega^{n-m}T^{m+1}\cga=R^i_{m,n-m}
\label{12.894}
\end{equation}
where the leading part of $R^i_{m,n-m}$ is:
\begin{equation}
-S_0^i lv\Omega^{n-m}T^{m+1}\chf+\Omega^{n-m}T^{m+1}\chdl^i
\label{12.895}
\end{equation}
(compare with \ref{12.543}, \ref{12.544}). Recalling the decomposition \ref{12.550}, 
we decompose equations \ref{12.894} into their 
$S_0$ and $\Omega^\prime_0$ components. In view of \ref{12.547} these are:
\begin{eqnarray}
&&\left(\frac{k}{3}+e_\bot\right)\Omega^{n-m}T^{m+1}\cv+h_\bot\Omega^{n-m}T^{m+1}\cga=R_{m,n-m\bot}
\nonumber\\
&&e_{||}\Omega^{n-m}T^{m+1}\cv+(1+h_{||})\Omega^{n-m}T^{m+1}\cga=R_{m,n-m||}
\label{12.896}
\end{eqnarray}
(compare with \ref{12.553}, \ref{12.554}). Solving for 
$(\Omega^{n-m}T^{m+1}\cv,\Omega^{n-m}T^{m+1}\cga)$ we obtain:
\begin{eqnarray}
&&\Omega^{n-m}T^{m+1}\cv=\frac{R_{m,n-m\bot}-(h_{\bot}/(1+h_{||}))R_{m,n-m||}}
{(k/3)+e_\bot-((e_{||}h_\bot)/(1+h_{||}))} \nonumber\\
&&\Omega^{n-m}T^{m+1}\cga=\frac{R_{m,n-m||}-(e_{||}/((k/3)+e_\bot))R_{m,n-m\bot}}
{1+h_{||}-((e_{||}h_\bot)/((k/3)+e\bot))}
\label{12.897}
\end{eqnarray}
In view of \ref{12.552} these imply that, pointwise:
\begin{eqnarray}
&&|\Omega^{n-m}T^{m+1}\cv|\leq C\sum_i|R^i_{m,n-m}|\nonumber\\
&&|\Omega^{n-m}T^{m+1}\cga|\leq C\sum_i|R^i_{m,n-m}|
\label{12.898}
\end{eqnarray}
Hence, taking $L^2$ norms on ${\cal K}^{\tau}$:
\begin{eqnarray}
&&\|\Omega^{n-m}T^{m+1}\cv\|_{L^2({\cal K}^{\tau})}\leq C\sum_i\|R^i_{m,n-m}\|_{L^2({\cal K}^{\tau})}
\nonumber\\
&&\|\Omega^{n-m}T^{m+1}\cga\|_{L^2({\cal K}^{\tau})}\leq C\sum_i\|R^i_{m,n-m}\|_{L^2({\cal K}^{\tau})}
\label{12.899}
\end{eqnarray}
Now, by \ref{12.895} and the definitions \ref{12.805} and \ref{12.892} we have, to leading terms,
\begin{eqnarray}
&&\sum_i\|R^i_{m,n-m}\|_{L^2({\cal K}^{\tau})}\leq C\left(
\|\Omega^{n-m}T^{m+1}\chf\|_{L^2({\cal K}^{\tau})}+
\sum_i\|\Omega^{n-m}T^{m+1}\chdl^i\|_{L^2({\cal K}^{\tau})}\right) \nonumber\\
&&\hspace{32mm}\leq C\tau^{c_m-2}\left(\s^{(m+1,n-m)}P(\tau_1)+\s^{(m+1,n-m)}Q(\tau_1)\right)
\label{12.900}
\end{eqnarray}
for all $\tau\in(0,\tau_1]$, $\tau_1\in(0,\delta]$. Substituting in \ref{12.899}, multiplying by 
$\tau^{-c_m+2}$ and taking the supremum over $\tau\in(0,\tau_1]$ then yields, in view of the 
definitions \ref{12.806}, \ref{12.807}, 
\begin{eqnarray}
&&\s^{(m+1,n-m)}V(\tau_1)\leq C\left(\s^{(m+1,n-m)}P(\tau_1)+\s^{(m+1,n-m)}Q(\tau_1)\right)\nonumber\\
&&\s^{(m+1,n-m)}\Ga(\tau_1)\leq C\left(\s^{(m+1,n-m)}P(\tau_1)+\s^{(m+1,n-m)}Q(\tau_1)\right)\nonumber\\
&&\label{12.901}
\end{eqnarray}
Adding these inequalities and substituting on the right the estimates \ref{12.860} and \ref{12.893} 
for $\s^{(m+1,n-m)}P(\tau_1)$ and $\s^{(m+1,n-m)}Q(\tau_1)$ respectively, we obtain:
\begin{eqnarray}
&&\s^{(m+1,n-m)}V(\tau_1)+\s^{(m+1,n-m)}\Ga(\tau_1)\leq \nonumber\\
&&\hspace{35mm}C\left\{\frac{\s^{(m+1,n-m)}V(\tau_1)}{\sqrt{2c_m}}
+\frac{\tau_1^2\s^{(m+1,n-m)}\Ga(\tau_1)}{\sqrt{2c_m+4}}\right.\nonumber\\
&&\hspace{45mm}+\frac{\sqrt{\s^{(Y;m,n-m)}\cA(\tau_1)}}{\sqrt{2c_m}}
+\tau_1^{1/2}\frac{\sqrt{\s^{(E;m,n-m)}\cA(\tau_1)}}{\sqrt{2c_m+1}}\nonumber\\
&&\hspace{80mm}\left.+\frac{\sqrt{\s^{(Y;m,n-m)}\cBb(\tau_1,\tau_1)}}{(a_m+2)}\right\}\nonumber\\
&&\hspace{25mm}+C\tau_1^{1/2}\sqrt{\max\{\s^{(m,n-m)}\cB(\tau_1,\tau_1),
\s^{(m,n-m)}\cBb(\tau_1,\tau_1)\}}\nonumber\\
&&\hspace{25mm}+C\left\{\tau_1\s^{(m,n-m)}L(\tau_1)+\tau_1\s^{(m,n-m)}\underline{L}(\tau_1)
+\tau_1\s^{(m,n-m)}M(\tau_1)
\right\}\nonumber\\
\label{12.902}
\end{eqnarray}
Choosing $c_m$ large enough so that in regard to the coefficient of $\s^{(m+1,n-m)}V(\tau_1)$ on the 
right we have:
\begin{equation}
\frac{C}{\sqrt{2c_m}}\leq \frac{1}{2}
\label{12.903}
\end{equation}
the coefficient of $\s^{(m+1,n-m)}\Ga(\tau_1)$ on the right is a fortiori not greater than 1/2 and the 
inequality \ref{12.902} then implies:
\begin{eqnarray}
&&\s^{(m+1,n-m)}V(\tau_1)+\s^{(m+1,n-m)}\Ga(\tau_1)\leq \nonumber\\
&&\hspace{38mm}C\left\{\frac{\sqrt{\s^{(Y;m,n-m)}\cA(\tau_1)}}{\sqrt{2c_m}}
+\tau_1^{1/2}\frac{\sqrt{\s^{(E;m,n-m)}\cA(\tau_1)}}{\sqrt{2c_m+1}}\right.\nonumber\\
&&\hspace{80mm}\left.+\frac{\sqrt{\s^{(Y;m,n-m)}\cBb(\tau_1,\tau_1)}}{(a_m+2)}\right\}\nonumber\\
&&\hspace{25mm}+C\tau_1^{1/2}\sqrt{\max\{\s^{(m,n-m)}\cB(\tau_1,\tau_1),
\s^{(m,n-m)}\cBb(\tau_1,\tau_1)\}}\nonumber\\
&&\hspace{25mm}+C\left\{\tau_1\s^{(m,n-m)}L(\tau_1)+\tau_1\s^{(m,n-m)}\underline{L}(\tau_1)
+\tau_1\s^{(m,n-m)}M(\tau_1)
\right\}\nonumber\\
\label{12.904}
\end{eqnarray}
(for new constants $C$). Substituting this in \ref{12.860} we arrive at the following 
proposition. 

\vspace{2.5mm}

\noindent{\bf Proposition 12.6} \ \ Recalling the definitions \ref{12.805} - \ref{12.807}, the 
transformation function differences satisfy to principal terms the following inequalities:
\begin{eqnarray*}
&&\s^{(m+1,n-m)}P(\tau_1), \s^{(m+1,n-m)}V(\tau_1), \s^{(m+1,n-m)}\Ga(\tau_1) \ \leq \\
&&\hspace{38mm}C\left\{\frac{\sqrt{\s^{(Y;m,n-m)}\cA(\tau_1)}}{\sqrt{2c_m}}
+\tau_1^{1/2}\frac{\sqrt{\s^{(E;m,n-m)}\cA(\tau_1)}}{\sqrt{2c_m+1}}\right.\\
&&\hspace{80mm}\left.+\frac{\sqrt{\s^{(Y;m,n-m)}\cBb(\tau_1,\tau_1)}}{(a_m+2)}\right\}\\
&&\hspace{25mm}+C\tau_1^{1/2}\sqrt{\max\{\s^{(m,n-m)}\cB(\tau_1,\tau_1),
\s^{(m,n-m)}\cBb(\tau_1,\tau_1)\}}\\
&&\hspace{25mm}+C\tau_1\left\{\s^{(m,n-m)}L(\tau_1)+\s^{(m,n-m)}\underline{L}(\tau_1)
+\s^{(m,n-m)}M(\tau_1)\right\}
\end{eqnarray*}
for all $\tau_1\in(0,\delta]$ and $m=1,...,n$.

\vspace{5mm}

\section{Estimates for $(\Omega^{n+1}\chf,\Omega^{n+1}\cv,\Omega^{n+1}\cga)$}

Let us recall from \ref{4.53} that the function $f$ is simply the restriction of $x^0$ to ${\cal K}$ 
and from \ref{9.51} that the function $f_N$ is simply the restriction of $x^0_N$ to ${\cal K}$.  
In view of \ref{3.a19} and \ref{10.130} we then have:
\begin{equation}
Ef=\left.E^0\right|_{\cal K}=\left.\pi\right|_{\cal K}, \ \ \ 
E_N f_N=\left.E_N^0\right|_{\cal K}=\left.\pi_N\right|_{\cal K}
\label{12.905}
\end{equation}
Let us define the difference:
\begin{equation}
\s^{(n+1)}\cf=E^{n+1}f-E_N^{n+1}f_N
\label{12.906}
\end{equation}
By \ref{12.905} we have:
\begin{equation}
\s^{(n+1)}\cf=E^{n}\pi-E_N^n\pi_N=E^{n-1}(EE^0)-E_N^{n-1}(E_N E_N^0) \ \ : \ \mbox{on} \ {\cal K}
\label{12.907}
\end{equation}
We also recall from \ref{4.53} that the functions $g^i$ are simply the restrictions of the $x^i$ to 
${\cal K}$ and from \ref{9.51} that the functions $g_N^i$ are simply the restrictions of the $x_N^i$ 
to ${\cal K}$, hence:
\begin{equation}
Eg^i=\left.E^i\right|_{\cal K}, \ \ \ E_N g_N^i=\left.E_N^i\right|_{\cal K}
\label{12.908}
\end{equation}
Let us define the differences:
\begin{equation}
\s^{(n+1)}\check{g}^i=E^{n+1}g^i-E_N^{n+1}g_N^i
\label{12.909}
\end{equation}
By \ref{12.908} we have:
\begin{equation}
\s^{(n+1)}\check{g}^i=E^{n-1}(EE^i)-E_N^{n-1}(E_N E_N^i) \ \ : \ \mbox{on} \ {\cal K}
\label{12.910}
\end{equation}
In regard to $EE^\mu$ we have the expansion \ref{3.a31}, with the coefficients $\se$, $e$, $\oe$ 
given by \ref{3.a32}. We see that:
\begin{equation}
[\se]_{P.A.}=0, \ \ \ [e]_{P.A.}=\frac{1}{2c}\tchib, \ \ \ [\oe]_{P.A.}=\frac{1}{2c}\tchi
\label{12.911}
\end{equation}
hence:
\begin{equation}
[EE^\mu]_{P.A.}=\frac{1}{2c}(\tchib N^\mu+\tchi\Nb^\mu)
\label{12.912}
\end{equation}
In regard to $E_N E_N^\mu$ we have the analogous expansion \ref{10.182} and the coefficient $\se_N$ 
is given by \ref{10.196}. To derive the corresponding formulas for  $e_N$, $\oe_N$, we note that 
according to the expansion \ref{10.182} we have:
$$e_N=-\frac{1}{2c_N}h_{\mu\nu,N}(E_N E_N^\mu)\Nb_N^\nu
=\frac{1}{2c_N}h_{\mu\nu,N}E_N^\mu(E_N\Nb_N^\nu)+\frac{1}{2c_N}(E_Nh_{\mu\nu,N})E_N^\mu\Nb_N^\nu$$
while by the 2nd of \ref{10.117} and by \ref{10.142}:
$$h_{\mu\nu,N}E_N^\mu(E_N\Nb_N^\nu)=\skb_N=\tchib_N-H_N\sbeta_N\ss_{\Nb,N}$$
and by \ref{9.6}:
$$(E_N h_{\mu\nu,N})E_N^\mu\Nb_N^\nu=\sbeta_N\beta_{\Nb,N}E_NH_N+H_N(\beta_{\Nb,N}\sss_N+\sbeta_N
\ss_{\Nb,N})$$
Substituting, we then obtain:
\begin{equation}
e_N=\frac{1}{2c_N}\left(\tchib_N+\beta_{\Nb,N}(\sbeta_N E_N H_N+H_N\sss_N)\right)
\label{12.913}
\end{equation}
Similarly we derive the conjugate formula:
\begin{equation}
\oe_N=\frac{1}{2c_N}\left(\tchi_N+\beta_{N,N}(\sbeta_N E_N H_N+H_N\sss_N)\right)
\label{12.914}
\end{equation}
The above lead to the conclusion that $\s^{(n+1)}\cf$ and $\s^{(n+1)}\check{g}^i$ are given by:
\begin{eqnarray}
&&\s^{(n+1)}\cf=\frac{1}{2c}\left(\s^{(n-1)}\ctchi+\s^{(n-1)}\ctchib\right) \label{12.915}\\
&&\s^{(n+1)}\check{g}^i=\frac{1}{2c}\left(\Nb^i\s^{(n-1)}\ctchi+N^i\s^{(n-1)}\ctchib\right) \label{12.916}
\end{eqnarray}
up to terms of order $n$ with vanishing principal acoustical part, which, as elsewhere in the present 
chapter, can be ignored. 

Now, \ref{12.915} yields, to leading terms:
\begin{equation}
\left(\int_0^\tau\tau^{\prime 2}\|\s^{(n+1)}\cf\|^2_{L^2(S_{\tau^\prime,\tau^\prime})}d\tau^\prime
\right)^{1/2}\leq C\left(\|r\s^{(n-1)}\ctchib\|_{L^2({\cal K}^{\tau})}
+\tau\|\s^{(n-1)}\ctchi\|_{L^2({\cal K}^{\tau})}\right)
\label{12.917}
\end{equation}
hence in view of the definitions \ref{12.397}, \ref{12.398}:
\begin{equation}
\left(\int_0^\tau\tau^{\prime 2}\|\s^{(n+1)}\cf\|^2_{L^2(S_{\tau^\prime,\tau^\prime})}d\tau^\prime
\right)^{1/2}\leq C\tau^{c_0}\left(\s^{(n-1)}\oXb(\tau)
+\tau\s^{(n-1)}\oX(\tau)\right)
\label{12.918}
\end{equation}
Since
\begin{eqnarray*}
&&\|\s^{(n+1)}\cf\|^2_{L^2({\cal K}^{\tau})}=\int_0^{\tau}\|\s^{(n+1)}\cf\|^2_
{L^2(S_{\tau^\prime,\tau^\prime})}d\tau^\prime\\
&&\hspace{25mm}=\int_0^{\tau}\tau^{\prime -2}\frac{\partial}{\partial\tau^\prime}
\left(\int_0^{\tau^\prime}\tau^{\prime\prime 2}\|\s^{(n+1)}\cf\|^2_
{L^2(S_{\tau^{\prime\prime},\tau^{\prime\prime}})}
d\tau^{\prime\prime}\right)d\tau^\prime\\
&&\hspace{15mm}=\tau^{-2}\int_0^{\tau}\tau^{\prime 2}\|\s^{(n+1)}\cf\|^2_{L^2(S_{\tau^\prime,\tau^\prime})}
d\tau^\prime\\
&&\hspace{25mm}+2\int_0^{\tau}\tau^{\prime -3}\left(\int_0^{\tau^\prime}\tau^{\prime\prime 2}
\|\s^{(n+1)}\cf\|^2_{L^2(S_{\tau^{\prime\prime},\tau^{\prime\prime}})}d\tau^{\prime\prime}\right)
d\tau^\prime
\end{eqnarray*}
substituting the estimate \ref{12.918} yields, for all $\tau_1\in(0,\delta]$:
\begin{equation}
\sup_{\tau\in(0,\tau_1]}\left\{\tau^{-c_0+1}\|\s^{(n+1)}\cf\|_{L^2({\cal K}^{\tau})}\right\}
\leq C\left(\s^{(n-1)}\oXb(\tau_1)+\tau_1\s^{(n-1)}\oX(\tau_1)\right)
\label{12.919}
\end{equation}
Substituting the estimates of Proposition 12.3 into the estimates for $\s^{(n-1)}\oXb(\tau_1)$, 
$\s^{(n-1)}\oX(\tau_1)$ and the results in turn into the right hand side of \ref{12.919} we obtain:
\begin{eqnarray}
&&\sup_{\tau\in(0,\tau_1]}\left\{\tau^{-c_0+1}\|\s^{(n+1)}\cf\|_{L^2({\cal K}^{\tau})}\right\}
\leq \frac{C}{\sqrt{2c_0}}\left(\sqrt{\s^{(E;0,n)}\cA(\tau_1)}
+\frac{\sqrt{\s^{(Y;0,n)}\cA(\tau_1)}}{\sqrt{2c_0}}\right)\nonumber\\
&&\hspace{38mm}+\frac{C}{(a_0+2)}\left(\sqrt{\s^{(E;0,n)}\cBb(\tau_1,\tau_1)}
+\frac{\sqrt{\s^{(Y;0,n)}\cBb(\tau_1,\tau_1)}}{\sqrt{2c_0}}\right)\nonumber\\
&&\hspace{45mm}+C\tau_1^{1/2}\sqrt{\max\{\s^{(0,n)}\cB(\tau_1,\tau_1),\s^{(0,n)}\cBb(\tau_1,\tau_1)\}}
\nonumber\\
&&\label{12.920}
\end{eqnarray}

Let us define:
\begin{equation}
\s^{(n+1)}\chf=\tau^{-2}\s^{(n+1)}\cf
\label{12.921}
\end{equation}
We have: 
\begin{eqnarray}
&&\|\s^{(n+1)}\chf\|^2_{L^2({\cal K}^{\tau_1})}=\int_0^{\tau_1}\tau^{-4}\|\s^{(n+1)}\cf\|^2
_{L^2(S_{\tau,\tau})}d\tau\nonumber\\
&&\hspace{30mm}=\int_0^{\tau_1}\tau^{-4}\frac{\partial}{\partial\tau}\left(\|\s^{(n+1)}\cf\|^2
_{L^2({\cal K}^{\tau})}\right)d\tau\nonumber\\
&&\hspace{15mm}=\tau_1^{-4}\|\s^{(n+1)}\cf\|^2_{L^2({\cal K}^{\tau_1})}
+4\int_0^{\tau_1}\tau^{-5}\|\s^{(n+1)}\cf\|^2_{L^2({\cal K}^{\tau})}d\tau\nonumber\\
&&\hspace{15mm}\leq\left(\sup_{\tau\in(0,\tau_1]}\left\{\tau^{-c_0+1}\|\s^{(n+1)}\cf\|
_{L^2({\cal K}^{\tau})}\right\}\right)^2\cdot\left(1+\frac{3}{c_0-3}\right)\tau_1^{2c_0-6}\nonumber\\
&&\label{12.922}
\end{eqnarray}
Hence (recalling that $c_0\geq 4$) by the estimate \ref{12.920}: 
\begin{eqnarray}
&&\sup_{\tau\in(0,\tau_1]}\left\{\tau^{-c_0+3}\|\s^{(n+1)}\chf\|_{L^2({\cal K}^{\tau})}\right\}
\leq \frac{C}{\sqrt{2c_0}}\left(\sqrt{\s^{(E;0,n)}\cA(\tau_1)}
+\frac{\sqrt{\s^{(Y;0,n)}\cA(\tau_1)}}{\sqrt{2c_0}}\right)\nonumber\\
&&\hspace{38mm}+\frac{C}{(a_0+2)}\left(\sqrt{\s^{(E;0,n)}\cBb(\tau_1,\tau_1)}
+\frac{\sqrt{\s^{(Y;0,n)}\cBb(\tau_1,\tau_1)}}{\sqrt{2c_0}}\right)\nonumber\\
&&\hspace{45mm}+C\tau_1^{1/2}\sqrt{\max\{\s^{(0,n)}\cB(\tau_1,\tau_1),\s^{(0,n)}\cBb(\tau_1,\tau_1)\}}
\nonumber\\
&&\label{12.923}
\end{eqnarray}

Next, recalling the definition \ref{12.502} of $\cdl^i$ we define in relation to $\s^{(n+1)}\cf$, 
$\s^{(n+1)}\check{g}^i$ (see \ref{12.906}, \ref{12.909}):
\begin{equation}
\s^{(n+1)}\cdl^i=\s^{(n+1)}\check{g}^i-N_0^i\s^{(n+1)}\cf
\label{12.924}
\end{equation}
From \ref{12.915}, \ref{12.916} we obtain:
\begin{equation}
\s^{(n+1)}\cdl^i=\frac{1}{2c}\left\{(\Nb^i-N_0^i)\s^{(n-1)}\ctchi+(N^i-N_0^i)\s^{(n-1)}\ctchib\right\}
\label{12.925}
\end{equation}
up to terms which can be ignored. In view of \ref{12.515} this implies, to leading terms:
\begin{equation}
\|\s^{(n+1)}\cdl^i\|_{L^2({\cal K}^{\tau})}\leq C\left\{\|r\s^{(n-1)}\ctchib\|_{L^2({\cal K}^{\tau})}
+\|\s^{(n-1)}\ctchi\|_{L^2({\cal K}^{\tau})}\right\}
\label{12.926}
\end{equation}
hence in view of the definitions \ref{12.397}, \ref{12.398}:
\begin{equation}
\|\s^{(n+1)}\cdl^i\|_{L^2({\cal K}^{\tau})}\leq C\tau^{c_0}\left(\s^{(n-1)}\oXb(\tau)
+\s^{(n-1)}\oX(\tau)\right)
\label{12.927}
\end{equation}
Substituting the estimates of Proposition 12.3 into the estimates for $\s^{(n-1)}\oXb(\tau_1)$, 
$\s^{(n-1)}\oX(\tau_1)$ and the results in turn into the right hand side of \ref{12.927} we obtain:
\begin{eqnarray}
&&\sup_{\tau\in(0,\tau_1]}\left\{\tau^{-c_0}\sum_i\|\s^{(n+1)}\cdl^i\|_{L^2({\cal K}^{\tau})}\right\}
\leq \frac{C}{\sqrt{2c_0}}\left(\sqrt{\s^{(E;0,n)}\cA(\tau_1)}
+\frac{\sqrt{\s^{(Y;0,n)}\cA(\tau_1)}}{\sqrt{2c_0}}\right)\nonumber\\
&&\hspace{38mm}+\frac{C}{(a_0+2)}\left(\sqrt{\s^{(E;0,n)}\cBb(\tau_1,\tau_1)}
+\frac{\sqrt{\s^{(Y;0,n)}\cBb(\tau_1,\tau_1)}}{\sqrt{2c_0}}\right)\nonumber\\
&&\hspace{45mm}+C\tau_1^{1/2}\sqrt{\max\{\s^{(0,n)}\cB(\tau_1,\tau_1),\s^{(0,n)}\cBb(\tau_1,\tau_1)\}}
\nonumber\\
&&\label{12.928}
\end{eqnarray}

Let us define:
\begin{equation}
\s^{(n+1)}\chdl^i=\tau^{-3}\s^{(n+1)}\cdl^i
\label{12.929}
\end{equation}
We have: 
\begin{eqnarray}
&&\|\s^{(n+1)}\chdl^i\|^2_{L^2({\cal K}^{\tau_1})}=\int_0^{\tau_1}\tau^{-6}\|\s^{(n+1)}\cdl^i\|^2
_{L^2(S_{\tau,\tau})}d\tau\nonumber\\
&&\hspace{30mm}=\int_0^{\tau_1}\tau^{-6}\frac{\partial}{\partial\tau}\left(\|\s^{(n+1)}\cdl^i\|^2
_{L^2({\cal K}^{\tau})}\right)d\tau\nonumber\\
&&\hspace{15mm}=\tau_1^{-6}\|\s^{(n+1)}\cdl^i\|^2_{L^2({\cal K}^{\tau_1})}
+6\int_0^{\tau_1}\tau^{-7}\|\s^{(n+1)}\cdl^i\|^2_{L^2({\cal K}^{\tau})}d\tau\nonumber\\
&&\hspace{15mm}\leq\left(\sup_{\tau\in(0,\tau_1]}\left\{\tau^{-c_0}\|\s^{(n+1)}\cdl^i\|
_{L^2({\cal K}^{\tau})}\right\}\right)^2\cdot\left(1+\frac{3}{c_0-3}\right)\tau_1^{2c_0-6}\nonumber\\
&&\label{12.930}
\end{eqnarray}
Hence (recalling that $c_0\geq 4$) by the estimate \ref{12.928}: 
\begin{eqnarray}
&&\sup_{\tau\in(0,\tau_1]}\left\{\tau^{-c_0+3}\sum_i\|\s^{(n+1)}\chdl^i\|_{L^2({\cal K}^{\tau})}\right\}
\leq \frac{C}{\sqrt{2c_0}}\left(\sqrt{\s^{(E;0,n)}\cA(\tau_1)}
+\frac{\sqrt{\s^{(Y;0,n)}\cA(\tau_1)}}{\sqrt{2c_0}}\right)\nonumber\\
&&\hspace{38mm}+\frac{C}{(a_0+2)}\left(\sqrt{\s^{(E;0,n)}\cBb(\tau_1,\tau_1)}
+\frac{\sqrt{\s^{(Y;0,n)}\cBb(\tau_1,\tau_1)}}{\sqrt{2c_0}}\right)\nonumber\\
&&\hspace{45mm}+C\tau_1^{1/2}\sqrt{\max\{\s^{(0,n)}\cB(\tau_1,\tau_1),\s^{(0,n)}\cBb(\tau_1,\tau_1)\}}
\nonumber\\
&&\label{12.931}
\end{eqnarray}

Let us define (see \ref{10.802}):
\begin{equation}
\s^{(n+1)}\cv=E^{n+1}v-E_N^{n+1}v_N, \ \ \ \s^{(n+1)}\cga=E^{n+1}\gamma-E_N^{n+1}\gamma_N
\label{12.932}
\end{equation}
Differentiating the identification equations \ref{12.534} implicitly with respect to $\vartheta$ 
we obtain:
\begin{equation}
\frac{\partial\hat{F}^i}{\partial v}\frac{\partial v}{\partial\vartheta}
+\frac{\partial\hat{F}^i}{\partial\gamma}\frac{\partial\gamma}{\partial\vartheta}
+\frac{\partial\hat{F}^i}{\partial\vartheta}=0
\label{12.933}
\end{equation}
where the argument of the partial derivatives of $\hat{F}^i$ is 
$((\tau,\vartheta),(v(\tau,\vartheta),\gamma(\tau,\vartheta)))$. Also, differentiating the $N$th 
approximant identification equations \ref{12.535} implicitly with respect to $\vartheta$ we obtain:
\begin{equation}
\frac{\partial\hat{F}^i_N}{\partial v}\frac{\partial v_N}{\partial\vartheta}
+\frac{\partial\hat{F}^i_N}{\partial\gamma}\frac{\partial\gamma_N}{\partial\vartheta}
+\frac{\partial\hat{F}^i_N}{\partial\vartheta}=\frac{\partial D^i_N}{\partial\vartheta}
\label{12.934}
\end{equation}
where the argument of the partial derivatives of $\hat{F}^i_N$ is 
$((\tau,\vartheta),(v_N(\tau,\vartheta,\gamma_N(\tau,\vartheta)))$. Since 
$$E=\frac{1}{\sqrt{\sh}}\frac{\partial}{\partial\vartheta}$$
multiplying \ref{12.933} by $1/\sqrt{\sh}$ the equation takes the form: 
\begin{equation}
\frac{\partial\hat{F}^i}{\partial v}Ev+\frac{\partial\hat{F}^i}{\partial\gamma}E\gamma+E\hat{F}^i=0
\label{12.935}
\end{equation}
Similarly, since
$$E_N=\frac{1}{\sqrt{\sh_N}}\frac{\partial}{\partial\vartheta}$$
multiplying \ref{12.934} by $1/\sqrt{\sh_N}$ the equation takes the form:
\begin{equation}
\frac{\partial\hat{F}^i_N}{\partial v}E_N v_N
+\frac{\partial\hat{F}^i_N}{\partial\gamma}E_N\gamma_N+E_N\hat{F}^i_N=E_N D^i_N
\label{12.936}
\end{equation}
Applying then $E^n$ as an implicit differential operator to \ref{12.935} we obtain:
\begin{equation}
\frac{\partial\hat{F}^i}{\partial v}E^{n+1}v
+\frac{\partial\hat{F}^i}{\partial\gamma}E^{n+1}\gamma
+E^{n+1}\hat{F}^i=G^i_n
\label{12.937}
\end{equation}
where, while the left hand side is of top order, $n+1$, the right hand side, $G^i_n$, is only of 
order $n$. Similarly applying $E_N^n$ as an implicit differential operator to \ref{12.936} we obtain:
\begin{equation}
\frac{\partial\hat{F}^i_N}{\partial v}E_N^{n+1}v_N
+\frac{\partial\hat{F}^i_N}{\partial\gamma}E_N^{n+1}\gamma_N+E_N^{n+1}\hat{F}^i_N
=G^i_{n,N}+E_N^{n+1}D^i_N
\label{12.938}
\end{equation}
where $G^i_{n,N}$ corresponds to $G^i_n$ for the $N$th approximants. Subtracting \ref{12.938} from 
\ref{12.937} then yields the following equations for $(\s^{(n+1)}\cv,\s^{(n+1)}\cga)$:
\begin{equation}
\frac{\partial\hat{F}^i}{\partial v}\s^{(n+1)}\cv+\frac{\partial\hat{F}^i}{\partial\gamma}\s^{(n+1)}\cga
=I^i_n
\label{12.939}
\end{equation}
Here:
\begin{eqnarray}
&&I^i_n=-\left(\frac{\partial\hat{F}^i}{\partial v}-\frac{\partial\hat{F}^i_N}{\partial v}\right)
E_N^{n+1}v_N
-\left(\frac{\partial\hat{F}^i}{\partial\gamma}-\frac{\partial\hat{F}^i_N}{\partial\gamma}\right)
E_N^{n+1}\gamma_N\nonumber\\
&&\hspace{8mm}-(E^{n+1}\hat{F}^i-E_N^{n+1}\hat{F}^i_N)+(G^i_n-G^i_{n,N})-E_N^{n+1}D^i_N
\label{12.940}
\end{eqnarray}
The principal term in $I^i_n$ is the 3rd term. This term is of top order, $n+1$, while the other terms 
in $I^i_n$ are of lower order. From Proposition 4.5 and from \ref{9.117}, noting that in 
calculating $E^{n+1}\hat{F}^i$ and $E_N^{n+1}\hat{F}^i_N$, $(v,\gamma)$ are held fixed, 
we deduce that the principal part of this term is:
\begin{eqnarray}
&&-S_0^i lv (E^{n+1}\hat{f}-E_N^{n+1}\hat{f}_N)+E^{n+1}\hat{\delta}^i-E_N^{n+1}\hat{\delta}^i_N
-\tau(E^{n+1}E^i-E_N^{n+1} E_N^i)\nonumber\\
&&-v\hat{f}(E^{n+1}-E_N^{n+1})(S_0^i l)-\frac{1}{6}v^3(E^{n+1}-E_N^{n+1})(S_0^i k)
-\gamma(E^{n+1}-E_N^{n+1})\Omega^{\prime i}_0 \nonumber\\
&&\label{12.941}
\end{eqnarray}
Here $E^i$ and $E_N^i$ are the remainder functions in Proposition 4.5 and \ref{9.117} respectively, 
the last given by \ref{9.118}. In view of the factor $\tau$, the corresponding term can be absorbed 
by imposing a smallness condition on $\delta$. 

In regard to the first term in \ref{12.941}, we recall that $\hat{f}$ is defined through \ref{4.240} 
and \ref{4.222} as:
\begin{equation}
\hat{f}(\tau,\vartheta)=\tau^{-2}(f(\tau,\vartheta)-f(0,\vartheta))
\label{12.942}
\end{equation}
We also recall that $\hat{f}_N$ is defined through \ref{9.116} and \ref{9.114} as:
\begin{equation}
\hat{f}_N=\tau^{-2}(f_N(\tau,\vartheta)-f(0,\vartheta))
\label{12.943}
\end{equation}
Denoting then by $f_0$ the $\tau$-independent function:
\begin{equation}
f_0(\vartheta)=f(0,\vartheta)
\label{12.944}
\end{equation}
and comparing with \ref{12.906} and \ref{12.921} we obtain, in regard to the first term in \ref{12.941}, 
\begin{equation}
E^{n+1}\hat{f}-E_N^{n+1}\hat{f}_N=\s^{(n+1)}\chf-\tau^{-2}(E^{n+1}-E_N^{n+1})f_0
\label{12.945}
\end{equation}
Consider next the second term in \ref{12.941}, the difference 
$E^{n+1}\hat{\delta}^i-E_N^{n+1}\hat{\delta}^i_N$. We recall that $\hat{\delta}^i$ is defined through 
\ref{4.244}, \ref{4.223}, \ref{4.222} as:
\begin{equation}
\hat{\delta}^i=\tau^{-3}\left(g^i(\tau,\vartheta)-g^i(0,\vartheta)
-N_0^i(\vartheta)(f(\tau,\vartheta)-f(0,\vartheta))\right)
\label{12.946}
\end{equation}
We also recall that $\hat{\delta}^i_N$ is defined through \ref{9.116}, \ref{9.115}, \ref{9.114} as:
\begin{equation}
\hat{\delta}^i_N=\tau^{-3}\left(g^i_N(\tau,\vartheta)-g^i(0,\vartheta)
-N_0^i(\vartheta)(f_N(\tau,\vartheta)-f(0,\vartheta))\right)
\label{12.947}
\end{equation}
Denoting then by $g_0^i$ the $\tau$-independent functions:
\begin{equation}
g_0^i(\vartheta)=g^i(0,\vartheta)
\label{12.948}
\end{equation}
and comparing with \ref{12.906}, \ref{12.909}, \ref{12.924}, \ref{12.929} we obtain:
\begin{eqnarray}
&&E^{n+1}\hat{\delta}^i-E_N^{n+1}\hat{\delta}^i_N=\s^{(n+1)}\chdl^i\nonumber\\
&&\hspace{35mm}-\tau^{-3}\left((E^{n+1}-E_N^{n+1})g_0^i-N_0^i(E^{n+1}-E_N^{n+1})f_0\right)\nonumber\\
&&-\tau^{-3}\sum_{j=1}^{n+1}\left(\begin{array}{c}n+1\\j\end{array}\right)
\left[(E^j N_0^i)E^{n+1-j}(f-f_0)-(E_N^j N_0^i)E_N^{n+1-j}(f_N-f_0)\right]\nonumber\\
&&\label{12.949}
\end{eqnarray}

The 2nd terms in \ref{12.945} and \ref{12.949} as well as the last three terms in \ref{12.941} are of 
the form of a 0th order coefficient times:
\begin{equation}
(E^{n+1}-E_N^{n+1})\phi_0
\label{12.950}
\end{equation}
where $\phi_0$ is a given smooth $\tau$-independent function on ${\cal K}$ defined by the initial data 
on $\partial_-{\cal K}=S_{0,0}$. As for the sum in \ref{12.949}, its principal part is contained 
in the $j=n+1$ term 
$$(E^{n+1}N_0^i)(f-f_0)-(E_N^{n+1}N_0^i)(f_N-f_0)$$
which up to the 0th order term $\tau^2(\hat{f}-\hat{f}_N)E_N^{n+1}N_0^i$ is equal to:
\begin{equation}
\tau^2\hat{f}(E^{n+1}-E_N^{n+1})N_0^i
\label{12.951}
\end{equation}
which is also of the above form. We have:
\begin{equation}
(E^{n+1}-E_N^{n+1})\phi_0=
\left(\frac{1}{\sqrt{\sh}}\frac{\partial}{\partial\vartheta}\right)^{n+1}\phi_0-
\left(\frac{1}{\sqrt{\sh_N}}\frac{\partial}{\partial\vartheta}\right)^{n+1}\phi_0
\label{12.952}
\end{equation}
Hence:
\begin{equation}
[(E^{n+1}-E_N^{n+1})\phi_0]_{P.P.}=-\frac{1}{2}\frac{\partial\phi_0}{\partial\vartheta}
\left[\frac{1}{\sh^{(n+3)/2}}\frac{\partial^n\sh}{\partial\vartheta^n}-
\frac{1}{\sh_N^{(n+3)/2}}\frac{\partial^n\sh_N}{\partial\vartheta^n}\right]_{P.P.}
\label{12.953}
\end{equation}
Introducing the difference quantity: 
\begin{equation}
\s^{(n)}\check{\sh}=\sh^{-(n+1)/2}\Omega^n\log\sh-\sh_N^{-(n+1)/2}\Omega^n\log\sh_N
\label{12.954}
\end{equation}
we can equivalently express \ref{12.953} in the form:
\begin{equation}
[(E^{n+1}-E_N^{n+1})\phi_0]_{P.P.}=-\frac{1}{2}(\Omega\phi_0)[\s^{(n)}\check{\sh}]_{P.P.}
\label{12.955}
\end{equation}
The quantity $\s^{(n)}\check{\sh}$ must be treated as a top order quantity. The point is that 
$\Omega^n\log\sh$ can only be controlled by controlling $\Omega^n T\log\sh$, which is achieved 
through the equation:
\begin{equation}
T\log\sh=2(\chi+\chib)
\label{12.956}
\end{equation}
This equation follows from the first variational formulas of Proposition 3.1 and is also readily 
deduced from the commutation relation:
\begin{equation}
[T,E]=-(\chi+\chib)E
\label{12.957}
\end{equation}
which follows from the first two of the commutation relations \ref{3.a14}. Similarly for the $N$th 
approximants we have the equation:
\begin{equation}
T\log\sh_N=2(\chi_N+\chib_N)
\label{12.958}
\end{equation} 
deduced from the commutation relation:
\begin{equation}
[T,E_N]=-(\chi_N+\chib_N)E_N
\label{12.959}
\end{equation}
following from the first two of the commutation relations \ref{9.b1}.

In view of \ref{12.955}, from \ref{12.945} we have, to principal terms,
\begin{equation}
E^{n+1}\hat{f}-E_N^{n+1}\hat{f}_N=\s^{(n+1)}\chf+\frac{1}{2}\tau^{-2}(\Omega f_0)\s^{(n)}\check{\sh}
\label{12.960}
\end{equation}
and from \ref{12.949}, \ref{12.951} we have, to principal terms, 
\begin{eqnarray}
&&E^{n+1}\hat{\delta}^i-E_N^{n+1}\hat{\delta}^i_N=\s^{(n+1)}\chdl^i\nonumber\\
&&\hspace{32mm}+\frac{1}{2}\tau^{-3}(\Omega g_0^i-N_0^i\Omega f_0)\s^{(n)}\check{\sh}\nonumber\\
&&\hspace{32mm}+\frac{1}{2}\tau^{-1}\hat{f}(\Omega N_0^i)\s^{(n)}\check{\sh}
\label{12.961}
\end{eqnarray}
Moreover, the principal part of the last three terms in \ref{12.941} is contained in:
\begin{equation}
\frac{1}{2}\left\{v\hat{f}\Omega(S_0^i l)+\frac{1}{6}v^3\Omega(S_0^i k)
-\gamma\Omega\Omega_0^{\prime i}\right\}\s^{(n)}\check{\sh}
\label{12.962}
\end{equation}
In regard to the coefficient of the 2nd term on the right in \ref{12.961}, from Proposition 4.3, 
recalling the fact that $\left.\sh\right|_{S_{0,0}}=1$, we have:
\begin{equation}
\Omega g_0^i-N_0^i\Omega f_0=\left.(\Omega^i-N^i\Omega^0)\right|_{S_{0,0}}=\Omega_0^{\prime i}
\label{12.963}
\end{equation}

We conclude from the above that the leading part of $I^i_n$, the right hand side of \ref{12.939}, is:
\begin{equation}
-S_0^i lv\s^{(n+1)}\chf+\s^{(n+1)}\chdl^i
+\frac{1}{2}\left(\tau^{-3}\Omega_0^{\prime i}+\tau^{-2}V^i\right)\s^{(n)}\check{\sh}
\label{12.964}
\end{equation}
where $V^i$ can be assumed to satisfy:
\begin{equation}
\sum_i|V^i|\leq C \ : \ \mbox{in ${\cal K}^{\delta}$}
\label{12.965}
\end{equation}
Recalling the decomposition \ref{12.550}, we decompose equations \ref{12.939} into their $S_0$ and 
$\Omega_0^\prime$ components. In view of \ref{12.547} these are:
\begin{eqnarray}
&&\left(\frac{k}{3}+e_\bot\right)\s^{(n+1)}\cv+h_\bot\s^{(n+1)}\cga=I_{n\bot}\nonumber\\
&&e_{||}\s^{(n+1)}\cv+(1+h_{||})\s^{(n+1)}\cga=I_{n||} 
\label{12.966}
\end{eqnarray}
(compare with \ref{12.553}, \ref{12.554} and with \ref{12.896}). Solving for 
$(\s^{(n+1)}cv,\s^{(n+1)}\cga)$ we obtain:
\begin{eqnarray}
&&\s^{(n+1)}\cv=\frac{I_{n\bot}-(h_\bot/(1+h_{||}))I_{n||}}{(k/3)+e_\bot-((e_{||}h_\bot)/(1+h_{||}))}
\nonumber\\
&&\s^{(n+1)}\cga=\frac{I_{n||}-(e_{||}/((k/3)+e_\bot))I_{n\bot}}{1+h_{||}
-((e_{||}h_{\bot})((k/3)+e_\bot))}
\label{12.967}
\end{eqnarray}
Now, by \ref{12.964} the leading part of $I_{n\bot}$ is:
\begin{equation}
-lv\s^{(n+1)}\chf+\s^{(n+1)}\chdl_\bot+\frac{1}{2}\tau^{-2}V_\bot\s^{(n)}\check{\sh}
\label{12.968}
\end{equation}
and the leading part of $I_{n||}$ is:
\begin{equation}
\s^{(n+1)}\chdl_{||}+\frac{1}{2}(\tau^{-3}+\tau^{-2}V_{||})\s^{(n)}\check{\sh}
\label{12.969}
\end{equation}
Therefore, in view of \ref{12.552} and \ref{12.965}, the formulas \ref{12.967} yield, to leading 
terms:
\begin{eqnarray}
&&\|\s^{(n+1)}\cv\|_{L^2(S_{\tau,\tau})}\leq C\|\s^{(n+1)}\chf\|_{L^2(S_{\tau,\tau})}
+C\sum_i\|\s^{(n+1)}\chdl^i\|_{L^2(S_{\tau,\tau})}\nonumber\\
&&\hspace{60mm}+C\tau^{-2}\|\s^{(n)}\check{\sh}\|_{L^2(S_{\tau,\tau})}\nonumber\\
&&\|\s^{(n+1)}\cga\|_{L^2(S_{\tau,\tau})}\leq C\tau\|\s^{(n+1)}\chf\|_{L^2(S_{\tau,\tau})}
+C\sum_i\|\s^{(n+1)}\chdl^i\|_{L^2(S_{\tau,\tau})}\nonumber\\
&&\hspace{62mm}+C\tau^{-3}\|\s^{(n)}\check{\sh}\|_{L^2(S_{\tau,\tau})}\nonumber\\
&&\label{12.970}
\end{eqnarray}
In regard to $\s^{(n)}\check{\sh}$, let us introduce, for arbitrary $\tau_1\in(0,\delta]$, the quantity:
\begin{equation}
\s^{(n)}H(\tau_1)=\sup_{\tau\in(0,\tau_1]}\left\{\tau^{-c_0+1}
\|\s^{(n)}\check{\sh}\|_{L^2({\cal K}^{\tau})}\right\}
\label{12.971}
\end{equation}
We then have:
\begin{eqnarray}
&&\int_0^{\tau_1}\tau^{-4}\|\s^{(n)}\check{\sh}\|^2_{L^2(S_{\tau,\tau})}d\tau
=\int_0^{\tau_1}\tau^{-4}\frac{\partial}{\partial\tau}
\left(\|\s^{(n)}\check{\sh}\|^2_{L^2({\cal K}^{\tau})}\right)d\tau\nonumber\\
&&\hspace{15mm}=\tau_1^{-4}\|\s^{(n)}\check{\sh}\|^2_{L^2({\cal K}^{\tau_1})}
+4\int_0^{\tau_1}\tau^{-5}\|\s^{(n)}\check{\sh}\|^2_{L^2({\cal K}^{\tau})}d\tau\nonumber\\
&&\hspace{15mm}\leq\left(\s^{(n)}H(\tau_1)\right)^2 \left(1+\frac{2}{c_0-3}\right) \tau_1^{2c_0-6}
\label{12.972}
\end{eqnarray}
provided that $c_0>3$. We also have:
\begin{eqnarray}
&&\int_0^{\tau_1}\tau^{-6}\|\s^{(n)}\check{\sh}\|^2_{L^2(S_{\tau,\tau})}d\tau
=\int_0^{\tau_1}\tau^{-6}\frac{\partial}{\partial\tau}
\left(\|\s^{(n)}\check{\sh}\|^2_{L^2({\cal K}^{\tau})}\right)d\tau\nonumber\\
&&\hspace{15mm}=\tau_1^{-6}\|\s^{(n)}\check{\sh}\|^2_{L^2({\cal K}^{\tau_1})}
+6\int_0^{\tau_1}\tau^{-7}\|\s^{(n)}\check{\sh}\|^2_{L^2({\cal K}^{\tau})}d\tau\nonumber\\
&&\hspace{15mm}\leq\left(\s^{(n)}H(\tau_1)\right)^2 \left(1+\frac{3}{c_0-4}\right) \tau_1^{2c_0-8}
\label{12.973}
\end{eqnarray}
provided that $c_0>4$. We shall in fact take:
\begin{equation}
c_0\geq 5
\label{12.974}
\end{equation}
Replacing $\tau$ by $\tau^\prime\in(0,\tau]$ in \ref{12.970} and taking $L^2$ norms with respect 
to $\tau^\prime$ on $(0,\tau)$, then yields, taking account of \ref{12.972} and \ref{12.973} 
with $(\tau^\prime,\tau)$ in the role of $(\tau,\tau_1)$, 
\begin{eqnarray}
&&\|\s^{(n+1)}\cv\|_{L^2({\cal K}^{\tau})}\leq C\|\s^{(n+1)}\chf\|_{L^2({\cal K}^{\tau})}
+C\sum_i\|\s^{(n+1)}\chdl^i\|_{L^2({\cal K}^{\tau})}\nonumber\\
&&\hspace{56mm}+C\tau^{c_0-3}\s^{(n)}H(\tau)\nonumber\\
&&\|\s^{(n+1)}\cga\|_{L^2({\cal K}^{\tau})}\leq C\tau\|\s^{(n+1)}\chf\|_{L^2{(\cal K}^{\tau})}
+C\sum_i\|\s^{(n+1)}\chdl^i\|_{L^2({\cal K}^{\tau})}\nonumber\\
&&\hspace{58mm}+C\tau^{c_0-4}\s^{(n)}H(\tau) 
\label{12.975}
\end{eqnarray}
Substituting the estimates \ref{12.923} and \ref{12.931}, multiplying the first of the inequalities 
\ref{12.975} by $\tau^{-c_0+3}$ and the second by $\tau^{-c_0+4}$ and taking the supremum over 
$\tau\in(0,\tau_1]$, noting that $\s^{(n)}H(\tau)$ is a non-decreasing function of $\tau$, we arrive 
at the following lemma. 

\vspace{2.5mm}

\noindent{\bf Lemma 12.12} \ \ The quantities $\s^{(n+1)}\cv$, $\s^{(n+1)}\cga$ satisfy to principal 
terms the following inequalities:
\begin{eqnarray*}
&&\sup_{\tau\in(0,\tau_1]}\left\{\tau^{-c_0+3}\|\s^{(n+1)}\cv\|_{L^2({\cal K}^{\tau})}\right\}
\leq C\s^{(n)}H(\tau_1)\\
&&\hspace{50mm}+\frac{C}{\sqrt{2c_0}}\left(\sqrt{\s^{(E;0,n)}\cA(\tau_1)}
+\frac{\sqrt{\s^{(Y;0,n)}\cA(\tau_1)}}{\sqrt{2c_0}}\right)\\
&&\hspace{38mm}+\frac{C}{(a_0+2)}\left(\sqrt{\s^{(E;0,n)}\cBb(\tau_1,\tau_1)}
+\frac{\sqrt{\s^{(Y;0,n)}\cBb(\tau_1,\tau_1)}}{\sqrt{2c_0}}\right)\\
&&\hspace{45mm}+C\tau_1^{1/2}\sqrt{\max\{\s^{(0,n)}\cB(\tau_1,\tau_1),\s^{(0,n)}\cBb(\tau_1,\tau_1)\}}
\end{eqnarray*}
\begin{eqnarray*}
&&\sup_{\tau\in(0,\tau_1]}\left\{\tau^{-c_0+4}\|\s^{(n+1)}\cga\|_{L^2({\cal K}^{\tau})}\right\}
\leq C\s^{(n)}H(\tau_1)\\
&&\hspace{50mm}+\frac{C\tau_1}{\sqrt{2c_0}}\left(\sqrt{\s^{(E;0,n)}\cA(\tau_1)}
+\frac{\sqrt{\s^{(Y;0,n)}\cA(\tau_1)}}{\sqrt{2c_0}}\right)\\
&&\hspace{38mm}+\frac{C\tau_1}{(a_0+2)}\left(\sqrt{\s^{(E;0,n)}\cBb(\tau_1,\tau_1)}
+\frac{\sqrt{\s^{(Y;0,n)}\cBb(\tau_1,\tau_1)}}{\sqrt{2c_0}}\right)\\
&&\hspace{45mm}+C\tau_1^{3/2}\sqrt{\max\{\s^{(0,n)}\cB(\tau_1,\tau_1),\s^{(0,n)}\cBb(\tau_1,\tau_1)\}}
\end{eqnarray*}

\vspace{2.5mm}

The objective now is to derive an appropriate bound for $\s^{(n)}H(\tau_1)$. By \ref{12.956}, 
\ref{12.958},
\begin{equation}
T\Omega^n\log\sh=2(\Omega^n\chi+\Omega^n\chib), \ \ \ 
T\Omega^n\log\sh_N=2(\Omega^n\chi_N+\Omega^n\chib_N)
\label{12.976}
\end{equation}
and we have:
\begin{equation}
\left.\log\sh\right|_{S_{0,0}}=\left.\log\sh_N\right|_{S_{0,0}}=0
\label{12.977}
\end{equation}
hence, integrating, we obtain:
\begin{eqnarray}
&&(\Omega^n\log\sh)(\tau,\vartheta)=2\int_0^{\tau}(\Omega^n\chi+\Omega^n\chib)(\tau^\prime,\vartheta)
d\tau^\prime\nonumber\\
&&(\Omega^n\log\sh_N)(\tau,\vartheta)=2\int_0^{\tau}(\Omega^n\chib_N+\Omega^n\chib_N)
(\tau^\prime,\vartheta)d\tau^\prime
\label{12.978}
\end{eqnarray}
Since to principal terms:
$$\Omega^n\chi=\sh^{n/2} E^n\chi, \ \ \ \Omega^{n}\chib=\sh^{n/2} E^n\chib$$
and correspondingly for the $N$th approximants:
$$\Omega^n\chi_N=\sh_N^{n/2} E_N^n\chi_N, \ \ \ \Omega^{n}\chib_N=\sh_N^{n/2} E_N^n\chib_N$$
by \ref{12.978} we have, to principal terms:
\begin{eqnarray}
&&(\Omega^n\log\sh)(\tau,\vartheta)=2\int_0^{\tau}(\sh(\tau^\prime,\vartheta))^{n/2}
(E^n\chi+E^n\chib)(\tau^\prime,\vartheta)d\tau^\prime\nonumber\\
&&(\Omega^n\log\sh_N)(\tau,\vartheta)=2\int_0^{\tau}(\sh_N(\tau^\prime,\vartheta))^{n/2}
(E_N^n\chi_N+E_N^n\chib_N)(\tau^\prime,\vartheta)d\tau^\prime\nonumber\\
&&\label{12.979}
\end{eqnarray}
It follows that to principal terms it holds:
\begin{equation}
(\sh^{(n+1)/2}\s^{(n)}\check{\sh})(\tau,\vartheta)=
2\int_0^{\tau}(\sh^{n/2}(\s^{(n)}\check{\chi}+\s^{(n)}\check{\chib}))(\tau^\prime,\vartheta)
d\tau^\prime
\label{12.980}
\end{equation}
Here, we have defined:
\begin{equation}
\s^{(n)}\check{\chi}=E^n\chi-E_N^n\chi_N, \ \ \ \s^{(n)}\check{\chib}=E^n\chib-E_N^n\chib_N
\label{12.981}
\end{equation}
Applying inequality \ref{12.366} to \ref{12.980} we deduce:
\begin{equation}
\|\s^{(n)}\check{\sh}\|_{L^2(S_{\tau,\tau})}\leq 2k\tau^{1/2}
\left\{\|\s^{(n)}\check{\chi}\|_{L^2({\cal K}^{\tau})}
+\|\s^{(n)}\check{\chib}\|_{L^2({\cal K}^{\tau})}\right\}
\label{12.982}
\end{equation}
where $k$ is a constant greater than 1, but which can be chosen as close to 1 as we wish by suitably 
restricting $\delta$. Taking $L^2$ norms with respect to $\tau$ on $(0,\tau_1)$, \ref{12.982} implies:
\begin{eqnarray}
&&\|\s^{(n)}\check{\sh}\|_{L^2({\cal K}^{\tau_1})}\leq 
2k\left(\int_0^{\tau_1}\tau\|\s^{(n)}\check{\chi}\|^2_{L^2({\cal K}^{\tau})}d\tau\right)^{1/2}
\nonumber\\
&&\hspace{24mm}+2k\left(\int_0^{\tau_1}\tau\|\s^{(n)}\check{\chib}\|^2_{L^2({\cal K}^{\tau})}d\tau\right)^{1/2}
\label{12.983}
\end{eqnarray}

We recall the formulas \ref{3.a20}:
\begin{eqnarray}
&&\chi=\rho\tchi+\frac{1}{2}\sbeta^2 LH \nonumber\\
&&\chib=\rhob\tchib+\frac{1}{2}\sbeta^2\Lb H \label{12.984}
\end{eqnarray}
and the corresponding formulas \ref{10.138}, \ref{10.141} for the $N$th approximants:
\begin{eqnarray}
&&\chi_N=\rho_N\tchi_N+\frac{1}{2}\sbeta_N^2 L_N H_N +\vep_{\chi,N} \nonumber\\
&&\chib_N=\rhob_N\tchib_N+\frac{1}{2}\sbeta_N^2\Lb_N H_N +\vep_{\chib,N} \label{12.985}
\end{eqnarray}
The formulas \ref{12.984} imply:
\begin{eqnarray}
&&[E^n\chi]_{P.P.}=\rho E^n\tchi+\frac{1}{2}\sbeta^2[E^n LH]_{P.P.} \nonumber\\
&&[E^n\chib]_{P.P.}=\rhob E^n\tchib+\frac{1}{2}\sbeta^2[E^n\Lb H]_{P.P.} \label{12.986}
\end{eqnarray}
and similarly from the formulas \ref{12.985} for the $N$th approximants. It follows that, 
to principal terms:
\begin{eqnarray}
&&\s^{(n)}\check{\chi}=\rho\s^{(n)}\ctchi+\frac{1}{2}\sbeta^2(E^n LH-E_N^n L_N H_N) \nonumber\\
&&\s^{(n)}\check{\chib}=\rhob\s^{(n)}\ctchib+\frac{1}{2}\sbeta^2(E^n\Lb H-E_N^n\Lb_N H_N) \label{12.987}
\end{eqnarray}

We shall first deal with the non-acoustical principal terms, the 2nd terms in each of \ref{12.987}. 
Comparing with \ref{10.621}, \ref{10.622}, the principal part of the 2nd term in the first of 
\ref{12.987} is:
\begin{equation}
-\sbeta^2 H^\prime (g^{-1})^{\mu\nu}\beta_\nu L(E^n\beta_\mu-E_N^n\beta_{\mu,N})
\label{12.988}
\end{equation}
Comparing with \ref{10.364}, \ref{10.365}, the principal part of the 2nd term in the second of 
\ref{12.987} is:
\begin{equation}
-\sbeta^2 H^\prime (g^{-1})^{\mu\nu}\beta_\nu\Lb(E^n\beta_\mu-E_N^n\beta_{\mu,N})
\label{12.989}
\end{equation}
We apply formula \ref{10.368} to express \ref{12.988}, \ref{12.989} as the sum of $E$, $N$, and 
$\Nb$ components. The dominant contribution to the $L^2$ norm on ${\cal K}$ is, in each case, that 
of the $\Nb$ component, which is
\begin{equation}
(\eta^2\beta_N/2c)\sbeta^2 H^\prime\Nb^\mu L(E^n\beta_\mu-E_N^n\beta_{\mu,N})
\label{12.990}
\end{equation}
for \ref{12.988}, and 
\begin{equation}
(\eta^2\beta_N/2c)\sbeta^2 H^\prime\Nb^\mu\Lb(E^n\beta_\mu-E_N^n\beta_{\mu,N})
\label{12.991}
\end{equation}
for \ref{12.989}. In regard to \ref{12.991}, we write, as in \ref{10.373}:
\begin{equation}
N^\mu\Lb(E^n\beta_\mu-E_N^n\beta_{\mu,N})+\ogamma\Nb^\mu\Lb(E^n\beta_\mu-E_N^n\beta_{\mu,N})
=\s^{(Y;0,n)}\cxi_{\Lb}
\label{12.992}
\end{equation}
Since $\ogamma\sim\tau$ along ${\cal K}$, the contribution to 
$\|\Nb^\mu\Lb(E^n\beta_\mu-E_N^n\beta_{\mu,N})\|^2_{L^2({\cal K}^{\tau_1})}$ of the right hand side of 
\ref{12.992} is bounded by a constant multiple of:
\begin{eqnarray}
&&\int_0^{\tau_1}\tau^{-2}\|\s^{(Y;0,n)}\cxi_{\Lb}\|^2_{L^2(S_{\tau,\tau})}d\tau
=\int_0^{\tau_1}\tau^{-2}\frac{\partial}{\partial\tau}
\left(\|\s^{(Y;0,n)}\cxi_{\Lb}\|^2_{L^2({\cal K}^{\tau})}\right)d\tau\nonumber\\
&&=C\left\{\tau_1^{-2}\|\s^{(Y;0,n)}\cxi_{\Lb}\|^2_{L^2({\cal K}^{\tau_1})}
+2\int_0^{\tau_1}\tau^{-3}\|\s^{(Y;0,n)}\cxi_{\Lb}\|^2_{L^2({\cal K}^{\tau})}d\tau\right\}\nonumber\\
&&\hspace{30mm}\leq C\s^{(Y;0,n)}\cA(\tau_1)\cdot\tau_1^{2c_0-2}\left(1+\frac{1}{c_0-1}\right)
\label{12.993}
\end{eqnarray}
(provided that $c_0>1$) where we have appealed to \ref{9.a39} and \ref{9.299}. Hence, in view of 
\ref{12.974}, replacing $\tau_1$ by $\tau$, this dominant contribution to \\
$\|\Nb^\mu\Lb(E^n\beta_\mu-E_N^n\beta_{\mu,N})\|_{L^2({\cal K}^{\tau})}$ is bounded by:
\begin{equation}
C\tau^{c_0-1}\sqrt{\s^{(Y;0,n)}\cA(\tau)}
\label{12.994}
\end{equation}
To estimate the contribution to 
$\|\Nb^\mu\Lb(E^n\beta_\mu-E_N^n\beta_{\mu,N})\|^2_{L^2({\cal K}^{\tau_1})}$ of the first term in the 
left hand side of \ref{12.992}, we express this term as in \ref{10.377} as 
$\lambda\s^{(E;0,n)}\csxi$ up to lower order terms. Then since $\lambda\sim\tau^2$ while 
$a\sim\tau^3$ along ${\cal K}$ this contribution is bounded by a constant multiple of:
\begin{eqnarray}
&&\int_0^{\tau_1}\tau^{-1}\|\sqrt{a}\s^{(E;0,n)}\csxi\|^2_{L^2(S_{\tau,\tau})}d\tau
=\int_0^{\tau_1}\tau^{-1}\frac{\partial}{\partial\tau}
\left(\|\sqrt{a}\s^{(E;0,n)}\csxi\|^2_{L^2({\cal K}^{\tau})}\right)d\tau\nonumber\\
&&=C\left\{\tau_1^{-1}\|\s^{\sqrt{a}(E;0,n)}\csxi\|^2_{L^2({\cal K}^{\tau_1})}
+\int_0^{\tau_1}\tau^{-2}\|\sqrt{a}\s^{(E;0,n)}\csxi\|^2_{L^2({\cal K}^{\tau})}d\tau\right\}\nonumber\\
&&\hspace{30mm}\leq C\s^{(E;0,n)}\cA(\tau_1)\cdot\tau_1^{2c_0-1}\left(1+\frac{1}{2c_0-1}\right)
\label{12.995}
\end{eqnarray}
appealing again to \ref{9.a39} and \ref{9.299}. Hence, replacing $\tau_1$ by $\tau$, the 
corresponding contribution to 
$\|\Nb^\mu\Lb(E^n\beta_\mu-E_N^n\beta_{\mu,N})\|_{L^2({\cal K}^{\tau})}$ is bounded by:
\begin{equation}
C\tau^{c_0-\frac{1}{2}}\sqrt{\s^{(E;0,n)}\cA(\tau)}
\label{12.996}
\end{equation}
In regard to \ref{12.990} we express $\Nb^\mu L(E^n\beta_\mu-E_N^n\beta_{\mu,N})$ as in \ref{10.436} 
as $\lambdab\s^{(E;0,l)}\sxi$ up to lower order terms. Then since $\lambdab\sim\tau$ while 
$a\sim\tau^3$ along ${\cal K}$, the $L^2$ norm on ${\cal K}^{\tau_1}$ of \ref{12.990} is bounded by 
a constant multiple of 
$$\int_0^{\tau_1}\tau^{-1}\|\sqrt{a}\s^{(E;0,n)}\csxi\|^2_{L^2(S_{\tau,\tau})}d\tau$$
therefore by \ref{12.995}, \ref{12.996} the $L^2({\cal K}^{\tau})$ norm of \ref{12.990} is also 
bounded by 
\begin{equation}
C\tau^{c_0-\frac{1}{2}}\sqrt{\s^{(E;0,n)}\cA(\tau)}
\label{12.997}
\end{equation}
The above lead to the conclusion that the contribution of the non-acoustical principal terms 
in \ref{12.987} to the right hand side of \ref{12.983} is bounded by:
\begin{equation}
\tau_1^{c_0}\left\{\frac{C\sqrt{\s^{(Y;0,n)}\cA(\tau_1)}}{\sqrt{2c_0}}
+\frac{C\tau_1^{1/2}\sqrt{\s^{(E;0,n)}\cA(\tau_1)}}{\sqrt{2c_0+1}}\right\}
\label{12.998}
\end{equation}

We now turn to the principal acoustical terms, the 1st terms in each of \ref{12.987}. From the 
definitions \ref{10.234}, \ref{10.236}, \ref{10.238} we have, to principal terms, 
\begin{equation}
\lambda\s^{(n)}\ctchi=\cth_n-(E^n f-E_N^n f_N), \ \ \ 
\lambdab\s^{(n)}\ctchib=\cthb_n-(E^n\fb-E_N^n\fb_N)
\label{12.999}
\end{equation}
In regard to the first, $\cth_n$ is estimated in $L^2({\cal K}^{\tau})$ in Proposition 10.3: 
\begin{eqnarray}
&&\|\cth_n\|_{L^2({\cal K}^{\tau})}\leq \tau^{c_0-1}
\left\{C\sqrt{\frac{\s^{(Y;0,n)}\cBb(\tau,\tau)}{2a_0+1}}\right.\nonumber\\
&&\hspace{40mm}\left.+C\tau^{1/2}\sqrt{\max\{\s^{(0,n)}\cB(\tau,\tau),\s^{(0,n)}\cBb(\tau,\tau)\}}\right\}\nonumber\\
&&\label{12.1000}
\end{eqnarray}
to principal terms. 
The leading contribution to $E^n f-E_N^n f_N$ comes from \ref{10.364}:
\begin{equation}
-\frac{1}{2}\beta_N^2(E^n\Lb H-E_N^n\Lb_N H_N)
\label{12.1001}
\end{equation}
the principal part of which is \ref{10.365}:
\begin{equation}
\beta_N^2 H^\prime (g^{-1})^{\mu\nu}\beta_\nu\Lb(E^n\beta_\mu-E_N^n\beta_{\mu,N})
\label{12.1002}
\end{equation}
This is expressed as a sum of $E$, $N$, and $\Nb$ components by \ref{10.369}, and the dominant 
contribution is that of the $\Nb$ component:
\begin{equation}
-(\eta^2\beta_N/2c)\beta_N^2 H^\prime\Nb^\mu\Lb(E^n\beta_\mu-E_N^n\beta_{\mu,N})
\label{12.1003}
\end{equation}
This is similar to \ref{12.991} and is bounded in $L^2({\cal K}^{\tau})$ by (see \ref{12.994}, 
\ref{12.996})
\begin{equation}
\tau^{c_0-1}\left(C\sqrt{\s^{(Y;0,n)}\cA(\tau)}+C\tau^{1/2}\sqrt{\s^{(E;0,n)}\cA(\tau)}\right)
\label{12.1004}
\end{equation}
Since along ${\cal K}$ $\lambda\sim\tau^2$ while $\rho\sim\tau$, we conclude from the above that 
the contribution of the 1st term in the first of \ref{12.987} to the 1st term on the right in the 
inequality \ref{12.983} is bounded by:
\begin{eqnarray}
&&\tau_1^{c_0-1}\left\{\frac{C\sqrt{\s^{(Y;0,n)}\cA(\tau_1)}}{\sqrt{2c_0-2}}
+\tau_1^{1/2}\frac{C\sqrt{\s^{(E;0,n)}\cA(\tau_1)}}{\sqrt{2c_0-1}}\right.\nonumber\\
&&\hspace{15mm}+\frac{C}{\sqrt{2c_0-2}}\sqrt{\frac{\s^{(Y;0,n)}\cBb(\tau_1,\tau_1)}{2a_0+1}}\nonumber\\
&&\hspace{15mm}\left.+\tau_1^{1/2}\frac{C}{\sqrt{2c_0-1}}
\sqrt{\max\{\s^{(0,n)}\cB(\tau_1,\tau_1),\s^{(0,n)}\cBb(\tau_1,\tau_1)\}}\right\}\nonumber\\
&&\label{12.1005}
\end{eqnarray}

In regard to the second of \ref{12.999}, $\cthb_n$ is estimated in $L^2({\cal K}^{\tau})$ in 
Proposition 10.9 :
\begin{eqnarray}
&&\|\cthb_n\|_{L^2({\cal K}^{\tau})}\leq \tau^{c_0-3}\left\{C\sqrt{\s^{(Y;0,n)}\cA(\tau)}
+C\tau^{1/2}\sqrt{\s^{(0,n)}\cA(\tau)}\right.\nonumber\\
&&\hspace{40mm}+C\sqrt{\frac{\s^{(Y;0,n)}\cBb(\tau,\tau)}{2a_0+1}}\nonumber\\
&&\hspace{40mm}\left.+C\tau^{1/2}\sqrt{\max\{\s^{(0,n)}\cB(\tau,\tau),\s^{(0,n)}\cBb(\tau,\tau)\}}
\right\}\nonumber\\
&&\hspace{20mm}+C\tau\|\Omega^{n+1}\chf\|_{L^2({\cal K}^{\tau})}
+C\tau\|\Omega^{n+1}\cv\|_{L^2({\cal K}^{\tau})}+C\tau^2\|\Omega^{n+1}\cga\|_{L^2({\cal K}^{\tau})}
\nonumber\\
&&\label{12.1006}
\end{eqnarray}
In regard to the first of the last three terms here, we recall \ref{12.960}. The left hand 
side of \ref{12.960} is, to principal terms, 
\begin{eqnarray*}
&&\sh^{-(n+1)/2}\left[\Omega^{n+1}\hat{f}-\frac{1}{2}(\Omega\hat{f})\Omega^n\log\sh\right]\\
&&-\sh_N^{-(n+1)/2}\left[\Omega^{n+1}\hat{f}_N-\frac{1}{2}(\Omega\hat{f}_N)\Omega^n\log\sh_N\right]
\end{eqnarray*}
which in view of the definition \ref{12.954} is, to principal terms, 
\begin{equation}
\sh^{-(n+1)/2}\Omega^{n+1}\chf-\frac{1}{2}(\Omega\hat{f})\s^{(n)}\check{\sh}
\label{12.1007}
\end{equation}
Replacing the left hand side of \ref{12.960} by \ref{12.1007} and recalling \ref{12.942}, \ref{12.944},  
we conclude that, to leading terms,
\begin{equation}
\sh^{-(n+1)/2}\Omega^{n+1}\chf=\s^{(n+1)}\chf+\frac{1}{2}\tau^{-2}(\Omega f)\s^{(n)}\check{\sh}
\label{12.1008}
\end{equation}
It follows that:
\begin{equation}
\|\Omega^{n+1}\chf\|_{L^2({\cal K}^{\tau})}\leq k\|\s^{(n+1)}\chf\|_{L^2({\cal K}^{\tau})}
+C\left(\int_0^{\tau}\tau^{\prime-4}\|\s^{(n)}\check{\sh}\|^2_{L^2(S_{\tau^\prime,\tau^\prime})}
d\tau^\prime\right)^{1/2}
\label{12.1009}
\end{equation}
where $k$ is a constant greater than 1, but which can be chosen as close to 1 as we wish by suitably 
restricting $\delta$. Hence by \ref{12.972} and the estimate \ref{12.923}:
\begin{eqnarray}
&&\|\Omega^{n+1}\chf\|_{L^2({\cal K}^{\tau})}\leq\tau^{c_0-3}\left\{C\s^{(n)}H(\tau)\right.\nonumber\\
&&\hspace{38mm}+\frac{C}{\sqrt{2c_0}}\left(\sqrt{\s^{(E;0,n)}\cA(\tau)}
+\frac{\sqrt{\s^{(Y;0,n)}\cA(\tau)}}{\sqrt{2c_0}}\right)\nonumber\\
&&\hspace{38mm}+\frac{C}{(a_0+2)}\left(\sqrt{\s^{(E;0,n)}\cBb(\tau,\tau)}
+\frac{\sqrt{\s^{(Y;0,n)}\cBb(\tau,\tau)}}{\sqrt{2c_0}}\right)\nonumber\\
&&\left.\hspace{45mm}+C\tau^{1/2}\sqrt{\max\{\s^{(0,n)}\cB(\tau,\tau),\s^{(0,n)}\cBb(\tau,\tau)\}}
\right\}\nonumber\\
&&\label{12.1010}
\end{eqnarray}
In regard to the last two terms in \ref{12.1006}, as in showing the equality, to principal terms, of 
the left hand side of \ref{12.960} to \ref{12.1007}, we similarly show, recalling the definitions 
\ref{12.932}, that, to principal terms, 
\begin{eqnarray}
&&\sh^{-(n+1)/2}\Omega^{n+1}\cv=\s^{(n+1)}\cv+\frac{1}{2}(\Omega v)\s^{(n)}\check{\sh}\nonumber\\
&&\sh^{-(n+1)/2}\Omega^{n+1}\cga=\s^{(n+1)}\cga+\frac{1}{2}(\Omega\gamma)\s^{(n)}\check{\sh}
\label{12.1011}
\end{eqnarray}
It follows that:
\begin{eqnarray}
&&\|\Omega^{n+1}\cv\|_{L^2({\cal K}^{\tau})}\leq k\|\s^{(n+1)}\cv\|_{L^2({\cal K}^{\tau})}
+C\|\s^{(n)}\check{\sh}\|_{L^2({\cal K}^{\tau})} \label{12.1012}\\
&&\|\Omega^{n+1}\cga\|_{L^2({\cal K}^{\tau})}\leq k\|\s^{(n+1)}\cga\|_{L^2({\cal K}^{\tau})}
+C\|\s^{(n)}\sh\|_{L^2({\cal K}^{\tau})} \label{12.1013}
\end{eqnarray}
By Lemma 12.12 and the definition \ref{12.971} we then obtain:
\begin{eqnarray}
&&\|\Omega^{n+1}\cv\|_{L^2({\cal K}^{\tau})}\leq\tau^{c_0-3}\left\{C\s^{(n)}H(\tau)\right.\nonumber\\
&&\hspace{45mm}+\frac{C}{\sqrt{2c_0}}\left(\sqrt{\s^{(E;0,n)}\cA(\tau)}
+\frac{\sqrt{\s^{(Y;0,n)}\cA(\tau)}}{\sqrt{2c_0}}\right)\nonumber\\
&&\hspace{38mm}+\frac{C}{(a_0+2)}\left(\sqrt{\s^{(E;0,n)}\cBb(\tau,\tau)}
+\frac{\sqrt{\s^{(Y;0,n)}\cBb(\tau,\tau)}}{\sqrt{2c_0}}\right)\nonumber\\
&&\hspace{45mm}\left.+C\tau^{1/2}\sqrt{\max\{\s^{(0,n)}\cB(\tau,\tau),\s^{(0,n)}\cBb(\tau,\tau)\}}
\right\}\nonumber\\
&&\label{12.1014}
\end{eqnarray}
\begin{eqnarray}
&&\|\Omega^{n+1}\cga\|_{L^2({\cal K}^{\tau})}\leq\tau^{c_0-4}\left\{C\s^{(n)}H(\tau)\right.\nonumber\\
&&\hspace{45mm}+\frac{C\tau}{\sqrt{2c_0}}\left(\sqrt{\s^{(E;0,n)}\cA(\tau)}
+\frac{\sqrt{\s^{(Y;0,n)}\cA(\tau)}}{\sqrt{2c_0}}\right)\nonumber\\
&&\hspace{38mm}+\frac{C\tau}{(a_0+2)}\left(\sqrt{\s^{(E;0,n)}\cBb(\tau,\tau)}
+\frac{\sqrt{\s^{(Y;0,n)}\cBb(\tau,\tau)}}{\sqrt{2c_0}}\right)\nonumber\\
&&\hspace{45mm}\left.+C\tau^{3/2}\sqrt{\max\{\s^{(0,n)}\cB(\tau,\tau),\s^{(0,n)}\cBb(\tau,\tau)\}}
\right\}\nonumber\\
&&\label{12.1015}
\end{eqnarray}

In regard again to the second of \ref{12.999} the leading contribution to $E^n\fb-E_N^n\fb_N$ comes 
from \ref{10.631}:
\begin{equation}
\lambdab\beta_{\Nb}\sbeta(E^{n+1}H-E_N^{n+1}H_N)
\label{12.1016}
\end{equation}
the principal part of which is \ref{10.632}:
\begin{equation}
-2\lambdab\beta_{\Nb}\sbeta H^\prime(g^{-1})^{\mu\nu}\beta_\nu E(E^n\beta_\mu-E_N^n\beta_{\mu,N})
\label{12.1017}
\end{equation}
This is expressed as a sum of $E$, $N$ and $\Nb$ components by \ref{10.369}, and the dominant contribution is that of the $\Nb$ component:
\begin{equation}
(\eta^2\beta_N/c)\lambdab\beta_{\Nb}\sbeta H^\prime\Nb^\mu E(E^n\beta_\mu-E_N^n\beta_{\mu,N})
\label{12.1018}
\end{equation}
To estimate the $L^2$ norm of $\Nb^\mu E(E^n\beta_\mu-E_N^n\beta_{\mu,N})$ on ${\cal K}^{\tau_1}$, we write $\Nb^\mu=\rhob^{-1}\Lb^\mu$ 
and express $\Lb^\mu E(E^n\beta_\mu-E_N^n\beta_{\mu,N})$ up to lower order terms as 
as $\s^{(E;0,n)}\cxi_{\Lb}$. Then since along ${\cal K}$ $\rhob\sim\tau^2$, we have, to leading terms, 
\begin{eqnarray}
&&\|\Nb^\mu E(E^n\beta_\mu-E_N^n\beta_{\mu,N})\|^2_{L^2({\cal K}^{\tau_1})}\leq 
C\int_0^{\tau_1}\tau^{-4}\|\s^{(E;0,n)}\cxi_{\Lb}\|^2_{L^2(S_{\tau,\tau})}
d\tau\nonumber\\
&&\hspace{15mm}=C\int_0^{\tau_1}\tau^{-4}\frac{\partial}{\partial\tau}\left(
\|\s^{(E;0,n)}\cxi_{\Lb}\|^2_{L^2({\cal K}^{\tau})}\right)d\tau\nonumber\\
&&\hspace{5mm}=C\left\{\tau_1^{-4}\|\s^{(E;0,n)}\cxi_{\Lb}\|^2_{L^2({\cal K}^{\tau_1})}
+4\int_0^{\tau_1}\tau^{-5}\|\s^{(E;0,n)}\cxi_{\Lb}\|^2_{L^2({\cal K}^{\tau})}d\tau\right\}\nonumber\\
&&\hspace{25mm}\leq C\tau_1^{2c_0-4}\left(\s^{(E;0,n)}\cA(\tau_1)\right)^2\left(1+\frac{2}{c_0-2}\right)
\label{12.1019}
\end{eqnarray}
(appealing to \ref{9.a39} and \ref{9.299})  
hence, replacing $\tau_1$ by $\tau\in(0,\tau_1]$, the $L^2$ norm of \ref{12.1018} on ${\cal K}^{\tau}$ 
is bounded by:
\begin{equation}
C\tau^{c_0-2}\sqrt{\s^{(E;0,n)}\cA(\tau)}
\label{12.1020}
\end{equation}
Substituting \ref{12.1010}, \ref{12.1014}, \ref{12.1015} in \ref{12.1006} and taking also into account 
the bound \ref{12.1020}, we conclude, in view of the fact that along ${\cal K}$ $\lambdab\sim\tau$ 
while $\rhob\sim\tau^2$, that the contribution of the 1st term in the second of \ref{12.987} 
to the 2nd term on the right in the inequality \ref{12.983} is bounded by:
\begin{eqnarray}
&&\tau_1^{c_0-1}\left\{\frac{C\tau_1\s^{(n)}H(\tau_1)}{\sqrt{2c_0}}+\frac{C\sqrt{\s^{(Y;0,n)}\cA(\tau_1)}}{\sqrt{2c_0-2}}
+\tau_1^{1/2}\frac{C\sqrt{\s^{(E;0,n)}\cA(\tau_1)}}{\sqrt{2c_0-1}}\right.\nonumber\\
&&\hspace{15mm}+\frac{C}{\sqrt{2c_0-2}}\sqrt{\frac{\s^{(Y;0,n)}\cBb(\tau_1,\tau_1)}{2a_0+1}}\nonumber\\
&&\hspace{15mm}\left.+\tau_1^{1/2}\frac{C_{0,n}}{\sqrt{2c_0-1}}
\sqrt{\max\{\s^{(0,n)}\cB(\tau_1,\tau_1),\s^{(0,n)}\cBb(\tau_1,\tau_1)\}}\right\}\nonumber\\
&&\label{12.1021}
\end{eqnarray}
Combining this with \ref{12.1005} and with \ref{12.998} we conclude through \ref{12.983} that:
\begin{eqnarray}
&&\|\s^{(n)}\check{\sh}\|_{L^2({\cal K}^{\tau_1})}\leq 
\tau_1^{c_0-1}\left\{\frac{C\tau_1\s^{(n)}H(\tau_1)}{\sqrt{2c_0}}\right.\nonumber\\
&&\hspace{37mm}+\frac{C\sqrt{\s^{(Y;0,n)}\cA(\tau_1)}}{\sqrt{2c_0-2}}
+\tau_1^{1/2}\frac{C\sqrt{\s^{(E;0,n)}\cA(\tau_1)}}{\sqrt{2c_0-1}}\nonumber\\
&&\hspace{37mm}+\frac{C}{\sqrt{2c_0-2}}\sqrt{\frac{\s^{(Y;0,n)}\cBb(\tau_1,\tau_1)}{2a_0+1}}\nonumber\\
&&\hspace{33mm}\left.+\tau_1^{1/2}\frac{C}{\sqrt{2c_0-1}}
\sqrt{\max\{\s^{(0,n)}\cB(\tau_1,\tau_1),\s^{(0,n)}\cBb(\tau_1,\tau_1)\}}\right\}\nonumber\\
&&\label{12.1022}
\end{eqnarray}
Replacing $\tau_1$ by $\tau\in(0,\tau_1]$, multiplying by $\tau^{-c_0+1}$ and taking the 
supremum over $\tau\in(0,\tau_1]$ then yields, in view of the definition \ref{12.971}, 
\begin{eqnarray}
&&\s^{(n)}H(\tau_1)\leq \frac{C\tau_1\s^{(n)}H(\tau_1)}{\sqrt{2c_0}}\nonumber\\
&&\hspace{18mm}+\frac{C\sqrt{\s^{(Y;0,n)}\cA(\tau_1)}}{\sqrt{2c_0-2}}
+\tau_1^{1/2}\frac{C\sqrt{\s^{(E;0,n)}\cA(\tau_1)}}{\sqrt{2c_0-1}}\nonumber\\
&&\hspace{18mm}+\frac{C}{\sqrt{2c_0-2}}\sqrt{\frac{\s^{(Y;0,n)}\cBb(\tau_1,\tau_1)}{2a_0+1}}\nonumber\\
&&\hspace{18mm}+\tau_1^{1/2}\frac{C}{\sqrt{2c_0-1}}
\sqrt{\max\{\s^{(0,n)}\cB(\tau_1,\tau_1),\s^{(0,n)}\cBb(\tau_1,\tau_1)\}}\nonumber\\
&&\label{12.1023}
\end{eqnarray}
In reference to the coefficient of $\s^{(n)}H(\tau_1)$ on the right, 
imposing the condition 
\begin{equation}
\frac{C\delta}{\sqrt{2c_0}}\leq\frac{1}{2}
\label{12.1024}
\end{equation}
then yields the following lemma. 

\vspace{2.5mm}

\noindent{\bf Lemma 12.13} \ \ The quantity $\s^{(n)}H(\tau_1)$ satisfies the bound:
\begin{eqnarray*}
&&\s^{(n)}H(\tau_1)\leq \frac{C\sqrt{\s^{(Y;0,n)}\cA(\tau_1)}}{\sqrt{2c_0-2}}
+\tau_1^{1/2}\frac{C\sqrt{\s^{(E;0,n)}\cA(\tau_1)}}{\sqrt{2c_0-1}}\\
&&\hspace{18mm}+\frac{C}{\sqrt{2c_0-2}}\sqrt{\frac{\s^{(Y;0,n)}\cBb(\tau_1,\tau_1)}{2a_0+1}}\\
&&\hspace{18mm}+\tau_1^{1/2}\frac{C}{\sqrt{2c_0-1}}
\sqrt{\max\{\s^{(0,n)}\cB(\tau_1,\tau_1),\s^{(0,n)}\cBb(\tau_1,\tau_1)\}}
\end{eqnarray*}

\vspace{2.5mm}

Substituting finally the bound for $\s^{(n)}H(\tau_1)$ of Lemma 12.13 in the estimates \ref{12.1010}, 
\ref{12.1014} and \ref{12.1015} yields the following proposition.

\vspace{2.5mm} 

\noindent{\bf Proposition 12.7} \ \ The transformation function differences satisfy for 
all $\tau_1\in(0,\delta]$ the following estimates to principal terms:
\begin{eqnarray*}
&&\sup_{\tau\in(0,\tau_1]}\left\{\tau^{-c_0+3}\|\Omega^{n+1}\chf\|_{L^2({\cal K}^{\tau})}\right\}, \ 
\sup_{\tau\in(0,\tau_1]}\left\{\tau^{-c_0+3}\|\Omega^{n+1}\cv\|_{L^2({\cal K}^{\tau})}\right\} \\
&&\hspace{30mm}\leq \frac{C\sqrt{\s^{(Y;0,n)}\cA(\tau_1)}}{\sqrt{2c_0-2}}
+\frac{C\sqrt{\s^{(E;0,n)}\cA(\tau_1)}}{\sqrt{2c_0-1}}\\
&&\hspace{30mm}+\frac{C\sqrt{\s^{(Y;0,n)}\cBb(\tau_1,\tau_1)}}{\sqrt{(2c_0-2)(2a_0+1)}}
+\frac{C\sqrt{\s^{(E;0,n)}\cBb(\tau_1,\tau_1)}}{a_0+2}\\
&&\hspace{30mm}+C\tau_1^{1/2}\sqrt{\max\{\s^{(0,n)}\cB(\tau_1,\tau_1),\s^{(0,n)}\cBb(\tau_1,\tau_1)\}}
\end{eqnarray*}
\begin{eqnarray*}
&&\sup_{\tau\in(0,\tau_1]}\left\{\tau^{-c_0+4}\|\Omega^{n+1}\cga\|_{L^2({\cal K}^{\tau})}\right\}
\leq \frac{C\sqrt{\s^{(Y;0,n)}\cA(\tau_1)}}{\sqrt{2c_0-2}}
+\tau_1^{1/2}\frac{C\sqrt{\s^{(E;0,n)}\cA(\tau_1)}}{\sqrt{2c_0-1}}\\
&&\hspace{50mm}+\frac{C}{\sqrt{2c_0-2}}\sqrt{\frac{\s^{(Y;0,n)}\cBb(\tau_1,\tau_1)}{2a_0+1}}\\
&&\hspace{38mm}+\tau_1^{1/2}\frac{C}{\sqrt{2c_0-1}}
\sqrt{\max\{\s^{(0,n)}\cB(\tau_1,\tau_1),\s^{(0,n)}\cBb(\tau_1,\tau_1)\}}
\end{eqnarray*}

\pagebreak

\chapter{The Top Order Energy Estimates}

\section{Estimates for $\s^{(V;m,n-m)}\check{b}$}

Let us recall form Chapter 9 the functions $\s^{(V;m,l)}\check{b}$ representing according to 
\ref{9.a33} the boundary values on ${\cal K}$ of $\s^{(V;m,l)}\cxi(B)$:
\begin{equation}
\s^{(V;m,l)}\cxi(B)=\s^{(V;m,l)}\check{b} \ \mbox{: on ${\cal K}$}
\label{13.1}
\end{equation}
In view of the definition \ref{9.a34}, let us add the integral
\begin{equation}
2C^\prime\int_{{\cal K}^{\ub_1}}\Omega^{1/2}\left(\s^{(V;m,l)}\check{b}\right)^2
\label{13.2}
\end{equation}
to both sides of the $(m,l)$ difference energy identity \ref{9.296}. Then this identity takes the form:
\begin{eqnarray}
&&\s^{(V;m,l)}\check{{\cal E}}^{\ub_1}(u_1)+\s^{(V;m,l)}\underline{\check{{\cal E}}}^{u_1}(\ub_1)
+\s^{(V;m,l)}\check{{\cal F}}^{\prime\ub_1}=
\s^{(V;m,l)}\check{{\cal G}}^{\ub_1,u_1}\nonumber\\
&&\hspace{70mm}+2C^\prime\int_{{\cal K}^{\ub_1}}\Omega^{1/2}\left(\s^{(V;m,l)}\check{b}\right)^2
\nonumber\\
&&\label{13.3}
\end{eqnarray}
Estimating the integral in question amounts to estimating 
$\|\s^{(V;m,l)}\check{b}\|^2_{{\cal K}^{\ub_1}}$. 
We shall be concerned with the energy identities of top order $m+l=n$. 

To derive an expression for the functions $\s^{(V;m,l)}\check{b}$, we shall first derive an expression 
for the boundary values of $\s^{(V;m,l)}\xi(B)$, where the $\s^{(V;m,l)}\xi$ are the 1-forms:
\begin{equation}
\s^{(V;m,l)}\xi=V^\mu d(E^l T^m\beta_\mu)
\label{13.4}
\end{equation}
The starting point for this derivation is the boundary condition \ref{6.164}:
\begin{equation}
\s^{(V)}\xi(B)=\s^{(V)}b \ \mbox{: on ${\cal K}$}
\label{13.5}
\end{equation}

The vectorfield $T$ being tangential to ${\cal K}$, the boundary condition \ref{13.5} implies 
the derived boundary conditions:
\begin{equation}
T^m(\s^{(V)}\xi(B))=T^m\s^{(V)}b \ \mbox{: on ${\cal K}$}
\label{13.6}
\end{equation}
The vectorfield $E$ being also tangential to ${\cal K}$, these in turn imply the derived boundary 
conditions:
\begin{equation}
E^l T^m(\s^{(V)}\xi(B))=E^l T^m\s^{(V)}b \ \mbox{: on ${\cal K}$}
\label{13.7}
\end{equation}

Recalling that $\s^{(V)}\xi$ is the 1-form
$\s^{(V)}\xi=V^\mu d\beta_\mu$
we have 
$$\s^{(V)}\xi(B)=V^\mu B\beta_\mu$$  
Hence the left hand side in \ref{13.6} is:
\begin{eqnarray}
&&T^m(V^\mu B\beta_\mu)=V^\mu T^m B\beta_\mu+\sum_{k=0}^{m-1}
\left(\begin{array}{c}m\\k\end{array}\right)(T^{m-k}V^\mu)T^k B\beta_\mu \nonumber\\
&&\hspace{20mm}=V^\mu BT^m\beta_\mu+\s^{(V)}R_{m,0}
\label{13.8}
\end{eqnarray}
where
\begin{eqnarray}
&&\s^{(V)}R_{m,0}=V^\mu[T^m,B]\beta_\mu
+\sum_{k=0}^{m-1}\left(\begin{array}{c}m\\k\end{array}\right)(T^{m-k}V^\mu)BT^k\beta_\mu\nonumber\\
&&\hspace{23mm}+\sum_{k=1}^{m-1}\left(\begin{array}{c}m\\k\end{array}\right)(T^{m-k}V^\mu) 
[T^k,B]\beta_\mu 
\label{13.9}
\end{eqnarray}
Applying then $E^l$ to \ref{13.8} we conclude that the left hand side in \ref{13.7} is:
\begin{eqnarray}
&&E^l T^m(V^\mu B\beta_\mu)=V^\mu E^l BT^m\beta_\mu +\sum_{i=0}^{l-1}
\left(\begin{array}{c}l\\i\end{array}\right)(E^{l-i}V^\mu)E^i BT^m\beta_\mu \nonumber\\
&&\hspace{50mm}+E^l\s^{(V)}R_{m,0}\nonumber\\
&&\hspace{25mm}=V^\mu BE^l T^m\beta_\mu+\s^{(V)}R_{m,l} 
\label{13.10}
\end{eqnarray}
where
\begin{eqnarray}
&&\s^{(V)}R_{m,l}=V^\mu[E^l,B]T^m\beta_\mu+\sum_{i=0}^{l-1}
\left(\begin{array}{c}l\\i\end{array}\right)(E^{l-i}V^\mu)BE^i T^m\beta_\mu \nonumber\\
&&\hspace{17mm}+\sum_{i=1}^{l-1}\left(\begin{array}{c}l\\i\end{array}\right)(E^{l-i}V^\mu)
[E^i,B]T^m\beta_\mu+E^l\s^{(V)}R_{m,0} 
\label{13.11}
\end{eqnarray}
Since according to \ref{13.4} we have 
$$V^\mu BE^l T^m\beta_\mu=\s^{(V;m,l)}\xi(B)$$
by \ref{13.10} the boundary conditions \ref{13.7} take the form:
\begin{equation}
\s^{(V;m,l)}\xi(B)=E^l T^m\s^{(V)}b-\s^{(V)}R_{m,l} \ \ \mbox{: on ${\cal K}$}
\label{13.12}
\end{equation}

In regard to the commutator terms in \ref{13.9} and \ref{13.11} we have the operator identities 
\begin{eqnarray}
&&[T^m,B]=\sum_{k=0}^{m-1}T^{m-1-k}[T,B]T^k \nonumber\\
&&[E^l,B]=\sum_{i=0}^{l-1}E^{l-1-i}[E,B]E^i \label{13.13}
\end{eqnarray}
We recall from Chapter 6 that the vectorfield $B$ is expressed through the decomposion \ref{6.156} 
by \ref{6.157} and \ref{6.158}. Defining along ${\cal K}$ the function
\begin{equation}
\check{\Lambda}=\rho\hat{\Lambda}
\label{13.14}
\end{equation}
\ref{6.158} takes the form:
\begin{equation}
B_{||}=\check{\Lambda}\left(\frac{1}{4}(\beta_N+r\beta_{\Nb})T-\lambda\sbeta E\right)
\label{13.15}
\end{equation}
Recalling that $\hat{\Lambda}$ is defined by \ref{6.159} with $q$ defined by \ref{6.151}, we have: 
\begin{equation}
\check{\Lambda}=-\frac{\eta^2\Lambda}{G\epsilon}
\label{13.16}
\end{equation}
Now, $\Lambda$ is defined by \ref{6.166}. In view of \ref{4.14}, \ref{4.15}, \ref{4.31} and 
Proposition 4.1, following the analysis leading from \ref{6.166} to \ref{6.173} we conclude that 
\begin{equation}
\Lambda=\nu^2\tilde{\Lambda}
\label{13.17}
\end{equation}
where $\tilde{\Lambda}$ is a known smooth function of $(\sigma,\nu)$. Substituting for $\nu$ from 
the 1st of \ref{4.37} and for $\ep$ in terms of $\epb$ from \ref{4.36} we then obtain:
\begin{equation}
\check{\Lambda}=\frac{\eta^6(\beta_N-r\beta_{\Nb})^2}{4c^2 G}\frac{\epb}{r}\tilde{\Lambda}
\label{13.18}
\end{equation}
Moreover, by \ref{4.40} and Proposition 4.2:
\begin{equation}
\frac{\epb}{r}=\frac{1}{j(\kappa,\epb)}
\label{13.19}
\end{equation}
where $\kappa$ stands for the quadruplet \ref{4.43}. Hence \ref{13.18} takes the form:
\begin{equation}
\check{\Lambda}=\frac{\eta^6(\beta_N-r\beta_{\Nb})^2}{4c^2 G}\frac{\tilde{\Lambda}}{j}
\label{13.20}
\end{equation}
In view of the fact that (see \ref{4.39}):
\begin{equation}
\left.j\right|_{\partial_-{\cal B}}=-\left.\frac{\eta^2\beta_N^3}{8c^2}H^\prime
\right|_{\partial_-{\cal B}}
\label{13.21}
\end{equation}
Comparing with \ref{4.112} we see that the fact that $l$ is strictly negative implies that 
$j$ is bounded away from zero on $\partial_-{\cal B}=\partial_-{\cal K}$. Therefore we can 
introduce the assumption:
\begin{equation}
|j|\geq C^{-1} \ \mbox{on ${\cal K}^{\delta}$}
\label{13.22}
\end{equation}
Summarizing, we have:
\begin{equation}
B=B_{\bot}+B_{||}, \  \ \ B_{\bot}=-\frac{1}{2}M, \  \ B_{||}=B^T T+\sB E
\label{13.23}
\end{equation}
where
\begin{equation}
B^T=\frac{1}{4}(\beta_N+r\beta_{\Nb})\check{\Lambda}, \ \ 
\sB=-\lambda\sbeta\check{\Lambda}
\label{13.24}
\end{equation}
are order 0 functions on ${\cal K}$. The commutation relations \ref{3.a14} then give:
$$[T,B_{\bot}]=\zeta E, \ \ \ [T,B_{||}]=(TB^T)T+(T\sB-(\chi+\chib)\sB)E$$
$$[E,B_{\bot}]=\frac{1}{2}(\chi+\chib)E, \ \ [E,B_{||}]=(EB^T)T+(E\sB+(\chi+\chib)B^T)E$$
hence:
\begin{eqnarray}
&&[T,B]=(TB^T)T+\left(T\sB+\zeta-(\chi+\chib)\sB\right)E \nonumber\\
&&[E,B]=(EB^T)T+\left(E\sB+\frac{1}{2}(\chi+\chib)(1+2B^T)\right)E
\label{13.25}
\end{eqnarray}
We see that these commutators are linear combinations of $T$ and $E$ with coefficients which are 
order 1 functions on ${\cal K}$. It then follows through \ref{13.13} in connection with \ref{13.9} 
and \ref{13.11} that $[T^m,B]\beta_\mu$ is of order $m$ and $[E^l,B]T^m\beta_\mu$ is of order 
$m+l$. Since $V^\mu$ is either $E^\mu$ or $Y^\mu$ and, by \ref{6.43}, \ref{6.44}, \ref{6.50}:
\begin{equation}
Y^\mu=N^\mu+r\Nb^\mu \ \mbox{: along ${\cal K}$}
\label{13.26}
\end{equation}
so the $V^\mu$ are order 0 functions along ${\cal K}$, it follows from \ref{13.9} that 
$\s^{(V)}R_{m,0}$ is of order $m$, the 2nd sum in \ref{13.9} being of order $m-1$, and then 
from \ref{13.10} that $\s^{(V)}R_{m,l}$ is of order $m+l$, the 2nd sum in \ref{13.11} being 
of order $m+l-1$. 

We must now derive an analogue of the boundary conditions \ref{13.12} for the $N$th approximants. 
We first derive such an analogue of the basic boundary condition \ref{13.5}. To do this we revisit 
the proof of Proposition 6.2. Reviewing this proof, in particular equations \ref{6.110} and \ref{6.a15},  
we remark that what is directly proved is that:
\begin{equation}
A_+^\mu V\beta_{+\mu}=A_-^\mu V\beta_{-\mu}
\label{13.27}
\end{equation}
where $A_{\pm}^\mu$ are the rectangular components of the vectorfields $A_{\pm}$ along ${\cal K}$ 
defined in the statement of Proposition 6.2. The equation \ref{2.88} of the wave system is then  
appealed to to express the left hand side of \ref{13.27} as $\s^{(V)}\xi_+(A_+)$. On the other hand, 
the right hand side is left as it is, and, after rescaling $A_{\pm}$ by the factor $\kappa G_+^{-1}$ 
to define $B_{\pm}$ as in \ref{6.155}, \ref{6.160}, omitting the subscript +, the right hand side 
becomes 
\begin{equation}
\s^{(V)}b=B_-^\mu V\beta_{-\mu}=B_-^\mu\left.V\beta^\prime_{\mu}\right|_{{\cal K}}
\label{13.28}
\end{equation}
as in \ref{6.256} (see \ref{6.250}, \ref{6.251}). The left hand side after the rescaling becomes 
$\s^{(V)}\xi(B)$, omitting the subscript +, and the boundary condition \ref{13.5} results. 
Now, the proof of Proposition 6.2 holds with $V$ being a vectorfield along ${\cal K}$ which is 
tangential to ${\cal K}$, thus a differential operator interior to ${\cal K}$, not only an actual 
variation field. Moreover, the proof is accomplished by applying $V$ to the condition \ref{6.107}. In 
view of \ref{6.106} and \ref{6.108}, the left hand side of \ref{6.107} is $-\delta^{-1} Q$ where 
$Q$ is the boundary quantity \ref{4.17}. In fact, what is actually shown in the proof of 
Proposition 6.2 is that:
\begin{equation}
VQ=-\delta \s^{(V)}\varepsilon 
\label{13.29}
\end{equation}
where:
\begin{equation}
\s^{(V)}\varepsilon=A_+^\mu V\beta_{+\mu}-A_-^\mu V\beta_{-\mu}
\label{13.30}
\end{equation}

For the actual solution $Q$ vanishes, hence also $\s^{(V)}\varepsilon$ vanishes for any $V$ as above. 
We now define the corresponding boundary quantities for the $N$th approximate solution. 
As in Chapter 9 we define, in $(\tau,\sigma,\vartheta)$ coordinates where ${\cal K}$ corresponds to 
$\sigma=0$, 
\begin{eqnarray}
&&\beta_{\mu,N+}(\tau,\vartheta)=\beta_{\mu,N}(\tau,0,\vartheta) \label{13.31}\\
&&\beta_{\mu,N-}(\tau,\vartheta)=\beta^{\prime}_\mu
(f_N(\tau,\vartheta),w_N(\tau,\vartheta),\psi_N(\tau,\vartheta)) \label{13.32}
\end{eqnarray}
so that, in accordance with \ref{9.55}, 
\begin{equation}
\triangle_N\beta_{\mu}(\tau,\vartheta)=\beta_{\mu,N+}(\tau,\vartheta)-\beta_{\mu,N-}(\tau,\vartheta)
\label{13.33}
\end{equation}
We also define:
\begin{eqnarray}
&&\sigma_{N\pm}=-(g^{-1})^{\mu\nu}\beta_{\mu,N\pm}\beta_{\nu,N\pm} \label{13.34}\\
&&G_{N\pm}=G(\sigma_{N\pm}) \label{13.35}
\end{eqnarray}
$G$ being a given smooth function of $\sigma$. We denote:
\begin{eqnarray}
&&\triangle_N\sigma(\tau,\vartheta)=\sigma_{N+}(\tau,\vartheta)-\sigma_{N-}(\tau,\vartheta) 
\label{13.36}\\
&&\triangle_N G(\tau,\vartheta)=G_{N+}(\tau,\vartheta)-G_{N-}(\tau,\vartheta) \label{13.37}
\end{eqnarray}
Also, 
\begin{equation}
\triangle_N(G\beta^\mu)(\tau,\vartheta)=(G_{N+}\beta^\mu_{N+})(\tau,\vartheta)
-(G_{N-}\beta^\mu_{N-})(\tau,\vartheta)
\label{13.38}
\end{equation}
We then define:
\begin{equation}
Q_N=\triangle_N(G\beta^\mu)\triangle_N\beta_\mu 
\label{13.39}
\end{equation}
all quantities appearing here being known smooth functions of $(\tau,\vartheta)$. 

Now, by \ref{9.61} and \ref{9.65} we have:
\begin{eqnarray}
&&\left.\frac{\partial^n\beta_{\mu,N+}}{\partial\tau^n}\right|_{\tau=0}=
\left.\frac{\partial^n\beta_{\mu+}}{\partial\tau^n}\right|_{\tau=0} \ \ \mbox{: for $n=0,...,N$}
\nonumber\\
&&\left.\frac{\partial^n\beta_{\mu,N-}}{\partial\tau^n}\right|_{\tau=0}=
\left.\frac{\partial^n\beta_{\mu-}}{\partial\tau^n}\right|_{\tau=0} \ \ \mbox{: for $n=0,...,N-1$}
\label{13.40}
\end{eqnarray}
It then follows that:
\begin{equation}
\left.\frac{\partial^n\triangle_N\beta_\mu}{\partial\tau^n}\right|_{\tau=0} 
=\left.\frac{\partial^n\triangle\beta_\mu}{\partial\tau^n}\right|_{\tau_0} \ \ 
\mbox{: for $n=0,...,N-1$}
\label{13.41}
\end{equation}
in accordance with \ref{9.68}. Moreover:
\begin{eqnarray}
&&\left.\frac{\partial^n\sigma_{N+}}{\partial\tau^n}\right|_{\tau=0}=
\left.\frac{\partial^n\sigma_+}{\partial\tau^n}\right|_{\tau=0} \ \ \mbox{: for $n=0,...,N$}\nonumber\\
&&\left.\frac{\partial^n\sigma_{N-}}{\partial\tau^n}\right|_{\tau=0}=
\left.\frac{\partial^n\sigma_-}{\partial\tau^n}\right|_{\tau=0} \ \ \mbox{: for $n=0,...,N-1$} 
\label{13.42}
\end{eqnarray}
which imply:
\begin{eqnarray}
&&\left.\frac{\partial^n G_{N+}}{\partial\tau^n}\right|_{\tau=0}=
\left.\frac{\partial^n G_+}{\partial\tau^n}\right|_{\tau=0} \ \ \mbox{: for $n=0,...,N$}\nonumber\\
&&\left.\frac{\partial^n G_{N-}}{\partial\tau^n}\right|_{\tau=0}=
\left.\frac{\partial^n G_-}{\partial\tau^n}\right|_{\tau=0} \ \ \mbox{: for $n=0,...,N-1$}
\label{13.43}
\end{eqnarray}
It then follows that:
\begin{equation}
\left.\frac{\partial^n\triangle_N(G\beta^\mu)}{\partial\tau^n}\right|_{\tau=0}=
\left.\frac{\partial^n\triangle(G\beta^\mu)}{\partial\tau^n}\right|_{\tau=0} \ \ 
\mbox{: for $n=0,...,N-1$}
\label{13.44}
\end{equation}
In view of the definitions \ref{4.17} and \ref{13.39}, equations \ref{13.41} and \ref{13.44} imply:
\begin{equation}
\left.\frac{\partial^n Q_N}{\partial\tau^n}\right|_{\tau=0}=
\left.\frac{\partial^n Q}{\partial\tau^n}\right|_{\tau=0}=0 \ \ \mbox{: for $n=0,...,N-1$}
\label{13.45}
\end{equation}
Moreover, for $n=N$ we have:
\begin{eqnarray}
&&\left.\frac{\partial^N Q_N}{\partial\tau^N}\right|_{\tau=0}
-\left.\frac{\partial^N Q}{\partial\tau^N}\right|_{\tau=0}=\nonumber\\
&&\hspace{10mm}\left(\left.\frac{\partial^N\triangle_N(G\beta^\mu)}{\partial\tau^N}\right|_{\tau=0}
-\left.\frac{\partial^N\triangle(G\beta^\mu)}{\partial\tau^N}\right|_{\tau=0}\right)
\left.\triangle\beta_\mu\right|_{\tau=0}\nonumber\\
&&\hspace{10mm}+\left.\triangle(G\beta^\mu)\right|_{\tau=0}\left(\left.
\frac{\partial^N\triangle_N\beta_\mu}{\partial\tau^N}\right|_{\tau=0}
-\left.\frac{\partial^N\triangle\beta_\mu}{\partial\tau^N}\right|_{\tau=0}\right)=0\nonumber\\
&&\label{13.46}
\end{eqnarray}
in view of the fact that
$$\left.\triangle(G\beta^\mu)\right|_{\tau=0}=\left.\triangle\beta_\mu\right|_{\tau=0}=0$$
Therefore also:
\begin{equation}
\left.\frac{\partial^N Q_N}{\partial\tau^N}\right|_{\tau=0}=0
\label{13.47}
\end{equation}
It then follows that:
\begin{equation}
Q_N=O(\tau^{N+1})
\label{13.48}
\end{equation}
in the sense of Definition 9.1. 

Given now a vectorfield $V_N$ along ${\cal K}$ which is tangential to ${\cal K}$ so that it is expressed 
as a linear combination of $\Omega=\partial/\partial\vartheta$ and $T=\partial/\partial\tau$ with 
coefficients which are smooth functions of $(\tau,\vartheta)$, we seek the analogue of \ref{13.29}, 
\ref{13.30} for the $N$th approximants. From \ref{13.38}:
\begin{equation}
V_N\triangle_N(G\beta^\mu)=V_N(G_{N+}\beta^\mu_{N+})-V_N(G_{N-}\beta^\mu_{N-})
\label{13.49}
\end{equation}
and in view of \ref{13.34}, \ref{13.35}:
\begin{equation}
V_N G_{N\pm}=G^\prime(\sigma_{N\pm})V_N\sigma_{N\pm}, \ \ \ 
V_N\sigma_{N\pm}=-2\beta^\mu_{N\pm}V_N\beta_{\mu,N\pm}
\label{13.50}
\end{equation}
Since according to the first of \ref{4.20}:
$$G^\prime=\frac{1}{2}GF,$$
we have:
\begin{equation}
G^\prime(\sigma_{N\pm})=\frac{1}{2}G_{N\pm}F_{N\pm}
\label{13.51}
\end{equation}
where:
\begin{equation}
F_{N\pm}=F(\sigma_{N\pm})
\label{13.52}
\end{equation}
We then obtain:
\begin{equation}
V_N(G_{N\pm}\beta^\mu_{N\pm})=-G_{N\pm}F_{N\pm}\beta^\mu_{N\pm}\beta^\nu_{N\pm}
V_N\beta_{\nu,N\pm}+G_{N\pm}V_N\beta^\mu_{N\pm}
\label{13.53}
\end{equation}
Let us define:
\begin{equation}
(h_{N\pm}^{-1})^{\mu\nu}=(g^{-1})^{\mu\nu}-F_{N\pm}\beta^\mu_{N\pm}\beta^\nu_{N\pm}
\label{13.54}
\end{equation}
Here the $+$ subscript is \ref{9.b7} along ${\cal K}$ while the $-$ subscript is the analogous 
definition along ${\cal K}$ corresponding to the prior solution. In terms of \ref{13.54}, \ref{13.53} 
takes the form:
\begin{equation}
V_N(G_{N\pm}\beta^\mu_{N\pm})=G_{N\pm}(h_{N\pm}^{-1})^{\mu\nu}V_N\beta_{\nu,N\pm}
\label{13.55}
\end{equation}
In analogy with \ref{6.106} let us define the function $\delta_N$ by:
\begin{equation}
\delta_N^2=\triangle_N\beta^\mu\triangle_N\beta_\mu 
\label{13.56}
\end{equation}
together with the condition that it is continuous and that its sign is such that the vector 
with rectangular components $\delta_N^{-1}\triangle_N\beta^\mu$ points along $\partial_-{\cal K}$ to 
the interior of ${\cal K}$. Then $\delta_N$ is smooth on ${\cal K}^\delta$ for suitably small $\delta$ 
and we have, by \ref{13.41} we have:
\begin{equation}
\left.\frac{\partial^n\delta_N}{\partial\tau^n}\right|_{\tau=0}=
\left.\frac{\partial^n\delta}{\partial\tau^n}\right|_{\tau=0} \ \ \mbox{: for $n=0,...,N-1$}
\label{13.57}
\end{equation}
We then define, in analogy with \ref{6.113},
\begin{equation}
K^\mu_{N\pm}=G_{N\pm}(h_{N\pm}^{-1})^{\mu\nu}\delta_N^{-1}\triangle_N\beta_\nu 
\label{13.58}
\end{equation}
In terms of this definition we have, by \ref{13.49} and \ref{13.55}:
\begin{equation}
\triangle_N\beta_\mu V_N\triangle_N(G\beta^\mu)=\delta_N\left(K^\mu_{N+}V_N\beta_{\mu,N+}
-K^\mu_{N-}V_N\beta_{\mu,N-}\right)
\label{13.59}
\end{equation}
Also:
\begin{equation}
\triangle_N(G\beta^\mu)V_N\triangle_N\beta_\mu=\triangle_N(G\beta^\mu)\left(V_N\beta_{\mu,N+}
-V_N\beta_{\mu,N-}\right)
\label{13.60}
\end{equation}
Adding then \ref{13.59} and \ref{13.60}, the sum of the left hand sides is $V_N Q_N$. Noting that 
in view of \ref{6.108} we may denote:
\begin{equation}
\triangle_N(G\beta^\mu)=-\triangle_N I^\mu
\label{13.61}
\end{equation}
the sum of the right hand sides is 
$-\delta_N\left(A^\mu_{N+}V_N\beta_{N+}-A^\mu_{N-}V_N\beta_{N-}\right)$
where, in analogy with the definition in the statement of Proposition 6.2:
\begin{equation}
A^\mu_{N\pm}=\frac{\triangle_N I^\mu}{\delta_N}-K^\mu_{N\pm}
\label{13.62}
\end{equation}
and we obtain, in analogy with \ref{13.29}, \ref{13.30}, 
\begin{equation}
V_N Q_N=-\delta_N\s^{(V_N)}\vep_N
\label{13.63}
\end{equation}
where:
\begin{equation}
\s^{(V_N)}\vep_N=A^\mu_{N+}V_N\beta_{\mu,N+}-A^\mu_{N-}V_N\beta_{\mu,N-}
\label{13.64}
\end{equation}
Now $\delta_N$ is a smooth function on ${\cal K}^\delta$ of fixed sign and by \ref{6.222} 
$|\delta_N|/\tau$ is a smooth function on ${\cal K}^\delta$ which is bounded from below by 
a positive constant. Taking then $V_N$ to be the $N$th approximant variation fields along ${\cal K}$:
\begin{equation}
V_N=\left\{\begin{array}{l} E_N=\sh_N^{-1/2}\Omega\\
Y_N=\rho_N^{-1}T\end{array}\right.
\label{13.65}
\end{equation}
(see \ref{6.271}) \ref{13.48} implies:
\begin{equation}
\s^{(V_N)}\vep_N=\left\{\begin{array}{lll} O(\tau^N) & : & \mbox{for $V_N=E_N$}\\
O(\tau^{N-2}) & : & \mbox{for $V_N=Y_N$} \end{array}\right.
\label{13.66}
\end{equation}
Equation \ref{13.64} with $\s^{(V_N)}\vep_N$ satisfying the estimates \ref{13.66} is the boundary 
condition satisfied by the $N$th approximate solution. To bring this condition to the form \ref{13.5} 
of the boundary condition for the actual solution, we must define $\kappa_N$,  
the $N$th approximant version of the positive quantity $\kappa$, so that we can perform a rescaling 
analogous to that of \ref{6.155}, \ref{6.160} to define, omitting the subscript +, the $N$th 
approximant versions of $B^\mu$, $B_-^\mu$:
\begin{equation}
B_N^\mu=\kappa_N G_N^{-1} A_N^\mu, \ \ \ B_{N-}^\mu=\kappa_N G_N^{-1} A_{N-}^\mu
\label{13.67}
\end{equation}
Moreover, defining in connection with the $N$th approximate solution in ${\cal N}$, the vectorfield 
$B_N$ along ${\cal N}$ so that:
\begin{equation}
B_N x_N^\mu=B_N^\mu 
\label{13.68}
\end{equation}
we must express, in regard to the first term on the right in \ref{13.64}, $B_N^\mu V_N\beta_{\mu,N}$ 
in terms of $V_N^\mu B_N\beta_{\mu,N}$, where:
\begin{equation}
V_N^\mu=V_N x_N^\mu 
\label{13.69}
\end{equation}
Note that $V_N^\mu$ need only be defined on ${\cal N}$. 

To define $\kappa_N$, we first define $w_N$, the $N$th approximant version of the positive function 
$w$ defined by \ref{6.147}. Since $\zeta_\mu$ corresponds for the $N$th approximants to 
$\triangle_N\beta_\mu/\delta_N$, setting, in analogy with \ref{4.11},
\begin{equation}
\mu_N=(h_N^{-1})^{\mu\nu}\triangle_N\beta_\mu\triangle_N\beta_\nu
\label{13.70}
\end{equation}
we define:
\begin{equation}
w_N=\frac{\mu_N}{\delta_N^2}
\label{13.71}
\end{equation}
Note that by \ref{9.b7} along ${\cal K}$ we have:
\begin{equation}
\delta_N^2=\mu_N+F_N\nu_N^2
\label{13.72}
\end{equation}
where, in analogy with the 1st of \ref{4.8}, 
\begin{equation}
\nu_N=\beta_N^\mu\triangle_N\beta_\mu
\label{13.73}
\end{equation}
The relation \ref{13.72} is analogous to the relation \ref{4.13}. The function $w_N/\tau$ is smooth 
and bounded from below by a positive constant on ${\cal K}^{\delta}$ (see \ref{6.220}). We then define 
$\kappa_N$ in analogy to \ref{6.149} to be the positive root of:
\begin{equation}
\kappa_N^2=\frac{a_N}{w_N}
\label{13.74}
\end{equation}
Since $a_N/\tau^3$ is smooth and bounded from below by a positive constant on ${\cal K}^{\delta}$ 
(see \ref{6.187}), it follows that $\kappa_N/\tau$ is smooth and bounded from below by a positive 
constant on ${\cal K}^{\delta}$. 

To define the vectorfield $B_N$ along ${\cal K}$ so that \ref{13.68} holds, we expand $B^\mu_N$ 
in the $N_N^\mu$, $\Nb_N^\mu$, $E_N^\mu$ frame:
\begin{equation}
B_N^\mu=B_N^N N_N^\mu+B_N^{\Nb}\Nb_N^\mu+B_N^E E_N^\mu 
\label{13.75}
\end{equation}
We then define the vectorfield $B_N$ according to:
\begin{equation}
B_N=B_N^N N_N+B_N^{\Nb}\Nb_N+B_N^E E_N
\label{13.76}
\end{equation}
$N_N$, $\Nb_N$ being given by \ref{9.b13}. As for the functions $V_N^\mu$ to be defined on ${\cal K}$, 
since $T=L_N+\Lb_N$, by \ref{9.b10}:
\begin{equation}
Tx_N^\mu=\rho_N N_N^\mu+\rhob_N\Nb_N^\mu+\vep_N^\mu+\vepb_N^\mu 
\label{13.77}
\end{equation}
hence by \ref{13.65} and \ref{9.102} along ${\cal K}$ we have: 
\begin{equation}
V_N^\mu=\left\{\begin{array}{l} E_N^\mu\\ Y_N^\mu\end{array}\right.
\label{13.78}
\end{equation}
where:
\begin{equation}
Y_N^\mu=N_N^\mu +r_N\Nb_N^\mu+\frac{1}{\rho_N}\left(\frac{\hat{\nu}_N}{c_N}\Nb_N^\mu+\vep_N^\mu 
+\vepb_N^\mu\right)
\label{13.79}
\end{equation}
Here by Proposition 9.1 and the estimate \ref{9.106} the error term is $O(\tau^{N-1})$. 

By \ref{9.a5} we have:
\begin{equation}
B_N^\mu V_N\beta_{\mu,N}=V_N^\mu B_N\beta_{\mu,N}+\omega_N(B_N,,V_N)
\label{13.80}
\end{equation}
Using Proposition 9.2 we derive the following estimates:
\begin{equation}
\omega_N(B_N,V_N)=\left\{\begin{array}{lll} O(\tau^N) & : & \mbox{if $V_N=E_N$}\\
O(\tau^{N-1}) & : & \mbox{if $V_N=Y_N$} \end{array} \right. 
\label{13.81}
\end{equation}
In view of \ref{13.80} and \ref{13.67}, the boundary condition \ref{13.64} takes the form:
\begin{equation}
V_N^\mu B_N\beta_{\mu,N}=\s^{(V_N)}b_N+\s^{(V_N)}\tilde{\vep}_N
\label{13.82}
\end{equation}
where:
\begin{equation}
\s^{(V_N)}b_N=B^\mu_{N-}V_N\beta_{\mu,N-} 
\label{13.83}
\end{equation}
is the $N$th approximant version of \ref{13.28} and $\s^{(V_N)}\tilde{\vep}_N$ is the error term:
\begin{equation}
\s^{(V_N)}\tilde{\vep}_N=\kappa_N G_N^{-1}\s^{(V_N)}\vep_N-\omega_N(B_N,V_N)
\label{13.84}
\end{equation}
The estimates \ref{13.66} and \ref{13.81} imply:
\begin{equation}
\s^{(V_N)}\tilde{\vep}_N=\left\{\begin{array}{lll} O(\tau^N) & : & \mbox{if $V_N=E_N$}\\
O(\tau^{N-1}) & : & \mbox{if $V_N=Y_N$} \end{array} \right. 
\label{13.85}
\end{equation}
In conclusion, \ref{13.82} is the analogue of the boundary condition \ref{13.5} for the $N$th 
approximants. 

Following an argument analogous to that leading from \ref{13.5} to \ref{13.12}, 
with $B_N$ and $V_N^\mu$ as defined above in the roles of $B$ and $V^\mu$, we deduce from 
\ref{13.82} the analogue of \ref{13.12} for the $N$th approximants, in the form:
\begin{equation}
V_N^\mu B_N(E_N^l T^m\beta_{\mu,N})=E_N^l T^m\s^{(V_N)}b_N-\s^{(V_N)}R_{m,l,N}
+E_N^l T^m\s^{(V_N)}\tilde{\vep}_N
\label{13.86}
\end{equation}
Here $\s^{(V_N)}R_{m,l,N}$ is the $N$th approximant version of the remainder term $\s^{(V)}R_{m,l}$ 
defined by \ref{13.9}, \ref{13.11}, the commutation relations \ref{9.b1} playing the role of the 
commutation relations \ref{3.a14}. Subtracting then \ref{13.86} from \ref{13.12} results in 
the boundary condition \ref{13.1} with: 
\begin{eqnarray}
&&\s^{(V;m,l)}\check{b}=E^l T^m\s^{(V)}b-E_N^l T^m\s^{(V_N)}b_N \nonumber\\
&&\hspace{18mm}+(V^\mu B-V_N^\mu B_N)E_N^l T^m\beta_{\mu,N} \nonumber\\
&&\hspace{18mm}-\s^{(V)}R_{m,l}+\s^{(V_N)}R_{m,l,N}-E_N^l T^m\s^{(V_N)}\tilde{\vep}_N \label{13.87}
\end{eqnarray}
Here, by \ref{13.85} the error term satisfies the estimates:
\begin{equation}
E_N^l T^m\s^{(V_N)}\tilde{\vep}_N=\left\{\begin{array}{lll} O(\tau^{N-m}) & : & \mbox{if $V_N=E_N$}\\
O(\tau^{N-1-m}) & : & \mbox{if $V_N=Y_N$} \end{array} \right. 
\label{13.88}
\end{equation}

We shall presently estimate the contribution to $\|\s^{(V;m,l)}\check{b}\|_{L^2({\cal K}^{\tau_1})}$, 
at the top order $m+l=n$, of the difference:
\begin{equation}
E^l T^m\s^{(V)}b-E_N^l T^m\s^{(V_N)}b_N 
\label{13.89}
\end{equation}
which contains the principal part of the right hand side of \ref{13.87}. More precisely, the 
contribution of the actual principal part of \ref{13.89}, which in the case $V=Y$ is, by \ref{6.258} 
and \ref{6.260},
\begin{eqnarray}
&&\sh^{-l/2}\hat{B}_-^N N^\mu\left\{\left.\frac{\partial\beta^\prime_\mu}{\partial t}\right|_{\cal K}
\tau^2\Omega^l T^{m+1}\chf+\left.\frac{\partial\beta^\prime_\mu}{\partial u^\prime}\right|_{\cal K}
\tau\Omega^l T^{m+1}\cv+\left.\frac{\partial\beta^\prime_\mu}{\partial\vartheta^\prime}\right|_
{\cal K}\tau^3\Omega^l T^{m+1}\cga\right\}\nonumber\\
&&+\sh^{-l/2}\hat{B}_-^{\Nb}\Nb^\mu\left\{\left.\frac{\partial\beta^\prime_\mu}{\partial t}\right|_{\cal K}
\tau^2\Omega^l T^{m+1}\chf+\left.\frac{\partial\beta^\prime_\mu}{\partial u^\prime}\right|_{\cal K}
\tau\Omega^l T^{m+1}\cv+\left.\frac{\partial\beta^\prime_\mu}{\partial\vartheta^\prime}\right|_
{\cal K}\tau^3\Omega^l T^{m+1}\cga\right\}\nonumber\\
&&+\sh^{-l/2}\hat{B}_-^E E^\mu\left\{\left.\frac{\partial\beta^\prime_\mu}{\partial t}\right|_{\cal K}
\tau^2\Omega^l T^{m+1}\chf+\left.\frac{\partial\beta^\prime_\mu}{\partial u^\prime}\right|_{\cal K}
\tau\Omega^l T^{m+1}\cv+\left.\frac{\partial\beta^\prime_\mu}{\partial\vartheta^\prime}\right|_
{\cal K}\tau^3\Omega^l T^{m+1}\cga\right\}\nonumber\\
&&\label{13.90} 
\end{eqnarray}
and in the case $V=E$ is, by \ref{6.271} and \ref{6.269}, 
\begin{eqnarray}
&&\sh^{-(l+1)/2}\rho\hat{B}_-^N N^\mu\left\{\left.\frac{\partial\beta^\prime_\mu}{\partial t}\right|_{\cal K}
\tau^2\Omega^{l+1} T^m\chf+\left.\frac{\partial\beta^\prime_\mu}{\partial u^\prime}\right|_{\cal K}
\tau\Omega^{l+1} T^m\cv+\left.\frac{\partial\beta^\prime_\mu}{\partial\vartheta^\prime}\right|_
{\cal K}\tau^3\Omega^{l+1} T^m\cga\right\}\nonumber\\
&&+\sh^{-(l+1)/2}\rho\hat{B}_-^{\Nb}\Nb^\mu\left\{\left.\frac{\partial\beta^\prime_\mu}{\partial t}\right|_{\cal K}
\tau^2\Omega^{l+1} T^m\chf+\left.\frac{\partial\beta^\prime_\mu}{\partial u^\prime}\right|_{\cal K}
\tau\Omega^{l+1} T^m\cv+\left.\frac{\partial\beta^\prime_\mu}{\partial\vartheta^\prime}\right|_
{\cal K}\tau^3\Omega^{l+1} T^m\cga\right\}\nonumber\\
&&+\sh^{-(l+1)/2}\rho\hat{B}_-^E E^\mu\left\{\left.\frac{\partial\beta^\prime_\mu}{\partial t}\right|_{\cal K}
\tau^2\Omega^{l+1} T^m\chf+\left.\frac{\partial\beta^\prime_\mu}{\partial u^\prime}\right|_{\cal K}
\tau\Omega^{l+1} T^m\cv+\left.\frac{\partial\beta^\prime_\mu}{\partial\vartheta^\prime}\right|_
{\cal K}\tau^3\Omega^{l+1} T^m\cga\right\}\nonumber\\
&&\label{13.91}
\end{eqnarray}

We consider first the case $m=0$, $l=n$. In this case \ref{13.91} takes the form:
\begin{eqnarray}
&&\sh^{-(n+1)/2}\rho\hat{B}_-^N N^\mu\left\{\left.\frac{\partial\beta^\prime_\mu}{\partial t}\right|_{\cal K}
\tau^2\Omega^{n+1}\chf+\left.\frac{\partial\beta^\prime_\mu}{\partial u^\prime}\right|_{\cal K}
\tau\Omega^{n+1}\cv+\left.\frac{\partial\beta^\prime_\mu}{\partial\vartheta^\prime}\right|_
{\cal K}\tau^3\Omega^{n+1}\cga\right\}\nonumber\\
&&+\sh^{-(n+1)/2}\rho\hat{B}_-^{\Nb}\Nb^\mu\left\{\left.\frac{\partial\beta^\prime_\mu}{\partial t}\right|_{\cal K}
\tau^2\Omega^{n+1}\chf+\left.\frac{\partial\beta^\prime_\mu}{\partial u^\prime}\right|_{\cal K}
\tau\Omega^{n+1}\cv+\left.\frac{\partial\beta^\prime_\mu}{\partial\vartheta^\prime}\right|_
{\cal K}\tau^3\Omega^{n+1}\cga\right\}\nonumber\\
&&+\sh^{-(n+1)/2}\rho\hat{B}_-^E E^\mu\left\{\left.\frac{\partial\beta^\prime_\mu}{\partial t}\right|_{\cal K}
\tau^2\Omega^{n+1}\chf+\left.\frac{\partial\beta^\prime_\mu}{\partial u^\prime}\right|_{\cal K}
\tau\Omega^{n+1}\cv+\left.\frac{\partial\beta^\prime_\mu}{\partial\vartheta^\prime}\right|_
{\cal K}\tau^3\Omega^{n+1}\cga\right\}\nonumber\\
&&\label{13.92}
\end{eqnarray}
Here the estimates of Proposition 12.7 apply. In view of the results \ref{6.252} - \ref{6.254} for 
the components of $\hat{B}_-$, as well as \ref{6.262}, \ref{6.264}, each of the three terms in 
parenthesis in the 1st of \ref{13.92} makes a leading contribution, only the second term in parenthesis 
in the 2nd of \ref{13.92} makes a leading contribution, while the contribution of the 3rd of \ref{13.92} 
is depressed relative to the preceding by a factor of $\tau$. It follows that \ref{13.92} is bounded in 
$L^2({\cal K}^{\tau})$ by:
\begin{equation}
C\tau^3\|\Omega^{n+1}\chf\|_{L^2({\cal K}^{\tau})}+C\tau^3\|\Omega^{n+1}\cv\|_{L^2({\cal K}^{\tau})}
+C\tau^4\|\Omega^{n+1}\cga\|_{L^2({\cal K}^{\tau})}
\label{13.93}
\end{equation}
which by Proposition 12.7 is bounded by:
\begin{eqnarray}
&&C\tau^{c_0}\left\{\frac{\sqrt{\s^{(Y;0,n)}\cA(\tau)}}{\sqrt{2c_0-2}}
+\frac{\sqrt{\s^{(E;0,n)}\cA(\tau)}}{\sqrt{2c_0-1}}\right.\nonumber\\
&&\hspace{12mm}\left.+\frac{\sqrt{\s^{(Y;0,n)}\cBb(\tau,\tau)}}{\sqrt{(2c_0-2)(2a_0+1)}}
+\frac{\sqrt{\s^{(E;0,n)}\cBb(\tau,\tau)}}{a_0+2}\right\}\nonumber\\
&&+C\tau^{c_0+\frac{1}{2}}
\sqrt{\max\{\s^{(0,n)}\cB(\tau,\tau),\s^{(0,n)}\cBb(\tau,\tau)\}}
\label{13.94}
\end{eqnarray}
In the case $m=0$, $l=n$, \ref{13.90} takes the form: 
\begin{eqnarray}
&&\sh^{-n/2}\hat{B}_-^N N^\mu\left\{\left.\frac{\partial\beta^\prime_\mu}{\partial t}\right|_{\cal K}
\tau^2\Omega^n T\chf+\left.\frac{\partial\beta^\prime_\mu}{\partial u^\prime}\right|_{\cal K}
\tau\Omega^n T\cv+\left.\frac{\partial\beta^\prime_\mu}{\partial\vartheta^\prime}\right|_
{\cal K}\tau^3\Omega^n T\cga\right\}\nonumber\\
&&+\sh^{-n/2}\hat{B}_-^{\Nb}\Nb^\mu\left\{\left.\frac{\partial\beta^\prime_\mu}{\partial t}\right|_{\cal K}
\tau^2\Omega^n T\chf+\left.\frac{\partial\beta^\prime_\mu}{\partial u^\prime}\right|_{\cal K}
\tau\Omega^n T\cv+\left.\frac{\partial\beta^\prime_\mu}{\partial\vartheta^\prime}\right|_
{\cal K}\tau^3\Omega^n T\cga\right\}\nonumber\\
&&+\sh^{-n/2}\hat{B}_-^E E^\mu\left\{\left.\frac{\partial\beta^\prime_\mu}{\partial t}\right|_{\cal K}
\tau^2\Omega^n T\chf+\left.\frac{\partial\beta^\prime_\mu}{\partial u^\prime}\right|_{\cal K}
\tau\Omega^n T\cv+\left.\frac{\partial\beta^\prime_\mu}{\partial\vartheta^\prime}\right|_
{\cal K}\tau^3\Omega^n T\cga\right\}\nonumber\\
&&\label{13.95} 
\end{eqnarray}
Here the estimates of Proposition 12.3 apply. In view of the results \ref{6.252} - \ref{6.254} for 
the components of $\hat{B}_-$, as well as \ref{6.262}, \ref{6.264}, each of the three terms in 
parenthesis in the 1st of \ref{13.95} makes a leading contribution, only the second term in parenthesis 
in the 2nd of \ref{13.95} makes a leading contribution, while the contribution of the 3rd of \ref{13.95} 
is depressed relative to the preceding by a factor of $\tau$. It follows that \ref{13.95} is bounded in 
$L^2({\cal K}^{\tau})$ by:
\begin{equation}
C\tau^2\|\Omega^n T\chf\|_{L^2({\cal K}^{\tau})}+C\tau^2\|\Omega^n T\cv\|_{L^2({\cal K}^{\tau})}
+C\tau^3\|\Omega^n T\cga\|_{L^2({\cal K}^{\tau})}
\label{13.96}
\end{equation}
which by Proposition 12.3 is bounded by: 
\begin{eqnarray}
&&C\tau^{c_0}\left\{\frac{\left(\sqrt{\s^{(Y;0,n)}\cA(\tau)}+\sqrt{\s^{(E;0,n)}\cA(\tau)}
\right)}{\sqrt{2c_0}}\right.\nonumber\\
&&\hspace{12mm}\left.+\frac{\left(\sqrt{\s^{(Y;0,n)}\cBb(\tau,\tau)}
+\sqrt{\s^{(E;0,n)}\cBb(\tau,\tau)}\right)}{a_0+2}\right\}\nonumber\\
&&+C\tau^{c_0+\frac{1}{2}}\sqrt{\max\{\s^{(0,n)}\cB(\tau,\tau),\s^{(0,n)}\cBb(\tau,\tau)\}}
\label{13.97}
\end{eqnarray}
the contribution of $\|\Omega^n T\cga\|_{L^2({\cal K}^{\tau})}$ being depressed relative to the others 
by a factor of $\tau$. 

We turn to the remaining cases $m=1,...,n$, $l=m-n$. Then in the case $V=Y$, \ref{13.90} holds and 
the estimates of Proposition 12.6 apply. In view again of the results \ref{6.252} - \ref{6.254} for 
the components of $\hat{B}_-$, as well as \ref{6.262}, \ref{6.264}, each of the three terms in 
parenthesis in the 1st of \ref{13.90} makes a leading contribution, only the second term in parenthesis 
in the 2nd of \ref{13.90} makes a leading contribution, while the contribution of the 3rd of \ref{13.90} 
is depressed relative to the preceding by a factor of $\tau$. It follows that \ref{13.90} is bounded in 
$L^2({\cal K}^{\tau})$ by:
\begin{equation}
C\tau^2\|\Omega^{n-m}T^{m+1}\chf\|_{L^2({\cal K}^{\tau})}
+C\tau^2\|\Omega^{n-m}T^{m+1}\cv\|_{L^2({\cal K}^{\tau})}
+C\tau^3\|\Omega^{n-m}T^{m+1}\cga\|_{L^2({\cal K}^{\tau})}
\label{13.98}
\end{equation}
which by Proposition 12.6 is, to leading terms, bounded by: 
\begin{eqnarray}
&&C\tau^{c_m}\left\{\frac{\sqrt{\s^{(Y;m,n-m)}\cA(\tau)}}{\sqrt{2c_m}}
+\tau^{1/2}\frac{\sqrt{\s^{(E;m,n-m)}\cA(\tau)}}{\sqrt{2c_m+1}}
+\frac{\sqrt{\s^{(Y;m,n-m)}\cBb(\tau,\tau)}}{(a_m+2)}\right\}\nonumber\\
&&+C\tau^{c_m+\frac{1}{2}}\sqrt{\max\{\s^{(m,n-m)}\cB(\tau,\tau),
\s^{(m,n-m)}\cBb(\tau,\tau)\}} \label{13.99}
\end{eqnarray}
the contribution of $\|\Omega^{n-m}T^{m+1}\cga\|_{L^2({\cal K}^{\tau})}$ being depressed relative 
to the others by a factor of $\tau$. 

In the case $V=E$, for $m=1,...,n$, $l=m-n$ \ref{13.91} holds and for $m=1$ the estimates of 
Proposition 12.3 apply while for $m=2,...,n$ those of Proposition 12.6 
with $m-1$ in the role of $m$ apply. 
In view again of the results \ref{6.252} - \ref{6.254} for 
the components of $\hat{B}_-$, as well as \ref{6.262}, \ref{6.264}, each of the three terms in 
parenthesis in the 1st of \ref{13.91} makes a leading contribution, only the second term in parenthesis 
in the 2nd of \ref{13.91} makes a leading contribution, while the contribution of the 3rd of \ref{13.91} 
is depressed relative to the preceding by a factor of $\tau$. It follows that \ref{13.91} is bounded in 
$L^2({\cal K}^{\tau})$ by:
\begin{equation}
C\tau^3\|\Omega^{n-m+1}T^{m}\chf\|_{L^2({\cal K}^{\tau})}
+C\tau^3\|\Omega^{n-m+1}T^{m}\cv\|_{L^2({\cal K}^{\tau})}
+C\tau^4\|\Omega^{n-m+1}T^{m}\cga\|_{L^2({\cal K}^{\tau})}
\label{13.100}
\end{equation}
which, by Proposition 12.3 for $m=1$ and Proposition 12.6 with $m$ replaced by $m-1$ for $m=2,..,n$, 
is, in view of the fact that $c_{m}\leq c_{m-1}$, bounded with an extra factor of $\tau$ at least.

\vspace{5mm}

\section{The Borderline Error Integrals Contributed by $\s^{(V;m,n-m)}\check{Q}_1$, 
$\s^{(V;m,n-m)}\check{Q}_2$}

Let us recall the error terms $\s^{(V;m,l)}\check{Q}_1$, $\s^{(V;m,l)}\check{Q}_2$, and 
$\s^{(V;m,l)}\check{Q}_3$ defined by \ref{9.281}, \ref{9.282}, and \ref{9.283} respectively. 
Of these $\s^{(V;m,l)}\check{Q}_1$, $\s^{(V;m,l)}\check{Q}_2$ correspond exactly to 
$\s^{(V)}Q_1$, $\s^{(V)}Q_2$ as defined by \ref{6.64}, \ref{6.65}, with the variation differences 
$\s^{(m,l)}\check{\dot{\phi}}_\mu$ (see \ref{9.236}) in the role of the variations $\dot{\phi}_\mu$, 
and the difference 1-form $\s^{(V;m,l)}\cxi$ (see \ref{9.275}) in the role of the 1-form $\s^{(V)}\xi$, 
thus with $\s^{(V;m,l)}\check{S}$ (see \ref{9.276}) in the role of $\s^{(V)}S$ (see \ref{6.102}). 
Consequently the estimates of Sections 7.2, 7.3 translate to estimates for the error integrals 
associated to $\s^{(V;m,l)}\check{Q}_1$, $\s^{(V;m,l)}\check{Q}_2$, in particular at the top order 
$m+l=n$. We shall presently estimate the error integrals of top order corresponding to the borderline 
error integrals encountered in Sections 7.2 and 7.3. 

In Section 7.2 we encountered two borderline error integrals, those in \ref{7.122} and \ref{7.124}. 
The corresponding top order integrals 
$$\int_0^{u_1}\s^{(V;m,n-m)}\cE^{\ub_1}(u)\frac{du}{u} \ \mbox{and} \ 
\int_0^{\ub_1}\s^{(V;m,n-m)}\cEb^{u_1}(\ub)\frac{d\ub}{\ub}$$
are bounded according to:
\begin{eqnarray}
&&\int_0^{u_1}\s^{(V;m,n-m)}\cE^{\ub_1}(u)\frac{du}{u} \nonumber\\
&&\hspace{15mm}\leq\s^{(V;m,n-m)}\cB(\ub_1,u_1)\ub_1^{2a_m}\int_0^{u_1} u^{2b_m-1}du \nonumber\\
&&\hspace{15mm}=\frac{\s^{(V;m,n-m)}\cB(\ub_1,u_1)}{2b_m}\ub_1^{2a_m}u_1^{2b_m} 
\label{13.101}
\end{eqnarray}
and:
\begin{eqnarray}
&&\int_0^{\ub_1}\s^{(V;m,n-m)}\cEb^{u_1}(\ub)\frac{d\ub}{\ub} \nonumber\\
&&\hspace{15mm}\leq\s^{(V;m,n-m)}\cBb(\ub_1,u_1)u_1^{2b_m}\int_0^{\ub_1} \ub^{2a_m-1}d\ub \nonumber\\
&&\hspace{15mm}=\frac{\s^{(V;m,n-m)}\cBb(\ub_1,u_1)}{2a_m}\ub_1^{2a_m}u_1^{2b_m} 
\label{13.102}
\end{eqnarray}
by \ref{9.297} and \ref{9.298} respectively. 

In Section 7.3 we encountered one borderline error integral, in the right hand side of \ref{7.186}. 
This is similar to \ref{7.124}, so the corresponding top order integral is bounded according to 
\ref{13.102}. 

\vspace{5mm}

\section{The Borderline Error Integrals Associated to $\cth_n$ and to $\cnu_{m-1,n-m+1} : m=1,...,n$}

We now come to the main point, the treatment of the borderline integrals contained in the 
error integral (see \ref{9.289}):
\begin{equation}
\s^{(V;m,l)}\check{{\cal G}}_3^{\ub_1,u_1}=\int_{{\cal R}_{\ub_1,u_1}}2a\Omega^{3/2}
\s^{(V;m,l)}\check{Q}_3
\label{13.103}
\end{equation}
associated to the error term $\s^{(V;m,l)}\check{Q}_3$, which is given by \ref{9.183}. The multiplier 
field $X$ having been fixed according to \ref{7.67}, this error integral is:
\begin{equation}
\s^{(V;m,l)}\check{{\cal G}}_3^{\ub_1,u_1}=-2\int_{{\cal R}_{\ub_1,u_1}}\Omega^{1/2}
\left(3\s^{(V;m,l)}\cxi_L+\s^{(V;m,l)}\cxi_{\Lb}\right)V^\mu\s^{(m,l)}\check{\tilde{\rho}}_\mu
\label{13.104}
\end{equation}
The borderline integrals contained in \ref{13.104} are contributed by the principal acoustical part 
of $\s^{(m,l)}\check{\tilde{\rho}}_\mu$. This principal acoustical part is obtained from the formulas 
\ref{8.158}, \ref{8.159}. Here we are at the top order $m+l=n$. The formula \ref{8.158} reads:
\begin{equation}
[\s^{(0,n)}\tilde{\rho}_\mu]_{P.A.}=\frac{1}{2}\rhob E^n\tchib L\beta_\mu
+\frac{1}{2}\rho E^n\tchi\Lb\beta_\mu
\label{13.105}
\end{equation}
In the case $m=1$ the formula \ref{8.159} reads:
\begin{eqnarray}
&&[\s^{(1,n-1)}\tilde{\rho}_\mu]_{P.A.}=\rho E^{n+1}\lambda \Lb\beta_\mu
+\rhob E^{n+1}\lambdab L\beta_\mu\nonumber\\
&&\hspace{26mm}+aE^n\tchi(\pi\Lb\beta_\mu+\rho E\beta_\mu)\nonumber\\
&&\hspace{26mm}+aE^n\tchib(\pi L\beta_\mu+\rhob E\beta_\mu) 
\label{13.106}
\end{eqnarray}
In the case $m\geq 2$ we must appeal to \ref{8.137}, \ref{8.138} to express 
$[E^{n-m+1}T^{m-1}\tchi]_{P.A.}$, $[E^{n-m+1}T^{m-1}\tchib]_{P.A.}$. In terms of the 
definitions \ref{12.64}, equations \ref{8.137}, \ref{8.138} take the form:
\begin{equation}
[T\chi]_{P.A.}=2[E\mu]_{P.A.}, \ \ \ [T\chib]_{P.A.}=2[E\mub]_{P.A.}
\label{13.107}
\end{equation}
Then by \ref{12.107}, \ref{12.108} and \ref{12.110}, \ref{12.111} we have:
\begin{equation}
[E^{n-m+1}T^{m-1}\tchi]_{P.A.}=[\mu_{m-2,n-m+2}]_{P.A.}, \ \ \ [E^{n-m+1}T^{m-1}\tchib]_{P.A.}=[\mub_{m-2,n-m+2}]_{P.A.}
\label{13.108}
\end{equation}
Therefore for $m\geq 2$ the formula \ref{8.159} implies:
\begin{eqnarray}
&&[\s^{(m,n-m)}\tilde{\rho}_\mu]_{P.A.}=\rho E^{n-m+2}T^{m-1}\lambda \Lb\beta_\mu 
+\rhob E^{n-m+2}T^{m-1}\lambdab L\beta_\mu \nonumber\\
&&\hspace{28mm}+a[\mu_{m-2,n-m+2}]_{P.A.}(\pi\Lb\beta_\mu+\rho E\beta_\mu)\nonumber\\
&&\hspace{28mm}+a[\mub_{m-2,n-m+2}]_{P.A.}(\pi L\beta_\mu+\rhob E\beta_\mu) 
\label{13.109}
\end{eqnarray}
Since similar formulas hold for the $N$th approximants it follows that:
\begin{equation}
[\s^{(0,n)}\check{\tilde{\rho}}_\mu]_{P.A.}=\frac{1}{2}\rhob \s^{(n)}\ctchib L\beta_\mu
+\frac{1}{2}\rho \s^{(n)}\ctchi\Lb\beta_\mu
\label{13.110}
\end{equation}
\begin{eqnarray}
&&[\s^{(1,n-1)}\check{\tilde{\rho}}_\mu]_{P.A.}=\rho \s^{(0,n+1)}\cla\Lb\beta_\mu
+\rhob \s^{(0,n+1)}\clab L\beta_\mu\nonumber\\
&&\hspace{26mm}+a\s^{(n)}\ctchi(\pi\Lb\beta_\mu+\rho E\beta_\mu)\nonumber\\
&&\hspace{26mm}+a\s^{(n)}\ctchib(\pi L\beta_\mu+\rhob E\beta_\mu) 
\label{13.111}
\end{eqnarray}
and, in view of \ref{12.117}, \ref{12.118}, \ref{12.119},  for $m\geq 2$:
\begin{eqnarray}
&&[\s^{(m,n-m)}\check{\tilde{\rho}}_\mu]_{P.A.}=\rho \s^{(m-1,n-m+2)}\cla \Lb\beta_\mu 
+\rhob \s^{(m-1,n-m+2)}\clab L\beta_\mu \nonumber\\
&&+a(\pi\Lb\beta_\mu+\rho E\beta_\mu)
\left(\frac{1}{2}(2\pi\lambda)^{m-1}\s^{(n)}\ctchi
+\sum_{i=1}^{m-1}(2\pi\lambda)^{m-1-i} \ \s^{(i-1,n-i+2)}\cla\right)\nonumber\\
&&+a(\pi L\beta_\mu+\rhob E\beta_\mu) 
\left(\frac{1}{2}(2\pi\lambdab)^{m-1}\s^{(n)}\ctchib
+\sum_{i=1}^{m-1}(2\pi\lambdab)^{m-1-i} \ \s^{(m-1-i,n-i+2)}\clab\right)\nonumber\\
&&\label{13.112}
\end{eqnarray}
Comparing now the coefficients of $\s^{(n)}\ctchi$, $\s^{(n)}\ctchib$ 
in the last two terms in \ref{13.111} with those in \ref{13.110}, taking account of the 
fact that by \ref{10.440} $\rho^{-1}L\beta_\mu$, $\Lb\beta_\mu$, $E\beta_\mu$ are bounded, the 
coefficients in \ref{13.111} are bounded with at least one extra factor of $\rho$ or $\rhob$ 
relative to the corresponding coefficients in \ref{13.110}. Since also $a_1\leq a_0$, $b_1\leq b_0$,  
it follows that the last two terms in \ref{13.111} will not make borderline contributions. Similarly, 
comparing the coefficients of $\s^{(i-1,n-i+2)}\cnu$, $\s^{(i-1, n-i+2)}\cnub$ :  $i=1,...,m-1$   
in the last two terms in \ref{13.112} with those of the first two terms in \ref{13.111} and 
\ref{13.112} with $i$ in the role of $m$, we see that the former are bounded with at least one 
extra factor of $\rho$ or $\rhob$ relative to the latter. Since also $a_m\leq a_i$, $b_m\leq b_i$,  
it follows that the last two terms in \ref{13.112} will not make borderline contributions. 
Thus, in regard to \ref{13.111} and \ref{13.112} we can restrict attention  to the first two terms 
on the right. 

Now, the principal acoustical terms in \ref{13.104} enter through the factor 
$$V^\mu\s^{(m,n-m)}\check{\tilde{\rho}}_\mu$$
Substituting \ref{13.110} - \ref{13.112}, the coefficient of $(1/2)\rho\s^{(n)}\ctchi$ from 
\ref{13.110}, of $\rho\s^{(0,n+1)}\cla$ from \ref{13.111}, and of $\rho\s^{(m-1,n-m+2)}\cla$ from 
\ref{13.112}, is:
\begin{equation}
V^\mu\Lb\beta_\mu=\left\{\begin{array}{lll} \rhob\ss_{\Nb} & : & \mbox{for $V=E$}\\
c\rhob\sss+\ogamma s_{\Nb\Lb} & : & \mbox{for $V=Y$} \end{array}\right.
\label{13.113}
\end{equation} 
Since $\ogamma\sim\sqrt{\rhob}$, we see that the dominant coefficient occurs in the case $V=Y$, 
through the term $\ogamma s_{\Nb\Lb}$. Only this case makes a borderline contribution. 
On the other hand, when substituting \ref{13.110} - \ref{13.112}, the coefficient of 
$(1/2)\rhob\s^{(n)}\ctchib$ from \ref{13.110}, of $\rhob\s^{(0,n+1)}\clab$ from \ref{13.111}, and of 
$\rhob\s^{(m-1,n-m+2)}\clab$ from \ref{13.112}, is:
\begin{equation}
V^\mu L\beta_\mu=\left\{\begin{array}{lll} \rho\ss_N & : & \mbox{for $V=E$}\\
s_{NL}+c\rho\ogamma\sss & : & \mbox{for $V=Y$} \end{array}\right.
\label{13.114}
\end{equation}
Here by \ref{10.440} the coefficient in both cases is bounded by a constant multiple of $\rho$ and  
in both cases we have a borderline contribution.

We now treat the borderline error integrals involving $\s^{(n)}\ctchi$ and $\s^{(m-1,n-m+2)}\cla 
\ : \ m=1,...,n$. As we have seen above, these arise in the case $V=Y$ only and are given by: 
\begin{equation}
-\int_{{\cal R}_{\ub_1,u_1}}\Omega^{1/2}\left(3\s^{(Y;0,n)}\cxi_L+\s^{(Y;0,n)}\cxi_{\Lb}\right)
(Y^\mu\Lb\beta_\mu)\rho\s^{(n)}\ctchi
\label{13.115}
\end{equation}
(case $m=0$, $l=n$ in \ref{13.104}), and, for $m=1,...,n$:
\begin{equation}
-2\int_{{\cal R}_{\ub_1,u_1}}\Omega^{1/2}\left(3\s^{(Y;m,n-m)}\cxi_L+\s^{(Y;m,n-m)}\cxi_{\Lb}\right)
(Y^\mu\Lb\beta_\mu)\rho\s^{(m-1,n-m+2)}\cla 
\label{13.116}
\end{equation}
(cases $m=1,...,n$, $l=n-m$ in \ref{13.104}). 

Consider first the integral \ref{13.115}. From the first of each of the definitions \ref{10.234}, 
\ref{10.236}, \ref{10.238} we can express:
\begin{equation}
\lambda E^n\tchi-\lambda_N E_N^n\tchi_N=\cth_n-E^n f+E_N^n f_N
\label{13.117}
\end{equation}
By \ref{12.3} with $l=n+1$, up to the 0th order term $(\lambda-\lambda_N)E_N^n\tchi_N$, 
the left hand side here is $\lambda\s^{(n)}\ctchi$. We can therefore replace \ref{13.115} by:
\begin{equation}
-\int_{{\cal R}_{\ub_1,u_1}}\Omega^{1/2}\left(3\s^{(Y;0,n)}\cxi_L+\s^{(Y;0,n)}\cxi_{\Lb}\right)
(Y^\mu\Lb\beta_\mu)\rho\lambda^{-1}(\cth_n-E^n f+E_N^n f_N)
\label{13.118}
\end{equation}
From \ref{13.113} we have:
\begin{equation}
|Y^\mu\Lb\beta_\mu|\rho\lambda^{-1}\leq C\frac{\rho}{\sqrt{\rhob}}
\label{13.119}
\end{equation}
The last being $\sim \ub/u$, \ref{13.118} is bounded in absolute value by the sum of 
a constant multiple of:
\begin{equation}
\int_{{\cal R}_{\ub_1,u_1}}(\ub/u)|\s^{(Y;0,n)}\cxi_{\Lb}||\cth_n|+
\int_{{\cal R}_{\ub_1,u_1}}(\ub/u)|\s^{(Y;0,n)}\cxi_{\Lb}||E^n f-E_N^n f_N|
\label{13.120}
\end{equation}
and a constant multiple of: 
\begin{equation}
\int_{{\cal R}_{\ub_1,u_1}}(\ub/u)|\s^{(Y;0,n)}\cxi_L||\cth_n|+
\int_{{\cal R}_{\ub_1,u_1}}(\ub/u)|\s^{(Y;0,n)}\cxi_L||E^n f-E_N^n f_N|
\label{13.121}
\end{equation}

Consider first \ref{13.120}. Since $\ub/u\leq 1$, the first integral is bounded by:
\begin{equation}
\int_0^{\ub_1}\left\{\int_{\Cb_{\ub}^{u_1}}|\cth_n||\s^{(Y;0,n)}\cxi_{\Lb}|\right\}d\ub
\leq \int_0^{\ub_1}\|\cth_n\|_{L^2(\Cb_{\ub}^{u_1})}\|\s^{(Y;0,n)}\cxi_{\Lb}\|_{L^2(\Cb_{\ub}^{u_1})}
d\ub
\label{13.122}
\end{equation}
Here we substitute the estimate for $\|\cth_n\|_{L^2(\Cb_{\ub}^{u_1})}$ of Proposition 10.3 ($l=n$). 
The borderline integral will be contributed by the leading term on the right:
\begin{equation}
C\sqrt{\frac{\s^{(Y;0,n)}\cBb(\ub_1,u_1)}{2a_0+1}}\cdot\ub^{a_0+2}u_1^{b_0-3}
\label{13.123}
\end{equation}
Substituting also:
\begin{equation}
\|\s^{(Y;0,n)}\cxi_{\Lb}\|_{L^2(\Cb_{\ub}^{u_1})}\leq C\sqrt{\s^{(Y;0,n)}\cEb^{u_1}(\ub)}
\leq C\ub^{a_0}u_1^{b_0}\sqrt{\s^{(Y;0,n)}\cBb(\ub_1,u_1)}
\label{13.124}
\end{equation}
(by \ref{9.291} and \ref{9.298}) we obtain that \ref{13.122} is to leading terms bounded by:
\begin{equation}
C\frac{\s^{(Y;0,n)}\cBb(\ub_1,u_1)}{\sqrt{2a_0+1}}\cdot u_1^{2b_0-3}\cdot\int_0^{\ub_1}\ub^{2a_0+2}d\ub
=C\frac{\s^{(Y;0,n)}\cBb(\ub_1,u_1)}{\sqrt{2a_0+1}(2a_0+3)}\cdot \ub_1^{2a_0}u_1^{2b_0}\cdot
\left(\frac{\ub_1}{u_1}\right)^3
\label{13.125}
\end{equation} 

Since again $\ub/u\leq 1$, the second of the integrals \ref{13.120} is bounded by:
\begin{eqnarray}
&&\int_0^{\ub_1}\left\{\int_{\Cb_{\ub}^{u_1}}|E^n f-E_N^n f_N||\s^{(Y;0,n)}\cxi_{\Lb}|\right\}d\ub
\nonumber\\
&&\hspace{30mm}\leq\int_0^{\ub_1}\|E^n f-E_N^n f_N\|_{L^2(\Cb_{\ub}^{u_1})}
\|\s^{(Y;0,n)}\cxi_{\Lb}\|_{L^2(\Cb_{\ub}^{u_1})}d\ub\nonumber\\
&&\label{13.126}
\end{eqnarray}
From Section 10.5, the leading contribution to $E^n f-E_N^n f_N$ comes from the principal part of 
\ref{10.364}, that is from \ref{10.365}, more precisely from the $\Nb$ component of this, which is 
(see \ref{10.369}):
$$-(\beta_N^3/2c)H^\prime\eta^2\Nb^\mu\Lb(E^n\beta_\mu-E_N^n\beta_{\mu,N})$$
This is bounded in $L^2(S_{\ub,u})$ by a constant multiple of:
\begin{equation}
\|\Nb^\mu\Lb(E^n\beta_\mu-E_N^n\beta_{\mu,N})\|_{L^2(S_{\ub,u})}\leq g_0(\ub,u)+g_1(\ub,u)
\label{13.127}
\end{equation}
according to \ref{10.379}. Here $g_0$ and $g_1$ are non-negative functions of $(\ub,u)$ satisfying 
\ref{10.380} and \ref{10.381}. The leading contribution is that of $g_0$. This contribution to  
$\|E^n f-E_N^n f_N\|^2_{L^2(\Cb_{\ub}^{u_1})}$ is then bounded by a constant multiple of:
\begin{eqnarray}
&&\int_{\ub}^{u_1}g_0^2(\ub,u)du=\int_{\ub}^{u_1}u^{-2}\frac{\partial}{\partial u}
\left(\int_{\ub}^u u^{\prime 2}g_0^2(\ub,u^\prime)du^\prime\right)du\nonumber\\
&&=u_1^{-2}\int_{\ub}^{u_1}u^2 g_0^2(\ub,u)du+2\int_{\ub}^{u_1}u^{-3}
\left(\int_{\ub}^u u^{\prime 2}g_0^2(\ub,u^\prime)du^\prime\right)du\nonumber\\
&&\leq Cu_1^{-2}\s^{(Y;0,n)}\cEb^{u_1}(\ub)+2C\int_{\ub}^{u_1}u^{-3}\s^{(Y;0,n)}\cEb^u(\ub)du
\nonumber\\
&&\label{13.128}
\end{eqnarray}
by \ref{10.380}. Substituting, from \ref{9.298}, 
$$\s^{(Y;0,n)}\cEb^u(\ub)\leq \ub^{2a_0}u^{2b_0}\s^{(Y;0,n)}\cBb(\ub,u_1) \ \forall u\in[\ub,u_1]$$
we conclude (in view of \ref{10.391}) that \ref{13.128} is bounded by:
$$C\ub^{2a_0}u_1^{2b_0-2}\s^{(Y;0,n)}\cBb(\ub,u_1)$$
consequently $\|E^n f-E_N^n f_N\|_{L^2(\Cb_{\ub}^{u_1})}$ is to leading terms bounded by:
\begin{equation}
C\ub^{a_0}u_1^{b_0-1}\sqrt{\s^{(Y;0,n)}\cBb(\ub_1,u_1)}
\label{13.129}
\end{equation}
Substituting this together with \ref{13.124} in \ref{13.126} then yields that the last is bounded by:
\begin{equation}
C\s^{(Y;0,n)}\cBb(\ub_1,u_1)\cdot u_1^{2b_0-1}\cdot\int_0^{\ub_1}\ub^{2a_0}d\ub
=C\frac{\s^{(Y;0,n)}\cBb(\ub_1,u_1)}{(2a_0+1)}\cdot \ub_1^{2a_0}u_1^{2b_0}\cdot\frac{\ub_1}{u_1}
\label{13.130}
\end{equation}

Consider now \ref{13.121}. Since $\ub/u\leq \ub^{1/2} u^{-1/2}$, the first integral is bounded by:
\begin{eqnarray}
&&\left(\int_{{\cal R}_{\ub_1,u_1}}\ub|\cth_n|^2\right)^{1/2}
\left(\int_{{\cal R}_{\ub_1,u_1}}u^{-1}|\s^{(Y;0,n)}\cxi_L|^2\right)^{1/2}\nonumber\\
&&=\left(\int_0^{\ub_1}\ub\|\cth_n\|^2_{L^2(\Cb_{\ub}^{u_1})}d\ub\right)^{1/2}
\left(\int_0^{u_1}u^{-1}\|\s^{(Y;0,n)}\cxi_L\|^2_{L^2(C_u^{\ub_1})}du\right)^{1/2}\nonumber\\
&&\label{13.131}
\end{eqnarray}
In regard to the integral in the first factor, we substitute the estimate for 
$\|\cth_n\|_{L^2(\Cb_{\ub}^{u_1})}$ of Proposition 10.3 ($l=n$). 
The leading contribution is \ref{13.123}, which gives a bound to leading terms for the integral 
in question by:
\begin{equation}
\frac{\s^{(Y;0,n)}\cBb(\ub_1,u_1)}{(2a_0+1)}\cdot u_1^{2b_0-6}\cdot\int_0^{\ub_1}\ub^{2a_0+5}d\ub
=\frac{\s^{(Y;0,n)}\cBb(\ub_1,u_1)}{(2a_0+1)(2a_0+6)}\cdot\ub_1^{2a_0}u_1^{2b_0}\cdot
\left(\frac{\ub_1}{u_1}\right)^6
\label{13.132}
\end{equation}
In regard to the integral in the second factor we substitute:
\begin{equation}
\|\s^{(Y;0,n)}\cxi_L\|_{L^2(C_u^{\ub_1})}\leq C\sqrt{\s^{(Y;0,n)}\cE^{\ub_1}(u)}
\leq C\ub_1^{a_0}u^{b_0}\sqrt{\s^{(Y;0,n)}\cB(\ub_1,u_1)}
\label{13.133}
\end{equation}
(by \ref{9.290} and \ref{9.297}) to obtain that the integral in question is bounded by:
\begin{equation}
\s^{(Y;0,n)}\cB(\ub_1,u_1)\cdot\ub_1^{2a_0}\cdot\int_0^{u_1}u^{2b_0-1}du
=\frac{\s^{(Y;0,n)}\cB(\ub_1,u_1)}{2b_0}\cdot \ub_1^{2a_0}u_1^{2b_0}
\label{13.134}
\end{equation}
Substituting this together with \ref{13.132} in \ref{13.131} then yields that the last is bounded by:
\begin{equation}
C\frac{\sqrt{\s^{(Y;0,n)}\cBb(\ub_1,u_1)\s^{(Y;0,n)}\cB(\ub_1,u_1)}}{\sqrt{(2a_0+1)(2a_0+6)2b_0}}
\cdot\ub_1^{2a_0}u_1^{2b_0}\cdot \left(\frac{\ub_1}{u_1}\right)^3
\label{13.135}
\end{equation}

Since again $\ub/u\leq \ub^{1/2} u^{-1/2}$, the second of the integrals \ref{13.121} is bounded by:
\begin{eqnarray}
&&\left(\int_{{\cal R}_{\ub_1,u_1}}\ub|E^n f-E_N^n f_N|^2\right)^{1/2}
\left(\int_{{\cal R}_{\ub_1,u_1}}u^{-1}|\s^{(Y;0,n)}\cxi_L|^2\right)^{1/2}\nonumber\\
&&=\left(\int_0^{\ub_1}\ub\|E^n f-E_N^n f_N\|^2_{L^2(\Cb_{\ub}^{u_1})}d\ub\right)^{1/2}
\left(\int_0^{u_1}u^{-1}\|\s^{(Y;0,n)}\cxi_L\|^2_{L^2(C_u^{\ub_1})}du\right)^{1/2}\nonumber\\
&&\label{13.136}
\end{eqnarray}
Here, by \ref{13.129} the integral in the first factor is to leading terms bounded by:
\begin{equation}
C\s^{(Y;0,n)}\cBb(\ub_1,u_1)\cdot u_1^{2b_0-2}\cdot\int_0^{\ub_1}\ub^{2a_0+1}d\ub
=C\frac{\s^{(Y;0,n)}\cBb(\ub_1,u_1)}{(2a_0+2)}\cdot\ub_1^{2a_0}u_1^{2b_0}\cdot
\left(\frac{\ub_1}{u_1}\right)^2
\label{13.137}
\end{equation}
while the integral in the second factor is bounded by \ref{13.134}, hence \ref{13.136} is to leading 
terms bounded by:
\begin{equation}
C\frac{\sqrt{\s^{(Y;0,n)}\cBb(\ub_1,u_1)\s^{(Y;0,n)}\cB(\ub_1,u_1)}}{\sqrt{(2a_0+2)2b_0}}\cdot 
\ub_1^{2a_0}u_1^{2b_0}\cdot\frac{\ub_1}{u_1} 
\label{13.138}
\end{equation}

We turn to the integral \ref{13.116}. From the first of each of the definitions \ref{10.235}, 
\ref{10.237}, \ref{10.239} we can express:
\begin{equation}
\lambda E^{n-m}T^{m-1}E^2\lambda-\lambda_N E_N^{n-m}T^m E_N^2\lambda_N=\cnu_{m-1,n-m+1}
+E^{n-m}T^{m-1}j-E_N^{n-m}T^{m-1}j_N
\label{13.139}
\end{equation}
Moreover we can replace $\cnu_{m-1,n-m+1}$ by $\cnu_{m-1,n-m+1}-\check{\tau}_{m-1,n-m+1}$, the 
difference being a lower order term. 
By \ref{12.4} with $m-1$ and $n-m+2$ in the role of $m$ and $l$ respectively, the left hand side 
here is $\lambda\s^{(m-1,n-m+2)}\cla$ up to lower order terms. We can therefore replace \ref{13.116} by:
\begin{eqnarray}
&&-2\int_{{\cal R}_{\ub_1,u_1}}\Omega^{1/2}\left(3\s^{(Y;m,n-m)}\cxi_L+\s^{(Y;m,n-m)}\cxi_{\Lb}\right)
(Y^\mu\Lb\beta_\mu)\rho\lambda^{-1}\cdot\nonumber\\
&&\hspace{15mm}\cdot(\cnu_{m-1,n-m+1}-\check{\tau}_{m-1,n-m+1}+E^{n-m}T^{m-1}j-E_N^{n-m}T^{m-1}j_N)
\nonumber\\
&&\label{13.140}
\end{eqnarray}
In view of \ref{13.119}, \ref{13.140} is bounded in absolute value by the sum of 
a constant multiple of:
\begin{eqnarray}
&&\int_{{\cal R}_{\ub_1,u_1}}(\ub/u)|\s^{(Y;m,n-m)}\cxi_{\Lb}||\cnu_{m-1,n-m+1}-\check{\tau}_{m-1,n-m+1}|\nonumber\\
&&\hspace{10mm}+\int_{{\cal R}_{\ub_1,u_1}}(\ub/u)|\s^{(Y;m,n-m)}\cxi_{\Lb}|
|E^{n-m}T^{m-1}j-E_N^{n-m}T^{m-1}j_N|\nonumber\\
&&\label{13.141}
\end{eqnarray}
and a constant multiple of: 
\begin{eqnarray}
&&\int_{{\cal R}_{\ub_1,u_1}}(\ub/u)|\s^{(Y;m,n-m)}\cxi_L||\cnu_{m-1,n-m+1}-\check{\tau}_{m-1,n-m+1}|\nonumber\\
&&\hspace{10mm}+\int_{{\cal R}_{\ub_1,u_1}}(\ub/u)|\s^{(Y;m,n-m)}\cxi_L|
|E^{n-m}T^{m-1}j-E_N^{n-m}T^{m-1}j_N|\nonumber\\
&&\label{13.142}
\end{eqnarray}

Consider first \ref{13.141}. Since $\ub/u\leq 1$, the first integral is bounded by:
\begin{eqnarray}
&&\int_0^{\ub_1}\left\{\int_{\Cb_{\ub}^{u_1}}|\cnu_{m-1,n-m+1}-\check{\tau}_{m-1,n-m+1}||\s^{(Y;m,n-m)}\cxi_{\Lb}|\right\}d\ub 
\nonumber\\
&&\hspace{10mm}\leq \int_0^{\ub_1}\|\cnu_{m-1,n-m+1}-\check{\tau}_{m-1,n-m+1}\|_{L^2(\Cb_{\ub}^{u_1})}
\|\s^{(Y;m,n-m)}\cxi_{\Lb}\|_{L^2(\Cb_{\ub}^{u_1})}d\ub\nonumber\\
&&\label{13.143}
\end{eqnarray}
Here we substitute the estimate for $\|\cnu_{m-1,n-m+1}-\check{\tau}_{m-1,n-m+1}\|_{L^2(\Cb_{\ub}^{u_1})}$ of Proposition 10.5 ($l=n-m$). 
The borderline integral will be contributed by the leading term on the right:
\begin{equation}
C\sqrt{\frac{\s^{(Y;m,n-m)}\cBb(\ub_1,u_1)}{2a_m+1}}\cdot\ub^{a_m+2}u_1^{b_m-3}
\label{13.144}
\end{equation}
Substituting also:
\begin{eqnarray}
&&\|\s^{(Y;m,n-m)}\cxi_{\Lb}\|_{L^2(\Cb_{\ub}^{u_1})}\leq C\sqrt{\s^{(Y;m,n-m)}\cEb^{u_1}(\ub)}
\nonumber\\
&&\hspace{34mm}\leq C\ub^{a_m}u_1^{b_m}\sqrt{\s^{(Y;m,n-m)}\cBb(\ub_1,u_1)}
\label{13.145}
\end{eqnarray}
(by \ref{9.291} and \ref{9.298}) we obtain that \ref{13.143} is to leading terms bounded by:
\begin{eqnarray}
&&C\frac{\s^{(Y;m,n-m)}\cBb(\ub_1,u_1)}{\sqrt{2a_m+1}}\cdot u_1^{2b_m-3}
\cdot\int_0^{\ub_1}\ub^{2a_m+2}d\ub\nonumber\\
&&\hspace{20mm}=C\frac{\s^{(Y;m,n-m)}\cBb(\ub_1,u_1)}{\sqrt{2a_m+1}(2a_m+3)}\cdot 
\ub_1^{2a_m}u_1^{2b_m}\cdot\left(\frac{\ub_1}{u_1}\right)^3
\label{13.146}
\end{eqnarray} 

Since again $\ub/u\leq 1$, the second of the integrals \ref{13.141} is bounded by:
\begin{eqnarray}
&&\int_0^{\ub_1}\left\{\int_{\Cb_{\ub}^{u_1}}|E^{n-m}T^{m-1}j-E_N^{n-m}T^{m-1}j_N||\s^{(Y;m,n-m)}\cxi_{\Lb}|\right\}d\ub
\nonumber\\
&&\leq\int_0^{\ub_1}\|E^{n-m}T^{m-1}j-E_N^{n-m}T^{m-1}j_N\|_{L^2(\Cb_{\ub}^{u_1})}
\|\s^{(Y;m,n-m)}\cxi_{\Lb}\|_{L^2(\Cb_{\ub}^{u_1})}d\ub\nonumber\\
&&\label{13.147}
\end{eqnarray}
From Section 10.6, the leading contribution to $E^{n-m}T^{m-1}j-E_N^{n-m}T^{m-1}j_N$ comes from the principal part of 
\ref{10.483}, that is from \ref{10.484}, more precisely from the $\Nb$ component of \ref{10.486}, which is: 
$$(\beta_N^3/4c)H^\prime\eta^2\Nb^\mu\Lb(E^{n-m}T^m\beta_\mu-E_N^{n-m}T^m\beta_{\mu,N})$$
This is bounded in $L^2(S_{\ub,u})$ by a constant multiple of:
\begin{equation}
\|\Nb^\mu\Lb(E^{n-m}T^m\beta_\mu-E_N^{n-m}T^m\beta_{\mu,N})\|_{L^2(S_{\ub,u})}\leq g_0(\ub,u)+g_1(\ub,u)
\label{13.148}
\end{equation}
according to \ref{10.511}. Here $g_0$ and $g_1$ are non-negative functions of $(\ub,u)$ satisfying 
\ref{10.512} and \ref{10.513}. The leading contribution is that of $g_0$. This contribution to  
$\|E^{n-m}T^{m-1}j-E_N^{n-m}T^{m-1}j_N\|^2_{L^2(\Cb_{\ub}^{u_1})}$ is then bounded by a constant multiple of:
\begin{eqnarray}
&&\int_{\ub}^{u_1}g_0^2(\ub,u)du=\int_{\ub}^{u_1}u^{-2}\frac{\partial}{\partial u}
\left(\int_{\ub}^u u^{\prime 2}g_0^2(\ub,u^\prime)du^\prime\right)du\nonumber\\
&&=u_1^{-2}\int_{\ub}^{u_1}u^2 g_0^2(\ub,u)du+2\int_{\ub}^{u_1}u^{-3}
\left(\int_{\ub}^u u^{\prime 2}g_0^2(\ub,u^\prime)du^\prime\right)du\nonumber\\
&&\leq Cu_1^{-2}\s^{(Y;m,n-m)}\cEb^{u_1}(\ub)+2C\int_{\ub}^{u_1}u^{-3}\s^{(Y;m,n-m)}\cEb^u(\ub)du
\nonumber\\
&&\label{13.149}
\end{eqnarray}
by \ref{10.512}. Substituting, from \ref{9.298}, 
$$\s^{(Y;m,n-m)}\cEb^u(\ub)\leq \ub^{2a_m}u^{2b_m}\s^{(Y;m,n-m)}\cBb(\ub,u_1) \ \forall u\in[\ub,u_1]$$
we conclude (in view of \ref{10.516}) that \ref{13.149} is bounded by:
$$C\ub^{2a_m}u_1^{2b_m-2}\s^{(Y;m,n-m)}\cBb(\ub,u_1)$$
consequently $\|E^{n-m}T^{m-1}j-E_N^{n-m}T^{m-1}j_N\|_{L^2(\Cb_{\ub}^{u_1})}$ is to leading terms bounded by:
\begin{equation}
C\ub^{a_m}u_1^{b_m-1}\sqrt{\s^{(Y;m,n-m)}\cBb(\ub_1,u_1)}
\label{13.150}
\end{equation}
Substituting this together with \ref{13.145} in \ref{13.147} then yields that the last is bounded by:
\begin{equation}
C\s^{(Y;m,n-m)}\cBb(\ub_1,u_1)\cdot u_1^{2b_m-1}\cdot\int_0^{\ub_1}\ub^{2a_m}d\ub
=C\frac{\s^{(Y;m,n-m)}\cBb(\ub_1,u_1)}{(2a_m+1)}\cdot \ub_1^{2a_m}u_1^{2b_m}\cdot\frac{\ub_1}{u_1}
\label{13.151}
\end{equation}

Consider now \ref{13.142}. Since $\ub/u\leq \ub^{1/2} u^{-1/2}$, the first integral is bounded by:
\begin{eqnarray}
&&\left(\int_{{\cal R}_{\ub_1,u_1}}\ub|\cnu_{m-1,n-m+1}-\check{\tau}_{m-1,n-m+1}|^2\right)^{1/2}
\left(\int_{{\cal R}_{\ub_1,u_1}}u^{-1}|\s^{(Y;m,n-m)}\cxi_L|^2\right)^{1/2}=\nonumber\\
&&\left(\int_0^{\ub_1}\ub\|\cnu_{m-1,n-m+1}-\check{\tau}_{m-1,n-m+1}\|^2_{L^2(\Cb_{\ub}^{u_1})}d\ub\right)^{1/2}
\left(\int_0^{u_1}u^{-1}\|\s^{(Y;m,n-m)}\cxi_L\|^2_{L^2(C_u^{\ub_1})}du\right)^{1/2}\nonumber\\
&&\label{13.152}
\end{eqnarray}
In regard to the integral in the first factor, we substitute the estimate for 
$\|\cnu_{m-1,n-m+1}-\check{\tau}_{m-1,n-m+1}\|_{L^2(\Cb_{\ub}^{u_1})}$ of Proposition 10.5 ($l=n-m$). 
The leading contribution is \ref{13.144}, which gives a bound to leading terms for the integral 
in question by:
\begin{equation}
\frac{\s^{(Y;m,n-m)}\cBb(\ub_1,u_1)}{(2a_m+1)}\cdot u_1^{2b_m-6}\cdot\int_0^{\ub_1}\ub^{2a_m+5}d\ub
=\frac{\s^{(Y;m,n-m)}\cBb(\ub_1,u_1)}{(2a_m+1)(2a_m+6)}\cdot\ub_1^{2a_m}u_1^{2b_m}\cdot
\left(\frac{\ub_1}{u_1}\right)^6
\label{13.153}
\end{equation}
In regard to the integral in the second factor we substitute:
\begin{equation}
\|\s^{(Y;m,n-m)}\cxi_L\|_{L^2(C_u^{\ub_1})}\leq C\sqrt{\s^{(Y;m,n-m)}\cE^{\ub_1}(u)}
\leq C\ub_1^{a_m}u^{b_m}\sqrt{\s^{(Y;m,n-m)}\cB(\ub_1,u_1)}
\label{13.154}
\end{equation}
(by \ref{9.290} and \ref{9.297}) to obtain that the integral in question is bounded by:
\begin{equation}
\s^{(Y;m,n-m)}\cB(\ub_1,u_1)\cdot\ub_1^{2a_m}\cdot\int_0^{u_1}u^{2b_m-1}du
=\frac{\s^{(Y;m,n-m)}\cB(\ub_1,u_1)}{2b_m}\cdot \ub_1^{2a_m}u_1^{2b_m}
\label{13.155}
\end{equation}
Substituting this together with \ref{13.153} in \ref{13.152} then yields that the last is bounded by:
\begin{equation}
C\frac{\sqrt{\s^{(Y;m,n-m)}\cBb(\ub_1,u_1)\s^{(Y;m,n-m)}\cB(\ub_1,u_1)}}{\sqrt{(2a_m+1)(2a_m+6)2b_m}}
\cdot\ub_1^{2a_m}u_1^{2b_m}\cdot \left(\frac{\ub_1}{u_1}\right)^3
\label{13.156}
\end{equation}

Since again $\ub/u\leq \ub^{1/2} u^{-1/2}$, the second of the integrals \ref{13.142} is bounded by:
\begin{eqnarray}
&&\left(\int_{{\cal R}_{\ub_1,u_1}}\ub|E^{n-m}T^{m-1}j-E_N^{n-m}T^{m-1}j_N|^2\right)^{1/2}\cdot
\nonumber\\
&&\hspace{40mm}\cdot\left(\int_{{\cal R}_{\ub_1,u_1}}u^{-1}|\s^{(Y;m,n-m)}\cxi_L|^2\right)^{1/2}\nonumber\\
&&=\left(\int_0^{\ub_1}\ub\|E^{n-m}T^{m-1}j-E_N^{n-m}T^{m-1}j_N\|^2_{L^2(\Cb_{\ub}^{u_1})}d\ub\right)^{1/2}\cdot\nonumber\\
&&\hspace{40mm}\cdot\left(\int_0^{u_1}u^{-1}\|\s^{(Y;m,n-m)}\cxi_L\|^2_{L^2(C_u^{\ub_1})}du\right)^{1/2}\nonumber\\
&&\label{13.157}
\end{eqnarray}
Here, by \ref{13.150} the integral in the first factor is to leading terms bounded by:
\begin{equation}
C\s^{(Y;m,n-m)}\cBb(\ub_1,u_1)\cdot u_1^{2b_m-2}\cdot\int_0^{\ub_1}\ub^{2m_0+1}d\ub
=C\frac{\s^{(Y;m,n-m)}\cBb(\ub_1,u_1)}{(2a_m+2)}\cdot\ub_1^{2a_m}u_1^{2b_m}\cdot
\left(\frac{\ub_1}{u_1}\right)^2
\label{13.158}
\end{equation}
while the integral in the second factor is bounded by \ref{13.155}, hence \ref{13.157} is to leading 
terms bounded by:
\begin{equation}
C\frac{\sqrt{\s^{(Y;m,n-m)}\cBb(\ub_1,u_1)\s^{(Y;m,n-m)}\cB(\ub_1,u_1)}}{\sqrt{(2a_m+2)2b_m}}\cdot 
\ub_1^{2a_m}u_1^{2b_m}\cdot\frac{\ub_1}{u_1} 
\label{13.159}
\end{equation}

\vspace{5mm}

\section{The Borderline Error Integrals Associated to $\cthb_n$ and to $\cnub_{m-1,n-m+1} : m=1,...,n$}

We now treat the borderline error integrals involving $\s^{(n)}\ctchib$ and 
$\s^{(m-1,n-m+2)}\clab \ : \ m=1,...,n$. As we have seen in the previous section, these arise in both 
cases $V=E, Y$ and are given by:
\begin{equation}
-\int_{{\cal R}_{\ub_1,u_1}}\Omega^{1/2}\left(3\s^{(V;0,n)}\cxi_L+\s^{(V;0,n)}\cxi_{\Lb}\right)
(V^\mu L\beta_\mu)\rhob\s^{(n)}\ctchib 
\label{13.160}
\end{equation}
(case $m=0$, $l=n$ in \ref{13.104}), and, for $m=1,...,n$:
\begin{equation}
-2\int_{{\cal R}_{\ub_1,u_1}}\Omega^{1/2}\left(3\s^{(V;m,n-m)}\cxi_L+\s^{(V;m,n-m)}\cxi_{\Lb}\right)
(V^\mu L\beta_\mu)\rhob\s^{(m-1,n-m+2)}\clab 
\label{13.161}
\end{equation}
(cases $m=1,...,n$, $l=n-m$ in \ref{13.104}). 

Consider first the integral \ref{13.160}. From the second of each of the definitions \ref{10.234}, 
\ref{10.236}, \ref{10.238} we can express:
\begin{equation}
\lambdab E^n\tchib-\lambdab_N E_N^n\tchib_N=\cthb_n-E^n\fb+E_N^n\fb_N
\label{13.162}
\end{equation}
By \ref{12.3} with $l=n+1$, up to the 0th order term $(\lambdab-\lambdab_N)E_N^n\tchib_N$, the left 
hand side here is $\lambdab_N\s^{(n)}\ctchib$. We can therefore replace \ref{13.160} by:
\begin{equation}
-\int_{{\cal R}_{\ub_1,u_1}}\Omega^{1/2}\left(3\s^{(V;0,n)}\cxi_L+\s^{(V;0,n)}\cxi_{\Lb}\right)
(V^\mu L\beta_\mu)\rhob\lambdab^{-1}(\cthb_n-E^n\fb+E_N^n\fb_N)
\label{13.163}
\end{equation}
From \ref{13.114} we have:
\begin{equation}
|V^\mu L\beta_\mu|\rhob\lambdab^{-1}\leq C\rhob 
\label{13.164}
\end{equation}
The last being $\sim u^2$, \ref{13.163} is bounded in absolute value by the sum of a constant multiple 
of:
\begin{equation}
\int_{{\cal R}_{\ub_1,u_1}}u^2|\s^{(V;0,n)}\cxi_L||\cthb_n|+
\int_{{\cal R}_{\ub_1,u_1}}u^2|\s^{(V;0,n)}\cxi_L||E^n\fb-E_N^n\fb_N|
\label{13.165}
\end{equation}
and a constant multiple of:
\begin{equation}
\int_{{\cal R}_{\ub_1,u_1}}u^2|\s^{(V;0,n)}\cxi_{\Lb}||\cthb_n|+
\int_{{\cal R}_{\ub_1,u_1}}u^2|\s^{(V;0,n)}\cxi_{\Lb}||E^n\fb-E_N^n\fb_N|
\label{13.166}
\end{equation}

Consider first \ref{13.165}. The first integral is:
\begin{equation}
\int_0^{u_1}u^2\left\{\int_{C_u^{\ub_1}}|\cthb_n||\s^{(V;0,n)}\cxi_L|\right\}du\leq 
\int_0^{u_1}u^2\|\cthb_n\|_{L^2(C_u^{\ub_1})}\|\s^{(V;0,n)}\cxi_L\|_{L^2(C_u^{\ub_1})}du
\label{13.167}
\end{equation}
Here, in regard to 
$\|\cthb_n\|_{L^2(C_u^{\ub_1})}$ we substitute the estimate of Proposition 10.7 ($l=n$) with $u$ in 
the role of $u_1$, noting that for $u\leq \ub_1$, $C_u^{\ub_1}$ stands of $C_u^u$, therefore 
the proposition reads, for $u\leq\ub_1$:
\begin{equation}
\|\cthb_n\|_{L^2(C_u^u)}\leq k\|\cthb_n\|_{L^2({\cal K}^u)}+
\Cb\sqrt{\max\{\s^{(0,n)}\cB(u,u),\s^{(0,n)}\cBb(u,u)\}}\cdot u^{a_0+b_0}
\label{13.168}
\end{equation}
and for $u>\ub_1$:
\begin{equation}
\|\cthb_n\|_{L^2(C_u^{\ub_1})}\leq k\|\cthb_n\|_{L^2({\cal K}^{\ub_1})}+
\Cb\sqrt{\max\{\s^{(0,n)}\cB(\ub_1,u),\s^{(0,n)}\cBb(\ub_1,u)\}}\cdot\ub_1^{a_0}u^{b_0}
\label{13.169}
\end{equation}
In regard to $\|\s^{(V;0,n)}\cxi_L\|_{L^2(C_u^{\ub_1})}$, noting again that for $u\leq\ub_1$, 
$C_u^{\ub_1}$ stands for $C_u^u$, by \ref{9.290} and \ref{9.297} we have, for $u\leq\ub_1$:
\begin{equation}
\|\s^{(V;0,n)}\cxi_L\|_{L^2(C_u^u)}\leq C\sqrt{\s^{(V;0,n)}\cE^{u}(u)}
\leq Cu^{a_0+b_0}\sqrt{\s^{(V;0,n)}\cB(u,u)}
\label{13.170}
\end{equation}
and for $u>\ub_1$:
\begin{equation}
\|\s^{(V;0,n)}\cxi_L\|_{L^2(C_u^{\ub_1})}\leq C\sqrt{\s^{(V;0,n)}\cE^{\ub_1}(u)}
\leq C\ub_1^{a_0}u^{b_0}\sqrt{\s^{(V;0,n)}\cB(\ub_1,u)}
\label{13.171}
\end{equation}
It follows that the contribution of the second term on the right in \ref{13.168}, \ref{13.169} to 
\ref{13.167} is bounded by:
\begin{eqnarray}
&&C\Cb\max\{\s^{(0,n)}\cB(\ub_1,u_1),\s^{(0,n)}\cBb(\ub_1,u_1)\}
\left\{\int_0^{\ub_1}u^{2a_0+2b_0+2}du+\ub_1^{2a_0}\int_{\ub_1}^{u_1}u^{2b_0+2}du\right\}\nonumber\\
&&\hspace{10mm}\leq C\Cb\frac{\max\{\s^{(0,n)}\cB(\ub_1,u_1),\s^{(0,n)}\cBb(\ub_1,u_1)\}}{2b_0+3}
\cdot\ub_1^{2a_0}u_1^{2b_0}\cdot u_1^3
\label{13.172}
\end{eqnarray}
which is depressed relative to a borderline contribution by the factor $u_1^3$. The borderline 
contribution is that of the first term on the right in \ref{13.168}, \ref{13.169}. By \ref{13.170}, 
\ref{13.171} this contribution to \ref{13.167} is bounded by:
\begin{equation}
C\sqrt{\s^{(V;0,n)}\cB(\ub_1,u_1)}\cdot k\left\{\int_0^{\ub_1}u^{a_0+b_0+2}
\|\cthb_n\|_{L^2({\cal K}^u)}du+\ub_1^{a_0}\|\cthb_n\|_{L^2({\cal K}^{\ub_1})}
\int_{\ub_1}^{u_1}u^{b_0+2}du\right\}
\label{13.173}
\end{equation}
Here we appeal to the estimate for $\|\cthb_n\|_{L^2({\cal K}^{\tau_1})}$ of Proposition 10.9 ($l=n$). 
The last three terms on the right:
$$\|\tau\Omega^{n+1}\chf\|_{L^2({\cal K}^{\tau_1})}, \|\tau\Omega^{n+1}\cv\|_{L^2({\cal K}^{\tau_1})}, 
\|\tau^2\Omega^{n+1}\cga\|_{L^2({\cal K}^{\tau_1})}$$
are by Proposition 12.7 bounded proportionally to $\tau_1^{a_0+b_0-2}$. So these terms are depressed 
relative to the leading terms by a factor of $\tau_1$. The leading terms in the estimate for 
$\|\cthb_n\|_{L^2({\cal K}^{\tau_1})}$ are then:
\begin{equation}
C\left\{\sqrt{\frac{\s^{(Y;0,n)}\cBb(\tau_1,\tau_1)}{2a_0+1}}+\sqrt{\s^{(Y;0,n)}\cA(\tau_1)}\right\}
\tau_1^{a_0+b_0-3}
\label{13.174}
\end{equation}
Placing then $u$ or $\ub_1$ in the role of $\tau_1$, in regard to the factor in parenthesis in 
\ref{13.173} we have:
\begin{eqnarray}
&&\int_0^{\ub_1}u^{2a_0+2b_0-1}du+\ub_1^{2a_0+b_0-3}\int_{\ub_1}^{u_1}u^{b_0+2}du=
\nonumber\\
&&\frac{\ub_1^{2a_0+2b_0}}{2a_0+2b_0}+\ub_1^{2a_0+b_0-3}
\frac{(u_1^{b_0+3}-\ub_1^{b_0+3})}{b_0+3}\leq \frac{\ub_1^{2a_0+b_0-3}u_1^{b_0+3}}{b_0+3}
\label{13.175}
\end{eqnarray}
(as $2a_0+2b_0>b_0+3$). We conclude that \ref{13.173} is bounded by:
\begin{equation}
C\frac{\sqrt{\s^{(V;0,n)}B(\ub_1,u_1)}}{b_0+3}\left\{\sqrt{\frac{\s^{(Y;0,n)}\cBb(\ub_1,\ub_1)}{2a_0+1}}+\sqrt{\s^{(Y;0,n)}\cA(\ub_1)}\right\}\cdot\ub_1^{2a_0}u_1^{2b_0}
\cdot\left(\frac{\ub_1}{u_1}\right)^{b_0-3}
\label{13.176}
\end{equation}
(for a different constant $C$). 

The second of the integrals \ref{13.165} is:
\begin{eqnarray}
&&\int_0^{u_1}u^2\left\{\int_{C_u^{\ub_1}}|E^n\fb-E_N^n\fb_N||\s^{(V;0,n)}\cxi_L|\right\}du\nonumber\\
&&\hspace{30mm}\leq 
\int_0^{u_1}u^2\|E^n\fb-E_N^n\fb_N\|_{L^2(C_u^{\ub_1})}\|\s^{(V;0,n)}\cxi_L\|_{L^2(C_u^{\ub_1})}du
\nonumber\\
&&\label{13.177}
\end{eqnarray}
From Section 10.7, the leading contribution to $E^n\fb-E_N^n\fb_N$ comes from the principal part of 
\ref{10.631}, that is from \ref{10.632}, more precisely from the $\Nb$ component of this which is 
(see \ref{10.369}):
$$\beta_N\beta_{\Nb}\sbeta H^\prime\eta^2\rho\Nb^\mu E(E^n\beta_\mu-E_N^n\beta_{\mu,N})$$
To estimate this in $L^2(C_u^{\ub_1})$, we write: 
\begin{equation}
N^\mu E(E^n\beta_\mu-E_N^n\beta_{\mu,N})+\ogamma\Nb^\mu E(E^n\beta_\mu-E_N^n\beta_{\mu,N})=
\s^{(Y;0,n)}\csxi
\label{13.178}
\end{equation}
In regard to the 1st term on the left, we can express 
$$\rho N^\mu E(E^n\beta_\mu-E_N^n\beta_{\mu,N}=L^\mu E(E^n\beta_\mu-E_N^n\beta_{\mu,N})$$
as 
$$E^\mu L(E^n\beta_\mu-E_N^n\beta_{\mu,N})=\s^{(E;0,n)}\cxi_L$$
up to lower order terms. Therefore we can express:
\begin{equation} 
\rho\Nb^\mu E(E^n\beta_\mu-E_N^n\beta_{\mu,N})
\label{13.179}
\end{equation}
up to lower order terms as:
\begin{equation}
\ogamma^{-1}\left(\rho\s^{(Y;0,n)}\csxi-\s^{(E;0,n)}\cxi_L\right)
\label{13.180}
\end{equation}
Since $\ogamma\sim u$, $\rho\sim\ub$ while $\sqrt{a}\sim u\sqrt{\ub}$, {13.180} is bounded pointwise by:
$$C(u^{-3/2}|\sqrt{a}\s^{(Y;0,n)}\csxi|+u^{-1}|\s^{(E;0,n)}\cxi_L|$$
hence in $L^2(C_u^{\ub_1})$ by:
\begin{equation}
Cu^{-3/2}\sqrt{\s^{(Y;0,n)}\cE^{\ub_1}(u)}+Cu^{-1}\sqrt{\s^{(E;0,n)}\cE^{\ub_1}(u)}
\label{13.181}
\end{equation}
by \ref{9.290}. It follows that $\|E^n\fb-E_N^n\fb_N\|_{L^2(C_u^{\ub_1})}$ is to leading terms bounded by:
\begin{equation}
Cu^{-3/2}\sqrt{\s^{(Y;0,n)}\cE^{\ub_1}(u)}\leq C\ub_1^{a_0}u^{b_0-\frac{3}{2}}
\sqrt{\s^{(Y;0,n)}\cB(\ub_1,u)}
\label{13.182}
\end{equation}
by \ref{9.297}. Together with \ref{13.171}, which by \ref{13.170} holds, a fortiori, also for 
$u\leq \ub_1$, we conclude that \ref{13.177} is to leading terms bounded by:
\begin{equation}
C\frac{\sqrt{\s^{(Y;0,n)}\cB(\ub_1,u_1)\s^{(V;0,n)}\cB(\ub_1,u_1)}}{2b_0+\frac{3}{2}}\cdot 
\ub_1^{2a_0}u_1^{2b_0}\cdot u_1^{3/2}
\label{13.183}
\end{equation}
This is depressed relative to a borderline contribution by the factor $u_1^{3/2}$. 

Consider now \ref{13.166}. Since $\ub\leq u$, the first integral is bounded by:
\begin{eqnarray}
&&\left(\int_{{\cal R}_{\ub_1,u_1}}u^5|\cthb_n|^2\right)^{1/2}
\left(\int_{{\cal R}_{\ub_1,u_1}}\ub^{-1}|\s^{(V;0,n)}\cxi_{\Lb}|^2\right)^{1/2}\nonumber\\
&&=\left(\int_0^{u_1}u^5\|\cthb_n\|^2_{L^2(C_u^{\ub_1})}du\right)^{1/2}
\left(\int_0^{\ub_1}\ub^{-1}\|\s^{(V;0,n)}\cxi_{\Lb}\|^2_{L^2(\Cb_{\ub}^{u_1})}d\ub\right)^{1/2} \nonumber\\
&&\label{13.184}
\end{eqnarray}
In regard to the integral in the first factor, we substitute the estimate \ref{13.168}, \ref{13.169}. 
As we have shown above, the leading contribution is that of the first term on the right, which to 
leading terms is bounded by $k$ times \ref{13.174} with $u$ in the role of $\tau_1$ for $u\leq\ub_1$, 
$\ub_1$ in the role of $\tau_1$ for $u>\ub_1$. We then conclude that the integral in question is 
to leading terms bounded by:
\begin{equation}
C\left({\frac{\s^{(Y;0,n)}\cBb(\ub_1,\ub_1)}{2a_0+1}}+\s^{(Y;0,n)}\cA(\ub_1)\right)
\left\{\int_0^{\ub_1}u^{2a_0+2b_0-1}du+\int_{\ub_1}^{u_1}\ub_1^{2a_0+2b_0-6}u^5 du\right\}
\label{13.185}
\end{equation}
Here, in regard to the 2nd integral we have, since in that integral $u>\ub_1$, 
$$\int_{\ub_1}^{u_1}\ub_1^{2a_0+2b_0-6}u^5 du\leq \ub_1^{2a_0}\int_{\ub_1}^{u_1}u^{2b_0-1}du
=\ub_1^{2a_0}\frac{(u_1^{2b_0}-\ub_1^{2b_0})}{2b_0}$$
It then follows that \ref{13.185} is bounded by:
\begin{equation}
\frac{C}{2b_0}\left(\frac{\s^{(Y;0,n)}\cBb(\ub_1,\ub_1)}{2a_0+1}+\s^{(Y;0,n)}\cA(\ub_1)\right)
\cdot\ub_1^{2a_0}u_1^{2b_0}
\label{13.186}
\end{equation}
In regard to the integral in the second factor in \ref{13.184} we substitute:
\begin{equation}
\|\s^{(V;0,n)}\cxi_{\Lb}\|_{L^2(\Cb_{\ub}^{u_1})}\leq C\sqrt{\s^{(V;0,n)}\cEb^{u_1}(\ub)}
\leq C\ub^{a_0}u_1^{b_0}\sqrt{\s^{(V;0,n)}\cBb(\ub_1,u_1)}
\label{13.187}
\end{equation}
(by \ref{9.291} and \ref{9.298}) to obtain that the integral in question is bounded by:
\begin{equation}
\s^{(V;0,n)}\cBb(\ub_1,u_1)\cdot u_1^{2b_0}\cdot\int_0^{\ub_1}\ub^{2a_0-1}du
=\frac{\s^{(V;0,n)}\cBb(\ub_1,u_1)}{2a_0}\cdot\ub_1^{2a_0}u_1^{2b_0}
\label{13.188}
\end{equation}
Substituting this together with \ref{13.186} in \ref{13.184} then yields that the last is bounded by:
\begin{equation}
\frac{C\ub_1^{2a_0}u_1^{2b_0}}{\sqrt{2a_0}\sqrt{2b_0}}\sqrt{\s^{(V;0,n)}\cBb(\ub_1,u_1)}
\sqrt{\frac{\s^{(Y;0,n)}\cBb(\ub_1,\ub_1)}{2a_0+1}+\s^{(Y;0,n)}\cA(\ub_1)}
\label{13.189}
\end{equation}

Since again $\ub\leq u$, the second of the integrals \ref{13.166} is bounded by:
\begin{eqnarray}
&&\left(\int_{{\cal R}_{\ub_1,u_1}}u^5 |E^n\fb-E_N^n\fb_N|^2\right)^{1/2}
\left(\int_{{\cal R}_{\ub_1,u_1}}\ub^{-1}|\s^{(V;0,n)}\cxi_{\Lb}|^2\right)^{1/2}\nonumber\\
&&=\left(\int_0^{u_1}u^5\|E^n\fb-E_N^n\fb_N\|^2_{L^2(C_u^{\ub_1})}du\right)^{1/2}
\left(\int_0^{\ub_1}\ub^{-1}\|\s^{(V;0,n)}\cxi_{\Lb}\|^2_{L^2(\Cb_{\ub}^{u_1})}d\ub\right)^{1/2}\nonumber\\
&&\label{13.190}
\end{eqnarray}
Here, by \ref{13.182} the integral in the first factor is to leading terms bounded by:
\begin{equation}
C\s^{(Y;0,n)}\cB(\ub_1,u_1)\cdot\ub_1^{2a_0}\cdot\int_0^{u_1}u^{2b_0+2}du
=C\frac{\s^{(Y;0,n)}\cB(\ub_1,u_1)}{(2b_0+3)}\cdot\ub_1^{2a_0}u_1^{2b_0}\cdot u_1^3
\label{13.191}
\end{equation}
while the integral in the second factor is bounded by \ref{13.188}, hence \ref{13.190} is to 
leading terms bounded by:
\begin{equation}
C\frac{\sqrt{\s^{(V;0,n)}\cBb(\ub_1,u_1)}}{\sqrt{2a_0}}\frac{\sqrt{\s^{(Y;0,n)}\cB(\ub_1,u_1)}}
{\sqrt{2b_0+3}}\cdot\ub_1^{2a_0}u_1^{2b_0}\cdot u_1^{3/2}
\label{13.192}
\end{equation}
which is depressed relative to a borderline contribution by the factor $u_1^{3/2}$. 

We turn to the integral \ref{13.161}. From the second of each of the definitions \ref{10.235}, 
\ref{10.237}, \ref{10.239} we can express:
\begin{equation}
\lambdab E^{n-m}T^{m-1}E^2\lambdab-\lambdab_N E_N^{n-m}T^m E_N^2\lambdab_N=\cnub_{m-1,n-m+1}
+E^{n-m}T^{m-1}\jb-E_N^{n-m}T^{m-1}\jb_N
\label{13.193}
\end{equation}
Moreover we can replace $\cnub_{m-1,n-m+1}$ by $\cnub_{m-1,n-m+1}-\check{\taub}_{m-1,n-m+1}$, the 
difference being a lower order term. 
By \ref{12.4} with $m-1$ and $n-m+2$ in the role of $m$ and $l$ respectively, the left hand side 
here is $\lambdab\s^{(m-1,n-m+2)}\clab$ up to lower order terms. We can therefore replace \ref{13.161} by:
\begin{eqnarray}
&&-2\int_{{\cal R}_{\ub_1,u_1}}\Omega^{1/2}\left(3\s^{(V;m,n-m)}\cxi_L+\s^{(V;m,n-m)}\cxi_{\Lb}\right)
(V^\mu L\beta_\mu)\rhob\lambdab^{-1}\cdot\nonumber\\
&&\hspace{15mm}\cdot(\cnub_{m-1,n-m+1}-\check{\taub}_{m-1,n-m+1}+E^{n-m}T^{m-1}\jb
-E_N^{n-m}T^{m-1}\jb_N)
\nonumber\\
&&\label{13.194}
\end{eqnarray}
In view of \ref{13.164}, \ref{13.194} is bounded in absolute value by the sum of 
a constant multiple of:
\begin{eqnarray}
&&\int_{{\cal R}_{\ub_1,u_1}}u^2|\s^{(V;m,n-m)}\cxi_L||\cnub_{m-1,n-m+1}-\check{\taub}_{m-1,n-m+1}|\nonumber\\
&&\hspace{10mm}+\int_{{\cal R}_{\ub_1,u_1}}u^2|\s^{(V;m,n-m)}\cxi_L|
|E^{n-m}T^{m-1}\jb-E_N^{n-m}T^{m-1}\jb_N|\nonumber\\
&&\label{13.195}
\end{eqnarray}
and a constant multiple of: 
\begin{eqnarray}
&&\int_{{\cal R}_{\ub_1,u_1}}u^2|\s^{(Y;m,n-m)}\cxi_{\Lb}||\cnub_{m-1,n-m+1}-\check{\taub}_{m-1,n-m+1}|\nonumber\\
&&\hspace{10mm}+\int_{{\cal R}_{\ub_1,u_1}}u^2|\s^{(V;m,n-m)}\cxi_{\Lb}|
|E^{n-m}T^{m-1}\jb-E_N^{n-m}T^{m-1}\jb_N|\nonumber\\
&&\label{13.196}
\end{eqnarray}

Consider first \ref{13.195}. The first integral is:
\begin{eqnarray}
&&\int_0^{u_1}u^2\left\{\int_{C_u^{\ub_1}}|\cnub_{m-1,n-m+1}-\check{\taub}_{m-1,n-m+1}||\s^{(V;m,n-m)}\cxi_L|\right\}du
\nonumber\\
&&\hspace{5mm}\leq 
\int_0^{u_1}u^2\|\cnub_{m-1,n-m+1}-\check{\taub}_{m-1,n-m+1}\|_{L^2(C_u^{\ub_1})}\|\s^{(V;m,n-m)}\cxi_L\|_{L^2(C_u^{\ub_1})}du
\nonumber\\
&&\label{13.197}
\end{eqnarray}
Here, in regard to 
$\|\cnub_{m-1,n-m+1}-\check{\taub}_{m-1,n-m+1}\|_{L^2(C_u^{\ub_1})}$ we substitute the estimate of Proposition 10.8 ($l=n-m$) 
with $u$ in the role of $u_1$, noting that for $u\leq \ub_1$, $C_u^{\ub_1}$ stands of $C_u^u$, therefore 
the proposition reads, for $u\leq\ub_1$:
\begin{eqnarray}
&&\|\cnub_{m-1,n-m+1}-\check{\taub}_{m-1,n-m+1}\|_{L^2(C_u^u)}\leq 
k\|\cnub_{m,n-m+1}-\check{\taub}_{m-1,n-m+1}\|_{L^2({\cal K}^u)}\nonumber\\
&&\hspace{30mm}+C_{m-1}u^2(\Pib u)^{m-1}\|\cthb_n\|_{L^2({\cal K}^u)}\nonumber\\
&&\hspace{10mm}+C^\prime_{m-1}u^2\sum_{i=0}^{m-2}(\Pib u)^{m-2-i}
\|\cnub_{i,n-i}-\check{\taub}_{i,n-i}\|_{L^2({\cal K}^u)}\nonumber\\
&&\hspace{10mm}+\Cb u^{a_m+b_m+\frac{1}{2}}\sqrt{\max\{\s^{[m,n-m]}\cB(u,u),\s^{[m,n-m]}\cBb(u,u)\}}
\nonumber\\
&&\label{13.198}
\end{eqnarray}
and for $u>\ub_1$:
\begin{eqnarray}
&&\|\cnub_{m-1,n-m+1}-\check{\taub}_{m-1,n-m+1}\|_{L^2(C_u^{\ub_1})}\leq 
k\|\cnub_{m-1,n-m+1}-\check{\taub}_{m-1,n-m+1}\|_{L^2({\cal K}^{\ub_1})}\nonumber\\
&&\hspace{30mm}+C_{m-1}\ub_1 u(\Pib\ub_1)^{m-1}\|\cthb_n\|_{L^2({\cal K}^{\ub_1})}\nonumber\\
&&\hspace{10mm}+C^\prime_{m-1}\ub_1 u\sum_{i=0}^{m-2}(\Pib\ub_1)^{m-2-i}
\|\cnub_{i,n-i}-\check{\taub}_{i,n-i}\|_{L^2({\cal K}^{\ub_1})}\nonumber\\
&&\hspace{10mm}+\Cb\ub_1^{a_m}u^{b_m+\frac{1}{2}}
\sqrt{\max\{\s^{[m,n-m]}\cB(\ub_1,u),\s^{[m,n-m]}\cBb(\ub_1,u)\}}
\nonumber\\
&&\label{13.199}
\end{eqnarray}
In regard to $\|\s^{(V;m,n-m)}\cxi_L\|_{L^2(C_u^{\ub_1})}$, noting again that for $u\leq\ub_1$, 
$C_u^{\ub_1}$ stands for $C_u^u$, by \ref{9.290} and \ref{9.297} we have, for $u\leq\ub_1$:
\begin{equation}
\|\s^{(V;m,n-m)}\cxi_L\|_{L^2(C_u^u)}\leq C\sqrt{\s^{(V;m,n-m)}\cE^{u}(u)}
\leq Cu^{a_m+b_m}\sqrt{\s^{(V;m,n-m)}\cB(u,u)}
\label{13.200}
\end{equation}
and for $u>\ub_1$:
\begin{equation}
\|\s^{(V;m,n-m)}\cxi_L\|_{L^2(C_u^{\ub_1})}\leq C\sqrt{\s^{(V;m,n-m)}\cE^{\ub_1}(u)}
\leq C\ub_1^{a_m}u^{b_m}\sqrt{\s^{(V;m,n-m)}\cB(\ub_1,u)}
\label{13.201}
\end{equation}
It follows that the contribution of the last term on the right in \ref{13.198}, \ref{13.199} to 
\ref{13.197} is bounded by:
\begin{eqnarray}
&&C\Cb\max\{\s^{[m,n-m]}\cB(\ub_1,u_1),\s^{[m,n-m]}\cBb(\ub_1,u_1)\}\cdot\nonumber\\
&&\hspace{40mm}\cdot\left\{\int_0^{\ub_1}u^{2a_m+2b_m+\frac{5}{2}}du+\ub_1^{2a_m}\int_{\ub_1}^{u_1}u^{2b_m+\frac{5}{2}}du\right\}\nonumber\\
&&\hspace{10mm}\leq C\Cb\frac{\max\{\s^{[m,n-m]}\cB(\ub_1,u_1),\s^{[m,n-m]}\cBb(\ub_1,u_1)\}}{2b_m+\frac{7}{2}}
\cdot\ub_1^{2a_m}u_1^{2b_m}\cdot u_1^{7/2}\nonumber\\
&&\label{13.202}
\end{eqnarray}
which is depressed relative to a borderline contribution by the factor $u_1^{7/2}$. The remaining terms 
on the right in \ref{13.198}, \ref{13.199} refer to the boundary ${\cal K}$. In view of the 
non-increasing with $m$ property of the exponents $a_m$, $b_m$ these terms other than the first term 
are depressed relative to the first term by at a factor of $u^2$ at least in the case of \ref{13.198},  
a factor of $\ub_1 u$ at least in the case of \ref{13.199}. Thus, the only borderline contribution 
comes from the first term on the right. Moreover, we have:
$$\|\cnub_{m-1,n-m+1}-\check{\taub}_{m-1,n-m+1}\|_{L^2({\cal K}^{\tau_1})}\leq 
\|\cnub_{m-1,n-m+1}\|_{L^2({\cal K}^{\tau_1})}+\|\check{\taub}_{m-1,n-m+1}\|_{L^2({\cal K}^{\tau_1})}$$
and the last can be neglected, the quantity $\check{\taub}_{m-1,n-m+1}$ being of order $n$ only 
(see \ref{10.655}). By \ref{13.200}, \ref{13.201}, the dominant contribution to \ref{13.197} is 
bounded by:
\begin{eqnarray}
&&C\sqrt{\s^{(V;m,n-m)}B(\ub_1,u_1)}\cdot k\left\{\int_0^{\ub_1}u^{a_m+b_m+2}
\|\cnub_{m-1,n-m+1}\|_{L^2({\cal K}^u)}du\right.\nonumber\\
&&\hspace{45mm}\left.+\ub_1^{a_m}\|\cnub_{m-1,n-m+1}\|_{L^2({\cal K}^{\ub_1})}
\int_{\ub_1}^{u_1}u^{b_m+2}du\right\}\nonumber\\
&&\label{13.203}
\end{eqnarray}
Here we appeal to the estimate for $\|\cnub_{m-1,n-m+1}\|_{L^2({\cal K}^{\tau_1})}$ of Proposition 10.9 ($l=n-m$). In regard to the last three terms on the right, 
$$\|\tau^3\Omega^{n-m+2}T^{m-1}\chf\|_{L^2({\cal K}^{\tau_1})}, 
\|\tau^2\Omega^{n-m+2}T^{m-1}\cv\|_{L^2({\cal K}^{\tau_1})}, 
\|\tau^4\Omega^{n-m+2}T^{m-1}\cga\|_{L^2({\cal K}^{\tau_1})},$$
are for $m=1$ bounded by Proposition 12.7 proportionally to $\tau_1^{a_0+b_0}$, $\tau_1^{a_0+b_0-1}$, 
$\tau_1^{a_0+b_0}$ respectively, for $m=2$ by Proposition 12.3 proportionally to $\tau_1^{a_0+b_0+1}$, 
$\tau_1^{a_0+b_0}$, $\tau_1^{a_0+b_0+2}$ respectively, and for $m\geq 3$ by Proposition 12.6 
proportionally to $\tau_1^{a_{m-2}+b_{m-2}+1}$, $\tau_1^{a_{m-2}+b_{m-2}}$, 
$\tau_1^{a_{m-2}+b_{m-2}+2}$ respectively. In view of the fact that the exponents $a_m$, $b_m$ are 
non-increasing in $m$, these terms are depressed relative to the leading terms 
by a factor of $\tau_1^2$. Also, the first term on the right, which involves 
$\cnub_{m-2-i,n-m+2+i}$ for $i=0,...,m-2$ and $\cthb_n$, is depressed relative to the leading terms 
by a factor of $\tau_1$. The leading terms in the estimate for 
$\|\cnub_{m-1,n-m+1}\|_{L^2({\cal K}^{\tau_1})}$ are then:
\begin{equation}
C\left\{\sqrt{\frac{\s^{(Y;m,n-m)}\cBb(\tau_1,\tau_1)}{2a_m+1}}+\sqrt{\s^{(Y;m,n-m)}\cA(\tau_1)}\right\}
\tau_1^{a_m+b_m-3}
\label{13.204}
\end{equation}
Placing then $u$ or $\ub_1$ in the role of $\tau_1$, in regard to the factor in parenthesis in 
\ref{13.203} we have:
\begin{eqnarray}
&&\int_0^{\ub_1}u^{2a_m+2b_m-1}du+\ub_1^{2a_m+b_m-3}\int_{\ub_1}^{u_1}u^{b_m+2}=
\nonumber\\
&&\frac{\ub_1^{2a_m+2b_m}}{2a_m+2b_m}+\ub_1^{2a_m+b_m-3}
\frac{(u_1^{b_m+3}-\ub_1^{b_m+3})}{b_m+3}\leq \frac{\ub_1^{2a_m+b_m-3}u_1^{b_m+3}}{b_m+3}\nonumber\\
&&\label{13.205}
\end{eqnarray}
(as $2a_m+2b_m>b_m+3$). We conclude that \ref{13.203} is bounded by:
\begin{eqnarray}
&&C\frac{\sqrt{\s^{(V;m,n-m)}B(\ub_1,u_1)}}{b_m+3}\left\{\sqrt{\frac{\s^{(Y;m,n-m)}\cBb(\ub_1,\ub_1)}{2a_m+1}}+\sqrt{\s^{(Y;m,n-m)}\cA(\ub_1)}\right\}\cdot\nonumber\\
&&\hspace{60mm}\cdot\ub_1^{2a_m}u_1^{2b_m}\cdot\left(\frac{\ub_1}{u_1}\right)^{b_m-3}\label{13.206}
\end{eqnarray}
(for a different constant $C$).

The second of the integrals \ref{13.195} is:
\begin{eqnarray}
&&\int_0^{u_1}u^2\left\{\int_{C_u^{\ub_1}}|E^{n-m}T^{m-1}\jb-E_N^{n-m}T^{m-1}\jb_N|
|\s^{(V;m,n-m)}\cxi_L|\right\}du\nonumber\\
&&\leq 
\int_0^{u_1}u^2\|E^{n-m}T^{m-1}\jb-E_N^{n-m}T^{m-1}\jb_N\|_{L^2(C_u^{\ub_1})}\|\s^{(V;m,n-m)}\cxi_L\|_{L^2(C_u^{\ub_1})}du
\nonumber\\
&&\label{13.207}
\end{eqnarray}
From Section 10.7, the leading contribution to $E^{n-m}T^{m-1}\jb-E_N^{n-m}T^{m-1}\jb_N$ comes from the 
$\Nb$ component of \ref{10.679} which is (see \ref{10.369}):
$$(\beta_N\beta_{\Nb}^2/4c)H^\prime\eta^2\Nb^\mu L(E^{n-m}T^m\beta_\mu-E_N^n T^m\beta_{\mu,N})$$
Expressing 
$$\Nb^\mu L(E^{n-m}T^m\beta_\mu-E_N^{n-m}T^m\beta_{\mu,N})$$
up to lower order terms as 
$$\lambdab\s^{(E;m,n-m)}\csxi$$ 
and noting that $\lambdab\leq Cu^{-1/2}\sqrt{a}$, it follows that 
$\|E^{n-m}T^{m-1}\jb-E_N^{n-m}T^{m-1}\jb_N\|_{L^2(C_u^{\ub_1})}$ is to leading terms bounded by:
\begin{eqnarray}
&&Cu^{-1/2}\|\sqrt{a}\s^{(E;m,n-m)}\csxi\|_{L^2(C_u^{\ub_1})}\nonumber\\
&&\hspace{10mm}\leq Cu^{-1/2}\sqrt{\s^{(E;m,n-m)}\cE^{\ub_1}(u)}
\leq C\ub_1^{a_m}u^{b_m-\frac{1}{2}}\sqrt{\s^{(E;m,n-m)}\cB(\ub_1,u)}\nonumber\\
&&\label{13.208}
\end{eqnarray}
by \ref{9.290} and \ref{9.297}. Together with \ref{13.201}, which by \ref{13.200} holds, a fortiori, also for 
$u\leq \ub_1$, we conclude that \ref{13.207} is to leading terms bounded by:
\begin{equation}
C\frac{\sqrt{\s^{(E;m,n-m)}\cB(\ub_1,u_1)\s^{(V;m,n-m)}\cB(\ub_1,u_1)}}{2b_m+\frac{5}{2}}\cdot 
\ub_1^{2a_m}u_1^{2b_m}\cdot u_1^{5/2}
\label{13.209}
\end{equation}
This is depressed relative to a borderline contribution by the factor $u_1^{5/2}$. 

Consider now \ref{13.196}. Since $\ub\leq u$, the first integral is bounded by:
\begin{eqnarray}
&&\left(\int_{{\cal R}_{\ub_1,u_1}}u^5|\cnub_{m-1,n-m+1}-\check{\taub}_{m-1,n-m+1}|^2\right)^{1/2}
\left(\int_{{\cal R}_{\ub_1,u_1}}\ub^{-1}|\s^{(V;m,n-m)}\cxi_{\Lb}|^2\right)^{1/2}\nonumber\\
&&=\left(\int_0^{u_1}u^5\|\cnub_{m-1,n-m+1}-\check{\taub}_{m-1,n-m+1}\|^2_{L^2(C_u^{\ub_1})}du\right)^{1/2}\cdot\nonumber\\
&&\hspace{30mm}\cdot\left(\int_0^{\ub_1}\ub^{-1}\|\s^{(V;m,n-m)}\cxi_{\Lb}\|^2_{L^2(\Cb_{\ub}^{u_1})}d\ub\right)^{1/2} \nonumber\\
&&\label{13.210}
\end{eqnarray}
In regard to the integral in the first factor, we substitute the estimate \ref{13.198}, \ref{13.199}. 
As we have shown above, the leading contribution comes from the first term on the right and is to 
leading terms is bounded by $k$ times \ref{13.204} with $u$ in the role of $\tau_1$ for $u\leq\ub_1$, 
$\ub_1$ in the role of $\tau_1$ for $u>\ub_1$. We then conclude that the integral in question is 
to leading terms bounded by:
\begin{eqnarray}
&&C\left({\frac{\s^{(Y;m,n-m)}\cBb(\ub_1,\ub_1)}{2a_m+1}}+\s^{(Y;m,n-m)}\cA(\ub_1)\right)
\cdot\nonumber\\
&&\hspace{30mm}\cdot\left\{\int_0^{\ub_1}u^{2a_m+2b_m-1}du+\int_{\ub_1}^{u_1}\ub_1^{2a_m+2b_m-6}u^5 du\right\}\nonumber\\
&&\label{13.211}
\end{eqnarray}
Here, in regard to the 2nd integral we have, since in that integral $u>\ub_1$, 
$$\int_{\ub_1}^{u_1}\ub_1^{2a_m+2b_m-6}u^5 du\leq \ub_1^{2a_m}\int_{\ub_1}^{u_1}u^{2b_m-1}du
=\ub_1^{2a_m}\frac{(u_1^{2b_m}-\ub_1^{2b_m})}{2b_m}$$
It then follows that \ref{13.211} is bounded by:
\begin{equation}
\frac{C}{2b_m}\left(\frac{\s^{(Y;m,n-m)}\cBb(\ub_1,\ub_1)}{2a_m+1}+\s^{(Y;m,n-m)}\cA(\ub_1)\right)
\cdot\ub_1^{2a_m}u_1^{2b_m}
\label{13.212}
\end{equation}
In regard to the integral in the second factor in \ref{13.210} we substitute:
\begin{equation}
\|\s^{(V;m,n-m)}\cxi_{\Lb}\|_{L^2(\Cb_{\ub}^{u_1})}\leq C\sqrt{\s^{(V;m,n-m)}\cEb^{u_1}(\ub)}
\leq C\ub^{a_m}u_1^{b_m}\sqrt{\s^{(V;m,n-m)}\cBb(\ub_1,u_1)}
\label{13.213}
\end{equation}
(by \ref{9.291} and \ref{9.298}) to obtain that the integral in question is bounded by:
\begin{equation}
\s^{(V;m,n-m)}\cBb(\ub_1,u_1)\cdot u_1^{2b_m}\cdot\int_0^{\ub_1}\ub^{2a_m-1}du
=\frac{\s^{(V;m,n-m)}\cBb(\ub_1,u_1)}{2a_m}\cdot\ub_1^{2a_m}u_1^{2b_m}
\label{13.214}
\end{equation}
Substituting this together with \ref{13.212} in \ref{13.210} then yields that the last is bounded by:
\begin{equation}
\frac{C\ub_1^{2a_m}u_1^{2b_m}}{\sqrt{2a_m}\sqrt{2b_m}}\sqrt{\s^{(V;m,n-m)}\cBb(\ub_1,u_1)}
\sqrt{\frac{\s^{(Y;m,n-m)}\cBb(\ub_1,\ub_1)}{2a_m+1}+\s^{(Y;m,n-m)}\cA(\ub_1)}
\label{13.215}
\end{equation}

Since again $\ub\leq u$, the second of the integrals \ref{13.196} is bounded by:
\begin{eqnarray}
&&\left(\int_{{\cal R}_{\ub_1,u_1}}u^5 |E^{n-m}T^{m-1}\jb-E_N^{n-m}T^{m-1}\jb_N|^2\right)^{1/2}
\cdot\nonumber\\
&&\hspace{40mm}\cdot\left(\int_{{\cal R}_{\ub_1,u_1}}\ub^{-1}|\s^{(V;m,n-m)}\cxi_{\Lb}|^2\right)^{1/2}\nonumber\\
&&=\left(\int_0^{u_1}u^5\|E^{n-m}T^{m-1}\jb-E_N^{n-m}T^{m-1}\jb_N\|^2_{L^2(C_u^{\ub_1})}du\right)^{1/2}
\cdot\nonumber\\
&&\hspace{30mm}\cdot\left(\int_0^{\ub_1}\ub^{-1}\|\s^{(V;m,n-m)}\cxi_{\Lb}\|^2_{L^2(\Cb_{\ub}^{u_1})}d\ub\right)^{1/2}\nonumber\\
&&\label{13.216}
\end{eqnarray}
Here, by \ref{13.208} the integral in the first factor is to leading terms bounded by:
\begin{equation}
C\s^{(E;m,n-m)}\cB(\ub_1,u_1)\cdot\ub_1^{2a_m}\cdot\int_0^{u_1}u^{2b_m+4}du
=C\frac{\s^{(E;m,n-m)}\cB(\ub_1,u_1)}{(2b_m+5)}\cdot\ub_1^{2a_m}u_1^{2b_m}\cdot u_1^5
\label{13.217}
\end{equation}
while the integral in the second factor is bounded by \ref{13.214}, hence \ref{13.216} is to 
leading terms bounded by:
\begin{equation}
C\frac{\sqrt{\s^{(V;m,n-m)}\cBb(\ub_1,u_1)}}{\sqrt{2a_m}}\frac{\sqrt{\s^{(E;m,n-m)}\cB(\ub_1,u_1)}}
{\sqrt{2b_m+5}}\cdot\ub_1^{2a_m}u_1^{2b_m}\cdot u_1^{5/2}
\label{13.218}
\end{equation}
which is depressed relative to a borderline contribution by the factor $u_1^{5/2}$. 

\vspace{5mm}

\section{The Top Order Energy Estimates}

We proceed to the energy estimates at the top order $m+l=n$. We introduce the notation 
(see \ref{9.300}, \ref{9.301}):
\begin{equation}
\s^{(m,l)}{\cal M}=\max\{\s^{(m,l)}\cB,\s^{(m,l)}\cBb\} 
\label{13.219}
\end{equation}
\begin{equation}
\s^{[m,l]}{\cal M}=\max_{i=0,...,m}\s^{(i,l+m-i)}{\cal M}=\max\{\s^{[m,l]}\cB,\s^{[m,l]}\cBb\}
\label{13.220}
\end{equation}

We first treat the case $m=0$. Consider first the boundary integral \ref{13.2}. By \ref{13.94} 
for the case $V=E$, and \ref{13.97} for the case $V=Y$, this is bounded to leading terms by:
\begin{equation}
C\ub_1^{2c_0}\left\{\frac{\s^{(0,n)}\cA(\ub_1)}{c_0}+\frac{\s^{(0,n)}\cBb(\ub_1,\ub_1)}{a_0^2}\right\}
\label{13.221}
\end{equation}
for both $V=E,Y$. Moreover, there are non-leading terms which are homogeneous in the $(0,n)$ 
energies and are bounded by:
\begin{equation}
C\ub_1^{2c_0+1}\left\{\s^{(0,n)}\cA(\ub_1)+\s^{(0,n)}{\cal M}(\ub_1,\ub_1)\right\}
\label{13.222}
\end{equation}

In addition, there is an inhomogeneous term \ref{13.88} arising from the error committed by the 
$N$th approximation in satisfying the boundary conditions. Since for an arbitrary function $f$ 
on ${\cal K}$ we have:
\begin{equation}
\|f\|_{L^2({\cal K}^{\tau_1})}\leq C\tau_1^{1/2}\|f\|_{L^\infty({\cal K}^{\tau_1})}
\label{13.223}
\end{equation}
by \ref{13.88} the contribution of this term to $\|\s^{(V;m,l)}\check{b}\|_{L^2({\cal K}^{\ub_1})}$ 
is bounded by:
\begin{equation}
C_{m,l,N}\ub_1^{N-\frac{1}{2}-m} \ \ \mbox{for both $V=E,Y$ and for any $m=0,...,n$}
\label{13.224}
\end{equation}
The corresponding contribution to the boundary integral \ref{13.2} is then bounded by:
\begin{equation}
C_{m,l,N}\ub_1^{2N-1-2m} \ \ \mbox{for both $V=E,Y$ and for any $m=0,...,n$}
\label{13.225}
\end{equation}

Combining \ref{13.221}, \ref{13.222}, and \ref{13.225} in the case $(m,l)=(0,n)$, we obtain:
\begin{eqnarray}
&&2C^\prime\int_{{\cal K}^{\ub_1}}\Omega^{1/2}\left(\s^{(V;0,n)}\check{b}\right)^2\leq
\ub_1^{2c_0}\left\{C\left(\frac{\s^{(0,n)}\cA(\ub_1)}{c_0}+\frac{\s^{(0,n)}\cBb(\ub_1,\ub_1)}{a_0^2}
\right)\right.\nonumber\\
&&\hspace{32mm}\left.+C\ub_1\left(\s^{(0,n)}\cA(\ub_1)+\s^{(0,n)}{\cal M}(\ub_1,\ub_1)\right)
+C_{0,n,N}\right\}\nonumber\\
&&\mbox{: for both $V=E,Y$}\label{13.226}
\end{eqnarray}
under the assumption that:
\begin{equation}
2N-1\geq 2c_0
\label{13.227}
\end{equation}

In $\s^{(V;m,n-m)}\check{{\cal G}}_1^{\ub_1,u_1}$ the leading terms are bounded by \ref{13.101} and \ref{13.102}. 
Moreover there are non-leading terms, which are bounded by:
\begin{equation}
C\ub_1^{2a_m}u_1^{2b_m+\frac{1}{2}}\s^{(m,n-m)}{\cal M}(\ub_1,u_1)
\label{13.228}
\end{equation}
the dominant of these terms being \ref{7.118} and \ref{7.119}. In 
$\s^{(V;m,n-m)}\check{{\cal G}}_2^{\ub_1,u_1}$, all terms are bounded by \ref{13.228}, the dominant 
terms being \ref{7.172} and \ref{7.187}. These results hold for both $V=E,Y$ and for any $m=0,...,n$. 

Considering $\s^{(V;m,n-m)}\check{{\cal G}}_3^{\ub_1,u_1}$ we first restrict attention to the case 
$m=0$. According to the results of the previous two sections, the leading terms  
from Section 13.3 (case $V=Y$ ) are bounded by \ref{13.125}, \ref{13.130}, \ref{13.135}, \ref{13.138}, 
and the leading terms from Section 13.4 (both cases $V=E,Y$) are bounded by \ref{13.176}, \ref{13.189}. 
These results simplify to:
\begin{eqnarray}
&&C\ub_1^{2a_0}u_1^{2b_0}\left\{\frac{1}{\sqrt{b_0}}
\frac{\s^{(E;0,n)}\cBb(\ub_1,u_1)}{\sqrt{a_0}}+\frac{\s^{(E;0,n)}\cB(\ub_1,u_1)}{b_0}\right.\nonumber\\
&&\hspace{18mm}\left.+\frac{1}{\sqrt{b_0}}\left(\frac{1}{\sqrt{a_0}}+\frac{1}{\sqrt{b_0}}\right)
\left(\frac{\s^{(Y;0,n)}\cBb(\ub_1,u_1)}{a_0}+\s^{(Y;0,n)}\cA(\ub_1)\right)
\right\}\nonumber\\
&&\label{13.229}
\end{eqnarray}
for the bound on the leading terms in the case $V=E$, and:
\begin{eqnarray}
&&C\ub_1^{2a_0}u_1^{2b_0}\left\{\left(\frac{1}{\sqrt{a_0}}+\frac{1}{\sqrt{b_0}}\right)
\frac{\s^{(Y;0,n)}\cBb(\ub_1,u_1)}{\sqrt{a_0}}+\frac{\s^{(Y;0,n)}\cB(\ub_1,u_1)}{b_0}
\right.\nonumber\\
&&\hspace{18mm}\left.+\left(\frac{1}{\sqrt{a_0}}+\frac{1}{\sqrt{b_0}}\right)
\frac{\s^{(Y;0,n)}\cA(\ub_1)}{\sqrt{b_0}}\right\}\nonumber\\
&&\label{13.230}
\end{eqnarray}
for the bound on the leading terms in the case $V=Y$. We may simplify the above further to a bound 
for the leading terms by:
\begin{equation}
C\ub_1^{2a_0}u_1^{2b_0}\left\{\left(\frac{1}{a_0}+\frac{1}{b_0}\right)\s^{(0,n)}{\cal M}(\ub_1,u_1)
+\left(\frac{1}{\sqrt{a_0}}+\frac{1}{\sqrt{b_0}}\right)
\frac{\s^{(0,n)}\cA(\ub_1)}{\sqrt{b_0}}\right\}
\label{13.231}
\end{equation}
in both cases $V=E,Y$. Moreover, there are non-leading terms, which are homogeneous in the 
$(0,n)$ energies and are bounded by:
\begin{equation}
C\ub_1^{2a_0}u_1^{2b_0+\frac{1}{2}}\left\{\s^{(0,n)}{\cal M}(\ub_1,u_1)+\s^{(0,n)}\cA(\ub_1)\right\}
\label{13.232}
\end{equation}

In addition, there are, for each $m=0,...,n$ inhomogeneous terms arising from the error committed 
by the $N$th approximation. The leading such inhomogeneous term is contributed through the factor 
$V^\mu\s^{(m,l)}\check{\tilde{\rho}}_\mu$ in \ref{13.104} by 
$$-\hat{E}_N^l\hat{T}_N^m\tilde{\kappa}^\prime_{\mu,N}$$
the last term in the expression for $\s^{(m,l)}\check{\tilde{\rho}}_\mu$ of Proposition 9.11. 
This term satisfies \ref{9.274}, which implies :
\begin{equation}
\|u\ub^{2+m-N}\hat{E}_N^l \hat{T}_N^m\tilde{\kappa}^\prime_{\mu,N}\|_{L^\infty({\cal R}_{\ub_1,u_1})}
\leq D_{m,l,N}(u_1)
\label{13.233}
\end{equation}
where $D_{m,l,N}(u_1)$ denotes a constant depending on the initial data including the up to the $N$th 
order derived data (see Chapter 5) on $\Cb_0^{u_1}$. In reference to \ref{13.104} we can then 
estimate:
\begin{eqnarray}
&&\int_{{\cal R}_{\ub_1,u_1}}|\s^{(V;m,l)}\cxi_L||V^\mu\hat{E}_N^l\hat{T}_N^m\tilde{\kappa}^\prime_{\mu,N}|\nonumber\\
&&\leq\left(\int_{{\cal R}_{\ub_1,u_1}}u^{-1}|\s^{(V;m,l)}\cxi_L|^2\right)^{1/2}
\left(\int_{{\cal R}_{\ub_1,u_1}}u|V^\mu\hat{E}_N^l\hat{T}_N^m\tilde{\kappa}^\prime_{\mu,N}|^2\right)^{1/2}\nonumber\\
&&\label{13.234}
\end{eqnarray}
Here, the integral in the 1st factor is bounded by:
\begin{equation}
C\frac{\s^{(V;m,l)}\cB(\ub_1,u_1)}{2b_m}\cdot\ub_1^{2a_m}u_1^{2b_m}
\label{13.235}
\end{equation}
On the other hand by \ref{13.233} the integral in the 2nd factor is bounded by:
\begin{equation}
CD_{m,l,N}(u_1)\int_{R_{\ub_1,u_1}}\ub^{2N-2m-5}dud\ub
\leq CD_{m,l,N}(u_1)\frac{u_1\ub_1^{2N-2m-4}}{2(N-m-2)}
\label{13.236}
\end{equation}
Therefore \ref{13.234} is bounded by:
\begin{equation}
C\sqrt{\frac{\s^{(V;m,l)}\cB(\ub_1,u_1)}{2b_m}} \ \ub_1^{a_m}u_1^{b_m+\frac{1}{2}}\cdot
\sqrt{\frac{D_{m,l,N}(u_1)}{2(N-m-2)}} \ \ub_1^{N-m-2}
\label{13.237}
\end{equation}
which does not exceed:
\begin{equation}
C\s^{(m,l)}{\cal M}(\ub_1,u_1)\ub_1^{2a_m}u_1^{2b_m+1}+D_{m,l,N}(u_1)\ub_1^{2(N-m-2)}
\label{13.238}
\end{equation}
Again in reference to \ref{13.104} we can also
estimate:
\begin{eqnarray}
&&\int_{{\cal R}_{\ub_1,u_1}}|\s^{(V;m,l)}\cxi_{\Lb}||V^\mu\hat{E}_N^l\hat{T}_N^m\tilde{\kappa}^\prime_{\mu,N}|\nonumber\\
&&\leq\left(\int_{{\cal R}_{\ub_1,u_1}}\ub^{-1}|\s^{(V;m,l)}\cxi_L|^2\right)^{1/2}
\left(\int_{{\cal R}_{\ub_1,u_1}}\ub|V^\mu\hat{E}_N^l\hat{T}_N^m\tilde{\kappa}^\prime_{\mu,N}|^2\right)^{1/2}\nonumber\\
&&\label{13.239}
\end{eqnarray}
Here, the integral in the 1st factor is bounded by:
\begin{equation}
C\frac{\s^{(V;m,l)}\cBb(\ub_1,u_1)}{2a_m}\cdot\ub_1^{2a_m}u_1^{2b_m}
\label{13.240}
\end{equation}
On the other hand by \ref{13.233} the integral in the 2nd factor is again bounded by \ref{13.235}, 
therefore \ref{13.239} is bounded by:
\begin{equation}
C\sqrt{\frac{\s^{(V;m,l)}\cBb(\ub_1,u_1)}{2a_m}} \ \ub_1^{a_m}u_1^{b_m+\frac{1}{2}}\cdot
\sqrt{\frac{D_{m,l,N}(u_1)}{2(N-m-2)}} \ \ub_1^{N-m-2}
\label{13.241}
\end{equation}
which likewise does not exceed \ref{13.238}. Under the assumption: 
\begin{equation}
N\geq c_m+m+2
\label{13.242}
\end{equation}
(compare with \ref{13.227}) we have $\ub_1^{2(N-m-2)}\leq\ub_1^{2c_m}\leq\ub_1^{2a_m}u_1^{2b_m}$ therefore \ref{13.238} does not 
exceed:
\begin{equation}
\ub_1^{2a_m}u_1^{2b_m}\left\{Cu_1\s^{(m,l)}{\cal M}(\ub_1,u_1)+D_{m,l,N}(u_1)\right\}
\label{13.243}
\end{equation}

In addition to the leading inhomogeneous term just considered, there are for $m=0$ inhomogeneous 
terms arising from the error term $E_N^n\vep_{\theta,N}$ in the propagation equation \ref{10.256} 
for $\cth_n$ and from the error term $E_N^n\vep_{\thetab,N}$ in the propagation equation 
\ref{10.257} for $\cthb_n$. Also, for $m=1,...,n$, $l=n-m$, inhomogeneous terms arising from the 
error term $E_N^{n-m}T^{m-1}\vep_{\nu,N}$ in the propagation equation \ref{10.336} for 
$\cnu_{m-1,n-m+1}$ and from the error term $E_N^{n-m}T^{m-1}\vep_{\nub,N}$ in the propagation 
equation \ref{10.337} for $\cnub_{m-1,n-m+1}$. 

In regard to the case $m=0$, $\vep_{\theta,N}$ and $\vep_{\thetab,N}$ are given by \ref{10.a1} 
and \ref{10.a2} respectively. By \ref{10.189}, \ref{10.192} and Proposition 9.10 :
\begin{equation}
\vep_{\theta,N} \ \vep_{\thetab,N} \ ={\bf O}(\tau^{N-2})+{\bf O}(u^{-3}\tau^N)
\label{13.244}
\end{equation}
therefore we have:
\begin{equation}
\|u\ub^{2-N}E_N^n\vep_{\theta,N}\|_{L^\infty({\cal R}_{\ub_1,u_1})}, \ 
\|u\ub^{2-N}E_N^n\vep_{\thetab,N}\|_{L^\infty({\cal R}_{\ub_1,u_1})} \ 
\leq D_{0,n,N}(u_1)
\label{13.245}
\end{equation}
It follows that for $(\ub,u)\in R_{\ub_1,u_1}$ the contribution of the error term 
$E_N^n\vep_{\theta,N}$ to $\|\cth_n\|_{L^2(S_{\ub,u})}$ is bounded by:
\begin{equation}
\sqrt{D_{0,n,N}(u_1)}\cdot u^{-1}\ub^{N-1}
\label{13.246}
\end{equation}
and the contribution of the error term $E_N^n\vep_{\thetab,N}$ to $\|\cthb_n\|_{L^2(S_{\ub,u})}$ is 
bounded by:
\begin{equation}
\sqrt{D_{0,n,N}(u_1)}\cdot\ub^{N-2}
\label{13.247}
\end{equation}
The corresponding contributions to 
$$\|\cth_n\|_{L^2(\Cb_{\ub}^u)}, \ \|\cthb_n\|_{C_u^{\ub})}$$
are then both bounded by:
\begin{equation}
\sqrt{D_{0,n,N}(u_1)}\cdot\ub^{N-\frac{3}{2}}
\label{13.248}
\end{equation}
The contributions through $\|\cth_n\|_{L^2(\Cb_{\ub}^u)}$ to \ref{13.122} and to \ref{13.131} are then 
bounded by 
\begin{equation}
\sqrt{D_{0,n,N}(u_1)}\sqrt{\s^{(Y;0,n)}\cBb(\ub_1,u_1)} \ 
\frac{\ub_1^{N-\frac{1}{2}+a_0}u_1^{b_0}}{N-\frac{1}{2}+a_0}
\label{13.249}
\end{equation}
and by 
\begin{equation}
\sqrt{D_{0,n,N}(u_1)} \ \frac{\ub_1^{N-\frac{1}{2}}}{\sqrt{2N-1}}\cdot 
\sqrt{\s^{(Y;0,n)}\cB(\ub_1,u_1)} \ \frac{\ub_1^{a_0}u_1^{b_0}}{\sqrt{2b_0}}
\label{13.250}
\end{equation}
respectively. Therefore the sum of these contributions is bounded by:
\begin{equation}
C\s^{(0,n)}{\cal M}(\ub_1,u_1)\ub_1^{2a_0+1}u_1^{2b_0}+D_{0,n,N}(u_1)\ub_1^{2N-2}
\label{13.251}
\end{equation}
To be able to replace $\ub_1^{2N-2}$ in the 2nd term by $\ub_1^{2a_0}u_1^{2b_0}$ it suffices that 
\begin{equation}
N\geq c_0+1
\label{13.252}
\end{equation}
a condition weaker than \ref{13.227}, which is in turn weaker than \ref{13.242} for $m=0$. 
Also, the contributions through $\|\cthb_n\|_{L^2(C_u^{\ub})}$ to \ref{13.167} and to \ref{13.184} are 
bounded by 
\begin{equation}
\sqrt{D_{0,n,N}}\sqrt{\s^{(V;0,n)}\cB(\ub_1,u_1)} \ \frac{\ub_1^{N-\frac{3}{2}+a_0}u_1^{b_0+3}}{b_0+3}
\label{13.253}
\end{equation}
and by 
\begin{equation}
\sqrt{D_{0,n,N}(u_1)} \ \frac{\ub_1^{N-\frac{3}{2}}u_1^3}{\sqrt{6}}\cdot 
\sqrt{\s^{(V;0,n)}\cBb(\ub_1,u_1)} \ \frac{\ub_1^{a_0}u_1^{b_0}}{\sqrt{2a_0}}
\label{13.254}
\end{equation} 
Therefore the sum of these contributions is bounded by:
\begin{equation}
C\s^{(0,n)}{\cal M}(\ub_1,u_1)\ub_1^{2a_0+1}u_1^{2b_0}+D_{0,n,N}(u_1)\ub^{2N-4}u_1^6
\label{13.255}
\end{equation}
To be able to replace $\ub_1^{2N-4}u_1^6$ in the 2nd term by $\ub_1^{2a_0}u_1^{2b_0}$ it suffices that 
\begin{equation}
N\geq c_0-1
\label{13.256}
\end{equation}
a condition even weaker than \ref{13.252}. 

Combining the above results for the case $m=0$, that is \ref{13.101}, \ref{13.102}, \ref{13.228} 
with $m=0$, \ref{13.231}, \ref{13.232}, and \ref{13.243} with $(m,l)=(0,n)$, and \ref{13.251}, 
\ref{13.255}, we conclude that: 
\begin{eqnarray}
&&\s^{(V;0,n)}\check{{\cal G}}^{\ub_1,u_1}\leq \ub_1^{2a_0}u_1^{2b_0}\left\{
C\left[\left(\frac{1}{a_0}+\frac{1}{b_0}\right)
\s^{(0,n)}{\cal M}(\ub_1,u_1)\right.\right.
\nonumber\\
&&\hspace{45mm}\left.+\left(\frac{1}{\sqrt{a_0}}+\frac{1}{\sqrt{b_0}}\right)
\frac{\s^{(0,n)}\cA(\ub_1)}{\sqrt{b_0}}\right]\nonumber\\
&&\hspace{30mm}\left.+C u_1^{1/2}\left[\s^{(0,n)}{\cal M}(\ub_1,u_1)+\s^{(0,n)}\cA(\ub_1)\right]
+D_{0,n,N}(u_1)\right\}\nonumber\\
&&\label{13.257}
\end{eqnarray}
for both $V=E,Y$. We substitute this bound taking $u_1=\ub_1$ on the right in the energy identity 
\ref{13.3} with $u_1=\ub_1$ for $(m,l)=(0,n)$. We also substitute on the right the bound \ref{13.226}. 
Keeping only the term $\s^{(V;0,n)}\check{{\cal F}}^{\prime\ub_1}$ on the left we obtain:
\begin{eqnarray}
&&\s^{(V;0,n)}\check{{\cal F}}^{\prime\ub_1}\leq \ub_1^{2c_0}\left\{
C\left[\left(\frac{1}{a_0}+\frac{1}{b_0}\right)
\s^{(0,n)}{\cal M}(\ub_1,\ub_1)\right.\right.
\nonumber\\
&&\hspace{45mm}\left.+\left(\frac{1}{\sqrt{a_0}}+\frac{1}{\sqrt{b_0}}\right)
\frac{\s^{(0,n)}\cA(\ub_1)}{\sqrt{b_0}}\right]\nonumber\\
&&\hspace{30mm}\left.+C\ub_1^{1/2}\left[\s^{(0,n)}{\cal M}(\ub_1,\ub_1)
+\s^{(0,n)}\cA(\ub_1)\right]+D_{0,n,N}(\ub_1)\right\}\nonumber\\
&&\label{13.258}
\end{eqnarray}
for both $V=E,Y$. Here we denote by $D_{0,n,N}(\ub_1)$ what was hitherto denoted by 
$$D_{0,n,N}(\ub_1)+C_{0,n,N}$$
Replacing then $\ub_1$ by any $\ub\in(0,\ub_1]$, multiplying by $\ub^{-2c_0}$ and taking the 
supremum over $\ub\in(0,\ub_1]$ yields, in view of the definitions \ref{9.299}, \ref{9.302}, 
\begin{eqnarray}
&&\s^{(0,n)}\cA(\ub_1)\leq 
\left[\frac{C}{\sqrt{b_0}}\left(\frac{1}{\sqrt{a_0}}+\frac{1}{\sqrt{b_0}}\right)
+C\ub_1^{1/2}\right]\s^{(0,n)}\cA(\ub_1)\nonumber\\
&&\hspace{20mm}+C\left(\frac{1}{a_0}+\frac{1}{b_0}\right)\s^{(0,n)}{\cal M}(\ub_1,\ub_1)
\nonumber\\
&&\hspace{20mm}+C\ub_1^{1/2}\s^{(0,n)}{\cal M}(\ub_1,\ub_1)+D_{0,n,N}(\ub_1)\nonumber\\
&&\label{13.259}
\end{eqnarray}
Taking $a_0$ and $b_0$ sufficiently large and $\delta$ sufficiently small so that in regard to the 
coefficient of $\s^{(0,n)}\cA(\ub_1)$ on the right we have:
\begin{equation}
\frac{C}{\sqrt{b_0}}\left(\frac{1}{\sqrt{a_0}}+\frac{1}{\sqrt{b_0}}\right)+C\delta^{1/2}
\leq\frac{1}{2}
\label{13.260}
\end{equation}
this implies:
\begin{eqnarray}
&&\s^{(0,n)}\cA(\ub_1)\leq 
C\left(\frac{1}{a_0}+\frac{1}{b_0}\right)\s^{(0,n)}{\cal M}(\ub_1,\ub_1)
\nonumber\\
&&\hspace{20mm}+C\ub_1^{1/2}\s^{(0,n)}{\cal M}(\ub_1,\ub_1)+D_{0,n,N}(\ub_1)\nonumber\\
&&\label{13.261}
\end{eqnarray}
(for new constants $C$). Substituting the bound \ref{13.261} in \ref{13.257} and in 
\ref{13.226} and substituting the resulting bounds on the right in the energy identity \ref{13.3} 
with $(m,l)=(0,n)$ the identity gives, dropping the term $\s^{(V;m,l)}\check{{\cal F}}^{\prime\ub_1}$ 
on the left, 
\begin{eqnarray}
&&\s^{(V;0,n)}\cE^{\ub_1}(u_1)+\s^{(V;0,n)}\cEb^{u_1}(\ub_1) \nonumber\\
&&\hspace{5mm}\leq \ub_1^{2a_0}u_1^{2b_0}\left\{\left[C\left(\frac{1}{a_0}+\frac{1}{b_0}\right)
+C u_1^{1/2}\right]\s^{(0,n)}{\cal M}(\ub_1,u_1)+D_{0,n,N}(u_1)\right\}\nonumber\\
&&\label{13.262}
\end{eqnarray}
Keeping only one of $\s^{(V;0,n)}\cE^{\ub_1}(u_1)$, $\s^{(V;0,n)}\cEb^{u_1}(\ub_1)$ 
at a time on the left, 
replacing $(\ub_1,u_1)$ by any $(\ub,u)\in R_{\ub_1,u_1}$, multiplying by $\ub^{-2a_0}u^{-2b_0}$, 
and taking the supremum over $(\ub,u)\in R_{\ub_1,u_1}$,  
yields, in view of the definitions \ref{9.297}, \ref{9.298}, 
\begin{eqnarray*}
&&\s^{(V;0,n)}\cB(\ub_1,u_1), \ \s^{(V;0,n)}\cBb(\ub_1,u_1) \\
&&\hspace{20mm}\leq\left[C\left(\frac{1}{a_0}+\frac{1}{b_0}\right)
+C u_1^{1/2}\right]\s^{(0,n)}{\cal M}(\ub_1,u_1)+D_{0,n,N}(u_1) 
\end{eqnarray*}
for both $V=E,Y$, therefore:
\begin{equation}
\s^{(0,n)}{\cal M}(\ub_1,u_1)\leq\left[C\left(\frac{1}{a_0}+\frac{1}{b_0}\right)
+C u_1^{1/2}\right]\s^{(0,n)}{\cal M}(\ub_1,u_1)+D_{0,n,N}(u_1) 
\label{13.263}
\end{equation}
Taking $a_0$ and $b_0$ sufficiently large and $\delta$ sufficiently small so that in regard to the 
coefficient of $\s^{(0,n)}{\cal M}(\ub_1,u_1)$ on the right we have:
\begin{equation}
C\left(\frac{1}{a_0}+\frac{1}{b_0}\right)+C\delta^{1/2}\leq\frac{1}{2}
\label{13.264}
\end{equation}
\ref{13.263} implies the estimate:
\begin{equation}
\s^{(0,n)}{\cal M}(\ub_1,u_1)\leq D_{0,n,N}(u_1) 
\label{13.265}
\end{equation}
Moreover, taking $u_1=\ub_1$ and substituting in \ref{13.261} yields the estimate:
\begin{equation}
\s^{(0,n)}\cA(\ub_1)\leq D_{0,n,N}(\ub_1)
\label{13.266}
\end{equation}

We turn to the cases $m=1,...,n$. Consider first the boundary integral \ref{13.2}. For 
$m=1,...,n$ leading terms are present only in the case $V=Y$ and are by \ref{13.99} bounded by:
\begin{equation}
C\ub_1^{2c_m}\left\{\frac{\s^{(Y;m,n-m)}\cA(\ub_1)}{c_m}
+\frac{\s^{(Y;m,n-m)}\cBb(\ub_1,\ub_1)}{a_m^2}\right\}
\label{13.267}
\end{equation}
Moreover, there are non-leading terms, present in both cases $V=E,Y$, which are homogeneous in the 
$(m,n-m)$ energies and are bounded by:
\begin{equation}
C\ub_1^{2c_m+1}\left\{\s^{[m,n-m]}\cA(\ub_1)+\s^{[m,n-m]}{\cal M}(\ub_1,\ub_1)\right\}
\label{13.268}
\end{equation}
In addition, there is the inhomogeneous term \ref{13.225}. Combining \ref{13.267}, \ref{13.268} 
and \ref{13.225} with $l=n-m$, and simplifying, we obtain:
\begin{eqnarray}
&&2C^\prime\int_{{\cal K}^{\ub_1}}\Omega^{1/2}\left(\s^{(V;m,n-m)}\check{b}\right)^2\nonumber\\
&&\hspace{27mm}\leq \ub_1^{2c_m}\left\{C\left(\frac{\s^{(m,n-m)}\cA(\ub_1)}{c_m}
+\frac{\s^{(m,n-m)}\cBb(\ub_1,\ub_1)}{a_m^2}\right)\right.\nonumber\\
&&\hspace{17mm}\left.+C\ub_1\left(\s^{[m,n-m]}\cA(\ub_1)+\s^{[m,n-m]}{\cal M}(\ub_1,\ub_1)\right)
+C_{m,n-m,N}\right\}\nonumber\\
&&\mbox{: for both $V=E,Y$}\label{13.269}
\end{eqnarray}
under the assumption that:
\begin{equation}
2N-1-2m\geq 2c_m
\label{13.270}
\end{equation}
which is weaker than the assumption \ref{13.242}. 

Turning to $\s^{(V;m,n-m)}\check{{\cal G}}_3^{\ub_1,u_1}$ we now consider the cases $m=1,...,n$. 
According to the results of the previous two sections, the leading terms  
from Section 13.3 (case $V=Y$ ) are bounded by \ref{13.146}, \ref{13.151}, \ref{13.156}, \ref{13.159}, 
and the leading terms from Section 13.4 (both cases $V=E,Y$) are bounded by \ref{13.206}, \ref{13.215}. 
These results simplify to:
\begin{eqnarray}
&&C\ub_1^{2a_m}u_1^{2b_m}\left\{\frac{1}{\sqrt{b_m}}
\frac{\s^{(E;m,n-m)}\cBb(\ub_1,u_1)}{\sqrt{a_m}}
+\frac{\s^{(E;m,n-m)}\cB(\ub_1,u_1)}{b_m}\right.\nonumber\\
&&\hspace{10mm}\left.+\frac{1}{\sqrt{b_m}}\left(\frac{1}{\sqrt{a_m}}+\frac{1}{\sqrt{b_m}}\right)
\left(\frac{\s^{(Y;m,n-m)}\cBb(\ub_1,u_1)}{a_m}+\s^{(Y;m,n-m)}\cA(\ub_1)\right)
\right\}\nonumber\\
&&\label{13.271}
\end{eqnarray}
for the bound on the leading terms in the case $V=E$, and:
\begin{eqnarray}
&&C\ub_1^{2a_m}u_1^{2b_m}\left\{\left(\frac{1}{\sqrt{a_m}}+\frac{1}{\sqrt{b_m}}\right)
\frac{\s^{(Y;m,n-m)}\cBb(\ub_1,u_1)}{\sqrt{a_m}}+\frac{\s^{(Y;m,n-m)}\cB(\ub_1,u_1)}{b_m}
\right.\nonumber\\
&&\hspace{20mm}\left.+\left(\frac{1}{\sqrt{a_m}}+\frac{1}{\sqrt{b_m}}\right)
\frac{\s^{(Y;m,n-m)}\cA(\ub_1)}{\sqrt{b_m}}\right\}\nonumber\\
&&\label{13.272}
\end{eqnarray}
for the bound on the leading terms in the case $V=Y$. We may simplify the above further to a bound 
for the leading terms by:
\begin{equation}
C\ub_1^{2a_m}u_1^{2b_m}\left\{\left(\frac{1}{a_m}+\frac{1}{b_m}\right)\s^{(m,n-m)}{\cal M}(\ub_1,u_1)
+\left(\frac{1}{\sqrt{a_m}}+\frac{1}{\sqrt{b_m}}\right)
\frac{\s^{(m,n-m)}\cA(\ub_1)}{\sqrt{b_m}}\right\}
\label{13.273}
\end{equation}
in both cases $V=E,Y$. Moreover, there are non-leading terms, which are homogeneous in the 
$(m,n-m)$ energies and are bounded by:
\begin{equation}
C\ub_1^{2a_m}u_1^{2b_+\frac{1}{2}}\left\{\s^{[m,n-m]}{\cal M}(\ub_1,u_1)+\s^{[m,n-m]}\cA(\ub_1)\right\}
\label{13.274}
\end{equation}
In addition, we have the leading inhomogeneous term \ref{13.243} with $l=n-m$, as well as the 
inhomogeneous terms arising from the error term $E_N^{n-m}T^{m-1}\vep_{\nu,N}$ in the 
propagation equation \ref{10.336} for $\cnu_{m-1,n-m+1}$ and from the error term 
$E_N^{n-m}T^{m-1}\vep_{\nub,N}$ in the propagation equation \ref{10.337} for $\cnub_{m-1,n-m+1}$. 
By \ref{10.334} the contribution of the first to $\|\cnu_{m-1,n-m+1}\|_{L^2(S_{\ub,u})}$ is bounded by:
\begin{equation}
\sqrt{D_{m,n-m,N}(u_1)}\cdot u^{-1}\ub^{N-m}
\label{13.275}
\end{equation}
and the corresponding contribution to $\|\cnu_{m-1,n-m+1}\|_{L^2(\Cb_{\ub}^u)}$ by:
\begin{equation}
\sqrt{D_{m,n-m,N}(u_1)}\cdot\ub^{N-m-\frac{1}{2}}
\label{13.276}
\end{equation}
while the contribution of the second to $\|\cnub_{m-1,n-m+1}\|_{L^2(S_{\ub,u})}$ is bounded by:
\begin{equation}
\sqrt{D_{m,n-m,N}(u_1)}\cdot\ub^{N-2-m}
\label{13.277}
\end{equation}
and the corresponding contribution to $\|\cnub_{m-1,n-m+1}\|_{L^2(C_u^{\ub})}$ by:
\begin{equation}
\sqrt{D_{m,n-m,N}(u_1)}\cdot\ub^{N-m-\frac{3}{2}}
\label{13.278}
\end{equation}
The contributions through $\|\cnu_{m-1,n-m+1}\|_{L^2(\Cb_{\ub}^u)}$ to \ref{13.143} 
and to \ref{13.152} are then bounded by 
\begin{equation}
\sqrt{D_{m,n-m,N}(u_1)}\sqrt{\s^{(Y;m,n-m)}\cBb(\ub_1,u_1)} \ 
\frac{\ub_1^{N+\frac{1}{2}-m+a_m}u_1^{b_m}}{N+\frac{1}{2}-m+a_m}
\label{13.279}
\end{equation}
and by 
\begin{equation}
\sqrt{D_{m,n-m,N}(u_1)} \ \frac{\ub_1^{N+\frac{1}{2}-m}}{\sqrt{2N+1-2m}}\cdot 
\sqrt{\s^{(Y;m,n-m)}\cB(\ub_1,u_1)} \ \frac{\ub_1^{a_m}u_1^{b_m}}{\sqrt{2b_m}}
\label{13.280}
\end{equation}
respectively. Therefore the sum of these contributions is bounded by:
\begin{equation}
C\s^{(m,n-m)}{\cal M}(\ub_1,u_1)\ub_1^{2a_m+1}u_1^{2b_m}+D_{m,n-m,N}(u_1)\ub_1^{2N-2m}
\label{13.281}
\end{equation}
To be able to replace $\ub_1^{2N-2m}$ in the 2nd term by $\ub_1^{2a_m}u_1^{2b_m}$ it suffices that 
\begin{equation}
N\geq c_m+m
\label{13.282}
\end{equation}
a condition weaker than \ref{13.270}, which is in turn weaker than \ref{13.242}. 
Also, the contributions through $\|\cnub_{m-1,n-m+1}\|_{L^2(C_u^{\ub})}$ to \ref{13.197} 
and to \ref{13.210} are bounded by 
\begin{equation}
\sqrt{D_{m,n-m,N}}\sqrt{\s^{(V;m,n-m)}\cB(\ub_1,u_1)} \ \frac{\ub_1^{N-\frac{3}{2}-m+a_m}u_1^{b_m+3}}{b_m+3}
\label{13.283}
\end{equation}
and by 
\begin{equation}
\sqrt{D_{m,n-m,N}(u_1)} \ \frac{\ub_1^{N-\frac{3}{2}-m}u_1^3}{\sqrt{6}}\cdot 
\sqrt{\s^{(V;m,n-m)}\cBb(\ub_1,u_1)} \ \frac{\ub_1^{a_m}u_1^{b_m}}{\sqrt{2a_m}}
\label{13.284}
\end{equation} 
Therefore the sum of these contributions is bounded by:
\begin{equation}
C\s^{(m,n-m)}{\cal M}(\ub_1,u_1)\ub_1^{2a_m+1}u_1^{2b_m}+D_{m,n-m,N}(u_1)\ub^{2N-4-2m}u_1^6
\label{13.285}
\end{equation}
To be able to replace $\ub_1^{2N-4-2m}u_1^6$ in the 2nd term by $\ub_1^{2a_m}u_1^{2b_m}$ 
it suffices that 
\begin{equation}
N\geq c_m+m-1
\label{13.286}
\end{equation}
a condition even weaker than \ref{13.282}.

Combining the above results for $m=1,...,n$, that is \ref{13.101}, \ref{13.102}, \ref{13.228}, \ref{13.273}, \ref{13.274}, and \ref{13.243} with $l=n-m$, and \ref{13.281}, 
\ref{13.285}, we conclude that: 
\begin{eqnarray}
&&\s^{(V;m,n-m)}\check{{\cal G}}^{\ub_1,u_1}\leq \ub_1^{2a_m}u_1^{2b_m}\left\{
C\left[\left(\frac{1}{a_m}+\frac{1}{b_m}\right)
\s^{(m,n-m)}{\cal M}(\ub_1,u_1)\right.\right.
\nonumber\\
&&\hspace{48mm}\left.+\left(\frac{1}{\sqrt{a_m}}+\frac{1}{\sqrt{b_m}}\right)
\frac{\s^{(m,n-m)}\cA(\ub_1)}{\sqrt{b_m}}\right]\nonumber\\
&&\hspace{10mm}\left.+C u_1^{1/2}\left[\s^{[m,n-m]}{\cal M}(\ub_1,u_1)+\s^{[m,n-m]}\cA(\ub_1)\right]
+D_{m,n-m,N}(u_1)\right\}\nonumber\\
&&\label{13.287}
\end{eqnarray}
for both $V=E,Y$. We substitute this bound taking $u_1=\ub_1$ on the right in the energy identity 
\ref{13.3} with $l=n-m$ and $u_1=\ub_1$. We also substitute on the right the bound \ref{13.269}. 
Keeping only the term $\s^{(V;m,n-m)}\check{{\cal F}}^{\prime\ub_1}$ on the left we obtain:
\begin{eqnarray}
&&\s^{(V;m,n-m)}\check{{\cal F}}^{\prime\ub_1}\leq \ub_1^{2c_m}\left\{
C\left[\left(\frac{1}{a_m}+\frac{1}{b_m}\right)
\s^{(m,n-m)}{\cal M}(\ub_1,\ub_1)\right.\right.
\nonumber\\
&&\hspace{48mm}\left.+\left(\frac{1}{\sqrt{a_m}}+\frac{1}{\sqrt{b_m}}\right)
\frac{\s^{(m,n-m)}\cA(\ub_1)}{\sqrt{b_m}}\right]\nonumber\\
&&\hspace{10mm}\left.+C\ub_1^{1/2}\left[\s^{[m,n-m]}{\cal M}(\ub_1,\ub_1)
+\s^{[m,n-m]}\cA(\ub_1)\right]+D_{m,n-m,N}(\ub_1)\right\}\nonumber\\
&&\label{13.288}
\end{eqnarray}
for both $V=E,Y$. Here we denote by $D_{m,n-m,N}(\ub_1)$ what was hitherto denoted by 
$$D_{m,n-m,N}(\ub_1)+C_{m,n-m,N}$$
We substitute:
\begin{equation}
\s^{[m,n-m]}\cA=\max\{\s^{(m,n-m)}\cA,\s^{[m-1,n-m+1]}\cA\}\leq \s^{(m,n-m)}\cA+\s^{[m-1,n-m+1]}\cA
\label{13.289}
\end{equation}
Replacing then $\ub_1$ by any $\ub\in(0,\ub_1]$, multiplying by $\ub^{-2c_m}$ and taking the 
supremum over $\ub\in(0,\ub_1]$ yields, in view of the definitions \ref{9.299}, \ref{9.302}, 
\begin{eqnarray}
&&\s^{(m,n-m)}\cA(\ub_1)\leq 
\left[\frac{C}{\sqrt{b_m}}\left(\frac{1}{\sqrt{a_m}}+\frac{1}{\sqrt{b_m}}\right)
+C\ub_1^{1/2}\right]\s^{(m,n-m)}\cA(\ub_1)\nonumber\\
&&\hspace{26mm}+C\left(\frac{1}{a_m}+\frac{1}{b_m}\right)\s^{(m,n-m)}{\cal M}(\ub_1,\ub_1)
\nonumber\\
&&\hspace{26mm}+C\ub_1^{1/2}\left[\s^{[m,n-m]}{\cal M}(\ub_1,\ub_1)+\s^{[m-1,n-m+1]}\cA(\ub_1)\right]\nonumber\\
&&\hspace{26mm}+D_{m,n-m,N}(\ub_1)
\label{13.290}
\end{eqnarray}
Taking $a_m$ and $b_m$ sufficiently large and $\delta$ sufficiently small so that in regard to the 
coefficient of $\s^{(m,n-m)}\cA(\ub_1)$ on the right we have:
\begin{equation}
\frac{C}{\sqrt{b_m}}\left(\frac{1}{\sqrt{a_m}}+\frac{1}{\sqrt{b_m}}\right)+C\delta^{1/2}
\leq\frac{1}{2}
\label{13.291}
\end{equation}
this implies:
\begin{eqnarray}
&&\s^{(m,n-m)}\cA(\ub_1)\leq 
C\left(\frac{1}{a_m}+\frac{1}{b_m}\right)\s^{(m,n-m)}{\cal M}(\ub_1,\ub_1)
\nonumber\\
&&\hspace{26mm}+C\ub_1^{1/2}\left[\s^{[m,n-m]}{\cal M}(\ub_1,\ub_1)+\s^{[m-1,n-m+1]}\cA(\ub_1)\right]\nonumber\\
&&\hspace{26mm}+D_{m,n-m,N}(\ub_1)
\label{13.292}
\end{eqnarray}
(for new constants $C$). Substituting the bound \ref{13.292} in \ref{13.287} and in 
\ref{13.269} and substituting the resulting bounds on the right in the energy identity \ref{13.3} 
with $l=n-m$ and substituting also 
\begin{equation}
\s^{[m,n-m]}{\cal M}=\max\{\s^{(m,n-m)}{\cal M},\s^{[m-1,n-m+1]}{\cal M}\}\leq \s^{(m,n-m)}{\cal M}
+\s^{[m-1,n-m+1]}{\cal M}
\label{13.293}
\end{equation}
the identity gives, dropping the term $\s^{(V;m,n-m)}\check{{\cal F}}^{\prime\ub_1}$ 
on the left, 
\begin{eqnarray}
&&\s^{(V;m,n-m)}\cE^{\ub_1}(u_1)+\s^{(V;m,n-m)}\cEb^{u_1}(\ub_1) \nonumber\\
&&\hspace{10mm}\leq \ub_1^{2a_m}u_1^{2b_m}\left\{\left[C\left(\frac{1}{a_m}+\frac{1}{b_m}\right)
+C u_1^{1/2}\right]\s^{(m,n-m)}{\cal M}(\ub_1,u_1)\right.\nonumber\\
&&\left.+Cu_1^{1/2}\left[\s^{[m-1,n-m+1]}{\cal M}(\ub_1,u_1)
+\s^{[m-1,n-m+1]}\cA(\ub_1)\right]+D_{m,n-m,N}(u_1)\right\} \nonumber\\
&&\label{13.294}
\end{eqnarray}
Keeping only one of $\s^{(V;m,n-m)}\cE^{\ub_1}(u_1)$, $\s^{(V;m,n-m)}\cEb^{u_1}(\ub_1)$ 
at a time on the left, 
replacing $(\ub_1,u_1)$ by any $(\ub,u)\in R_{\ub_1,u_1}$, multiplying by $\ub^{-2a_m}u^{-2b_m}$, 
and taking the supremum over $(\ub,u)\in R_{\ub_1,u_1}$,  
yields, in view of the definitions \ref{9.297}, \ref{9.298}, 
\begin{eqnarray*}
&&\s^{(V;m,n-m)}\cB(\ub_1,u_1), \ \s^{(V;m,n-m)}\cBb(\ub_1,u_1) \\
&&\hspace{20mm}\leq\left[C\left(\frac{1}{a_m}+\frac{1}{b_m}\right)
+C u_1^{1/2}\right]\s^{(m,n-m)}{\cal M}(\ub_1,u_1)\nonumber\\
&&+Cu_1^{1/2}\left[\s^{[m-1,n-m+1]}{\cal M}(\ub_1,u_1)
+\s^{[m-1,n-m+1]}\cA(\ub_1)\right]+D_{m,n-m,N}(u_1)
\end{eqnarray*}
for both $V=E,Y$, therefore:
\begin{eqnarray}
&&\s^{(m,n-m)}{\cal M}(\ub_1,u_1)\leq\left[C\left(\frac{1}{a_m}+\frac{1}{b_m}\right)
+C u_1^{1/2}\right]\s^{(m,n-m)}{\cal M}(\ub_1,\ub_1)\nonumber\\
&&\hspace{22mm}+Cu_1^{1/2}\left[\s^{[m-1,n-m+1]}{\cal M}(\ub_1,u_1)
+\s^{[m-1,n-m+1]}\cA(\ub_1)\right]\nonumber\\
&&\hspace{33mm}+D_{m,n-m,N}(u_1)
\label{13.295}
\end{eqnarray}
Taking $a_m$ and $b_m$ sufficiently large and $\delta$ sufficiently small so that in regard to the 
coefficient of $\s^{(m,n-m)}{\cal M}(\ub_1,u_1)$ on the right we have:
\begin{equation}
C\left(\frac{1}{a_m}+\frac{1}{b_m}\right)+C\delta^{1/2}\leq\frac{1}{2}
\label{13.296}
\end{equation}
\ref{13.295} implies:
\begin{eqnarray}
&&\s^{(m,n-m)}{\cal M}(\ub_1,u_1)\leq Cu_1^{1/2}\left[\s^{[m-1,n-m+1]}{\cal M}(\ub_1,u_1)
+\s^{[m-1,n-m+1]}\cA(\ub_1)\right]\nonumber\\
&&\hspace{33mm}+D_{m,n-m,N}(u_1)
\label{13.297}
\end{eqnarray}
Substituting \ref{13.293} in \ref{13.292} and substituting for $\s^{(m,n-m)}{\cal M}(\ub_1,\ub_1)$ 
in the resulting inequality from \ref{13.297} we then obtain:
\begin{eqnarray}
&&\s^{(m,n-m)}\cA(\ub_1)\leq 
C\ub_1^{1/2}\left[\s^{[m-1,n-m+1]}{\cal M}(\ub_1,\ub_1)
+\s^{[m-1,n-m+1]}\cA(\ub_1)\right]\nonumber\\
&&\hspace{26mm}+D_{m,n-m,N}(\ub_1)
\label{13.298}
\end{eqnarray}
The inequalities \ref{13.297}, \ref{13.298} constitute a recursion in $m$ with starting point the 
estimates \ref{13.265}, \ref{13.266} and yield estimates for $\s^{(m,n-m)}{\cal M}(\ub_1,u_1)$, 
$\s^{(m,n-m)}\cA(\ub_1)$ for all $m=0,...,n$, that is a priori bounds for all energies at the top 
order. 

\vspace{2.5mm} 

We remark that the addition of a term of the form:
\begin{equation}
C_m\ub_1\left[\s^{[m-1,n-m+1]}{\cal A}(\ub_1)+\s^{[m-1,n-m+1]}{\cal M}(\ub_1,u_1)\right]
\label{13.299}
\end{equation}
to the right hand side of \ref{13.269}, and of a term of the form:
\begin{equation}
C_m u_1^{1/2}\left[\s^{[m-1,n-m+1]}{\cal M}(\ub_1,u_1)+\s^{[m-1,n-m+1]}{\cal A}(\ub_1)\right]
\label{13.300}
\end{equation}
to the right hand side of \ref{13.287} does not affect the argument, which, being recursive in $m$, 
depends only on the assumptions \ref{13.260}, \ref{13.264}, \ref{13.291}, \ref{13.296}. 

\vspace{2.5mm}

In the preceding we have ignored the lower order terms. In Sections 14.1 and 14.8 of the next chapter 
we shall show how the lower order terms are taken account of. 

\pagebreak

\chapter{Lower Order Estimates, Recovery of the Bootstrap Assumptions, and Completion of the Argument} 

\section{The Bootstrap Assumptions Needed} 

The estimates of the preceding chapters depend on certain bootstrap assumptions which we shall now 
discuss. The fundamental bootstrap assumptions concern the quantities of order 0. These are:
\begin{equation}
\mbox{ the positive quantities}  \ \sigma, \ \beta_0 ; \ c, \ \sh; \ \rho, \rhob ; \ 
\mbox{ and the quantity} \ \ \pi 
\label{14.a1}
\end{equation}
the fact that $\beta_0$ is positive being an immediate consequence of the definition \ref{1.53}. 
We shall show that the bootstrap assumptions on the quantities \ref{14.a1} control all the other 
0th order quantities, namely all components of $\beta_\mu$ and of the frame field $E^\mu$, $N^\mu$, 
$\Nb^\mu$.  

By the results of Chapter 9 we can chose $\delta_0>0$ and suitably small so that the quantities
\begin{equation}
\inf_{{\cal R}_{\delta_0,\delta_0}}\sigma_N, \ \inf_{{\cal R}_{\delta_0,\delta_0}}c_N, \ 
\inf_{{\cal R}_{\delta_0,\delta_0}}\sh_N, \ \inf_{{\cal R}_{\delta_0,\delta_0}}(\rho_N/\ub), 
\ \inf_{{\cal R}_{\delta_0,\delta_0}}(\rhob_N/u^2)
\label{14.a2}
\end{equation}
are all positive. We then set:
\begin{equation}
\Ab_*=\frac{1}{2}\min\left\{\inf_{{\cal R}_{\delta_0,\delta_0}}\sigma_N, \inf_{{\cal R}_{\delta_0,\delta_0}}c_N, \inf_{{\cal R}_{\delta_0,\delta_0}}\sh_N, 
\inf_{{\cal R}_{\delta_0,\delta_0}}(\rho_N/\ub), 
\inf_{{\cal R}_{\delta_0,\delta_0}}(\rhob_N/u^2)\right\}
\label{14.a3}
\end{equation}
and:
\begin{equation}
\oA_*=2\max\left\{\sup_{{\cal R}_{\delta_0,\delta_0}}\beta_{0,N}, 
\sup_{{\cal R}_{\delta_0,\delta_0}}|\pi_N|, 
\sup_{{\cal R}_{\delta_0,\delta_0}}\sh_N, \sup_{{\cal R}_{\delta_0,\delta_0}}(\rho_N/\ub), 
\sup_{{\cal R}_{\delta_0,\delta_0}}(\rhob_N/u^2)\right\}
\label{14.a4}
\end{equation}

In the following $\delta$ is a positive real number not exceeding $\delta_0$ which enters the 
argument.

\vspace{2.5mm}

\noindent{\large{\bf Fundamental Bootstrap Assumptions:}}

\vspace{2.5mm}

\noindent On ${\cal R}_{\delta,\delta}$ we have:
\begin{eqnarray*}
&&\sigma\geq \Ab_*, \ \ c\geq \Ab_*, \ \ \sh\geq \Ab_* \\
&&\frac{\rho}{\ub}\geq \Ab_*, \ \ \frac{\rhob}{u^2}\geq \Ab_*
\end{eqnarray*}
and:
\begin{eqnarray*}
&&\beta_0\leq \oA_*, \ \ |\pi|\leq \oA_*, \ \ \sh\leq \oA_*\\
&&\frac{\rho}{\ub}\leq \oA_*, \ \ \frac{\rhob}{u^2}\leq \oA_*
\end{eqnarray*}

\vspace{2.5mm}

Since 
$$\sigma=-(g^{-1})^{\mu\nu}\beta_\mu\beta_\nu=\beta_0^2-\sum_i\beta_i^2$$
the first of the upper bound bootstrap assumptions implies:
\begin{equation}
\sigma\leq \oA_*^2 \ \mbox{: on ${\cal R}_{\delta,\delta}$}
\label{14.a5}
\end{equation}
therefore the range of $\sigma$ in ${\cal R}_{\delta,\delta}$ is contained in the compact interval 
$[\Ab_*,\oA_*^2]$ of the positive real line. 
Now the functions $G$ and $F$ (see \ref{2.131}) are given smooth positive functions on the positive real line. 
It follows that the 
functions $G\circ\sigma$ and $F\circ\sigma$, which we have denoted by $G(\sigma)$ and $F(\sigma)$, 
have fixed positive lower and upper bounds 
in ${\cal R}_{\delta,\delta}$. Hence (see \ref{2.132}) 
$$\eta(\sigma)=\frac{1}{\sqrt{1+\sigma F(\sigma)}}$$
is in ${\cal R}_{\delta,\delta}$ bounded below by a fixed positive constant and above by a constant strictly less than 1. Also, $H(\sigma)=\eta^{2}(\sigma)F(\sigma)$ is bounded above and below 
by fixed positive constants in ${\cal R}_{\delta,\delta}$ and all the derivatives of $H$ composed 
with $\sigma$ are  in ${\cal R}_{\delta,\delta}$ bounded in absolute value by fixed positive constants. 
Moreover by \ref{2.108} the conformal factor $\Omega(\sigma)$ has fixed positive lower and upper 
bounds in ${\cal R}_{\delta,\delta}$.

Since $\sigma>0$, the first of the upper bound bootstrap assumptions also implies:
\begin{equation}
\sqrt{\sum_i\beta_i^2}<\beta_0\leq\oA_*
\label{14.a6}
\end{equation}

Consider next the second of the upper bound bootstrap assumptions. Since 
$$h_{\mu\nu}E^\mu E^\nu=1 \ \ \mbox{and} \ \ h_{\mu\nu}=g_{\mu\nu}+H\beta_\mu\beta_\nu$$
we have:
\begin{equation}
g_{\mu\nu}E^\mu E^\nu=1-H\sbeta^2\leq 1
\label{14.a7}
\end{equation}
Since
$$g_{\mu\nu}E^\mu E^\nu=-(E^0)^2+\sum_i(E^i)^2 \ \ \mbox{and} \ \ E^0=\pi$$
it follows that:
\begin{equation}
\sqrt{\sum_i(E^i)^2}\leq \sqrt{1+\pi^2}\leq \sqrt{1+\oA_*^2}
\label{14.a8}
\end{equation}

We turn to the vectors $N^\mu$, $\Nb^\mu$. Since these vectors are null future-directed relative 
to $h_{\mu\nu}$ and the future null cone of $h_{\mu\nu}$ (the sound cone) is contained within the 
future null cone of $g_{\mu\nu}$ (the light cone), these vectors are timelike future-directed 
relative to $g_{\mu\nu}$. Since $N^0=\Nb^0=1$, we then have:
\begin{equation}
\sqrt{\sum_i(N^i)^2}<1, \ \ \ \sqrt{\sum_i(\Nb^i)^2}<1
\label{14.a9}
\end{equation}
Furthermore since by the definition \ref{1.53} the 1-form $\beta$ has positive evaluation in 
the interior of the future null cone of $g_{\mu\nu}$ (and on the boundary), we have:
\begin{equation}
\beta_N>0, \ \ \ \beta_{\Nb}>0
\label{14.a10}
\end{equation}
In fact \ref{14.a9} and \ref{14.a10} can be made more precise using the first of the lower bound 
and the first of the upper bound bootstrap assumptions. We have:
\begin{equation}
1-\sum_i(N^i)^2=-g_{\mu\nu}N^\mu N^\nu=H\beta_N^2, \ \ \ 
1-\sum_i(\Nb^i)^2=-g_{\mu\nu}\Nb^\mu\Nb^\nu=H\beta_{\Nb}^2
\label{14.a11}
\end{equation}
Also, 
\begin{eqnarray*}
&&\beta_N=\beta_0+\sum_i\beta_i N^i\geq\beta_0-\sqrt{\sum_i\beta_i^2}\sqrt{\sum_i(N^i)^2}
\geq \beta_0-\sqrt{\sum_i\beta_i^2}\\
&&\beta_{\Nb}=\beta_0+\sum_i\beta_i \Nb^i\geq\beta_0-\sqrt{\sum_i\beta_i^2}\sqrt{\sum_i(\Nb^i)^2}
\geq \beta_0-\sqrt{\sum_i\beta_i^2}
\end{eqnarray*}
by \ref{14.a9}. It follows that:
\begin{equation}
\beta_N, \beta_{\Nb}\geq \frac{\sigma}{2\beta_0}
\label{14.a12}
\end{equation}
Thus, the first of the lower bound together with the first of the upper bound bootstrap assumptions 
give a positive lower bound on $\beta_N$, $\beta_{\Nb}$, which together with the positive lower 
bound on $H$ provide a positive lower bound for the right hand sides in \ref{14.a11}.  

We turn to $c$, which satisfies the second of the lower bound bootstrap assumptions. An upper bound 
for $c$ by a constant strictly less than 1 is obtained using the first lower and upper bound 
bootstrap assumptions, as follows. We have:
\begin{eqnarray*}
&&c=-\frac{1}{2}h_{\mu\nu}N^\mu\Nb^\nu=-\frac{1}{2}(g_{\mu\nu}N^\mu\Nb^\nu+H\beta_N\beta_{\Nb})\\
&&\hspace{5mm}=\frac{1}{2}\left(1-\sum_i N^i\Nb^i-H\beta_N\beta_{\Nb}\right) 
\end{eqnarray*}
and, by \ref{14.a11}, 
$$-\sum_i N^i\Nb^i\leq\frac{1}{2}\sum_i(N^i)^2+\frac{1}{2}(\Nb^i)^2
\leq 1-\frac{1}{2}H\beta_N^2-\frac{1}{2}H\beta_{\Nb}^2$$
we obtain:
\begin{equation}
c\leq 1-\frac{1}{4}H(\beta_N+\beta_{\Nb})^2
\label{14.a13}
\end{equation}
Then the positive lower bound on $\beta_N$, $\beta_{\Nb}$ and on $H$ provides an upper bound 
on $c$ by a constant less than 1. 

\vspace{2.5mm}

We remark that the fundamental bootstrap assumptions imply 
a bound on the coercivity constant $C^\prime$ in Proposition 7.1. This is obtained through 
\ref{7.29} - \ref{7.35}, after suitably restricting $\delta_0$.

\vspace{2.5mm}

In addition to the above fundamental bootstrap assumptions, the fundamental assumptions \ref{10.775} 
and \ref{10.793}, which concern 0th order quantities defined on the boundary ${\cal K}$, have been used. 
In regard to the assumption \ref{10.775}, by virtue of \ref{4.215} we have, 
restricting suitably $\delta_0$, 
\begin{equation}
2\Bb\leq\frac{r_N}{\tau}\leq\frac{\oB}{2}\ \ \mbox{: in ${\cal K}^{\delta_0}$}
\label{14.a14}
\end{equation}
where:
\begin{equation}
\Bb=\inf_{S_{0,0}}\left(-\frac{l}{8c_0}\right), \ \ \ \oB=\sup_{S_{0,0}}\left(-2\frac{l}{c_0}\right)
\label{14.a15}
\end{equation}
In regard to the assumption \ref{10.793}, setting 
\begin{equation}
j_N=j(\kappa_N,\epb_N)
\label{14.a16}
\end{equation}
we define: 
\begin{equation}
B=2\sup_{{\cal K}^{\delta_0}}|j_N|
\label{14.a17}
\end{equation}

There is an additional fundamental bootstrap assumption needed to guarantee that the 
requirements of the shock development problem are fulfilled. By \ref{4.218}, \ref{4.219} 
we have: 
\begin{eqnarray}
&&\left. (t_*\circ(w,\psi)-f)\right|_{S_{0,0}}=
\left.\frac{\partial}{\partial\tau}(t_*\circ(w,\psi)-f)\right|_{S_{0,0}}=0, \nonumber\\
&&\left.\frac{\partial^2}{\partial\tau^2}(t_*\circ(w,\psi)-f)\right|_{S_{0,0}}=-\frac{2k}{3l} 
\label{14.b1}
\end{eqnarray}
Therefore setting:
\begin{equation}
\Kb=\inf_{S_{0,0}}\left(-\frac{k}{6l}\right)
\label{14.b2}
\end{equation}
we have, restricting suitably $\delta_0$, 
\begin{equation}
\tau^{-2}\left(t_*\circ(w_N,\psi_N)-f_N\right)\geq 2\Kb \ : \ \mbox{in ${\cal K}^{\delta_0}$} 
\label{14.b3}
\end{equation}

\vspace{2.5mm}

\noindent{\large{\bf Fundamental Boundary Bootstrap Assumptions:}}

\vspace{2.5mm}

\noindent On ${\cal K}^{\delta}$ we have:
\begin{eqnarray*}
&&\Bb\leq\frac{r}{\tau}\leq\oB \\
&&|j|\leq B
\end{eqnarray*}
and:
$$\tau^{-2}\left(t_*\circ(w,\psi)-f\right)\geq \Kb$$

\vspace{2.5mm}

In Chapter 10 we encountered the assumptions \ref{10.355}, \ref{10.439}, \ref{10.440}, \ref{10.546}, 
\ref{10.610}, \ref{10.701}. These assumptions where used in Chapter 10 in the derivation of the 
top order acoustical estimates and where also used in Chapter 12 in the derivation of the next to 
to top order acoustical estimates. However to control all the lower order terms appearing, which 
we have called ignorable, we must introduce additional bootstrap assumptions. Let us denote by 
$n_*$ the integral part of $n/2$: 
\begin{equation}
n_*=\left\{\begin{array}{lll} n/2 &:&\mbox{if $n$ is even}\\
(n-1)/2 &:&\mbox{if $n$ is odd}\end{array}\right. 
\ \mbox{, thus} \ \ \ 
n=\left\{\begin{array}{lll}2n_* &:&\mbox{if $n$ is even}\\
2n_*+1 &:&\mbox{if $n$ is odd}\end{array}\right.
\label{14.a18}
\end{equation}
Let us set:
\begin{eqnarray}
&&\overline{C}=2\max_{i+j\leq n_*+2}\max_{\mu}\sup_{{\cal R}_{\delta_0,\delta_0}}
|E_N^i T^j\beta_{\mu,N}|\nonumber\\
&&C^{\Lb}=2\max_{i+j\leq n_*+1}\max_{\mu}\sup_{{\cal R}_{\delta_0,\delta_0}}
|E_N^i T^j\Lb_N\beta_{\mu,N}|\nonumber\\
&&C^L=2\max_{i\leq n_*+1}\max_{\mu}\sup_{{\cal R}_{\delta_0,\delta_0}}\ub^{-1}|E_N^i L_N\beta_{\mu,N}|
\nonumber\\
&&C_*=\max\{\overline{C},C^{\Lb},C^L\} \label{14.a19}
\end{eqnarray}

\vspace{2.5mm}

\noindent{\large{\bf Bootstrap Assumptions on the $\beta_\mu$:}}

\vspace{2.5mm}

\noindent On ${\cal R}_{\delta,\delta}$ we have:
\begin{eqnarray*}
&&|E^i T^j\beta_\mu|\leq C_* \ : \ \mbox{for all $i+j\leq n_*+2$}\\
&&|E^i T^j\Lb\beta_\mu|\leq C_* \ : \ \mbox{for all $i+j\leq n_*+1$}\\
&&|E^i T^j L\beta_\mu|\leq \ub C_* \ : \ \mbox{for all $i+j\leq n_*+1$}
\end{eqnarray*}

\vspace{2.5mm}

Let us set:
\begin{eqnarray}
&&C^{\chi}=2\max_{i\leq n_*}\sup_{{\cal R}_{\delta_0,\delta_0}}|E_N^i\tchi_N|\nonumber\\
&&C^{\chib}=2\max_{i\leq n_*}\sup_{{\cal R}_{\delta_0,\delta_0}}|E_N^i\tchib_N|\nonumber\\
&&C^{\lambda}=2\max\{\max_{i\leq n_*+1}\sup_{{\cal R}_{\delta_0,\delta_0}}u^{-2}|E_N^i\lambda_N|, 
\max_{i\leq n_*}\sup_{{\cal R}_{\delta_0,\delta_0}}u^{-1}|E_N^i T\lambda_N|, \nonumber\\
&&\hspace{55mm}\max_{i+j\leq n_*-1}\sup_{{\cal R}_{\delta_0,\delta_0}}|E_N^i T^{j+2}\lambda_N|\}\nonumber\\
&&C^{\lambdab}=2\max\{\max_{i\leq n_*+1}\sup_{{\cal R}_{\delta_0,\delta_0}}\ub^{-1}|E_N^i\lambdab_N|, 
\max_{i+j\leq n_*}\sup_{{\cal R}_{\delta_0,\delta_0}}|E_N^i T^{j+1}\lambdab_N|\} \nonumber\\
&&C^\prime_*=\max\{C^{\chi},C^{\chib},C^{\lambda},C^{\lambdab}\}\label{14.a20}
\end{eqnarray}
(see \ref{9.a13}, \ref{9.a21}). 

\vspace{2.5mm}

\noindent{\large{\bf Acoustical Bootstrap Assumptions:}}

\vspace{2.5mm}

\noindent On ${\cal R}_{\delta,\delta}$ we have:
\begin{eqnarray*}
&&|E^i\tchi|\leq C^\prime_* \ : \ \mbox{for all $i\leq n_*$}\\
&&|E^i\tchib|\leq C^\prime_* \ : \ \mbox{for all $i\leq n_*$}\\
&&|E^i\lambda|\leq u^2 C^\prime_* \ : \ \mbox{for all $i\leq n_*+1$}\\
&&|E^i T\lambda|\leq u C^\prime_* \ : \ \mbox{for all $i\leq n_*$}\\
&&|E^i T^{j+2}\lambda|\leq C^\prime_* \ : \ \mbox{for all $i+j\leq n_*-1$}\\
&&|E^i\lambdab|\leq \ub C^\prime_* \ : \ \mbox{for all $i\leq n_*+1$}\\
&&|E^i T^{j+1}\lambdab|\leq C^\prime_* \ : \ \mbox{for all $i+j\leq n_*$}
\end{eqnarray*}

\vspace{2.5mm}

Furthermore, to control the lower order terms appearing in the boundary conditions which involve the transformation functions (see proof of Proposition 10.9 and of Lemmas 12.7, 12.8, and 12.11) and 
in the differential consequences of the identification equations (see proof of Propositions 12.3, 12.6 
and of Lemma 12.12) bootstrap assumptions in regard to the transformation functions $\hat{f}$, $v$, 
$\gamma$ are needed. Let us set:
\begin{equation}
C^{\prime\prime}_*=2\max_{i+j\leq n_*+1}\max\left\{\sup_{{\cal K}^{\delta_0}}|\Omega^i T^j\hat{f}_N|, 
\sup_{{\cal K}^{\delta_0}}|\Omega^i T^j v_N|, 
\sup_{{\cal K}^{\delta_0}}|\Omega^i T^j\gamma_N|\right\}
\label{14.a21}
\end{equation}

\vspace{2.5mm}

\noindent{\large{\bf Boundary Bootstrap Assumptions:}}

\vspace{2.5mm}

\noindent On ${\cal K}^{\delta}$ we have:
\begin{eqnarray*}
&&|\Omega^i T^j\hat{f}|\leq C^{\prime\prime}_* \ : \ \mbox{for all $i+j\leq n_*+1$}\\
&&|\Omega^i T^j v|\leq C^{\prime\prime}_* \ : \ \mbox{for all $i+j\leq n_*+1$}\\
&&|\Omega^i T^j\gamma|\leq C^{\prime\prime}_* \ : \ \mbox{for all $i+j\leq n_*+1$}
\end{eqnarray*}

\vspace{2.5mm}

We shall show that the above bootstrap assumptions suffice to control all the lower order terms,  
in particular those appearing in the top order source differences $\s^{(m,n-m)}\check{\tilde{\rho}}_\mu 
\ : \ m=0,...,n$ which enter the top order energy estimates (see Section 9.6). 

We first remark that by equations \ref{3.a15}, \ref{3.a31}, \ref{3.a33}, \ref{3.a38} 
$E$ and $T$ applied to the frame field components $E^\mu$, $N^\mu$, $\Nb^\mu$ give:
\begin{eqnarray}
&&EE^\mu=\se E^\mu+eN^\mu+\oe\Nb^\mu \nonumber\\
&&TE^\mu=(\sf+\sfb)E^\mu+(f+\fb)N^\mu+(\of+\ofb)\Nb^\mu \nonumber\\
&&EN^\mu=\sk E^\mu+k N^\mu+\ok\Nb^\mu \nonumber\\
&&TN^\mu=(\sm+\sn)E^\mu+(m+n)N^\mu+(\om+\on)\Nb^\mu \nonumber\\
&&E\Nb^\mu=\skb E^\mu+\kb N^\mu+\okb\Nb^\mu \nonumber\\
&&T\Nb^\mu=(\smb+\snb)E^\mu+(\mb+\nb)N^\mu+(\omb+\onb)\Nb^\mu 
\label{14.a22}
\end{eqnarray} 
where the coefficients depend on the 1st order quantities $\tchi, \tchib, E\lambda, E\lambdab$ 
and $E\beta_\mu, L\beta_\mu, \Lb\beta_\mu$ as well as on the 0th order quantities $c, \lambda,\lambdab$, 
the frame field components themselves, and the $\beta_\mu$. Also, by equations \ref{3.90}, \ref{3.91} 
and \ref{3.a23}, 
\begin{eqnarray}
&&Ec=-\frac{1}{2}\left(\beta_N\beta_{\Nb}EH+H(\beta_N\ss_{\Nb}+\beta_{\Nb}\ss_N)\right)+c(k+\okb)
\nonumber\\
&&Tc=-\frac{1}{2}\left\{\beta_N\beta_{\Nb}TH+H\left((\beta_N\lambdab+\beta_{\Nb}\lambda)\sss
+\beta_{\Nb}s_{NL}+\beta_Ns_{\Nb\Lb}\right)\right\} \nonumber\\
&&\hspace{9mm}+c(m+n+\omb+\onb) \label{14.a23}
\end{eqnarray}
are likewise expressed. Now the acoustical bootstrap assumptions directly control $E^i\tchi, E^i\tchib$ for 
$i\leq n_*$ and $E^i T^j\lambda, E^i T^j\lambdab$ for $i+j\leq n_*+1$. The acoustical estimates 
themselves of the next to the top order directly control $E^{n-1}\tchi$, $E^{n-1}\tchib$ and 
$E^{n-m}T^m\lambda, E^{n-m}T^m\lambdab \ : \ m=0,...,n$. As we shall see in Sections 14.5 and 14.8, 
lower order acoustical estimates are deduced which directly control $E^i\tchi, E^i\tchib$ for 
$i\leq n-2$ and $E^i T^j\lambda, E^i T^j\lambdab$ for $i+j\leq n-1$. Thus the quantities 
not directly controlled are $E^i T^j\tchi, E^i T^j\tchib$ for $j\geq 1$ and $i+j\leq n_*$ in regard 
to the acoustical bootstrap assumptions, $i+j\leq n-1$ in regard to the acoustical estimates of up to 
the next to the top order. Applying though  $E^i T^{j-1}$ to the formula \ref{14.227} of Section 14.5 
below and to its conjugate, and using \ref{14.a22}, \ref{14.a23}, recalling also that by the first two of the commutation relations 
\ref{3.a14}
$$[E,T]=(\chi+\chib)E,$$
we arrive at recursion formulas for $E^i T^j\tchi$, $E^i T^j\tchib$ of the form:
\begin{eqnarray}
&&E^i T^j \tchi=2\pi\lambda E^{i+1}T^{j-1}\tchi+\sum_{k=0}^{j-1}\sum_{l=0}^i(a_{k,l}E^l T^k\tchi
+\overline{a}_{k,l}E^l T^k\tchib)+b_{j,i}\nonumber\\
&&E^i T^j \tchib=2\pi\lambdab E^{i+1}T^{j-1}\tchib+\sum_{k=0}^{j-1}\sum_{l=0}^i (\underline{a}_{k,l}
E^l T^k\tchi+\overline{\underline{a}}_{k,l} E^l T^k\tchib)+\underline{b}_{j,i} \nonumber\\
&&\label{14.a24}
\end{eqnarray}
where the coefficients $a_{k,l}$, $\overline{a}_{k,l}$, $\underline{a}_{k,l}$, 
$\overline{\underline{a}}_{k,l}$, of order $i+j-k-l$, and the inhomogeneous terms 
$b_{j,i}$, $\overline{b}_{j,i}$, of order $i+j$, depend only on quantities up to the respective order, 
which are under direct control. The last include $E^l T^k\beta_\mu$ for $k+l\leq i+j+1$, 
and $E^l T^k\Lb\beta_\mu$, $E^l T^k L\beta_\mu$ for $k+l\leq i+j$, the third being expressible in 
terms of the previous two since $L=T-\Lb$. In this way also $E^i T^j\tchi$, $E^i T^j\tchib$ are 
indirectly controlled, by the bootstrap assumptions for $i+j\leq n_*$ and by the estimates of 
Sections 14.5 and 14.8 for $i+j\leq n-2$ and by the next to the top order acoustical estimates and 
the estimates of Section 14.4 for $i+j=n-1$. 

In showing how the lower order terms are treated we shall refer to the notion of {\em total order}. 
The total order of a term is simply the sum of the orders of its factors. The total order of each term 
in a sum appearing in an equation satisfied by the actual solution is the same. Therefore each side 
of the equation has the same total order and thus the notion of total order can be assigned to 
the equation. The total order of the propagation equations for $E^{n-1}\tchi$ and $E^{n-1}\tchib$, 
\ref{12.27} and \ref{12.30} respectively, is $n+1$. The total order of the propagation equations 
for $E^n\lambda$ and $E^n\lambdab$, \ref{12.99} and \ref{12.102}, is also $n+1$, 
and so is the total order of the propagation equations for $E^{n-m}T^m\lambda$ and $E^{n-m}T^m\lambdab$ 
for $m\geq 1$, equation \ref{12.112} and its conjugate. In these equations the order of the right 
hand side is also $n+1$. On the other hand, in the case of the propagation equations for 
$\theta_n$ and $\thetab_n$, equations \ref{10.249} and \ref{10.252} respectively 
with $n$ in the role of $l$, 
while the order of the right hand sides is $n+1$, the total order is $n+2$. The same is true of the 
propagation equations for $\nu_{m-1,n-m+1}$ and $\nub_{m-1,n-m+1}$, $m=1,...,n$, equations \ref{10.320} 
and \ref{10.323} with $(m-1,n-m+1)$ in the role of $(m,l)$. 

In regard to the rescaled sources $\s^{(m,n-m)}\tilde{\rho}_\mu$, recall the definition of 
$\s^{(C,m,l)}\sigma_\mu$, $C=E,T$, in the paragraph preceding the recursion formulas \ref{8.28}, \ref{8.29}. Then by \ref{8.42}, \ref{8.45} - \ref{8.47}, and \ref{8.50}, \ref{8.51}, 
$\s^{(C,m,l)}\sigma_\mu$, $C=E,T$, is of order $m+l+2$ and of total order $m+l+3$. Now the formula 
\ref{9.256} for the rescaled sources with $n-m$ in the role of $l$ reads:
\begin{equation}
\s^{(m,n-m)}\tilde{\rho}_\mu=\sum_{j=0}^{m-1}\hat{E}^{n-m}\hat{T}^{m-1-j}\s^{(T,j,0)}\sigma_\mu 
+\sum_{k=0}^{n-m-1}\hat{E}^{n-m-1-k}\s^{(E,m,k)}\sigma_\mu 
\label{14.a25}
\end{equation}
We see then that each term in each of the two sums is of order $n+1$ and of total order $n+2$. 

Consider a lower order term in $\s^{(m,n-m)}\tilde{\rho}_\mu$ or in the propagation equation 
for one of the top order acoustical quantities $\theta_n$, $\thetab_n$, 
$\nu_{m-1,n-m+1}$, $\nub_{m-1,n-m+1}$. Such a term is of total order 
$n+2$ and order at most $n$. It is then expressed as the product of at most $n+2$ factors each of which 
is of order at least 1 and at most $n+1$. Now there can be at most one factor of order exceeding 
$n_*+1$ otherwise the total order would be at least $2(n_*+2)\geq n+3$. When there is such a factor 
present, this factor is to be estimated either in $L^2(S_{\ub,u})$, as in the case where 
the factor is of order $n-1$ or less or of order $n$ but with vanishing P.A. part, 
or after integration with respect to either $\ub$ or $u$ 
in $L^2(C_u^{\ub})$ or $L^2(\Cb_{\ub}^u)$ as appropriate when the factor in question 
is an acoustical quantity of order $n$. The remaining factors are of order $n_*+1$ or less hence 
can all be estimated pointwise by appealing to the bootstrap assumptions. Denoting by $h$ the 
factor of highest order and by $l_1,...,l_k$ the remaining factors, the term is:
$$l_1\cdot\cdot\cdot l_k h$$
The corresponding term in the associated difference source or propagation equation for the 
associated top order acoustical difference quantity is:
$$l_1\cdot\cdot\cdot l_k h-l_{1,N}\cdot\cdot\cdot l_{k,N}h_N$$
which is expressed in the form:
\begin{eqnarray*}
&&l_1\cdot\cdot\cdot l_k (h-h_N)+(l_1\cdot\cdot\cdot l_k-l_{1,N}\cdot\cdot\cdot l_{k,N})h_N = \\
&&l_1\cdot\cdot\cdot l_k (h-h_N)+
\left(\sum_{i=1}^k l_{1,N}\cdot\cdot\cdot l_{i-1,N}(l_i-l_{i,N})l_{i+1}\cdot\cdot\cdot l_k\right)h_N 
\end{eqnarray*}
In the first term,  the factor $h-h_N$ is estimated as discussed above and the other factors 
are estimated pointwise. In the $i$th term in the sum, the difference $l_i-l_{i,N}$ is estimated 
in $L^2(S_{\ub,u})$ and the other factors are estimated pointwise. When there is no factor of order 
exceeding $n_*+1$, that is the term in question is of order $n_*+1$ or less, it can be treated 
in the same way as the above sum. The lower order terms in the propagation equations for 
$E^{n-1}\tchi$, $E^{n-1}\tchib$, $E^{n-m}T^m\lambda$, $E^{n-m}T^m\lambdab$ are treated in a similar 
manner. 

The lower order terms appearing in the $n+1$th order differential consequences of the identification
equations involved in estimating 
$$(\Omega^{n+1-m}T^m\chf, \Omega^{n+1-m}T^m\cv, \Omega^{n+1-m}T^m\cga) \ : \ m=0,...,n+1$$ 
(see proof of Propositions 12.3, 12.6 and of Lemma 12.12) are treated in a similar manner, 
appealing to the boundary bootstrap assumptions. So are the lower order terms appearing in 
the boundary conditions (see proof of Proposition 10.9 and of Lemmas 12.7, 12.8, and 12.11) and 
in the estimation of $\s^{(V;m,l)}\check{b}$ in Section 13.1.

\vspace{5mm}

\section{$L^2(S_{\ub,u})$ Estimates for $\s^{(n-1)}\ctchi$ and for $\s^{(m,n-m)}\cla \ : \ 
m=0,...,n-1$}

Our aim in Section 14.4 is to derive $L^2(S_{\ub,u})$ estimates for the $n$th order 
variation differences \ref{9.236}. This will require $L^2(C_u^{\ub})$ estimates for their 
$L$ derivatives in terms of the top order energies. In these estimates the $n$th order acoustical 
quantities $\s^{(n-1)}\ctchi$ and $\s^{(m,n-m)}\cla \ : \ m=0,...,n-1$ enter. Now, these have all been 
estimated in $L^2(\Cb_{\ub}^u)$, see \ref{12.411}, \ref{12.412} and Proposition 12.2, and \ref{12.665} 
and Proposition 12.5. However what is needed in the estimates just mentioned is $L^2(C_u^{\ub})$ 
estimates for these quantities. We shall in fact presently derive $L^2(S_{\ub,u})$ estimates for 
these quantities, which, although weaker in respect to powers of growth in $\ub$ and $u$, suffice 
for the purpose establishing the required estimates for the $L$ derivatives of the $n$th order variation 
differences. 

For simplicity of notation we denote in the following:
\begin{equation}
\max_{m=0,...,n}D_{m,n-m,N}(u_1)=D
\label{14.1}
\end{equation}
it been understood that for estimates in regard to the spacetime region ${\cal R}_{\ub_1,u_1}$ the 
argument of $D$ is $u_1$. 

To derive an estimate for $\s^{(n-1)}\ctchi$ in $L^2(S_{\ub,u})$ we revisit the proof of Lemma 12.2. 
In reference to inequality \ref{12.128}, we must now estimate the integral:
\begin{equation}
\int_0^{\ub_1}\|\s^{(n-1)}\chF\|_{L^2(S_{\ub,u})}d\ub
\label{14.2}
\end{equation}
pointwise relative to $u$. In regard to \ref{12.131}, as previously noted the 1st can be absorbed in the 
inequality, while the 2nd is bounded by:
\begin{eqnarray}
&&C\int_0^{\ub_1}\|\rho\s^{(n-1)}\ctchib\|_{L^2(S_{\ub,u})}d\ub
\leq C\ub_1^{1/2}\|\rho\s^{(n-1)}\ctchib\|_{L^2(C_u^{\ub_1})}\nonumber\\
&&\hspace{10mm}\leq C\ub_1^{a_0+\frac{3}{2}}u^{b_0-1}\s^{(n-1)}\Xb(\ub_1,u) 
\leq C\sqrt{D}\ub_1^{a_0+\frac{3}{2}}u^{b_0-1}\label{14.3}
\end{eqnarray}
by definition \ref{12.460}, Propositions 12.2 and 12.3, and the top order energy estimates. 
We turn to the principal difference terms \ref{12.132}. Consider first the 2nd of these difference 
terms, which gives the dominant contribution through its $\Nb$ component, last of \ref{12.134}.  
To estimate this in $L^2(C_u^{\ub_1})$ we use the fact that
\begin{equation}
N^\mu E(E^n\beta_\mu-E_N^n\beta_{\mu,N})+\ogamma\Nb^\mu E(E^n\beta_\mu-E_N^n\beta_{\mu,N})
=\s^{(Y;0,n)}\cxi
\label{14.4}
\end{equation}
This requires considering first the $N$ component, second of \ref{12.134}. we write: 
\begin{eqnarray}
&&\rho N^\mu E^{n+1}\beta_\mu=L^\mu E^{n+1}\beta_\mu 
=E^n(L^\mu E\beta_\mu)-\sum_{i=1}^n\left(\begin{array}{c} n\\i \end{array}\right) 
E^i(\rho N^\mu)E^{n+1-i}\beta_\mu\nonumber\\
&&\hspace{8mm}=E^n(E^\mu L\beta_\mu)-E^n(\rho N^\mu) E\beta_\mu \ \mbox{: up to ignorable terms}\nonumber\\
&&\hspace{8mm}=E^\mu E^n L\beta_\mu+(E^n E^\mu)L\beta_\mu-E^n(\rho N^\mu) E\beta_\mu \ 
\mbox{: up to ignorable terms}\nonumber\\
&&\label{14.5}
\end{eqnarray}
Using \ref{12.7}, \ref{12.8}, \ref{12.10}, and the fact that by the 1st of the commutation relations 
\ref{3.a14}:
\begin{equation}
[E^n,L]\beta_\mu=\sum_{i=0}^{n-1}E^{n-1-i}(\chi E^{i+1}\beta_\mu)=
\rho\beta_\mu E^{n-1}\tchi \ \mbox{: up to ignorable terms}
\label{14.6}
\end{equation}
we then obtain:
\begin{eqnarray}
&&\rho N^\mu E^{n+1}\beta_\mu=E^\mu LE^n\beta_\mu+\frac{3}{2}\sss\rho E^{n-1}\tchi
+\frac{1}{2c}(s_{NL}-2\pi\lambdab\ss_N)E^{n-1}\tchib\nonumber\\
&&\hspace{25mm}-\frac{\ss_N}{c}E^n\lambdab 
\  \ \ \mbox{: up to ignorable terms} \label{14.7}
\end{eqnarray}
Similarly for the $N$th approximants. Taking the difference then yields: 
\begin{eqnarray}
&&\rho N^\mu E^{n+1}\beta_\mu-\rho_N N_N^\mu E_N^{n+1}\beta_{\mu,N}=\s^{(E;0,n)}\cxi_L\nonumber\\
&&\hspace{20mm}+\frac{3}{2}\sss\rho \s^{(n-1)}\ctchi
+\frac{1}{2c}(s_{NL}-2\pi\lambdab\ss_N)\s^{n-1}\ctchib-\frac{\ss_N}{c}\s^{(0,n)}\lambdab \nonumber\\
&&\hspace{50mm} \mbox{: up to ignorable terms} \label{14.8}
\end{eqnarray}
In deriving an estimate for the contribution of the $\Nb$ component through \ref{14.4}, this must be 
multiplied by $\ogamma^{-1}\sim u^{-1}$. The contribution of the term in $\s^{(n-1)}\ctchi$ to 
\ref{14.2} can be absorbed, while the contribution of the term in $\s^{(n-1)}\ctchib$ is bounded as 
in \ref{14.3} but with an extra $u^{-1}$ factor, that is by:
\begin{equation}
C\sqrt{D}\ub_1^{a_0+\frac{3}{2}}u^{b_0-2}
\label{14.9}
\end{equation}
The contribution of the term in $\s^{(0,n)}\clab$ is bounded by:
\begin{eqnarray}
&&Cu^{-1}\int_0^{\ub_1}\|\s^{(0,n)}\clab\|_{L^2(S_{\ub,u})}d\ub
\leq Cu^{-1}\ub_1^{1/2}\|\s^{(0,n)}\clab\|_{L^2(C_u^{\ub_1})}\nonumber\\
&&\hspace{10mm}\leq C\ub_1^{a_0+\frac{1}{2}}u^{b_0-1}\s^{(0,n)}\Lab(\ub_1,u)
\leq C\sqrt{D}\ub_1^{a_0+\frac{1}{2}}u^{b_0-1}\label{14.10}
\end{eqnarray}
by definition \ref{12.461}, Propositions 12.2 and 12.3, and the top order energy estimates. 
As for the contribution of $\s^{(E;0,n)}\cxi_L$, it is bounded by:
\begin{eqnarray}
&&Cu^{-1}\int_0^{\ub_1}\|\s^{(E;0,n)}\cxi_L\|_{L^2(S_{\ub,u})}d\ub
\leq Cu^{-1}\ub_1^{1/2}\|\s^{(E;0,n)}\cxi_L\|_{L^2(C_u^{\ub_1})}\nonumber\\
&&\hspace{10mm}\leq Cu^{-1}\ub_1^{1/2}\sqrt{\s^{(E;0,n)}\cE^{\ub_1}(u)}\leq 
C\sqrt{D}\ub_1^{a_0+\frac{1}{2}}u^{b_0-1} \label{14.11}
\end{eqnarray}
What remains is then the contribution of the right hand side of \ref{14.4}, which, as it is to be 
multiplied by $\rho\ogamma^{-1}\sim\ub^{1/2}u^{-2}a^{1/2}$ is bounded by:
\begin{eqnarray}
&&Cu^{-2}\int_0^{\ub_1}\ub^{1/2}\|a^{1/2}\s^{(Y;0,n)}\csxi\|_{L^2(S_{\ub,u})}d\ub\nonumber\\
&&\leq C u^{-2}\left(\int_0^{\ub_1}\ub d\ub\right)^{1/2}\|a^{1/2}\s^{(Y;0,n)}\csxi\|_{L^2(C_u^{\ub_1})}
\nonumber\\
&&\leq Cu^{-2}\ub_1\sqrt{\s^{(Y;0,n)}\cE^{\ub_1}(u)}\leq C\sqrt{D}\ub_1^{a_0+1}u^{b_0-2}\label{14.12}
\end{eqnarray}
This is the dominant contribution to the integral \ref{14.2}. The contribution of the $E$ component 
in \ref{12.134} is bounded by, since $\rho\sim\ub^{1/2}u^{-1}a^{1/2}$, 
\begin{eqnarray}
&&Cu^{-1}\int_0^{\ub_1}\ub^{1/2}\|a^{1/2}\s^{(E;0,n)}\csxi\|_{L^2(S_{\ub,u})}d\ub\nonumber\\
&&\leq Cu^{-1}\left(\int_0^{\ub_1}\ub d\ub\right)^{1/2}\|a^{1/2}\s^{(E;0,n)}\csxi\|_{L^2(C_u^{\ub_1})}
\nonumber\\
&&\leq Cu^{-1}\ub_1\sqrt{\s^{(E;0,n)}\cE^{\ub_1}(u)}\leq C\sqrt{D}\ub_1^{a_0+1}u^{b_0-1}\label{14.13}
\end{eqnarray}

As for the  1st of the principal difference terms \ref{12.132}, expanded in \ref{12.164} in its $E$, 
$N$, and $\Nb$ components, to estimate the $\Nb$ component we write: 
\begin{eqnarray}
&&\Nb^\mu E^n L\beta_\mu=E^n(\Nb^\mu L\beta_\mu)-\sum_{i=1}^n\left(\begin{array}{c} n\\i \end{array}
\right)(E^i\Nb^\mu)E^{n-i}L\beta_\mu\nonumber\\
&&\hspace{8mm}=E^n(\lambdab E^\mu E\beta_\mu)-\rho(\ss_N-c\pi\sss)E^{n-1}\tchib 
\ \mbox{: up to ignorable terms}\nonumber\\
&&\hspace{8mm}=\lambdab E^\mu E^{n+1}\beta_\mu 
+\frac{\rho}{2}\ss_{\Nb}E^{n-1}\tchi+\frac{\rho}{2}(2c\pi-\ss_N)E^{n-1}\tchib 
+\sss E^n\lambdab\nonumber\\
&&\hspace{8mm} \ \mbox{: up to ignorable terms} \label{14.14}
\end{eqnarray}
Similarly for the $N$th approximants. Taking the difference then yields:
\begin{eqnarray}
&&\Nb^\mu E^n L\beta_\mu -\Nb_N^\mu E_N^n L_N\beta_{\mu,N}=\lambdab\s^{(E;0,n)}\sxi
\nonumber\\
&&\hspace{20mm}+\frac{\rho}{2}\ss_{\Nb}\s^{(n-1)}\ctchi+\frac{\rho}{2}(2c\pi-\ss_N)\s^{(n-1)}\ctchib 
+\sss \s^{(0,n)}\clab\nonumber\\
&&\hspace{46mm} \ \mbox{: up to ignorable terms} \label{14.15}
\end{eqnarray}
Here, the contribution to \ref{14.2} of the term in $\s^{(n-1)}\ctchib$ is bounded as in \ref{14.3}, 
that of the term in $\s^{(n-1)}\clab$ as in \ref{14.10} with an extra $u$ factor, while that of the term 
in $\s^{(n-1)}\ctchi$ is absorbed. On the other hand, the contribution of the 1st term on the right 
is bounded as in \ref{14.13}. To estimate the contributions of the $E$ and $N$ components in 
\ref{12.164} we write
$$E^n L\beta_\mu=LE^n\beta_\mu+\rho(E^{(n-1)}\tchi E\beta_\mu \ \mbox{: up to ignorable terms}$$
and similarly for the $N$th approximants, hence taking the difference we obtain:
\begin{equation}
E^n L\beta_\mu-E_N^n L_N\beta_{\mu,N}=L(E^n\beta_\mu-E_N^n\beta_{\mu,N})
+\rho\s^{(n-1)}\ctchi E\beta_\mu  \ \mbox{: up to ignorable terms} 
\label{14.16}
\end{equation}
Here the contribution of the second term is absorbed. The first term contributes 
\begin{equation}
\s^{(E;0,n)}\cxi_L \ \mbox{: in the case of the $E$ component}
\label{14.17}
\end{equation}
and 
\begin{equation}
N^\mu L(E^n\beta_\mu-E_N^n\beta_{\mu,N}) \ \mbox{: in the case of the $N$ component}
\label{14.18}
\end{equation}
The contribution to \ref{14.2} of \ref{14.17} is bounded as in \ref{14.11} but with an extra $u$ 
factor. We estimate the contribution of \ref{14.18} using
\begin{equation}
N^\mu L(E^n\beta_\mu-E_N^n\beta_{\mu,N})+\ogamma\Nb^\mu L(E^n\beta_\mu -E_N^n\beta_{\mu,N})
=\s^{(Y;0,n)}\cxi_L
\label{14.19}
\end{equation}
The contribution of the right hand side is similar to that of \ref{14.17}, while the contribution of 
the 2nd term on the left is bounded by $Cu$ times the bound for the contribution of \ref{14.15}. 

The above results lead to the conclusion that $\s^{(n-1)}\ctchi$ satisfies the $L^2(S_{\ub,u})$ 
estimate:
\begin{equation}
\|\s^{(n-1)}\ctchi\|_{L^2(S_{\ub,u})}\leq C\sqrt{D}\ub^{a_0+\frac{1}{2}}u^{b_0-\frac{3}{2}}
\label{14.20}
\end{equation}

To derive an estimate for $\s^{(0,n)}\cla$ in $L^2(S_{\ub,u})$ we revisit the proof of Lemma 12.3. 
In reference to inequality \ref{12.196}, we must now estimate the integral: 
\begin{equation}
\int_0^{\ub_1}\|\s^{(0,n)}\chG\|_{L^2(S_{\ub,u})}d\ub
\label{14.21}
\end{equation}
pointwise relative to $u$. In regard to \ref{12.199}, as previously noted the 1st can be absorbed in the 
inequality, while the 2nd is estimated by \ref{14.10} with an extra $u$ factor and the 4th is estimated 
by \ref{14.3}. As for the 3rd, by the estimate \ref{14.20} it is bounded by: 
\begin{equation}
C\sqrt{D}\int_0^{\ub_1}\ub^{a_0+\frac{3}{2}}u^{b_0-\frac{3}{2}}d\ub
=\frac{C\sqrt{D}}{\left(a_0+\frac{5}{2}\right)}\ub_1^{a_0+\frac{5}{2}}u^{b_0-\frac{3}{2}}
\label{14.22}
\end{equation}
We turn to the principal difference terms \ref{12.200}. The dominant contribution comes from the 
1st difference term, which is expanded in \ref{12.202} in its $E$, $N$, and $\Nb$ components. The 
dominant contribution is that of the $\Nb$ component. This can only be estimated in terms of 
the $(0,n)$ energies through \ref{12.203}. However this involves $\s^{(Y;0,n)}\cxi_{\Lb}$ 
which does not have integrability on the $C_u^{\ub}$. Thus the only way to proceed is to replace 
the 1st difference term in \ref{12.200} by:
\begin{equation}
\frac{\rho}{4}\beta_N^2 E^n TH-\frac{\rho_N}{4}\beta_{N,N}^2 E_N^n TH_N
-\frac{\rho}{4}\beta_N^2 E^n LH+\frac{\rho_N}{4}\beta_{N,N}^2 E_N^n L_N H_N
\label{14.23}
\end{equation}
These two differences can be expanded in their $E$, $N$, and $\Nb$ components. The dominant 
contribution is that of the $\Nb$ component of the 1st difference, which can be expressed in terms of:
\begin{equation}
\rho\Nb^\mu E^n T\beta_\mu-\rho_N\Nb_N^\mu E_N^n T\beta_{\mu,N}
\label{14.24}
\end{equation}
Here we use the fact that:
\begin{eqnarray}
&&\rho N^\mu E(E^{n-1}T\beta_\mu-E_N^{n-1}T\beta_{\mu,N})
+\rho\ogamma\Nb^\mu E(E^{n-1}T\beta_\mu-E_N^{n-1}T\beta_{\mu,N})\nonumber\\
&&\hspace{50mm}=\rho\s^{(Y;1,n-1)}\csxi \label{14.25}
\end{eqnarray}
The 1st term on the left is what results when we express the $N$ component of the same difference. 
We have:
\begin{eqnarray}
&&\rho N^\mu E^n T\beta_\mu=L^\mu E^n T\beta_\mu 
=L^\mu\left\{E^{n-1}TE\beta_\mu+E^{n-1}\left((\chi+\chib)E\beta_\mu\right)\right\}\nonumber\\
&&=L^\mu E^{n-1}TE\beta_\mu+\rho\ss_N(\rho E^{n-1}\tchi+\rhob E^{n-1}\tchib) \ \mbox{: up to ignorable 
terms} \nonumber\\
\label{14.26}
\end{eqnarray}
Here, the 1st term on the right is:
\begin{equation}
E^{n-1}T(L^\mu E\beta_\mu)-(E\beta_\mu)E^{n-1}T(\rho N^\mu)  \ \mbox{: up to ignorable terms} 
\label{14.27}
\end{equation}
and we have:
\begin{eqnarray}
&&E^{n-1}T(L^\mu E\beta_\mu)=E^{n-1}T(E^\mu L\beta_\mu)\nonumber\\
&&=E^\mu E^{n-1}TL\beta_\mu+(L\beta_\mu)E^{n-1}TE^\mu  \ \mbox{: up to ignorable terms} 
\label{14.28}
\end{eqnarray}
and:
\begin{eqnarray}
&&E^\mu E^{n-1}TL\beta_\mu=E^\mu E^{n-1}LT\beta_\mu+\sss E^{n-1}\zeta \ \mbox{: up to ignorable terms} 
\nonumber\\
&&\hspace{22mm}=E^\mu LE^{n-1}T\beta_\mu+\sss E^{n-1}\zeta \ \mbox{: up to ignorable terms} \nonumber\\ 
&&\label{14.29}
\end{eqnarray}
Now, since (see \ref{3.a14}) $\zeta=2(\eta-\etab)$, by \ref{3.122}, \ref{3.a25} and \ref{12.64} we have:
\begin{equation}
[\zeta]_{P.A.}=2\rho\mu-2\rhob\mub 
\label{14.30}
\end{equation}
Using also \ref{12.74} and the fact that by \ref{12.66} and \ref{12.70}:
\begin{equation}
[T(\rho N^\mu)]_{P.A.}=c^{-1}N^\mu T\lambdab+2E^\mu\rho\mu+2\pi N^\mu\rho\mub
\label{14.31}
\end{equation}
we deduce from the above (see \ref{12.107}, \ref{12.108}, \ref{12.110}, \ref{12.111}):
\begin{eqnarray}
&&\rho N^\mu E^n T\beta_\mu=E^\mu LE^{n-1}T\beta_\mu+\sss\mu_{0,n-1}
+(c^{-1}s_{NL}-2\pi\ss_N\rho-2\sss\rhob)\mub_{0,n-1}\nonumber\\
&&\hspace{20mm}+\ss_N\rho^2 E^{n-1}\tchi+\ss_N\rho\rhob E^{n-1}\tchib-c^{-1}\ss_N E^{n-1}T\lambdab\nonumber\\
&&\hspace{20mm} \ \mbox{: up to ignorable terms} \label{14.32}
\end{eqnarray}
A similar formula holds for the $N$th approximants. Taking the difference we conclude that 
the 1st term on the left in \ref{14.25} is (see \ref{12.117}):
\begin{eqnarray}
&&\rho N^\mu E(E^{n-1} T\beta_\mu-E_N^{n-1}T\beta_{\mu,N})=\s^{(E;1,n-1)}\cxi_L\nonumber\\
&&\hspace{20mm}+\sss\cmu_{0,n-1}+(c^{-1}s_{NL}-2\pi\ss_N\rho-2\sss\rhob)\cmub_{0,n-1}\nonumber\\
&&\hspace{20mm}+\ss_N\rho^2 \s^{(n-1)}\ctchi+\ss_N\rho\rhob E^{(n-1)}\ctchib-c^{-1}\ss_N 
\s^{(1,n-1)}\clab\nonumber\\
&&\hspace{20mm} \ \mbox{: up to ignorable terms} \label{14.33}
\end{eqnarray}
We shall presently estimate, in connection with \ref{14.21}, the integral:
\begin{equation}
\int_0^{\ub_1}\|\rho N^\mu E(E^{n-1}T\beta_\mu-E_N^{n-1}T\beta_{\mu,N})\|_{L^2(S_{\ub,u})}d\ub
\label{14.34}
\end{equation}
On the right hand side of \ref{14.33} we have 6 terms. The contribution to \ref{14.34} of the 1st term 
is bounded by:
\begin{equation}
\ub_1^{1/2}\|\s^{(E;1,n-1)}\cxi_L\|_{L^2(C_u^{\ub_1})}\leq C\ub_1^{1/2}
\sqrt{\s^{(E;1,n-1)}\cE^{\ub_1}(u)}\leq C\sqrt{D}\ub_1^{a_1+\frac{1}{2}}u^{b_1}
\label{14.35}
\end{equation}
Recalling from \ref{12.118} that
$$\cmu_{0,n-1}=\s^{(0,n)}\cla+\pi\lambda\s^{(n-1)}\ctchi  \ \mbox{: up to ignorable terms}$$
the partial contribution of $\s^{(0,n)}\cla$ through the 2nd term can be absorbed 
while by \ref{14.20} the contribution of the remainder is bounded by:
\begin{equation}
C\sqrt{D}\ub_1^{a_0+\frac{5}{2}}u^{b_0+\frac{1}{2}}
\label{14.36}
\end{equation}
Recalling from \ref{12.119} that
$$\cmub_{0,n-1}=\s^{(0,n)}\clab+\pi\lambdab\s^{(n-1)}\ctchib \ \mbox{: up to ignorable terms}$$
the contribution to \ref{14.34} of the 3rd term on the right in \ref{14.33} is by \ref{14.3} and 
\ref{14.10} bounded by:
\begin{equation}
C\sqrt{D}\ub_1^{a_0+\frac{1}{2}}u^{b_0+1}
\label{14.37}
\end{equation}
By \ref{14.20} the contribution of the 4th term is bounded by:
\begin{equation}
C\sqrt{D}\ub_1^{a_0+\frac{7}{2}}u^{b_0-\frac{3}{2}}
\label{14.38}
\end{equation}
while by \ref{14.3} the contribution of the 5th term is bounded by:
\begin{equation}
C\sqrt{D}\ub_1^{a_0+\frac{3}{2}}u^{b_0+1}
\label{14.39}
\end{equation}
Finally, the contribution to \ref{14.34} of the 6th term on the right in \ref{14.33} is bounded by:
\begin{equation}
\ub_1^{1/2}\|\s^{(1,n-1)}\clab\|_{L^2(C_u^{\ub_1})}\leq\ub_1^{a_1+\frac{1}{2}}u^{b_1}
\s^{(1,n-1)}\Lab(\ub_1,u)\leq C\sqrt{D}\ub_1^{a_1+\frac{1}{2}}u^{b_1}
\label{14.40}
\end{equation}
by definition \ref{12.666}, Propositions 12.5 and 12.6, and the top order energy estimates. 
Combining the above results we conclude that \ref{14.34} is bounded by:
\begin{equation}
C\sqrt{D}\ub_1^{a_1+\frac{1}{2}}u^{b_1}
\label{14.41}
\end{equation}
up to terms which can be absorbed. 

Now, the contribution of \ref{14.24} to the integral \ref{14.21} is bounded in terms of the integral:
\begin{equation}
\int_0^{\ub_1}\|\rho\Nb^\mu E(E^{n-1}T\beta_\mu-E_N^{n-1}T\beta_{\mu,N})\|_{L^2(S_{\ub,u})}d\ub
\label{14.42}
\end{equation}
Since $\ogamma^{-1}\sim u^{-1}$ the contribution of the 1st term on the left in \ref{14.25} to this 
integral is bounded by $Cu^{-1}$ times the integral \ref{14.34}, that is by:
\begin{equation}
C\sqrt{D}\ub_1^{a_1+\frac{1}{2}}u^{b_1-1}
\label{14.43}
\end{equation}
On the other hand, the contribution of the right hand side of \ref{14.25} to \ref{14.42} is,  
in view of the fact that $\rho\ogamma^{-1}\sim\ub^{1/2}u^{-2}a^{1/2}$, bounded by: 
\begin{eqnarray}
&&Cu^{-2}\int_0^{\ub_1}\ub^{1/2}\|a^{1/2}\s^{(Y;1,n-1)}\csxi\|_{L^2(S_{\ub,u})}d\ub\nonumber\\
&&\leq C u^{-2}\left(\int_0^{\ub_1}\ub d\ub\right)^{1/2}\|a^{1/2}\s^{(Y;1,n-1)}\csxi\|_{L^2(C_u^{\ub_1})}
\nonumber\\
&&\leq Cu^{-2}\ub_1\sqrt{\s^{(Y;1,n-1)}\cE^{\ub_1}(u)}\leq C\sqrt{D}\ub_1^{a_1+1}u^{b_1-2}
\label{14.44}
\end{eqnarray}
This is the dominant contribution to the integral \ref{14.21}. 

Let us also consider the leading, $\Nb$, component of the 2nd of the differences \ref{14.23}. 
This can be expressed in terms of: 
\begin{equation}
\rho\Nb^\mu E^n L\beta_\mu-\rho_N\Nb_N^\mu E_N^n L_N\beta_{\mu,N}
\label{14.45}
\end{equation}
To estimate the contribution of this to \ref{14.21} we write:
\begin{eqnarray}
&&\Nb^\mu E^n L\beta_\mu=E^n(\Nb^\mu L\beta_\mu)-(L\beta_\mu)E^n\Nb^\mu \ \mbox{: up to ignorable terms}
\nonumber\\
&&\hspace{10mm}=E^n(\lambdab E^\mu E\beta_\mu)-(\ss_N-c\pi\sss)\rho E^{n-1}\tchib 
\ \mbox{up to ignorable terms} \nonumber\\
&&\hspace{10mm}=\lambdab E^\mu E^{n+1}\beta_\mu+\sss E^n\lambdab+\frac{1}{2}\ss_{\Nb}\rho E^{n-1}\tchi
+\frac{1}{2}(2c\pi-\ss_N)\rho E^{n-1}\tchib\nonumber\\
&&\hspace{10mm} \ \mbox{: up to ignorable terms} \label{14.46}
\end{eqnarray}
Similarly for the $N$th approximants. Taking the difference we then obtain that \ref{14.45} is given by:
\begin{eqnarray}
&&c\rho^2\s^{(E;0,n)}\csxi+\sss\rho\s^{(0,n)}\clab
+\frac{1}{2}\ss_{\Nb}\rho^2\s^{(n-1)}\ctchi+\frac{1}{2}(2c\pi-\ss_N)\rho^2 \s^{(n-1)}\ctchib\nonumber\\
&&\hspace{10mm} \ \mbox{: up to ignorable terms} \label{14.47}
\end{eqnarray}
It follows that the contribution of \ref{14.45} to \ref{14.21} is bounded by:
\begin{eqnarray}
&&Cu^{-1}\int_0^{\ub_1}\ub^{3/2}\|a^{1/2}\s^{(E;0,n)}\sxi\|_{L^2(S_{\ub,u})}d\ub
+C\int_0^{\ub_1}\ub\|\s^{(0,n)}\clab\|_{L^2(S_{\ub,u})}d\ub\nonumber\\
&&+C\int_0^{\ub_1}\ub^2\|\s^{(n-1)}\ctchi\|_{L^2(S_{\ub,u})}d\ub
+C\int_0^{\ub_1}\ub\|\rho\s^{(n-1)}\ctchib\|_{L^2(S_{\ub,u})}d\ub\nonumber\\
&&\leq C\sqrt{D}u_1^{a_0+\frac{3}{2}}u^{b_0-\frac{1}{2}} \label{14.48}
\end{eqnarray}
This is dominated by \ref{14.44} by a factor of $\ub_1^{1/2}u^{3/2}$ at least. 

The above results lead to the conclusion that $\s^{(0,n)}\cla$ satisfies the $L^2(S_{\ub,u})$ 
estimate:
\begin{equation}
\|\s^{(0,n)}\cla\|_{L^2(S_{\ub,u})}\leq C\sqrt{D}\ub^{a_1+\frac{1}{2}}u^{b_1-\frac{3}{2}}
\label{14.49}
\end{equation}

We proceed to prove by induction that 
\begin{equation}
\|\s^{(m,n-m)}\cla\|_{L^2(S_{\ub,u})}\leq C\sqrt{D}\ub^{a_{m+1}+\frac{1}{2}}u^{b_{m+1}-\frac{3}{2}} 
\ \mbox{: for $m=0,...,n-1$}
\label{14.50}
\end{equation}
the case $m=0$ being the estimate \ref{14.49}. The inductive hypothesis is then that the estimate holds 
with $m$ replaced by $j=0,...,m-1$. To derive an estimate for $\s^{(m,n-m)}\cla$ in $L^2(S_{\ub,u})$ 
we revisit the proof of Lemma 12.9. In reference to inequality \ref{12.564} we must now estimate the 
integral 
\begin{equation}
\int_0^{\ub_1}\|\s^{(m,n-m)}\chG\|_{L^2(S_{\ub,u})}d\ub 
\label{14.51}
\end{equation}
pointwise relative to $u$. In regard to \ref{12.566}, as previously noted the 1st can be absorbed in the 
inequality, while the 2nd is bounded by:
\begin{eqnarray}
&&\ub_1^{1/2}\|\s^{(m,n-m)}\clab\|_{(C_u^{\ub_1})}\leq 
\ub_1^{a_m+\frac{1}{2}}u^{b_m}\s^{(m,n-m)}\Lab(\ub_1,u)\nonumber\\
&&\hspace{33mm}\leq C\sqrt{D}\ub_1^{a_m+\frac{1}{2}}u^{b_m}
\label{14.52}
\end{eqnarray}
by definition \ref{12.666}, Propositions 12.4, 12.5 and 12.6, and the top order energy estimates. 
By the inductive hypothesis the 3rd of \ref{12.566} is bounded by:
\begin{equation}
C\sqrt{D}\int_0^{\ub_1}\ub^{a_m+\frac{3}{2}}u^{b_m-\frac{3}{2}}d\ub
=\frac{C\sqrt{D}}{\left(a_m+\frac{5}{2}\right)}\ub_1^{a_m+\frac{5}{2}}u^{b_m-\frac{3}{2}}
\label{14.53}
\end{equation}
By \ref{12.567}, the inductive hypothesis and \ref{14.20}, the 4th of \ref{12.566} is bounded by:
\begin{eqnarray}
&&C\left\{\sum_{j=0}^{m-1}u^{2j}\int_0^{\ub-1}\ub\|\s^{(m-1-j,n-m+j+1)}\cla\|_{L^2(S_{\ub,u})}d\ub
\right.\nonumber\\
&&\hspace{45mm}\left.+u^{2m}\int_0^{\ub_1}\ub\|\s^{(n-1)}\ctchi\|_{L^2(S_{\ub,u})}d\ub\right\}
\nonumber\\
&&\leq C\sqrt{D}\left\{\sum_{j=0}^{m-1}u^{b_{m-j}-\frac{3}{2}+2j}
\int_0^{\ub_1}\ub^{a_{m-j}+\frac{3}{2}}d\ub
+u^{b_0-\frac{3}{2}+2m}\int_0^{\ub_1}\ub^{a_0+\frac{3}{2}}d\ub\right\}\nonumber\\
&&\leq\frac{C\sqrt{D}}{\left(a_m+\frac{5}{2}\right)}\ub_1^{a_m+\frac{5}{2}}u^{b_m-\frac{3}{2}}
\label{14.54}
\end{eqnarray}
in view of the non-increasing with $m$ property of the exponents $a_m$, $b_m$. Also, by \ref{12.568} 
the 5th of \ref{12.566} is bounded by:
\begin{eqnarray}
&&C\left\{\sum_{j=0}^{m-1}\ub_1^{j+\frac{3}{2}}\|\s^{(m-1-j,n-m+j+1)}\clab\|_{L^2(C_u^{\ub_1})}
+\ub_1^{m+\frac{1}{2}}\|\rho\s^{(n-1)}\ctchib\|_{C_u^{\ub_1})}\right\}\nonumber\\
&&\leq C\left\{\sum_{j=0}^{m-1}\ub_1^{a_{m-1-j}+j+\frac{3}{2}}u^{b_{m-1-j}}
\s^{(m-1-j,n-m+j+1)}\Lab(\ub_1,u)\right.\nonumber\\
&&\hspace{45mm}\left.+\ub_1^{a_0+m+\frac{3}{2}}u^{b_0-1}\s^{(n-1)}\Xb(\ub_1,u)\right\}\nonumber\\
&&\leq C\sqrt{D}\ub_1^{a_{m-1}+\frac{3}{2}}u^{b_{m-1}}
\label{14.55}
\end{eqnarray}
by definitions \ref{12.460}, \ref{12.666}, Propositions 12.2, 12.3, 12.5, 12.6, and the top order 
energy estimates. 

We turn to the principal difference terms \ref{12.571}. The dominant contribution comes from the 
1st difference term, which is expanded in \ref{12.573} in its $E$, $N$, and $\Nb$ components. The 
dominant contribution is that of the $\Nb$ component. This can only be estimated in terms of 
the $(m,n-m)$ energies through \ref{12.574}. However this involves $\s^{(Y;m,n-m)}\cxi_{\Lb}$ 
which does not have integrability on the $C_u^{\ub}$. Thus the only way to proceed is to replace 
the 1st difference term in \ref{12.571} by:
\begin{eqnarray}
&&\frac{\rho}{4}\beta_N^2 E^{n-m}T^{m+1}H-\frac{\rho_N}{4}\beta_{N,N}^2 E_N^{n-m}T^{m+1}H_N \nonumber\\
&&-\frac{\rho}{4}\beta_N^2 E^{n-m}T^m LH+\frac{\rho_N}{4}\beta_{N,N}^2 E_N^{n-m}T^m L_N H_N
\label{14.56}
\end{eqnarray}
These two differences can be expanded in their $E$, $N$, and $\Nb$ components. The dominant 
contribution is that of the $\Nb$ component of the 1st difference, which can be expressed in terms of:
\begin{equation}
\rho\Nb^\mu E^{n-m}T^{m+1}\beta_\mu-\rho_N\Nb_N^\mu E_N^{n-m}T^{m+1}\beta_{\mu,N}
\label{14.57}
\end{equation}
Here we use the fact that:
\begin{eqnarray}
&&\rho N^\mu E(E^{n-1-m}T^{m+1}\beta_\mu-E_N^{n-1-m}T^{m+1}\beta_{\mu,N}) \nonumber\\
&&\hspace{25mm}+\rho\ogamma\Nb^\mu E(E^{n-1-m}T^{m+1}\beta_\mu-E_N^{n-1-m}T^{m+1}\beta_{\mu,N})
\nonumber\\
&&\hspace{50mm}=\rho\s^{(Y;m+1,n-1-m)}\csxi \label{14.58}
\end{eqnarray}
The 1st term on the left is what results when we express the $N$ component of the same difference. 
We have:
\begin{eqnarray}
&&\rho N^\mu E^{n-m} T^{m+1}\beta_\mu
=L^\mu\left\{E^{n-1-m}T^{m+1}E\beta_\mu+E^{n-1-m}\left(T^m(\chi+\chib)\cdot E\beta_\mu\right)\right\} 
\nonumber\\
&& \ \mbox{: up to ignorable terms}\nonumber\\
&&=L^\mu E^{n-1-m}T^{m+1}E\beta_\mu+2\rho\ss_N(\rho\mu_{m-1,n-m}+\rhob\mub_{m-1,n-m}) \nonumber\\
&& \ \mbox{: up to ignorable terms} 
\label{14.59}
\end{eqnarray}
Here, the 1st term on the right is:
\begin{equation}
E^{n-1-m}T^{m+1}(L^\mu E\beta_\mu)-(E\beta_\mu)E^{n-1-m}T^{m+1}(\rho N^\mu)  
\ \mbox{: up to ignorable terms} 
\label{14.60}
\end{equation}
and we have:
\begin{eqnarray}
&&E^{n-1-m}T^{m+1}(L^\mu E\beta_\mu)=E^{n-1-m}T^{m+1}(E^\mu L\beta_\mu)\nonumber\\
&&=E^\mu E^{n-1-m}T^{m+1}L\beta_\mu+(L\beta_\mu)E^{n-1-m}T^{m+1}E^\mu  \nonumber\\
&&\ \mbox{: up to ignorable terms}
\label{14.61}
\end{eqnarray}
and:
\begin{eqnarray}
&&E^\mu E^{n-1-m}T^{m+1}L\beta_\mu=E^\mu LE^{n-1-m}T^{m+1}\beta_\mu+\sss E^{n-1-m}T^m\zeta \nonumber\\
&& \hspace{33mm}\ \mbox{: up to ignorable terms} 
\label{14.62}
\end{eqnarray}
Using \ref{14.31}, \ref{12.74}, and the fact that by \ref{14.30} (see \ref{12.107}, \ref{12.108}, 
\ref{12.110}, \ref{12.111}):
\begin{equation}
[E^{n-1-m}T^m\zeta]_{P.A.}=2\rho\mu_{m,n-1-m}-2\rhob\mub_{m,n-1-m}
\label{14.63}
\end{equation}
we deduce from the above:
\begin{eqnarray}
&&\rho N^\mu E^{n-m}T^{m+1}\beta_\mu=E^\mu LE^{n-1-m}T^{m+1}\beta_\mu 
-c^{-1}\ss_N E^{n-1-m}T^{m+1}\lambdab\nonumber\\
&&\hspace{30mm}+\sss\rho\mu_{m,n-1-m}+(c^{-1}s_{NL}-2\pi\ss_N\rho-2\rhob\sss)\mub_{m,n-1-m}\nonumber\\
&&\hspace{30mm}+2\ss_N\rho(\rho\mu_{m-1,n-m}+\rhob\mub_{m-1,n-m}) \nonumber\\
&&\hspace{30mm} \ \mbox{: up to ignorable terms}
\label{14.64}
\end{eqnarray}
A similar formula holds for the $N$th approximants. Taking the difference we conclude that 
the 1st term on the left in \ref{14.58} is (see \ref{12.117}):
\begin{eqnarray}
&&\rho N^\mu E(E^{n-1-m}T^{m+1}\beta_\mu-E_N^{n-1-m}T^{m+1}\beta_{\mu,N}) \nonumber\\
&&\hspace{20mm}=\s^{(E;m+1,n-1-m)}\cxi_L-c^{-1}\ss_N \s^{(m+1,n-1-m)}\clab\nonumber\\
&&\hspace{20mm}+\sss\rho\cmu_{m,n-1-m}+(c^{-1}s_{NL}-2\pi\ss_N\rho-2\rhob\sss)\cmub_{m,n-1-m}\nonumber\\
&&\hspace{20mm}+2\ss_N\rho(\rho\cmu_{m-1,n-m}+\rhob\cmub_{m-1,n-m}) \nonumber\\
&&\hspace{20mm} \ \mbox{: up to ignorable terms} 
\label{14.65}
\end{eqnarray}
We shall presently estimate, in connection with \ref{14.51}, the integral:
\begin{equation}
\int_0^{\ub_1}\|\rho N^\mu E(E^{n-1-m}T^{m+1}\beta_\mu-E_N^{n-1-m}T^{m+1}\beta_{\mu,N})\|_{L^2(S_{\ub,u})}d\ub
\label{14.66}
\end{equation}
On the right hand side of \ref{14.65} we have 6 terms. The contribution to \ref{14.66} of the 1st term 
is bounded by:
\begin{eqnarray}
&&\ub_1^{1/2}\|\s^{(E;m+1,n-1-m)}\cxi_L\|_{L^2(C_u^{\ub_1})}\leq C\ub_1^{1/2}
\sqrt{\s^{(E;m+1,n-1-m)}\cE^{\ub_1}(u)}\nonumber\\
&&\hspace{48mm}\leq C\sqrt{D}\ub_1^{a_{m+1}+\frac{1}{2}}u^{b_{m+1}}
\label{14.67}
\end{eqnarray}
while the contribution of the 2nd term is bounded by:
\begin{eqnarray}
&&\ub_1^{1/2}\|\s^{(m,n-m)}\clab\|_{(C_u^{\ub_1})}\leq 
\ub_1^{a_m+\frac{1}{2}}u^{b_m}\s^{(m,n-m)}\Lab(\ub_1,u)\nonumber\\
&&\hspace{33mm}\leq C\sqrt{D}\ub_1^{a_m+\frac{1}{2}}u^{b_m}
\label{14.68}
\end{eqnarray}
by definition \ref{12.666}, Propositions 12.4, 12.5 and 12.6, and the top order energy estimates.
The contribution of the 3rd term is bounded by:
\begin{equation}
C\int_0^{\ub_1}\ub\|\cmu_{m,n-1-m}\|_{L^2(S_{\ub,u})}d\ub
\label{14.69}
\end{equation}
Now, $\|\cmu_{m,n-1-m}\|_{L^2(S_{\ub,u})}$ being bounded by \ref{12.567} with $m+1$ in the role of $m$, 
the $j=0$ term in the sum, that is $\|\s^{(m,n-m)}\cla\|_{L^2(S_{\ub,u})}$, can be absorbed in the 
integral inequality \ref{12.564}. The contribution of the remainder to \ref{14.69} is bounded  by:
\begin{eqnarray}
&&C\int_0^{\ub_1}\ub\left\{\sum_{j=1}^{m}u^{2j}\|\s^{(m-j,n-m+j)}\cla\|_{L^2(S_{\ub,u})}\right.
\nonumber\\
&&\hspace{50mm}\left. +u^{2m}\|\s^{(n-1)}\ctchi\|_{L^2(S_{\ub,u})}\right\}d\ub\nonumber\\
&&\leq C\sqrt{D}\int_0^{\ub_1}\ub\left\{\sum_{j=1}^m u^{2j}\ub^{a_{m+1-j}+\frac{1}{2}}
u^{b_{m+1-j}-\frac{3}{2}}+u^{2m}\ub^{a_0+\frac{1}{2}}u^{b_0-\frac{3}{2}}\right\}d\ub\nonumber\\
&&\leq\frac{C\sqrt{D}}{\left(a_m+\frac{5}{2}\right)}\ub_1^{a_m+\frac{5}{2}}u^{b_m+\frac{1}{2}}
\label{14.70}
\end{eqnarray}
by virtue of the inductive hypothesis and \ref{14.20}. Since $|s_{NL}|\leq C\rho$, the contribution to 
\ref{14.66} of the 4th term on the right in \ref{14.65} is bounded by:
\begin{equation}
Cu\int_0^{\ub_1}\|\cmub_{m,n-1-m}\|_{L^2(S_{\ub,u})}d\ub
\label{14.71}
\end{equation}
which by \ref{12.568} with $m+1$ in the role of $m$ is bounded by:
\begin{eqnarray}
&&Cu\left\{\sum_{j=0}^m \ub_1^{j+\frac{1}{2}}\|\s^{(m-j,n-m+j)}\clab\|_{L^2(C_u^{\ub_1})}
+\ub_1^{m-\frac{1}{2}}\|\rho\s^{(n-1)}\ctchib\|_{C_u^{\ub_1})}\right\}\nonumber\\
&&\leq C\sqrt{D}\ub_1^{a_m+\frac{1}{2}}u^{b_m+1}
\label{14.72}
\end{eqnarray}
Finally, the contribution of the 5th term is bounded by:
\begin{equation}
\frac{C\sqrt{D}}{\left(a_m+\frac{7}{2}\right)}\ub_1^{a_m+\frac{7}{2}}u^{b_m-\frac{3}{2}}
\label{14.73}
\end{equation}
(see \ref{14.54}) and the contribution of the 6th term is bounded by:
\begin{equation}
C\sqrt{D}\ub_1^{a_{m-1}+\frac{3}{2}}u^{b_{m-1}+2}
\label{14.74}
\end{equation}
(see \ref{14.55}). Combining the above results we conclude that \ref{14.66} is bounded by:
\begin{equation}
C\sqrt{D}\ub_1^{a_{m+1}+\frac{1}{2}}u^{b_{m+1}}
\label{14.75}
\end{equation}
up to terms which can be absorbed. 

Now, the contribution of \ref{14.57} to the integral \ref{14.51} is bounded in terms of the integral:
\begin{equation}
\int_0^{\ub_1}\|\rho\Nb^\mu E(E^{n-1-m}T^{m+1}\beta_\mu-E_N^{n-1-m}T^{m+1}\beta_{\mu,N})\|_{L^2(S_{\ub,u})}d\ub
\label{14.76}
\end{equation}
Since $\ogamma^{-1}\sim u^{-1}$ the contribution of the 1st term on the left in \ref{14.58} to this 
integral is bounded by $Cu^{-1}$ times the integral \ref{14.66}, that is by:
\begin{equation}
C\sqrt{D}\ub_1^{a_{m+1}+\frac{1}{2}}u^{b_{m+1}-1}
\label{14.77}
\end{equation}
On the other hand, the contribution of the right hand side of \ref{14.58} to \ref{14.76} is,  
in view of the fact that $\rho\ogamma^{-1}\sim\ub^{1/2}u^{-2}a^{1/2}$, bounded by: 
\begin{eqnarray}
&&Cu^{-2}\int_0^{\ub_1}\ub^{1/2}\|a^{1/2}\s^{(Y;m+1,n-1-m)}\csxi\|_{L^2(S_{\ub,u})}d\ub\nonumber\\
&&\leq C u^{-2}\left(\int_0^{\ub_1}\ub d\ub\right)^{1/2}\|a^{1/2}\s^{(Y;m+1,n-1-m)}\csxi\|_{L^2(C_u^{\ub_1})}
\nonumber\\
&&\leq Cu^{-2}\ub_1\sqrt{\s^{(Y;m+1,n-1-m)}\cE^{\ub_1}(u)}\leq C\sqrt{D}\ub_1^{a_{m+1}+1}u^{b_{m+1}-2}
\label{14.78}
\end{eqnarray}
This is the dominant contribution to the integral \ref{14.51}. The above lead to the conclusion that 
\begin{equation}
\|\s^{(m,n-m)}\cla\|_{L^2(S_{\ub,u})}\leq C\sqrt{D}\ub^{a_{m+1}+\frac{1}{2}}u^{b_{m+1}-\frac{3}{2}}
\label{14.79}
\end{equation}
and the inductive step is complete. 

We summarize the results of this section in the following proposition. 

\vspace{2.5mm}

\noindent{\bf Proposition 14.1} \ \ \  The next to top order acoustical quantities $\s^{(n-1)}\ctchi$ and 
$\s^{(m,n-m)}\cla \ : \ m=0,...,n-1$ satisfy the following $L^2$ estimates on the $S_{\ub,u}$:
$$\|\s^{(n-1)}\ctchi\|_{L^2(S_{\ub,u})}\leq C\sqrt{D}\ub^{a_0+\frac{1}{2}}u^{b_0-\frac{3}{2}}$$
and for $m=0,...,n-1$:
$$\|\s^{(m,n-m)}\cla\|_{L^2(S_{\ub,u})}\leq C\sqrt{D}\ub^{a_{m+1}+\frac{1}{2}}u^{b_{m+1}-\frac{3}{2}}$$

\vspace{5mm}

\section{$L^2({\cal K}_\sigma^\tau)$ Estimates for the $n$th Order Acoustical Differences}

Let us denote by ${\cal K}_\sigma^\tau$ the hypersurfaces, in terms of $(\ub,u,\vartheta)$ 
coordinates (see \ref{9.1}), 
\begin{equation}
{\cal K}_\sigma^\tau=\{(\ub,u,\vartheta) \ : \ u-\ub=\sigma, \ \ub\in[0,\tau], \ \vartheta\in S^1\}
\label{14.80}
\end{equation}
The ${\cal K}_\sigma^\tau$ are generated by integral curves of $T$ and we have 
${\cal K}_0^\tau={\cal K}^\tau$. In the present section we shall derive $L^2$ estimates for the 
$n$th order acoustical quantities $\s^{(n-1)}\ctchi$, $\s^{(m,n-m)}\cla \ : \ m=0,...,n$, and 
$\s^{(n-1)}\ctchib$, $\s^{(m,n-m)}\clab \ : \ m=0,...,n$, on the ${\cal K}_\sigma^\tau$. 

To derive an estimate for $\s^{(n-1)}\ctchi$ in $L^2({\cal K}_\sigma^{\tau_1})$, we revisit the 
proof of Lemma 12.2 starting with the inequality \ref{12.128}. Placing $\tau$ in the role of $\ub_1$ 
and $\sigma+\tau$ in the role of $u$ this inequality reads:
\begin{equation}
\|\s^{(n-1)}\ctchi\|_{L^2(S_{\tau,\sigma+\tau})}\leq k\int_0^\tau 
\|\s^{(n-1)}\chF\|_{L^2(S_{\ub,\sigma+\tau})}d\ub
\label{14.81}
\end{equation}
In regard to \ref{12.131}, as previously noted the 1st can be absorbed in the inequality while the 
2nd is bounded by (see \ref{14.3}):
\begin{equation}
C\sqrt{D}\tau^{a_0+\frac{3}{2}}(\sigma+\tau)^{b_0-1}
\label{14.82}
\end{equation}
The contribution of this to: 
\begin{equation}
\|\s^{(n-1)}\ctchi\|_{L^2({\cal K}_\sigma^{\tau_1})}=
\left\{\int_0^{\tau_1}\|\s^{(n-1)}\ctchi\|^2_{L^2(S_{\tau,\sigma+\tau})}d\tau\right\}^{1/2}
\label{14.83}
\end{equation}
is bounded by:
\begin{equation}
\frac{C\sqrt{D}}{\sqrt{2(a_0+2)}}\tau_1^{a_0+2}u_1^{b_0-1} \ \ \mbox{where} \ \ u_1=\sigma+\tau_1
\label{14.84}
\end{equation}
Proceeding to the principal difference terms \ref{12.132}, we consider the 2nd of these difference 
terms, which gives the dominant contribution through its $\Nb$ component, last of \ref{12.134}.  
In estimating this we shall use an appropriate version of Lemma 12.1b. The point is that we have 
integrals of the form:
\begin{equation}
g(u)=\int_0^\tau f(\ub,u)d\ub, \ \ u=\sigma+\tau, \ \ \mbox{$\sigma$ is fixed}
\label{14.85}
\end{equation}
$f$ being a non-negative function on $R_{\delta,\delta}$, and we are considering:
\begin{equation}
\|g\|_{L^2(\sigma,u_1)}=\left(\int_\sigma^{u_1}g^2(u)du\right)^{1/2}  
\ \ \mbox{where} \ \ u_1=\sigma+\tau_1
\label{14.86}
\end{equation}
Let us define:
\begin{equation}
\overline{f}(\ub,u)=\left\{\begin{array} {lll} f(\ub,u) & : & \mbox{if $\ub\in[0,u-\sigma]$}\\
0 & : & \mbox{if $\ub\in(u-\sigma,\tau_1]$} \end{array}\right.
\label{14.87}
\end{equation}
We then have:
\begin{equation}
g(u)=\int_0^{\tau_1}f(\ub,u)d\ub
\label{14.88}
\end{equation}
It follows that:
\begin{equation}
\|g\|_{L^2(\sigma,u_1)}\leq\int_0^{\tau_1}\|\overline{f}(\ub,\cdot)\|_{L^2(\sigma,u_1)}d\ub
\label{14.89}
\end{equation}
Since:
\begin{eqnarray}
&&\|\overline{f}(\ub,\cdot)\|^2_{L^2(\sigma,u_1)}=\int_\sigma^{u_1}\overline{f}^2(\ub,u)du\nonumber\\
&&\hspace{25mm}=\int_{\sigma+\ub}^{u_1}f^2(\ub,u)du\leq\int_{\ub}^{u_1}f^2(\ub,u)du \label{14.90}
\end{eqnarray}
we conclude that:
\begin{equation}
\|g\|_{L^2(\sigma,u_1)}\leq \int_0^{\tau_1}\left(\int_{\ub}^{u_1}f^2(\ub,u)du\right)^{1/2}d\ub
\label{14.91}
\end{equation}
We apply this to the 1st of \ref{12.142} with $\tau$ in the role of $\ub_1$, setting:
\begin{equation}
f(\ub,u)=\ub u^{-2}\|\s^{(E;0,n)}\cxi_{\Lb}\|_{L^2(S_{\ub,u})}
\label{14.92}
\end{equation}
as in \ref{12.143}. By \ref{12.144} and the top order energy estimates:
\begin{eqnarray}
&&\int_{\ub}^{u_1}f^2(\ub,u)du\leq C\ub^{2a_0+2}u_1^{2b_0-4}\s^{(E;0,n)}\cBb(\ub_1,u_1)\nonumber\\
&&\hspace{25mm}\leq CD\ub^{2a_0+2}u_1^{2b_0-4}
\label{14.93}
\end{eqnarray}
Substituting in \ref{14.91} we then obtain:
\begin{eqnarray}
&&\|g\|_{L^2(\sigma,u_1)}\leq C\sqrt{D}\int_0^{\tau_1}\ub^{a_0+1}d\ub \cdot u_1^{b_0-2}\nonumber\\
&&\hspace{18mm}=\frac{C\sqrt{D}}{(a_0+2)}\tau_1^{a_0+2}u_1^{b_0-2} \label{14.94}
\end{eqnarray}
as the bound for the corresponding contribution to 
$\|\s^{(n-1)}\ctchi\|_{L^2({\cal K}_\sigma^{\tau_1})}$. In connection with the 3rd of \ref{12.142}, 
we set:
\begin{equation}
f(\ub,u)=\ub u^{-2}\|\s^{(n-1)}\ctchi\|_{L^2(S_{\ub,u})}
\label{14.95}
\end{equation}
Since by definition \ref{12.411}, Propositions 12.2 and 12.3, and the top order energy estimates, 
\begin{equation}
\|\s^{(n-1)}\ctchi\|_{L^2(\Cb_{\ub}^u)}\leq\ub^{a_0}u^{b_0}\s^{(n-1)}X(\ub_1,u_1)
\leq C\sqrt{D}\ub^{a_0}u^{b_0}
\label{14.96}
\end{equation}
we similarly obtain:
\begin{equation}
\|g\|_{L^2(\sigma,u_1)}\leq\frac{C\sqrt{D}}{(a_0+2)}\tau_1^{a_0+2}u_1^{b_0-2}
\label{14.97}
\end{equation}
as the bound for the corresponding contribution to 
$\|\s^{(n-1)}\ctchi\|_{L^2({\cal K}_\sigma^{\tau_1})}$. In connection with the 4th of \ref{12.142}, 
we set:
\begin{equation}
f(\ub,u)=\ub u^{-2}\|\s^{(0,n)}\cla\|_{L^2(S_{\ub,u})}
\label{14.98}
\end{equation}
Since by definition \ref{12.412}, Propositions 12.2 and 12.3, and the top order energy estimates,
\begin{equation}
\|\s^{(0,n)}\cla\|_{L^2(\Cb_{\ub}^{u})}\leq\ub^{a_0+\frac{1}{2}}u^{b_0+\frac{1}{2}}
\s^{(0,n)}\La(\ub_1,u_1)\leq C\sqrt{D}\ub^{a_0+\frac{1}{2}}u^{b_0+\frac{1}{2}}
\label{14.99}
\end{equation}
we obtain:
\begin{equation}
\|g\|_{L^2(\sigma,u_1)}\leq\frac{C\sqrt{D}}{\left(a_0+\frac{5}{2}\right)}
\tau_1^{a_0+\frac{1}{2}}u_1^{b_0-\frac{3}{2}}
\label{14.100}
\end{equation}
the corresponding contribution to $\|\s^{(n-1)}\ctchi\|_{L^2({\cal K}_\sigma^{\tau_1})}$. 
The contributions of the remaining terms are straightforward to estimate. We deduce in this manner 
an estimate for $\s^{(n-1)}\ctchi$ in $L^2({\cal K}_\sigma^{\tau_1})$ of the form:
\begin{equation}
\|\s^{(n-1)}\ctchi\|_{L^2({\cal K}_\sigma^{\tau_1})}\leq C\sqrt{D}\tau_1^{a_0}u_1^{b_0}
\label{14.101}
\end{equation}

To derive an estimate for $\s^{(0,n)}\cla$ in $L^2({\cal K}_\sigma^{\tau_1})$ we revisit the proof of 
Lemma 12.3 stating with the inequality \ref{12.196}, which we write in a form similar to \ref{14.81}: 
\begin{equation}
\|\s^{(0,n)}\cla\|_{L^2(S_{\tau,\sigma+\tau})}\leq k\int_0^{\tau}
\|\s^{(0,n)}\chG\|_{L^2(S_{\ub,\sigma+\tau})}d\ub
\label{14.102}
\end{equation}
In regard to \ref{12.199} the 1st can be absorbed in the inequality while the 2nd is bounded by 
\begin{eqnarray}
&&\tau^{1/2}\|\s^{(0,n)}\clab\|_{L^2(C_{\sigma+\tau}^{\tau})}\leq 
\tau^{a_0+\frac{1}{2}}(\sigma+\tau)^{b_0}\s^{(0,n)}\Lab(\tau,\sigma+\tau)\nonumber\\
&&\hspace{32mm}\leq C\sqrt{D}\tau^{a_0+\frac{1}{2}}(\sigma+\tau)^{b_0}
\label{14.103}
\end{eqnarray}
by definition\ref{12.461}, Propositions 12.2 and 12.3, and the top order energy estimates. The 
contribution of this to:
\begin{equation}
\|\s^{(0,n)}\cla\|_{L^2({\cal K}_\sigma^{\tau_1})}=\left\{\int_0^{\tau_1}\|\s^{(0,n)}\cla\|^2_{L^2
(S_{\tau,\sigma+\tau})}d\tau\right\}^{1/2}
\label{14.104}
\end{equation}
is bounded by:
\begin{equation}
\frac{C\sqrt{D}}{\sqrt{2(a_0+1)}}\tau_1^{a_0+1}u_1^{b_0}
\label{14.105}
\end{equation}
The 4th of \ref{12.199} coincides with the 2nd of \ref{12.131} and its contribution has been 
estimated by \ref{14.84}. To estimate the contribution of the 3rd of \ref{12.199} we apply \ref{14.91} 
taking:
\begin{equation}
f(\ub,u)=\ub\|\s^{(n-1)}\ctchi\|_{L^2(S_{\ub,u})}
\label{14.106}
\end{equation}
We then have:
\begin{eqnarray}
&&\int_{\ub}^{u_1}f^2(\ub,u)du=\ub^2\|\s^{(n-1)}\ctchi\|^2_{L^2(\Cb_{\ub}^{u_1})}\nonumber\\
&&\hspace{20mm}\leq\ub^{2a_0+2}u_1^{2b_0}\left(\s^{(n-1)}X(\ub,u_1)\right)^2
\leq CD\ub^{2a_0+2}u_1^{2b_0}\nonumber\\
&&\label{14.107}
\end{eqnarray}
by the definition \ref{12.411}, Propositions 12.2 and 12.3, and the top order energy estimates. 
Substituting in \ref{14.91} we then obtain:
\begin{equation}
\|g\|_{L^2(\sigma,u_1)}\leq \frac{C\sqrt{D}}{(a_0+2)}\tau_1^{a_0+2}u_1^{b_0}
\label{14.108}
\end{equation}
as the bound for the corresponding contribution to $\|\s^{(0,n)}\cla\|_{L^2({\cal K}_\sigma^{\tau_1})}$. 
The dominant contribution of the principal differences \ref{12.200} is through the $\Nb$ 
component of the first difference, expressed by the last of \ref{12.202}. More precisely, the 
dominant contribution is that the 1st of \ref{12.217}. To estimate its contribution to 
$\|\s^{(0,n)}\cla\|_{L^2({\cal K}_\sigma^{\tau_1})}$ we apply \ref{14.91} taking:
\begin{equation}
f(\ub,u)=\ub u^{-1}\|\s^{(Y;0,n)}\cxi_{\Lb}\|_{L^2(S_{\ub,u})}
\label{14.109}
\end{equation}
as in \ref{12.221}. By \ref{12.222}:
\begin{eqnarray}
&&\int_{\ub}^{u_1}f^2(\ub,u)du\leq C\ub^{2a_0+2}u_1^{2b_0-2}\s^{(Y;0,n)}\cBb(\ub_1,u_1)\nonumber\\
&&\hspace{22mm}\leq CD\ub^{2a_0+2}u_1^{2b_0-2}
\label{14.110}
\end{eqnarray}
by the top order energy estimates. Substituting in \ref{14.91} we then obtain:
\begin{equation}
\|g\|_{L^2(\sigma,u_1)}\leq \frac{C\sqrt{D}}{(a_0+2)}\tau_1^{a_0+2}u_1^{b_0-1}
\label{14.111}
\end{equation}
as the bound for the corresponding contribution to $\|\s^{(0,n)}\cla\|_{L^2({\cal K}_\sigma^{\tau_1})}$. 
The contributions of the remaining terms are estimated in a similar manner. We deduce in this way an 
estimate for $\s^{(0,n)}\cla$ in $L^2({\cal K}_\sigma^{\tau_1})$ of the form:
\begin{equation}
\|\s^{(0,n)}\cla\|_{L^2({\cal K}_\sigma^{\tau_1})}\leq C\sqrt{D}
\tau_1^{a_0+\frac{1}{2}}u_1^{b_0+\frac{1}{2}}
\label{14.112}
\end{equation}

To derive an estimate for $\s^{(n-1)}\ctchib$ in $L^2({\cal K}_{\sigma_1}^{\tau_1})$ 
we revisit the proof of Lemma 12.4. We start with the inequality \ref{12.263}, 
placing $\tau$ in the role of $\ub$, 
$\sigma+\tau$ in the role of $u$ and $\sigma_1+\tau$ in the role of $u_1$:
\begin{equation}
\|\s^{(n-1)}\ctchib\|_{L^2(S_{\tau,\sigma_1+\tau})}\leq k\|\s^{(n-1)}\ctchib\|_{L^2(S_{\tau,\tau})}
+k\int_0^{\sigma_1}\|\s^{(n-1)}\chFb\|_{L^2(S_{\tau,\sigma+\tau})}d\sigma
\label{14.113}
\end{equation}
In regard to \ref{12.268}, as previously noted the 1st can be absorbed in the inequality. The 2nd is:
\begin{eqnarray}
&&\int_0^{\sigma_1}\|\s^{(n-1)}\ctchi\|_{L^2(S_{\tau,\sigma+\tau})}d\sigma
\leq\sigma_1^{1/2}\|\s^{(n-1)}\ctchi\|_{L^2(\Cb_{\tau}^{\sigma_1+\tau})}\nonumber\\
&&\hspace{40mm}\leq\tau^{a_0}(\sigma_1+\tau)^{b_0+\frac{1}{2}}\s^{(n-1)}X(\ub_1,u_1) \nonumber\\
&&\hspace{30mm}\mbox{where:} \ \ \ub_1=\tau_1, \ u_1=\sigma_1+\tau_1 \label{14.114}
\end{eqnarray}
by the definition \ref{12.411}. By Propositions 12.2 and 12.3 and the top order energy estimates 
this is bounded by:
\begin{equation}
C\sqrt{D}\tau^{a_0}(\sigma_1+\tau)^{b_0+\frac{1}{2}}
\label{14.115}
\end{equation}
The corresponding contribution to:
\begin{equation}
\|\s^{(n-1)}\ctchib\|_{L^2({\cal K}_{\sigma_1}^{\tau_1})}=
\left\{\int_0^{\tau_1}\|\s^{(n-1)}\ctchib\|^2_{L^2(S_{\tau,\sigma_1+\tau})}d\tau\right\}^{1/2}
\label{14.116}
\end{equation}
is then bounded by:
\begin{equation}
\frac{C\sqrt{D}}{\sqrt{2a_0+1}}\tau_1^{a_0+\frac{1}{2}}u_1^{b_0+\frac{1}{2}}
\label{14.117}
\end{equation}
The dominant contribution of the principal differences \ref{12.270} is through the $\Nb$ 
component of the first difference, expressed by the last of \ref{12.272}. More precisely, the 
dominant contribution is that the 1st of \ref{12.280}:
\begin{eqnarray}
&&\int_0^{\sigma_1}(\sigma+\tau)^{-1}\|\s^{(Y;0,n)}\cxi_{\Lb}\|_{L^2(S_{\tau,\sigma+\tau})}d\sigma
\nonumber\\
&&\leq\left(\int_0^{\sigma_1}(\sigma+\tau)^{-2}d\sigma\right)^{1/2}\left(\int_0^{\sigma_1}
\|\s^{(Y;0,n)}\cxi_{\Lb}\|^2_{L^2(S_{\tau,\sigma+\tau})}d\sigma\right)^{1/2}\nonumber\\
&&\leq\tau^{-1/2}\|\s^{(Y;0,n)}\cxi_{\Lb}\|_{L^2(\Cb_\tau^{\sigma_1+\tau})}
\leq C\tau^{-1/2}\sqrt{\s^{(Y;0,n)}\cEb^{\sigma_1+\tau}(\tau)}\nonumber\\
&&\leq C\tau^{a_0-\frac{1}{2}}(\sigma_1+\tau)^{b_0}\sqrt{\s^{(Y;0,n)}\cBb(\ub_1,u_1)}
\leq C\sqrt{D}\tau^{a_0-\frac{1}{2}}u_1^{b_0}
\label{14.118}
\end{eqnarray}
by the top order energy estimates. The corresponding contribution to 
$\|\s^{(n-1)}\ctchib\|_{L^2({\cal K}_{\sigma_1}^{\tau_1})}$ is then bounded by:
\begin{equation}
\frac{C\sqrt{D}}{\sqrt{2a_0}}\tau_1^{a_0}u_1^{b_0}
\label{14.119}
\end{equation}
The remaining contributions from the integral on the right in \ref{14.113} are estimated in a similar 
manner. We then deduce from \ref{14.113} the inequality:
\begin{equation}
\|\s^{(n-1)}\ctchib\|_{L^2({\cal K}_{\sigma_1}^{\tau_1})}\leq 
k\|\s^{(n-1)}\ctchib\|_{L^2({\cal K}^{\tau_1})}+C\sqrt{D}\tau_1^{a_0}u_1^{b_0}
\label{14.120}
\end{equation}
We have:
\begin{eqnarray}
&&\|\s^{(n-1)}\ctchib\|^2_{L^2({\cal K}^{\tau_1})}=
\int_0^{\tau_1}\|\s^{(n-1)}\ctchib\|^2_{L^2(S_{\tau,\tau})}d\tau\nonumber\\
&&=\int_0^{\tau_1}\tau^{-2}\frac{\partial}{\partial\tau}\left(\int_0^{\tau}\tau^{\prime 2}
\|\s^{(n-1)}\ctchib\|^2_{L^2(S_{\tau^\prime,\tau^\prime})}d\tau^\prime\right)d\tau\nonumber\\
&&=\tau_1^{-2}\int_0^{\tau_1}\tau^2\|\s^{(n-1)}\ctchib\|^2_{L^2(S_{\tau,\tau})}d\tau\nonumber\\
&&\hspace{15mm}+2\int_0^{\tau_1}\tau^{-3}\left(\int_0^{\tau}\tau^{\prime 2}
\|\s^{(n-1)}\ctchib\|^2_{L^2(S_{\tau^\prime,\tau^\prime})}d\tau^\prime\right)d\tau
\label{14.121}
\end{eqnarray}
Since $r\sim\tau$ the last is bounded by a constant multiple of: 
\begin{eqnarray}
&&\tau_1^{-2}\int_0^{\tau_1}\|r\s^{(n-1)}\ctchib\|^2_{L^2(S_{\tau,\tau})}d\tau\nonumber\\
&&\hspace{15mm}+2\int_0^{\tau_1}\tau^{-3}\left(\int_0^{\tau}
\|r\s^{(n-1)}\ctchib\|^2_{L^2(S_{\tau^\prime,\tau^\prime})}d\tau^\prime\right)d\tau\nonumber\\
&&=\tau_1^{-2}\|r\s^{(n-1)}\ctchib\|^2_{L^2({\cal K}^{\tau_1})}
+2\int_0^{\tau_1}\tau^{-3}\|r\s^{(n-1)}\ctchib\|^2_{L^2({\cal K}^{\tau})}d\tau\nonumber\\
&&\leq \tau_1^{2c_0-2}\left(1+\frac{2}{2c_0-2}\right)\left(\s^{(n-1)}\oXb(\tau_1)\right)^2
\nonumber\\
&&\leq CD\tau_1^{2c_0-2} \label{14.122}
\end{eqnarray}
by the definition \ref{12.398}, Propositions 12.1 and 12.3, and the top order energy estimates. 
Substituting then in \ref{14.120} yields the estimate:
\begin{equation}
\|\s^{(n-1)}\ctchib\|_{L^2({\cal K}_{\sigma_1}^{\tau_1})}\leq C\sqrt{D}\tau_1^{a_0}u_1^{b_0-1}
\label{14.123}
\end{equation}

To derive an estimate for $\s^{(0,n)}\clab$ in $L^2({\cal K}_{\sigma_1}^{\tau_1})$ 
we revisit the proof of Lemma 12.6. We start with the inequality \ref{12.306} in the form similar to 
\ref{14.113}:
\begin{equation}
\|\s^{(0,n)}\clab\|_{L^2(S_{\tau,\sigma_1+\tau})}\leq k\|\s^{(0,n)}\clab\|_{L^2(S_{\tau,\tau})}
+k\int_0^{\sigma_1}\|\s^{(0,n)}\chGb\|_{L^2(S_{\tau,\sigma+\tau})}d\sigma
\label{14.124}
\end{equation}
In regard to \ref{12.312}, as previously noted the 1st can be absorbed in the inequality. The 2nd is:
\begin{eqnarray}
&&\tau\int_0^{\sigma_1}\|\s^{(0,n)}\cla\|_{L^2(S_{\tau,\sigma+\tau})}d\sigma
\leq\tau\sigma_1^{1/2}\|\s^{(0,n)}\cla\|_{L^2(\Cb_{\tau}^{\sigma_1+\tau})}\nonumber\\
&&\hspace{40mm}\leq\tau^{a_0+\frac{3}{2}}(\sigma_1+\tau)^{b_0+1}\s^{(0,n)}\La(\ub_1,u_1) \nonumber\\
&&\hspace{30mm}\mbox{where:} \ \ \ub_1=\tau_1, \ u_1=\sigma_1+\tau_1 \label{14.125}
\end{eqnarray}
by the definition \ref{12.412}. By Propositions 12.2 and 12.3 and the top order energy estimates 
this is bounded by:
\begin{equation}
C\sqrt{D}\tau^{a_0+\frac{3}{2}}(\sigma_1+\tau)^{b_0+1}
\label{14.126}
\end{equation}
The corresponding contribution to:
\begin{equation}
\|\s^{(0,n)}\clab\|_{L^2({\cal K}_{\sigma_1}^{\tau_1})}=
\left\{\int_0^{\tau_1}\|\s^{(0,n)}\clab\|^2_{L^2(S_{\tau,\sigma_1+\tau})}d\tau\right\}^{1/2}
\label{14.127}
\end{equation}
is then bounded by:
\begin{equation}
\frac{C\sqrt{D}}{\sqrt{2(a_0+2)}}\tau_1^{a_0+2}u_1^{b_0+1}
\label{14.128}
\end{equation}
Similarly (see \ref{14.114}, \ref{14.115}) the contribution of the 4th of \ref{12.312} is bounded by:
\begin{equation}
\frac{C\sqrt{D}}{\sqrt{2a_0+3}}\tau_1^{a_0+\frac{3}{2}}u_1^{b_0+\frac{1}{2}}
\label{14.129}
\end{equation}

To estimate the contribution of the 3rd of \ref{12.312}, we apply the following corollary of 
Lemma 12.5. Let $f$ be a non-negative function on $R_{\delta,\delta}$. For 
$(\ub_1,u_1)\in R_{\delta,\delta}$, where $u_1=\sigma_1+\ub_1$, we define on $[0,\ub_1]$ the function:
\begin{equation}
g(\ub)=\int_{\ub}^{\sigma_1+\ub}f(\ub,u)du
\label{14.130}
\end{equation}
We then have:
\begin{eqnarray}
&&\|g\|_{L^2(0,\ub_1)}\leq\int_0^{\ub_1}\|f(\cdot,u)\|_{L^2(0,u)}du\nonumber\\
&&\hspace{20mm}+\int_{\ub_1}^{u_1}\|f(\cdot,u)\|_{L^2(0,\ub_1)}du
\label{14.131}
\end{eqnarray}
This follows by applying Lemma 12.5 with $f\chi_{\sigma_1}$ in the role of $f$, where $\chi_{\sigma_1}$ 
is the characteristic function of the domain: 
$$\{(\ub,u)\in R_{\delta,\delta} \ : \ u-\ub\leq\sigma_1\}$$

To estimate then the contribution of the 3rd of \ref{12.312}, we set:
\begin{equation}
f(\ub,u)=\ub\|\s^{(n-1)}\ctchib\|_{L^2(S_{\ub,u})}
\label{14.132}
\end{equation}
We then have, by the definition \ref{12.460}, 
\begin{eqnarray*}
&&\|f(\cdot,u)\|^2_{L^2(0,u)}=\int_0^u \ub^2\|\s^{(n-1)}\ctchib\|^2_{L^2(S_{\ub,u})}d\ub\\
&&\hspace{20mm}\leq C\|\rho\s^{(n-1)}\ctchib\|^2_{L^2(C_u^u)}\leq Cu^{2a_0+2b_0}
\left(\s^{(n-1)}\Xb(\ub_1,u_1)\right)^2\\
&&\hspace{20mm}\mbox{: for $u\in[0,\ub_1)$}
\end{eqnarray*}
\begin{eqnarray*}
&&\|f(\cdot,u)\|^2_{L^2(0,\ub_1)}=\int_0^{\ub_1} \ub^2\|\s^{(n-1)}\ctchib\|^2_{L^2(S_{\ub,u})}d\ub\\
&&\hspace{20mm}\leq C\|\rho\s^{(n-1)}\ctchib\|^2_{L^2(C_u^{\ub_1})}\leq C\ub_1^{2(a_0+1)}u^{2(b_0-1)}\left(\s^{(n-1)}\Xb(\ub_1,u_1)\right)^2\\
&&\hspace{20mm}\mbox{: for $u\in[\ub_1,u_1]$}
\end{eqnarray*}
hence by Propositions 12.2 and 12.3, and the top order energy estimates:
\begin{eqnarray}
&&\|f(\cdot,u)\|_{L^2(0,u)}\leq C\sqrt{D}u^{a_0+b_0} \ \ \mbox{: for $u\in[0,\ub_1)$} \nonumber\\
&&\|f(\cdot,u)\|_{L^2(0,\ub_1)}\leq C\sqrt{D}\ub_1^{a_0+1}u^{b_0-1} \ \ \mbox{: for $u\in[\ub_1,u_1]$} 
\label{14.133}
\end{eqnarray}
Substituting in \ref{14.131} we obtain:
\begin{eqnarray}
&&\|g\|_{L^2(0,\ub_1)}\leq C\sqrt{D}\left\{\int_0^{\ub_1}u^{a_0+b_0}du
+\int_{\ub_1}^{u_1}\ub_1^{a_0+1}u^{b_0-1}du\right\}\nonumber\\
&&\hspace{20mm}\leq\frac{C\sqrt{D}}{b_0}\ub_1^{a_0+1}u_1^{b_0}
\label{14.134}
\end{eqnarray}
as the bound for the contribution of the 3rd of \ref{12.312} to 
$\|\s^{(0,n)}\clab\|_{L^2({\cal K}_{\sigma_1}^{\tau_1})}$ (recall that $\ub_1=\tau_1$, 
$u_1=\sigma_1+\tau_1$). The contributions of the principal difference terms \ref{12.317} are bounded 
in a similar manner. We then deduce from \ref{14.124} the inequality:
\begin{equation}
\|\s^{(0,n)}\clab\|_{L^2({\cal K}_{\sigma_1}^{\tau_1})}
\leq k\|\s^{(0,n)}\clab\|_{L^2({\cal K}^{\tau_1})}
+C\sqrt{D}\tau_1^{a_0}u_1^{b_0+1}
\label{14.135}
\end{equation}
By the definition \ref{12.399}, Propositions 12.1 and 12.3, and the top order energy estimates:
\begin{equation}
\|\s^{(0,n)}\clab\|_{L^2({\cal K}^{\tau_1})}\leq\tau_1^{c_0}\s^{(0,n)}\oLab(\tau_1)
\leq C\sqrt{D}\tau_1^{c_0}
\label{14.136}
\end{equation}
Substituting we conclude that:
\begin{equation}
\|\s^{(0,n)}\clab\|_{L^2({\cal K}_{\sigma_1}^{\tau_1})}\leq C\sqrt{D}\tau_1^{a_0}u_1^{b_0}
\label{14.137}
\end{equation}

Revisiting the proof of Lemma 12.9 and proceeding in the same manner as in the derivation of the 
estimate for $\s^{(0,n)}\cla$ in $L^2({\cal K}_{\sigma}^{\tau_1})$ we deduce:
\begin{equation}
\|\s^{(m,n-m)}\cla\|_{L^2({\cal K}_{\sigma}^{\tau_1})}\leq C\sqrt{D}
\tau_1^{a_m+\frac{1}{2}}u_1^{b_m+\frac{1}{2}} \ \ \mbox{: for $m=1,...,n$  }
\label{14.138}
\end{equation}
(where $u_1=\sigma+\tau_1$). Also, revisiting the proof of Lemma 12.10 and proceeding in the same manner 
as in the derivation of the estimate for $\s^{(0,n)}\clab$ in $L^2({\cal K}_{\sigma_1}^{\tau_1})$ 
we deduce:
\begin{equation}
\|\s^{(m,n-m)}\clab\|_{L^2({\cal K}_{\sigma_1}^{\tau_1})}\leq C\sqrt{D}\tau_1^{a_m}u_1^{b_m} 
\ \ \ \mbox{: for $m=1,...,n$}
\label{14.139}
\end{equation}

We summarize the results of this section in the following proposition. 

\vspace{2.5mm} 

\noindent{\bf Proposition 14.2} \ \ \ The $n$th order acoustical difference quantities satisfy the 
following estimates in $L^2({\cal K}_{\sigma}^{\tau})$:
\begin{eqnarray*}
&&\|\s^{(n-1)}\ctchi\|_{L^2({\cal K}_{\sigma}^{\tau})}\leq C\sqrt{D}\tau^{a_0}u^{b_0}\\
&&\|\s^{(n-1)}\ctchib\|_{L^2({\cal K}_{\sigma}^{\tau})}\leq C\sqrt{D}\tau^{a_0}u^{b_0-1}\\
&&\|\s^{(m,n-m)}\cla\|_{L^2({\cal K}_{\sigma}^{\tau})}\leq C\sqrt{D}
\tau^{a_m+\frac{1}{2}}u^{b_m+\frac{1}{2}} \ \ \ \mbox{: for $m=0,...,n$}\\
&&\|\s^{(m,n-m)}\clab\|_{L^2({\cal K}_{\sigma}^{\tau})}\leq C\sqrt{D}\tau^{a_m}u^{b_m} 
\ \ \ \mbox{: for $m=0,...,n$}
\end{eqnarray*}
In the above, $u=\sigma+\tau$. 

\vspace{5mm}

\section{$L^2(S_{\ub,u})$ Estimates for the $n$th Order Variation Differences}

In the present section we shall derive $L^2$ estimates on the $S_{\ub,u}$ for the $n$th order variation 
differences (see \ref{9.236}): 
\begin{equation}
\s^{(m,n-m)}\check{\dot{\phi}}_\mu=E^{n-m}T^m\beta_\mu-E_N^{n-m}T^m\beta_{\mu,N}
\label{14.140}
\end{equation}

Let $f$ be a function on ${\cal N}$ which vanishes on $\Cb_0$. In terms of $(\ub,u,\vartheta^\prime)$ 
coordinates on ${\cal N}$ which are adapted to the flow of $L$ (see \ref{10.342}, \ref{10.343}) 
we have:
\begin{equation}
f(\ub_1,u,\vartheta^\prime)=\int_0^{\ub_1}(Lf)(\ub,u,\vartheta^\prime)d\ub
\label{14.141}
\end{equation}
In view of the relation \ref{10.352}, we write this in the form:
\begin{equation}
\left(f\sh^{\prime 1/4}\right)(\ub_1,u,\vartheta^\prime)=\int_0^{\ub_1}
\left(\frac{\sh^\prime(\ub_1,u,\vartheta^\prime)}{\sh^\prime(\ub,u,\vartheta^\prime)}\right)^{1/4}
\left((Lf)\sh^{\prime 1/4}\right)(\ub,u,\vartheta^\prime)d\ub
\label{14.142}
\end{equation}
By \ref{10.356} this implies:
\begin{equation}
\left|\left(f\sh^{\prime 1/4}\right)(\ub_1,u,\vartheta^\prime)\right|\leq k\int_0^{\ub_1}
\left|\left((Lf)\sh^{\prime 1/4}\right)(\ub,u,\vartheta^\prime)\right|d\ub
\label{14.143}
\end{equation}
where $k$ is a constant greater than 1, but which can be chosen as close to 1 as we wish by suitably 
restricting $\delta$. Taking $L^2$ norms with respect to $\vartheta^\prime\in S^1$ \ref{14.143} implies:
\begin{equation}
\left\|\left(f\sh^{\prime 1/4}\right)(\ub_1,u)\right\|_{L^2(S^1)}\leq 
k\int_0^{\ub_1}\left\|\left((Lf)\sh^{\prime 1/4}\right)(\ub,u)\right\|_{L^2(S^1)}d\ub
\label{14.144}
\end{equation}
which in view of the relation \ref{10.352} is simply:
\begin{equation}
\|f\|_{L^2(S_{\ub_1,u})}\leq k\int_0^{\ub_1}\|Lf\|_{L^2(S_{\ub,u})}d\ub
\label{14.145}
\end{equation}

We shall apply \ref{14.145} with the $\s^{(m,n-m)}\check{\dot{\phi}}_\mu$ in the role of $f$. 
Now, in view of the expansion:
\begin{equation}
(h^{-1})^{\lambda\nu}=E^\lambda E^\nu-\frac{1}{2c}\Nb^\lambda N^\nu-\frac{1}{2c}N^\lambda\Nb^\nu
\label{14.146}
\end{equation}
we can express:
\begin{equation}
\delta_\mu^\nu=h_{\mu\lambda}(h^{-1})^{\lambda\nu}=c_{E,\mu}E^\nu
+c_{N,\mu}N^\nu+c_{\Nb,\mu}\Nb^\nu
\label{14.147}
\end{equation}
where:
\begin{equation}
c_{E,\mu}=h_{\mu\nu}E^\nu, \ \ c_{N,\mu}=-\frac{h_{\mu\nu}\Nb^\nu}{2c}, \ \ 
c_{\Nb,\mu}=-\frac{h_{\mu\nu}N^\nu}{2c}
\label{14.148}
\end{equation}
are bounded coefficients. Therefore for any $\theta_\mu$ we can express:
\begin{equation}
\theta_\mu=c_{E,\mu}E^\nu\theta_\nu+c_{N,\mu}N^\nu\theta_\nu+c_{\Nb,\mu}\Nb^\nu\theta_\nu
\label{14.149}
\end{equation}
It follows that there is a constant $C$ such that:
\begin{equation}
|\theta_\mu|\leq C\left\{|E^\nu\theta_\nu|+|N^\nu\theta_\nu|+|\Nb^\nu\theta_\nu|\right\}
\label{14.150}
\end{equation}
We apply this to $\theta_\mu=L\s^{(m,n-m)}\check{\dot{\phi}}_\mu$. Taking then $L^2$ norms on 
$S_{\ub,u}$ we obtain:
\begin{eqnarray}
&&\|L\s^{(m,n-m)}\check{\dot{\phi}}_\mu\|_{L^2(S_{\ub,u})}\leq \nonumber\\
&&\hspace{20mm}C\left\{\|E^\nu L\s^{(m,n-m)}\check{\dot{\phi}}_\nu\|_{L^2(S_{\ub,u})}
+\|N^\nu L\s^{(m,n-m)}\check{\dot{\phi}}_\nu\|_{L^2(S_{\ub,u})}\right.\nonumber\\
&&\hspace{65mm}\left.+\|\Nb^\nu L\s^{(m,n-m)}\check{\dot{\phi}}_\nu\|_{L^2(S_{\ub,u})}\right\} 
\nonumber\\
&&\label{14.151}
\end{eqnarray}
Substituting in \ref{14.145} with $\s^{(m,n-m)}\check{\dot{\phi}}_\mu$ in the role of $f$ we 
then obtain:
\begin{eqnarray}
&&\|\s^{(m,n-m)}\check{\dot{\phi}}_\mu\|_{L^2(S_{\ub_1,u})}\leq C\int_0^{\ub_1}\|E^\nu L\s^{(m,n-m)}\check{\dot{\phi}}_\nu\|_{L^2(S_{\ub,u})}d\ub
\nonumber\\
&&\hspace{35mm}+C\int_0^{\ub_1}\|N^\nu L\s^{(m,n-m)}\check{\dot{\phi}}_\nu\|_{L^2(S_{\ub,u})}d\ub 
\nonumber\\
&&\hspace{35mm}+C\int_0^{\ub_1}\|\Nb^\nu L\s^{(m,n-m)}\check{\dot{\phi}}_\nu\|_{L^2(S_{\ub,u})}d\ub
\nonumber\\
&&\label{14.152}
\end{eqnarray}
Recalling \ref{9.275} the 1st integral on the right is:
\begin{eqnarray}
&&\int_0^{\ub_1}\|\s^{(E;m,n-m)}\cxi_L\|_{L^2(S_{\ub,u})}d\ub
\leq \ub_1^{1/2}\|\s^{(E;m,n-m)}\cxi_L\|_{L^2(C_u^{\ub_1})}\nonumber\\
&&\hspace{10mm}\leq C\ub_1^{1/2}\sqrt{\s^{(E;m,n-m)}\cE^{\ub_1}(u)}\leq C\ub_1^{a_m+\frac{1}{2}}u^{b_m}
\sqrt{\s^{(E;m,n-m)}\cB(\ub_1,u)}\nonumber\\
&&\hspace{40mm}\leq C\sqrt{D}\ub_1^{a_m+\frac{1}{2}}u^{b_m}
\label{14.153}
\end{eqnarray}
by the top order energy estimates. 

Consider next the 3rd integral on the right \ref{14.152}. Recalling that 
\begin{equation}
\Nb^\mu L\s^{(m,n-m)}\check{\dot{\phi}}_\mu=\Nb^\mu L(E^{n-m}T^m\beta_\mu-E_N^{n-m}T^m\beta_{\mu,N})
\label{14.154}
\end{equation}
to estimate this integral we must use the fact that $\Nb^\mu L\beta_{\mu}=\lambdab E^\mu E\beta_\mu$. 
We must therefore first commute $L$ with $E^{n-m}T^m$, then express $\Nb^\mu E^{n-m}T^m L\beta_\mu$ 
in terms of $E^{n-m}T^m(\Nb^\mu L\beta_\mu)$, then use the fact just mentioned, then express 
$E^{n-m}T^m (\lambdab E^\mu E\beta_\mu)$ in terms of $\lambdab E^\mu E^{n-m}T^m E\beta_\mu$, and 
finally commute $T^m$ with $E$ in the last to obtain $\lambdab E^\mu E^{n-m+1}T^m \beta_\mu$. 
A similar procedure applied to the $N$th approximants will give us an expression for 
$\Nb_N^\mu L_N E_N^{n-m}T^m\beta_{\mu,N}$ in terms of $\lambdab_N E_N^\mu E_N^{n-m+1}T^m\beta_{\mu,N}$. 
Subtracting will then yield an expression for \ref{14.154} in terms of $\lambdab \s^{(E;m,n-m)}\csxi$ 
up to 0th order difference terms. 

Consider first the case $m=0$. Then the first commutation consists of commuting $L$ with $E^n$. 
Since by the 1st of the commutation relations \ref{3.a14}:
\begin{equation}
[L,E^n]=-\sum_{i=0}^{n-1}E^{n-1-i}\chi E^{i+1}
\label{14.155}
\end{equation}
this will leave an $n$th order remainder with a principal acoustical part of:
\begin{equation}
-\ss_{\Nb}\rho E^{n-1}\tchi
\label{14.156}
\end{equation}
contributed by the $i=0$ term in the above sum. Expressing $\Nb^\mu E^n L\beta_\mu$ in terms of 
$E^n (\Nb^\mu L\beta_\mu)$ will leave an $n$th order remainder with a principal acoustical part of:
\begin{equation}
-(\ss_N-c\pi\sss)\rho E^{n-1}\tchib
\label{14.157}
\end{equation} 
contributed by $-(E^n\Nb^\mu)L\beta_\mu$ (see \ref{12.7}). Expressing $E^n(\lambdab E^\mu E\beta_\mu)$ 
in terms of $\lambdab E^\mu E^{n+1}\beta_\mu$ will leave an $n$th order remainder with a principal 
acoustical part of:
\begin{equation}
\sss E^n\lambdab+\frac{\rho}{2}(\ss_N E^{n-1}\tchib+\ss_{\Nb}E^{n-1}\tchi)
\label{14.158}
\end{equation}
contributed by the terms $(E^n\lambdab)E^\mu E\beta_\mu$ and $\lambdab(E^n E^\mu)E\beta_\mu$ in 
$(E^n(\lambdab E^\mu))E\beta_\mu$ (see \ref{12.8}). Combining and subtracting the analogous result for the $N$th 
approximants we conclude that in the case $m=0$ \ref{14.154} is given up to ignorable terms by:
\begin{equation}
\lambdab\s^{(E;0,n)}\csxi+\s^{(0,n)}\clab+\left(c\pi\sss-\frac{1}{2}\ss_N\right)\rho\s^{(n-1)}\ctchib
-\frac{1}{2}\ss_{\Nb}\rho\s^{(n-1)}\ctchi
\label{14.159}
\end{equation}
Therefore in the case $m=0$ the 3rd of the integrals on the right in \ref{14.152} is bounded by:
\begin{eqnarray}
&&\int_0^{\ub_1}\|\lambdab\s^{(E;0,n)}\csxi\|_{L^2(S_{\ub,u})}d\ub
+C\int_0^{\ub_1}\|\s^{(0,n)}\clab\|_{L^2(S_{\ub,u})}d\ub\nonumber\\
&&+C\int_0^{\ub_1}\|\rho\s^{(n-1)}\ctchib\|_{L^2(S_{\ub,u})}d\ub
+C\int_0^{\ub_1}\|\rho\s^{(n-1)}\ctchi\|_{L^2(S_{\ub,u})}d\ub \nonumber\\
&&\label{14.160}
\end{eqnarray}
up to the contribution of what we have called ignorable terms. These are terms of order $n$ with 
vanishing principal acoustical part. The contribution of these ignorable terms can be absorbed. 
Since $\lambdab\sim \ub$ while $a^{1/2}\sim \ub^{1/2}u$, the 1st of \ref{14.160} is bounded by: 
\begin{eqnarray}
&&Cu^{-1}\int_0^{\ub_1}\ub^{1/2}\|a^{1/2}\s^{(E;0,n)}\csxi\|_{L^2(S_{\ub,u})}d\ub
\leq Cu^{-1}\ub_1\|a^{1/2}\s^{(E;0,n)}\csxi\|_{L^2(C_u^{\ub_1})}\nonumber\\
&&\hspace{10mm}\leq Cu^{-1}\ub_1\sqrt{\s^{(E;0,n)}\cE^{\ub_1}(u)}\leq C\ub_1^{a_0+1}u^{b_0-1}
\sqrt{\s^{(E;0,n)}\cB(\ub_1,u)}\nonumber\\
&&\hspace{40mm}\leq C\sqrt{D}\ub_1^{a_0+1}u^{b_0-1}
\label{14.161}
\end{eqnarray}
by the top order energy estimates. The 2nd of \ref{14.160} is bounded by:
\begin{eqnarray}
&&C\ub_1^{1/2}\|\s^{(0,n)}\clab\|_{L^2(C_u^{\ub_1})}\leq 
C\ub_1^{a_0+\frac{1}{2}}u^{b_0}\s^{(0,n)}\Lab(\ub_1,u)\nonumber\\
&&\hspace{20mm}\leq C\sqrt{D}\ub_1^{a_0+\frac{1}{2}}u^{b_0}\label{14.162}
\end{eqnarray}
(see \ref{14.10}), while the 3rd of \ref{14.160} is bounded according to \ref{14.3}. As for the 4th of 
\ref{14.160}, this is bounded using the estimate for $\s^{(n-1)}\ctchi$ in $L^2(S_{\ub,u})$ of 
Proposition 14.1 by:
\begin{equation}
C\sqrt{D}\int_0^{\ub_1}\ub^{a_0+\frac{3}{2}}u^{b_0-\frac{3}{2}}d\ub
=\frac{C\sqrt{D}}{\left(a_0+\frac{5}{2}\right)}\ub_1^{a_0+\frac{5}{2}}u^{b_0-\frac{3}{2}}
\label{14.163}
\end{equation}
Combining the above we conclude that in the case $m=0$ the 3rd of the integrals on the right in 
\ref{14.152} is bounded by:
\begin{equation}
C\sqrt{D}\ub_1^{a_0+\frac{1}{2}}u^{b_0-\frac{1}{2}}
\label{14.164}
\end{equation}

Consider next the cases $m=1,...,n$. Then in the first commutation we must first commute $L$ with 
$E^{n-m}$. This leaves a ignorable remainder. We must then commute $L$ with $T^m$ and then apply 
$E^{n-m}$. Since by the 3rd of the commutation relations \ref{3.a14}:
\begin{equation}
[L,T^m]=-\sum_{j=0}^{m-1}T^{m-1-j}\zeta ET^j
\label{14.165}
\end{equation}
this will leave an $n$th order remainder with a principal acoustical part of, by \ref{14.63}, 
\begin{equation}
-2\ss_{\Nb}(\rho\mu_{m-1,n-m}-\rhob\mub_{m-1,n-m})
\label{14.166}
\end{equation}
contributed by the $j=0$ term in the above sum. Expressing $\Nb^\mu E^{n-m}T^m L\beta_\mu$ in terms of 
$E^{n-m}T^m(\Nb^\mu L\beta_\mu)$ will leave an $n$th order remainder with a principal acoustical part of:
\begin{equation}
-2(\ss_N-c\pi\sss)\rho\mub_{m-1,n-m}
\label{14.167}
\end{equation} 
contributed by $-(E^{n-m}T^m\Nb^\mu)L\beta_\mu$ (see \ref{12.66}). Expressing 
$E^{n-m}T^m(\lambdab E^\mu E\beta_\mu)$ 
in terms of $\lambdab E^\mu E^{n-m+1}T^m\beta_\mu$ will leave an $n$th order remainder with a principal 
acoustical part of:
\begin{equation}
\sss E^{n-m}T^m\lambdab+\rho(\ss_N\mub_{m-1,n-m}+\ss_{\Nb}\mu_{m-1,n-m})
\label{14.168}
\end{equation}
contributed by the terms $(E^{n-m}T^m\lambdab)E^\mu E\beta_\mu$ and 
$\lambdab(E^{n-m}T^m E^\mu)E\beta_\mu$ in \\$(E^{n-m}T^m(\lambdab E^\mu))E\beta_\mu$ (see \ref{12.74}).
Finally, in the last commutation we must commute $T^m$ with $E$ and then apply $E^{n-m}$. 
Since by the first two of the commutation relations \ref{3.a14}:
\begin{equation}
[T^m,E]=-\sum_{j=0}^{m-1}T^{m-1-j}(\chi+\chib)ET^j
\label{14.169}
\end{equation}
this will leave an $n$th order remainder with a principal acoustical part of:
\begin{equation}
-\lambdab\sss\cdot\left\{ \begin{array}{lll} 
(\rho E^{n-1}\tchi+\rhob E^{n-1}\tchib) & : & \mbox{for $m=1$}\\
2(\rho\mu_{m-2,n-m+1}+\rhob\mub_{m-2,n-m+1}) & : & \mbox{for $m=2,...,n$}
\end{array}\right.
\label{14.170}
\end{equation}
contributed by the $j=0$ term in the above sum. The expression for $m=2,...,n$ results in view of 
\ref{12.64} by \ref{12.105} - \ref{12.108} and \ref{12.110}, \ref{12.111}. Combining and subtracting 
the analogous result for the $N$th 
approximants we conclude that for $m=1,...,n$ \ref{14.154} is given up to ignorable terms by:
\begin{eqnarray}
&&\lambdab\s^{(E;m,n-m)}\csxi+\sss\s^{(m,n-m)}\clab\nonumber\\
&&+\left[\left(2c\pi\sss-\ss_N\right)\rho+2\ss_{\Nb}\rhob\right]\cmub_{m-1,n-m}
-\ss_{\Nb}\rho\cmu_{m-1,n-m}\nonumber\\
&&-\lambdab\sss\cdot\left\{ \begin{array}{lll} 
(\rho \s^{(n-1)}\ctchi+\rhob \s^{(n-1)}\ctchib) & : & \mbox{for $m=1$}\\
2(\rho\cmu_{m-2,n-m+1}+\rhob\cmub_{m-2,n-m+1}) & : & \mbox{for $m=2,...,n$}
\end{array}\right. \nonumber\\
&&\label{14.171}
\end{eqnarray}
It follows that for $m=1,...,n$ the 3rd of the integrals \ref{14.152} is bounded, 
up to ignorable terms which can be absorbed, by the sum of:
\begin{eqnarray}
&&\int_0^{\ub_1}\|\lambdab\s^{(E;m,n-m)}\csxi\|_{L^2(S_{\ub,u})}d\ub
+C\int_0^{\ub_1}\|\s^{(m,n-m)}\clab\|_{L^2(S_{\ub,u})}d\ub\nonumber\\
&&+Cu\int_0^{\ub_1}\|\cmub_{m-1,n-m}\|_{L^2(S_{\ub,u})}d\ub
+C\int_0^{\ub_1}\ub\|\cmu_{m-1,n-m}\|_{L^2(S_{\ub,u})}d\ub\nonumber\\
&&\label{14.172}
\end{eqnarray}
and:
\begin{equation}
Cu^2\int_0^{\ub_1}\|\rho \s^{(n-1)}\ctchib\|_{L^2(S_{\ub,u})}d\ub
+C\int_0^{\ub_1}\ub^2\|\s^{(n-1)}\ctchi\|_{L^2(S_{\ub,u})}d\ub
\label{14.173}
\end{equation}
in the case $m=1$, 
\begin{equation}
Cu^2\int_0^{\ub_1}\ub\|\cmub_{m-2,n-m+1}\|_{L^2(S_{\ub,u})}d\ub
+C\int_0^{\ub_1}\ub^2\|\cmu_{m-2,n-m+1}\|_{L^2(S_{\ub,u})}d\ub
\label{14.174}
\end{equation}
for $m=2,...,n$. Since $\lambdab\sim \ub$ while $a^{1/2}\sim \ub^{1/2}u$, the 1st of \ref{14.172} is bounded by: 
\begin{eqnarray}
&&Cu^{-1}\int_0^{\ub_1}\ub^{1/2}\|a^{1/2}\s^{(E;m,n-m)}\csxi\|_{L^2(S_{\ub,u})}d\ub
\leq Cu^{-1}\ub_1\|a^{1/2}\s^{(E;m,n-m)}\csxi\|_{L^2(C_u^{\ub_1})}\nonumber\\
&&\hspace{10mm}\leq Cu^{-1}\ub_1\sqrt{\s^{(E;m,n-m)}\cE^{\ub_1}(u)}\leq C\ub_1^{a_m+1}u^{b_m-1}
\sqrt{\s^{(E;m,n-m)}\cB(\ub_1,u)}\nonumber\\
&&\hspace{40mm}\leq C\sqrt{D}\ub_1^{a_m+1}u^{b_m-1}
\label{14.175}
\end{eqnarray}
by the top order energy estimates. The 2nd of \ref{14.172} is bounded by: 
\begin{eqnarray}
&&C\ub_1^{1/2}\|\s^{(m,n-m)}\clab\|_{(C_u^{\ub_1})}\leq 
C\ub_1^{a_m+\frac{1}{2}}u^{b_m}\s^{(m,n-m)}\Lab(\ub_1,u)\nonumber\\
&&\hspace{33mm}\leq C\sqrt{D}\ub_1^{a_m+\frac{1}{2}}u^{b_m}
\label{14.176}
\end{eqnarray}
(see \ref{14.52}), while, comparing with \ref{14.71}, \ref{14.72} with $m-1$ in the role of $m$, 
the 3rd of \ref{14.172} is bounded by:
\begin{equation}
C\sqrt{D}\ub_1^{a_m+\frac{1}{2}}u^{b_m+1}
\label{14.177}
\end{equation}
As for the 4th of \ref{12.172}, we compare with \ref{14.69}, \ref{14.70} with $m-1$ in the role of $m$, 
to which the $j=0$ term in the sum (see \ref{12.567}) is added. The contribution of this term 
$$C\int_0^{\ub_1}\ub\|\s^{(m-1,n-m+1)}\cla\|_{L^2(S_{\ub,u})}d\ub$$
is by Proposition 14.1 bounded by:
$$C\sqrt{D}\int_0^{\ub_1}\ub^{a_m+\frac{3}{2}}u^{b_m-\frac{3}{2}}d\ub
=\frac{C\sqrt{D}}{\left(a_m+\frac{5}{2}\right)}\ub_1^{a_m+\frac{5}{2}}u^{b_m-\frac{3}{2}}$$
The 4th of \ref{14.172} is then seen to be bounded by:
\begin{equation}
\frac{C\sqrt{D}}{\left(a_m+\frac{5}{2}\right)}\ub_1^{a_m+\frac{5}{2}}u^{b_m-\frac{3}{2}}
\label{14.178}
\end{equation}
Turning to \ref{14.173} (case $m=1$), the 1st is bounded by \ref{14.3} by:
\begin{equation}
C\sqrt{D}\ub_1^{a_0+\frac{3}{2}}u^{b_0+1}
\label{14.179}
\end{equation}
while the 2nd is bounded by Proposition 14.1 by:
\begin{equation}
\frac{C\sqrt{D}}{\left(a_0+\frac{7}{2}\right)}\ub_1^{a_0+\frac{7}{2}}u^{b_0-\frac{3}{2}}
\label{14.180}
\end{equation}
Finally in regard to \ref{14.174} (cases $m=2,...,n$), comparing with the 3rd of \ref{14.172} with 
$m-1$ in the role of $m$ 
(see \ref{14.177}), the 1st is bounded by:
\begin{equation}
C\sqrt{D}\ub_1^{a_{m-1}+\frac{3}{2}}u^{b_{m-1}+1}
\label{14.181}
\end{equation}
and comparing with the 4th of \ref{14.172} with $m-1$ in the role of $m$ (see \ref{14.178}), the 2nd 
is bounded by:
\begin{equation}
\frac{C\sqrt{D}}{\left(a_{m-1}+\frac{7}{2}\right)}\ub_1^{a_{m-1}+\frac{7}{2}}u^{b_{m-1}-\frac{3}{2}}
\label{14.182}
\end{equation}
an extra $\ub$ factor in the integrant being present here. 
Combining the above we conclude that for $m=1,...,n$ the 3rd of the integrals on the right in 
\ref{14.152} is bounded by:
\begin{equation}
C\sqrt{D}\ub_1^{a_m+\frac{1}{2}}u^{b_m-\frac{1}{2}}
\label{14.183}
\end{equation}
This extends the bound \ref{14.164} to all $m=0,...,n$. 

To estimate the 2nd integral on the right in \ref{14.152} we use the fact that:
\begin{equation}
N^\mu L\s^{(m,n-m)}\check{\dot{\phi}}_\mu+\ogamma\Nb^\mu L\s^{(m,n-m)}\check{\dot{\phi}}_\mu 
=\s^{(Y;m,n-m)}\cxi_L
\label{14.184}
\end{equation}
Since $\ogamma\sim u$ the contribution to this integral of the 2nd term on the left is bounded by 
$Cu$ times the bound for the 3rd integral, that is, by:
\begin{equation}
C\sqrt{D}\ub_1^{a_m+\frac{1}{2}}u^{b_m+\frac{1}{2}}
\label{14.185}
\end{equation}
while the contribution of the right hand side is bounded by:
\begin{eqnarray}
&&\int_0^{\ub_1}\|\s^{(Y;m,n-m)}\cxi_L\|_{L^2(S_{\ub,u})}d\ub
\leq \ub_1^{1/2}\|\s^{(Y;m,n-m)}\cxi_L\|_{L^2(C_u^{\ub_1})}\nonumber\\
&&\hspace{10mm}\leq C\ub_1^{1/2}\sqrt{\s^{(Y;m,n-m)}\cE^{\ub_1}(u)}\leq C\ub_1^{a_m+\frac{1}{2}}u^{b_m}
\sqrt{\s^{(Y;m,n-m)}\cB(\ub_1,u)}\nonumber\\
&&\hspace{40mm}\leq C\sqrt{D}\ub_1^{a_m+\frac{1}{2}}u^{b_m}
\label{14.186}
\end{eqnarray}
by the top order energy estimates. We conclude that the 2nd integral on the right in \ref{14.152} 
is bounded by:
\begin{equation}
C\sqrt{D}\ub_1^{a_m+\frac{1}{2}}u^{b_m}
\label{14.187}
\end{equation}

The bounds \ref{14.153}, \ref{14.187}, \ref{14.183} (including, for $m=0$, \ref{14.164}), for the 
three terms on the right in \ref{14.152} yield the following proposition. 

\vspace{2.5mm}

\noindent{\bf Proposition 14.3} \ \ \ 
The $n$th order variation differences $\s^{(m,n-m)}\check{\dot{\phi}}_\mu$ 
satisfy the following $L^2$ estimates on the $S_{\ub,u}$:
$$\|\s^{(m,n-m)}\check{\dot{\phi}}_\mu\|_{L^2(S_{\ub,u})}\leq 
C\sqrt{D}\ub^{a_m+\frac{1}{2}}u^{b_m-\frac{1}{2}} \  \ \mbox{: for $m=0,...,n$}$$

\vspace{2.5mm}

In addition to the above  estimates for $n$th order variation differences, we shall need 
$L^2$ estimates for the $n$th order quantities:
\begin{equation}
q_{1,m,\mu}=E^{n-m}T^{m-1}\Lb\beta_\mu-E_N^{n-m}T^{m-1}\Lb_N\beta_{\mu,N} \ \ : \ m=1,...,n
\label{14.188}
\end{equation}
Since $L=T-\Lb$, $L_N=T-\Lb_N$, we can express:
\begin{equation}
E^{n-m}T^{m-1}L\beta_\mu-E_N^{n-m}T^{m-1}L_N\beta_{\mu,N}=\s^{(m,n-m)}\check{\dot{\phi}}_\mu 
-q_{1,m,\mu} \ \ : \ m=1,...,n 
\label{14.189}
\end{equation}
therefore $L^2(S_{\ub,u})$ estimates for the $n$th order differences 
\begin{equation}
\qb_{1,m,\mu}=E^{n-m}T^{m-1}L\beta_\mu-E_N^{n-m}T^{m-1}L_N\beta_{\mu,N} \ \ : \ m=1,...,n
\label{14.190}
\end{equation}
will immediately follow. 

To derive the $L^(S_{\ub,u})$ estimates for the $q_{1,m,\mu}$ we apply \ref{14.145} with $q_{1,m,\mu}$ 
in the role of $f$. We then apply \ref{14.150} with $Lq_{1,m,\mu}$ in the role of $\theta_\mu$ to 
arrive at the following analogue of \ref{14.152}:
\begin{eqnarray}
&&\|q_{1,m,\mu}\|_{L^2(S_{\ub_1,u})}\leq C\int_0^{\ub_1}\|E^\nu Lq_{1,m,\nu}\|_{L^2(S_{\ub,u})}d\ub
\nonumber\\
&&\hspace{29mm}+C\int_0^{\ub_1}\|N^\nu Lq_{1,m,\nu}\|_{L^2(S_{\ub,u})}d\ub 
\nonumber\\
&&\hspace{29mm}+C\int_0^{\ub_1}\|\Nb^\nu Lq_{1,m,\nu}\|_{L^2(S_{\ub,u})}d\ub
\nonumber\\
&&\label{14.191}
\end{eqnarray}
Consider then $Lq_{1,m,\mu}$. By the definition \ref{14.188} up to a 0th order term this is:
\begin{equation}
LE^{n-m}T^{m-1}\Lb\beta_\mu-L_N E_N^{n-m}T^{m-1}\Lb_N\beta_{\mu,N}
\label{14.192}
\end{equation}
Commuting here $L$ with $E^{n-m}T^{m-1}$ leaves only an order $n-1$ remainder, which is ignorable. 
Hence \ref{14.192} can be replaced by:
\begin{equation}
E^{n-m}T^{m-1}L\Lb\beta_\mu-E_N^{n-m}T^{m-1}L_N\Lb_N\beta_{\mu,N}
\label{14.193}
\end{equation}
By \ref{9.205} and the 1st of \ref{9.230} with $\beta_\mu$ in the role of $f$ we have, in view of the 
fact that $\Omega a\square_{\tilde{h}}\beta_\mu=0$,
\begin{eqnarray}
&&L\Lb\beta_\mu=aE^2\beta_\mu-\frac{1}{2}\chib L\beta_\mu-\frac{1}{2}\chi\Lb\beta_\mu+2\etab E\beta_\mu 
\nonumber\\
&&\hspace{12mm}-\frac{1}{4}((L\log\Omega)\Lb\beta_\mu+(\Lb\log\Omega)L\beta_\mu)
+\frac{1}{2}a(E\log\Omega)(E\beta_\mu) \nonumber\\
&&\label{14.194}
\end{eqnarray}
Similarly, by \ref{9.210} and the 2nd of \ref{9.230} with $\beta_{\mu,N}$ in the role of $f$ we have, 
in view of Proposition 9.10, 
\begin{eqnarray}
&&L_N\Lb_N\beta_{\mu,N}=a_N E_N^2\beta_{\mu,N}-\frac{1}{2}\chib_N L_N\beta_{\mu,N}
-\frac{1}{2}\chi_N\Lb_N\beta_{\mu,N}+2\etab_N E_N\beta_{\mu,N} 
\nonumber\\
&&\hspace{21mm}-\frac{1}{4}((L_N\log\Omega_N)\Lb_N\beta_{\mu,N}
+(\Lb_N\log\Omega_N)L_N\beta_{\mu,N})\nonumber\\
&&\hspace{21mm}+\frac{1}{2}a_N(E_N\log\Omega_N)(E_N\beta_{\mu,N})
-\tilde{\kappa}^\prime_{\mu,N} \label{14.195}
\end{eqnarray}
Taking into account the fact that by \ref{3.a22}, \ref{3.a25} and \ref{12.64} 
\begin{equation}
[\etab]_{P.A.}=\rhob\mub
\label{14.196}
\end{equation}
it follows that \ref{14.193} can be replaced up to ignorable terms by:
\begin{eqnarray}
&&a E^{n+1}\beta_\mu-a_N E_N^{n+1}\beta_{\mu,N}\nonumber\\
&&-\frac{1}{2}(\rhob L\beta_\mu)\s^{(n-1)}\ctchib
-\frac{1}{2}(\rho\Lb\beta_\mu)\s^{(n-1)}\ctchi \nonumber\\
&&+2\rhob(E\beta_\mu)\cmub_{0,n-1} \label{14.197}
\end{eqnarray}
for $m=1$, and by:
\begin{eqnarray}
&&a E^{n-m}T^{m-1}E^2\beta_\mu-a_N E_N^{n-m}T^{m-1}E_N^2\beta_{\mu,N}\nonumber\\
&&-(\rhob L\beta_\mu)\cmub_{m-2,n-m+1}
-(\rho\Lb\beta_\mu)\cmu_{m-2,n-m+1}\nonumber\\
&&+2\rhob(E\beta_\mu)\cmub_{m-1,n-m} \label{14.198}
\end{eqnarray}
for $m=2,...,n$. Moreover for $m=2,...,n$ we must commute $T^{m-1}$ with $E^2$ in the 1st term 
and with $E_N^2$ in the 2nd. By the commutation relations \ref{3.a14} and \ref{9.a1} the first 
two terms in \ref{14.198} then become:
\begin{eqnarray}
&&a E^{n-m+2}T^{m-1}\beta_\mu-a_N E_N^{n-m+2}T^{m-1}\beta_{\mu,N}\nonumber\\
&&-a(E\beta_\mu)\cdot\left\{\begin{array}{lll}(\rho\s^{(n-1)}\ctchi+\rhob\s^{(n-1)}\ctchib) &:& 
\mbox{for $m=2$}\\2(\rho\cmu_{m-3,n-m+2}+\rhob\cmub_{m-3,n-m+2}) &:& \mbox{for $m=3,...,n$}
\end{array}\right.\nonumber\\
&&\label{14.199}
\end{eqnarray}
up to ignorable terms. We shall first estimate the contribution of the remainders, that is what 
remains in each of \ref{14.197}, \ref{14.198}, \ref{14.199} after the principal difference 
constituted by the first two terms has been subtracted, to the right hand side of \ref{14.191} or 
equivalently to
\begin{equation}
\int_0^{\ub_1}\|Lq_{1,m,\mu}\|_{L^2(S_{\ub,u})}d\ub
\label{14.200}
\end{equation}
The contribution to \ref{14.200} of the remainder in \ref{14.197} is bounded by:
\begin{eqnarray}
&&Cu^2\int_0^{\ub_1}\|\rho\s^{(n-1)}\ctchib\|_{L^2(S_{\ub,u})}d\ub
+C\int_0^{\ub_1}\|\s^{(n-1)}\ctchi\|_{L^2(S_{\ub,u})}\nonumber\\
&&+Cu^2\int_0^{\ub_1}\|\cmub_{0,n-1}\|_{L^2(S_{\ub,u})}d\ub
\leq C\sqrt{D}\ub_1^{a_0+\frac{1}{2}}u^{b_0+1}\label{14.201}
\end{eqnarray}
by \ref{14.3}, Proposition 14.1, and \ref{14.71}, \ref{14.72} with $m=0$. 
The contribution to \ref{14.200} of the remainder in \ref{14.198} is bounded by:
\begin{eqnarray}
&&Cu^2\int_0^{\ub_1}\ub\|\cmub_{m-2,n-m+1}\|_{L^2(S_{\ub,u})}d\ub
+Cu^2\int_0^{\ub_1}\|\cmub_{m-1,n-m}\|_{L^2(S_{\ub,u})}d\ub \nonumber\\
&&+C\int_0^{\ub_1}\ub\|\cmu_{m-2,n-m+1}\|_{L^2(S_{\ub,u})}d\ub
\leq C\sqrt{D}\ub_1^{a_{m-1}+\frac{1}{2}}u^{b_{m-1}+\frac{1}{2}}\label{14.202}
\end{eqnarray}
by \ref{14.71}, \ref{14.72} with $m-2$ and $m-1$ in the role of $m$ and the fact that the last 
integral is similar to \ref{14.172}, which is bounded by \ref{14.178}, with $m-1$ in the role of $m$. 
The contribution to \ref{14.200} of the remainder in \ref{14.199} is bounded in the case $m=2$ by:
\begin{eqnarray}
&&Cu^2\int_0^{\ub_1}\ub^2\|\s^{(n-1)}\ctchi\|_{L^2(S_{\ub,u})}d\ub
+Cu^4\int_0^{\ub_1}\|\rho\s^{(n-1)}\ctchib\|_{L^2(S_{\ub,u})}d\ub\nonumber\\
&&\hspace{48mm}\leq C\ub_1^{a_0+\frac{3}{2}}u^{b_0+\frac{5}{2}}
\label{14.203}
\end{eqnarray}
by Proposition 14.1 and \ref{14.3}, and in the cases $m=3,...,n$ by:
\begin{eqnarray}
&&Cu^2\int_0^{\ub_1}\ub^2\|\cmu_{m-3,n-m+2}\|_{L^2(S_{\ub,u})}d\ub
+Cu^4\int_0^{\ub_1}\ub\|\cmub_{m-3,n-m+2}\|_{L^2(S_{\ub,u})}d\ub\nonumber\\
&&\hspace{48mm}\leq C\sqrt{D}\ub_1^{a_{m-2}+\frac{3}{2}}u^{b_{m-2}+\frac{5}{2}}
\label{14.204}
\end{eqnarray}
by \ref{14.178} with $m-2$ in the role of $m$ and \ref{14.72} with $m-3$ in the role of $m$. 

We reduce in this way the problem of estimating the right hand side of \ref{14.191} to that 
of estimating the three integrals:
\begin{eqnarray}
&&\int_0^{\ub_1}\|aE^\mu E(E^{n+1-m}T^{m-1}\beta_\mu-E_N^{n+1-m}T^{m-1}\beta_{\mu,N})\|_{L^2(S_{\ub,u})}
d\ub, \nonumber\\
&&\int_0^{\ub_1}\|aN^\mu E(E^{n+1-m}T^{m-1}\beta_\mu-E_N^{n+1-m}T^{m-1}\beta_{\mu,N})\|_{L^2(S_{\ub,u})}
d\ub, \nonumber\\
&&\int_0^{\ub_1}\|a\Nb^\mu E(E^{n+1-m}T^{m-1}\beta_\mu-E_N^{n+1-m}T^{m-1}\beta_{\mu,N})\|_{L^2(S_{\ub,u})}
d\ub \nonumber\\
&&\label{14.205}
\end{eqnarray}
corresponding to the difference constituted by the first two terms in \ref{14.197} and \ref{14.199}. 
The first integral is:
\begin{eqnarray}
&&\int_0^{\ub_1}\|a\s^{(E;m-1,n-m+1)}\csxi\|_{L^2(S_{\ub,u})}d\ub \nonumber\\
&&\leq Cu\int_0^{\ub_1}\ub^{1/2}\|a^{1/2}\s^{(E;m-1,n-m+1)}\csxi\|_{L^2(S_{\ub,u})}d\ub\nonumber\\
&&\leq Cu\ub_1\|a^{1/2}\s^{E;m-1,n-m+1)}\csxi\|_{L^2(C_u^{\ub_1})}\leq Cu\ub_1\sqrt{\s^{(E;m-1,n-m+1)}\cE^{\ub_1}(u)}\nonumber\\
&&\leq C\sqrt{D}\ub_1^{a_{m-1}+1}u^{b_{m-1}+1} \label{14.206}
\end{eqnarray}
by the top order energy estimates. 

To estimate the second integral we write $a N^\mu=\lambda L^\mu$ and use the fact that 
$L^\mu E\beta_\mu=E^\mu L\beta_\mu$. This requires that we first commute one $E$ factor with $T^{m-1}$, 
then express $\lambda L^\mu E^{n+1-m}T^{m-1}E\beta_\mu$ in terms of 
$\lambda E^{n+1-m}T^{m-1}(L^\mu E\beta_\mu)$, then use the fact just mentioned, then express 
$\lambda E^{n+1-m}T^{m-1}(E^\mu L\beta_\mu)$ in terms of $\lambda E^\mu E^{n+1-m}T^{m-1}L\beta_\mu$, 
and finally commute $E^{n+1-m}T^{m-1}$ with $L$ in the last to obtain 
$\lambda E^\mu LE^{n+1-m}T^{m-1}\beta_\mu$. A similar procedure applied to the $N$th approximants 
will give us an expression for $a_N N_N^\mu E_N^{n+2-m}T^{m-1}\beta_{\mu,N}$ in terms of 
$\lambda_N E_N^\mu L_N E_N^{n+1-m}T^{m-1}\beta_{\mu,N}$. Subtracting will then give an expression 
for 
$$aN^\mu E(E^{n+1-m}T^{m-1}\beta_\mu-E_N^{n+1-m}T^{m-1}\beta_{\mu,N}) \ \mbox{in terms of} \ 
\lambda \s^{(E;m-1,n-m+1)}\cxi_L$$ 
up to 0th order difference terms. 

The first commutation is only present for $m=2,...,n$ and leaves then an $n$th order remainder 
with a principal acoustical part of 
\begin{eqnarray}
&&a\ss_N[E^{n-m+1}T^{m-2}(\chi+\chib)]_{P.A.}\nonumber\\
&&=a\ss_N\left\{\begin{array}{lll}\rho\s^{(n-1)}\ctchi+\rhob\s^{(n-1)}\ctchib &:& \mbox{for $m=2$}\\
2\rho\cmu_{m-3,n-m+2}+2\rhob\cmub_{m-3,n-m+2} &:& \mbox{for $m=3,...,n$}\end{array}\right. \nonumber\\
&&\label{14.207}
\end{eqnarray}
once differences from the $N$th approximants are taken. 
This is similar to the remainder in \ref{14.199} hence its contribution to the second integral 
\ref{14.205} is bounded by:
\begin{equation}
C\sqrt{D}\ub_1^{a_{m-2}+\frac{3}{2}}u^{b_{m-2}+\frac{5}{2}} \ \mbox{: for all $m=2,...,n$}
\label{14.208}
\end{equation}
(see \ref{14.203}, \ref{14.204}). Expressing $\lambda L^\mu E^{n+1-m}T^{m-1}E\beta_\mu$ in terms of \\
$\lambda E^{n+1-m}T^{m-1}(L^\mu E\beta_\mu)$ will leave an $n$th order remainder with a principal acoustical part of:
\begin{eqnarray}
&&\lambda(E\beta_\mu)[E^{n-m+1}T^{m-1}(\rho N^\mu)]_{P.A.}\nonumber\\
&&=\rhob\left\{\begin{array}{lll}\ss_N\cmub_{0,n-1}+\sss\lambdab\s^{(n-1)}\ctchi &:& \mbox{for $m=1$}\\
\ss_N\s^{(m-1,n-m+1)}\clab &\s & \s \\
+2\lambdab(\sss\cmu_{m-2,n-m+1}+\pi\ss_N\cmub_{m-2,n-m+1}) &:& 
\mbox{for $m=2,...,n$}\end{array}\right.\nonumber\\
&&\label{14.209}
\end{eqnarray}
once differences from the $N$th approximants are taken (see \ref{12.7}, \ref{12.10}, \ref{14.31}). 
The contribution of this remainder to the second integral \ref{14.205} is bounded by:
\begin{equation}
C\sqrt{D}\ub_1^{a_{m-1}+\frac{1}{2}}u^{b_{m-1}+2} \ \mbox{: for all $m=1,...,n$}
\label{14.210}
\end{equation}
by \ref{14.72} with $m=0$, Proposition 14.1, \ref{14.52} with $m-1$ in the role of $m$,  
\ref{14.178} with $m-1$ in the role of $m$, and \ref{14.72} with $m-2$ in the role of $m$. Expressing  
$\lambda E^{n+1-m}T^{m-1}(E^\mu L\beta_\mu)$ in terms of $\lambda E^\mu E^{n+1-m}T^{m-1}L\beta_\mu$ 
will leave an $n$th order remainder with a principal acoustical part of:
\begin{eqnarray}
&&\lambda(L\beta_\mu)[E^{n-m+1}T^{m-1}E^\mu]_{P.A.}\nonumber\\
&&=\frac{1}{2}\rhob\left\{\begin{array}{lll}s_{NL}\s^{(n-1)}\ctchib+\sss\lambdab\s^{(n-1)}\ctchi &:& \mbox{for $m=1$}\\2s_{NL}\cmub_{m-2,n-m+1}+2\sss\lambdab\cmu_{m-2,n-m+1} &:& \mbox{for $m=2,...,n$}
\end{array}\right.\nonumber\\
&&\label{14.211}
\end{eqnarray}
once differences from the $N$th approximants are taken (see \ref{12.8}, \ref{12.74}). The contribution 
of this remainder to the second integral \ref{14.205} is bounded by:
\begin{equation}
C\sqrt{D}\cdot\left\{\begin{array}{lll}\ub_1^{a_0+\frac{3}{2}}u^{b_0+1} &:& \mbox{for $m=1$}\\
\ub_1^{a_{m-1}+\frac{3}{2}}u^{b_{m-1}+\frac{3}{2}} &:& \mbox{for $m=2,...,n$}\end{array}\right.
\label{14.212}
\end{equation}
by \ref{14.3} (recall that $|s_{NL}|\leq C\rho$), Proposition 14.1, \ref{14.178} with $m-1$ in the 
role of $m$, and \ref{14.72} with $m-2$ in the role of $m$. Finally, the second commutation leaves 
an $n$th order remainder, namely
\begin{equation}
\lambda E^\mu[E^{n-m+1}T^{m-1},L]\beta_\mu 
\label{14.213}
\end{equation}
The principal acoustical part of this is:
\begin{equation}
\lambda\sss\cdot\left\{\begin{array}{lll}\rho\s^{(n-1)}\ctchi &:& \mbox{for $m=1$}\\
2\rho\cmu_{m-2,n-m+1}-2\rhob\cmub_{m-2,n-m+1} &:& \mbox{for $m=2,...,n$}\end{array}\right.
\label{14.214}
\end{equation}
once differences from the $N$th approximants are taken. The contribution of this remainder to 
the second integral \ref{14.205} is bounded by:
\begin{equation}
C\sqrt{D}\cdot\left\{\begin{array}{lll}\ub_1^{a_0+\frac{5}{2}}u^{b_0+\frac{1}{2}} &:& \mbox{for $m=1$}\\
\ub_1^{a_{m-1}+\frac{1}{2}}u^{b_{m-1}+\frac{5}{2}} &:& \mbox{for $m=2,...,n$}\end{array}\right.
\label{14.215}
\end{equation}
by Proposition 14.1, \ref{14.178} with $m-1$ in the role of $m$, and \ref{14.72} with $m-2$ in the 
role of $m$. Finally, the principal contribution to the second of the integrals \ref{14.205} that of 
$\lambda\s^{(E;m-1,n-m+1)}\cxi_L$ is:
\begin{eqnarray}
&&\int_0^{\ub_1}\|\lambda\s^{(E;m-1,n-m+1)}\cxi_L\|_{L^2(S_{\ub,u})}d\ub \nonumber\\
&&\leq Cu^2\int_0^{\ub_1}\|\s^{(E;m-1,n-m+1)}\cxi_L\|_{L^2(S_{\ub,u})}d\ub\nonumber\\
&&\leq Cu^2\ub_1^{1/2}\|\s^{E;m-1,n-m+1)}\cxi_L\|_{L^2(C_u^{\ub_1})}\leq Cu^2\ub_1^{1/2}\sqrt{\s^{(E;m-1,n-m+1)}\cE^{\ub_1}(u)}\nonumber\\
&&\leq C\sqrt{D}\ub_1^{a_{m-1}+\frac{1}{2}}u^{b_{m-1}+2} \label{14.216}
\end{eqnarray}
by the top order energy estimates. Combining with the estimates \ref{14.208}, \ref{14.210}, 
\ref{14.212}, \ref{14.215} for the remainders we conclude that the second of the integrals \ref{14.205} 
is bounded by:
\begin{equation}
C\sqrt{D}\ub_1^{a_{m-1}+\frac{1}{2}}u^{b_{m-1}+2}
\label{14.217}
\end{equation}

To estimate the third of the integrals \ref{14.205} we note that this integral is:
$$\int_0^{\ub_1}\|a\Nb^\mu E\s^{(m-1,n-m+1}\check{\dot{\phi}}_\mu\|_{L^2(S_{\ub,u})}d\ub$$
while the second integral is:
$$\int_0^{\ub_1}\|a N^\mu E\s^{(m-1,n-m+1)}\check{\dot{\phi}}_\mu\|_{L^2(S_{\ub,u})}d\ub$$
and use the fact that:
\begin{equation}
N^\mu E\s^{(m-1,n-m+1)}\check{\dot{\phi}}_\mu+\ogamma\Nb^\mu E\s^{(m-1,n-m+1)}\check{\dot{\phi}}_\mu
=\s^{(Y;m-1,n-m+1)}\sxi
\label{14.218}
\end{equation}
Since $\ogamma\sim u$, the contribution to the third integral of the 1st term on the left is bounded 
by $Cu^{-1}$ times the bound for the second integral, that is by: 
\begin{equation}
C\sqrt{D}\ub_1^{a_{m-1}+\frac{1}{2}}u^{b_{m-1}+1}
\label{14.219}
\end{equation} 
while the contribution of the right hand side is bounded by:
\begin{eqnarray}
&&Cu^{-1}\int_0^{\ub_1}\|a\s^{(Y;m-1,n-m+1)}\csxi\|_{L^2(S_{\ub,u})}d\ub \nonumber\\
&&\leq C\int_0^{\ub_1}\ub^{1/2}\|a^{1/2}\s^{(Y;m-1,n-m+1)}\csxi\|_{L^2(S_{\ub,u})}d\ub\nonumber\\
&&\leq C\ub_1\|a^{1/2}\s^{Y;m-1,n-m+1)}\csxi\|_{L^2(C_u^{\ub_1})}\leq C\ub_1\sqrt{\s^{(Y;m-1,n-m+1)}\cE^{\ub_1}(u)}\nonumber\\
&&\leq C\sqrt{D}\ub_1^{a_{m-1}+1}u^{b_{m-1}} \label{14.220}
\end{eqnarray}
by the top order energy estimates. 

Combining the results \ref{14.206}, \ref{14.217}, \ref{14.220} for the three integrals on the 
right in \ref{14.205} we arrive at the following lemma. 

\vspace{2.5mm}

\noindent{\bf Lemma 14.1} \ \ \ The $n$th order quantities $q_{1,m,\mu}$ satisfy the following 
$L^2$ estimates on the $S_{\ub,u}$:
$$\|q_{1,m,\mu}\|_{L^2(S_{\ub,u})}\leq C\sqrt{D} \ub^{a_{m-1}+\frac{1}{2}}u^{b_{m-1}+\frac{1}{2}}
\  \ \mbox{: for $m=1,...,n$}$$ 

\vspace{5mm}

\section{$L^2(S_{\ub,u})$ Estimates for the $n-1$th Order Acoustical Differences}

In the present section we shall derive $L^2(S_{\ub,u})$ estimates for the $n-1$th order 
acoustical differences  on the basis of the $L^2({\cal K}_\sigma^\tau)$ estimates for the $n$th order 
acoustical differences of Section 14.3. 

We begin with $\s^{(n-2)}\ctchi$ and $\s^{(n-2)}\ctchib$. In the following, as in Section 12.1, 
we shall denote by $R_k$ a generic term of order $k$ with vanishing order $k$ acoustical part. 
We recall the 2nd variation equation \ref{12.13}:
\begin{equation}
L\tchi=A\tchi-\rho\tchi^2+\oA\tchib-\beta_N\sbeta ELH+\frac{1}{2}\rho\beta_N^2 E^2 H+R_1
\label{14.221}
\end{equation} 
where the coefficients $A$ and $\oA$ are given by \ref{12.14}, \ref{12.15}. We turn to the cross 
variation equation \ref{10.199}. Substituting for $\sk$ in terms of $\tchi$ 
from \ref{3.a21}, the left hand side becomes:
\begin{eqnarray}
&&\Lb\sk=\Lb\tchi-H\sbeta N^\mu E\Lb\beta_\mu-H\left[(\beta_{\Nb}-2\pi c\sbeta)\frac{\ss_N}{c}
+2\sbeta\sss\right](E\lambda+\pi\lambda\tchi) \nonumber\\
&&\hspace{8mm}+\rhob H\sbeta\ss_N\tchib+R_1
\label{14.222}
\end{eqnarray}
Substituting also for the structure functions on the right hand side from Chapter 3, 
we bring the cross variation equation \ref{10.199} to the form:
\begin{eqnarray}
&&\Lb\tchi=2E^2\lambda+2\pi \lambda E\tchi
-\beta_N\sbeta E\Lb H+\rhob\frac{\beta_N}{2}(\beta_N+2\beta_{\Nb})E^2 H\nonumber\\
&&\hspace{8mm}+H\sbeta N^\mu E\Lb\beta_\mu+\rhob H(\beta_N-2\pi c\sbeta)N^\mu E^2\beta_\mu\nonumber\\
&&\hspace{8mm}+\rhob\tchi^2+B E\lambda+\tilde{A}\tchi+\tilde{\oA}\tchib+R_1 
\label{14.223}
\end{eqnarray}
where:
\begin{equation}
B=\left[\frac{\beta_N}{c}(\beta_N+\beta_{\Nb})+\sbeta^2\right] EH+4H\sbeta\sss
+\frac{H}{c}\left(\frac{3}{2}\beta_N+\beta_{\Nb}-3\pi c\sbeta\right)\ss_N 
\label{14.224}
\end{equation} 
\begin{eqnarray}
&&\tilde{A}=\frac{\beta_N}{4c}(4\pi c\sbeta-\beta_N-2\beta_{\Nb})\Lb H
+\frac{\rhob}{2}\left[4(\beta_N+\beta_{\Nb})\sbeta+\pi\beta_N(\beta_N-2\beta_{\Nb})\right] EH
\nonumber\\
&&\hspace{6mm}+\frac{\rhob}{2} H(-2\pi c\sbeta+3\beta_N+2\beta_{\Nb})\sss
+\frac{\pi\rhob}{2}H(\beta_N-2\pi c\sbeta)\ss_N\nonumber\\
&&\hspace{6mm}+\rhob H\left[(1-4\pi^2 c)\sbeta+\pi(2\beta_{\Nb}-\beta_N)\right]\ss_{\Nb}
\label{14.225}
\end{eqnarray}
\begin{equation}
\tilde{\oA}=-\frac{\beta_N^2}{4c}\Lb H+\frac{\rhob\beta_N}{2}(2\sbeta+\pi\beta_N)EH
+\frac{\rhob}{2}H\left[-2(1+3\pi^2 c)\sbeta+3\pi\beta_N\right]\ss_N 
\label{14.226}
\end{equation}
Adding equations \ref{14.221} and \ref{14.223} yields:
\begin{eqnarray}
&&T\tchi=2E^2\lambda+2\pi\lambda E\tchi
-\beta_N\sbeta ETH+\left[\frac{1}{2}(\rho+\rhob)\beta_N^2+\rhob\beta_N\beta_{\Nb}\right]E^2 H
\nonumber\\
&&\hspace{8mm}+H\sbeta N^\mu E\Lb\beta_\mu+\rhob H(\beta_N-2\pi c\sbeta)N^\mu E^2\beta_\mu \nonumber\\
&&\hspace{8mm}+(\rhob-\rho)\tchi^2+B E\lambda+A^\prime\tchi+\oA^\prime\tchib+R_1 
\label{14.227}
\end{eqnarray}
where:
\begin{equation}
A^\prime=A+\tilde{A}, \ \ \ \oA^\prime=\oA+\tilde{\oA}
\label{14.228}
\end{equation}
We recall from \ref{12.16} that an equation analogous to \ref{14.221} is satisfied by the $N$th 
approximants, with a error term $\vep_{L\tchi,N}$ satisfying \ref{12.19}. Similarly from \ref{10.208} 
we obtain an equation analogous to \ref{14.223} for the $N$th approximants, with an error term 
$\vep_{\Lb\tchi,N}=O(\tau^{N+1})$. Therefore an equation analogous to \ref{14.227} is satisfied by 
the $N$th approximants, with an error term:
\begin{equation}
\vep_{T\chi,N}=\vep_{L\chi,N}+\vep_{\Lb\chi,N}=O(\tau^N)
\label{14.229}
\end{equation}
Applying $E^{n-2}$ to \ref{14.227} and using the fact that 
\begin{eqnarray}
&&[T, E^{n-2}]\tchi=-\sum_{i=0}^{n-3}E^{n-3-i}\left((\chi+\chib)E^{i+1}\tchi\right)\nonumber\\
&&\hspace{17mm}=-(n-2)(\chi+\chib)E^{n-2}\tchi+\mbox{terms of order $n-2$}\nonumber\\
&&\label{14.230}
\end{eqnarray}
we obtain:
\begin{eqnarray}
&&TE^{n-2}\tchi=2E^n\lambda+2\pi\lambda E^{n-1}\tchi\nonumber\\
&&\hspace{16mm}-\beta_N\sbeta E^{n-1}TH
+\left[\frac{1}{2}(\rho+\rhob)\beta_N^2+\rhob\beta_N\beta_{\Nb}\right]E^n H\nonumber\\
&&\hspace{16mm}+H\sbeta N^\mu E^{n-1}\Lb\beta_\mu+\rhob H(\beta_N-2\pi c\sbeta)N^\mu E^n\beta_\mu \nonumber\\
&&\hspace{16mm}+B E^{n-1}\lambda
+\left[A^\prime+2(\rhob-\rho)\tchi-(n-2)(\chi+\chib)\right]E^{n-2}\tchi\nonumber\\
&&\hspace{16mm}+\oA^\prime E^{n-2}\tchib+R_{n-1}
\label{14.231}
\end{eqnarray} 
Applying $E_N^{n-2}$ to the $N$th approximant version of \ref{14.227} an analogous equation is deduced 
for the $N$th approximants, with an error term:
\begin{equation}
E_N^{n-2}\vep_{T\chi,N}=O(\tau^N)
\label{14.232}
\end{equation}
Subtracting the two equations then yields:
\begin{eqnarray}
&&T\s^{(n-2)}\ctchi=2\s^{(0,n)}\cla+2\pi\lambda \s^{(n-1)}\ctchi\nonumber\\
&&\hspace{16mm}+2\beta_N\sbeta H^\prime\beta^\mu\s^{(1,n-1)}\check{\dot{\phi}}_\mu
-\left[(\rho+\rhob)\beta_N^2+2\rhob\beta_N\beta_{\Nb}\right]H^\prime\beta^\mu
\s^{(0,n)}\check{\dot{\phi}}_\mu\nonumber\\
&&\hspace{16mm}+H\sbeta N^\mu q_{1,1,\mu}+\rhob H(\beta_N-2\pi c\sbeta)N^\mu
\s^{(0,n)}\check{\dot{\phi}}_\mu \nonumber\\
&&\hspace{16mm}+B\s^{(0,n-1)}\cla
+\left[A^\prime+2(\rhob-\rho)\tchi-(n-2)(\chi+\chib)\right]\s^{(n-2)}\ctchi\nonumber\\
&&\hspace{16mm}+\oA^\prime\s^{(n-2)}\ctchib\label{14.233}
\end{eqnarray}
up to ignorable terms. In a similar manner we deduce the conjugate equation:
\begin{eqnarray}
&&T\s^{(n-2)}\ctchib=2\s^{(0,n)}\clab+2\pi\lambdab \s^{(n-1)}\ctchib\nonumber\\
&&\hspace{16mm}+2\beta_{\Nb}\sbeta H^\prime\beta^\mu\s^{(1,n-1)}\check{\dot{\phi}}_\mu
-\left[(\rho+\rhob)\beta_{\Nb}^2+2\rho\beta_N\beta_{\Nb}\right]H^\prime\beta^\mu
\s^{(0,n)}\check{\dot{\phi}}_\mu\nonumber\\
&&\hspace{16mm}+H\sbeta\Nb^\mu\qb_{1,1,\mu}+\rho H(\beta_{\Nb}-2\pi c\sbeta)\Nb^\mu
\s^{(0,n)}\check{\dot{\phi}}_\mu \nonumber\\
&&\hspace{16mm}+\Bb\s^{(0,n-1)}\clab
+\left[\oAb^\prime+2(\rho-\rhob)\tchib-(n-2)(\chi+\chib)\right]\s^{(n-2)}\ctchib\nonumber\\
&&\hspace{16mm}+\Ab^\prime\s^{(n-2)}\ctchi\label{14.234}
\end{eqnarray}
up to ignorable terms. 
Here, the coefficients $\Bb$, $\oAb^\prime$, $\Ab^\prime$ are the conjugates of the coefficients 
$B$, $A^\prime$, $\oA^\prime$ respectively.

Consider next $\s^{(0,n-1)}\cla$ and $\s^{(0,n-1)}\clab$. We have:
\begin{equation}
TE^{n-1}\lambda=E^{n-1}T\lambda-\sum_{i=0}^{n-2}E^{n-2-i}\left((\chi+\chib)E^{i+1}\lambda\right)
\label{14.235}
\end{equation}
Now the highest order, $n-1$, terms in the sum are:
\begin{eqnarray*}
&&\mbox{from $i=0$:} \ \ \ (E\lambda)E^{n-2}(\chi+\chib)=(E\lambda)(\rho E^{n-2}\tchi+\rhob E^{n-2}\tchib)
+R_{n-1}\\
&&\mbox{from all $i=0,...,n-2$:} \ \ \ (n-1)(\chi+\chib)E^{n-1}\lambda
\end{eqnarray*}
Therefore \ref{14.235} is:
\begin{equation}
TE^{n-1}\lambda=E^{n-1}T\lambda-(E\lambda)(\rho E^{n-2}\tchi+\rhob E^{n-2}\tchib)
-(n-1)(\chi+\chib)E^{n-1}\lambda+R_{n-1}
\label{14.236}
\end{equation}
A similar formula holds for the $N$th approximants. Subtracting the two then yields:
\begin{equation}
T\s^{(0,n-1)}\cla=\s^{(1,n-1)}\cla-(E\lambda)(\rho\s^{(n-2)}\ctchi+\rhob\s^{(n-2)}\tchib)
-(n-1)(\chi+\chib)\s^{(0,n-1)}\cla
\label{14.237}
\end{equation}
up to ignorable terms. In a similar manner we deduce the conjugate equation:
\begin{equation}
T\s^{(0,n-1)}\clab=\s^{(1,n-1)}\clab-(E\lambdab)(\rho\s^{(n-2)}\ctchi+\rhob\s^{(n-2)}\tchib)
-(n-1)(\chi+\chib)\s^{(0,n-1)}\clab
\label{14.238}
\end{equation}
up to ignorable terms. 

Now, an analogue of the inequality \ref{12.366} holds with any ${\cal K}_\sigma$ in the role of 
${\cal K}$. That is, of $f$ is a function on ${\cal N}$ vanishing on $\Cb_0$, so that in terms of 
$(\tau,\sigma,\vartheta)$ (see \ref{9.1}) coordinates we have:
\begin{equation}
f(\tau,\sigma,\vartheta)=\int_0^{\tau}(Tf)(\tau^\prime,\sigma,\vartheta)d\tau^\prime
\label{14.239}
\end{equation}
we have:
\begin{eqnarray}
&&\|f\|_{L^2(S_{\tau,\sigma+\tau})}\leq k\int_0^{\tau}\|Tf\|_{L^2(S_{\tau^\prime,\sigma+\tau^\prime})}
d\tau^\prime\nonumber\\
&&\hspace{20mm}\leq k\tau^{1/2}\|Tf\|_{L^2({\cal K}_\sigma^\tau)}
\label{14.240}
\end{eqnarray}
Applying this with $\s^{(n-2)}\ctchi$, $\s^{(n-2)}\ctchib$ in the role of $f$, by virtue of equations 
\ref{14.233}, \ref{14.234} we obtain, in view of the fact that by the bootstrap assumptions the 
coefficients $B$, $\Bb$, $A^\prime$, $\oA^\prime$, $\oAb^\prime$, $\Ab^\prime$ are all bounded by 
a fixed constant, 
\begin{eqnarray}
&&\|\s^{(n-2)}\ctchi\|_{L^2(S_{\tau,\sigma+\tau})}\leq 
C\tau^{1/2}\|\s^{(0,n)}\cla\|_{L^2({\cal K}_\sigma^\tau)}
+C\tau^{1/2}(\sigma+\tau)^2\|\s^{(n-1)}\ctchi\|_{L^2({\cal K}_\sigma^\tau)}\nonumber\\
&&\hspace{35mm}+C\sum_{\mu}\int_0^{\tau}\|\s^{(1,n-1)}\check{\dot{\phi}}_\mu\|_{L^2
(S_{\tau^\prime,\sigma+\tau^\prime})}d\tau^\prime\nonumber\\
&&\hspace{35mm}+C\sum_{\mu}\int_0^{\tau}(\sigma+\tau^\prime)\|\s^{(0,n)}\check{\dot{\phi}}_\mu\|_{L^2
(S_{\tau^\prime,\sigma+\tau^\prime})}d\tau^\prime\nonumber\\
&&\hspace{35mm}+C\sum_{\mu}\int_0^{\tau}\|q_{1,1,\mu}\|_{L^2(S_{\tau^\prime,\sigma+\tau^\prime})}
d\tau^\prime 
\nonumber\\
&&\hspace{35mm}+C\int_0^{\tau}\|\s^{(0,n-1)}\cla\|_{L^2(S_{\tau^\prime,\sigma+\tau^\prime})}d\tau^\prime
\nonumber\\
&&\hspace{20mm}+C\int_0^{\tau}\|\s^{(n-2)}\ctchi\|_{L^2(S_{\tau^\prime,\sigma+\tau^\prime})}d\tau^\prime
+C\int_0^{\tau}\|\s^{(n-2)}\ctchib\|_{L^2(S_{\tau^\prime,\sigma+\tau^\prime})}d\tau^\prime\nonumber\\
&&\label{14.241}
\end{eqnarray}
\begin{eqnarray}
&&\|\s^{(n-2)}\ctchib\|_{L^2(S_{\tau,\sigma+\tau})}\leq 
C\tau^{1/2}\|\s^{(0,n)}\clab\|_{L^2({\cal K}_\sigma^\tau)}
+C\tau^{3/2}\|\s^{(n-1)}\ctchib\|_{L^2({\cal K}_\sigma^\tau)}\nonumber\\
&&\hspace{35mm}+C\sum_{\mu}\int_0^{\tau}\|\s^{(1,n-1)}\check{\dot{\phi}}_\mu\|_{L^2
(S_{\tau^\prime,\sigma+\tau^\prime})}d\tau^\prime\nonumber\\
&&\hspace{35mm}+C\sum_{\mu}\int_0^{\tau}(\sigma+\tau^\prime)\|\s^{(0,n)}\check{\dot{\phi}}_\mu\|_{L^2
(S_{\tau^\prime,\sigma+\tau^\prime})}d\tau^\prime\nonumber\\
&&\hspace{35mm}+C\sum_{\mu}\int_0^{\tau}\|\qb_{1,1,\mu}\|_{L^2(S_{\tau^\prime,\sigma+\tau^\prime})}
d\tau^\prime 
\nonumber\\
&&\hspace{35mm}+C\int_0^{\tau}\|\s^{(0,n-1)}\clab\|_{L^2(S_{\tau^\prime,\sigma+\tau^\prime})}
d\tau^\prime\nonumber\\
&&\hspace{20mm}+C\int_0^{\tau}\|\s^{(n-2)}\ctchi\|_{L^2(S_{\tau^\prime,\sigma+\tau^\prime})}d\tau^\prime
+C\int_0^{\tau}\|\s^{(n-2)}\ctchib\|_{L^2(S_{\tau^\prime,\sigma+\tau^\prime})}d\tau^\prime\nonumber\\
&&\label{14.242}
\end{eqnarray}
Also, applying \ref{14.240} with $\s^{(0,n-1)}\cla$, $\s^{(0,n-1)}\clab$ we obtain, by virtue of 
equations \ref{14.237}, \ref{14.238}, 
\begin{eqnarray}
&&\|\s^{(0,n-1)}\cla\|_{L^2(S_{\tau,\sigma+\tau})}\leq 
C\tau^{1/2}\|\s^{(1,n-1)}\cla\|_{L^2({\cal K}_\sigma^\tau)}\nonumber\\
&&\hspace{35mm}+C\int_0^{\tau}\|\s^{(0,n-1)}\cla\|_{L^2
(S_{\tau^\prime,\sigma+\tau^\prime})}d\tau^\prime\nonumber\\
&&\hspace{35mm}+C\int_0^{\tau}\tau^\prime(\sigma+\tau^\prime)^2\|\s^{(n-2)}\ctchi\|_{L^2
(S_{\tau^\prime,\sigma+\tau^\prime})}d\tau^\prime\nonumber\\
&&\hspace{35mm}+C\int_0^{\tau}(\sigma+\tau^\prime)^4\|\s^{(n-2)}\ctchib\|_{L^2
(S_{\tau^\prime,\sigma+\tau^\prime})}d\tau^\prime \nonumber\\
&&\label{14.243}
\end{eqnarray}
\begin{eqnarray}
&&\|\s^{(0,n-1)}\clab\|_{L^2(S_{\tau,\sigma+\tau})}\leq 
C\tau^{1/2}\|\s^{(1,n-1)}\clab\|_{L^2({\cal K}_\sigma^\tau)}\nonumber\\
&&\hspace{35mm}+C\int_0^{\tau}\|\s^{(0,n-1)}\clab\|_{L^2
(S_{\tau^\prime,\sigma+\tau^\prime})}d\tau^\prime\nonumber\\
&&\hspace{35mm}+C\int_0^{\tau}\tau^{\prime 2}\|\s^{(n-2)}\ctchi\|_{L^2
(S_{\tau^\prime,\sigma+\tau^\prime})}d\tau^\prime\nonumber\\
&&\hspace{35mm}+C\int_0^{\tau}\tau^\prime(\sigma+\tau^\prime)^2\|\s^{(n-2)}\ctchib\|_{L^2
(S_{\tau^\prime,\sigma+\tau^\prime})}d\tau^\prime \nonumber\\
&&\label{14.244}
\end{eqnarray}
Substituting in the above the $L^2({\cal K}_\sigma^\tau)$ estimates of Proposition 14.2 
for the $n$th order acoustical quantities, the $L^2(S_{\ub,u})$ estimates of Proposition 14.3 
for the $n$th order variation differences, and the $L^2(S_{\ub,u})$ estimates of Lemma 14.1 
for the quantities $q_{1,m,\mu}$, we obtain:
\begin{eqnarray}
&&\|\s^{(n-2)}\ctchi\|_{L^2(S_{\tau,\sigma+\tau})}\leq 
C\sqrt{D}\tau^{a_1+\frac{1}{2}}(\sigma+\tau)^{b_1+\frac{1}{2}}
+C\int_0^{\tau}\|\s^{(0,n-1)}\cla\|_{L^2(S_{\tau^\prime,\sigma+\tau^\prime})}d\tau^\prime
\nonumber\\
&&\hspace{24mm}+C\int_0^{\tau}\|\s^{(n-2)}\ctchi\|_{L^2(S_{\tau^\prime,\sigma+\tau^\prime})}d\tau^\prime
+C\int_0^{\tau}\|\s^{(n-2)}\ctchib\|_{L^2(S_{\tau^\prime,\sigma+\tau^\prime})}d\tau^\prime\nonumber\\
&&\label{14.245}
\end{eqnarray} 
\begin{eqnarray}
&&\|\s^{(n-2)}\ctchib\|_{L^2(S_{\tau,\sigma+\tau})}\leq 
C\sqrt{D}\tau^{a_1+\frac{1}{2}}(\sigma+\tau)^{b_1}+C\int_0^{\tau}\|\s^{(0,n-1)}\clab\|_{L^2(S_{\tau^\prime,\sigma+\tau^\prime})}
d\tau^\prime\nonumber\\
&&\hspace{24mm}+C\int_0^{\tau}\|\s^{(n-2)}\ctchi\|_{L^2(S_{\tau^\prime,\sigma+\tau^\prime})}d\tau^\prime
+C\int_0^{\tau}\|\s^{(n-2)}\ctchib\|_{L^2(S_{\tau^\prime,\sigma+\tau^\prime})}d\tau^\prime\nonumber\\
&&\label{14.246}
\end{eqnarray}
\begin{eqnarray}
&&\|\s^{(0,n-1)}\cla\|_{L^2(S_{\tau,\sigma+\tau})}\leq 
C\sqrt{D}\tau^{a_1+1}(\sigma+\tau)^{b_1+\frac{1}{2}}\nonumber\\
&&\hspace{35mm}+C\int_0^{\tau}\|\s^{(0,n-1)}\cla\|_{L^2
(S_{\tau^\prime,\sigma+\tau^\prime})}d\tau^\prime\nonumber\\
&&\hspace{35mm}+C\int_0^{\tau}\tau^\prime(\sigma+\tau^\prime)^2\|\s^{(n-2)}\ctchi\|_{L^2
(S_{\tau^\prime,\sigma+\tau^\prime})}d\tau^\prime\nonumber\\
&&\hspace{35mm}+C\int_0^{\tau}(\sigma+\tau^\prime)^4\|\s^{(n-2)}\ctchib\|_{L^2
(S_{\tau^\prime,\sigma+\tau^\prime})}d\tau^\prime \nonumber\\
&&\label{14.247}
\end{eqnarray}
\begin{eqnarray}
&&\|\s^{(0,n-1)}\clab\|_{L^2(S_{\tau,\sigma+\tau})}\leq 
C\sqrt{D}\tau^{a_1+\frac{1}{2}}(\sigma+\tau)^{b_1}\nonumber\\
&&\hspace{35mm}+C\int_0^{\tau}\|\s^{(0,n-1)}\clab\|_{L^2
(S_{\tau^\prime,\sigma+\tau^\prime})}d\tau^\prime\nonumber\\
&&\hspace{35mm}+C\int_0^{\tau}\tau^{\prime 2}\|\s^{(n-2)}\ctchi\|_{L^2
(S_{\tau^\prime,\sigma+\tau^\prime})}d\tau^\prime\nonumber\\
&&\hspace{35mm}+C\int_0^{\tau}\tau^\prime(\sigma+\tau^\prime)^2\|\s^{(n-2)}\ctchib\|_{L^2
(S_{\tau^\prime,\sigma+\tau^\prime})}d\tau^\prime \nonumber\\
&&\label{14.248}
\end{eqnarray}
With $\sigma$ fixed, let us define the quantities: 
\begin{eqnarray}
&&x(\tau)=\tau^{-a_1-\frac{1}{2}}(\sigma+\tau)^{-b_1-\frac{1}{2}}\|\s^{(n-2)}\ctchi\|_{L^2
(S_{\tau,\sigma+\tau})} \nonumber\\
&&\xb(\tau)=\tau^{-a_1-\frac{1}{2}}(\sigma+\tau)^{-b_1}\|\s^{(n-2)}\ctchib\|_{L^2(S_{\tau,\sigma+\tau})}
\nonumber\\
&&y(\tau)=\tau^{-a_1-1}(\sigma+\tau)^{-b_1-\frac{1}{2}}\|\s^{(0,n-1)}\cla\|_{L^2(S_{\tau,\sigma+\tau})}
\nonumber\\
&&\yb(\tau)=\tau^{-a_1-\frac{1}{2}}(\sigma+\tau)^{-b_1}\|\s^{(0,n-1)}\clab\|_{L^2(S_{\tau,\sigma+\tau})}
\nonumber\\
&&\label{14.249}
\end{eqnarray}
We then have:
\begin{eqnarray}
&&\int_0^{\tau}\|\s^{(n-2)}\ctchi\|_{L^2(S_{\tau^\prime,\sigma+\tau^\prime})}d\tau^\prime
=\int_0^{\tau}\tau^{\prime a_1+\frac{1}{2}}(\sigma+\tau^\prime)^{b_1+\frac{1}{2}}
x(\tau^\prime)d\tau^\prime\nonumber\\
&&\hspace{10mm}\leq\tau^{a_1+\frac{1}{2}}(\sigma+\tau)^{b_1+\frac{1}{2}}\int_0^{\tau}x(\tau^\prime)d\tau^\prime 
\label{14.250}
\end{eqnarray}
We also have:
\begin{equation}
\int_0^{\tau}\|\s^{(n-1)}\ctchib\|_{L^2(S_{\tau^\prime,\sigma+\tau^\prime})}d\tau^\prime
=\int_0^{\tau}\tau^{\prime a_1+\frac{1}{2}}(\sigma+\tau^\prime)^{b_1}\xb(\tau^\prime)d\tau^\prime
\label{14.251}
\end{equation}
In connection with \ref{14.245} we bound the last by:
\begin{equation}
\int_0^{\tau}\tau^{\prime a_1+\frac{1}{2}}(\sigma+\tau^\prime)^{b_1+\frac{1}{2}}
\xb(\tau^\prime)\frac{d\tau^\prime}
{\sqrt{\tau^\prime}}\leq \tau^{a_1+\frac{1}{2}}(\sigma+\tau)^{b_1+\frac{1}{2}}
\int_0^{\tau}\xb(\tau^\prime)\frac{d\tau^\prime}{\sqrt{\tau^\prime}}
\label{14.252}
\end{equation}
and in connection with \ref{14.246} simply by:
\begin{equation}
\tau^{a_1+\frac{1}{2}}(\sigma+\tau)^{b_1}\int_0^{\tau}\xb(\tau^\prime)d\tau^\prime
\label{14.253}
\end{equation}
In connection with \ref{14.247}, \ref{14.250} implies:
\begin{equation}
\int_0^{\tau}\tau^\prime(\sigma+\tau^\prime)^2\|\s^{(n-2)}\ctchi\|_{L^2
(S_{\tau^\prime,\sigma+\tau^\prime})}d\tau^\prime\leq 
\tau^{a_1+\frac{3}{2}}(\sigma+\tau)^{b_1+\frac{5}{2}}\int_0^{\tau}x(\tau^\prime)d\tau^\prime
\label{14.254}
\end{equation}
and \ref{14.252} implies:
\begin{equation}
\int_0^{\tau}(\sigma+\tau^\prime)^4\|\s^{(n-2)}\ctchib\|_{L^2(S_{\tau^\prime,\sigma+\tau^\prime})}
d\tau^\prime\leq\tau^{a_1+1}(\sigma+\tau)^{b_1+4}\int_0^{\tau}\xb(\tau^\prime)\frac{d\tau^\prime}
{\sqrt{\tau^\prime}}
\label{14.255}
\end{equation}
In connection with \ref{14.248}, \ref{14.250} implies:
\begin{equation}
\int_0^{\tau}\tau^{\prime 2}\|\s^{(n-2)}\ctchi\|_{L^2(S_{\tau^\prime,\sigma+\tau^\prime})}d\tau^\prime
\leq \tau^{a_1+\frac{5}{2}}(\sigma+\tau)^{b_1+\frac{1}{2}}\int_0^{\tau}x(\tau^\prime)d\tau^\prime
\label{14.256}
\end{equation}
and \ref{14.253} implies:
\begin{equation}
\int_0^{\tau}\tau^\prime(\sigma+\tau^\prime)^2\|\s^{(n-2)}\ctchib\|_{L^2
(S_{\tau^\prime,\sigma+\tau^\prime})}d\tau^\prime\leq \tau^{a_1+\frac{3}{2}}(\sigma+\tau)^{b_1+2}
\int_0^\tau\xb(\tau^\prime)d\tau^\prime
\label{14.257}
\end{equation}
Moreover, we have:
\begin{eqnarray}
&&\int_0^{\tau}\|\s^{(0,n-1)}\cla\|_{L^2(S_{\tau^\prime,\sigma+\tau^\prime})}d\tau^\prime
=\int_0^{\tau}\tau^{\prime a_1+1}(\sigma+\tau^\prime)^{b_1+\frac{1}{2}}y(\tau^\prime)d\tau^\prime
\nonumber\\
&&\hspace{10mm}\leq \tau^{a_1+1}(\sigma+\tau)^{b_1+\frac{1}{2}}\int_0^{\tau}y(\tau^\prime)d\tau^\prime
\label{14.258}
\end{eqnarray}
and:
\begin{eqnarray}
&&\int_0^{\tau}\|\s^{(0,n-1)}\clab\|_{L^2(S_{\tau^\prime,\sigma+\tau^\prime})}d\tau^\prime
=\int_0^{\tau}\tau^{\prime a_1+\frac{1}{2}}(\sigma+\tau^\prime)^{b_1}d\tau^\prime\nonumber\\
&&\hspace{10mm}\leq \tau^{a_1+\frac{1}{2}}(\sigma+\tau)^{b_1}\int_0^{\tau}\yb(\tau^\prime)d\tau^\prime
\label{14.259}
\end{eqnarray}
Substituting the above in \ref{14.245} - \ref{14.248} yields the following system of integral 
inequalities:
\begin{eqnarray}
&&x(\tau)\leq C\sqrt{D}+C\tau^{1/2}\int_0^{\tau}y(\tau^\prime)d\tau^\prime \nonumber\\
&&\hspace{10mm}+C\int_0^{\tau}x(\tau^\prime)d\tau^\prime
+C\int_0^{\tau}\xb(\tau^\prime)\frac{d\tau^\prime}{\sqrt{\tau^\prime}} \label{14.260}
\end{eqnarray}
\begin{eqnarray}
&&\xb(\tau)\leq C\sqrt{D}+C\int_0^{\tau}\yb(\tau^\prime)d\tau^\prime \nonumber\\
&&\hspace{10mm}+C(\sigma+\tau)^{1/2}\int_0^{\tau}x(\tau^\prime)d\tau^\prime
+C\int_0^{\tau}\xb(\tau^\prime)d\tau^\prime \label{14.261}
\end{eqnarray}
\begin{eqnarray}
&&y(\tau)\leq C\sqrt{D}+C\int_0^{\tau}y(\tau^\prime)d\tau^\prime 
+C\tau^{1/2}(\sigma+\tau)^2\int_0^{\tau}x(\tau^\prime)d\tau^\prime\nonumber\\
&&\hspace{10mm}
+C(\sigma+\tau)^{7/2}\int_0^{\tau}\xb(\tau^\prime)\frac{d\tau^\prime}{\sqrt{\tau^\prime}} \label{14.262}
\end{eqnarray}
\begin{eqnarray}
&&\yb(\tau)\leq C\sqrt{D}+C\int_0^{\tau}\yb(\tau^\prime)d\tau^\prime 
+C\tau^2(\sigma+\tau)^{1/2}\int_0^{\tau}x(\tau^\prime)d\tau^\prime \nonumber\\
&&\hspace{10mm}+C\tau(\sigma+\tau)^2\int_0^{\tau}\xb(\tau^\prime)d\tau^\prime \label{14.263}
\end{eqnarray}
which implies:
\begin{equation}
x(\tau), \ \xb(\tau), \ y(\tau), \ \yb(\tau) \ \leq \ C\sqrt{D}
\label{14.264}
\end{equation}

For $m=2,...,n$ we have:
\begin{eqnarray} 
&&TE^{n-m}T^{m-1}\lambda-E^{n-m}T^m\lambda=-(n-m)(\chi+\chib)E^{n-m}T^{m-1}\lambda \nonumber\\
&&TE^{n-m}T^{m-1}\lambdab-E^{n-m}T^m\lambdab=-(n-m)(\chi+\chib)E^{n-m}T^{m-1}\lambdab \nonumber\\
&&\label{14.265}
\end{eqnarray}
up to terms of order $n-2$, and similarly for the $N$th approximants, hence, taking difference 
we obtain:
\begin{eqnarray}
&&T\s^{(m-1,n-m)}\cla=\s^{(m,n-m)}\cla -(n-m)(\chi+\chib)\s^{(m-1,n-m)}\cla \nonumber\\
&&T\s^{(m-1,n-m)}\clab=\s^{(m,n-m)}\clab-(n-m)(\chi+\chib)\s^{(m-1,n-m)}\clab \nonumber\\
&&\label{14.266}
\end{eqnarray}
up to ignorable terms. Applying then \ref{14.240} with $\s^{(m-1,n-m)}\cla$, $\s^{(m-1,n-m)}\clab$ 
in the role of $f$ yields the inequalities:
\begin{eqnarray}
&&\|\s^{(m-1,n-m)}\cla\|_{L^2(S_{\tau,\sigma+\tau})}\leq 
C\tau^{1/2}\|\s^{(m,n-m)}\cla\|_{L^2({\cal K}_\sigma^\tau)}\nonumber\\
&&\hspace{40mm}+C\int_0^\tau\|\s^{(m-1,n-m)}\cla\|_{L^2(S_{\tau^\prime,\sigma+\tau^\prime})}d\tau^\prime\nonumber\\
&&\|\s^{(m-1,n-m)}\clab\|_{L^2(S_{\tau,\sigma+\tau})}\leq 
C\tau^{1/2}\|\s^{(m,n-m)}\clab\|_{L^2({\cal K}_\sigma^\tau)}\nonumber\\
&&\hspace{40mm}+C\int_0^\tau\|\s^{(m-1,n-m)}\clab\|_{L^2(S_{\tau^\prime,\sigma+\tau^\prime})}d\tau^\prime\nonumber\\
\label{14.267}
\end{eqnarray}
Substituting the $L^2({\cal K}_\sigma^\tau)$ estimates for $\s^{(m,n-m)}\cla$, $\s^{(m,n-m)}\clab$ of 
Proposition 14.2 we then deduce:
\begin{eqnarray}
&&\|\s^{(m-1,n-m)}\cla\|_{L^2(S_{\tau,\sigma+\tau})}\leq C\sqrt{D}
\tau^{a_m+1}(\sigma+\tau)^{b_m+\frac{1}{2}}\nonumber\\
&&\|\s^{(m-1,n-m)}\clab\|_{L^2(S_{\tau,\sigma+\tau})}\leq C\sqrt{D}
\tau^{a_m+\frac{1}{2}}(\sigma+\tau)^{b_m}\nonumber\\
&&\hspace{25mm}\mbox{: for $m=2,...,n$} \label{14.268}
\end{eqnarray}

We summarize the results of this section in the following proposition.  

\vspace{2.5mm}

\noindent{\bf Proposition 14.4} \ \ \ The $n-1$th order acoustical difference quantities satisfy the 
following $L^2$ estimates on the $S_{\ub,u}$:
\begin{eqnarray*}
&&\|\s^{(n-2)}\ctchi\|_{L^2(S_{\ub,u})}\leq C\sqrt{D}\ub^{a_1+\frac{1}{2}}u^{b_1+\frac{1}{2}}\\
&&\|\s^{(n-2)}\ctchib\|_{L^2(S_{\ub,u})}\leq C\sqrt{D}\ub^{a_1+\frac{1}{2}}u^{b_1}\\
&&\|\s^{(m-1,n-m)}\cla\|_{L^2(S_{\ub,u})}\leq C\sqrt{D}\ub^{a_m+1}u^{b_m+\frac{1}{2}} 
\ \ \ \mbox{: for $m=1,...,n$}\\
&&\|\s^{(m-1,n-m)}\clab\|_{L^2(S_{\ub,u})}\leq C\sqrt{D}\ub^{a_m+\frac{1}{2}}u^{b_m} 
\ \ \ \mbox{: for $m=1,...,n$}
\end{eqnarray*} 

\vspace{2.5mm}

In addition, in regard to the transformation functions, Propositions 12.3 and 12.6 together with 
the top order energy estimates imply:
\begin{eqnarray}
&&\|\Omega^{n-m}T^{m+1}\chf\|_{L^2({\cal K}^{\tau})}\leq C\sqrt{D}\tau^{c_m-2} \nonumber\\
&&\|\Omega^{n-m}T^{m+1}\cv\|_{L^2({\cal K}^{\tau})}\leq C\sqrt{D}\tau^{c_m-2} \nonumber\\
&&\|\Omega^{n-m}T^{m+1}\cga\|_{L^2({\cal K}^{\tau})}\leq C\sqrt{D}\tau^{c_m-2} 
\label{14.a26}
\end{eqnarray}
Applying then inequality \ref{12.366} we obtain:
\begin{eqnarray}
&&\|\Omega^{n-m}T^m\chf\|_{L^2(S_{\tau,\tau})}\leq C\sqrt{D}\tau^{c_m-\frac{3}{2}} \nonumber\\
&&\|\Omega^{n-m}T^m\cv\|_{L^2(S_{\tau,\tau})}\leq C\sqrt{D}\tau^{c_m-\frac{3}{2}} \nonumber\\
&&\|\Omega^{n-m}T^m\cga\|_{L^2(S_{\tau,\tau})}\leq C\sqrt{D}\tau^{c_m-\frac{3}{2}} 
\label{14.a27}
\end{eqnarray}

\vspace{5mm}

\section{$L^2(S_{\ub,u})$ Estimates for All $n$th Order Derivatives of the $\beta_\mu$}

We now define for each $k=0,...,n$ and each $m=k,...,n$ the $n$th order quantities:
\begin{equation}
q_{k,m,\mu}=E^{n-m}T^{m-k}\Lb^k\beta_\mu-E_N^{n-m}T^{m-k}\Lb_N^k\beta_{\mu,N}
\label{14.269}
\end{equation}
These quantities express all $n$th order derivatives of the $\beta_\mu$. For $k=0$ we have, simply:
\begin{equation}
q_{0,m,\mu}=\s^{(m,n-m)}\check{\dot{\phi}}_\mu 
\label{14.270}
\end{equation}
while for $k=1$ the definition \ref{14.269} reduces to the definition \ref{14.188}. Thus in the cases 
$k=0$ and $k=1$ $L^2(S_{\ub,u})$ estimates for $q_{k,m,\mu}$ are provided by Proposition 14.3 and 
Lemma 14.1 respectively. We now proceed to derive $L^2(S_{\ub,u})$ estimates for $q_{k,m,\mu}$ 
for $k=2,...,n$. This shall be accomplished through a recursive inequality for the $L^2(S_{\ub,u})$ 
norms to be presently established. 

We write:
\begin{equation}
\Lb^k\beta_\mu=\Lb^{k-2}\Lb^2\beta_\mu=\Lb^{k-2}T\Lb\beta_\mu-\Lb^{k-2}L\Lb\beta_\mu 
\label{14.271}
\end{equation}
and similarly for the $N$th approximants. Applying $E^{n-m}T^{m-k}$ to \ref{14.271} 
and $E_N^{n-m}T^{m-k}$ to its 
analogue for the $N$th approximants, and subtracting, we obtain: 
\begin{eqnarray}
&&q_{k,m,\mu}=E^{n-m}T^{m-k}\Lb^{k-2}T\Lb\beta_\mu-E_N^{n-m}T^{m-k}\Lb_N^{k-2}T\Lb_N\beta_{\mu,N} 
\nonumber\\
&&\hspace{13mm}-E^{n-m}T^{m-k}\Lb^{k-2}L\Lb\beta_\mu 
+E_N^{n-m}T^{m-k}\Lb_N^{k-2}L_N\Lb_N\beta_{\mu,N} \nonumber\\
\label{14.272}
\end{eqnarray}
The first difference on the right is:
\begin{equation}
q_{k-1,m,\mu}+E^{n-m}T^{m-k}[\Lb^{k-2},T]\Lb\beta_\mu-E_N^{n-m}T^{m-k}[\Lb_N^{k-2},T]\Lb_N\beta_{\mu,N}
\label{14.273}
\end{equation}
Here the difference involving the commutators is:
\begin{equation}
\sum_{j=0}^{k-3}\left\{E^{n-m}T^{m-k}\Lb^{k-3-j}(\zeta E\Lb^{j+1}\beta_\mu)
-E_N^{n-m}T^{m-k}\Lb_N^{k-3-j}(\zeta_N E_N\Lb_N^{j+1}\beta_{\mu,N})\right\}
\label{14.274}
\end{equation}
The highest order term contributed by this sum is
$$(k-2)\zeta E^{n-m+1}T^{m-k}\Lb^{k-2}\beta_\mu$$
which is of order $n-1$ with vanishing acoustical part of order $n-1$, the highest order acoustical 
term being 
$$(E\Lb\beta_\mu)E^{n-m}T^{m-k}\Lb^{k-3}\zeta$$
which is of order $n-2$. In the present context we can consider terms of order $n-1$ with vanishing 
acoustical part of order $n-1$ as ignorable. Thus the difference involving the commutators in 
\ref{14.273} can be ignored. 

We turn to the second difference on the right in \ref{14.272}. Here we substitute \ref{14.194}, 
\ref{14.195}. In view of \ref{14.196}, this gives the following expression for this difference, 
up to terms which can be ignored:
\begin{eqnarray}
&&-a E^{n-m}T^{m-k}\Lb^{k-2}E^2\beta_\mu+a_N E_N^{n-m}T^{m-k}\Lb_N^{k-2}E_N^2\beta_{\mu,N}\nonumber\\
&&+\frac{1}{2}E^{n-m}T^{m-k}\left\{\rhob(L\beta_\mu)\Lb^{k-2}\tchib 
+\rho(\Lb\beta_\mu)\Lb^{k-2}\tchi-4\rhob(E\beta_\mu)\Lb^{k-2}\mub\right\}\nonumber\\
&&-\frac{1}{2}E_N^{n-m}T^{m-k}\left\{\rhob_N(L_N\beta_{\mu,N})\Lb_N^{k-2}\tchib_N 
+\rho_N(\Lb_N\beta_{\mu,N})\Lb_N^{k-2}\tchi_N\right.\nonumber\\
&&\hspace{59mm}\left.-4\rhob_N(E_N\beta_{\mu,N})\Lb_N^{k-2}\mub_N\right\}\nonumber\\
&&\label{14.275}
\end{eqnarray}
Now since $[\Lb\tchib]_{P.A.}=[\Lb\mub]_{P.A.}=0$, the 1st and 3rd term in the parentheses 
contribute non-ignorable terms only for $k=2$, in which case \ref{14.275} reduces to:
\begin{eqnarray}
&&-a E^{n-m}T^{m-2}E^2\beta_\mu+a_N E_N^{n-m}T^{m-2}E_N^2\beta_{\mu,N}\nonumber\\
&&+\frac{1}{2}(L\beta_\mu)\rhob\left\{\begin{array}{lll}
\s^{(n-2)}\ctchib &:& \mbox{for $m=2$}\\
\cmub_{m-3,n-m+1} &:& \mbox{for $m=3,...,n$}
\end{array}\right.\nonumber\\
&&+\frac{1}{2}(\Lb\beta_\mu)\rho\left\{\begin{array}{lll}
\s^{(n-2)}\ctchi &:& \mbox{for $m=2$}\\
\cmu_{m-3,n-m+1} &:& \mbox{for $m=3,...,n$}
\end{array}\right.\nonumber\\
&&-2(E\beta_\mu)\rhob\cmub_{m-2,n-m} \label{14.276}
\end{eqnarray}
up to terms which can be ignored. On the other hand, for $k=3,...,n$, in regard to the 2nd term 
in the parentheses, writing 
$$\Lb^{k-2}\tchi=(T-L)^{k-2}\tchi$$
and expanding the binomial only the term $T^{k-2}\tchi$ contains a $k-1$th order acoustical part 
in $2T^{k-3}E\mu$. Therefore for $k=3,....,n$ \ref{14.275} reduces to:
\begin{eqnarray}
&&-a E^{n-m}T^{m-k}\Lb^{k-2}E^2\beta_\mu+a_N E_N^{n-m}T^{m-k}\Lb_N^{k-2}E_N^2\beta_{\mu,N}\nonumber\\
&&+\rho(\Lb\beta_\mu)\cmu_{m-3,n-m+1} \label{14.277}
\end{eqnarray}
up to terms which can be ignored. In regard to the first two terms in \ref{14.276}, we have:
\begin{eqnarray}
&&-a E^{n-m}T^{m-2}E^2\beta_\mu+a_N E_N^{n-m}T^{m-2}E_N^2\beta_{\mu,N}\nonumber\\
&&=-a E^{n-m+2}T^{m-2}\beta_\mu+a_N E_N^{n-m+2}T^{m-2}\beta_{\mu,N}\nonumber\\
&&+(E\beta_\mu)E^{n-m+1}T^{m-3}(\chi+\chib)-(E_N\beta_{\mu,N})E_N^{n-m+1}T^{m-3}(\chi_N+\chib_N)
\nonumber\\
&&\mbox{: up to terms which can be ignored} \nonumber\\
&&=-aq_{0,m-2,\mu}\nonumber\\
&&+a(E\beta_\mu)\left\{\begin{array}{lll}
\rho\s^{(n-2)}\ctchi+\rhob\s^{(n-2)}\ctchib &:& \mbox{for $m=3$}\\
2\rho\cmu_{m-4,n-m+2}+2\rhob\cmub_{m-4,n-m+2} &:& \mbox{for $m=4,...,n$}
\end{array}\right.\nonumber\\
&&\mbox{: up to terms which can be ignored} \label{14.278}
\end{eqnarray}
In regard to the first two terms in \ref{14.277}, noting that the principal acoustical part of 
the commutator $[\Lb^{k-2},E^2]\beta_\mu$ is contained in the term 
$-(E\beta_\mu)E\Lb^{k-3}\chib$, hence is equal to $-(E\beta_\mu)\rhob E\tchib$ for $k=3$ and 
vanishes for $k\geq 4$, we conclude that for $k=3$ we have:
\begin{eqnarray}
&&-a E^{n-m}T^{m-3}\Lb E^2\beta_\mu+a_N E_N^{n-m}T^{m-3}\Lb_N E_N^2\beta_{\mu,N}\nonumber\\
&&=-aq_{1,m-2,\mu}+a(E\beta_\mu)\left\{\begin{array}{lll}\rhob\s^{(n-2)}\ctchib &:& \mbox{for $m=3$}\\
2\rhob\cmub_{m-4,n-m+2} &:& \mbox{for $m=4,...,n$}\end{array}\right. \nonumber\\
&&\mbox{: up to terms which can be ignored}\label{14.279}
\end{eqnarray}
while for $k\geq 4$ we have simply:
\begin{eqnarray}
&&-a E^{n-m}T^{m-k}\Lb^{k-2}E^2\beta_\mu+a_N E_N^{n-m}T^{m-k}\Lb_N^{k-2}E_N^2\beta_{\mu,N}\nonumber\\
&&=-aq_{k-2,m-2,\mu} \ \ \mbox{: up to terms which can be ignored} \label{14.280}
\end{eqnarray}

We now use Proposition 14.4 to estimate in $L^2(S_{\ub,u})$ the remainders in 
\ref{14.276} - \ref{14.279}.
From \ref{12.567} with $(m-1,n-1)$ in the role of $(m,n)$ we obtain, by virtue of Proposition 14.4, 
\begin{equation}
\|\cmu_{m-2,n-m}\|_{L^2(S_{\ub,u})}\leq C\sqrt{D}\ub^{a_{m-1}+\frac{1}{2}}u^{b_{m-1}+1}
\label{14.281}
\end{equation}
From \ref{12.568} with $(m-1,n-1)$ in the role of $(m,n)$ we obtain, by virtue of Proposition 14.4,
\begin{equation}
\|\cmub_{m-2,n-m}\|_{L^2(S_{\ub,u})}\leq C\sqrt{D}\ub^{a_{m-1}+\frac{1}{2}}u^{b_{m-1}}
\label{14.282}
\end{equation}
Using \ref{14.281}, \ref{14.282}, and the results of Proposition 14.4 concerning $\s^{(n-2)}\ctchi$, 
$\s^{(n-2)}\ctchib$, we conclude that the remainder in \ref{14.276} is bounded in $L^2(S_{\ub,u})$ by:
\begin{equation}
C\sqrt{D}\left\{\begin{array}{lll} \ub^{a_1+\frac{1}{2}}u^{b_1+\frac{3}{2}} &:& \mbox{for $m=2$}\\
\ub^{a_{m-1}+\frac{1}{2}}u^{b_{m-1}+2} &:& \mbox{for $m=3,...,n$} \end{array}\right.
\label{14.283}
\end{equation}
the remainder in \ref{14.277} is bounded in $L^2(S_{\ub,u})$ by:
\begin{equation}
C\sqrt{D}\ub^{a_{m-2}+\frac{3}{2}}u^{b_{m-2}+1}
\label{14.284}
\end{equation}
the remainder in \ref{14.278} is bounded in $L^2(S_{\ub,u})$ by:
\begin{equation}
C\sqrt{D}\left\{\begin{array}{lll} \ub^{a_1+\frac{3}{2}}u^{b_1+\frac{7}{2}} &:& \mbox{for $m=3$}\\
\ub^{a_{m-3}+\frac{3}{2}}u^{b_{m-3}+4} &:& \mbox{for $m=4,...,n$} \end{array}\right.
\label{14.285}
\end{equation}
and the remainder in \ref{14.279} is bounded in $L^2(S_{\ub,u})$ by:
\begin{equation}
C\sqrt{D}\left\{\begin{array}{lll} \ub^{a_1+\frac{3}{2}}u^{b_1+4} &:& \mbox{for $m=3$}\\
\ub^{a_{m-3}+\frac{3}{2}}u^{b_{m-3}+4} &:& \mbox{for $m=4,...,n$} \end{array}\right.
\label{14.286}
\end{equation}

In conclusion, \ref{14.272} - \ref{14.273} and \ref{14.276} - \ref{14.279} yield for $k\geq 2$ 
the recursive inequality (recall that $a\sim\ub u^2$):
\begin{equation}
\|q_{k,m,\mu}\|_{L^2(S_{\ub,u})}\leq \|q_{k-1,m,\mu}\|_{L^2(S_{\ub,u})}
+C\ub u^2\|q_{k-2,m-2,\mu}\|_{L^2(S_{\ub,u})}+r_{k,m}
\label{14.287}
\end{equation}
where, by the estimates \ref{14.283} - \ref{14.286}:
\begin{equation}
r_{k,m}=C\sqrt{D}\left\{\begin{array}{lll} \ub^{a_1+\frac{1}{2}}u^{b_1+\frac{3}{2}} &:& 
\mbox{for $k=m=2$}\\
\ub^{a_{m-1}+\frac{1}{2}}u^{b_{m-1}+2} &:& \mbox{for $k=2$, $m\geq 3$}\\
\ub^{a_{m-2}+\frac{3}{2}}u^{b_{m-2}+1} &:& \mbox{for $k\geq 3$} \end{array}\right.
\label{14.288}
\end{equation}

\vspace{2.5mm}

\noindent{\bf Proposition 14.5} \ \ \ The $n$th order quantities $q_{k,m,\mu}$ satisfy the following 
$L^2$ estimates on the $S_{\ub,u}$:
\begin{eqnarray*}
&&\|q_{k,m,\mu}\|_{L^2(S_{\ub,u})}\leq C\sqrt{D}\ub^{a_{m-1}+\frac{1}{2}}u^{b_{m-1}+\frac{1}{2}}\\
&&\hspace{5mm}\mbox{: for all $k=1,...,n$ ; $m=k,...,n$}
\end{eqnarray*}

\vspace{2.5mm}

\noindent {\em Proof:} The proposition is true for $k=1$ by Lemma 14.1. We now take $k\geq 2$ and 
assume as the inductive hypothesis that the statement holds with $k$ replaced by any 
$k^\prime=1,...,k-1$. We then have:
$$\|q_{k-1,m,\mu}\|_{L^2(S_{\ub,u})}\leq C\sqrt{D}\ub^{a_{m-1}+\frac{1}{2}}u^{b_{m-1}+\frac{1}{2}}$$
and:
$$\|q_{k-2,m-2,\mu}\|_{L^2(S_{\ub,u})}\leq C\sqrt{D}\left\{
\begin{array}{lll} \ub^{a_{m-2}+\frac{1}{2}}u^{b_{m-2}-\frac{1}{2}} &:& \mbox{for $k=2$}\\
\ub^{a_{m-3}+\frac{1}{2}}u^{b_{m-3}+\frac{1}{2}} &:& \mbox{for $k\geq 3$} \end{array}\right.$$
Substituting these in \ref{14.287} and taking taking into account \ref{14.288} yields the 
statement for $k^\prime=k$ and the inductive step is complete. 

\vspace{5mm}

\section{$L^2(S_{\ub,u})$ Estimates for $\Omega^{n-1}\log\sh$ and $\Omega^{n-1}b$} 

Let us denote:
\begin{eqnarray}
&&\s^{(n-1)}\check{\sh}=\Omega^{n-1}\log\left(\frac{\sh}{\sh_N}\right)\nonumber\\
&&\hspace{12mm}=\Omega^{n-1}\log\sh-\Omega^{n-1}\log\sh_N \label{14.289}
\end{eqnarray}
Note that this is analogous, but does not quite coincide with , \ref{12.954} with $n-1$ in the role 
of $n$. According to \ref{12.956}, \ref{12.958} we have:
\begin{equation}
T\log\sh=2(\chi+\chib), \ \ \ 
T\log\sh_N=2(\chi_N+\chib_N) 
\label{14.290} 
\end{equation}
Applying $\Omega^{n-1}$ to these and taking the difference gives:
\begin{equation}
T\s^{(n-1)}\check{\sh}=2(\Omega^{n-1}\chi-\Omega^{n-1}\chi_N)+2(\Omega^{n-1}\chib-\Omega^{n-1}\chib_N)
\label{14.291}
\end{equation}

Now, we have:
\begin{equation}
\Omega=\sqrt{\sh}E=\sqrt{\sh_N}E_N
\label{14.292}
\end{equation}
It follows that if $f$ is an arbitrary function on ${\cal N}$ it holds:
\begin{eqnarray}
&&\Omega^k f=\sum_{i=1}^k P_{i,k}\sh^{i/2}E^i f \nonumber\\
&&\hspace{8mm}=\sum_{i=1}^k P_{i,k,N}\sh_N^{i/2}E_N^i f \label{14.293}
\end{eqnarray}
Here 
\begin{equation}
P_{k,k}=1
\label{14.294}
\end{equation}
and for each $i=1,...,k-1$, $P_{i,k}$ is a polynomial in 
$\Omega\log\sh, . . . ,\Omega^{k-i}\log\sh$, the total number of $\Omega$ occurring in each term 
of the polynomial being equal to $k-i$. Thus the highest order term in $P_{i,k}$ is the 
linear term:
\begin{equation}
c_{i,k}\Omega^{k-i}\log\sh
\label{14.295}
\end{equation}
and we have 
\begin{equation}
c_{1,k}=\frac{1}{2}
\label{14.296}
\end{equation}
Also $P_{i,k,N}$ is the same as $P_{i,k}$ but with $\log\sh_N$ in the role of 
$\log\sh$. 

Applying \ref{14.293} to the functions $\chi$, $\chi_N$ and to the functions $\chib$, $\chib_N$ 
we obtain:
\begin{eqnarray}
&&\Omega^{n-1}\chi-\Omega^{n-1}\chi_N=\sum_{i=1}^{n-1}\left\{P_{i,n-1}\sh^{i/2}E^i\chi
-P_{i,n-1,N}\sh_N^{i/2}E_N^i\chi_N\right\}\nonumber\\
&&\Omega^{n-1}\chib-\Omega^{n-1}\chib_N=\sum_{i=1}^{n-1}\left\{P_{i,n-1}\sh^{i/2}E^i\chib 
-P_{i,n-1,N}\sh_N^{i/2}E_N^i\chib_N\right\}\nonumber\\
&&\label{14.297}
\end{eqnarray}
The highest order terms on the right in each of the above are:
\begin{eqnarray*}
&&\sh^{(n-1)/2}E^{n-1}\chi-\sh_N^{(n-1)/2}E_N^{n-1}\chi_N \\
&&+\frac{1}{2}(\Omega^{n-2}\log\sh)\sh^{1/2}E\chi
-\frac{1}{2}(\Omega^{n-2}\log\sh_N)\sh_N^{1/2}E_N\chi_N
\end{eqnarray*}
or, up to 2nd order terms, 
\begin{equation}
\sh^{(n-1)/2}(E^{n-1}\chi-E_N^{n-1}\chi_N)+\frac{1}{2}(E\chi)\s^{(n-2)}\check{\sh}
\label{14.298}
\end{equation}
and:
\begin{eqnarray*}
&&\sh^{(n-1)/2}E^{n-1}\chib-\sh_N^{(n-1)/2}E_N^{n-1}\chib_N \\
&&+\frac{1}{2}(\Omega^{n-2}\log\sh)\sh^{1/2}E\chib
-\frac{1}{2}(\Omega^{n-2}\log\sh_N)\sh_N^{1/2}E_N\chib_N
\end{eqnarray*}
or, up to 2nd order terms, 
\begin{equation}
\sh^{(n-1)/2}(E^{n-1}\chib-E_N^{n-1}\chib_N)+\frac{1}{2}(E\chib)\s^{(n-2)}\check{\sh}
\label{14.299}
\end{equation}
Actually the contribution of the 2nd terms involving $\s^{(n-2)}\check{\sh}$, 
in \ref{14.298}, \ref{14.299}, can be absorbed in the estimate for $\s^{(n-2)}\check{\sh}$ being 
of 1 order lower than the quantity under estimation. In fact, 
$$T\s^{(n-2)}\check{\sh}=2\Omega^{n-2}(\chi-\chi_N)+2\Omega^{n-2}(\chib-\chib_N)$$
and this can itself be estimated in $L^2(S_{\ub,u})$. We conclude that, up to terms which can be 
absorbed, 
\begin{eqnarray}
&&\|T\s^{(n-1)}\check{\sh}\|_{L^2(S_{\ub,u})}\leq 
C\left\{\|E^{n-1}\chi-E_N^{n-1}\chi_N\|_{L^2(S_{\ub,u})}\right.\nonumber\\
&&\hspace{40mm}\left.+\|E^{n-1}\chib-E_N^{n-1}\chib_N\|_{L^2(S_{\ub,u})}\right\}
\label{14.300}
\end{eqnarray}
We now apply the inequality \ref{14.240} with $\s^{(n-1)}\check{\sh}$ in the role of the function $f$ 
to obtain, by \ref{14.300}, 
\begin{eqnarray}
&&\|\s^{(n-1)}\check{\sh}\|_{L^2(S_{\tau,\sigma+\tau})}\leq 
C\int_0^{\tau}\|E^{n-1}\chi-E_N^{n-1}\chi_N\|_{L^2(S_{\tau^\prime,\sigma+\tau^\prime})}d\tau^\prime 
\nonumber\\
&&\hspace{31mm}+C\int_0^{\tau}
\|E^{n-1}\chib-E_N^{n-1}\chib_N\|_{L^2(S_{\tau^\prime,\sigma+\tau^\prime})}
d\tau^\prime \nonumber\\
&&\label{14.301}
\end{eqnarray}
By \ref{3.a20} and \ref{10.138}, \ref{10.141} we have, to principal terms, 
\begin{eqnarray*}
&&E^{n-1}\chi-E_N^{n-1}\chi_N=\rho\s^{(n-1)}\ctchi+\frac{1}{2}\sbeta^2(E^{n-1}LH-E_N^{n-1}L_N H_N)\\
&&E^{n-1}\chib-E_N^{n-1}\chib_N=\rhob\s^{(n-1)}\ctchib+\frac{1}{2}\sbeta^2(E^{n-1}\Lb H
-E_N^{n-1}\Lb_N H_N)
\end{eqnarray*}
Moreover, to principal terms, 
\begin{eqnarray*}
&&E^{n-1}LH-E_N^{n-1}L_N H_N=-2H^\prime\beta^\mu E^{n-1}L\beta_\mu+2H^\prime_N\beta^\mu_{\s,N}E_N^{n-1}L_N\beta_{\mu,N}\\
&&\hspace{38mm}=-2H^\prime\beta^\mu\qb_{1,1,\mu}
\end{eqnarray*}
\begin{eqnarray*}
&&E^{n-1}\Lb H-E_N^{n-1}\Lb_N H_N=-2H^\prime\beta^\mu E^{n-1}\Lb\beta_\mu 
+2H^\prime_N\beta^\mu_{\s,N}E_N^{n-1}\Lb_N\beta_{\mu,N}\\
&&\hspace{38mm}=-2H^\prime\beta^\mu q_{1,1,\mu}
\end{eqnarray*}
by the definitions \ref{14.188}, \ref{14.190}. Therefore, to principal terms we have:
\begin{eqnarray}
&&E^{n-1}\chi-E_N^{n-1}\chi_N=\rho\s^{(n-1)}\ctchi-\sbeta^2 H^\prime\beta^\mu\qb_{1,1,\mu}\nonumber\\
&&E^{n-1}\chib-E_N^{n-1}\chib_N=\rhob\s^{(n-1)}\ctchib-\sbeta^2 H^\prime\beta^\mu q_{1,1,\mu}\nonumber\\
&&\label{14.302}
\end{eqnarray}
It follows that in regard to the two integrals on the right in \ref{14.301} we can estimate:
\begin{eqnarray}
&&\int_0^{\tau}\|E^{n-1}\chi-E_N^{n-1}\chi_N\|_{L^2(S_{\tau^\prime,\sigma+\tau^\prime})}d\tau^\prime 
\leq C\tau^{3/2}\|\s^{(n-1)}\ctchi\|_{L^2({\cal K}_\sigma^\tau)}\nonumber\\
&&\hspace{60mm}+C\int_0^{\tau}\sum_\mu\|\qb_{1,1,\mu}\|_{L^2(S_{\tau^\prime,\sigma+\tau^\prime})}d\tau^\prime\nonumber\\
&&\label{14.303}\\
&&\int_0^{\tau}\|E^{n-1}\chib-E_N^{n-1}\chib_N\|_{L^2(S_{\tau^\prime,\sigma+\tau^\prime})}d\tau^\prime 
\leq C\tau^{1/2}(\sigma+\tau)^2\|\s^{(n-1)}\ctchib\|_{L^2({\cal K}_\sigma^\tau)}\nonumber\\
&&\hspace{60mm}+C\int_0^{\tau}\sum_\mu\|q_{1,1,\mu}\|_{L^2(S_{\tau^\prime,\sigma+\tau^\prime})}d\tau^\prime\nonumber\\
&&\label{14.304}
\end{eqnarray}
By Proposition 14.2 the 1st term on the right in \ref{14.303} is bounded by 
$C\sqrt{D}\ub^{a_0+\frac{3}{2}}u^{b_0}$ while by Proposition 14.3 and Lemma 14.1 (see \ref{14.189}, 
\ref{14.190}) the 2nd term is bounded by $C\sqrt{D}\ub^{a_1+\frac{3}{2}}u^{b_1-\frac{1}{2}}$ 
(recall that $\ub=\tau$, $u=\sigma+\tau$). Also, by Proposition 14.2 the 1st term on the right in 
\ref{14.304} is bounded by $C\sqrt{D}\ub^{a_0+\frac{1}{2}}u^{b_0+1}$ while by Lemma 14.1 the 2nd term 
is bounded by $C\sqrt{D}\ub^{a_0+\frac{3}{2}}u^{b_0+\frac{1}{2}}$. It then follows through \ref{14.301} 
that:
\begin{equation}
\|\s^{(n-1)}\check{\sh}\|_{L^2(S_{\ub,u})}\leq C\sqrt{D}\ub^{a_1+\frac{1}{2}}u^{b_1+\frac{1}{2}}
\label{14.305}
\end{equation}

We proceed to derive an $L^2(S_{\ub,u})$ estimate for $\Omega^{n-1}\check{b}$, where we denote:
\begin{equation}
\check{b}=b-b_N
\label{14.306}
\end{equation}
By \ref{3.a8},
\begin{eqnarray*}
&&[\Lb,L]=\left[\frac{\partial}{\partial u}+b\frac{\partial}{\partial\vartheta},
\frac{\partial}{\partial\ub}-b\frac{\partial}{\partial\vartheta}\right]\\
&&\hspace{10mm}=-\left(\frac{\partial b}{\partial\ub}+\frac{\partial b}{\partial u}\right)
\frac{\partial}{\partial\vartheta}=-(Tb)\Omega
\end{eqnarray*}
Comparing with the 3rd of \ref{3.a14} we conclude that:
\begin{equation}
Tb=-\sh^{-1/2}\zeta
\label{14.307}
\end{equation}
From \ref{9.11} and the 3rd of \ref{9.b1} we similarly deduce:
\begin{equation}
Tb_N=-\sh_N^{-1/2}\zeta_N
\label{14.308}
\end{equation}
Taking the difference and applying $\Omega^{n-1}$ we then obtain:
\begin{equation}
T\Omega^{n-1}\check{b}=-\Omega^{n-1}(\sh^{-1/2}\zeta)+\Omega^{n-1}(\sh_N^{-1/2}\zeta_N)
\label{14.309}
\end{equation}
Applying then \ref{14.293} to the functions $\sh^{-1/2}\zeta$ and $\sh_N^{-1/2}\zeta_N$ gives:
\begin{equation}
T\Omega^{n-1}\check{b}=-\sum_{i=1}^{n-1}\left\{P_{i,n-1}\sh^{i/2}E^i(\sh^{-1/2}\zeta)
-P_{i,n-1,N}\sh_N^{i/2}E_N^i(\sh_N^{-1/2}\zeta_N)\right\}
\label{14.310}
\end{equation}
The highest order terms on the right are:
\begin{eqnarray*}
&&-\sh^{(n-1)/2}E^{n-1}(\sh^{-1/2}\zeta)+\sh_N^{(n-1)/2}E_N^{n-1}(\sh_N^{-1/2}\zeta_N)\\
&&-\frac{1}{2}(\Omega^{n-2}\log\sh)\sh^{1/2}E(\sh^{-1/2}\zeta)
+\frac{1}{2}(\Omega_N^{n-2}\log\sh_N)\sh_N^{1/2}E_N(\sh_N^{-1/2}\zeta_N)
\end{eqnarray*}
Here, as in \ref{14.298}, \ref{14.299}, the 2nd difference can be absorbed, 
and the above can be replaced by:
\begin{eqnarray*}
&&-\sh^{(n-2)/2}E^{n-1}\zeta+\sh_N^{(n-2)/2}E_N^{n-1}\zeta_N\\
&&+\frac{1}{2}\sh^{-1/2}(\Omega^{n-1}\log\sh)\zeta-\frac{1}{2}\sh_N^{-1/2}(\Omega^{n-1}\log\sh_N)\zeta_N
\end{eqnarray*}
which can in turn be replaced by:
\begin{equation}
-\sh^{(n-2)/2}(E^{n-1}\zeta-E_N^{n-1}\zeta_N)+\frac{1}{2}\sh^{-1/2}\s^{(n-1)}\check{\sh}\zeta
\label{14.311}
\end{equation}
We conclude that, up to terms which can be 
absorbed, 
\begin{equation}
\|T\Omega^{n-1}\check{b}\|_{L^2(S_{\ub,u})}\leq C\|E^{n-1}\zeta-E_N^{n-1}\zeta_N\|_{L^2(S_{\ub,u})}
+C\ub\|\s^{(n-1)}\check{\sh}\|_{L^2(S_{\ub,u})}
\label{14.312}
\end{equation}
(recall that $|\zeta|\leq C\ub$). We now apply inequality \ref{14.240} with $\Omega^{n-1}\check{b}$ 
in the role of the function $f$ to obtain, by \ref{14.312},
\begin{eqnarray}
&&\|\Omega^{n-1}\check{b}\|_{L^2(S_{\tau,\sigma+\tau})}\leq 
C\int_0^{\tau}\|E^{n-1}\zeta-E_N^{n-1}\zeta_N\|_{L^2(S_{\tau^\prime,\sigma+\tau^\prime})}d\tau^\prime 
\nonumber\\
&&\hspace{31mm}+C\int_0^{\tau}
\tau^\prime\|\s^{(n-1)}\check{\sh}\|_{L^2(S_{\tau^\prime,\sigma+\tau^\prime})}
d\tau^\prime \nonumber\\
&&\label{14.313}
\end{eqnarray}
From \ref{3.a22}, \ref{3.a25} and the corresponding equations for the $N$th approximants we deduce:
\begin{equation}
[E^{n-1}\zeta-E_N^{n-1}\zeta_N]_{P.A.}=2\rho\s^{(0,n)}\cla-2\rhob\s^{(0,n)}\clab
+\pi a\s^{(n-1)}\ctchi-\pi a\s^{(n-1)}\ctchib
\label{14.314}
\end{equation}
and:
\begin{eqnarray}
&&[E^{n-1}\zeta-E_N^{n-1}\zeta_N]_{P.P.}-[E^{n-1}\zeta-E_N^{n-1}\zeta_N]_{P.A.}\nonumber\\
&&\hspace{30mm}=\sbeta H^\prime \beta^\mu(-\rhob\beta_{\Nb} \qb_{1,1,\mu}
+\rho\beta_N q_{1,1,\mu})\nonumber\\
&&+\frac{a}{c}\left[-(\beta_N-\beta_{\Nb})H^\prime \beta^\mu 
+H\left(-c\pi\sbeta+\beta_N-\frac{\beta_{\Nb}}{2}\right)N^\mu\right.\nonumber\\
&&\hspace{36mm}\left.+H\left(c\pi\sbeta-\beta_{\Nb}+\frac{\beta_N}{2}\right)\Nb^\mu\right]\s^{(0,n)}\check{\dot{\phi}}_\mu 
\nonumber\\
&&\label{14.315}
\end{eqnarray}
We bound the contribution of \ref{14.314} to the first integral on the right in \ref{14.314} by:
$$\tau^{1/2}\|[E^{n-1}\zeta-E_N^{n-1}\zeta_N]_{P.A.}\|_{L^2({\cal K}_\sigma^\tau)}$$
and then applying Proposition 14.2. This gives a bound by:
\begin{equation}
C\sqrt{D}\ub^{a_0+\frac{1}{2}}u^{b_0+2}
\label{14.316}
\end{equation}
(recall again that $\ub=\tau$, $u=\sigma+\tau$). We bound the contribution of \ref{14.315} to the 
same integral using Proposition 14.3 and Lemma 14.1 (see \ref{14.189}, \ref{14.190}). 
This gives a bound by:
\begin{equation}
C\sqrt{D}\ub^{a_1+\frac{3}{2}}u^{b_1+\frac{3}{2}}
\label{14.317}
\end{equation}
Finally by \ref{14.305} the second integral on the right in \ref{14.314} is bounded by:
\begin{equation}
C\sqrt{D}\ub^{a_1+\frac{5}{2}}u^{b_1+\frac{1}{2}}
\label{14.318}
\end{equation}
Combining we then conclude through \ref{14.314} that:
\begin{equation}
\|\Omega^{n-1}\check{b}\|_{L^2(S_{\ub,u})}\leq C\sqrt{D}\ub^{a_1+\frac{1}{2}}u^{b_1+2}
\label{14.319}
\end{equation}

We summarize the results of this section in the following proposition.  

\vspace{2.5mm}

\noindent{\bf Proposition 14.6} \ \ \ The $n-1$th order quantities $\s^{(n-1)}\check{\sh}$, 
$\Omega^{n-1}\check{b}$ satisfy the following $L^2$ estimates on the $S_{\ub,u}$:
\begin{eqnarray*}
&&\|\s^{(n-1)}\check{\sh}\|_{L^2(S_{\ub,u})}\leq C\sqrt{D}\ub^{a_1+\frac{1}{2}}u^{b_1+\frac{1}{2}}\\
&&\|\Omega^{n-1}\check{b}\|_{L^2(S_{\ub,u})}\leq C\sqrt{D}\ub^{a_1+\frac{1}{2}}u^{b_1+2}
\end{eqnarray*}

\vspace{5mm}

\section{Lower Order $L^2(S_{\ub,u})$ Estimates}

Let us consider the $L^2(S_{\ub,u})$ estimates for the $n$th order variation differences 
$\s^{(m,n-m)}\check{\dot{\phi}}_\mu : m=0,...,n$ of Proposition 14.3, 
for the $n$th order quantities $q_{1,m,\mu} : m=1,...,n$ of Lemma 14.1, 
for the $n$th order quantities 
$\qb_{1,m,\mu} : m=1,...,n$ which follow from the preceding through \ref{14.189}, \ref{14.190}, 
as well as the $L^2(S_{\ub,u})$ estimates for the $n-1$th order acoustical difference quantities 
$\s^{(n-2)}\ctchi$, $\s^{(n-2)}\tchib$, $\s^{(m-1,n-m)}\cla : m=1,...,n$, 
$\s^{(m-1,n-m)}\clab : m=1,...,n$ of Proposition 14.4, 
recalling also \ref{14.233}, \ref{14.234}. To deduce from these $L^2(S_{\ub,u})$ estimates of 
one order lower, that is $L^2(S_{\ub,u})$ estimates for the $n-1$th order variation differences 
$\s^{(m-1,n-m)}\check{\dot{\phi}}_\mu : m=1,...,n$, for the $n-1$th order quantities 
\begin{equation}
\stackrel{1}{q}_{1,m-1,\mu}=E^{n-m}T^{m-2}\Lb\beta_\mu-E_N^{n-m}T^{m-2}\Lb_N\beta_{\mu,N} : m=2,...,n
\label{14.320}
\end{equation}
which correspond to the $q_{1,m-1,\mu} : m=2,...,n$ at one order lower, 
for the $n-1$th order quantities 
\begin{equation}
\stackrel{1}{\qb}_{1,m-1,\mu}=E^{n-m}T^{m-2}L\beta_\mu-E_N^{n-m}T^{m-2}L_N\beta_{\mu,N} : m=2,...,n
\label{14.321}
\end{equation}
which correspond to the $\qb_{1,m,\mu}$ at one order lower, and for the $n-2$th order acoustical 
quantities $\s^{(n-3)}\ctchi$, $\s^{(n-3)}\ctchib$, $\s^{(m-2,n-m)}\cla : m=2,...,n$, 
$\s^{(m-2,n-m)}\clab : m=2,...,n$, we must estimate the $L^2(S_{\ub,u})$ norm of $T$ applied to each 
of the lower order quantities in terms of the $L^2(S_{\ub,u})$ norm of a corresponding quantity of 
one order higher, that is of $\s^{(m,n-m)}\check{\dot{\phi}}_\mu$ in the case of 
$\s^{(m-1,n-m)}\check{\dot{\phi}}_\mu$, $q_{1,m,\mu}$ in the 
case of $\stackrel{1}{q}_{1,m-1,\mu}$, $\qb_{1,m,\mu}$ in the case of $\stackrel{1}{\qb}_{1,m-1,\mu}$, 
$\s^{(m-1,n-m)}\cla$ in the case of $\s^{(m-2,n-m)}\cla$, 
$\s^{(m-1,n-m)}\clab$ in the case of $\s^{(m-2,n-m)}\clab$. In the case of 
$\s^{(n-3)}\ctchi$, $\s^{(n-3)}\ctchib$, by \ref{14.233}, \ref{14.234} with $n-1$ in the role of $n$ 
the estimate involves 
$\s^{(1,n-2)}\check{\dot{\phi}}_\mu$. We then apply the inequality \ref{14.240} to obtain the 
desired $L(S_{\ub,u})$ estimates of one order lower. This procedure involves lowering the index $m$ 
of the quantity being estimated by one when the quantity considered is of one order lower. 
Since in \ref{14.240} there is an integration with respect to $\tau$, the estimate is by one 
power of $\ub$ higher than that of the corresponding quantity of one higher order, and there is 
a factor equal to the reciprocal of that power in the coefficient. Thus the $L^2(S_{\ub,u})$ 
estimates of one order lower read: 
\begin{eqnarray}
&&\|\s^{(m-1,n-m)}\check{\dot{\phi}}_\mu\|_{L^2(S_{\ub,u})}\leq 
\frac{C\sqrt{D}}{\left(a_m+\frac{3}{2}\right)}\ub^{a_m+\frac{3}{2}}u^{b_m-\frac{1}{2}} 
\ \ : \ m=1,...n\nonumber\\
&&\|\stackrel{1}{q}_{1,m-1,\mu}\|_{L^2(S_{\ub,u})}\leq 
\frac{C\sqrt{D}}{\left(a_m+\frac{3}{2}\right)}\ub^{a_m+\frac{3}{2}}u^{b_m+\frac{1}{2}} 
\ \ : \ m=2,...,n\nonumber\\
&&\|\stackrel{1}{\qb}_{1,m-1,\mu}\|_{L^2(S_{\ub,u})}\leq 
\frac{C\sqrt{D}}{\left(a_m+\frac{3}{2}\right)}\ub^{a_m+\frac{3}{2}}u^{b_m-\frac{1}{2}} 
\ \ : \ m=2,...,n\nonumber\\
&&\label{14.322}
\end{eqnarray}
\begin{eqnarray}
&&\|\s^{(n-3)}\ctchi\|_{L^2(S_{\ub,u})}\leq 
\frac{C\sqrt{D}}{\left(a_2+\frac{3}{2}\right)}\ub^{a_2+\frac{3}{2}}u^{b_2+\frac{1}{2}}\nonumber\\
&&\|\s^{(n-3)}\ctchib\|_{L^2(S_{\ub,u})}\leq 
\frac{C\sqrt{D}}{\left(a_2+\frac{3}{2}\right)}\ub^{a_2+\frac{3}{2}}u^{b_2}\nonumber\\
&&\|\s^{(m-2,n-m)}\cla\|_{L^2(S_{\ub,u})}\leq 
\frac{C\sqrt{D}}{\left(a_m+2\right)}\ub^{a_m+2}u^{b_m+\frac{1}{2}} \ \ : \ m=2,...,n\nonumber\\
&&\|\s^{(m-2,n-m)}\clab\|_{L^2(S_{\ub,u})}\leq 
\frac{C\sqrt{D}}{\left(a_m+\frac{3}{2}\right)}\ub^{a_m+\frac{3}{2}}u^{b_m} \ \ : \ m=2,...,n\nonumber\\
&&\label{14.323}
\end{eqnarray} 

We see that if to a quantity of the higher order is associated the pair of indices $(m,n-m)$, the 
indices $m$ and $n-m$ corresponding to the number of $T$ and $E$ derivatives respectively, and 
the pair of exponents $(a_m,b_m)$, indicating the growth of the quantity, the exponents $a_m$ and 
$b_m$ corresponding to the powers of $\ub$ and $u$ respectively, then to the corresponding quantity 
of one order lower but with the same number of $T$ derivatives is associated the index pair 
$(m,n-1-m)$ and, according to the above, the exponent pair $(a_{m+1}+1,b_{m+1})$. We remark that in 
general for positive exponents $a$, $b$ and $a^\prime$, $b^\prime$ we have 
$$\ub^{a^\prime}u^{b^\prime}\leq \ub^a u^b \ \ \mbox{: for all $(\ub,u)\in R_{1,1}$}$$
if and only if $a^\prime\geq a$ and $c^\prime\geq c$ where we denote $c=a+b$, 
$c^\prime=a^\prime+b^\prime$. Now, in connection with the pair of quantities just discussed, 
the condition for the quantity of one order lower to be able to be absorbed in the estimates involving 
the higher order quantity is that its growth is strictly dominated by the growth of the higher order 
quantity. Therefore, in view of the remark, the necessary and sufficient condition for the quantity 
to which the index pair $(m,n-1-m)$ is associated to be able to be absorbed in the estimates involving 
the corresponding quantity with associated index pair $(m,n-m)$ is that:
\begin{equation}
a_{m+1}+1>a_m \ \ \mbox{and} \ \ c_{m+1}+1>c_m
\label{14.324}
\end{equation}
So this condition on the exponents $a_m$, $b_m$ must be imposed in addition to the condition that 
they are non-increasing with $m$. With this condition satisfied, all the terms which we have 
called ignorable in the preceding chapters can indeed be absorbed in the estimates. The choice of 
relation between the exponent pairs $(a_m,b_m)$ for different $m$ which ensures that the lower order 
terms are maximally depressed relative to the higher order terms is simply to choose 
all the $a_m$ equal and all the $b_m$ equal:
\begin{equation}
a_m=a, \ \ \ b_m=b \ \ \ : \ m=0,...,n
\label{14.325}
\end{equation}
So we take this optimal choice from now on. 

Applying then repeatedly the preceding argument we deduce 
the $L^2(S_{\ub,u})$ estimates of $i$ orders lower in the form:
\begin{eqnarray}
&&\|\s^{(m-i,n-m)}\check{\dot{\phi}}_\mu\|_{L^2(S_{\ub,u})}\leq 
\frac{C\sqrt{D}}{\left(a+\frac{3}{2}\right)\cdot\cdot\cdot\left(a+\frac{1}{2}+i\right)}
\ub^{a+\frac{1}{2}+i}u^{b-\frac{1}{2}} 
\ \ : \ m=i,...n\nonumber\\
&&\|\stackrel{i}{q}_{1,m-i,\mu}\|_{L^2(S_{\ub,u})}\leq 
\frac{C\sqrt{D}}{\left(a+\frac{3}{2}\right)\cdot\cdot\cdot\left(a+\frac{1}{2}+i\right)}
\ub^{a+\frac{1}{2}+i}u^{b+\frac{1}{2}} 
\ \ : \ m=i+1,...,n\nonumber\\
&&\|\stackrel{i}{\qb}_{1,m-i,\mu}\|_{L^2(S_{\ub,u})}\leq 
\frac{C\sqrt{D}}{\left(a+\frac{3}{2}\right)\cdot\cdot\cdot\left(a+\frac{1}{2}+i\right)}
\ub^{a+\frac{1}{2}+i}u^{b-\frac{1}{2}} 
\ \ : \ m=i+1,...,n\nonumber\\
&&\label{14.326}
\end{eqnarray} 
\begin{eqnarray}
&&\|\s^{(n-2-i)}\ctchi\|_{L^2(S_{\ub,u})}\leq 
\frac{C\sqrt{D}}{\left(a+\frac{3}{2}\right)\cdot\cdot\cdot\left(a+\frac{1}{2}+i\right)}
\ub^{a+\frac{1}{2}+i}u^{b+\frac{1}{2}}\nonumber\\
&&\|\s^{(n-2-i)}\ctchib\|_{L^2(S_{\ub,u})}\leq 
\frac{C\sqrt{D}}{\left(a+\frac{3}{2}\right)\cdot\cdot\cdot\left(a+\frac{1}{2}+i\right)}
\ub^{a+\frac{1}{2}+i}u^{b}\nonumber\\
&&\|\s^{(m-1-i,n-m)}\cla\|_{L^2(S_{\ub,u})}\leq 
\frac{C\sqrt{D}}{\left(a+2\right)\cdot\cdot\cdot\left(a+1+i\right)}\ub^{a+1+i}u^{b+\frac{1}{2}} \ \ : \ m=i+1,...,n\nonumber\\
&&\|\s^{(m-1-i,n-m)}\clab\|_{L^2(S_{\ub,u})}\leq 
\frac{C\sqrt{D}}{\left(a+\frac{3}{2}\right)\cdot\cdot\cdot\left(a+\frac{1}{2}+i\right)}
\ub^{a+\frac{1}{2}+i}u^{b} \ \ : \ m=i+1,...,n\nonumber\\
&&\label{14.327}
\end{eqnarray} 
where, in \ref{14.326}, we denote:
\begin{eqnarray}
&&\stackrel{i}{q}_{1,m-i,\mu}=E^{n-m}T^{m-1-i}\Lb\beta_\mu-E_N^{n-m}T^{m-1-i}\Lb_N\beta_{\mu,N} 
\nonumber\\
&&\stackrel{i}{\qb}_{1,m-i,\mu}=E^{n-m}T^{m-1-i}L\beta_\mu-E_N^{n-m}T^{m-1-i}L_N\beta_{\mu,N}
\nonumber\\
&&\label{14.328}
\end{eqnarray}

In addition, in regard to the transformation functions, we deduce from the estimates \ref{14.a27} 
the following lower order estimates:
\begin{eqnarray}
&&\|\Omega^{n-m}T^{m-i}\chf\|_{L^2(S_{\tau,\tau})}\leq \frac{C\sqrt{D}}{\left(c-\frac{1}{2}\right)
\cdot\cdot\cdot\left(c-\frac{3}{2}+i\right)}\tau^{c-\frac{3}{2}+i} \nonumber\\
&&\|\Omega^{n-m}T^{m-i}\cv\|_{L^2(S_{\tau,\tau})}\leq \frac{C\sqrt{D}}{\left(c-\frac{1}{2}\right)
\cdot\cdot\cdot\left(c-\frac{3}{2}+i\right)}\tau^{c-\frac{3}{2}+i} \nonumber\\
&&\|\Omega^{n-m}T^{m-i}\cga\|_{L^2(S_{\tau,\tau})}\leq \frac{C\sqrt{D}}{\left(c-\frac{1}{2}\right)
\cdot\cdot\cdot\left(c-\frac{3}{2}+i\right)}\tau^{c-\frac{3}{2}+i} \nonumber\\
&&\mbox{: for $i=1,...,n-1$, $m=i,...,n$} \label{14.a28}
\end{eqnarray}

\vspace{2.5mm}

In assessing the contribution of the lower order terms to the top order energy estimates these 
terms are of course not estimated in terms $D$ but in terms of $\s^{[n,n]}{\cal M}$ and 
$\s^{[n,n]}{\cal A}$. Since the lower order terms contribute with progressively higher powers of 
$\delta$ and inverse powers of $a$, they can all be absorbed in the estimates.

\vspace{5mm}

\section{Pointwise Estimates and Recovery of the Bootstrap Assumptions}

We shall now derive pointwise estimates for all the lower order quantities appearing above, that is 
for $\s^{(m-i,n-m)}\check{\dot{\phi}}_\mu : m=i,...,n; \ i=1,...,n$, for 
$\stackrel{i}{q}_{1,m-i,\mu}$, $\stackrel{i}{\qb}_{1,m-i,\mu}$ : $m=i+1,...,n; \ i=1,...,n-1$, and for  
$\s^{(n-2-i)}\ctchi$, $\s^{(n-2-i)}\ctchib$ : $i=1,...,n-2$, and for 
$\s^{(m-1-i,n-m)}\cla$, $\s^{(m-1-i,n-m)}\clab$ : $m=i+1,...,n; \ i=1,...,n-1$. To derive the 
pointwise estimates we shall use the following version of the Sobolev inequality on $S^1$. 

\vspace{2.5mm}

\noindent{\bf Lemma 14.2} \ \ \ Let $f\in H_1(S^1)$. Then $f\in C^{0,1/2}(S^1)$ and we have:
$$\sup_{S^1}|f-\overline{f}|\leq\sqrt{2}\|f-\overline{f}\|^{1/2}_{L^2(S^1)}
\|\Omega f\|^{1/2}_{L^2(S^1)}$$
where $\overline{f}$ is the mean value of $f$ on $S^1$. 

\vspace{2.5mm}

\noindent {\em Proof:} This is standard. For the inequality, consider the function 
$\check{f}=f-\overline{f}$ . This function has vanishing mean hence by continuity there is a 
$\vartheta_0\in S^1$ where $\check{f}$ vanishes. Then for an arbitrary $\vartheta\in S^1$ 
we have:
\begin{eqnarray*}
&&\check{f}^2(\vartheta)=\int_{\vartheta_0}^{\vartheta}(\Omega(\check{f}^2))(\vartheta^\prime)
d\vartheta^\prime
=2\int_{\vartheta_0}^{\vartheta}(\check{f}\Omega\check{f})(\vartheta^\prime)d\vartheta^\prime\\
&&\hspace{10mm}\leq\int_{S^1}|\check{f}||\Omega\check{f}|
\leq\|\check{f}\|_{L^2(S^1)}\|\Omega\check{f}\|_{L^2(S^1)}
\end{eqnarray*}
Since $\Omega\check{f}=\Omega f$ the inequality follows. 

\vspace{2.5mm}

Since 
$$\|f-\overline{f}\|^2_{L^2(S^1)}=\int_{S^1}(f-\overline{f})^2=\int_{S^1}f^2-\int_{S^1}\overline{f}^2
\leq \|f\|^2_{L^2(S^1)}$$
and 
$$|\overline{f}|\leq\frac{1}{2\pi}\|f\|_{L^1(S^1)}\leq\frac{1}{\sqrt{2\pi}}\|f\|_{L^2(S^1)}$$
Lemma 14.2 implies:
\begin{equation}
\sup_{S^1}|f|\leq\frac{1}{\sqrt{2\pi}}\|f\|_{L^2(S^1)}
+\sqrt{2}\|f\|^{1/2}_{L^2(S^1)}\|\Omega f\|^{1/2}_{L^2(S^1)}
\label{14.329}
\end{equation}

Let us first revisit Proposition 14.6. Following the same argument with $n-i$ in the role of $n$ and 
using the estimates \ref{14.326} and \ref{14.327} we deduce the lower order $L^2(S_{\ub,u})$ 
estimates:
\begin{eqnarray}
&&\|\s^{(n-1-i))}\check{\sh}\|_{L^2(S_{\ub,u})}\leq 
\frac{C\sqrt{D}}{\left(a+\frac{3}{2}\right)\cdot\cdot\cdot\left(a+\frac{1}{2}+i\right)}
\ub^{a+\frac{1}{2}+i}u^{b+\frac{1}{2}}\nonumber\\
&&\|\Omega^{n-1-i}\check{b}\|_{L^2(S_{\ub,u})}\leq 
\frac{C\sqrt{D}}{\left(a+\frac{3}{2}\right)\cdot\cdot\cdot\left(a+\frac{1}{2}+i\right)}
\ub^{a+\frac{1}{2}+i}u^{b+2}\nonumber\\
&&\hspace{18mm} : \ i=1,...,n-1\label{14.330}
\end{eqnarray}
In view of the definition \ref{14.289} we have:
\begin{equation}
\Omega\s^{(n-2-i)}\check{\sh}=\s^{(n-1-i)}\check{\sh} \ \ : \ i=0,...,n-2
\label{14.331}
\end{equation}
Applying then Lemma 14.2 with the restrictions of $\s^{(n-2-i)}\check{\sh}$ and 
$\Omega^{n-2-i}\check{b}$ to $S_{\ub,u}$ in the role of the function $f$ and appealing 
to Proposition 14.6 for the case $i=0$ and to the estimates \ref{14.330} for the cases $i=1,..,n-2$ 
we deduce the pointwise estimates:
\begin{eqnarray}
&&\sup_{S_{\ub,u}}|\s^{(n-2-i)}\check{\sh}|\leq 
\frac{C\sqrt{D}}{\left(a+\frac{3}{2}\right)\cdot\cdot\cdot\left(a+\frac{1}{2}+i\right)}
\frac{\ub^{a+1+i}u^{b+\frac{1}{2}}}{\sqrt{a+\frac{3}{2}+i}}\nonumber\\
&&\sup_{S_{\ub,u}}|\Omega^{(n-2-i)}\check{b}|\leq 
\frac{C\sqrt{D}}{\left(a+\frac{3}{2}\right)\cdot\cdot\cdot\left(a+\frac{1}{2}+i\right)}
\frac{\ub^{a+1+i}u^{b+2}}{\sqrt{a+\frac{3}{2}+i}}\nonumber\\
&&\hspace{18mm} : \ i=0,...,n-2\label{14.332}
\end{eqnarray}

Consider next the lower order variation differences. In view of \ref{14.292} we have, 
for $i=0,...,n-1$ and $m=i+1,...,n$:
\begin{eqnarray}
&&\Omega\s^{(m-1-i,n-m)}\check{\dot{\phi}}_\mu=\Omega(E^{n-m}T^{m-1-i}\beta_\mu 
-E_N^{n-m}T^{m-1-i}\beta_{\mu,N})\nonumber\\
&&=\sh^{1/2}E^{n-m+1}T^{m-1-i}\beta_\mu-\sh_N^{1/2}E_N^{n-m+1}T^{m-1-i}\beta_{\mu,N}\nonumber\\
&&=\sh^{1/2}\s^{(m-1-i,n-m+1)}\check{\dot{\phi}}
+(\sh^{1/2}-\sh_N^{1/2})E_N^{n-m+1}T^{m-1-i}\beta_{\mu,N}\nonumber\\
&&\label{14.333}
\end{eqnarray}
By Proposition 14.3 for the case $i=0$ and the first of the estimates \ref{14.326} for the cases 
$i=1,...,n-1$ the first term on the right in \ref{14.333} is bounded in $L^2(S_{\ub,u})$ by:
\begin{equation}
\frac{C\sqrt{D}}{\left(a+\frac{3}{2}\right)\cdot\cdot\cdot\left(a+\frac{1}{2}+i\right)}
\ub^{a+\frac{1}{2}+i}u^{b-\frac{1}{2}} 
\label{14.334}
\end{equation}
while by the first of \ref{14.330} with $i=n-1$ the second term on the right in \ref{14.333} is bounded 
in $L^2(S_{\ub,u})$ by:
\begin{equation}
\frac{C\sqrt{D}}{\left(a+\frac{3}{2}\right)\cdot\cdot\cdot\left(a-\frac{1}{2}+n\right)}
\ub^{a-\frac{1}{2}+n}u^{b+\frac{1}{2}}
\label{14.335}
\end{equation}
Therefore for all $i=0,...,n-1$ and $m=i+1,...,n$ we have:
\begin{equation}
\|\Omega\s^{(m-1-i,n-m)}\check{\dot{\phi}}_\mu\|_{L^2(S_{\ub,u})}\leq 
\frac{C\sqrt{D}}{\left(a+\frac{3}{2}\right)\cdot\cdot\cdot\left(a+\frac{1}{2}+i\right)}
\ub^{a+\frac{1}{2}+i}u^{b-\frac{1}{2}} 
\label{14.336}
\end{equation} 
Applying then Lemma 14.2 with 
$\left.\s^{(m-1-i,n-m)}\check{\dot{\phi}}_\mu\right|_{S_{\ub,u}}$ 
in the role of $f$ and using also the 1st of \ref{14.326} with $i+1$ in the role of $i$ we obtain:
\begin{equation}
\sup_{S_{\ub,u}}|\s^{(m-1-i,n-m)}\check{\dot{\phi}}_\mu|
\leq\frac{C\sqrt{D}}{\left(a+\frac{3}{2}\right)\cdot\cdot\cdot\left(a+\frac{1}{2}+i\right)}
\frac{\ub^{a+1+i}u^{b-\frac{1}{2}}}{\sqrt{a+\frac{3}{2}+i}}
\label{14.337}
\end{equation}
for all $i=0,...,n-1$ and $m=i+1,...,n$. 

Next, in view of \ref{14.328} we have, for $i=1,...,n-1$ and $m=i+1,...,n$:
\begin{eqnarray}
&&\Omega\stackrel{i}{q}_{1,m-i,\mu}=\Omega(E^{n-m}T^{m-1-i}\Lb\beta_\mu 
-E_N^{n-m}T^{m-1-i}\Lb_N\beta_{\mu,N})\nonumber\\
&&=\sh^{1/2}E^{n-m+1}T^{m-1-i}\Lb\beta_\mu-\sh_N^{1/2}E_N^{n-m+1}T^{m-1-i}\Lb_N\beta_{\mu,N}\nonumber\\
&&=\sh^{1/2}\stackrel{i-1}{q}_{1,m-i,\mu}
+(\sh^{1/2}-\sh_N^{1/2})E_N^{n-m+1}T^{m-1-i}\Lb_N\beta_{\mu,N}\nonumber\\
&&\label{14.338}
\end{eqnarray}
By Lemma 14.1 for the case $i=1$ and the second of the estimates \ref{14.326} for the cases 
$i=2,...,n-1$ the first term on the right in \ref{14.338} is bounded in $L^2(S_{\ub,u})$ by:
\begin{equation}
\frac{C\sqrt{D}}{\left(a+\frac{3}{2}\right)\cdot\cdot\cdot\left(a-\frac{1}{2}+i\right)}
\ub^{a-\frac{1}{2}+i}u^{b+\frac{1}{2}}
\label{14.339}
\end{equation}
while by the first of \ref{14.330} with $i=n-1$ the second term on the right in \ref{14.338} is 
bounded in $L^2(S_{\ub,u})$ by: 
\begin{equation}
\frac{C\sqrt{D}}{\left(a+\frac{3}{2}\right)\cdot\cdot\cdot\left(a-\frac{1}{2}+n\right)}
\ub^{a-\frac{1}{2}+n}u^{b+\frac{1}{2}}
\label{14.340}
\end{equation}
Therefore for all $i=1,...,n-1$ and $m=i+1,...,n$ we have:
\begin{equation}
\|\Omega\stackrel{i}{q}_{1,m-i,\mu}\|_{L^2(S_{\ub,u})}\leq 
\frac{C\sqrt{D}}{\left(a+\frac{3}{2}\right)\cdot\cdot\cdot\left(a-\frac{1}{2}+i\right)}
\ub^{a-\frac{1}{2}+i}u^{b+\frac{1}{2}}
\label{14.341}
\end{equation}
Using also the 2nd of \ref{14.326} in applying Lemma 14.2 with 
$\left.\stackrel{i}{q}_{1,m-i,\mu}\right|_{S_{\ub,u}}$ in the role of $f$ we obtain:
\begin{equation}
\sup_{S_{\ub,u}}|\stackrel{i}{q}_{1,m-i,\mu}|\leq 
\frac{C\sqrt{D}}{\left(a+\frac{3}{2}\right)\cdot\cdot\cdot\left(a-\frac{1}{2}+i\right)}
\frac{\ub^{a+i}u^{b+\frac{1}{2}}}{\sqrt{a+\frac{1}{2}+i}}
\label{14.342}
\end{equation}
for all $i=1,...,n-1$ and $m=i+1,...,n$. Since (see \ref{14.328}):
\begin{equation}
\stackrel{i}{q}_{1,m-i,\mu}+\stackrel{i}{\qb}_{1,m-i,\mu}=\s^{(m-i,n-m)}\check{\dot{\phi}}_\mu 
\label{14.343}
\end{equation}
the estimates \ref{14.337} with $i-1$ in the role of $i$ and \ref{14.342} imply:
\begin{equation}
\sup_{S_{\ub,u}}|\stackrel{i}{\qb}_{1,m-i,\mu}|\leq 
\frac{C\sqrt{D}}{\left(a+\frac{3}{2}\right)\cdot\cdot\cdot\left(a-\frac{1}{2}+i\right)}
\frac{\ub^{a+i}u^{b-\frac{1}{2}}}{\sqrt{a+\frac{1}{2}+i}}
\label{14.344}
\end{equation}
for all $i=1,...,n-1$ and $m=i+1,...,n$. 

We turn to the lower order acoustical differences. In view of \ref{14.292} we have, for 
$i=1,...,n-2$:
\begin{eqnarray}
&&\Omega\s^{(n-2-i)}\ctchi=\Omega(E^{n-2-i}\tchi-E_N^{n-2-i}\tchi_N)\nonumber\\
&&=\sh^{1/2}E^{n-1-i}\tchi-\sh_N^{1/2}E_N^{n-1-i}\tchi_N\nonumber\\
&&=\sh^{1/2}\s^{(n-1-i)}\ctchi+(\sh^{1/2}-\sh_N^{1/2})E_N^{n-1-i}\tchi_N
\label{14.345}
\end{eqnarray}
Similarly, 
\begin{equation}
\Omega\s^{(n-2-i)}\ctchib=\sh^{1/2}\s^{(n-1-i)}\ctchib+(\sh^{1/2}-\sh_N^{1/2})E_N^{n-1-i}\ctchib_N
\label{14.346}
\end{equation}
The first terms on the right in each of \ref{14.345}, \ref{14.346} are bounded in $L^2(S_{\ub,u})$ 
by Proposition 14.4 in the case $i=1$ and by the first two of the estimates \ref{14.327} with $i-1$ in the role of $i$ 
in the cases $i=2,...,n-2$, 
while the second terms on the right in each of \ref{14.345}, \ref{14.346} are bounded in 
$L^2(S_{\ub,u})$ using the first of \ref{14.330} with $i=n-1$. We obtain:
\begin{eqnarray}
&&\|\Omega\s^{(n-2-i)}\ctchi\|_{L^2(S_{\ub,u})}\leq 
\frac{C\sqrt{D}}{\left(a+\frac{3}{2}\right)\cdot\cdot\cdot\left(a-\frac{1}{2}+i\right)}
\ub^{a-\frac{1}{2}+i}u^{b+\frac{1}{2}}\nonumber\\
&&\|\Omega\s^{(n-2-i)}\ctchib\|_{L^2(S_{\ub,u})}\leq 
\frac{C\sqrt{D}}{\left(a+\frac{3}{2}\right)\cdot\cdot\cdot\left(a-\frac{1}{2}+i\right)}
\ub^{a-\frac{1}{2}+i}u^{b}\nonumber\\
&&\label{14.347}
\end{eqnarray}
Using also the first two of the estimates \ref{14.327} in applying Lemma 14.2 with 
$\left.\s^{(n-2-i)}\ctchi\right|_{S_{\ub,u}}$ and $\left.\s^{(n-2-i)}\ctchib\right|_{S_{\ub,u}}$ 
in the role of $f$ we deduce:
\begin{eqnarray}
&&\sup_{S_{\ub,u}}|\s^{(n-2-i)}\ctchi|\leq  
\frac{C\sqrt{D}}{\left(a+\frac{3}{2}\right)\cdot\cdot\cdot\left(a-\frac{1}{2}+i\right)}
\frac{\ub^{a+i}u^{b+\frac{1}{2}}}{\sqrt{a+\frac{1}{2}+i}}\nonumber\\
&&\sup_{S_{\ub,u}}|\s^{(n-2-i)}\ctchib|\leq  
\frac{C\sqrt{D}}{\left(a+\frac{3}{2}\right)\cdot\cdot\cdot\left(a-\frac{1}{2}+i\right)}
\frac{\ub^{a+i}u^{b}}{\sqrt{a+\frac{1}{2}+i}}\nonumber\\
&&\label{14.348}
\end{eqnarray}
for all $i=1,...,n-2$. Taking account of \ref{14.292} we have, for $i=1,...,n-1$, $m=i+1,...,n$:
\begin{eqnarray}
&&\Omega\s^{(m-1-i,n-m)}\cla=\sh^{1/2}\s^{(m-1-i,n-m+1)}\cla
+(\sh^{1/2}-\sh_N^{1/2})E_N^{n-m+1}T^{m-1-i}\lambda_N \nonumber\\
&&\label{14.349}\\
&&\Omega\s^{(m-1-i,n-m)}\clab=\sh^{1/2}\s^{(m-1-i,n-m+1)}\clab
+(\sh^{1/2}-\sh_N^{1/2})E_N^{n-m+1}T^{m-1-i}\lambdab_N \nonumber\\
&&\label{14.350}
\end{eqnarray}
The first term on the right in each of \ref{14.349}, \ref{14.350} 
is bounded in $L^2(S_{\ub,u})$ according to 
the estimates of Proposition 14.4 with $m-1$ in the role of $m$ in the case $i=1$ and according to 
the last two of the estimates \ref{14.327} with $(m-1,i-1)$ in the role of $(m,i)$ 
in the cases $i=2,...,n-1$. By the first of \ref{14.330} with $i=n-1$ the second term on the right in 
\ref{14.349} is bounded in $L^2(S_{\ub,u})$ by:
\begin{equation}
\frac{C\sqrt{D}}{\left(a+\frac{3}{2}\right)\cdot\cdot\cdot\left(a-\frac{1}{2}+n\right)}
\ub^{a-\frac{1}{2}+n}u^{b+\frac{1}{2}}\cdot\left\{\begin{array}{lll}u^2&:&\mbox{if $m=i+1$}\\
u&:&\mbox{if $m=i+2$}\\1&:&\mbox{if $m\geq i+3$}\end{array}\right.
\label{14.351}
\end{equation}
(see \ref{9.a21}), and the second term on the right in \ref{14.350} is bounded in $L^2(S_{\ub,u})$ by:
\begin{equation}
\frac{C\sqrt{D}}{\left(a+\frac{3}{2}\right)\cdot\cdot\cdot\left(a-\frac{1}{2}+n\right)}
\ub^{a-\frac{1}{2}+n}u^{b+\frac{1}{2}}\cdot\left\{\begin{array}{lll}\ub&:&\mbox{if $m=i+1$}\\
1&:&\mbox{if $m\geq i+2$}\end{array}\right.
\label{14.352}
\end{equation}
(see \ref{9.a13}). For both of \ref{14.349}, \ref{14.350} the bound for the first term dominates 
and we obtain:
\begin{eqnarray}
&&\|\Omega\s^{(m-1-i,n-m)}\cla\|_{L^2(S_{\ub,u})}\leq 
\frac{C\sqrt{D}}{(a+2)\cdot\cdot\cdot(a+i)}\ub^{a+i}u^{b+\frac{1}{2}} \nonumber\\
&&\label{14.353}\\
&&\|\Omega\s^{(m-1-i,n-m)}\clab\|_{L^2(S_{\ub,u})}\leq 
\frac{C\sqrt{D}}{\left(a+\frac{3}{2}\right)\cdot\cdot\cdot\left(a-\frac{1}{2}+i\right)}
\ub^{a-\frac{1}{2}+i}u^b \nonumber\\
&&\label{14.354}
\end{eqnarray}
Using also the last two of the estimates \ref{14.327} in applying Lemma 14.2 with 
$\left.\s^{(m-1-i,n-m)}\cla\right|_{S_{\ub,u}}$ and $\left.\s^{(m-1-i,n-m)}\clab\right|_{S_{\ub,u}}$ 
in the role of $f$ we deduce:
\begin{eqnarray}
&&\sup_{S_{\ub,u}}|\s^{(m-1-i,n-m)}\cla|\leq  
\frac{C\sqrt{D}}{(a+2)\cdot\cdot\cdot(a+i)}
\frac{\ub^{a+\frac{1}{2}+i}u^{b+\frac{1}{2}}}{\sqrt{a+1+i}}\nonumber\\
&&\sup_{S_{\ub,u}}|\s^{(m-1-i,n-m)}\clab|\leq  
\frac{C\sqrt{D}}{\left(a+\frac{3}{2}\right)\cdot\cdot\cdot\left(a-\frac{1}{2}+i\right)}
\frac{\ub^{a+i}u^{b}}{\sqrt{a+\frac{1}{2}+i}}\nonumber\\
&&\label{14.355}
\end{eqnarray}
for all $i=1,...,n-1$ and $m=i+1,...,n$. 

In addition, in regard to the transformation functions, we similarly deduce using the estimates 
\ref{14.a27} and \ref{14.a28}, the pointwise estimates:
\begin{eqnarray}
&&\sup_{S_{\tau,\tau}}|\Omega^{n-m}T^{m-i}\chf|\leq\frac{C\sqrt{D}}{\left(c-\frac{1}{2}\right)
\cdot\cdot\cdot\left(c-\frac{5}{2}+i\right)}\frac{\tau^{c-2+i}}{\sqrt{c-\frac{3}{2}+i}}
\nonumber\\
&&\sup_{S_{\tau,\tau}}|\Omega^{n-m}T^{m-i}\cv|\leq\frac{C\sqrt{D}}{\left(c-\frac{1}{2}\right)
\cdot\cdot\cdot\left(c-\frac{5}{2}+i\right)}\frac{\tau^{c-2+i}}{\sqrt{c-\frac{3}{2}+i}}
\nonumber\\
&&\sup_{S_{\tau,\tau}}|\Omega^{n-m}T^{m-i}\cga|\leq\frac{C\sqrt{D}}{\left(c-\frac{1}{2}\right)
\cdot\cdot\cdot\left(c-\frac{5}{2}+i\right)}\frac{\tau^{c-2+i}}{\sqrt{c-\frac{3}{2}+i}}
\nonumber\\
&&\label{14.a29}
\end{eqnarray}
for all $i=1,...,n-1$ and $m=i,...,n$. 

The above estimates \ref{14.337}, \ref{14.342}, \ref{14.344}, \ref{14.348}, \ref{14.355}, \ref{14.a29} 
allow us to recover all the bootstrap assumptions as {\em strict inequalities}, if 
\begin{equation}
n_*+2\leq n-1
\label{14.b4}
\end{equation}
that is (see \ref{14.a18}) if we choose: 
\begin{equation}
n\geq 5
\label{14.b5}
\end{equation} 
provided that $\delta$ is suitably small depending on:
$$\frac{\Ab_*}{\sqrt{D}}, \frac{\oA_*}{\sqrt{D}}, \frac{\Bb}{\sqrt{D}}, \frac{\oB}{\sqrt{D}}, 
\frac{B}{\sqrt{D}}, \frac{\Kb}{\sqrt{D}}, 
\frac{C_*}{\sqrt{D}}, \frac{C^\prime_*}{\sqrt{D}}, \frac{C^{\prime\prime}_*}{\sqrt{D}}$$

\vspace{5mm}

\section{Completion of the Argument}

Now, the results of Chapters 10, 12, 13, in particular the top order energy estimates, as well as the 
results of the present chapter up to this point, are in the nature of {\em a priori bounds}, and 
hold for a solution defined in a domain ${\cal R}_{\delta,\delta}$, provided that $\delta$ satisfies 
the smallness conditions stated in the preceding chapers, as long as the bootstrap assumptions hold 
on this domain. We shall presently show that once $\delta$ is subjected to the further smallness 
condition stated at the end of the previous section, allowing the recovery of the bootstrap assumptions,  
the whole argument can be completed to an actual proof of an existence theorem of a solution defined 
in ${\cal R}_{\delta,\delta}$ and satisfying all the stated bounds. 

To accomplish this we start not at $\Cb_0$, as the data on $S_{0,0}$ is singular and moreover $\lambdab$ 
vanishes everywhere on $\Cb_0$ (see Chapter 5), but at $\Cb_{\tau_0}$ for some $\tau_0>0$ and 
much smaller than the minimum of the upper bounds placed on $\delta$ by the smallness conditions 
stated above, including the smallness condition stated at the end of the preceding section which 
guarantees the recovery of the bootstrap assumptions as strict inequalities. 
In fact, denoting this minimal upper bound by $\delta_*$, we set, successively :
\begin{equation}
\tau_0=\tau_{0,m}:=2^{-m}\delta_* \ : \ m=1,2,...  
\label{14.a30}
\end{equation}
thus $\tau_0\rightarrow 0$ as $m\rightarrow\infty$. To provide the initial data on $\Cb_{\tau_0}$, we 
shall use the $N$th approximate solution $(x_N^\mu, b_N, \beta_{\mu,N})$ as a point of reference. 
We consider the restriction of this to $\Cb_{\tau_0}$. We first perform a transformation as in Section 
10.7, equations \ref{10.598} - \ref{10.601} to adapt the coordinate $\vartheta$ to the integral curves 
of $\Lb_N$ on $\Cb_{\tau_0}$. This sets:
\begin{equation}
b_N=0, \ \ \ \Lb_N=\frac{\partial}{\partial u} \ \ \ \mbox{: \ on $\Cb_{\tau_0}$}
\label{14.a31}
\end{equation}
The initial data to be provided on $\Cb_{\tau_0}$ consist of the $x^\mu$ and the $\beta_\mu$ on 
$\Cb_{\tau_0}$ as we are thinking of the coordinate $\vartheta$ on $\Cb_{\tau_0}$ to be adapted to 
the generators of $\Cb_{\tau_0}$ thus in correspondence with \ref{14.a31} we have:
\begin{equation}
b=0, \ \ \ \Lb=\frac{\partial}{\partial u} \ \ \ \mbox{: \ on $\Cb_{\tau_0}$}
\label{14.a32}
\end{equation}
Since we are to confine ourselves to values of $u$ not exceeding $\delta_*$, the initial data 
need only be provided on $\Cb_{\tau_0}^{\delta_*}$ which in the following we denote simply by 
$\Cb_{\tau_0}$. 

The setup of characteristic initial data has already been discussed in the first paragraph of 
Section 5.1. The initial data on $\Cb_{\tau_0}$ consist of the 6 functions 
$(x^\mu : \mu=0,1,2  \ ; \ \beta_\mu : \mu=0,1,2)$ as smooth functions on $\Cb_{\tau_0}$ which is 
identified with: 
\begin{equation}
\{(u,\vartheta) \ : \ u\in[\tau_0,\delta_*], \ \vartheta\in S^1\}
\label{14.a33}
\end{equation}
The 2nd member of the characteristic system \ref{2.62} constitutes on $\Cb_{\tau_0}$ 2 constraints on 
the initial data:
\begin{equation}
\frac{\partial x^i}{\partial u}=\Nb^i\frac{\partial t}{\partial u} \ \ : \ i=1,2
\label{14.a34}
\end{equation}
the $\Omega^\mu=\partial x^\mu/\partial\vartheta$ being defined on $\Cb_{\tau_0}$ by the data for the 
$x^\mu$, the $h_{\mu\nu}$ being defined on $\Cb_{\tau_0}$ by the data for the $\beta_\mu$, and the 
$N^i$, $\Nb^i$ being determined through these as in Section 2.2. Equation \ref{2.88}, the 2nd member 
of the wave system, induced on $\Cb_{\tau_0}$, namely the 3rd of \ref{2.92}:  
$$(\Lb x^\mu)\Omega\beta_\mu-(\Omega x^\mu)\Lb\beta_\mu=0$$ 
(here the spatial dimension is 2), constitutes on $\Cb_{\tau_0}$ the 3rd constraint on the initial data:
\begin{equation}
\frac{\partial x^\mu}{\partial u}\frac{\partial\beta_\mu}{\partial\vartheta}
-\frac{\partial x^\mu}{\partial\vartheta}\frac{\partial\beta_\mu}{\partial u}=0
\label{14.a35}
\end{equation}
Equation \ref{2.86}, the 1st member of the wave system, or \ref{2.91} in the form:  
$$N^\mu\Lb\beta_\mu-\lambda\sh^{-1}\Omega^\mu\Omega\beta_\mu=0$$
which it takes by virtue of \ref{2.72}, \ref{2.74} and \ref{2.93}, constitutes on $\Cb_{\tau_0}$ 
the 4th constraint on the initial data:
\begin{equation}
N^\mu\frac{\partial\beta_\mu}{\partial u}-\lambda\sh^{-1}\Omega^\mu
\frac{\partial\beta_\mu}{\partial\vartheta}=0
\label{14.a36}
\end{equation}
$\sh$ and $\lambda$ being determined on $\Cb_{\tau_0}$ by the initial data:
\begin{equation}
\sh=h_{\mu\nu}\Omega^\mu\Omega^\nu \ ; \ \ \lambda=c\rhob \ ; \ \  c=-\frac{1}{2}h_{\mu\nu}N^\mu\Nb^\nu 
\ , \ \rhob=\frac{\partial t}{\partial u}
\label{14.a37}
\end{equation}
Summarizing, the 6 functions $(x^\mu : \mu=0,1,2  \ ; \ \beta_\mu : \mu=0,1,2)$ on \ref{14.a33} 
representing the initial data on $\Cb_{\tau_0}$ are subject to the 4 constraints, equations 
\ref{14.a34}, \ref{14.a35}, \ref{14.a36}. Consequently there are 2 functional degrees of freedom 
on $\Cb_{\tau_0}$ in the initial data. 

In addition to the constraints on $\Cb_{\tau_0}$ there are constraints on 
$S_{\tau_0,\tau_0}$, its past boundary, which according to \ref{14.a33} is identified with: 
\begin{equation}
\{(\tau_0,\vartheta) \ : \ \vartheta\in S^1\}
\label{14.a38}
\end{equation}
The data for the $x^\mu$ restricted to $S_{\tau_0,\tau_0}$ define $\overline{S}_{\tau_0}$, a curve in 
$\mathbb{M}^2$, Minkowski spacetime of 2-spatial dimensions, parametrized on $S^1$, by:
\begin{equation}
\overline{S}_{\tau_0} \ : \ \vartheta\in S^1\mapsto (x^\mu(\tau_0,\vartheta) : \mu=0,1,2)\in 
\mathbb{M}^2 
\label{14.a39}
\end{equation}
the $x^\mu$ being rectangular coordinates in $\mathbb{M}^2$. The prior solution 
defines the $\beta_\mu$ as functions of the parameter $\vartheta$ along $\overline{S}_{\tau_0}$. 
We denote these functions by:
\begin{equation}
\left.\beta^\prime_\mu\right|_{\overline{S}_{\tau_0}} 
\label{14.a40}
\end{equation}
We also denote:
\begin{equation}
\left.\beta^{\prime\mu}\right|_{\overline{S}_{\tau_0}}=(g^{-1})^{\mu\nu}
\left.\beta^\prime_\nu\right|_{\overline{S}_{\tau_0}}, \ \ \ 
\left.\sigma^\prime\right|_{\overline{S}_{\tau_0}}=
-\left.\beta^{\prime\mu}\right|_{\overline{S}_{\tau_0}}
\left.\beta^\prime_\mu\right|_{\overline{S}_{\tau_0}}, \ \ \ 
\left. G^\prime\right|_{\overline{S}_{\tau_0}}=
G\circ\left.\sigma^\prime\right|_{\overline{S}_{\tau_0}}
\label{14.a41}
\end{equation}
(in the last $G$ stands for the function on the positive real line). 
Then the linear jump condition \ref{4.1}:
$$\Omega^\mu\triangle\beta_\mu=0$$
constitutes at $S_{\tau_0,\tau_0}$ a boundary constraint on the initial data:
\begin{equation}
\Omega^\mu(\tau_0,\vartheta)\left(\beta_\mu(\tau_0,\vartheta)
-\left.\beta^\prime_\mu\right|_{\overline{S}_{\tau_0}}(\vartheta)\right)=0 \ \ : \ 
\forall\vartheta\in S^1
\label{14.a42}
\end{equation}
Moreover, the nonlinear jump condition \ref{1.329}:
$$(\triangle\beta_\mu)\triangle(G\beta_\mu)=0$$
constitutes at $S_{\tau_0,\tau_0}$ another boundary constraint on the initial data: 
\begin{equation}
\left(\beta_\mu(\tau_0,\vartheta)
-\left.\beta^\prime_\mu\right|_{\overline{S}_{\tau_0}}(\vartheta)\right)
\left((G\beta^\mu)(\tau_0,\vartheta)-\left(\left.G^\prime\right|_{\overline{S}_{\tau_0}}
\left.\beta^{\prime\mu}\right|_{\overline{S}_{\tau_0}}\right)(\vartheta)\right)=0 \ \ : \ 
\forall\vartheta\in S^1
\label{14.a43}
\end{equation}
We remark that the linear jump condition \ref{4.2} constitutes the boundary condition for the 
1st derived datum $\lambdab$ (see Sections 4.2 and 5.1) in terms of $\lambda$, so this condition 
does not constrain the $\left.\beta_\mu\right|_{S_{\tau_0,\tau_0}}$. Summarizing, the 3 functions on 
$S^1$ representing the $\left.\beta_\mu\right|_{S_{\tau_0,\tau_0}}$ are subject to the 2 boundary 
constraints, equations \ref{14.a42}, \ref{14.a43}. 

Now the $N$th approximants $(x_N^\mu, \beta_{\mu,N})$ on $\Cb_{\tau_0}$, after the transformation 
in connection with achieving \ref{14.a31}, satisfy approximately the constraint equations 
\ref{14.a34}, \ref{14.a35}, \ref{14.a36} and the boundary constraints \ref{14.a42}, \ref{14.a43}. 
We shall presently estimate the errors committed, that is the quantities by which the $N$th 
approximants fail to satisfy these constraints. By the 2nd of \ref{9.15} and by Proposition 9.1 the 
error committed in regard to \ref{14.a34} is:
\begin{equation}
\vepb_N^i=O(\tau_0^{N+1})
\label{14.a44}
\end{equation}
By the 3rd of \ref{9.44} the error committed in regard to \ref{14.a35} is:
\begin{equation}
\omega_{\Lb\Omega,N}=O(\tau_0^{N+1})
\label{14.a45}
\end{equation}
Noting that by \ref{9.b6}, \ref{9.b7} the left hand side of \ref{14.a36} is:
\begin{equation}
-\frac{1}{\rho}\left\{a h^{-1}(dx^\mu,d\beta_\mu)+\frac{1}{2}(dx^\mu\wedge d\beta_\mu)(\Lb,L)\right\}
\label{14.a46}
\end{equation}
on $\Cb_{\tau_0}$, by \ref{9.a6}, the 1st of \ref{9.44} and by Propositions 9.2, 9.3 the error 
committed in regard to \ref{14.a36} is:
\begin{equation}
-\frac{1}{\rho_N}(\delta_N-\omega_{L\Lb,N})=O(\tau_0^{N-1})+
\left.\rhob_N^{-1}\right|_{\Cb_{\tau_0}}O(\tau_0^{N+1})
\label{14.a47}
\end{equation}

We turn to the errors committed in regard to the boundary constraints \ref{14.a42}, \ref{14.a43}. 
To assess these we define, in analogy with \ref{14.a39}, the curve $\overline{S}_{\tau_0,N}$ in 
$\mathbb{M}^2$ by:
\begin{equation}
\overline{S}_{\tau_0,N} \ : \ \vartheta\in S^1\mapsto (x_N^\mu(\tau_0,\vartheta) : \mu=0,1,2)\in 
\mathbb{M}^2 
\label{14.a48}
\end{equation}
The prior solution defines the $\beta_\mu$ as functions of the parameter $\vartheta$ along 
$\overline{S}_{\tau_0,N}$. In analogy with \ref{14.a40} we denote these functions by:
\begin{equation}
\left.\beta^\prime_\mu\right|_{\overline{S}_{\tau_0,N}} 
\label{14.a49}
\end{equation}
The error committed in regard to \ref{14.a42} is then:
\begin{equation}
\Omega_N^\mu(\tau_0,\vartheta)\left(\beta_{\mu,N}(\tau_0,\vartheta)
-\left.\beta^\prime_\mu\right|_{\overline{S}_{\tau_0,N}}(\vartheta)\right):=
\iota^\prime_{\Omega,\tau_0,N}(\vartheta)
\label{14.a50}
\end{equation}
a function on $S^1$. To estimate this error we remark that according to the definitions \ref{9.55} and 
\ref{9.59} we have:
\begin{equation}
\Omega_N^\mu(\tau_0,\vartheta)\left(\beta_{\mu,N}(\tau_0,\vartheta)
-\beta^\prime_\mu(f_N(\tau_0,\vartheta),w_N(\tau_0,\vartheta),\psi_N(\tau_0,\vartheta))\right)
=\iota_{\Omega,N}(\tau_0,\vartheta)
\label{14.a51}
\end{equation}
the prior solution being expressed here in terms of $(t,u^\prime,\vartheta^\prime)$ coordinates as in 
\ref{9.55}. Subtracting \ref{14.a51} from \ref{14.a50} we obtain:
\begin{eqnarray}
&&\iota^\prime_{\Omega,\tau_0,N}(\vartheta)=\iota_{\Omega,N}(\tau_0,\vartheta)+\nonumber\\
&&\hspace{18mm}\Omega_N^\mu(\tau_0,\vartheta)\left(\beta^\prime_\mu
(f_N(\tau_0,\vartheta),w_N(\tau_0,\vartheta),\psi_N(\tau_0,\vartheta))
-\left.\beta^\prime_\mu\right|_{\overline{S}_{\tau_0,N}}(\vartheta)\right)\nonumber\\
&&\label{14.a52}
\end{eqnarray}
In regard to the factor in parenthesis in the 2nd term, while according to the definitions \ref{9.51}:
\begin{equation}
f_N(\tau_0,\vartheta)=x_N^0(\tau_0,\vartheta), \ \ \ 
g_N^i(\tau_0,\vartheta)=x_N^i(\tau_0,\vartheta) \ : \ i=1,2
\label{14.a53}
\end{equation}
are, in view of \ref{14.a48} the rectangular coordinates of the point on 
$\overline{S}_{\tau_0,N}\subset\mathbb{M}^2$ corresponding to the parameter $\vartheta$,
$$x^{\prime i}(f_N(\tau_0,\vartheta),w_N(\tau_0,\vartheta),\psi_N(\tau_0,\vartheta))$$
the spatial rectangular coordinates of the point in the prior solution corresponding to 
$$t=f_N(\tau_0,\vartheta), \ u^\prime=w_N(\tau_0,\vartheta), \ 
\vartheta^\prime=\psi_N(\tau_0,\vartheta)$$
(see \ref{4.56}), differ from $g_N^i(\tau_0,\vartheta)$ by the quantity 
\begin{equation}
x^{\prime i}(f_N(\tau_0,\vartheta),w_N(\tau_0,\vartheta),\psi_N(\tau_N,\vartheta))
-g_N^i(\tau_0,\vartheta)=\tau_0^3 D_N^i(\tau_0,\vartheta)
\label{14.a54}
\end{equation}
by which, with the definitions \ref{9.53}, \ref{9.54} the $N$th approximants fail to satisfy the 
identification equations \ref{4.57}. Recall here Proposition 4.5 and definition \ref{9.119}. 
As a consequence of this error the factor in parenthesis in the 2nd term on the right in 
\ref{14.a52} fails to vanish. By Proposition 9.7, 
\begin{equation}
D_N^i(\tau_0,\vartheta)=O(\tau_0^{N-1})
\label{14.a55}
\end{equation}
 
Let us define $(w_{\tau_0}(\vartheta),\psi_{\tau_0}(\vartheta))$ so as to satisfy:
\begin{equation}
x^{\prime i}(f_N(\tau_0,\vartheta),w_{\tau_0}(\vartheta),\psi_{\tau_0}(\vartheta))-g_N^i(\tau_0,\vartheta)=0
\label{14.a56}
\end{equation} 
while being suitably close to $(w_N(\tau_0,\vartheta),\psi_N(\tau_0,\vartheta))$. Setting, in 
analogy with \ref{9.54}, 
\begin{equation}
w_{\tau_0}(\vartheta)=\tau_0 v_{\tau_0}(\vartheta), \ \ \ 
\psi_{\tau_0}(\vartheta)=\vartheta+\tau_0^3\gamma_{\tau_0}(\vartheta)
\label{14.a57}
\end{equation}
this means that the regularized identification equations of Proposition 4.5 must be satisfied 
with $\tau=\tau_0$ and $(v_{\tau_0}(\vartheta),\gamma_{\tau_0}(\vartheta))$ in the role of 
$(v,\gamma)$: 
\begin{equation}
\hat{F}_N^i((\tau_0,\vartheta),(v_{\tau_0}(\vartheta),\gamma_{\tau_0}(\vartheta))=0
\label{14.a58}
\end{equation}
if we replace in the definition of $\hat{F}^i$ at $\tau=\tau_0$, $\hat{f}(\tau_0,\vartheta)$ by:
$$\hat{f}_N(\tau_0,\vartheta)=\tau_0^{-2}\left(f_N(\tau_0,\vartheta)-f(0,\vartheta)\right)$$
and $\hat{\delta}^i(\tau_0,\vartheta)$ by:
$$\hat{\delta}_N^i(\tau_0,\vartheta)=\tau_0^{-3}\left[g_N^i(\tau_0,\vartheta)-g^i(0,\vartheta)
-N_0^i(\vartheta)(f_N(\tau_0,\vartheta)-f(0,\vartheta))\right]$$
(see \ref{4.222}, \ref{4.223}, \ref{4.240}, \ref{4.244}), so that $\hat{F}^i$ at $\tau=\tau_0$ 
is replaced, according to \ref{9.117}, \ref{9.118}, by $\hat{F}_N^i$ at $\tau=\tau_0$. Recall that 
according to \ref{9.119}: 
\begin{equation}
\hat{F}_N^i((\tau_0,\vartheta),(v_N(\tau_0,\vartheta),\gamma_N(\tau_0,\vartheta))
=D_N^i(\tau_0,\vartheta)
\label{14.a59}
\end{equation}
Moreover, $(v_{\tau_0}(\vartheta),\gamma_{\tau_0}(\vartheta))$ should be close to 
$(v_N(\tau_0,\vartheta),\gamma_N(\tau_0,\vartheta))$. In view of \ref{12.545}, 
provided that $\tau_0$ is suitably small, we can apply the implicit function theorem 
to conclude that there is a unique solution 
$(v_{\tau_0}(\vartheta),\gamma_{\tau_0}(\vartheta))$ of \ref{14.a58} in a neighborhood of 
$(v_N(\tau_0,\vartheta),\gamma_N(\tau_0,\vartheta))$. To ensure the required smallness of 
$\tau_0$ we start the sequence \ref{14.a30} from $m=M$ for some sufficiently large $M$. 
We can then furthermore show that by virtue of 
\ref{14.a55}:
\begin{equation}
v_{\tau_0}(\vartheta)-v_N(\tau_0,\vartheta)=O(\tau_0^{N-1}), \ \ \ 
\gamma_{\tau_0}(\vartheta)-\gamma_N(\tau_0,\vartheta)=O(\tau_0^{N-1})
\label{14.a60}
\end{equation}
therefore by \ref{14.a57} and \ref{9.54}:
\begin{equation}
w_{\tau_0}(\vartheta)-w_N(\tau_0,\vartheta)=O(\tau_0^{N}), \ \ \ 
\psi_{\tau_0}(\vartheta)-\psi_N(\tau_0,\vartheta)=O(\tau_0^{N+2})
\label{14.a61}
\end{equation}
Now, by \ref{14.a53} and \ref{14.a56} we have:
\begin{equation}
\left.\beta^\prime_\mu\right|_{\overline{S}_{\tau_0,N}}(\vartheta)=
\beta^\prime_\mu(f_N(\tau_0,\vartheta),w_{\tau_0}(\vartheta),\psi_{\tau_0}(\vartheta))
\label{14.a62}
\end{equation}
the prior solution being again expressed in terms of $(t,u^\prime,\vartheta^\prime)$ coordinates as in 
\ref{14.a51}. Therefore the factor in parenthesis in the 2nd term on the right in \ref{14.a52} is: 
\begin{eqnarray}
&&\beta^\prime_\mu(f_N(\tau_0,\vartheta),w_N(\tau_0,\vartheta),\psi_N(\tau_0,\vartheta))-
\beta^\prime_\mu(f_N(\tau_0,\vartheta),w_{\tau_0}(\vartheta),\psi_{\tau_0}(\vartheta))\nonumber\\
&&\hspace{52mm}=O(\tau_0^{N}) \label{14.a63}
\end{eqnarray}
A more precise estimate for the 2nd term on the right in \ref{14.a52} can in fact be obtained. 
This is because by \ref{4.102}, \ref{4.103} and Proposition 4.3:
\begin{eqnarray}
&&\Omega_0^\mu\left.\frac{\partial\beta^\prime_\mu}{\partial u^\prime}\right|_{\partial_-{\cal B}}
=\left.\left(\Omega^{\prime\mu}+\pi L^{\prime\mu}\right)\left(T^\prime\beta^\prime_\mu\right)
\right|_{\partial_-{\cal B}}\nonumber\\
&&\hspace{20mm}=\left.\left[\frac{\mu^\prime}{\alpha}\hat{T}^{\prime i}
\left(\frac{\partial\beta^\prime_i}{\partial\vartheta^\prime}
+\pi\frac{\partial\beta^\prime_i}{\partial t}\right)\right]\right|_{\partial_-{\cal B}}=0 \label{14.a64}
\end{eqnarray}
It then follows that the 2nd term on the right in \ref{14.a52} is in fact:
\begin{equation}
O(\tau_0^{N+1})
\label{14.a65}
\end{equation}
Combining then with the first result of Proposition 9.4 we conclude through \ref{14.a52} that:
\begin{equation}
\iota^\prime_{\Omega,\tau_0,N}=O(\tau_0^{N+1})
\label{14.a66}
\end{equation}

As for the error committed in regard to \ref{14.a43}, this being the nonlinear jump condition 
\ref{1.329} at $S_{\tau_0,\tau_0}$, we recall that the last is equivalent to \ref{9.107} at 
$S_{\tau_0,\tau_0}$, that is:
\begin{equation}
\ep(\tau_0,\vartheta)+j(\kappa(\tau_0,\vartheta),\epb(\tau_0,\vartheta))\epb^2(\tau_0,\vartheta)=0
\label{14.a67}
\end{equation}
where (see \ref{4.3}):
\begin{eqnarray}
&&\ep(\tau_0,\vartheta)=N^\mu(\tau_0,\vartheta)\left(\beta_\mu(\tau_0,\vartheta)
-\left.\beta^\prime_\mu\right|_{\overline{S}_{\tau_0}}(\vartheta)\right) \nonumber\\
&&\epb(\tau_0,\vartheta)=\Nb^\mu(\tau_0,\vartheta)\left(\beta_\mu(\tau_0,\vartheta)
-\left.\beta^\prime_\mu\right|_{\overline{S}_{\tau_0}}(\vartheta)\right) \label{14.a68}
\end{eqnarray}
$\kappa(\tau_0,\vartheta)$ stands for the quadruplet (see \ref{4.43}):
\begin{equation}
\kappa(\tau_0,\vartheta)=\left(\sigma(\tau_0,\vartheta),c(\tau_0,\vartheta),
(N^\mu\beta_\mu)(\tau_0,\vartheta),(\Nb^\mu\beta_\mu)(\tau_0,\vartheta)\right)
\label{14.a69}
\end{equation}
and $j$ stands for the function defined by Proposition 4.2. Recalling the definition \ref{9.55}, 
we have, by \ref{14.a62}, \ref{14.a63}:
\begin{equation}
\beta_{\mu,N}(\tau_0,\vartheta)-\left.\beta^\prime_\mu\right|_{\overline{S}_{\tau_0,N}}(\vartheta)
-\triangle_N\beta_\mu(\tau_0,\vartheta)=O(\tau_0^{N})
\label{14.a70}
\end{equation}
Hence, denoting:
\begin{eqnarray}
&&\ep^\prime_{\tau_0,N}(\vartheta)=N_N^\mu(\tau_0,\vartheta)\left(\beta_{\mu,N}(\tau_0,\vartheta)
-\left.\beta^\prime_\mu\right|_{\overline{S}_{\tau_0,N}}(\vartheta)\right) \nonumber\\
&&\epb^\prime_{\tau_0,N}(\vartheta)=\Nb_N^\mu(\tau_0,\vartheta)\left(\beta_{\mu,N}(\tau_0,\vartheta)
-\left.\beta^\prime_\mu\right|_{\overline{S}_{\tau_0,N}}(\vartheta)\right) \label{14.a71}
\end{eqnarray}
we have:
\begin{equation}
\ep^\prime_{\tau_0,N}(\vartheta)-\ep_N(\tau_0,\vartheta)=O(\tau_0^{N}), \ \ \ 
\epb^\prime_{\tau_0,N}(\vartheta)-\epb_N(\tau_0,\vartheta)=O(\tau_0^{N})
\label{14.a72}
\end{equation}
the quantities $\ep_N$, $\epb_N$ being defined by \ref{9.79}. A more precise estimate for the 1st 
difference can in fact be derived. This is because by \ref{4.102}, \ref{4.103} and Proposition 4.3:
\begin{eqnarray}
&&N_0^\mu\left.\frac{\partial\beta^\prime_\mu}{\partial u^\prime}\right|_{\partial_-{\cal B}}
=\left.\left(L^{\prime\mu}T^\prime\beta^\prime_\mu\right)
\right|_{\partial_-{\cal B}}\nonumber\\
&&\hspace{20mm}=\left.\left(\frac{\mu^\prime}{\alpha}\hat{T}^{\prime i}
\frac{\partial\beta^\prime_i}{\partial t}\right)\right|_{\partial_-{\cal B}}=0 \label{14.a73}
\end{eqnarray}
It follows that the 1st of \ref{14.a72} can be precised to:
\begin{equation}
\ep^\prime_{\tau_0,N}(\vartheta)-\ep_N(\tau_0,\vartheta)=O(\tau_0^{N+1})
\label{14.a74}
\end{equation}
The quantity by which \ref{14.a67} fails to be satisfied by the $N$th approximants is:
\begin{equation}
\ep^\prime_{\tau_0,N}(\vartheta)+j(\kappa_N(\tau_0,\vartheta),\epb^\prime_{\tau_0,N}(\vartheta)) \
\epb^{\prime 2}_{\tau_0,N}(\vartheta):=d^\prime_{\tau_0,N}(\vartheta)
\label{14.a75}
\end{equation}
Comparing with the definition \ref{9.111} we see that by virtue of \ref{14.a74} and the 2nd of 
\ref{14.a72}:
\begin{equation}
d^\prime_{\tau_0,N}(\vartheta)-d_N(\tau_0,\vartheta)=O(\tau_0^{N+1})
\label{14.a76}
\end{equation}
Therefore by Proposition 9.6:
\begin{equation}
d^\prime_{\tau_0,N}=O(\tau_0^{N+1})
\label{14.a77}
\end{equation}

%Let us denote by ${\cal D}_*$ the domain:
%\begin{equation}
%{\cal D}_*=\left\{(\tau_0,u,\vartheta) \ : \ \tau_0\in[0,\tau_{0,M}], \ u\in [\tau_0,\delta_*], \ 
%\vartheta\in S^1\right\}
%\label{14.a78}
%\end{equation}

By virtue of the estimates \ref{14.a44}, \ref{14.a45}, \ref{14.a47}, and 
\ref{14.a66}, \ref{14.a77}, for the quantities by which by the $N$th approximants fail to satisfy 
the constraint equations \ref{14.a34}, \ref{14.a35}, \ref{14.a36}, on $\Cb_{\tau_0}$, and the 
boundary constraints \ref{14.a42}, \ref{14.a67} (which is equivalent to \ref{14.a43}), on 
$S_{\tau_0,\tau_0}$, we can show that we can find smooth initial data $x^\mu$, $\beta_\mu$ on 
$\Cb_{\tau_0}$,  satisfying the constraint equations on 
$\Cb_{\tau_0}$ and the boundary constraints on $S_{\tau_0,\tau_0}$, and verifying the following 
closeness conditions to the corresponding $N$th approximants:
\begin{eqnarray}
&&\left.x^\mu\right|_{\Cb_{\tau_0}}-\left.x_N^\mu\right|_{\Cb_{\tau_0}}=
O(\tau_0^{N-1})+\left.\rhob_N^{-1}\right|_{\Cb_{\tau_0}}O(\tau_0^{N+1}) \nonumber\\
&&\left.\beta_\mu\right|_{\Cb_{\tau_0}}-\left.\beta_{\mu,N}\right|_{\Cb_{\tau_0}}=
O(\tau_0^{N-1})+\left.\rhob_N^{-1}\right|_{\Cb_{\tau_0}}O(\tau_0^{N+1}) \label{14.a79}
\end{eqnarray}
Moreover, this is true for each $\tau_0\in[0,\tau_{0,M}]$, where $\tau_{0,M}$ is the starting point of 
the sequence \ref{14.a30} as discussed above in connection with equation \ref{14.a58}. To simplify the 
notation we introduce the following. 

{\bf Definition 14.1:} \ \ Let 
$\{\left.f\right|_{\Cb_{\tau_0}}  :  \tau_0\in[0,\tau_{0,M}]\}$ be a 1-parameter family of functions 
such that for each $\tau_0\in[0,\tau_{0,M}]$, $\left.f\right|_{\Cb_{\tau_0}}$ is a smooth function on 
$\Cb_{\tau_0}$. Let $p$ be a positive integer and $q$ a non-negative integer, $q<p$. We write:
$$\left.f\right|_{\Cb_{\tau_0}}=\cO(u^{-q}\tau_0^p)$$
if for each pair of non-negative integers $k$, $l$ there is a constant $C_{k,l}$ independent of 
$\tau_0$ such that:
$$\left|\frac{\partial^k}{\partial u^k}\frac{\partial^l}{\partial\vartheta^l}
\left.f\right|_{\Cb_{\tau_0}}\right|\leq C_{k,l}u^{-q-k}\tau_0^{p}$$
(compare with Definitions 9.1 and 9.2). 

It follows from \ref{14.a79} that in the sense of Definition 14.1:
\begin{eqnarray}
&&\left(\left.x^\mu\right|_{\Cb_{\tau_0}}-\left.x_N^\mu\right|_{\Cb_{\tau_0}}\right)
=\cO(\tau_0^{N-1})\nonumber\\
&&\left(\left.\beta_\mu\right|_{\Cb_{\tau_0}}-\left.\beta_{\mu,N}\right|_{\Cb_{\tau_0}}\right)
=\cO(\tau_0^{N-1}) \label{14.a92}
\end{eqnarray}
The first of \ref{14.a92} implies:
\begin{equation}
\left(\left.\Omega^\mu\right|_{\Cb_{\tau_0}}-\left.\Omega_N^\mu\right|_{\Cb_{\tau_0}}\right)
=\cO(\tau_0^{N-1})
\label{14.a93}
\end{equation}
which, combined with the second of \ref{14.a92}, in turn implies:
\begin{eqnarray}
&&\left(\left.N^\mu\right|_{\Cb_{\tau_0}}-\left.N_N^\mu\right|_{\Cb_{\tau_0}}\right)
=\cO(\tau_0^{N-1}) \nonumber\\
&&\left(\left.\Nb^\mu\right|_{\Cb_{\tau_0}}-\left.\Nb_N^\mu\right|_{\Cb_{\tau_0}}\right)
=\cO(\tau_0^{N-1}) \label{14.a94}
\end{eqnarray}
and:
\begin{equation}
\left(\left.\sh\right|_{\Cb_{\tau_0}}-\left.\sh_N\right|_{\Cb_{\tau_0}}\right)=\cO(\tau_0^{N-1})  \label{14.a95}
\end{equation}
Moreover, \ref{14.a94} together with the second of \ref{14.a92} imply:
\begin{equation}
\left(\left.c\right|_{\Cb_{\tau_0}}-\left.c_N\right|_{\Cb_{\tau_0}}\right)=\cO(\tau_0^{N-1}) \label{14.a96}
\end{equation}
Also, \ref{14.a95} together with the first of \ref{14.a92} implies:
\begin{equation}
\left(\left.\pi\right|_{\Cb_{\tau_0}}-\left.\pi_N\right|_{\Cb_{\tau_0}}\right)=\cO(\tau_0^{N-1}) \label{14.a97}
\end{equation}

To simplify the notation in regard to boundary data we introduce the following.

{\bf Definition 14.2:} \ \ Let 
$\{\left.f\right|_{S_{\tau_0,\tau_0}}  :  \tau_0\in[0,\tau_{0,M}]\}$ be a 1-parameter family of functions 
such that for each $\tau_0\in[0,\tau_{0,M}]$, $\left.f\right|_{S_{\tau_0,\tau_0}}$ is a smooth function on 
$S_{\tau_0,\tau_0}$, and let $p$ be a positive integer. We write:
$$\left.f\right|_{S_{\tau_0,\tau_0}}=\ocO(\tau_0^p)$$
if for each non-negative integer $l$ there is a constant $C_l$ independent of 
$\tau_0$ such that:
$$\left|\frac{\partial^l}{\partial\vartheta^l}
\left.f\right|_{S_{\tau_0,\tau_0}}\right|\leq C_l\tau_0^{p}$$

Note that if $\{\left.f\right|_{\Cb_{\tau_0}}  :  \tau_0\in[0,\tau_{0,M}]\}$ is a 1 parameter family 
of functions as in Definition 14.1 such that $\left.f\right|_{\Cb_{\tau_0}}=\cO(u^{-q}\tau_0^p)$ then the 
corresponding family $\{\left.f\right|_{S_{\tau_0,\tau_0}}  :  \tau_0\in[0,\tau_{0,M}]\}$ is a 
1 parameter family of functions as in Definition 14.2 such that 
$\left.f\right|_{S_{\tau_0,\tau_0}}=\ocO(\tau_0^{p-q})$. Thus the 1st of \ref{14.a92} implies:
\begin{equation}
\left.x^\mu\right|_{S_{\tau_0,\tau_0}}-\left.x_N^\mu\right|_{S_{\tau_0,\tau_0}}=\ocO(\tau_0^{N-1})
\label{14.a80}
\end{equation}
that is, the parametric equations of the curves $\overline{S}_{\tau_0}$ (see \ref{14.a40}) and 
$\overline{S}_{\tau_0,N}$ (see \ref{14.a48}) in $\mathbb{M}^2$ differ by $\ocO(\tau_0^{N-1})$. 
Setting, in accordance with \ref{4.53},
\begin{equation}
f(\tau_0,\vartheta)=x^0(\tau_0,\vartheta), \ \ \ g^i(\tau_0,\vartheta)=x^i(\tau_0,\vartheta) 
\ : \ i=1,2
\label{14.a81}
\end{equation}
these being the rectangular coordinates of the point on $\overline{S}_{\tau_0}$ corresponding to the 
parameter $\vartheta$, from the point of view of the prior solution the 
$(t,u^\prime,\vartheta^\prime)$ coordinates of this point are:
\begin{equation}
t=f(\tau_0,\vartheta), \ \ u^\prime=w(\tau_0,\vartheta), \ \ \vartheta^\prime=\psi(\tau_0,\vartheta)
\label{14.a82}
\end{equation}
(see \ref{4.54}, \ref{4.56}). The pair $(w(\tau_0,\vartheta),\psi(\tau_0,\vartheta)$ is determined 
by the identification equations \ref{4.57} at $S_{\tau_0,\tau_0}$:
\begin{equation}
x^{\prime i}(f(\tau_0,\vartheta),w(\tau_0,\vartheta),\psi(\tau_0,\vartheta))=g^i(\tau_0,\vartheta) \ \ : \ i=1,2
\label{14.a83}
\end{equation}
or (see \ref{4.60} - \ref{4.63}):
\begin{equation}
F^i((\tau_0,\vartheta), (w(\tau_0,\vartheta),\psi(\tau_0,\vartheta)-\vartheta))=0 \ : \ i=1,2
\label{14.a84}
\end{equation}
where:
\begin{equation}
F^i((\tau_0,\vartheta),(u^\prime,\varphi^\prime))=
x^{\prime i}(f(\tau_0,\vartheta),u^\prime,\vartheta+\varphi^\prime)-g^i(\tau_0,\vartheta)
\label{14.a85}
\end{equation}
According to Proposition 4.5 at $\tau=\tau_0$, setting 
\begin{equation}
w(\tau_0,\vartheta)=\tau_0 v(\tau_0,\vartheta), \ \ \ 
\psi(\tau_0,\vartheta)=\vartheta+\tau_0^3\gamma(\tau_0,\vartheta)
\label{14.a86}
\end{equation}
we have:
\begin{equation}
F^i((\tau_0,\vartheta),(\tau_0 v,\tau_0^3\gamma))=\tau_0^3\hat{F}^i((\tau_0,\vartheta),(v,\gamma))
\label{14.a87}
\end{equation}
and the pair $(v(\tau_0,\vartheta),\gamma(\tau_0,\vartheta))$ is determined as the solution of the 
regularized identification equations:
\begin{equation}
\hat{F}^i((\tau_0,\vartheta),(v(\tau_0,\vartheta),\gamma(\tau_0,\vartheta))=0
\label{14.a88}
\end{equation}
the functions $F^i((\tau_0,\vartheta),(v,\gamma)) : i=1,2$ being given by Proposition 4.5 at 
$\tau=\tau_0$ in terms of:
\begin{eqnarray}
&&\hat{f}(\tau_0,\vartheta)=\tau_0^{-2}(f(\tau_0,\vartheta)-f(0,\vartheta))\nonumber\\
&&\hat{\delta}^i(\tau_0,\vartheta)=\tau_0^{-3}\left[g^i(\tau_0,\vartheta)-g^i(0,\vartheta)
-N_0^i(\vartheta)(f(\tau_0,\vartheta)-f(0,\vartheta))\right]\nonumber\\
&&\label{14.a89}
\end{eqnarray}
(see \ref{4.222}, \ref{4.223}, \ref{4.240}, \ref{4.244}). On the other hand the pair 
$(v_N(\tau_0,\vartheta),\gamma_N(\tau_0,\vartheta))$ satisfies \ref{14.a59} with the right hand side 
satisfying the estimate \ref{14.a55}. Moreover, according to \ref{9.117} the functions 
$\hat{F}^i_N((\tau_0,\vartheta),(v,\gamma))$ are obtained from the functions 
$\hat{F}^i((\tau_0,\vartheta),(v,\gamma))$ by replacing $\hat{f}(\tau_0,\vartheta)$, 
$\hat{\delta}^i(\tau_0,\vartheta)$ by $\hat{f}_N(\tau_0,\vartheta)$, 
$\hat{\delta}^i_N(\tau_0,\vartheta)$ respectively. Using the fact that by \ref{14.a80} and 
\ref{14.a81}, \ref{9.114} - \ref{9.116}:
\begin{eqnarray}
&&\hat{f}(\tau_0,\vartheta)-\hat{f}_N(\tau_0,\vartheta)=\ocO(\tau_0^{N-3}) \nonumber\\
&&\hat{\delta}^i(\tau_0,\vartheta)-\hat{\delta}_N^i(\tau_0,\vartheta)=\ocO(\tau_0^{N-4}) \label{14.a90}
\end{eqnarray}
we deduce (see \ref{12.1}):
\begin{equation}
\left.\chf\right|_{S_{\tau_0,\tau_0}}=\ocO(\tau_0^{N-3}), \ \ \ 
\left.\cv\right|_{S_{\tau_0,\tau_0}}, \left.\cga\right|_{S_{\tau_0,\tau_0}} =\ocO(\tau_0^{N-4}) 
\label{14.a91}
\end{equation}
(the 1st being a restatement of the 1st of \ref{14.a90}). 

The 1st of \ref{14.a92} implies:
\begin{equation}
\left.\rhob\right|_{\Cb_{\tau_0}}-\left.\rhob_N\right|_{\Cb_{\tau_0}}=\cO(u^{-1}\tau_0^{N-1})
\label{14.a98}
\end{equation}
which, together with \ref{14.a96} implies:
\begin{equation}
\left.\lambda\right|_{\Cb_{\tau_0}}-\left.\lambda_N\right|_{\Cb_{\tau_0}}=\cO(u^{-1}\tau_0^{N-1})
\label{14.a99}
\end{equation}
The 2nd of \ref{14.a92} together with \ref{14.a93} - \ref{14.a95} implies:
\begin{eqnarray}
&&\left.\left(\sbeta-\sbeta_N\right)\right|_{\Cb_{\tau_0}}, 
\left.\left(\beta_N-\beta_{N,N}\right)\right|_{\Cb_{\tau_0}}, 
\left.\left(\beta_{\Nb}-\beta_{\Nb,N}\right)\right|_{\Cb_{\tau_0}} = \cO(\tau_0^{N-1})\nonumber\\
&&\left.\left(\sss-\sss_N\right)\right|_{\Cb_{\tau_0}}, 
\left.\left(\ss_N-\ss_{N,N}\right)\right|_{\Cb_{\tau_0}}, 
\left.\left(\ss_{\Nb}-\ss_{\Nb,N}\right)\right|_{\Cb_{\tau_0}} = \cO(\tau_0^{N-1})\nonumber\\
&&\left.\left(s_{\Nb\Lb}-s_{\Nb\Lb,N}\right)\right|_{\Cb_{\tau_0}} = \cO(u^{-1}\tau_0^{N-1}) 
\label{14.a100}
\end{eqnarray}
Given that $\sk$, $\skb$ are determined by $EN^\mu$, $E\Nb^\mu$ through the first two 
of the expansions \ref{3.a15} and that
$\tchi$, $\tchib$ are determined by $\sk$, $\skb$ according to \ref{3.a21}, and also similarly that 
$\sk_N$, $\skb_N$ are determined through $E_N N_N^\mu$, $E_N\Nb_N^\mu$ through the expansions 
\ref{10.117} and that $\tchi_N$, $\tchib_N$ are determined by $\sk_N$, $\skb_N$ according to 
\ref{10.133}, \ref{10.142}, the estimates \ref{14.a93} - \ref{14.a95} together with \ref{14.a100} 
imply:
\begin{equation}
\left.(\tchi-\tchi_N)\right|_{\Cb_{\tau_0}}, 
\left.(\tchib-\tchib_N)\right|_{\Cb_{\tau_0}} = \cO(\tau_0^{N-1})
\label{14.a101}
\end{equation}

We proceed to the 1st derived data $\lambdab$, $s_{NL}$ on $\Cb_{\tau_0}$. Here we follow the 
argument of Chapter 5. The pair $(\lambdab,s_{NL})$ satisfy along the generators of $\Cb_{\tau_0}$, the 
lines of constant $\vartheta$, the linear inhomogeneous system of ordinary differential equations 
of Proposition 5.1. Revisiting the argument leading to this proposition in Section 5.1, we see that 
the 1st equation of the system is simply the propagation equation for $\lambdab$ of Proposition 3.3, 
while the 2nd equation with the right hand side brought to the left is:
\begin{equation}
-\Omega a N^\mu\square_{\tilde{h}}\beta_\mu=0 \ \mbox{: on $\Cb_{\tau_0}$} 
\label{14.a102}
\end{equation}
In regard to the coefficients in the 1st equation, we find:
\begin{eqnarray}
&&a_{11}=\pb+\frac{\eta^2}{4c^2}\frac{dH}{d\sigma}\beta_{\Nb}^2(\beta_N\sss-\sbeta\ss_N)\lambda
\nonumber\\
&&a_{12}=\frac{\eta^2}{4c^2}\frac{dH}{d\sigma}\beta_{\Nb}^3 \label{14.a103}
\end{eqnarray}
The coefficients in the 2nd equation are expressed by \ref{5.33} in terms of $f_\mu$, $e_\mu$, $g$. 
We find:
\begin{equation}
f_\mu=f^{\Lb}\Lb\beta_\mu+f^E E\beta_\mu-\rhob E(E\beta_\mu)
\label{14.a104}
\end{equation}
where:
\begin{eqnarray}
&&f^{\Lb}=\frac{1}{2c}\left\{\tchi+\frac{1}{2}\left[\frac{dH}{d\sigma}(\eta^2\sbeta^2+\sigma)+2H\right]
(\beta_N\sss-2\sbeta\ss_N)\right\} \nonumber\\
&&f^E=-\frac{1}{c}\left\{2\teta+\frac{1}{2}\rhob\left[\eta^2\frac{dH}{d\sigma}\sbeta\beta_{\Nb}
(\beta_N\sss-2\sbeta\ss_N)+\frac{c}{\eta^2}\left(\sigma\frac{dH}{d\sigma}+2H\right)E\sigma\right]
\right\}\nonumber\\
&&\label{14.a105}
\end{eqnarray}
with $\teta$ given by the 1st of \ref{3.a25}. We also find:
\begin{equation}
e_\mu=e^{\Lb}\Lb\beta_\mu+e^E E\beta_\mu 
\label{14.a106}
\end{equation}
where:
\begin{eqnarray}
&&e^{\Lb}=\frac{\beta_{\Nb}}{4c}\left[\frac{dH}{d\sigma}(\eta^2\sbeta^2+\sigma)+2H\right] \nonumber\\
&&e^E=-\frac{\eta^2}{2c}\rhob\sbeta\beta_{\Nb}^2\frac{dH}{d\sigma}
\label{14.a107}
\end{eqnarray}
Moreover, the function $g$ is given by:
\begin{equation}
g=\frac{1}{2}\left(\chib+\frac{1}{2}\Lb\log\Omega\right)
\label{14.a108}
\end{equation}
Recall here that $\chib$ is given in terms of $\tchib$ by \ref{3.a20}. Substituting \ref{14.a104} 
and \ref{14.a106} in \ref{5.33} we obtain:
\begin{eqnarray}
&&a_{21}=\ss_N(c^{-1}\sn-f^E)+\sss(\on-\lambda f^{\Lb})-\rhob N^\mu E(E\beta_\mu)\nonumber\\
&&a_{22}=n-g-\ss_N e^E-\sss\lambda e^{\Lb} \label{14.a109}
\end{eqnarray}

The quantities $\lambdab_N$, $s_{NL,N}$ defined by the $N$th approximate solution satisfy along 
$Cb_{\tau_0}$ a linear homogeneous system of ordinary differential equations:
\begin{eqnarray}
&&\frac{\partial\lambdab_N}{\partial u}=a_{11,N}\lambdab_N+a_{12,N}s_{NL,N}+\vep_{1,N}\nonumber\\
&&\frac{\partial s_{NL,N}}{\partial u}=a_{21,N}\lambdab_N+a_{22,N}s_{NL,N}+\vep_{2,N} 
\label{14.a110}
\end{eqnarray}
the 1st equation being the propagation equation \ref{10.164} along $\Cb_{\tau_0}$, so the 
inhomogeneous term is:
\begin{equation}
\vep_{1,N}=\vep_{\lambda,N}=\cO(\tau_0^N)
\label{14.a111}
\end{equation}
the last following from the estimate \ref{10.167}. The 2nd of the equations \ref{14.a110}
corresponds to the equation of Proposition 9.10 multiplied by $N_N^\mu$, that is 
\begin{equation}
-\Omega_N a_N N_N^\mu\square_{\tilde{h}^\prime_N}\beta_{\mu,N}=-N_N^\mu\tilde{\kappa}^\prime_{\mu,N}
\label{14.a112}
\end{equation}
(compare with \ref{14.a102}) along $\Cb_{\tau_0}$, so the inhomogeneous term is:
\begin{equation}
\vep_{2,N}=-N_N^\mu\tilde{\kappa}^\prime_{\mu,N}=\cO(u^{-1}\tau_0^{N-2})
\label{14.a113}
\end{equation}
the last following from the estimate of Proposition 9.10. The coefficients $a_{11,N}$, $a_{12,N}$, 
$a_{21,N}$, $a_{22,N}$ are the $N$th approximant analogues of the corresponding coefficients $a_{11}$, 
$a_{12}$, $a_{21}$, $a_{22}$ of the system of Proposition 5.1. Subtracting the system \ref{14.110} 
from the system of Proposition 5.1 we obtain the following linear inhomogeneous system for the 
differences
\begin{equation}
\clab=\lambdab-\lambdab_N, \ \ \check{s}_{NL}=s_{NL,N} 
\label{14.a114}
\end{equation}
along $\Cb_{\tau_0}$:
\begin{eqnarray}
&&\frac{\partial\clab}{\partial u}=a_{11}\clab+a_{12}\check{s}_{NL}+b_1\nonumber\\
&&\frac{\partial\check{s}_{NL}}{\partial u}=a_{21}\clab+a_{22}\check{s}_{NL}+b_2
\label{14.a115}
\end{eqnarray}
where the inhomogeneous terms are:
\begin{eqnarray}
&&b_1=(a_{11}-a_{11,N})\lambdab_N+(a_{12}-a_{12,N})s_{NL,N}-\vep_{1,N} \nonumber\\
&&b_2=(a_{21}-a_{21,N})\lambdab_N+(a_{22}-a_{22,N})s_{NL,N}-\vep_{2,N} 
\label{14.a116}
\end{eqnarray}
Using \ref{14.a100}, \ref{14.a101}, as well as \ref{14.a98}, \ref{14.a99} together with \ref{14.93} - 
\ref{14.a97} and the 2nd of \ref{14.a92}, we deduce:
\begin{eqnarray}
&&\left.(a_{11}-a_{11,N})\right|_{\Cb_{\tau_0}}, \left.(a_{12}-a_{12,N})\right|_{\Cb_{\tau_0}} 
= \cO(u^{-1}\tau_0^{N-1})\nonumber\\
&&\left.(a_{21}-a_{21,N})\right|_{\Cb_{\tau_0}}, \left.(a_{22}-a_{22,N})\right|_{\Cb_{\tau_0}} 
= \cO(u^{-1}\tau_0^{N-1})
\label{14.a117}
\end{eqnarray}
Then in view of the estimates \ref{14.a111}, \ref{14.a113} we obtain:
\begin{eqnarray}
&&b_1=\cO(u^{-1}\tau_0^N)\nonumber\\
&&b_2=\cO(u^{-1}\tau_0^{N-2}) \label{14.a118}
\end{eqnarray}

We proceed to analyze the boundary conditions for the system \ref{14.a115} on $S_{\tau_0,\tau_0}$. 
We first consider the boundary condition for $\clab$. Subtracting \ref{9.102} from \ref{4.6}, the 
last in the form 
$$r\lambdab=\lambda,$$
we obtain:
\begin{equation}
r\clab+\lambdab_N\check{r}=\cla+\hat{\nu}_N \ \ \mbox{: on $S_{\tau_0,\tau_0}$}
\label{14.a119}
\end{equation}
where we denote:
\begin{equation}
\cla=\lambda-\lambda_N, \ \ \check{r}=r-r_N
\label{14.a120}
\end{equation}
Hence the boundary condition for $\clab$ on $S_{\tau_0,\tau_0}$ is:
\begin{equation}
\clab=\frac{\cla-\lambdab_N\check{r}+\hat{\nu}_N}{r_N+\check{r}}
\label{14.a121}
\end{equation}
By \ref{14.a99}:
\begin{equation}
\left.\cla\right|_{S_{\tau_0,\tau_0}}=\ocO(\tau_0^{N-2})
\label{14.a122}
\end{equation}
Also, by \ref{9.106}:
\begin{equation}
\left.\hat{\nu}_N\right|_{S_{\tau_0,\tau_0}}=\ocO(\tau_0^{N+1})
\label{14.a123}
\end{equation}
To estimate $\check{r}$ we note that by \ref{4.6} and Proposition 4.2:
\begin{equation}
r=-\frac{\ep}{\epb}=j(\kappa,\epb)\epb
\label{14.a124}
\end{equation}
while by \ref{9.99} and \ref{9.111}:
\begin{equation}
r_N=-\frac{\ep_N}{\epb_N}=j(\kappa_N,\epb_N)\epb_N)-\frac{d_N}{\epb_N}
\label{14.a125}
\end{equation}
Subtracting \ref{14.a125} from \ref{14.a124} we obtain:
\begin{equation}
\check{r}=j(\kappa,\epb)\epb-j(\kappa_N,\epb_N)\epb_N+\frac{d_N}{\epb_N}
\label{14.a126}
\end{equation}
By the first triplet of the estimates \ref{14.a100} and by \ref{14.a96} and the 2nd of \ref{14.a92} 
(see\ref{9.108}, \ref{9.109}):
\begin{equation}
\left.(\kappa-\kappa_N)\right|_{S_{\tau_0,\tau_0}}=\ocO(\tau_0^{N-1})
\label{14.a127}
\end{equation}
Also, by \ref{9.105} and Proposition 9.6:
\begin{equation}
\left.\frac{d_N}{\epb_N}\right|_{S_{\tau_0,\tau_0}}=\ocO(\tau_0^N)
\label{14.a128}
\end{equation}
Thus, what remains to be estimated in order to estimate $\check{r}$ is the difference 
$\epb-\epb_N$ on $S_{\tau_0,\tau_0}$. According to the 2nd of \ref{4.3} and \ref{9.79} 
this difference is:
\begin{equation}
\epb-\epb_N=\Nb^\mu\triangle\beta_\mu-\Nb_N^\mu\triangle_N\beta_\mu \ \ \mbox{: at $S_{\tau_0,\tau_0}$}
\label{14.a129}
\end{equation}
Moreover, by \ref{4.150} and \ref{9.55}:
\begin{eqnarray}
&&\left.\triangle\beta_\mu\right|_{S_{\tau_0,\tau_0}}=\left.\beta_\mu\right|_{S_{\tau_0,\tau_0}}
-\left.\beta^\prime_{\mu}\right|_{\overline{S}_{\tau_0}}\nonumber\\
&&\left.\triangle_N\beta_\mu\right|_{S_{\tau_0,\tau_0}}=\left.\beta_{\mu,N}\right|_{S_{\tau_0,\tau_0}}
-\left.\beta^\prime_{\mu}\right|_{\overline{S}_{\tau_0,N}} \label{14.a130}
\end{eqnarray}
and:
\begin{eqnarray}
&&\left.\beta^\prime_{\mu}\right|_{\overline{S}_{\tau_0}}(\vartheta)
=\beta^\prime_\mu(f(\tau_0,\vartheta),w(\tau_0,\vartheta,\psi(\tau_0,\vartheta))\nonumber\\
&&\left.\beta^\prime_{\mu}\right|_{\overline{S}_{\tau_0,N}}(\vartheta)
=\beta^\prime_\mu(f_N(\tau_0,\vartheta),w_N(\tau_0,\vartheta,\psi_N(\tau_0,\vartheta)) \label{14.a131}
\end{eqnarray}
The estimates \ref{14.a91} imply (see \ref{14.a86} and \ref{9.54}):
\begin{eqnarray}
&&\left.(f-f_N)\right|_{S_{\tau_0,\tau_0}}=\ocO(\tau_0^{N-1})\nonumber\\
&&\left.(w-w_N)\right|_{S_{\tau_0,\tau_0}}=\ocO(\tau_0^{N-3})\nonumber\\
&&\left.(\psi-\psi_N)\right|_{S_{\tau_0,\tau_0}}=\ocO(\tau_0^{N-1}) \label{14.a132}
\end{eqnarray}
It follows that:
\begin{equation}
\left.\beta^\prime_\mu\right|_{\overline{S}_{\tau_0}}
-\left.\beta^\prime_{\mu}\right|_{\overline{S}_{\tau_0,N}}=\ocO(\tau_0^{N-3})
\label{14.a133}
\end{equation}
which, together with the 2nd of \ref{14.a92}, implies:
\begin{equation}
\left.\triangle\beta_\mu\right|_{S_{\tau_0,\tau_0}}-
\left.\triangle_N\beta_\mu\right|_{S_{\tau_0,\tau_0}}=\ocO(\tau_0^{N-3})
\label{14.a134}
\end{equation}
This, together with the 2nd of \ref{14.a94}, implies through \ref{14.a129}:
\begin{equation}
\left.(\epb-\epb_N)\right|_{S_{\tau_0,\tau_0}}=\ocO(\tau_0^{N-3})
\label{14.a135}
\end{equation}
This estimate combined with \ref{14.a127} and \ref{14.a128} yields through \ref{14.a126}:
\begin{equation}
\left.\check{r}\right|_{S_{\tau_0,\tau_0}}=\ocO(\tau_0^{N-3})
\label{14.a136}
\end{equation}
Using the estimates \ref{14.a122}, \ref{14.a123} and \ref{14.a136}, we obtain through \ref{14.a121} 
the desired boundary estimate for $\clab$:
\begin{equation}
\left.\clab\right|_{S_{\tau_0,\tau_0}}=\ocO(\tau_0^{N-3})
\label{14.a137}
\end{equation}

We turn to the boundary condition for $s_{NL}$. To derive this condition we start from the remark that 
any local solution containing a future neighborhood of $S_{\tau_0,\tau_0}$ in ${\cal K}$ must satisfy 
the boundary condition \ref{14.a67} with $\tau_0$ replaced by any $\tau\in[\tau_0,\tau_0+\varepsilon)$ 
for some $\varepsilon>0$. Differentiating then with respect to $\tau$ at $\tau=\tau_0$, that is 
applying $T$ to 
$$\ep+j(\kappa,\epb)\epb^2=0$$
and evaluating the result at $S_{\tau_0,\tau_0}$ we obtain:
\begin{equation}
T\ep+\left(2j(\kappa,\epb)+\epb\frac{\partial j}{\partial\epb}(\kappa,\epb)\right)\epb T\epb 
+\epb^2\frac{\partial j}{\partial\kappa}(\kappa,\epb)\cdot T\kappa=0 \ \ \mbox{: at $S_{\tau_0,\tau_0}$}
\label{14.a138}
\end{equation}
From \ref{4.3}, to determine $T\ep$ and $T\epb$ at $S_{\tau_0,\tau_0}$ we must determine 
$\left.(T\triangle\beta_\mu)\right|_{S_{\tau_0,\tau_0}}$, which in view of the 1st of \ref{14.a130},  
\ref{14.a131} involves the derivatives of the transformation functions $v$, $\gamma$ with respect to 
$\tau$ at $\tau=\tau_0$. In this regard, recall from \ref{4.119} that:
\begin{equation}
\left.\frac{\partial f}{\partial\tau}\right|_{S_{\tau_0,\tau_0}}
=\left.(\rho+\rhob)\right|_{S_{\tau_0,\tau_0}}, \ \ 
\left.\frac{\partial g^i}{\partial\tau}\right|_{S_{\tau_0,\tau_0}}
=\left.(\rho N^i+\rhob\Nb^i)\right|_{S_{\tau_0,\tau_0}} 
\label{14.a139}
\end{equation}
are already determined, $\rho=c^{-1}\lambdab$ on $S_{\tau_0,\tau_0}$ having just been determined. 

From the 1st of \ref{14.a130} we have:
\begin{equation}
\left.T\triangle\beta_\mu\right|_{S_{\tau_0,\tau_0}}=\left.T\beta_\mu\right|_{S_{\tau_0,\tau_0}}
-\left.T\beta^\prime_\mu\right|_{\overline{S}_{\tau_0}}
\label{14.a140}
\end{equation}
hence, since by  \ref{4.3},
\begin{eqnarray}
&&T\ep=N^\mu T\triangle\beta_\mu+(TN^\mu)\triangle\beta_\mu \nonumber\\
&&T\epb=\Nb^\mu T\triangle\beta_\mu+(T\Nb^\mu)\triangle\beta_\mu 
\label{14.a78}
\end{eqnarray}
we obtain:
\begin{eqnarray}
&&\left.T\ep\right|_{S_{\tau_0,\tau_0}}=\left.(s_{NL}+\lambda\sss)\right|_{S_{\tau_0,\tau_0}}
-\left.N^\mu T\beta^\prime_\mu\right|_{\overline{S}_{\tau_0}}+
\left.(TN^\mu)\right|_{S_{\tau_0,\tau_0}}\left.\triangle\beta_\mu\right|_{S_{\tau_0,\tau_0}}\nonumber\\
&&\left.T\epb\right|_{S_{\tau_0,\tau_0}}=\left.(s_{\Nb\Lb}+\lambdab\sss)\right|_{S_{\tau_0,\tau_0}}
-\left.\Nb^\mu T\beta^\prime_\mu\right|_{\overline{S}_{\tau_0}}+
\left.(T\Nb^\mu)\right|_{S_{\tau_0,\tau_0}}\left.\triangle\beta_\mu\right|_{S_{\tau_0,\tau_0}}
\nonumber\\
&&\label{14.a141}
\end{eqnarray}
In regard to the last term on the right in each of the above, since according to \ref{4.4},
\begin{equation}
\triangle\beta_\mu=-\frac{1}{2c}h_{\mu\nu}(\epb N^\nu+\ep\Nb^\mu)
\label{14.a142}
\end{equation}
by the last four of the expansions \ref{3.a15} we have: 
\begin{eqnarray}
&&(TN^\mu)\triangle\beta_\mu=(m+n)\ep+(\om+\on)\epb \nonumber\\
&&(T\Nb^{\mu})\triangle\beta_\mu=(\omb+\onb)\epb+(\mb+\nb)\ep \label{14.a143}
\end{eqnarray}
Now, the terms in $n$, $\on$ are contributed by the expansion of $\Lb N^\mu$ and the terms in 
$\mb$, $\omb$ are contributed by the expansion of $\Lb\Nb^\mu$, hence do not depend on the 1st derived 
data $(\lambdab,s_{NL})$ being directly defined in terms of the initial data. On the other hand, 
\begin{eqnarray}
&&\om=\frac{1}{4c}(\beta_N^2 LH+2H\beta_N s_{NL})\nonumber\\
&&m=-\pi\sm-\om, \ \ \sm=-\sbeta\beta_N LH+\frac{\beta_N^2}{2c}\lambdab EH-H\sbeta s_{NL} 
\label{14.a144}
\end{eqnarray}
(see \ref{3.57}, \ref{3.68} and 1st of \ref{3.47}) and:
\begin{eqnarray}
&&LH=\frac{dH}{d\sigma}L\sigma\nonumber\\
&&L\sigma=\frac{\eta^2}{c}\left[(\beta_N\sss-2\sbeta\ss_N)\lambdab+\beta_{\Nb}s_{NL}\right]
\label{14.a145}
\end{eqnarray}
(see \ref{5.14}, \ref{5.15}), so $\om$ and $m$ depend, linearly, on the 1st derived data. 
Also, 
\begin{equation}
\nb=\frac{1}{4c}(\beta_{\Nb}^2 LH+2H\beta_{\Nb}\sss\lambdab)
\label{14.a146}
\end{equation}
(see \ref{3.62}) so this likewise depends, linearly, on the 1st derived data. 
Moreover, by the 2nd of \ref{3.48}:
\begin{equation}
\onb=-\pi\snb-\nb
\label{14.a147}
\end{equation}
and by the 2nd of \ref{3.a24}:
\begin{equation}
\snb=E\lambdab+(c^{-1}\beta_N\beta_{\Nb}EH-2\okb)\lambdab-\sbeta\beta_{\Nb}LH
\label{14.a148}
\end{equation}
This also depends linearly on the first derived data, but depends on $E\lambdab$ as well as $\lambdab$. 
However $\lambdab$ on $S_{\tau_0,\tau_0}$ with its derivatives tangentially to $S_{\tau_0,\tau_0}$ 
has already been estimated by \ref{14.a137}. Note that by \ref{14.a144} - \ref{14.a148} the coefficient 
of $s_{NL}$ in each of \ref{14.a143} is $\ocO(\tau_0)$ on $S_{\tau_0,\tau_0}$ (see \ref{14.a135} and 
\ref{9.89}). 

Consider next $T\kappa$ on $S_{\tau_0,\tau_0}$ which enters the 3rd term on the left in \ref{14.a138}. 
Decomposing $T\kappa$ into the sum of $L\kappa$ and $\Lb\kappa$, the second does not depend on the 
1st derived data $(\lambdab,s_{NL})$ being directly defined in terms of the initial data. 
As for 
\begin{equation}
L\kappa=(L\sigma, Lc, L\beta_N, L\beta_{\Nb})
\label{14.a149}
\end{equation}
$L\sigma$ is given by \ref{14.a145}, $Lc$ is given by \ref{3.85}:
\begin{equation}
Lc=-\frac{1}{2}\beta_N\beta_{\Nb}LH-\frac{1}{2}H(\beta_N\lambdab\sss+\beta_{\Nb}s_{NL})
+c(m+\onb)
\label{14.a150}
\end{equation}
and we have:
\begin{eqnarray}
&&L\beta_N=s_{NL}+m\beta_N+\om\beta_{\Nb}+\sm\sb\nonumber\\
&&L\beta_{\Nb}=\lambdab\sss+\nb\beta_N+\onb\beta_{\Nb}+\snb\sbeta
\label{14.a151}
\end{eqnarray}
By the above $L\kappa$ depends linearly on the 1st derived data with the coefficient of $s_{NL}$ 
being $\ocO(1)$ on $S_{\tau_0,\tau_0}$. 

We now turn to the derivation of the estimates for $(Tv,T\gamma)$ on $S_{\tau_0,\tau_0}$ which are 
needed in order to estimate $\left.T\beta^\prime_\mu\right|_{\overline{S}_{\tau_0}}$, thus the 
2nd term on the right in each of \ref{14.a141}. The required estimates are to be obtained by 
differentiating the identification equations (see Proposition 4.5) implicitly with respect to $\tau$ 
at $\tau=\tau_0$:
\begin{equation}
\frac{\partial\hat{F}^i}{\partial v}\frac{\partial v}{\partial\tau}+
\frac{\partial\hat{F}^i}{\partial\gamma}\frac{\partial\gamma}{\partial\tau}
+\frac{\partial\hat{F}^i}{\partial\tau}=0
\label{14.a152}
\end{equation}
the arguments of the partial derivatives of $\hat{F}^i$ being $((\tau_0,\vartheta),(v,\gamma))$ 
with $(v,\gamma)=(v(\tau_0,\vartheta),\gamma(\tau_0,\vartheta))$, and the arguments of 
$(\partial v/\partial\tau,\partial\gamma/\partial\tau)$ being $(\tau_0,\vartheta)$. From \ref{14.a152} 
the corresponding equation for the $N$th approximants is to be subtracted, namely the derivative of 
\ref{9.119} with respect to $\tau$ at $\tau=\tau_0$:
\begin{equation}
\frac{\partial\hat{F}_N^i}{\partial v}\frac{\partial v_N}{\partial\tau}+
\frac{\partial\hat{F}_N^i}{\partial\gamma}\frac{\partial\gamma_N}{\partial\tau}
+\frac{\partial\hat{F}_N^i}{\partial\tau}=\frac{\partial D_N^i}{\partial\tau}
\label{14.a153}
\end{equation}
the arguments of the partial derivatives of $\hat{F}_N^i$ being $((\tau_0,\vartheta),(v,\gamma))$ 
with $(v,\gamma)=(v_N(\tau_0,\vartheta),\gamma_N(\tau_0,\vartheta))$, and the arguments of 
$(\partial v_N/\partial\tau,\partial\gamma_N/\partial\tau)$ and of $\partial D_N^i/\partial\tau$ 
being $(\tau_0,\vartheta)$. Subtracting \ref{14.a153} from \ref{14.a152} gives:
\begin{eqnarray}
&&\frac{\partial\hat{F}^i}{\partial v}\frac{\partial\cv}{\partial\tau}
+\frac{\partial\hat{F}^i}{\partial\gamma}\frac{\partial\cga}{\partial\tau}=
-\frac{\partial\hat{F}^i}{\partial\tau}+\frac{\partial\hat{F}^i}{\partial\tau}\nonumber\\
&&\hspace{32mm}-\left(\frac{\partial\hat{F}^i}{\partial v}-\frac{\partial\hat{F}_N^i}{\partial v}\right)
\frac{\partial v_N}{\partial\tau}
-\left(\frac{\partial\hat{F}^i}{\partial\gamma}-\frac{\partial\hat{F}_N^i}{\partial\gamma}\right)\frac{\partial\gamma_N}{\partial\tau}\nonumber\\
&&\hspace{32mm}-\frac{\partial D_N^i}{\partial\tau} \label{14.a154}
\end{eqnarray}
The principal part of the right hand side is the difference constituted by the first two terms:
\begin{eqnarray}
&&-\frac{\partial\hat{F}^i}{\partial\tau}+\frac{\partial\hat{F}^i}{\partial\tau}=
-S_0^i(\vartheta)v_N(\tau_0,\vartheta)l(\vartheta)
\frac{\partial\chf}{\partial\tau}(\tau_0,\vartheta)
+\frac{\partial\chdl^i}{\partial\tau}(\tau_0,\vartheta)\nonumber\\
&&\hspace{32mm}\mbox{: to leading terms}
\label{14.a155}
\end{eqnarray}
by Proposition 4.5 and \ref{9.117}, \ref{9.118}. Now, the analogues of \ref{14.a139} for the $N$th 
approximants are, by the definitions \ref{9.51} and \ref{9.15} (see \ref{9.14}):
\begin{equation}
\left.\frac{\partial f_N}{\partial\tau}\right|_{S_{\tau_0,\tau_0}}=
\left.(\rho_N+\rhob_N)\right|_{S_{\tau_0,\tau_0}}, \ \ 
\left.\frac{\partial g_N^i}{\partial\tau}\right|_{S_{\tau_0,\tau_0}}=
\left.(\rho_N N_N^i+\rhob_N\Nb_N^i+\vep_N^i+\vepb_N^i)\right|_{S_{\tau_0,\tau_0}}
\label{14.a156}
\end{equation}
Subtracting \ref{14.156} from \ref{14.a139} we obtain:
\begin{eqnarray}
&&\left.\frac{\partial\check{f}}{\partial\tau}\right|_{S_{\tau_0,\tau_0}}=
\left.\left(\rho-\rho_N+\rhob-\rhob_N\right)\right|_{S_{\tau_0,\tau_0}}\nonumber\\
&&\left.\frac{\partial\check{g}^i}{\partial\tau}\right|_{S_{\tau_0,\tau_0}}=
\left.\left(\rho N^i-\rho_N N_N^i+\rhob\Nb^i-\rhob_N\Nb_N^i-\vep_N^i-\vepb_N^i\right)
\right|_{S_{\tau_0,\tau_0}} \nonumber\\
&&\label{14.a157}
\end{eqnarray}
In view of the estimate \ref{14.a137} together with \ref{14.a84}, and \ref{14.a98}, \ref{14.a92} and 
Proposition 9.1 we conclude that:
\begin{equation}
\left.\frac{\partial\check{f}}{\partial\tau}\right|_{S_{\tau_0,\tau_0}}=\ocO(\tau_0^{N-3}), \ \ 
\left.\frac{\partial\check{g}^i}{\partial\tau}\right|_{S_{\tau_0,\tau_0}}=\ocO(\tau_0^{N-3})
\label{14.a158}
\end{equation}
It follows that:
\begin{equation}
\left.\frac{\partial\chf}{\partial\tau}\right|_{S_{\tau_0,\tau_0}}=
\tau_0^{-2}\left.\frac{\partial\check{f}}{\partial\tau}\right|_{S_{\tau_0,\tau_0}}
-2\tau_0^{-1}\left.\chf\right|_{S_{\tau_0,\tau_0}}=\ocO(\tau_0^{N-5})
\label{14.a159}
\end{equation}
(see \ref{14.a91}) and: 
\begin{equation}
\left.\frac{\partial\chdl^i}{\partial\tau}\right|_{S_{\tau_0,\tau_0}}=
\tau_0^{-3}\left.\left(\frac{\partial\check{g}^i}{\partial\tau}
-N_0^i\frac{\partial\check{f}}{\partial\tau}\right)\right|_{S_{\tau_0,\tau_0}}
-3\tau_0^{-1}\left.\chdl^i\right|_{S_{\tau_0,\tau_0}}=\ocO(\tau_0^{N-5})
\label{14.a160}
\end{equation}
(see \ref{14.a90}) the last because from \ref{14.a157} the leading part of 
$$\left.\left(\frac{\partial\check{g}^i}{\partial\tau}
-N_0^i\frac{\partial\check{f}}{\partial\tau}\right)\right|_{S_{\tau_0,\tau_0}}$$ 
is:
$$\left.(N^i-N_0^i)(\rho-\rho_N)\right|_{S_{\tau_0,\tau_0}}=\ocO(\tau_0)\ocO(\tau_0^{N-3})
=\ocO(\tau_0^{N-2})$$
The difference constituted by the first two terms on the right in \ref{14.a154} is then 
$\ocO(\tau_0^{N-5})$. By \ref{14.a90}, \ref{14.a91} the difference constituted by the next two 
terms on the right in \ref{14.a154} is $\ocO(\tau_0^{N-4})$. As for the last term on the right in 
\ref{14.a154} it is $\ocO(\tau_0^{N-2})$ by Proposition 9.7. We then conclude through \ref{14.a154} 
that:
\begin{equation}
\left.\left(\frac{\partial\cv}{\partial\tau},\frac{\partial\cga}{\partial\tau}\right)
\right|_{S_{\tau_0,\tau_0}}=\ocO(\tau_0^{N-5})
\label{14.a161}
\end{equation}

We now return to \ref{14.a130}. From the 2nd of \ref{14.a130} we have:
\begin{equation}
\left.T\triangle_N\beta_\mu\right|_{S_{\tau_0,\tau_0}}=\left.T\beta_{\mu,N}\right|_{S_{\tau_0,\tau_0}}
-\left.T\beta^\prime_\mu\right|_{\overline{S}_{\tau_0,N}}
\label{14.a162}
\end{equation}
Subtracting \ref{14.a162} from \ref{14.a140} we obtain:
\begin{eqnarray}
&&\left.(T\triangle\beta_\mu-T\triangle_N\beta_\mu)\right|_{S_{\tau_0,\tau_0}}=
\left.(T\beta_\mu-T\beta_{\mu,N})\right|_{S_{\tau_0,\tau_0}}\nonumber\\
&&\hspace{25mm}-\left.\frac{\partial\beta^\prime}{\partial t}\right|_{\overline{S}_{\tau_0}}
\left.\frac{\partial(f-f_N)}{\partial\tau}\right|_{S_{\tau_0,\tau_0}}
-\left.\frac{\partial\beta^\prime_\mu}{\partial u^\prime}\right|_{\overline{S}_{\tau_0}}
\left.\frac{\partial(w-w_N)}{\partial\tau}\right|_{S_{\tau_0,\tau_0}}\nonumber\\
&&\hspace{45mm}
-\left.\frac{\partial\beta^\prime}{\partial\vartheta^\prime}\right|_{\overline{S}_{\tau_0}}
\left.\frac{\partial(\psi-\psi_N)}{\partial\tau}\right|_{S_{\tau_0,\tau_0}}\nonumber\\
&&\hspace{40mm}-\left(\left.\frac{\partial\beta^\prime_\mu}{\partial t}\right|_{\overline{S}_{\tau_0}}
-\left.\frac{\partial\beta^\prime_\mu}{\partial t}\right|_{\overline{S}_{\tau_0,N}}\right)
\left.\frac{\partial f_N}{\partial\tau}\right|_{S_{\tau_0,\tau_0}}\nonumber\\
&&\hspace{40mm}-\left(\left.\frac{\partial\beta^\prime_\mu}{\partial u^\prime}\right|_{\overline{S}_{\tau_0}}
-\left.\frac{\partial\beta^\prime_\mu}{\partial u^\prime}\right|_{\overline{S}_{\tau_0,N}}\right)
\left.\frac{\partial w_N}{\partial\tau}\right|_{S_{\tau_0,\tau_0}}\nonumber\\
&&\hspace{40mm}-\left(\left.\frac{\partial\beta^\prime_\mu}{\partial\vartheta^\prime}\right|_{\overline{S}_{\tau_0}}
-\left.\frac{\partial\beta^\prime_\mu}{\partial\vartheta^\prime}\right|_{\overline{S}_{\tau_0,N}}\right)
\left.\frac{\partial\psi_N}{\partial\tau}\right|_{S_{\tau_0,\tau_0}}\nonumber\\
\label{14.a163}
\end{eqnarray}
Here, we recall that for any function $\phi^\prime=\phi^\prime(t,u^\prime,\vartheta^\prime)$ defined 
by the prior solution, such as the partial derivatives of $\beta^\prime_\mu$ with respect to the 
$t, u^\prime, \vartheta^\prime$ coordinates, we denote:
\begin{eqnarray}
&&\left.\phi^\prime\right|_{\overline{S}_{\tau_0}}(\vartheta)=
\phi^\prime(f(\tau_0,\vartheta),w(\tau_0,\vartheta,\psi(\tau_0,\vartheta))\nonumber\\
&&\left.\phi^\prime\right|_{\overline{S}_{\tau_0,N}}(\vartheta)=
\phi^\prime(f_N(\tau_0,\vartheta),w_N(\tau_0,\vartheta),\psi_N(\tau_0,\vartheta))\label{14.a164}
\end{eqnarray}
Since on $S_{\tau_0,\tau_0}$:
\begin{equation}
f-f_N=\tau_0^2\chf, \ \ w-w_N=\tau_0\cv, \ \ \psi-\psi_N=\tau_0^3\cga
\label{14.a165}
\end{equation}
the estimates \ref{14.a91} imply:
\begin{equation}
f-f_N=\ocO(\tau_0^{N-1}), \ \ w-w_N=\ocO(\tau_0^{N-3}), \ \ \psi-\psi_N=\ocO(\tau_0^{N-1})
\label{14.a166}
\end{equation}
Also, the estimates \ref{14.a159}, \ref{14.a161} imply:
\begin{eqnarray}
&&\left.\frac{\partial(f-f_N)}{\partial\tau}\right|_{S_{\tau_0,\tau_0}}=
\tau_0^2\left.\frac{\partial\chf}{\partial\tau}\right|_{S_{\tau_0,\tau_0}}
+2\tau_0\left.\chf\right|_{S_{\tau_0,\tau_0}}=\ocO(\tau_0^{N-3})\nonumber\\
&&\left.\frac{\partial(w-w_N)}{\partial\tau}\right|_{S_{\tau_0,\tau_0}}=
\tau_0\left.\frac{\partial\cv}{\partial\tau}\right|_{S_{\tau_0,\tau_0}}
+\left.\cv\right|_{S_{\tau_0,\tau_0}}=\ocO(\tau_0^{N-4})\nonumber\\
&&\left.\frac{\partial(\psi-\psi_N)}{\partial\tau}\right|_{S_{\tau_0,\tau_0}}=
\tau_0^3\left.\frac{\partial\cga}{\partial\tau}\right|_{S_{\tau_0,\tau_0}}
+3\tau_0^2\left.\cga\right|_{S_{\tau_0,\tau_0}}=\ocO(\tau_0^{N-2})\nonumber\\
&&\label{14.a167}
\end{eqnarray}
By virtue of the above estimates, \ref{14.a163} implies:
\begin{eqnarray}
&&\left.(T\triangle\beta_\mu-T\triangle_N\beta_\mu)\right|_{S_{\tau_0,\tau_0}}=
\left.(T\beta_\mu-T\beta_{\mu,N})\right|_{S_{\tau_0,\tau_0}}\nonumber\\
&&\hspace{35mm}+\ocO(\tau_0^{N-4})\left.\frac{\partial\beta^\prime_\mu}{\partial u^\prime}
\right|_{\overline{S}_{\tau_0,N}}+\ocO(\tau_0^{N-3})\nonumber\\
&&\label{14.a168}
\end{eqnarray}

We now return to the boundary condition for $s_{NL}$ on $S_{\tau_0,\tau_0}$, revisiting the argument 
which started below \ref{14.a137} and ended with the analysis of $L\kappa$ (see \ref{14.a149} - 
\ref{14.a151}). The analogue of \ref{14.a138} for the $N$th approximants, obtained by applying $T$ to 
\ref{9.111} and evaluating the result at $S_{\tau_0,\tau_0}$, is: 
\begin{eqnarray}
&&T\ep_N+\left(2j(\kappa_N,\epb_N)+
\epb_N\frac{\partial j}{\partial\epb}(\kappa_N,\epb_N)\right)\epb_N T\epb_N \nonumber\\
&&\hspace{7mm}+\epb_N^2\frac{\partial j}{\partial\kappa}(\kappa_N,\epb_N)\cdot T\kappa_N=Td_N \ \ \ \mbox{: at $S_{\tau_0,\tau_0}$}
\label{14.a169}
\end{eqnarray}
By Proposition 9.6:
\begin{equation}
\left.Td_N\right|_{S_{\tau_0,\tau_0}}=\ocO(\tau_0^N)
\label{14.a170}
\end{equation}
Subtracting \ref{14.a169} from \ref{14.a138} we obtain: 
\begin{eqnarray}
&&(T\ep-T\ep_N)+\left(2j(\kappa,\epb)+\epb\frac{\partial j}{\partial\epb}(\kappa,\epb)\right)
(\epb T\epb-\epb_N T\epb_N)\nonumber\\
&&+\left(2j(\kappa,\epb)-2j(\kappa_N,\epb_N)+\epb\frac{\partial j}{\partial\epb}(\kappa,\epb)
-\epb_N\frac{\partial j}{\partial\epb}(\kappa_N,\epb_N)\right)\epb_N T\epb_N\nonumber\\
&&+(\epb^2-\epb_N^2)\frac{\partial j}{\partial\kappa}(\kappa,\epb)\cdot T\kappa
+\epb_N^2\left(\frac{\partial j}{\partial\kappa}(\kappa,\epb)\cdot T\kappa
-\frac{\partial j}{\partial\kappa}(\kappa_N,\epb_N)\cdot T\kappa_N\right)\nonumber\\
&&=-Td_N \label{14.a171}
\end{eqnarray}
By \ref{9.79},
\begin{eqnarray}
&&T\ep_N=N_N^\mu T\triangle_N\beta_\mu+(TN_N^\mu)\triangle_N\beta_\mu \nonumber\\
&&T\epb_N=\Nb_N^\mu T\triangle_N\beta_\mu+(T\Nb_N^\mu)\triangle_N\beta_\mu 
\label{14.a172}
\end{eqnarray}
hence, subtracting from \ref{14.a79},
\begin{eqnarray}
&&\left.T\ep\right|_{S_{\tau_0,\tau_0}}-\left.T\ep_N\right|_{S_{\tau_0,\tau_0}}=
\left.N^\mu T(\triangle\beta_\mu-\triangle_N\beta_\mu)\right|_{S_{\tau_0,\tau_0}}\nonumber\\
&&\hspace{35mm}+\left.(N^\mu-N_N^\mu)T\triangle_N\beta_\mu\right|_{S_{\tau_0,\tau_0}}\nonumber\\
&&\hspace{20mm}+\left.(TN^\mu)\triangle\beta_\mu\right|_{S_{\tau_0,\tau_0}}
-\left.(TN_N^\mu)\triangle_N\beta_\mu\right|_{S_{\tau_0,\tau_0}}
\label{14.a173}
\end{eqnarray}
and:
\begin{eqnarray}
&&\left.T\epb\right|_{S_{\tau_0,\tau_0}}-\left.T\epb_N\right|_{S_{\tau_0,\tau_0}}=
\left.\Nb^\mu T(\triangle\beta_\mu-\triangle_N\beta_\mu)\right|_{S_{\tau_0,\tau_0}}\nonumber\\
&&\hspace{35mm}+\left.(\Nb^\mu-\Nb_N^\mu)T\triangle_N\beta_\mu\right|_{S_{\tau_0,\tau_0}}\nonumber\\
&&\hspace{20mm}+\left.(T\Nb^\mu)\triangle\beta_\mu\right|_{S_{\tau_0,\tau_0}}
-\left.(T\Nb_N^\mu)\triangle_N\beta_\mu\right|_{S_{\tau_0,\tau_0}}
\label{14.a174}
\end{eqnarray} 
In view of the estimate \ref{14.a137} and the ensuing discussion ending with the analysis of $L\kappa$, 
the difference which the last two terms in \ref{14.a173} constitute is:
\begin{equation}
\ocO(\tau_0)(s_{NL}-s_{NL,N})+\ocO(\tau_0^{N-3})
\label{14.a175}
\end{equation}
Also, the 2nd term on the right in \ref{14.a173} is:
\begin{equation}
\ocO(\tau_0^{N-1})
\label{14.a176}
\end{equation}
On the other hand, by \ref{14.a168} the 1st term on the right in \ref{14.a173} is:
\begin{eqnarray}
&&\left.N^\mu(T\triangle\beta_\mu-T\triangle_N\beta_\mu)\right|_{S_{\tau_0,\tau_0}}=
\left.(N^\mu T\beta_\mu-N_N^\mu T\beta_{\mu,N})\right|_{S_{\tau_0,\tau_0}}\nonumber\\
&&\hspace{35mm}+\ocO(\tau_0^{N-4})\left.N^\mu\frac{\partial\beta^\prime_\mu}{\partial u^\prime}
\right|_{\overline{S}_{\tau_0,N}}+\ocO(\tau_0^{N-3})\nonumber\\
&&\label{14.a177}
\end{eqnarray}
Because of the fact that, according to \ref{9.92}, 
\begin{equation} 
\left. \left(N^\mu\frac{\partial\beta^\prime_\mu}{\partial u^\prime}\right)\right|_{\partial_-{\cal B}}
=0
\label{14.a178}
\end{equation}
the 2nd term on the right in \ref{14.a177} is actually $\ocO(\tau_0^{N-3})$. As for the 1st term on the  
right in \ref{14.a177} it is:
\begin{eqnarray}
&&\left.(N^\mu T\beta_\mu-N_N^\mu T\beta_{\mu,N})\right|_{S_{\tau_0,\tau_0}}=
\left.\left(s_{NL}-s_{NL,N}+\lambda\sss-\lambda_N\sss_N\right)\right|_{S_{\tau_0,\tau_0}}
+\ocO(\tau_0^{N-1})\nonumber\\
&&\hspace{44mm}=\left.\check{s}_{NL}\right|_{S_{\tau_0,\tau_0}}+\ocO(\tau_0^{N-2})
\label{14.a179}
\end{eqnarray} 
Here we have used the fact that by the 1st of \ref{9.b10}, the 1st of \ref{9.44}, and by \ref{9.a11} 
(see also \ref{9.42} and \ref{9.101}):
\begin{eqnarray}
&&N_N^\mu\Lb_N\beta_{\mu,N}-\lambda_N\sss_N=\frac{1}{\rho_N}\left(\rho_N N_N^\mu\Lb_N\beta_{\mu,N}
-a_N\sss_N\right)\nonumber\\
&&\hspace{20mm}=\frac{1}{\rho_N}\left[(L_N x_N^\mu)\Lb_N\beta_{\mu,N}-\vep_N^\mu\Lb_N\beta_{\mu,N}
-a_N(E_N x_N^\mu)E_N\beta_{\mu,N}\right]\nonumber\\
&&\hspace{6mm}=\frac{1}{\rho_N}\left[\frac{1}{2}(L_N x_N^\mu)\Lb_N\beta_{\mu,N}
+\frac{1}{2}(\Lb_N x_N^\mu)L_N\beta_{\mu,N}-a_N(E_N x_N^\mu)E_N\beta_{\mu,N}\right.\nonumber\\
&&\hspace{45mm}\left.+\frac{1}{2}\omega_{L\Lb,N}-\vep_N^\mu\Lb_N\beta_{\mu,N}\right]\nonumber\\
&&\hspace{20mm}=\frac{1}{\rho_N}\left(-\delta_N^\prime+\frac{1}{2}\omega_{L\Lb,N}-\vep_N^\mu\Lb_N\beta_{\mu,N}\right)=O(\tau^{N-1}) \label{14.a180}
\end{eqnarray}
the last by virtue of Propositions 9.1 and 9.2 and Lemma 9.2. 

In view of the estimate \ref{14.a137} and the ensuing discussion ending with the analysis of $L\kappa$, 
the difference which the last two terms in \ref{14.a174} constitute is:
\begin{equation}
\ocO(\tau_0)(s_{NL}-s_{NL,N})+\ocO(\tau_0^{N-3})
\label{14.a181}
\end{equation}
Also, the 2nd term on the right in \ref{14.a174} is:
\begin{equation}
\ocO(\tau_0^{N-1})
\label{14.a182}
\end{equation}
On the other hand, by \ref{14.a168} the 1st term on the right in \ref{14.a173} is:
\begin{equation}
\left.N^\mu(T\triangle\beta_\mu-T\triangle_N\beta_\mu)\right|_{S_{\tau_0,\tau_0}}=
\left.(N^\mu T\beta_\mu-N_N^\mu T\beta_{\mu,N})\right|_{S_{\tau_0,\tau_0}}+\ocO(\tau_0^{N-4})
\label{14.a183}
\end{equation}
and we have:
\begin{eqnarray}
&&\left.(\Nb^\mu T\beta_\mu-\Nb_N^\mu T\beta_{\mu,N})\right|_{S_{\tau_0,\tau_0}}=
\left.\left(s_{\Nb\Lb}-s_{\Nb\Lb,N}+\lambdab\sss-\lambdab_N\sss_N\right)\right|_{S_{\tau_0,\tau_0}}
+\ocO(\tau_0^{N-2})\nonumber\\
&&\hspace{44mm}=\ocO(\tau_0^{N-3})
\label{14.a184}
\end{eqnarray} 
Here we have used the estimate \ref{14.a137} and the conjugate of the expression \ref{14.a180}, 
which gives:
\begin{eqnarray}
&&\Nb_N^\mu L_N\beta_{\mu,N}-\lambdab_N\sss_N=\frac{1}{\rhob_N}\left(\rhob_N \Nb_N^\mu L_N\beta_{\mu,N}
-a_N\sss_N\right)\nonumber\\
&&\hspace{20mm}=\frac{1}{\rhob_N}\left(-\delta_N^\prime-\frac{1}{2}\omega_{L\Lb,N}-\vepb_N^\mu L_N\beta_{\mu,N}\right)=\rhob_N^{-1}O(\tau^{N}) \nonumber\\
&&\label{14.a185}
\end{eqnarray}
by virtue of Propositions 9.1 and 9.2 and Lemma 9.2. 

Finally using the analysis of $L\kappa$, in particular \ref{14.a149} - \ref{14.a151}, we deduce:
\begin{equation} 
\left.(T\kappa-T\kappa_N)\right|_{S_{\tau_0,\tau_0}}=
\ocO(1)\left.\check{s}_{NL}\right|_{S_{\tau_0,\tau_0}}+\ocO(\tau_0^{N-3})
\label{14.a186}
\end{equation}

The results \ref{14.a175} - \ref{14.a179} in regard to \ref{14.a173}, the results \ref{14.a181} - 
\ref{14.a184} in regard to \ref{14.a174}, together with the result \ref{14.a186}, taking into account 
of the fact that $\epb_N=O(\tau)$ (see \ref{9.89}), reduce equation \ref{14.a171} at 
$S_{\tau_0,\tau_0}$ to the form:
$$(1+\ocO(\tau_0))\left.\check{s}_{NL}\right|_{S_{\tau_0,\tau_0}}=\ocO(\tau_0^{N-3})$$
which, taking $M$ suitably large so that $\tau_{0,M}$ is suitably small, implies: 
\begin{equation}
\left.\check{s}_{NL}\right|_{S_{\tau_0,\tau_0}}=\ocO(\tau_0^{N-3})
\label{14.a187}
\end{equation}

We now revisit the linear inhomogeneous system of ordinary differential equations \ref{14.a115} 
along the generators of $\Cb_{\tau_0}$. The initial conditions on $S_{\tau_0,\tau_0}$ for this system 
satisfy the estimates \ref{14.a137}, \ref{14.a187}. The coefficients satisfy the estimates 
\ref{14.a117} and the inhomogeneous terms the estimates \ref{14.a118}. Integrating then yields 
the following estimate for the 1st derived data on $\Cb_{\tau_0}$:
\begin{equation}
\left.\clab\right|_{\Cb_{\tau_0}}, \left.\check{s}_{NL}\right|_{\Cb_{\tau_0}} = \cO(\tau_0^{N-3})
\label{14.a188}
\end{equation} 
The system \ref{14.a115} yields moreover:
\begin{equation}
\left.\frac{\partial\clab}{\partial u}\right|_{\Cb_{\tau_0}}, 
\left.\frac{\partial\check{s}_{NL}}{\partial u}\right|_{\Cb_{\tau_0}} = \cO(\tau_0^{N-3})
\label{14.a189}
\end{equation}

Proceeding in the same way, following the approach of Chapter 5, we deduce estimates 
for the higher order derived data which extend \ref{14.a188} and \ref{14.a189} to:
\begin{equation}
\left.T^m\clab\right|_{\Cb_{\tau_0}}, \left.T^m\check{s}_{NL}\right|_{\Cb_{\tau_0}} = 
\cO(\tau_0^{N-3-m}) \ \mbox{: for all $m=0,1,2,...$}
\label{14.a190}
\end{equation}
and:
\begin{equation}
\left.\frac{\partial T^m\clab}{\partial u}\right|_{\Cb_{\tau_0}}, 
\left.\frac{\partial T^m\check{s}_{NL}}{\partial u}\right|_{\Cb_{\tau_0}} = \cO(\tau_0^{N-3-m}) 
\ \mbox{: for all $m=0,1,2...$}
\label{14.a191}
\end{equation}
Moreover, in determining the boundary condition for $T^m s_{NL}$ on $S_{\tau_0,\tau_0}$ 
we derive at the same time estimates for the $m+1$th order $T$ - derivatives of the transformation 
functions, which extend \ref{14.a159}, \ref{14.a161} to: 
\begin{equation}
\left.\left(T^{m+1}\chf, \ T^{m+1}\cv, \ T^{m+1}\cga\right)\right|_{S_{\tau_0,\tau_0}}=\ocO(\tau_0^{N-5-m}) 
 \ \mbox{: for all $m=0,1,2,...$}
\label{14.a192}
\end{equation}

Let now $\alpha_\mu$, $\alphab_\mu$, $\salpha_\mu$ be the rectangular components of the co-frame field 
dual to the frame field with rectangular components $N^\mu$, $\Nb^\mu$, $E^\mu$:
\begin{eqnarray}
&&\alpha_\mu N^\mu=1, \ \ \alpha_\mu\Nb^\mu=0, \ \ \alpha_\mu E^\mu=0\nonumber\\
&&\alphab_\mu N^\mu=0, \ \ \alphab_\mu\Nb_\mu=1, \ \ \alphab_\mu E^\mu=0\nonumber\\
&&\salpha_\mu N^\mu=0, \ \ \salpha_\mu\Nb^\mu=0, \ \ \salpha_\mu E^\mu=1 \label{14.a193}
\end{eqnarray}
We then have:
\begin{equation}
\alpha_\mu=-\frac{1}{2c}h_{\mu\nu}\Nb^\nu, \ \ \alphab_\mu=-\frac{1}{2c}h_{\mu\nu}N^\nu, \ \ 
\salpha_\mu=h_{\mu\nu}E^\nu
\label{14.a194}
\end{equation}
Then any 1-form with rectangular components $\theta_\mu$ can be expanded as:
\begin{equation}
\theta_\mu=\alpha_\mu\theta_N+\alphab_\mu\theta_{\Nb}+\salpha_\mu\stheta
\label{14.a195}
\end{equation}
where we denote:
\begin{equation}
\theta_N=\theta_\mu N^\mu, \ \ \theta_{\Nb}=\theta_\mu\Nb^\mu, \ \ \stheta=\theta_\mu E^\mu 
\label{14.a196}
\end{equation}

Taking then $\theta_\mu=L\beta_\mu$, we have: 
\begin{equation}
\theta_N=s_{NL}, \ \ \theta_{\Nb}=\lambdab\sss, \ \ \stheta=\rho\ss_N
\label{14.a197}
\end{equation}
hence the estimates \ref{14.a188} imply:
\begin{equation}
\left.\left(L\beta_\mu-L_N\beta_{\mu,N}\right)\right|_{\Cb_{\tau_0}}=\cO(\tau_0^{N-3})
\label{14.a198}
\end{equation}
Moreover, by \ref{14.a189} also:
\begin{equation}
\left.\frac{\partial}{\partial u}\left(L\beta_\mu-L_N\beta_{\mu,N}\right)\right|_{\Cb_{\tau_0}}
=\cO(\tau_0^{N-3})
\label{14.a199}
\end{equation}
Together with the 2nd of \ref{14.a92} these give:
\begin{equation}
\left.(T\beta_\mu-T\beta_{\mu,N})\right|_{\Cb_{\tau_0}}=\cO(\tau_0^{N-3})
\label{14.a200}
\end{equation}
\begin{equation}
\left.\frac{\partial}{\partial u}(T\beta_\mu-T\beta_{\mu,N})\right|_{\Cb_{\tau_0}}=\cO(\tau_0^{N-3})
\label{14.a201}
\end{equation}

Recalling from \ref{3.a22} that $\zeta$ is given by:
\begin{equation}
\zeta=2(\rho\teta-\rhob\tetab)+\frac{1}{2}\sbeta(\beta_{\Nb}\rhob LH-\beta_N\rho\Lb H)
\label{14.a202}
\end{equation}
with $\teta$, $\tetab$ given by \ref{3.a25}, and recalling also from Section 10.3 the 
corresponding formulas for the $N$th approximants, taking into account that: 
\begin{equation}
\left.(\teta-\teta_N)\right|_{\Cb_{\tau_0}}=\cO(u^{-1}\tau_0^{N-1}), \ \ 
\left.(\tetab-\tetab_N)\right|_{\Cb_{\tau_0}}=\cO((\tau_0^{N-3})
\label{14.a203}
\end{equation}
the last by \ref{14.a188} and \ref{14.a198}, we conclude that:
\begin{equation}
\left.(\zeta-\zeta_N)\right|_{\Cb_{\tau_0}}=\cO(\tau_0^{N-3})
\label{14.a204}
\end{equation} 
Also, from the propagation equation for $\lambda$ of Proposition 3.3 and the analogous 
propagation equation for $\lambda_N$ from Section 10.3, together with the estimates \ref{14.a99}, 
\ref{14.a188}, \ref{14.a198} we deduce:
\begin{equation}
\left.T\cla\right|_{\Cb_{\tau_0}}=\cO(\tau_0^{N-3})
\label{14.a205}
\end{equation}

Proceeding by induction using the estimates \ref{14.a190}, \ref{14.a191} we deduce:
\begin{eqnarray}
&&\left.\left(T^m\beta_\mu-T^m\beta_{\mu,N}\right)\right|_{\Cb_{\tau_0}}, 
\left.\frac{\partial}{\partial u}\left(T^m\beta_\mu-T^m\beta_{\mu,N}\right)\right|_{\Cb_{\tau_0}} =\cO(\tau_0^{N-2-m}) \nonumber\\
&&\hspace{22mm}\mbox{ : for $m=1,2,...$}
\label{14.a206}
\end{eqnarray}
and:
\begin{equation}
\left.(\zeta-\zeta_N)\right|_{\Cb_{\tau_0}}=\cO(\tau_0^{N-3-m}) \ \ \mbox{: for all $m=0,1,2,...$}
\label{14.a207}
\end{equation}
and:
\begin{equation}
\left.T^m\cla\right|_{\Cb_{\tau_0}}=\cO(\tau_0^{N-2-m}) \ \ \mbox{: for $m=1,2,...$}
\label{14.a208}
\end{equation}

The initial data being now given on $\Cb_{\tau_0}$ we denote by ${\cal N}_{\tau_0}$ the 
part of ${\cal N}$ corresponding to $\ub\geq \tau_0$:
\begin{equation}
{\cal N}_{\tau_0}=\{(\ub,u,\vartheta) \ : \ \ub\geq\tau_0, \ u\geq\ub, \ \vartheta\in S^1\}
\label{14.a209}
\end{equation}
Similarly we denote by ${\cal K}_{\tau_0}$ the part of ${\cal K}$ corresponding to $\tau\geq\tau_0$:
\begin{equation}
{\cal K}_{\tau_0}=\{(\tau,\tau,\vartheta) \ : \ \tau\geq\tau_0, \ \vartheta\in S^1\}
\label{14.a210}
\end{equation}
Recalling the definitions \ref{6.80}, \ref{6.83}, given $u_1\geq\ub_1\geq\tau_0$ we now denote by 
$R_{\ub_1,u_1,\tau_0}$ the trapezoid: 
\begin{equation}
R_{\ub_1,u_1,\tau_0}=\{(\ub,u) \ : \ u\in[\ub,u_1], \ \ub\in[\tau_0,\ub_1]\}
\label{14.a211}
\end{equation}
and by ${\cal R}_{\ub_1,u_1,\tau_0}$ the region in ${\cal N}_{\tau_0}$ corresponding to 
${\cal R}_{\ub_1,u_1,\tau_0}$:
\begin{equation}
{\cal R}_{\ub_1,u_1,\tau_0}=\bigcup_{(\ub,u)\in R_{\ub_1,u_1,\tau_0}}S_{\ub,u}
\label{14.a212}
\end{equation}
Also, given $\tau_1\geq\tau_0$ we denote:
\begin{equation}
{\cal K}_{\tau_0}^{\tau_1}=\{(\tau,\tau,\vartheta) \ : \ \tau\in[\tau_0,\tau_1], \ \vartheta\in S^1\}
\label{14.a213}
\end{equation}
In the following we suppose that we have a smooth solution of the problem defined on 
${\cal R}_{\ub_0,u_0,\tau_0}$ for some $u_0\leq\delta$, $u_0\geq\ub_0>\tau_0$ and satisfying there 
the bootstrap assumptions of Section 14.2 with ${\cal R}_{\ub_0,u_0,\tau_0}$ in the role of 
${\cal R}_{\delta,\delta}$. We shall presently review Chapters 10, 12, 13, as well as the results of 
the preceding sections of the present chapter, marking the changes caused by giving the initial 
data on $\Cb_{\tau_0}$ rather than $\Cb_0$. 

We begin with the {\bf changes in connection with the propagation equations for $\s^{(n-1)}\ctchi$ and 
$\s^{(m,n-m)}\cla$ : $m=0,...,n$} of Chapter 12. 

In inequality \ref{12.128} the integral with respect 
to $\ub$ is now on interval $[\tau_0,\ub_1]$ and there is an additional term
\begin{equation}
k\|\s^{(n-1)}\ctchi\|_{L^2(S_{\tau_0,u})}
\label{14.a214}
\end{equation}
on the right. Taking the $L^2$ norm with respect to $u$ on $[\ub_1,u_1]$ to obtain the estimate for 
$\|\s^{(n-1)}\ctchi\|_{L^2(\Cb_{\ub_1}^{u_1})}$, \ref{14.a214} contributes $k$ times:
\begin{eqnarray}
&&\left\{\int_{\ub_1}^{u_1}\|\s^{(n-1)}\ctchi\|^2_{L^2(S_{\tau_0,u})}du\right\}^{1/2}
\leq\|\s^{(n-1)}\ctchi\|_{L^2(\Cb_{\tau_0}^{u_1})} \nonumber\\
&&\hspace{48mm}\leq C\tau_0^{N-1}
\label{14.a215}
\end{eqnarray}
the last by virtue of \ref{14.a101}. Similarly in obtaining the estimate for 
$\|\s^{(n-1)}\ctchi\|_{L^2({\cal K}_{\tau_0}^{u_1})}$ (see Lemma 12.2). 

In inequality \ref{12.196} the integral with respect 
to $\ub$ is now on interval $[\tau_0,\ub_1]$ and there is an additional term
\begin{equation}
k\|\s^{(0,n)}\cla\|_{L^2(S_{\tau_0,u})}
\label{14.a216}
\end{equation}
on the right. Taking the $L^2$ norm with respect to $u$ on $[\ub_1,u_1]$ to obtain the estimate for 
$\|\s^{(0,n)}\cla\|_{L^2(\Cb_{\ub_1}^{u_1})}$, \ref{14.a216} contributes $k$ times:
\begin{eqnarray}
&&\left\{\int_{\ub_1}^{u_1}\|\s^{(0,n)}\cla\|^2_{L^2(S_{\tau_0,u})}du\right\}^{1/2}
\leq\|\s^{(0,n)}\cla\|_{L^2(\Cb_{\tau_0}^{u_1})} \nonumber\\
&&\hspace{48mm}\leq C\tau_0^{N-\frac{3}{2}}
\label{14.a217}
\end{eqnarray}
the last by virtue of \ref{14.a99}. Similarly in obtaining the estimate for 
$\|\s^{(0,n)}\cla\|_{L^2({\cal K}_{\tau_0}^{u_1})}$ (see Lemma 12.3). 

In inequality \ref{12.564} the integral with respect 
to $\ub$ is now on interval $[\tau_0,\ub_1]$ and there is an additional term
\begin{equation}
k\|\s^{(m,n-m)}\cla\|_{L^2(S_{\tau_0,u})}
\label{14.a218}
\end{equation}
on the right. Here $m=1,...,n$. Taking the $L^2$ norm with respect to $u$ on $[\ub_1,u_1]$ to obtain the estimate for 
$\|\s^{(m,n-m)}\cla\|_{L^2(\Cb_{\ub_1}^{u_1})}$, \ref{14.a218} contributes $k$ times:
\begin{eqnarray}
&&\left\{\int_{\ub_1}^{u_1}\|\s^{(m,n-m)}\cla\|^2_{L^2(S_{\tau_0,u})}du\right\}^{1/2}
\leq\|\s^{(m,n-m)}\cla\|_{L^2(\Cb_{\tau_0}^{u_1})} \nonumber\\
&&\hspace{48mm}\leq Cu_1\tau_0^{N-2-m}
\label{14.a219}
\end{eqnarray}
the last by virtue of \ref{14.a208}. Similarly in obtaining the estimate for 
$\|\s^{(m,n-m)}\cla\|_{L^2({\cal K}_{\tau_0}^{u_1})}$ for $m=1,...,n$ (see Lemma 12.9). 

We turn to the {\bf changes in connection with the propagation equations for $\cth_n$ and 
$\cnu_{m-1,n-m+1}$ : $m=1,...,n$} of Chapter 10. 

In inequality \ref{10.361} with $l=n$ the integral with respect 
to $\ub$ is now on interval $[\tau_0,\ub_1$ and there is an additional term
\begin{equation}
k\|\cth_n\|_{L^2(S_{\tau_0,u})}
\label{14.a220}
\end{equation}
on the right. Taking the $L^2$ norm with respect to $u$ on $[\ub_1,u_1]$ to obtain the estimate for 
$\|\cth_n\|_{L^2(\Cb_{\ub_1}^{u_1})}$, \ref{14.a220} contributes $k$ times:
\begin{eqnarray}
&&\left\{\int_{\ub_1}^{u_1}\|\cth_n\|^2_{L^2(S_{\tau_0,u})}du\right\}^{1/2}
\leq\|\cth_n\|_{L^2(\Cb_{\tau_0}^{u_1})} \nonumber\\
&&\hspace{48mm}\leq C\tau_0^{N-\frac{3}{2}}
\label{14.a221}
\end{eqnarray}
the last by virtue of \ref{14.a92}, \ref{14.a95}, \ref{14.a99}, \ref{14.a101}. Similarly in obtaining the estimate for 
$\|\cth_n\|_{L^2({\cal K}_{\tau_0}^{u_1})}$ (see Propositions 10.3 and 10.4). 

In inequality \ref{10.480} with $l=n-m$ the integral with respect 
to $\ub$ is now on interval $[\tau_0,\ub_1]$ and there is an additional term
\begin{equation}
k\|\cnu_{m-1,n-m+1}-\check{\tau}_{m-1,n-m+1}\|_{L^2(S_{\tau_0,u})}
\label{14.a222}
\end{equation}
on the right. Here $m=1,...,n$. Taking the $L^2$ norm with respect to $u$ on $[\ub_1,u_1]$ to obtain the estimate for 
$\|\cnu_{m-1,n-m+1}-\check{\tau}_{m-1,n-m+1}\|_{L^2(\Cb_{\ub_1}^{u_1})}$, \ref{14.a222} contributes $k$ times:
\begin{eqnarray}
&&\left\{\int_{\ub_1}^{u_1}\|\cnu_{m-1,n-m+1}-\check{\tau}_{m-1,n-m+1}\|^2_{L^2(S_{\tau_0,u})}du\right\}^{1/2}
\label{14.a223}\\
&&\hspace{10mm}\leq\|\cnu_{m-1,n-m+1}-\check{\tau}_{m-1,n-m+1}\|_{L^2(\Cb_{\tau_0}^{u_1})}\leq C\tau_0^{N-\frac{3}{2}-m}\nonumber
\end{eqnarray}
the last by virtue of \ref{14.a99}, \ref{14.a191}, \ref{14.a208}. Similarly in obtaining the estimate for 
$\|\cnu_{m-1,n-m+1}\|_{L^2({\cal K}_{\tau_0}^{u_1})}$ (see Propositions 10.5 and 10.6). 

Next we consider {\bf the change to the top order energy identities \ref{9.285} with $(m,l)=(m,n-m)$ 
due to the fact that the initial data is now given on $\Cb_{\tau_0}$ rather than on $\Cb_0$}. 
This change is simply that the term $\s^{(V;m,n-m)}\cEb^{u_1}(0)$, which vanishes in accordance with 
\ref{9.295}, is replaced by (see \ref{9.291}):
\begin{equation}
\s^{(V;m,n-m)}\cEb^{u_1}(\tau_0)=\int_{\Cb_{\tau_0}}\Omega^{1/2}(3a(\s^{(V;m,n-m)}\csxi)^2
+(\s^{(V;m,n-m)}\cxi_{\Lb})^2)
\label{14.a224}
\end{equation}
Recalling that:
\begin{eqnarray}
&&\s^{(V;m,n-m)}\csxi=V^\mu E(E^{n-m}T^m\beta_\mu-E_N^{n-m}T^m\beta_{\mu,N}) \nonumber\\
&&\s^{(V;m,n-m)}\cxi_{\Lb}=V^\mu\Lb(E^{n-m}T^m\beta_\mu-E_N^{n-m}T^m\beta_{\mu,N}) \label{14.a225}
\end{eqnarray}
For $m=0$ the 2nd of \ref{14.a92} together with \ref{14.a95} imply:
\begin{equation}
\left.\s^{(V;0,n)}\csxi\right|_{\Cb_{\tau_0}}=\cO(\tau_0^{N-1}), \ \ 
\left.\s^{(V;0,n)}\cxi_{\Lb}\right|_{\Cb_{\tau_0}}=\cO(u^{-1}\tau_0^{N-1})
\label{14.a226}
\end{equation}
hence:
\begin{equation}
\s^{(V;0,n)}\cEb^{u_1}(\tau_0)\leq C\tau_0^{2N-3}
\label{14.a227}
\end{equation}
For $m=1,...,n$ \ref{14.a206} together with \ref{14.a95} imply:
\begin{equation}
\left.\s^{(V;m,n-m)}\csxi\right|_{\Cb_{\tau_0}}, \left.\s^{(V;m,n-m)}\cxi_{\Lb}\right|_{\Cb_{\tau_0}} 
 = \cO(\tau_0^{N-2-m})
\label{14.a228}
\end{equation}
hence:
\begin{equation}
\s^{(V;m,n-m)}\cEb^{u_1}(\tau_0)\leq C u_1\tau_0^{2N-4-2m}
\label{14.a229}
\end{equation}

In addition to the above change, the integral in the definition \ref{9.287} of 
$\s^{(V;m,n-m)}\cE^{\ub_1}(u)$ and in \ref{9.290} is now on $C_{u,\tau_0}^{\ub_1}$, the part of $C_u^{\ub_1}$ lying in 
${\cal N}_{\tau_0}$, that is corresponding to $\ub\geq\tau_0$. Also, the integral 
in the definition \ref{9.288} of $\s^{(V;m,n-m)}\cF^{\ub_1}$, in \ref{9.292}, and the integrals 
in \ref{9.a34}, \ref{9.a35}, \ref{9.a39} are now on ${\cal K}_{\tau_0}^{\ub_1}$ and 
${\cal K}_{\tau_0}^{\tau_1}$ respectively. Also, the error integral \ref{9.289} is now on 
${\cal R}_{\ub_1,u_1,\tau_0}$. Moreover, the supremum in the definitions \ref{9.297}, \ref{9.298} 
of $\s^{(V;m,n-m)}\cB(\ub_1,u_1)$, $\s^{(V;m,n-m)}\cBb(\ub_1,u_1)$ is now over $R_{\ub_1,u_1,\tau_0}$, 
and the supremum in the definition \ref{9.299} of $\s^{(V;m,n-m)}\cA(\tau_1)$ is now over 
$\tau\in[\tau_0,\tau_1]$. These changes affect the succeeding chapters. 
In particular, in Chapter 10, the $L^2$ norms 
on ${\cal K}^\tau$ in Propositions 10.4, 10.6 are now $L^2$ norms on ${\cal K}_{\tau_0}^\tau$, 
and the $L^2$ norms on $C_{u_1}^{\ub_1}$ and on ${\cal K}^{\ub_1}$ in Propositions 10.7 and 10.8 
are now $L^2$ norms on $C_{u_1,\tau_0}^{\ub_1}$ and on ${\cal K}_{\tau_0}^{\ub_1}$ respectively. 
Also, the $L^2$ norms on ${\cal K}^{\tau_1}$ in Proposition 10.9 are $L^2$ norms on 
${\cal K}_{\tau_0}^{\tau_1}$. In Chapter 12, the $L^2$ norms on ${\cal K}^{\tau}$ in Lemmas 12.2, 
12.3, 12.7, 12.8, 12.9 are now $L^2$ norms on ${\cal K}_{\tau_0}^{\tau}$, and the $L^2$ norms on 
$C_{u_1,\tau_0}^{\ub_1}$ and on ${\cal K}^{\ub_1}$ in Lemmas 12.4, 12.6, 12.10 are now $L^2$ norms 
on $C_{u_1,\tau_0}^{\ub_1}$ and on ${\cal K}_{\tau_0}^{\ub_1}$. Moreover, the supremum in the 
definitions \ref{12.394} - \ref{12.396}, \ref{12.397} - \ref{12.400}, \ref{12.532}, \ref{12.795} - 
\ref{12.796}, \ref{12.805} - \ref{12.807}, \ref{12.892}, \ref{12.971} is now over 
$\tau\in[\tau_0,\tau_1]$. Also, the supremum in the definitions \ref{12.411} - \ref{12.412}, 
\ref{12.460} - \ref{12.461}, \ref{12.665} - \ref{12.666} is now over $(\ub,u)\in R_{\ub_1,u_1,\tau_0}$.  
In addition, the suprema in Lemma 12.12 and in Proposition 12.7 are now over $\tau\in[\tau_0,\tau_1]$. 
In Chapter 13 the integral \ref{13.2} is now on ${\cal K}_{\tau_0}^{\ub_1}$. Also, in 
the energy identity \ref{13.3} the initial data term $\s^{(V;m,n-m)}\cEb^{u_1}(\tau_0)$, estimated 
above, now appears. The remaining changes are obvious. 

However, {\bf a number of results rely on the use of the inequality \ref{12.366}, which changes 
as a consequence of the initial data being now given on $\Cb_{\tau_0}$.} We must presently assess 
these changes. Since in \ref{12.365} the integration is to be on $[\tau_0,\tau]$, and the function 
$f$ does not vanish on $S_{\tau_0,\tau_0}$, \ref{12.365} is replaced by:
\begin{equation}
f(\tau,\vartheta)=f(\tau_0,\vartheta)+\int_{\tau_0}^{\tau}(Tf)(\tau^\prime,\vartheta)d\tau^\prime
\label{14.a230}
\end{equation}
hence the inequality \ref{12.366} is replaced by:
\begin{eqnarray}
&&\|f\|_{L^2(S_{\tau,\tau})}\leq k\|f\|_{L^2(S_{\tau_0,\tau_0})}+k\int_{\tau_0}^{\tau}
\|Tf\|_{L^2(S_{\tau^\prime,\tau^\prime})}d\tau^\prime\nonumber\\
&&\hspace{17mm}\leq k\|f\|_{L^2(S_{\tau_0,\tau_0})}+k\tau^{1/2}\|Tf\|_{L^2({\cal K}_{\tau_0}^{\tau})}
\label{14.a231}
\end{eqnarray}

As a consequence, the inequality \ref{12.487} changes to:
\begin{eqnarray}
&&\|\Omega^n\check{f}\|_{L^2(S_{\tau,\tau})}\leq k\|\Omega\check{f}\|_{L^2(S_{\tau_0,\tau_0})}
+k\tau^{1/2}\|\Omega^n Tf\|_{L^2({\cal K}_{\tau_0}^{\tau})}\nonumber\\
&&\hspace{17mm}\leq k\|\Omega^n\check{f}\|_{L^2(S_{\tau_0,\tau_0})}
+k\tau^{c_0+\frac{1}{2}}\sup_{\tau\in[\tau_0,\tau_1]}
\left\{\tau^{-c_0}\|\Omega^n T\check{f}\|_{L^2({\cal K}_{\tau_0}^{\tau})}\right\}\nonumber\\
&&\label{14.a232}
\end{eqnarray}
Then in \ref{12.488}, ${\cal K}^\tau$ being replaced by ${\cal K}_{\tau_0}^\tau$ and the 
supremum on the right being now on $[\tau_0,\tau_1]$, an additional term of the form 
\begin{equation}
C\tau_0^{-5/2}\|\Omega^n\check{f}\|_{L^2(S_{\tau_0,\tau_0})}
\label{14.a233}
\end{equation}
appears on the right hand side, which by the 1st of \ref{14.a92} is bounded by:
\begin{equation}
C\tau_0^{N-\frac{7}{2}}
\label{14.a234}
\end{equation}
The effect of this change is that to the right hand side of \ref{12.490}  a term of the form 
\ref{14.a234} is added. Similarly, the inequality \ref{12.527} changes by the addition of the term 
\begin{equation}
k\|\Omega^n\check{\delta}^i\|_{L^2(S_{\tau_0,\tau_0})}
\label{14.a235}
\end{equation}
apart from the obvious change of the replacement of ${\cal K}^\tau$ by ${\cal K}_{\tau_0}^\tau$ and 
of the supremum on $[0,\tau_1]$ by the supremum on $[\tau_0,\tau_1]$. An additional term of the 
form 
\begin{equation}
C\tau_0^{-7/2}\|\Omega^n\check{\delta}^i\|_{L^2(S_{\tau_0,\tau_0})}
\label{14.a236}
\end{equation}
then appears on the right hand side of \ref{12.528}, which by the 1st of \ref{14.a92} is bounded by 
(see also \ref{14.a90}):
\begin{equation}
C\tau_0^{N-\frac{9}{2}}
\label{14.a237}
\end{equation}
The effect of this change is that to the right hand side of \ref{12.530}  a term of the form 
\ref{14.a237} is added. 

Next, inequality \ref{12.784} changes in accordance with inequality \ref{14.a231} to:
\begin{equation}
\|f\|_{L^2({\cal K}_{\tau_0}^{\tau})}\leq k\tau^{1/2}\|f\|_{L^2(S_{\tau_0,\tau_0})}+
k\left(\int_{\tau_0}^{\tau}\tau^\prime\|Tf\|^2_{L^2({\cal K}_{\tau_0}^{\tau^\prime})}d\tau^\prime
\right)^{1/2}
\label{14.a238}
\end{equation}
As a consequence, inequality \ref{12.785} changes by the addition to the right hand side of the terms 
\begin{equation}
k\|f\|_{L^2(S_{\tau_0,\tau_0})}+k^2\tau\|Tf\|_{L^2(S_{\tau_0,\tau_0})}
\label{14.a239}
\end{equation}
inequalities \ref{12.787} and \ref{12.789} change by the addition to the right hand side 
of terms of the form 
\begin{equation}
\sum_{i=0}^{j}k^{i+1}\tau^i\|T^i f\|_{L^2(S_{\tau_0,\tau_0})}
\label{14.a240}
\end{equation}
with $j=2$ and $j=3$ respectively. Therefore terms of the form 
\begin{eqnarray}
&&C\left\{\sum_{i=0}^1\tau^i\|\Omega^{n-m}T^{m+i-1}\cv\|_{L^2(S_{\tau_0,\tau_0})}
\right\}\leq C\tau\tau_0^{N-4-m}\label{14.a241}\\
&&C\left\{\sum_{i-0}^2\tau^i\|\Omega^{n-m}T^{m+i-2}\chf\|_{L^2(S_{\tau_0,\tau_0})}
\right\}\leq C\tau^2\tau_0^{N-4-m}\label{14.a242}\\
&&C\left\{\sum_{i-0}^3\tau^i\|\Omega^{n-m}T^{m+i-3}\cga\|_{L^2(S_{\tau_0,\tau_0})}
\right\}\leq C\tau^3\tau_0^{N-4-m}\label{14.a243}
\end{eqnarray}
are to be added to the right hand sides of \ref{12.786}, \ref{12.788}, \ref{12.790} respectively, 
the bounds holding by virtue of the estimates \ref{14.a192}. This affects the bounds \ref{12.809}, 
\ref{12.808}, \ref{12.810}, by the addition of the terms: 
\begin{eqnarray}
&&C\tau\tau_0^{N-4-m}\leq C\tau^{c_m-\frac{1}{2}}\cdot\tau_0^{N-\frac{5}{2}-c_m-m} \label{14.a244}\\
&&C\tau^2\tau_0^{N-4-m}\leq C\tau^{c_m+\frac{1}{2}}\cdot\tau_0^{N-\frac{5}{2}-c_m-m} \label{14.a245}\\
&&C\tau^3\tau_0^{N-4-m}\leq C\tau^{c_m+\frac{3}{2}}\cdot\tau_0^{N-\frac{5}{2}-c_m-m} \label{14.a246}
\end{eqnarray}
respectively. Note that he extra factor $\tau_0^{N-\frac{5}{2}-c_m-m}$ tends to 0 as 
$\tau_0\rightarrow 0$, if we choose:
\begin{equation}
N>c_m+m+\frac{5}{2}
\label{14.a247}
\end{equation}
a condition which implies \ref{13.242}. 

Next, as a consequence of the replacement of inequality \ref{12.366} by inequality \ref{14.a231}, 
in \ref{12.847} the term 
\begin{equation}
k\|\Omega^{n-m}T^m\check{f}\|_{L^2(S_{\tau_0,\tau_0})}
\label{14.a248}
\end{equation}
must be added to the right hand side. It then follows by \ref{14.a192} that an additional term of the form
\begin{equation}
C\tau_0^{N-\frac{9}{2}-m}
\label{14.a249}
\end{equation}
appears on the right hand side of \ref{12.848}. Similarly, in \ref{12.849} a sum of the form 
\begin{equation}
\sum_{j=1}^i C_{i,j}\tau^{j-1}\|\Omega^{n-m}T^{m-i+j}\check{f}\|_{L^2(S_{\tau_0,\tau_0})}
\leq C_i\tau^{i-1}\tau_0^{N-2-m}
\label{14.a250}
\end{equation}
must be added to the right hand side, the bound by virtue of \ref{14.a192}. The effect of this in 
\ref{12.854} is the appearance on the right hand side of an additional term 
\begin{equation}
C\tau_0^{N-\frac{9}{2}-m}\leq C\tau_1^{c_m-2}\cdot\tau_0^{N-\frac{5}{2}-c_m-m}
\label{14.a251}
\end{equation}
the extra factor tending to 0 as $\tau_0\rightarrow 0$ by the choice \ref{14.a247}. The argument 
in regard to $\chdl^i$ from the formula \ref{12.878} to the estimate \ref{12.889} parallels 
the argument just discussed in regard to $\chf$ from the formula \ref{12.844} to the estimate \ref{12.855}. 

The replacement of \ref{12.365}, \ref{12.366} by \ref{14.a230}, \ref{14.a231} causes an additional 
term of the form
\begin{equation}
\|\s^{(n)}\check{\sh}\|_{L^2(S_{\tau_0,\tau_0})}\leq C\tau_0^{N-1}
\label{14.a252}
\end{equation}
to appear on the right hand side of \ref{12.982}, the bound by virtue of \ref{14.a95}. The 
contribution of this term to $\s^{(n)}H(\tau_1)$, defined by \ref{12.971}, is bounded by:
\begin{equation}
C\tau_0^{N+\frac{1}{2}-c_0}
\label{14.a253}
\end{equation}

We finally consider the {\bf changes in the results of the preceding sections of the present chapter 
caused by the initial data being given on $\Cb_{\tau_0}$}. There are no changes in Sections 14.1 and 
14.6. The changes in regard to Sections 14.2 and 14.3 are similar to those discussed above 
in connection with the propagation equations for $\s^{(n-1)}\ctchi$ and 
$\s^{(m,n-m)}\cla$ : $m=0,...,n$ of Chapter 12. Therefore we focus on the changes in Sections 14.4 - 
14.5 and 14.7 - 14.9. 

In regard to Section 14.4, the integral in \ref{14.141} is now on $[\tau_0,\ub_1]$ and the initial 
data term $f(\tau_0,u,\vartheta)$ is added on the right (the coordinates while being adapted to the 
flow of $L$ in ${\cal N}_{\tau_0}$ they are also adapted to the flow of $\Lb$ on $\Cb_{\tau_0}$. 
As a consequence, inequality \ref{14.145} is replaced by:
\begin{equation}
\|f\|_{L^2(S_{\ub_1,u})}\leq k\|f\|_{L^2(S_{\tau_0,u_1})}
+k\int_{\tau_0}^{\ub_1}\|Lf\|_{L^2(S_{\ub,u})}d\ub
\label{14.a254}
\end{equation}
Then in \ref{14.152} the term 
\begin{equation}
k\|\s^{(m,n-m)}\check{\dot{\phi}}_\mu\|_{L^2(S_{\tau_0,u})}\leq C\left\{\begin{array}{lll}
\tau_0^{N-1}  &:& \mbox{for $m=0$}\\
\tau_0^{N-2-m} &:& \mbox{for $m=1,...,n$}\end{array}\right.
\label{14.a255}
\end{equation}
must be added to the right hand side, the bound by virtue of \ref{14.a206} and the 2nd of \ref{14.a92}.  
As a result, the right hand side of the inequality of Proposition 14.3 changes to:
\begin{equation}
\left(C\sqrt{D}+C\left\{\begin{array}{lll}\tau_0^{N-1-c_0} &:& \mbox{for $m=0$}\\
\tau_0^{N-2-m-c_m} &:& \mbox{for $m=1,...,n$}\end{array}\right.\right)
\ub^{a_m+\frac{1}{2}}u^{b_m-\frac{1}{2}}
\label{14.a256}
\end{equation}
the 2nd term in parenthesis tending to 0 as $\tau_0\rightarrow 0$ by the choice \ref{14.a247}. 

Similarly, in \ref{14.191} the term 
\begin{equation}
k\|q_{1,m,\mu}\|_{L^2(S_{\tau_0,u})}
\label{14.a257}
\end{equation}
must be added to the right hand side. From the definition \ref{14.188} we have:
\begin{eqnarray}
&&q_{1,m,\mu}=\Lb E^{n-m}T^{m-1}\beta_\mu-\Lb_N E_N^{n-m}T^{m-1}\beta_{\mu,N}\nonumber\\
&&\hspace{12mm}-[\Lb,E^{n-m}T^{m-1}]\beta_\mu+[\Lb_N,E_N^{n-m}T^{m-1}]\beta_{\mu,N}
\label{14.a258}
\end{eqnarray}
Evaluating this on $\Cb_{\tau_0}$ the first difference is 
\begin{equation}
\left\{\begin{array}{lll}\cO(u^{-1}\tau_0^{N-1}) &:& \mbox{for $m=1$}\\
\cO(\tau_0^{N-1-m}) &:& \mbox{for $m=2,...,n$}\end{array}\right.
\label{14.a259}
\end{equation}
As for the commutator difference, using \ref{14.a204} and the fact that by \ref{14.a101}:
\begin{equation}
\left.(\chib-\chib_N)\right|_{\Cb_{\tau_0}}=\cO(u^{-1}\tau_0^{N-1})
\label{14.a260}
\end{equation}
a similar estimate is obtained. We then conclude that \ref{14.a257} is bounded by:
\begin{equation}
C\left\{\begin{array}{lll}u^{-1}\tau_0^{N-1} &:& \mbox{for $m=1$}\\
\tau_0^{N-1-m} &:& \mbox{for $m=2,...,n$}\end{array}\right.
\label{14.a261}
\end{equation}
As a result, the right hand side of the inequality of Lemma 14.1 changes to:
\begin{equation}
\left(C\sqrt{D}+C\tau_0^{N-2-c_{m-1}-m}\right)
\ub^{a_{m-1}+\frac{1}{2}}u^{b_{m-1}+\frac{1}{2}} \ \mbox{: for $m=1,...,n$}
\label{14.a262}
\end{equation}
the 2nd term in parenthesis tending to 0 as $\tau_0\rightarrow 0$ by the choice \ref{14.a247}. 

In regard to Section 14.5, the integral in \ref{14.239} is now on $[\tau_0,\ub_1]$ and the initial 
data term $f(\tau_0,\sigma,\vartheta)$ is added on the right. 
As a consequence, inequality \ref{14.240} is replaced by:
\begin{eqnarray}
&&\|f\|_{L^2(S_{\tau,\sigma+\tau})}\leq k\|f\|_{L^2(S_{\tau_0,\sigma+\tau_0})}
+k\int_{\tau_0}^{\tau}\|Tf\|_{L^2(S_{\tau^\prime,\sigma+\tau^\prime})}d\tau^\prime \nonumber\\
&&\hspace{20mm}\leq k\|f\|_{L^2(S_{\tau_0,\sigma+\tau_0})}
+k\tau^{1/2}\|Tf\|_{L^2({\cal K}_{\sigma,\tau_0}^{\tau})}
\label{14.a265}
\end{eqnarray}
where we denote by ${\cal K}_{\sigma,\tau_0}^\tau$ the part of ${\cal K}_{\sigma}^\tau$ lying in 
${\cal N}_{\tau_0}$. Then in each of \ref{14.241} - \ref{14.244} there is an additional term on 
the right:
\begin{eqnarray}
&&k\|\s^{(n-2)}\ctchi\|_{L^2(S_{\tau_0,\sigma+\tau_0})}\leq C\tau_0^{N-1} \label{14.a266}\\
&&k\|\s^{(n-2)}\ctchib\|_{L^2(S_{\tau_0,\sigma+\tau_0})}\leq C\tau_0^{N-1} \label{14.a267}\\
&&k\|\s^{(0,n-1}\cla\|_{L^2(S_{\tau_0,\sigma+\tau_0})}\leq C(\sigma+\tau_0)^{-1}\tau_0^{N-1} 
\label{14.a268}\\
&&k\|\s^{(0,n-1)}\clab\|_{L^2(S_{\tau_0,\sigma+\tau_0})}\leq C\tau_0^{N-3} \label{14.a269}
\end{eqnarray}
respectively, the bounds by virtue of \ref{14.a101}, \ref{14.a99}, \ref{14.a188}. Now, the substitutions 
stated just below \ref{14.244}, result in terms due to the initial data on $\Cb_{\tau_0}$ in addition 
to the terms \ref{14.a266} - \ref{14.a269}. The substitutions from Proposition 14.3 and from Lemma 14.1 
result in the additional terms discussed above. There are also initial data contributions 
associated to the $L^2({\cal K}_{\sigma,\tau_0}^{\tau})$ estimates of Proposition 14.2 for 
$\s^{(n-1)}\ctchi$ and for $\s^{(m,n-m)}\cla$, the case $m=0$ pertaining to \ref{14.241} - \ref{14.244}. 
These are similar to the initial data contributions, discussed above, which are associated to 
the propagation equations for $\s^{(n-1)}\ctchi$ and $\s^{(m,n-m)}\cla$ : $m=0,...,n$ of Chapter 12. 
That is, the contribution, for arbitrary $f$, of the initial data term $k\|f\|_{L^2(S_{\tau_0,u})}$ 
to the bound for $\|f\|_{L^2(S_{\ub,u})}$, 
setting $\ub=\tau^\prime$, $u=\sigma+\tau^\prime$ and taking the 
$L^2$ norm with respect to $\tau^\prime$ on $[\tau_0,\tau]$, gives a contribution of:
\begin{equation}
\left\{\int_{\sigma+\tau_0}^{\sigma+\tau}\|f\|^2_{L^2(S_{\tau_0,u})}du\right\}^{1/2}
=\|f\|_{L^2(\Cb_{\tau_0,\tau,\sigma})}
\label{14.a270}
\end{equation}
to $\|f\|_{L^2({\cal K}_{\sigma,\tau_0}^{\tau})}$. Here we denote by $\Cb_{\tau_0,\tau,\sigma}$ the 
part of $\Cb_{\tau_0}$ corresponding to $\sigma+\tau_0\leq u\leq\sigma+\tau$. Thus, in regard to the 
1st term on the right in \ref{14.241}, the initial data contribution to 
$\|\s^{(0,n)}\cla\|_{L^2({\cal K}_{\sigma,\tau_0}^{\tau})}$ is:
\begin{equation}
\|\s^{(0,n)}\cla\|_{L^2(\Cb_{\tau_0,\tau,\sigma})}\leq C(\sigma+\tau_0)^{-1}\tau_0^{N-1}
\label{14.a271}
\end{equation}
by \ref{14.a99}. In regard to the 2nd term on the right in \ref{14.241}, the initial data 
contribution to $\|\s^{(n-1)}\ctchi\|_{L^2({\cal K}_{\sigma,\tau_0}^{\tau})}$ is:
\begin{equation}
\|\s^{(n-1)}\ctchi\|_{L^2(\Cb_{\tau_0,\tau,\sigma})}\leq C\tau\tau_0^{N-1}
\label{14.a272}
\end{equation}
by \ref{14.a101}. Also, in regard to the 1st term on the right in \ref{14.243}, the initial data 
contribution to $\|\s^{(1,n-1)}\cla\|_{L^2({\cal K}_{\sigma,\tau_0}^{\tau})}$ is:
\begin{equation}
\|\s^{(1,n-1)}\cla\|_{L^2(\Cb_{\tau_0,\tau,\sigma})}\leq C\tau^{1/2}\tau_0^{N-3}
\label{14.a273}
\end{equation}
by \ref{14.a205}. The above considerations lead to the conclusion that the initial data 
on $\Cb_{\tau_0}$ contributions to the inequalities \ref{14.245} - \ref{14.248} amount to a replacement 
of the $C\sqrt{D}$ factor in the 1st term in each by:
\begin{equation}
C\sqrt{D}+C\tau_0^{N-\frac{7}{2}-c_1}
\label{14.a274}
\end{equation}
The result \ref{14.264} is then accordingly modified to:
\begin{equation}
x(\tau), \ \xb(\tau), \ y(\tau), \ \yb(\tau) \ \leq \ C\sqrt{D}+C\tau_0^{N-\frac{7}{2}-c_1}
\label{14.a275}
\end{equation}
Note that with the choice \ref{14.a247} the 2nd term on the right tends to 0 as $\tau_0\rightarrow 0$. 

We turn to the inequalities \ref{14.267}. Here $m=2,...,n$. In each of these two inequalities 
there is an additional term on the right:
\begin{eqnarray}
&&k\|\s^{(m-1,n-m)}\cla\|_{L^2(S_{\tau_0,\sigma+\tau_0})}\leq C\tau_0^{N-1-m} \nonumber\\
&&k\|\s^{(m-1,n-m)}\clab\|_{L^2(S_{\tau_0,\sigma+\tau_0})}\leq C\tau_0^{N-2-m} \label{14.a276}
\end{eqnarray}
in the first and second inequality respectively, the bounds by virtue of \ref{14.a208} and 
\ref{14.a190}. Moreover, in regard to the 1st term on the right in the first inequality, there is an 
initial data contribution to $\|\s^{(m,n-m)}\cla\|_{L^2({\cal K}_{\sigma,\tau_0}^{\tau})}$, which is:
\begin{equation}
\|\s^{(m,n-m)}\cla\|_{L^2(\Cb_{\tau_0,\tau,\sigma})}\leq C\tau^{1/2}\tau_0^{N-2-m}
\label{14.a277}
\end{equation}
by \ref{14.a208}. The effect of these changes is to replace the $C\sqrt{D}$ factor in the last two 
estimates of Proposition 14.4 by:
\begin{equation}
C\sqrt{D}+C\tau_0^{N-\frac{5}{2}-m-c_m}
\label{14.a278}
\end{equation}
the first two estimates being modified by the replacement of $C\sqrt{D}$ by \ref{14.a274} according 
to \ref{14.a275}. Note again that with the choice \ref{14.a247} the 2nd term in \ref{14.a278} tends 
to 0 as $\tau_0\rightarrow 0$. 

In regard to the inequalities \ref{14.270}, since these rely on \ref{12.366} which has 
been replaced by \ref{14.a231}, there is an additional term in each:
\begin{eqnarray}
&&k\|\Omega^{n-m}T^m\chf\|_{L^2(S_{\tau_0,\tau_0})}\leq C\tau_0^{N-4-m} \nonumber\\
&&k\|\Omega^{n-m}T^m\cv\|_{L^2(S_{\tau_0,\tau_0})}\leq C\tau_0^{N-4-m} \nonumber\\
&&k\|\Omega^{n-m}T^m\cga\|_{L^2(S_{\tau_0,\tau_0})}\leq C\tau_0^{N-4-m} \label{14.a279}
\end{eqnarray}
the bounds by virtue of \ref{14.a192}. The effect is to replace $C\sqrt{D}$ in \ref{14.270} by:
\begin{equation}
C\sqrt{D}+C\tau_0^{N-2-m-c_m}
\label{14.a280}
\end{equation}

In regard to Section 14.7, as a result of the replacement of inequality \ref{14.240} by inequality 
\ref{14.a265}, the additional term 
\begin{equation}
k\|\s^{(n-1)}\check{\sh}\|_{L^2(S_{\tau_0,\tau_0})}\leq C\tau_0^{N-1}
\label{14.a281}
\end{equation}
now appears on the right hand side of \ref{14.301}, the bound by virtue of \ref{14.a95}.  
In addition, in regard to the 1st term on the right in \ref{14.303} there is the initial data 
contribution \ref{14.a272} to $\|\s^{(n-1)}\ctchi\|_{L^2({\cal K}_{\sigma,\tau_0}^{\tau})}$. Also,  
in \ref{14.303} and \ref{14.304} there are contributions from the additional terms in Proposition 14.3 
and in Lemma 14.1 previously discussed. The effect of these changes is to replace $C\sqrt{D}$ in 
\ref{14.305} and in the corresponding estimate in Proposition 14.6 by:
\begin{equation}
C\sqrt{D}+C\tau_0^{N-3-c_1}
\label{14.a282}
\end{equation}
Because of the vanishing of $\check{b}$ on $\Cb_{\tau_0}$, the replacement 
of inequality \ref{14.240} by inequality \ref{14.a265} does not cause an additional term to 
appear on the right hand side of \ref{14.313}. However there are initial data contributions 
from the first integral on the right in \ref{14.313} (the integrals are now on $[\tau_0,\tau]$) 
corresponding to the initial data contributions to 
$\|\s^{(0,n)}\cla\|_{L^2({\cal K}_{\sigma,\tau_0}^{\tau})}$ and 
$\|\s^{(n-1)}\ctchi\|_{L^2({\cal K}_{\sigma,\tau_0}^{\tau})}$ discussed above (see \ref{14.314}), 
besides the contributions from the additional terms in Proposition 14.3 
and in Lemma 14.1 (see \ref{14.314}). The effect of these changes is to replace $\sqrt{D}$ in 
\ref{14.319} and in the corresponding estimate of Proposition 14.6 by \ref{14.a282}. 

In regard to Section 14.8, we consider first the changes to the inequalities \ref{14.322} and 
\ref{14.323}. The argument preceding these inequalities applies the inequality \ref{14.240} 
and relies on the estimates of Proposition 14.3, Lemma 14.1, and Proposition 14.4. As a result of 
the inequality \ref{14.240} having been replaced by the inequality \ref{14.a265} and of the 
modifications, discussed above, to the estimates of Proposition 14.3, Lemma 14.1, and Proposition 14.4, 
caused by the initial data being given on $\Cb_{\tau_0}$ rather than $\Cb_0$, the inequalities 
\ref{14.322} and \ref{14.323} are modified to:
\begin{eqnarray}
&&\|\s^{(m-1,n-m)}\check{\dot{\phi}}_\mu\|_{L^2(S_{\ub,u})}\leq 
\ub^{a_m+\frac{3}{2}}u^{b_m-\frac{1}{2}}
\left\{\frac{C\sqrt{D}}{\left(a_m+\frac{3}{2}\right)}+C\tau_0^{N-2-m-c_m}\right\} \nonumber\\
&&\hspace{37mm} \ : \ m=1,...n\nonumber\\
&&\|\stackrel{1}{q}_{1,m-1,\mu}\|_{L^2(S_{\ub,u})}\leq 
\ub^{a_m+\frac{3}{2}}u^{b_m+\frac{1}{2}}
\left\{\frac{C\sqrt{D}}{\left(a_m+\frac{3}{2}\right)}+C\tau_0^{N-2-m-c_m}\right\} \nonumber\\
&&\hspace{31mm} \ : \ m=2,...,n\nonumber\\
&&\|\stackrel{1}{\qb}_{1,m-1,\mu}\|_{L^2(S_{\ub,u})}\leq 
\ub^{a_m+\frac{3}{2}}u^{b_m-\frac{1}{2}}
\left\{\frac{C\sqrt{D}}{\left(a_m+\frac{3}{2}\right)}+C\tau_0^{N-2-m-c_m}\right\} \nonumber\\
&&\hspace{31mm} \ : \ m=2,...,n\nonumber\\
&&\label{14.a283}
\end{eqnarray}
and:
\begin{eqnarray}
&&\|\s^{(n-3)}\ctchi\|_{L^2(S_{\ub,u})}\leq 
\ub^{a_2+\frac{3}{2}}u^{b_2+\frac{1}{2}}
\left\{\frac{C\sqrt{D}}{\left(a_2+\frac{3}{2}\right)}+C\tau_0^{N-4-c_2}\right\} \nonumber\\
&&\|\s^{(n-3)}\ctchib\|_{L^2(S_{\ub,u})}\leq 
\ub^{a_2+\frac{3}{2}}u^{b_2}
\left\{\frac{C\sqrt{D}}{\left(a_2+\frac{3}{2}\right)}+C\tau_0^{N-4-c_2}\right\} \nonumber\\
&&\|\s^{(m-2,n-m)}\cla\|_{L^2(S_{\ub,u})}\leq 
\ub^{a_m+2}u^{b_m+\frac{1}{2}}
\left\{\frac{C\sqrt{D}}{\left(a_m+2\right)}+C\tau_0^{N-\frac{5}{2}-m-c_m}\right\} \nonumber\\
&&\hspace{35mm} \ : \ m=2,...,n\nonumber\\
&&\|\s^{(m-2,n-m)}\clab\|_{L^2(S_{\ub,u})}\leq 
\ub^{a_m+\frac{3}{2}}u^{b_m}
\left\{\frac{C\sqrt{D}}{\left(a_m+\frac{3}{2}\right)}+C\tau_0^{N-\frac{5}{2}-m-c_m}\right\} \nonumber\\
&&\hspace{35mm} \ : \ m=2,...,n\nonumber\\
&&\label{14.a284}
\end{eqnarray} 
Note that with the choice \ref{14.a247} in each case the additional term in parenthesis tends to 0 as 
$\tau_0\rightarrow 0$. 

We turn to the inequalities \ref{14.326}. Here and in the remainder the exponents are chosen 
to satisfy \ref{14.325}. We now repeatedly apply the inequality \ref{14.a265} rather than
the inequality \ref{14.240}. This results in the appearance of additional terms which depend on 
the initial data on $\Cb_{\tau_0}$. The modified inequalities take the form:
\begin{eqnarray}
&&\|\s^{(m-i,n-m)}\check{\dot{\phi}}_\mu\|_{L^2(S_{\ub,u})}\leq 
\ub^{a+\frac{1}{2}+i}u^{b-\frac{1}{2}}
\left\{\frac{C\sqrt{D}}{\left(a+\frac{3}{2}\right)\cdot\cdot\cdot\left(a+\frac{1}{2}+i\right)}
+A_i\tau_0^{N-2-m-c}\right\}\nonumber\\ 
&&\hspace{36mm} \ : \ m=i,...n\nonumber\\
&&\|\stackrel{i}{q}_{1,m-i,\mu}\|_{L^2(S_{\ub,u})}\leq 
\ub^{a+\frac{1}{2}+i}u^{b+\frac{1}{2}}
\left\{\frac{C\sqrt{D}}{\left(a+\frac{3}{2}\right)\cdot\cdot\cdot\left(a+\frac{1}{2}+i\right)}
+A_i\tau_0^{N-2-m-c}\right\} \nonumber\\ 
&&\hspace{31mm} \ : \ m=i+1,...,n\nonumber\\
&&\|\stackrel{i}{\qb}_{1,m-i,\mu}\|_{L^2(S_{\ub,u})}\leq 
\ub^{a+\frac{1}{2}+i}u^{b-\frac{1}{2}}
\left\{\frac{C\sqrt{D}}{\left(a+\frac{3}{2}\right)\cdot\cdot\cdot\left(a+\frac{1}{2}+i\right)}
+A_i\tau_0^{N-2-m-c}\right\} \nonumber\\ 
&&\hspace{31mm} \ : \ m=i+1,...,n\nonumber\\
&&\label{14.a285}
\end{eqnarray}
Here we denote by $A_i$ positive constants depending on $i$ as well as on the exponent $a$ which are 
non-increasing in $a$.  Similarly, the inequalities \ref{14.327} are modified to:
\begin{eqnarray}
&&\|\s^{(n-2-i)}\ctchi\|_{L^2(S_{\ub,u})}\leq \ub^{a+\frac{1}{2}+i}u^{b+\frac{1}{2}}
\left\{\frac{C\sqrt{D}}{\left(a+\frac{3}{2}\right)\cdot\cdot\cdot\left(a+\frac{1}{2}+i\right)}
+A_i\tau_0^{N-3-i-c}\right\}
\nonumber\\
&&\|\s^{(n-2-i)}\ctchib\|_{L^2(S_{\ub,u})}\leq \ub^{a+\frac{1}{2}+i}u^{b}
\left\{\frac{C\sqrt{D}}{\left(a+\frac{3}{2}\right)\cdot\cdot\cdot\left(a+\frac{1}{2}+i\right)}
+A_i\tau_0^{N-3-i-c}\right\}\nonumber\\
&&\|\s^{(m-1-i,n-m)}\cla\|_{L^2(S_{\ub,u})}\leq \ub^{a+1+i}u^{b+\frac{1}{2}}
\left\{\frac{C\sqrt{D}}{\left(a+2\right)\cdot\cdot\cdot\left(a+1+i\right)} 
+A_i\tau_0^{N-\frac{5}{2}-m-c}\right\} \nonumber\\
&&\hspace{30mm} \ : \ m=i+1,...,n\nonumber\\
&&\|\s^{(m-1-i,n-m)}\clab\|_{L^2(S_{\ub,u})}\leq \ub^{a+\frac{1}{2}+i}u^{b}
\left\{\frac{C\sqrt{D}}{\left(a+\frac{3}{2}\right)\cdot\cdot\cdot\left(a+\frac{1}{2}+i\right)}
+A_i\tau_0^{N-\frac{5}{2}-m-c}\right\} \nonumber\\
&&\hspace{30mm} \ : \ m=i+1,...,n\nonumber\\
&&\label{14.a286}
\end{eqnarray} 
The inequalities \ref{14.a28} are likewise modified to:
\begin{eqnarray}
&&\|\Omega^{n-m}T^{m-i}\chf\|_{L^2(S_{\tau,\tau})}\leq \tau^{c-\frac{3}{2}+i}
\left\{\frac{C\sqrt{D}}{\left(c-\frac{1}{2}\right)
\cdot\cdot\cdot\left(c-\frac{3}{2}+i\right)}+B_i\tau_0^{N-\frac{5}{2}-m-c}\right\} \nonumber\\
&&\|\Omega^{n-m}T^{m-i}\cv\|_{L^2(S_{\tau,\tau})}\leq \tau^{c-\frac{3}{2}+i}
\left\{\frac{C\sqrt{D}}{\left(c-\frac{1}{2}\right)
\cdot\cdot\cdot\left(c-\frac{3}{2}+i\right)}+B_i\tau_0^{N-\frac{5}{2}-m-c}\right\}\nonumber\\
&&\|\Omega^{n-m}T^{m-i}\cga\|_{L^2(S_{\tau,\tau})}\leq \tau^{c-\frac{3}{2}+i}
\left\{\frac{C\sqrt{D}}{\left(c-\frac{1}{2}\right)\cdot\cdot\cdot\left(c-\frac{3}{2}+i\right)} 
+B_i\tau_0^{N-\frac{5}{2}-m-c}\right\} \nonumber\\
&&\mbox{: for $i=1,...,n-1$, $m=i+1,...,n$} \label{14.a287}
\end{eqnarray}
Here we denote by $B_i$ positive constants depending on $i$ as well as on the exponent $c$ which are 
non-increasing in $c$. 

Proceeding to the estimates \ref{14.330}, these are accordingly modified to:
\begin{eqnarray}
&&\|\s^{(n-1-i))}\check{\sh}\|_{L^2(S_{\ub,u})}\leq \ub^{a+\frac{1}{2}+i}u^{b+\frac{1}{2}}
\left\{\frac{C\sqrt{D}}{\left(a+\frac{3}{2}\right)\cdot\cdot\cdot\left(a+\frac{1}{2}+i\right)}
+A_i\tau_0^{N-3-i-c}\right\}\nonumber\\
&&\|\Omega^{n-1-i}\check{b}\|_{L^2(S_{\ub,u})}\leq \ub^{a+\frac{1}{2}+i}u^{b+2}
\left\{\frac{C\sqrt{D}}{\left(a+\frac{3}{2}\right)\cdot\cdot\cdot\left(a+\frac{1}{2}+i\right)}
+A_i\tau_0^{N-3-i-c}\right\}\nonumber\\
&&\hspace{18mm} : \ i=1,...,n-1\label{14.a288}
\end{eqnarray}
The estimates \ref{14.332} are then modified to: 
\begin{eqnarray}
&&\sup_{S_{\ub,u}}|\s^{(n-2-i)}\check{\sh}|\leq 
\left\{\frac{C\sqrt{D}}{\left(a+\frac{3}{2}\right)\cdot\cdot\cdot\left(a+\frac{1}{2}+i\right)}
+A_i\tau_0^{N-3-i-c}\right\}
\frac{\ub^{a+1+i}u^{b+\frac{1}{2}}}{\sqrt{a+\frac{3}{2}+i}}\nonumber\\
&&\sup_{S_{\ub,u}}|\Omega^{(n-2-i)}\check{b}|\leq 
\left\{\frac{C\sqrt{D}}{\left(a+\frac{3}{2}\right)\cdot\cdot\cdot\left(a+\frac{1}{2}+i\right)}
A_i\tau_0^{N-3-i-c}\right\}
\frac{\ub^{a+1+i}u^{b+2}}{\sqrt{a+\frac{3}{2}+i}}\nonumber\\
&&\hspace{18mm} : \ i=0,...,n-2\label{14.a289}
\end{eqnarray}

In the same way, the estimates \ref{14.337}, \ref{14.342}, \ref{14.344} are modified to:
\begin{eqnarray}
&&\sup_{S_{\ub,u}}|\s^{(m-1-i,n-m)}\check{\dot{\phi}}_\mu|
\leq\left\{\frac{C\sqrt{D}}{\left(a+\frac{3}{2}\right)\cdot\cdot\cdot\left(a+\frac{1}{2}+i\right)}
+A_i\tau_0^{N-2-m-c}\right\}
\frac{\ub^{a+1+i}u^{b-\frac{1}{2}}}{\sqrt{a+\frac{3}{2}+i}}\nonumber\\
&&\ : \ i=0,...,n-1 \ ; \ m=i+1,...,n \label{14.a290}\\
&&\sup_{S_{\ub,u}}|\stackrel{i}{q}_{1,m-i,\mu}|\leq 
\left\{\frac{C\sqrt{D}}{\left(a+\frac{3}{2}\right)\cdot\cdot\cdot\left(a-\frac{1}{2}+i\right)}
+A_i\tau_0^{N-2-m-c}\right\}
\frac{\ub^{a+i}u^{b+\frac{1}{2}}}{\sqrt{a+\frac{1}{2}+i}}\nonumber\\
&&\ : \ i=1,,...,n-1 \ : \ m=i+1,....,n \label{14.a291}\\
&&\sup_{S_{\ub,u}}|\stackrel{i}{\qb}_{1,m-i,\mu}|\leq 
\left\{\frac{C\sqrt{D}}{\left(a+\frac{3}{2}\right)\cdot\cdot\cdot\left(a-\frac{1}{2}+i\right)}
+A_i\tau_0^{N-2-m-c}\right\}
\frac{\ub^{a+i}u^{b-\frac{1}{2}}}{\sqrt{a+\frac{1}{2}+i}}\nonumber\\
&&\ : \ i=1,...,n-1 \ ; \ m=i+1,...,n \label{14.a292}
\end{eqnarray}
Moreover, the estimates \ref{14.348}, \ref{14.355} are modified to: 
\begin{eqnarray}
&&\sup_{S_{\ub,u}}|\s^{(n-2-i)}\ctchi|\leq  
\left\{\frac{C\sqrt{D}}{\left(a+\frac{3}{2}\right)\cdot\cdot\cdot\left(a-\frac{1}{2}+i\right)}
+A_i\tau_0^{N-3-i-c}\right\}
\frac{\ub^{a+i}u^{b+\frac{1}{2}}}{\sqrt{a+\frac{1}{2}+i}}\nonumber\\
&&\sup_{S_{\ub,u}}|\s^{(n-2-i)}\ctchib|\leq  
\left\{\frac{C\sqrt{D}}{\left(a+\frac{3}{2}\right)\cdot\cdot\cdot\left(a-\frac{1}{2}+i\right)}
+A_i\tau_0^{N-3-i-c}\right\}
\frac{\ub^{a+i}u^{b}}{\sqrt{a+\frac{1}{2}+i}}\nonumber\\
&&\ : \ i=1,...,n-2 \label{14.a293}\\
&&\sup_{S_{\ub,u}}|\s^{(m-1-i,n-m)}\cla|\leq  
\left\{\frac{C\sqrt{D}}{(a+2)\cdot\cdot\cdot(a+i)}+A_i\tau_0^{N-\frac{5}{2}-m-c}\right\}
\frac{\ub^{a+\frac{1}{2}+i}u^{b+\frac{1}{2}}}{\sqrt{a+1+i}}\nonumber\\
&&\sup_{S_{\ub,u}}|\s^{(m-1-i,n-m)}\clab|\leq  
\left\{\frac{C\sqrt{D}}{\left(a+\frac{3}{2}\right)\cdot\cdot\cdot\left(a-\frac{1}{2}+i\right)}
+A_i\tau_0^{N-\frac{5}{2}-m-c}\right\}
\frac{\ub^{a+i}u^{b}}{\sqrt{a+\frac{1}{2}+i}}\nonumber\\
&&\ : i=1,...,n-1 \ : \ m=i+1,...,n \label{14.a294}
\end{eqnarray}

Finally, the estimates \ref{14.a29} are modified to: 
\begin{eqnarray}
&&\sup_{S_{\tau,\tau}}|\Omega^{n-m}T^{m-i}\chf|\leq\left\{\frac{C\sqrt{D}}{\left(c-\frac{1}{2}\right)
\cdot\cdot\cdot\left(c-\frac{5}{2}+i\right)}+B_i\tau_0^{N-\frac{5}{2}-m-c}\right\}
\frac{\tau^{c-2+i}}{\sqrt{c-\frac{3}{2}+i}}
\nonumber\\
&&\sup_{S_{\tau,\tau}}|\Omega^{n-m}T^{m-i}\cv|\leq\left\{\frac{C\sqrt{D}}{\left(c-\frac{1}{2}\right)
\cdot\cdot\cdot\left(c-\frac{5}{2}+i\right)}+B_i\tau_0^{N-\frac{5}{2}-m-c}\right\}
\frac{\tau^{c-2+i}}{\sqrt{c-\frac{3}{2}+i}}
\nonumber\\
&&\sup_{S_{\tau,\tau}}|\Omega^{n-m}T^{m-i}\cga|\leq\left\{\frac{C\sqrt{D}}{\left(c-\frac{1}{2}\right)
\cdot\cdot\cdot\left(c-\frac{5}{2}+i\right)}+B_i\tau_0^{N-\frac{5}{2}-m-c}\right\}
\frac{\tau^{c-2+i}}{\sqrt{c-\frac{3}{2}+i}}
\nonumber\\
&& \ : i=1,...,n-1 \ : \ m=i,...,n \label{14.a295}
\end{eqnarray}

\vspace{5mm}

To fix the notation, we denote by $(x_{\tau_0}^\mu, b_{\tau_0}, \beta_{\tau_0,\mu})$ 
a solution arising from initial data on $\Cb_{\tau_0}$ as in the 
present section, and defined on a domain ${\cal R}_{\ub_1,u_1,\tau_0}$ for some $\ub_1>\tau_0$ and 
$u_1\geq \ub_1$ (see \ref{14.a211}, \ref{14.a212}). The positive real number $\tau_0$ shall be taken from the 
sequence \ref{14.a30}, the starting point of which is $\tau_{0,M}$ for sufficiently large $M$, in 
accordance with the statement just above \ref{14.a60} and the statement just above \ref{14.a187}. 
The bootstrap assumptions of Section 14.1 refer now to the solution 
$(x_{\tau_0}^\mu, b_{\tau_0}, \beta_{\tau_0,\mu})$ and to the domain ${\cal R}_{\ub_1,u_1,\tau_0}$, 
the quantities $\Ab_*$, $\oA*$, $\Bb$, $\oB$, $B$, $C_*$, $C^\prime_*$, $C^{\prime\prime}_*$,  
which refer to the $N$th approximate solution on ${\cal R}_{\delta_0,\delta_0}$, being as before. 
Given any such $\tau_0$, the local shock continuation theorem gives us a solution 
$(x_{\tau_0}^\mu, b_{\tau_0}, \beta_{\tau_0,\mu})$ defined on a domain 
${\cal R}_{\tau_0+\vep,\tau_0+\vep,\tau_0}$ for some $\vep>0$. Moreover, by 
continuity, suitably restricting if necessary $\vep$, the bootstrap assumptions can be 
assumed to hold on ${\cal R}_{\tau_0+\vep,\tau_0+\vep,\tau_0}$. Note that this $\vep$ may be much 
smaller than $\tau_0$. Fixing then any $\delta<\delta_*$ (see \ref{14.a30}), let us denote by 
$\ou_*$ the least upper bound of the set of real numbers $\ou$ greater than $\tau_0$ but not exceeding 
$\delta$ such that we have a solution $(x_{\tau_0}^\mu, b_{\tau_0}, \beta_{\tau_0,\mu})$ defined 
on ${\cal R}_{\ou,\ou,\tau_0}$ and satisfying the bootstrap assumptions there. We then have 
$\ou_*\geq\tau_0+\vep$ and $\ou_*\leq\delta$. We shall show that in fact $\ou_*=\delta$. The argument 
is by contradiction. Suppose that $\ou_*<\delta$. Since for any $\ou<\ou_*$ there is a solution 
$(x_{\tau_0}^\mu, b_{\tau_0}, \beta_{\tau_0,\mu})$ defined on ${\cal R}_{\ou,\ou,\tau_0}$ and 
satisfying the bootstrap assumptions there, all the estimates of the present section, including the 
top order energy estimates as modified according to the discussion following \ref{14.a229} hold 
on ${\cal R}_{\ou,\ou,\tau_0}$, for all $\ou<\ou_*$. The bounds being independent of $\ou$ 
the solution $(x_{\tau_0}^\mu, b_{\tau_0}, \beta_{\tau_0,\mu})$ extends continuously with all its 
partial derivatives of order up to $n-1$ to the future boundary $C_{\ou_*,\tau_0}^{\ou_*}$ of 
${\cal R}_{\ou_*,\ou_*,\tau_0}$. Then with  $M$ taken suitably large the estimates 
\ref{14.a289} - \ref{14.a295} imply that the bootstrap assumptions hold as strict inequalities on 
${\cal R}_{\ou_*,\ou_*,\tau_0}$, including therefore $C_{\ou_*,\tau_0}^{\ou_*}$. Moreover, the 
$L^2(S_{\ub,u})$ estimates hold also for $u=\ou_*$ and all $\ub\in[\tau_0,\ub_*]$ and the 
$L^2(C_{\ou_*,\tau_0}^{\ou_*})$ hold as well. Thus we have characteristic initial data on 
$C_{\ou_*,\tau_0}^{\ou_*}$ as well as on the part of $\Cb_{\tau_0}$ corresponding to 
$\ou_*\leq u\leq\delta_*$. Then by the local existence theorem for the standard characteristic 
initial value problem, we have a solution in the domain in ${\cal N}_{\tau_0}$ corresponding 
to the rectangle $[\tau_0,\tau_0+\vep^\prime]\times[\ou_*,\ou_*+\vep^\prime]$ 
in the $(\ub,u)$ plane, for some $\vep^\prime>0$, 
and, restricting $\vep^\prime$ if necessary, satisfying the bootstrap assumptions on this domain. 
The solution $(x_{\tau_0}^\mu, b_{\tau_0}, \beta_{\tau_0,\mu})$ is then extended to 
${\cal R}_{\ou_*,\ou_*,\tau_0}\bigcup {\cal R}_{\tau_0+\vep^\prime,\ou_*+\vep^\prime,\tau_0}$ the 
aforementioned domain corresponding to ${\cal R}_{\tau_0+\vep^\prime,\ou_*+\vep^\prime,\tau_0}
\setminus {\cal R}_{\ou_*,\ou_*,\tau_0}$. Let us denote by 
$\oub_*$ the least upper bound of the set of real numbers $\oub$ greater than $\tau_0$ but not exceeding 
$\ou_*$ for each of which the solution extends to the domain corresponding to a maximal rectangle 
in the $(\ub,u)$ plane of the form $[\tau_0,\oub]\times[\ou_*,\ou_*+\vepb(\oub)]$, 
$0<\vepb(\oub)\leq\delta$, with the bootstrap assumptions holding there. We shall presently 
show that $\oub_*=\ou_*$, that is $S_{\oub_*,\ou_*}\subset{\cal K}_{\tau_0}$. For, suppose that on the 
contrary, $\oub_*<\ou_*$. Since for any $\oub<\oub_*$, the solution extends to 
${\cal R}_{\ou_*,\ou_*,\tau_0}\bigcup {\cal R}_{\oub,\ou_*+\vepb(\oub),\tau_0}$ with the bootstrap 
assumptions holding there, all the estimates also hold in this domain. But the local existence 
theorem for the standard characteristic initial value problem provides a solution the domain 
corresponding to a rectangle $[\ub_1,\ub_2]\times[u_1,u_2]$ in the $(\ub,u)$ plane given initial data 
on the characteristic hypersurfaces corresponding to the sides $\{\ub_1\}\times[u_1,u_2]$, 
$[\ub_1,\ub_2]\times\{u_1\}$ provided that these sides are suitably small, that is the positive real 
numbers $\ub_2-\ub_1$, $u_2-u_1$ are suitably small, {\em depending only on the norm of the data}. 
Here, the norm of the data on the part of $C_{\ou_*,\tau_0}^{\ou_*}$ corresponding to 
$\oub\leq \ub\leq \oub_*$, in fact on all of $C_{\ou_*,\tau_0}^{\ou_*}$, as well as the norm of the 
data on the part of $\Cb_{\oub}^{\ou_*+\vepb(\oub)}$ corresponding to $u\geq\ou_*$, 
in fact in all of $\Cb_{\oub}^{\ou_*+\vepb(\oub)}$, satisfy fixed bounds by virtue of the aforementioned 
estimates. Therefore, taking $\oub$ suitably close to $\oub_*$ we have a solution on the 
domain corresponding to a rectangle in the $(\ub,u)$ plane of the form 
$[\oub,\oub_*]\times[\ou_*,\ou_*+\vepb(\oub_*)]$ for some $\vepb(\oub_*)>0$ and not exceeding 
$\vepb(\oub)$. Since the bootstrap assumptions are satisfied as strict inequalities on all of 
$C_{\ou_*,\tau_0}^{\ou_*}$, suitably restricting $\vepb(\oub_*)$ if necessary, the bootstrap 
assumptions also hold in the domain corresponding to the rectangle 
$[\oub,\oub_*]\times[\ou_*,\ou_*+\vepb(\oub_*)]$. Therefore the solution 
$(x_{\tau_0}^\mu, b_{\tau_0}, \beta_{\tau_0,\mu})$ extends to the domain 
${\cal R}_{\ou_*,\ou_*,\tau_0}\bigcup {\cal R}_{\oub_*,\ou_*+\vepb(\oub_*),\tau_0}$ 
and the bootstrap assumptions hold there, hence so do all the estimates, and the bootstrap 
assumptions then hold as strict inequalities in this closed domain. Then if $\oub_*< \ou_*$ the local 
existence theorem for the standard characteristic initial value problem with initial data 
on the part of $\Cb^{\ou_*+\vepb(\oub_*)}_{\oub_*}$ corresponding to $u\geq\ou_*$ and 
on the part of $C_{\ou_*,\tau_0}^{\ou_*}$ corresponding to $\ub\geq\oub_*$ give an extension 
to the domain corresponding to a rectangle in the $(\ub,u)$ plane of the form 
$[\oub_*,\oub]\times [\ou_*,\ou_*+\vepb(\oub)]$ for some $\oub>\oub_*$ and $\vepb(\oub)>0$, 
contradicting the definition of $\oub_*$. We conclude that $\oub_*=\ou_*$. The same argument 
then applies with $\ou_*$ in the role of $\oub_*$ taking $\oub$ suitably close to $\ou_*$. 
This shows that the solution $(x_{\tau_0}^\mu, b_{\tau_0}, \beta_{\tau_0,\mu})$  extends to  
${\cal R}_{\ou_*,\ou_*+\vep,\tau_0}$  and the bootstrap assumptions hold there, hence so do all 
the estimates, and then the bootstrap assumptions hold as strict inequalities in this closed domain. 
We are then in a position to apply the local shock continuation theorem to obtain an extension of the solution 
to a domain corresponding in the $(\ub,u)$ plane to the triangle: 
$$\{(\ub,u) \ : \ u\in[\ub,\ub+\vep_*], \ \ub\in[\ou_*,\ou_*+\vep_*]\}$$
for some $\vep_*>0$, not exceeding $\vep$, with the bootstrap assumptions holding on this domain. 
Since the bootstrap assumptions then hold in the union with 
${\cal R}_{\ou_*,\ou_*+\vep,\tau_0}$, which contains ${\cal R}_{\ou_*+\vep_*,\ou_*+\vep_*,\tau_0}$, 
we arrive at a contradiction to the definition of $\ou_*$ above. We have thus demonstrated that 
$\ou_*=\delta$. 

In conclusion, for any $\tau_0\in(0,\tau_{0,M}]$, in particular any member of the sequence
\begin{equation}
\left(\tau_{0,m} \ : \ m=M,M+1,M+2, \ . \ . \ .\right)
\label{14.a296}
\end{equation}
we have a solution 
$(x_{\tau_0}^\mu, b_{\tau_0}, \beta_{\tau_0,\mu})$ defined on ${\cal R}_{\delta,\delta,\tau_0}$ 
and satisfying there all the estimates of the present section. In particular the pointwise estimates 
\ref{14.a289} - \ref{14.a294} hold for the differences from the $N$th approximate solution 
$(x_N^\mu, b_N, \beta_{\mu,N})$, the last being of course defined on all of ${\cal R}_{\delta,\delta}$. 
These give uniform bounds on all partial derivatives of the difference 
$(x_{\tau_0}^{\mu}-x_N^{\mu}, b_{\tau_0}-b_N, \beta_{\tau_0,\mu}-\beta_{N,\mu})$ of order up to 
$n-1$ on ${\cal R}_{\delta,\delta,\tau_0}$, tending rapidly as $\tau_0\rightarrow 0$ 
to the corresponding bounds of 
Section 14.9. In fact, by the same methods uniform H\"{o}lder $C^{0,1/2}$ bounds for the 
$n-1$th order derivatives can also be deduced. It then follows by the Ascoli-Arzel\`{a} theorem 
that there is a subsequence of the sequence \ref{14.a296} and $C^{n-1}$ functions 
$(x^\mu-x_N^\mu, b-b_N, \beta_\mu-\beta_{\mu,N})$ on the closed domain ${\cal R}_{\delta,\delta}$, 
vanishing on $\Cb_0^{\delta}$ together with the partial derivatives with respect to $\tau$ of 
order up to $n-1$, such that for any $\eta>0$ the corresponding subsequence of 
$(x_{\tau_{0,m}}^\mu-x_N^\mu, b_{\tau_{0,m}}-b_N, \beta_{\mu,\tau_{0,m}}-\beta_{\mu,N})$ converges 
uniformly on ${\cal R}_{\delta,\delta,\eta}$ with all its derivatives of order up to $n-1$ to 
$(x^\mu-x_N^\mu, b-b_N, \beta_\mu-\beta_{\mu,N})$. Moreover, in view of the fact that the bounds 
\ref{14.a289} - \ref{14.a294} all contain a factor of a power at least $a+1$ of $\ub$, therefore 
are bounded in the part of ${\cal R}_{\delta,\delta}$ corresponding to $\ub\leq\eta$ by a constant 
multiple of $\eta^{a+1}$, the convergence is in fact uniform on all of ${\cal R}_{\delta,\delta}$. 
Moreover, the pointwise estimates \ref{14.a295} hold for the differences 
$(\chf_{\tau_0}, \cv_{\tau_0}, \cga_{\tau_0})$, defined on ${\cal K}_{\tau_0}^\delta$, of the transformation functions associated to 
the solution $(x_{\tau_0}^\mu, b_{\tau_0}, \beta_{\tau_0,\mu})$ from the corresponding transformation 
functions associated to the $N$th approximate solution. These give uniform bounds on the partial 
derivatives with respect to $\vartheta$ and $\tau$ of order up to $n-1$  of 
$(\chf_{\tau_0}, \cv_{\tau_0}, \cga_{\tau_0})$ on ${\cal K}_{\tau_0}^\delta$, tending rapidly 
as $\tau_0\rightarrow 0$ to the corresponding bounds of Section 14.9. 
In fact, by the same methods uniform H\"{o}lder $C^{0,1/2}$ bounds for the 
$n-1$th order derivatives can also be deduced. Then again by the Ascoli-Arzel\`{a} theorem 
there is a subsequence of previous subsequence and $C^{n-1}$ functions $(\chf, \cv, \cga)$ on 
${\cal K}^{\delta}$, vanishing at $S_{0,0}$ together with 
the partial derivatives with respect to $\tau$ of order up to $n-1$, such that for any $\eta>0$ 
the corresponding subsequence of $(\chf_{\tau_0,m}, \cv_{\tau_0,m}, \cga_{\tau_0,m})$ converges 
uniformly on ${\cal K}_\eta^\delta$ with all its derivatives of order up to $n-1$ to 
$(\chf, \cv, \cga)$. Moreover, in view of the fact that the bounds 
\ref{14.a295} contain a factor of a power at least $c-1$ of $\tau$, therefore 
are bounded in ${\cal K}^{\eta}$  by a constant 
multiple of $\eta^{c-1}$, the convergence is in fact uniform on all of ${\cal K}^{\delta}$. 
It follows that $(x^\mu, b, \beta_\mu)$ is a solution of the characteristic and wave systems 
in ${\cal R}_{\delta,\delta}$, satisfying the initial conditions on $\Cb_0^{\delta}$, and together 
with the associated transformation functions 
$(\hat{f}=\chf+\hat{f}_N, v=\cv+v_N, \gamma=\cga+\gamma_N)$ the identification equations and the 
boundary conditions on ${\cal K}^{\delta}$. {\bf We therefore have a $C^{n-1}$ solution of the 
shock development problem on ${\cal R}_{\delta,\delta}$.}

Furthermore, given any $S_{\ub,u}\subset {\cal R}_{\delta,\delta}$, by weak compactness of closed balls 
in $H_k(S^1)$, for a subsequence of the aforementioned subsequence the restrictions to 
$S_{\ub,u}$ of derived functions associated to 
$(x_{\tau_{0,m}}^\mu, b_{\tau_{0,m}}, \beta_{\tau_{0,m}})$ and 
$(\hat{f}_{\tau_{0,m}}, v_{\tau_{0,m}}, \gamma_{\tau_{0,m}})$ converge weakly in $H_k(S^1)$, for the 
appropriate $k$. It then follows that the restriction to $S_{\ub,u}$ of the corresponding derived 
functions associated to $(x^\mu, b, \beta_\mu)$ and $(\hat{f}, \cv, \cga)$ belong to $H_k(S^1)$ 
and satisfy the appropriate $H_k(S^1)$ bounds. As a consequence the solution possess weak 
derivatives of order $n$ in $L^2(S_{\ub,u})$ satisfying the bounds of Proposition 14.3, Lemma 14.1, 
Proposition 14.4, and the estimates \ref{14.a27}. Similarly, given any 
${\cal K}_{\sigma}^{\delta-\sigma}$, $\sigma\in [0,\delta)$ (see Section 14.4), 
by weak compactness of closed balls in $H_k([0,\delta-\sigma]\times S^1)$, for a subsequence of the 
aforementioned subsequence, the restrictions to ${\cal K}_{\sigma}^{\delta-\sigma}$ of 
the corresponding $\tchi_{\tau_{0,m}}$, $\tchib_{\tau_{0,m}}$, $\lambda_{\tau_{0,m}}$, 
$\lambdab_{\tau_{0,m}}$ converge weakly in $H^k([0,\delta-\sigma]\times S^1)$ for $k=n-1$ in the case of 
the first pair, $k=n$ in the case of the second pair. As a consequence the functions 
$\tchi$, $\tchib$, $\lambda$, $\lambdab$ associated to the solution belong to 
$H^k([0,\delta-\sigma]\times S^1)$ for $k=n-1$, $k=n$ accordingly, 
and satisfy the bounds of Proposition 14.2. 
The same argument in connection with ${\cal K}^{\delta}={\cal K}_0^{\delta}$, 
shows that the transformation functions associated to the solution possess 
weak derivatives of order $n+1$ in $L^2({\cal K}^{\delta})$ satisfying the estimates \ref{14.a26}. 
Finally a similar argument in connection with the $C_u^u$, $u\in(0,\delta]$, and with the 
$\Cb_{\ub}^{\delta}$, $\ub\in [0,\delta)$ shows that the solution possesses well defined top order 
energies and fluxes satisfying the top order energy estimates of Section 13.5. 

Now, the smallness restrictions on 
$\delta$ of Section 13.5, that is the restrictions on $\delta$ implied by \ref{13.260}, \ref{13.264}, 
\ref{13.291}, \ref{13.296}, are, like the restrictions on $\delta$ from Chapter 12, 
independent of $n$. On the other hand the 
smallness restriction on $\delta$ required to recover the bootstrap assumptions 
(see end of Section 14.9) depends on the quantity denoted by $D$ which is defined by \ref{14.1} and 
does depend on $n$ as well as $N$. We therefore let $n_0$ be the smallest positive integer $n$ 
satisfying \ref{14.b4} (appropriate to $d=2$ spatial dimensions) that is $n_0=5$, then choose the exponents $a$, $b$ 
large enough as required. We then let $N_0$ be the smallest positive integer $N$ verifying (see 
\ref{14.a247})
\begin{equation}
N>c+n+\frac{5}{2}
\label{14.a297}
\end{equation}
with $n=n_0$. 
The power series approximation can then be truncated at order $N_0$. 
Let $D_0$ be the quantity $D$ corresponding to the choices $n_0$, $N_0$. Given then $n>n_0$, 
the argument having closed at order $n_0$, the bootstrap assumptions are no longer needed, therefore 
no new smallness restrictions on $\delta$ are required to proceed inductively to orders 
$n_0+1,...,n$. However since $c$ depends on $n$ and the condition \ref{14.a297} is required  
the power series approximation is to be truncated at successively higher orders. 
This argument establishes that, 
given that the data of the problem are $C^\infty$, the solution is also $C^\infty$. The 
uniqueness of the solution readily follows by considering the $(m,l)$ energies corresponding 
to their difference for $m+l=1$. That is, in the definitions of Section 9.7 we place one of the 
solutions in the role of the $N$th approximate solution. 

For $d>2$ spatial dimensions, after the essential change in regard to the top order acoustical 
estimates discussed in Chapter 11, the argument is carried out in exactly the same manner, with 
one exception: Lemma 14.2. That is, the imbedding of $H_1(S^1)$ into $C^{0,1/2}(S^1)$ and the 
Sobolev inequality on $S^1$. In $d$ spatial dimensions we define $k_0$ to be the smallest integer 
greater than $(d-1)/2$. Then in the role of Lemma 14.2 we place the imbedding of $H_{k_0}(S_{\ub,u})$ 
into $C^{0,1/2}(S_{\ub,u})$ and the Sobolev inequality:
\begin{equation}
\sup_{S_{\ub,u}}|f|\leq C_S(\ub,u)\|f\|_{H_{k_0}(S_{\ub,u})}
\label{14.a298}
\end{equation}
where $C_S(\ub,u)$, the Sobolev constant of $S_{\ub,u}$, can be bounded in terms of $C_S(0,u-\ub)$, 
the Sobolev constant of $S_{0,u-\ub}$ because the pointwise bounds on $\chi$, $\chib$ contained in the 
bootstrap assumptions allow us to derive upper and lower bounds for the eigenvalues of 
$F_{\ub}^*\left.\sh\right|_{S_{\ub,u}}$, the pullback to $S_{0,u-\ub}$ by the diffeomorphism 
$F_{\ub} \ : \ S_{0,u-\ub}\rightarrow S_{\ub,u}$ (see Section 2.2) of the induced metric on $S_{\ub,u}$, 
with respect to $\left.\sh\right|_{S_{0,u-\ub}}$, the induced metric on $S_{0,u-\ub}$.  The 
condition \ref{14.b4} is replaced by the condition:
\begin{equation}
n_*+2\leq n-k_0
\label{14.a299}
\end{equation}
that is (see \ref{14.a18}):
\begin{equation}
n\geq n_0:=2k_0+2
\label{14.a300}
\end{equation}

We are now ready to state the theorem which concludes this monograph. 
Recall that $\rho$ and $\rhob$ corresponding to the solution 
are smooth functions of $(\ub,u,\vartheta)$ which are approximated by $\rho_{N_0}$, $\rhob_{N_0}$ according to 
the pointwise estimates in connection with the recovery of the bootstrap assumptions. It follows that  
(see \ref{9.a13}, \ref{9.a14}, \ref{9.a19}):
\begin{equation}
\rho=-\frac{k}{3l}\ub+O(\ub u), \ \ \rhob=\frac{k}{12c_0}(3u^2-\ub^2)+O(u^3)
\label{14.a301}
\end{equation}
The first implies that (see \ref{2.69}, \ref{2.70}) along a generator of a given $C_u$, denoting 
by $t_0$ the value of $t$ at the point where the generator intersects $\Cb_0$, $t-t_0$ is 
a smooth function of $\ub$ whose Taylor expansion begins with quadratic terms. Recall also that 
by Proposition the vectorfield $N$ along $\Cb_0$ is equal to the vectorfield $L^\prime$ defined by 
the prior solution. 

\vspace{5mm}

\noindent{\bf {\large{Theorem:}} \ In any spatial dimension $d\geq 2$, given a prior maximal classical solution as discussed in Section 1.5 
there is a $\delta>0$ and a unique solution of the restricted shock development problem stated in Section 1.6, defined on ${\cal R}_{\delta,\delta}$
by the triplet $(x^\mu,b,\beta_\mu)$ where $x^\mu : \mu=0,...,d$ and $\beta_\mu : \mu=0,...,d$ 
are smooth functions and $b$ is a smooth mapping of $R_{\delta,\delta}$ into the space of smooth 
vectorfields on $S^{d-1}$. The smooth mapping ${\cal R}_{\delta,\delta}\rightarrow\mathbb{M}^{d}$ by 
$(\ub,u,\vartheta)\mapsto (x^\mu(\ub,u,\vartheta) :  \mu=0,...,d)$ has negative Jacobian 
in ${\cal R}_{\delta,\delta}$ except on $\Cb_0$ where the Jacobian vanishes. The boundaries $\Cb^\delta_0$ 
and ${\cal K}^\delta$ are mapped onto smooth hypersurfaces in $\mathbb{M}^d$, an acoustically null 
hypersurface $\underline{{\cal C}}$ and the shock hypersurface respectively, the latter being acoustically timelike except 
along its past boundary where it is acoustically null and transverse to the latter. The 
$\beta_\mu$ expressed through the inverse mapping in terms of rectangular coordinates in $\mathbb{M}^d$ 
are smooth functions on the image of ${\cal R}_{\delta,\delta}$ in $\mathbb{M}^d$ except on 
$\underline{{\cal C}}$, the image of $\Cb_0$. The new solution $\beta_\mu$ so expressed extends the prior solution $\beta^\prime_\mu$ in a $C^1$ manner across $\underline{{\cal C}}$. However, 
along a generator of the image of a given $C_u$, a transversal to $\underline{{\cal C}}$ acoustically null hypersurface, $\beta_{\mu}-\beta_{\mu,0}$ expands in half integral powers of $t-t_0$ beginning with the power 2, $t_0$ and $\beta_{\mu,0}$ being the values of $t$ and $\beta_\mu$ respectively at the point where the generator intersects $\underline{{\cal C}}$}.

\pagebreak

\section*{\huge{Bibliography}} 

\vspace{5mm}

\noindent [Be1] \ \ Bernstein, S. ``Sur la g\'{e}n\'{e}ralization du probl\`{e}me de Dirichlet II", 
Math. Ann. {\bf 69}, 82-136 (1910). 

\noindent [Be2] \ \ Bernstein, S. ``Sur les surfaces d\'{e}finies au moyen de leur coubure moyenne 
et totale", Ann. Sci. \'{E}cole Norm. Sup. {\bf 27}, 233-256 (1910). 

\noindent [Be3] \ \ Bernstein, S. ``Sur les \'{e}quations du calcul des variations", Ann. Sci. 
\'{E}cole Norm. Sup. {\bf 29}, 431-485 (1912). 

\noindent [CB] \ \ Choquet-Bruhat, Y. ``Th\'{e}or\`{e}me d'\'{e}xistence pour certain syst\`{e}mes 
d'\'{e}quations aux d\'{e}rive\'{e}s partielles nonlin\'{e}aires", Acta Mathematica {\bf 88}, 141-225 
(1952). 

\noindent [Ch-A] \ \ Christodoulou, D. {\em The Action Principle and Partial Differential Equations}, 
Ann. Math. Stu. {\bf 146}, Princeton University Press, 2000. 

\noindent [Ch-Kl] \ \ Christodoulou, D. and Klainerman, S. {\em The Global Nonlinear Stability of the 
Minkowski Space}, Princeton Mathematical Series {\bf 41}, Princeton University Press, 1993. 

\noindent [Ch-Li] \ \ Christodoulou, D. and Lisibach, A. ``Shock development in spherical symmetry", 
Ann. PDE {\bf 2}, no. 1, Art. 3, 246 pp. (2016). 

\noindent [Ch-Mi] \ \ Christodoulou, D. and Miao, S. {\em Compressible Flow and Euler's Equations}, 
Surveys of Modern Mathematics {\bf 9}, International Press \& Higher Education Press, 2014.

\noindent [Ch-S] \ \ Christodoulou, D. {\em The Formation of Shocks in 3-Dimensional Fluids}, EMS 
Monographs in Mathematics, EMS Publishing House, 2007. 

\noindent [DG] \ \ De Giorgi, E. ``Sula differentiabilit\`{a} e l'analiticit`{a} delle estemali degli 
integrali multipli regolari", Mem. Accad. Sci. Torino Cl. Fis. Mat. Natur. (3) {\bf 3}, 25-43 (1957). 

\noindent [Ga] \ \ G$\stackrel{\circ}{\mbox{a}}$rding, L. ``Le probl\`{e}me de la deriv\'{e}e oblique 
pour l'equation des ondes", C. R. Acad. Sci. Paris, S\'{e}r.  A {\bf 285}, 773-775 (1977). 

\noindent [Je-Se] Jenkins, H. and Serrin, J. ``The Dirichlet problem for the minimal surface equation 
in higher dimensions", J. Reine Angew. Math. {\bf 229}, 170-187 (1968). 

\noindent [Kl] \ \ Klainerman, S. ``Uniform decay estimates and the Lorentz invariance of the classical 
wave equation", Comm. Pure Appl. Math. {\bf 38}, 321-332 (1985). 

\noindent [Ma1] \ \ Majda, A. ``The stability of multidimensional shock fronts", Mem. Amer. Math. Soc. 
{\bf no. 275}, 1983. 

\noindent [Ma2] \ \ Majda, A. ``The existence of multidimensional shock fronts", Mem. Amer. Math. Soc. 
{\bf no. 281}, 1983. 

\noindent [Ma-Th] \ \ Majda, A. and Thomann, E. ``Multidimensional shock fronts for second order wave 
equations", Comm. in Partial Differential Equations {\bf 12}(7), 777-828 (1987). 

\noindent [Na] Nash, J. ``Continuity of solutions of parabolic and elliptic equations", Amer. J. Math. 
{\bf 80}, 931-954 (1958). 

\noindent [Ne] Newton, I. {\em Philosophiae Naturalis Principia Mathematica}, Motte's translation of the 3rd edition (1725) (1st edition 1686). 

\noindent [No] Noether, E. ``Invariante Variationsprobleme", Nachr. Ges. Wiss. G\"{o}tingen, Math-Phys. Kl. {\bf 1918}, 235-257 (1918). 

\noindent [Si-Mo] Siegel, C.L. and Moser J. {\em Lectures on Celestial Mechanics}, Grundlehren 
der Mathematishen Wissenschaften {\bf 187}, Springer 1971.

\end{document}